\documentclass[10pt]{book}
\usepackage{graphics}
\usepackage{graphicx}
\usepackage{amsmath,amscd}
\usepackage{amssymb}
\usepackage{theorem}
\usepackage{gastex}

\usepackage{tikz}
\usepackage{kbordermatrix}
\usetikzlibrary{positioning}
\usepackage{url}
\usepackage[round]{natbib}
\usepackage{multind}
\usepackage{hyperref}
\newtheorem{theorem}{Theorem}[section]
\newtheorem{proposition}[theorem]{Proposition}
\newtheorem{corollary}[theorem]{Corollary}
\newtheorem{lemma}[theorem]{Lemma}

{\theorembodyfont{\rmfamily}%
  \newtheorem{example}[theorem]{Example}
  \newtheorem{remark}[theorem]{Remark} 
\newtheorem{exercise}{\unskip\hspace{3pt}}[chapter]}
\numberwithin{figure}{section}
\newenvironment{solution}[1]{%
  \medskip\par\noindent{\bf{#1}}\ \ }{%
  \vskip\theorempostskipamount\par}
\newenvironment{proof}{\noindent\textit{Proof.}}
{\QED\vskip\theorempostskipamount} 
\newenvironment{proofof}[1]{\noindent\textit{Proof
    \protect{#1}.}}
                       {\QED\vskip\theorempostskipamount}
\def\petitcarre{\vrule height4pt width 4pt depth0pt}
\def\QED{\relax\ifmmode\eqno{\hbox{\petitcarre}}\else{%
  \unskip\nobreak\hfil\penalty50\hskip2em\hbox{}\nobreak\hfil
  \petitcarre
  \parfillskip=0pt \finalhyphendemerits=0\par\smallskip}
  \fi}

\def\exosection#1{\subsection*{#1}}
\newcommand{\R}{\mathbb{R}}
\newcommand{\Z}{\mathbb{Z}}
\newcommand{\N}{\mathbb{N}}
\newcommand{\Q}{\mathbb{Q}}
\newcommand{\T}{\mathbb{T}}
\newcommand{\C}{\mathbb{C}}
\newcommand{\charac}{\chi}

\DeclareMathOperator{\Card}{Card}
\DeclareMathOperator{\Div}{Div}
\DeclareMathOperator{\rank}{rank}
\DeclareMathOperator{\Sep}{Sep}
\DeclareMathOperator{\Hom}{Hom}
\DeclareMathOperator{\im}{Im}
\DeclareMathOperator{\Ad}{Ad}
\DeclareMathOperator{\Pal}{Pal}
\DeclareMathOperator{\Osc}{\rm Osc}
\DeclareMathOperator{\Pref}{\rm Pref}
\def\u(#1){\underline{#1\!}\,}
\DeclareMathOperator{\id}{id}
\DeclareMathOperator{\Inf}{Inf}
\DeclareMathOperator{\Fac}{Fac}
\DeclareMathOperator{\Tr}{Tr}

\DeclareMathOperator{\diam}{diam}
\DeclareMathOperator{\interior}{int}
\DeclareMathOperator{\spr}{spr}
\DeclareMathOperator{\DG}{DG}

\def\1{\mathbf{1}}
\newcommand{\edge}[1]{\stackrel{#1}{\rightarrow}}
\newcommand{\ledge}[1]{\stackrel{#1}{\longrightarrow}}
\newcommand\norm[1]{\lVert#1\rVert}

\newcommand{\Pg}{\mathfrak P}
\newcommand{\Qg}{\mathfrak Q}
\newcommand{\Sg}{\mathfrak S}
\newcommand{\Ga}{\mathfrak A}
\newcommand{\Gb}{\mathfrak B}
\newcommand{\Gk}{\mathfrak K}

\newcommand{\A}{\mathcal A}
\newcommand{\B}{\mathcal B}
\newcommand{\La}{\mathcal L}
\newcommand{\G}{\mathcal{G}}
\newcommand{\HH}{\mathcal{H}}
\newcommand{\I}{\mathcal{I}}
\newcommand{\RR}{\mathcal R}
\newcommand{\CR}{\mathcal CR}
\newcommand{\MR}{\mathcal MR}
\newcommand{\cL}{\mathcal L}
\newcommand{\KK}{\mathcal K}
\newcommand{\M}{\mathcal M}
\newcommand{\Sa}{\mathcal S}
\newcommand{\E}{\mathcal E}
\newcommand{\D}{\mathcal D}
\newcommand{\cT}{\mathcal T}
\newcommand{\F}{\mathcal F}
\newcommand{\Id}{\mathcal J}

\newcommand{\0}{\mathbf{0}}

\numberwithin{equation}{section}

\title{Dimension groups and dynamical systems}
\author{Fabien Durand, Dominique Perrin}
\makeindex{subject}
\makeindex{names}
\makeindex{symbols}
\begin{document}
\maketitle
\tableofcontents
\chapter{Introduction}
In this monograph, we introduce the reader to the connection between
topological dynamical systems and dimension groups.
In this way, we will be able to distinguish topological
dynamical systems, which can appear in many
different forms by comparing their dimension groups,
which are easier to handle.

Dimension groups are ordered abelian groups associated
with a family of associative algebras
called approximately finite,
or AF-algebras. 

These algebras are themselves
a class of $C^*$-algebras which are direct limits
of finite dimensional algebras and were introduced by 
\cite{Bratteli1972}. The algebra is
build from a special kind of graph called
a \emph{Bratteli diagram}\index{subject}{Bratteli diagram}.
\index{names}{Bratteli, Ola}%

Dimension groups
where introduced by \cite{Elliott1976}
\index{names}{Elliott, George A.}% as a
tool for classifying AF-algebras and he proved
that the dimension group (together with an additional
information called the scale) provides a complete
algebraic invariant for these algebras.

The connection of these ideas with dynamical systems
was first done by Krieger \citep{Krieger1977}
(see also~\cite{Krieger1980})
who defined a dimension group for every shift of finite type.
The link with Bratteli diagrams
was done by \cite{Versik1982}\index{names}{Vershik, Anatol M.} who used a lexicographic
order on paths of the Bratteli diagrams to define a topological
dynamical system on the set of infinite paths of the graph.
Later, Herman, Putnam and Skau showed that every minimal
system on a Cantor space is isomorphic to such system.
As a consequence, a dimension group is attached to
any minimal Cantor system and subsequent work
by Giordano, Putnam and Skau~\citep{GiordanoPutnamSkau1995}
\index{names}{Giordano, Thierry}\index{names}{Putnam, Ian F.}%
\index{names}{Skau, Christian F.} showed that
this group is related to the orbit structure of the system.

In this expository presentation, written after the unpublished notes
by Bernard Host~\citep{Host1995}
\index{names}{Host, Bernard}%
 (see also~\citep{Host2000}), we present the basic elements
of this theory, insisting  
on the computational
and algorithmic aspects allowing one to effectively compute
the dimension groups. The computation applies in particular
to the case of substitution shifts, explicitly
presented previously in~\cite{Durand&Host&Skau:1999}.
\index{names}{Durand, Fabien}%

In the first chapter (Chapter~\ref{chapterTopologicalDynamicalSystems})
we present the basic notions of topological dynamical systems.
We restrict our attention to Cantor systems on which acts the group $\Z$
or the semigroup $\N$.
We define recurrent systems and minimal dynamical systems (Section~\ref{sectionRecurrentMinimal}). 
Next, we introduce in Section~\ref{sectionSymbolicSystems}
subshift dynamical systems, which are the
basic systems we are interested in. We define return words
and higher block shifts.
In Section~\ref{sectionSubstitutionSystems}, we introduce 
substitution shifts. 
We define the notion of recognizable
substitution and we state the Theorem of Moss\'e
(Theorem~\ref{theoremMosse}) asserting that any
aperiodic primitive substitution is recognizable.

In the second chapter (Chapter~\ref{chapterDirectLimitsOrderedGroups}),
we shift to an algebraic and combinatorial environment.
We first introduce, in Section~\ref{sectionOrderedGroups},
ordered groups (considering only abelian groups).
We define several notions, as that of order unit
and order ideal. We also define a simple ordered group
as one with no nontrivial ideals.
In Section~\ref{sectionDirectLimits} we define direct limits of ordered groups
and
we give examples of the computation of these ordered groups.
In the last part of this
section (Section~\ref{sectionDimensionGroups}), we finally define
dimension groups. These groups are defined as direct limits
of groups $\Z^n$ with the usual ordering. We prove
the abstract characterization by Effros, Handelman and Shen~\cite{EffrosHandelmanShen1980}
\index{names}{Effros, Edward G.}\index{names}{Handelman, David}%
\index{names}{Shen, Chao-Liang}%
 using the property of Riesz interpolation.

In Chapter~\ref{chapterOrderedCohomology}, we come to notions
of cohomology defined in a Cantor system. We first introduce the
notion of coboundary (Section~\ref{sectionCoboundaries})
and prove in Section~\ref{sectionGH}
the Gottshalk Hedlund Theorem (Proposition~\ref{theoremGH})
\index{names}{Gottschalk, Walter H.}\index{names}{Hedlund, Gustav A.}%
characterizing the continuous functions on a Cantor set
which are coboundaries. We next define the ordered cohomology
group $K^0(X,T)$ of a recurrent system $(X,T)$ as the quotient of the
group of integer valued continuous functions on $X$
by the subgroup formed by coboundaries.
In the next two sections (Sections~\ref{sectionFactorMaps}
and~\ref{sectionInduced}), we consider the effect
on the ordered cohomology group
of applying a factor map or taking the system induced on a clopen
set.
In a second part of this chapter, beginning with
Section~\ref{sectionInvariant}, we define invariant probability
measures on a Cantor system and recall that a substitutive
shift defined by a primitive substitution has a unique invariant
probability measure. We indicate a method to compute this measure.
We show in
Section~\ref{sectionInvariantStates} that there
is a close connection between the cohomology group and the
cone of invariant measures (Proposition~\ref{propositionKerov}).
We use this connection to give a description of the dimension
groups of Sturmian shifts (Theorem~\ref{theoremDGSturm}).

In Chapter~\ref{chapterDimensionGroupsPartitions}, we introduce the fundamental
tool of partitions in towers, or Kakutani-Rokhlin partitions.
We prove the theorem of Herman, Putnam and Skau which shows that any
minimal Cantor system can be represented as the limit of a sequence
of partitions in towers (Theorem~\ref{theoremKRPartitions}). 
In  Chapter~\ref{chapterDimensionGroupsPartitions}
we come back to partition in towers.
We first show how to associate an ordered group
to a partition in towers (Section~\ref{sectionOrderedGroupPartition}).
Next, in
Section~\ref{sectionOrderedGroupSequences}, we
use a sequence of partitions in towers to prove
that the  group $K^0(X,T)$ is, for any minimal dynamical system $(X,T)$,
a simple dimension group (Theorem~\ref{theoremDimensionGroup}) .
In the next sections, we present explicit methods to
compute the dimension group of a minimal shift space.
 In Section
\ref{chapterReturnWords}, we use return words and in Section
\ref{chapterDimensionGroupsRauzyGraphs}, we use Rauzy graphs.
\index{names}{Rauzy, G\'erard}%
Finally, in Section~\ref{chapterFixedPointMorphism}, we show
how to compute the dimension group of a substitutive shift,
as exposed in~\cite{Durand&Host&Skau:1999}.

We introduce Bratteli diagrams in Chapter~\ref{chapterBratteliDiagrams}.
We define the telescoping of a diagram.
We define the dimension group of a Bratteli diagram
and prove that it is a complete invariant for telescoping
equivalence (Theorem~\ref{theoremDGBratteli}). We next
introduce ordered Bratteli diagams and show that one may associate
a dynamical system to every properly ordered Bratteli diagram.
We prove the Bratteli-Versik model Theorem (Theorem~\ref{ch5:theo:BVmodel})
showing the completeness of the model for minimal Cantor systems.
We next prove the Strong Orbit Equivalence Theorem (Theorem \ref{ch5:th:GPS})
showing that dimension groups are a complete invariant for
strong orbit equivalence and the related Orbit Equivalence
Theorem (Theorem~\ref{theoremOrbitEquivalence}).

In Chapter~\ref{ch5:sec:examples}, we focus on substitution shifts
and their representations. We begin by considering odometers,
which have BV-representations close to substitution shifts.
We characterise, as a main result, the family of BV-systems
associated with stationary Bratteli diagrams as the
disjoint union of stationary odometers and substitution
minimal systems (Theorem~\ref{ch5:subsec:Bratteli-substitution}).
We develop next the description of linearly recurrent shifts,
which are characterized by their BV-representation (Theorem~\ref{ch5:theorem:LRBVrepresentation}).
We introduce in Section~\ref{sectionSadicShifts} the notion
of an $\Sa$-adic representation. The main result is an explicit
description of the dimension group of a unimodular $\Sa$-adic shift
(Theorem~\ref{theo:dg}). In the last section (Section~\ref{sectionDerivatives}),
we consider the family of substitutive shifts, a natural
generalization of substitution shifts. The main result is
a characterization by a finiteness condition of substitutive sequences
(Theorem~\ref{theoremCharacterisationSubstitutive}).

Chapter~\ref{chapterDendricShifts} describes the class of dendric shifts, defined
by a restrictive condition on the possible extensions of a word.
This class is a simultaneous generalisations of several other
classes of interest, such as Sturmian shifts or interval exchange shifts
(introduced in the next chapter). The main result is the Return Theorem
(Theorem~\ref{theoremReturn}) which states that the set of return words
in a minimal dendric shift is a basis of the free goup on the alphabet.
We use this result to describe the $\Sa$-adic representation of dendric shifts
and show that it can be defined using elementary automorphisms
of the free group (Theorem~\ref{theoremSadicDendric}).
We illustrate these results by considering the class
of Sturmian shifts (Section~\ref{sectionChapter6Sturmian}).
The last part of the chapter is devoted to specular shifts,
a class of eventually dendric shifts which plays
a role in the next chapter, when we introduce linear involutions.
The main result is a description of the dimension group
of a specular shift (Theorem~\ref{theoremDGSpecular}).

In Chapter \ref{chapterIET}, we introduce the notion of interval
exchange transformation. We prove Keane's Theorem charcterizing
minimal interval exchanges (Theorem~\ref{theoremIDOC}).
We develop the notion of Rauzy induction and characterize the
subintervals reached by iterating the transformation (Theorem~\ref{theo:birauzy1}).
We generalize Rauzy induction to a two-sided version and
characterize the intervals reached by this more genreral transformation
(Theorem~\ref{theo:birauzy2}). We link these transformations
with automorphisms of the free group (Theorem~\ref{theo:inductionbi}).
We also relate these results with the theorem of Boshernizan and Carroll
giving a finiteness condition on the systems induced by an interval
exchange when the lengths of the intervals belong to a quadratic field
(Theorem \ref{theo:quadratic}). In the last section (Section~\ref{sectionInvolutions})
we define linear involutions and show that the natural coding of
a linear involution without connexions is a specular shift (Theorem~\ref{theoremInvolutionSpecular}).

In the last chapter (Chapter~\ref{chapterBratteli}) we give 
a brief introduction to the link between Bratteli diagrams and the the vast subject 
of $C^*$-algebras. We define approximately finite algebras
and show their relation Bratteli diagrams. We relate simple Bratteli diagrams
and simple AF algebras (Theorem~\ref{theoremSimpleAF}). We prove Elliott's Theorem
showing that AF algebras are characterized by their dimension groups
(Theorem~\ref{theoremElliott}).

The book ends with three appendices, to be used as a reference for notions
from several domains of mathematics used in this book . 
\paragraph{Acknowledgements}
This book started as notes following  a seminar held in Marne-la-Vall\'ee
in june 2016 and gathering Marie-Pierre B\'eal,
\index{names}{B\'eal, Marie-Pierre} Val\'erie Berth\'e,
\index{names}{Berth\'e, Val\'erie}
Francesco
Dolce,\index{names}{Dolce, Francesco} Pavel Heller,
\index{names}{Heller, Pavel} Revekka Kyriakoglou,\index{names}{Kyriakoglou, Revekka}
 Julien Leroy,\index{names}{Leroy, Julien} Dominique Perrin\index{names}{Perrin, Dominique} and
Giuseppina Rindone\index{names}{Rindone, Giuseppina}. The text follows,
at least in an initial version,
closely the notes of Bernard Host~\citep{Host1995},
\index{names}{Host, Bernard}
trying to develop more explicitly some arguments.
We have also chosen to follow a slightly different presentation,
not assuming systematically that the dynamical systems are  minimal.
We wish to thank Brian Marcus
\index{names}{Marcus, Brian H.} and Mike Boyle\index{Boyle, Michael}%
 for fruitful discussions
during the preparation of this text. We are also grateful
to Paulina Cecchi,\index{names}{Cecchi Bernales, Paulina}
Francesco Dolce, Pavel Heller,
Maryam Hosseini,\index{names}{Hosseini, Maryam}
Amir Khodayan Karim,\index{names}{Khodayan Karim, Amir}
Revekka Kyriakoglou,
 Christophe Reutenauer, \index{names}{Reutenauer, Christophe}
 and Gw\'ena\"el Richomme \index{names}{Richomme, Gw\'ena\"el}
 for reading the manuscript and founding many errors.
 Special thanks are due to Soren Eilers \index{names}{Eilers, Soren} for reading
 closely Chapter~\ref{chapterBratteli}
  and
 to Christian Choffrut~\index{names}{Choffrut, Christian}
 for his careful reading of most chapters.

%%%%%%%%%%%%%%%%%%%%%%
% chapter Topological dynamical systems
%%%%%%%%%%%%%%%%%
\chapter{Topological dynamical systems}
\label{chapterTopologicalDynamicalSystems}
We present in this chapter some definitions and
basic properties concerning
topological dynamical systems and
symbolic dynamical systems. We define
shift spaces and the important particular
case of substitution shifts, obtained by iterating
a substitution. We also prove some difficult results,
including Mosse recognizability theorem and the
basic results on Sturmian sequences and their generalizations.

This long chapter aims at serving both as an introduction to the
subject and also as a reference when reading forthcoming
chapters. It should thus be considered both a tutorial and a memento.

We begin, in Section~\ref{sectionRecurrentMinimal}, with some
general definitions concerning topological dynamical
systems. The adjective topological is used to distinguish these
systems which are based on a topological space from dynamical
systems based on a measurable space. We define the notion
of recurrent and of minimal system.

In Section~\ref{sectionSymbolicSystems}, we consider shift spaces
and their language. In Sections~\ref{sectionSFT}, 
\ref{sectionSubstitutionSystems}, \ref{sectionSturmianShifts}
and \ref{sectionToeplitz},
we present several particular types of symbolic systems, namely
shifts of finite type, Sturmian shifts, substitution shifts
and finally Toeplitz shifts. 

%%%%%%%%%%%%%%%%%%%%%%%%%
\section{Recurrent and minimal dynamical systems}
\label{sectionRecurrentMinimal}
A \emph{topological dynamical system}
\index{subject}{topological!dynamical system|see{dynamical system}}
\index{subject}{dynamical system!topological}  
is a pair $(X,T)$ where $X$ is a \emph{compact
metric space}
\index{subject}{compact!space}%
\index{subject}{space!compact metric}%
 and $T:X\rightarrow X$ a continuous map.

\begin{example}\label{exampleInit}
As a simple example, consider $X=[0,1]$, which is metric and compact
as a closed interval of the real line $\R$ and the transformation
$T:x\mapsto (x+\alpha)\mod 1$ for some $\alpha\in\R$, which 
is a continuous map from $X$ into $X$.
\end{example}

Thus, in such a system, to each point $x$
in the space $X$ is associated
 a sequence $(x,T(x),T^{2}(x),\ldots)$ of points. It is
convenient to imagine the action of $T$ as the sequence
of positions of the point $x$ in the space $X$ at discrete times $0,1,2,\ldots$.
The effect of the hypothesis that $X$ is compact
is to guarantee that the sequence of these points
will remain at bounded distance of $x$.

When $T$ is a homeomorphism, we say that the system $(X,T)$
is  \emph{invertible}.
\index{subject}{invertible dynamical system}%
\index{subject}{dynamical system!invertible}%
Although
we will meet most of the time invertible dynamical systems,
we do not make this hypothesis systematically, mentioning
each time when it is necessary. Note that, since $X$
is assumed to be compact, if $T$ is invertible, its
inverse is continuous and 
thus $T$ is a homeomorphism (Exercise~\ref{exerciseInverseContinuous}).

\begin{example}\label{exampleInit2}
The system of Example~\ref{exampleInit} is not invertible
since $1$ is not in the image of $T$. If we consider,
instead of $[0,1]$, the \emph{torus} $\T=\R/\Z$
in which $0$ and $1$ are identified, the  transformation $T$
is simply
$x\mapsto x+\alpha$
and becomes a homeomorphism.
\end{example}

For $x\in X$, we often denote $Tx$ instead of $T(x)$. We also
denote $T^0$ for the identity on $X$ and for $n\ge 0$,
$T^{-n}(x)=\{y\in X\mid T^n(y)=x\}$.

An important example is when $X$ is a {\em Cantor space}\index{subject}{Cantor!space}\index{subject}{space!Cantor}\index{names}{Cantor, Georg}, that is, 
a \emph{totally disconnected}\index{subject}{totally disconnected}
\index{subject}{disconnected!totally}%
\index{subject}{space!totally disconnected}%
 compact metric space without \emph{isolated points}.
\index{isolated point}
 We say then that
$(X,T)$ is a \emph{Cantor dynamical system}\index{subject}{Cantor!dynamical system}
\index{subject}{dynamical system!Cantor}
or \emph{Cantor system} (we shall come back shortly to Cantor spaces).

In particular, let $A$ be a finite set called an \emph{alphabet}.
The set $A^\Z$ of all bi-infinite sequences endowed with the
product topology is a compact space. Let $d$ be the distance
 on $A^\Z$ defined for $x\ne y$ by 
\begin{displaymath}
d(x,y)=1/\max\{n\ge 0\mid x_i=y_i, |i|\le n\}
\end{displaymath}
The topology defined by this distance is the same as the product
topology and thus $A^\Z$ is a compact metric space.
It is actually a Cantor space (see below). 

The \emph{shift}\index{subject}{shift transformation} transformation $S :  A^\Z \to A^\Z$ is defined
for $x=(x_n)_{n\in \Z}$ by $y=Sx$ where
\begin{equation}
y_n=x_{n+1},\label{eqShift}
\end{equation}
for all $n\in\Z$. The shift is obviously continuous and thus $(A^\Z,S)$ is a 
topological dynamical
system. 

As a variant, the set $A^\N$ of one-sided infinite words is also a
topological space for the product topology and this topology
is defined by the metric analogous to the one above. It is also
a Cantor space. The one-sided shift
transformation is defined by Equation \eqref{eqShift} for $n\in\N$. It is not invertible
as soon as $\Card(A)\ge 2$. Thus $(A^\N,S)$ is an example
of a non invertible dynamical system on a Cantor space.

\subsection{Recurrent dynamical systems}
A point $x\in X$ in a topological dynamical system $(X,T)$ is \emph{recurrent}
\index{subject}{recurrent!point} if for every open set $U$ containing $x$
there is an $n\ge 1$ such that $T^n(x)\in U$
(and, in fact, then an infinity of such $n$, 
see Exercise~\ref{exerciseRecurrentPoint}).

A system $(X,T)$ is \emph{recurrent}%
\index{subject}{recurrent!dynamical system}\index{subject}{dynamical system!recurrent} 
(or \emph{topologically transitive}%
\index{subject}{topologically transitive}%
\index{subject}{dynamical system!topologically transitive}) if for every pair of
nonempty open sets $U,V$ in $X$, there is an integer $n\geq 0$
such that $U\cap T^{-n}V\ne\emptyset$
(or equivalently $T^{n}U\cap V\ne\emptyset$, see Exercise~\ref{exerciseRecurrent}).

The \emph{positive orbit}\index{subject}{orbit!positive} of a point $x\in X$
in a dynamical system $(X,T)$ is the set
$\mathcal{O}_+ (x) = \{ T^nx \mid n\ge 0 \}$. Its \emph{orbit}
\index{subject}{orbit}%
\index{symbols}{O(x)@$\mathcal{O}_+(x)$}%
is the set $\mathcal{O}(x)=\cup_{n\ge 0}T^{-n}(\mathcal{O}_+(x))$.
\index{symbols}{O@$\mathcal{O}(x)$}% 
Thus, 
we have also $\mathcal{O}(x)=\{y\in X\mid T^nx=T^my \mbox{ for some } m,n\ge 0\}$
and when $T$ is invertible,
$\mathcal{O}(x)=\{T^nx\mid n\in\Z\}$.
If a point has a dense positive orbit, it is recurrent.

\begin{proposition}\label{propositionRecurrentSystem}
The following conditions are equivalent for a topological dynamical
system.
\begin{enumerate}
\item[\rm(i)] $(X,T)$ is recurrent.
%\item[\rm(ii)] The set of recurrent points is dense.
\item[\rm(ii)] There is a  point in $X$ with a dense positive orbit.
\end{enumerate}
\end{proposition}
The proof is left as Exercise~\ref{exerciseCondRecurrent}.

A \emph{morphism}\index{subject}{morphism!of dynamical systems}
\index{subject}{dynamical system!morphism of}
of dynamical systems from $(X,T)$ to $(X',T')$ is a continuous map
$\phi:X\rightarrow X'$ such that $\phi\circ T=T'\circ \phi$
(see the diagram below).
\begin{displaymath}
\begin{CD}
X@>{T}>>X\\
@VV{\phi}V @VV{\phi}V\\
X'@>{T'}>>X'
\end{CD}
\end{displaymath}
If $\phi$ is onto, it is
called a \emph{factor map}\index{subject}{factor!map} 
and $(X',T')$ is called a \emph{factor system}\index{subject}{factor!system} of $(X,T)$.

A factor of a topological dynamical system inherits some of its dynamical properties.
For instance, a factor of a recurrent system is recurrent.

When it is bijective, the morphism $\phi$
is called an \emph{isomorphism} of dynamical systems,
\index{subject}{isomorphism!of dynamical systems}%
or also a \emph{topological conjugacy}
\index{subject}{topological!conjugacy}%
(or simply a \emph{conjugacy}).
\index{subject}{conjugacy!of dynamical systems}%
  Since $X$ is compact, the inverse is also continuous
(Exercise~\ref{exerciseInverseContinuous})
and thus $\phi$ is a homeomorphism.

Conjugate systems are indistinguishable concerning their dynamical properties.
It can be very difficult to exhibit a conjugacy between dynamical systems
and much of what will follow in this book addresses this problem. In particular,
we will be looking for \emph{invariants}\index{subject}{invariant!under conjugacy},
that is properties shared by conjugate systems and easier to determine 
than the conjugacy itself. Even more interesting are \emph{complete invariants}\index{subject}{complete!invariant}\index{subject}{invariant!complete},
which characterize the conjugacy class. This applies to other equivalences
than conjugacy, as we shall see. Formally, let $\theta$
be an equivalence on the class of dynamical systems, such as conjugacy.
A map $f$ assigning
an element $f(X,T)$ to a dynamical system $(X,T)$ is an invariant
for $\theta$ if $(X,T)\equiv(X',T')\bmod\theta$ implies $f(X,T)=f(X',T')$.
It is a complete invariant if the converse holds.

We will also be interested in the effective
verification of properties of dynamical systems (such as conjugacy for example).
We say that a property is \emph{decidable}\index{subject}{decidable property}
\index{subject}{property!decidable} if there
is an algorithm which allows to verify it (see the Notes for a reference).
Otherwise, it is called \emph{undecidable}\index{subject}{undecidable property}
\index{subject}{property!undecidable}
Thus decidability is desirable, although not the end of things since
decidable properties can be very difficult to verify.

\subsection{Minimal dynamical system}

Given a topological dynamical system $(X,T)$,
a subset $Y$ of $X$ is \emph{stable}\index{subject}{stable! subset}
if $TY\subset Y$. The empty set and $X$ are always (trivial) stable sets.
As a stronger condition, the
 set $Y\subset X$ is \emph{invariant}\index{subject}{invariant!set}
if $T^{-1}Y=Y$.

A topological dynamical system is \emph{minimal}\index{subject}{minimal!dynamical system}
\index{subject}{dynamical system!minimal} if $X$ is the only closed and stable nonempty set. 
 A factor of a minimal dynamical system is minimal (Exercise~\ref{exerciseFactorMinimal}).

The following
characterization of minimal systems is sometimes used for definition.

\begin{proposition}\label{propositionMinimal}
A topological dynamical system is minimal if and only if
the positive orbit of every point
is dense in $X$.
\end{proposition}
\begin{proof}
Assume first that $(X,T)$ is minimal. For every $x\in X$, the
closure of the positive orbit of $x$ is a closed stable nonempty set
and thus it is equal to $X$. Conversely, let $Y\subset X$
be a closed stable nonempty set. For any $y\in Y$, the
closure of the positive orbit of $y$ is contained in $Y$. Since
it is equal to $X$, we conclude that $Y=X$.
\end{proof}
It follows directly from Propositions~\ref{propositionMinimal}
and \ref{propositionRecurrentSystem} that a minimal system
is recurrent.

The simplest example of a minimal dynamical system is a finite system
(also called a \emph{periodic system}\index{subject}{periodic!dynamical system}),
\index{subject}{finite!dynamical system}\index{subject}{dynamical system!finite}
\index{subject}{dynamical system!periodic} formed 
of a finite set $X=\{1,2,\ldots,n\}$ on which acts
a circular permutation $T$. A \emph{periodic point}
\index{subject}{periodic!point}\index{subject}{point!periodic}%
in  a dynamical system $(X,T)$ is a point $x\in X$ such that $T^nx=x$ for some
$n\ge 1$. A dynamical system $(X,T)$
is  called \emph{aperiodic}\index{subject}{dynamical system!aperiodic}
\index{subject}{aperiodic!dynamical system}%
if it does not contain any periodic point. 
A minimal system is aperiodic if and only if it is infinite.

As a second example, we find the rotations
of the circle.
\begin{example}\label{exampleRotations}
consider the \emph{unit circle}
\index{subject}{unit!circle} $\mathbb{S}^1 = \{ z \in \mathbb{C} \mid |z|=1 \}$
\index{symbols}{S@$\mathbb{S}^1$} and fix some $\lambda \in \mathbb{S}^1$. 
Let $R_\lambda$ be the map defined by  $R_\lambda ( z ) =  \lambda z$. 
Then $(\mathbb{S}^1 , R_\lambda )$ is a dynamical system.
It is minimal if and only if $\lambda$ is not a root of unity
(Exercise~\ref{exerciseIrrationalRotation}).
Otherwise it is a disjoint union of periodic systems. 

We have introduced above (Example~\ref{exampleInit2}) the
 \emph{torus}\index{subject}{torus} $\mathbb{T}$
\index{symbols}{T@$\mathbb{T}$} as
 the topological space $\mathbb{R} / \mathbb{Z}$. For $\alpha\in \R$,
let $T_\alpha:\T\to \T$ be the map $T_\alpha(x)=x+\alpha$
and set $\lambda=\exp(2i\pi\alpha)$. The
map $\phi:x\mapsto e^{2i\pi x}$ is a homeomorphism from
$\T$ onto $\mathbb{S}^1$ and $R_\lambda\circ\phi=\phi\circ T_\alpha$
because
\begin{displaymath}
R_\lambda\circ \phi(x)=R_\lambda e^{2i\pi x}=e^{2i\pi x}e^{2i\pi\alpha}=\phi\circ T_\alpha(x).
\end{displaymath}
  Consequently the topological dynamical systems $(\mathbb{S}^1  , R_\lambda )$ and $(\mathbb{T}  , T_\alpha )$ are  isomorphic,  The transformation $R_\alpha$ is called a \emph{rotation}
of angle $\alpha$. \index{rotation of the circle}.
The isomorphism is the map $\phi:x\mapsto \exp(2i\pi x)$ (see the diagram below).
\begin{displaymath}
\begin{CD}
\T@>{T_\alpha}>>\T\\
@VV{\phi} V   @VV{\phi} V\\
\mathbb{S}^1@>{R_\lambda}>>\mathbb{S}^1
\end{CD}
\end{displaymath}
\end{example}
\subsection{Induced systems}\label{sec:inducedsystems}
Let $(X,T)$ be a  minimal topological dynamical system and let $U$
be a nonempty clopen subset of $X$.  
Since $X$ is minimal, for every $x\in X$, there is an $n>0$
such that $T^nx\in U$ and we can define the integer
\begin{displaymath}
n(x)=\inf\{n>0\mid T^nx\in U\}
\end{displaymath}
called the \emph{entrance time}\index{subject}{entrance time} of $x$ in $U$.

Since $U$ is clopen, the function $x\mapsto n(x)$ is continuous.
Indeed, for each $n\ge 1$, the set of $x\in X$ such that $n(x)=n$
is $T^{-n}(U)\setminus\cup_{i=1}^{n-1}T^{-i}(U)$ which is open.
Thus, the map $x\mapsto n(x)$ is locally constant.
Consequently the map $T_U:U\rightarrow U$ defined by
\begin{displaymath}
T_U(x)=T^{n(x)}(x)
\end{displaymath}
\index{symbols}{T@$T_U$}%
is continuous being locally a power of $T$. It is called the \emph{induced transformation}
\index{subject}{induced!transformation} on $U$
and $(U,T_U)$ is called an \emph{induced system}
\index{subject}{induced!system}%
of $(X,T)$ or also a \emph{derivative} of $(X,T)$.
\index{subject}{derivative!of dynamical system}
\index{subject}{dynamical system!derivative}%
This system
 is minimal. Indeed, the orbit under $T$ of every point $x$ of $U$
is dense in $X$ and thus its orbit under $T_U$ is dense in $U$.

For $x\in U$, the integer $n(x)$ is called the \emph{return time}
\index{subject}{return!time}%
to $U$. For $x\in X$, the function
\begin{displaymath}
m(x)=\begin{cases}n(x)&\mbox{ if $x\notin U$}\\0&\mbox{otherwise}\end{cases}
\end{displaymath}
is called the \emph{waiting time}\index{subject}{waiting time}
to access $U$.

Note that the induced system can be defined even if $(X,T)$
is not minimal, provided the clopen set $U$ is such
that the return time $n(x)$ is  bounded on $U$.

The inverse operation can be described as follows. Let
$(X,T)$ be a topological dynamical system and let $f:X\to \N$
be a continuous function. Set
\begin{equation}
\hat{X}=\{(x,i)\mid x\in X, 0\le i<f(x)\}\label{equationTower}
\end{equation}
and define  a map $\hat{T}:\hat{X}\to\hat{X}$ by
\begin{displaymath}
\hat{T}(x,i)=\begin{cases}(x,i+1)&\mbox{ if $i+1<f(x)$}\\(Tx,0)\mbox{ otherwise}
\end{cases}\end{displaymath}
Then $(\hat{X},\hat{T})$ is a  topological dynamical system,
which is minimal if $(X,T)$ is minimal (see Exercise~\ref{exerciseTowerConstruction}).
The map $x\mapsto (x,0)$ identifies $X$ to the system induced
by $\hat{X}$ on $X\times\{0\}$.
The space $\hat{X}$ defined above is
said to be obtained from $X$ by the \emph{tower construction}\index{subject}{tower!construction}  relative to $f$
and the system $(\hat{X},\hat{T})$ is called a \emph{primitive}\index{subject}{primitive}
of $X$.

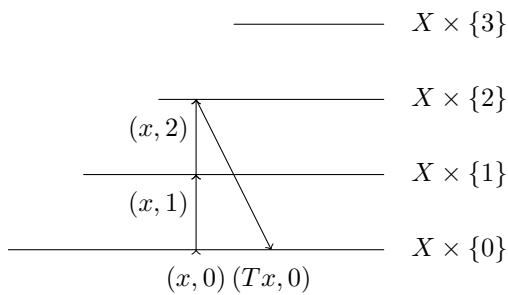
\begin{figure}[hbt]
  \centering
  \tikzset{node/.style={draw,minimum size=0.4cm,inner sep=0pt}}
	\tikzset{title/.style={minimum size=0.5cm,inner sep=0pt}}
  \begin{tikzpicture}
\draw(0,0)--(5,0);\node[title]at(6,0){$X\times\{0\}$};
\draw(1,1)--(5,1);\node[title]at(6,1){$X\times\{1\}$};
\draw(2,2)--(5,2);\node[title]at(6,2){$X\times\{2\}$};
\draw(3,3)--(5,3);\node[title]at(6,3){$X\times\{3\}$};
\node[title](x0)at(2.5,-.4){$(x,0)$};
\node[title](x1)at(2,.6){$(x,1)$};
\draw[->](2.5,0)edge node{}(2.5,1);
\node[title](x2)at(2,1.6){$(x,2)$};
\draw[->](2.5,1)edge node{}(2.5,2);
\draw[->](2.5,2)edge node{}(3.5,0);\node[title]at(3.5,-.4){$(Tx,0)$};
  \end{tikzpicture}
  \caption{The tower construction}\label{figureTowerConstruction}
\end{figure}
One can usefully figure out the result from the tower construction
as a superposition of `floors' above $X$. The number of floors
above $x\in X$ is $f(x)$. The transformation consists in accessing
the point above until there is no floor above and then using
the transformation defined on $X$(see Figure~\ref{figureTowerConstruction}).

%Again, the primitive of $(X,T)$ can be defined even if $X$ is not
%minimal, provided $f$ is continuous and bounded.

\begin{example}
Consider the system $(X,T)$ with $X=[0,1]$ and $T(x)=x+\alpha\bmod 1$
of Example~\ref{exampleInit} for $1/2<\alpha<1$. The
transformation $T_U$ induced by $T$ on $U=[0,\alpha [$ is 
given by 
\begin{displaymath}
T_U (x)=\begin{cases} x+2\alpha-1&\mbox{ if $0\leq x<1-\alpha$}\\x+\alpha-1&\mbox{ otherwise }\end{cases}\end{displaymath}
Thus, $(U,T_U)$ is isomorphic to the system $(X,T_{2-(1/\alpha)})$
via the conjugacy $x\mapsto x/\alpha$.
\end{example}
%%%%%%%%%%%%%%%%%%%%%%%%
\subsection{Dynamical systems on Cantor spaces}\label{sectionCantorSpaces}
%%%%%%%%%%%%%%%%%%%%%%%%%%

We have defined a Cantor space as a totally disconnected
compact metric space without isolated points.
Recall that a topological space is called totally disconnected
if every connected component is reduced to a point
(see Appendix~\ref{appendixTopo}).

Let us give classical examples of Cantor spaces. 
There is first the original one, as defined originally by Cantor.  
The \emph{Cantor set}\index{subject}{Cantor!set}
 is the subset $C$ of $[0,1]$ defined by $C = \cap_n C_n$ where $C_0 = [0,1]$ and 
$$
C_{n+1} = \frac{C_n}{3} \cup \left( \frac{2}{3} + \frac{C_n}{3} \right) .
$$

It is a Cantor space for the induced topology.

There is the abstract topological one, that we used for the definition.
It is a compact metric space which is totally disconnected without isolated point or, equivalently, a topological space with a countable basis of 
the topology consisting of clopen sets
 and without isolated points (see Appendix~\ref{appendixTopo}).

There is the algebraic one, the ring $\Z_p$ of $p$-\emph{adic integers}%
\index{subject}{p-adic@$p$-adic!integer}, where $p$ is a prime number, 
for the $p$-\emph{adic topology}%
\index{subject}{p-adic@$p$-adic!topology} which is induced by the distance 
defined for $x\ne y$ by 
$d(x,y)=p^{-n}$ if $p^n$ is the higher power of $p$ dividing
$x-y$. 
An important example of a minimal Cantor system
 is the \emph{odometer}\index{subject}{odometer} $(\Z_p , T)$
\index{symbols}{Z@$\Z_p$} where $T$ is the 
transformation defined by  $T ( x) = x+1$.

Finally, there is the symbolic one that we have already seen and
consists in considering the set $A^\Z$ of two-sided infinite sequences
on a finite alphabet $A$.

From a topological point of view all these topological spaces are the same as they are homeomorphic (see the notes section for a reference to a proof).

A \emph{Cantor system}\index{subject}{dynamical system!Cantor}\index{subject}{Cantor!dynamical system}  is a dynamical system $(X , T)$ where $X$ is a $T$-stable Cantor set. 
The odometers are Cantor systems.
A \emph{symbolic system}\index{subject}{dynamical system!symbolic}\index{subject}{symbolic!dynamical system} is a dynamical system $(X , T)$ where $X$ is a $T$-stable closed subset of $A^\Z$. The transformation $T$
need not be the shift (see the example of the odometer).
Observe that symbolic systems are not necessarily Cantor systems since
 closed subsets could have isolated points.
But infinite minimal symbolic systems are Cantor systems
(Exercise~\ref{exerciseMinimalSymbolicCantor}).

The pair $(A^\Z,S)$ is called the \emph{full shift}
\index{subject}{full shift}%
 (on the alphabet $A$).
If $X$ is a closed shift-invariant subset of $A^\Z$, then the topological dynamical system $(X,T)$, where $T$ is the restriction of $S$ to $X$, 
is called a \emph{subshift}
\index{subject}{subshift}\index{subject}{shift space!subshift of}%
 of the full shift on the alphabet $A$ or a {\em shift space}
\index{subject}{shift space}%
or a \emph{two-sided shift space}.\index{subject}{two-sided!shift space}
Thus, a shift space is a symbolic system.
For the full shifts on any finite alphabet and in general for shift spaces, 
the transformation will usually be denoted by $S$. Thus, we will often
use the notation $X$ instead of $(X,S)$ for a shift space.

Similarly one can define the \emph{one-sided full shift}
\index{subject}{full shift!one-sided}%
\index{subject}{one-sided!full shift} as follows.
We still denote by $S$ the one-sided transformation,
called the \emph{one-sided shift}\index{subject}{shift transformation!one-sided}
 which is defined for $x\in A^\mathbb{N}$ by $y=Sx$
if $y_n=x_{n+1}$, as in Equation~\eqref{eqShift}.
The pair $(A^\N , S)$ is called the one-sided full shift on the alphabet $A$.
If $X$ is a closed $S$-invariant subset of $A^\N$, the pair $(X,S)$ is called a \emph{one-sided shift space}
\index{subject}{one-sided!shift space}\index{subject}{shift space!one-sided}%
 or \emph{one-sided subshift}.
\index{subject}{one-sided!subshift}\index{subject}{subshift!one-sided}%

To every shift space $(X,S)$, one may associate a one-sided
shift space $(Y,S)$ called its \emph{associated one-sided shift space}%
\index{subject}{one-sided!shift!associated}
by considering the set $Y$
of $y_0y_1\cdots $ such that $y_i=x_i$ ($i\ge 0$) for some $x=(x_n)_{n\in\Z}$.
The map $\theta:X\to Y$ defined by $\theta(x)=y$ is a surjective morphism.

Conversely, for every one-sided shift space $(Y,S)$ there
is a unique two-sided shift space $(X,S)$ such that $(Y,S)$ is associated
to $(X,S)$. Indeed, the set $X$ of sequences $x\in A^\Z$ such that
$x_nx_{n+1}\cdots$ belongs to $Y$ for every $n\in\Z$
is closed and shift-invariant. 
It is clear that $(X,S)$ is the unique shift space
such that $Y=\theta(X)$. We also say that
$X$ is the two-sided shift space \emph{associated}
\index{subject}{two-sided!shift space!associated} to $Y$.

This shows that one-sided and two-sided shift spaces are closely related
objects. In general, in this book, we consider two-sided shift
spaces rather than one-sided shift spaces because it is often convenient
to have a transformation which is invertible. In general,
by a shift space, we mean a two-sided shift space.

\begin{example}\label{exampleGolden}
The \emph{golden mean shift}
\index{subject}{golden mean!shift}%
\index{subject}{shift space!golden mean} is the set $X$
of two-sided sequences on $A=\{a,b\}$ with no consecutive $b$.
Thus $X$ is the set of labels of two-sided infinite paths in the graph of Figure
\ref{figureGolden}. It is recurrent, as
one may easily verify. It is not minimal since it contains
the one-point set $\{a^\Z\}$ which is closed and shift invariant.
\begin{figure}[hbt]
\centering
%\gasset{Nadjust=wh}
%\begin{picture}(20,20)(0,-3)
%\node(1)(0,0){$1$}\node(2)(20,0){$2$}

%\drawloop[loopangle=180](1){$a$}
%\drawedge[curvedepth=3](1,2){$b$}
%\drawedge[curvedepth=3](2,1){$a$}
%\end{picture}
%Figure en tikz
\tikzset{node/.style={circle,draw,minimum size=0.4cm,inner sep=0pt}}
	\tikzset{title/.style={minimum size=0.5cm,inner sep=0pt}}
        \tikzstyle{every loop}=[->,shorten >=1pt,looseness=12]
        \tikzstyle{loop left}=[in=130,out=220,loop]
        %\tikzset{state/ .style={minimum size=0.5cm,inner sep=0pt}}
	\begin{tikzpicture}
	\node[node](1) {$1$};
        %\node(E)[title]($1$)
	\node[node](2) [right = 2cm of 1] {$2$};

	\draw[bend left, ->,>=stealth] (1) edge node {} (2);
        \node[title](B)[above right = 0.3cm and 0.9cm of 1]{$b$};
        \draw[bend left, ->,>=stealth] (2) edge node {} (1);
        \node[title](A)[below right = 0.2cm and 0.9cm of 1]{$a$};
        %\draw [loop left] (1)  node {$a$} (1);
        \draw (1) edge [loop left,>=stealth] node {} (1);
        \node [title][left = 0.6cm of 1]{$a$};
	\end{tikzpicture}
\caption{The golden mean shift.}\label{figureGolden}
\end{figure}
\end{example}
\begin{example}\label{exampleFibonacci00}
Let $\varphi:a\mapsto ab,b\mapsto a$ be the \emph{Fibonacci morphism}.
\index{subject}{Fibonacci!substitution}%
\index{subject}{substitution!Fibonacci}%
Since $\varphi(a)$ begins with $a$, any $\varphi^n(a)$
is a prefix of $\varphi^{n+1}(a)$. Let $x\in \{a,b\}^\N$
be the sequence $x$ having all $\varphi^n(a)$ as prefixes.
Thus
\begin{displaymath}
x=abaababa\cdots
\end{displaymath}
It is known as the \emph{Fibonacci word}
(we prefer to keep the name `Fibonacci sequence' for the
well-known sequence $F_{n+1}=|\varphi^n(b)|$).
\index{subject}{Fibonacci!sequence} The subshift
of $\{a,b\}^\N$ which is the closure of the orbit $x$ is the \emph{one-sided Fibonacci
shift}.\index{Fibonacci!one-sided shift}%
We will see that it is minimal (Example~\ref{exampleFibonacci2}).
\end{example}
%%%%%%%%%%%%%%%%%%%%%%%%
\section{More on shift spaces}
%%%%%%%%%%%%%%%%%%%%%%%%%%
\label{sectionSymbolicSystems}
In this section, we develop in more detail the notions related
to shift spaces and their language. We will see how the notions
of recurrence and minimality can be expressed adequately for shift spaces.
We will also introduce important notions like return words or Rauzy graphs.
%%%%%%%%%%%%%%%%%%%%%%%%
\subsection{Some combinatorics on words}\label{sectionCombinatoricsonWords}

Let $A$ be a nonempty set called an \emph{alphabet}\index{subject}{alphabet}.
We will generally assume that the alphabet is finite.
A {\em word}\index{subject}{word}
 over $A$ is an element of the \emph{free monoid}\index{subject}{monoid!free}
\index{subject}{free!monoid} generated by $A$, denoted by $A^*$. 
\index{symbols}{A@$A^*$}%
If $u = u_0u_1 \cdots u_{n-1}$ (with $u_i\in A$, $0\leq i\leq n-1$) is a word,
 its {\em length}\index{subject}{word!length}\index{subject}{length!of word}
 is $n$ and is denoted by $|u|$. 
 \index{symbols}{u@$\lvert u\rvert$}%
 For $a\in A$, we denote by $|u|_a$ the number of occurrences
 of the letter $a$ in $u$.\index{symbols}{u@$\lvert u\rvert_a$}
 
The {\em empty word}\index{subject}{word!empty}
 is denoted by $\varepsilon$.
\index{symbols}{e@$\varepsilon$} It is the unique word of length $0$. 
The set of nonempty words over $A$,
called the \emph{free semigroup}\index{subject}{semigroup!free}
on $A$ is denoted by $A^+$.

We consider the free monoid $A^*$ as embedded in the free
group on $A$ (see Appendix~\ref{appendixGroups}). Consequently,
when $u=vw$, we also write $v=uw^{-1}$ and $w=v^{-1}u$.

A   {\em factor}\index{subject}{factor!of word}\index{subject}{word!factor of}
(also called a {\em subword}\index{subject}{subword}\index{word!subword of}
or a \emph{block}\index{subject}{block})\index{word!block of}
 of a word $u$ is a finite word $y$ such that there exist two words $v$ and $w$ satisfying $u = vyw$.
 When $v$ (resp. $w$) is the empty word, we say that $y$ is a {\em prefix} \index{subject}{prefix!of word}\index{subject}{word!prefix of}(resp. {\em suffix})\index{subject}{suffix!of word}\index{subject}{word!suffix of} of $u$.
 A factor $y$ (resp. prefix, resp. suffix) of a word $u$ is
 \emph{proper}\index{subject}{proper!factor}\index{subject}{proper!prefix}
 \index{subject}{proper!suffix}\index{subject}{factor!of word!proper}%
 \index{subject}{prefix!of word!proper}\index{subject}{suffix!of word!proper}%
 if $y\ne u$.

Two words $u,v$ are \emph{conjugate}\index{subject}{conjugate!words}%
\index{subject}{word!conjugate}\index{subject}{conjugacy!of words} if $u=rs$ and $v=sr$ for some
words $r,s$ or, equivalently, if $v$ is obtained from $u$ by a circular
permutation of its letters. Conjugacy is an equivalence relation on words.

A word $w$ is \emph{primitive} if it is not a power of another word.
Formally, $w$ is primitive if $w=u^n$ implies $n=1$. A primitive word
of length $n$ has $n$ distinct conjugates (Exercise~\ref{exercisePrimitive}).
Any nonempty word $w$ can be written uniquely $w=u^n$ with $u$ primitive and $n\ge 1$. The integer $n$ is called the
\emph{exponent}\index{subject}{exponent of word} of $w$.

A word $w=a_1a_2\cdots a_n$ with $a_i\in A$ has \emph{period}
\index{subject}{period!of word} $p$ if $a_i=a_{i+p}$ for $1\le i\le n-p$.
In this case, $w$ has a prefix of exponent $\lfloor n/p\rfloor$.
Indeed, one has $w=u^qv$ where $n=pq+r$ with $0\le r<p$, $|u|=p$
and $|v|=r$.

\begin{figure}[hbt]
  \centering
\tikzset{node/.style={circle,draw,minimum size=0.1cm,inner sep=0pt}}
\tikzset{title/.style={minimum size=0.5cm,inner sep=0pt}}
\begin{tikzpicture}
  \node[node](r)at(0,1){};\node[node](uh)at(1,1){};\node[node](s)at(3,1){};
  \node[node](ss)at(4,1){};
  \node[node](r')at(0,0){};\node[node](ub)at(1.5,0){};\node[node](s')at(3.5,0){};
  \node[node](ss')at(4.5,0){};

  \draw[above](r)edge node{$r$}(uh);\draw[above](uh)edge node{$u$}(s);
  \draw[above](s)edge node{$s$}(ss);
  \draw[above](r')edge node{$r'$}(ub);\draw[above](ub)edge node{$u$}(s');
  \draw[above](s')edge node{$s'$}(ss');
\end{tikzpicture}
\caption{Two overlapping occurrences of $u$.}\label{figureOverlap}
  \end{figure}
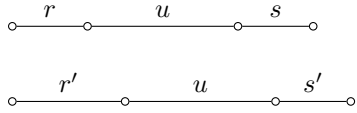

There is an important connection between the periods of a word
and the overlap of its factors.
\begin{proposition}\label{propositionOverlappingFactor}
  If a word $w$
has two overlapping occurrences in $u$, that is, if
$w=rus=r'us'$ with $|r|\le|r'|\le|r|+|u|$ Then  $u$ has period $|r|-|r'|$.
\end{proposition}
The proof is straightfoward (see Figure~\ref{figureOverlap}).

 A set of words
 on the alphabet $A$ is also called a \emph{language}\index{subject}{language}
 on $A$. If $U$ is a language on $A$, we denote by
 $U^*$\index{symbols}{U@$U^*$}%  
  the submonoid of $A^*$ generated by $U$, that is
\begin{equation}
   U^*=\{w\in A^*\mid w=u_1u_2\cdots u_n, u_i\in U, n\ge 0\}.\label{equationStar}
   \end{equation}
When $U=\{u\}$, we denote simply $u^*$ instead of $\{u\}^*$.
Thus $u^*=\{u^n\mid n\ge 0\}$. For $U,V\subset A^*$, we also denote
 \begin{equation}
   UV=\{uv\mid u\in U,v\in V\}.\label{equationProduct}
   \end{equation}
 
If $k,l$ are integers such that $0 \leq k \leq l < |u|$, we let $u_{[k,l]}$ denote the subword $u_k u_{k+1} \cdots u_l$ of $u$.
We define $u_{[k,l+1)}$ to be $u_{[k,l]}$.
If $l < k$, then $u_{[k,l]}$ is the empty word.
If $y$ is a factor of $u$, the {\em occurrences}\index{subject}{occurrence of factor}
\index{subject}{factor!occurrence of} of $y$ in $u$ are the integers $i$ such that $u_{[i,i + |y| - 1]}= y$. 
If $y$ has an occurrence in $u$, we also say that $y$ {\em occurs} in $u$.

The \emph{reversal}
\index{subject}{reversal!of word}\index{subject}{word!reversal}%
of a word $u=u_0u_1\cdots u_n$, with $u_i\in A$,
is the word $\tilde{u}=u_n\cdots u_1u_0$.
\index{symbols}{u@$\tilde{u}$}%
 The reversal of a set $U$ of words is the
set $\tilde{U}=\{\tilde{u}\mid u\in U\}$.

The elements of $A^{\mathbb{K}}$, where $\mathbb{K}$ is equal to $\mathbb{N}$,
$-\N$ or $\mathbb{Z}$, are called {\em sequences} or \emph{infinite words}.
\index{subject}{word!infinite}%
When we need to precise which kind of sequences we are dealing with we sometimes say {\em right infinite sequence},\index{subject}{right!infinite sequence}
\index{subject}{sequence!right-infinite}%
or \emph{one-sided sequence}\index{subject}{sequence!one-sided}
when $K=\N$,
\emph{left-infinite sequence}\index{subject}{sequence!left-infinite}
when $K=-\N$  and {\em two-sided sequence}
(or also \emph{bi-infinite sequence})
\index{subject}{bi-infinite!sequence}\index{subject}{sequence!bi-infinite}%
\index{subject}{two-sided!sequence}\index{subject}{sequence!two-sided}%
 in the last case. 
For $x = (x_n)_{n \in \mathbb{Z}} \in A^\mathbb{Z}$, we let $x^+$ and $x^-$ respectively denote the sequences $(x_n)_{n \geq 0}$ and $(x_n)_{n < 0}$.
\index{symbols}{x@$x^+$}\index{symbols}{x@$x^-$}%
For $x\in A^{-\N}$ and $y\in A^{\N}$, we denote $z=x\cdot y$ the
two sided sequence $z$ such that $x=z^-$ and $y=z^+$.
The notion of \emph{factor}\index{subject}{factor!of sequence}
 is naturally extended to sequences, as well as the notion of \emph{prefix}
 \index{subject}{prefix!of sequence} when $\mathbb{K} = \N$.
 
The set of subwords of length $n$ of $x$ is written $\mathcal{L}_n(x)$ and the set of subwords of $x$, called the {\em language}
\index{subject}{language!of a sequence}%
 of $x$, is denoted by $\mathcal{L}(x)$. We also denote
 by $\cL_{\le n}(x)$ the set of words of length at most $n$ in $\cL(x)$.

 The \emph{reversal}\index{subject}{reversal!of sequence}
 of a right-infinite sequence $x=x_0x_1\cdots$ is the
 left-infinite sequence $\tilde{x}=\cdots x_1x_0$.
 \index{symbols}{x@$\tilde{x}$}%
 The language of $\tilde{x}$ is the reversal of $\cL(x)$.

For a finite word $u\in A^+$, we denote by $u^\omega$
\index{symbols}{u@$u^\omega$}
the right infinite sequence $uuu\cdots$ and by $u^\infty$
\index{symbols}{u@$u^\infty$} the two-sided
infinite sequence $x=\cdots uuu\cdot uuu\cdots$, where
the dot is placed to the left of $x_0$. In this way, we have 
$x^+=u^\omega$.

An integer $p\ge 1$ is a \emph{period} of a sequence $x\in A^\N$ if
$x_i=x_{i+p}$ for all $i\ge 0$.\index{subject}{period!of sequence}
Clearly, the sequence $x$ has period $p$ if and only
if $x=u^\omega$ with $p=|u|$.
\begin{proposition}\label{propositionMinimalPeriod}
Every sequence $x$ has a unique minimal period
\index{subject}{minimal!period of sequence} and its multiples
form the set of periods of $x$. Moreover, $x$ has  minimal period $p$
if and only if $x=u^\omega$ with $p=|u|$ and $u$ primitive.
\end{proposition}
\begin{proof}
The set of $p\in\Z$
such that $|p|$ is a period of $x$ is a subgroup of $\Z$.
Indeed, if $p,q$ are periods of $x$ with $p\le q$, then
$x_{i+q-p}=x_{(i+q-p)+p}=x_{i+q}=x_i$ for every $i\ge 0$ and
thus $q-p$ is also a period of $x$.
This shows that the set of periods of $x$ coincides with
the set of multiples of the minimal period. Next, the map
$u\mapsto |u|$ sends the set of words $u$ such that $x=u^\omega$
onto the set of periods of $x$. This proves the second statement.
\end{proof}
A refinement of this statement, called Fine-Wilf Theorem,
is given in Exercise~\ref{exerciseFineWilf}.
These notions carry easily to two-sided sequences.

The sequence $(p_n(x))_{n\ge 0}$\index{symbols}{pn(x)@$p_n(x)$}
defined by  $p_n(x)= \Card(\mathcal{L}_n(x))$
is  the {\em factor complexity}
\index{subject}{factor!complexity}%
(or \emph{word complexity}\index{subject}{word!complexity}
or simply the \emph{complexity}\index{subject}{complexity of sequence}) of $x$. 
Note that $p_0(x)=1$, that
$p_n(n)\le p_{n+1}(x)$ and that $p_{n+m}(x)\le p_n(x)p_m(x)$ for all $n,m\ge 0$.

\begin{example}\label{exampleFibonacciWordComplexity}
  Let $x$ be the Fibonacci word (Example~\ref{exampleFibonacci00}). We have $p_1(x)=2$
  since every letter $a,b$ appears in $x$. Next
  $p_2(x)=3$ since $\cL_2(x)=\{aa,ab,ba\}$ as one may verify. We will see that
  actually, one has $p_n(x)=n+1$ for all $n\ge 1$ (see Section~\ref{sectionSturmianShifts}).
\end{example}

The sequence $x\in A^\mathbb{N}$ is {\em eventually periodic}\index{subject}{eventually!periodic sequence}
 if there exist a word $u$ and a nonempty word $v$ such that $x=uv^{\omega}$, where $v^{\omega}= vvv\cdots$.
 A sequence that is not eventually periodic is called {\em aperiodic}.
 \index{subject}{aperiodic!sequence}\index{subject}{sequence!aperiodic}
%In this case $|v|$ is called a {\em length period} of $x$. 
 It is {\em periodic}\index{subject}{periodic!sequence}
 \index{subject}{sequence!periodic} if $u$ is the empty word.

The following result is classical.\index{subject}{Theorem!Morse-Hedlund}
\index{subject}{Morse-Hedlund!Theorem}%
\begin{theorem}[Morse, Hedlund]\label{theoremCovenHedlundSequence}
\index{names}{Morse, Marston}\index{names}{Hedlund, Gustav A.}
Let $x$ be a two-sided sequence. The following conditions are equivalent.
\begin{enumerate}
\item[\rm(i)] For some $n\ge 1$, one has $p_n(x)\le n$.
\item[\rm(ii)] For some $n\ge 1$, one has $p_n(x)=p_{n+1}(x)$.
\item[\rm(iii)] $x$ is periodic. 
\end{enumerate}
Morover, in this case, the least period of $x$ is $\max p_n(x)$.
\end{theorem}
\begin{proof}
(i) $\Rightarrow$ (ii).
Since $p_n(x)\le p_{n+1}(x)$ for all $n\ge 0$, the hypothesis implies
that $p_n(x)=p_{n+1}(x)$ for some $n\ge 0$. 

(ii) $\Rightarrow$ (iii). For every
$w\in\cL_n(x)$, there is a unique letter $a\in A$ such that
$wa\in\cL_{n+1}(x)$. This implies that two consecutive
occurrences of a word $u$ of length $n$ in $x$ are separated
by a fixed word depending only on $u$ and thus that $x$ is periodic.

(iii) $\Rightarrow$ (i) is obvious.

Let $n$ be the least period of $x$.
Since a primitive word of length $n$ has $n$ distinct conjugates, we have
$p_n(x)=n$ and $p_m(x)=n$ for all $m\ge n$.
\end{proof}
Thus, by Proposition~\ref{theoremCovenHedlundSequence},
 either $p_n(x)\ge n+1$ for all $n\ge 1$ or $p_n(x)$
is eventually constant. The case $p_n(x)=n+1$ for all $n\ge 1$
corresponds to the Sturmian sequences (see below).

Note that for a one-sided sequence, the same result holds
with condition (iii) replaced by the condition that $x$ is
eventually periodic. The proof is the same.

\begin{proposition}\label{propositionRecurrentSequence}
The following conditions are equivalent for $x\in A^\N$.
\begin{enumerate}
\item[\rm(i)] Every $u\in\cL(x)$ has at least two occurrences in $x$.
\item[\rm(ii)] Every $u\in\cL(x)$ has an infinite number of occurrences in $x$.
\item[\rm(iii)] For every $u,w\in\cL(x)$ there is $v\in\cL(x)$
such that $uvw\in\cL(x)$.
\end{enumerate}
\end{proposition}
\begin{proof}
(i) $\Rightarrow$ (ii). Assume that $u$ has a finite number of occurrences
in $x$. Let $v$ be the shortest prefix of $x$ containing all these occurrences.
Since $v$ has a second occurrence in $x$, we have a contradiction.

(ii) $\Rightarrow$ (iii). Assume that $u=x_{[i,j)}$.
Since $w$ has an infinite number of occurrences
in $x$, there is an index $k$ larger than $j$ such that
$w=x_{[k,\ell)}$. Set
$v=x_{[j,k)}$. Then $uvw=x_{[i,\ell)}$.

(iii) $\Rightarrow$ (i) is clear considering $u=w$.
\end{proof}
A word $u$ is {\em recurrent}\index{subject}{recurrent!word in sequence}
 in $x\in A^\N$ if condition (ii) above is satisfied. The sequence
$x$ itself is called \emph{recurrent}\index{subject}{recurrent!sequence}
if one of the conditions is satisfied. Thus a sequence
$x\in A^\N$ is recurrent if and only if $x$ is a recurrent
point of the full shift. We could of course use
Proposition~\ref{propositionRecurrentSystem} to prove Proposition~\ref{propositionRecurrentSequence}.
We also say that
the language $\cL(x)$ is %\emph{recurrent}\index{subject}{recurrent!language}
%\index{subject}{language!recurrent}%
%or
\emph{irreducible}\index{subject}{irreducible!language}
\index{subject}{language!irreducible}%
if condition (iii) is satisfied.

\begin{example}\label{exampleChampernowne}
Let $A=\{a,b\}$ and let $x=abaaabbabb\cdots$ be the sequence formed
of all words on $A$ in radix order (that is, ordered first by length,
then lexicographically). It is a recurrent sequence in which all
words on $A$ appear, that is, such that $\cL(x)=A^*$.
As a variant of this example, the \emph{Champernowne sequence}
\index{subject}{Champernowne!sequence}
\index{names}{Champernowne, David G.} is the sequence
$x=01234567891011121314151617181920\cdots$ formed of the decimal
representation of all numbers
in increasing order.
\end{example}

\begin{proposition}\label{propositionUniformlyRecurrentSequence}
The following conditions are equivalent for a sequence $x\in A^\N$.
\begin{enumerate}
\item[\rm(i)] Every $u\in\cL(x)$ occurs infinitely often in $x$ and the greatest
difference of two successive occurrences of $u$ is bounded.
\item[\rm(ii)] For every $u\in\cL(x)$, there is an $n\ge 1$
such that $u$ occurs in every word of $\cL_n(x)$.
\end{enumerate}
\end{proposition}
\begin{proof}
(i) $\Rightarrow$ (ii). Let $k$ be the maximum of the differences between
successive occurrences of $u$. Then $u$ appears in every word of $\cL_n(x)$
for $n=|u|+k$.

(ii) $\Rightarrow$ (i) is clear.
\end{proof}
A sequence $x\in A^\N$ is {\em uniformly recurrent}
\index{subject}{uniformly!recurrent!sequence}\index{subject}{recurrent!uniformly!sequence}%
if one of these conditions hold. We also say in this case that $\cL(x)$ is 
\emph{uniformly recurrent}.\index{subject}{uniformly!recurrent!language}\index{subject}{recurrent!uniformly!language}
\index{subject}{language!uniformly recurrent}%

\begin{example}
The Fibonacci word $x$ (Example~\ref{exampleFibonacciWordComplexity})
is uniformly recurrent. Indeed, let $u\in\cL(x)$ and let $n$
be the minimal integer such that $u$ is a factor of $\varphi^n(a)$.
Then $\varphi^{n+2}(a)=\varphi^n(a)\varphi^n(b)\varphi^n(a)$
and thus $u$ has a second occurrence in $\varphi^{n+2}(a)$
at bouded distance of the first one, which implies that it
occurs infinitely often at bounded distance.
\end{example}
 
A sequence $x\in A^\N$ is {\em linearly recurrent}
\index{subject}{linearly!recurrent!sequence}\index{subject}{recurrent!linearly!sequence}
with \emph{constant} $K$
\index{subject}{constant!of linear recurrence} if it is recurrent and  the greatest difference between successive occurrences
of $u$ is bounded by $K|u|$.

Most of this terminology extends naturally to a two-sided infinite sequence.
In particular, $x\in A^\Z$ is \emph{periodic} 
\index{subject}{periodic!two-sided sequence} if $x=v^\infty$ 
for some $v\in A^+$.
In this case, $x^+$ is periodic. Similarly, $x\in A^\Z$ is 
\emph{recurrent}
(resp. \emph{uniformly recurrent})
\index{subject}{uniformly!recurrent!two-sided sequence}\index{subject}{recurrent!uniformly!two-sided sequence}%
if $\cL(x)$ is  recurrent (resp. uniformly recurrent).
The same extension holds for linearly recurrent sequences.

%Let $A$ and $B$ be finite or infinite alphabets. By a {\em morphism}\index{subject}{morphism} from $A^*$ to $B^*$ we mean a homomorphism of free monoids.
\begin{proposition}\label{propositionLRhasLinearComplexity}
Let $x$ be  a two-sided sequence which is linearly recurrent with constant $K$. Then
\begin{enumerate}
\item Every word of $\cL_n(x)$ appears in every word of 
$\cL_{(K+1)n-1}(x)$.
\item The factor complexity of $x$ is at most $Kn$.
\end{enumerate}
Moreover, if $x$ is not periodic, it is $(K+1)$-power free, that is
for every nonempty word $u\in\cL(x)$, $u^n\in\cL(x)$ implies
$n\le K$.
\end{proposition}
\begin{proof}
1. Let $u\in\cL_n(x)$. Since two successive occurrences of $u$
differ by at most $Kn$, every word of $\cL_{(K+1)n-1}$ contains $u$
as a factor.

2. Set $p_n(x)=\Card(\cL_n(x))$. A word of length $(K+1)n-1$
has at most $Kn$ factors of length $n$. Thus, by Assertion 1,
$p_n(x)\le Kn$.

Assume now that $u^{K+1}\in\cL(x)$. Set $n=|u|$. The word
$u^{K+1}$ has length $(K+1)n$ and at most $n$ factors of length $n$.
By Assertion 1, this implies $p_n(x)\le n$. 
By Theorem~\ref{theoremCovenHedlundSequence}, this implies that $x$ is periodic.
\end{proof}
Clearly, a linearly recurrent sequence is uniformly recurrent
but the converse is not true (see Exercise~\ref{exerciseSturmLR}).

\begin{example}
The Fibonacci word is linearly recurrent (Exercise~\ref{exerciseFibonacciLR}).
\end{example}
%%%%%%%%%%%%%%%%%%%%%%%%%%%%%
\subsection{The language of a shift space}

%Let $A$ be a finite alphabet.
%The set $A^\Z$  is
%a Cantor space for the product topology which is metrizable for the distance
%$d(x,y)=2^{-r(x,y)}$ where $r(x,y)=\inf\{|n| \mid x_n \not =y_n \}$
%(we set $d(x,y)=0$ if $x=y$).

Let $X$ be a shift space.
The {\em language} of $X$
\index{subject}{language!of a shift space}%
\index{subject}{shift space!language of}%
 is the set $\mathcal{L} (X) $ \index{symbols}{l@$\cL(X)$}
 of subwords of elements belonging to $X$.
The set $\cL(X)$ is of course the union of the languages $\cL(x)$
for $x\in X$. We also denote by $\cL_n(X)$\index{symbols}{l@$\cL_n(X)$}
 the set of words
of length $n$ in $\cL(X)$ and by $\cL_{\le n}(X)$ the set of those
of length at most $n$.

The same notation can be used for the language of a one-sided shift
space. The languages of a two-sided shift space and of its
associated one-sided shift space are actually the same.

A set $L$ of words on the alphabet $A$ is \emph{factorial}
\index{subject}{factorial!language}%
if it contains the factors of its elements. A word
$u\in L$ is \emph{extendable}\index{subject}{extendable!word}
\index{subject}{word!extendable}%
 in $L$ if there
are letters $a,b\in A$ such that $aub$ belongs to  $L$.
The language $L$ is said to be \emph{extendable}
\index{subject}{extendable!set of words} if for every $u\in L$
is extendable. The language of a shift
space is factorial and extendable and, conversely, for every factorial
extendable set $L$, there is a unique shift space $X$ such that
$\cL(X)=L$ (Exercise~\ref{exerciseFactorialExtendable}).

For two words $u,v$ such that $uv\in \mathcal{L}(X)$, the set 
\begin{displaymath}
[u\cdot v]_X=\{x\in X\mid x_{[-|u|,|v|-1]}=uv\}
\end{displaymath}
\index{symbols}{u@$[u\cdot v]$}%
is nonempty. It is called the \emph{cylinder}\index{subject}{cylinder}
with basis $(u,v)$. 
We set $[v]_X= [\varepsilon\cdot v]_X$,
where $\varepsilon$ is the empty word,
or equivalently $[v]_X=\{x\in X\mid x_{[0,|v|-1]}=v\}$. 
\index{symbols}{u@$[u]$}%
Any cylinder is open and every 
open set in a shift space is a union of cylinders.
The clopen sets in $X$ are the finite
unions of cylinders. 

For any sequence $x\in A^\Z$ there is a smallest shift space
containing $x$ called the \emph{subshift generated}
\index{subject}{subshift!generated!by a sequence}%
\index{subject}{shift space!generated!by a sequence}%
 by $x$ and denoted $\Omega(x)$.
\index{symbols}{Omega@$\Omega(x)$} It is the closure of the orbit of $x$.
%The same notation can be used for the one-sided shift space $\Omega(x)$
%generated by a one-sided sequence $x$.
For example, if $x=u^\infty$ is a periodic sequence,
the shift $\Omega(x)$ is periodic (that is, it is a periodic dynamical system).
\index{subject}{periodic!shift}\index{subject}{shift space!periodic}%

The following property is sometimes  taken
for definition of shift spaces.
\begin{proposition} A set $X\subset A^\Z$ is a shift space if and only if there
is a set $F\subset A^*$ of finite words such that $X$ is the set
$X_F$ of infinite words without factor in $F$.
\end{proposition}
\begin{proof} Indeed, such a set
is clearly closed and invariant by the shift. Conversely, let $X$
be a shift space. Since $X$ is closed, its complement
$Y=A^\Z\setminus X$ is open. Thus for every $y\in Y$, there
is a cylinder $[u_y\cdot v_y]_{A^\Z}$ containing $y$
and contained in $Y$. Set $F=\{u_yv_y\mid y\in Y\}$.
Then $X_F\subset X$ since $X_F\cap Y=\emptyset$. Conversely, if $x\in X$
has a factor  $u_yv_y$ in $F$, then $T^nx$ is in $[u_y\cdot v_y]_{A^\Z}$ for some
$n\in \Z$ and this is a contradiction since $X$ is shift invariant.
Thus $X=X_F$.
\end{proof}
Let $X$ and $Y$ be shift spaces on alphabets $A,B$ respectively.
Given nonnegative integers $m,n$ such that $-m\le n$, let $f:\cL_{m+n+1}(X)\to B$
be a map called a \emph{block map}.
\index{subject}{block!map}%
We call $m$ the \emph{memory}\index{subject}{memory of block map}
and $n$ the \emph{anticipation}\index{subject}{anticipation of block map}
of the block map $f$.

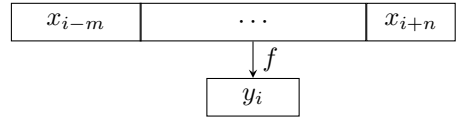
\begin{figure}[hbt]
\centering
%\gasset{Nmr=0,Nw=10,Nh=4}
%\begin{picture}(50,10)

%\node(k)(20,10){$x_{n}$}
%\node(k+1)(30,10){$x_{n+1}$}\node(k+n)(52,10){$x_{n+k-1}$}
%\node[Nw=12](d)(41,10){$\cdots$}
%\node(y)(20,0){$y_n$}
%\drawedge(k,y){$f$}
%\end{picture}
\tikzset{node/.style={rectangle,draw,minimum width=.6cm,minimum height=.5cm}}
\tikzstyle{node2}=[rectangle,draw,minimum width=3cm,minimum height=.5cm]
	\tikzset{title/.style={minimum size=0.5cm,inner sep=0pt}}
\begin{tikzpicture}
  \node[node](k) at (.65,1){$\quad x_{i-m}\quad $};
  %\node[node](k+1) at (1,1){$\quad x_{i-m+1}\ \ $};
  \node[node2](m) at (3,1){$\ \cdots\ $};
  \node[node](n+k-1) at (5.1,1){$\ x_{i+n}\ $};
  \node[node](y) at (3,0){$\quad y_i\quad$};

  \draw[ ->,>=stealth,right] (m) edge node {$f$} (y);
  %\node[title] at (0,0.5){$f$};
\end{tikzpicture}
\caption{The sliding block code.}\label{figureSlidingBlockCode}
\end{figure}

The \emph{sliding block code}\index{subject}{sliding block code}
\index{subject}{block!code}%
induced by $f$ is the map $\varphi:X\to B^\Z$ defined, for all $x\in X$, by 
$y=\varphi(x)$ if (see Figure~\ref{figureSlidingBlockCode})
\begin{displaymath}
y_i=f(x_{i-m}\cdots x_{i+n}) \quad (i\in\Z).
\end{displaymath}
When $ \varphi(X)$ is included in $ Y$, we write $\varphi:X\to Y$.

\begin{theorem}[Curtis, Hedlund, Lyndon]\label{theoremCurtisHedlundLyndon}
\index{subject}{Curtis-Hedlund-Lyndon Theorem}%
\index{subject}{Theorem!Curtis-Hedlund-Lyndon}%
\index{names}{Curtis, Morton L.}%
\index{names}{Lyndon, Roger C.}%
\index{names}{Hedlund, Gustav A.}
Let $X$ and $Y$ be shift spaces.
A map $\varphi:X\to Y$ is a factor map if and only if it is
a sliding block code.
\end{theorem} 
\begin{proof}
A sliding block code is clearly continuous and commutes with the shifts,
that is $\varphi\circ T=S\circ \varphi$. Thus it is a morphism of dynamical systems.
Conversely, let $\varphi:X\to Y$ be a morphism. 
For every $b\in B$ the set $[b]$ is clopen.
The map $\varphi$ being continuous, it is also the case for $\varphi^{-1}([b])$.
The family  $\{ \varphi^{-1}([b]) | b\in B \}$ being a partition, the set $B$ being finite and each  $\varphi^{-1}([b])$ being a finite union of cylinders,
there is an integer $n$ such that $\varphi(x)_0$ depends only
on $x_{[-n,n]}$. Set $f(x_{-n}\cdots x_n)=\varphi(x)_0$.
Then $\varphi$ is the sliding block code associated
to $f$.
\end{proof}

The \emph{factor complexity}\index{subject}{factor!complexity}
(or \emph{word complexity}
\index{subject}{word!complexity} 
or simply the \emph{complexity}\index{subject}{shift space!complexity})
  of the shift space $X$ is the
sequence 
\begin{displaymath}
p_n(X)=\Card(\mathcal{L}_n(X)).\index{symbols}{pn(X)@$p_n(X)$}
\end{displaymath}
Observe that $p_0(X)=1$ and that $p_n(X)\le p_{n+1}(X)$. Indeed,
for every $w\in\cL_n(X)$, there is a letter $a\in A$ such that $wa\in \cL(X)$.
\begin{example}
Let $X$ be the golden mean shift (Example~\ref{exampleGolden}).
The factor complexity of $X$ is given by $p_n(X)=F_{n+1}$ where
$F_n$ is the \emph{Fibonacci sequence}\index{subject}{Fibonacci!sequence of numbers}
defined by $F_0=0, F_1=1$ and $F_{n+1}=F_n+F_{n-1}$ for $n\ge 1$.
Indeed,  the number of words in $\cL_{n+1}(X)$ ending with $a$
is equal to $p_n(X)$ and the number of those ending with $b$
is $p_{n-1}(X)$ since the $b$ has to be preceded by $a$.
Thus $p_{n+1}(X)=p_n(X)+p_{n-1}(X)$.
\end{example}
The following result is the counterpart for shift spaces
of Theorem~\ref{theoremCovenHedlundSequence}.
\index{subject}{Theorem!Morse-Hedlund}\index{subject}{Morse-Hedlund!Theorem}%
\begin{theorem}[Morse, Hedlund]\label{theoremCovenHedlund}
\index{names}{Morse, Marston}\index{names}{Hedlund, Gustav A.}
Let $X$ be a shift space. The following conditions are equivalent.
\begin{enumerate}
\item[\rm(i)] For some $n\ge 1$, one has $p_n(X)\le n$.
\item[\rm(ii)] For some $n\ge 1$, one has $p_n(X)=p_{n+1}(X)$.
\item[\rm(iii)] $X$ is finite.
\end{enumerate}
\end{theorem}
\begin{proof}
(i) $\Rightarrow$ (ii).
Since $p_n(X)\le p_{n+1}(X)$ for all $n\ge 0$, the hypothesis implies
that $p_n(X)=p_{n+1}(X)$ for some $n\ge 0$. 

(ii) $\Rightarrow$ (iii). By Theorem~\ref{theoremCovenHedlundSequence},
every $x\in X$ is periodic and its period is bounded by
$\max p_n(X)$.

(iii) $\Rightarrow$ (i). If $X$ is finite, it is a finite
union of orbits of periodic points of the form $u^\infty$ for some word $u$.
Since the
set $\cL_n(u^\infty)$ has at most $|u|$ elements, $p_n(u^\infty)$
is bounded and thus also $p_n(X)$.
\end{proof}
Thus, by Theorem~\ref{theoremCovenHedlund},
 either $p_n(X)\ge n+1$ for all $n\ge 1$ or $p_n(X)$
is eventually constant. The case $p_n(X)=n+1$ for all $n\ge 1$
corresponds to the Sturmian shifts (see below).

%We say that a shift space $(X,S)$ is \emph{aperiodic}
%\index{subject}{aperiodic!shift space}\index{subject}{shift space!aperiodic}%
%if it does not contain any periodic point.

A shift space $X$ is \emph{irreducible}\index{subject}{irreducible!shift space}
\index{subject}{shift space!irreducible} if the language $\cL(X)$
is irreducible, that is, if for every $u,v\in \mathcal{L}(X)$
there is a $w$ such that $uwv\in \mathcal{L}(X)$.
\begin{proposition}\label{propositionRecurrentShiftSpace}
A shift space  is recurrent if and only if it is irreducible.
\end{proposition}
\begin{proof}
Assume first that $X$ is recurrent and consider $u,v\in\cL(X)$.
Let $n>|u|$ be such that $S^{-n}([v]_X)\cap[u]_X\ne\emptyset$. Then
every $x\in S^{-n}([v]_X)\cap[u]_X$ is in $[uwv]_X$ for some $w$. Thus
$uwv$ belongs to $\cL(X)$, showing that $X$ is irreducible.

Conversely, assume that $X$ is irreducible. Let $U,V$
 be open sets in $X$. We can find cylinder sets $[u]_X$ and  $[v]_X$ that are respectively included in $\subset U$
and $V$.  Since $X$ irreducible there is some word
$w$ such that $uwv\in\cL(X)$. Then $[uwv]_X$ is nonempty and
in $S^{-n}V\cap U$ for $n=|uw|$. Thus $X$ is recurrent.
\end{proof}
A shift space is of course irreducible if and only if it is recurrent
as a topological dynamical system. Thus we could also have used Proposition
\ref{propositionRecurrentSystem} to prove Proposition~\ref{propositionRecurrentShiftSpace}.

Thus, for a shift space,  the property of being recurrent
can be expressed by a property of its language $\cL(X)$.
Likewise, we can translate the property of being minimal.
A shift space $X$ is \emph{uniformly recurrent}
\index{subject}{uniformly!recurrent!shift}\index{subject}{recurrent!uniformly!shift}%
\index{subject}{shift space!uniformly recurrent}%
 if the language $\cL(X)$ is uniformly recurrent, that is,
if for every $u\in \mathcal{L}(X)$
there is an $n\ge 1$ such that $u$
is a factor of every word in $\mathcal{L}_n(X)$.
\begin{proposition}\label{propositionMinimalUR}
The following conditions are equivalent
for a shift space $X$.
\begin{enumerate}
\item[\rm (i)] The shift space $X$ is minimal.
\item[\rm (ii)] The shift space $X$ is uniformly recurrent.
  \item[\rm (iii)] Every $x\in X$  is uniformly recurrent and $\cL(x)=\cL(X)$.
\end{enumerate}
\end{proposition}
\begin{proof}
  (i) $\Rightarrow$ (ii). Let $u\in \cL(X)$. Since $[u]_X$
  is clopen and $X$ is minimal, for every $x\in X$,
  the entrance time
  \begin{displaymath}
    n(x)=\min\{n>0\mid S^nx\in[u]_X\}
    \end{displaymath}
  exists and is bounded. Set $n=\max n(x)$.
  Then  every word in $\cL(X)$
  of length $n+|u|$ has a factor $u$, showing that $X$
  is uniformly recurrent.

  (ii) $\Rightarrow$ (iii) is clear since every $u\in\cL(X)$ appears in every
  word of  $\cL_n(x)$.

(iii) $\Rightarrow$ (i) is clear since the orbit of every $x\in X$
is dense.
\end{proof}
Again, we could have used Proposition~\ref{propositionMinimal}
to prove Proposition~\ref{propositionMinimalUR}.

\subsection{Special words}
For $w\in \mathcal{L}(X)$, there
is at least one letter $a\in A$ such that $wa$ belongs to $\mathcal{L}(X)$ and
symmetrically, at least one letter $a\in A$ such that $aw$ belongs to $\mathcal{L}(X)$. The
word $w$ is called \emph{right-special}\index{subject}{right!special word}
\index{subject}{word!right-special} if there is at least two letters
$a\in A$ such that $wa$ belongs to $\mathcal{L}(X)$.
Symmetrically, $w$ is \emph{left-special}\index{subject}{left!special word}
\index{subject}{word!left-special} if there is at least two letters
$a\in A$ such that $aw$ belongs to $\mathcal{L}(X)$.
It is \emph{bispecial}
\index{subject}{bispecial!word}\index{subject}{word!bispecial}%
if it is both left and right-special.

Special words are closely linked with the factor complexity of a shift space.
Let indeed $p_n(X)=\Card(\cL_n(X))$ be the factor complexity of the shift space $X$
on the alphabet $A$. Set for $n\ge 0$ 
\begin{eqnarray*}
s_n(X) & = & p_{n+1}(X)-p_n(X), \\ \index{symbols}{sn(X)@$s_n(X)$}
b_n(X) & = & s_{n+1}(X)-s_n(X).\index{symbols}{bn(X)@$b_n(X)$}
\end{eqnarray*}
For a word $w\in\cL(X)$, let
\begin{eqnarray*}
\ell_X(w)&=&\Card\{a\in A\mid aw\in\cL(X)\}\\
r_X(w)&=&\Card\{b\in A\mid wb\in\cL(X)\}\\
e_X(w)&=&\Card\{(a,b)\in A\times A\mid awb\in\cL(X)\}
\end{eqnarray*}
\index{symbols}{l@$\ell_X(w)$}\index{symbols}{r@$r_X(w)$}%
\index{symbols}{e@$e_X(w)$}%
Thus $\ell_X(w)>1$ (resp. $r_X(w)>1$) if and only if $w$
is left-special (resp. right-special). Define
also the \emph{multiplicity}\index{subject}{multiplicity of word} of $w\in\cL(X)$ as
\begin{displaymath}
m_X(w)=e_X(w)-\ell_X(w)-r_X(w)+1.
\end{displaymath}
\index{symbols}{mX(w)@$m_X(w)$}
The word $w$ is called \emph{neutral}
\index{subject}{neutral!word}\index{subject}{word!neutral}%
if $m_X(w)=0$.

\begin{proposition}
\label{propositionCassaigne}
We have for all $n \ge 0$,
\begin{equation}
s_n(X) = \sum_{w \in \cL_n(X)}(\ell_X(w)-1) = \sum_{w \in \cL_n(X)}(r_X(w)-1)
\label{eqs_n}
\end{equation}
and
\begin{equation}
b_n(X) = \sum_{w \in \cL_n(X)} m_X(w).
\label{eqb_n}
\end{equation}
In particular, the number of left-special (resp. right-special) words of length $n$ is bounded by $s_n(X)$.
\end{proposition}
\begin{proof}
We have
\begin{eqnarray*}
\sum_{w \in \cL_n(X)}(\ell_X(w)-1) & = & \sum_{w \in \cL_n(X)}\ell_X(w) - \Card(\cL_n(X)) \\
 & = & \Card(\cL_{n+1}(X)) - \Card(\cL_n(X)) = p_{n+1}(X)-p_n(X) \\
 & = & s_n(X)
\end{eqnarray*}
with the same result for $\sum_{w \in \cL_n(X)}(r_X(w)-1)$. Next,

\begin{eqnarray*}
\sum_{w \in \cL_n(X)}m_X(w) & = & \sum_{w \in \cL_n(X)}(e_X(w)-\ell_X(w)-r_X(w)+1) \\
& = & p_{n+2}(X)-2p_{n+1}(X)+p_n(X) = s_{n+1}(X)-s_n(X)\\
&=& b_n(X).
\end{eqnarray*}
\end{proof}

\subsection{Return words}\label{sectionReturnWords}
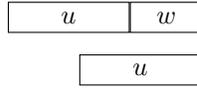
\begin{figure}[hbt]
\centering
\tikzset{node/.style={draw,minimum size=0.4cm,inner sep=0pt}}
	\tikzset{title/.style={minimum size=0.5cm,inner sep=0pt}}
\begin{tikzpicture}
  \node[node](wh) at (0,0.7){$\qquad u\qquad$};
  \node[node](u) at (1.3,0.7){$\quad w\quad$};
  \node[node](wb) at (0.95,0) {$\qquad u\qquad$};
\end{tikzpicture}
\caption{A right return word}\label{figureRightReturn}
\end{figure}

Let $X$ be a shift space.
For $u\in \mathcal{L}(X)$ a \emph{right return word}
\index{subject}{right!return word}\index{subject}{return!word}%
\index{subject}{return!word!right} to $u$ is a nonempty word $w$ such that
$uw$ is in $\mathcal{L}(X)$, $uw$ has $u$ as a proper suffix
and $uw$ has no factor $u$ which is not a prefix or a suffix.
Thus a right return word is a word $w$ such that, reading
from left to right and having already seen $u$, after reading $w$, one sees
again the word $u$ for the first time (see Figure~\ref{figureRightReturn}).
%\begin{figure}[hbt]
%\centering
%\gasset{Nmr=0,Nadjust=h}
%\begin{picture}(50,10)
%\node[Nw=30](w)(0,5){$w$}\node[Nw=20](u)(25,5){$u$}
%\node[Nw=30](w2)(20,0){$w$}
%%\end{picture}
%\caption{A right return word}\label{figureRightReturn}
%\end{figure}

For example, in the golden mean shift, the word $w=aab$
is a right return word to $u=b$.

Symmetrically, a \emph{left return word}\index{subject}{left!return word}
\index{subject}{word!left return}
\index{subject}{return!word!left} to $u$ is a word $w$ such that
$wu$ is in $\mathcal{L}(X)$, $wu$ has $u$ as a proper prefix and $wu$ has
no other factor $u$.

We denote by $\RR_X(u)$ (resp. $\RR'_X(u)$)
\index{symbols}{R@$\RR_X(u)$}\index{symbols}{R@$\RR'_X(u)$}%
 the set of
right (resp. left) return words to $u$. The set
$U=\RR_X(u)$ is a \emph{prefix code},\index{subject}{prefix!code}
\index{subject}{code!prefix}%
that is, no word in $U$ is a prefix of another one.
Symmetrically, the set $U'=\RR'_X(u)$ is a \emph{suffix code},
\index{subject}{suffix!code}\index{subject}{code!suffix}%
that is no word of $U'$ is a suffix of another one.

For any prefix code $U$, the set $U^*$ satisfies
\begin{equation}
  r,rs\in U^*\Rightarrow s\in U^*.\label{equationUnitary}
  \end{equation}
The following property of return words is easy to verify.
\begin{proposition}
Every word $w\in \cL(X)$  which begins and ends with $u$ has a unique
factorisation $w=uw_1w_2\cdots w_n$ with $w_i$ in $\RR_X(u)$ and $n\ge 0$.
\end{proposition}

For example, if $X$ is the golden mean shift, we have
\begin{displaymath}
\RR_X(b)=\{ab,aab,aaab,\ldots\}\mbox{ and }
\RR_X'(b)=\{ba,baa,\ldots\}.
\end{displaymath}

Clearly a recurrent shift space $(X,S)$ is minimal if and only if
$\RR_X(w)$ is finite for every $w\in \cL(X)$.
\subsection{Rauzy graphs}\label{sectionChapter2RauzyGraphs}
Let $(X,S)$ be a shift space on the alphabet $A$.
The \emph{Rauzy graph}\index{subject}{Rauzy!graph}\index{subject}{graph!Rauzy}
\index{names}{Rauzy, G\'erard} of $X$ of order $n$,
denoted $\Gamma_n(X)$,\index{symbols}{G@$\Gamma_n(X)$} is
the following labeled graph.  The set of vertices of $\Gamma_n(X)$
is the set $\cL_{n-1} (X)$ and the set of edges is $\cL_n(X)$.
The origin and end of the edge $w$ are the words $u,v$
 such that $w=ua=bv$ with $a,b$ in $A$. The label of the
edge $w$ is $a$.

\begin{example}
Let $X$ be the Fibonacci shift.
\index{subject}{shift space!Fibonacci}%
 The Rauzy graphs of order
$n=1,2,3$ are represented in Figure~\ref{figureRauzy}
(with the edge from $w=ua$ labeled $a$).
\begin{figure}[hbt]
\centering
\tikzset{node/.style={circle,draw,minimum size=0.4cm,inner sep=0.2pt}}
\tikzstyle{every loop}=[->,shorten >=1pt,looseness=12]
\tikzstyle{loop left}=[in=130,out=220,loop]
\tikzstyle{loop right}=[in=330,out=50,loop]
\begin{tikzpicture}

\node[node](11) at (0,0){$\varepsilon$};

\draw[left](11) edge[loop left]node {$a$}(11);
\draw[right](11) edge[loop right]node {$b$}(11);

\node[node](21)at(3,0){$a$};\node[node](22)at(5,0){$b$};

\draw[left](21)edge[loop left]node{$a$}(21);
\draw[above, bend left, ->](21) edge node{$b$}(22);
\draw[below, bend left, ->](22)edge node{$a$}(21);

\node[node](31)at (6,0){$aa$};
\node[node](32)at (8,1){$ab$};
\node[node](33)at(8,-1){$ba$};

\draw[above, bend left, ->](31)edge node{$b$}(32);
\draw[right,bend left, ->](32)edge node{$a$}(33);
\draw[below, bend left, ->](33)edge node{$a$}(31);
\draw[left, bend left, ->](33)edge node{$b$}(32);

\end{tikzpicture}
\caption{The Rauzy graphs of order $n=1,2,3$ of the Fibonacci shift.}
\label{figureRauzy}
\end{figure}
\end{example}
Every infinite path $\cdots \edge{a_{i-1}}p_i\edge{a_i}p_{i+1}\edge{a_{i+1}}\cdots$ in $\Gamma_n(X)$
has a label, which is the sequence $(a_i)_{i\in \Z}$.
The set of these labels
is a shift space $X_n$. We have
\begin{displaymath}
X_1\supset X_2\supset\cdots\supset X_n\supset\cdots\supset X
\end{displaymath}
and $X=\cap_ {n\ge 0}X_n$. Thus the sequence $X_n$ approximates $X$
from above. A graph is \emph{strongly connected}\index{subject}{strongly connected graph}
\index{subject}{graph!strongly connected} if for every pair of vertices
$v,w$ there is a path
from $v$ to $w$.
\begin{proposition}\label{propositionFausse}
If a shift space $X$ is recurrent, all graphs $\Gamma_n(X)$
are strongly connected.
\end{proposition}
\begin{proof}
Assume that $X$ is recurrent (or, equivalently, irreducible). If $u,v$ belong to $\cL_{n-1}(X)$,
there is some $w\in\cL(X)$ such that $uwv$ belongs to $\cL(X)$. But then there
is a path labeled $wv$ from $u$ to $v$ in $\Gamma_n(X)$. Thus $\Gamma_n(X)$
is strongly connected.

\end{proof}
The converse of Proposition~\ref{propositionFausse} is not true,
as shown by the example of the shift space $X$ such that $\cL(X)=a^*b^*\cup b^*a^*$.
\subsection{Higher block shifts}
Let $X$ be a shift space on the alphabet $A$ and let $k\ge 1$
be an integer. Let $f:\mathcal{L}_k(X)\rightarrow A_k$ be a bijection
from  the set $\mathcal{L}_k(X)$ of blocks of  length $k$ of $X$
onto an alphabet $A_k$.
The map $\gamma_k:X\rightarrow A_k^\Z$ defined for $x\in X$ by $y=\gamma_k(x)$
if for every $n\in\Z$
\begin{displaymath}
y_n=f(x_n\cdots x_{n+k-1})
\end{displaymath}
is the $k$-th \emph{higher block code}\index{subject}{higher block!code}
\index{subject}{block!code} on $X$.
The set $X^{(k)}=\gamma_k(X)$
\index{symbols}{X@$X^{(k)}$}  is a shift space on $A_k$, called the $k$-th
\emph{higher block presentation}
\index{subject}{higher block!presentation!of a shift}%
\index{subject}{block!presentation}\index{subject}{presentation!higher block!of shift}%
of $X$ (one also uses the term of coding by \emph{overlapping blocks} of length $k$).
\index{subject}{coding!by overlapping blocks}

The higher block code is an isomorphism of dynamical systems and the inverse
of $\gamma_k$ is the map $y\mapsto x$ such that
$x_n$ is the first letter of $f^{-1}(y_n)$ for all $n$.
\begin{figure}[hbt]
\centering
%\gasset{Nmr=0,Nw=10,Nh=4}
%\begin{picture}(50,10)

%\node(k)(20,10){$x_{n}$}
%\node(k+1)(30,10){$x_{n+1}$}\node(k+n)(52,10){$x_{n+k-1}$}
%\node[Nw=12](d)(41,10){$\cdots$}
%\node(y)(20,0){$y_n$}
%\drawedge(k,y){$f$}
%\end{picture}
\tikzset{node/.style={draw,minimum size=0.4cm,inner sep=0pt}}
	\tikzset{title/.style={minimum size=0.5cm,inner sep=0pt}}
\begin{tikzpicture}
  \node[node](k) at (-0.2,1){$\quad x_n\quad $};
  \node[node](k+1) at (1,1){$\quad x_{n+1}\ \ $};
  \node[node](...) at (2,1){$\ \cdots\ $};
  \node[node](n+k-1) at (3,1){$\ x_{n+k-1}\ $};
  \node[node](y) at (-0.2,0){$\quad y_n\quad$};
  \draw[ ->,>=stealth] (k) edge node {} (y);
  \node[title] at (0,0.5){$f$};
\end{tikzpicture}
\caption{The $k$-th higher block code.}\label{figureSlidingBlock}
\end{figure}
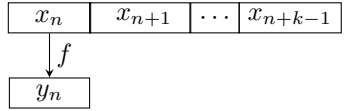

We sometimes, when no confusion arises, identify $A_k$ and $\cL_k(X)$
and write simply $y_0y_1\cdots=(x_0x_1\cdots x_{k-1})(x_1x_2\cdots x_k)\cdots$.
\begin{example}
Consider again the golden mean shift $X$ (Example~\ref{exampleGolden}).
We have $\cL_3(X)=\{aaa,aab,aba,baa,bab\}$. Set $f:aaa\mapsto x,
 aab\mapsto y,aba\mapsto z, baa\mapsto t,bab\mapsto u$.
The third higher block shift $X^{(3)}$ of $X$ is the set of two-sided infinite
paths in the graph of Figure~\ref{figureHigher} on the left
(this graph is, up to the labeling, the Rauzy graph $\Gamma_3(X)$,
see below on the right).
%\begin{figure}[hbt]

%\centering
%\gasset{Nadjust=wh}
%\begin{picture}(20,20)
%\node(1)(0,10){$1$}\node(2)(20,20){$2$}\node(3)(20,0){$3$}

%\drawloop[loopangle=180](1){$x$}
%\drawedge[curvedepth=3](1,2){$y$}
%\drawedge[curvedepth=3](2,3){$z$}
%\drawedge[curvedepth=3](3,2){$u$}
%\drawedge[curvedepth=3](3,1){$t$}
%\end{picture}
%\caption{The third higher block coding of the golden mean shift}\label{figureHigher}

%\end{figure}
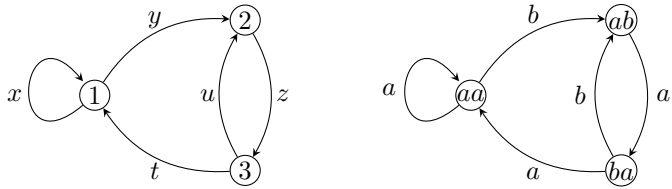
\begin{figure}[hbt]
\centering
\tikzset{node/.style={circle,draw,minimum size=0.4cm,inner sep=0pt}}
	\tikzset{title/.style={minimum size=0.5cm,inner sep=0pt}}
        \tikzstyle{every loop}=[->,shorten >=1pt,looseness=12,>=stealth]
        \tikzstyle{loop left}=[in=130,out=220,loop]
        \tikzstyle{edge}=[->,>=stealth]
\begin{tikzpicture}
  \node[node](1) at (0,1) {$1$};
  \node[node](2) at (2,2) {$2$};
  \node[node](3) at (2,0) {$3$};

  \draw[bend left, ->,>=stealth] (1) edge node {} (2);\node[title] at (0.8,2){$y$};
  \draw[bend left, ->,>=stealth] (2) edge node {} (3);\node[title] at (2.5,1){$z$};
  \draw[bend left, ->,>=stealth] (3) edge node {} (2);\node[title] at (1.5,1){$u$};
  \draw[bend left, ->,>=stealth] (3) edge node {} (1);\node[title] at (0.8,0){$t$};
  \draw (1) edge [loop left] node {} (1);\node [title][left = 0.6cm of 1]{$x$};

 \node[node](1) at (5,1) {$aa$};
  \node[node](2) at (7,2) {$ab$};
  \node[node](3) at (7,0) {$ba$};

  \draw[bend left, ->,>=stealth,above] (1) edge node {$b$} (2);
  \draw[bend left, ->,>=stealth,right] (2) edge node {$a$} (3);
  \draw[bend left, ->,>=stealth,left] (3) edge node {$b$} (2);
  \draw[bend left, ->,>=stealth,below] (3) edge node {$a$} (1);
  \draw[left ] (1) edge [loop left] node {$a$} (1);
\end{tikzpicture}

\caption{The third higher block coding of the golden mean shift $X$
and the graph $\Gamma_3(X)$}\label{figureHigher}
\end{figure}
\end{example}
Note that for $n\ge k$, the Rauzy graph $\Gamma_n(X^{(k)})$ is the same,
up to the labels, as the graph $\Gamma_{n+k-1}(X)$. As an example,
the graph of Figure~\ref{figureHigher} on the left can be identified either
with $\Gamma_2(X^{(2)})$ or with $\Gamma_3(X)$
(see Figure~\ref{figureHigher} on the right).
%%%%%%%%%%%%%%%%%%
\section{Shifts of finite type}\label{sectionSFT}
A shift space $X$ is  \emph{of finite type}
\index{subject}{shift space!of finite type}%
\index{subject}{subshift!of finite type}%
if $\mathcal{L}(X)=A^*\setminus A^*IA^*$ where $I\subset A^*$ is finite set. 
In other words, a  two-sided infinite
sequence $x$ is in $X$ if and only if
it has no factor in the finite set $I$.
The elements of $I$ are called the \emph{forbidden blocks}
\index{subject}{forbidden blocks} of $X$. 

The class of shifts of finite type is closed under conjugacy (Exercise~\ref{exerciseConjugacySFT}).

A well-know example is the golden mean shift
\index{subject}{golden mean!shift}
\index{subject}{shift space!golden mean} $X$  where $X$ is the set
of two-sided sequences on $A=\{a,b\}$ with no consecutive $b$
(see Example~\ref{exampleGolden}).
Thus the set of forbidden blocks is $I=\{bb\}$ and $X$ is the set of labels of two-sided infinite paths in the graph of Figure
\ref{figureGolden}.

As a more general example, for every shift space $X$ and $n\ge 1$,
the set $X_n$ of labels of bi-infinite  paths
\index{subject}{bi-infinite!path} in the Rauzy graph
$\Gamma_n(X)$ is a shift of finite type. Indeed, it is defined
by the finite set of forbidden blocks which is the
set of words of length $n$ which are not in $\cL_n(X)$.

Given a finite graph $G=(V,E)$, the \emph{edge shift}
\index{subject}{edge!shift}%
on $G$ is the shift space $X$ where $X\subset E^\Z$
is the set of bi-infinite paths in $G$ and $S$ is the shift
on $E^\Z$. An edge shift is a shift of finite type, since
it is defined by forbidden blocks of length $2$. Moreover,
if the graph $G$ is strongly connected, the edge shift
on $G$ is irreducible (and thus recurrent).
\begin{proposition}\label{propositionSFTEdgeShift}
Any shift
of finite type is conjugate to an edge shift on some graph. The shift is recurrent
if and only if the graph can be chosen strongly connected.
\end{proposition}
The proof is left as an exercise (Exercise~\ref{exerciseEdgeShifts}).
\begin{example}
The edge shift on the graph $G$ represented in Figure~\ref{figureEdgeShift}
is conjugate to the golden mean shift by the $1$-block map
$e\mapsto a,f\mapsto b,g\mapsto a$.
\begin{figure}[hbt]
\centering
\tikzset{node/.style={circle,draw,minimum size=0.4cm,inner sep=0pt}}
	\tikzset{title/.style={minimum size=0.5cm,inner sep=0pt}}
        \tikzstyle{every loop}=[->,shorten >=1pt,looseness=12]
        \tikzstyle{loop left}=[in=130,out=220,loop]
        %\tikzset{state/ .style={minimum size=0.5cm,inner sep=0pt}}
	\begin{tikzpicture}
	\node[node](1) {$1$};
	\node[node](2) [right = 2cm of 1] {$2$};
	\draw[bend left, ->,>=stealth,above] (1) edge node {$f$} (2);
         \draw[bend left, ->,>=stealth,below] (2) edge node {$g$} (1);
        \draw[left] (1) edge [loop left,>=stealth] node {$e$} (1);
        \end{tikzpicture}
\caption{An edge shift}\label{figureEdgeShift}
\end{figure}
\end{example}

%%%%%%%%%%%%%%
\section{Substitution shifts}
\label{sectionSubstitutionSystems}

Let $A$ and $B$ be finite  alphabets. 
By a {\em morphism}\index{subject}{morphism!of monoids}
\index{subject}{monoid!morphism}from $A^*$ to $B^*$ we mean 
a morphism of  monoids, that is a map $\sigma:A^*\to B^*$ such that 
$\sigma(\varepsilon)=\varepsilon$ and $\sigma(uv)=\sigma(u)\sigma(v)$ for all $u,v\in A^*$.
When $\sigma (A) = B$, we say $\sigma$ is a {\em letter-to-letter morphism}\index{subject}{letter-to-letter morphism}.
\index{subject}{morphism!letter-to-letter}%
Thus, letter-to-letter morphisms are onto.

We set $|\sigma | = \max_{a\in A} |\sigma (a)|$ and $\langle \sigma \rangle = \min_{a\in A} |\sigma (a)|$.
\index{symbols}{sigma@$\lvert\sigma\rvert$}%
\index{symbols}{sigma@$\langle \sigma \rangle$}%
The morphism $\sigma$ is of {\em constant length}
\index{subject}{constant!length morphism}\index{subject}{morphism!constant length} if $\langle \sigma \rangle = |\sigma|$. It is
is {\em growing}\index{subject}{growing morphism}\index{subject}{morphism!growing} if $\lim_{n\to\infty} \langle \sigma^n \rangle = +\infty$
or, equivalently, if $|\sigma^n(a)|\to\infty$ for every $a\in A$.

We say that a morphism
$\sigma$ is {\em erasing}
\index{subject}{erasing morphism}\index{subject}{morphism!erasing}%
 if there exists $b\in A$ such that $\sigma (b)$ is the empty word
and \emph{nonerasing}\index{subject}{nonerasing morphism}
\index{subject}{morphism!nonerasing} otherwise. 
A growing morphism is nonerasing.
If $\sigma$ is nonerasing, it induces, by infinite concatenation, a map from $A^{\mathbb{N}}$ to $B^{\mathbb{N}}$ defined for $x\in A^\N$ by
\begin{displaymath}
\sigma(x)=\sigma(x_0)\sigma(x_1)\cdots
\end{displaymath}
and a map from $A^{\mathbb{Z}}$ to $B^{\mathbb{Z}}$ defined for $x\in A^\Z$ by
\begin{displaymath}
\sigma(x)=\cdots\sigma(x_{-1})\cdot\sigma(x_0)\sigma(x_1)\cdots
\end{displaymath} 
These maps are also denoted by $\sigma$. Both are continuous because
if $x,y$ coincide on the first $n$ letters, then so do
$\sigma(x),\sigma(y)$.

The \emph{language}\index{subject}{language!of a morphism}
 of the morphism $\sigma : A^* \to A^*$ is the set $\mathcal{L} (\sigma )$
 \index{symbols}{l@$\cL(\sigma)$} of words occurring in some $\sigma^n (a)$, $a\in A$, $n \in \mathbb{N}$ and $X (\sigma )$\index{symbols}{X(sigma)@$X(\sigma)$}
 will denote the set of sequences $y \in A^\Z$ whose subwords belong to $ \mathcal{L}(\sigma)$. The
set $X(\sigma)$ is clearly a shift space since it is
closed and invariant. It is called
the \emph{shift generated}
\index{subject}{subshift!generated!by a morphism}%
\index{subject}{shift space!generated!by a morphism}  by $\sigma$.

By definition $\cL(\sigma^n)$ is included in $\cL(\sigma)$ and
$X(\sigma^n)$ in $X(\sigma)$
for all $n\ge 1$ (the equality does not always hold, see Exercise~\ref{exercisesigma^n}).
Observe that $\mathcal{L} (\sigma )$ is finite if and only if $X (\sigma )$ is empty. Observe also that $\cL(X(\sigma))$ is included in  $\cL(\sigma)$
but that the equality might not hold in general. For example,
if $\sigma$ is the morphism defined by $a\mapsto a$, $b\mapsto ba$, then $X(\sigma)$ is reduced to $a^\infty$
and thus $b$ is not a factor of $X(\sigma)$. We will soon
introduce  a class of morphisms (the substitutions) where this
kind of phenomenon cannot happen.

%We say that a morphism $\sigma:A^*\to A^*$ is a \emph{substitution}
%\index{subject}{substitution}%
%if $\cL(\sigma)$ is extendable (or equivalently if every
%word in $\cL(\sigma)$ is a factor of some $x\in X(\sigma)$,
%that is, if $\cL(X(\sigma))=\cL(\sigma)$).
%Note that in this case $\cL(\sigma)$ is infinite or that, equivalently,
%$\lim_{n\to\infty}|\sigma^n(a)|=\infty$ for some $a\in A$.
%Note also that every growing morphism is a substitution.
%When $\sigma$ is a substitution,
%the shift space $X(\sigma)$ is  called a \emph{substitution shift}.
%\index{subject}{substitution!shift}%

%Observe also that, under the mild assumption that every letter $a\in A$ appears
%in some $\sigma(b)$ for $b\in A$, one may replace $\sigma$ by one of its powers
%since then $\cL(\sigma^n)=\cL(\sigma)$ and thus $X(\sigma^n)=X(\sigma)$.

\begin{example}\label{exampleFibonacci0}
The \emph{Fibonacci morphism}\index{subject}{Fibonacci!morphism}
\index{names}{Fibonacci, Leonardo}%
\index{subject}{morphism!Fibonacci}%
(Example~\ref{exampleFibonacci00})
is   defined by $\varphi:a\mapsto ab$, $b\mapsto a$.
The shift space associated to $\varphi$ is called the
\emph{Fibonacci shift}.\index{subject}{Fibonacci!shift}
\index{subject}{shift space!Fibonacci}%
%The Fibonacci shift is Sturmian (see Exercise~\ref{exerciseFiboisSturmian})
An example of a two-sided infinite sequence in $X(\varphi)$ is
$z=\tilde{x}\cdot x$ where $x$ is the Fibonacci word. Indeed,
$\varphi^2(a)=aba$ ends with $a$ and thus there is a unique left-infinite
word $y$ having all $\varphi^{2n}(a)$ as suffixes. Since $\varphi^2(a)$
is a palindrome, we have $y=\tilde{x}$. Then all factors of $z$
are factors of some $\varphi^{2n}(aa)$ and thus, as $aa$ is in $\cL (\varphi)$, the sequence $z$ belongs to $X(\varphi)$.
\end{example}

\begin{example}
The \emph{Thue-Morse morphism}
\index{subject}{Thue-Morse! morphism}\index{subject}{morphism!Thue-Morse}%
\index{names}{Morse, Marston}\index{names}{Thue, Axel}%
 is the morphism defined by
$\tau:a\mapsto ab,b\mapsto ba$. The associated shift space
is the \emph{Thue-Morse shift}\index{subject}{Thue-Morse!shift space}\index{subject}{shift space!Thue-Morse}.
\end{example}
A morphism $\sigma:A^*\to A^*$ is called \emph{periodic}\index{subject}{periodic!morphism}
\index{subject}{morphism!periodic} if its associated shift space
is periodic and \emph{aperiodic}\index{subject}{aperiodic!morphism}
otherwise. When $\sigma$ is periodic, its period can be bounded
a priori and thus this property is decidable for a morphism
(Exercise~\ref{exerciseFixedPointPeriodic2}).
\begin{example}
The shift space associated to the morphism $\varphi:a\mapsto aba,b\mapsto b$ is periodic. Indeed,
$\varphi(ab)=(ab)^2$ and thus the set $X(\varphi)$ is
formed of two points.
\end{example}

Let $\sigma$ be a morphism and let $X=X(\sigma)$ be the two-sided
shift defined by $\sigma$. A nonerasing morphism $\sigma$
is a \emph{substitution}\index{subject}{substitution}
if one has
\begin{displaymath}
  \cL(X)=\cL(\sigma).
\end{displaymath}
Substitutions form a fairly
general class of morphisms and,
for example, they can fail to be growing (see the example
of the Chacon binary substitution of Exercise~\ref{exerciseChaconMinimal}). Their
definition expresses the possibility to extend
the finite words which belong
to $\cL(\sigma)$ into bi-infinite words in $X$. 
The property of being a substitution is easy to verify
(Exercise~\ref{exerciseDefSubstitution}). The examples
of morphisms seen before, as the Fibonacci or the Thue-Morse
morphisms are all substitutions.

Note that if $\sigma:A^*\to A^*$ is a substitution, then $X(\sigma)$
is not empty since $\cL(\sigma)$ always contains the alphabet.
Consequently, there is a letter $a\in A$ such that $\lim_{n\to\infty}|\sigma^n(a)|=\infty$.

\begin{example}\label{exampleb,ab}
  The morphism $\sigma:a\to a,b\to b$ is not a substitution
  since $\cL(\sigma)=\{a,b\}$ while $X=\{a,b\}^\Z$.
  As a less trivial example,
  let $\sigma:a\to b,b\to ab$.
  Then $\cL(\sigma)=ab^*\cup b^*$,
  $X=b^\infty$ and $\cL(X)=b^*$. Thus $\sigma$ is not a substitution.
  \end{example}

%We now introduce a family of morphisms which define
%in a satisfactory way infinite sequences.
%We will use term \emph{substitution}\index{subject}{substitution}
% for a morphism $\sigma:A^*\to A^*$
% such that for some letter $a\in A$,
% \begin{enumerate}
% \item $\lim_{n\to\infty}|\sigma^n(a)|=\infty$.
% \item for every letter  $b\in A$ there is $n\ge 1$
%   and $u\in A^*$ such that $aub$ is a factor of $\sigma^n(a)$.
% \end{enumerate}
% The first condition is equivalent to $X(\sigma)\ne\emptyset$
% and is therefore a minimal requirement. The second condition
% implies
% \begin{enumerate}
% \item[3] For every $b\in A$ there is $n\ge 1$ such that
%   $b$ appears in $\sigma^n(a)$
% \end{enumerate}
% which is a mild requirement since a letter $b$ which does not satisfy
% it could be removed.

% If  conditions 1 and 3 are satisfied,  condition 2
% is equivalent to $\cL(X(\sigma))=\cL(\sigma)$ (Exercise~\ref{exerciseDefSubstitution}).

 %If $\sigma$ is a substitution, then $X(\sigma)=X(\sigma^n)$ for all $n\ge 1$.
%Moreover, some power has an admissible fixed-point
% $x$. Additionnally, the  subshift
% generated by $x$ (that is the shift space $X$
% such that $\cL(X)=cL(x)$) coincides with $X(\sigma)$
% (Exercise~\ref{exerciseDefSubstitution2}).
%Thus the shift space generated by $\sigma$
%and the subshift generated by $x$ are fortunately the same thing.

\subsection{Fixed points}
Let $\sigma:A^* \to A^*$ be a morphism.
If there exists a letter $a\in A$ such that $\sigma(a)$ begins with $a$
 and if, moreover,  $\lim_{n\to+\infty}|\sigma^n(a)|=+\infty$, then $\sigma$ is said to be {\em right-prolongable on $a$}.
\index{subject}{right!prolongable}%
\index{subject}{prolongable!right}\index{subject}{morphism!right prolongable}%

Suppose that $\sigma$ is right-prolongable on $a \in A$. 
Since for all $n\in\mathbb{N}$, the word $\sigma^n(a)$ is a prefix of $\sigma^{n+1}(a)$ and because $|\sigma^n(a)|$ tends to infinity with $n$, there
is a unique right infinite word  denoted $\sigma^\omega(a)$
which has all $\sigma^n(a)$ as prefixes. Indeed, for every $i\ge 0$,
the $i$-th
letter of this  infinite word
is the common $i$-th letter of all words $\sigma^n(a)$, $n\in \N$,
longer than $i$. Then
$x=\sigma^\omega(a)$  is a \emph{fixed point}\index{subject}{fixed point!one-sided}
of $\sigma$, which means by definition that
 $\sigma(x)=x$.

By an \emph{admissible one-sided fixed point}
\index{subject}{admissible!one-sided fixed point}%
 of $\sigma$, we  mean a one-sided infinite sequence $x=\sigma^\omega(a)$
where $\sigma$ is right-prolongable on $a$.
\index{subject}{fixed point!admissible!one-sided}%

Observe that a morphism can have other fixed points, either
finite or infinite. In fact, if $\sigma$
is the morphism $a\mapsto a,b\mapsto ba$, then $a$ is a finite fixed point
since $\sigma(a)=a$ and $a^\omega$ is also a fixed point
since $\sigma(a^\omega)=a^\omega$, but it is not admissible.

\begin{proposition}\label{propositionGrowing}
Every  growing morphism  has a power which has an admissible
 one-sided fixed point.
\end{proposition}
\begin{proof}
Let $\sigma:A^*\to A^*$ be a growing morphism.
Let $a\in A$. Since $\sigma$ is growing, it is nonerasing
and thus all $\sigma^n(a)$
are nonempty.
Since $A$ is finite, there are $n,p\ge 1$  such that
$\sigma^n(a)$ and $\sigma^{n+p}(a)$ begin with the same letter, say $b$.
Since $\sigma$ is growing, $\lim_{k\to\infty} |\sigma^{kp}(b)|=\infty$.
Thus $\sigma^p$ is right-prolongable on $b$ and  has an
admissible one-sided fixed point.
\end{proof}

A sequence $x$ which is an admissible
 one-sided fixed point of a substitution
$\sigma$ is said to be {\em purely substitutive}
\index{subject}{purely!substitutive!sequence!one-sided}\index{subject}{substitutive!purely}%
  (with respect to $\sigma$). 
The subshift generated by $x$ is then called a \emph{substitution shift}.

%By Proposition~\ref{propositionGrowing},
%every growing morphism has a power which is a substitution
%(and it generates the same shift space).

\begin{example}\label{exampleThueMorseSeqeunce}
The Thue-Morse morphism $\tau:a\mapsto ab,b\mapsto ba$ is prolongable on $a$
and $b$.
The fixed point $x=\lim\tau^n(a)$ is called the \emph{Thue-Morse sequence}.
\index{subject}{Thue-Morse!sequence}%
One has
\begin{displaymath}
x=abbabaab\cdots
\end{displaymath}
One may show that $x_n=a$ if and only if the number of $1$ in the
binary expansion of $n$ is even (Exercise~\ref{exerciseThueMorse}).
\end{example}

If $x\in A^\mathbb{N}$ is purely substitutive 
(with respect to $\sigma$) and $\phi:A^*\to B^*$ is a 
letter-to-letter morphism, then $y=\phi (x)$ is said to be {\em substitutive}\index{subject}{substitutive} (with respect
to ($\sigma  ,\phi$)) and the shift space it generates, denoted
$X(\sigma,\phi)$, is called 
a {\em substitutive shift}\index{subject}{shift space!substitutive}\index{subject}{substitutive!shift space}.

Two-sided fixed points \index{subject}{fixed point!two-sided}%
of $\sigma$ can be similarly defined.
Assume that $\sigma$ is right-prolongable on $a \in A$ and 
{\em left-prolongable}\index{subject}{left!prolongable}\index{subject}{prolongable!left}\index{subject}{morphism!left prolongable} on $b\in A$, that is, $\sigma(b)$ ends with $b$ and $\lim_{n \to +\infty} |\sigma^n(b)|=+\infty$.
Let $x=\sigma^\omega(a)$ and let $y$ be the left infinite sequence
having all $\sigma^n(b)$ as suffixes. Let $z\in A^\Z$ be the two-sided
sequence such that $x=z^+$ and $y=z^-$. Then $z$
 is a fixed point of $\sigma$ denoted $\sigma^\omega(b\cdot a)$.
\index{symbols}{sigma@$\sigma^\omega(b\cdot a)$}%
It can happen that $\sigma^{ \omega}(b \cdot a)$ does not belong to $X (\sigma )$. 
In fact, $\sigma^{ \omega}(b \cdot a)$ belongs to $X (\sigma )$ if and only if $ba$ belongs to $\mathcal{L}(\sigma)$.
In this case, we say that $\sigma^{ \omega}(b \cdot a)$ is an {\em admissible}\index{subject}{admissible!fixed point}\index{subject}{fixed point!admissible!two-sided} fixed point of $\sigma$.
Equivalently, $z\in A^\Z$ is an admissible fixed point of $\sigma$ if, and only if $z^+,z^-$ are admissible one-sided
fixed points and $z$ belongs to $X(\sigma )$. 

When $z\in A^\Z$ is a two-sided admissible
fixed-point of $\sigma$, we say, as in the
one-sided case, that $z$ is a \emph{purely substitutive}
\index{subject}{purely!substitutive!sequence} two-sided sequence
and if $\phi$ is a letter-to-letter morphism, we say that
$x=\phi(z)$ is a \emph{substitutive} two-sided sequence. Likewise, the
shift generated by $x$, denoted $X(\sigma,\phi)$ is called
a \emph{substitutive shift}.\index{subject}{substitutive!shift}

\begin{proposition}
Every growing morphism has a power with an admissible 
two-sided fixed point.
\end{proposition}
\begin{proof}
Let $\sigma:A^*\to A^*$ be a growing morphism.
Let $a,b\in A$ be such that $ab\in\cL(\sigma)$. There are integers
$n,p$ such that simultaneously 
$\sigma^n(a),\sigma^{n+p}(a)$ end with the same letter $c$
and $\sigma^n(b),\sigma^{n+p}(b)$ begin with the same letter $d$.
This will be also true, with the same $p$
but possibly different letters $c,d$, for $n+1,n+2,\ldots$
and thus
we may also assume that $p$ divides $n$.
Thus $\tau^\omega(c.d)$ is a two-sided fixed point of $\tau=\sigma^p$.
But $cd$ is a factor of $\sigma^{n+p}(ab)$ and thus $cd\in\cL(\tau)$,
showing that $\tau^\omega(c.d)$ is an admissible two-sided fixed point of $\tau$.
\end{proof}
\begin{example}
Let  $\varphi:a\mapsto ab,b\mapsto a$ be the  Fibonacci morphism.
Then $\varphi^2:a\mapsto aba,b\mapsto ab$ is right prolongable 
on $a$ and left prolongable on $a$ and $b$. Since moreover $ba$ is in $\cL(\varphi)$,
the two sided infinite sequences $\varphi^\omega(a\cdot a)$
and $\varphi^\omega(b\cdot a)$ are
 admissible fixed points of $\varphi^2$.
\end{example}
Observe that a morphism $\sigma$ could have non admissible fixed points. 
\begin{example}
Let $\sigma : A^* \to A^*$, with $A= \{a,b,c\}$ be the morphism defined by $a\mapsto ab, b\mapsto ac, c\mapsto aa$. 
It can be checked that $\sigma$ has a unique fixed point in $A^\N$ but no admissible fixed point in $A^\Z$ whereas $\sigma^3$ has three admissible fixed points  : $\sigma^\omega (a.a), \sigma^\omega (b.a)$ and $\sigma^\omega (c.a)$.
\end{example}

\subsection{Primitive morphisms}

A morphism $\varphi:A^*\rightarrow A^*$ is said
to be \emph{primitive}\index{subject}{primitive!morphism}\index{subject}{morphism!primitive} if there is an $n\ge 1$ such that
for every $a,b\in A$,  the letter $b$ appears
in the word $\varphi^n(a)$.

In the following statement, we exclude the case of a one-letter
alphabet for which $\sigma:a\to a$ is a primitive
morphism which is neither growing nor a substitution.

\begin{proposition}\label{propositionPrimitiveGrowing}
If $\sigma:A^*\to A^*$ is a primitive morphism and $A$
has at least two letters, then $\sigma$ is a growing substitution. 
\end{proposition}
\begin{proof}
Let $n\ge 1$ be such that every letter $b\in A$ occurs in 
$\sigma^n (a)$ for every $a\in A$. Then $|\sigma^n(a)|\ge\Card(A)$
and thus $|\sigma^{n+m}(a)|\ge\Card(A)^m$ for all $m\ge 1$.
This shows that $\sigma$ is growing.

Since $\sigma$ is growing, there are letters $a,b,c$ such that $abc\in\cL(x)$.
For large enough $n$, each of the three words
$\sigma^n(a),\sigma^n(b),\sigma^n(c)$ contains every letter of $A$.
This shows that every letter is extendable in $\cL(\sigma)$
and thus, using Exercise~\ref{exerciseDefSubstitution}, that
$\sigma$ is a substitution.
\end{proof}

The converse of Proposition~\ref{propositionPrimitiveGrowing}
is not true since, for example, $a\mapsto aba,b\mapsto bb$
is a growing substitution which is not primitive.

As a consequence of Proposition~\ref{propositionPrimitiveGrowing},
the shift space $X(\sigma)$ generated by a primitive morphism $\sigma$
on at least two letters, is a substitution
shift. 

The following result is well known. 
\begin{proposition}\label{propositionPrimitiveMinimal}
A primitive substitution shift 
is minimal.
\end{proposition}
\begin{proof}
Let $\sigma:A^*\to A^*$ be a primitive morphism.
Let $n\ge 1$ be such that every $b\in A$ occurs
in every  $\sigma^n(a)$ for $a\in A$.
By Proposition~\ref{propositionMinimalUR}, it is enough to prove
that $\cL(\sigma)$ is uniformly recurrent. Let $u\in \cL(\sigma)$.
Let $a\in A$ and $m\ge 1$ be such that $u$ is a factor of $\sigma^m(a)$.
Then $u$ is a factor of every $\sigma^{n+m}(b)$ and thus
a factor of every word of $\cL(\sigma)$ of length $2|\sigma|^{n+m}$.
Therefore $\cL(\sigma)$ is uniformly recurrent.
\end{proof}
The converse of Proposition~\ref{propositionPrimitiveMinimal}
is not  true.  A well-known example
is the \emph{binary
Chacon substitution} $\sigma:0\mapsto 0010,1\mapsto 1$.
\index{subject}{Chacon!binary!substitution}%
\index{subject}{substitution!binary Chacon}%
\index{names}{Chacon, Rafael V.}%
This substitution is not primitive but the corresponding shift space
is minimal (see Exercise~\ref{exerciseChaconMinimal}). 
The converse is true however for  a substitution 
which is growing 
(Exercise~\ref{exerciseGrowingMinimal}).

When $\sigma$ is primitive,
 we easily check that, for all $k\geq 1$, $\mathcal{L} (\sigma ) = \mathcal{L} (\sigma^k )$ and that $\cL(\sigma) = \mathcal{L} (\sigma^\omega (a))$ for any letter $a \in A$ on which $\sigma$ is right-prolongable. 
In particular, we also have that $\Omega(x) 
%= \Omega(\sigma^\omega(a) 
= X (\sigma ) = X (\sigma^k )$ for all $x \in X (\sigma )$ and $k\geq 1$. 
We say that $X (\sigma )$ is a {\em primitive substitution shift}.
\index{subject}{primitive!substitution!shift}%
\index{subject}{shift space!primitive substitutive}%
A sequence $x$ is {\em primitive substitutive}
\index{subject}{primitive!substitutive!sequence}%
 if it is substitutive with
respect to a primitive substitution.
The subshift $\Omega (x)$ that it generates is then
primitive substitutive.
Since a letter-to-letter morphism is a morphism of dynamical systems, we have the
following corollary of Proposition~\ref{propositionPrimitiveMinimal}.
\begin{corollary}
A primitive substitutive shift is minimal.
\end{corollary}  
%It is clearly minimal.
%We say $\sigma $ is {\em aperiodic}\index{subject}{aperiodic!substitution}\index{subject}{substitution!aperiodic} whenever $(X (\sigma ), S)$ is aperiodic.

The \emph{composition matrix}\index{subject}{composition!matrix}\index{subject}{morphism!composition matrix}\index{subject}{matrix!composition} of 
a morphism $\sigma:A^*\rightarrow B^*$
is the $B\times A$-matrix $M(\sigma)$ defined for every $a\in A$ and $b\in B$ by
\begin{displaymath}
M(\sigma)_{b,a}=|\sigma(a)|_b
\end{displaymath}
where $|\sigma(a)|_b$ denotes the number of occurrences of the letter $b$ 
in the word $\sigma(a)$. 
The transposed matrix of $M(\sigma)$\index{symbols}{M@$M(\sigma)$} is called the \emph{incidence matrix}
\index{subject}{incidence!matrix!of morphism}%
\index{subject}{morphism!incidence matrix}%
\index{subject}{matrix!incidence}%
 of $\sigma$. Thus, on the composition matrix, the columns correspond to
the words $\sigma(a)$ (in the sense that the column of index
$a$ is the vector $(|\sigma(a)|_b)_{b\in B}$) while, in the incidence matrix,
this role is played by the rows.

Such a matrix gives, as we shall see, important information on the morphism.
(although distinct morphisms may well have the same composition matrix).
An important property is that  $\sigma:A^*\to B^*$ and
$\tau:B^*\to C^*$ are morphisms, then
\begin{equation}
M(\tau\circ\sigma)=M(\tau) M(\sigma).\label{eqMatrix}
\end{equation}
Note that, with incidence matrices, the relative order of $\sigma,\tau$
is reversed since $M(\tau\circ\sigma)^t=M(\sigma)^t M(\tau)^t$.

The following mnemonic can be useful since there is a risk
of confusion between the composition matrix and the incidence
matrix. The morphism is read on the Columns
of the Composition matrix (note the initial letter C). Next, Composition
matrices muliply in the same order as the Composition of the morphisms.

When $A=B$, the composition matrix is a square $A\times A$-matrix.
It follows for \eqref{eqMatrix} that for every $k\ge 0$, one has
\begin{equation}
M(\sigma^k)=M(\sigma)^k\label{eqMatrix2}
\end{equation}

For example, the composition matrix of the  substitution 
$\sigma:0\mapsto 01,1\mapsto 00$ is

\begin{displaymath}
M(\sigma)=\begin{bmatrix}1&2\\1&0\end{bmatrix}.
\end{displaymath}

We easily derive from \eqref{eqMatrix2} the following statement.
\begin{proposition}
The morphism $\sigma$ is primitive if and only if the 
matrix $M(\sigma)$ is primitive.
\end{proposition}

\begin{example}\label{exampleFibonacci2}
The Fibonacci morphism is primitive. Indeed, $\varphi^2(a)=aba$
and $\varphi^2(b)=ab$ both contain $a$ and $b$. Accordingly
\begin{displaymath}
M(\varphi)=\begin{bmatrix}1&1\\1&0\end{bmatrix},\quad
M(\varphi)^2=\begin{bmatrix}2&1\\1&1\end{bmatrix}.
\end{displaymath}
Thus the Fibonacci shift is minimal.
\end{example}

Let $\varphi:A^*\rightarrow A^*$ be a primitive morphism. 
Then, the matrix $M$ is a primitive nonnegative
matrix. Thus, we may use the Perron Frobenius Theorem
concerning the dominant eigenvalue $\lambda_M$ of $M$.
(see Appendix~\ref{appendixLinearAlgebra}).

We illustrate the use of this theorem with the following example.
\begin{example}
Let $M=\begin{bmatrix}1&1\\1&0\end{bmatrix}$. The eigenvalues of
$M$ are $\lambda=(1+\sqrt{5})/2$
(called the \emph{golden mean}\index{subject}{golden mean!number})
and $\hat{\lambda}=(1-\sqrt{5})/2$. 
Since $\hat{\lambda}<\lambda$, we have $\lambda_M=\lambda$.
As $\lambda^2=\lambda+1$, a left eigenvector corresponding
to $\lambda$ is $\begin{bmatrix}\lambda& 1\end{bmatrix}$.
The sequence $(M^n/\lambda_M^n)$ tends to the matrix
\begin{displaymath}
\frac{1}{1+\lambda^2}\begin{bmatrix}\lambda^2&\lambda\\\lambda&1\end{bmatrix}
=\frac{1}{1+\lambda^2}\begin{bmatrix}\lambda\\1\end{bmatrix}
\begin{bmatrix}\lambda&1\end{bmatrix}.
\end{displaymath}
\end{example}
\begin{example}
Let $\varphi$ be the Thue-Morse morphism. Then
\begin{displaymath}
M=\begin{bmatrix}1&1\\1&1\end{bmatrix}.
\end{displaymath}
The dominant eigenvalue is equal to $2$ and the corresponding row eigenvector is
$v=\begin{bmatrix}1/2 ,&1/2\end{bmatrix}$.
\end{example}
Thus, if $M$ is the incidence matrix of a primitive morphism, we have
for any $a\in A$,
\begin{equation}
\lim_{n\rightarrow \infty} |\varphi^{n+1}(a)|/|\varphi^{n}(a)|=\lambda_M.
\label{equationLambda}
\end{equation}
since  $|\varphi^n(a)|/\lambda_M^n$ is the sum of
coefficients of the row of index $a$
of $M^n/\lambda_M^n$.

We now prove the following important property of primitive morphisms.
For a morphism $\varphi$, we denote $\langle\varphi\rangle=\min_{a\in A}|\varphi(a)|$ and $|\varphi|=\max_{a\in A}|\varphi(a)|$.

\begin{proposition}\label{propositionBalance}
  For every primitive morphism $\varphi$, the sequence of quotients $(|\varphi^n|/\langle\varphi^n\rangle)_n$
  is bounded.
  \end{proposition}
For this we need the following lemma (which will also be used again later).
\begin{lemma}\label{lemmaFabien}
  Let $M$ be a primitive matrix with dominant eigenvalue $\lambda_M$
  and positive left and right eigenvectors $x=(x_a)$, $y=(y_b)$
  such that $\sum_{a\in A}x_a=1$ and $\sum_{a\in A}x_ay_a=1$.
 There are $c>0$ and $\tau<\lambda_M$
such that for every $a,b\in A$ and $n\ge 1$, we have
\begin{eqnarray}
\bigl| |\varphi^n(a)|_b-\lambda_M^ny_ax_b\bigr|&\le& c\tau^n\label{eqFabien1},\\ 
\bigl||\varphi^n(a)|-\lambda_M^ny_a\bigr|&\le& c\Card(A)\tau^n\label{eqFabien2},\\ 
\bigl||\varphi^n(a)|_b-|\varphi^n(a)|x_b\bigr|&\le& c(1+\Card(A))\tau^n \label{eqFabien3}.
\end{eqnarray}
\end{lemma}
\begin{proof}
The first equation results directly from \eqref{eqGeometricRate}.
For the second one, we write, using the triangular inequality
and \eqref{eqFabien1}
\begin{eqnarray*}
\bigl||\varphi^n(a)|-\lambda_M^n y_a\bigr|&=&|\sum_{b\in A}|\varphi^n(a)|_b-\lambda_M^n\sum_{b\in A}y_ax_b|\\
&\le&\sum_{b\in A}||\varphi^n(a)|_b-\lambda_M^n y_ax_b|\\
&\le&c\Card(A)\tau^n
\end{eqnarray*}
which is \eqref{eqFabien2}. Finally, using \eqref{eqFabien1} and \eqref{eqFabien2}, we obtain, since $x_b\le 1$,
\begin{eqnarray*}
\bigl||\varphi^n(a)|_b-|\varphi^n(a)|x_b\bigr|&=&\bigl||\varphi^n(a)|_b-\lambda_M^ny_ax_b+
\lambda_M^ny_ax_b-|\varphi^n(a)|x_b\bigr|\\
&\le& 
||\varphi^n(a)|_b-\lambda_M^ny_ax_b|+x_b||\varphi^n(a)|-\lambda_M^ny_a|\\
&\le & c(1+\Card(A))\tau^n
\end{eqnarray*}
which is \eqref{eqFabien3}.
\end{proof}

\begin{proofof}{of Proposition~\ref{propositionBalance}}
  Let $M$ be the incidence matrix of $\varphi$.
  Equation \eqref{eqFabien2} shows that for every $a\in A$, the sequence 
  $(|\varphi^n(a)|/\lambda_M^n)_n$ converges to $y_a>0$. This implies
  that for every $a,b\in A$, the sequence of quotients $(|\varphi^n(a)|/|\varphi^n(b)|)_n$
  converges to $y_a/y_b$ and thus proves the statement.
\end{proofof}
 We say that the shift space $X$ is {\em linearly
recurrent}\index{subject}{linearly!recurrent!subshift}\index{subject}{recurrent!linearly!sequence} (LR) (with constant $K\ge 0$) if it is minimal
and if for all $u$ in  $\cL(X)$ and for all right return words $w$
to $u$ in $X$ we have $|w| \leq K|u|$. 
%\end{definition}

\begin{proposition}\label{propositionPrimitiveSubstitutiveIsLR}
A primitive substitution shift is linearly recurrent.
\end{proposition}
\begin{proof}
Let $\sigma:A^*\to A^*$ be a primitive substitution and let $X=X(\sigma)$
be the corresponding shift space. Set as usual $|\sigma|=\max_{a\in A}|\sigma(a)|$
and $\langle \sigma\rangle=\min_{a\in A}|\sigma(a)|$. 

Since $\sigma$ is
primitive, it follows from Proposition~\ref{propositionBalance}
that there is a constant $k$ such that 
\begin{equation}
|\sigma^n|\le k\langle \sigma^n\rangle\label{eqBalance}
\end{equation}
for all $n\ge 1$. 

The substitution $\sigma$ being primitive, the sequence $(\langle\sigma^{n}\rangle)_n$ is strictly increasing. 
Thus, for any $w\in\cL(X)$ there is an integer $n$ such that 
$\langle\sigma^{n-1}\rangle < |w|\le\langle\sigma^{n}\rangle$. 
The right inequality implies that $w$ is a factor of $\sigma^n(ab)$ for some $a,b\in A$. 
%Therefore $w$ appears in all words of length $2|\sigma^{n+1}|$.
%If $n$ is chosen minimal, we have $\langle\sigma^{n-1}\rangle< |w|$.
Every return word to $w$ is then a factor 
of a return word to  $\sigma^n(ab)$. Let
$R$ be the maximal length of return words to words of length $2$.
Then, for every return word $u$ to $w$, we have
\begin{eqnarray*}
|u|&\le&R|\sigma^{n}|\le R k\langle\sigma^{n}\rangle \leq R k |\sigma | \langle\sigma^{n-1}\rangle\\
&<&R k |\sigma||w|.
\end{eqnarray*}
 This shows that $X$ is linearly recurrent
with constant $K=kR|\sigma|$.
\end{proof}
We have seen an illustration of this result for the Fibonacci shift
(Exercise~\ref{exerciseFibonacciLR}).
We add an illustration  on the Thue-Morse shift.
\index{subject}{Thue-Morse!shift space}\index{subject}{shift space!Thue-Morse}
\begin{example}\label{exampleThueMorseLR}
Let $\sigma:a\mapsto ab,b\mapsto ba$ be the Thue-Morse morphism.
According to the above, since $k=1$ for a constant length
morphism, and since the maximal length $R$
of return words of length $2$
is $6$,  the Thue-Morse shift is LR with $K=12$.
\end{example}

\subsection{Circular codes}
We say that a set $U\subset A^+$ of nonempty
words is a {\em code}\index{subject}{code} if every word $w\in A^*$ admits at most one
decomposition in a concatenation of elements of $U$. Thus $U$ is a code
if the submonoid $U^*$ generated by $U$ is isomorphic to the free monoid on $U$.
A prefix code is obviously a code.

A \emph{coding morphism}\index{subject}{coding!morphism} for a code $U\subset A^+$ is a morphism $\varphi:B^*\to A^*$ whose restriction to $B$
is a bijection onto $U$. Thus, in particular $\varphi(B)=U$.
Since $U$ is a code, such a morphism is injective.

The set $U\subset A^+$ is a {\em circular code}\index{subject}{circular code}\index{subject}{code!circular}
whenever it is a code and moreover
if for $p,q\in A^*$ one has 
\begin{equation}
pq,qp\in U^*\Rightarrow p,q\in U^*.\label{eqVeryPure}
\end{equation}

One may visualize this definition as follows. Imagine the word $pq$ written
on a circle (or infinitely repeated). Then $pq$ and $qp$ are two decompositions in words of $U$.
Thus the circular codes are such that the decomposition is unique on a circle
(or, equivalently, on the infinite repetition $\cdots pqpq\cdots$).

\begin{example}
The set $U = \{ 0, 01\}$ is a circular code but $U' = \{010,101\}$ is a code which is not circular. 
\end{example}

\begin{proposition}\label{propositionReturnCircular}
Let $X$ be a shift space. For every $u\in \mathcal{L} (X)$,
the set $\RR_X (u)$ is a circular code. 
\end{proposition}
\begin{proof}
The set $U=\RR_X(u)$ is a prefix code and thus it is a code. Next, assume that
$pq,qp$ belong to $U^*$. Set $p=u_1u_2\cdots u_ns$ with $u_i\in U$
and $s$ a proper prefix of some element of $U$. 
Let $k\ge 1$ be such that $(pq)^k$
is longer than $u$. Since $(pq)^k$ is in $U^*$, the word $u(pq)^k$
ends with $u$ and thus $(pq)^k$  ends with $u$.
Next, since $(qp)^{k+1}=q(pq)^kp$ also ends with $u$ for the same reason, we obtain
that $up$ ends with $u$. 
Hence $s$ belongs to $U^*$ by~\eqref{equationUnitary}.
This forces $s=\varepsilon$ because $U$ is a prefix code.
Thus $p$ belongs to $U^*$, which implies that $q$ is also in $U^*$.
\end{proof}

Circular codes have a property of unique decomposition of
sequences.

\begin{proposition}\label{propositionCircularCodesInjective}
Let $U\subset A^+$ be a finite  code and let $\varphi:B^*\to A^*$ be a coding
morphism for $U$. The following conditions are equivalent.
\begin{itemize}
\item[\rm(i)] The code $U$ is circular.
  \item[\rm(ii)] For every $x\in A^\Z$
there is at most one pair $(k,y)$ such that $x=S^k\varphi(y)$
with $0\le k<|\varphi(y_0)|$.
\end{itemize}
In particular, in this case, the map $\varphi:B^\Z\to A^\Z$
is injective.
\end{proposition}
The proof is left as an exercise (Exercise~\ref{exerciseUniformSync}).
Note that in the case where $U=\RR_X(u)$, the proof is immediate
since the occurrences of $u$ in a sequence determine the factorisation
in words of $U$.

\subsection{Recognizable morphisms}\label{sectionRecognizable}

Let $\varphi:A^*\to B^*$ be a nonerasing morphism.
Let $X$ be a shift space
on the alphabet $A$ and let $Y$ be the  closure
under the shift of $\varphi(X)$.
Then  every $y\in Y$
has a representation as
$y=S^k\varphi(x)$ with $x\in X$ and $0\le k<|\varphi(x_0)|$.
Thus $Y=\{S^k\varphi(x)\mid x\in X, 0\le k <|\varphi(x_0)|\}$.
We say that $\varphi$ is \emph{recognizable}
\index{subject}{recognizable!morphism}\index{subject}{morphism!recognizable}%
on $X$
if every $y\in Y$ has only one such representation.

As an equivalent definition, consider 
the system obtained from $X$
by the tower construction\index{subject}{tower!construction}
 using to the function $x\mapsto |\varphi(x_0)|$,
 already introduced in Section
\ref{sectionRecurrentMinimal}. It
 is the dynamical system $(X^\varphi,T)$ where
\begin{displaymath}
X^\varphi=\{(x,i)\mid x\in X,\ 0\le i<|\varphi(x_0)|\}
\end{displaymath}
\index{symbols}{X@$X^{\tau}$}%
and 
\begin{displaymath}
  T(x,i)=\begin{cases}(x,i+1)&\mbox{ if $i+1<|\varphi(x_0)|$}\\
  (Sx,0)&\mbox{ otherwise.}\end{cases}
\end{displaymath}
\index{symbols}{T@$T(x,i)$}%
The map $\widehat{\varphi}:(x,i)\mapsto S^i\varphi(x)$
\index{symbols}{phi@$\widehat{\varphi}$}% is a morphism 
of dynamical systems from $(X^\varphi,T)$
onto $(Y,S)$. The morphism $\varphi$ is recognizable on $X$
if $\widehat{\varphi}^{-1}(\{y\} )$ has only one element for every $y\in Y$.
Thus $\varphi$ is recognizable on $X$ if and only if $\widehat{\varphi}$
is injective (note that $\varphi$ can be injective
although $\widehat{\varphi}$ is not, see Exercise~\ref{exerciseFiniteToOne}). Consequently  $\hat{\varphi}$ is a homeomorphism. Since $X$ is isomorphic to the system induced
by $(X^\varphi , T)$ on $X\times\{0\}$, we have proved the following useful
statement.
\begin{proposition}\label{propositionRecognizableInduced}
  If $\varphi$ is recognizable on $X$, then
  \begin{enumerate}
  \item The map $\hat{\varphi}:X^\varphi\to Y$ is a homeomorphism.
    \item $\varphi(X)$ is a clopen subset of $Y$.
    \item $(X,T)$ is isomorphic
      to the shift space induced by $Y$ on $\varphi(X)$.
      \end{enumerate}
\end{proposition}
Note that we may consider $X^\varphi$ as a shift space on the alphabet
\begin{displaymath}
  A^\varphi=\{(a,i)\mid a\in A,\ 0\le i<|\varphi(a)|\}.
\end{displaymath}
Indeed, there is a unique embedding $\alpha$ from $(A^\Z)^\varphi$ into
the full shift on $A^\varphi$ such that $(\alpha(x,i))_0=(x_0,i)$.
The image of $(A^\Z)^\varphi$ by $\alpha$ is actually a shift of finite type,
as illustrated in the next example.
\begin{example} Let $X$ be the full shift on $\{u,v,w\}$.
  Consider the morphism $\varphi:u\mapsto a,v\mapsto ba,w\mapsto bba$.
  It is clearly recognizable on $\{ u,v,w \}^\Z$.
  The shift space $\alpha(X^\varphi)$ is the edge shift on the graph $G$
  represented in Figure~\ref{figureXvarphi}. 
  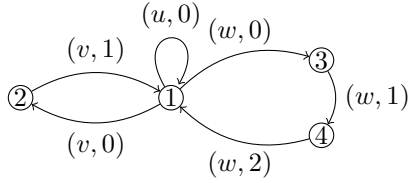
\begin{figure}[hbt]
\centering
\tikzset{node/.style={circle,draw,minimum size=0.1cm,inner sep=0.4pt}}
\tikzset{title/.style={minimum size=0.5cm,inner sep=0pt}}
\tikzstyle{every loop}=[->,shorten >=1pt,looseness=15]
        \tikzstyle{loop left}=[in=60,out=120,loop]
        \begin{tikzpicture}
          \node[node](1)at(0,0){$1$};
          \node[node](2)at(-2,0){$2$};
          \node[node](3)at(2,.5){$3$};
          \node[node](4)at(2,-.5){$4$};

          \draw[above,->](1)edge[loop left]node{$(u,0)$}(1);
          \draw[below,->,bend left](1)edge node{$(v,0)$}(2);
          \draw[above,->,bend left](2)edge node{$(v,1)$}(1);
          \draw[above,->,bend left](1)edge node{$(w,0)$}(3);
          \draw[right,->,bend left](3)edge node{$(w,1)$}(4);
          \draw[below,->,bend left](4)edge node{$(w,2)$}(1);
        \end{tikzpicture}
        \caption{The shift space $\alpha(X^\varphi)$.}\label{figureXvarphi}
\end{figure}

  \end{example}

The following is a reformulation of Proposition~\ref{propositionCircularCodesInjective}.
\begin{proposition}\label{propositionCircularRecognizable}
A nonerasing morphism $\varphi:B^*\to A^*$ is recognizable on 
the full shift  $B^\Z$
if and only if $\varphi$ is injective on $B$ and $\varphi(B)$
is a circular code.
\end{proposition}
As a positive example of application of Proposition~\ref{propositionCircularRecognizable}, we have the case of the Fibonacci morphism.
\begin{example}
The Fibonacci morphism $\varphi:a\mapsto ab,b\mapsto a$ is recognizable 
on the full shift. Indeed, $\{a,ab\}$ is a circular code.
\end{example}
The negative example of the Thue-Morse morphism shows the necessity of the
hypothesis that $\varphi(B)$ is circular in Proposition~\ref{propositionCircularCodesInjective}.

\begin{example}
Let $\tau:a\mapsto ab,b\mapsto ba$ be the Thue-Morse morphism. 
We have $(ab)^\infty=\tau (a^\infty)=S\tau (b^\infty)$.
Thus $\tau$ is not recognizable on $A^\Z$.
But we will see below (Example~\ref{exampleThueMorseRecognizable}) that $\tau$ is recognizable on $X(\tau )$.
\end{example}

The following important theorem will be used
several times.
In particular in the representation
of substitution shifts by sequences of partitions in towers in 
Chapter~\ref{chapterDimensionGroupsPartitions}.
\begin{theorem}[Moss\'e]\label{theoremMosse}
Every primitive aperiodic morphism $\varphi$
 is recognizable on $X(\varphi)$.
\end{theorem}
The proof is given below. 

\begin{example}\label{exampleThueMorseRecognizable}
The Thue-Morse morphism $\tau:a\mapsto ab,b\mapsto ba$ is recognizable
on the Thue-Morse shift $X=X(\tau)$. Indeed, let $y=S^k\tau(x)$
with $0\le k<2$. Every word in $\cL(\tau)$ of length at least
$5$ contains $aa$ or $bb$. Thus $y$ has an infinite
number of occurrences of both $aa$ and $bb$. Consider  an $i>1$ such that
$y_i=y_{i+1}=a$. Then $y_{-k}\cdots y_{i}$ is in $\tau(A^*)$
and thus has even length. Consequently, $i$ and $k$ have opposite parities.
This shows that $k$ is unique and thus also $x$ since $\tau$
has constant length and $ba\not = ab$.
\end{example}
The next example shows that the hypothesis that $\varphi$ is aperiodic is
necessary in Theorem~\ref{theoremMosse}.
\begin{example}
Let $\sigma:a\mapsto ab,b\mapsto ca,c\mapsto bc$. We have $\sigma(abc)=(abc)^2$
and thus $X(\sigma)$ is formed of the three shifts of $(abc)^\infty$.
Therefore  $\sigma$ is primitive but periodic. Since
\begin{displaymath}
(abc)^\infty=\sigma(abc)^\infty=S\sigma(bca)^\infty
\end{displaymath}
the morphism $\sigma$ is not recognizable 
\end{example}

The proof of Theorem~\ref{theoremMosse} uses the following concept.
Let $\varphi:A^*\to B^*$ be a nonerasing morphism, let $X$ be a shift space on the alphabet $A$, let $Y$ be the orbit closure of $\varphi(X)$
under the shift. Let $(p,q)$ be a pair of words
on the alphabet $B$ and $(u,v)$ be a pair of words on the alphabet $A$
such that $uv\in\cL(X)$.
We say that $(p,q)$  is \emph{parsable}\index{subject}{parsable pair}
for $\varphi$ if there exists a pair $(u,v)$ of words on $A$
with $uv\in\cL(X)$ such that

\begin{enumerate}
\item[(i)] $p$ is a suffix of $\varphi(u)$ and
\item[(ii)] $q$ is a prefix of $\varphi(v)$.
\end{enumerate}
We also say  in this case that $(p,q)$ is $(u,v)$-parsable for $\varphi$.

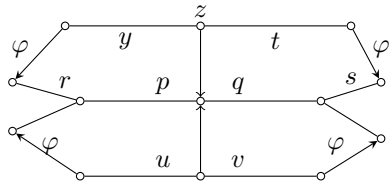
\begin{figure}[hbt]
\centering
\tikzset{node/.style={circle,draw,minimum size=0.1cm,inner sep=0pt}}
	\tikzset{title/.style={minimum size=0.5cm,inner sep=0pt}}
\begin{tikzpicture}
  \node[node](y) at (0.2,3) {};\node[node](z) at (2,3){};\node[title](Z) at (2,3.2){$z$};\node[node](t) at (4,3){};
  \node[title] at (1,2.8){$y$};\node[title] at (3,2.8){$t$};
  \node[node](r) at (-0.5,2.25){};\node[node](ru) at (-0.5,1.6){};
  \node[node](p) at (0.4,2){};\node[node](pq) at (2,2){};\node[node](q)at(3.6,2){};
  \node[node](s) at (4.4,2.25){};
  \node[node](u) at (0.4,1){};\node[node](uv) at (2,1){};\node[node](v) at (3.6,1){};
  \node[node](sv) at (4.4,1.5){};

  \draw (y) edge node {} (z);\draw (z) edge node {} (t);
  \draw[->,>=stealth] (y) edge node{} (r);\node[title] at (-0.4,2.7){$\varphi$};
  \draw (r) edge node {} (p);\node[title] at (0.2,2.25){$r$};
  \draw (p) edge node {} (pq);\node[title] at (1.5,2.2){$p$};
  \draw (pq) edge node {} (q);\node[title] at (2.5,2.2){$q$};
  \draw[->,>=stealth] (t) edge node{} (s);\node[title] at (4.4,2.7){$\varphi$};
  \draw (q) edge node{} (s);\node[title] at (4,2.3){$s$};
  \draw(ru) edge node{} (p);
  \draw(q) edge node {} (sv);
  \draw[->,>=stealth] (u) edge node{} (ru);\node[title] at (0,1.4){$\varphi$};
  \draw (u) edge node {} (uv);\node[title] at (1.5,1.2){$u$};
  \draw (uv) edge node {} (v);\node[title] at (2.5,1.2){$v$};
  \draw[->,>=stealth] (v) edge node{} (sv);\node[title] at (3.8,1.4){$\varphi$};
  \draw[->](z)edge node{}(pq);\draw[->](uv)edge node{}(pq);
\end{tikzpicture}
\caption{A synchronizing pair}\label{figureSynchronizingPair}
\end{figure}

Given a  pair $(p,q)$ which is $(u,v)$-parsable for $\varphi$,
a word $z$ on $A$ such that $pq$ is a factor of $\varphi(z)$ is \emph{synchronized}
\index{subject}{synchronized word} with $(p,q)$ if
there is a  factorization
$z=yt$ such that (see Figure~\ref{figureSynchronizingPair})
\begin{enumerate}
  \item[(i)] $p$ is a suffix of $\varphi(y)$,
\item[(ii)] $q$ is a prefix of $\varphi(t)$,
\item[(iii)] the first letters of $v,t$ are equal.
  \end{enumerate}
Finally, a pair $(p,q)$ which is $(u,v)$-parsable for $\varphi$
is \emph{synchronizing}\index{subject}{synchronizing pair}
for $\varphi$ if  every $z\in\cL(X)$
such that  $pq$ is a factor of $\varphi(z)$ is synchronized with $(p,q)$.
%\begin{figure}[hbt]
%\centering
%\gasset{AHnb=0,Nw=1,Nh=1}
%\begin{picture}(80,25)
%\node(y)(5,20){}\node[ExtNL=y,NLdist=1](yt)(35,20){$z$}\node(t)(65,20){}
%\node(r)(0,12){}\node[ExtNL=y,NLdist=1](ru)(5,8){}\node(p)(10,10){}\node(pq)(35,10){}\node(q)(60,10){}\node[ExtNL=y,NLdist=1](sv)(65,8){}\node(s)(70,12){}
%\node(u)(10,0){}\node(uv)(35,0){}\node(v)(60,0){}

%\drawedge[ELside=r](y,yt){$y$}\drawedge[ELside=r](yt,t){$t$}
%\drawedge(r,p){$r$}\drawedge(ru,p){}\drawedge(p,pq){$p$}\drawedge(pq,q){$q$}
%\drawedge(q,sv){}\drawedge(q,s){$s$}
%\drawedge[AHnb=1,ELside=r](y,r){$\varphi$}\drawedge[AHnb=1](yt,pq){$\varphi$}\drawedge[AHnb=1](t,s){$\varphi$}
%\drawedge[AHnb=1](u,ru){$\varphi$}\drawedge[AHnb=1](uv,pq){$\varphi$}
%\drawedge[AHnb=1](v,sv){$\varphi$}
%\drawedge(u,uv){$u$}\drawedge(uv,v){$v$}
%\end{picture}
%\caption{A synchronizing pair}\label{figureSynchronizingPair}
%\end{figure}
Let $\varphi:A^*\to B^*$ be a morphism, let $X$ be a shift space on $A$
and let $Y$ be the orbit closure of $\varphi(X)$ under the shift.
\begin{proposition}\label{propositionRecognizableInjective}
  The morphism $\varphi$ is recognizable on $X$ if and only
  if there is an integer $L$ such that every pair of words of length
  $L$, which is parsable for $\varphi$, is synchronizing for $\varphi$.
\end{proposition}
\begin{proof}
  Let us first show that the condition is sufficient.

We have to prove 
that for $x,x'\in X$ such that  $S^j\varphi(x)= S^k\varphi(x')$
with $0\le j<|\varphi(x_0)|$ and
$0\le k<|\varphi(x'_0)|$, one has $j=k$ and $x_0=x'_0$. We may suppose that $j\le k$. 

By the hypothesis, for $N$
large enough, the pair $(q,r)$
with $q=\varphi(x_{[-N,0)})$ and $r=\varphi(x_{[0,N)})$ is synchronizing
for $\varphi$.
\begin{figure}[hbt]
\centering
\tikzset{node/.style={circle,draw,minimum size=0.1cm,inner sep=0pt}}
\begin{tikzpicture}
\node[node](x0)at(2,3){};
\node[node](x1)at(4,3){};
\node[node](x2)at(6,3){};
\node[node](y0)at(0,1.5){};
\node[node](y1)at(1,1.5){};
\node[node](y2)at(3.5,1.5){};
\node[node](y3)at(4,1.5){};\node[node](y4)at(7,1.5){};\node[node](y5)at(8,1.5){};
\node[node](z0)at(1,0){};\node[node](z1)at(3.7,0){};
\node[node](x'0)at(4.3,0){};\node[node](z2)at(6,0){};

\draw[above](x0)edge node{$x_{[-N,0)}$}(x1);
\draw[above](x1)edge node{$x_{[0,N)}$}(x2);
\draw[above](y0)edge node{$p$}(y1);
\draw[above](y1)edge node{$q'$}(y2);\draw[above](y2)edge node{$q''$}(y3);
\draw[above](y3)edge node{$r$}(y4);\draw[above](y4)edge node{$s$}(y5);
\draw[left,->](x0)edge node{$\varphi$}(y1);\draw[left,->](x1)edge node{$\varphi$}(y3);
\draw[right,->](x2)edge node{$\varphi$}(y4);
\draw[above](z0)edge node{$x'_{[-M,0)}$}(z1);
\draw[above](z1)edge node{$x'_0$}(x'0);
\draw[above](x'0)edge node{$x'_{[1,M)}$}(z2);
\draw[left,->](z0)edge node{$\varphi$}(y0);
\draw[left,->](z1)edge node{$\varphi$}(y2);\draw[->](x'0)edge node{}(4.5,1.5);
\draw[right,->](z2)edge node{$\varphi$}(y5);
\end{tikzpicture}
\caption{The intersection $\varphi(x)\cap S^{k-j}\varphi(x')$.}
\label{figureIntersection}
\end{figure}
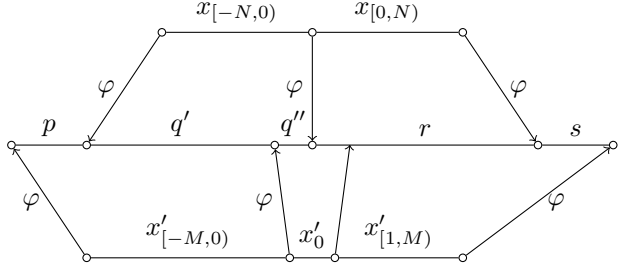

We can also assume that $|\varphi(x_{[-N,0)})|\ge k-j$ and thus
that $q=q'q''$ with $|q''|=k-j$ (see Figure~\ref{figureIntersection}).
For $M\ge 1$ large enough,
we have  $\varphi(x'_{[-M,0)})=pq'$ and
$\varphi(x'_{[0,M)})=q''rs$ for some words $p$ and $s$
(see Figure~\ref{figureIntersection}). Since $(q,r)$ is synchronizing
for $\varphi$, there
is a prefix $t$ of $x'_{[0,M)}$ such that $\varphi(t)=q''$
  and $x_0$ is the first letter of $tx'_{[0,M)}$.  Since
$|q''|=k-j\le k<|\varphi(x'_0)|$, this
forces $|t|=|q''|=0$, $j=k$ and  $x_0=x'_0$.

Conversely, assume that $\varphi$ is recognizable on $X$.
Set $\psi=\hat{\varphi}^{-1}$. Since $\psi:Y\to X^\varphi$ is
a conjugacy between shift spaces (we consider $X^\varphi$ as a shift
space on the alphabet $A^\varphi$),
it is a sliding block code. Thus
there is an integer $L\ge 1$ such that the symbol $(\psi (y))_0$
depends only on $y_{[-L,L)}$.

\begin{figure}[hbt]
\centering
\tikzset{node/.style={circle,draw,minimum size=0.1cm,inner sep=0pt}}
	\tikzset{title/.style={minimum size=0.5cm,inner sep=0pt}}
\begin{tikzpicture}
  \node[node](m) at (0.2,3) {};
  \node[node](a) at (1.8,3){};
  \node[node](n)at(2.2,3){};
  
  \node[node](t) at (4,3){};
  
  \node[node](r) at (-0.5,2.5){};\node[node](ru) at (-0.5,1.6){};
  \node[node](p') at (0.4,2){};\node[node](p'')at(1.4,2){};
  \node[node](pq) at (2,2){};
  \node[node](q'')at(2.4,2){};
  \node[node](q)at(3.6,2){};
  \node[node](s) at (4.4,2.5){};
  \node[node](u) at (0.4,1){};\node[node](uv) at (2,1){};\node[node](v) at (3.6,1){};
  \node[node](sv) at (4.4,1.5){};
  \node[title]at(-2,1){$x$};
  \node[title]at(-2,1.6){$y$};
  \node[title]at(-2,2.5){$y'$};
  \node[title]at(-2,3){$x'$};
  
  \draw(-1.5,1)edge node{}(u);
  \draw(-1.5,3)edge node{}(m);
  \draw(t)edge node{}(5.4,3);
  \draw(v)edge node{}(5.4,1);
  \draw[above](m) edge node {$m$} (a);
  \draw[->](a)edge node{}(p'');
  \draw[above](a) edge node {$a$} (n);
  \draw[->](n)edge node{}(q'');
  \draw[above](n) edge node {$n$} (t);
  \draw[->,>=stealth,left] (m) edge node{$\varphi$} (r);
  \draw[color=blue,line width=1,above](r) edge node {$r$} (p');
  \draw[line width=1,color=violet,above](p') edge node {$p'$} (p'');
  \draw[line width=1,color=violet,above](p'') edge node {$p''$} (pq);
  \draw[line width=1,color=violet,above](pq) edge node {$q'$} (q'');
  \draw[line width=1,color=violet,above](q'') edge node {$q''$} (q);
  \draw[->,>=stealth,right] (t) edge node{$\varphi$} (s);
  \draw [color=blue,line width =1,above](q) edge node{$s$} (s);
  \draw[color=red,line width=1](ru) edge node{} (p);
  \draw[line width=1,color=red](q) edge node {} (sv);
  \draw[->,>=stealth] (u) edge node{} (ru);\node[title] at (0,1.4){$\varphi$};
  \draw (u) edge node {} (uv);\node[title] at (1.5,1.2){$u$};
  \draw (uv) edge node {} (v);\node[title] at (2.5,1.2){$v$};
  \draw[->,>=stealth] (v) edge node{} (sv);\node[title] at (3.8,1.4){$\varphi$};
  \draw[color=red,line width=1](-1.5,1.6)edge node{}(ru);
  \draw[color=red,line width=1](sv)edge node{}(5.4,1.5);
  \draw[color=blue,line width=1](-1.5,2.5)edge node{}(r);
  \draw[color=blue,line width=1](s)edge node{}(5.4,2.5);
  \draw[->](pq)edge node{}(uv);
\end{tikzpicture}
\caption{Proof that the condition is necessary.}\label{figureSufficiencyRecognizability}
\end{figure}
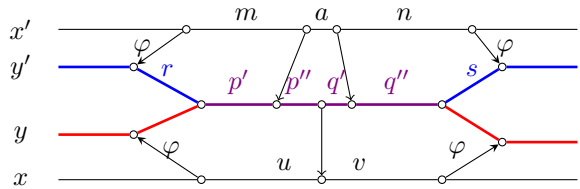

Let $(p,q)$ be pair of words of length $L$
which is  $(u,v)$-parsable for $\varphi$, with $uv\in\cL(X)$.
Let $x\in X$ be such that $u$ is a suffix of $x^-$ and $v$ is a prefix
of $x^+$. Set $y=\varphi(x)$. Let $z\in\cL(X)$ be such that $\varphi(z)=rpqs$.
Let $m,n,p',p'',q',q''$ be words and $a\in A$ be a letter such that $z=man$
with (see Figure~\ref{figureSufficiencyRecognizability})
\begin{displaymath}
  \varphi(m)=rp',\ \varphi(a)=p''q',\ \varphi(n)=q''s,\ p=p'p'',\ q=q'q''
\end{displaymath}
and $q'$ not empty. 

Let $x'\in X$ be such that $x'^-$ ends with $m$ and $x'^+$ begins with $an$.
Set $y'=S^{|p''|}\varphi(x')$ (in such a way that $\psi(y')_0=(a,|p''|)$).
Since $y_{[-L,L)}=y'_{[-L,L)}=pq$, it
follows from the definition of $L$ that $\psi (y)_0=\psi (y')_0$.
But $\psi (y)_0=(b,0)$
where $b$ is the first letter of $v$ while $\psi (y')_0=(a,|p''|)$.
We conclude that $p''$ is empty and $a=b$, which shows that
the pair $(p,q)$ is synchronizing for $\varphi$.
\end{proof}

\begin{proofof}{of Theorem~\ref{theoremMosse}}
  Assume, by contradiction, that $\varphi:A^*\to A^*$ is not
  recognizable on $X=X(\varphi)$. By Proposition~\ref{propositionRecognizableInjective}, this implies that for every $\ell$, there is a pair $(p,q)$
  of words of length $\ell$ which is parsable for $\varphi$ in $\cL(X)$ but not
  synchronizing for $\varphi$.

  Fix an integer $k\ge 1$
  which will be chosen later. By the hypothesis, there
  is for every $n\ge 1$ a word $r_n$   in $\cL_{2k}(X)$
  and a factorization $\varphi^{n-1}(r_n)=u_nv_n$
  such that the pair $(p_n,q_n)$ with $p_n=\varphi(u_n)$ and
  $q_n=\varphi(v_n)$ is not synchronizing for $\varphi$. Thus, there is for every $n\ge 1$
  a word $z_n\in\cL(X)$  such that $pq$ is factor of $\varphi(z_n)$
 and
  \begin{enumerate}
  \item[(i)] $\varphi^{n-1}(z_n)$ is not synchronized with $(p_n,q_n)$.
  \end{enumerate}
  Moreover, choosing $z_n$ of minimal length, we have
  \begin{enumerate}
  \item[(ii)] $z_n=a_nw_nb_n$ with $a_n,b_n\in A$ and $p_nq_n=s_n\varphi^n(w_n)t_n$.   \end{enumerate}
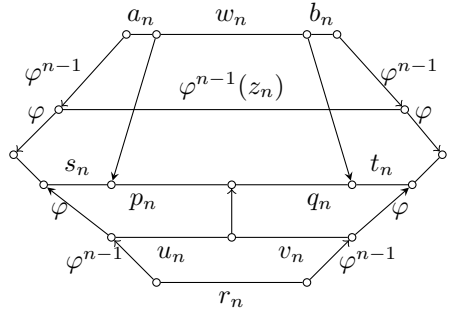
\begin{figure}[hbt]
\centering
\tikzset{node/.style={circle,draw,minimum size=0.1cm,inner sep=0pt}}
	\tikzset{title/.style={minimum size=0.5cm,inner sep=0pt}}
        \begin{tikzpicture}
          \node[node](an)at(.6,4){};
  \node[node](y) at (1,4) {};\node[node](t) at (3,4){};
  \node[node](bn)at(3.4,4){};
  \node[node](phin-1l)at(-.3,3){};\node[node](phin-1r)at(4.3,3){};
  \node[node](r) at (-0.5,2){};\node[node](rr)at(-.9,2.4){};
  \node[node](p) at (0.4,2){};\node[node](pq) at (2,2){};\node[node](q)at(3.6,2){};
  \node[node](ss)at(4.8,2.4){};
  \node[node](s) at (4.4,2){};
  \node[node](u) at (.4,1.3){};\node[node](uv) at (2,1.3){};
  \node[node](v) at (3.6,1.3){};
  %\node[title](zn)at(2,2.7){$z_n$};
  \node[node](tl)at(1,.7){};\node[node](tr)at(3,.7){};

\draw[above](an)edge node{$a_n$}(y);
\draw [above](y) edge node {$w_n$} (t);
\draw[above](t)edge node{$b_n$}(bn);
\draw[left,->](an)edge node{$\varphi^{n-1}$}(phin-1l);
\draw[right,->](bn)edge node{$\varphi^{n-1}$}(phin-1r);
\draw[above](phin-1l)edge node{$\varphi^{n-1}(z_n)$}(phin-1r);
\draw[->,>=stealth,left] (y) edge node{} (p);
\draw[above,->](phin-1l)edge node{$\varphi$}(rr);
\draw(rr)edge node{}(r);\draw(s)edge node{}(ss);
\draw[above,->](phin-1r)edge node{$\varphi$}(ss);
  \draw[above] (r) edge node {$s_n$} (p);
  \draw [below,near start](p) edge node {$p_n$} (pq);
  \draw[below,near end] (pq) edge node {$q_n$} (q);
  \draw[->,>=stealth] (t) edge node{} (q);
  \draw[above] (q) edge node{$t_n$} (s);

  \draw[->,>=stealth,left] (u) edge node{$\varphi$} (r);
  \draw[below] (u) edge node {$u_n$} (uv);
  \draw[below] (uv) edge node {$v_n$} (v);
  \draw[->,>=stealth,right] (v) edge node{$\varphi$} (s);
  \draw[->](uv)edge node{}(pq);
  \draw[left,->](tl)edge node{$\varphi^{n-1}$}(u);
  \draw[right,->](tr)edge node{$\varphi^{n-1}$}(v);
  \draw[below](tl)edge node{$r_n$}(tr);
\end{tikzpicture}

\caption{Proof of Mosse's Theorem.}\label{figureProofMosse}
\end{figure}
Since $\varphi^n(r_n)=p_nq_n$ and $|r_n|=2k$, we have
\begin{displaymath}
  |p_nq_n|\le 2k|\varphi^n| .
\end{displaymath}
Moreover, since $p_nq_n=s_n\varphi^n(w_n)t_n$, we have
\begin{displaymath}
  |p_nq_n|\ge |w_n|\langle \varphi^n\rangle.
\end{displaymath}
Thus, we obtain
\begin{equation}
  |w_n|\le 2k|\varphi^n|/\langle \varphi^n\rangle.\label{eqw_n} 
  \end{equation}
But, the morphism $\varphi$ being primitive, by Proposition \ref{propositionBalance},
the right-hand side of \eqref{eqw_n} is bounded independently of $n$.
Thus, there
is an infinite set $E$ of integers $n$ such that $r_n=r$, $z_n=z$, $a_n=a$,
$b_n=b$ and $w_n=w$
for every $n\in E$. Hence we have the equalities
\begin{eqnarray}
  z&=&awb,\\
\varphi^n(r)&=&s_n\varphi^n(w)t_n. \label{eqphinuv}
  \end{eqnarray}

Consider $n,m\in E$ with $n<m$. We have
\begin{eqnarray*}
  s_m\varphi^m(w)t_m&=&\varphi^m(r)=\varphi^{m-n}(\varphi^n(r))\\
  &=&\varphi^{m-n}(s_n)\varphi^m(w)\varphi^{m-n}(t_n).
\end{eqnarray*}
Suppose that $s_m\ne\varphi^{m-n}(s_n)$ for arbitrary large values of $k$. Assume first that
$|s_m|>|\varphi^{m-n}(s_n)|$.
Then, by Proposition~\ref{propositionOverlappingFactor},
    the word $\varphi^m(w)$ is periodic
of period $|s_m|-|\varphi^{m-n}(s_n)|\le|s_m|$. This implies that $\cL(X)$
contains words of exponent larger than $\lfloor e\rfloor$ with
\begin{displaymath}
  e=|\varphi^m(w)|/|s_m|\ge |p_mq_m|/|s_m|
  \ge 2k\langle\varphi^m\rangle/|\varphi^m|
\end{displaymath}
which tends to infinity with $k$.
Thus we may  choose $k$ to contradict Proposition~\ref{propositionLRhasLinearComplexity} since, by Proposition~\ref{propositionPrimitiveSubstitutiveIsLR}, the shift space $X$ is linearly recurrent.

In the case where $|s_m|<|\varphi^{m-n}(s_n)|$, we find that
$\varphi^m(w)$ is periodic of period $|\varphi^{m-n}(s_n)|-|s_m|<|\varphi^{m-n}(s_n)|$ and thus that $\cL(X)$ contains word of exponent larger than
$\lfloor e\rfloor$ with
\begin{displaymath}
  e=|\varphi^m(w)|/|\varphi^{m-n}(s_n)|\ge |p_mq_m|/|\varphi^m(a)|
  \ge 2k\langle\varphi^m\rangle/|\varphi^m|
\end{displaymath}
whence a contradiction again.
\begin{figure}[hbt]
\centering
\tikzset{node/.style={circle,draw,minimum size=0.1cm,inner sep=0pt}}
	\tikzset{title/.style={minimum size=0.5cm,inner sep=0pt}}
        \begin{tikzpicture}
          \node[node](r)at(2,2){};\node[node](rr)at(3,2){};
          \node[node](pn)at(1,1){};%\node[node](pnqn)at(2.5,1){};
          \node[node](qn)at(4,1){};
          \node[node](um)at(.3,0){};\node[node](umvm)at(2.5,0){};
          \node[node](vm)at(4.7,0){};
          \node[node](pm)at(0,-.5){};\node[node](pmqm)at(2.5,-.5){};
          \node[node](qm)at(5,-.5){};
          \draw[above](r)edge node{$r$}(rr);
          \draw[left,->](r)edge node{$\varphi^n$}(pn);
          \draw[right,->](rr)edge node{$\varphi^n$}(qn);
          \draw[above,color=red](pn)edge node{$\varphi^n(r)$}(qn);
          %\draw[above,color=red](pnqn)edge node{$q_n$}(qn);
          \draw[left,->](pn)edge node{$\varphi^{m-n-1}$}(um);
          %\draw[left,->](pnqn)edge node{}(umvm);
          \draw[right,->](qn)edge node{}(vm);
          \draw[above,color=green](um)edge node{$u_m$}(umvm);
          \draw[above,color=green](umvm)edge node{$v_m$}(vm);
          \draw[left,->](um)edge node{$\varphi$}(pm);
          \draw[left,->](umvm)edge node{}(pmqm);
          \draw[right,->](vm)edge node{$\varphi$}(qm);
          \draw[below](pm)edge node{$p_m$}(pmqm);
          \draw[below](pmqm)edge node{$q_m$}(qm);

          \node[node](a)at(8,2){};\node[node](w)at(8.5,2){};
          \node[node](b)at(10.5,2){};\node[node](bb)at(11,2){};
          \node[node](phia)at(7,1.3){};\node[node](sn)at(7.5,1){};
          \node[node](phinw)at(8,1){};\node[node](tn)at(10.5,1){};
          \node[node](phib)at(11,1){};\node[node](phibb)at(11.5,1.3){};
          \node[node](0)at(6,.3){};\node[node](sm)at(7,0){};
          \node[node](1)at(8,0){};\node[node](tm)at(10.5,0){};
          \node[node](2)at(11.5,0){};\node[node](3)at(12.3,.3){};

          \draw[above](a)edge node{$a$}(w);\draw[above](w)edge node{$w$}(b);
          \draw[above](b)edge node{$b$}(bb);
          \draw[left,->](a)edge node{$\varphi^n$}(phia);
          \draw[left,->](w)edge node{}(phinw);
          \draw[left,->](b)edge node{}(tn);
          \draw[right,->](bb)edge node{$\varphi^n$}(phibb);
          \draw[above](phia)edge node{}(sn);
          \draw[above,color=red](sn)edge node{$s_n$}(phinw);
          \draw[above,color=red](phinw)edge node{$\varphi^n(w)$}(tn);
          \draw[above,color=red](tn)edge node{$t_n$}(phib);
          \draw[above](phib)edge node{}(phibb);
          \draw[left,->](phia)edge node{}(0);
          \draw[left,->](sn)edge node{}(sm);
          \draw[left,->](phinw)edge node{}(1);
          \draw[left,->](tn)edge node{}(tm);
          \draw[left,->](phib)edge node{}(2);
          \draw[right,->](phibb)edge node{$\varphi^{m-n-1}$}(3);
          \draw[below](0)edge node{}(sm);
          \draw[below,color=green](sm)edge node{$\varphi^{m-n-1}(s_m)$}(1);
          \draw[below,color=green](1)edge node{$\varphi^{m-1}(w)$}(tm);
          \draw[below,color=green](tm)edge node{$\varphi^{m-n-1}(t_m)$}(2);
          \draw[below](2)edge node{}(3);
        \end{tikzpicture}
        \caption{The case $s_m=\varphi^{m-n}(s_n)$.}\label{figuresm=phin-m}
  \end{figure}
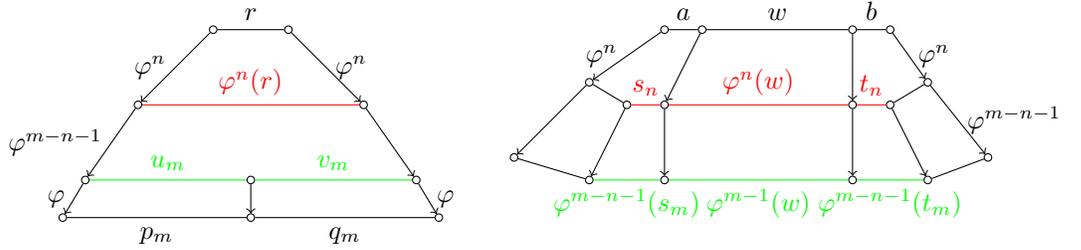

Thus $s_m=\varphi^{m-n}(s_n)$. Now we have by Equation~\eqref{eqphinuv}
$s_n\varphi^n(w)t_n=\varphi^n(r)$ and
thus, applying $\varphi^{m-n-1}$ to both sides (see Figure~\ref{figuresm=phin-m}),
\begin{eqnarray*}
\varphi^{m-n-1}(s_n)\varphi^{m-1}(w)\varphi^{m-n-1}(t_n)&=&\varphi^{m-1}(r)=u_mv_m
\end{eqnarray*}
Let $y,t$ be such that $\varphi^{m-1}(w)=yt$ with
\begin{equation}
  u_m=\varphi^{m-n-1}(s_n)y,\quad
  v_m=t\varphi^{m-n-1}(t_n).\label{eqContrad}
\end{equation}
Applying $\varphi$ on both sides of \eqref{eqContrad}, we obtain
\begin{displaymath}
  p_m=\varphi^{m-n}(s_n)\varphi(y),\quad q_m=\varphi(t)\varphi^{m-n}(t_n).
\end{displaymath}

Since  $\varphi^{m-1}(z)=\varphi^{m-1}(a)yt\varphi^{m-1}(b)$
we find that
\begin{enumerate}
\item[(i)] $p_m$ is a suffix of $\varphi^{m-1}(a)\varphi(y)$,
\item[(ii)] $q_m$ is a prefix of  $\varphi(t)\varphi^{m-1}(b)$,
\item[(iii)] the first letters of $t$ and $v_m$ are equal.
  \end{enumerate}
  This shows that $\varphi^{m-1}(z)$ is synchronized with $(p_m,q_m)$,
a contradiction.

  \end{proofof}
The definition of recognizability for a primitive morphism
can be formulated in terms of a fixed point of the morphism
(see Exercise~\ref{exerciseMosseOriginal}).

A striking feature of Mosse Theorem is that it holds for non injective
morphisms. We illustrate this in the following example.

\begin{example}
  The morphism $\varphi:a\mapsto abc$, $b\mapsto a$, $c\mapsto d$, $d\mapsto bcd$
  is primitive and aperiodic. Thus it is recognizable on  $X(\varphi)$.
  It is not injective since $\varphi(ac)=\varphi(bd)=abcd$. However,
  one can check that the restriction of $\varphi$ to $\cL(X)$
  is injective.
  \end{example}
\subsection{Block presentations}
We will need later to associate to a primitive morphism $\varphi$
and an integer $k\ge 1$ the $k$-block \emph{presentation}\index{subject}{higher block!presentation!of a morphism}\index{subject}{presentation!higher block!of morphism}
 of $\varphi$
denoted 
$\varphi_k$.\index{symbols}{phik@$\varphi_k$}

Let $\varphi:A^*\rightarrow A^*$ be a primitive morphism
and let $(X,S)$ be its associated shift space. For
$k\ge 1$, consider an alphabet $A_k$
in one-to-one correspondence by $f: \cL_k(X)\rightarrow A_k$ with the set $\cL_k(X)$
of factors (or blocks) of length $k$ of $X$. The map
$f$ extends naturally  to a map, still denoted $f$,
from $\cL_{k+n}(X)$ to $\cL_{n+1}(X^{(k)})$
 defined for $n\ge 0$ by 
\begin{equation}
  f(a_1a_2\cdots a_{k+n})=f(a_1\cdots a_k)f(a_2\cdots a_{k+1})\cdots f(a_{n+1}\cdots a_{k+n}).\label{equationfn+k}
\end{equation}

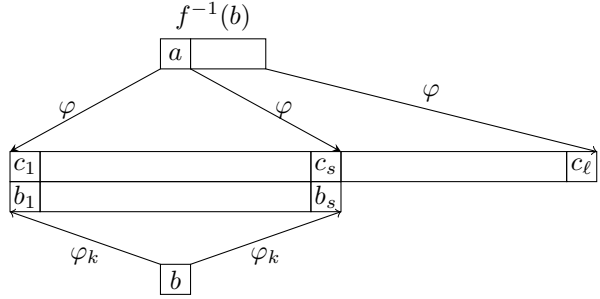
\begin{figure}[hbt]
\centering
\tikzset{node/.style={draw,minimum size=0.4cm,inner sep=0pt}}
	\tikzset{title/.style={minimum size=0.5cm,inner sep=0pt}}
\begin{tikzpicture}
\node[title](fb)at (2.5,2.5){$f^{-1}(b)$};
\node[node](a) at (2,2){$a$};
\node[node,text width=1cm](f) at (2.7,2){};
\node[node](c1)at (0,.5){$c_1$};
\node[node,text width=3.6cm](c)at (2,.5){};
\node[node](cs)at (4,.5){$c_s$};
\node[node,text width=3cm](d)at (5.7,.5){};
\node[node](cl)at (7.4,.5){$c_\ell$};
\node[node](b1)at(0,.1){$b_1$};\node[node,text width=3.6cm](bh)at(2,.1){};
\node[node](bs)at(4,.1){$b_s$};
\node[node](b)at(2,-1){$b$};

\draw[->,above,left,>=stealth](1.8,1.8)-- node{$\varphi$}(-.2,.7);
\draw[->,above,right,>=stealth](2.2,1.8)-- node{$\varphi$}(4.2,.7);
\draw[->,above](3.2,1.8)-- node{$\varphi$}(7.6,.7);
\draw[->,below](1.8,-.8)-- node{$\varphi_k$}(-.2,-.1);
\draw[->,below](2.2,-.8)-- node{$\varphi_k$}(4.2,-.1);
\end{tikzpicture}
\caption{The morphism $\varphi_k$.}\label{figureBlockPresentation}
\end{figure}

We define a morphism
$\varphi_k:A_k^*\rightarrow A_k^*$ as follows. Let
$b\in A_k$ and let $a$ be the first letter of $f^{-1}(b)$ (see Figure~\ref{figureBlockPresentation}).
Set $s=|\varphi(a)|$.
To compute $\varphi_k(b)$, we first compute
the word $\varphi(f^{-1}(b))=c_1c_2\cdots c_\ell$.
Note that $\ell\ge |\varphi(a)|+k-1=s+k-1$.
 We set 
 $$
 \varphi_k(b)=b_1b_2\cdots b_s
 $$ 
 where
\begin{displaymath}
b_1=f(c_1c_2\cdots c_k),\ b_2=f(c_2c_3\cdots c_{k+1}),\ldots,\
b_s=f(c_s\cdots c_{s+k-1}).
\end{displaymath}
In other terms, $\varphi_k(b)$ is the prefix of length $s=|\varphi(a)|$
of $f\circ\varphi\circ f^{-1}(b)$ where $f$ is the map defined by \eqref{equationfn+k} (see Figure~\ref{figureBlockPresentation}).

\begin{example}\label{exampleBlockPresentation}
Let $\varphi:a\mapsto ab,b\mapsto a$ be the Fibonacci morphism. We have
$\cL_2(X)=\{aa,ab,ba\}$. Set $A_2=\{x,y,z\}$ and let
$f:aa\mapsto x,ab\mapsto y,ba\mapsto z$. Since $f\circ\varphi\circ f^{-1}(x)=f(\varphi(aa))=f(abab)=yzy$, 
keeping the prefix of length $|\varphi(a)|=2$, we have
$\varphi_2(x)=yz$. Similarly, we have $\varphi_2(y)=yz$ and
$\varphi_2(z)=x$.
\end{example}

Let $\pi:A_k^*\rightarrow A^*$ be the morphism defined by
$\pi(b)=a$ where
$a$ is  the first letter of $f^{-1}(b)$.

Then we have for each $n\ge 1$ the following commutative diagram
which expresses the fact that $\varphi_k^n$ is the
counterpart of $\varphi^n$ for $k$-blocks. 

\begin{equation}
\begin{CD}
A_k^+ @>{\varphi_k^n}>> A_k^+\\
@VV{\pi}V               @VV{\pi}V\\
A^+ @>{\varphi^n}>> A^+
\end{CD}\label{DiagramExtension}
\end{equation}

Indeed, we have $\varphi \circ \pi(b)=\pi \circ \varphi_k(b)$
for every $b\in A_k$ by definition of $\varphi_k$.
Since $\varphi\pi$ and $\pi\varphi_k$ are morphisms,
this implies $\varphi\circ\pi=\pi\circ\varphi_k$
and thus $\varphi^n\circ\pi=\pi\circ\varphi_k^n$ for all $n\ge 1$.
This proves~\eqref{DiagramExtension}.
In particular, since $\pi$ is length preserving, we have
\begin{equation}
|\varphi_k^n(b)|=|\varphi^n(a)|\label{eqLgphikn}
\end{equation}
for $n\ge 1$ and $a=\pi(b)$.

\begin{figure}[hbt]
\centering
\tikzset{node/.style={draw,minimum size=0.4cm,inner sep=0pt}}
	\tikzset{title/.style={minimum size=0.5cm,inner sep=0pt}}
\begin{tikzpicture}

\node[node](a) at (2,2){$a$};
\node[node,text width=1cm](f) at (2.7,2){$\ \rho(u)$};\node[node,text width=.6cm](u)at(3.5,2){};

\node[node,text width=3.6cm](c)at (2,.5){$\quad\varphi(a)$};
\node[node,text width=3cm](d)at (5.3,.5){};
\node[node,text width=1.6cm](cm)at (7.6,.5){};
\node[node,text width=3.6cm](bh)at(2,.1){};
\node[node,text width=3cm](bt)at(5.3,.1){};
\node[node](b)at(2,-1.5){$b$};\node[node,text width=1cm](v)at(2.7,-1.5){$\quad v$};

\draw[->,above,left](1.8,1.8)-- node{$\varphi$}(.2,.7);
\draw[->,above,right](2.2,1.8)-- node{$\varphi$}(3.8,.7);
\draw[->,above,right](3.2,1.8)-- node{$\varphi$}(6.8,.7);
\draw[->,above](3.8,1.8)-- node{$\varphi$}(8.4,.7);
\draw[->,below](1.8,-1.3)-- node{$\varphi_k$}(.2,-.1);
\draw[->,below](2.2,-1.3)-- node{$\varphi_k$}(3.8,-.1);
\draw[->,below](3.2,-1.3)-- node{$\varphi_k$}(6.8,-.1);
\end{tikzpicture}
\caption{Comparing $\varphi_k(f(au))$ and $\varphi(au)$.}\label{figureBlockPresentation2}
\end{figure}
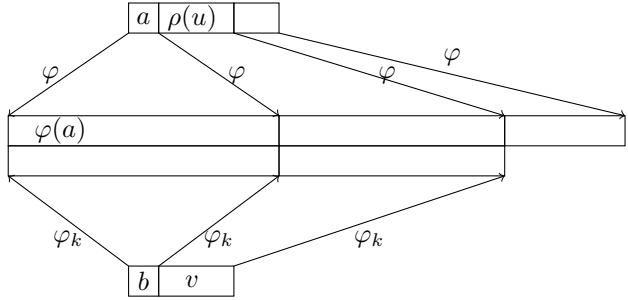

In the following, we denote $u\le v$ to express that the word $u$ is a prefix of $v$.
\begin{proposition}\label{propositionPresentation}
We have for every $u\in \cL(X)$ of length at least $k$,
\begin{equation}
\varphi_k(f(u))\leq f(\varphi(u)).\label{inequationphik}
\end{equation}
\end{proposition}
\begin{proof}
 For  a word $w$
of length at least $n$, we denote by $\Pref_n(w)$ its prefix of length
$n$ and we set $\rho(u)=\Pref_{|u|-k+1}(u)$ for $u$ of length at least $k$.
For $u\in\cL(X)$, set $\ell(u)=|\varphi(\rho(u))|$.
By definition, $\varphi_k(f(u))$ is, for every $u\in\cL_k(X)$,
the prefix of length $|\varphi(a)|$ of $f(\varphi(u))$
where $a$ is the first letter of $u$. We prove that
this property extends to longer words $u$ and that,
for every $u\in\cL_{\ge k}(X)$, one has
\begin{equation}
\varphi_k(f(u))=\Pref_{\ell(u)}f(\varphi(u)).\label{equationphik}
\end{equation}

Indeed, arguing by induction on the length of $u$,
consider  $a\in A$ and $u\in\cL_{\ge k}(X)$. Set $f(au)=bv$
with $b\in A_k$ (see Figure~\ref{figureBlockPresentation2}). Then, since $\rho(au)=a\rho(u)$, since $f(\varphi(au))=\varphi(a)f(\varphi(u))$
and since $f^{-1}(b)=\Pref_k(au)$,
\begin{eqnarray*}
\Pref_{\ell(au)}f(\varphi(au))&=&\Pref_{|\varphi(a)|+\ell(u)}f(\varphi(au))\\
&=&\Pref_{|\varphi(a)|}f(\varphi(au))
\Pref_{\ell(u)}f(\varphi(u))\\
&=&\Pref_{|\varphi(a)|}f(\varphi(f^{(-1)}(b))\Pref_{\ell(u)}f(\varphi(u))\\
&=&\varphi_k(b)\varphi_k(v)=\varphi_k(bv)
\end{eqnarray*}
proving~\eqref{equationphik} and thus \eqref{inequationphik}.

\end{proof}

\begin{proposition}\label{propositionphi_kPrimitive}
When $\varphi$ is primitive, then for $k\ge 1$,
\begin{enumerate}
\item the morphism $\varphi_k$ is primitive,
\item for every $b\in A_k$, $u=f^{(-1)}(b)$ and $n\ge 1$, $\varphi_k^n(b)\le f(\varphi^n(u))$,
\item the shift space associated to $\varphi_k$ is
the $k$-th higher block presentation of $X$.
\end{enumerate}
\end{proposition}
\begin{proof}
1. Since $\varphi$ is primitive, there is an integer
$n$ such that for every $a\in A$, the word $\varphi^n(a)$
contains all the words of $\cL_k(\varphi)$ as factors.
Then for every $b\in A_k$, the word $\varphi_k^n(b)$ contains
all letters of $A_k$. Thus $\varphi_k$ is primitive.

2. Since $\varphi_k(f(u))\le f(\varphi(u))$ for every $u\in\cL_{\ge k}(X)$,
we have also $\varphi_k^n(f(u))\le f(\varphi^n(u))$ for every
$n\ge 1$. 

3. This shows that for every $b\in A_k$ and $n\ge 1$, the word
$\varphi_k^n(b)$ is in $\cL(X^{(k)})$, which implies the conclusion
since $X^{(k)}$ is minimal.
\end{proof}
We denote by $M_k$ the incidence matrix of the morphism $\varphi_k$.
\index{symbols}{Mk@$M_k$}%
Thus, we have $M_k = M(\varphi_k)^t$.
\begin{proposition}\label{propositionMk}
All matrices $M_k$, for $k\ge 1$, have the same dominant
eigenvalue.
\end{proposition} 
\begin{proof} 
Since  $|\varphi_k^n(b)|=|\varphi^n(a)|$ for all $n\ge 1$
and $b\in A_k$ with $a=\pi(b)$ 
by \eqref{eqLgphikn}, this follows from Equation
\eqref{equationLambda}.
\end{proof}
Actually, one can prove that all matrices $M_k$ for $k\ge 2$
have the same nonzero eigenvalues (Exercise~\ref{exerciseShiftEquivM_k}).
 \begin{example}\label{exampleFibonacci3}
Let $\varphi:a\mapsto ab,b\mapsto a$ be the Fibonacci morphism
as in Example~\ref{exampleBlockPresentation}.
Set $A_2=\{x,y,z\}$ and $f:x\mapsto aa,y\mapsto ab,z\mapsto ba$.
Then $\varphi_2:x\mapsto yz,y\mapsto yz, z\mapsto x$
as we have already seen. Thus
\begin{displaymath}
M_2=\begin{bmatrix}0&1&1\\0&1&1\\1&0&0\end{bmatrix}
\end{displaymath}
\end{example}
%%%%%%%%%%%%%%%%%%%%%%
 \section{Sturmian and Arnoux-Rauzy shifts}\label{sectionSturmianShifts}
 In this section, we introduce an important class of minimal shift spaces.
 We begin with the classical Sturmian shifts, which
 are on a binary alphabet, and we present next the
 episturmian and Arnoux Rauzy shifts which are a generalization
 of the former to an arbitrary finite alphabet.
\subsection{Sturmian shifts}
A one-sided sequence $s$ on the binary alphabet $\{0,1\}$ is
\emph{Sturmian}\index{subject}{one-sided!sequence!Sturmian}
\index{subject}{Sturmian!one-sided sequence} if its
word complexity is given by $p_n (s) = n+1$
(we will say often, by a slight abuse of language,
that its complexity is $n+1$). Thus, by Proposition~\ref{theoremCovenHedlund}, Sturmian sequences
have the minimal possible nonconstant word complexity.

Given two words $u,v\in\{0,1\}^*$ their \emph{balance}\index{subject}{balance!of two words}\index{subject}{word!balance of} is
\begin{displaymath}
  \delta(u,v)=\lvert |u|_1-|v|_1\rvert.
\end{displaymath}
A set of words $U\subset\{0,1\}^*$ is \emph{balanced}\index{subject}{balanced!set of words} if for all $u,v\in U$
\begin{displaymath}
  |u|=|v|\Rightarrow \delta(u,v)\le 1.
\end{displaymath}
A sequence $s$ is \emph{balanced}\index{subject}{balanced!sequence}
if $\cL(x)$ is balanced.
The following statement gives an alternative definition of Sturmian sequences.
\begin{proposition}\label{propositionEquivalentDefSturm}
A sequence $s$ is Sturmian if and only if it is balanced and aperiodic.
\end{proposition}
The proof is left as Exercise~\ref{exerciseEquivalentDefSturm}.

A  one-sided shift space $X$ on a binary alphabet
is called \emph{Sturmian}
\index{subject}{Sturmian!shift space}\index{subject}{shift space!Sturmian} if
it is generated by a Sturmian sequence. It can be shown that a Sturmian
shift is minimal (we will prove this below, see Equation~\eqref{equationReturnPal}).

A two-sided shift
is Sturmian if its language is the language of a one-sided Sturmian sequence.
Thus, the complexity of a Sturmian shift is $n+1$ and conversely
a minimal shift of complexity $n+1$ is Sturmian. However, it should
be noted that there exist non minimal shifts of complexity $n+1$
which are not Sturmian. An example is the shift space $X$ such that
$\cL(X)=0^*1^*$ (recall from \eqref{equationStar} and
\eqref{equationProduct} that $0^*1^*=\{0^n1^m\mid n,m\ge 0\}$).

An element of a Sturmian shift is a \emph{Sturmian sequence}.\index{subject}{Sturmian!sequence}
\index{subject}{sequence!Sturmian}%

We shall see 
shortly an example of a Sturmian shift (Example~\ref{exampleFibonacci0bis}).

%The language of  Sturmian subshifts can be characterized by
%several equivalent properties. One of them is
%the \emph{balance property}.

The definition implies that
 Sturmian shifts are such that for every $n\ge 1$ there is
exactly one right special word (and one left special word). 
Since a prefix of a left special word is left special, this
implies that the left special words in a Sturmian shift are the prefixes
of one right infinite word.

\begin{example}\label{exampleFibonacci0bis}
The Fibonacci shift (see Example~\ref{exampleFibonacci0})
is Sturmian (see Exercise~\ref{exerciseFiboisSturmian}). 
\end{example}

As well known, Sturmian shifts correspond to the coding of irrational rotations
on the circle. Actually, for every irrational real number $\alpha$
with $0<\alpha\le 1$ and $0\le\rho< 1$, let $s_{\alpha,\rho}=(s_n)_{n\in \Z}$ be the biinfinite word 
defined by
\begin{equation}
s_n=\lfloor (n+1)\alpha+\rho\rfloor -\lfloor n\alpha+\rho\rfloor. \label{eqRotation}
\end{equation}
It can be shown that  $s^+_{\alpha,\rho}$ is a Sturmian sequence (Exercise~\ref{MechanicalIsSturm}.
Conversely every
Sturmian shift is the closure of the orbit of a sequence
$s_{\alpha,\rho}$ with $\rho=0$ (we actually prove this  statement below
in Proposition~\ref{propositionRotationIsSturm}). 

For example, if $\alpha=(3-\sqrt{5})/2$, then $s_0s_1\cdots=001001\cdots$.
The closure of the orbit of $s$ is the Fibonacci shift.
\index{subject}{Fibonacci!shift}\index{subject}{shift space!Fibonacci}%

\subsection{Episturmian shifts}
Sturmian shifts can be generalized to arbitrary alphabets
as follows. A shift space $X$ on an alphabet $A$ is called
\emph{episturmian}\index{subject}{episturmian!shift}\index{subject}{shift space!episturmian}
 if $\cL(X)$ is closed under reversal and for every $n\ge 1$, there
 is at most one right-special word of length $n$.

 It is
called \emph{strict episturmian}
\index{subject}{strict episturmian!shift}%
\index{subject}{episturmian!shift!strict}%
(or also \emph{Arnoux-Rauzy})
\index{subject}{Arnoux-Rauzy!shift}%
\index{subject}{shift space!Arnoux Rauzy}%
\index{names}{Arnoux, Pierre}\index{names}{Rauzy, G\'erard} if for each $n\ge 1$, there is a unique right-special word $w$ of length $n$ which is moreover such that
$wa\in \cL(X)$ for every $a\in A$. Again, the definition applies
to two-sided or one-sided shifts.

As an equivalent definition, a minimal shift $X$ on the alphabet
$A$ is an Arnoux-Rauzy shift
if for every $n$ there is a unique left-special (resp. right-special)
factor $w$ of length $n$ which is moreover such that $wa\in\cL(X)$
(resp. $aw\in\cL(X)$) for every $a\in A$
(Exercise~\ref{exerciseEquivalentDefAR}).

A Sturmian shift is strict episturmian. Indeed, if $X$ is Sturmian,
it can be shown that $\cL(X)$ is closed under reversal (Exercise~\ref{exerciseSturmReversal}).

A one-sided infinite sequence $s$ is called {\em episturmian}
\index{subject}{episturmian!sequence}\index{subject}{sequence!episturmian} (resp. strict episturmian
or Arnoux-Rauzy)\index{subject}{strict episturmian!sequence}
\index{subject}{episturmian!sequence!strict}%
\index{subject}{sequence!strict episturmian}\index{subject}{Arnoux-Rauzy!sequence}%
if the subshift generated by $s$ is episturmian (resp. strict episturmian). 
It is called \emph{standard}\index{subject}{standard!episturmian sequence}
\index{subject}{episturmian!sequence!standard} if its left special factors are
prefixes of $s$. For every  episturmian one-sided shift $X$,
there is a unique standard infinite word $s$ in $X$.
%Accordingly, when $X$ is a strict episturmian two-sided shift,
%there is a unique standard infinite word $y$ which belongs
%to the one-sided shift associated to $X$, that is, such that
%$y=x^+$ for some $x\in X$.

\begin{example}\label{exampleTribonacci}
The morphism $\varphi:a\mapsto ab,b\mapsto ac,c\mapsto a$
is called the \emph{Tribonacci morphism}.
\index{subject}{Tribonacci!morphism}\index{subject}{morphism!Tribonacci}%
The fixed point $s=\varphi^\omega(a)$ is called the \emph{Tribonacci sequence}.
\index{subject}{Tribonacci!sequence}\index{subject}{sequence!Tribonacci}%
It is a strict and standard episturmian sequence (Exercise~\ref{exerciseTriboisSturmian}).
There are three two sided infinite words $z$ such that
$z^+=s$, namely $\varphi^{3\omega}(a\cdot a)$, $\varphi^{3\omega}(b\cdot a)$
and $\varphi^{3\omega}(c\cdot a)$.
\end{example}

For $a\in A$, denote by $L_a$ the map
\index{subject}{elementary!automorphism}%
\index{subject}{automorphism!elementary}%
\index{symbols}{L@$L_a$}%
 defined for every $b\in A$  by
\begin{displaymath}
L_a(b)=\begin{cases}ab&\mbox{if $b\ne a$}\\a&\mbox{otherwise.}\end{cases}
\end{displaymath}
Each $L_a$ is an automorphism of the free group on $A$ since
$L_a^{-1}(b)=a^{-1}b$ for $b\ne a$ and $L_a^{-1}(a)=a$ and
these maps are called the \emph{elementary automorphisms}.
We define $L_u$ for $u\in A^*$ by extending the map $a\mapsto L_a$
to a morphism $u\mapsto L_u$
from $A^*$ into the group of automorphisms of the free group
by $L_{ua}=L_u\circ L_a$.

A \emph{palindrome}\index{subject}{palindrome} is a word equal to its reversal.
It happens that palindromes play an essential role in Sturmian sequences,
due to the fact that a bispecial factor of a Sturmian word,
being the only one of its length, has to be a palindrome.

For a word $w$, the \emph{palindromic closure}\index{subject}{palindromic closure}
of $w$, denoted $w^{(+)}$,\index{symbols}{w@$w^{(+)}$} is the shortest palindrome which has $w$ as a prefix.
The \emph{iterated palindromic closure}\index{subject}{iterated palindromic closure}
\index{subject}{palindromic closure!iterated} of a word $w$, denoted $\Pal(w)$,
 is defined by $\Pal(\varepsilon)=\varepsilon$  and
$\Pal(ua)=(\Pal(u)a)^{(+)}$ for $u\in A^*$ and $a\in A$.
 \index{symbols}{Pal@$\Pal(u)$}%
 The map $w\mapsto\Pal(w)$ is called the \emph{palindromization map}.
 \index{subject}{palindromization map}%

 We may extend the palindromization map to
 one-sided infinite sequences.
 Indeed, since $\Pal(u)$ is a prefix of $\Pal(uv)$, there is for every $x\in A^\N$,
a unique right
infinite word $s$ such that $\Pal(u)$ is a prefix of $s$ for every
prefix $u$ of $x$.
We set $\Pal(x) = s$. The map $x\mapsto \Pal(x)$ is continuous.

An important property of the palindromization
map is \emph{Justin Formula}. This formula expresses
that, for every $u,v\in A^*$, one has
\index{subject}{Justin Formula}\index{subject}{Formula!Justin}%
\index{names}{Justin, Jacques}%
\begin{equation}
\Pal(uv)=L_u(\Pal(v))\Pal(u)\label{eqJustin0}
\end{equation}
The proof of Formula~\ref{eqJustin0} is given as Exercise~\ref{exerciseJustin4}.
One may verify that
the palindromization map is the unique function $f:A^*\to A^*$
such that $f(uv)=L_u(f(v))f(u)$ for every $u,v\in A^*$.

As a consequence, for every $u\in A^*$ and $x\in A^\N$, one has
\begin{equation}
\Pal(ux)=L_u(\Pal(x)).\label{eqJustin}
\end{equation}
Moreover, one may verify that
the palindromization map is the unique continuous map $f:A^\N\to A^\N$
such that $f(ux)=L_u(f(x))$ for every $u\in A^*$ and $x\in A^\N$.

We will use the following result in Chapter \ref{chapterDimensionGroupsPartitions}.

\begin{theorem}\label{standardEpisturmianTheorem}
  A one-sided infinite sequence $s$ on the alphabet $A$
  is  standard episturmian  if and only
if there exists a one-sided sequence $x$ on $A$ such that
 $s=\Pal(x)$.
Moreover, the  episturmian sequence $s$ is strict if and only if
every letter of $A$ occurs infinitely often in $x$.
\end{theorem}
The proof is given below. 

The one-sided infinite word $x$ is called the \emph{directive sequence}
\index{subject}{directive!sequence}%
of the standard sequence $s$.

Note that since a Sturmian sequence
is strict episturmian, the directive sequence of a Sturmian sequence
has an infinite number of occurrences of $0$ and $1$.

The sequence $s$ is actually in the limit set $\cap_{n\ge 0} L_{a_0\cdots a_{n-1}}(A^\N)$ obtained using the infinite sequence of morphisms
\begin{equation}
  \ldots A^*\edge{L_{a_1}}A^*\edge{L_{a_0}}A^*.\label{eqEpisturm}
\end{equation}
In such a sequence, the infinite iteration of one morphism is replaced by
the composition of an infinite sequence of morphisms, a feature
that we will introduce later with $\Sa$-adic representations of shift spaces
(see Section~\ref{sectionSadicShifts}).

If $x=x_0x_1\cdots$, the words $u_n=\Pal(x_0\cdots x_{n-1})$ are
the palindromic prefixes of $s$. It can be shown 
that  the set of left return words to $u_n$ satisfies
\begin{equation}
\RR'_X(u_n)\subset\{L_{x_0\cdots x_{n-1}}(a)\mid a\in A\}\label{equationReturnPal}
\end{equation}
with equality when $s$ is strict and that the sets of return words
to other words have the same form up to conjugacy (Exercise~\ref{exerciseReturnSturm}).
This formula implies in particular that the
set of return words is finite and thus that
an episturman shift is minimal.

Formula \eqref{equationReturnPal} also shows the remarkable fact that in a strict episturmian
shift $X$, the set of return words to the words $u_n$ 
(and, as a consequence to every word in $\cL(X)$) has a
constant cardinality. We shall have more to say about this
later, when we introduce dendric shifts (Chapter~\ref{chapterDendricShifts}).
\begin{example}\label{exampleTribonacci2}
The directive sequence of the Tribonacci sequence $s$ (see Example~\ref{exampleTribonacci})
is $(abc)^\omega$. Indeed,
we have by Justin's Formula $\Pal((abc)^\omega)=L_{abc}(\Pal(abc)^\omega)$
whence the result since $L_{abc}=\varphi^3$. The palindromic
prefixes of $s$ are 
\begin{displaymath}
a,aba,abacaba,\ldots.
\end{displaymath}
 We have for example
\begin{eqnarray*}
\RR'_X(abacaba)&=&\{L_{abc}(a),L_{abc}(b),L_{abc}(c)\}\\
&=&\{abacaba,abacab,abac\}.
\end{eqnarray*}
\end{example}

We now give an example of an episturmian shift which is not strict.
\begin{example}\label{exampleEpistutmianNotStrict}
We have $\Pal(ab^\omega)=(ab)^\omega$, which is not Sturmian.
  \end{example}

The proof of Theorem~\ref{standardEpisturmianTheorem}
uses the following notion. A one-sided infinite sequence $s$ is \emph{palindrome
  closed} \index{subject}{palindrome!closed sequence}
\index{subject}{sequence!palindrome closed} if for every prefix $u$
of $s$, the word $u^{(+)}$ is also a prefix of $s$.
The following property is not surprising.
\begin{proposition}\label{proposition(iii)->(i)}
  If $s=\Pal(x)$ for some infinite sequence $x$, then $s$ is palindrome closed.
\end{proposition}
\begin{proof}
  Set $x=x_0x_1\cdots$ and $u_n=\Pal(x_0\cdots x_{n-1})$. Let $u$ be a prefix of $s$. Then
  $|u_n|<u\le|u_{n+1}|$ for some $n$. Since
  $u^{(+)}$ is a palindrome with $u_n$ as a prefix, we have $u^{(+)}=u_{n+1}$.
  Thus $s$ is palindrome closed.
\end{proof}

%\begin{proposition}\label{propositionAl}
%  Let  $s$ be a palindrome closed  sequence. Let
%  $v$ be a palindrome factor of $s$ and let $s=uvs'$ with $|u|$ minimal.
%  Then $v$ is the longest palindrome suffix of $uv$ and
%  $uv\tilde{u}$ is a prefix of $s$.
%\end{proposition}
%\begin{proof}
%  Assume that $v$ is not the longest palindrome suffix of $uv$. Then $uv=u'wv$ with $wv$ palindrome and $w$ nonempty.
%  Thus
%  $uv=u'v\tilde{w}$, a contradiction since $|u'|<|u|$. Then
%  $(uv)^{(+)}=uv\tilde{u}$ (see Exercise \ref{exerciseJustin0}).
 % Since $s$ is palindrome closed, $(uv)^{(+)}$ is a prefix of $s$.
%\end{proof}
A factor  of a word $w$ is \emph{unioccurrent}\index{subject}{unioccurrent factor} if it has exactly one occurrence in $w$.
A word $w$ is a \emph{Justin word} if it has
a palindromic suffix which is unioccurrent. This suffix is then
also the longest palindromic suffix of $w$. As an example,
every palindrome is a Justin word but $w=abca$ is not a Justin word.

A one-sided sequence $s$ is \emph{weakly palindrome closed}
if for every prefix $u$ of $s$ which is a Justin word, the word 
$u^{(+)}$ is a prefix of $s$.
\begin{proposition}\label{lemma1DJP}
  Let $s$ be a weakly palindrome closed sequence. Then
  every prefix of $s$ is a Justin word. Consequently $s$ is palindrome closed.
\end{proposition}
\begin{proof}
  Let $w$ be a prefix of $s$. We argue by induction on $n=|w|$. The
  property is true if $n=1$ since every letter is a Justin word.
  Next, consider $w=va$ with $a\in A$.
  By induction hypothesis, $v$ is a Justin word. Set $v=v_1v_2$ with $v_2$ palindrome unioccurrent in $v$. Then $v_2$ is also the longest
  palindromic suffix of $v$ and $v^{(+)}=v_1v_2\tilde{v_1}$
  (see Exercise~\ref{exerciseJustin0}). If $v_1$ is not empty,
  since $s$ is weakly palindrome closed, $v^{(+)}=v_1v_2\tilde{v}_1$ is a prefix of $s$
  and thus $v_1$ ends with $a$.
  Then $av_2a$ is an unioccurrent palindrome suffix of $w$
  and the property is true. Otherwise, $v$ is a palindrome. If $a$ does not
  occur in $v$, then $a$ is a palindromic suffix of $w$ unioccurrent in $w$.
  Otherwise, let $ua$ be the prefix of $v$ ending with $a$ with $|u|$ minimal.
  
  Claim 1: the word $u$ is a palindrome. By induction
  hypothesis, we have $u=u_1u_2$ with $u_2$ a palindrome unioccurent in $u$.
  Again if $u_1$ is not empty, since $u^{(+)}=u_1u_2\tilde{u_1}$
  is a prefix of $s$, the word $u_1$ ends with $a$, a contradiction with
  the definition of $u$. This proves Claim 1.

  Let $qa$ be the longest prefix of $v$ with $q$ a palindrome. Since $v$
  is a palindrome, the word
  $aqa$ is a palindromic suffix of $w$.

  Claim 2:
  the palindrome $aqa$ is  unioccurrent in $w$. Otherwise, set $v=xaqay$
  with $aqa$ unioccurrent in $xaqa$.

  Let us first prove that $z=xaq$
  is a palindrome. Otherwise, we have by induction hypothesis
  $z=z_1z_2$ with $z_2$ a palindrome unioccurrent in $z$ and $z_1$ nonempty. Since
  $s$ is weakly palindrome closed, $z^{(+)}=z_1z_2\tilde{z_1}$ is a prefix
  of $s$ and thus $a$ is the last letter of $z_1$. But then $az_2a$
  and $aqa$ are two palindromic unioccurrent suffixes of $za$,
  which forces $z_2=q$. Since $q$ is a prefix
  of $z$ this contradicts the unioccurrence of $z_2$.
  Thus $z$ is a palindrome.

  But $za$ is a prefix of $v$ and $|z|>|q|$,
  which contradicts the definition of $q$. This proves Claim 2 and thus the
  proposition.
  
\end{proof}
The following result is the key of the proof of Theorem~\ref{standardEpisturmianTheorem}.
\begin{proposition}\label{proposition3cond}
  The following conditions are equivalent for a one-sided sequence $s$.
  \begin{enumerate}
  \item[\rm(i)] $s$ is palindrome closed.
  \item[\rm(ii)] $s$ is weakly palindrome closed.
  \item[\rm(iii)] $s=\Pal(x)$ for some one-sided sequence $x$.
  \end{enumerate}
\end{proposition}
\begin{proof}
  (i)$\Leftrightarrow$ (ii) is clear by Proposition~\ref{lemma1DJP}.

  (ii)$\Rightarrow$ (iii) Let $x_1$ be the first letter of $s$
  and set $u_1=x_1$. Assume that $x_0,\ldots,x_{n-1}$ and
  $u_i=\Pal(x_0\cdots x_{i-1})$ for $1\le i\le n$ are already defined.
  Let $x_n$ be the letter following $u_n$ in $s$. By Proposition
  \ref{lemma1DJP}, $u=u_nx_n$ is a Justin word. Thus, since $s$
  is palindrome closed, $u^{(+)}$ is a prefix of $s$. This
  allows us to define $u_{n+1}=u^{(+)}=\Pal(x_0\cdots x_{n})$.
  The word $x=x_0x_1\cdots$ build in this way is clearly
  such that $s=\Pal(x)$.

  (iii)$\Rightarrow$ (i) is Proposition~\ref{proposition(iii)->(i)}.
  \end{proof}
\begin{proposition}\label{proposition5DJP}
  If $s=\Pal(x)$ for some infinite word $x$, then all its
  left-special factors are prefixes of $s$.
\end{proposition}
\begin{proof}
  Let $u$ be the shortest left-special factor of $s$ which is not
  a prefix of $s$. Set $u=va$ with $a\in A$. Then $v$ is left-special
  and shorter than $u$. Thus $v$ is a prefix of $s$. Let $c\in A$
  be such that $vc$ is a prefix of $s$. Since $c\ne a$, $v$ is right-special.
  Thus $\tilde{v}$ is a prefix of $s$ and consequently $v=\tilde{v}$.
  Since $u$ is left-special, we have $bva,b'va\in\cL(s)$ for
  two distinct letters $b,b'$. 

  Set $x=x_0x_1\cdots$ and $u_n=\Pal(x_0\cdots x_{n-1})$. Let
  $n$ be such that $u$ is a factor of $u_{n+1}$ but not of $u_n$.
  Then $bva$ or $avb$ is a prefix of $x_nu_n$ (Exercise~\ref{exerciselemma6}).
  Since $|u_n|>|v|$, $vc$ is a prefix of $u_n$ whence $c=a$ or $c=b$.
  Since $c\ne a$, we have $c=b$. In the same way, $c=b'$,
  whence $b=b'$, a contradiction.
\end{proof}

\begin{proofof}{of Theorem~\ref{standardEpisturmianTheorem}}
  Assume first that $s=\Pal(x)$. By Proposition~\ref{proposition5DJP}
  its left-special factors are prefixes of $s$ and thus there
  is at most one of each length. Moreover, every factor $u$ of $s$
  is a factor of some palindromic prefix and thus $\tilde{u}$
  is also a factor of $s$. Thus $s$ is standard episturmian.

  Conversely, consider a standard episturmian sequence $s$.
  Let us prove that $s$ is weakly palindrome closed. Otherwise,
  there is a leftmost occurrence $s=uavw\tilde{v}bs'$ of a palindrome $w$
  with $a\ne b$. Then $vw\tilde{v}$ is left-special and thus a prefix
  of $s$. This contradicts the hypothesis that  the
  leftmost occurrence of $w$ appears after $uav$ and proves the claim.
  By Proposition~\ref{proposition3cond}, this implies that
  $s=\Pal(x)$ for some sequence $x$.

  Let us prove now the last assertion. For this, it is enough to prove that
  for $a\in A$ and $n\ge 1$, one has $u_na\in\cL(s)$
  if and only if $a$ appears in $x_nx_{n+1}\cdots$.

  If $a=x_m$
  for some $m\ge n$, then $u_ma$ is a prefix of $s$ and thus
  $u_na\in\cL(s)$ since $u_n$ is a suffix of $u_m$. Thus $u_na\in\cL(s)$.

  Conversely, if $u_na\in\cL(s)$, we may assume that $a\ne x_n$
  and thus that $u_na$ is not a prefix of $s$.
  Let $b\in A$ be such
  that $bu_na\in\cL(s)$ and let $m$ be such that
  $bu_na$ is a factor of $u_{m+1}$ but not of $u_m$.
  Then, by Exercise~\ref{exerciselemma6},
  either $bu_na$ or  is a suffix of $u_mx_m$ or it
  is a prefix of $x_mu_m$. The first case implies
  $a=x_m$ and the conclusion.
  The second case is not possible. Indeed, since
  $|bu_na|\le |u_mx_m|$, we have $n<m$ and thus
  $u_na$ being a prefix of $u_m$, this forces $a=x_n$.
  \end{proofof}
\subsection{Sturmian sequences, rotations and continued fractions}
Let us come back to binary Sturmian sequences. We will use the traditional
alphabet $A=\{0,1\}$.
To every standard Sturmian sequence $s\in \{0,1\}^\N$, we associate a real number
$\alpha$, with $0<\alpha\le 1$, called its slope
\index{subject}{slope of Sturmian sequence}%
and defined as follows. 

Let $x$ be the directive sequence of $s$. Since $x$ has an infinite
number of occurrences of $0$ and $1$, we can write
\begin{displaymath}x=0^{d_1}1^{d_2}0^{d_3}\cdots
\end{displaymath}
 Set
\begin{displaymath}
\pi(x)=[0;1+d_1,d_2,d_3,\ldots].
\end{displaymath}
where $[a_0;a_1,a_2,\ldots]$ denotes the continued fraction
with coefficients $a_0,a_1,\ldots$ (see Appendix~\ref{appendixAlgebraicNumberTheory}). The \emph{slope} of $s$ is the real number $\pi(x)$.

Recall that we denote $\T=\R/\Z$ the one-dimensional torus.
For $0<\alpha\le 1$, let $R_\alpha:\T\to\T$ be the transformation defined by
$R_\alpha(z)=z+\alpha$.
\index{symbols}{R@$R_\alpha$}%
 The pair $(\T,R_\alpha)$ is a topological dynamical
system called the \emph{rotation}\index{subject}{rotation} of \emph{angle}
\index{subject}{angle of rotation}\index{subject}{rotation!angle}%
$\alpha$.
\begin{proposition}\label{propositionRotationIsSturm}
Let $s$ be a standard Sturmian sequence, let $X$ be the subshift generate by $s$ and let $\alpha$ be the slope of $s$.
The map $\gamma_\alpha:\T\to X$,
\index{symbols}{gamma@$\gamma_\alpha$} defined, for $z\in\T$, by $\gamma_\alpha(z) = y$ where
\begin{displaymath}
y_n=\begin{cases}0&\mbox{if \ $R_\alpha^nz\in[0,1-\alpha)$}\\1&\mbox{otherwise,}
\end{cases}\end{displaymath}
\begin{itemize}
\item
is such that $s=\gamma_\alpha(\alpha)$,
\item
is an injective map and
\item
satisfies $\gamma_\alpha\circ R_\alpha=S\circ\gamma_\alpha$ (see Figure
\ref{figureNaturalRepresentation}).
\end{itemize}
\end{proposition}
\begin{figure}[hbt]
\centering
\tikzset{node/.style={circle,draw,minimum size=0.1cm,inner sep=0pt}}
\begin{tikzpicture}
\node(Tg)at(0,1.5){$\T$};\node(Td)at(2,1.5){$\T$};
\node(Xg)at(0,0){$X$};\node(Xd)at(2,0){$X$};

\draw[above,->](Tg)edge node{$R_\alpha$}(Td);
\draw[left,->](Tg)edge node{$\gamma_\alpha$}(Xg);
\draw[right,->](Td)edge node{$\gamma_\alpha$}(Xd);
\draw[above,->](Xg)edge node{$S$}(Xd);
\end{tikzpicture}
\caption{The map $\gamma_\alpha$.}\label{figureNaturalRepresentation}
\end{figure}
The map $\gamma_\alpha$ is called the \emph{natural coding}
\index{subject}{natural!coding!of rotation}%
\index{subject}{rotation!natural coding of}\index{subject}{coding!natural}%
 of $(\T,R_\alpha)$.
 Note first that $\gamma_\alpha(z)=s_{\alpha,\rho}$ for
 $z=\lfloor \alpha+\rho\rfloor$ where $s_{\alpha,\rho}$
 is the sequence defined by Equation~\ref{eqRotation}.
 This shows that every Sturmian shift if generated by
 such a sequence. It gives also the following corollary.
 The \emph{slope}\index{subject}{slope!of word} of a word $u\in\{0,1\}^*$
 is the real number $\pi(u)=|u|_1/|u|$.
 \begin{corollary}\label{corollarySlope}
For every factor
$u$ of a standard Sturmian sequence $s$ of slope $\alpha$, the slope of $u$
tends to $\alpha$ when $|u|\to\infty$.
 \end{corollary}
 Indeed, by Proposition~\ref{propositionRotationIsSturm}, we have $s=s_{\alpha,0}$ and thus the statement
 follows by Exercise~\ref{exerciseMechanicalIsSturm}.
 
Observe also that $\gamma_\alpha$ is not a conjugacy
since it is neither continuous nor surjective. 
Indeed, it is not continuous since
$\gamma_\alpha(0)_0=0$ while  $\gamma_\alpha(z)_0=1$ for
all $z\in(\alpha,1)$ and thus for values of $z$ arbitrary close
to $0$. Thus $\gamma_\alpha$ is not continuous at zero (it is however
continuous at all other points of $\T$).

Next,
the sequence $t=\lim_{z\rightarrow 1_-}\gamma_\alpha(z)$  is not in $\gamma_\alpha(\T)$. 
Indeed, let 
$$c_\alpha=\gamma_\alpha(\alpha)^+ .$$
By definition of $c_\alpha$ and $t$ we have $t^+=1 c_\alpha$.
Moreover, if there exists some $y$ with 
$t=\gamma_\alpha(y)$ then we should have $y=0$, hence $\gamma_\alpha(0)_0=0$,
a contradiction.

The sequence $c_\alpha$\index{symbols}{c@$c_\alpha$} is called the \emph{characteristic sequence}
of slope $\alpha$. \index{subject}{characteristic!sequence}
Note that
the left-special words in $\cL(X)$ are the prefixes
of $c_\alpha$. Thus $c_\alpha$ is a standard sequence.

Moreover as in a Sturmian shift $X$ there is a unique left special word of length $n$ for all $n\in \N$, this implies that $c_\alpha$ is unique: if $x^+$ is standard for some $x\in X$ then $x=c_\alpha$.
We thus say that $X$ is the Sturmian shift of slope $\alpha$.

To prove Proposition~\ref{propositionRotationIsSturm}, we first prove the following lemma.
\begin{lemma}
For every irrational $\alpha$ with
$0<\alpha\le 1$, the sequence $c_\alpha=\gamma_\alpha(\alpha)^+$
satisfies
\begin{equation}
c_\alpha=\begin{cases}L_0(c_{\alpha/(1-\alpha)})&\mbox{ if $\alpha<1/2$}\\
L_1(c_{(1-\alpha)/\alpha})&\mbox{ otherwise.}
\end{cases}\label{eqc_alpha}
\end{equation}
\end{lemma}
\begin{proof}
  Let $\alpha=[0;1+d_1,d_2,\ldots]$ be the continued fraction expansion of
  $\alpha$.
Assume that $d_1>0$ (or equivalently $\alpha<1/2$). We have to prove that
\begin{equation}
c_\alpha=L_0(c_{\alpha/(1-\alpha)}).\label{eqc_alpha}
\end{equation}
We consider the transformation $R_\alpha$ as a map defined on $[0,1)$
translating the semi intervals $[0,1-\alpha)$ and $[1-\alpha,1)$
as indicated in Figure~\ref{figureRotation} on the left.
Consider the transformation $R'$ induced by $R=R_\alpha$ on the semi-interval
$[0,1-\alpha)$. It is defined by
\begin{displaymath}
R'(z)=\begin{cases}R(z)&\mbox{ if $z\in [0,1-2\alpha)$}\\
R^2(z)&\mbox{ otherwise}\end{cases}
\end{displaymath}
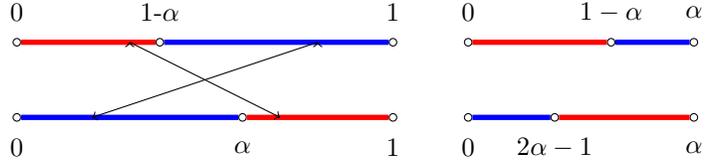
\begin{figure}[hbt]
\centering
\tikzset{node/.style={circle,draw,minimum size=0.1cm,inner sep=0pt}}
\tikzset{title/.style={minimum size=0.5cm,inner sep=0pt}}
\begin{tikzpicture}
%original
\node[node](0h)at(0,1){};\node[title]at(0,1.4){$0$};
\node[node](alphah)at(3,1){};\node[title]at(3,1.4){$1-\alpha$};
\node[node](1h)at(5,1){};\node[title]at(5,1.4){$1$};
\node[node](0b)at(0,0){};\node[title]at(0,-0.4){$0$};
\node[node](alphab)at(1.9,0){};\node[title]at(1.9,-.4){$\alpha$};
\node[node](1b)at(5,0){};\node[title]at(5,-.4){$1$};

\draw[line width=2pt,red](0h)edge node{}(alphah);
\draw[line width=2pt,blue](alphah)edge node{}(1h);
\draw[line width=2pt,blue](0b)edge node{}(alphab);
\draw[line width=2pt,red](alphab)edge node{}(1b);
\draw[->](1.5,1)edge node{}(3.5,0);\draw[->](4,1)edge node{}(1,0);
%transformed
\node[node](0h)at(6,1){};\node[title]at(6,1.4){$0$};
\node[node](alphah)at(7.15,1){};\node[title]at(7.15,1.4){$1-2\alpha$};
\node[node](1h)at(9,1){};\node[title]at(9,1.4){$1-\alpha$};
\node[node](0b)at(6,0){};\node[title]at(6,-0.4){$0$};
\node[node](alphab)at(7.9,0){};\node[title]at(7.9,-.4){$\alpha$};
\node[node](1b)at(9,0){};\node[title]at(9,-.4){$1-\alpha$};

\draw[line width=2pt,red](0h)edge node{}(alphah);
\draw[line width=2pt,blue](alphah)edge node{}(1h);
\draw[line width=2pt,blue](0b)edge node{}(alphab);
\draw[line width=2pt,red](alphab)edge node{}(1b);

\end{tikzpicture}
\caption{The action of $R_\alpha$ when $\alpha<1/2$.}\label{figureRotation}
\end{figure}
The transformation $R'$ is pictured on the right side of Figure  \ref{figureRotation}.

Let $\rho (z)\in \{ 0,1\}^\Z$ be the natural coding of the $R'$-orbit of $z\in [0,1-\alpha)$ with the intervals $[0,1-2\alpha )$, coded by $0$, and $[1-2\alpha , 1-\alpha)$ coded by $1$ as we did to define $\gamma_\alpha$.
Then $\gamma_\alpha (z) = L_0 (\rho (z))$.
Moreover, normalizing the map $R'$ to $[0,1)$ it is the rotation of angle $\alpha / (1-\alpha )$ and $\rho (z)$ is $\gamma_{\alpha /(1-\alpha )} (z/(1-\alpha ))$.
Hence we have the relation
$$\gamma_\alpha(z)=L_0(\gamma_{\alpha /(1-\alpha )} (z/(1-\alpha )))$$
which proves  Equation~\eqref{eqc_alpha} taking $z=\alpha$.

If $\alpha>1/2$, we consider instead the transformation $R'$
  induced by $R_\alpha$ on $[0,\alpha)$ (see Figure~\ref{figureRotation2}).
Let $\rho (z)\in \{ 0,1\}^\Z$ be the natural coding of the $R'$-orbit of $z\in [0,\alpha)$ with the intervals $[0,1-\alpha )$, coded by $0$, and $[1-\alpha , \alpha)$ coded by $1$.
Then $\gamma_\alpha (z) = L_1 (\rho (z))$.
Moreover, normalizing the map $R'$ to $[0,1)$ it is the rotation of angle $\alpha / (1-\alpha )$ and $\rho (z)$ is $\gamma_{\alpha /(1-\alpha )} (z/(1-\alpha ))$.
Hence we have the relation
$$\gamma_\alpha(z)=L_0(\gamma_{\alpha /(1-\alpha )} (z/(1-\alpha )))$$
which proves  Equation~\eqref{eqc_alpha} taking $z=\alpha$.
  
    \begin{figure}[hbt]
\centering
\tikzset{node/.style={circle,draw,minimum size=0.1cm,inner sep=0pt}}
\tikzset{title/.style={minimum size=0.5cm,inner sep=0pt}}
\begin{tikzpicture}
%original
\node[node](0h)at(0,0){};\node[title]at(0,-.4){$0$};
\node[node](alphah)at(3,0){};\node[title]at(3,-.4){$\alpha$};
\node[node](1h)at(5,0){};\node[title]at(5,-.4){$1$};
\node[node](0b)at(0,1){};\node[title]at(0,1.4){$0$};
\node[node](alphab)at(1.9,1){};\node[title]at(1.9,1.4){1-$\alpha$};
\node[node](1b)at(5,1){};\node[title]at(5,1.4){$1$};

\draw[line width=2pt,blue](0h)edge node{}(alphah);
\draw[line width=2pt,red](alphah)edge node{}(1h);
\draw[line width=2pt,red](0b)edge node{}(alphab);
\draw[line width=2pt,blue](alphab)edge node{}(1b);
\draw[->](1.5,1)edge node{}(3.5,0);\draw[->](4,1)edge node{}(1,0);
%transformed
\node[node](0h)at(6,0){};\node[title]at(6,-.4){$0$};
\node[node](alphah)at(7.15,0){};\node[title]at(7.15,-.4){$2\alpha -1$};
\node[node](1h)at(9,0){};\node[title]at(9,-.4){$\alpha$};
\node[node](0b)at(6,1){};\node[title]at(6,1.4){$0$};
\node[node](alphab)at(7.9,1){};\node[title]at(7.9,1.4){$1-\alpha$};
\node[node](1b)at(9,1){};\node[title]at(9,1.4){$\alpha$};

\draw[line width=2pt,blue](0h)edge node{}(alphah);
\draw[line width=2pt,red](alphah)edge node{}(1h);
\draw[line width=2pt,red](0b)edge node{}(alphab);
\draw[line width=2pt,blue](alphab)edge node{}(1b);

\end{tikzpicture}
\caption{The action of $R_\alpha$ when $\alpha>1/2$.}\label{figureRotation2}
    \end{figure}

    We obtain this time a new rotation of angle $(1-\alpha)/\alpha$
    and for every $z\in[0,\alpha)$, we have $\gamma_\alpha(z)=L_1(\gamma_\alpha (R'(z)))$.
      Since $R'(z)=R_{(1-\alpha)/\alpha}(z/\alpha)$, the result follows
      also in this case.
\end{proof}

\begin{proofof}{of Proposition \ref{propositionRotationIsSturm}}
 We have to
prove that for all $u\in\{0,1\}^*$ and $y\in\{0,1\}^\N$, we have
\begin{equation}
c_{\pi(uy)}=L_u(c_{\pi(y)}).\label{eqJustinBis}
\end{equation}
Indeed, by Justin Formula~\eqref{eqJustin}, since the
map $x\mapsto c_{\pi(x)}$ is continuous, this implies that 
$c_{\pi(x)}=\Pal(x)$ and thus that $c_\alpha=s$.

Equation~\eqref{eqJustinBis} is true for $u=\varepsilon$. 
Note that if $\alpha=\pi(x)$ and $x=0x'$, then $\pi(x')=\alpha/(\alpha-1)$.
 Arguing
by induction on $|u|$, consider $u=0v$. Then, by~\eqref{eqc_alpha} and
using the induction hypothesis,
we have
\begin{eqnarray*}
c_{\pi(uy)}&=&L_0(c_{\pi(vy)})=L_0(L_v(c_{\pi(y)}))\\
&=&L_u(c_{\pi(y)}).
\end{eqnarray*}
The case $u=1v$ is similar.
\end{proofof}
The transformation used in the proof is called a \emph{Rauzy induction}
\index{subject}{Rauzy!induction}%
(we shall meet again this notion in Chapter~\ref{ch5:sec:examples}
when we consider interval exchange transformations).
\begin{example}\label{exampleSlopeFibonacci}
The slope of the Fibonacci word is $\alpha=\frac{3-\sqrt{5}}{2}$.
Indeed, its directive word is $0101\ldots$ and thus $\alpha=[0,2,1,1,\ldots]$.
The Fibonacci word is actually an approximation of the line
of equation $y=\alpha x+\alpha$ (see Figure~\ref{figureSlope}).
\begin{figure}[hbt]
\centering
\tikzset{node/.style={circle,draw,minimum size=0.1cm,inner sep=0pt}}
\tikzset{title/.style={minimum size=0.5cm,inner sep=0pt}}
\begin{tikzpicture}
%horiz
\node[title]at(-.4,-.4){$0$};
\draw[->](0,0)edge node{}(6,0);\node[title]at(6,-.4){$x$};
\draw(0,1)edge node{}(5,1);
\draw(0,2)edge node{}(5,2);
%\draw(0,3)edge node{}(5,3);
%vert
\draw[->](0,0)edge node{}(0,3);\node[title]at(-.4,3){$y$};
\draw(1,0)edge node{}(1,2.5);
\draw(2,0)edge node{}(2,2.5);\draw(3,0)edge node{}(3,2.5);
\draw(4,0)edge node{}(4,2.5);\draw(5,0)edge node{}(5,2.5);
%line
\draw(0,.381)edge node{}(5,2.286);\node[title]at(-.4,.381){$\alpha$};
%mecanical
\draw[color=red,above,line width=2pt](0,0)edge node{$0$}(1,0);
\draw[color=red,above,line width=2pt](1,0)edge node{$1$}(2,1);
\draw[color=red,above,line width=2pt](2,1)edge node{$0$}(3,1);
\draw[color=red,above,line width=2pt](3,1)edge node{$0$}(4,1);
\draw[color=red,above,line width=2pt](4,1)edge node{$1$}(5,2);
\end{tikzpicture}
\caption{The Fibonacci word as an approximation of the line $y=\alpha(x+1)$.}\label{figureSlope}
\end{figure}
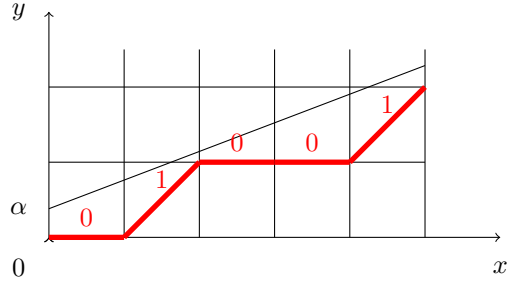
\end{example}

\begin{theorem}
Let $X$  be a Sturmian shift and let $\alpha$ be its slope. 
There exists a unique factor map $f:X\to \T$ such that its restriction to $\gamma_\alpha(\T)$
is the inverse of $\gamma_\alpha$. 
Moreover, one has
\begin{itemize}
\item
$f$ is one-to-one on $X\setminus \gamma_\alpha^{-1}( \mathcal{O}_{R_\alpha} (0)) $,
\item
$\Card (f^{-1} (\{ n\alpha \}) = 2$  for all $n\in \mathbb{Z}$.
\end{itemize}
\end{theorem}

\begin{proof}
\end{proof}
%%%%%%%%%%%%%%%%
\section{Toeplitz shifts}\label{sectionToeplitz}

A sequence $x = (x_n)_{n \in \mathbb{K}}$, with $\mathbb{K} = \mathbb{Z}$ or $\mathbb{N}$, on the alphabet $A$ satisfying 

\[
\forall n \in \mathbb{K}, \exists p\in \mathbb{N}, \forall k\in \mathbb{N}, x_n = x_{n+kp} 
\]
is called a {\em Toeplitz sequence}.\index{subject}{Toeplitz!sequence}%
\index{subject}{sequence!Toeplitz}\index{names}{Toeplitz, Otto}

For such a sequence $x$ and $p \geq 1$, we let 
\[
{\rm Per}_p(x)= \{ n \in \mathbb{Z}; x_{n}=x_{n+k p}\text{ for all } k\in \mathbb{K} \} .
\] 

It is clear that $x$ is a Toeplitz sequence if there exists a sequence $(p_n)_{n\geq 1}$ in $\mathbb{N}\setminus\{0\}$ such that 
$\mathbb{Z}=\bigcup_{n\geq 1} {\rm Per}_{p_n}(x)$. 
Equivalently, $x$ is a Toeplitz sequence if all finite blocks in $x$ appear periodically. 
Hence, Toeplitz sequences are uniformly recurrent.
We say that $p_n$ is an  \emph{essential period} if for any $1 \leq p <p_n$ the sets ${\rm Per}_{p}(x)$ and ${\rm Per}_{p_n}(x)$ do not coincide. If the sequence $(p_n)_{n\geq 1}$ is formed by essential periods and $p_n$ divides $p_{n+1}$, we call it a \emph{periodic structure} of $x$. 
Clearly, if $(p_n)_{n\geq 1}$ is a periodic structure, then 
$(p_{i_n})_{n\geq 1}$ is also a periodic structure for any strictly increasing sequence of positive integers $(i_n)_{n\geq 1}$. 

A shift space $(X,S)$ is a \emph{Toeplitz shift}
\index{subject}{Toeplitz!shift}\index{subject}{shift space!Toeplitz} if $X$ is generated by a Toeplitz sequence $x \in X$.
Note that, unless it is periodic,  a Toeplitz shift contains elements
which are not Toeplitz sequences.\marginpar{Argument?}

\begin{example}\label{exampleToeplitz}
Let $\sigma$ be the substitution $0\mapsto 01,1\mapsto 00$.
 The substitution shift
$X(\sigma)$ is a Toeplitz shift. Indeed, let us show that the admissible
fixed point  $x=\sigma^{2\omega}(1\cdot 0)$
is a Toeplitz sequence
(sometimes known as the \emph{period-doubling sequence}).
\index{subject}{period-doubling sequence}\index{subject}{sequence!period-doubling}% 
First all symbols
$x_{2n}$ of even index are equal to $0$. This implies
that for $k\in\Z$ the blocks $x_{[2k2^N,2k2^N+2^N)}$ of length
$2^N$ are all equal. Thus,
for any $n\ge 1$, let $N$ be such that $n<2^N$. Then $x_n=x_{n+kp}$
for all $k\in\Z$ with $p=2^{N+1}$. The period structure
is $(2^n)_{n\ge 1}$.
\end{example}

More generally, a substitution $\sigma$ of constant length $n$ is said 
to have a \emph{coincidence}\index{subject}{coincidence} at index $k$ for $1\le k\le n$ if the
$k$-th letter of every $\sigma(a)$ is the same. 
\begin{proposition}\label{propositionCoicidence}
A fixed
point of a constant length substitution having a coincidence
at index $1$
is a Toeplitz sequence.
\end{proposition} 
The
proof is similar to the one above. The general case of a coincidence at index $k$ is proposed
as Exercise~\ref{exerciseShiftToeplitz}.

There is a constructive way to obtain all Toeplitz sequences.
Let $A$ be a finite alphabet and $?$ a letter not in $A$ (usually the symbol $?$ is referred as a ``hole'').
Let $x$ $\in (A\cup \{ \text{?} \})^{\mathbb{Z}}$. 
Given $x, y \in (A\cup \{\text{?}\})^{\mathbb{Z}}$, define 
$F_x(y)$ as the sequence obtained from $x$ replacing consecutively all the $?$ by the symbols of $y$, where $y_0$ is placed in the first $?$ to the right of coordinate $0$. In particular, if 
$x$ has no holes, $F_{x}(y)=x$ for every $y \in (A\cup \{\text{?}\})^{\mathbb{Z}}$. 
In addition, observe that:
\begin{equation}\label{eq:magicformula}
\text{ if } \ z=F_{x}(y) \text{ then } \ F_{z}=F_{x}\circ F_{y}.
\end{equation}

Now, consider a sequence of finite words $(w(n))_{n\in \mathbb{N}}$ in $A\cup \{ \text{?}\}$.  
For each $n$, let $y(n)$ be the periodic sequence $w(n)^\infty =\cdots w(n)w(n)w(n).w(n)w(n)w(n)\cdots \in (A\cup \{ \text{?} \})^{\mathbb{Z}}$, where the central dot indicates the position to the left of coordinate $0$. 

We define the sequence $(z (n))_{n\geq 1}$ by: $z(1)=y(1)$ and, for every $n\geq 1$,
\begin{equation}
\label{eq:holetoeplitz}
z(n+1)=F_{z(n)}(y(n+1)) .
\end{equation}

It is not complicated to see that 
$z(n) = u_{n}^{\infty}$  for some word $u_{n}$ of length $|w(1)||w(2)|\cdots |w(n)|$ and that the limit $z=\lim_{n\to \infty} z(n)$ is well defined as a sequence in $(A\cup \{\text{?}\})^{\mathbb{Z}}$. 
Moreover, if the $w(n)$ does not start or finish with a hole for infinitely many $n$, then the limit sequence belongs to $A^{\mathbb{Z}}$, {\em i.e.}, $z$ has no holes. 
It is clear that $z$ is a Toeplitz sequence. 

On the other way round, let $(p_n)_{n\geq 1}$ be a  periodic structure of the Toeplitz sequence $x$.
Then for every $n \geq 1$ we can define the \emph{skeleton} of $x$ at scale $p_n$ by $z(n)_{m}$  equal to $y_m$ if $m \in {\rm Per}_{p_n}(y)$ and to $?$ otherwise. 
Since $z(n)$ has period $p_n$ and $p_n$ divides $p_{n+1}$, there exists a periodic sequence $y(n)$ such that \eqref{eq:holetoeplitz} holds.

We illustrate this constuction on the period doubling sequence.
\begin{example}
  Let $x$ be the period-doubling sequence.
  Consider the periodic sequence $y=(10?0)^\infty$. Then $Sx$ is the result of
  inserting $y$ into itself indefinitely as  shown below. Since
  $\sigma^2:0\mapsto 0100,1\mapsto 0101$, the sequence $Sx$ is the fixed point of
  $0\mapsto 1000,1\mapsto 1010$ and is
  \begin{displaymath}
0100\ 1010\ 1000\ 1010\cdot 1000\ 1010\ 1000\ 1000
  \end{displaymath}
  and the iterated insertion of $y$ in itself gives successively
  \newcommand{\red}[1]{\textcolor{red}{#1}}
  \begin{eqnarray*}
    10?0\       10?0\        10?0       \cdot 10?0\       10?0\       10?0\       10?0\       10?0\\
    10\red{0}0\ 10?0\        10\red{0}0 \cdot 10\red{1}0\ 10\red{0}0\ 10?0\       10\red{0}0\ 10\red{1}0\\
    1000 \      10\red{0 }0\ 1000       \cdot 1010\       1000\       10\red{1}0\ 1000\       1010
    \end{eqnarray*}
  \end{example}

%It is well known that a minimal shift space $(X,S)$ is a Toeplitz shift if it is 
%an almost 
%one-to-one extension of an odometer. 
%The odometer is given by $\mathbb{Z}_{(p_n)}$, where $(p_n)_{n\geq 1}$ is a periodic structure of a Toeplitz point $x$ generating $X$: $X = \overline{ \{ S^n x | n\in \mathbb{Z}\}}$.

\section{Exercises}
\exosection{Section \protect{\ref{sectionRecurrentMinimal}}}

\begin{exercise}\label{exerciseInverseContinuous}
Let $\phi:X\to Y$ be a continuous map between compact metric spaces $X,Y$.
Show that if $\phi$  is bijective, 
its inverse is continuous.
\end{exercise}
\begin{exercise}\label{exerciseRecurrent}
Let $(X,T)$ be a topological dynamical system. Show
that for two nonempty  sets $U,V\subset X$ and $n\ge 0$, one has
$U\cap T^{-n}V\ne\emptyset$ if and only if $T^nU\cap V\ne\emptyset$.
\end{exercise}
\begin{exercise}\label{exerciseRecurrentPoint}
Let $(X,T)$ be a topological dynamical system.
%A point $x\in X$ is \emph{recurrent}\index{subject}{recurrent!point} if
%for every open set $U$ containing $x$, there is an $n> 0$ such that
%$T^n(x)\in U$.
Show that if $x\in X$ is recurrent, for every open
set $U$ containing $x$, there are infinitely many  $n>0$
such that $T^nx$ belongs to $U$.
\end{exercise}

\begin{exercise}\label{exerciseCondRecurrent}
  Prove Proposition~\ref{propositionRecurrentSystem}.
  Hint: For (i)$\Rightarrow$ (ii), use Baire category Theorem
  (Theorem~\ref{theoremBaireCategory}).
\end{exercise}

\begin{exercise}\label{exerciseFactorMinimal}
Show that a factor of a minimal system is minimal.
\end{exercise}
\begin{exercise}\label{exerciseIrrationalRotation}
Let  $R_\alpha$ be the transformation $x\mapsto x+\alpha$
on the torus $\T=\R/\Z$. Show that the system $(X,T)$ is minimal
if and only if $\alpha$
is irrational .
%Hint: use the continued fraction expansion of $\alpha$.
\end{exercise}
\begin{exercise}\label{exerciseTowerConstruction}
Let $(X,T)$ be a topological dynamical system and let $h:X\to\N$
be a continuous function. The set
\begin{displaymath}
X^{h}=\{(x,i)\mid 0\le i< h(x)\}
\end{displaymath}
\index{symbols}{X@$X^h$}%
is a compact metric space as a closed subset of the product of finitely many
copies of $X$ (indeed, since $h$ is continuous on the
compact space $X$, it is bounded). Moreover, the map
\begin{displaymath}
T^{h}(x,i)=\begin{cases}(x,i+1)&\mbox{ if $i+1<h(x)$}\\
(Tx,0)&\mbox{ otherwise}\end{cases}
\end{displaymath}
is continuous. Thus $(X^h,T^h)$ is a topological dynamical system
said to be obtained from $X$ by
the \emph{tower construction}\index{subject}{tower!construction}
relative to $h$.
Show that:
\begin{enumerate}
\item
$T^h$ is invertible if, and only if, $T$ is invertible,
\item
$(X^h , T^h)$ is minimal if, and only if, $(X,T)$ is minimal,
\item
the system induced on $X\times\{0\}$ is isomorphic to $(X,T)$.
\end{enumerate}
\end{exercise}
%\exosection{Section~\ref{sectionCantorSpaces}}
\begin{exercise}\label{exerciseMorphismCantor}
Let $\phi: \{0,1\}^\N \to [0,1]$ be the map defined by $y=\phi(x)$ 
where
\begin{displaymath}
y=\sum_{n\ge 1}x_n2^{-n}.
\end{displaymath}
Show that $\phi$ is a morphism from the full shift $(\{0,1\}^\N ,S)$ to the dynamical system $([0,1],T)$
with $T(x)=2x\bmod 1$.
\end{exercise}
\begin{exercise}\label{exerciseMinimalSymbolicCantor}
Show that a minimal infinite symbolic system is a Cantor system.
\end{exercise}

\exosection{Section~\ref{sectionSymbolicSystems}}
\begin{exercise}\label{exercisePrimitive}
Show that a primitive word of length $n$ has $n$ distinct conjugates.
  \end{exercise}
\begin{exercise}\label{exerciseFineWilf}
Prove the following property of periods of words,
known as \emph{Fine-Wilf Theorem}.
\index{subject}{Fine-Wilf Theorem}\index{subject}{Theorem!Fine-Wilf}%
\index{names}{Fine, Nathan J.}\index{names}{Wilf, Herbert S.}%
Let $p,q\ge 1$ be integers and $d=\gcd(p,q)$.
If a word $w$ has periods $p$ and $q$ with $|w|\ge p+q-d$,
it has period $d$.
\end{exercise}
\begin{exercise}\label{exerciseProposition1.3.14}
  Let $x$ be a one-sided infinite sequence.
  A factor of  $x$ is \emph{conservative}
  \index{subject}{conservative factor}\index{subject}{factor!conservative}%
  if it is not right-special. Let $n\ge 1$ and let $c$
  be the number of conservative factors of $x$ of length $n$.
  Show that if $x$ has a factor of length $n+c$ whose
  factors of length $n$ are all conservative, then $x$ is eventually
  periodic.
  \end{exercise}
\begin{exercise}\label{exerciseFactorialExtendable}
Show that for every factorial extendable set $L$, there is a unique
shift space $X$ such that $\cL(X)=L$.
\end{exercise}
\begin{exercise}\label{exerciseFibonacciLR}
Show that the Fibonacci word is linearly recurrent.
\end{exercise}
\exosection{Section~\ref{sectionSFT}}
\begin{exercise}\label{exerciseMagicWords}
Show that a shift space $X$ is of finite type
\index{subject}{shift space!of finite type}%
 if and only if
there is an $n\ge 1$ such that every $v\in\cL_n(X)$ satisfies
\begin{equation}
uv,vw\in\cL(X)\Rightarrow uvw\in\cL(X)\label{eqMagicWords}
\end{equation}
for every $u,w\in\cL(X)$.
\end{exercise}
\begin{exercise}\label{exerciseConjugacySFT}
Show that the class of shifts of finite type is closed under conjugacy.
\index{subject}{shift space!of finite type}%
\end{exercise}
\begin{exercise}\label{exerciseEdgeShifts}
Prove Proposition~\ref{propositionSFTEdgeShift}.
\end{exercise}
\begin{exercise}\label{exerciseSoficShift}
Given a finite graph $G$ with edges labeled by letters of
an alphabet $A$, we denote by
$X_G$ the  set of labels of infinite paths in $G$.
A shift space $X$ on the alphabet $A$
is \emph{sofic}\index{subject}{sofic shift}\index{subject}{shift space!sofic}
if there is a finite labeled graph $G$, called a \emph{presentation}
\index{subject}{sofic shift!presentation of}\index{subject}{presentation!of sofic shift}%
of $X$, such that $X=X_G$.
Show that sofic shifts are the factors of shifts of finite
type.
\end{exercise}
\begin{exercise}\label{exerciseMinimalCover}
A labeled graph $G$ is \emph{right-resolving}
\index{subject}{right!resolving!graph}%
if the edges going out of of the same vertex have different labels.
A \emph{right-resoving presentation}
\index{subject}{right!resolving!presentation}\index{subject}{sofic shift!presentation of!right-resolving}%
of a sofic shift is a right-resolving graph $G$ such
that $X=X_G$. Show that 
\begin{enumerate}
\item every sofic shift has a right-resolving
presentation 
\item  every irreducible
sofic shift has  a unique minimal strongly connected 
right-resolving presentation.%\index{subject}{sofic shift!presentation of!right-resolving!minimal}%
\end{enumerate}
Hint: consider the follower sets $F(u)=\{v\in\cL(X)\mid uv\in \cL(X)\}$
for $u\in\cL(X)$. If $G=(Q,E)$ is a right-resolving presentation
of $X$, denote $I(u)=\{q\in Q\mid\mbox{ there is a path ending at $q$
labeled $u$}\}$ for $u\in\cL(X)$. Show that $I(u)=I(v)$
implies $F(u)=F(v)$.
\end{exercise}
\begin{exercise}\label{exerciseStrongShiftEquivalence}
Two nonnegative integral square matrices  $M,N$ are
\emph{elementary equivalent}\index{subject}{elementary!equivalence}
if there are nonnegative integral matrices $U,V$ such that
\begin{displaymath}
M=UV,\quad N=VU
\end{displaymath}
The matrices are \emph{strong shift equvalent}\index{subject}{strong!shift equivalence}\index{subject}{shift!equivalence!strong} if there is a sequence $(M_1,M_2,\ldots,M_k)$ of matrices
such that $M_1=M$, $M_k=N$ and $M_i$ elementary equivalent
to $M_{i+1}$ for $1\le i\le k-1$.

Let $M$ be a nonnegative $n\times n$-matrix. Denote by $X_M$
\index{symbols}{X@$X_M$}%
the edge shift on a graph with adjacency matrix $M$.
Show that if $M,N$ are strong shift equivalent,
then $X_M,X_N$ are conjugate.
\end{exercise}
\exosection{Section \protect{\ref{sectionSubstitutionSystems}}}
\begin{exercise}\label{exerciseDefSubstitution}
  Show that the following conditions are equivalent for a
  nonerasing morphism
  $\sigma:A^*\to A^*$.
  \begin{enumerate}
  \item[(i)] $\sigma$ is a substitution.
  \item[(ii)] The language $\cL(\sigma)$ is extendable.
  \item[(iii)] Every letter $a\in A$ is extendable in $L$.
    \end{enumerate}
\end{exercise}

\begin{exercise}\label{exercisesigma^n}
  Let $\sigma:A^*\to A^*$ be a morphism.
  Show that if $\sigma$ is primitive, then $\cL(\sigma^n)=\cL(\sigma)$
  for every $n\ge 1$. Give an example of a morphism such that $\cL(\sigma^2)$
  is strictly contained in $\cL(\sigma)$.
  \end{exercise}
\begin{exercise}\label{exerciseThueMorse}
Let $x$ be the Thue-Morse sequence.
Prove that $x_n=a$ if and only if the number of $1$ in the binary expansion
of $n$ is even.
\end{exercise}
\begin{exercise}\label{exerciseRudinShapiro}
Let $\varepsilon_n\in\{-1,1\}$ be the parity of the number
of (possibly overlapping) factors $11$ in the binary representation
of $n$. The sequence $\varepsilon_0\varepsilon_1\varepsilon_2\cdots$
is the \emph{Rudin-Shapiro sequence}.\index{subject}{Rudin-Shapiro sequence}
\index{sequence!Rudin-Shapiro} Show that it is the image under
the morphism $\phi:a\mapsto 1,b\mapsto 1,c\mapsto-1,d\mapsto -1$ of the fixed point
$x$
beginning with $a$ of the substitution $\sigma:a\mapsto ab,b\mapsto ac,c\mapsto db,d\mapsto dc$.
\end{exercise}
\begin{exercise}\label{exerciseChaconMinimal}
  Let $\sigma:0\mapsto 0010,1\mapsto 1$ be
  the Chacon binary
  substitution. Show that $\sigma$ is actually a substitution
  and  that the \emph{Chacon binary shift}
\index{subject}{Chacon!binary!shift}%
generated by $\sigma$  is minimal.
\end{exercise}
\begin{exercise}\label{exerciseGrowingMinimal}
Let $\sigma:A^*\to A^*$ be a growing morphism.
Show that if the shift space generated by $\sigma$ is minimal, then
$\sigma$ is primitive.
\end{exercise}

\begin{exercise}\label{exerciseElementarySubstitution}
A  morphism $\varphi:A^*\to B^*$ is 
\emph{elementary}
\index{subject}{elementary!morphism}\index{subject}{morphism!elementary}%
(not to be confused with the elementary automorphisms of Section~\ref{sectionSturmianShifts})
if it cannot 
be written $\varphi=\alpha\circ\beta$ with $\beta:A^*\to C^*$,
$\alpha:C^*\to A^*$ and $\Card(B)<\Card(A)$. 
Prove that an elementary nonerasing morphism defines an injective
 map from $A^\N$ to $A^\N$.
\end{exercise}
\begin{exercise}\label{exerciseFixedPointPeriodic1}
Prove that if $x$
is a periodic fixed point of a primitive elementary morphism, then every
letter $a\in A$ can be followed by at most one letter in $x$
and thus that the period of $x$ is at most $\Card(A)$.
\end{exercise}
\begin{exercise}\label{exerciseFixedPointPeriodic2}
Prove that if a fixed point $x$ of a primitive morphism 
$\varphi$ is periodic, then the period of $x$ is at most
$|\varphi|^{\Card(A)-1}$ where $|\varphi|=\max\{|\varphi(a)|\mid a\in A\}$
(Hint: Use Exercise~\ref{exerciseFixedPointPeriodic1} and
the fact that if $\varphi=\alpha\circ\beta$ and if $x$
is a periodic fixed point of $\varphi$, then $y=\beta(x)$ is a periodic
fixed point of $\psi=\beta\circ\alpha$ and that
the period of $x$ is at most the period of $y$ times 
$|\alpha|\le |\varphi|$).
\end{exercise}
\begin{exercise}\label{exerciseStableSubmonoid}
Prove that a submonoid $M$ of $A^*$ is generated by a code
if and only if it satisfies
\begin{equation}
u,uv,vw,w\in M\Rightarrow v\in M.\label{eqStability}
\end{equation}
for all $u,v,w\in A^*$.
\end{exercise}
\begin{exercise}\label{exerciseFlowerAutomaton}
Let $U\subset A^+$ be finite set of words. The
\emph{flower automaton}
\index{subject}{flower automaton}%
\index{subject}{automaton!flower}%
of $U$ is the following labeled graph $\A(U)=(Q,E)$.
\index{symbols}{A@$\A(U)$} The set $Q$ of vertices 
of $\A(U)$ is the set of pairs $(u,v)$ of nonempty words such that
$uv\in U$ plus the special vertex $\omega$. For every $u,a,v$
with $a\in A$ and $u,v\in A^*$ such that $uav\in U$ there
is an edge $e$ labeled $a$ with
\begin{enumerate}
\item[(i)] $e:(u,av)\mapsto(ua,v)$ if $u,v\ne\varepsilon$
\item[(ii)] $e:(\omega,av)\mapsto (a,v)$ if $u=\varepsilon$, $v\ne\varepsilon$,
\item[(iii)] $e:(u,a)\mapsto(ua,\omega)$ if $u\ne\varepsilon$, $v=\varepsilon$,
\item[(iv)] and finally $e:\omega\mapsto\omega$ if $u=v=\varepsilon$.
\end{enumerate}
Show that the number of paths labeled $w$ from $\omega$
to $\omega$ is equal to the number of factorizations of
$w$ in words of $U$. Show that $U$ is
\begin{enumerate}
\item  a code if and only if there is a unique path 
from $p$ to $q$ labeled $w$ for every $p,q\in Q$ and $w\in A^*$,
\item a circular code if for every nonempty word $w$ there is at most one
$p\in Q$ such that there is a
cycle labeled $w$ from $p$ to $p$.
\end{enumerate}
\end{exercise}
\begin{exercise}\label{exerciseRankOne}
Let $U$ be a  code.  A pair $(x,y)$ of words in $U^*$
is \emph{synchronizing}\index{subject}{synchronizing pair}
 if for every $u,v\in U^*$, one has
\begin{displaymath}
uxyv\in U^*\Rightarrow ux,yv\in U^*.
\end{displaymath}
Note that this definition is coherent with the definition of
synchronizing pair given in Section~\ref{sectionSubstitutionSystems}.
Consider the flower automaton $\A(U)=(Q,E)$.
Let $\mu$ be the morphism from $A^*$ into the monoid of $Q\times Q$-matrices
with integer elements defined by
\begin{displaymath}
\mu(w)_{p,q}=\begin{cases}1&\mbox{ if $p\edge{w}q$}\\
0&\mbox{ otherwise}\end{cases}
\end{displaymath}
Show that the following conditions are equivalent for $x,y\in U^*$
\begin{enumerate}
\item[(i)] $(x,y)$ is synchronizing.
\item[(ii)] $\mu(xy)_{p,q}=\mu(x)_{p,\omega}\mu(y)_{\omega,q}$ for all $p,q\in Q$.
and thus $\mu(xy)$ has rank one.
\end{enumerate}
Moreover, if $\mu(x),\mu(y)$ have rank one, then $(x,y)$ is synchronizing.
\end{exercise}
\begin{exercise}\label{exerciseUniformSync}
A code $U\subset A^+$ is said to have \emph{finite synchronization delay}
\index{subject}{synchronization delay} $n$
(or to be \emph{uniformly synchronized})
\index{subject}{uniformly!synchronized code}%
\index{subject}{code!uniformly synchronized}%
if  there is an integer $n$ such that every pair $x,y$ of words in $U^*$ of length at least $n$
is synchronizing.

Show that the following conditions are equivalent
for a finite code $U$ on the alphabet $A$.
\begin{enumerate}
\item[(i)] $U$ is a circular code.
\item[(ii)] $U$ has finite synchronization delay.
\item[(iii)] every sequence in $A^\Z$ has at most one factorization
in words of $U$, that is for every $x\in A^\Z$ there is
at most one pair $(k,y)$ with $y\in B^\Z$ and $0\le k<|\varphi(y_0)|$
such that $x=S^k\varphi(y)$.
\end{enumerate}
Hint: for (i) $\Rightarrow$ (ii), use Exercise~\ref{exerciseRankOne}.
\end{exercise}
\begin{exercise}\label{exerciseFiniteToOne}
  Let $\varphi:A^*\to B^*$ be a nonerasing morphism. Let $X=A^\Z$,
  $Y=B^\Z$
  and
  $X^\varphi=\{(x,i)\mid x\in A^\Z, 0\le i<|\varphi(x_0)\}$ as in
  Section \ref{sectionRecognizable} and let $\varphi:X^\varphi\to B^\Z$
  be defined by $\hat{\varphi}(x,i)=S^i\varphi(x)$.
  Show that $\varphi$ is injective if and only if $\hat{\varphi}$
  is finite-to-one, that is, for every $y\in Y$, the set $\hat{\varphi}^{-1}(y)$
  is finite.\index{subject}{finite-to-one map}
 
\end{exercise}

\begin{exercise}\label{exerciseMosseOriginal}
  Let $\varphi:A^*\to A^*$ be a nonerasing morphism with
  an admissible fixed point $x\in A^\N$. Let $f:\N\to\N$
  be defined by $f(n)=|\varphi(x_{[0,n)})|$. Show
    that $\varphi$ is recognizable if and only if the
    following condition is satisfied. There is an integer
    $\ell>0$ such that whenever $x_{[i-\ell,i+\ell)}=x_{[j-\ell,j+\ell)}$
        and $i\in f(\N)$, then $j\in f(\N)$ and $x_{f^{-1}(i)}=x_{f^{-1}(j)}$.
\end{exercise}
\exosection{Section \ref{sectionSturmianShifts}}
\begin{exercise}\label{exerciseFiboisSturmian}
Show that the Fibonacci shift is Sturmian. Hint: show that the
left special words are the prefixes of the $\varphi^n(a)$.
\index{subject}{Fibonacci!shift}\index{subject}{shift space!Fibonacci}%
\end{exercise}
\begin{exercise}\label{exerciseBalanced1}
  Let $U$ be a factorial set over $\{0,1\}$.
  Show that if $U$ is balanced,
then $\Card(U\cap\{0,1\}^n)\le n+1$.
\end{exercise}
\begin{exercise}\label{exerciseBalanced2}
  Show that if $x$ is unbalanced, there is a palindrome $w$
  such that $0w0$ and $1w1$ belong to $\cL(x)$.
  \end{exercise}
\begin{exercise}\label{exerciseEquivalentDefSturm}
  Prove Proposition~\ref{propositionEquivalentDefSturm}.
  Hint: use Exercises~\ref{exerciseBalanced1}
  and \ref{exerciseBalanced2}.
\end{exercise}
\begin{exercise}\label{exerciseSturmReversal}
  Show that the set of factors of a Sturmian sequence is closed under
  reversal.
  \end{exercise}
\begin{exercise}\label{exerciseTriboisSturmian}
  Prove that the Tribonacci sequence is a strict standard episturmian
  sequence.
\index{subject}{Tribonacci!shift}\index{subject}{shift space!Tribonacci}%
\end{exercise}

\begin{exercise}\label{exerciseMechanicalIsSturm}
  Show that for $\alpha,\rho$ with $0<\alpha< 1$ and $0\le \rho<1$
  with $\alpha$ irrational, the sequence $s=(s_n)_{n\ge 0}$ defined
  by $s_n=\lfloor(n+1)\alpha+\rho\rfloor - \lfloor n\alpha+\rho\rfloor$
  is Sturmian. Moreover, define the
  \emph{slope}\index{subject}{slope!of word} of a word $u$ on $\{0,1\}$ as the
  real number $\pi(u)=|u|_1/|u|$. Show that
  for every factor $u$ of $s$, the slope of $u$ tends to $\alpha$ when $|u|\to\infty$
  (justifying the name of slope for $\alpha$, see Figure~\ref{figureSlope}).
  \end{exercise}
\begin{exercise}\label{exerciseJustin0}
  Show that for every word $u$, the palindromic closure of $u$
  is $u^{(+)}=uv^{-1}\tilde{u}$ where $v$ is the longest palindrome suffix of $u$.
  \end{exercise}
\begin{exercise}\label{exerciseJustin1}
  Prove that for every $w\in A^*$ and $a\in A$, one has
  \begin{displaymath}
    \Pal(wa)=\begin{cases}\Pal(w)a\Pal(w)&\mbox{ if $|w|_a=0$}\\
    \Pal(w)\Pal(w_1)^{-1}\Pal(w)&\mbox{ if $w=w_1aw_2$ with $|w_2|_a=0$}
    \end{cases}
    \end{displaymath}
\end{exercise}
\begin{exercise}\label{exerciseJustin2}
  Prove that if $p$ is a palindrome and $a$ a letter, then $L_a(p)a$
  is a palindrome.
\end{exercise}
\begin{exercise}\label{exerciseJustin3}
  Show that for $a\in A$ and $w\in A^*$, one has
  \begin{equation}
    \Pal(aw)=L_a(\Pal(w))a.\label{equationJustin3}
    \end{equation}
  \end{exercise}
\begin{exercise}\label{exerciseJustin4}
  Prove Justin Formula $\Pal(vw)=L_v(\Pal(w))\Pal(v)$.
  \end{exercise}
\begin{exercise}\label{exerciseReturnSturm}
Let $s=\Pal(x)$ be a standard  episturmian sequence
with directive sequence $x=x_0x_1\cdots$ and let $X$  be the
shift genrerated by $s$.
Let  $u_n=\Pal(x_0x_1\cdots x_{n-1})$ for $n\ge 1$
be the palindromic prefixes of $s$.
Show that one has
\begin{equation}
\RR'_X(u_n)\subset\{L_{x_0\cdots x_{n-1}}(a)\mid a\in A\}
\end{equation}
with equality if $s$ is strict. Show that for $w\in\cL(s)$, one
has $\RR'_X(w)=y^{-1}\RR'_X(u_n)y$ where $y$ is the shortest
word such that $yw$ is a prefix of $u_n$ for some $n\ge 1$.
\end{exercise}
\begin{exercise}\label{exerciseEquivalentDefAR}
  Show that a minimal shift $X$ on the alphabet $A$ is strict episturmian (that is, Arnoux-Rauzy)
  if and only if, for every $n\ge 1$, there is a unique right-special
  (resp. left-special) factor $w$ of length $n$ such that
  $wa\in\cL(X)$ (resp. $aw\in\cL(X)$) for every $a\in A$.
  \end{exercise}
\begin{exercise}\label{exerciselemma6}
  Let $s=\Pal(x)$ for some infinite sequence $x$. Set $x=x_0x_1\cdots$
  and let $u_n=\Pal(x_0\cdots x_{n-1})$. Let $u$ be a factor of $s$
  and let $n\ge 0$ be such that $u$ is a factor of $u_{n+1}$
  but not of $u_n$. Set $u_{n+1}=tw\tilde{t}$ where $w$ is the unioccurrent
  palindrome suffix of $u_nx_n$. Show that
  \begin{enumerate}
  \item The leftmost occurrence of $u$ in $u_{n+1}$ is $u_{n+1}=yuy'=yzwz'y'$
    for some words $y,y',z,z'$.
  \item If $u=avb$ with $v$ a palindrome prefix of $s$, then $u$
    is either a suffix of $u_nx_n$ or a prefix of $x_nu_n$.
    \end{enumerate}
  \end{exercise}
\begin{exercise}\label{exerciseStandard}
Let $(d_1,d_2,\ldots)$ be a sequence of integers with $d_1\ge 0$
and $d_n>0$ for $n\ge 2$.
The \emph{standard sequence}\index{subject}{standard!sequence}
 with \emph{directive sequence}\index{subject}{directive!sequence!of
standard sequence} $(d_1,d_2,\ldots)$
is the sequence $(s_n)$ of words defined by $s_0=0$, $s_1=0^{d_1}1$
and $s_n=s_{n-1}^{d_n}s_{n-2}$ for $n\ge 2$. Show that each
$s_n$ is a primitive word which is prefix of the characteristic sequence of slope $\alpha=[0,1+d_1,d_2,\ldots]$.
\end{exercise}
\begin{exercise}\label{cfrac1}
Let $(s_n)$ be the standard sequence with directive
sequence $(d_1,d_2,\ldots)$
and let $\alpha=[0,1+d_1,d_2,\ldots]$.
 Show that for $n\ge 3$, the
word $s_n^{1+d_{n+1}}$ is a prefix of the characteristic word $c_\alpha$ but not the
word $s_n^{2+d_{n+1}}$.
\end{exercise}
\begin{exercise}\label{exercisecfrac2}
A sequence $x$ is said to be $d$-\emph{power free}\index{subject}{d-power free@$d$-power free sequence}\index{subject}{sequence!$d$-power free} if for every nonempty
word $w$, $w^n\in\cL(x)$ implies $n<d$. Show that if a Sturmian
sequence of slope $\alpha=[a_0,a_1,a_2,\ldots]$ is $d$-power free 
for some $d$,
the $a_i$ are bounded (note that the converse
is also true, see the Notes). Hint : use Exercise~\ref{cfrac1}.
\end{exercise}
\begin{exercise}\label{exerciseSturmLR}
Show that if a Sturmian sequence of slope $\alpha=[a_0,a_1,\ldots]$
is linearly recurrent, the $a_i$ are bounded.
\end{exercise}
\exosection{Section~\ref{sectionToeplitz}}
\begin{exercise}\label{exerciseDoublingSequence}
Let $x\in\{0,1\}^\N$ be the period-doubling sequence,
\index{subject}{period-doubling sequence} which is the fixed
point of the substitution $\sigma:0\mapsto 01,1\mapsto 00$. Show that
\begin{displaymath}
x_n=\nu_2(n+1)\bmod 2
\end{displaymath}
where $\nu_2(m)$ is the number of $0$ ending the binary representation 
of $m$.
\end{exercise}

\begin{exercise}
\label{exerciseShiftToeplitz}
Let $\sigma $ be a primitive substitution of constant length $n$ having a coincidence at index $k$.
Then $(X_\sigma , S)$ is a Toeplitz shift.\index{subject}{coincidence}
\end{exercise}

%%%%%%%%%%%%%%%%%%%%%%%
\section{Solutions}
\exosection{Section \protect{\ref{sectionRecurrentMinimal}}}

\begin{solution}{\ref{exerciseInverseContinuous}}
Let $F$ be a closed set of $X$.
Then, it is compact.
Thus, the set $(\phi^{-1})^{-1} (F) = \phi (F)$ is compact and hence closed proving that $\phi^{-1}$ is continuous. 
\end{solution}

\begin{solution}{\protect{\ref{exerciseRecurrent}}}
Assume that $U\cap T^{-n}V$ is nonempty. 
Let $x$ in $U\cap T^{-n}V\ne\emptyset$. 
Then $T^nx$ belongs to $T^nU\cap V$. 
Conversely, if $y$ is in $T^nU\cap V$, there
is some $x\in U$ such that $T^n(x)=y$.
Thus $x$ is an element of $U\cap T^{-n}V$.
\end{solution}

\begin{solution}{\protect{\ref{exerciseRecurrentPoint}}}
If $T^nx$ is in $U$, then $x$ belongs to $U\cap T^{-n}U$. 
Since the latter is an open set, there is an $m\geq 0$ such that $T^mx$ belongs to $U\cap T^{-n}U$
and thus $T^{n+m}x$ is in $U$. 
The same argument can be repeated to obtain the conclusion.
\end{solution}

\begin{solution}{\protect{\ref{exerciseCondRecurrent}}}
(i)$\Rightarrow$(ii) Let $(U_n)_{n\ge 0}$ be a countable basis
  of open sets (this exists for any  metric space). For every
  $n\ge 0$, the open set $V_n=\cup_{m\ge 0}T^{-m}U_n$ is
  dense since $(X,T)$ is recurrent. The set
\begin{displaymath}
V=\cap_{n\ge 0}V_n.
\end{displaymath}
is formed of the points with a dense positive orbit.
By Baire Category Theorem (Theorem~\ref{theoremBaireCategory}),
the set $V$ is dense in $X$.
%(ii)$\Rightarrow$ (iii) \marginpar{A COMPLETER}

(ii)$\Rightarrow$ (i) Let $x_0\in X$ be a point with a dense positive orbit.
For every pair $U,V$ of nonempty open sets, there are, by Exercise~\ref{exerciseRecurrentPoint}, arbitrary large integers $n,m$ such that $T^nx$ is in $U$ and $T^mx$ in $V$.
Choosing $n<m$, we obtain the set $U\cap T^{m-n}V$ is nonempty.
\end{solution}

\begin{solution}{\protect{\ref{exerciseFactorMinimal}}}
Let $\phi:(X,T)\to(X',T')$ be a factor map from a minimal system $(X,T)$ to
$(X',T')$. Let $Y'$ be a closed stable nonempty subset of $X'$. Then
$Y=f^{-1}(Y')$ is nonempty and 
closed. It is also stable because for $y\in Y$, 
we have $\phi(T(y))=T'(\phi(y))\in T'Y'\subset Y'$
and thus $T(y)$is in $Y$. Hence $Y=X$ which implies $Y'=X'$.
\end{solution}

\begin{solution}{\ref{exerciseIrrationalRotation}}
Assume that $\alpha$ is irrational.
For every $q\ge 1$, there is $p\in \N$ such that $\alpha $ is in $(p/q, (p+1)/q)$.
Then every $x\in [0,1)$ belongs to an interval $[np/q,(n+1)p/q)$, for some $n$,
and thus $|x-n\alpha|\le 1/q$.

Suppose $\alpha $ is rational with $\alpha =p/q$. 
Then the orbit $(n\alpha )_n $ of $0$ in $\mathcal{T}$ will visit finitely many points and thus cannot be dense in $\mathcal{T}$.

%Let 
%\begin{displaymath}
%\alpha=d_0+\cfrac{1}{d_1+\cfrac{1}{d_2+\ddots}}
%\end{displaymath}
%% where $d_0\ge 0$ and $d_i\ge 1$ are integers,
%be the continued fraction expansion
%of $\alpha$. We denote $\alpha=[d_0,d_1,\ldots]$.
%Set $p_{-1}=1$, $q_{-1}=0$, $p_0=d_0$, $q_0=1$ and inductively
%\begin{displaymath}
%p_n=d_np_{n-1}+p_{n-2},\quad q_n=d_nq_{n-1}+q_{n-2}
%\end{displaymath}
%in such a way that
%$p_n/q_n=[d_0,d_1,\ldots,d_n]$. Then
%\begin{displaymath}
%\alpha-\frac{p_n}{q_n}=\frac{(-1)^n}{q_nq'_{n+1}}
%\end{displaymath}
%where $q'_{n+1}=d'_{n+1}q_n+q_{n-1}$ with $d'_{n+1}=[d_{n+1},d_{n+2},\ldots]$.
%Since $q'_{n+1}>d_{n+1}q_n+q_{n-1}=q_{n+1}$, we obtain
%\begin{displaymath}
%|\alpha-\frac{p_n}{q_n}|\le \frac{1}{q_nq_{n+1}}.
%\end{displaymath}
%This shows that $\alpha$ is arbitrary close to a rational
%and thus that the numbers $\lfloor n\alpha\rfloor$ are dense
%in $[0,1]$.
\end{solution}
\begin{solution}{\ref{exerciseTowerConstruction}}
If $T$ is invertible, the inverse of $T^h$ is the map which sends $(x,i)$ to $(x,i-1)$
if $i>1$ and to $(T^{-1}x,1)$ otherwise.

Conversely if $T^h$ is invertible the inverse of $T$ is the map which sends $x$ to the unique $(y,h(y)-1)$ such that $T^h (y,h(y)-1)=(x,0)$.

If $(X,T)$ is minimal, the orbit of every $(x,i)\in X^h$ is clearly dense in $X^h$.
Thus $(X^h,T^h)$ is minimal.

If $(X^h , T^h)$ is minimal then each orbit is dense in $X\times \{ 0\}$ and thus $(X,T)$ is minimal. 

The last assertion result from the fact that $f(x)$ is the return time to $U=X\times\{0\}$
and thus
\begin{displaymath}
(T^h)_U(x,0)=(T^h)^{h(x)}(x,0)=(Tx,0).
\end{displaymath}
\end{solution}

%\exosection{Section~\ref{sectionCantorSpaces}}
\begin{solution}{\ref{exerciseMorphismCantor}}
The map $\phi$ is clearly continuous and satisfies $\phi\circ T=S\circ \phi$.
\end{solution}

\begin{solution}{\ref{exerciseMinimalSymbolicCantor}}
Let $(X,T)$ be an infinite minimal symbolic system. Since $X$
is a closed subset of $A^\Z$, it is compact and totally disconnected. 
Suppose $x$ is an isolated point: There is a open set $U$ such that $U\cap X = \{ x\}$.
As each orbit $(T^n (y))_n$ is dense there is some $n<m$ such that $T^n (y)$ and $T^m (y)$ belongs to $U$ and consequently $T^n (y)=T^m (y) = x$.
Thus the orbit of $x$ should be $\{ x , T(x) , \dots , T^{m-n-1} (x)\}$.
By minimality this would imply $X = \{ x , T(x) , \dots , T^{m-n-1} (x)\}$, a contradiction.
\end{solution}
\exosection{Section~\ref{sectionSymbolicSystems}}
\begin{solution}{\ref{exercisePrimitive}}
  Assume that $w$ is a primitive word of length $n$ and that
  $w=uv=vu$. Then for every $k\ge 1$ we have $w^k=u^kv$,
  as shown easily by induction on $k$. This implies that
  $w^\omega=u^\omega$ and thus $w=u$.
  \end{solution}

\begin{solution}{\ref{exerciseFineWilf}}
Assume $p<q$. We use induction on $p+q$. The result is trivial
if $p=d$. Otherwise, let $u$ be the prefix of
length $q-d$ of $w$. Then, for $1\le i\le p-d$, we have
$u_i=u_{i+p}=u_{i+p-q}$. This shows that $u$ has period $q-p$ and $p$.
By induction hypothesis, $u$ has period $d$. Since $|u|\ge p$,
the word $w$ has also period $d$.
\end{solution}
\begin{solution}{\ref{exerciseProposition1.3.14}}
  Let $u=a_0a_1\cdots a_{n+c-1}$ be a factor of $x$ of length $n+c$
  whose factors of length $n$ are all conservative. Set $x=vuy$ and
  for $0\le i\le c$, set $p_i=a_ia_{i+1}\cdots a_{i+n-1}$. Since there
  are only $c$ conservative factors of $x$ of length $n$, there
  are indices $i,j$ with $0\le i<j\le c$ such that $p_i=p_j$.
  Then $x=va_0\cdots a_{i-1}(p_ia_{i+n}\cdots a_{j+n-1})^\omega$.
  \end{solution}
\begin{solution}{\ref{exerciseFactorialExtendable}}
Let $X$ be the set of $x\in A^\Z$ with all its factors in $L$. Clearly
$X$ is the largest shift space such that
$\cL(X)\subset L$. For $u\in L$, the language $L$ being extendable, there are sequences $(a_n)_n A^\N$ and $(b_n)_n \in A^\N$
such that $a_n\cdots a_1ub_1b_2\cdots b_n\in L$ for every $n\ge 0$. Then
$\cdots a_2a_1\cdot ub_1b_2\cdots$ is in $X$ and thus
$u\in\cL(X)$.
\end{solution}
\begin{solution}{\ref{exerciseFibonacciLR}}
Let $\varphi:a\mapsto ab,b\mapsto a$ be the Fibonacci morphism and let $x$
be the Fibonacci word. Let $F_n$
be the Fibonacci sequence defined by $F_0=0$, $F_1=1$ and $F_{n+1}=F_n+F_{n-1}$
for $n\ge 1$. Note that $F_n\le F_{n+1}$ implies $F_{n+1}\le 2F_n$ and
$F_{n+2}\le 3F_n$  for $n\ge 1$.

For every $n\ge 1$, we have $\varphi^{n+1}(a)=\varphi^n(a)\varphi^n(b)=\varphi^n(a)\varphi^{n-1}(a)$.
Thus $|\varphi^n(a)|=F_{n+2}$ and $|\varphi^n(b)|=F_{n+1}$.

Let $w\in\cL(x)$ and let $n$ be the least integer such that
$|w|< F_{n+1}$. By the choice of $n$, we have $F_n\le |w|$. 

Since $w$ has
no factor in $\varphi^{n}(A)$, it is a factor of $\varphi^{n}(A^2)$.
But the largest difference between the occurrences of a word of
length $2$ in the Fibonacci word is $5$ (this bound is
reached by $aa$ in $aababaa$). Thus two occurrences of $w$
in $x$ are separated by at most $5F_{n+2}\le 15F_{n}\le 15|w|$.
This shows that $x$ is linearly recurrent.
\end{solution}
\exosection{Section~\ref{sectionSFT}}
\begin{solution}{\ref{exerciseMagicWords}}
Assume first that $(X,S)$ is a shift of finite type defined
by a finite set $I$ of forbidden blocks. Let $n$ be the
maximal length of the words of $I$. It is clear that every
$v\in\cL_n(X)$ satisfies~\eqref{eqMagicWords}.

Conversely, consider the Rauzy graph $G=\Gamma_n(X)$. By condition
\eqref{eqMagicWords}, the label of every every infinite path in $G$ 
is in $X$. Thus the edge shift
on $G$ can be identified with $(X,S)$.
\end{solution}
\begin{solution}{\ref{exerciseConjugacySFT}}
Let $\varphi:X\to Y$ be a sliding block code from $(X,S)$
to a shift of finite type $(Y,S)$ which is a conjugacy. We may suppose that
$\varphi,\varphi^{-1}$ are defined by  block maps $f,g$ with memory 
and anticipation $\ell$. Let $k$ be the integer such
that every word in $\cL_k(Y)$ satisfies \eqref{eqMagicWords}.
Set $m=k+4\ell$.
Then, one may verify that every word in $\cL_m(X)$
satisfies \eqref{eqMagicWords}. Thus $(X,S)$ is a shift
of finite type.
\end{solution}
\begin{solution}{\ref{exerciseEdgeShifts}}
Let $I\subset A^*$ be a finite set of forbidden blocks
and let $(X,S)$ be the shift of finite type corresponding to $I$.
Let $n$ be the maximal length of the words in $I$ and consider
the Rauzy graph $\Gamma_n(X)$.
Then the edge shift on $\Gamma_n(X)$ is the $n$-th higher block presentation
of $(X,S)$. If $(X,S)$ is recurrent, then $\Gamma_n(X)$ is strongly connected.
\end{solution}
\begin{solution}{\ref{exerciseSoficShift}}
Suppose first that $X$ is sofic and let $G=(V,E)$ be a finite labeled graph
such that $X=X_G$. The one-block map assigning to an edge its label
is a factor map from the edge shift on $G$ onto $X$. Thus $X$
is a factor of a shift of finite type.

Conversely, let $Y$ be a shift of finite type, which may be assumed
to be an edge shift.
 Let$X$ be the image
of $Y$ by a sliding block code $f:\cL_{m+n-1}(X)$
with memory $m$ and anticipation $n$. Consider the
$n+m-1$-th higher block presentation $Z$ of $Y$. Then $Z$ is conjugate
to $Y$ and thus is a shift of finite type by 
Exercise~\ref{exerciseConjugacySFT}. Let $G$ be the Rauzy
graph $\Gamma_{m+n-1}(Y)$. Then $Y$ can be identified with
the edge shift on $G$. Let $H$ be the graph which the same
as $G$ but with the labeling defined by the block map $f$.
Clearly $X=X_H$ and thus $X$ is a sofic shift.
\end{solution}
\begin{solution}{\ref{exerciseMinimalCover}}
Let $X=X_G$ be a sofic shift where $G=(V,E)$ is a finite graph with labels
in $A$. Let $H=(W,F)$ be the following graph. The set
$W$ is the set $\Pg(V)$ of subsets of $V$. For $P,Q\subset V$, 
there is an edge
from $P$ to $Q$ labeled $a$ if 
\begin{displaymath}
Q=\{q\in V\mid \mbox{there is an edge $p\edge{a}q$ with $p\in P$}\}.
\end{displaymath}
Then $H$ is right-resolving and it is easy to see that $X=X_H$.

Assume now that $X$ is an irreducible sofic shift. Any
minimal right-resolving presentation of $X$ is clearly
strongly connected. 

Set $L=\cL(X)$.
Let $M$ be the following labeled graph. Its vertices 
are the \emph{follower sets}\index{subject}{follower set}
 \begin{displaymath}
F(u)=\{v\in L\mid uv\in L\}
\end{displaymath}
for $u\in L$. There is an
edge from $p$ to $q$ labeled $a$ if $p=F(u)$ and $q=F(ua)$.
Let us show that any strongly connected right-resolving 
minimal presentation $G=(Q,E)$
of $X$
can be identified with $M$.  For every $p\in Q$, let
$F(p)$ be the set of all $u\in L$ such that there
is a path with label $u$ starting at $p$. 
 Define an equivalence on $Q$ by $p\sim q$
if $F(p)=F(q)$. The quotient $Q/\sim$ is clearly again
a right-resolving presentation of $X$. Since $G$
is minimal, the equivalence $\sim$ is the equality.

Consider, for $u\in L$, the set $I(u)$
of all $q\in Q$ such that there is a path in $G$ with label $u$ 
ending at $q$. Clearly, if $I(u)=I(v)$,
then $F(u)=F(v)$. Take $u\in L$ such that $I(u)$ is of minimal
cardinality. For every $p,q\in I(u)$, we have $p\sim q$.
Thus $I(u)$ has only one element. Since $G$ is strongly connected,
every element of $Q$ appears in this way. This allows us
to associate to every  $p\in Q$  the follower set $F(u)$
where $u\in L$ is such that $I(u)=\{p\}$. This identifies
$G$ with $M$.
\end{solution}
\begin{solution}{\ref{exerciseStrongShiftEquivalence}}
It is enough to consider the case of two elementary
equivalent matrices $M=UV$ and $N=VU$. Let $G_M=(V_M,E_M)$
and $G_N=(V_N,E_N)$ be graphs with matrices $M,N$ respectively.
Let also $G_U$ be the graph on $V_M\cup V_N$ having $U_{xy}$ edges from 
$x\in V_M$ to $y\in V_N$ and similarly for $G_V$. 

Since $M=UV$, there is a bijection $e\mapsto (u(e),v(e))$
from $E_M$ onto the paths of length $2$ made of an edge of $G_U$
followed by an edge of $G_V$. We denote $e(u,v)$ the inverse map.
Similarly, we have a bijection $f\mapsto (v(f),u(f))$ from
$E_N$ onto the paths formed of an edge of $G_V$ followed
by an edge of $G_U$ with an inverse map denoted $f(v,u)$.

We define a $2$-block map $s$ (with memory $0$) from $\cL_2(X_M)$ to $E_N$ by
$s(e_0e_1)=f(v(e_0)u(e_1))$. Let $\sigma:X_M\to X_N$ be the
sliding block code defined by $s$. Similarly, let
$t$ be the $2$-block map (with memory $0$) from $\cL_2(X_N)$ to $E_M$ defined by
$t(f_0f_1)=e(u(f_0)v(f_1))$. Let $\tau:X_N\to X_M$ be the corresponding sliding
block code. We have (see Figure~\ref{figureElementaryShift})
\begin{displaymath}
\tau\circ \sigma=S_M
\end{displaymath}
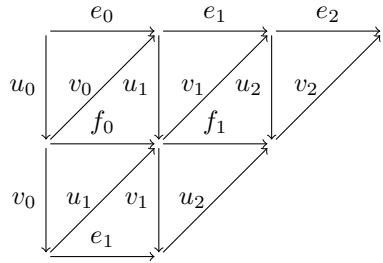
\begin{figure}[hbt]
\centering
\tikzset{node/.style={circle,draw,minimum size=0.1cm,inner sep=0pt}}
\tikzset{title/.style={minimum size=0.1cm,inner sep=0pt}}
\begin{tikzpicture}
\node[title](03)at(0,4.5){};\node[title](13)at(1.5,4.5){};\node[title](23)at(3,4.5){};
\node[title](33)at(4.5,4.5){};
\node[title](02)at(0,3){};\node[title](12)at(1.5,3){};\node[title](22)at(3,3){};
\node[title](01)at(0,1.5){};\node[title](11)at(1.5,1.5){};
%\node[title](00)at(0,0){};
\draw[->,above](03)edge node{$e_0$}(13);
\draw[->,above](13)edge node{$e_1$}(23);
\draw[->,above](23)edge node{$e_2$}(33);
\draw[->,left](03)edge node{$u_0$}(02);
\draw[->,left](13)edge node{$u_1$}(12);
\draw[->,left](23)edge node{$u_2$}(22);
\draw[->,left](02)edge node{$v_0$}(13);%diag
\draw[->,left](12)edge node{$v_1$}(23);%diag
\draw[->,left](22)edge node{$v_2$}(33);%diag
\draw[->,above](02)edge node{$f_0$}(12);
\draw[->,above](12)edge node{$f_1$}(22);
\draw[->,left](02)edge node{$v_0$}(01);
\draw[->,left](12)edge node{$v_1$}(11);
\draw[->,left](01)edge node{$u_1$}(12);%diag
\draw[->,left](11)edge node{$u_2$}(22);%diag
\draw[->,above](01)edge node{$e_1$}(11);
%\draw[->,left](01)edge node{$u_1$}(00);
%\draw[->,left](00)edge node{$v_1$}(11);%diag
\end{tikzpicture}
\caption{The conjugacies $\sigma$ and $\tau$.}\label{figureElementaryShift}
\end{figure}
where $S_M$ is the shift transformation on $X_M$.
\end{solution}
\exosection{Section \protect{\ref{sectionSubstitutionSystems}}}
\begin{solution}{\ref{exerciseDefSubstitution}}
  (i) $\Rightarrow$ (ii) is clear since $\cL(\sigma)=\cL(X)$
  implies that $\cL(\sigma)$ is extendable.

  (ii) $\Rightarrow$ (iii) is obvious.

  (iii) $\Rightarrow$ (ii). Consider $u\in\cL(\sigma)$. Let
  $a\in A$ and $n\ge 1$ be such that $u$ is a factor of $\sigma^n(a)$.
  Set $\sigma^n(a)=pus$ with $p,s\in A^*$. Let $b,c\in A$ be such that
  $bac\in\cL(X)$ and let $m\ge 1$, $d\in A$ be such that $bac$
  is a factor of $\sigma^m(d)$. Then $\sigma^n(b)pus\sigma^n(c)$
  is a factor of $\sigma^{n+m}(d)$. Since $\sigma$ is nonerasing,
  we have $euf\in\cL(\sigma)$ where $e$ is the last
  letter of $\sigma^n(b)p$ and $f$ is the first letter of $s\sigma^n(c)$.
  This shows that $\cL(\sigma)$ is extendable.

  (ii) $\Rightarrow$ (i). Let $u\in\cL(\sigma)$.
  An induction on $n$ proves that there exists a sequence
  $(a_n,b_n)$ of letters such that $a_n\cdots a_0ub_0\cdots b_n\in\cL(\sigma)$
  for every $n\ge 0$. Then
  $x=\cdots a_n\cdots a_0ub_0\cdots b_n\cdots$
  is in $X$ and thus $\cL(\sigma)\subset\cL(X)$. The other
  inclusion holds for any morphism.
\end{solution}

\begin{solution}{\ref{exercisesigma^n}}
   Assume that $\sigma$ is primitive.
  Let $m\ge 1$ be such that $|\sigma^m(b)|_a>0$ for every $a,b\in A$.
  Consider $w\in\cL(\sigma)$. Let $p\ge 1$ and $a\in A$ be such that
  $w$ is a factor of $\sigma^p(a)$. Let $q\ge m$ be such that $p+q$ is
  a multiple of $n$. Then $w$ is a factor of $\sigma^{p+q}(a)$ and
  thus $w$ is in $\cL(\sigma^n)$.

  Let $\sigma:a\mapsto bb,b\mapsto c,c\mapsto c$. Then $bb$ is in $\cL(\sigma)$
  although it is not in $\cL(\sigma^2)$.
  \end{solution}
\begin{solution}{\ref{exerciseThueMorse}}
Let $\bar{a}=b$ and $\bar{b}=a$. 
With this notation, the Thue-Morse morphism is defined by $\tau(x)=x\bar{x}$
for every $x\in \{a,b\}$. The
property is true for $n=0,1$. Next, since $\tau(x)=x$, and $\tau(x_0\cdots x_{n-1})=x_0\cdots x_{2n-1}$,
we have $\tau(x_n)=x_{2n}x_{2n+1}$. 
Thus $x_{2n}=x_n$ and $x_{2n+1}=\bar{x_n}$. This proves the property
by induction on $n$. 
\end{solution}
\begin{solution}{\ref{exerciseRudinShapiro}}
Consider the labeled graph represented in Figure~\ref{figureRudinShapiro}.
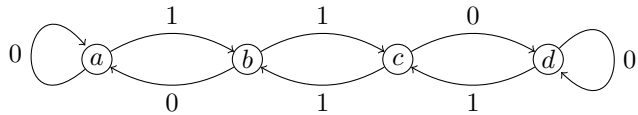
\begin{figure}[hbt]
\centering
\tikzset{node/.style={circle,draw,minimum size=0.4cm,inner sep=0pt}}
\tikzset{title/.style={minimum size=0.1cm,inner sep=0pt}}
\tikzstyle{every loop}=[->,shorten >=1pt,looseness=12]
        \tikzstyle{loop left}=[in=130,out=220,loop]
        \tikzstyle{loop right}=[in=-45,out=45,loop]
\begin{tikzpicture}
\node[node](a)at(0,0){$a$};\node[node](b)at(2,0){$b$};
\node[node](c)at(4,0){$c$};\node[node](d)at(6,0){$d$};

\draw[left,->](a)edge[loop left] node{$0$}(a);
\draw[->,above,bend left](a)edge node{$1$}(b);
\draw[->,below,bend left](b)edge node{$0$}(a);
\draw[->,above,bend left](b)edge node{$1$}(c);
\draw[->,below,bend left](c)edge node{$1$}(b);
\draw[->,above,bend left](c)edge node{$0$}(d);
\draw[->,below,bend left](d)edge node{$1$}(c);
\draw[right,->](d)edge[loop right] node{$0$}(d);
\end{tikzpicture}
\caption{The automaton of the Rudin-Shapiro sequence.}\label{figureRudinShapiro}
\end{figure}

Let $q_n$ be the vertex reached from $a$ following a path
labeled by the binary representation $b(n)$ of $n$. It is easy
to verify by induction on the length of $b(n)$ that
\begin{enumerate}
\item One has $x_n=q_n$.
\item One has 
\begin{displaymath}
\varepsilon_n=\begin{cases}1&\mbox{ if $q_n\in\{a,b\}$}\\
-1&\mbox{ otherwise}
\end{cases}
\end{displaymath}
Thus we conclude that $\varepsilon_n=\phi(x_n)$.
\end{enumerate}
\end{solution}
\begin{solution}{\ref{exerciseChaconMinimal}}
  Since $000,010\in\cL(\sigma)$, the letters
  $0,1$ are extendable in $\cL(\sigma)$. Thus it follows from Exercise~\ref{exerciseDefSubstitution} that $\sigma$ is a substitution.
Set $v_n=\sigma^n(0)$. Then $v_{n+1}=v_nv_n1v_n$. Any $u\in\cL(X(\sigma))$ is
a factor of some $v_n$. Thus for every $n\ge 1$ there exists $N$
such that every word of $\cL_n(X)$ is a factor of $v_N$
and thus an integer $M$ such that every word in $\cL_n(X(\sigma))$
is a factor of every word in $\cL_M(X(\sigma))$.
\end{solution}
\begin{solution}{\ref{exerciseGrowingMinimal}}
  Since $\sigma$ is growing, some power $\sigma^n$
  of $\sigma$ is prolongable on some $a\in A$. Set $x=\sigma^n(a)$.
  Since $X=X(\sigma)$ is minimal, $X(\sigma)$ is generated by $x$.
  Since $\sigma$ is a substitution, $\cL(\sigma)=\cL(X)=\cL(x)$.
Thus every $b\in A$ appears in $x$.
Since $X$ is minimal, $x$ is uniformly recurrent.
Since $\sigma$ is growing, for every $b\in A$ there is some $n\ge 1$
such that $a$ appears in $\sigma^n(b)$. Since
$\sigma$ is prolongable on $a$, the letter $a$ also appears in all $\sigma^m(b)$
for $m\ge n$. Thus, there is an $N$ such that $a$ appears
in all $\sigma^N(b)$ for $b\in A$. By minimality again,
there is a $k\ge 1$ such that every $c\in A$ appears in $\sigma^k(a)$.
Then every $c$ appears in every $\sigma^{N+k}(b)$.
\end{solution}
\begin{solution}{\protect{\ref{exerciseElementarySubstitution}}}
Let $\varphi:A^*\to B^*$ be a nonerasing morphism which is not injective
as a map from $A^\N$ to $B^\N$. Set $U=\varphi(A)$.
Let $Y$ be the basis of the
intersection of all free submonoids containing $U^*$
and let $\beta:B\to Y$ be a bijection from an alphabet $B$
with $Y$.
Then there is a morphism $\alpha:A^*\to C^*$
such that $\varphi=\alpha\circ\beta$.
For a word $x\in U$, let $\lambda(x)\in Y$ be the first 
symbol in its decomposition in words of $Y$, that is
$x\in\lambda(x)Y^*$.
If some $y\in Y$ does  not appear as an initial symbol in the words
of $U$, set $Z=(Y\setminus y)y^*$. Then $Z^*$ is free
and $U^*\subset Z^*\subset Y^*$. Thus $Y=Z$, a contradiction.
Since $\varphi$ is not injective on $A^\N$, the map $\lambda$
is not injective and thus $\Card(Y)<\Card(U)$,
showing that $\varphi$ is not elementary.
\end{solution}
\begin{solution}{\protect{\ref{exerciseFixedPointPeriodic1}}}
Let $p(n)$ be the number of factors of length $n$ of $x$.
Let us show that if $p(n)<p(n+1)$, then there is an $m>n$ such
that $p(m)<p(m+1)$. Indeed if $p(n)<p(n+1)$ there is
a factor $u$ of length $n$ of $x$ which is right-special, that
is, there are two distinct letters
$a,b$ such that $ua,ub\in F(x)$. Since $\varphi$ is elementary,
it is
injective on $A^\N$ (Exercise~\ref{exerciseElementarySubstitution}).
Thus there are some $v,w$ such that $\varphi(av)\ne\varphi(bw)$.
 If $|\varphi(u)|>|u|$
or if $\varphi(a),\varphi(b)$ begin by the same letter,
the longest common prefix of $\varphi(uav),\varphi(ubw)$ is a right-special
word of length $m>n$. Otherwise we may replace $u$ by $\varphi(u)$.
Since $\varphi$ is primitive, the second case can happen only a bounded
number of times.
If $x$ is periodic, then $p(n)$ is bounded and by the preceding argument,
 no letter
can be right-special, which implies that the period of $x$
is at most $\Card(A)$.
\end{solution}

\begin{solution}{\protect{\ref{exerciseFixedPointPeriodic2}}}
We use an induction on $\Card(A)$. If $\Card(A)=1$, the result is true
since the period is $1$. Otherwise, either $\varphi$ is elementary
and thus injective by Exercise~\ref{exerciseElementarySubstitution}.
Then, by Exercise~\ref{exerciseFixedPointPeriodic1}, the period
of $x$ is at most equal to $\Card(A)$. Finally, assume that
$\varphi=\alpha\circ\beta$ with $\beta:A^*\to B^*$ and
$\alpha:B^*\to A^*$ and $\Card(B)<\Card(A)$. Set $y=\beta(x)$ and $\psi=\beta\circ\alpha$.
Then $\psi(y)=\beta\circ\alpha\circ\beta(x)=\beta(x)=y$. Thus
$y$ is a fixed point of $\psi$. Since $\psi$ is primitive, we may apply
the induction hypothesis. Since the period of $y$ is at most
$|\psi|^{\Card(B)-1}$, the period of $x$ is at most
$|\alpha|\,|\psi|^{\Card(B)-1}\le|\varphi|^{\Card(A)-1}$.
\end{solution}
\begin{solution}{\ref{exerciseStableSubmonoid}}
Assume first that $M=U^*$ where $U$ is a code. Then, the unique
factorisation of $u(vw)=(uv)w$ forces $v\in M$. Conversely,
let $U$ be the the set of nonempty words in $M$ which cannot be written
as a product of strictly shorter words in $M$.
Take $x_1x_2\cdots x_n=y_1y_2\cdots y_m$ with $x_i,y_j\in U$,
and $n,m\ge 0$ minimal. By definition of $U$, we have $n,m\ge 1$.
Assume $|x_1|\le |y_1|$ and
set $y_1=x_1v$. Then, $x_2\cdots x_n=vy_1\cdots y_m$. By \eqref{eqStability},
we have $v\in M$, which forces $v=\varepsilon$ and $x_1=y_1$,
a contradiction with the minimality of $n,m$.

\end{solution}
\begin{solution}{\ref{exerciseFlowerAutomaton}}
 The first assertion is clear since, by construction of $\A(U)$,
the nonempty simple paths from $\omega$ to $\omega$
are in bijection with $U$.

 1. Suppose first that $U$ is a code. By definition there
can be at most one path from $\omega$ to $\omega$ labeled
$w$ since otherwise $w$ would have two factorisations in words of $U$.
Next, if there are two distinct paths labeled $w$ from $p$ to $q$,
there would be two distinct paths labeled $uwv$ from
$\omega$ to $\omega$ for words $u,v$ labeling paths
from $\omega$ to $p$ and $q$ to $\omega$ respectively,
a contradiction. The converse is obvious.

2. The number of paths labeled $w$ from $(u,v)$ to itself is the 
number of factorisations $w=vxu$ with $x\in U^*$.
Assume that $U$ is circular. There cannot be cycles labeled $w$ from $\omega$ to
$\omega$ and also from $(u,v)$ to $(u,v)$ by definition of a circular
code. Next, suppose that there are cycles labeled $w$ around $(u,v)$ 
and around $(u',v')$. Set $w=vxu=v'x'u'$. Assuming $|v|<|v'|$,
set $v'=vz$ and $t=x'u'$. Then $xu=zt$ implies $ztv=xuv$, which
shows that $ztv\in U^*$. But $vzt,ztv\in U^*$ implies
$zt,v\in U^*$ by \eqref{eqVeryPure}. Thus $w=vxu=vzt$
is in $U^*$, a contradiction.
The converse is obvious.
\end{solution}
\begin{solution}{\ref{exerciseRankOne}}
(i) $\Rightarrow$ (ii) the first
assertion results from the definition of a synchronizing pair.
The second one is clear since $\mu(xy)$ is the
product of the column of index $\omega$ of $\mu(x)$
by the row of index $\omega$
of $\mu(y)$.
 
If $\mu(x),\mu(y)$ have rank one, consider $u,v\in A^*$ such that
$uxyv\in U^*$. Since $\mu(x)$ has rank one, it follows from $uxyv,x\in U^*$
that $ux,xyv\in U^*$. Similarly, $uxy,yv\in U^*$. This shows that
$(x,y)$ is synchronizing.
\end{solution}
\begin{solution}{\ref{exerciseUniformSync}}
(i) $\Rightarrow$ (ii).
Let $\A(U)=(Q,E)$ be the flower automaton of $U$ and let $\mu$ be as in Exercise
\ref{exerciseRankOne}.
Let $S$ be the finite semigroup $\mu(A^+)$. For $m\in S$e write
$p\edge{m}q$ for $m_{p,q}=1$. Since $S$ is finite, there
is for every $m\in S$ an integer $k\ge 1$ such that $m^k$ is idempotent,
that is $m^{2k}=m^k$. Let $e$ be such an idempotent. For every $p,q\in Q$
there is, since $e=e^3$  some $r,s\in Q$ such that $p\edge{e}r\edge{e}s\edge{e}q$. By unambiguity, $r=s$, that is $r$ is a fixed point of $m$.
But an element of $S$ cannot have more than one fixed point by Exercise~\ref{exerciseFlowerAutomaton}
and thus $e$ has rank one.

Let $J$ be the set of elements of $S$ of rank one. Since all idempotents
of $S$ are in $J$, we have $S^n\subset J$ for $n=\Card(S)+1$. 
Indeed, if $m\in S^n$ with $n\ge \Card(S)+1$, there
is a factorisation $m=uvw$ with $uv=u$. Then $uv^k=u$ for all $k\ge 1$.
When $v^k$ is idempotent, we have $v^k\in J$ and thus $m\in J$.
This
shows that for every word $x$ of length $n$, $\mu(x)$ has rank one.
We conclude using Exercise~\ref{exerciseRankOne} that $U$
has finite synchronization delay.

(ii) $\Rightarrow$ (iii) is clear by definition of a synchronizing pair.

(iii) $\Rightarrow$ (i). Suppose that $p,q\in A^*$ are such 
that $pq,qp\in U^*$. Then the sequence $\cdots qp\cdot pqpq\cdots$
has two factorizations unless $p,q\in U^*$. Thus $U$ is circular.
\end{solution}

\begin{solution}{\ref{exerciseFiniteToOne}}
  As seen in Section~\ref{sectionRecognizable}, we may consider $X^\varphi$ as a shift on the alphabet
  $A^\varphi=\{(a,i)\mid 0\le i<|\varphi(a)|\}$. Assume first that
  $\varphi$ is injective. For every $y\in\hat{\varphi}(X^\varphi)$, the set
  $\hat{\varphi}^{-1}(y)$ has at most $(\Card(A^\varphi))^2$ elements.
  Indeed, assume that $\hat{\varphi}(x)=y$. For $n\ge 1$, set $x_{-n}=(a,i)$ and $x_n=(b,j)$.
  Then, there are unique $k\ge 0$ and $a_1,\ldots ,a_k\in A$ such that
  \begin{displaymath}
\varphi(aa_1\ldots a_kb)=y_{-n-i}\ldots y_{n+|\varphi(b)|-j}.
    \end{displaymath}
  Thus, the symbols $x_{-n}$ and $x_n$ and the sequence $y$ determine
  the symbols $x_i$ for $-n\le i\le n$. Conversely, if there
  are distinct $w,w'$ such that $\varphi(w)=\varphi(w')$, then
  the set $\hat{\varphi}^{-1}(\varphi(w)^\infty)$ is infinite.
  
  \end{solution}
\begin{solution}{\ref{exerciseMosseOriginal}}
  Assume first that $\varphi$ is recognizable
  and suppose that $x_{[i-\ell,i+\ell)}=x_{[j-\ell,j+\ell)}$.
      Set $p=x_{[i-\ell,i)}$
     and $q=x_{[i,i+\ell)}$. If $i\in f(\N)$,
   the pair $(p,q)$   is parsable.
       For $\ell$ large enough, it is synchronizing.
       Let $r<s$ be such that $f(r)\le j-\ell<j+\ell\le f(s)$.
       Then $pq$ is a factor of $\varphi(x_{[r,s)})$. Since
      $(p,q)$ is synchronizing, we have $j\in f(\N)$
      and $x_{f^{-1}(i)}=x_{f^{-1}(j)}$.

            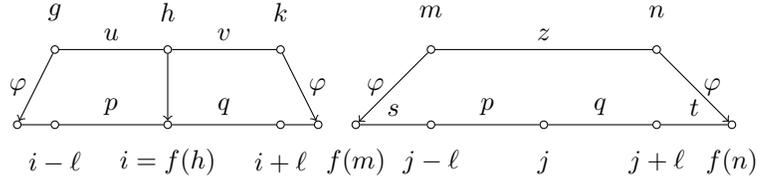
\begin{figure}[hbt]
\centering
\tikzset{node/.style={circle,draw,minimum size=0.1cm,inner sep=0pt}}
\tikzset{title/.style={minimum size=0.1cm,inner sep=0pt}}

\begin{tikzpicture}
  \node[node](u)at(0,1){};\node[node](uv)at(1.5,1){};\node[node](v)at(3,1){};
  \node[title](g)at(0,1.5){$g$};\node[title](h)at(1.5,1.5){$h$};
  \node[title](k)at(3,1.5){$k$};
  \node[node](fg)at(-.5,0){};\node[node](i-l)at(0,0){};
  \node[node](i)at(1.5,0){};\node[node](i+l)at(3,0){};\node[node](fk)at(3.5,0){};
  \node[title]at(0,-.5){$i-\ell$};\node[title]at(1.5,-.5){$i=f(h)$};
  \node[title]at(3,-.5){$i+\ell$};

  \draw[above](u)edge node{$u$}(uv);\draw[above](uv)edge node{$v$}(v);
  \draw[left,->](u)edge node{$\varphi$}(fg);
  \draw[left,->](uv)edge node{}(i);
  \draw[right,->](v)edge node{$\varphi$}(fk);
  \draw[below](fg)edge node{}(i-l);\draw[above](i-l)edge node{$p$}(i);
  \draw[above](i)edge node{$q$}(i+l);\draw[below](i+l)edge node{}(fk);

  \node[node](m)at(5,1){};\node[node](n)at(8,1){};
  \node[title]at(5,1.5){$m$};\node[title]at(8,1.5){$n$};
  \node[node](fm)at(4,0){};\node[node](j-l)at(5,0){};
  \node[node](j)at(6.5,0){};\node[node](j+l)at(8,0){};
  \node[node](fn)at(9,0){};
  \node[title]at(4,-.5){$f(m)$};\node[title]at(5,-.5){$j-\ell$};
  \node[title]at(6.5,-.5){$j$};\node[title]at(8,-.5){$j+\ell$};
  \node[title]at(9,-.5){$f(n)$};
  \draw[above](m)edge node{$z$}(n);
  \draw[left,->](m)edge node{$\varphi$}(fm);
  \draw[right,->](n)edge node{$\varphi$}(fn);
  \draw[above](fm)edge node{$s$}(j-l);
  \draw[above](j-l)edge node{$p$}(j);\draw[above](j)edge node{$q$}(j+l);
  \draw[above](j+l)edge node{$t$}(fn);
\end{tikzpicture}
\caption{Alternative definition of recognizability}\label{figureAltDefRecognizability}
       \end{figure}
                  Conversely, let $p,q\in\cL_\ell(X)$ be such that the pair $(p,q)$ is $(u,v)$-parsable with $uv\in\cL(X)$.
                  Let $g<h<k$ be such that $u=x_{[g,h)}$ and $v=x_{[h,k)}$
                      (see Figure~\ref{figureAltDefRecognizability}
                      on the left). We have $pq=x_{[i-\ell,i+\ell)}$
                          and $i=f(h)$.
          Let $z\in\cL(X)$ be such that $pq$ is a factor of $\varphi(z)$.
          Set $\varphi(z)=spqt$ and $z=x_{[m,n)}$
            (see Figure~\ref{figureAltDefRecognizability} on the right). Then
            $pq=x_{[j-\ell,j+\ell)}$ with $j=f(m)+|s|+\ell$.
              By the hypothesis, we have $j\in f(\N)$ and $x_h=x_{f^{-1}(j)}$,
              which implies that $(p,q)$ is synchronizing.
\end{solution}
\exosection{Section \ref{sectionSturmianShifts}}
\begin{solution}{\protect{\ref{exerciseFiboisSturmian}}}
The words $F_n=\varphi^n(a)$ are left-special. Indeed, this is
true for $n=0$ since $aa,ba\in\cL(\varphi)$ and 
\begin{displaymath}
aF_{n+1}=\varphi(bF_n),\quad abF_{n+1}=\varphi(aF_n)
\end{displaymath}
shows the claim by induction on $n$. It is easy to see
(again by induction) that conversely every left-special
word is a prefix of some $F_n$. This implies that there
is exactly one left-special word of each length which is moreover
extendable by every letter and is additionally a prefix of $x$
and thus the conclusion that $x$ is standard episturmian.
\end{solution}
\begin{solution}{\ref{exerciseBalanced1}}
  Set $p_n(U)=\Card(U\cap \{0,1\}^n)$.
  Assume by contradiction that $p_n(U)\ge n+2$ for some $n\ge 1$.
  If $n$ is chosen minimal, we have $p_{n-1}(U)\le n$. Thus
  there exist two distinct left-special words $u,v$ of length $n-1$. Since
  $u\ne v$, there is a word $w$ such that $w0$ is a prefix
  of, say $u$, and $w1$ a prefix of $v$. Then
  $0w0,1w1$ are in $U$, a contradiction.
\end{solution}
\begin{solution}{\ref{exerciseBalanced2}}
  Let $u,v$ be words of minimal length  such that $\delta(u,v)\ge 2$.
  The first letters of $u,v$ are distinct. Assuming that $u$ starts with $0$,
  we have factorizations $u=0wau'$ and $v=1wbv'$ where $a,b$ are distinct
  letters. If $a=1$ and $b=0$, then $\delta(u',v')=\delta(u,v)$,
  a contradiction with the minimality of $n$. Thus $a=0$ and $b=1$.
  By minimality again, we found $u=0w0$ and $v=1w1$. If $w$
  is not palindrome, there is a word $z$ such that $za$ is a prefix of $w$
  and $b\tilde{z}$ is a suffix of $w$ with $a,b$ distinct letters.
  If $a=0$, then $0z0$ and $1\tilde{z}1$ are in $\cL(s)$, a contradiction.
  Thus $a=1$ and $b=0$. Set $u=0z1u''$ and $v=v''1\tilde{z}0$.
  Then $\delta(u'',v'')=\delta(u,v)$, a contradiction again.
  We conclude that $w$ is a palindrome.
  \end{solution}
\begin{solution}{\ref{exerciseEquivalentDefSturm}}
  Let $p_n(x)=\Card(\cL_n(x))$ be the complexity of $x$.
  If $x$ is aperiodic, then $p_n(x)\ge n+1$ by the variant
  of Theorem~\ref{theoremCovenHedlundSequence} for one-sided sequences.
  Since $x$ is balanced, $p_n(x)\le n+1$ by Exercise~\ref{exerciseBalanced1}.
  Thus $p_n(x)=n+1$ showing that $x$ is Sturmian.

  Conversely, assume that $x$ is Sturmian and unbalanced. We
  prove that $x$ is eventually periodic. By
  Exercise \ref{exerciseBalanced2}, there is a palindrome $w$
  such that $0w0,1w1$ are factors of $x$. Thus $w$ is right-special.
  Set $n=|w|+1$. Since $x$ is Sturmian there is a unique right-special
  factor of length $n+1$, which must be either $0w$ or $1w$.
  Assume that $0w$ is right-special. Let $v$ be a word of length $n-1$
  such that $u=1w1v$ is a factor of $x$. In view of using Exercise \ref{exerciseProposition1.3.14}, let us show that all factors of length $n$
  of $u$ are conservative.
  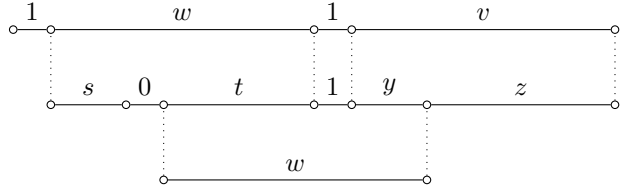
\begin{figure}[hbt]
    \centering
\tikzset{node/.style={circle,draw,minimum size=0.1cm,inner sep=0pt}}
\tikzset{title/.style={minimum size=0.1cm,inner sep=0pt}}
    \begin{tikzpicture}
      \node[node](u0)at(0,1){};\node[node](u1)at(.5,1){};
      \node[node](u2)at(4,1){};\node[node](u3)at(4.5,1){};
      \node[node](u4)at(8,1){};

      \node[node](1)at(.5,0){};\node[node](2)at(1.5,0){};\node[node](3)at(2,0){};
      \node[node](4)at(4,0){};\node[node](5)at(4.5,0){};
      \node[node](6)at(5.5,0){};\node[node](7)at(8,0){};

      \node[node](w3)at(2,-1){};\node[node](w6)at(5.5,-1){};

      \draw[above](u0)edge node{$1$}(u1);\draw[above](u1)edge node{$w$}(u2);
      \draw[above](u2)edge node{$1$}(u3);\draw[above](u3)edge node{$v$}(u4);

      \draw[above](1)edge node{$s$}(2);\draw[above](2)edge node{$0$}(3);
      \draw[above](3)edge node{$t$}(4);\draw[above](4)edge node{$1$}(5);
      \draw[above](5)edge node{$y$}(6);\draw[above](6)edge node{$z$}(7);

      \draw[above](w3)edge node{$w$}(w6);

      \draw[dotted](1)edge node{}(u1);\draw[dotted](4)edge node{}(u2);
      \draw[dotted](5)edge node{}(u3);\draw[dotted](7)edge node{}(u4);
      \draw[dotted](w3)edge node{}(3);\draw[dotted](w6)edge node{}(6);
    \end{tikzpicture}
    \caption{The words $u,v$ and $w$}\label{figureEquivalentDefSturm}
  \end{figure}
  It is enough to show that $0w$ is not
  a factor of $u$. Assume the contrary. Then we have $w=s0t$,
  $v=yz$ and $w=t1y$ (see Figure~\ref{figureEquivalentDefSturm}). Since $w$ is a palindrome, the first equality
  implies $w=\tilde{t}0\tilde{s}$. 
  The words $t$ and $\tilde{t}$ having the same length, this implies $t1 = \tilde{t}0$, a contradiction.
  We conclude that $u$ is a factor of $x$ of length $2n$ in
  which all factors of length $n$ are conservative. Since $x$
  is Sturmian, there are $n$ conservative factors of length $n$.
  Thus $x$ is eventually periodic by Exercise~\ref{exerciseProposition1.3.14}.
\end{solution}
\begin{solution}{\ref{exerciseSturmReversal}}
  Let $x$ be a Sturmian sequence. Set $U=\cL(x)\cup\widetilde{\cL(x)}$.
  Since $U$ is balanced, we have $\Card(U\cap \{0,1\}^n)\le n+1$
  by Exercise~\ref{exerciseBalanced1}.
  Since $\Card(\cL(x)\cap \{0,1\}^n)=n+1$, this forces
  $\cL(x)=\widetilde{\cL(x)}$
  \end{solution}
\begin{solution}{\ref{exerciseTriboisSturmian}}
  Let $X$ be the shift generated by $\varphi$. Let us first
  verify that $\cL(X)$ is closed under reversal. Denote
  $\rho(w)=\tilde{w}$ the reversal of $w$.  Let $\psi:A^*\to A^*$
  be defined by $\psi(w)=\varphi(w)a$. Observe
  that the map $\tilde{\varphi}$ defined by $\tilde{\varphi}(w)=a^{-1}\varphi(w)a$
  is a morphism and that $\rho(\varphi(w))= \widetilde{\varphi}(\tilde{w})$.
  This implies that $\psi^n(a)$ is a palindrome for every $n\ge 1$ because
  \begin{eqnarray*}
    \rho(\psi^{n+1}(a))&=&\rho(\varphi(\psi^n(a)a)=a\tilde{\varphi}(\rho(\psi^n(a)))\\
    &=&a\tilde{\varphi}(\psi^n(a))=\varphi(\psi^n(a))a=\psi^{n+1}(a).
  \end{eqnarray*}
  Next, we claim that the prefixes of the Tribonacci word are the left
  special words of $\cL(X)$ and that they have three extensions.
  Indeed, let us verify by induction on $n$ that
  the words $T_n=\varphi^n(a)$ are left-special with $3$ extensions. 
This is
true for $n=0$ since $aa,ba,ca\in\cL(X)$. Next the equalities
\begin{displaymath}
aT_{n+1}=\varphi(cT_n),\quad abT_{n+1}=\varphi(aT_n),\quad acT_{n+1}=\varphi(bT_n)
\end{displaymath}
prove the property by induction. Conversely, any left-special word is a prefix
of some $T_n$, whence the assertion. 
\end{solution}
\begin{solution}{\ref{exerciseMechanicalIsSturm}}
  By Proposition~\ref{propositionEquivalentDefSturm}, it is equivalent to
  prove that $s$ is balanced and aperiodic. Let
  $u=s_{[n,n+p)}$ and set $h(u)=|u|_1$. Then
    $h(u)=\lfloor\alpha(n+p)+\rho\rfloor-\lfloor \alpha n+\rho\rfloor$. Thus
    \begin{equation}
      \alpha|u|< h(u)<\alpha|u|+1.\label{ineqHeight}
    \end{equation}
    This implies that $\lfloor \alpha|u|\rfloor \le h(u)\le\lfloor \alpha|u|\rfloor+1$
    and shows that $h(u)$ take only two consecutive values when $u$
    range over the factors of fixed length of $s$. Thus $s$ is balanced.
    Moreover, by \eqref{ineqHeight}, the real number $\pi(u)=h(u)/|u|$ satisfies
    \begin{displaymath}
      |\pi(u)-\alpha|<\frac{1}{|u|}.
    \end{displaymath}
    This proves  that $\pi(u)$ tends to $\alpha$ when $|u|\to\infty$
    (and thus the last statement) and that $s$ is aperiodic, since otherwise the limit of
    $\pi(u)$ would be rational.
  \end{solution}
\begin{solution}{\ref{exerciseJustin0}}
  Any word $uv^{-1}\tilde{u}$ is palindrome and has $u$ as a prefix.
  Let us show that the palindromic closure is of this form.
  Set $u^{(+)}=ur=\tilde{r}\tilde{u}$. We have $|r|<|u|$ since we can
  choose the last letter of $u$ as palindrome suffix. Set $u=\tilde{r}s$.
  Then $ur=\tilde{r}\tilde{u}$ implies $s=\tilde{s}$. This
  proves that $u^{(+)}$ is of the form above and thus corresponds
  to the longest possible $v$.
  \end{solution}
\begin{solution}{\ref{exerciseJustin1}}
  The formula is clear when $|w|_a=0$. Otherwise, 
  the longest palindromic suffix of $\Pal(w)a$ is $a\Pal(w_1)a$.
  By Exercise~\ref{exerciseJustin0}, we have
  \begin{eqnarray*}
    \Pal(wa)&=&(\Pal(w)a)^{(+)}=\Pal(w)a(a\Pal(w_1)a)^{(-1)}a\Pal(w)\\
    &=&\Pal(w)\Pal(w_1)^{-1}\Pal(w).
    \end{eqnarray*}
\end{solution}
\begin{solution}{\ref{exerciseJustin2}}
  This will follow from the fact that for every word $u$ and every letter $a$,
  one has the formula $L_a(\tilde{u})a=\widetilde{L_a(u)a}$ easily
  established by induction on the length of $u$.
\end{solution}
\begin{solution}{\ref{exerciseJustin3}}
  We use induction on the length of $w$. 
  The formula is true if $w$ is empty. Otherwise, set $w=vb$.
  Assume first that $|v|_b=0$. Then $\Pal(w)=\Pal(v)b\Pal(v)$. Next,
  using the induction hypothesis,
  \begin{eqnarray*}
    L_a(\Pal(w))a&=&L_a(\Pal(v))b\Pal(v)a=L_a(\Pal(v))aa^{-1}L_a(b)L_a(\Pal(v))a\\
    &=&\Pal(av)a^{-1}L_a(b)\Pal(av)\\
    &=&\begin{cases}\Pal(av)b\Pal(av)=\Pal(avb)=\Pal(aw)&\mbox{ if $a\ne b$}\\
      \Pal(av)a\Pal(av)=\Pal(ava)=\Pal(aw)&\mbox{ otherwise}
      \end{cases}
  \end{eqnarray*}
  Now, if $|v|_b>0$, set $v=v_bv_2$ with $|v_2|_b=0$. By Exercise \ref{exerciseJustin1}, we have $\Pal(w)=\Pal(v)\Pal(v_1)^{-1}\Pal(v)$. Thus
  \begin{eqnarray*}
    L_a(\Pal(w))a&=&L_a(\Pal(v))aa^{-1}L_a(\Pal(v_1))^{-1}L_a(\Pal(v))a\\
    &=&\Pal(av)\Pal(av_1)^{-1}\Pal(av)\\
    &=&\Pal(aw).
    \end{eqnarray*}
  \end{solution}
\begin{solution}{\ref{exerciseJustin4}}
  We use induction on the length of $v$. The formula is true for $|v|=0$.
  Otherwise, set $v=av_1$. 
  Using Exercise~\ref{exerciseJustin3} and the induction hypothesis, we have
  \begin{eqnarray*}
    \Pal(vw)&=&\Pal(av_1w)=L_a(\Pal(v_1w))a=L_a(L_{v_1}(\Pal(w))\Pal(v_1))a\\
    &=&L_{av_1}(\Pal(w))L_a(\Pal(v_1))a=L_v(\Pal(w))\Pal(v).
  \end{eqnarray*}
  \end{solution}
\begin{solution}{\ref{exerciseReturnSturm}}
  We first observe that for every $a\in A$, there is at most
  one $w\in \RR'_X(u_n)$ which ends with $a$. Indeed, assume
  that $w,w'$ are distinct words in $\RR_X'(u_n)$ ending
  with the same letter. Set
  $w=w_1w_2$ and $w'=w'_1w_2$ with $w'_1,w_1$ ending with different
  letters. Then $w_2u_n$ is left-special and thus begins with $u_n$,
  a contradiction.

  Consider now $w\in\RR'_X(u_n)$. Let $a$ be the last letter of $w$.
  Since $u_n$ can be followed by $a$, there  is an $m\ge 0$
such that $x_{n+m}=a$. Then $u_{n+m+1}=(u_{n+m}a)^{(+)}$ has a factor
$(u_na)^{(+)}=\Pal(x_0x_1\cdots x_{n-1}a)$.
By Justin's Formula~\eqref{eqJustin0}, we have for every $a\in A$
\begin{displaymath}
L_{x_0x_1\cdots x_{n-1}}(a)u_n=\Pal(x_0x_1\cdots x_{n-1}a).
\end{displaymath}
Since $\Pal(x_0x_1\cdots x_{n-1}a)=(u_na)^{(+)}$ is the shortest palindrome with
prefix $u_na$, it begins and ends with $u_n$ with no other
occurrence of $u_n$. Thus the left return word to $u_n$
ending with $a$ is $L_{x_0x_1\cdots x_{n-1}}(a)$. This proves the
inclusion $\RR'_X(u_n)\subset\{L_{x_0x_1\cdots x_n}(a)\mid a\in A\}$.

Assume now that $s$ is strict and consider $a\in A$.
Since
every letter appears infinitely often in $x$, there is an $m\ge 0$
such that $x_{n+m}=a$. Then $u_{n+m+1}=(u_{n+m}a)^{(+)}$ has a factor
$(u_na)^{(+)}=\Pal(x_0x_1\cdots x_{n-1}a)$. This word
is clearly a left return word to $u_n$. This shows that
$L_{x_0x_1\cdots x_{n-1}}(a)$ is in $\RR'_X(u_n)$.

The statement concerning $\RR'_X(w)$ is clear since
$u_n$ is the shortest left special word having $w$ as a suffix.
\end{solution}
\begin{solution}{\ref{exerciseEquivalentDefAR}}
  The condition is necessary since a strict episturmian shift is 
   closed under reversal. Conversely, let $X$ be a minimal shift
  space satisfying the condition. Since $X$ is minimal, every $u\in\cL(X)$
  is a prefix of a right-special word. This implies that every word
  in $\cL(X)$ is a factor of some bispecial word. But the condition
  implies that bispecial words are palindrome. Indeed, if $w$ is bispecial,
  then $\tilde{w}$ is also bispecial and thus $w=\tilde{w}$. This
  shows that $\cL(X)$ is closed by reversal.
  \end{solution}
\begin{solution}{\ref{exerciselemma6}}
 
  \begin{enumerate}
  \item Set $u_{n+1}=yuy'$ with $|y|$ minimal. Since $u_{n+1}=u_nx_n\tilde{t}$,
    and since $|yu|\ge |u_nx_n|$, we have $|y'|\le|t|$. Similarly,
    since $\tilde{u}$ is not a factor of $u_n$, we have $|y'u|\ge|u_nx_n|$
    and thus
    $|y|\le|t|$. Set $t=yz$ and $\tilde{t}=z'y'$. Then
    \begin{displaymath}
      u_{n+1}=yuy'=yzwz'y'.
      \end{displaymath}
  \item We have $u=avb=zwz'$. If $z,z'$ are nonempty, then $w$ is a factor
    of $v$. Since $v$ is a palindrome prefix of $s$ shorter than $u_{n+1}$,
    it is a prefix of $u_n$ and thus $w$ is a factor of $u_n$. This contradicts the fact that $w$ is
    unioccurrent in $u_nx_n$
    \end{enumerate}
  \end{solution}
\begin{solution}{\ref{exerciseStandard}}
We first prove that the words $s_n$ are primitive.
For two words $x,y$ on the alphabet $\{0,1\}$, denote
\begin{displaymath}
M(x,y)=\begin{bmatrix}|x|_0&|x|_1\\|y|_0& |y|_1\end{bmatrix}
\end{displaymath}
We prove by induction on $n\ge 1$ that
$\det M(s_n,s_{n-1})=1$. This implies that 
$|s_n|_0,|s_n|_1$ are relatively prime and thus every $s_n$ is primitive.
The equality is true for $n=1$. Next, for $i\ge 1$, we have
\begin{displaymath}
M(s_n,s_{n-1})=\begin{bmatrix}d_n&1\\1&0\end{bmatrix}M(s_{n-1},s_{n-2})
\end{displaymath}
whence the conclusion. 

Let us now show that that every $s_n$ is a prefix
of $c_\alpha$. For this, let 
\begin{displaymath}
h_n=\begin{cases}L_{0^{d_1}1^{d_2}\cdots1^{d_n}}&\mbox{ if $n$ is even}\\
L_{0^{d_1}1^{d_2}\cdots0^{d_n}}&\mbox{ if $n$ is odd}.
\end{cases}
\end{displaymath}
One can easily verify by induction on $n\ge 1$ that
\begin{eqnarray*}
  h_{2n-1}(1)&=&s_{2n-1}=h_{2n}(1)\\
  h_{2n}(0)&=&s_{2n}=h_{2n+1}(0)
  \end{eqnarray*}
By Justin Formula, this implies that $s_n$
is a prefix of $c_\alpha$.
\end{solution}
\begin{solution}{\ref{cfrac1}}
We show that for $n\ge 3$, one has
\begin{displaymath}
s_{n-1}s_n=s_nt_{n-1}\mbox{ with } t_n=s_{n-1}^{d_{n-1}}s_{n-2}s_{n-1}.
\end{displaymath}
Indeed,
\begin{eqnarray*}
s_{n-1}s_n&=&s_{n-1}s_{n-1}^{d_{n}}s_{n-2}=s_{n-1}^{d_{n}}s_{n-2}^{d_{n-1}}s_{n-3}s_{n-2}\\
&=&s_{n-1}^{d_{n}}s_{n-2}s_{n-2}^{d_{n-1}-1}s_{n-3}s_{n-2}=s_nt_{n-1}.
\end{eqnarray*}
Observe that $t_{n-1}$ is not a prefix of $s_n$ since otherwise
$s_n=t_{n-1}u$ for some word $u$ and $s_{n-1}s_nu=s_n^2$,
a contradiction since $s_n$ is primitive by Exercise~\ref{exerciseStandard}.
Clearly $s_{n+1}s_n$ is a prefix of the characteristic word $c_\alpha$. Since 
\begin{displaymath}
s_{n+1}s_n=s_n^{d_{n+1}}s_{n-1}s_n=s_n^{d_{n+1}}s_nt_{n-1}=s_n^{d_{n+1}+1}t_{n-1}
\end{displaymath}
the word $s_n^{d_{n+1}+1}$ is a prefix of $c_\alpha$ and since
$t_{n-1}$ is not a prefix of $s_n$, the word $s_n^{d_{n+1}+2}$ is
not a prefix of $c_\alpha$.
\end{solution}
\begin{solution}{\ref{exercisecfrac2}}
It suffices to consider the characteristic Sturmian sequence
$c_\alpha$ with $\alpha=[0,1+d_1,d_2,\ldots]$.
If the sequence of $a_i$ is unbounded, then 
$s_n^{d_{n+1}}$ is a prefix of $c_\alpha$ and consequently $c_\alpha$
is not $d$-power free for any $d$.
\end{solution}
\begin{solution}{\ref{exerciseSturmLR}}
Let $x$ be a two-sided Surmian sequence of slope $\alpha=[a_0,a_1,\ldots]$.
Assume that $x$ is linearly recurrent with constant $K$.
Since $x$
is Sturmian it is not periodic. Thus, by
 Proposition~\ref{propositionLRhasLinearComplexity}, it is
$(K+1)$-power free. By Exercise \ref{exercisecfrac2}, the coefficients
$a_i$ are bounded.
\end{solution}
\exosection{Section~\ref{sectionToeplitz}}
\begin{solution}{\ref{exerciseDoublingSequence}}
The formula is true for $n=0$. Next, $\nu_2(2n+2)\equiv\nu_2(n+1)+1\bmod 2$
and $\nu_2(2n+1)=0$ for all $n\ge 0$.
\end{solution}

\begin{solution}{\ref{exerciseShiftToeplitz}}
  We treat the case of the one-sided shift $X$ generated by $\sigma$.
  The case of the two-sided shift is similar.
Observe that the powers of $\sigma$ also have coincidences.
Thus we can suppose $\sigma $ has an admissible one-sided fixed point $x$.  

There is a letter $a$ such that for all letters $b$ one has $\sigma (b) = u(b) a v(b)$ where $b\mapsto |u(b)|$ is constant equal to $k$. 
Let $k(N) = \sum_{j=0}^{N} k n^{j}$. 
The sequence $x$ being a concatenation of images of $\sigma$ one has $x_k = x_{k+in}$ for all $i$. 
As a consequence, for all $N$, the word $\sigma^N (a)$ appears in $x$ at the occurrence $k(N)+n^{N+1}i$ for all $i\ge 0$.
One can check that $\sigma^N (a)$ is a prefix of $S^{k(N)} x$. 
Let $y$ be an accumulation point of $(S^{k(N)} x)_N$. 
From the observation above it is a Toeplitz sequence
and thus $X$ is a Toeplitz shift.
\end{solution}
%%%%%%%%%%%%%%%%%%%%%%
\section{Notes}

Topological dynamical systems are usually presented together
with measure-theoretic ones (in which there is a measure
on the space and the transformation preserves the measure,
see Chapter~\ref{chapterOrderedCohomology}).
 For a more
detailed introduction to topological dynamical systems,
see for example~\cite{Brown1976} or \cite{Petersen1983}.
\index{names}{Brown, James R.}\index{names}{Petersen, Karl}%
\subsection{Topological dynamical systems}
Cantor spaces (Section~\ref{sectionCantorSpaces})
 are a classical object in topology.
See~\cite{Willard2004}\index{names}{Willard, Stephen} for more details and, in particular
\cite[Theorem 30.3]{Willard2004} for a proof that all Cantor spaces are homeomorphic.

The mathematical definition of a decidable property,
using Turing machines and due to Turing, can be found in any textbook
on Theory of Computation (see \citep{AhoHopcroftUllman1974} for example).
Its adequacy to describe effectively computable properties
is known as the \emph{Church-Turing thesis}.
\index{names}{Aho, Alfred V.}\index{names}{Hopcroft, Jeffrey D.}\index{names}{Ullman, Jeffrey D.}%
\index{names}{Turing, Alan}\index{names}{Church, Alonzo}%
It is worth noting that a famous undecidable property is closely related to
the subject of this book. Indeed, given two morphism $\alpha,\beta:A^*\to B^*$,
the existence of a word $w$ such that $\alpha(w)=\beta(w)$ is an
undecidable property (called \emph{Post correspondance problem}, see \citep{AhoHopcroftUllman1974}).
\index{subject}{Post correspondance problem}\index{names}{Post, Emil}

The minimality
of irrational rotations (Exercise~\ref{exerciseIrrationalRotation})
 is known as the
(one-dimensional) Kronecker Theorem (see~\cite{HardyWright2008}).
\index{subject}{Kronecker Thorem}\index{subject}{Theorem!Kronecker}%
\index{names}{Kronecker, Leopold}%
\index{names}{Hardy, Godfrey}\index{names}{Wright, Edward M.}%

\subsection{Shift spaces}
Shift spaces are the basic object of symbolic dynamics. The
classical reference on symbolic dynamics is~\cite{LindMarcus1995}.
\index{names}{Lind, Douglas}\index{names}{Marcus, Brian H.}%
Many classes of shift spaces (such as shifts of finite type)
are described there in much more detail than we do here.
Shift spaces have also been called \emph{symbolic flows}
\index{subject}{symbolic!flow} (see \citep{Furstenberg1981} for example).
\index{names}{Furstenberg, Harry}

Section \ref{sectionCombinatoricsonWords} is a brief introduction
to combinatorics on words. For a more detailed exposition,
see \cite{Lothaire1983}. The combinatorial properties
of words are a source for many interesting algorithmic
problems (see the recent \cite{CrochemoreLecroqRytter2021} for an appealing presentation
of these problems).\index{names}{Crochemore, Maxime}\index{names}{Lecroq, Thierry}\index{names}{Rytter, Wojciech}
The original reference to the Curtis-Hedlund-Lyndon Theorem
is \cite{Hedlund1969}.

The classical
reference for the Morse Hedlund Theorem (Theorem~\ref{theoremCovenHedlund})
is~\cite{MorseHedlund1938}.
\index{names}{Morse, Marston}%
 The case of one-sided shifts
is considered in \citep{CovenHedlund1973}.
\index{names}{Coven, Ethan M.}\index{names}{Hedlund, Gustav A.}%
See also ~\cite[Theorem 1.3.13]{Lothaire2002} for example.
Proposition \ref{propositionCassaigne} is from \cite{Cassaigne1997}
 (see also~\cite[Theorem 4.5.4]{BertheRigo2010}).

The Champernowne sequence (Example~\ref{exampleChampernowne})
is from \cite{Champernowne1933}. It is usually
referred to as the \emph{Champernowne constant},\index{subject}{Champernowne!constant}
 which
is the real number $0.123456789101112\ldots$.

The bound $K\le 15$ given in Exercise~\ref{exerciseFibonacciLR} for the 
 linear
recurrence of the Fibonacci word is far from optimal. It
is shown in \cite{ChenFeiDuMousaviSchaefferShallit2014},
using the software Walnut \citep{mousavi2016automatic}
based on decidable properties of substitutive sequences,
that $K\le 3$. Actually, one has
more precisely $K\le (3+\sqrt{5})/2$,
this bound being the best possible~\citep{Shallit2020}.
\index{names}{Du, Chen Fei}\index{names}{Moussavi, Hamoon}%
\index{names}{Schaeffer, Luke}\index{names}{Shallit, Jeffrey O.}%

Linearly recurrent sequences form an important class of sequences
introduced in~\cite{Durand&Host&Skau:1999}.
We shall study these sequences in more detail in Chapter~\ref{chapterSubstitutionShifts}.

\subsection{Shifts of finite type}

Our brief introduction to shifts of finite type
(also called \emph{topological Markov chains})
\index{subject}{topological!Markov chain} follows \cite{LindMarcus1995}.
Sofic shifts (Exercise~\ref{exerciseSoficShift}) were originally
introduced by Benjamin Weiss in \cite{Weiss1973}.
\index{names}{Weiss, Benjamin}%

The notion of \emph{right-resolving presentation} 
(Exercise~\ref{exerciseMinimalCover}) is close to the notion of \emph{deterministic automaton}
which is classical in automata theory (see~\cite{BerstelPerrinReutenauer2009}).
The unique
right-resolving minimal presentation of a sofic shift 
(Exercise~\ref{exerciseMinimalCover}) is due 
to \cite{Fischer1975} and it is often
called its \emph{Fischer cover}.
\index{names}{Fischer, Roland}%
It is closely related with the notion of \emph{minimal
automaton} of a language (see~\cite{BerstelPerrinReutenauer2009} for example).

The notion of strong shift equivalence 
(Exercise~\ref{exerciseStrongShiftEquivalence}) was introduced
by Robert Williams in \citep{Williams1973}.
\index{names}{Williams, Robert F.}%
The converse of the statement of Exercise 
\ref{exerciseStrongShiftEquivalence} is also true,
and the equivalence is \emph{Williams Classification Theorem}.
\index{subject}{Williams Classification Theorem}%
\index{subject}{Theorem!Williams Classification}%
see~\cite[Theorem 7.2.7]{LindMarcus1995}.
The proof of the converse uses the \emph{Decomposition Theorem}
\index{subject}{Decomposition Theorem} which asserts that
every conjugacy between shifts of finite type can be decomposed
in elementary conjugacies called \emph{input splits}
and \emph{output splits}
\index{subject}{input!split}\index{subject}{output!split}%
and their inverses \emph{input merges}
and \emph{output merges}
\index{subject}{input!merge}\index{subject}{output!merge}%
\cite[Theorem 7.1.2]{LindMarcus1995}. We will define
ouput splits in Chapter~\ref{chapterSubstitutionShifts}.
There is a close connection between strong shift equivalence
and another notion called \emph{shift equivalence}
\index{subject}{shift!equivalence} that we will
present in Chapter~\ref{chapterDirectLimitsOrderedGroups}
(see Exercise~\ref{exerciseShiftEquivalence1}).
\subsection{Substitution shifts}

The substitution shifts of Section~\ref{sectionSubstitutionSystems} form an important class of shifts
and a good part of this book is focused on this
class. For a more detailed treatment, see the
classical reference~\citep{Queffelec2010}\index{names}{Queff\'elec, Martine} or \citep{BertheRigo2010},
where, in particular, Theorem~\ref{propositionPrimitiveMinimal} appears.

The terminology concerning substitution shifts admits
some variations. What we call purely substitutive
sequence is also called a \emph{purely morphic sequence}
\index{subject}{purely!morphic sequence}%
 (see~\cite{AlloucheShallit2003} for example)
or a substitutive sequence (as in~\cite{Queffelec2010}).
Similarly, what we call a substitutive sequence is called
a \emph{morphic sequence} in \cite{Queffelec2010}, \cite{Rigo2014}
and \cite{PytheasFogg2002}.

We warn the reader that our definition of a substitution and of a substitution shift
differs from than that used in \cite{Queffelec2010},
where a substitution is a morphism assumed to be growing and
right-prolongable on some letter.

The property of the Thue-Morse sequence stated in Exercise~\ref{exerciseThueMorse} shows that the Thue-Morse sequence is an \emph{automatic sequence}.
\index{subject}{automatic sequence} Automatic
sequences form a class of substitutive sequences
intoduced by \cite{Cobham1972}\index{names}{Cobham, Alan} under the name
of \emph{uniform tag sequences}.\index{subject}{tag sequence}
 This theory is developped in~\cite{AlloucheShallit2003}.
\index{names}{Allouche, Jean-Paul}\index{names}{Shallit, Jeffrey O.}%
Another example of automatic sequence is the period-doubling sequence of 
Example~\ref{exampleToeplitz} (Exercise~\ref{exerciseDoublingSequence}).
Still another example is the Rudin-Shapiro sequence,
also called the \emph{Golay-Rudin-Shapiro sequence}
 (Exercise~\ref{exerciseRudinShapiro}). It is named after its independent invention
by Golay and Shapiro in connection with extremal problems
in analysis and  in physics, and later
by Rudin~\citep{Rudin1959} (see~\cite{AlloucheShallit2003}).
\index{names}{Golay, Marcel J.E.}%
\index{names}{Rudin, Walter}\index{names}{Shapiro, Harold S.}%

The Chacon binary shift (Exercise~\ref{exerciseChaconMinimal})
is a classical symbolic system (see~\cite{PytheasFogg2002}).

The definition of an elementary morphism and the property
stated in Exercise~\ref{exerciseElementarySubstitution} is 
due to \cite{EhrenfeuchtRozenberg1978}. 
\index{names}{Ehrenfeucht, Andrew}\index{names}{Rozenberg, Gregorz}%
The decidability \index{decidability!of periodicity of morphisms} of periodicity of fixed points
of morphisms (Exercises \ref{exerciseFixedPointPeriodic1} and \ref{exerciseFixedPointPeriodic2}) is due to \cite{Pansiot1986}
\index{names}{Pansiot, Jean-Jacques}%
 and \cite{HarjuLinna1986}
\index{names}{Harju, Tero}\index{names}{Linna, Matti}%
independently.
(see also~\cite{Kurka2003}).
\index{names}{K{\r{u}}rka, Petr}%
 It was extended by \cite{Durand2013}
to the decidability of a more general
question, with periodicity replaced by
eventual periodicity and purely substitutive
sequences replaced by  substitutive sequences.
See \cite[Problem 89]{CrochemoreLecroqRytter2021} for an
algorithmic description of this question.

The constant $K=12$ for the linear recurrrence of the Thue-Morse
shift given in Example \ref{exampleThueMorseLR} is not
 optimal. It is shown in \cite{SchaefferShallit2012}
\index{names}{Schaeffer, Luke}\index{names}{Shallit, Jeffrey O.}%
that the optimal bound is a computable rational number
for every constant length substitution. Actually,
the optimal bound for the Thue-Morse shift, computed using the sotware
Walnut \citep{mousavi2016automatic}
\index{names}{Moussavi, Hamoon} is $K=10$
\citep{Shallit2020}.

Codes and circular codes are described in detail in~\cite{BerstelPerrinReutenauer2009}.
\index{names}{Berstel, Jean}\index{names}{Perrin, Dominique}%
\index{names}{Reutenauer, Christophe}%
 A submonoid satisfying condition \eqref{eqStability}
is called \emph{stable} \index{subject}{stable!submonoid}
\index{subject}{submonoid!stable}%
and a submonoid satisfying \eqref{eqVeryPure}
is called \emph{very pure}.\index{subject}{very pure submonoid}
\index{subject}{submonoid!very pure}%
These notions were originally introduced by Sch\"utzenberger
\citep{Schutzenberger1955}.
\index{names}{Sch\"utzenberger, Marcel-Paul}%
The flower automaton (Exercise~\ref{exerciseFlowerAutomaton})
is a particular case of \emph{finite automaton}.
\index{subject}{finite!automaton}%
The uniqueness of paths with given origin, end and label
defines the so-called \emph{unambigous automata}.
\index{subject}{unambiguous automaton}%
The property of uniform synchronization of finite circular codes 
(Exercise~\ref{exerciseUniformSync}) is due
to \cite{Restivo1975}
\index{names}{Restivo, Antonio}
(see also~\cite[Theorem 10.2.7]{BerstelPerrinReutenauer2009}).

The notion of recognizability for morphisms was introduced initially
by \cite{Martin1973}\index{names}{Martin, John C.} and its status remained uncertain during
many years. 
The definition of recognizability given here 
is from \cite{BezuglyiKwiatkowskiMedynets2009}
who have proved that any (primitive or not) substitution $\sigma$
is recognizable on $X(\sigma)$ for aperiodic points.
\index{names}{Bezuglyi, Sergey}\index{names}{Kwiatkowski, Jan}\index{names}{Medynets, Konstantin}%
The basic result on the subject is Moss\'e's Theorem
(Theorem~\ref{theoremMosse}), which is
from  ~\citep{Mosse1992,Mosse1996}.
\index{names}{Moss\'e, Brigitte}%
Our presentation follows~\cite{Kyriakoglou2019}\index{names}{Kyriakoglou, Revekka}
where Proposition~\ref{propositionRecognizableInjective} is from.
The idea of introducing synchronizing pairs to formulate
recognizability was initiated by \cite{Cassaigne1994}
\index{names}{Cassaigne, Julien}%
and pushed forward by \cite{MignosiSeebold1993}
\index{names}{Mignosi, Filippo}%
and \cite{KloudaStarosta2019}.
\index{names}{Klouda, Karel}\index{names}{Starosta, Stepan}%
The use of synchronizing pairs avoids to rely heavily
on a fixed point of the morphism, as the original formulation
of the Theorem and subsequent proofs do.
The proof of Theorem~\ref{theoremMosse}
presented here follows however essentially
\cite{Kurka2003}.

The finite-to-one maps between shift spaces (and especially shifts of finite
type) introduced in Exercise \ref{exerciseFiniteToOne} are
studied in \cite{LindMarcus1995} (see also \cite{AshleyMarcusPerrinTuncel1993}).
\index{names}{Ashley, John}\index{names}{Marcus, Brian H.}%
\index{names}{Perrin, Dominique}\index{names}{Tuncel, Selim}%
\subsection{Sturmian shifts}
The notion of Sturmian shift (and many ideas of symbolic dynamics
including the term `symbolic dynamics' itself)
was introduced by  \cite{MorseHedlund1938,MorseHedlund1940}.
An introduction  can be found
in~\cite{PytheasFogg2002}\index{names}{Fogg N. Pytheas},
\cite{Lothaire2002}\index{names}{Lothaire, M.} or~\cite{BertheRigo2010}.
\index{names}{Berth\'e, Val\'erie}%
The proof of the equivalent definition of Sturmian sequences
using balanced words (Exercise~\ref{exerciseEquivalentDefSturm})
is taken from \cite{Lothaire2002}.

Arnoux-Rauzy words
are named after the paper~\citep{ArnouxRauzy1991}\index{names}{Arnoux, Pierre}
\index{names}{Rauzy, G\'erard}
in which they are introduced as a generalization
on more than two letters of Sturmian sequences. A reference for
episturmian sequences is~\cite{DroubayJustinPirillo2001}
\index{names}{Droubay, Xavier}\index{names}{Justin, Jacques}%
\index{names}{Pirillo, Giuseppe}%
where Theorem~\ref{standardEpisturmianTheorem} is proved.
We follow here the proof of \cite{DroubayJustinPirillo2001}
using palindrome closure and the notion of Justin word
(called property \emph{Ju} in \citep{DroubayJustinPirillo2001}).
The function $\Pal$ has been introduced by \cite{deLuca1997}
(see also \cite{Reutenauer2019}).
\index{names}{de Luca, Aldo}\index{names}{Reutenauer, Christophe}%
 Justin Formula (Equation~\ref{eqJustin} is from \cite{JustinVuillon2000}\index{names}{Vuillon, Laurent}.
The directive word of a standard episturmian sequence is
 called in \cite{PytheasFogg2002}
the \emph{additive coding sequence}\index{subject}{additive coding sequence}.

Standard sequences (Exercise~\ref{exerciseStandard})
are defined in \cite{Lothaire2002} p. 75.

The statement of \ref{exercisecfrac2} (with a converse) 
is  \cite[Theorem 2.2.31]{Lothaire2002} (see also \cite{Berstel1999}
\index{names}{Berstel, Jean} where
it is credited to \cite{Mignosi1991}.\index{names}{Mignosi, Filippo} Powers in Sturmian
sequences have been extensively studied (see~\cite{DamanikLenz2002,DamanikLenz2003}).
\index{names}{Damanik, David}\index{names}{Lenz, Daniel H.}%
We shall see in Chapter~\ref{chapterDendricShifts}
a closely related statement concerning linearly
recurrent sequences (Corollary~\ref{corollarySturmianLR}).

\subsection{Toeplitz shifts}
Toeplitz sequences were introduced by \cite{JacobsKeane1969}
based on construction of \cite{Toeplitz1928}.
We refer to \cite{Downarowicz:2005}
\index{names}{Downarowicz, Tomasz} for a survey on Toeplitz shifts
(see also \cite{Williams:1984,JacobsKeane1969}).
\index{names}{Williams, Susan G.}
\index{names}{Jacobs, Konrad}\index{names}{Keane, Michael}
The coincidences in substitutions of constant length
have been introduced in \cite{Dekking1978}.
\index{names}{Dekking F. Michel}

%%%%%%%%%%%%%%%%%%%%%%%%%
%  chapter Direct Limits of Ordered Groups
%%%%%%%%%%%%%%%%%%%%%%%%%%%
\chapter{Ordered groups}
\label{chapterDirectLimitsOrderedGroups}
We now introduce notions concerning abelian
groups: ordered abelian groups 
and direct limits of abelian groups. This
will allow us to define 
dimension groups,
which are our main object of interest in this book.

In Section~\ref{sectionOrderedGroups}, we define abelian
ordered groups and the notions of positive morphism or unit
of an ordered group. In Section~\ref{sectionDirectLimits},
we introduce direct limits of ordered groups, an essential notion
for the following chapters. We come in Section~\ref{sectionDimensionGroups}
to the main focus of this book, that is, dimension groups,
as direct limits of groups $\Z^d$ with the usual order. We
 prove the important theorem of Effros, Handelman
and Shen characterizing dimension groups among abelian ordered
groups (Theorem~\ref{theoremEffrosHandelmanShen}). The use
of the term `dimension' in the name of dimension groups
will be explained in the last chapter (Chapter~\ref{chapterBratteli})
where the dimension groups are related to the dimensions of some algebras.
%%%%%%%%%%%%%%%%%%%%%
\section{Ordered abelian groups}\label{sectionOrderedGroups}
By an \emph{ordered group}\index{subject}{ordered!group}\index{subject}{group!ordered }
we mean an abelian group $G$\index{subject}{abelian!group}
\index{subject}{group!abelian}%
with a partial order $<$ which is compatible
with the group operations, that is, such that for all $g,h\in G$
with $g< h$,
one has $g+k< h+k$ for every $k\in G$. 
As usual the partial order $\leq $ is defined by $g\leq h$ whenever $g<h$ or $g=h$, and $g>h$ means $h<g$.

In an ordered group $G$, the \emph{positive cone}\index{subject}{positive!cone}
\index{subject}{cone!positive} is the set 
$$
G^+=\{g\in G\mid g> 0\}\cup \{ 0 \} = \{g\in G\mid g\ge 0\}.
$$
%Note that the term 'positive cone' can be slightly misleading since it contains $0$ and
%is actually the set of nonnegative elements of the group.
\index{symbols}{G@$G^+$} It is a \emph{submonoid} of $G$,
\index{subject}{submonoid!of abelian group}%
\index{subject}{abelian!monoid}\index{subject}{monoid!abelian} that is,
it contains $0$ and satisfies
$G^++G^+\subset G^+$. Moreover $G^+\cap(-G^+)=\{0\}$. Indeed
let $g\in G^+$. If $-g\in G^+$, then $g\ge 0$ and $0\ge g$  which implies
$g=0$ since $\le$ is an order relation.

The set $G^+$ is not a subgroup but, since $G$ is abelian,
the set $G^+-G^+$ is a subgroup which is itself an ordered
group with the same positive cone as $G$. 
%We will always
%assume that in an ordered group the posive cone generates $G$,
%that is, $G^+$ is such that $G^+-G^+=G$.

\begin{proposition}
For any pair $(G,G^+)$ formed of an abelian group $G$
and and a subset $G^+$ of $G$ such that
\begin{displaymath}
G^++G^+\subset G^+,\ G^+\cap (-G^+)=\{0\},
\end{displaymath}
the relation $g< h$ if $h-g\in G^+\setminus \{ 0 \}$ is a partial order
compatible with the group operation  and such that
$G^+$ is the positive cone.
\end{proposition}
\begin{proof}
 The first  condition on $G^+$ 
implies that the partial order $<$ is transitive and the second one that
it is antisymmetric. For $g,h\in G$
such that $g< h$ and $k\in G$, one has $(g+k)-(h+k)=g-h$
and thus $g+k< h+k$. 
\end{proof}
We will often denote an ordered group as a pair $(G,G^+)$
where $G^+$ is the positive cone of $G$.
\begin{example}
The sets $\R^d$ and $\Z^d$ are, for $d\ge 1$, abelian groups
for the componentwise addition
\begin{displaymath}
(x_1,x_2,\ldots,x_d)+(y_1,y_2,\ldots,y_d)=(x_1+y_1,x_2+y_2,\ldots,x_d+y_d)
\end{displaymath}
with $\0=(0,\ldots,0)$ as neutral element.
The pairs $(\R^d,\R_+^d)$ and  $(\Z^d,\Z_+^d)$ with
$\R_+$\index{symbols}{R@$\R_+$} (resp. $\Z_+$)
\index{symbols}{Z@$\Z_+$} formed of the nonnegative reals (resp. integers),
are  ordered groups. The corresponding partial order
on $\R^d$ is called the \emph{natural order}.\index{subject}{natural!order}
\index{subject}{order!natural}%
It is defined by 
\begin{displaymath}
(x_1,x_2,\ldots,x_d)< (y_1,y_2,\ldots,y_d)
\end{displaymath}
if $x_i< y_i$ for $1\le i\le d$.
This partial order is a \emph{lattice order},\index{subject}{lattice order}
\index{subject}{order!lattice} which means that every
pair $x,y$ has a least upper bound, that is, an element $z$ satisfying
$x,y\le z$ and such that $z\le z'$ whenever $x,y\le z'$.
\end{example}

In the next example, the order is a total order.
\begin{example}\label{exampleLexicographic}
Let $G=\Z^d$ ordered by the \emph{lexicographic order}
\index{subject}{lexicographic!order}\index{subject}{order!lexicographic}%
defined by $(x_1,x_2,\ldots,x_d)< (y_1,y_2,\ldots,y_d)$ if there is
an index $i$ with $1\le i\le d$ such that $x_1=y_1$,\ldots,$x_{i-1}=y_{i-1}$
and $x_i<y_i$.
Then $G^+$ is the set
\begin{displaymath}
\{(x_1,\ldots,x_d )\mid x_1=\ldots=x_{i-1}=0 \mbox{ and  $x_i>0$
for some $i$ with $1\le i\le d$}\}\cup\{\0\}.
\end{displaymath}
\end{example}
A \emph{subgroup}\index{subject}{subgroup!of ordered group} of an ordered group $(G,G^+)$ is a
pair $(H,H^+)$ where $H$ is a subgroup of $G$ and $H^+=H\cap G^+$.
Such a subgroup is itself an ordered group. Indeed, $H^+$
is clearly a submonoid and $H^+\cap(-H^+)\subset G^+\cap(-G^+))=\{0\}$.
In this way, the order on $H$ is the restriction
to $H\times H$ of the order on $G$.

An ordered group $G$ is \emph{directed}\index{subject}{directed!ordered group}
\index{subject}{ordered!group!directed} if for every $x,y\in G$
there is some $z\in G$ such that $x,y\le z$. In other terms,
$G$ is directed if every pair of elements has a common upper bound.
\begin{example}\label{exampleFirstComp1}
Let $G=\Z^2$ with the positive cone $G^+=\{(x_1,x_2)\mid x_1>0\}\cup\{(0,0)\}$.
It is a directed group.
\end{example}
In the next statement, we use the fact that if $S$ is a submonoid of
an abelian group $G$, then the set $S-S=\{s-t\mid s,t\in S\}$
is the subgroup generated by $S$ (Exercise~\ref{exerciseSubmonoid}).

\begin{proposition}\label{propositionClifford}
An ordered group $G$ is directed if and only if it is generated
by the positive cone, that is, if $G=G^+-G^+$.
\end{proposition}
\begin{proof}
Assume first that $G^+$ generates $G$. For $x,y\in G$, there
exist $z,t,u,v\in G^+$ such that $x=z-t$ and $y=u-v$. Then
$w=z+u$ is such that $x,y\le w$.

Conversely, for any $x\in G$, considering the pair $0,x$, there is some $y\in G$
such that $0,x\le y$. Then $x=y-(y-x)$ belongs to $G^+-G^+$.
\end{proof}
Note that a subgroup of a directed group need not be directed, as shown in
the next example.
\begin{example}\label{exampleFirstComp}
Let $G=\Z^2$ with the positive cone $G^+=\{(x_1,x_2)\mid x_1>0\}\cup\{(0,0)\}$
as in Example~\ref{exampleFirstComp1}.
Then $H=\{0\}\times\Z$ is a subgroup of a directed group.
But is not directed since $H^+=\{0\}$.
\end{example}

Let $(G,G^+)$ and $(H,H^+)$ be ordered groups.
A morphism $\varphi:G\rightarrow H$ is \emph{positive}\index{subject}{positive! morphism}
\index{subject}{morphism!positive} if $\varphi(G^+)$ is a subset of $H^+$. Note that
a morphism is positive if and only if it preserves the orders
on $G,H$, that is $g\le g'$ implies $\varphi(g)\le\varphi(g')$.

\subsection{Ideals and simple ordered groups}
An \emph{order ideal}\index{subject}{order!ideal}\index{subject}{ideal!of ordered group} $J$ of an ordered group $(G,G^+)$ 
is a subgroup $J$ of $G$ such that $J=J^+-J^+$ (with $J^+=J\cap G^+$)
and such that $0\le a\le b$ with $b\in J^+$ implies $a\in J$.

A \emph{face}\index{subject}{face in ordered group} in $G$ is a subset $F$ of $G^+$ which is a submonoid
and such that $0\le a\le b$ with $b\in F$ implies $a\in F$.

\begin{proposition}\label{propositionIdeal}
Let $\G=(G,G^+)$ be an ordered group. 
\begin{enumerate}
\item For every $g\in G^+$, the set
\begin{displaymath}
[g]=\{h\in G\mid 0\le h\le ng \mbox{ for some $n\ge 0$}\}
\end{displaymath}
is a face. 
\item For every face $F$, the subgroup $J=F-F$ satisfies $J\cap G^+=F$
and is the smallest order ideal of $\G$
containing $F$.
\end{enumerate}
\end{proposition}
\begin{proof}
1. If $h,k$ are in $[g]$ then $h\le ng$ and $k\le mg$ for some $n,m\ge 0$.
Thus $h+k\le(n+m)g$ showing that $h+k$ belongs to $[g]$. Thus $[g]$
is a submonoid.  If $0\le h\le k$ with $n\ge 0$ such that
$k\le ng$, then $0\le h\le ng$
and thus $h$ is in $[g]$. This shows that $[g]$ is a face.

2. The set $J=F-F$ is clearly a subgroup. The set $J^+=J\cap G^+$
is equal to $F$. In fact $F\subset J^+$ by definition. Conversely
if $h\in J^+$, set $h=a-b$ with $a,b\in F$. Then $h\le a$
implies $h\in F$ since $F$ is a face. This shows that
$J$ is an ideal. Finally, if $K$ is an ideal containing $F$,
then  $J\subset K$ since $K$ is a subgroup.
\end{proof}
An ordered group is \emph{simple}\index{subject}{simple!ordered group}
\index{subject}{ordered!group!simple} if it has no nonzero
proper order ideals. Note that a simple group is directed since
$G^+-G^+$ is an ideal of $(G,G^+)$, and thus $G^+-G^+ = G$.
\begin{proposition}
A subgroup of a simple ordered group is simple.
\end{proposition}
\begin{proof} Let $(H,H^+)$
be a subgroup of the simple group $(G,G^+)$. Let $J$
be a nonzero order ideal of $H$. Let $k$ be a nonzero
element of $J^+$. Then $F=\{h\in G\mid
h\le nk\mbox{ for some } n\ge 0\}$ is a face by Proposition
\ref{propositionIdeal}. Thus, by Proposition~\ref{propositionIdeal} again,
 $K=F-F$ is an ideal and $K\cap G^+=F$.
Since $G$
is simple, we have $K=G$ and also $F=G^+$. Thus every $h$ belonging to $H^+$
is in $F$ and consequently in $J$, which shows that $J=H$. 
Therefore $H$ is simple.
\end{proof}
%If $H$ is subgroup of $G$ and and if $(G,G^+)$
%is simple group, the  is simple. Indeed, let $H$ be a 
%subgroup of $(G,G^+,1)$.
\begin{example}\label{exampleSimple}
The ordered groups $\Z,\Q$ and $\R$ (with the natural order) are simple.
On the contrary, the ordered group   $(\Z^d,\Z_+^d)$ for $d\ge 2$
 is not simple. Indeed, for $d=2$,
the set $\Z\times\{0\}$ is an order ideal of $\Z^2$.

In contrast, for every irrational $\alpha$, the group
 $\Z+\alpha\Z$, with the order induced by $\R$, is simple
since it is a subgroup of $\R$. 
We shall meet this simple ordered group several times.
It is contained in the additive subgroup $\Z[\alpha]$ of $\R$
generated by the powers of $\alpha$. When $\alpha$ is
an \emph{algebraic integer},\index{subject}{algebraic!integer} that is, such that $p(\alpha)=0$ for some
polynomial $p(x)=x^{k+1}+a_{k}x^{k}+\ldots+a_1x+a_0$
with $a_i$ in $\Z$, the group $\Z[\alpha]$ is generated 
by $1,\alpha,\ldots,\alpha^k$ and is thus finitely generated.
For $k=2$, the groups $\Z+\alpha\Z$ and $\Z[\alpha]$ coincide
(see Appendix~\ref{appendixAlgebraicNumberTheory}).
\end{example}
An \emph{order unit}\index{subject}{order! unit} of the ordered group $G$
is a positive element $u$ such that for every $g\in G^+$
there is an integer $n>0$ such that $g<nu$.
Equivalently, $u$ is an order unit if the set $[u]$
defined in Proposition~\ref{propositionIdeal} is equal to $G^+$.

A \emph{unital ordered group}\index{subject}{unital!ordered group}
\index{subject}{ordered!group!unital} is a triple $(G,G^+,\mathbf{1}_G)$
\index{symbols}{G@$(G,G^+,u)$}%
formed of an ordered group $(G,G^+)$ and an order unit $\mathbf{1}_G$.

Note that if $u$ is a unit of a directed ordered group $G$, then
for every $g\in G$ there is an $n>0$ such that $g<nu$. Indeed,
since $G$ is directed there is $h$ such that $g,u\le h$. Then
$u\le h$ implies $h\in G^+$. If $n>0$ is such that $h<nu$, we have
$g\le h< nu$.

\begin{proposition}\label{propositionSimpleUnit}
A directed ordered group $(G,G^+ )$ is simple if and only if
every nonzero element of $G^+$ is an order unit.
\end{proposition}

The proof is left as an exercise (Exercise~\ref{exerciseSimple}).

\begin{example}
The triples $(\R^d,\R_+^d,\mathbf{1})$ and  $(\Z^d,\Z_+^d,\mathbf{1})$ 
are unital ordered groups
with order unit ${\mathbf 1}=(1,1,\ldots,1)$.
\end{example}

A \emph{morphism}\index{subject}{morphism!of unital ordered groups} of 
unital ordered groups from $\G=(G,G^+,\mathbf{1}_G)$
to ${\cal H}=(H,H^+,\mathbf{1}_H)$ is a group morphism $\varphi:G\rightarrow H$
which is positive and
such that $\varphi(\mathbf{1}_G)=\mathbf{1}_H$.

A \emph{subgroup}\index{subject}{subgroup!of unital ordered group}
 of a unital ordered group $\G=(G,G^+,u)$ is a unital ordered group
${\cal H}=(H,H^+,u)$ such that $H$ is a subgroup of $G$ containing $u$
and $H^+=H\cap G^+$. 
\subsection{Unperforated ordered groups}
An ordered group $(G,G^+)$ is \emph{unperforated}\index{subject}{unperforated group}
\index{subject}{group!unperforated}
if for every $g\in G$, and $n>0$, if $ng$ belongs to $G^+$ then $g$ is in $G^+$.
Otherwise, the group is \emph{perforated}.\index{subject}{perforated group}
\index{subject}{group!perforated}

For example, $\Z^d$ with the natural order is unperforated.
\begin{example}
The group $G=\Z$ with positive cone the submonoid $G^+$ of $\Z$ generated
by the set 
$\{2,5\}$ is perforated since $9=3\times 3$ belongs to $G^+$ whereas $3$ does not.
\end{example}
A group $G$ is \emph{torsion-free}\index{subject}{torsion-free group}\index{subject}{group!torsion-free} if for every $g\in G$ and $n>0$,  $ng=0$ implies $g=0$.
\begin{proposition}
An unperforated group 
is torsion-free.
\end{proposition}
\begin{proof}
Suppose that $ng=0$ for $g\in G$ and $n>0$. Then $g\in G^+$ since
$G$ is unperforated. But $-g=(n-1)g$ implies $-g\in G^+\cap(-G^+)$
and thus $g=0$.
\end{proof}
As an example, the group $G=\Z\times \Z/2\Z$ with 
$G^+=\{(\alpha,\beta)\in G\mid \alpha>0\}\cup\{(0,0)\}$
is a perforated group. Indeed, $2(0,1)=(0,0)$
and thus $G$ has torsion.
%%%%%%%%%%%%%%%%%%
\section{States}\label{sectionStates}
Let $\G=(G,G^+,\mathbf{1})$ be a unital ordered group.
A \emph{state}
\index{subject}{state!of unital ordered group}%
\index{subject}{ordered!group!state of}%
%(also called a \emph{trace})
%\index{subject}{trace of unital ordered group}%
%\index{subject}{unital!ordered group!trace on}%
%\index{subject}{ordered!group!trace}%
on $\G$ 
is a morphism of unital ordered groups 
from $\G$ to the unital ordered group $(\R,\R_+,1)$.
Thus a group morphism $p:G\rightarrow \R$ is a state if $p(g)\ge 0$
for every $g\in G^+$ and $p(\mathbf{1})=1$.

Let $S(\G)$\index{symbols}{S@$S(G)$} denote the set of states of $\G$. It is a convex set.
Indeed, let $p,q\in S(G)$ and let $t\in[0,1]$. Then $r=tp+(1-t)q$
is a morphism since the set of morphisms from
a group to $\R$ forms a vector space. It is
positive because for every $g\in G^+$, 
we have $r(g)\ge\min\{p(g),q(g)\}\ge 0$. Moreover $r(\1)=tp(\1)+(1-t)q(\1)=1$. Thus
$tp+(1-t)q$ belongs to $S(\G)$ which shows that $S(\G)$ is convex.
\begin{example}
Let $\G$ be the group $\Z^d$ with the usual order and $u=(1,\ldots,1)$. The set
$S(\G)$ is the $d-1$ simplex (see Appendix \ref{appendixMeasureIntegration}
for the definition) formed of the maps
$(\alpha_1,\ldots,\alpha_d)\mapsto p_1\alpha_1+\ldots+p_d\alpha_d$
for $p_1,\ldots,p_d\ge 0$ of sum $1$.
\end{example}

\begin{example}
Let $\lambda$ be irrational. The unital ordered group $G=\Z+\lambda\Z$
(with the order induced by the order on $\R$, and, order unit $1$) is simple,
as we have already seen (Example~\ref{exampleSimple}).
There is a unique state which is the identity. Indeed, let $p$
be a state on $G$. 
For every $x,y\in\Z$, as $p(1) = 1$, one classically obtains that $p(x)=x$, $p(y)=y$, and, thus
$p(x+\lambda y ) = x +p(\lambda )y$. 
Moreover, if $x+\lambda y\ge 0$, then $x+p(\lambda ) y\ge 0$. 
This implies $p (\lambda )=\lambda$ and show the unique state is the identity. 
\end{example}

\begin{proposition}
\label{prop:statescompact}
Let $\G$ be a directed unital ordered group.
The set $S(\G)$ is convex and compact for the product topology on $\R^G$. 
\end{proposition}
\begin{proof}
For every $g\in G$, we have
$g=g'-g''$ with $g',g''\in G^+$ since $G$ is directed (Proposition \ref{propositionClifford}). 
Let $n\ge 1$ be such that
$g',g''\le n\1$. Then $|p(g)|\le|p(g')|+|p(g'')|\le2n\1$ 
for every $p\in S(\G)$.
This implies that $p(g)\in[-n,n]$ for every $p\in S(\G)$
and thus that $S(\G)$ is compact for the product topology.
%For $g\in G$, we denote $\hat{g}:S(\G)\rightarrow \R$ the map defined by 
%$\hat{g}(p)=p(g)$ for every $p\in S(\G)$.
\end{proof}

We will now prove the following important result.

\begin{theorem}\label{theoremGoodearlHandelman}
For every directed unital ordered group $\G$, the set of states on $\G$ is nonempty.
\end{theorem}

We first prove the following lemmas.

\begin{lemma}\label{lemma3.1}
Let $\G=(G,G^+,u)$ be a unital ordered group
and let ${\cal H}=(H,H^+,u)$
be  a unital ordered subgroup of $\G$.
 Let $p$ be a state on $\cal H$.
 For $g\in G^+$,
let
\begin{eqnarray}
\alpha(g)&=&\sup\{p(x)/m\mid x\in H,m>0,x\le mg\},\label{eqalpha}\\
\beta(g)&=&\inf\{p(y)/n\mid y\in H,n>0,ng\le y\}.\label{eqbeta}
\end{eqnarray}
Then
\begin{enumerate}
\item $0\le\alpha(g)\le\beta(g)<\infty$.
\item If $q$ is a state on $H+\Z g$ which extends $p$, then
$\alpha(g)\le q(g)\le\beta(g)$.
\item If $\alpha(g)\le\gamma\le\beta(g)$, there is a state $q$ on 
$(H+\Z g,H\cap (G^++\Z_+ g),u)$
such that $q$ extends $p$ and $q(g)=\gamma$.
\end{enumerate}

\end{lemma}
\begin{proof}
1. Since $0\le g$, we have by \eqref{eqalpha} with $x=0$ and $m=1$, that
$\alpha(g)\ge p(0)/1=0$. Next, since $u$ is an order unit, we have $g\le ku$
for some $k>0$. Thus, using \eqref{eqbeta} with $y=ku$ and
$n=1$, we obtain $\beta(g)\le p(ku)/1=k<\infty$.

Consider $x,y\in H$ and $m,n>0$ such that $x\le mg$ and $ng\le y$. Then
$nx\le mng\le my$, and consequently $p(x)/m\le p(y)/n$. Therefore
$\alpha(g)\le\beta(g)$.

2. Let $x\in H$ and $m>0$ be such that $x\le mg$. Then
$p(x)=q(x)\le m q(g)$ and thus $\alpha(g)\le q(g)$. The proof
that $q(g)\le\beta(g)$ is similar.

3. We first claim that if $z+kg\ge 0$ for some $z\in H$ and $k\in\Z$, then
$p(z)+k\gamma\ge 0$. Indeed, if $k=0$, then $z\ge 0$ 
and $p(z)+k\gamma=p(z)\ge 0$. If $k>0$, then we have $-z\le kg$
with $-z\in H$ whence, by \eqref{eqalpha}, $p(-z)/k\le\alpha(g) \le\gamma$
and so $p(z)+k\gamma\ge 0$. Finally, if $k<0$, we have $-kg\le z$ 
with $z\in H$ and $-k>0$, hence $\gamma\le\beta(g)\le p(z)/(-k)$
and so $p(z)+k\gamma\ge 0$. This proves the claim.

As a consequence of the claim, we obtain that if $z+kg=0$ for some
$z\in H$ and $k\in\Z$, then $p(z)+k\gamma=0$. Indeed,
we have both $z+kg\ge 0$ which implies $z+k\gamma\ge 0$
and $-z+(-k)g\ge 0$ which implies $-z+(-k)\gamma\ge 0$.
As a consequence,
the map $z+kg\mapsto p(z)+k\gamma$ induces a morphism $q$ from 
$H+\Z g$ to $\R$. Moreover, the claim shows that $q$ is positive.
Since $q(g)=\gamma$ and $q(u)=p(u)=1$, the proof is complete.
\end{proof}

\begin{lemma}\label {lemmaTheorem3.2}
Let $\G=(G,G^+,u)$ be a directed unital ordered group
and let ${\cal H}=(H,H^+,u)$
be a subgroup of $\G$.
Every state on $\cal H$ extends to
a state on $\G$.
\end{lemma}
\begin{proof}
Consider the family of pairs $({\cal K},q)$ of a unital ordered
group ${\cal K}=(K,K^+,u)$ such that $K$
is a subgroup of $G$ which contains $H$
with $K^+=K\cap G^+$ and a state $q$ on ${\cal K}$
which extends $p$.
By Zorn's Lemma, this family has a maximal element $(K,q)\in{\cal K}$.
Suppose that $K\ne G$. Since $\G$ is directed,
$G^+$ generates $G$, which implies that there is some 
$g\in G^+\setminus K$. By Lemma \ref{lemma3.1} there
is a state $q$ on $K+\Z g$ which extends $p$. But then $(K+\Z g,q)$
is in ${\cal K}$, a contradiction. 
We conclude that $K=G$ and that $q$ is a state on $\G$ which extends $p$.
\end{proof}

We are now ready for the proof of Theorem~\ref{theoremGoodearlHandelman}.

\begin{proofof}{of Theorem~\ref{theoremGoodearlHandelman}}
Set $\G=(G,G^+,u)$. The map $p:nu\mapsto n$ is a state
on the unital ordered group ${\cal H}=(\Z u,\Z_+ u,u)$
which is a unital ordered subgroup of $\G$. By Lemma~\ref{lemmaTheorem3.2},
$p$ extends to a state of $G$.
\end{proofof}

\begin{example}\label{exampleState}
Let $\G=(G,G^+,u)$ with $G=\Z^2$, $G^+=\{(x_1,x_2)\in G\mid x_1>0\}\cup\{(0,0)\}$
and $u=(1,0)$. There is a unique state which is the projection on the first
component.
\end{example}

%One has in ordered groups the following minimax principle,
%based on Lemma~\ref{lemma3.1}.

%  \begin{proposition}\label{propositionMiniMax}
%Let $\G=(G,G^+,u)$ be a unital ordered group. For every $g\in G^+$,
%one has 
%\begin{displaymath}
%\inf\{p(g)\mid p\in S(\G)\}=\sup\{\varepsilon\in\Q_+\mid \varepsilon u\le g\}.
%\end{displaymath}
%  \end{proposition}
%  The proof is left as Exercise~\ref{exerciseMiniMax}.
The following result
shows that, for an unperforated simple ordered group,
 the order is determined by the set of states.
\begin{proposition}\label{propositionEffros}
If $\G=(G,G^+,u)$ is an unperforated simple unital  group, then
$G^+=\{g\in G\mid p(g)>0\mbox{ for every }p\in S(\G)\}\cup \{0\}$.
\end{proposition}
\begin{proof}
Since $G$ is simple, every nonzero element 
of $G^+$ is an order unit (Proposition~\ref{propositionSimpleUnit}).
Thus if $g$ belongs to $G^+\setminus \{0\}$, 
there is $n\ge 1$ such that $ng\ge u$. Then, for any
$p\in S(\G)$ (which is nonempty by 
Theorem~\ref{theoremGoodearlHandelman}), we have $p(ng)\ge p(u)=1$ and thus $p(g)>0$
since $p(ng)=np(g)$. Conversely, if $p(g)>0$ for every
$p\in S(\G)$, then since $S(\G)$ is closed by Proposition \ref{prop:statescompact}, there is $\varepsilon >0$
such that $p(g)>\varepsilon$ for every $p\in S(\G)$.
Let $r,s>0$ be such that $p(g)>r/s$ for every $p\in S(\G)$.
Since $p(sg-ru)\ge 0$ for every $p\in S(\G)$ and since $G$
is unperforated, we have $(sg-ru)+u\ge 0$ (Exercise~\ref{exerciseLemma6.1}).
Thus $sg\ge (r-1)u\ge 0$, which implies $g\in G^+$ since $G$ is
unperforated.
\end{proof}

\begin{example}\label{exampleStates}
Let $G=\Z^2$ with $G^+=\{(x_1,x_2)\mid x_1>0\}\cup\{(0,0)\}$
as in Example~\ref{exampleState}. There is a unique
state $p:(x_1,x_2)\to x_1$. Thus Proposition~\ref{propositionEffros}
is satisfied. In contrast, let $G=\Z^2$ with the lexicographic
order, that is with $G^+=\{(x_1,x_2)\mid x_1>0 \mbox{ or $x_1=0$
and $x_2\ge 0$}\}$. There is only one state which is the
projection on the first component. Thus Proposition
\ref{propositionEffros} does not hold. There is no contradiction since
$G$ is not simple.
\end{example}

\subsection{Infinitesimals}
Let $(G,G^+,u)$ be a unital unperforated ordered group. We say that an element 
$g\in G$ is \emph{infinitesimal}
\index{subject}{infinitesimal!element} 
if $ng\le u$ for every $n\in\Z$.
It is easy to see that the definition does not depend 
on the choice of the order unit $u$ (Exercise~\ref{exerciseInfinitesimals1}).

If $G=\Z^d$ with the usual order and unit
$\1=(1,1,\ldots,1)$, there are no nonzero infinitesimals. On the
contrary the following example exhibits nonzero infinitesimals.
\begin{example}\label{exampleInfinitesimals}
Let $G=\Z^2$ with $G^+=\{(\alpha,\beta)\in G\mid \alpha>0\}\cup\{(0,0)\}$
and $u=(1,0)$ as in Example~\ref{exampleState}. Any element $(0,\beta)$ with $\beta\in\Z$
is infinitesimal.
\end{example}
Let us introduce the following useful notation. For
$\varepsilon\in\Q$, 
the inequality  $g\leq \varepsilon u$ means that $qg \leq pu$
for some integers $p,q\ge 1$ such that $\varepsilon=p/q$. 

 Using this notation, one can give as an equivalent
definition  that $g$ is infinitesimal
if $-\varepsilon u \leq g\leq \varepsilon u$ 
for all $\varepsilon >0 $ in $\mathbb{Q}$ 
(Exercise~\ref{exerciseInfinitesimals2}).

Another equivalent definition is the following.
\begin{proposition}\label{propositionInfinitesimals}
Let $\G=(G,G^+,u)$ be an unperforated simple  unital group.
An element $g\in G$ is infinitesimal if and only 
if $p(g)=0$ for all $p \in  S(\G)$. 
\end{proposition}
\begin{proof}
By Theorem~\ref{theoremGoodearlHandelman}, the set $S(\G)$
of states is nonempty.
Assume first that $g$ is infinitesimal and let $p\in S(\G)$.
Since $-\varepsilon\le p(g)\le\varepsilon$ for every $\varepsilon\in\Q_+$,
we conclude that $p(g)=0$. 

Conversely, suppose that $|p(g)|\ge 1/n$ for some trace $p$.
Then $u-ng$ and $u+ng$ cannot be both in $G^+$. Thus either
$ng\le u$ or $-ng\le u$ is false, a contradiction.
\end{proof}

The collection of infinitesimal elements of $G$ forms a  subgroup, 
called
the \emph{infinitesimal subgroup of $G$}\index{subject}{infinitesimal!subgroup}, 
which we denote by $\Inf(G)$.\index{symbols}{Inf@$\Inf(G)$}%

The quotient group $G/\Inf(G)$ of a simple ordered group $G$
is also a simple ordered group for the induced order, 
and the infinitesimal subgroup of 
$G/\Inf(G)$ is trivial (see Exercise~\ref{exerciseInfGroup}). 
Furthermore, an order unit for $G$ maps 
to an order unit for  $G/\Inf(G)$.
Moreover the traces space of $G$ and $G/\Inf (G)$ are isomorphic.

%When  $S(G,G^+,u)$ consists of a unique trace,  notice that 
%$G/{\rm Inf}(G)$ is isomorphic to $(I(G,G^+,u), I(G,G^+,u)\cap
%\RR^+,1)$, as ordered groups with unit.

It may be interesting to summarize the properties of some
of the various orders on $\Z^2$ that we have considered
in the examples.
\begin{figure}[hbt]
\centering
\begin{tabular}{|l|l|l|l|l|}
\hline
$G^+$            &directed&total&simple\\\hline
$\Z_+\times\{0\}$&no      &no  &no    \\ \hline
natural          &yes     &no  &no    \\ \hline
lexicographic    &yes      &yes&no    \\ \hline
$(\Z_+\setminus\{0\})\times\Z\cup\{(0,0)\}$    &yes      &no &yes  \\ \hline
\end{tabular}
\caption{The properties of various orders on $\Z^2$.}
\end{figure}
The first row concerns the group $\Z^2$ ordered by the
first component. It is not directed since $G^+$ generates
$\Z\times\{0\}$.
The third row concerns the group of Example~\ref{exampleLexicographic}.
 The last one is the group of Example~\ref{exampleFirstComp}.

\subsection{Image subgroup}

Let $\G=(G,G^+,u)$ be a unital ordered group. The \emph{image subgroup} 
associated to $\G$
\index{subject}{image subgroup} is the subgroup of $(\R,\R_+,1)$ defined as
\begin{equation}
I(\G)=\bigcap_{p\in S(\G)}p(G)\label{eqImageSubgroup}
\end{equation}
\index{symbols}{I@$I(G)$}%
When $S(\G)$ consists of a unique trace, the group $I(\G)$
is isomorphic to $\G/\Inf(\G)$.
\begin{example}
Let $G=\Z^2$ with $G^+=\{(\alpha,\beta)\in G\mid\alpha>0\}\cup\{(0,0\}$
as in Example~\ref{exampleInfinitesimals}. Then
the image subgroup is $I(\G)=(\R,\R_+,1)$.
\end{example}
%%%%%%%%%%%%%%%%%%%%%%%
\section{Direct limits}\label{sectionDirectLimits}
We now introduce the important notion of direct limit,
which is central in this book. We will first
formulate it for ordinary abelian groups and next for ordered ones.
\subsection{Direct limits of abelian groups}
Let $G_n$ be for each $n\ge 0$
an   abelian group
 and let $i_{n+1,n}:G_n\rightarrow G_{n+1}$ be for every $n\ge 0$
a  morphism. The sets
\begin{displaymath}
\Delta=\{(g_n)_{n\ge 0}\mid g_n\in G_n, g_{n+1}=i_{n+1,n}(g_n) \mbox{ for every $n$ large enough}\}
\end{displaymath}
\index{symbols}{Delta@$\Delta$}%
and
\begin{displaymath}
\Delta^0=\{(g_n)_{n\ge 0}\mid g_n\in G_n, g_{n}=0 \mbox{ for every $n$ large enough}\}
\end{displaymath}
\index{symbols}{Delta@$\Delta^0$}%
are subgroups of the direct product $\Pi_{n\ge 0}G_n$ 
and $\Delta^0\subset \Delta$.
Let $G$ be the quotient group $G=\Delta/\Delta^0$ and 
$\pi : \Delta \to G$ be the natural projection.
The
group $G$, denoted $G=\displaystyle{\lim\limits_\rightarrow G_n}$,
\index{symbols}{lim@$\lim\limits_\rightarrow G_n$}%
 is called the \emph{direct limit}
 \index{subject}{direct limit!of groups}
 (or \emph{inductive limit}\index{subject}{inductive limit!of groups})
of the sequence $(G_n)_{n\ge 0}$ with the maps $i_{n+1,n}$.
The maps $i_{n+1,n}$ and more generally the maps $i_{m,n}=i_{m,m-1}\circ\cdots\circ
i_{n+1,n}$ for $n<m$ are called the 
\emph{connecting morphisms}\index{subject}{connecting!morphism}\index{subject}{morphism!connecting}.

See Exercise~\ref{exerciseAltDefDirectLimits} for an alternative
definition of the direct limit as a quotient of the union
of the $G_n$.

Given $g\in G_n$, all the sequences $(g_k)_{k\ge 0}\in\Pi_{k\ge 0} G_k$ such that
\begin{displaymath}
g_n=g\mbox{ and }g_{m+1}=i_{m+1,m} (g_m )\mbox{ for all } m\ge n
\end{displaymath}
belong to $\Delta$ and have the same projection, denoted $i_n(g)$ in $G$.
One easily checks that this defines a morphism $i_n$ from $G_n$ to $G$. 
It is called the 
\emph{natural morphism}\index{subject}{natural!morphism}
\index{subject}{morphism!natural into direct limit} from $G_n$ into $G$.

Note that $i_n:G_n\to G$ is such that
$i_n=i_{n+1}\circ i_{n+1,n}$ for every $n\ge 0$.
Its kernel is
\begin{displaymath}
\ker(i_n)=\cup_{m\ge n}\ker(i_{m,m-1}\circ\cdots\circ i_{n+1,n}).
\end{displaymath}
The group $G$ is the union of  the ranges of the $i_n$.
Thus, for every $g\in G$, there exist an integer $n\ge 0$
and an element $g_n\in G_n$ such that $g=i_n(g_n)$.

\begin{example}\label{exampleDyadic}
Consider the case of the sequence $\Z\edge{2}\Z\edge{2}\Z\ldots$ 
where each map
is the multiplication by $2$. The direct limit $G$ of this sequence
can be identified with the group $\Z[1/2]$ of \emph{dyadic rationals},
formed of all
rational numbers $p/q$ with $q$ a power of $2$.
\index{subject}{dyadic rationals}%
\index{symbols}{Z@$\Z[1/2]$}%
Indeed, $\Delta$ is formed of the sequences $(g_n)_{n\ge 0}$ such that
for some $k\ge 1$, one has $g_{n+1}=2g_n$ for every $n\ge k$.
Consider the map $\pi:\Delta\rightarrow \Z[1/2]$ sending
such a sequence on $2^{-k}g_k$. This map is a well defined
group morphism and its kernel is $\Delta^0$. Thus
it induces an isomorphism from $G$ onto $\Z[1/2]$.
The natural morphism from $G_n=\Z$ to $G$
is $i_n(g)=2^{-n}g$.
\end{example}

Direct limits have the following property (called a
\emph{universal property}\index{subject}{universal!property of direct limits})
expressing that the direct limit is, in some sense, the largest
possible abelian group.
\begin{figure}[hbt]
\centering
\begin{tikzpicture}
\node(G_n)at(0,1){$G_n$};\node(G_{n+1})at(2,1){$G_{n+1}$};
\node(H)at(1,0){$H$};
\draw[->,above](G_n)edge node{$i_{n+1,n}$}(G_{n+1});
\draw[->,left](G_n)edge node{$\alpha_n$}(H);
\draw[->,right](G_{n+1})edge node{$\alpha_{n+1}$}(H);

\node(G_n)at(5,1){$G_n$};\node(G)at(7,1){$G$};
\node(H)at(6,0){$H$};
\draw[->,above](G_n)edge node{$i_n$}(G);
\draw[->,left](G_n)edge node{$\alpha_n$}(H);
\draw[->,right](G)edge node{$\varphi$}(H);
\end{tikzpicture}
\caption{The universal property of the direct limit}
\label{figureUniversal}
\end{figure}
\begin{proposition}\label{propositionUniversalProperty}
Let $(G_n)$ be a sequence of abelian groups with
connecting morphisms $i_{n+1,n}:G_n\to G_{n+1}$. For every abelian
group $H$ and every
sequence $(\alpha_n)$ of morphisms from $G_n$ to $H$
such that $\alpha_n=\alpha_{n+1}\circ i_{n+1,n}$
(see Figure~\ref{figureUniversal} on the left),
there is a unique morphism $\varphi$ from the direct limit
$G=\lim\limits_\rightarrow G_n$ to $H$ such that
\begin{equation}
  \alpha_n=\varphi\circ i_n\label{equationUniversalProperty}
  \end{equation}
for all $n\ge 0$ (see Figure~\ref{figureUniversal} on the right).
\end{proposition}
The verification is left as Exercise~\ref{exerciseUniversal}.
\subsection{Direct limits of ordered groups}

Let now $\G_n=(G_n,G_n^+,\mathbf{1}_n)$ be for each $n\ge 0$
a  directed unital ordered group
 and let $i_{n+1,n}:G_n\rightarrow G_{n+1}$ be for every $n\ge 0$
a  morphism of unital ordered groups.
Let $G$ be the direct limit of the sequence $(G_n)_{n\geq 0}$, let
$G^+$ be the projection in $G$ of the set
\begin{displaymath}
\Delta^+=\{(g_n)_{n\ge 0}\mid g_{n}\in G_n^+ \mbox{ for every  large enough $n$}\}
\end{displaymath}
\index{symbols}{Delta@$\Delta^+$}%
and let $\mathbf{1}$ be the projection in $G$ of the sequence $(\mathbf{1}_n)_{n\ge 0}\in\Delta$.
\begin{proposition}\label{propositionDirectLimit}
The triple $\G=(G,G^+,\mathbf{1})$ is a directed unital ordered group
and every natural morphism $i_n:\G_n\to\G$, $n\in \N$,
is a morphism of unital ordered groups and $G^+$
is the union of the $i_n(G^+)$.
\end{proposition}
\begin{proof}
We first verify that $G^+$ satisfies the two conditions defining an ordered group.
First, if $g,g'$ belong to  $\Delta^+$, then $g_n,g'_n$ belong to $G_n^+$ for every 
large enough $n$ and thus $g_n+g'_n$ also belong to $G_n^+$ for every large enough $n$.
Thus $g+g'$ belongs to $\Delta^+$. This shows that $\Delta^++\Delta^+$ is included in $\Delta^+$
and it implies that $G^++G^+$ is a subset of $G^+$. Similarly, the facts that
$G=G^+-G^+$ and $G^+\cap(-G^+)=\{0\}$ follow from
$\Delta=\Delta^+-\Delta^+$ and $\Delta^+\cap(-\Delta^+)=\Delta^0$.
Thus $\G$ is directed.

Finally, let us show that $\mathbf{1}$ is an order unit. 
Let $g=(g_n)_{n\ge 0}\in\Delta$. If $g_{n+1}=i_n(g_n)$ for some $n\ge 0$,
and if  $g_{n}\le k\mathbf{1}_n$ for some $k\ge 1$,
then, since
$i_{n+1,n}$ is a morphism of ordered groups with order unit, we
have $g_{n+1}\le k\mathbf{1}_{n+1}$. Thus, there is a $k\ge 1$
such that $g_n\le k\mathbf{1}_n$ for every $n$ large enough. This implies
that the projection of $g$ in $G$ is bounded by $k\mathbf{1}$.
\end{proof}
The triple $\G$ is called the \emph{direct limit}
\index{subject}{direct limit!of ordered groups} (or \emph{inductive limit}\index{subject}{inductive limit!of ordered groups})
of the sequence $(\G_n)_{n\ge 0}$
(see Exercise~\ref{exerciseAltDefDirectLimits2} for an alternative
definition).

The universal property of direct limits (Proposition~\ref{propositionUniversalProperty})
holds for direct limits of ordered groups.

\begin{example}\label{exampleDyadic2}
Consider again the case of the sequence $\Z\longrightarrow^{{\hskip -14 pt {\small \times 2}}} \Z\longrightarrow^{{\hskip -14 pt {\small \times 2}}}\cdots$
of Example~\ref{exampleDyadic} with $\Z$ ordered as usual.
The direct limit of the corresponding sequence of ordered groups
is the ordered group $(\Z[1/2],\Z_+[1/2],1)$
where $\Z_+[1/2]$ is the set of non negative dyadic rationals.
\end{example}

\subsection{Ordered group of a matrix}
We can generalize Example~\ref{exampleDyadic} by considering
$G_n=\Z^d$ for some integer $d\ge 1$,
an integer $d\times d$-matrix $M$ and the sequence of morphisms
$i_n$ being the multiplication by $M$ on the elements of $\Z^d$
considered as column vectors.

Thus we address the description of the direct limit of a sequence
\begin{displaymath}
  \Z^d\edge{M}\Z^d\edge{M}\Z^d\ldots
\end{displaymath}
of groups all equal to $\Z^d$ with the same connecting morphisms.

Define the \emph{eventual range}\index{subject}{eventual! range} of $M$ as 
\begin{displaymath}
\RR_M=\cap_{k\ge1}M^k\Q^d
\end{displaymath}
\index{symbols}{R@$\RR_M$}%
and the \emph{eventual kernel}\index{subject}{eventual! kernel} of $M$ as
\begin{displaymath}
\KK_M=\cup_{k\ge1}\ker(M^k).
\end{displaymath}
\index{symbols}{K@$\KK_M$}%
Note that 
\begin{equation}
\Q^d=\RR_M\oplus \KK_M\label{eqRK}
\end{equation}
 and that the multiplication by $M$
defines an automorphism of
$\RR_M$.  
Indeed, since 
\begin{displaymath}
\ldots \subset M^2\Q^d\subset M\Q^d\subset \R^d\mbox{ and } \ker M\subset \ker M^2\subset\ldots
\end{displaymath}
there is some $h\ge 0$ such that $\RR_M=M^h\Q^d$ and $\KK_M=\ker M^h$.
Then $M\RR_M=M^{h+1}\Q^d=M^h\Q^d=\RR_M$ and thus the multiplication by $M$ is 
an automorphism of $\RR_M$.
If $x$ is in $\RR_M$, then $Mx=0$ implies $x=0$. Thus $\RR_M\cap \KK_M=\{0\}$.
Next, for every $x\in \Q^d$ there is, since $M^hx$ belongs to $\RR_M$,
some $y\in\RR_M$ such that $M^hx=M^hy$. Then
$x=y+(x-y)$ belongs to $\RR_M+\KK_M$. This proves Equation~\eqref{eqRK}

Let also
\begin{equation}
\Delta_M=\{v\in \RR_M\mid M^kv\in \Z^d\mbox{ for some $k\ge 0$}\}
\label{eqDeltaM}
\end{equation}
\index{symbols}{Delta@$\Delta_M$}%
The following result describes  direct limits with identical connecting
morphisms $i_{n+1,n}$ for all $n\ge 0$. 
\begin{proposition}\label{propositionDeltaM}
 For every integer $d\times d$-matrix, the direct limit $G$ of the sequence
\begin{displaymath}
\Z^d\longrightarrow^{{\hskip -12 pt M}}  \hskip 3pt \Z^d \longrightarrow^{{\hskip -12 pt M}}\hskip 3pt \Z^d\ldots ,
\end{displaymath} where each map is the
multiplication by $M$, is isomorphic to $\Delta_M$.
If moreover the matrix $M$ is nonnegative,
the triple $(\Delta_M,\Delta_M^+,\1_M)$, where
\begin{displaymath}
\Delta_M^+=\{v\in \RR_M\mid M^kv\in \Z_+^d\mbox{ for every large enough $k\ge 0$}\}
\end{displaymath}
and $\1_M$ is the projection on $\RR_M$ along $\KK_M$
of the vector $[1\ 1\ \ldots 1]^t$,
is  a unital ordered group.
If $M$ is primitive, the group $(\Delta_M,\Delta_M^+,\1_M)$ is simple.
\end{proposition}
\begin{proof}
Let
\begin{displaymath}
\Delta=\{(x_n)_{n\ge 0}\mid x_n\in\Z^d \mbox{ for all $n\ge 0$ }, x_{n+1}=Mx_n\mbox{ for every $n$ large enough}\}
\end{displaymath}
and
\begin{displaymath}
\Delta^0=\{(x_n)_{n\ge 0}\in\Delta\mid x_n=0\mbox{ for every $n$ large enough}\}.
\end{displaymath}
We have by definition $G=\Delta/\Delta^0$.
Let $x\in\Delta$. We may assume, by choosing $k$
large enough, that $x_k$ is in $\RR_M$ and $x_{n+1}=Mx_n$
for every $n\ge k$. Since the multiplication by
$M$ is an automorphism of $\RR_M$ there is a unique $y\in \RR_M$
such that $M^ky=x_k$. The map $\pi:x\in\Delta\mapsto y\in\Delta_M$ is a well-defined
group morphism and its kernel is $\Delta^0$.
Thus $\pi$ induces an isomorphism from $G$
onto $\Delta_M$. This proves the first statement.

Assume now that $M$ is nonnegative. Let
\begin{displaymath}
\Delta^+=\{(x_n)_{n\ge 0}\mid x_n\in\Z_+^d\mbox{ for every $n$ large enough}\}.
\end{displaymath}
Since $(G,G^+)$ is an ordered group by Proposition~\ref{propositionDirectLimit}
and since $\pi(\Delta^+)=\Delta_M^+$, the group $(\Delta_M,\Delta_M^+)$
is an ordered group.

Set $u=[1\ 1\ \ldots 1]^t$  and let $v$ be an element of $\Delta_M^+$.
Set $u=u'+w$  with $u'\in\RR_M$ and $w\in \KK_M$.
There is  an $n\ge 1$ such that $nu-v\in \Q^d_+$. Then
$M^k(nu-v)=M^k(nu'-v)$ is  in $\Z_+^d$ for $k$ large enough. Consequently
$nu'-v$ belongs to $\Delta_M^+$. 
Thus $u'$ is an order unit.

Finally, assume that $M$ is primitive. Let $u,v\in\Delta_M^+$ with $u$ nonzero.
There is an integer $k$ such that
$M^k$ is strictly positive and $M^ku,M^kv$ are also strictly positive.
The same argument as above shows that $u$ is an order unit,
proving  the last
statement by Proposition~\ref{propositionSimpleUnit}.
\end{proof}

 The ordered group $(\Delta_M,\Delta_M^+,\1_M)$  is called the
\emph{ordered group of the matrix}\index{subject}{ordered!group!of a matrix} $M$. 

The following result is very useful.
\begin{proposition}\label{propositionDGMatPrim}
  Let  $M$ be a primitive matrix. Then
  \begin{equation}
\label{equation:conepositif}
\Delta^+_M = \{ v \in \Delta_M \mid z\cdot v > 0 \} \cup \{0\}
  \end{equation}
  where $z$ is a positive left eigenvector of $M$ for the dominant eigenvalue.
  \end{proposition}
\begin{proof}
Let $\lambda$ be the dominant eigenvalue of $M$.
  Assume first that $v$ is in $\Delta_M^+\setminus\{0\}$.
  Then $M^kv$ belongs to $\Z_+^d$ for some $k\ge 0$.
  Since $M$ is primitive, we may assume that all
  entries of $M^kv$ are positive.
  By assertion (iii) of Perron-Frobenius (Theorem \ref{theorem:perron}), the vector
  $\lim_{n\to +\infty }\lambda^{-n}M^nv=(tz)v=t(z\cdot v)$ where $t$ and $z$ are respectively a positive right and left eigenvectors
  of $M$ relative to $\lambda$ such that $z\cdot t=1$. This implies that $z\cdot v$ is strictly positive.
  Conversely, if $z$ is strictly positive, the vector $\lim_{n\to +\infty }\lambda^{-n}M^nv$ has all its components
  positive and thus there exists $n\ge 0$ such that $M^kv$ is an element of $\Z_+^d$. This shows that
  $v$ is in $\Delta_M^+$.
  \end{proof}

\begin{example}\label{exampleGolden1}
Consider the primitive matrix $M={\scriptscriptstyle\begin{bmatrix}1&1\\1&0\end{bmatrix}}$.
Since $M$ is invertible, we have $\RR_M=\Q^2$
and $\Delta_M=\Z^2$. Next,  the
dominant eigenvalue of $M$ is $\lambda=(1+\sqrt{5})/2$,
and a corresponding row eigenvector is 
\begin{displaymath}
z=\begin{bmatrix}\lambda&1\end{bmatrix}.
\end{displaymath}
Thus, one has $ \Delta_M^+=\{v\in\Z^2\mid z\cdot v\ge 0\}$.
We recall that $\lambda^2 = \lambda +1$
The map $(\alpha,\beta)\mapsto \lambda\alpha+\beta$
is a positive isomorphism from $(\Delta_M,\Delta_M^+)$
to the group of algebraic integers $\Z[\lambda]=\Z+\lambda\Z$.
The
order unit $\1_M=(1,1)$ is mapped to $\lambda+1$.

This shows that the direct limit of the sequence
$\Z\edge{M}\Z\edge{M}\cdots$ is isomorphic to the group of algebraic integers
$\Z+\lambda \Z$ with the order induced by the reals
and $1+\lambda$ as ordered unit. One can normalize the order
unit to be $1$ as follows. 
The map $x+\lambda y\to (x+\lambda y)/(1+\lambda )$ is an isomorphism of groups which sends
$\Z[\lambda ]$ to itself (because $\lambda^2 = \lambda +1$ is invertible), and 
sends $1+\lambda$ to $1$.  
Thus we find
that the group $(\Delta_M,\Delta_M^+,1_M)$ is isomorphic
to $(\Z[\lambda] , \Z[\lambda] \cap \R_+ , 1)$.
\end{example}
We next give an example with a non primitive matrix.
\begin{example}\label{exampleTriangular}
Consider now $M={\scriptscriptstyle\begin{bmatrix}1&1\\0&1\end{bmatrix}}$.
We have again $\Delta_M=\Z^2$ but this time 
$\Delta_M^+=\{(x,y)\mid y>0\}\cup\{(x,0)\mid x\ge 0\}$.
Thus we find that ordered group of the matrix $M$ is
$\Z^2$ with the lexicographic order.
\end{example}

Note that there can be nonzero vectors $v$ in $\Delta_M$ such that $z\cdot v=0$
and thus that such a vector cannot be in $\Delta_M^+$ (see Example\ref{examplePrimitiveUnimodular}).

Let $(H,H^+,\mathbf{1})$ be a unital ordered group and for every $n\ge 1$,
let
\begin{displaymath}
j_n:(G_n,G_n^+,\mathbf{1}_n)\rightarrow (H,H^+,\mathbf{1})
\end{displaymath}
be a morphism such that $j_{n+1}\circ i_{n+1,n}=j_n$ for every $n$.

The following result will be used later (see Lemma~\ref{lemma1}).

\begin{proposition}\label{propositionMorphismDirectLimit}
There exists a unique morphism of unital ordered groups $j:(G,G^+,\mathbf{1}_G)\rightarrow
(H,H^+,\mathbf{1}_H)$ such that $j\circ i_n=j_n$, for every $n$,
where $i_n:G_n\to G$
is the  natural morphism into $G$ and and $j_n:G_n\to H$ is defined above.
It is surjective if $\cup_n j_n(G_n)=H$ and it is injective if
$\ker(j_n)\subset \ker(i_n)$ for every $n$.
\end{proposition}

\begin{proof}
  Let $g=i_n(g_n)\in G$. Set $j(g)=j_n(g_n)$. Then,
  by the universal property of direct limits, $j$ is a well-defined
morphism from $G$ into $H$. It is a morphism of unital
ordered groups satisfying $j\circ i_n=j_n$. 

Conversely, suppose $j:(G,G^+,\mathbf{1}_G)\rightarrow
(H,H^+,\mathbf{1}_H)$ is a morphism of unital ordered groups such that $j\circ i_n=j_n$
for every $n$.
Let $g$ be an element of $G$. 
If $g=i_n(g_n)$, then $j(g)=j\circ i_n(g_n)=j_n(g_n)$.
This shows the uniqueness.

The last assertions follow easily.
\end{proof}

%%%%%%%%%%%%%%
\section{Dimension groups}\label{sectionDimensionGroups}
A \emph{dimension group}\index{subject}{dimension group} is
a direct limit 
\begin{displaymath}
\Z^{k_1}\edge{M_1}\Z^{k_2}\edge{M_2}\Z^{k_3}\cdots
\end{displaymath}
of groups $\Z^{k_i}$, with $k_i\ge 1$, ordered in the usual way
and with order unit $(1,\ldots,1)$,
with the morphisms defined by nonnegative matrices $M_i$.
Thus a dimension group is a unital ordered group.

The definition of dimension groups  implies some properties
of these groups which hold in any group $\Z^d$ with the natural order.

First of all, a dimension group $(G,G^+ )$ 
is unperforated. Let indeed $g$ be a nonzero element
of $G=\lim\limits_\rightarrow\Z^{k_n}$ and assume that
$ng$ belongs to $G^+$ for some $n>0$. Let $m\ge 1$ be such that
$g=i_m(x)$ for $x\in \Z_+^{k_m}$. Then $nx>0$ implies $x>0$
and thus $g$ belongs to $G^+$. Next, dimension groups satisfy the Riesz interpolation property
that we introduce now.

\subsection{Riesz groups}
An ordered group $G$ satisfies the \emph{Riesz interpolation property}
\index{subject}{Riesz! interpolation property} if for any 
$x_1$, $x_2$ , $y_1$, $y_2$ in $G$
such that $x_1\le y_1,y_2$ and $x_2\le y_1,y_2$, there exists
some $z\in G$ such that $x_1,x_2\le z\le y_1,y_2$.

This property
is equivalent to the \emph{Riesz decomposition property}\index{subject}{Riesz! decomposition property}, requiring that given $x_1,x_2,y_1,y_2\in G^+$
 if $x_1+x_2=y_1+y_2$, then there
exists $z_{ij}\in G^+$ with $1\le i,j\le 2$ such that
$x_i=\sum _jz_{ij}$ and $y_j=\sum_i z_{ij}$ (Exercise~\ref{exerciseRiesz}).

As a variant of the interpolation property, we have that for
every $x,y_1,\ldots,y_k$ in $G^+$ such that
$x\le y_1+\ldots+y_k$, there are $x_1,\ldots,x_k$ in $G^+$
such that $x=x_1+\ldots+x_k$ and $x_i\le y_i$ for $1\le i\le k$
(see Exercise~\ref{exerciseRiesz2}).

An ordered group $G$ is said to be  a \emph{Riesz group}
\index{subject}{Riesz!group}\index{subject}{group!Riesz}%
if it  satisfies the Riesz interpolation
property.

The groups $\Z$ and $\R$ clearly have the Riesz interpolation property
as any totally ordered group. More generally, any lattice ordered
group is a Riesz group.
Next, we have the following more subtle example.

\begin{example}\label{exampleIsRiesz} Any dense subgroup of $\R^2$ is a Riesz group.

\begin{figure}[hbt]
\centering
\tikzset{node/.style={draw,minimum size=0.4cm,inner sep=0pt,fill=gray!10}}
	\tikzset{title/.style={minimum size=0.5cm,inner sep=0pt}}
\begin{tikzpicture}(50,20)
  \node[title]at(-1,1){$x_1$};
  \node[title]at(1,0){$x_2$};
  \draw(0,1)--(6,1);
  \draw(0,2)--(6,2);
  \draw(0,3)--(4,3);
\fill[red](2,1)--(2,2)--(4,2)--(4,1);
\draw(2,0)--(6,0);
\draw(0,1)--(0,3);
\draw(2,0)--(2,3);
\draw(4,0)--(4,3);
\draw(6,0)--(6,2);
\node[title]at(5,3){$y_1$};
\node[title]at(7,2){$y_2$};
%\node[node]() at (15,7.5){};
\end{tikzpicture}
\caption{The Riesz interpolation property}\label{figureRiesz}
\end{figure}
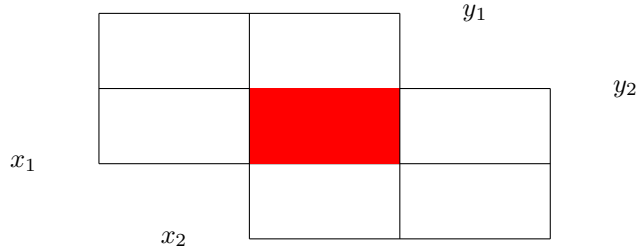

Indeed, let $x_1,x_2$, $y_1,y_2\in \R^2$ with $x_1\le y_1, y_2$
and $x_2\le y_1, y_2$ as in Figure~\ref{figureRiesz}. 
The set of points $z$
such that $x_1,x_2\le z\le y_1,y_2$ is the central rectangle.
%\begin{figure}[hbt]
%\centering
%\begin{picture}(50,20)
%\put(-5,3){$x_1$}\put(5,-2){$x_2$}
%\put(0,10){\framebox(10,5){}}\put(10,10){\framebox(10,5){}}
%\put(0,5){\framebox(10,5){}}\put(10,5){\framebox(10,5){}}\put(20,5){\framebox(10,5){}}
%\put(10,0){\framebox(10,5){}}\put(20,0){\framebox(10,5){}}
%\put(23,17){$y_1$}\put(33,12){$y_2$}
%\node[Nmr=0,Nh=5,Nw=10,fillcolor=red](c)(15,7.5){}
%\end{picture}
%\caption{The Riesz interpolation property}\label{figureRiesz}
%\end{figure}

\end{example}
Every group $\Z^k$ with the natural order is a Riesz group.
This implies that a direct limit of such groups is a Riesz group
and therefore that a dimension group is a Riesz group.

\begin{example}
Let $G$ be the quotient of $\Z^4$ by the subgroup generated by $(1,1,-1,-1)$
and the order induced by the natural order. Denote by $[x]$
the projection on $G$ of $x\in\Z^4$. The group $G$
is not a Riesz group. Indeed, one has, with $x=(1,0,0,0)$,
$y=(0,1,0,0)$, $z=(0,0,1,0)$ and $t=(0,0,0,1)$ the inequality
$[x]\le [z]+[t]$ since $[x]+[y]=[z]+[t]$. However, $[x]$
cannot be written as a sum of two positive smaller elements
and thus the decomposition property fails to hold.
\end{example}
\subsection{The Effros-Handelman-Shen Theorem}
\index{subject}{Effros-Handelman-Shen Theorem}%
\index{subject}{Theorem!Effros-Handelman-Shen}%
The following important theorem
characterizes dimension groups among countable ordered groups.
\begin{theorem}[Effros, Handelman, Shen]\label{theoremEffrosHandelmanShen}
 A countable ordered group is
a dimension group if and only if it is an unperforated directed 
Riesz group.
\end{theorem}
Theorem~\ref{theoremEffrosHandelmanShen} gives
a much easier way to verify that an ordered group is  a dimension group
than using the definition, since it does not require to find
an infinite sequence of morphisms. For example, it shows
directly that any countable dense subroup of $\R^2$
is a dimension group since it is an unperforated Riesz group
(see Example~\ref{exampleIsRiesz}).

The essential step of the proof is the following lemma.
In the proof, we will find it convenient to identify
the group $\Z^n$ with the \emph{free abelian group}
\index{subject}{free!abelian group}\index{subject}{group!free!abelian}%
on a set $A$ with $n$ elements, denoted $\Z(A)$. \index{symbols}{Z@$\Z(A)$}%
The elements of $\Z(A)$ have the form
$\sum_{a\in A}x_a a$ with $x_a\in\Z$.
This
amounts to identify the set $A$ with the canonical basis of $\Z^n$.

The ordered group $(\Z^n,\Z_+^n)$ is then identified
with $(\Z(A),\Z(A)^+)$ where $\Z(A)^+$ is the set
of sums $\sum_{a\in A} x_a a$ with $x_a\ge 0$.

\begin{lemma}\label{lemma1EffrosHandelmanShen}
Let $(G,G^+)$ be an unperforated Riesz group and let $n\ge 1$.
Let $\alpha:(\Z^n,\Z_+^n)\to (G,G^+)$ be a morphism
and let $x\in \Z^n$ be such that $\alpha(x)=0$. There is
an integer $m\ge 1$, a surjective morphism $\eta:(\Z^n,\Z_+^n)\to (\Z^m,\Z_+^m)$
and  a morphism $\beta:(\Z^m,\Z_+^m)\to (G,G^+)$ such that
$\eta(x)=0$ and $\beta\circ \eta=\alpha$ (see the diagram
below).
\begin{figure}[hbt]
\centering
\tikzset{node/.style={minimum size=0.4cm,inner sep=0pt}}
\begin{tikzpicture}
\node[node](A)at(0,2){$(\Z^n,\Z_+^n)$};\node[node](B)at(2,2){$(\Z^m,\Z_+^m)$};
\node[node](G)at(1,0){$(G,G^+)$};

\draw[->,above](A)edge node{$\eta$}(B);
\draw[->,left](A)edge node{$\alpha$}(G);
\draw[->,right](B)edge node{$\beta$}(G);
\end{tikzpicture}
\caption{The diagram of Lemma \ref{lemma1EffrosHandelmanShen}.}\label{diagramShenLemma}
\end{figure}
\end{lemma}
\begin{proof}
We first remark that if we can prove the statement with a morphism
$\eta:(\Z^n,\Z_+^n)\to(\Z^m,\Z_+^m)$ which is not surjective,
we may replace $(\Z^m,\Z_+^m)$ by the ordered group $\eta(\Z^m,\Z_+^m)$,
which is isomorphic to some $(\Z^{m'},\Z_+^{m'})$ (by 
Exercise~\ref{exerciseSubgroupsZ^n}) and thus $\eta$ becomes surjective.

We identify as above $\Z^n$ and $\Z(A)$.
For $x=\sum_{a\in A}x_aa\in\Z(A)$, set $\|x\|=\max\{|x_a|\mid a\in A\}$
and $m(x)=\Card\{a\in A\mid |x_a|=\|x\|\}$.
Let 
\begin{displaymath}
A_+=\{a\in A\mid x_a>0\} \mbox{ and } A_-=\{a\in A\mid x_a<0\}.
\end{displaymath}
We will use an induction on the pairs $(\|x\|,m(x))$ ordered
lexicographically.
Suppose first that $\|x\|=0$. Then $x=0$ and there is nothing to prove.

Assume next that $\|x\|>0$. %Suppose first  that  $A_-$ is empty
%(the case where $A_+$ is empty is solved by changing $x$ into $-x$).
% Let $a_0$ be an element
%of $A_+$. We can write
%\begin{displaymath}
%-\alpha(a_0)=(x_{a_0}-1)\alpha(a_0)+\sum_{a\in A_+,a\ne a_0}x_a\alpha(a).
%\end{displaymath}
%This shows that $-\alpha(a_0)$ is in $G^+$. Since $\alpha$ is
%a positive morphism, we have also $\alpha(a_0)\in G^+$. Thus
%$\alpha(a_0)=0$. We then define a new set $B$ as $B=A\setminus A_+$
%and a morphism $\eta:\Z(A)\to\Z(B)$ by $\eta(\sum_{a\in A}y_aa)=\sum_{a\in B}y_aa$.
%The restriction $\beta$ of $\alpha$ to $\Z(B)$ is a solution.
%We now come to the case where
Since $\alpha(x)=0$, we may assume that $A_+$ and $A_-$ are both nonempty.
Changing if necessary $x$ into $-x$, we may suppose that
$\|x\|=\max \{x_a\mid a\in A_+\}$. Choose $a_0\in A_+$ such
that $x_{a_0}=\|x\|$. Since $\alpha(x)=0$, we have
\begin{displaymath}
x_{a_0}\alpha(a_0)\le\sum_{a\in A_+}x_a\alpha(a)=\sum_{a\in A_-}(-x_a)\alpha(a)\le
x_{a_0}\sum_{a\in A_-}\alpha(a).
\end{displaymath}
Since $G$ is unperforated, we derive that 
$\alpha(a_0)\le \sum_{a\in A_-}\alpha(a)$. Since $G$ is a Riesz group,
there are some $g_a\in G^+$, for each $a\in A_-$,
 such that $\alpha(a_0)=\sum_{a\in A_-}g_a$
with $g_a\le \alpha(a)$ for all $a\in A_-$.

Consider the set $B=(A\setminus \{a_0\})\cup C$ where $C=\{a'\mid a\in A_-\}$
is a copy of $A_-$. We define two positive morphisms $\eta:\Z(A)\to \Z(B)$
 and $\beta:\Z(B)\to G$ by
\begin{displaymath}
\eta(a)=\begin{cases}\sum_{a\in A_-}a'&\mbox{ if $a=a_0$}\\
a&\mbox{ if $a\in A\setminus(A_-\cup\{a_0\})$}\\
a+a'&\mbox{ if $a\in A_-$}
\end{cases}
\end{displaymath}
and
\begin{displaymath}
\beta(b)=\begin{cases}\alpha(b)&\mbox{ if $b\in A\setminus(A_-\cup\{a_0\})$}\\
\alpha(b)-g_b&\mbox{ if $b\in A_-$}\\
g_a&\mbox{ if $b=a'\in C$}.
\end{cases}
\end{displaymath}
It is easy to verify  that $\alpha=\beta\circ\eta$. Next, we claim that
$y=\eta(x)$ is such that
\begin{displaymath}
(\|y\|,m(y))<(\|x\|,m(x)).
\end{displaymath} 
Indeed, we have
\begin{equation}
y=\sum_{a\ne a_0}x_a a+\sum_{a\in A_-}(x_a+x_{a_0})a'.\label{eqShen}
\end{equation}
For every $a\in A_-$, we have $-x_{a_0}\le x_a<0$ and thus
$0\le x_a+x_{a_0}<x_{a_0}$. This shows, by inspection
of the right-hand side of Equation~\eqref{eqShen} that $\|\eta(x)\|\le\|x\|$.
In the case of equality, we clearly have less terms with maximal
absolute value since there is no term $a_0$. Thus $m(y)<m(x)$.

This allows us to apply the induction hypothesis to the
morphism $\beta:(\Z(B),\Z(B)^+)\to (G,G^+)$ and $y\in \Z(B)$.
The solution is a pair of morphisms $\eta':\Z(B)\to\Z(B')$
and $\beta':\Z(B')\to G$ such that $\eta'(y)=0$
 with the diagram of Figure~\ref{diagram2ShenLemma} being commutative. Since $\eta'\circ\eta(x)=0$,
the pair $(\eta'\circ\eta,\beta')$ is a solution. 
\begin{figure}[hbt]
\centering
\tikzset{node/.style={minimum size=0.4cm,inner sep=0.8pt}}
\begin{tikzpicture}
\node[node](A)at(0,2){$(\Z(A),\Z(A)^+)$};\node[node](B)at(3,2){$(\Z(B),\Z(B)^+)$};
\node[node](B')at(6,2){$(\Z(B'),\Z(B')^+)$};
\node[node](G)at(3,0){$(G,G^+)$};

\draw[->,above](A)edge node{$\eta$}(B);\draw[->,above](B)edge node{$\eta'$}(B');
\draw[->,left](A)edge node{$\alpha$}(G);\draw[->,left](B)edge node{$\beta$}(G);
\draw[->,right](B')edge node{$\beta'$}(G);
\end{tikzpicture}
\caption{The induction step in Lemmas~\ref{lemma1EffrosHandelmanShen}
and \ref{lemma2EffrosHandelmanShen}.}\label{diagram2ShenLemma}
\end{figure}
\end{proof}

We prove a second lemma using iteratively the first one.

\begin{lemma}\label{lemma2EffrosHandelmanShen}
Let $(G,G^+)$ be an unperforated Riesz group and let $n\ge 1$.
Let $\alpha:(\Z^n,\Z_+^n)\to (G,G^+)$ be a morphism. 
Then there is $m\ge 1$, a  surjective morphism $\eta:(\Z^n,\Z_+^n)\to (\Z^m,\Z_+^m)$
and a morphism $\beta:(\Z^m,\Z_+^m)\to(G,G^+)$ such that $\ker \eta=\ker\alpha$
and $\beta\circ\eta=\alpha$ (see the diagram~\ref{diagramShenLemma}).
\end{lemma}

\begin{proof}
Since $\ker\alpha$ is a subgroup of $\Z^n$, it is finitely
generated (see Exercise~\ref{exerciseSubgroupsZ^n}).
Let $x_1,x_2,\ldots,x_k$ be a set of generators of $\ker\alpha$.
We proceed by induction on $k$. If $k=0$, there is nothing to prove.
Otherwise, by Lemma~\ref{lemma1EffrosHandelmanShen}, we find an
integer $m\ge 1$, a  surjective morphism $\eta:(\Z^n,\Z_+^n)\to(\Z^m,\Z_+^m)$
and a morphism $\beta:(\Z^m,\Z_+^m)\to (G,G^+)$ such that $\eta(x_1)=0$ and 
$\alpha=\beta\circ\eta$. Then $\ker\beta$ is generated by 
$\eta(x_2),\ldots,\eta(x_k)$. By induction hypothesis, there
is an integer $m'$, a  surjective morphism 
$\eta':(\Z^m,\Z_+^m)\to (\Z^{m'},\Z_+^{m'})$ and a morphism 
$\beta':(\Z^{m'},\Z_+^{m'})\to (G,G^+)$ such that
$\ker\eta'=\ker\beta$. Then $\alpha=\beta'\circ \eta'\circ\eta$
(see Figure~\ref{diagram2ShenLemma}) and $\ker\alpha=\ker\eta'\circ\eta$,
whence the conclusion.
\end{proof}
We prove a third lemma. We will only need one direction
of the equivalence but we state the full result.
\begin{lemma}[Shen]\label{lemma3EffrosHandelmanShen}
A directed  ordered countable group $(G,G^+)$ is a dimension group if and
only if for every morphism $\alpha:(\Z^n,\Z_+^n)\to G$
there is an integer $m\ge 1$ and two morphisms $\eta:(\Z^n,\Z_+^n)\to (\Z^m,\Z_+^m)$
and $\beta:(\Z^m,\Z_+^m)\to (G,G^+)$ 
with $\eta$ surjective such that $\alpha=\beta\circ\eta$
with $\ker\alpha=\ker\eta$.
\end{lemma}
\begin{proof}
Assume first that $(G,G^+)=\lim\limits_\rightarrow(\Z^{k_n},\Z_+^{k_n})$.
Let $\alpha:(\Z(A),\Z_+(A))\to (G,G^+)$ be a morphism.
Choosing $n$ large enough, we can find in $\Z_+^{k_n}$
elements $x_a$ such that $i_n(x_a)=\alpha(a)$ for all $a\in A$. 
Set $\eta(a)=x_a$, $a\in A$. 
Let $u_1,\ldots,u_k$ be a set of generators of $\ker\alpha$. 
Choosing $n$ large enough,
we will have $\eta(u_1)=\ldots=\eta(u_k)=0$
and thus $\ker\alpha=\ker\eta$. Thus the morphisms
$\eta$  and $i_n$ are a solution.

Let us prove the converse.
Since $G$ is countable, $G^+$ is also countable. Let
$S=\{g_0,g_1,\ldots\}$ with $g_n\in G^+$ be a set of generators of $G$ 
and consider a set $A=\{a_0,a_1,\ldots\}$ in bijection with $S$. 
\begin{figure}[hbt]
\centering
\tikzset{node/.style={minimum size=0.4cm,inner sep=0.8pt}}
\begin{tikzpicture}
\node[node](A0)at(0,2){$\Z(A_0)$};\node[node](B0)at(2,2){$\Z(B_0)$};
\node[node](A1)at(4,2){$\Z(A_1)$};\node[node](B1)at(6,2){$\Z(B_1)$};
\node[node](etc)at(7.5,2){$\ \ldots$};
\node[node](G)at(3,0){$G$};

\draw[->,above](A0)edge node{$\eta_0$}(B0);
\draw[->,above](B0)edge node{$\theta_0$}(A1);
\draw[->,left](A0)edge node{$\alpha_0$}(G);
\draw[->,left](B0)edge node{$\beta_0$}(G);
\draw[->,above](A1)edge node{$\eta_1$}(B1);
\draw[->,left](A1)edge node{$\alpha_1$}(G);
\draw[->,left](B1)edge node{$\beta_1$}(G);
\draw[->,above](B1)edge node{}(etc);
\end{tikzpicture}
\caption{The construction of $(A_i,B_i)$.}\label{diagram3ShenLemma}
\end{figure}
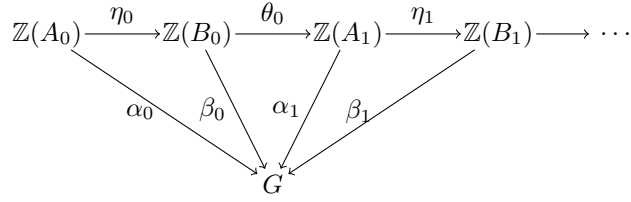

We are going to build a sequence of finite sets $(A_0,B_0,A_1,B_1,\ldots)$
and morphisms $\alpha_n$, $\beta_n$, $\theta_n$ and $\eta_n$ with
$\ker\alpha_n=\ker \eta_n$ and $\eta_n$ surjective,
as in Figure~\ref{diagram3ShenLemma}.
Set $A_0=B_0=\{a_0\}$ and define $\alpha_0(a_0)=\beta_0(a_0)=g_0$
while $\eta_0$ is the identity.

Assume that $A_n,B_n,\beta_n,\eta_n$ are already defined. Set $A_{n+1}=B_n\cup\{a_{n+1}\}$ and let $\theta_n$ be the natural inclusion of
 $\Z(B_n)$ into $\Z(A_{n+1})$. Define $\alpha_{n+1}:\Z(A_{n+1})\to G$ by
\begin{displaymath}
\alpha_{n+1}(a)=\begin{cases}g_{n+1}&\mbox{ if $a=a_{n+1}$}\\
\beta_n(a)&\mbox{ otherwise}
\end{cases} 
\end{displaymath}
From Lemma \ref{lemma2EffrosHandelmanShen} applied to the morphism $\alpha_{n+1}$,
we obtain $B_{n+1}$, $\eta_{n+1}$ and $\beta_{n+1}$. Thus
the iteration can continue indefinitely (unless $S$ is finite,
in which case we stop).

Consider the sequence 
\begin{displaymath}
\Z(B_0)\edge{\gamma_0}\Z(B_1)\edge{\gamma_1}\Z(B_2)\cdots
\end{displaymath}
with $\gamma_n=\eta_{n+1}\circ\theta_n $, which is
obtained by telescoping the top line of Figure~\ref{diagram3ShenLemma}.
Let $(H,H^+)=\lim\limits_\rightarrow (\Z(B_n) , \Z_+ (B_n) ) $ be its direct limit.
In case $S$ is finite, we take for $H$ the last $\Z(B_n)$ instead of the
direct limit.

Let $i_n:\Z(B_n)\to H$ be the natural morphism.
By the universal property of direct
limits, there is a morphism  $h:(H,H^+)\to (G,G^+)$  such that $h\circ i_n=\beta_n$.

The morphism $h$ is injective. 
Indeed, since $\alpha_n=\beta_n\circ\eta_n$ with $\ker\alpha_n=\ker\eta_n$,
$\beta_n$ is injective.

Finally, $h$ is surjective. Indeed, let $g\in S$.
Then $g=g_n$ for some $n$ and thus $g=\alpha_n(a_n)$,
which implies that $g$ is in the image of $h$.
Since $G$ is directed, it is generated by $G^+$ and thus
by $S$. Therefore the image
of $h$ is $G$.
\end{proof}

The proof of Theorem~\ref{theoremEffrosHandelmanShen}
is now reduced to a concluding sentence.

\begin{proofof}{of Theorem~\ref{theoremEffrosHandelmanShen}}
We have already seen that a dimension group is a countable
directed Riesz group. Conversely,  let $G$ be a countable
unperforated directed Riesz group. By 
Lemma~\ref{lemma2EffrosHandelmanShen} 
the condition of Shen's Lemma (Lemma~\ref{lemma3EffrosHandelmanShen})
 is satisfied.
Thus $G$ is a dimension group.
\end{proofof}
Let us illustrate the proof on the example of the
group $G=\Z[1/2]$ of dyadic rationals. The
submonoid $G^+$ is generated by $S=\{1,1/2,1/4,\ldots\}$. We start
with $A_0=B_0=\{a_0\}$ and $\alpha_0(a_0)=1$. Next, we find $A_1=\{a_0,a_1\}$
with $\alpha_1(a_1)=1/2$. Three iterations of the proof of 
Lemma~\ref{lemma2EffrosHandelmanShen} give $B_1=\{b_0\}$ with
$\eta(a_0)=2b_0$ and $\eta(a_1)=b_0$. Continuing in this way
(Exercise~\ref{exerciseChristian}), we obtain $G$ as the direct
limit of the sequence
\begin{displaymath}
\Z\edge{2}\Z\edge{2}\Z\edge{2}\ldots
\end{displaymath}
which was expected (see Example~\ref{exampleDyadic}).

%%%%%%%%%%%%%%%%%%%%%%%%%%
\section{Stationary systems}
Let us consider in more detail the dimension groups obtained
as the direct limit of a sequence
\begin{equation}
\Z^d\edge{M}\Z^d\edge{M}\Z^d\cdots\label{eqStationary}
\end{equation}
where at each step the matrix $M$ is a fixed nonnegative 
integer matrix.
Such a sequence is called a \emph{stationary system}.
\index{subject}{stationary!system}%

Recall that the direct limit of a stationary system has been
characterized in Proposition~\ref{propositionDeltaM} as a triple
$(\Delta_M,\Delta_M^+,\1_M)$ called the ordered group of $M$.
This unital ordered group  is a dimension group
which has properties closely related to algebraic number theory.
Indeed, the largest eigenvalue $\lambda$ of $M$ is an algebraic integer
since it is a root of the polynomial $\det(xI-M)$.
Consequently, all components of the corresponding eigenvector are
in the algebraic field $\Q[\lambda]$.

We already know by Proposition~\ref{propositionDeltaM}
that, wen $M$ is primitive, the group $\G=(\Delta_M,\Delta_M^+,\1_M)$ is simple.
We prove the following additional property.
A subset $C$ of $\Q^n$ is a \emph{cone}\index{subject}{cone}
if $\alpha x\in C$ for every $x\in C$
and $\alpha\in\Q_+$.
\begin{proposition}[Elliott]\label{propositionUniqueState}
A dimension group $\G=(\Delta_M,\Delta_M^+,\1_M)$
with $M$ primitive  has a unique state.
\end{proposition}
\begin{proof}
  We have seen in Proposition~\ref{propositionDGMatPrim}
  that when $M$ is primitive, we have
$\Delta_M^+=\{x\in\Delta_M\mid v\cdot x> 0\}\cup\{0\}$
where $v$ is a row eigenvector of $M$
relative to the maximal eigenvalue. 
We may further assume
that $v\cdot\1_M=1$.  This implies,
by Proposition~\ref{propositionEffros}
that $p:x\to v\cdot x$ is the unique state of $\G$.
Indeed, let $q$ be another state of $\G$. Since every
vector of $\RR_M$ is colinear to a vector of $\Delta_M$,
$p$ and $q$ can be uniquely extended to linear forms on $\RR_M$.
The set $C$ of $x\in\RR_M$ such that $p(x),q(x)$ have opposite
signs is a nonempty cone.
 Every nonempty cone
in $\Q^n$ contains integer points. Thus $C$ contains integer points,
which are consequently in $\Delta_M$, a contradiction.
\end{proof}
\begin{example}\label{exampleChristian}
  Let $M$ be the primitive matrix
  \begin{displaymath}
M=\begin{bmatrix}0&1&1&0\\0&1&0&1\\1&0&1&0\\0&1&1&0\end{bmatrix}.\quad
  \end{displaymath}
  Its eigenvalues are $-1,0,1,2$ and the corresponding eigenvectors are
  \begin{displaymath}
    x=\begin{bmatrix}2\\-1\\-1\\2\end{bmatrix},
      y=\begin{bmatrix}1\\1\\-1\\-1\end{bmatrix},
        z=\begin{bmatrix}0\\1\\-1\\0\end{bmatrix},
          t=\begin{bmatrix}1\\1\\1\\1\end{bmatrix}.
  \end{displaymath}
  The eventual range $\RR_M$ is the space generated by $x,z,t$.
  Thus we have
  \begin{displaymath}
    \Delta_M=\{\alpha x+\beta t+\gamma z
  \mid \alpha,\beta,\gamma\in \Q, M^k(\alpha x+\beta t+\gamma z)\in\Z^4
  \mbox{ for $k$ large enough}\}.
  \end{displaymath}
  The vector
  \begin{displaymath}
    v=\begin{bmatrix}1/6&1/3&1/3&1/6\end{bmatrix}
  \end{displaymath}
  is a left eigenvector corresponding to the eigenvalue $2$ such that
  $v\cdot x=v\cdot z=0$ and
  $v\cdot t=1$. Thus the unique state is $\alpha x+\beta t+\gamma z\mapsto\beta$.

  \end{example}
We will now suppose that $M$ is \emph{unimodular}.
\index{subject}{unimodular!matrix}\index{subject}{matrix!unimodular} This means
that $M$ has determinant $\pm 1$ or, equivalently
that $M$ is an element of the group $GL(n,\Z)$ of integer matrices
with integer inverse. 
The dominating eigenvalue $\lambda$ of $M$
 is then a unit of the ring
$\Z[\lambda]$. Indeed, set $\det(xI-M)=x^d+a_{d-1}x^{d-1}+\ldots+a_1x+a_0$.
Then $\lambda(\lambda^{d-1}+a_{d-1}\lambda^{d-2}+\cdots+a_1)=-a_0$. Since $a_0=\det(M)=\pm 1$,
we conclude that $\lambda$ is invertible.
\begin{proposition}\label{propositionDGPrimitiveUnimodular}
 Let $M\in GL(n,\Z)$ be a primitive  unimodular matrix with dominating eigenvalue $\lambda$.
 The group $\G_M/\Inf(\G_M)$ is isomorphic
to $\Z[\lambda]$ 
with order unit the projection of $u=\begin{bmatrix}1&1&\ldots&1\end{bmatrix}^t$
  in the quotient.
\end{proposition}
\begin{proof}
Let $p(x)$ be the minimal polynomial of $\lambda$. Set
$\det(xI-M)=p(x)q(x)$ and let $E=\ker p(M)$, $F=\ker q(M)$.
The subgroups $E$ and $F$ are invariant by $M$.
Since $p$ is irreducible and $\lambda$ has multiplicity $1$,
$p$ and $q$ are relatively prime.
 This implies that $\Z^n=E\oplus F$. 
Indeed, let $a(x),b(x)$ be such that $a(x)p(x)+b(x)q(x)=1$.
For every $w\in\Z^n$, we have $w=a(M)p(M)w+b(M)q(M)w$.
Since $p(M)q(M)=0$, the first term is in $F$ and the second one in $E$.
Let $v$ be a row eigenvector of $M$ corresponding
to $\lambda$. We have for every $w\in F$
\begin{eqnarray*}
  v\cdot w&=& v\cdot a(M)p(M)w\\
  &=&a(\lambda)p(\lambda)v\cdot w\\
  &=&0 .
\end{eqnarray*}
On the other hand, for every nonzero vector $w\in E$, we have $v\cdot w\ne 0$ since
otherwise $E$ would contain an invariant subgroup, in contradiction with the
fact that $p$ is irreducible.
Thus $F$ is the infinitesimal subgroup and $E$ is isomorphic with $\Z[\lambda]$
via the map $w\mapsto v\cdot w$. 
\end{proof}
\begin{example}\label{examplePrimitiveUnimodular}
  Let $M$ be the primitive matrix
  \begin{displaymath}
    M=\begin{bmatrix}1&1&0\\2&1&1\\0&1&0\end{bmatrix}
  \end{displaymath}
  The eigenvalues of $M$ are $\lambda=(3+\sqrt{5})/2$, $(3-\sqrt{5})/2$ and $-1$
  and thus $M$ is unimodular with dominating eigenvalue $\lambda$.
  We have $p(x)=x^2-3x+1$ and $q(x)=x+1$ and thus $\G_M/\Inf(\G_M)$ is isomorphic
  to $\Z[\lambda]$. The infinitesimal group $\Inf(\G_M)$ is generated by  $[1,-2,2]^t$
  which is an eigenvector for the value $-1$.
  A left eigenvector of $M$ for $\lambda$
  is $v=[2\lambda-2,\lambda,1]$. %In the basis formed of the three vectors
  %\begin{displaymath}
    %\begin{bmatrix}1\\1\\0\end{bmatrix},\begin{bmatrix}0\\1\\1\end{bmatrix},
        %\begin{bmatrix}1\\-2\\2\end{bmatrix},
  %\end{displaymath}
  %the matrix $M$ takes the form
  %\begin{displaymath}
  %\begin{bmatrix}
    %2&1&0\\1&1&0\\0&0&-1\end{bmatrix}
  %\end{displaymath}
  
    Since $v\cdot u=3\lambda-1=\lambda^2$ is a unit of $\Z[\lambda]$, the map $w\mapsto (1/\lambda^2)v\cdot w$
    induces an isomorphism from $\G_M/\Inf(\G_M)$ onto $\Z[\lambda]$ with unit $1$.
  \end{example}
 It is a natural question to ask when
 the  groups $\G_M/\Inf(\G_M)$ is isomorphic, as a unital group,
to $\Z[\lambda]$ with unit $1$.
Actually, this happens if and only if $v\cdot u$ is a unit,
as in the above example.
%%%%%%%%%%%%%%%%%%%%ùùù
\section{Exercises}
\exosection{Section \protect{\ref{sectionOrderedGroups}}}
\begin{exercise}\label{exerciseSubmonoid}
Let $S$ be a submonoid of an abelian group $G$.
Show that the subgroup generated by $S$ is the set
$S-S=\{s-t\mid s,t\in S\}$.
\end{exercise}
\begin{exercise}\label{exerciseSimple}
Show that an ordered group $(G,G^+)$ is simple if and only if 
every nonzero element of $G^+$ is an order unit.
\end{exercise}
\begin{exercise}\label{exerciseSubgroupsZ^n}
Show that a subgroup of $\Z^n$ can be generated by at most $n$ elements.
Hint: use induction on $n$.
\end{exercise}
%\begin{exercise}\label{exerciseMorphismOrderedGroups}
%Let $\G=(G,G_+)$ be an ordered group and let $\varphi:G\to H$
%be a morphism from $G$ onto an abelain group $G$. Set $H_+=\varphi(G_+)$.
%Show that $(H,H_+)$ is an ordered group if and only if
%$\Ker\varphi\cap G_+=\{0\}$.
%\end{exercise}
\exosection{Section \protect{\ref{sectionStates}}}
\begin{exercise}\label{exerciseInfinitesimals1}
  Show that the definition of an infinitesimal element
  in an unperforated ordered group does not depend
on the choice of the order unit.
\end{exercise}
\begin{exercise}\label{exerciseInfinitesimals2}
  Let $G=(G,G^+,u)$ be a unital unperforated ordered group.
Show that $g$ is infinitesimal if and only if $-\varepsilon u\le g\le\varepsilon u$
for every $\varepsilon\in\Q_+$.
\end{exercise}
\begin{exercise}\label{exerciseInfGroup}
Let $\G=(G,G^+,u)$ be a simple unital group.
Denote by $\dot{g}$ the image of $g$ in $G/\Inf(G)$.
Show that $G/\Inf(G)$ has a natural ordering defined by $\dot{g}\ge 0$
if $g+h\ge 0$ for some infinitesimal $h$, that
$G/\Inf(G)$ is simple  and that the infinitesimal subgroup
of $G/\Inf(G)$ is trivial.
\end{exercise}
\begin{exercise}\label{exercisef*}
 Let $\G=(G,G^+,u)$ be a unital ordered group. For $g\in G^+$, set 
\begin{eqnarray*}
f_*(g)&=&\sup\{\varepsilon\in\Q_+\mid\varepsilon u\le g\},\\\label{eqalf_*}
f^*(g)&=&\inf\{\varepsilon\in\Q_+\mid g\le\varepsilon u\}.\label{eqf^*}
\end{eqnarray*}
Show that $\alpha(g)=f_*(g)$ and $\beta(g)=f^*(g)$ where
$\alpha,\beta$ are as in  Lemma~\ref{lemma3.1} with $H=\Z u$.
  \end{exercise}
\begin{exercise}\label{exerciseMiniMax}
  Let $\G=(G,G^+,u)$ be a directed
  unital ordered group. Show that for every $g\in G^+$,
one has the following minimax principle:
\begin{displaymath}
\inf\{p(g)\mid p\in S(\G)\}=\sup\{\varepsilon\in\Q_+\mid \varepsilon u\le g\}.
\end{displaymath}
Hint: use Exercise~\ref{exercisef*}.
\end{exercise}
\begin{exercise}\label{exerciseLemma4.1}
  Let $\G=(G,G^+,u)$ be a directed unital ordered group. Show that
  for every $g\in G^+$ there is a state $p$ on $\G$ such that
  $p(g)=\inf\{\varepsilon\in\Q_+\mid g\le\varepsilon u\}$.
\end{exercise}
\begin{exercise}\label{exerciseLemma6.1}
  Let $\G=(G,G^+,u)$ be an unperforated directed unital group. Show that
  for  $z\in G$, if one has $p(z)\ge 0$ for every
  state $p$ of $\G$ then $z+u\ge 0$.
  \end{exercise}
\exosection{Section \protect{\ref{sectionDirectLimits}}}
\begin{exercise}\label{exerciseAltDefDirectLimits}
Let $(G_n)_{n\ge 0}$ be a sequence of abelian groups with
connecting morphisms $i_{n+1,n}:G_n\to G_{n+1}$. For $m\le n$,
set $i_{n,m}=i_{n,n-1}\circ\cdots\circ i_{m+1,m}$,
with convention that $i_{n,n}$ is the identity.
Let $G$ be the quotient of the disjoint union of the $G_n$ by the equivalence
generated by the pairs $(g,i_{n+1,n}(g))$ for all $g\in G_n$
and all $n\ge 0$. Denote by $[g]$ the class of $g\in \cup G_n$.

Show that $G$ is a group for the
operation
\begin{displaymath}
[g+h]=[g+i_{n,m}(h)]
\end{displaymath}
where $g\in G_n$ and $h\in G_m$ with $m<n$.

Show that $G$ is isomorphic with the direct limit of the 
sequence $(G_n)_{n\ge 0}$.
\end{exercise}
\begin{exercise}\label{exerciseUniversal}
Prove the universal property of direct limits (Equation~\eqref{equationUniversalProperty}).
\end{exercise}
\begin{exercise}\label{exerciseAltDefDirectLimits2}
Let $(\G_n)_{n\ge 0}$ be a sequence of unital ordered groups
$\G_n=(G_n,G_n^+,\1_n)$ with connecting morphisms $i_{n+1,n}$.
Let $G$ be the quotient of the disjoint union
of the union of the $G_n$ by the equivalence
generated by the pairs $(g,i_{n+1,n}(g))$ for all $g\in G_n$
and all $n\ge 0$, as in Exercise~\ref{exerciseAltDefDirectLimits}.
Let $G^+$ be the set of classes of the elements of $\cup G_n^+$
and let $\1$ be the class of $\1_0$.
Show that $(G,G^+,\1)$ is isomorphic to the direct limit of the sequence 
$\G_n$.
\end{exercise}
\begin{exercise}\label{exerciseDeltaM}
Let
\begin{displaymath}
M=\begin{bmatrix}1&1\\1&1\end{bmatrix}
\end{displaymath}
Show that $\Delta_M\sim \Z[1/2]$ and $\Delta_M^+\sim \Z_+[1/2]$.
\end{exercise}
\begin{exercise}\label{exerciseShiftEquivalence1}
Two  integral square matrices $M,N$ are \emph{shift equivalent over $\Z$
with lag $\ell$},%
\index{subject}{shift!equivalence!over $\Z$}
written $M\sim_{\Z} N$, if there are rectangular integral matrices
$R,S$ such that
\begin{equation}
 MR=RN,\qquad SM=NS,\label{eqShiftEquiv1}
\end{equation}
and
\begin{equation}
M^\ell=RS,\qquad N^\ell=SR.\label{eqShiftEquiv2}
\end{equation}
We denote this situation by $(R,S):M\sim_\Z N$ (lag $\ell$).
When $M,N$ are nonnegative, the matrices are \emph{shift equivalent},%
\index{subject}{shift!equivalence}
written $M\sim N$,\index{symbols}{M@$M\sim N$}%
if $R,S$ can be chosen to be nonnegative integral matrices.
We then denote $(R,S):M\sim N$ (lag $\ell$).
\index{symbols}{R@$(R,S):M\sim N$ (lag $\ell$)}%
Show that shift equivalence over $\Z$ (and shift equivalence)
is an equivalence relation.
\end{exercise}
\begin{exercise}\label{exerciseShiftEquivalenceSpectrum}
Show that two matrices which are shift equivalent over $\Z$ have
the same nonzero eigenvalues.
\end{exercise}
\begin{exercise}\label{exerciseShiftEquivalence2}
For a square integral matrix $M$, denote by $\delta_M$ the restriction
of $M$ to $\Delta_M$.
Show that $M\sim_\Z N$ if and only if $(\Delta_M,\delta_M)\simeq (\Delta_N,\delta_N)$ (the latter means that there is a linear isomorphism
$\theta:\Delta_M\to\Delta_N$
 such that $\delta_N\circ\theta=\theta\circ\delta_M$).
\end{exercise}
\begin{exercise}\label{exerciseShiftEquivalence3}
Show that two nonnegative integral matrices are
shift equivalent if and only if $(\Delta_M,\Delta_M^+,\delta_M)\simeq
(\Delta_N,\Delta_N^+,\delta_N)$ (in the sense that there is a linear
isomorphism $\theta:\Delta_M\to \Delta_N$ such that
$\theta(\Delta_M^+)=\Delta_N^+$ and $\theta:\Delta_M\to\delta_N$
is such that $\delta_N\circ\theta=\theta\circ\delta_M$).
\end{exercise}

\exosection{Section \protect{\ref{sectionDimensionGroups}}}
\begin{exercise}\label{exerciseSimplicial}
A submonoid of an abelian group is \emph{simplicial}
\index{subject}{simplicial submonoid}%
\index{subject}{submonoid!simplicial}%
if its positive cone is generated,
as a monoid,
by a finite independent set. Show that
an ordered group is isomorphic to $\Z^n$ with the
natural order if and only if its positive cone is 
simplicial.
\end{exercise}
\begin{exercise}\label{exerciseRiesz}
Show that an ordered abelian group  satisfies the Riesz interpolation property
if and only if it satisfies the Riesz decomposition property.
\end{exercise}
\begin{exercise}\label{exerciseRiesz2}
Let $(G,G^+)$ be an ordered group which satisfies the Riesz 
interpolation property.
Show that for all $x,y_1,\ldots,y_k\in G^+$ such that
$x\le y_1+y_2+\ldots+y_k$, there exist $x_1,x_2,\ldots,x_k\in G^+$
such that $x=x_1+x_2+\ldots+x_k$ and $x_i\le y_i$ for $1\le i\le k$.
\end{exercise}
\begin{exercise}\label{exerciseChristian}
Let $G=\Z[1/2]$ be the group of dyadic rationals. Show that
the proof of Lemma~\ref{lemma3EffrosHandelmanShen}
gives successively the following values for $A_n,B_n$ and the morphisms
$\alpha_n,\beta_n,\eta_n$.
\begin{displaymath}
A_{0}=\{a_{0}\}. \ A_{i}= \{b_{i-1}, a_{i}\} \text{ with } i>0 , B_{i}=\{b_{i}\} 
\end{displaymath}
\begin{displaymath}
\alpha_{i}(b_{i-1})=\frac{1}{2^{i-1}}, \alpha_{i}(a_{i})=\frac{1}{2^{i}}
\end{displaymath}
\begin{displaymath}
\eta_{i}(b_{i-1})=2 b_{i}, \eta_{i}(a_{i})=b_{i}, \quad\beta_{i}(b_{i})=\frac{1}{2^{i}}
\end{displaymath}

\end{exercise}
%%%%%%%%%%%%%%%%%%%%%%%%%
\section{Solutions}
\exosection{Section \protect{\ref{sectionOrderedGroups}}}
\begin{solution}{\ref{exerciseSubmonoid}}
The set $S-S$ contains $0$ and is closed under addition because
\begin{displaymath}
(s-t)+(s'-t')=(s+s')-(t+t').
\end{displaymath}
It is also closed by taking inverses because $-(s-t)=t-s$.
Thus it is a subgroup of $G$. Since any subgroup containing $S$
contains $S-S$, the statement is proved.
\end{solution}
\begin{solution}{\ref{exerciseSimple}}
If $G$ is simple, consider $u\in G^+$.
By Proposition~\ref{propositionIdeal}, the set
$J=[u]-[u]$ is an ideal such that $J\cap G^+=[u]$
Thus $J=G$ and $[u]=G^+$, which implies that $u$ a an order unit.
Conversely, if every nonzero element of $G^+$ is an order unit, 
consider an ideal
$J$ of $G$ not reduced to $0$. Let $u\in J^+$ with $u\ne 0$. Since $u$ 
is an order unit, we have $[u]=G^+$ and thus  $J=J^+-J^+=G^+-G^+=G$.
\end{solution}
\begin{solution}{\ref{exerciseSubgroupsZ^n}}
We use an induction on $n$. For $n=1$, the result is true
since a subgroup of $\Z$ is cyclic. Next, assume the result
true for $n-1$ and consider a subgroup $H$ of $\Z^n$. Let
$\pi:\Z^n\to\Z^{n-1}$ be the projection on the first $n-1$
components. Thus $\pi(x_1,x_2,\ldots,x_n)=(x_1,x_2,\ldots,x_{n-1})$.
By induction hypothesis, the group $\pi(H)$ is generated by
$k$ elements $\pi(h_1),\ldots,\pi(h_k)$ with $k\le n-1$.
On the other hand, $\ker(\pi)$ is isomorphic to $\Z$ and thus $H\cap\ker(\pi)$
is cyclic. Let $h_{k+1}$ be a generator of $\ker(\pi)$. For
every $h\in H$, we have $\pi(h)=\sum_{i=1}^kn_i\pi(h_i)$
for some $n_i\in\Z$.
Since $\pi(h-\sum_{i=1}^kn_i h_i)=0$, we have $h-\sum_{i=1}^kn_i h_i=n_{k+1}h_{k+1}$
and thus $h=\sum_{i=1}^{k+1}n_i h_i$. This shows that $H$ is generated by
$h_1,h_2,\ldots,h_{k+1}$.
\end{solution}
%\begin{solution}{\protect{\ref{exerciseMorphismOrderedGroups}}}
%Assume first that $\Ker\varphi\cap G_+=\{0\}$.
%The set $H_+$ is a submonoid and it generates $H$ since $\varphi$ is surjective%.
%Let $h\in H_+\cap (-H_+)$. Then $h=\varphi(g)=\varphi(-g')$
%for some $g,g'\in G_+$. Since $g+g'\in \Ker\varphi$, we have $g+g'=0$
%and thus $g=g'=0$ which implies $h=0$.

%Conversely, assume that $(H,H_+)$ is an ordered group. If $g\in\Ker\varphi\cap %G_+$, then 
%\end{solution}
\exosection{Section \protect{\ref{sectionStates}}}
\begin{solution}{\protect{\ref{exerciseInfinitesimals1}}}
Let $g$ be an infinitesimal element.
Let $v\in G^+$ be another order unit. By definition, there
is an integer $m$ such that $u\le mv$. Since
$g$ is infinitesimal, we have $mng\le u$ for any $n\in \Z$
and thus $mng\le mv$ which implies $ng\le v$.
\end{solution}
\begin{solution}{\protect{\ref{exerciseInfinitesimals2}}}
Assume first that $g$ is infinitesimal and consider $\varepsilon\in\Q_+$.
Set $\varepsilon=p/q$ with $p,q\ge 0$. We have $qg\le u\le pu$ and thus
$g\le\varepsilon u$. Similarly, $qg\le u$
implies $-pu\le -u\le -qg$ and thus $-\varepsilon u\le g$.

Conversely if $n\ge 0$, $g\le (1/n)u$ implies $ng\le u$ and
if $n\le 0$, $-(1/n)u\le g$ implies also $ng\le u$.
\end{solution}
\begin{solution}{\protect{\ref{exerciseInfGroup}}}
Set $H=G/\Inf(G)$ and $H^+=\{\dot{g}\mid g+h\ge 0\mbox{ for some }h\in\Inf (G)\}$.
If $\dot{g}\in H^+\cap (-H^+)$ we have $g+h\ge 0$ and $ -g-h'\ge 0$
for some $h,h'\in\Inf(G)$. Then $h-h'\ge 0$. Since
$\G$ is simple, $h-h'$ is an order unit although an infinitesimal
and thus $h=h'$. We conclude that $g+h=0$ and thus that $\dot{g}=0$.

If $g+h> 0$, there is since $\G$ is simple an integer $n>0$
such that $u\le n(g+h)$. Then $\dot{u}\le n\dot{g}$ since $nh\in\Inf(G)$.
Thus $(H,H^+)$ is simple.

Let $g\in G$ be such that $\dot{g}\in \Inf(H)$. Then for every $\varepsilon\in \Q_+$,
we have $\dot{g}\le \varepsilon\dot{u}$
and thus $g+h\le \varepsilon u$ for some $h\in\Inf(G)$.
Since $g-\varepsilon u\le g+h\le\varepsilon u$ it follows that $g\le 2\varepsilon u$.
A similar argument shows that $-2\varepsilon u\le g$. Thus
$g\in\Inf(G)$ and finally $\dot{g}=0$.
\end{solution}
\begin{solution}{\ref{exercisef*}}
 Note first that since $0u\le 1g$, we have $f_*(g)\ge 0$. Next,
for every $n\ge 0$ and $m>0$ such that $nu\le mg$, we
have $n/m=p(nu)/m\le\alpha(g)$ showing that $f_*(g)\le\alpha(g)$.
Next, consider $x\in H$ and $m>0$ such that $x\le mg$. Set $x=nu$
with $n\in \Z$. If $n<0$, then $p(x)/m=n/m<0\le f_*(g)$.
Next if $n\ge 0$, we have $p(x)/m=n/m\le f_*(g)$. Thus $\alpha(g)\le f_*(g)$.
We conclude
that $\alpha(g)=f_*(g)$. The proof that $f^*(g)=\beta(g)$ is similar.
  \end{solution}
\begin{solution}{\protect{\ref{exerciseMiniMax}}}
  By Exercise~\ref{exercisef*}, we have $f_*(g)=\alpha(g)$ and $f^*(g)=\beta(g)$
  where $\alpha,\beta$ are as in Lemma \ref{lemma3.1} with $H=\Z u$.
By Lemma~\ref{lemma3.1}, this shows that
\begin{displaymath}
0\le f_*(g)\le f^*(g)<\infty
\end{displaymath}
and that for every state $p$ on $H+\Z g$ one has $f_*(g)\le p(g)\le f^*(g)$.
By Lemma~\ref{lemmaTheorem3.2}, this inequality holds for every
state $p$ on $\G$.

We finally claim that if $f_*(g)\le \gamma\le f^*(g)$ there
is a state $p\in S(\G)$ such that $p(g)=\gamma$. Indeed,
by Lemma~\ref{lemma3.1} (3), there is a state $r$ on $H+\Z g$
such that $r(g)=\gamma$. By Lemma~\ref{lemmaTheorem3.2},
$r$ extends to a state $p$ on $\G$. This proves the claim.

This shows that $\inf\{p(g)\mid p\in S(\G)\}=f_*(g)$.
\end{solution}

\begin{solution}{\ref{exerciseLemma4.1}}
  Let $H=\Z u$. 
  By Lemma~\ref{lemma3.1}, there is a state $q$ on $H+\Z g$
  such that $q(g)=\beta(g)$. By Exercise \ref{exercisef*},
  we have $\beta(g)=f^*(g)=\inf\{\varepsilon\in \Q_+\mid g\le\varepsilon u\}$.
  By Lemma~\ref{lemmaTheorem3.2}, there is a state $p$ of $\G$ which
  extends $q$. Thus $p$ is a state of $\G$ such that $p(g)=\inf\{\varepsilon\in \Q_+\mid g\le\varepsilon u\}$.
\end{solution}
\begin{solution}{\ref{exerciseLemma6.1}}
  Since $G$ is directed, there is $j>0$ such that $z\le ju$. Set $t=ju-z$
  and $f^*(t)=\inf\{\varepsilon\in\Q_+\mid t\le\varepsilon u\}$.
  Since $t\in G^+$, by Exercise~\ref{exerciseLemma4.1}, there
  is a state $p$ on $\G$ such that $p(t)=f^*(t)$.
  Then $j-f^*(t)=p(ju-t)=p(z)\ge 0$ implies $f^*(t)\le j$.
  But since $f^*(t)<j+1$, we can find $k,n>0$ such that $nt\le ku$
  with $k/n<j+1$ and thus $kt\le (j+1)nt$. Finally, we obtain
  \begin{displaymath}
    kju=kz+kt\le kz+(j+1)nt\le kz+(j+1)ku
  \end{displaymath}
  whence $k(z+u)\ge 0$ and the conclusion $z+u\ge 0$ since
  $\G$ is unperforated.
  \end{solution}
\exosection{Section \protect{\ref{sectionDirectLimits}}}
\begin{solution}{\protect{\ref{exerciseAltDefDirectLimits}}}
It is easy to verify that $G$ is a group with neutral element $[0]$
since the operation is the unique operation
on $G$ which extends the operations of the $G_n$.
Consider the map $\pi:\cup_{n\ge 0}G_n\to \Delta$ which sends $g\in G_n$
to $\pi(g)=(g,i_{n+1,n}(g),\ldots)$. Since $\pi^{-1}(\Delta^0)=[0]$,
the map $\pi$
it induces an isomorphism from $G$ onto $\lim\limits_\rightarrow G_n$.
\end{solution}
\begin{solution}{\ref{exerciseUniversal}}
Let $(G_n)$ be a sequence of abelian groups with connecting
morphisms $i_{n+1,n}:G_n\to G_{n+1}$ and let $\alpha_n:G_n\to H$
be for each $n\ge 0$ a morphism such that $\alpha_n=\alpha_{n+1}\circ i_{n+1,n}$.
Let $(g_n)\in\Delta$ be such that $g_{m+1}=i_{m+1,m}(g_m)$
for every $m\ge n$. Then $\alpha_{m+1}(g_{m+1})=\alpha_{m+1}\circ i_{m+1,m}(g_m)=
\alpha_n(g_m)$. Thus there is a morphism $h:\Delta\to H$
such that $h(g_0,g_1,\ldots)=\alpha_m(g_m)$ for all $m\ge n$.
Since $\Delta^0\subset \ker h$, the morphism $h$ induces
a morphism $\varphi:G\to H$ such that $\alpha_n=\varphi\circ i_n$.
\end{solution}
\begin{solution}{\protect{\ref{exerciseAltDefDirectLimits2}}}
We have seen in Exercise~\ref{exerciseAltDefDirectLimits}
that the map $\pi$ which sends $g\in G_n$ to $(g,i_{n+1,n}(g),\ldots)$
induces an isomorphism from $G$ onto the direct limit of the $G_n$.
Since $\pi(g)\in \Delta^+$ if and only if $g\in \cup G_n^+$ 
 and since $\pi(1_0)=(1_n)_{n\ge 0}$, the triple $(G,G^+,\1)$
is a unital ordered group isomorphic to the direct
limit of the $\G_n$.
\end{solution}
\begin{solution}{\protect{\ref{exerciseDeltaM}}}
This follows from
\begin{displaymath}
\RR_M=\{\begin{bmatrix}x\\x\end{bmatrix}\mid x\in\R\}.
\end{displaymath}
\end{solution}
\begin{solution}{\protect{\ref{exerciseShiftEquivalence1}}}
We have to prove the transitivity. Assume that
$(R,S):M\sim_\Z N$ (lag $\ell$) and $(T,U):N\sim_\Z P$ (lag $k$).
Then $(RT,US):M\sim_\Z P$ (lag $k+\ell$). Indeed,
$MRT=RNT=RTP$, $USM=UNS=PUS$ and
\begin{displaymath}
M^{\ell+k}=RSM^k=RN^kS=RTUS,\quad P^{\ell+k}=P^\ell UT=UN^\ell T=USRT
\end{displaymath}
\end{solution}
\begin{solution}{\ref{exerciseShiftEquivalenceSpectrum}}
Assume that $(R,S):M\sim_\Z N$ (lag $\ell$). Let $\lambda$ be a nonzero
eigenvalue of $M$ and let $v\ne 0$ be a corresponding eigenvector.
Set $w=Sv$.
We have
\begin{displaymath}
Rw=RSv=M^\ell v=\lambda^\ell v
\end{displaymath}
and thus $w\ne 0$. Next
\begin{displaymath}
Nw=NSv=SMv=\lambda Sv=\lambda w.
\end{displaymath}
Thus $\lambda$ is an eigenvalue of $N$. The proof that every
nonzero eigenvalue of $N$ is an eigenvalue of $M$ is similar.
\end{solution}
\begin{solution}{\protect{\ref{exerciseShiftEquivalence2}}}
First suppose that $(R,S):M\sim_\Z N$ (lag $\ell$).
 Denote by $m,n$ the sizes of $M,N$. Let $\tilde{R}$ and $\tilde{S}$
be the maps defined 
respectively on $\RR_N$ and $\RR_M$ by $\tilde{R}(v)=Rv$ and $\tilde{S}(w)=Sw$.
It follows from \eqref{eqShiftEquiv1} and \eqref{eqShiftEquiv2}
that $\tilde{R}$ and $\tilde{S}$ are mutually inverse linear isomorphisms
between $\RR_M$ and $\RR_N$. Suppose that $w\in\Delta_N$
and let $k\ge 1$ be such that $N^kw\in\Z^n$. Then, since $R$ is integral,
\begin{displaymath}
M^k(Rw)=RN^k w\in \Z^m
\end{displaymath}
showing that $\tilde{R}(w)\in\Delta_M$. Thus
$\tilde{R}(\Delta_N)\subset \Delta_M$. Similarly $\tilde{S}(\Delta_M)\subset \Delta_N$. By~\eqref{eqShiftEquiv2},
we have the commutative diagrams below and this shows that
$(\Delta_M,\delta_M)\simeq(\Delta_N,\delta_N)$.

\begin{equation}
\begin{CD}
\Delta_N @>{\tilde{R}}>> \Delta_M@>\tilde{S}>>\Delta_N\\
@VV{\delta_N}V               @VV{\delta_M}V @VV{\delta_N}V\\
\Delta_N @>{\tilde{R}}>> \Delta_M@>\tilde{S}>>\Delta_N
\end{CD}\label{DiagramIntertwin}
\end{equation}
Conversely, suppose that $\theta:\Delta_M\to\Delta_N$
is a linear isomorphism such that $\delta_N\circ\theta=\theta\circ\delta_M$.
Since $\theta(\Delta_M)\subset\Delta_N$, there is,  for every $v\in\Z^m$,
a $k\ge 1$
such that $N^k\theta (M^mv)\in \Z^m$. Let
$S$ be the matrix of the linear map $v\mapsto N^k\theta (M^mv)$
from $\R^m$ to $\R^n$. Then for every $v\in \R^m$
\begin{displaymath}
SMv=N^k\theta (M^{m+1}v)=N^{k+1}\theta (M^mv)=NSv
\end{displaymath}
 and thus $SM=NS$. Similarly, there
is an $\ell$ such that $M^\ell\theta^{-1} (N^nw)\in\Z^n$ for every
$w\in \Z^n$. Let $R$ be the matrix of the linear map
$w\mapsto M^\ell\theta^{-1} (N^n w)$. Then we have as above $MR=RN$. 
Since finally $RS=M^\ell\theta^{-1}N^nN^k\theta M^m=M^{k+\ell+m+n}$
we conclude that $M$ and $N$ and shift equivalent over $\Z$.
\end{solution}

\begin{solution}{\protect{\ref{exerciseShiftEquivalence3}}}
Let $(R,S):M\sim N$ (lag $\ell$) with $R,S$ nonnegative. Then the
isomorphism $\tilde{S}:\Delta_M\to \Delta_N$ defined
in the solution of Exercise~\ref{exerciseShiftEquivalence2}
maps $\Delta_M^+$ into $\Delta_n^+$.
and similarly $\tilde{R}(\Delta_N^+)\subset \Delta_M^+$.
But $\tilde{R}\circ\tilde{S}=\delta_M^\ell$ maps $\Delta_M^+$
onto itself, so that $\tilde{S}$ maps
$\Delta_M^+$ onto $\Delta_N^+$. Hence $(\Delta_M,\Delta_M^+,\delta_M)\simeq
(\Delta_N,\Delta_N^+,\delta_N)$.

Conversely, if $(\Delta_M,\Delta_M^+,\delta_M)\simeq
(\Delta_N,\Delta_N^+,\delta_N)$, it is easy to verify that
the matrices $R,S$ defined in the solution of Exercise~\ref{exerciseShiftEquivalence2} are nonnegative.

\end{solution}

\exosection{Section \protect{\ref{sectionDimensionGroups}}}
\begin{solution}{\protect{\ref{exerciseSimplicial}}}
Assume first that $(G,G^+)$ is isomorphic to $\Z^n$ with the usual
order. Let $\alpha:\Z^n\to G$ be an isomorphism such that
$\alpha(\Z_+^n)=G^+$. Then $G^+$ is generated by the
images of the elementary basis vectors, which form an independent set.

Conversely, assume that $S\subset G^+$ is a finite independent set
which generates $G^+$ as a semigroup. Set $S=\{s_1,s_2,\ldots,s_n\}$ and let
$\alpha:\Z^n\to S$ be the linear map sending the  $i$th elementary basis 
vector
to $s_i$. Since $G=G^+-G^+$,
the map $\alpha$ is surjective. It is injective since $S$
is independent. Finally $\alpha(\Z_+^n)=G^+$ and thus $(G,G^+)$
is isomorphic to $(\Z^n,\Z_+^n)$.
\end{solution}
\begin{solution}{\protect{\ref{exerciseRiesz}}}
Assume first that $G$ satisfies the Riesz interpolation property and
consider $x_1,x_2,y_1,y_2\ge 0$ such that $x_1+x_2=y_1+y_2$. Since 
$0\le x_1\le y_1+y_2$, we have $0,x_1-y_2\le x_1,y_1$ and thus
by the interpolation property, there is some $z_{11}$ such that
$0,x_1-y_2\le z_{11}\le x_1,y_1$. Set 
\begin{displaymath}
z_{12}=x_1-z_{11},\ z_{21}=y_1-z_{11},\ z_{22}=y_2-z_{12}.
\end{displaymath}
These elements are all positive and $x_1=z_{11}+z_{12}$, $y_1=z_{11}+z_{21}$, $y_2=z_{12}+z_{22}$.
Finally $z_{21}+z_{22}=y_1-z_{11}+y_2-z_{12}=y_1+y_2-x_1=x_2$.

Conversely, assume that $G$ satisfies the decomposition property. By substraction, it is enough to prove the interpolation property for $0,x\le y_1,y_2$.
We  then have $0\le y_1\le y_1+(y_2-x)=y_2+(y_1-x)$. Let $z\ge 0$ be such that 
$y_1+z=y_2+(y_1-x)$. By the decomposition property, there
are $z_{ij}\ge 0$ such that $y_1=z_{11}+z_{12}$, $y_2=z_{11}+z_{21}$ and
$y_1-x=z_{12}+z_{22}$.
Then $x=z_{11}-z_{22}$ and thus $0,x\le z_{11}\le y_1,y_2$.
\end{solution}
\begin{solution}{\ref{exerciseRiesz2}}
We use an induction on $k$. The property holds trivially for $k=1$.
Next, assume that $x,y_1,y_2,\ldots,y_{k+1}\in G^+$ satisfy the hypothesis.
Since $0,x-y_{k+1}\le x,y_1+\ldots+y_k$ there is by Riesz
interpolation some $z\in G$ such that $0,x-y_{k+1}\le z\le x,y_1+\ldots+y_k$.
By induction hypothesis, there are $x_1,\ldots,x_k$ such that
$z=x_1+\ldots+x_k$ with $x_i\le y_i$ for $1\le i\le k$. But then
$x_1,\ldots,x_k,x_{k+1}$ with $x_{k+1}=x-z$ are a solution.
\end{solution}
\begin{solution}{\ref{exerciseChristian}}

Indeed, we have
\begin{displaymath}
\beta_{i}(\eta_{i}(b_{i-1}))=\alpha_{i}(b_{i-1})=\frac{1}{2^{i-1}}, \quad  \beta_{i}(\eta_{i}(a_{i}))=\alpha_{i}(a_{i})=\frac{1}{2^{i}}
\end{displaymath}
and
\begin{displaymath}
\text{ker}(\eta_{i})= \text{ker}(\alpha_{i})=  \Z( b_{i-1} -2 a_{i})
\end{displaymath}
\end{solution}
%%%%%%%%%%%%%%%%%%%%%%%%%%%%%
\section{Notes}
Ordered algebraic structures are a classical subject of which
ordered groups (and a fortiori ordered abelian groups)
are a particular case.
See~\citep{Fuchs1963} for a general introduction to ordered groups
and \cite{Goodearl1986} for a more detailed presentation
of the ordered groups described in this chapter.
Many authors assume ordered groups to be directed (see for example
\cite{Putnam2018}). This simplifies the presentation but
has the drawback of complicating the definition of subgroups of ordered groups.
We follow the choice of~\cite{GoodearlHandelman1976}.
This does not make any difference in the sequel since dimension groups
are directed.

Note that, by a result of \cite{GlassMadden1984}, the isomorphism
of finitely presented lattice ordered groups is undecidable (see also the notes of Chapter~\ref{chapterBratteliDiagrams}).\index{subject}{undecidability!of isomorphism of ordered groups}
\subsection{States}
The notion of state is closely related with the notion
of trace in an algebra (see~\citep{Davidson1996}).
Theorem~\ref{theoremGoodearlHandelman} is a Hahn-Banach 
\index{subject}{Hahn-Banach Theorem}%
\index{names}{Hahn, Hans}\index{names}{Banach, Stefan}%
type
existence theorem due
to \cite{GoodearlHandelman1976}.
\index{names}{Goodearl, Kenneth R.}%
\index{names}{Handelman, David}%

Proposition~\ref{propositionEffros} is Corollary 4.2
in \citep{Effros1981}
(it is stated there for a simple dimension group but actually
holds in this slightly more general case).
The result appears already in \cite{EffrosHandelmanShen1980}
and variants of it appear in \cite{Goodearl1986}.
%The
%proof uses an argument provided by \cite{Hosseini2020}.
%\index{names}{Hosseini, Maryam}%

The definition of the infinitesimal 
group is from \cite{GiordanoPutnamSkau1995}.
\subsection{Direct limits}
The notion of direct limit is classical and can be formulated
for other categories than groups, in particular
for algebras, as we shall see in Chapter~\ref{chapterBratteli}.

The group $\Delta_M$ defined by Equation~\eqref{eqDeltaM}
is called in~\citep{LindMarcus1995} the dimension group 
 of $M$. It is actually a dimension
group, in the sense of the definition given
in Section~\ref{sectionDimensionGroups}. 
Proposition~\ref{propositionDeltaM} is essentially Theorem 7.5.13 in~\citep{LindMarcus1995}.
\subsection{Dimension groups}
The definition of dimension groups (Section~\ref{sectionDirectLimits}) was introduced by 
G. Elliott in~\citep{Elliott1976}.
Simplicial semigroups 
(Exercise~\ref{exerciseSimplicial})
are taken from~\citep{Elliott1978}.

A group with the Riesz (interpolation or decomposition)
property was called a \emph{Riesz group}\index{subject}{Riesz!group}
 in~\citep{Fuchs1965}. Some authors add the requirement
that the group is unperforated (see \citep{Davidson1996} for example).

The Effros, Handelman and Shen Theorem 
(Theorem \ref{theoremEffrosHandelmanShen})
is from~\citep{EffrosHandelmanShen1980},
(see also the expositions in~\cite[Theorem 3.1]{Effros1981}, \cite[Section IV.7]{Davidson1996} or~\cite[Chapter 8]{Putnam2018}).
Part of the proof is already in~\cite{Shen1979}. In particular,
Lemma~\ref{lemma2EffrosHandelmanShen} is \cite[Theorem 3.1]{Shen1979}
and the argument of Lemma~\ref{lemma3EffrosHandelmanShen}
is already in~\citep{Elliott1978}.

The \emph{unimodular conjecture}\index{subject}{unimodular!conjecture}
proposed by Effros and Shen \citep{EffrosShen1979}
asks whether any dimension group $G$ can be obtained as a direct limit
\begin{displaymath}
\Z^{k}\edge{M_1}\Z^{k}\edge{M_2}\Z^{k}\ldots
\end{displaymath}
where $k\ge 1$ and all $M_n$ are unimodular matrices.
It was proved to hold when $G$ is simple and has one state
\cite{Riedel1981a} but disproved in the general case
\cite{Riedel1981b}.
\index{names}{Riedel, Norbert}%
\subsection{Exercises}
The minimax principle on ordered groups (Exercise~\ref{exerciseMiniMax})
is due to \cite{GoodearlHandelman1976}.

Our treatment of shift equivalence (Exercises~\ref{exerciseShiftEquivalence1},
\ref{exerciseShiftEquivalence2} and \ref{exerciseShiftEquivalence3})
follows~\cite{LindMarcus1995} where the pair $(\Delta_M,\delta_M)$
is called the \emph{dimension pair}\index{subject}{dimension!pair} of $M$
and the triple $(\Delta_M,\Delta^+_M,\delta)$ is called the
\emph{dimension triple}\index{subject}{dimension!triple}
 of $M$. The statement of
Exercise~\ref{exerciseShiftEquivalence3}
can be expressed by
the property that the dimension pair is a complete invariant of shift
equivalence over $\Z$ and the dimension triple a complete invariant
of shift equivalence, a result of \cite{Krieger1980}
(see also~\citep{Krieger1977}). 
\index{names}{Krieger, Wolfgang}%
Since the dimension triple $(\Delta_M,\Delta^+_M,\delta)$
is invariant by shift equivalence,
it is well defined on a shift of finite type $X$
by choosing an edge shift $X_M$ on a graph with adjacency matrix $M$
conjugate to $X$ (see Proposition~\ref{propositionSFTEdgeShift}).
%We call the pair $(\Delta_M,\Delta^+_M)$ the \emph{Krieger group}
%\index{subject}{Krieger group}
%of the shift of finite type.
The group $(\Delta_M,\Delta_M^+)$ can actually by defined
on more general shift spaces (see \citep{Putnam2014}
for  the definition on \emph{Smale spaces}).\index{subject}{Smale space}
\index{names}{Smale, Stephen J.}%

Note that shift equivalence is decidable, \index{subject}{decidability!of shift equivalence}
that is, given two matrices $M,N$, using the notation of Exercise~\ref{exerciseShiftEquivalence1}, one can effectively decide
whether there are 
$U,V$ and $\ell$ such that $(U,V):M\equiv N$ (lag $\ell$)
\citep{KimRoush1988} (see also \citep{KimRoush1979}).
Strong shift equivalence obviously implies shift equivalence.
After remaining many years as a conjecture, known as Williams
Conjecture, the converse was disproved by
\cite{KimRoush1992}, even for primitive matrices \cite{KimRoush1998}.
The proof uses notions introduced in \cite{BoyleKrieger1987}.
Strong shift equivalence is not known to be decidable or not
(see~\citep{Boyle2008}).

%%%%%%%%%%%%%%%%
%  chapter Ordered Cohomology
%%%%%%%%%%%%%%%%%%%%%
\chapter{Ordered cohomology}\label{chapterOrderedCohomology}
In this chapter, we define coboundaries and cohomologous functions
in the space of continuous integer valued functions on a topological
dynamical system.
These terms are used in homological algebra and are the dual
notion of boundaries and homologous elements. We briefly explain these terms
in the elementary setting of graphs.

Assume that $G=(V,E)$ is a graph on a set $V$ of vertices with
a set $E$ of edges. Each edge $e$ has its origin $\alpha(e)$
and its end $\omega(e)$. One defines the \emph{boundary}\index{subject}{boundary}
 operator
$\partial:\Z(E)\rightarrow\Z(V)$ from the free abelian group on
$E$ to the free abelian group on $V$ by $\partial(e)=\omega(e)-\alpha(e)$.
The elements of $\im(\partial)$ are called \emph{boundaries}
\index{subject}{boundary!operator}\index{subject}{operator!boundary}
and the elements of $\ker(\partial)$ are the \emph{cycles}\index{subject}{cycle}.
Elements of $\Z(E)$ equivalent modulo $\ker(\partial)$ are said
to be \emph{homologous}\index{subject}{homologous}.

Identifying $\Z^V$ to ${\rm Hom } (\Z (V) , \Z )$, by duality,
we have a \emph{coboundary operator}\index{subject}{coboundary!operator}
\index{subject}{operator!coboundary}%
\begin{displaymath}
\partial^t:\Z^V\rightarrow \Z^E
\end{displaymath}
(note that it operates in the reverse way). It is such that for $\phi\in\Z^V$

\begin{displaymath}
\partial^t\phi(e)=\phi(\omega(e))-\phi(\alpha(e)).
\end{displaymath}
Indeed, by definition of the dual operator, we have
\begin{displaymath}
\partial^t\phi(e)=\phi(\partial(e))=\phi(\omega(e)-\alpha(e))=\phi(\omega(e))-\phi(\alpha(e)).
\end{displaymath}
The elements of $\im(\partial^t)$ are called \emph{coboundaries}
\index{subject}{coboundary!element}
and the elements of $\ker(\partial^t)$ \emph{cocycles}\index{subject}{cocycle}.
In what follows, we will choose for $V$ a topological space $X$,
use the graph with edges $(x,Tx)$ for a transformation $T$ on $X$
 and use
the topological dual $C(X,\Z)$ instead of $\Z^V$.

The chapter is organized as follows.
In Section~\ref{sectionCoboundaries}, we define the coboundary
operator on continuous functions from a space $X$ to $\R$.
We prove in the next section an important result, due to
Gottschalk and Hedlund, characterizing coboundaries in a minimal
system (Theorem~\ref{theoremGH}). It is used in the next section
(Section~\ref{sectionOrderedGroupDynamicalSystem}) to define
the ordered cohomology group $K^0(X,T)$ of a system $(X,T)$.
In Sections~\ref{sectionFactorMaps}
and \ref{sectionInduced} we relate the group $K^0(X,T)$
to the operations of conjugacy and induction.

In Section~\ref{sectionInvariant}, we introduce
invariant Borel probability measures on topological dynamical systems.
We recall the basic notions of invariant measure and
of ergodic measure. We prove the unique ergodicity
of primitive substitution shifts (Theorem \ref{theoremMichel}).
We relate in Section~\ref{sectionInvariantStates}
coboundaries with invariant Borel probability measures
and  the notion of state. We prove the important result
stating that, for a minimal Cantor system, the states
of the dimension group are in one-to-one correspondence
with invariant Borel probability measures (Theorem~\ref{propositionKerov}).
We finally use this result to give a description
of the dimension groups of Sturmian shifts (Theorem~\ref{theoremDGSturm}).
%%%%%%%%%%%%%%%%%%%%%%%%%%%%ù
\section{Coboundaries}\label{sectionCoboundaries}
Let $(X,T)$ be a topological dynamical system.
We denote by
$C(X,\R)$\index{symbols}{C@$C(X,\R)$}  the group of  real valued  continuous functions on $X$
 and by $C(X,\Z)$ the group of integer valued continuous functions. We
denote by $C(X,\R_+)$ and $C(X,\Z_+)$ the corresponding 
sets of non-negative functions.

As any function with values in a discrete space, an integer valued function $f$
on $X$ is  continuous if and only if it is \emph{locally constant},
\index{subject}{locally constant function}%
\index{subject}{function!locally constant}%
 that is,
for every $x\in X$, there is a neighborhood of $x$ on which
$f$ is constant. When $X$ is a Cantor space, this neighborhood can be chosen
clopen.
When $X$ is a space without non trivial clopen sets, like the torus or any non trivial closed interval of $\R$, then $C(X,\Z)$ consists of constant functions.

Since $X$ is compact, a function $f\in C(X,\Z)$ takes only a finite
number of values. Indeed, the family $(f^{-1}(\{n\} ))_{n\in\Z}$ is a
covering of $X$ by open sets, which has a finite subcover.

For every $f\in C(X,\R)$, we define the \emph{coboundary}
\index{subject}{coboundary!function}\index{subject}{function!coboundary}%
of $f$ as the function
\begin{displaymath}
\partial_T f=f\circ T-f.
\end{displaymath}
\index{symbols}{delta@$\partial_T$}%
Clearly, the map $f\mapsto \partial_T (f)$ is an endomorphism
of the group of both $C(X,\Z)$ and $C(X,\R)$.
%We denote here $\partial_T$ the coboundary operator as
%in \cite{DurandHostSkau1999} instead of $\beta$
%as in~\cite{Host1995}.
 Note that the operator $\partial_T$
is the coboundary operator (as introduced above) related to
the graph with $X$ as set of vertices and edges
from each $x\in X$ to $Tx$.

A function $f\in C(X,\R)$ is a \emph{coboundary}\index{subject}{coboundary} if there is
a function $g\in C(X,\R)$ such that $f=\partial_T g$.
Two functions $f,f'$ are \emph{cohomologous}
\index{subject}{cohomologous functions}\index{subject}{function!cohomologous}%
if $f-f'$ is a coboundary.

Note that if $f\in C(X,\R)$ is a coboundary, then $f\circ T$ is
also a coboundary. Indeed, if $f=g\circ T-g$, then $f\circ T=g\circ T^2-g\circ T=\partial_T(g\circ T)$.

\begin{example}
Let $A=\{a,b\}$, $\delta\in\R$ and $(A^\mathbb{Z} , S)$ be the full shift on $A$. 
The continuous function $f$ defined, for  $x\in A^\Z$ by
\begin{displaymath}
f(x)=\begin{cases}\delta&\text{if $x\in [ab]$,}\\
-\delta&\text{if $x\in [ba]$,}\\
0&\text{otherwise,}\end{cases}
\end{displaymath}
is a coboundary.
 Indeed, it is the coboundary of any function $g$ defined for $x\in A^\Z$ by 
\begin{displaymath}
g(x)=\begin{cases}\alpha&\text{if $x\in [a]$,}\\
\beta&\text{if $x\in [b]$,}\end{cases}
\end{displaymath}
for $\beta-\alpha=\delta$.
\end{example}
We now give a natural example of a real valued continuous function on
a Cantor space.
\begin{example}
Let $X=\{0,1\}^\N$ be the one-sided full shift on $\{0,1 \}$. To every $x\in X$,
we associate the real number
\begin{displaymath}
f(x)=\sum_{n\ge 0}x_n2^{-n-1}
\end{displaymath}
which is the value of $x$ considered as an expansion 
$0.x_0x_1\cdots$ in base $2$. The map $f$ is continuous and
we have
\begin{displaymath}
\partial_S f(x)=\begin{cases}f(x)&\mbox{ if $x\le 1/2$}\\
f(x)-1&\mbox{otherwise}.
\end{cases}
\end{displaymath}
\end{example}
\begin{proposition}\label{propositionRemark2.1}
Let $(X,T)$ be a minimal dynamical system and $f\in C (X, \mathbb{R} )$.
Then

\begin{enumerate}
\item One has $\partial_Tf=0$ if and only if $f$ is constant.
\item If $f\in C(X,\Z)$ is a coboundary, it is the coboundary of some $h\in C(X,\Z)$.
\end{enumerate}
\end{proposition}

\begin{proof}
1. Suppose $\partial_Tf=0$.
For $c\in\R$, the set $Y= f^{-1} (\{ c\})$ is closed. 
Assume that $Y$ is nonempty. 
Since $\partial_T f=0$, the set $Y$ is invariant by $T$.
Hence, $(X, T)$ being minimal, this forces $Y=X$.

2. Assume that $f=\partial_T g$ with $g\in C(X,\R)$. 
Let  $\tau : \R\rightarrow \R/\Z$ be  the natural projection 
onto the torus $\T=\R / \Z$.
 Since $g\circ T-g$ belongs to $C(X,\Z)$, we have $\tau \circ (g\circ T-g)=0$
and thus $\partial_T(\tau\circ g)=\tau\circ g\circ T-\tau\circ g=0$. 
Since $\tau$ is continuous, the same argument as above
 implies that $\tau\circ g$ is constant. Thus there
exists $c\in \R$ such that $h(x)=g(x)-c$ is an integer for all $x\in X$.
Since $\partial_T h=\partial_T g$ we obtain the conclusion.
\end{proof}

For $n> 0$, we set
\begin{displaymath}
f^{(n)}=f+f\circ T+\cdots+f\circ T^{n-1}
\end{displaymath}
with $f^{(0)}=0$.\index{symbols}{f@$f^{(n)}$}
The family $f^{(n)}$ for $n\ge 0$ is called the \emph{cocycle}\index{subject}{cocycle}
associated to $f$. One has for all $n,m\ge 0$, the relation
\begin{equation}
f^{(m+n)}=f^{(m)}+f^{(n)}\circ T^{m}\label{eqcocycle}
\end{equation}
called the \emph{cocycle relation}\index{subject}{cocycle! relation}.
Indeed, we have
\begin{eqnarray*}
f^{(m+n)}&=&(f+f\circ T+\cdots+f\circ T^{m-1})+(f\circ T^m+\cdots+f\circ T^{m+n-1})\\
&=&f^{(m)}+(f+f\circ T+\cdots+f\circ T^{n-1})\circ T^m\\
&=&f^{(m)}+f^{(n)}\circ T^m.
\end{eqnarray*}
The following formula will be used often.

\begin{proposition}
\label{prop:cocyclerelation}
Let $(X,T)$ be a dynamical system and $g\in C (X,\R )$.
If $f=\partial_T g$, then we have for all $n\ge 0$,
\begin{equation}
f^{(n)}=g\circ T^n-g.\label{eq1}
\end{equation}
\end{proposition}
\begin{proof}
We have
\begin{eqnarray*}
f^{(n)}&=&g\circ T-g+g\circ T^2-g\circ T+\cdots+g\circ T^n-g\circ T^{n-1}\\
&=&g\circ T^n-g.
\end{eqnarray*}
\end{proof}
%%%%%%%%%%%%%%%%%%%%%%%%%
\section{Gottschalk and Hedlund Theorem}
\label{sectionGH}
\index{names}{Gottschalk, Walter H.}
\index{names}{Hedlund, Gustav A.}
We first prove the following simple property.
\begin{proposition}
\label{prop:cobordborne}
Let $(X,T)$ be a topological dynamical system.
If $f \in C(X, \R)$ is a coboundary, then the sequence
$(f^{(n)})$ is bounded uniformly.
\end{proposition}
\begin{proof}
Let $f=\partial_T g$ for some $g\in C (X,\R )$, then, using Proposition \ref{prop:cocyclerelation},
one gets $f^{(n)}=g\circ T^n-g$ and thus all  $|f^{(n)}|$ are bounded by $2\sup |g|$.
\end{proof}
Proposition~\ref{prop:cobordborne} is useful to prove that
a continuous function is not a coboundary.
\begin{example}
  Let $X=\{0,1\}^\N$ and let $f$ be the function defined by $f(x)=x_0$.
  Then $f^{(n)}(x)$ is the number of symbols $1$ in $x_{[0,n-1]}$. Since $f$
  is not bounded on $X$, it is not a coboundary.
  \end{example}

The following  consequence will be used in the next section.

\begin{corollary}\label{corollary1}
Let $(X,T)$ be a recurrent topological dynamical system.
If $f\in C(X,\R)$ is a non-negative coboundary, it is identically zero.
\end{corollary}
\begin{proof}
Let $x_0\in X$ be a recurrent point. 
By Proposition~\ref{prop:cobordborne}, 
 the sequence $f^{(n)}(x_0)$ is bounded. Since $x_0$ is recurrent,
 one has $f(T^nx_0)\geq 1/2\sup f$ for infinitely many values of $n$. Thus
$f^{(n)}(x_0)\rightarrow \infty$ as $n\rightarrow\infty$
unless $\sup f=0$.
\end{proof}
The statement is not true if the system is not recurrent, as shown by
the following example.

\begin{example}\label{examplea*b*}
  Let $X$ be the shift space such that $\cL(X)=a^*b^*$.
\begin{figure}[hbt]
    \centering
\tikzset{node/.style={circle,draw,minimum size=0.4cm,inner sep=0pt}}
	\tikzset{title/.style={minimum size=0.5cm,inner sep=0pt}}
        \tikzstyle{every loop}=[->,shorten >=1pt,looseness=12]
        \tikzstyle{loop left}=[in=130,out=220,loop]
        \tikzstyle{loop right}=[in=40,out=-40,loop]
	\begin{tikzpicture}
          \node[node](1)at(0,0){$1$};\node[node](2)at(2,0){$1$};

          \draw[left,->](1)edge[loop left]node{$a$}(1);
          \draw[above,->](1)edge node{$b$}(2);
          \draw[right,->](2)edge[loop right]node{$b$}(2);
        \end{tikzpicture}
        \caption{The shift $X$.}\label{figurea*b*}
    \end{figure}

  Thus $X$
  is the set of labels of two-sided infinite paths in the graph
  of Figure~\ref{figurea*b*}. Let $f$ be the characteristic function
  of the set of sequences $x\in X$ such that $x_0=a$. Then
  $f\circ S$ is the characteristic function of the set of sequences $x\in X$
  such that $x_1=a$. Thus $-\partial f$ is the characteristic
  function of the set reduced to the sequence $x\in X$ such that $x_0=a$
  and $x_1=b$. This function is therefore a non-negative coboundary
  which is not zero.

  \end{example}

In the minimal case there is  a converse of Proposition \ref{prop:cobordborne}
 by the following result.
We will use it several times.
\begin{theorem}[Gottschalk, Hedlund]\label{theoremGH}
Let $(X,T)$ be a  minimal topological dynamical system and $f$ be in $C(X,\R)$.
The following are equivalent.
\begin{enumerate}
\item
\label{enum:cob}
$f$ is a coboundary.
\item
\label{enum:cobbound}
The sequence $(f^{(n)})$ is bounded uniformly.
\item
\label{enum:cobboundpoint}
There exists $x_0\in X$ such that the sequence
$(f^{(n)}(x_0))_{n\ge 0}$ is bounded.
\end{enumerate}
\end{theorem}
\begin{proof}
\noindent 
We have already seen that \eqref{enum:cob} implies \eqref{enum:cobbound}.
Assertion \eqref{enum:cobbound} clearly implies \eqref{enum:cobboundpoint}.

%To prove that \eqref{enum:cobboundpoint} implies \eqref{enum:cob}, first
%observe that, by minimality,  $f^{(n)}(x)$ is bounded for all $x$. 
%Indeed, let $M=\sup_{n\ge 0}|f^{(n)(x_0)|$.

%set
%\begin{equation}
%g(x)=\sup_{n\ge 0}f^{(n)}(x)\label{eqKatok}
%\end{equation}

The proof that \eqref{enum:cobboundpoint} implies \eqref{enum:cob} is in three steps.
\paragraph{Step 1} For each clopen neighborhood $U$ of $x_0$, we set
\[
\Lambda(U)
=
\overline{ \{ f^{(n)}(x_0) \mid n\ge 0 , T^nx_0\in U \} } \:  \hbox{\rm and }
\Lambda = \bigcap_U \Lambda (U)
\]
$U$ running over all the clopen neighborhoods of $x_0$.
These sets are bounded and contain $0$. 
We claim that $\Lambda=\{0\}$. To prove it, we show that
for any $a\in\Lambda$, we have $2a\in\Lambda$, which will imply $a=0$
since $\Lambda$ is bounded. 
Suppose indeed that
$a$ belongs to $\Lambda$. 
Let $U$ be a clopen neighborhood of $x_0$ and $\varepsilon>0$.
Since $a$ belongs to $\Lambda(U)$, there exists $n\ge 0$ such
that $T^nx_0\in U$ and $|f^{(n)}(x_0)-a|<\varepsilon$.
As the maps $T^n$ and $f^{(n)}$ are continuous, there is
a clopen neighborhood $V\subset U$ of $x_0$ such that
$T^ny$ is in $U$ and $|f^{(n)}(y)-a|<2\varepsilon$ for all $y\in V$.
Since $a$ belongs to $\Lambda(V)$, there exists $m\ge 0$ with
$T^mx_0\in V$
and $|f^{(m)}(x_0)-a|<\varepsilon$.
As $T^mx_0$ belongs to $V$, we have $T^{n+m}x_0\in U$, $f^{(n+m)}(x_0)\in \Lambda (U)$  and
$|f^{(n)}T^m(x_0)-a|<2\varepsilon$.
By the cocycle relation~\eqref{eqcocycle}, we obtain
\begin{eqnarray*}
|f^{(n+m)}(x_0)-2a|&\le&|f^{(n)}(T^mx_0)-a|+|f^{(m)}(x_0)-a|\\
&<&3\varepsilon.
\end{eqnarray*}
Since $\varepsilon$ is arbitrary, this implies that $2a$ belongs to $\Lambda(U)$
and $U$ being an arbitrary neighborhood of $x_0$, it
implies that $2a$ also belongs to $\Lambda$, which proves the claim.
\paragraph{Step 2}
 Because each $\Lambda(U)$ is compact and $\Lambda(U\cap U')
\subset \Lambda(U)\cap\Lambda(U')$, for every $\varepsilon >0$, there 
exists a neighborhood $U_\varepsilon$ of $x_0$ such that
$\Lambda(U_\varepsilon)\subset[-\varepsilon,\varepsilon]$.

 We claim that there exists a function $g\in C(X,\R)$ such that
$g(T^n(x_0))=f^{(n)}(x_0)$ for all $n\ge 0$. For this,
it is enough to prove that for every $x\in X$
and every sequence $n_i\rightarrow \infty$
of integers such that $T^{n_i}x_0\rightarrow x$, the sequence
$f^{(n_i)}(x_0)$ converges. We then define $g(x)$ as this limit.
Fix $\varepsilon>0$. By minimality, there exists $n\ge 0$ such that
$T^nx_0$ is in $U_\varepsilon$. 

Let $W$ be a neighborhood of
$x$ such that $T^ny\in U_\varepsilon$ and $|f^{(n)}(y)-f^{(n)}(x)|<\varepsilon$
for all $y\in W$. For $i$ large enough, $y=T^{n_i}x_0$ is in $W$.
Then $T^ny=T^{n_i+n}x_0$ belongs to $U_\varepsilon$ and consequently
$|f^{(n_i+n)}(x_0)|\le\varepsilon$. Moreover, we
have $|f^{(n)}(y)-f^{(n)}(x)|=|f^{(n)}(T^{n_i}x_0)-f^{(n)}(x)|<\varepsilon$.
Thus
\begin{displaymath}
|f^{(n_i)}(x_0)+f^{(n)}(x)|\le \varepsilon+|f^{(n_i)}(x_0)+f^{(n)}(T^{n_i}x_0)|=\varepsilon+|f^{(n_i+n)}(x_0)|\le 2\varepsilon.
\end{displaymath}
For large enough $n_i,n_j$, we obtain
\begin{displaymath}
|f^{(n_i)}(x_0)-f^{(n_j)}(x_0)|\le |f^{(n_i)}(x_0)-f^{(n)}(x)|+|f^{(n)}(x)-f^{(n_j)}(x_0)|\le 4\varepsilon.
\end{displaymath}
It follows that $(f^{(n_i)}(x_0))$ is a Cauchy sequence and converges.
\paragraph{Step 3} Since $(X,T)$ is minimal, the orbit of $x_0$
is dense. For any $x=T^nx_0$ in this set, we have by construction
\begin{eqnarray*}
\partial_T g(x)&=&g\circ T^{n+1}(x_0)-g\circ T^n(x_0)=f^{(n+1)}(x_0)-f^{(n)}(x_0)\\
&=&f\circ T^n(x_0)=f(x)
\end{eqnarray*}
By continuity, this extends to any $x\in X$ and this proves that
$f$ is a coboundary.
\end{proof}
It can be shown that the function $g$ such that $f=\partial g$ is
determined uniquely up to a constant (Exercise~\ref{exerciseGH}).

The following example illustrates the fact that Theorem~\ref{theoremGH}
is false without the hypothesis of minimality.

\begin{example}\label{examplea*b*2}
  Let $X$ be the golden mean shift, which is the set of labels
  of infinite paths in the graph of Figure~\ref{figureGoldenMeanShift}.
  It is recurrent but not minimal since it contains $a^\Z$.
\begin{figure}[hbt]
    \centering
\tikzset{node/.style={circle,draw,minimum size=0.4cm,inner sep=0pt}}
	\tikzset{title/.style={minimum size=0.5cm,inner sep=0pt}}
        \tikzstyle{every loop}=[->,shorten >=1pt,looseness=12]
        \tikzstyle{loop left}=[in=130,out=220,loop]
        \tikzstyle{loop right}=[in=40,out=-40,loop]
	\begin{tikzpicture}
          \node[node](1)at(0,0){$1$};\node[node](2)at(2,0){$2$};

          \draw[left,->](1)edge[loop left]node{$a$}(1);
          \draw[above,->,bend left](1)edge node{$b$}(2);
          \draw[below,->,bend left](2)edge node{$a$}(1);
        \end{tikzpicture}
        \caption{The golden mean shift.}\label{figureGoldenMeanShift}
    \end{figure}

  Let $f$ be the characteristic
  function of the set of points $x\in X$ such that $x_0=b$ and let
  $y=a^\infty$. Then $f^{(n)}(y)=0$ for all $n\ge 0$ although
  $f$ is not a coboundary.
  \end{example}

We prove the following additional result.
The proof uses the Baire Category Theorem
\index{subject}{Baire Category Theorem}%
\index{subject}{Theorem!Baire Category}%
(see Appendix~\ref{appendixTopo}).
\begin{proposition}\label{lemma2}
 Let $(X,T)$ be a minimal topological dynamical system.
The following conditions are equivalent for $f\in C(X,\Z)$.
\begin{itemize}
\item[\rm (i)] There exists $g\in C(X,\Z)$ such that $f+\partial_T g\ge 0$.
\item[\rm (ii)] There exists $g\in C(X,\R)$ such that $f+\partial_T g\ge 0$.
\item[\rm (iii)] The family of functions $(f^{(n)})_{n\ge 0}$ is uniformly
bounded from below.
\item[\rm (iv)] For every $x\in X$,
the family of numbers $(f^{(n)}(x))_{n\ge 0}$ is 
bounded from below.
\end{itemize}
\end{proposition}
\begin{proof}
(i) $\Rightarrow$ (ii) is obvious.

(ii) $\Rightarrow$ (iii) If $f+\partial_T g\ge 0$, then
$(f+\partial_T g)^{(n)}=f^{(n)}+(\partial_T g)^{(n)}\ge 0$.
Hence, by Equation~\eqref{eq1},
  $f^{(n)}+g\circ T^n-g\ge 0$ and thus
$f^{(n)}\ge g-g\circ T^n$ is bounded from below.

(iii) $\Rightarrow$ (iv) is trivial.

(iv) $\Rightarrow$ (i).
For each $n\ge 0$, let $g_n$ and $g$ be defined by $g_n (x)  =\inf\{f^{(k)} (x)\mid 0\le k\le n\}$ and $g (x)=\inf\{f^{(n)} (x) \mid n\ge 0\}=\inf\{g_n (x)\mid n\ge 0\}$.
For each $n\ge 0$ and $k\ge 1$, we obtain $g_{n+k} (x)=\inf\{g_{k-1} (x),f^{(k)} (x)+g_n \circ T^k (x)\}$ and $g (x)=\inf\{g_{k-1} (x),f^{(k)}(x)+g\circ T^k (x)\}$. In particular, $g(x)=\inf\{0,f (x)+g\circ T (x)\}\le f (x)+g\circ T(x)$, which implies
$f+\partial_Tg\ge 0$. Because each $g_n$ belongs to $C(X,\Z)$, it is sufficient
to prove that $g=g_n$ for some $n\ge 0$.

For each $n\ge 0$, let $K_n=\{x\in X\mid g(x)=g_n(x)\}$. For all $x\in X$,
as each $g_n$ takes only integer values, there exists an $n\ge 0$ such
that $g(x)=g_n(x)$ and thus $x\in K_n$. Consequently,
$X=\cup_{n\ge 0}K_n$. Since each $K_n$ is closed, by 
Baire Category Theorem,
there exists $m\ge 0$ such that $K_m$ has a nonempty interior. By minimality
and compactness, there exists $p\ge 1$ such that $\cup_{1\le k\le p}T^{-k}K_m=X$.
Let $x\in X$ and $1\le k\le p$ be such that $T^kx\in K_m$. We obtain
\begin{displaymath}
g(x)=\inf\{g_{k-1}(x),f^{(k)}+g\circ T^k(x)\}=
\inf\{g_{k-1},f^{(k)}(x)+g_m\circ T^k(x)\}=g_{k+m}(x)
\end{displaymath}
and thus $x\in K_{k+m}\subset K_{p+m}$. We conclude that $K_{p+m}=X$
and $g=g_{p+m}$.
\end{proof}

%%%%%%%%%%%%%%%%%%%%%
%%%%%%%%%%%%
\section{Ordered cohomology group of a dynamical system}\label{sectionOrderedGroupDynamicalSystem}
Let $(X,T)$ be a topological dynamical system. 
Since the coboundary operator on $C(X,\Z)$ is a group morphism,
its image $\partial_T C(X,\Z)$ is a subgroup.
We denote by $H(X,T,\Z)$ the quotient group
\begin{displaymath}
H(X,T,\Z)=C(X,\Z)/\partial_T C(X,\Z).
\end{displaymath}
\index{symbols}{H@$H(X,T,\Z)$}%
and by $H^+(X,T,\Z)$ the image of $C(X,\Z_+)$ in this quotient.

We denote by $K^0(X,T)$\index{symbols}{K@$K^0(X,T)$} the triple
\begin{displaymath}
K^0(X,T)=(H(X,T,\Z),H^+(X,T,\Z), \mathbf{1}_X).
\end{displaymath}
where $\mathbf{1}_X$ is the image in $H(X,T,\Z)$ of the constant
function with value $1$ on $X$.

The proof of the following result uses Gottschalk and Hedlund Theorem 
(more precisely Corollary~\ref{corollary1}).
\begin{proposition}\label{propositionOrderedCohomologyGroup}
For any recurrent topological dynamical system $(X,T)$,
the triple $K^0(X,T)$ is a unital ordered group.
\end{proposition}
\begin{proof}
%Clearly $H^+(X,T,\Z)$ is a submonoid which generates
%$H(X,T,\Z)$ since both properties hold for $C(X,\Z_+)$.
Assume that $f,f'\in C(X,\Z_+)$ are such that
$f$ is cohomologous to $-f'$. Then $f+f'$ is a non-negative
coboundary.
By Corollary \ref{corollary1}, we have $f+f'=0$,
which implies $f=f'=0$. Thus $H^+(X,T,\Z)\cap -H^+(X,T,\Z)=\{0\}$.
Consequently $K^0(X,T)$ is an ordered group and $\1_X$
is clearly  a unit.
\end{proof}
Following the term
introduced by Boyle and Handelman,
the group $K^0(X,T)$ is called the \emph{ordered cohomology group} of the 
topological dynamical system
\index{subject}{ordered!cohomology group}
\index{subject}{dynamical system!ordered cohomology group of}
$(X,T)$.

Proposition~\ref{propositionOrderedCohomologyGroup} is not
true without the hypothesis that $X$ is recurrent, as illustrated
by the following example.
\begin{example}\label{examplea*b*3}
  Let $X$ be the shift space such that $\cL(X)=a^*b^*c^*$. It is
  the set of labels of infinite paths in the graph
  of Figure~\ref{figurea*b*c*}.
  \begin{figure}[hbt]
    \centering
\tikzset{node/.style={circle,draw,minimum size=0.4cm,inner sep=0pt}}
	\tikzset{title/.style={minimum size=0.5cm,inner sep=0pt}}
        \tikzstyle{every loop}=[->,shorten >=1pt,looseness=12]
        \tikzstyle{loop left}=[in=130,out=220,loop]
        \tikzstyle{loop right}=[in=40,out=-40,loop]
        \tikzstyle{loop above}=[in=60,out=120,loop]
	\begin{tikzpicture}
          \node[node](1)at(0,0){$1$};\node[node](2)at(2,0){$2$};
          \node[node](3)at(4,0){$3$};

          \draw[left,->](1)edge[loop left]node{$a$}(1);
          \draw[above,->](1)edge node{$b$}(2);
          \draw[above,->](2)edge[loop above] node{$b$}(2);
          \draw[above,->](2)edge node{$c$}(3);
          \draw[right,->](3)edge[loop right] node{$c$}(3);
          \draw[below,->,bend right](1)edge node{$c$}(3);
        \end{tikzpicture}
        \caption{The shift $X$.}\label{figurea*b*c*}
  \end{figure}
  
  We have, denoting by $\charac_U$ the characteristic function of
  the set $U$ and by $[u]$ the cylinder $[u]=\{x\in X\mid u=x_{[0,|u|-1]}\}$,
  \begin{displaymath}
    \partial\charac_{[a]}=\charac_{[aa]}-\charac_{[a]}=-\charac_{[ab]}-\charac_{[ac]}.
  \end{displaymath}
  Thus $\charac_{[ab]}+\charac_{[ac]}$ is a  coboundary. It can
  be verified that $\charac_{[ab]}$ is not a coboundary (Exercise~\ref{exercise[ab]}).
  Thus, the image in $G=H(X,S,\Z)$ of $\charac_{[ab]}$ is a nonzero element of
   $G^+\cap -G^+$ and thus $G$ is not an ordered group.
\end{example}

We give now a very simple example of computation of $K^0(X,T)$.
\begin{example}
Let $(X,S)$ be the shift space formed of the two
infinite sequences $x=(\cdots abab.abab \cdots )$ and $y=(\cdots baba.baba \cdots )$. Obviously,
$Sx=y$ and $Sy=x$. 
The characteristic functions of $x$
and $y$ are also exchanged by $S$. Thus $H(X,S,\Z)=\Z$
and $K^0(X,S)=(\Z,\Z_+,1)$.
\end{example}
We can generalize the last example by considering the case of a finite
set $X$ and a permutation $T$ on $X$. Then $H(X,T,\Z)=\Z^d$
where $d$ is the number of orbits of the permutation $T$.

We now give an elementary argument to show that  $H(X,S,\Z)$ is isomorphic to  $\Z^2$ when $(X,S)$ is a Sturmian shift
like the Fibonacci shift. 
\index{subject}{Sturmian!shift space}\index{subject}{shift space!Sturmian}
%For the positive cone,
We will see 
%soon (Section~\ref{sectionInvariant}) how to compute the
%positive cone $H(X,T,\Z_+)$ and 
later (using return words, in Section~\ref{chapterReturnWords})
how  this can be done by more general methods.
\begin{example}\label{exampleSturmian1}
Let $X$ be a Sturmian shift on $A=\{a,b\}$.
We denote by $\charac_{[w]}$
\index{symbols}{char@$\charac_{[w]}$}%
the characteristic function
\index{subject}{characteristic!function} of the cylinder set $[w]$.

Set $A=\{a,b\}$. We show by induction on $|w|$ that 
$\charac_{[w]}$ is cohomologous to an element of
 the subgroup $G$ generated by $\charac_{[a]}$ and $\charac_{[b]}$.
This is true if $|w|=1$. Assume that it holds
for words of length $n$ and consider a word $w$ of length $n+1$.
Then $w=ux$  for some nonempty word $u$ of length $n$
and some letter $x\in\{a,b\}$.
If $u$ is not right-special, we have $uA^\N=uxA^\N$
and thus $\charac_{[ux]}$ is cohomologous to an element
of $G$ by induction hypothesis. Otherwise, we have
\begin{equation}
\charac_{[u]}=\charac_{[ua]}+\charac_{[ub]}.
\label{eqRightSpecial}
\end{equation} 
As $w$ is either $ua$ or $ub$, it suffices to show that $\charac_{[ua]}$ and
 $\charac_{[ub]}$ are cohomologuous to some elements in $G$.
Set $u=yv$
with $y$ a letter and $v$ a word. Then $va$ and $vb$ cannot
be both left-special. Assume that $va$ is not left-special.
Then $\charac_{[ua]} = \charac_{[yva]} = \charac_{[va]}\circ T$.
By the induction hypothesis, $\charac_{[va]}$ is
cohomologous to an element of $G$
and thus also the map $\charac_{[ua]}$.
By Equation~\eqref{eqRightSpecial}, this implies that $\charac_{[ub]}$
is also cohomologous to an element of $G$.

Since $C(X,\Z)$ is generated by the functions
$\charac_{[w]}$ (see Proposition~\ref{lemma4} below), this shows that $H(X,S , \Z)$ is generated by the projections
of  $\charac_{[a]}$ and $\charac_{[b]}$. 
Thus the morphism sending $(\alpha,\beta)$
to the class of $\alpha \charac_{[a]}+\beta \charac_{[b]}$ is surjective from $\Z^2$
to $H(X,S, \Z)$. It is injective because if $f=\alpha \charac_{[a]}+\beta \charac_{[b]}$
with $\alpha$ or $\beta$ nonzero,
then $f^{(n)}(x)$ is not bounded because otherwise
the slope of $x$ would be a rational number (by Corollary~\ref{corollarySlope}).
Thus $f$ is not a coboundary
by Theorem~\ref{theoremGH}.
\end{example}
We may also consider the group $H(X,T,\R)=C(X,\R)/\partial_TC(X,\R)$
and define $H^+(X,T,\R)$ as the image of $C(X,\R)$ in $H(X,T,\R)$.
By Corollary~\ref{corollary1}, the pair $(H,X,T,\R),H^+(X,T,\R))$ is
an ordered group. We observe the following relation between $H^+(X,T,\Z)$
and $H^+(X,T,\R)$.
\begin{corollary}\label{corollary2}
$H(X,T,\Z)$ is a subgroup of $H(X,T,\R)$ and
\begin{equation}
H^+(X,T,\Z)=H(X,T,\Z)\cap H^+(X,T,\R).\label{equationRealValued}
\end{equation}
\end{corollary}
\begin{proof}
The group $C(X,\Z)$ is included in $C(X,\R)$ and if
$f\in C(X,\Z)$ is in $\partial_TC(X,\R)$, then
by Proposition~\ref{propositionRemark2.1}, it is in $\partial_TC(X,\Z)$.
Thus the inclusion of $C(X,\Z)$ in $C(X,\R)$ defines an injection
of $H(X,T,\Z)$ in $H(X,T,\R)$.

Next the left side of~\eqref{equationRealValued} is clearly
included in the right side.
Suppose conversely
 that $f=\partial_Tg+h$ with $f\in C(X,\Z)$, $g\in C(X,\R)$ and
$h\in C(X,\R_+)$. For every $n\ge 0$, we have by Equation~\eqref{eq1}
\begin{displaymath}
f^{(n)}(x)=g(T^nx)-g(x)+h^{(n)}(x)
\end{displaymath}
and the family $(f^{(n)})_{n\ge 0}$ is uniformly bounded from below. By
Proposition~\ref{lemma2}, $f$ is cohomologous in $C(X,\Z)$ to some
$f'\in C(X,\Z_+)$, which proves~\eqref{equationRealValued}.
\end{proof}
The next proposition characterizes functions which are cohomologous
to integer valued functions.
Before we need a classical lemma.
We provide a proof for sake of completeness. 
 Let $\tau:\R\rightarrow \R / \Z$ be the projection from $\R$
onto the torus $\mathbb{T}=\R/\Z$.

\begin{lemma}\label{lemmaTorus}
Let $X$ be a Cantor space and $f : X \to \mathbb{T}$ 
%\marginpar{FD : Tore et S1} 
be a continuous map. 
For every $\epsilon > 0$, there exists a
continuous function $h_\epsilon : X \to [ -\epsilon ,  1]$ such that $f = \tau\circ h_\epsilon$.
\end{lemma}

\begin{proof}
Let $\epsilon >0$.
For each $x \in X$, there exist a neighborhood $U_x$ and a continuous map   $ E_x:
U_x\rightarrow \R$ such that $E_x (U_x)$ is included in $[E_x (x)-\epsilon,E_x (x) + \epsilon]$ and $f(y)=\tau\circ E_x(y)$,  for all $y\in U_x$. 
Since $X$ is a Cantor space, we can suppose $U_x$ is a clopen set. 
Compactness of $X$ yields points $x_1, \dots , x_n$ such that $U_{x_1} , \dots , U_{x_n}$ is a finite clopen covering of $X$.
Then, for each $i\in [ 1,n]$ there exists a (possibly empty) clopen set $V_i$ included in $U_{x_i}$ such that $V_1 , \dots , V_n$ is a partition of $X$.
Let $I = \{ i \mid E_{x_i}( x_i)+\epsilon >1\}$.
Then the map $F_{\epsilon} = \sum_{i\not \in I} {\bf 1}_{V_i} E_{x(i)} + \sum_{i\in I} {\bf 1}_{V_i} (E_{x(i)} - 1)$ fulfills the requirement.
\end{proof}

\begin{proposition}\label{lemma3}
Let $U$ be a non empty clopen set in $X$ and $f$ belonging to $C(X,\R)$.
\begin{enumerate}
\item If $f^{(n)}(x)$ belongs to $\Z$ for every $x\in U$ such that $T^nx\in U$,
then $f$ is cohomologous to some $g\in C(X,\Z)$.
\item If $f^{(n)}(x)$ belongs to $\Z_+$ for every $x\in U$ such that $T^nx\in U$,
then $f$ is cohomologous to some $g\in C(X,\Z_+)$.
\end{enumerate}
\end{proposition}
\begin{proof}
1. Using the same method as the step 2 in the proof of Theorem~\ref{theoremGH},
we get that the sequence $(\tau(f^{(n_i)}(x_0))$ converges in $\T$
whenever $n_i\rightarrow\infty$ and $T^{n_i}x_0$ converges.
Thus there exists a continuous function $u : X\rightarrow \T$ such that $u(T^nx_0)=\tau (f^{(n)}(x_0))$
for every $n\ge 0$. If $x=T^nx_0$ for some $n\ge 0$, then
\begin{eqnarray*}
u(Tx)-u(x)&=&u(T^{n+1}x_0)-u(T^nx_0)\\
&=&\tau(f^{(n+1)}(x_0)-f^{(n)}(x_0))=\tau(f(T^nx_0))=
\tau(f(x)).
\end{eqnarray*}
By density, the same is true for every $x\in X$. By Lemma~\ref{lemmaTorus}
%\marginpar{FD:modif}
there exists $h\in C(X,\R)$ such that $\tau\circ h=u$. The function
$g=f-\partial_Th$ belongs to $C(X,\Z)$ and is cohomologous to $f$.

2. The family $(f^{(n)})$ is uniformly bounded from below and the result follows from Proposition~\ref{lemma2}.
\end{proof}
%%%%%%%%%%%%
\section{Cylinder functions}\label{sectionCylinder}
We consider in this section the case of shift spaces,
in which the cylinder functions play an important role.

Let $X\subset A^\Z$ be a shift space. 
Recall that we denote by $\cL_n(X)$ the set words of length $n$ in $\cL(X)$.
We denote
by $R_n(X)$\index{symbols}{R@$R_n(X)$} the group of maps from $\cL_n(X)$ into $\R$, by
$Z_n(X)$\index{symbols}{Z@$Z_n(X)$} the group of maps from $\cL_n(X)$ into $\Z$
and by $Z_n^+(X)$ the corresponding subset of non-negative maps.

For $\phi\in R_n(X)$, the function $\u(\phi):X\rightarrow \R$ given by
\begin{displaymath}
\u(\phi)(x)=\phi(x_{[0,n-1]})
\end{displaymath}
\index{symbols}{phi@$\u(\phi)$}%
is called the \emph{cylinder function}\index{subject}{cylinder! function}
associated to $\phi$. It belongs to $C(X,\Z)$ when $\phi$ belongs to $Z_n(X)$
and to $C(X,\Z_+)$ when $\phi$ belongs to $Z_n^+(X)$.

A function $f\in C(X,\Z)$ is a cylinder function if and only if
there exists $n\ge 1$ such that $f(x)$ depends only on $x_{[0,n-1]}$.

\begin{proposition}\label{lemma4}
Let $X$ be a shift space.
Every function in $C(X,\Z)$ (resp. $C(X,\Z_+)$) is cohomologous
to some cylinder function (resp. non-negative cylinder function).
\end{proposition}
\begin{proof}
Let $f\in C(X,\Z)$. Since $f$ is locally constant, there exists $n$
such that $f(x)$ depends only on $x_{[-n,n]}$. Then $f(S^nx)$ depends
only on $x_{[0,2n]}$ and $f\circ S^n$ is a cylinder function.
Since $f\circ S^n-f=(\partial_Sf)^{(n)}$ by Equation~\eqref{eq1}, $f$
is cohomologous to $f\circ S^n$
and the conclusion follows. Finally, if $f$ belongs to $C(X,\Z_+)$, then $f\circ S^n$ is non-negative.
\end{proof}
\begin{proposition}\label{lemma4ii}
Let $X$ be a shift space.
If a cylinder function with integer values is a coboundary, it is
the coboundary of some cylinder function.
\end{proposition}
\begin{proof}
Let $g\in C(X,\Z)$ and suppose that $f=\partial_Sg$ is a cylinder function.
We may choose $n$ large enough so that simultaneously $g(x)$
depends only on $x_{[-n,n]}$ and $f(x)$ depends only on $x_{[0,n]}$.
Assume that $y_{[0,2n]}=z_{[0,2n]}$. Since $f(x)$ depends only on
$x_{[0,n]}$, the value $f^{(n)}(x)$ depends only on $x_{[0,2n]}$
and thus  $f^{(n)}(y)= f^{(n)}(z)$. Similarly, since
$(S^ny)_{[-n,n]}=(S^nz)_{[-n,n]}$, we have $g(S^ny)=g(S^nz)$.
Thus, by Equation~\eqref{eq1}, we have 
\begin{eqnarray*}
g(y)&=&g(S^ny)-f^{(n)}(y)\\
&=&g(S^nz)-f^{(n)}(z)=g(z)
\end{eqnarray*}
and thus $g$ is a cylinder function.
\end{proof}

\section{Ordered group of a recurrent shift space}\label{chapterOrderedGroupRecurrent}
In this section, we show how to compute the ordered cohomology group
of a recurrent shift space using the cylinder functions
introduced in Section~\ref{sectionCylinder}. This description
will be used in the next chapters. We will first define
an ordered group $\G_n(X)$ associated with cylinder functions
corresponding to words of length $n$ (Proposition~\ref{propositionG_n}).
In a second part we show that the direct limit of these
groups is the cohomology group $K^0(X,S)$ of the shift space
$X$ (Proposition~\ref{propositionK0Recurrent}).

\subsection{Groups associated with cylinder functions}\label{sectionGroupCylinder}

For a word $w=a_1a_2\cdots a_n$ of length $n\ge 1$ with $a_i\in A$, we set
$p_n(w)=a_1\cdots a_{n-1}$ and $s_n(w)=a_2\cdots a_n$. Next, given
a shift space $X$, recall that $R_n(X)$ denotes the group
of maps from the set $\cL_n(X)$ of words of length $n$ in $\cL(X)$
into $\R$. 
We define, for $n\ge 1$,
three group morphisms
\begin{displaymath}
p_{n-1}^*,s_{n-1}^*,\partial_{n-1}:R_{n-1}(X)\rightarrow R_n(X)
\end{displaymath}
\index{symbols}{delta@$\partial_n$}%
by
\begin{displaymath}
p_{n-1}^*  ( \phi )=\phi\circ p_n,\ s_{n-1}^* (\phi )=\phi\circ s_n,
\end{displaymath}
and
\begin{displaymath}
 \partial_{n-1} (\phi )=s_{n-1}^* (\phi )-p_{n-1}^* (\phi )=\phi\circ s_n-\phi\circ p_n
\end{displaymath}
for every $\phi\in R_{n-1}(X)$. These morphisms map $Z_{n-1}(X)$ into $Z_n(X)$
(recall that $Z_n(X)$ denotes the group of functions from $\cL_n(X)$
into $\Z$).
Moreover, $p_{n-1}^*$ and $s_{n-1}^*$ are injective and positive.

Note that one has the equality
\begin{equation}
  s_n\circ p_{n+1}=p_n\circ s_{n+1}\label{eqChristian}
\end{equation}
since both sides send $a_1\ldots a_{n+1}$ to $a_2\cdots a_n$. It
follows from \eqref{eqChristian} that
\begin{equation}
  \partial_n\circ p_{n-1}^*=p_n^*\circ \partial_{n-1},\label{eqChristian2}
\end{equation}
as one may verify.

Note also that for every $\phi\in R_{n-1}(X)$, the cylinder functions
associated to $p_{n-1}^* (\phi )$ and to $\phi$ are the same, that is
\begin{equation}
\underline{p_{n-1}^* (\phi )}=\underline{\phi}.\label{eqCylinderFctsSame}
\end{equation}
Note also that for 
$\phi\in R_{n-1}(X)$ and $\psi\in R_{n}(X)$ one has
\begin{equation}
\psi=\partial_{n-1}(\phi)\Leftrightarrow \u(\psi)=\partial_S\u(\phi)
\label{eqCylinder}
\end{equation}
and thus the cylinder function $\u(\psi)$ 
associated to $\psi\in\partial_{n-1}R_{n-1}(X)$ is a coboundary.
We denote
\begin{equation}
\label{def:gnX}
G_n(X)=Z_n(X)/\partial_{n-1}Z_{n-1}(X)
\end{equation}
\index{symbols}{G@$G_n(X)$}%
the quotient  of the group $Z_n(X)$ by its subgroup $\partial_{n-1}Z_{n-1}(X)$,
we denote by $G_n^+(X)$ the image in $G_n(X)$ of $Z_n^+(X)$
and by $\mathbf{1}_n(X)$ the image in $G_n(X)$ of the constant function
$\mathbf{1}\in Z_n(X)$.
\begin{proposition}\label{propositionG_n}
For every recurrent  shift space $X$, the triple 
\begin{displaymath}
\G_n(X)=(G_n(X),G_n^+(X),\mathbf{1}_n(X))
\end{displaymath}
is a directed unital ordered group.
\end{proposition}
\begin{proof}
The set $G_n^+(X)$ is a submonoid of $G_n(X)$ because $Z_n^+(X)$ is a submonoid
of $Z_n(X)$ and it generates $G_n(X)$ because $Z_n^+(X)$ generates $Z_n(X)$.
%To prove that $G_n^+(X)\cap (-G_n^+(X))=\{0\}$, we have to prove that
%if $\phi\in Z_{n-1}(X)$ is such that $\partial_{n-1}\phi\in Z_n^+(X)$, then
%$\partial_{n-1}\phi=0$. Let $u,v\in L_{n-1}(X)$. Since $X$
%is recurrent, there
%exists $m\ge n$ and $w\in L_m(X)$ such that $u$ is a prefix of $w$
%and $v$ is a suffix of $w$. Then, since $u=w_{[1,n-1]}$ and
%$v=w_{[m-n+2,m]}$, we have, because the first sum is telescopic,
%\begin{eqnarray*}
%\phi(v)-\phi(u)&=&\sum_{i=1}^{m-n+1}(\phi(w_{[i+1,i+n-1]})-\phi(w_{[i,i+n-2]}))\\
%&=&\sum_{i=1}^{m-n+1}(\partial_{n-1}\phi)(w_{[i,i+n-1]})\ge 0.
%\end{eqnarray*}
% Since this is true for every $u,v$,
%it implies that $\phi$ is constant and thus that $\partial_{n-1}\phi=0$.
Let $\alpha,\beta\in Z_n^+(X)$ be such that $\alpha+\beta=\partial_{n-1}\phi$
for some $\phi\in Z_{n-1}(X)$. Then $\partial\u(\phi)$ is a positive
coboundary. By Corollary~\ref{corollary1}, this implies $\partial_S\u(\phi)=0$
and thus $\partial_{n-1}\phi=0$. Thus $\alpha+\beta=0$, which
implies $\alpha=\beta=0$. This shows that $G_n^+(X)\cap (-G_n^+(X))=\{0\}$
and thus $\G_n(X)$ is a directed ordered group. It is unital
because $\1$ is a unit of $Z_n(X)$.
\end{proof}
A direct proof of Proposition~\ref{propositionG_n},
that is without using Corollary~\ref{corollary1},
is given in Exercise~\ref{exerciseG_n}.

Note that, since $\partial_n\circ p_{n-1}^*=p_n^*\circ\partial_{n-1}$
by \eqref{eqChristian2}, the morphism
$p_n^*$ induces a morphism of the quotient groups
\begin{displaymath}
i_{n+1,n}:G_n(X)\rightarrow G_{n+1}(X)
\end{displaymath}
which is morphism of unital ordered groups.
\begin{example}
Let $X$ be the Fibonacci shift.
\index{subject}{Fibonacci!shift}\index{subject}{shift space!Fibonacci}%
 We have
 $\cL_1(X)=\{a,b\}$ and $\cL_2(X)=\{aa,ab,ba\}$. For $u\in\cL_n(X)$,
 denote by $\phi_u$ the element of $Z_n(X)$ which takes the value $1$
 on $u$ and $0$ elsewhere, that is, such that
 $\charac_{[u]}=\u(\phi_u)$. Since 
$\partial_1(\phi_{a})=\phi_{ba}+\phi_{aa}-(\phi_{aa}+\phi_{ab})=\phi_{ba}-\phi_{ab}$, the projections
of $\phi_{ba}$ and $\phi_{ab}$ in $G_2(X)$ are equal
and $G_2(X)\simeq\Z^2$. The matrix of the projection
of $Z_2(X)$ on $G_2(X)$ is
\begin{displaymath}
P=\begin{bmatrix}1&0&0\\0&1&1\end{bmatrix}
\end{displaymath}
The matrix of the morphism
$i_{2,1}$ is
\begin{displaymath}
\begin{bmatrix}1&0\\1&1\end{bmatrix}
\end{displaymath}
since for example
\begin{displaymath}
i_{2,1}\begin{bmatrix}1\\0\end{bmatrix}=P(p_1^*\charac_{[a]})=P\begin{bmatrix}1\\1\\0\end{bmatrix}=\begin{bmatrix}1\\1\end{bmatrix}.
\end{displaymath}
\end{example}
Note that Example~\ref{examplea*b*3} shows that Proposition
\ref{propositionG_n} does not hold without the hypothesis that
$X$ is recurrent.
\subsection{Ordered cohomology group}
We now prove the following result which describes the
ordered cohomology group $K^0(X,S)$ of a recurrent shift
space $X$.
\begin{proposition}\label{propositionK0Recurrent}
For every recurrent shift space  $X$,  the unital ordered group
$K^0(X,S)$ is the inductive limit
of the family $\G_n(X)$ with the morphisms $i_{n+1,n}$.
\end{proposition}
\begin{proof}
To every $\phi\in Z_n(X)$ we can associate the corresponding cylinder function
$\u(\phi)$ and its projection $\pi(\u(\phi))$ in $H(X,S,\Z)$. When
$\phi$ is in  $\partial_{n-1}(Z_{n-1}(X))$, the cylinder function $\u(\phi)$ is a coboundary and $\pi(\u(\phi))=0$. If $\phi$ belongs to $Z_n^+(X)$, then  $\u(\phi)$ is in $C(X,\Z_+)$
and the cylinder function associated with the constant function
equal to $1$ on $\cL_n(X)$ is the constant function equal to $1$ on $X$.

Thus we have defined a morphism of unital
ordered groups
\begin{displaymath}
j_n:\G_n (X)\rightarrow K^0(X,S)
\end{displaymath}
and clearly $j_{n+1}\circ i_{n+1,n}=j_n$. By the universal property
of direct limits (Proposition~\ref{propositionUniversalProperty}) the sequence $(j_n)_{n\ge 1}$
induces a morphism $j$ from the inductive limit of the $\G_n (X)$ to the 
group $K^0(X,S)$. Let us show that $j$ is an isomorphism.

Let $\phi\in Z_n(X)$, let $\alpha$ be its projection in $\G_n(X)$
and suppose that $j_n (\alpha )=0$. Then $\u(\phi)$ is a coboundary
and, by Proposition~\ref{lemma4ii}, it is the coboundary of some
cylinder function $\u(\psi)$. Let $m\ge n$ be such that
 $\psi$ belongs to $Z_m(X)$. 
Since  $p_m^*\circ\cdots\circ p_n^* (\phi )$ belongs to $Z_{m+1}(X)$,
we have by Equation~\eqref{eqCylinderFctsSame}
\begin{displaymath}
\underline{\phi}=\underline{p_m^*\circ\cdots\circ p_n^* (\phi )}.
\end{displaymath}
Now, since $\underline{\phi}=\partial_T (\underline{\psi})$,
we have by \eqref{eqCylinder}
\begin{displaymath}
\underline{p_m^*\circ\cdots\circ p_n^* (\phi )}=\partial_T (\underline{\psi})\quad
\Rightarrow \quad p_m^*\circ\cdots\circ p_n^* (\phi )=\partial_m (\psi )
\end{displaymath}
which implies $p_m^*\circ\cdots\circ p_n^* (\phi )\in\partial_m(Z_m(X))$
and $i_{m,m-1}\circ\cdots\circ i_{n+1,n} (\alpha )=0$.
Thus the image of $\alpha$ in the inductive limit is $0$. This
shows that $j$ is injective.

Moreover, every $f\in C(X,\Z)$ (resp. $C(X,\Z_+)$) is cohomologous
to a cylinder function (resp. to a non-negative cylinder function).
It follows that $j$ is onto and maps the positive cone of the inductive limit
onto $H^+(X,S,\Z)$.
\end{proof}
We now present an example 
which shows that the ordered cohomology group of a shift
of finite type may fail to be a dimension group.
This contrasts with the fact that, as we will see in Chapter~\ref{chapterDimensionGroupsPartitions}, the ordered cohomology group of a minimal
Cantor system is a simple dimension group (Theorem~\ref{theoremDimensionGroup}).
\begin{example}\label{exampleKimRoushWilliams}
Let $X$ be the set of two-sided infinite paths in the graph
represented in Figure~\ref{figureExampleKimRoushWilliams}.
\begin{figure}[hbt]
\centering
\tikzset{node/.style={circle,draw,minimum size=0.4cm,inner sep=0pt}}
\tikzstyle{every loop}=[->,shorten >=1pt,looseness=12]
\tikzstyle{loop left}=[in=130,out=220,loop]
\tikzstyle{loop right}=[in=330,out=50,loop]
\begin{tikzpicture}
\node[node](1) at(0,0){$1$};
\node[node](2) at (2,0){$2$};

\draw[left](1) edge [loop left]node {$a$}(1);
%\draw[below,bend left, ->](1) edge node{$c$}(2);
\draw[below,bend left, ->](2) edge node{$d,e$}(1);
\draw[above,bend left, ->](1) edge node{$b,c$}(2);
%\draw[below,bend left, ->](2) edge node{$e$}(1);
\draw[right](2)edge[loop right]node{$f$}(2);
\end{tikzpicture}
\caption{A shift of finite type}\label{figureExampleKimRoushWilliams}
\end{figure}
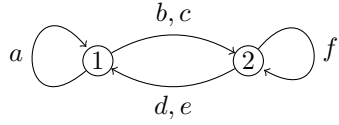
Since $a$ can be preceded by $a,d$ or $e$, we have
\begin{displaymath}
\chi_{[a]}\circ T=\chi_{[aa]}+\chi_{[da]}+\chi_{[ea]}.
\end{displaymath}
Similarly, since $a$ can be followed by $a,b$ or $c$, we have
\begin{displaymath}
\chi_{[a]}=\chi_{[aa)}+\chi_{[ab]}+\chi_{[ac]}.
\end{displaymath}

Thus $\partial \charac_{[a]}=\charac_{[da]}+\charac_{[ea]}-\charac_{[ab]}-\charac_{[ac]}$. We have then in $H(X,S,\Z)$ the equality
\begin{eqnarray*}
\charac_{[d]}+\charac_{[e]}&=&\charac_{[da]}+\charac_{[ea]}+\sum_{x\in\{d,e\},y\in\{b,c\}}\chi_{[xy]}\\
&\simeq&\chi_{[ab]}+\chi_{[ac]}+\sum_{x\in\{e,d\},y\in\{b,c\}}\chi_{[xy]}\\
&\simeq& \charac_{[b]}+\charac_{[c]}.
\end{eqnarray*}
But 
$\charac_{[b]}+\charac_{[c]}\simeq\charac_{[d]}+\charac_{[e]}$ implies
\begin{displaymath}
0\le \charac_{[b]}\le \charac_{[d]}+\charac_{[e]}.
\end{displaymath}

It can be shown, using an invariant Markov measure on $X$,
that there is no decomposition $\charac_{[b]}=\phi_1+\phi_2$
such that $\phi_1\le\charac_{[d]}$ and $\phi_2\le\charac_{[e]}$
(Exercise~\ref{exerciseKimRoushWilliams}).
Thus $H(X,S,\Z)$ is not a Riesz group.

\end{example}
%%%%%%%%%%%%%%%%%%%%%%%
\section{Factor maps and conjugacy}\label{sectionFactorMaps}
Let $(X,T)$ and $(X',T')$ be two topological dynamical systems
and let $\phi$ be a
factor map from $(X,T)$ to $(X',T')$. If $f\in C(X',\Z)$ is the
coboundary of $g$, then $f\circ\phi$ is the coboundary of
$g\circ\phi$. Therefore, the group homomorphism $f\mapsto f\circ\phi$
from $C(X',\Z)$ to $C(X,\Z)$ induces a group homomorphism
\begin{displaymath}
\phi^*:H(X',T',\Z)\rightarrow H(X,T,\Z).
\end{displaymath}
Clearly, the map $\phi^*$ is a morphism of ordered groups with order
units since $\phi^*(H^+(X',T',\Z))$ is a subset of $H^+(X,T,\Z)$
and $\phi^*(\mathbf{1}_{X'})=\mathbf{1}_X$.

\begin{proposition}\label{proposition4.7.1}
If $(X,T)$ is minimal, the morphism $\phi^*$ is injective
and 
\begin{displaymath}
\phi^*(H^+(X',T',\Z))=\phi^*(H(X',T',\Z))\cap H^+(X,T,\Z).
\end{displaymath}
\index{symbols}{phi@$\phi^*$}
\end{proposition}
\begin{proof}
If $f\in C(X',\Z)$ is such that $f\circ\phi$ is a coboundary, then
by Theorem~\ref{theoremGH} the sequence of functions
$((f\circ\phi)^{(n)})_{n\ge 0}$ is uniformly bounded. But, for all $n\ge 0$,
$(f\circ\phi)^{(n)}=f^{(n)}\circ\phi$ and since, by minimality, the map $\phi$ is onto,
$||f^{(n)}\circ\phi||_\infty=||f^{(n)}||_\infty$. Thus the sequence
$(f^{(n)})_{n\ge 0}$ is uniformly bounded and, by Theorem~\ref{theoremGH}
again, $f$ is a coboundary in $C(X',\Z)$. This shows that
$\phi^*$ is injective.

To prove the second assertion, we first note that
the inclusion from left to right is clear since
$\phi^*(H^+(X',T',\Z))$ is included in $H^+(X,T,\Z)$. Conversely, consider
$f\in C(X',\Z)$ and $g\in C(X,\Z)$
such that $f\circ\phi+\partial_T (g)\ge 0$. By Proposition~\ref{lemma2},
the sequence of functions $((f\circ\phi)^{(n)})_{n\ge 0}$ is bounded
from below. By the same argument as in  the proof of the first
assertion, the sequence $(f^{(n)})_{n\ge 0}$ is bounded from below.
Thus, by Proposition~\ref{lemma2} again, there is a function
$g'\in C(X',\Z)$ such that $f+\delta_{T'} (g')\ge 0$. Thus the class
of $f$ is in $H^+(X',T',\Z)$ and its image by $\phi^*$ is in
$\phi^*(H^+(X',T',\Z)$.
\end{proof}
We deduce from Proposition~\ref{proposition4.7.1} that
the ordered cohomology group is  invariant under conjugacy.
\begin{corollary}
If $\phi$ is a conjugacy from $(X,T)$ onto $(X',T')$, then the map
$\phi^*$ from $H(X',T',\Z)$ to $H(X,T,\Z)$ is an isomorphism.
\end{corollary}

Note that Proposition~\ref{lemma4ii} can  be deduced
from Proposition~\ref{proposition4.7.1}.
Indeed, let $(Y,S)$
be the one-sided shift space associated with the minimal shift space $(X,S)$
and $\theta:(X,S)\to(Y,S)$ be the natural morphism.
Assume that $f\in C(X,\Z)$ is a cylinder function with integer
values. Then $f=k\circ\theta$ for some
$k\in C(Y,\Z)$. If $f$ is a coboundary, by Proposition~\ref{proposition4.7.1},
$k$ is a coboundary. Thus $k=\partial_S (h)$
for some $h\in C(Y,S)$. Then
$f=\partial_S(h)\circ\theta=h\circ S\circ\theta-h\circ\theta=h\circ \theta\circ S-h\circ\theta$. Thus $f$ is the coboundary of $h\circ\theta$
which is a cylinder function.

Note also that this argument shows that 
\begin{displaymath}
\theta^*:(H(Y,S,\Z),H_+(Y,S,\Z),\mathbf{1}_Y)\rightarrow
(H(X,S,\Z),H_+(X,S,\Z),\mathbf{1}_X)
\end{displaymath}
is an isomorphism of unital ordered groups.
%%%%%%%%%%%%%%
\section{Groups of induced systems}\label{sectionInduced}
Let $(X,T)$ be a minimal topological dynamical system and let
$U$ be a nonempty clopen subset of $X$. We have already
seen the notion of induced system $(U,T_U)$ in Section \ref{sec:inducedsystems}. 
Recall that
$T_U(x)=T^{n(x)}(x)$ where $n(x)=\inf\{n>0\mid T^nx\in U\}$.
The system $(U,T_U)$ is again minimal.

\begin{figure}[hbt]
\centering
\tikzset{node/.style={circle,draw,minimum size=0.1cm,inner sep=0pt}}
\tikzset{title/.style={minimum size=0.5cm,inner sep=0pt}}
\begin{tikzpicture}
  \node[title]at(-.5,0){$x$};
  \node[node](y0)at(0,1){};\node[node](y1)at(1,1){};\node[node](y2)at(2,1){};
  \node[title]at(2.5,1){$\cdots$};
  \node[node](x0)at(0,0){};\node[node](x1)at(2,0){};\node[node](x2)at(3,0){};
  \node[node](x3)at(5,0){};\node[node](x4)at(6,0){};
  \node[title]at(6.5,0){$\cdots$};

\draw[above](y0)edge node{$y_0$}(y1);\draw[above](y1)edge node{$y_1$}(y2);
\draw[left,->](y0)edge node{$\varphi$}(x0);
\draw[left,->](y1)edge node{}(x2);
\draw[above,->](y2)edge node{$\varphi$}(x4);
\draw[below](x0)edge node{$u$}(x1);\draw(x1) edge node{}(x2);
\draw[below](x2)edge node{$u$}(x3);\draw(x3) edge node{}(x4);
\end{tikzpicture}
\caption{The unique $y$ such that $\varphi(y)=x$.}\label{figureInducedSystem}
  \end{figure}
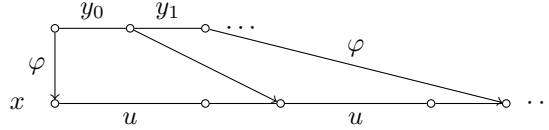

As an example, which will be considered often in the next chapters,
 let us consider a minimal shift space $X$, 
a word $u\in\cL(X)$
and the clopen set $U=[u]$.
Let $\RR'_X(u)$ be the set of left return words to $u$ and let
$\varphi:B\to \RR'_X(u)$ be a bijection extended to a morphism
from $B^*$ to $\RR'_X(u)^*$. Every $x\in U$ can be written
in a unique way as $x=\varphi(y)$ for some $y\in B^\Z$
(see Figure~\ref{figureInducedSystem}
or Proposition~\ref{propositionCircularCodesInjective} if you are not entirely convinced).

Then
$Y=\varphi^{-1}(U)$ is a shift space on $B$ and
$(U,S_U)$ is isomorphic to $(Y,T)$ where $T$ is the shift
of $B^\Z$. Indeed, by definition of return words, for every $x\in[u]$, 
the integer
$n(x)=|\varphi(y_0)|$ is the length of the unique prefix of $x^+$ in
$\RR_X'(u)$. Thus
\begin{displaymath}
S_U\circ\varphi(y)=S^{n(x)}(\varphi(y))=\varphi\circ T(y)
\end{displaymath}
showing that $\varphi$ is a conjugacy from $(Y,T)$
onto $(U,S_U)$. The shift $(Y,T)$ is called the \emph{derivative shift}
\index{subject}{derivative!shift}\index{subject}{shift space!derivative}%
of $X$ with respect to $[u]$.
\begin{example}
Let $X$ be the Fibonacci shift,
\index{subject}{Fibonacci!shift}\index{subject}{shift space!Fibonacci}%
 which is the
substitution shift generated by $\varphi:a\to ab,b\to a$. The system induced on $U=[a]$ is isomorphic
to $X$. Indeed, the set of left return words to $a$ is $\RR'(a)=\{ab,a\}$
and thus $\varphi$ is an isomorphism from $X$ onto $(U,T_U)$.
\end{example}

Let $I_U:C(X,\Z)\rightarrow C(U,\Z)$ and $R_U:C(U,\Z)\rightarrow C(X,\Z)$
\index{symbols}{I@$I_U$}\index{symbols}{R@$R_U$}%
be the morphisms of ordered groups defined by
\begin{displaymath}
(I_U(f))(x)=f^{(n(x))}(x)\mbox{ for all $x\in U$},
\end{displaymath}
and let $R_U(f)$ be the map equal to $f$ on $U$
and equal to $0$ elsewhere.

We will now show that an induction does not modify the cohomology group,
as an ordered group (although not as a unitary
ordered group since the unit is not preserved, see the remark after
the proof).

\begin{proposition}\label{propositionIuRu}
  The maps $I_U,R_U$ are reciprocal
  isomorphisms between $C(X,\Z)$
  and $C(U,\Z)$. They induce reciprocal isomorphisms of ordered groups from
$K^0(X,T)$ onto $K^0(U,T_U)$.
\end{proposition}
\begin{proof}
Let $f\in C(X,\Z)$ be the coboundary of $g\in C(X,\Z)$.
Then for every $x\in U$, by Equation~\eqref{eq1},
we have $f^{(n(x))}=g\circ T^{n(x)}-g$ and thus $I_U(f)$
is the coboundary of the restriction of $g$ to $U$.
This shows that $I_U$ induces a morphism 
\begin{displaymath}
i_U:H(X,T,\Z)\rightarrow H(U,T_U,\Z).
\end{displaymath}
\index{symbols}{i@$i_U$}
 Since $I_U$ maps $C(X,\Z_+)$ to
$C(U,\Z_+)$, it is a morphism of ordered groups
from $K^0(X,T)$ to $K^0(U,T_U)$.

Let us now show that $R_U$ induces a map of ordered groups
from $K^0(U,T_U)$ to $K^0(X,T)$.
For every $x\in X$, let $m(x)=\inf\{m\ge 0\mid T^m(x)\in U\}$.
Then the function $x\mapsto m(x)$ is continuous, coincides
with $n(x)$ outside $U$ and is zero on $U$. Suppose that
$f\in C(U,\Z)$ is the coboundary in $(U,T_U)$ of some $g\in C(U,\Z)$.
Let $h(x)=g(T^{m(x)}(x))$. We have
\begin{displaymath}
h(Tx)-h(x)=\begin{cases}
f(x)&\mbox{ if $x\in U$}\\
0&\mbox{ if $x\notin U$}.
\end{cases}
\end{displaymath}
Indeed, if $x\in U$, then $m(x)=0$, $m(Tx)=n(x)-1$ and
\begin{displaymath}
h(T(x))-h(x)=g(T^{n(x)} (x))-g(x)=g(T_U(x))-g(x)=f(x)
\end{displaymath}
and if $x\notin U$, then $m(T(x))=m(x)-1$ and $h(T(x))=h(x)$.
This shows that $\partial_T (h)=R_U (f)$ and thus that $R_U$
maps coboundaries of $(U,T_U)$ to coboundaries of $(X,T)$
and thus induces a map 
\begin{displaymath}
r_U:H(U,T_U,\Z)\rightarrow H(X,T,\Z).
\end{displaymath}

Finally, 
for any $f\in C(U,\Z)$ and $x\in U$, we have
\begin{eqnarray*}
I_U\circ R_U(f)(x) = (R_U (f))^{(n(x))}(x)=f(x)+R_U(f)(T(x))\ldots+R_U(f)(T^{n(x)-1}(x)) =f(x).
\end{eqnarray*}
Thus $I_U\circ R_U(f)=f$
and $I_U\circ R_U$ is the identity. Next for every
$f\in C(X,\Z)$, we have
\begin{equation}
R_U\circ I_U(f)-f=\partial_T (k)\label{eqInduction}
\end{equation}
where $k(x)=f^{(m(x))}(x)$. Indeed, for every $x\in X$, the above equation can
be rewritten
\begin{displaymath}
(R_U(f))^{(n(x))}(x)-f(x)=f^{(m(T(x)))}(T(x))-f^{(m(x))}(x).
\end{displaymath}
If $x\in U$,  the value of the right-hand side is $f^{(m(Tx))}(T(x))-f^{(m(x))}(x)=
f^{(n(x)-1)}(T(x))-f^{(0)}(x)=f^{(n(x)-1)}(T(x))=f^{(n(x))}(x)-f(x)$.
Next, if $x\notin U$, the value of
the left-hand side is $-f(x)$  and the value of the right-hand side
is $f^{(n(x)-1)}(T(x))-f^{(n(x))}(x)=-f(x)$.
Equation~\eqref{eqInduction} shows that $R_U\circ I_Uf-f$ is the
coboundary of $k$ and thus that $R_U\circ I_U$ is the identity
on $H(X,T,\Z)$.
\end{proof}
Note that the isomorphism from $K^0(X,T)$ onto $K^0(U, T_U)$
is not an isomorphism of unital ordered groups because
the units are not the same. Indeed, the image by $R_U$
of $\charac_U$ is the characteristic function of $U$ and not
the constant function $\charac_X$.
This phenomenon will play an important role when
we compute the dimension group of an induced system
(see Proposition~\ref{propositionBVPrimitive}).

%%%%%%%%%%%%%%%%%%%
\section{Invariant Borel probability measures}\label{sectionInvariant}
Let $(X,T)$ be a topological dynamical system. 
There is an important connection between the ordered cohomology group
$H(X,T,\Z)$ 
of  $(X,T)$ and the invariant Borel probability
measures on $(X,T)$.

\subsection{Invariant measures}
Recall (see Appendix~\ref{appendixMeasureIntegration}) that
a \emph{Borel probability measure}\index{subject}{Borel!probability measure}
\index{subject}{probability measure!Borel} on a topological space $X$
is a  Borel measure $\mu$ such that $\mu(X)=1$.

It is useful to see how such a measure
is defined in a shift space. Consider  a shift space $X$
on the alphabet $A$. A Borel probability measure $\mu$ on $X$
 determines  a map $\pi:\cL(X)\to[0,1]$ by 
\begin{equation}\pi(u)=\mu([u]_X)\label{eqProbaWords}
\end{equation}
where $[u]=\{x\in X\mid x_{[0,|u|-1]}=u\}$ is the cylinder defined by $u$.
This map satisfies, as a consequence of the equality
$[u]=\cup_{a\in A}[ua]$, the \emph{compatibility conditions}
\index{subject}{compatibility conditions}%
$\pi(\varepsilon)=1$ and
\begin{equation}
\sum_{a\in A, ua\in \cL (X)}\pi(ua)=\pi(u)\label{eqCompatibility}
\end{equation}
for every $u\in\cL(X)$. 

Conversely, any map satisfying these compatibility
conditions defines a unique Borel probability measure
satisfying~\eqref{eqProbaWords}
by the Caratheodory extension theorem (see Appendix~\ref{appendixMeasureIntegration}).
\index{subject}{Carath\'eodory Extension Theorem}%
\index{subject}{Theorem!Carath\'eodory Extension}%
\index{names}{Carath\'eodory, Constantin}%

The following example gives the simplest possible measure on
the full shift. It corresponds to a sequence of successive independent
and identically distributed choices of the letters forming a sequence.
\begin{example}\label{exampleBernoulliMeasure}
  Let $X=A^\Z$ be the full shift on $A$. Let $\pi:A\to[0,1]$ be
    a map such that $\sum_{a\in A}\pi(a)=1$. We extend $\pi$ to
    multiplicative morphism from $A^*$ into $[0,1]$. Since
    the compatibility conditions \eqref{eqCompatibility} are
    satisfied, there is a unique Borel probability measure $\mu$ on $X$
    satisfying \eqref{eqProbaWords}. It is called a \emph{Bernoulli measure}.
    \index{subject}{Bernoulli measure}\index{names}{Bernoulli, Jacob}%
\end{example}
As a more general example, where the successive choices are not independent,
but depend on the previous one, we have the following important class of measures.
\begin{example}\label{exampleMarkovChain}
  Let $P$ be an $A\times A$-stochastic matrix and let $v$ be a stochastic
  row $A$-vector, that is, such that $\sum_{a\in A}v_a=1$.
  Let $\pi:A^*\to[0,1]$ be defined by $\pi(\varepsilon)=1$ and 
  \begin{equation}
    \pi(a_1a_2\cdots a_n)=v_{a_1}P_{a_1,a_2}\ldots P_{a_{n-1},a_n}\label{eqMarkovChain}
  \end{equation}
  for $n\ge 1$ and $a_i\in A$. Since $v$ and $P$ are stochastic, the map
  $\pi$ satisfies the compatibility conditions \eqref{eqCompatibility}. Indeed
  $\sum_{a\in A}\pi(a)=\sum_{a\in A}v_a=1$ and next
  \begin{displaymath}
    \sum_{a\in A}\pi(a_1\ldots a_na)=\sum_{a\in A}\pi(a_1\ldots a_n)P_{a_n,a}=\pi(a_1\cdots a_n).
  \end{displaymath}
  The probability measure on $A^\Z$ defined by $\pi$ is called a
  \emph{Markov measure} and the pair $(v,P)$  a \emph{Markov chain}.
  \index{subject}{Markov chain}
\end{example}

Let now $(X,T)$ be a topological dynamical system.
A  Borel measure $\mu$ on $X$ is said to be  \emph{invariant}
\index{subject}{invariant! measure}
\index{subject}{measure!invariant} if one has
$\mu(T^{-1}(U))=\mu(U)$ for every Borel subset $U$ of $X$.

Let $\mu$ be an invariant Borel probability measure  on a shift
space $X$.
The associated map $\pi:\cL(X)\to[0,1]$ such that
$\pi(u)=\mu([u])$ satisfies in addition to the compatibility
conditions \eqref{eqCompatibility}, the symmetric conditions
\begin{equation}
\pi(u)=\sum_{a\in A, au\in \cL (X)}\pi(au).\label{eqCompatibilityBis}
\end{equation}
Indeed, $\sum_{a\in A, au\in \cL (X)}\pi(au)=\sum_{a\in A}\mu([au])=\mu(T^{-1}([u]))=\mu([u])=\pi(u)$. Conversely, for every map $\pi:\cL(X)\to[0,1]$ satisfying
the compatibility conditions \eqref{eqCompatibility}
and \eqref{eqCompatibilityBis}, there is by Carath\'eodory Extension Theorem
(see Appendix~\ref{appendixMeasureIntegration})
a unique invariant Borel probability measure $\mu$
such that $\mu([u])=\pi(u)$.

It may be puzzling that the compatibility conditions
\eqref{eqCompatibility} and \eqref{eqCompatibilityBis} are left-right symmetric
of each other,
although the notion of invariant measure is by no means the symmetric
of the notion of measure. This apparent contradiction comes from
the asymmetry of the definition of $\pi$ which is 'future oriented'
since $x\in[u]$ if $x^+$ begins by $u$ (and not $x^-$).

For example, a Bernoulli measure is an invariant Borel probability
measure since \eqref{eqCompatibilityBis} is obviously satisfied.
\begin{example}
  A Markov chain $(v,P)$ is called \emph{stationary}\index{subject}{stationary Markov chain}
  \index{subject}{Markov chain!stationary} if $v$
  is a left eigenvector of $P$ for the eigenvalue $1$.
  A Markov measure on $A^\Z$ defined by a Markov chain $(v,P)$
  is invariant if and only if $(v,P)$ is stationary. Indeed, if $vP=v$,
  the map $\pi$ defined by \eqref{eqMarkovChain} satisfies
  for every $a\in A$ and $w=a_1\cdots a_n$
  \begin{eqnarray*}
    \sum_{a\in A}\pi(aa_1\ldots a_n)&=&\sum_{a\in A}(v_aP_{a,a_1})P_{a_1,a_2}\cdots P_{a_{n-1},a_n}\\
      &=&\pi(a_1\cdots a_n).
  \end{eqnarray*}
  Conversely, if $\mu$ is invariant, then the associated map $\pi$ satisfies
  $v_b=\pi(b)=\sum_{a\in A}\pi(ab)=\sum_{a\in A}(v_aP)_b=(vP)_b$ and thus $vP=v$.
  \end{example}

The set of invariant Borel probability
measures on $X$ is convex since for invariant Borel probability measures
$\mu,\nu$ and $\alpha,\beta$ with
$\alpha+\beta=1$, $\alpha\mu+\beta\nu$ is an invariant Borel probability measure.

\begin{theorem}[Krylov, Bogolyubov]\label{theoremKrylovBogolioubov}
\index{names}{Krylov, Nicolai M.}\index{names}{Bogolyubov, Nikolai N.}
\index{subject}{Krylov-Bogolyubov Theorem}
\index{subject}{Theorem!Krylov-Bogolyubov}
Any topological dynamical system has at least one invariant Borel probability
measure.
\end{theorem}
\begin{proof}
Let $(X,T)$ be a topological dynamical system and
let $x\in X$. For every $n\ge 0$ we consider the probability measure
$\delta_{T^nx}$, where $\delta_y$ denotes 
the Dirac measure at the point $y$. 
Consequently
\begin{displaymath}
\mu_N=\frac{1}{N}\sum_{n<N}\delta_{T^n(x)}
\end{displaymath}
is also a Borel probability measure. 

By Theorem~\ref{theoremBanachAlaoglu},
 the  sequence $(\mu_N)$
has a cluster point $\mu$
for the weak-star
topology.
Let $(\mu_{N_i})$ be a the subsequence converging to $\mu$. For
every Borel subset $U$ of $X$,
\begin{eqnarray*}
\mu(U)&=&\lim\frac{1}{N_i}\Card\{n<N_i\mid T^n(x)\in U\}\\
&=&\lim\frac{1}{N_i}\Card\{n<N_i\mid T^{n+1}(x)\in U\}=\mu(T^{-1}(U)).
\end{eqnarray*}
Thus $\mu$ is an invariant Borel probability measure.
\end{proof}
A different proof uses the \emph{Markov-Kakutani fixed point Theorem}
\index{subject}{Markov-Kakutani fixed point Theorem}%
\index{subject}{Theorem!Markov-Kakutani}%
\index{names}{Markov, Andrei A.}%
\index{names}{Kakutani, Shizuo}%
(see Exercise~\ref{exerciseMarkovKakutani}).

A \emph{measure-theoretic dynamical system}%
\index{subject}{measure!theoretic dynamical system}%
\index{subject}{dynamical system!measure-theoretic} is
a system $(X,T,\mu)$
\index{symbols}{XT@$(X,T,\mu)$} where $X$ is a compact
metric space, $\mu$ is a Borel probability
measure on $X$ and $T:X\to X$ is a measurable map that preserves
the measure $\mu$, that is, such that 
$T^{-1}(U)$ is a Borel set with measure $\mu(T^{-1}(U))=\mu(U)$
for every Borel set $U\subset X$.

Thus, Theorem~\ref{theoremKrylovBogolioubov} shows that
every topological dynamical system may be considered as
a measure-theoretic one. Although this book is devoted
to topological dynamical systems, we will have occasions
to meet measure-theoretic ones.

A basic result concerning measure-theoretic dynamical systems
is \emph{Poincar\'e Recurrence Theorem},
\index{subject}{Poincar\'e Recurrence Theorem}%
\index{subject}{Theorem!Poincar\'e Recurrence}%
\index{names}{Poincar\'e, Henry}%
 which we will use
in Chapter~\ref{ch5:sec:examples}.
\begin{theorem}[Poincar\'e]\label{theoremPoincare}
Let $(X,T,\mu)$ be a measure-theoretic dynamical system.
 For every Borel set $U\subset X$
such that $\mu(U)>0$, the set of points $x\in U$
such that  for some $N\ge 1$ one has
$T^n (x)\notin U$ for all $n\ge N$ has measure $0$.
\end{theorem}
\begin{proof}
Let $V_N$ be the set of $x\in U$ such that $T^n(x)\notin U$ for all
$n\ge N$. It is a Borel set because 
\begin{displaymath}
V=U\cap\left(\cap_{n\ge N}T^{-n}(X\setminus U)\right)
\end{displaymath}
Next, by definition of $V_N$, we have  $T^{-n}(V_N)\cap V_N=\emptyset$ for all $n\ge N$,
which implies that the sets $T^{-N}(V_N),T^{-N-1}(V_N),\ldots$ are disjoint.
Therefore
\begin{displaymath}
\mu(X)\ge\mu(\bigcup_{n\ge N}T^{-n}(V_N))=\sum_{n\ge N}\mu(T^{-n}(V_N))=\sum_{n\ge 0}\mu(V_N).
\end{displaymath}
Since $\mu(X)=1$, this implies that $\mu(V_N)=0$.
Consequently, the set $V=\cup_{N\ge 1}V_N$ has measure $\mu(V)=\sum_{N\ge 1}\mu(V_N)=0$.
\end{proof}

\subsection{Ergodic measures}
Recall that, for a topological dynamical system $(X,T)$, a subset $U$ of $X$
is said to be \emph{invariant}\index{subject}{invariant!set}
if $T^{-1}(U)=U$.

An invariant Borel probability measure on $(X,T)$
 is \emph{ergodic}\index{subject}{ergodic!probability measure}\index{subject}{probability measure!ergodic} whenever $\mu (U)$ equals $0$ or $1$
for every invariant Borel subset $U$ of $X$. One also
says that the transformation $T$ is ergodic
\index{subject}{ergodic!transformation} with respect to $\mu$
or that the triple $(X,T,\mu)$ is ergodic.
\index{subject}{ergodic!system}%

As an example, a Bernoulli measure is ergodic (Exercise~\ref{exerciseBernoulliErgodic})
and a Markov measure defined by $(v,P)$ with $v>0$ and $vP=v$
is ergodic if and only if $P$ is irreducible (Exercise~\ref{exerciseMarkovErgodic}).

A basic result, that we state without proof,
is the Birkhoff's \emph{ergodic Theorem}.
\index{subject}{Birkhoff!ergodic Theorem}
\index{subject}{ergodic!theorem}\index{names}{Birkhoff, George D.}%
\begin{theorem}[Birkhoff]\label{theoremErgodic}
  Let $\mu$ be an ergodic measure on $(X,T)$.
  For every integrable function $f$ on $X$, the sequence
$f_n=\frac{1}{n}f^{(n)}(x)$ converges $\mu$-almost everywhere to $\int fd\mu$.
\end{theorem}
The functions $f_n$ are sometimes called the \emph{Birkhoff averages}.
\index{subject}{Birkhoff!averages}

A real valued measurable
function $f$ on $X$ is \emph{invariant}\index{subject}{invariant!function}
if $f=f\circ T$. Thus a set $U$ is invariant if and only if
its characteristic function is invariant. 

For two sets $U,V$ we write $U=V\bmod\mu$ if $U,V$ differ
by sets of measure $0$. The following statement gives a useful
variant of the definition of an ergodic measure.
\begin{proposition}\label{propositionCondEquiv}
The following conditions are equivalent for an invariant Borel
probability measure $\mu$ on $(X,T)$.
\begin{enumerate}
\item[\rm(i)]
 $\mu$ is ergodic.
\item[\rm(ii)] Every Borel set $U$
such that $T^{-1}(U)=U\bmod \mu$ is such that $\mu(U)=0$ or $\mu(U)=1$.
%\item[\rm(iii)]
%Every measurable invariant function is constant almost everywhere.
%\item[\rm(iv)] Every integrable function $f$ such that
%$f\circ T=f$ almost everywhere is constant almost everywhere.
\end{enumerate}
\end{proposition}
The proof is left as Exercise~\ref{exerciseInvariantFunction}.

Ergodicity is the measure-theoretic counterpart of minimality.
Both notions are different. For example a Bernoulli measure
on the full shift is ergodic although the full shift is not minimal.
Conversely, an invariant Borel probability measure on a minimal shift
need not be ergodic as we shall see.
We  note the following relation between ergodicity and recurrence.
\begin{proposition}\label{propositionErgodicRecurrent}
Let $(X,T)$ be  a topological dynamical system and let
$\mu$ be an invariant Borel probability on $(X,T)$. Assume
that $\mu(U)>0$ for every nonempty open set $U$. If $\mu$
is ergodic, then $(X,T)$ is recurrent.
\end{proposition}
The proof is left as an exercise (Exercise~\ref{exerciseErgodicRecurrent}).

Recall from Appendix~\ref{appendixMeasureIntegration}
that a Borel probability measure $\nu$
is \emph{absolutely continuous}\index{subject}{absolutely continuous measure}
\index{subject}{measure!absolutely continuous} with respect to $\mu$,
denoted $\nu\ll\mu$,\index{symbols}{nu@$\nu\ll\mu$} if for every Borel set $U\subset X$
such that $\mu(U)=0$ one has $\nu(U)=0$.

The following result shows that the ergodic measures
are the minimal invariant measures with respect to the preorder $\ll$.
\begin{proposition}\label{propositionAbsoluteCont}
Let $\mu,\nu$ be invariant Borel probability measures on $(X,T)$.
If $\mu$ is ergodic and $\nu\ll\mu$, then $\mu=\nu$.
\end{proposition}
\begin{proof}
Since $\nu\ll\mu$, by the Radon-Nikodym Theorem
(see Appendix \ref{appendixMeasureIntegration}), there is
a non-negative $\mu$-integrable function $f$ such that $\nu(U)=\int_U f d\mu$
for every Borel set $U\subset X$.

Consider the Borel set $B=\{x\in X\mid f(x)>1\}$. We will prove
that $B$ is invariant modulo a set of $\mu$-measure zero.
Note first that $\mu(T^{-1}(B)\setminus B)=\mu(B\setminus T(B))$
and $\nu(T^{-1}(B)\setminus B)=\nu(B\setminus T(B))$
since $\mu,\nu$ are invariant measures.
Assume that $\mu(T^{-1}(B)\setminus B)=\mu(B\setminus T(B))>0$.
 Then we have
\begin{eqnarray*} 
\mu(B\setminus T^{-1}(B))&<&\int_{B\setminus T^{-1}(B)} f d\mu=\nu(B\setminus T(B))
=\nu(T^{-1}B\setminus B)\\
&=&\int_{T^{-1}(B)\setminus B}f d\mu\le\mu((T^{-1}(B))\setminus B)=\mu(B)\setminus T^{-1}B)
\end{eqnarray*}
which is absurd
and thus $\mu(T^{-1}(B)\setminus B)=\mu(B\setminus T(B))=0$.
This in turn implies that $B=T^{-1}(B)$ modulo a set of $\mu$-measure
zero. By Proposition~\ref{propositionCondEquiv}, since $\mu$
is ergodic, this implies that $\mu(B)=0$ or $1$. If
$\mu(B)=1$, then $\nu(X)>1$ which is absurd. Thus $\mu(B)=0$.
Since $\nu(X)=\int_X f d\mu=1$, this
implies that $f=1$ almost everywhere or equivalently $\mu=\nu$.
\end{proof}

An \emph{extreme point}\index{subject}{extreme point} of a convex set
$K$ is a point which does not belong to any open line segment in $K$.
By the Krein-Millman Theorem (see Appendix~\ref{appendixMeasureIntegration}),
\index{subject}{Krein-Millman Theorem}%
\index{subject}{Theorem!Krein-Millman}%
the set ${\cal M}(X,T)$\index{symbols}{M@$\M(X,T)$} of invariant Borel probability measures on $(X,T)$
is the convex hull of its extreme points.
\index{names}{Krein, Mark}\index{names}{Millman, David}

\begin{proposition}
The ergodic measures are the extreme points of the set of invariant
probability measures.
\end{proposition}
\begin{proof}
Let $\mu$ be an extreme point of ${\cal M}(X,T)$. If $\mu$ is not ergodic,
there is an invariant Borel set $U$  such that $0<\mu(U)<1$.
The complement $V$ of $U$ is also invariant.
This defines two invariant Borel probability measures $\mu_1$ and $\mu_2$ by
\begin{displaymath}
\mu_1(W)=\frac {1}{\mu(U)}\mu(W\cap U),\quad
\mu_2(W)=\frac {1}{\mu(U)}\mu(W\cap V) .
\end{displaymath}
But then
\begin{displaymath}
\mu=\mu(U)\mu_1+(1-\mu(U))\mu_2
\end{displaymath}
shows that $\mu$ is not an extreme point.

Let now $\mu$ be ergodic.
Then $\mu=\alpha\mu_1+(1-\alpha)\mu_2$ for  $0\le\alpha\le 1$
and some
extreme points $\mu_1,\mu_2\in{\cal M}(X,T)$. Assume $\alpha>0$.
For any Borel set $U\subset X$, one has that $\mu(U)=0$
implies $\mu_1(U)=0$, that is $\mu_1$ is absolutely continuous
with respect to $\mu$. By Proposition~\ref{propositionAbsoluteCont},
this implies that $\mu=\mu_1$. Thus $\mu$ is an extreme point.
\end{proof}
In particular, we have the following important case.
\begin{corollary}\label{corollaryUniqueisErgodic}
If there is a unique invariant Borel probability measure, it  is ergodic. 
\end{corollary}

The following statement gives additional properties of the
set of ergodic measures.
\begin{proposition}\label{propositionMutuallySingular}
  Let $(X,T)$ be a topological dynamical system.
  \begin{enumerate}
    \item Two distinct
      ergodic measures are mutually singular.
    \item Every set of $n$ distinct ergodic measures forms a linearly
      independent set.
      \end{enumerate}
\end{proposition}
\begin{proof}
  Let $\mu,\nu$ be two distinct ergodic measures on $(X,T)$. Since
  $\mu,\nu$ are distinct, there exists an integrable function
  $f$ such that $\int fd\mu\ne\int fd\nu$. Let $U$ 
  be the set of points $x$ such that $f_n(x)$ converges to $\int fd\nu$.
  Then, by Birkhoff Theorem, $\mu(U)=\nu(X\setminus U)=0$. Thus $\mu,\nu$
  are mutually singular.

  Let $\mu_1,\ldots,\mu_n$ be $n$ distinct ergodic measures. Assume
  that $\sum_{k=1}^n\alpha_k\mu_k=0$. Let $1\le i\le n$ be an integer.
  For every $j\ne i$, since $\mu_j,\mu_i$ are mutually singular,
  there is a Borel set $G_J$ such that $\mu_j(G_j)=\mu_i(X\setminus G_j)=0$
  and consequently $\mu_i(G_j)=1$.
  Then  we have $\sum_{k=1}^n\alpha_k\mu_k(\cap_{j\ne i}G_j)=\alpha_i$
  and thus $\alpha_i=0$. This shows that $\mu_1,\ldots,\mu_n$
  are linearly independent.
  \end{proof}
\subsection{Unique ergodicity}
In view of Corollary~\ref{corollaryUniqueisErgodic}, a system with a unique
invariant Borel probability measure  is called \emph{uniquely ergodic}. 
\index{subject}{uniquely ergodic}%
It is called \emph{strictly ergodic}\index{subject}{strictly ergodic}
if it is minimal and uniquely ergodic.

It is easy to give examples of a system which is not uniquely
ergodic when it is not minimal. For example, 
 $\{0^\infty\}\cup \{1^\infty\}$ has clearly two ergodic measures. An example
of a minimal system which is not uniquely ergodic is given
in Exercise~\ref{exerciseNonUniquelyErgodic}.
\begin{theorem}[Oxtoby]\label{theoremOxtoby}
\index{names}{Oxtoby, John C.}
Let $(X,T)$ be a topological dynamical system and $\mu$ be an
invariant Borel probability measure on $(X,T)$. The following conditions are
equivalent.
\begin{enumerate}
\item[\rm(i)] $(X,T)$ is uniquely ergodic.
\item[\rm(ii)] $f_n(x)=\frac{1}{n}f^{(n)}(x)$ converges uniformly on $X$
to $\int fd\mu$ for every $f\in C(X,\R)$.
\item[\rm(iii)] $f_n(x)=\frac{1}{n}f^{(n)}(x)$ converges pointwise to
$\int f d\mu$ for every $f\in C(X,\R)$.
\end{enumerate}
\end{theorem}
\begin{proof}
  (i)$\Rightarrow$(ii). Suppose that (ii) does not hold. We use a diagonal argument
  to find a contradiction. We can
find $\varepsilon>0$, a map $g\in C(X,\R)$ and a sequence $(x_n)$
of points of $X$ such that
\begin{displaymath}
\left|g_n(x_n)-\int g d\mu\right|>\varepsilon.
\end{displaymath}
For every $n\ge 0$, the sum $\delta_{x,n}=\frac{1}{n}\sum_{i=0}^{n-1}\delta_{T^ix}$,
%\marginpar{F:$\delta^{(n)}$ pas defini}
where $\delta_x$ is the Dirac measure,
is a Borel probability measure for which 
\begin{displaymath}
\int fd\delta_{x,n}=f_n(x)
\end{displaymath}
for every $f\in C(X,\R)$ and $x\in X$.
Let $\nu$ be a cluster point of the sequence
$\delta_{x,n}$ for the weak-star topology. Refining $(x_n)$
%\marginpar{F: by the Cantor's diagonal argument}
if necessary, we have
\begin{displaymath}
\lim f_n(x_n)
=\int f d\nu
\end{displaymath}
for every $f\in C(X,\R)$. Hence 
\begin{displaymath}
\left|\int g d\mu - \int g d\nu \right|>\varepsilon
\end{displaymath}
showing that $\mu,\nu$ are distinct invariant probability measures. Thus
(i) does not hold.
%\marginpar{F: il faudrait faire un corollaire a Caratheodory pour les integrales}

(ii)$\Rightarrow$ (iii) is obvious.
%\marginpar{F: rien a montrer dans cette direction} We have
%\begin{displaymath}
%\int f_n d\mu=\frac{1}{n}\sum_{k<n}\int f(T^kx)d\mu=\frac{1}{n}\sum_{k<n}\int f d\mu=\int f d\mu.
%\end{displaymath}
%Thus $\lim f_n(x)=\int f d\mu$.

(iii)$\Rightarrow$ (i). For any invariant Borel probability measure $\nu$, we have
\begin{equation}
\int f_nd\nu=\frac{1}{n}\sum_{i<n}\int f\circ T^id\nu=\int fd\nu.\label{eqCC1}
\end{equation}
On the other hand, by Equations \eqref{eqCC1} and the Dominated Convergence Theorem, we have
\begin{equation}
\int fd\nu = \lim \int f_n d\nu = \int \lim f_n d\nu=\int \left(\int fd\mu \right) d\nu = \int f d\mu\label{eqCC2}
\end{equation}

Thus, by Equations \eqref{eqCC1} and \eqref{eqCC2},
we conclude that $\mu=\nu$.
\end{proof}
Theorem~\ref{theoremOxtoby} is a refinement for uniquely ergodic measures of Birkhoff's ergodic Theorem.
\marginpar{F: Weyl theorem ? }

\subsection{Unique ergodicity of primitive substitution shifts}
\index{subject}{unique ergodicity!of primitive substitution shifts}

In the case of substitution shifts, we have the following important result.
\begin{theorem}\label{theoremMichel}
Primitive substitution shifts are uniquely ergodic.
\end{theorem}
For two words $u,v$, we denote $|u|_v$ the number
of factors of $u$ equal to $v$. By convention,
$|u|_\varepsilon=|u|$.
We say that an infinite word $x\in A^\Z$ has \emph{uniform frequencies}
\index{subject}{uniform!frequencies}
if for every factor $v$ of $x$, all sequences
\begin{displaymath}
(f_{k,v}(x))_n=\frac{|x_k\cdots x_{k+n-1}|_v}{n}
\end{displaymath} 
\index{symbols}{f@$f_{k,v}(x)$}%
converge and tend to the same limit $f_v(x)$ when $n\to\infty$, uniformly in $k$.

Note that $x$ has uniform frequencies if and only if
for every $v\in\cL(x)$ there is an $f_v$
such that for every $\varepsilon$ there is an $N\ge 1$
such that for every $u\in\cL(x)$ of length at least $N$
we have $||u|_v/|u|-f_v|\le\varepsilon$. We express this
property by saying that the frequency of $v$ in $u$
tends to $f_v$ when $|u|\to\infty$ uniformly in $u$.

\begin{proposition}\label{propositionUniformFrequencies}
The following conditions are equivalent
for a minimal shift space $X$.
\begin{enumerate}
\item[\rm(i)] Every $x\in X$ has uniform frequencies.
\item[\rm(ii)] Some $x\in X$ has uniform frequencies.
\item[\rm(iii)] $X$ is uniquely ergodic.
\end{enumerate}
Moreover, in this case the unique invariant measure
is given by $\mu([v])=f_v(x)$ for every $x\in X$.
\end{proposition}
\begin{proof}
(i) implies (ii) is clear. Assume (ii).
For $v\in\cL(X)$, set $\pi(v)=f_v(x)$. This map satisfies
clearly the compatibility conditions~\eqref{eqCompatibility}
and \eqref{eqCompatibilityBis}. Thus there is a unique invariant
Borel probability measure $\mu$ such that $\mu([v])=\pi(v)$. 
Set $g_{k,v}=\charac_{[v]}\circ S^k$. Since
$|x_k\ldots x_{k+n}|_v=\sum_{j<n-|v|}\charac_{[v]}(S^{j+k}x)=(g_{k,v})^{(n-|v|)}(x)$, we have
\begin{displaymath}
(g_{k,v})^{(n)}\rightarrow_{n\to +\infty }\mu([v])=\int \charac_{[v]}d\mu
\end{displaymath}
uniformly in $k$.
But the family of linear combinations all $g_{k,v}$ is dense in $C(X,\R)$. 
Hence 
\begin{displaymath}
\frac{1}{n}g^{(n)}\rightarrow_{n\to +\infty } \int gd\mu
\end{displaymath}
for every continuous function $g$. By Oxtoby Theorem \ref{theoremOxtoby},
this implies that $(X,S)$ is uniquely ergodic.

Finally, by Oxtoby Theorem again, (iii) implies (i).
\end{proof}
Let $\varphi:A\to A^*$ be a primitive substitution. Let $M$ be the
incidence matrix of $\varphi$. Since $M$
is a primitive matrix, by Perron-Frobenius Theorem
(Theorem~\ref{theorem:perron}), the matrix $M$ has a dominant
eigenvalue $\lambda_M$   and positive left and right
eigenvectors  $x=(x_a),y=(y_b)$  relative to $\lambda_M$.
We may assume that $\sum_{a\in A}x_a=1$ and that
$xy=\sum_{a\in A}x_ay_a=1$. With this notation, we prove two lemmas.

\begin{lemma}\label{propositionPhin1}
 For 
every $a,b\in A$, the sequence
\begin{displaymath}
\frac{|\varphi^n(a)|_b}{|\varphi^n(a)|}
\end{displaymath}
converges at geometric rate to $x_b$.
\end{lemma}
\begin{proof}
 We have
\begin{displaymath}
\frac{|\varphi^n(a)|_b}{|\varphi^n(a)|}=
\frac{M_{a,b}^n}{\sum_{b\in A}M_{a,b}^n}\to
\frac{y_ax_b}{\sum_{b\in A}y_ax_b}=x_b
\end{displaymath}
and the result is proved since the 
convergence is at geometric rate by Theorem~\ref{theorem:perron} (iii).
\end{proof}
We will now extend Lemma~\ref{propositionPhin1} to arbitrary words
$u,v\in\cL(X)$. For this, set  $\ell=|v|$. We consider the alphabet $A_\ell$ in one-to-one correspondance with $\cL_\ell(X)$ via $f:\cL_\ell(X)\to A_\ell$
and the $\ell$-block presentation $\varphi_\ell$ of $\varphi$
(see Section~\ref{sectionSubstitutionSystems}).
Recall that, by Proposition~\ref{propositionphi_kPrimitive},
 $\varphi_\ell$ is primitive. 

Let $M_\ell$ be the incidence matrix of $\varphi_\ell$.
Recall that $M_\ell$ has dominant eigenvalue $\lambda_M$
(Proposition~\ref{propositionMk}).
Let $x^{(\ell)}$ be the left eigenvector of $M_\ell$
relative to $\lambda_M$ with coefficients of sum $1$.
\begin{lemma}\label{propositionPhin}
 For 
every $u,v\in\cL(X)$, the sequence
\begin{displaymath}
\frac{|\varphi^n(u)|_v}{|\varphi^n(u)|}
\end{displaymath}
converges at geometric rate to $\pi(v)=x_b^{(\ell)}$ where $b=f(v)$.
\end{lemma}
\begin{proof}
Suppose first that $u=a\in A$.
Let 
$c\in A_\ell$ be such that  $w=f^{-1}(c)$ begins with $a$.
We have by~\eqref{eqLgphikn}
$|\varphi^n(a)=|\varphi_\ell^n(c)|$
for all $n\ge 1$. We extend $f$ naturally
to a map from $\cL_{\ge \ell}(X)$ to $\cL(X^{(\ell)})$.
By Proposition~\ref{propositionphi_kPrimitive}, $\varphi_\ell^n(c)$ is the prefix of length
$|\varphi^n(a)|$ of $f(\varphi^n(w))$. Thus there is a word
$r$ of length $\ell-1$ such that
$\varphi_\ell^n(c)=f(\varphi^n(a)r)$. This shows
that  $|\varphi^n(a)|_v$ and
$|\varphi_\ell^n(c)|_b$ with $b=f(v)$ differ by a constant. Thus
\begin{displaymath}
\frac{|\varphi^n(a)|_v}{|\varphi^n(a)|}\sim\frac{|\varphi_\ell^n(c)|_b}{|\varphi_\ell^n(c)|}.
\end{displaymath}
By Lemma~\ref{propositionPhin1}, the right hand side tends at geometric rate to 
$x_b^{(\ell)}$, whence the conclusion in this case.

Consider finally arbitrary $u,v\in\cL(X)$. We use induction on $|u|$.
The result is true for $|u|=1$ by the previous case. Consider
now $u=u'a$ with $u'\in\cL(X)$ and $a\in A$. 
By induction hypothesis, we have 
$\bigl||\varphi^n(u')|_v-\pi(v)|\varphi^n(u')|\bigr|\le
cr^n|\varphi^n(u')|$ and $\bigl|\varphi^n(a)|_v-\pi(v)|\varphi^n(a)|\bigr|\le cr^n|\varphi^n(a)|$ for some $c>0$ and $r<1$.

Every occurrence of $v$ in $\varphi^n(u)$ is either an occurrence
in $\varphi^n(u')$ or in $\varphi^n(a)$ except $\theta\le|v|$
occurrences which begin in $\varphi^n(a)$ and end in $\varphi^n(a)$.
Thus
\begin{displaymath}
|\varphi^n(u)|_v=|\varphi^n(u')|_v+\theta+|\varphi^n(a)|_v
\end{displaymath}
with $\theta\le|v|$. This implies
\begin{eqnarray*}
\bigl||\varphi^n(u)|_v-\pi(v)|\varphi^n(u)|\bigr|&\le&\bigl||\varphi^n(u')|_v-\pi(v)|\varphi^n(u')|\bigr|+\theta+\\
&&\qquad\quad\bigl||\varphi^n(a)|_v-\pi(v)|\varphi^n(a)|\bigr|\\
&\le&cr^n|\varphi^n(u')|+cr^n|\varphi^n(a)|+|v|\\
&\le&cr^n|\varphi^n(u)|+|v|
\end{eqnarray*}
and thus the conclusion since $|v|/|\varphi^n(u)|$ tends to $0$ at
geometric rate.
\end{proof}
Recall from Chapter~\ref{chapterTopologicalDynamicalSystems} that,
for a morphism $\varphi:A\to A^*$, we denote $|\varphi|=\max_{a\in A}|\varphi(a)|$.
The following lemma is interesting in itself.
\begin{lemma}\label{lemmaQueffelec}
 For every nonempty word $u\in \cL(X)$, there is some $m\ge 0$
and words $v_i,w_i\in \cL(X)$ for $0\le i\le m$, with $v_m$ nonempty such that
\begin{equation}
u=v_0\varphi(v_1)\cdots\varphi^{m-1}(v_{m-1})\varphi^m(v_m)\varphi^{m-1}(w_{m-1})\cdots\varphi(w_1)w_0,
\label{eqQueffelec}
\end{equation}
with $|v_i|,|w_i|\le |\varphi|$.
\end{lemma}
\begin{proof}
We use induction on $|u|$. The result is true if $|u|<|\varphi|$
choosing $m=0$ and $v_0=u$.
Otherwise, by definition of $\cL(X)$, there exists a nonempty
word $u'\in \cL(X)$ such that
$u=v_0\varphi(u')w_0$. Choosing $u'$ of maximal length, we have moreover 
$|v_0|,|w_0|\le|\varphi|$. By induction hypothesis, we have
a decomposition \eqref{eqQueffelec} for $u'$, that is
\begin{displaymath}
u'=v'_0\varphi(v'_1)\cdots\varphi^{m-1}(v'_{m-1})\varphi^m(v'_m)\varphi^{m-1}(w_{m-1}')\cdots\varphi(w_1')w_0'.
\end{displaymath}
In this way, we obtain
\begin{displaymath}
u=v_0\varphi(u')w_0=v_0\varphi(v'_0)\cdots\varphi^{m}(v'_{m-1})\varphi^{m+1}(v'_m)\varphi^{m}(w'_{m-1})\cdots\varphi(w_0')w_0.
\end{displaymath}
which is of the form \eqref{eqQueffelec}.
\end{proof}
\begin{proofof}{of Theorem~\ref{theoremMichel}}
Let $x\in X$.
By Proposition~\ref{propositionUniformFrequencies} it is enough to prove
that $x$ has uniform frequencies. We will prove that $|x_k\cdots x_{k+n-1}|_v/n$
converges to $\pi(v)$ uniformly in $k$, where $\pi(v)$ is as
in Lemma~\ref{propositionPhin}.
Set $u=x_k\cdots x_{k+n-1}$. By Lemma \ref{lemmaQueffelec}, we have
\begin{displaymath}
|u|_v=\sum_{i\le m}|\varphi^i(v_i)|_v+\sum_{i< m}|\varphi^i(w_i)|_v+\theta
\end{displaymath}
where $\theta$ is the number of occurrences of $v$ which overlap
more than one $v_i,w_i$. Thus $\theta\le 2m|v|$.
By Lemma~\ref{propositionPhin} there are $c>0$ and $r<1$
such that 
\begin{equation}
||\varphi^i(v_i)|_v-\pi(v)|\varphi^n(v_i)||\le cr^i\label{eqMajv_i}
\end{equation}
for $0\le i\le m$ (and the analogue inequality for the $w_i$).
Since
\begin{displaymath}
n=\sum_{i\le m}|\varphi^i(v_i)|+\sum_{i<m}|\varphi^i(w_i)|
\end{displaymath}
we obtain
\begin{eqnarray*}
||u|_v-\pi(v)n|&=&||u|_v-\pi(v)(\sum_{i\le m}|\varphi^i(v_i)|+\sum_{i<m}|\varphi^i(w_i))|)\\
&\le&\sum_{i\le m}||\varphi^i(v_i)|_v-\pi(v)|\varphi^i(v_i)||\\
&&\quad
+\sum_{i<m}||\varphi^i(w_i)|_v-\pi(v)|\varphi^i(w_i)||+2m|v|\\
&\le&2c\sum_{i\le m}r^i+2m|v|\\
&\le&dr^m+2m|v|.
\end{eqnarray*}
But $m/n$ tends to $0$ when $n$ tends to infinity. Indeed, by
\eqref{eqFabien2}, there is a constant $e>0$ such that
$|\varphi^n(u)|\ge e|u|\lambda_M^n$ for $n$ large enough. This implies
that $n\ge|\varphi^m(v_m)|\ge e|v_m|\lambda_M^m$. Hence
$m/n\le m/e\lambda_M^m|v_m|$ which tends to $0$.

This shows that $|u|_v/n$ converges to $\pi(v)$ when $n$ tends to infinity 
independently of $u$, which concludes the proof.
\end{proofof}
\subsection{Computation of the unique invariant Borel probability  measure}
The computation of the unique invariant Borel probability measure can be done using the matrices $M_\ell$
of $\ell$-block presentations of the substitution, as we shall see below.

Let $\varphi:A\rightarrow A^*$ be a primitive substitution
and let $X$ be the associated shift space.
Let $\mu$ be the unique invariant Borel probability measure
 on $X$. For $\ell\ge 1$, let $A_\ell$ be an alphabet in bijection
with $\cL_\ell(X)$ via $f:\cL_\ell(X)\to A_\ell$
and let $x^{(\ell)}$ be the unique positive
row eigenvector of the matrix $M_\ell$ corresponding to the maximal
eigenvalue $\lambda$ and such that the sum of its
components is $1$. 
\begin{proposition}
For every $v\in \cL_\ell(X)$,
\begin{equation}
\mu([v])=x^{(\ell)}_c \label{eqInvariantMeasure}
\end{equation}
with $c=f(v)$.
\end{proposition}
\begin{proof}
By Theorem~\ref{propositionUniformFrequencies}, $\mu([v])$ is
equal to the frequency $f_v(x)$ of the word $v$ in any $x\in X$.
By Lemma \ref{propositionPhin}, the frequency of $v$
in $\varphi^n(a)$ tends to $x^{(\ell)}_c$ when
$n\to\infty$ and thus we have $f_v(x)=x^{(\ell)}_c$.
\end{proof}

We develop below the cases of the Fibonacci and of the Thue-Morse morphisms.
For the Fibonacci shift, we already now that it is uniquely
ergodic since this is true of every Sturmian shift (Exercise~\ref{exerciseErgodicityRotations}).
\begin{example}
Let $X$ be the  Fibonacci shift.
\index{subject}{Fibonacci!shift}\index{subject}{shift space!Fibonacci}%
 Since the Fibonacci
substitution is primitive, there is a unique invariant
  probability measure on $X$. Its values on the cylinder
$[w]$
defined by words $w$ of length
  at most $4$ are shown on Figure~\ref{figProbaFibo} with
  $\rho=(\sqrt{5}-1)/2$.

This is consistent with the value of the eigenvector 
of $M_2$
\begin{displaymath}
v^{(2)}=\begin{bmatrix}2\rho-1&1-\rho&1-\rho\end{bmatrix}
\end{displaymath}
Note that the Fibonacci shift is Sturmian of slope $\alpha=(3-\sqrt{5})/2$.
The unique invariant Borel probability measure $\mu$ on the Fibonacci shift
can also be computed using the natural representation $\gamma_\alpha$
(Exercise~\ref{exerciseErgodicityRotations})
giving $\mu([a])=1-\alpha$. Consistently with the above,
$\rho=1-\alpha$.
\end{example}
\begin{figure}[hbt]
\centering
\tikzset{node/.style={minimum size=0.1cm,inner sep=0pt}}
	\tikzset{title/.style={minimum size=0.5cm,inner sep=0pt}}
\begin{tikzpicture}
\node[node](1) at (0,2){$1$};
\node[node](a) at (2,3){$\rho$};\node[node](b) at (2,1){$1-\rho$};
\node[node](aa) at (4,4){$2\rho-1$};\node[node](ab) at (4,2.5){$1-\rho$};\node[node](ba) at (4,0.5){$1-\rho$};
\node[node](aab) at (6,4){$2\rho-1$};\node[node](aba) at (6,2.5){$1-\rho$};
\node[node](baa) at (6,1){$2\rho-1$};\node[node](bab) at (6,0){$2-3\rho$};
\node[node](aaba) at (8,4){$2\rho-1$};\node[node](abaa) at (8,3){$2\rho-1$};\node[node](abab) at (8,2){$2-3\rho$};
\node[node](baab) at (8,1){$2\rho-1$};\node[node](baba) at (8,0){$2-3\rho$};

\draw[above](1) edge node{$a$} (a);\draw[above](1) edge node{$b$} (b);
\draw[above](a) edge node{$a$} (aa);\draw[above](a) edge node{$b$} (ab);\draw[above](b) edge node{$a$} (ba);
\draw[above](aa) edge node{$b$} (aab);\draw[above](ab) edge node{$a$} (aba);
\draw[above](ba) edge node{$a$} (baa);\draw[above](ba) edge node{$b$} (bab);
\draw[above](aab) edge node{$a$} (aaba);\draw[above](aba) edge node{$a$} (abaa);\draw[above](aba) edge node{$b$} (abab);
\draw[above](baa) edge node{$b$} (baab);\draw[above](bab) edge node{$a$} (baba);

\end{tikzpicture}
\caption{The invariant Borel probability measure on the 
  Fibonacci set.}\label{figProbaFibo}
\end{figure}
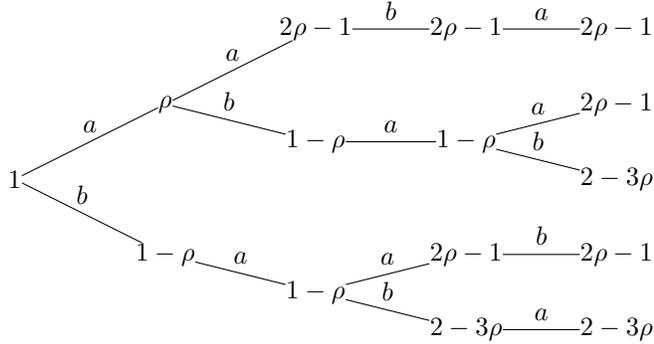
\begin{example}\label{exampleMeasureMorse}
Consider the Thue-Morse morphism. The matrix $M_2$
is
\begin{displaymath}
M_2=\begin{bmatrix}0&1&1&0\\0&1&0&1\\1&0&1&0\\0&1&1&0\end{bmatrix}
\end{displaymath}
The unique invariant Borel probability
measure on the Thue-Morse shift is shown in Figure~\ref{figProbaMorse}.
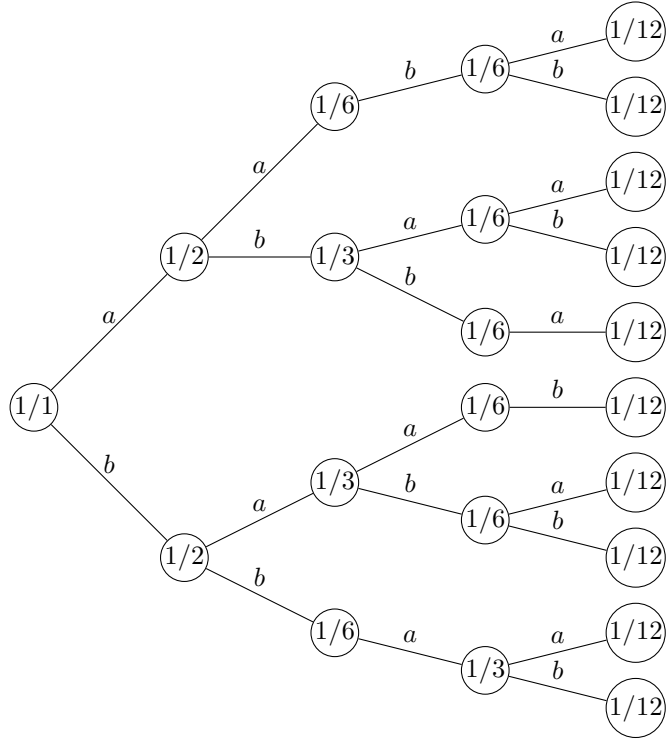
\begin{figure}[hbt]
\centering
\tikzset{node/.style={circle,draw,minimum size=0.1cm,inner sep=0pt}}
	\tikzset{title/.style={minimum size=0.5cm,inner sep=0pt}}
\begin{tikzpicture}
\node[node](1) at (0,5){$1/1$};
\node[node](a) at (2,7){$1/2$};\node[node](b) at (2,3){$1/2$};
\node[node](aa) at (4,9){$1/6$};\node[node](ab) at (4,7){$1/3$};\node[node](ba) at (4,4){$1/3$};\node[node](bb) at (4,2){$1/6$};
\node[node](aab) at (6,9.5){$1/6$};\node[node](aba) at (6,7.5){$1/6$};\node[node](abb) at (6,6){$1/6$};
\node[node](baa) at (6,5){$1/6$};\node[node](bab) at (6,3.5){$1/6$};\node[node](bba) at (6,1.5){$1/3$};
\node[node](aaba) at (8,10){$1/12$};\node[node](aabb) at (8,9){$1/12$};\node[node](abaa) at (8,8){$1/12$};\node[node](abab) at (8,7){$1/12$};
\node[node](abba) at (8,6){$1/12$};
\node[node](baab) at (8,5){$1/12$};\node[node](baba) at (8,4){$1/12$};\node[node](babb) at (8,3){$1/12$};\node[node](bbaa) at (8,2){$1/12$};
\node[node](bbab) at (8,1){$1/12$};

\draw[above](1) edge node{$a$} (a);\draw[above](1) edge node{$b$} (b);
\draw[above](a) edge node{$a$} (aa);\draw[above](a) edge node{$b$} (ab);\draw[above](b) edge node{$a$} (ba);\draw[above](b) edge node{$b$} (bb);
\draw[above](aa) edge node{$b$} (aab);\draw[above](ab) edge node{$a$} (aba);\draw[above](ab) edge node{$b$} (abb);
\draw[above](ba) edge node{$a$} (baa);\draw[above](ba) edge node{$b$} (bab);\draw[above](bb) edge node{$a$} (bba);
\draw[above](aab) edge node{$a$} (aaba);\draw[above](aab) edge node{$b$} (aabb);\draw[above](aba) edge node{$a$} (abaa);\draw[above](aba) edge node{$b$} (abab);\draw[above](abb) edge node{$a$} (abba);
\draw[above](baa) edge node{$b$} (baab);\draw[above](bab) edge node{$a$} (baba);\draw[above](bab) edge node{$b$} (babb);\draw[above](bba) edge node{$a$} (bbaa);\draw[above](bba) edge node{$b$} (bbab);

\end{tikzpicture}
\caption{The invariant Borel probability measure on the 
  Thue-Morse shift.}\label{figProbaMorse}
\end{figure}
The values on the words of length $2$ is consistent with
the value of the eigenvector of $M_2$ which is
\begin{displaymath}
v^{(2)}=\begin{bmatrix}1/6&1/3&1/3&1/6\end{bmatrix}
\end{displaymath}
\end{example}

%%%%%%%%%%%%%%%%%%%%%%%
%\section{Invariant measures and states}\label{sectionInvariantStates}
%There is an important connexion between the dimension group $H (X,T,\Z )$
%of a topological dynamical system $(X,T)$ and the invariant probability
%measures on $(X,T)$.
\section{Invariant measures and states}\label{sectionInvariantStates}

We now come to an essential point and relate invariant measures and coboundaries.

\begin{proposition}\label{propositionCoboundaryInt}
Let $\mu$ be an invariant measure on the dynamical system
$(X,T)$. For every coboundary $f\in \partial_TC(X,\R)$, one has
$\int fd\mu=0$.
\end{proposition}
\begin{proof}
When $\mu$ is a $T$-invariant Boel probability
measure we have $\mu\circ T^{-1}=\mu$.
Thus, for $f=\partial_T g$, by the change of variable
formula (see Appendix~\ref{appendixMeasureIntegration}), 
$\int\partial_T g d\mu=\int g\circ T d\mu-\int gd\mu=0$.
\end{proof}
Proposition~\ref{propositionCoboundaryInt} implies that for all $T$-invariant Borel probability measures $\mu$ on $(X,T)$ the map $f\mapsto\int fd\mu$ defines a group
homomorphism
\begin{displaymath}
\alpha_\mu:H(X,T,\Z)\rightarrow\R
\end{displaymath}
\index{symbols}{alpha@$\alpha_\mu$}%
\begin{proposition}
  Let $(X,T)$ be a recurrent dynamical system.
The map $\alpha_\mu$ is a morphism of unital ordered groups 
from $K^0(X,T)$ to $(\R,\R_+,1)$.
\end{proposition}
\begin{proof}
By definition of $\int fd\mu$, we have $\alpha_\mu(H(X,T,\Z_+))\subset \R_+$.
Moreover,
$\alpha_\mu(\mathbf{1}_X)=\mu(X)=1$ since $\mu$ is a probability measure. 
\end{proof}
 The proof of the following statement
uses  the Carath\'eodory extension Theorem
\index{subject}{Carath\'eodory Extension Theorem}%
\index{subject}{Theorem!Carath\'eodory Extension}%
\index{names}{Carath\'eodory, Constantin}%
(see Appendix~\ref{appendixMeasureIntegration}).

\begin{theorem}\label{propositionKerov}
Let $(X,T)$ be a recurrent dynamical system.
The map $\mu\mapsto\alpha_\mu$ is a bijection from
the space of $T$-invariant Borel probability measures on $(X,T)$
onto the set of states of $K^0(X,T)$.
\end{theorem}
\begin{proof}
Let $\alpha$ be a state on $K^0(X,T)$. For a clopen set $U\subset X$, we set
$\phi(U)=\alpha(\charac_U)$ where $\charac_U$ is the image
in $H(X,T,\Z)$ of the characteristic function of $U$. Thus
$\phi(U)\ge 0$ for every clopen set $U$, and $\phi(U\cup V)=\phi(U)+\phi(V)$ if
$U,V$ are disjoint clopen sets. Since $X$  is a Cantor
space, its topology is generated by the clopen sets.
Thus, there is, by the Carath\'eodory Theorem, a unique 
probability measure $\mu$ on $X$ such that
$\mu(U)=\phi(U)$ for every clopen set $U$.
This already shows that $\mu \mapsto \alpha_\mu$ is one-to-one.
For every clopen set $U$, the difference between the characteristic
functions of $U$ and $T^{-1}U$ is a coboundary. Thus these
functions have the same image in $H(X,T,\Z)$ and $\phi(U)=\phi(T^{-1}U)$.
It follows that $\mu$ is $T$-invariant. Moreover,
by construction $\alpha=\alpha_\mu$ since 
$\alpha_\mu(\charac_U)=\int\charac_Ud\mu=\mu(U)=\alpha(\charac_U)$
for every clopen set $U$. This shows that the map
$\mu\mapsto\alpha_\mu$ is  onto.
\end{proof}

Observe that Theorem~\ref{propositionKerov} allows us to give a
very simple
proof of Theorem~\ref{theoremMichel}. Indeed, if $X$ is a 
primitive substitution shift, then $K^0(X,S)$ is the limit
of a stationary system defined by a primitive matrix and thus
has a unique state by Proposition~\ref{propositionUniqueState}.
We conclude by Theorem~\ref{propositionKerov} that
$X$ is uniquely ergodic.\index{subject}{unique ergodicity!of primitive substitution shifts}

Observe also that Theorem~\ref{propositionKerov} implies
that there is a bijection between the invariant Borel probability
measures on a minimal system $(X,T)$ and an induced system
$(U,T_U)$ for $U\subset X$ clopen. Indeed, by Proposition
\ref{propositionIuRu}, the ordered groups $K^0(X,T)$
and $K^0(U,T_U)$ are isomorphic. A direct proof
of the bijection between invariant Borel probability measures
is given in Exercise~\ref{exerciseMeasureInduced}.

We note another important consequence of Theorem~\ref{propositionKerov}.

\begin{proposition}\label{nbErgodicMeasures}
  Let $(X,T)$ be a recurrent dynamical system. If the group $K^0(X,T)$ has
  rank $n$, there are at most $n$ ergodic measures on $(X,T)$.
  If additionnally $K^0(X,T)$ is a simple dimension group, this number is at most $n-1$.
\end{proposition}
\begin{proof}
  Let $m$ be the number of ergodic measures.
  Let us first show  that $m\le n$.
  Indeed, let $V$ be the vector space generated by the
  ergodic measures on $(X,T)$ and $W$ the vector
  space generated by the states on $K^0(X,T)$.
  By Proposition~\ref{propositionMutuallySingular},
  the $m$ ergodic measures form a basis
  of $V$ and thus $V$ has dimension $m$.
  Since $W$ is included is formed of linear
  forms on $K^0(X,T)$, it has dimension at most $n$.
  Since the map $\mu\to \alpha_\mu$ defines  a linear injective map from
  $V$ to $W$, the conclusion follows.

  Suppose now that $K^0(X,T)$ is a simple dimension group.
  Assume that $m=n$. Identify the group $H(X,T,\Z)$ with $\Z^n$. Let
  $f_1,f_2,\ldots,f_n$ be $n$ states generating a space of dimension $n$.
  Let $v_1,v_2,\ldots,v_n\in \R^n$ be a dual basis. Fix $g\in H^+(X,T,\Z)$
  and set $g=\sum_{i=1}\alpha_i v_i$. Let $D=\{h\in G\mid 0\le h\le g\}$.
  By Proposition~\ref{propositionEffros}, we have
  \begin{eqnarray*}
    D=\{h\in\Z^n\mid 0\le f_i(h)\le \alpha_i\mbox{ for $1\le i\le n$}\}
        \end{eqnarray*}
  Since the last set is contained in the compact set
  \begin{displaymath}
    K=\{\sum_{i=1}^n\beta_iv_i\mid 0\le\beta_i\le\alpha_i\mbox{ for $1\le i\le n$}\}\end{displaymath}
    it is finite. This implies that $D$ contains minimal nonzero elements.
    Such an element generates a nonzero ideal, contradicting the
    hypthesis that $G$ is simple.
\end{proof}
We will see later many examples where the above hypotheses are
satisfied.

\subsection{Dimension groups of Sturmian shifts}
As an illustration of Theorem~\ref{propositionKerov}, let us prove the following statement,
which gives the form of the cohomology group of Sturmian shifts.
\index{subject}{Sturmian!shift space!dimension group}%
\index{subject}{dimension group!of Sturmian shift}%
In the next statement the group $\Z+\alpha\Z$ of real numbers
of the form
$x+\alpha y$ is considered as a unital ordered group for the
order induced by $\R$ and the order unit $1$.
\begin{theorem}\label{theoremDGSturm}
Let $X$ be a Sturmian shift of slope $\alpha$. Then
$K^0(X,S)=\Z+\alpha\Z$.
\end{theorem}
\begin{proof}
We have already seen (Example~\ref{exampleSturmian1}) that
$H(X,S,\Z)=\Z^2$ for every Sturmian shift. Let
$\{0,1\}$ be the alphabet of $X$. We can identify the group $H(X,S,\Z)$ with
the pairs $(x,y)\in\Z^2$ via the map $(x,y)\mapsto x[0]+y[1]$.
We have seen (Exercise~\ref{exerciseErgodicityRotations})
that $(X,S)$ is uniquely ergodic and that the unique
invariant Borel probability measure $\mu$ is such that
$\mu([0])=1-\alpha$. By Theorem~\ref{propositionKerov},
the ordered group $K^0(X,S)$ has a unique state
which is $\alpha_\mu$. Thus, by Proposition~\ref{propositionEffros},
$H^+(X,S,\Z)=\{(x,y)\in\Z^2\mid (1-\alpha)x+\alpha y>0\}\cup\{0\}$,
the order unit being $(1,1)$.
Thus, the map $(x,y)\to (x,y-x)$ identifies $K^0(X,T)$ and $\Z+\alpha\Z$.
\end{proof}

\begin{example}\label{exampleK0Fibonacci}
Let $X$ be the Fibonacci shift.
\index{subject}{Fibonacci!shift}%
It is a Sturmian shift of slope $\alpha=(3-\sqrt{5})/2$
(Example~\ref{exampleSlopeFibonacci}). Thus
$K^0(X,S)=\Z[\alpha]$.
\end{example}

As an application of Theorem~\ref{theoremDGSturm}, we prove the following
result.
\begin{theorem}\label{theoremConjugacySturmian}
Two Sturmian shifts are never conjugate unless they differ
by a permutation of the two letters.
\end{theorem}
\begin{proof}
Let $X$ and $Y$ be Sturmian shifts of slopes $\alpha,\beta$
respectively. If they are conjugate, their
cohomology groups are isomorphic
and thus, by Theorem~\ref{theoremDGSturm}, the groups $\Z+\alpha\Z$
and $\Z+\beta\Z$ are isomorphic. Thus there are $a,b,c,d\in\Z$
such that $\beta=a+b\alpha$ and $\alpha=c+d\beta$. Then
$\alpha=c+ad+bd\alpha$ implies $bd=1$. Thus either $\alpha=\beta$
or $\alpha=1-\beta$.
\end{proof}
This statement is a striking illustration
of the  power of the cohomology
group. It would probably be very difficult to prove directly
the above statement by a direct argument.

We can also use Theorem~\ref{propositionKerov} to give
the form of the image subgroup  and the infinitesimal
subgroup  for a minimal
Cantor system $(X,T)$. 

The image subgroup of $K^0(X,T)$, denoted $I(X,T)$,
\index{subject}{image subgroup}\index{symbols}{I@$I(X,T)$}%
 is given by
\begin{equation}
I(X,T)=\cap_{\mu\in\M(X,T)}\left\{\int fd\mu\mid f\in C(X,\Z)\right\}.
\label{eqImageSubgroupDynam}
\end{equation}
where $\M(X,T)$\index{symbols}{M@$\M(X,T)$} denotes the set of invariant Borel probability measures
on $(X,T)$.
This follows directly from the definition of the image subgroup
given in Equation~\eqref{eqImageSubgroup} 
by Theorem \ref{propositionKerov}.

The infinitesimal subgroup of $K^0(X,T)$, denoted $\Inf(X,T)$
\index{symbols}{Inf@$\Inf(X,T)$} is given by
\begin{equation}
\Inf(X,T)=\{[f]\in H(X,T,\Z)\mid \int f d\mu=0 \mbox{ for all $\mu\in\M(X,T)$}\}.\label{eqInfSubgroup}
\end{equation}
If $(X,T)$ is uniquely ergodic, we have $K^0(X,T)/\Inf(K^0(X,T))=I(X,T)$.
\begin{example}
Let $X$ be the Fibonacci shift as in the previous example.
Then $\Inf(X,S)=\{0\}$ and $I(X,S)=\Z+\alpha\Z$.
\end{example}
%%%%%%%%%%%%%%%%%%%%

%%%%%%%%%%%%%%%%%%%
\section{Exercises}
\exosection{Section~\protect{\ref{sectionCoboundaries}}}

\begin{exercise}\label{exerciseNegativeCocycles}
Let $(X,T)$ be a topological dynamical system. For $f\in C(X,\R)$
extend the definition of $f^{(n)}(x)$ to negative indexes
by defining for $n\ge 1$ and $x\in X$
\begin{equation}
f^{(-n)}(x)=-f^{(n)}(T^{-n}x).\label{eqNegativeCocycle}
\end{equation}
Show that, with this definition, the \emph{cohomological
equation}\index{subject}{cohomological equation}
\begin{equation}
f^{(n+m)}(x)=f^{(n)}(x)+f^{(m)}(T^nx)\label{eqCohomologicalEquation}
\end{equation}
holds for every $n\in\Z$ and $x\in X$. 
\end{exercise}

\begin{exercise}\label{exerciseHomotopy}
A \emph{homotopy}\index{subject}{homotopy} between topological spaces
$X,Y$ is a family $f_t:X\to Y$ ($t\in[0,1]$) of maps 
such that the associated map $F:X\times [0,1]$
given by $F(x,t)=f_t(x)$ is continuous. One says that
 $f_0,f_1:X\to Y$ 
are \emph{homotopic}
\index{subject}{homotopic functions}\index{subject}{function!homotopic}%
if there is a homotopy $f_t$ connecting them.
Show that homotopy is an equivalence
relation compatible with composition.
\end{exercise}
\begin{exercise}\label{exerciseHomotopy2}
A continous function 
is \emph{nullhomotopic}\index{subject}{nullhomotopic}
 if it is homotopy equivalent to a constant function.
Show that  any function $f:X\to [a,b]$
from $X$ to a closed interval $[a,b]$ of $\R$ 
is
nullhomotopic.
\end{exercise}
\begin{exercise}\label{exerciseSuspension}
A \emph{continuous flow}\index{subject}{continuous!flow}
\index{subject}{flow!continuous}%
is a pair $(X,(T_t)_{t\in\R})$ of a compact metric space $X$
and a family $(T_t)_{t\in T}$ of homeomorphisms $T_t:X\to X$
such that 
\begin{enumerate}
\item[(i)] the map $(x,t)\mapsto T_t(x)$ is continuous from
$X\times \R$ to $X$.
\item[(ii)] $T_{t+s}=T_t\circ T_s$ for all $s,t\in\R$.
\end{enumerate}
For every topological dynamical system $(X,T)$, one can build a continuous
flow called the \emph{suspension flow}
\index{subject}{suspension flow}\index{subject}{flow!suspension}%
 over $(X,T)$ as follows.
Consider the quotient $\tilde{X}$ of $X\times \R$ by the
equivalence relation  which identifies
$(x,s+1)$ and $(Tx,s)$ for all $x\in X$ and $s\in\R$. Denote
by $[(x,s)]$ the equivalence class of $(x,s)$. Then define $T_t$ on the quotient by
\begin{displaymath}
T_t([(x,s)])=[(x,s+t)].
\end{displaymath}
Show that we obtain in this way a continuous flow.
\end{exercise}

\begin{exercise}\label{exerciseSuspensionPeriodic}
Show that the suspension flow of any periodic system
can be identified with the torus $\T=\R/\Z$.
\end{exercise}

\begin{exercise}\label{exerciseFlowEquivalence}
An \emph{equivalence}\index{subject}{equivalence!of flows}
between two continuous flows $(X,(T_t)_{t\in\R})$ and $(X',(T'_t)_{t\in \R})$
is a homeomorphism $\pi:X\to X'$ which maps orbits of $T_t$
to orbits of $T'_t$ in an orientation preserving way, that is, for all $x\in X$,
$\pi(T_t(x))=T'_{f_x(t)}(\pi(x))$ for some monotonically increasing
$f_x:\R\to\R$. Two flows
are \emph{equivalent}\index{subject}{equivalent!flows}
if there is an equivalence from one to the other.

Two (ordinary) topological dynamical systems are \emph{flow equivalent}
\index{subject}{flow!equivalence} 
if their suspension flows are flow equivalent.
Show that two equivalent dynamical systems are flow equivalent
but that the converse is false.
\end{exercise}

\begin{exercise}\label{exerciseCechCohomology}
Let $(X,T)$ be a topological dynamical system and let
$(\tilde{X},(T_t)_{t\in\R})$ be its suspension flow (as defined
in Exercise~\ref{exerciseSuspension}). The first \emph{\v{C}ech cohomology
group}\index{subject}{Cech@\v{C}ech cohomology group}
\index{names}{Cech@\v{C}ech, Eduard}
 of $\tilde{X}$, denoted $H^1(\tilde{X},\Z)$,
 is the group of continuous maps from $\tilde{X}$
to the torus $\T=\R/\Z$ modulo the group of nullhomotopic maps. 

Show that $H^1(\tilde{X},\Z)$ is
isomorphic to $H(X,T,\Z)$ (hint: 
Consider the map $\pi:C(X,\Z)\to C(\tilde{X},\T)$ which associates to 
$f\in C(X,\Z)$ the map $\pi(f)\in C(\tilde{X},\T)$ defined by
$\pi(f)(x,s))=\tau(f(x)s)$. Show that $\pi$ induces
an isomorphism from $H(X,T,\Z)$ onto $H^1(\tilde{X},\Z)$).
\end{exercise}

\exosection{Section~\protect{\ref{sectionGH}}}
\begin{exercise}\label{exerciseGH}
Show that if $(X,T)$ is a transitive system, then
for any $f\in C(X,\R)$ two solutions of the equation
$\partial g=f$ differ by a constant.
\end{exercise}
\begin{exercise}\label{exerciseProofKH1}
Let $(X,T)$ be a topological dynamical system. For $f\in C(X,\R)$,
assume that
the sequence $(f^{(n)}(x))_{n\ge 0}$ is uniformly bounded. Set
\begin{displaymath}
g(x)=\sup_{n\in\Z}f^{(n)}(x)
\end{displaymath}
where $f^{(n)}$ is defined for all $n\in\Z$ as in Exercise~\ref{exerciseNegativeCocycles}. Show that $\partial g=f$.
\end{exercise}
\begin{exercise}\label{exerciseProofKH2}
Let $(X,T)$ be a minimal system. Let $f\in C(X,\R)$ and $x_0\in X$
be such that $(f^{(n)}(x_0))_{n\ge 0}$ is bounded. Show that $f^{(n)}(x)$
is bounded for all $x\in X$ and
that $g(x)=\sup_{n\in\Z}f^{(n)}(x)$ is a continuous function such that
$\partial g=f$ (hint: use Exercise~\ref{exerciseProofKH1}).
\end{exercise}
%\begin{exercise}\label{exerciseGH2}
%Let $(X,T)$ be a dynamical system and let $f\in C(X,T)$
%be such that $f^{(n)}(x)$ is bounded uniformly in $n$ and $x$.
%Then $f=\partial g$ for
%\begin{displaymath}
%g(x)=\sup_{n\ge 1}(-f^{(n)}(x)).
%\end{displaymath}
%\end{exercise}
\begin{exercise}\label{exerciseProofHerman}
Let $(X,T)$ be a minimal dynamical system and $f\in C(X,\R)$ be such that
$(f^{(n)}(x_0))_{n\ge 0}$ is bounded for some $x_0\in X$. 

Let $(Y,S)$ be the dynamical system formed of $Y=X\times\R$
with $S(x,t)=(Tx,t+f(x))$. In this way, $S^n(x,t)=(T^nx,t+f^{(n)}(x))$
for all $n\ge 0$. 
\begin{enumerate}
\item Show that the closure $E$ of the set $\{S^n(x_0,0)\mid n\ge 0\}$
is compact and $S$-invariant.
\item Let $K$ be a minimal closed $S$-invariant nonempty subset of $E$.
Show that $K=\{(x,g(x))\mid x\in X\}$ for some 
$g\in C(X,\R)$.
\item Show that $\partial g=f$.
\end{enumerate}
\end{exercise}

%\exosection{Section~\protect{\ref{sectionInduced}}

\exosection{Section~\protect{\ref{chapterOrderedGroupRecurrent}}}
\begin{exercise}\label{exercise[ab]}
  Show that for the shift space $X$ of Example~\ref{examplea*b*3},
  $\charac_{[ab]}$ is not a coboundary. 
  \end{exercise}
\begin{exercise}\label{exerciseG_n}
  Prove directly (that is, without using Corollary~\ref{corollary1})
  that if $X$ is recurrent, $\G_n(X)$ is an ordered group.
  \end{exercise}

\exosection{Section~\protect{\ref{sectionInvariant}}}

\begin{exercise}\label{exerciseMarkovKakutani}
The \emph{Markov-Kakutani fixed point theorem}
\index{subject}{Markov-Kakutani fixed point Theorem}%
\index{subject}{Theorem!Markov-Kakutani}%
states that if $V$ is a  topological vector space
and $K$ is a convex and compact subset of $V$, then every
continuous linear map $T$ 
mapping $K$ into itself has a fixed point in $K$.

Use this theorem to give a proof of Theorem~\ref{theoremKrylovBogolioubov}.
Hint: Consider the map $T^*:\M_X\to \M_X$ defined
by $T^*\mu=\mu\circ T^{-1}$.
\end{exercise}

\begin{exercise}\label{exerciseInvariantFunction}
  Prove Proposition~\ref{propositionCondEquiv}. Hint: set $V=\cap_{n\ge 1}\cup_{i\ge n}T^{-i}U$.
  Show that $T^{-1}V=V$ and $U=V\bmod\mu$.
\end{exercise}

\begin{exercise}\label{exerciseBernoulliErgodic}
  Show that if an invariant measure $\mu$ on $(X,T)$ is ergodic, then
  \begin{equation}
    \lim_{n\to\infty}\frac{1}{n}\sum_{i=0}^{n-1}\mu(U\cap T^{-i}V)=\mu(U)\mu(V)\label{eqCriteriumErgodicity}
  \end{equation}
  for all Borel sets $U,V$ and
  that the converse holds provided \eqref{eqCriteriumErgodicity} is true
  for all clopen sets $U,V$.
  Use this to prove that Bernoulli measures are ergodic.
\end{exercise}

\begin{exercise}\label{exerciseMarkovErgodic}
  Let $\mu$ be an invariant measure defined by a Markov chain $(v,P)$
  with $vP=v$ and $v>0$. Show that $\mu$ is ergodic if and only if the matrix $P$ is 
irreducible.
  \end{exercise}

\begin{exercise}\label{exerciseErgodicRecurrent}
Prove Proposition~\ref{propositionErgodicRecurrent}
(hint: prove that the set of recurrent points has measure $1$).
\end{exercise}

\begin{exercise}\label{exerciseErgodicityRotations}
Use Oxtoby's Theorem (Theorem~\ref{theoremOxtoby})
to prove that irrational rotations are uniquely ergodic.
\index{subject}{unique ergodicity!of irrational rotations}%
 Conclude that Sturmian shifts are uniquely ergodic
and that the unique invariant Borel probability measure $\mu$
on a Sturmian shift of slope $\alpha$ satisfies $\mu([0])=1-\alpha$ .
Hint: use the fact that trigonometric polynomials are
dense in the space $\C(S^1,\C)$. 
\end{exercise}

\begin{exercise}\label{exerciseNonUniquelyErgodic}
Let $(k_i)_{i\ge 0}$ be a sequence of positive integers
such that $k_i$ divides $k_{i+1}$.
Let $x\in\{0,1\}^\Z$ be the Toeplitz sequence 
\index{subject}{Toeplitz!sequence} with periodic structure
$(k_i)_{i\ge 0}$ defined as follows. For $i\ge 1$, let
\begin{displaymath}
E_i=\cup_{m\in\Z}\{n\in\Z\mid |n-mk_i|\le k_{i-1}\}
\end{displaymath}
For $n\in\Z$, let $p(n)$ be the least integer $p$ such that
$n\in E_{p}$. Set $x_n\equiv p(n)\bmod 2$. Show that if 
\begin{equation}
\sum_{i\ge 1}\frac{k_{i-1}}{k_i}\le\frac{1}{12}\label{eqOxtoby}
\end{equation}
then $x$ is not uniquely ergodic.
\end{exercise}

\begin{exercise}\label{exerciseUnitary}

Show that if
  $(X,T)$ is a topological dynamical system and
$\mu$ is an invariant Borel probability measure,
then $(X,T,\mu)$ is a measure-theoretic dynamical system.
\end{exercise}
\begin{exercise}\label{exerciseUnitary2}
Let $(X,T,\mu)$ be a measure-theoretic dynamical system
with $T$ invertible. Let $H= {L}^2(X)$ be the Hilbert space of 
real valued square
integrable functions on $X$ (modulo a.e. vanishing functions).

Show
that the operator $U$ defined by $Uf=f\circ T$ is a unitary
operator from $H$  to itself.
\end{exercise}
\begin{exercise}\label{exerciseGH3}
Let $(X,T,\mu)$ be a measure-theoretic
dynamical system. Show that an
element $f$ of the Hilbert space $H=L^2(X)$ is a coboundary
of the form
 $f=Ug-g$ for some $g\in H$, with $U$ as in Exercise~\ref{exerciseUnitary2},
if and only if  $||f^{(n)}||$
is bounded (hint: Assume $||f^{(n)}||\le k$ and consider
the closure of the set of convex linear combinations of the $U^nf$.
Apply the Schauder-Tychonov fixed point theorem to
the map $h\mapsto f+Uh$).
\end{exercise}

\begin{exercise}\label{exoQueffelec}
Let  $\varphi:A^*\rightarrow A^*$ be a primitive morphism
and let $(X,S)$ be the associated shift space.
We indicate here a method to compute the frequency of the factors of length 
$k$ in $\cL(X)$
by a faster method
than by using Formula~\eqref{eqInvariantMeasure}.
 
Let $M_{k}$ be the incidence
matrix of $\varphi_k$. Let $p$ be an integer such that 
$|\varphi^p(a)|>k-2$ for all $a\in A$. Let $U$ be the 
$\cL_2(X)\times \cL_k(X)$--matrix defined as
follows. For $a,b\in A$ such that $ab\in \cL_2(X)$ and $y\in \cL_k(X)$, 
$U_{ab,y}$ is the number
 of occurrences of $y$ in $\varphi^p(ab)$ 
 that begin inside the prefix $\varphi^p(a)$. Show that
$$UM_k=M_2U,$$
and that  if $v_2$
is a row eigenvector of $M_2$ corresponding to the common dominant
eigenvalue $\rho$ of $M_2$ and $M_k$, 
then $v_k=v_2U$ is an eigenvector of
$M_k$ corresponding to $\rho$.
\end{exercise}
\begin{exercise}\label{exerciseQueffelec2}
Let $\mu:a\rightarrow ab,b\rightarrow ba$ be the morphism with fixpoint the Thue-Morse word.
\index{names}{Thue, Axel}
Show that for $k=5$, $p=3$, the matrix $U$ of the previous problem (with the 12 factors
of length 5 of the Thue-Morse word listed in alphabetic order) is
$$U=\left[
\begin{array}{cccccccccccc}
1&0&1&1&0&1&1&0&1&1&0&1\\
0&1&1&0&1&1&1&1&0&1&1&0\\
1&1&0&1&1&0&0&1&1&0&1&1\\
1&0&1&1&0&1&1&0&1&1&0&1
\end{array}
\right]$$
and that the vector $v_2U$ with $v_2=\left[\begin{array}{cccc}1&2&2&1\end{array}\right]$
is the vector with all components equal to 4. Deduce that  the 12 factors
of length 5 of the Thue-Morse word have the same frequency (see Example~\ref{exampleMeasureMorse}).
\end{exercise}
\begin{exercise}\label{exerciseShiftEquivM_k}
Let $k,p$ be as in the previous exercise. Let $V$ be the matrix
of the map $\pi:\cL_k(X)\to \cL_2(X)$ which sends $a_0a_1\cdots a_{k-1}$
to $a_0a_1$. Show that $M_2$ is shift equivalent
over $\Z$ to $M_k$ (see Exercise~\ref{exerciseShiftEquivalence1})
and more precisely that $(U,V):M_2\sim M_k$ (lag $k$).
Conclude that $M_2$ and $M_k$ have the same nonzero eigenvalues.
Hint: use Exercise~\ref{exerciseShiftEquivalenceSpectrum}.
\end{exercise}
\exosection{Section~\protect{\ref{sectionInvariantStates}}}

\begin{exercise}\label{exerciseKimRoushWilliams}
  Let $X$ be the shift of Example~\ref{exampleKimRoushWilliams}.
  Show that there is no decomposition $\charac_{[b]}=\phi_1+\phi_2$
  such that $\phi_1\le\charac_{[d]}$ and $\phi_2\le\charac_{[e]}$.
  Hint:
  To prove that $\charac_{[b]}\le\phi_1\le\charac_{[d]}$,
  find an invariant Markov measure such that $\mu([a^{n-1}b])>0$
  for all $n\ge 1$
  while $\mu([d])=0$ and use Theorem~\ref{propositionKerov}.
  \end{exercise}

  \begin{exercise}\label{exerciseMeasureInduced}
    Let $(X,T)$ be a minimal dynamical system and let $U\subset X$ be clopen.
    Let $n(x)=\inf\{n>0\mid T^nx\in U\}$ be the entrance time in $U$ and $T_U(x)=T^{n(x)}(x)$ be the
    transformation induced by $T$ on $U$ (see Section~\ref{sectionInduced}).
    Let also $X_n=\{x\in X\mid n(x)=n\}$.
    Show that for every invariant Borel probability measure $\nu$ on $(U,T_U)$,
    the measure defined on $X$ by
    \begin{equation}
\hat{\nu}(V)=\frac{1}{\lambda}\sum_{n\ge 1}\nu(T^n(V\cap X_n))\label{equationInducedMeasure}
    \end{equation}
    with $\lambda=\sum_{n\ge 1}n\nu(U\cap X_n)$ is an invariant Borel probability
    measure on $(X,T)$. Show that the map $\nu\mapsto\hat{\nu}$ is a bijection from
    the set of invariant Borel probability measures of $(U,T_U)$
    onto the set of invariant Borel probability measures on $(X,T)$.
  \end{exercise}

  \begin{exercise}\label{exerciseMeasureInducedShift}
    Let $\varphi:A^*\to B^*$ be an injective nonerasing morphism. Let $X$ be a shift space
    on $A$ and let $Y$ be the closure under the shift of $\varphi(X)$.
    For a nonempty word $w\in B^*$, a \emph{context} of $w$ is a pair $(u,v)$
    of words in $B^*$ such that $uwv=\varphi(a_1\cdots a_n)$ with
    $n\ge 0$, $a_i\in A$ and $u$ (resp. $v$) a proper prefix of $\varphi(a_1)$
    (resp. a proper suffix of $\varphi(a_n)$). Denote by $C(w)$ the set
    of contexts of $w$.
    
    Let $\nu$ be an invariant Borel probability measure on $X$. Set
    $\lambda=\sum_{a\in A}|\varphi(a)|\nu([a])$. For $w\in B^*$, let
    \begin{equation}
      \mu([w])=\frac{1}{\lambda}\sum_{(u,v)\in C(w)}\nu(\varphi^{-1}[uwv]).\label{eqContextualProba}
    \end{equation}
    Show that $\mu$ defines an invariant Borel probability measure on $Y$.
    \end{exercise}

%%%%%%%%%%%%%%%%%%%%%%%%%%%%
\section{Solutions}
\exosection{Section~\protect{\ref{sectionCoboundaries}}}
\begin{solution}{\protect{\ref{exerciseNegativeCocycles}}}
We have, extending the composition with $T$ to rational
fractions in $T$,
\begin{displaymath}
f^{(n)}=f\circ \frac{1-T^n}{1-T}.
\end{displaymath}
This holds also for $n<0$ since 
\begin{displaymath}
f\circ\frac{1-T^{-n}}{1-T}=-f\circ T^{-n}\circ\frac{1-T^n}{1-T}=-f^{(n)}(T^{-n}).
\end{displaymath}

The cohomlogical equation is then easy to verify since
\begin{eqnarray*}
f^{(n)}+f^{(m)}\circ T^n&=&f\circ\frac{1-T^n}{1-T}+f\circ\frac{1-T^m}{1-T}\circ T^n\\
&=&f\circ\frac{1-T^n+T^n-T^{n+m}}{1-T}
=f\circ\frac{1-T^{n+m}}{1-T}\\
&=&f^{(n+m)}.
\end{eqnarray*}

\end{solution}
\begin{solution}{\protect{\ref{exerciseHomotopy}}}
We have to verify the transitivity. Assume that $f_0,f_1$
are connected by $f_t$ and that $g_0=f_1,g_1$ are connected by $g_t$.
Then $f_0,g_1$ are connected by $h_t$ where
\begin{displaymath}
h_t(x)=\begin{cases}f_{2t}(x)&\mbox{ if $0\le t<1/2$}\\
g_{2t-1}(x)&\mbox{ if $1/2<t\le 1$}\end{cases}
\end{displaymath}

Let $f_0,f_1:X\to Y$ be connected by a homotopy $f_t$. For $g:Y\to Z$,
the maps $h_0=g\circ f_0$ and $h_1=g\circ f_1$ are connected by $h_t=g\circ f_t$.
Similarly, if $g_0,g_1:X\to Z$ are connected by $g_t$, then for $f:X\to Y$,
the maps $g_0\circ f$ and $g_1\circ f$ are connected by $g_t\circ f$.
This proves the compatibility with composition.
\end{solution}
\begin{solution}{\protect{\ref{exerciseHomotopy2}}}
Set $f_t(x)=ta+(1-t)f(x)$.
\end{solution}
\begin{solution}{\protect{\ref{exerciseSuspension}}}
First $\tilde{X}$ is compact because it can be identified
with the quotient of $X\times[0,1]$ by the equivalence
which identifies $(x,1)$ and $(Tx,0)$ for all $x\in X$.
As a continuous image of a compact metric space, it is
metrizable (see~\cite{Willard2004}). The map
$(y,t)\mapsto T_t(y)$ is well defined and continuous
since $(y,s)\mapsto (y,s+t)$ is continuous from
$X\times \R$ to itself.
\end{solution}
\begin{solution}{\protect{\ref{exerciseSuspensionPeriodic}}}
Set $X=\{0,1,\ldots,n-1\}$ with $T(i)=i+1\bmod n$.
Then $\tilde{X}$ can be identified with the torus $\T$
 by the map $(i,t)\in X\times[0,1[\mapsto (i+t)/n$.
\end{solution}
\begin{solution}{\protect{\ref{exerciseFlowEquivalence}}}
Assume that $\varphi:(X,S)\to(Y,T)$ is an equivalence. Then
the map $(x,t)\to(\varphi(x),t)$ induces
an equivalence from the suspension flow over $(X,S)$
onto the suspension flow over $(Y,T)$.
All periodic systems are flow equivalent by Exercise~\ref{exerciseSuspensionPeriodic}. Thus flow equivalence is weaker than equivalence.
\end{solution}
\begin{solution}{\protect{\ref{exerciseCechCohomology}}}
Consider the map $\pi:C(X,\Z)\to C(\tilde{X},\T)$ which associates to 
$f\in C(X,\Z)$ the map $\pi(f)\in C(\tilde{X},\T)$ defined by
$\pi(f)(x,s))=\tau(f(x)s)$. 

We first show that $\pi$
induces a surjective map from $C(X,\Z)$ to
$H^1(\tilde{X},\Z)$.
Let $\tilde{f}:\tilde{X}\to\T$ be continuous.
By Lemma~\ref{lemmaTorus} there is  a continuous function
$f:X\to\R$ such that $\tilde{f}(x,0)=\tau(f(x))$. Set for $0\le s<1$,
\begin{displaymath}
h(x,s)=\tau(f(x)(1-s)+f(Tx)s).
\end{displaymath}
Then $h$ is nullhomotopic (by Exercise~\ref{exerciseHomotopy2}
because $f(x)(1-s)+f(Tx)s$ extends to a continuous
map from $X\times[0,1]$ to $\R$)
and $h(x,0)=\tilde{f}(x,0)$.
Therefore the map $g_1(x,s)=\tilde{f}(x,s)/h(x,s)$ is 
homotopic to $\tilde{f}$ and
$g_1(x,0)=g_1(Tx,1)=0$ for all $x\in X$. For $x\in X$, let $r(x)$
be the number of times the loop $g_1(x,s)$ wraps around $\T$ 
as $s$ increases from $0$ to $1$. Then $r$ is continuous and,
for each $x\in X$, $g(x,s)=\tau(r(x)s)$ wraps around $\T$ the
the same number of times as $g_1(x,s)$. Hence $g_1$ and $g$ are
holomorphic. This shows that $\pi$ induces a surjective map.

Next, let $f,g\in C(X,\Z)$ be such that $\pi(f)$ and $\pi(g)$
are in the same homotopy class. We may write for $0\le s<1$,
$\tau((f(x)-g(x))s)=\tau(r(x,s))$ for some continuous function $r$
from $\tilde{X}$ to $\R$ (because $\tau(r(x,s))$ is nullhomotopic).
 Then $(f(x)-g(x))s=r(x,s)+P(x,s)$
where $P(x,s)\in\Z$. But since $P(x,s)$ is continuous in $s$,
this forces $P(x,s)=P(x,0)$ for all $s$. Set $p(x)=P(x,0)$. Then
\begin{displaymath}
p(Tx)=-r(Tx,0)=-r(x,1)=p(x)-(f(x)-g(x))
\end{displaymath}
which shows that $f(x)-g(x)$ is a coboundary. Thus $\pi$
induces an injective map from $H(X,T,\Z)$ to $H^1(\tilde{X},\Z)$.

\end{solution}

\exosection{Section~\protect{\ref{sectionGH}}}
\begin{solution}{\protect{\ref{exerciseGH}}}
Let $g,g'$ be two solutions.
Let $x$ be a recurrent point in $X$. Then for any
$y=T^nx$, we have by~\eqref{eq1}
\begin{displaymath}
g(y)=f^{(n)}(x)-f(x)+g(x)
\end{displaymath}
and thus
\begin{displaymath}
g(y)-g'(y)=g(x)-g'(x).
\end{displaymath}
Since the positive orbit of $x$ is dense, this shows
that $g-g'$ differ by a constant.
\end{solution}

\begin{solution}{\protect{\ref{exerciseProofKH1}}}
Since $f^{(n+1)}(x)=f(x)+f^{(n)}(Tx)$ for all $n\in\Z$,
by Exercise \ref{exerciseNegativeCocycles}, we have 
\begin{eqnarray*}
g(x)&=&\sup_{n\in\Z}f^{(n)}(x)=\sup_{n\in\Z}f^{((n+1)}(x)=\sup_{n\in\Z}(f(x)+f^{(n)}(Tx))\\
&=& f(x)+g(Tx)
\end{eqnarray*}
whence the result.
\end{solution}
\begin{solution}{\protect{\ref{exerciseProofKH2}}}
We first note that  $|f^{(n)}(x)|$ is bounded for all $x\in X$.
Indeed, set $M=\sup_{n\ge 0}|f^{(n)}(x_0)|$. If $|f^{(n)}(y)|>2M$, the same inequality holds
for any $z$ sufficiently close to $y$, in particular
for some iterate $T^m(x_0)$. But then $2M<|f^{(n)}(T^mx_0)|\le
|f^{(n+m)}(x_0)|+|f^{(m)}(x_0)|$ contrary to the definition of $M$.

Thus $f^{(n)}(x)$ is bounded for all $n\in \Z$ and,
by Exercise~\ref{exerciseProofKH1}, the map $g(x)=\sup_{n\in_Z}f^{(n)}(x)$
is such that $\partial g=f$.

Define the \emph{oscillation} 
\index{subject}{oscillation}%
 of a real valued function $h:X\to \R$ defined on a metric space
$X$ at a point $x$ as
\begin{displaymath}
\Osc_h(x)=\lim_{\delta\to 0}(\sup\{h(y)\mid d(x,y)<\delta\}
-\inf\{h(y)\mid d(x,y)<\delta\})
\end{displaymath}
Note that the oscillation of a function $h$ at $x$ vanishes
if and only if $h$ is continuous at $x$. Since $f$
is continuous, we have $\Osc_h(x)=0$ for $h=\partial g$. This
implies $\Osc_{g\circ T}(x)=\Osc_g(x)$ or $\Osc_g\circ T=\Osc_g$
and thus that the function
$x\mapsto \Osc_g(x)$ is invariant. For $\varepsilon>0$, let $O_{\varepsilon,n}$
be the set of $x\in X$ such that $f(x)-f^{(n)}(x)\le\varepsilon/2$.
Since $O_{\varepsilon,n}$ is closed, the set 
\begin{displaymath}
\{x\in O_{\varepsilon,n}\mid \Osc_g(x)\le\varepsilon\}
\end{displaymath}
 is closed invariant and nonempty it is equal to $X$.
Thus $g$ is continuous.
\end{solution}
\begin{solution}{\protect{\ref{exerciseProofHerman}}}
\begin{enumerate}
\item follows from the hypothesis that the sequence $f^{(n)}(x_0)$ is bounded.
\item Assume that $(x,u),(x,v)\in K$ for some $x\in X$ and $u\ne v$.
The map $S$ commutes with  the vertical translations $T_u:(x,t)\mapsto (x,t+u)$.
Thus $T_{u-v}K=\{(x,t+u-v)\mid (x,t)\in K\}$ is also an
$S$-invariant and minimal compact set.
 It intersects $K$ since $T_{u-v}(x,v)=(x,u)$.
By minimality of $K$,
this implies $T_{u-v}K=K$ which contradicts the fact that $K$ is bounded.
This shows that $K$ is the graph of a function $g$. Its
domain is $X$ since it is closed and $T$-invariant. The
function $g$ is continuous since $K$ is compact.
\item Since $K$ is $S$-invariant, we have $S(x,g(x))=(Tx,g(x)+f(x))\in K$
for every $x\in X$. Thus $g(Tx)=g(x)+f(x)$.
\end{enumerate}
\end{solution}

\exosection{Section~\protect{\ref{chapterOrderedGroupRecurrent}}}
\begin{solution}{\ref{exercise[ab]}}
  Suppose that $\charac_{[ab]}=\partial \phi$
  with $\phi\in C(X,\Z)$. Let $n$ be such that $\phi(x)$
  depends only on $x_{[-n,n)}$. Let $\alpha$ be the value
    of $\phi$ on $[a^n\cdot a^n]$. We proceed in four steps.
    \begin{enumerate}
      \item For $2\le i\le n$,
    we have $\phi\circ T(x)=\phi(x)$ and thus $\phi$
    is constant with value $\alpha$ on $[a^n\cdot a]$.
    \item Next, since $1=\phi\circ T(x)-\phi(x)$ for every
      $x\in[a^n\cdot ab]$, $\phi$ is constant with value $\alpha+1$
      on $[a^n\cdot b]$.
    \item Since $\phi\circ T(x)=\phi(x)$ for every $x\in [a^n\cdot ac]$,
      we have that $\phi$ is constant with value $\alpha$ on $[a^n\cdot c]$.
    \item Finally, we have $\phi\circ T(x)=\phi(x)$ for every $x\in [a^n\cdot bc]$, which implies that $\phi$ is constant with value $\alpha+1$
      on $[a^n\cdot c]$, a contradiction.
      \end{enumerate}
  \end{solution}
\begin{solution}{\ref{exerciseG_n}}
  To prove that $G_n^+(X)\cap (-G_n^+(X))=\{0\}$, it is
  enough to prove that
if $\phi\in Z_{n-1}(X)$ is such that $\partial_{n-1}\phi\in Z_n^+(X)$, then
$\partial_{n-1}\phi=0$. Let $u,v\in L_{n-1}(X)$. Since $X$
is recurrent, there
exists $m\ge n$ and $w\in L_m(X)$ such that $u$ is a prefix of $w$
and $v$ is a suffix of $w$. Then, since $u=w_{[1,n-1]}$ and
$v=w_{[m-n+2,m]}$, we have, because the first sum is telescopic,
\begin{eqnarray*}
\phi(v)-\phi(u)&=&\sum_{i=1}^{m-n+1}(\phi(w_{[i+1,i+n-1]})-\phi(w_{[i,i+n-2]}))\\
&=&\sum_{i=1}^{m-n+1}(\partial_{n-1}\phi)(w_{[i,i+n-1]})\ge 0.
\end{eqnarray*}
 Since this is true for every $u,v$,
it implies that $\phi$ is constant and thus that $\partial_{n-1}\phi=0$.
  \end{solution}

\exosection{Section \protect{\ref{sectionInvariant}}}
\begin{solution}{\protect{\ref{exerciseMarkovKakutani}}}
Let $(X,T)$ be a topological dynamical system.
The  set $\M_X$
is convex and compact by 
 Theorem~\ref{theoremBanachAlaoglu}.
The map $T^*:\M_X\to\M_X$ defined by $T^*\mu=\mu\circ T^{-1}$
is a continuous linear map. Thus, by the Markov-Kakutani Theorem,
it has a fixed point $\mu$, which is obviously an invariant
measure.
\end{solution}
\begin{solution}{\protect{\ref{exerciseInvariantFunction}}}
%Assertions (1) and (3) are equivalent. Indeed, if $\mu$ is ergodic
%and $f$ is invariant and not constant almost everywhere, then
%there is a set $U$ in the range of $f$ such that $0<\mu(f^{-1}(U)<1$.
%Since this set is invariant, we have a contradiction. Thus
%(1) implies (3). The converse is clear since the characteristic function
%of an invariant set is invariant. The same arguments show that
%(3) and (4) are equivalent.
%Finally (4) implies clearly (3). Conversely, suppose that
%$f\in L^1$ and $f\circ T=f$ a.e. Let $g(x)=\limsup f(T^nx)$.
%Then $g$ is $T$-invariant. Let us show that $f=g$ a.e.
%The set $U=\{x\in X\mid f(Tx)=f\}$ has measure $1$
%and thus the set $V=\cap_{n\ge 0}T^{-n}(U)$ has also measure $1$.
%But for every $x\in V$, we have $f(T^nx)=f(x)$ for all $n\ge 0$
%and thus $f(x)=g(x)$. This proves the claim.
Assume that $\mu$ is ergodic. Let $U$ be a Borel set such
that $T^{-1}(U)=U\bmod\mu$. Set 
\begin{displaymath}
V=\cap_{n\ge 0}\cup_{i\ge n}T^{-i}(U).
\end{displaymath}
It is clear that $T^{-1}V=V$. We claim that  $U=V\bmod\mu$.

To prove the claim, note first that for two sets $U,V$, one has
$U=V\bmod\mu$ if and only if $\mu(U\Delta V)=0$ where
$U\Delta V=(U\cup V)\setminus (U\cap V)$ is the
\emph{symmetric difference} \index{subject}{symmetric!difference} \index{symbols}{Delta@$\Delta$}
of $U,V$.
We will use the fact that
\begin{equation}
  U\Delta W\subset (U\Delta V)\cup(V\Delta W).\label{eqDelta}
  \end{equation}

 We first observe that
\begin{equation}
  \mu(U\Delta(\cup_{i\le n-1} T^{-i}U))=0.\label{eqSubClaim}
\end{equation}
Indeed, using \eqref{eqDelta} to telescope the union on the right hand side, we have
\begin{eqnarray*}
  T^{-n}U\Delta U&\subset& \cup_{i=0}^{n-1}(T^{-i-1}U\Delta T^{-i}U)\\
  &=&\cup_{i=0}^{n-1}T^{-i}(T^{-1}U\Delta U)
\end{eqnarray*}
and thus $\mu(T^{-n}U\Delta U)\le n\mu(T^{-1}U\Delta U)=0$. This implies
that
\begin{displaymath}
  \mu(U\Delta(\cup_{i\ge n}T^{-i}U))\le \sum_{i\ge n}\mu(U\Delta T^{-i}U)=0
\end{displaymath}
which proves \eqref{eqSubClaim}.
Now \eqref{eqSubClaim} implies that  $\mu(U\Delta V)=0$ and
thus that $\mu(U)=\mu(V)$, proving the claim.

Since $\mu$ is ergodic, we have $\mu(V)=0$ or $1$. This shows
that (i) implies (ii). The converse is obvious.
\end{solution}
\begin{solution}{\ref{exerciseBernoulliErgodic}}
  Denote $S(U,V)=\frac{1}{n}\sum_{i=0}^{n-1}\mu(U\cap T^{-i}V)$.
  Suppose first that \eqref{eqCriteriumErgodicity} holds for all clopen sets $U,V$. Then
  it holds for all Borel sets $U,V$. Indeed, for $\varepsilon>0$, let $U_0,V_0$ be clopen sets
  such that $\mu(U\Delta U_0)<\varepsilon$ and $\mu(V\Delta V_0)<\varepsilon$
  (they exist by Theorem~\ref{theoremCaratheodoryExtension}). Then
  $|S(U_0,V_0)-\mu(U_0)\mu(V_0)|<\varepsilon$
  for $n$ large enough. Since
  $(U\cap T^{-i}V)\Delta(U_0\cap T^{-i}V_0)\subset (U\Delta U_0)\cup (T^{-i}V\Delta V_0)$,
  we have $\mu((U\cap T^{-i}V)\Delta(U_0\cap T^{-i}V_0))<2\varepsilon$. Thus we obtain
  \begin{eqnarray*}
    |S(U,V)-\mu(U)\mu(V)|&\le&|S(U,V)-S(U_0,V_0)|+|S(U_0,V_0)-\mu(U_0)\mu(V_0)|\\
    &&+|\mu(U_0)\mu(V_0)-\mu(U)\mu(V_0)|\\&&+|\mu(U)\mu(V_0)-\mu(U)\mu(V)|\\
    &\le&5\varepsilon.
    \end{eqnarray*}

  Then
  for every invariant Borel set $U$, we have $\mu(U)=\frac{1}{n}\sum_{i=0}^{n-1}\mu(U)=
  \frac{1}{n}\sum_{i=0}^{n-1}\mu(U\cap T^{-i}U)=\mu(U)^2$
  and thus $\mu(U)$ is $0$ or $1$.

  Conversely, 
  by the ergodic Theorem
  (Theorem~\ref{theoremErgodic}) applied to the characteristic function $\charac_V$
  of the Borel set $V$, we have
  \begin{displaymath}
    \lim_{n\to\infty}\frac{1}{n}\charac^{(n)}_V(x)=\int\charac_Vd\mu=\mu(V)
  \end{displaymath}
  almost everywhere. Multiplying both sides by $\charac_U(x)$ and integrating,
  we obtain \eqref{eqCriteriumErgodicity}.

  If $\mu$ is a Bernoulli measure, then \eqref{eqCriteriumErgodicity}
  holds for all clopen sets $U,V$ and thus $\mu$ is ergodic.
\end{solution}
\begin{solution}{\ref{exerciseMarkovErgodic}}
  Assume first that $P$ is not irreducible. Let $B\subset A$ be a subset of
  the alphabet corresponding to a strongly connected component of the
  graph associated to $P$ (for $a,b\in A$, there is an edge from $a$ to $b$
  if and only if $P_{a,b}>0$). Then $U=B^\Z$ is an invariant subset of
  $A^\Z$ such that $0<\mu(U)<1$.

  Conversely, by Exercise~\ref{exerciseBernoulliErgodic}, it is enough to
  prove that $\mu(U\cap T^{-n}(V))$ converges in mean to $\mu(U)\mu(V)$
  for $U,V$ clopen. Since the family of sets
  $U,V$ for which this holds is closed under union and translation,
  it is enough to consider $U=[u]$ and $V=[w]$.
  Set $u=a_1\cdots a_s$ and $w=b_1\cdots b_t$. Then
  \begin{displaymath}
    \mu([u]\cap T^{-n}[w])=v_{a_1}P_{a_1,a_2}\cdots P_{a_{s-1},a_s}(P^{n-s})_{a_s,b_1}P_{b_1,b_2}\ldots P_{b_{t-1},b_t}.
  \end{displaymath}
  By Theorem~\ref{theorem:perron}, the sequence $P^n$ converges in mean to the matrix with all rows
  equal to $v$. Thus $\mu([u]\cap T^{-n}[w])$ converges in mean to
  \begin{displaymath}
    v_{a_1}P_{a_1,a_2}\cdots P_{a_{s-1},a_s}v_{b_1}P_{b_1,b_2}\ldots P_{b_{t-1},b_t}=\mu([u])\mu([w]).
    \end{displaymath}
  
  \end{solution}

\begin{solution}{\protect{\ref{exerciseErgodicRecurrent}}}
Let $R$ be the set of recurrent points in $(X,T)$ 
(see Exercise \ref{exerciseRecurrentPoint}).
We show that $\mu(R)=1$. This proves that
$(X,T)$ is recurrent by Exercise~\ref{exerciseCondRecurrent}.

Let $(U_n)_{n\ge 1}$ be a countable basis of open sets of $X$ (this
exists because $X$ is a metric space). Set
\begin{displaymath}
V_n=\cap_{k\ge0}T^{-k}(X\setminus U_n).
\end{displaymath}
We have $X\setminus R=\cup_{n\ge 1}V_n$. Since $T^{-1}V_n\subset V_n$
and $\mu(T^{-1}V_n)=\mu(V_n)$, it follows that $T^{-1}V_n=V_n\bmod\mu$
and thus $\mu(V_n)=0$ or $1$
by Proposition~\ref{propositionCondEquiv} since $\mu$ is ergodic.
 But $U_n\subset X\setminus V_n$
implies $\mu(U_n)\le\mu(X\setminus V_n)$. By our hypothesis,
this implies $\mu(X\setminus V_n)>0$ and thus $\mu(V_n)=0$.
This gives finally $\mu(R)=1$.
\end{solution}
\begin{solution}{\ref{exerciseErgodicityRotations}}
By Weierstrass Theorem,
\index{subject}{Weierstrass Approximation Theorem}%
\index{subject}{Theorem!Weierstrass Approximation}%
\index{names}{Weierstrass, Karl Theodor}%
the linear combinations of the functions $\chi_m(z)=z^m$
are dense in $C(S^1,\C)$. Thus, by Theorem~\ref{theoremOxtoby},
it is enough to prove that for every $m\ge 0$, the averages
$\chi_m^{(n)}$ converge uniformly to a constant.This is
trivially true if $m=0$.  Otherwise, the
rotation $R_\alpha(z)=\lambda z$ with $\lambda=e^{2i\pi\alpha}$ satisfies
\begin{eqnarray*}
\left|\frac{1}{n}\sum_{k=0}^{n-1}\chi_m(R^k_\alpha(z))\right|&=&
\left|\frac{1}{n}\sum_{k=0}^{n-1}e^{2i\pi km\alpha }\right|\\
&=&\frac{|1-e^{2i\pi mn\alpha }|}{n|1-e^{2i\pi m\alpha }|}\le
\frac{2}{n|1-e^{2i\pi m\alpha}|}\to 0.
\end{eqnarray*}
Every Sturmian shift $(X,S)$ of slope $\alpha$
is, by Proposition~\ref{propositionRotationIsSturm},
the image by the natural coding $\gamma_\alpha$ of the rotation
$\T,R_\alpha)$ 
of angle $\alpha$. Since $\gamma_\alpha$ is one-to-one except on
a denumerable set (the orbit of $0$ under $R_\alpha$), the
invariant Borel probability measures on $(\T,R_\alpha)$ and on $(X,S)$
are exchanged by $\gamma_\alpha$. Thus $(X,S)$ is uniquely ergodic
and its unique invariant measure $\mu$ is such that $\mu([0])=1-\alpha$,
$\mu((1])=\alpha$.
\end{solution}

\begin{solution}{\ref{exerciseNonUniquelyErgodic}}
It follows from the definition that for $1\le j\le i$,
the number of elements of $E_j$ in the interval
$0<n\le k_i$ is exactly $(k_i/k_j)(2k_{j-1}+1)$.
Hence an upper bound to the number of elements
of $E_1\cup E_2\cup\cdots\cup E_i$ in the interval
$0<n\le k_i$ is 
\begin{displaymath}
\sum_{j=1}^i\frac{3k_ik_{j-1}}{k_j}\le3k_i\sum_{j\ge 1}\frac{k_{j-1}}{k_j}<\frac{1}{4}k_i.
\end{displaymath}
It follows that $p(n)=i+1$ for at least $3/4$ of the numbers $n$ in
the interval $0<n\le k_i$. This implies that
\begin{displaymath}
\left|\frac{1}{k_i}\sum_{n=1}^{k_i}x(n)-\frac{1}{k_{i+1}}\sum_{n=1}^{k_{i+1}}x(n)\right|\ge\frac{1}{2}.
\end{displaymath}
and thus the frequency of $1$ is not defined.
\end{solution}

\begin{solution}{\protect{\ref{exerciseUnitary}}}
Since $T$ is continuous, it is measurable and
since $\mu$ is invariant, we have $\mu(T^{-1}V)=\mu(V)$
for every borel subset $V$ of $X$. Thus $T$ preserves $\mu$.
\end{solution}
\begin{solution}{\protect{\ref{exerciseUnitary2}}}
Set $\langle f,g\rangle=\int fg d\mu$.
We have, by change of variable 
(Equation~\eqref{eqChangeVariables}), for every $f,g\in{L}^2(X)$
\begin{displaymath}
\langle U^*Uf,g\rangle=\langle Uf,Ug\rangle=
\int (Uf)(Ug)d\mu=\int U(fg)d\mu=\int fg d\mu=\langle f,g\rangle.
\end{displaymath}
Thus $U^*U=I$. Since $T$ is invertible, $U$ is surjective. Let
$V$ be its right inverse, that is, such that $UV=I$. Then
$V=U^*UV=U^*$. Thus $UU^*=I$. This shows that $U$ is unitary.
\end{solution}
%\begin{solution}{\protect{\ref{exerciseGH2}}}
%We have
%\begin{eqnarray*}
%\partial g(x)&=&\sup_{n\ge 1}(-f^{(n)}(Tx)) -\sup_{n\ge 1}(-f^{(n)}(x))\\
%&=&\sup_{n\ge 1}(-f^{(n)}(Tx))-\sup_{n\ge 1} (-f(x)-f^{(n)}(Tx))\\
%&=&f(x).
%\end{eqnarray*}
%\end{solution}
\begin{solution}{\protect{\ref{exerciseGH3}}}
If $f=Ug-g$ with $g\in H$, then $||f^{(n)}||=||U^ng-g||\le2||g||$
and thus the sequence is bounded.

Conversely, assume that $||f^{(n)}||\le k$ for $n\ge 1$.
Let $S$ be the set of all convex linear combinations of $U^nf$
and let $\bar{S}$ be its closure in the weak topology of $H$.
Then $\bar{S}$ is convex and as it is
contained in the weakly compact set $\{h\in H\mid ||h||\le k\}$
it is weakly compact. Moreover $\bar{S}$
is invariant under the continuous affine map $h\mapsto f+Uh$.
By the Schauder-Tychonoff Theorem~\index{subject}{Schauder-Tychonoff Theorem}%
\index{subject}{Theorem!Schauder-Tychonov}%
this map has fixed point in $\bar{S}$, that is there is $g\in\bar{S}$
such that $g=x+Ug$ and thus $x=Ug-g$ is the coboundary of $g\in H$.
\end{solution}
%\exosection{ Section \protect{\ref{sectionInvariant}}}
\begin{solution}{\protect{\ref{exoQueffelec}}}
By the choice of $p$, the value of $\varphi_k^p$ is determined by the
two first letters $a_0a_1$ of $x=a_0a_1\cdots a_{k-1}$.
Let $\alpha:\cL_2(X)^*\to\cL_k(X)^*$ be the morphism defined by
$\alpha(a_0a_1)=\varphi_k^p(x)$. Then $U$ is the incidence
matrix of $\alpha$. Since obviously
$\varphi_k\circ\alpha=\alpha\circ\varphi_2$, we obtain
$UM_k=M_2U$. Since
\begin{displaymath}
v_kM_k=v_2UM_k=v_2M_2U=\rho v_2U=\rho v_k
\end{displaymath}
the last assertion follows.
\end{solution}
\begin{solution}{\protect{\ref{exerciseQueffelec2}}}
We have $\mu^3(aa)=abbabaab\cdot abbabaab$ and thus
\begin{displaymath}
\mu_5^3(aa)=(abbab)(bbaba)(babaa)(abaab)\cdots
\end{displaymath}
 whence
the value of the first row of $U$ and similarly for the others.
\end{solution}
\begin{solution}{\ref{exerciseShiftEquivM_k}}
The equalities
\begin{displaymath}
M_2U=UM_k,\quad VM_2=M_kV,
\end{displaymath}
and \begin{displaymath}
M_2^p=UV,\quad M_k^p=VU
\end{displaymath}
result of the  commutative diagram of
Figure~\ref{figureDiagramExtension2}.
The diagram on the left gives the
equality $VM_2=M_kV$ (see the diagram \eqref{DiagramExtension}).
The other three equalities result from the diagram on the right.
\begin{figure}[hbt]
\centering
\tikzset{node/.style={minimum size=0.1cm,inner sep=0pt}}
	\tikzset{title/.style={minimum size=0.5cm,inner sep=0pt}}
\begin{tikzpicture}
%diagram1
\node(h0)at(-5,2){$\cL_k(X)^*$};\node(h1)at(-3,2){$\cL_k(X)^*$};
\node(b0)at(-5,0){$\cL_2(X)^*$};\node(b1)at(-3,0){$\cL_2(X)^*$};
\draw[above,->](h0)edge node{$\varphi_k$}(h1);
\draw[below,->](b0)edge node{$\varphi_2$}(b1);
\draw[left,->](h0)edge node{$\pi$}(b0);
\draw[right,->](h1)edge node{$\pi$}(b1);
%diagram2
\node[title](h0)at(0,2){$\cL_k(X)^*$};\node[title](h1)at(2,2){$\cL_k(X)^*$};
\node[title](h2)at(4,2){$\cL_k(X)^*$};
\node[title](b0)at(0,0){$\cL_2(X)^*$};\node[title](b1)at(2,0){$\cL_2(X)^*$};
\node[title](b2)at(4,0){$\cL_2(X)^*$};

\draw[above,->](h0)edge node{$\varphi_k^p$}(h1);
\draw[above,->](h1)edge node{$\varphi_k$}(h2);
\draw[left,->](h0)edge node{$\pi$}(b0);%\draw[right,->](h1)edge node{$\pi$}(b1);
\draw[right,->](h2)edge node{$\pi$}(b2);
\draw[below,->](b0)edge node{$\varphi_2$}(b1);
\draw[below,->](b1)edge node{$\varphi_2^p$}(b2);
\draw[below,->](b0)edge node{$\alpha$}(h1);\draw[below,->](b1)edge node{$\alpha$}(h2);
\end{tikzpicture}\label{figureDiagramExtension2}
\caption{A commutative diagram}
\end{figure}
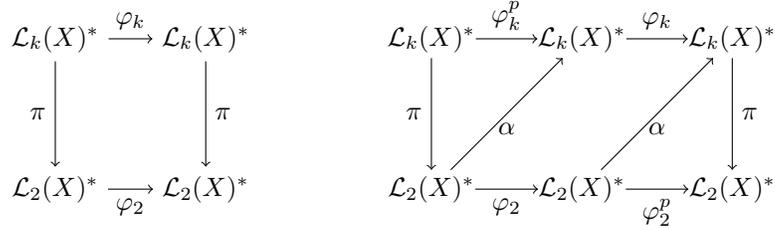
\end{solution}

\exosection{Section~\protect{\ref{sectionInvariantStates}}}

\begin{solution}{\ref{exerciseKimRoushWilliams}}
Let us show that there is no decomposition $\charac_{[b]}\simeq\phi_1+\phi_2$
such that $\phi_1\le \charac_{[d]}$ and $\phi_2\le \charac_{[e]}$. This
will show that $H(X,S,\Z)$ is not a Riesz group and thus not a dimension
group. Assume the contrary. Then there is an integer $n\ge 1$ 
 such that
$\phi_1,\phi_2$ can be defined in $Z_n(X)$. 
Since $\charac_{[a^{n-1}b]}\le\charac_{[b]}$, we must have
$\charac_{[a^{n-1}b]}\le \charac_{[d]}$ or $\charac_{[a^{n-1}b]}\le\charac_{[e]}$. It is easy
to see that this is not possible by assigning real values
to $a,b,c,d,e,f$ on the graph of Figure~\ref{figureExampleKimRoushWilliams}
defining an invariant Borel probability measure on $X$.
For example, using the stationary Markov chain  defined by
\begin{displaymath}
  v=\kbordermatrix{&a&b&c&d&e&f\\&\frac{1}{4}&\frac{1}{4}&0&0&\frac{1}{4}&\frac{1}{4}},\quad
  P=\kbordermatrix{
     &a          &b          &c&d&e          &f\\
    a&\frac{1}{2}&\frac{1}{2}&0&0&0          &0\\
    b&0          &0          &0&0&\frac{1}{2}&\frac{1}{2}\\
    c&0          &0          &1&0&0          &0 \\
    d&0          &0          &0&1&0          &0\\
    e&\frac{1}{2}&\frac{1}{2}&0&0&0          &0\\
    f&0          &0          &0&0&\frac{1}{2}&\frac{1}{2}
    }
\end{displaymath}
we obtain  an invariant
probability measure $\mu$ on $X$ (in fact, the nonzero entries of $P$
are all edges of the Rauzy graph $\Gamma_2(X)$).

By Theorem~\ref{propositionKerov},
we have for every cylinders $[u],[v]\in X$
\begin{displaymath}
\charac_{[u]}\le \charac_{[v]}\Rightarrow \mu([u])\le\mu([v]).
\end{displaymath}
This shows that $\charac_{[a^{n-1}b]}\le \charac_{[d]}$ is impossible since the
measure of the cylinder $[a^{n-1}b]$ is $1/2^{n+1}$ while the
second one has measure $0$. A similar argument shows that
$\charac_{[a^{n-1}b]}\le\charac_{[e]}$ is also impossible and we therefore
obtain a contradiction.

  \end{solution}
  \begin{solution}{\ref{exerciseMeasureInduced}}
    Formula~\ref{equationInducedMeasure} defines clearly a measure
    on $X$. Since $T^n(X_n)=\cup_{m\ge n}U\cap X_m$,
    one has $\hat{\nu}(X)=\frac{1}{\lambda}\sum_{n\ge 1}\nu(T^n(X_n))=1$
    and thus $\hat{\nu}$ is a probability measure.
    Next, for $V\subset U$, one has with $W_n=V\cap X_n$,
    \begin{eqnarray*}
      \hat{\nu}(V)&=&\frac{1}{\lambda}\sum_{n\ge 1}\nu(T^nW_n)
      =\frac{1}{\lambda}\sum_{n\ge 1}\nu(T_UW_n)\\
      &=&\frac{1}{\lambda}\sum_{n\ge 1}\nu(W_n)=\frac{1}{\lambda}\nu(V)
      \end{eqnarray*}
    showing (together wih $\hat{\nu}(U)\lambda=1$)
    that $\hat{\nu}$ determines $\nu$. Finally, let $V\subset X_n$.
    If $n>1$, one has $TV\subset X_{n-1}$ and thus $\hat{\nu}(TV)=\hat{\nu}(V)$
    by definition. If $n=1$, then since $TV\subset U$, we have
    \begin{eqnarray*}
      \hat{\nu}(TV)&=&\frac{1}{\lambda}\nu(TV)=\hat{\nu}(V).
      \end{eqnarray*}
    This shows that $\hat{\nu}$ is invariant.
    \end{solution}
\begin{solution}{\ref{exerciseMeasureInducedShift}}
  Let $\hat{\nu}$ be defined as in Exercise~\ref{exerciseMeasureInduced}.
  Then $\mu=\hat{\nu}\circ \hat{\varphi}^{-1}$
  and thus $\mu$ is an invariant Borel probability measure.
    \end{solution}
%%%%%%%%%%%%%%%%%%%%%%%%%%
\section{Notes}
For a general introduction to cohomology, see for example~\citep{Hatcher2001}.
\index{names}{Hatcher, Allen, E.}%
\subsection{Gottschalk and Hedlund Theorem}
Gottschalk and Hedlund Theorem (Theorem~\ref{theoremGH})
 is from \citep{GottschalkHedlund1955}.
A somewhat simpler, although incorrect, proof appears in
 \cite[Theorem 2.9.4]{KatokHasselblatt1995}.
It is reproduced as Exercise~\ref{exerciseProofKH2}
using a corrected version due to \cite{Petite2019}. 
\index{names}{Petite, Samuel}%
Yet another
proof (perhaps the most elegant one) is given in Exercise~\ref{exerciseProofHerman}. According to~\cite{Petite2019}) it should be credited to Michael Herman
 (through Sylvain Crovisier).
\index{names}{Herman, Michael-R.}%
\index{names}{Crovisier, Sylvain}%
Given $f\in C(X,T)$, one may consider $\partial g=f$ as a functional
equation with $g$ as unknown. More generally, an equation
of the form
\begin{equation}
f(x)=\lambda g(Tx)-g(x)\label{cohomologicalEquation}
\end{equation}
where $T:X\to X$ is a given map, $f$ is a given scalar
function on $X$, $\lambda$ is a given constant and
$g$ is an unknown scalar function
is called in \cite{KatokHasselblatt1995} a \emph{cohomological equation}.
\index{subject}{cohomological equation}%
 In this context, given
a dynamical system $(X,T)$, a map
$\alpha:\Z\times X\to\R$ is called a \emph{one-cocycle}\index{subject}{one-cocycle}
if it satisfies the identity
\begin{displaymath}
\alpha(n+m,x)=\alpha(n,T^mx)+\alpha(m,x).
\end{displaymath}
Thus, for every map $f\in C(X,\R)$,
the map $(n,x)\mapsto f^{(n)}(x)$ is a one-cocycle (see Exercise~\ref{exerciseNegativeCocycles}).

The analogue of the ordered
cohomology group for real valued functions
has been studied by \cite{Ormes2000}.

%The Baire Category Theorem, used in the proof of Proposition~\ref{lemma2},
%can be found for example in \cite[Theorem 25.3]{Willard2004}.
\subsection{Ordered cohomology group}

\cite{BoyleHandelman1996} have introduced 
 the term \emph{ordered cohomology group} of a
topological dynamical system $(X,T)$ 
for the  group $K^0(X,T)$, which had before been only defined
for minimal systems. They showed that this group is
a complete invariant for flow equivalence of irreducible shifts of finite type.
The group is actually defined for a class of systems which is
larger than  recurrent dynamical systems, namely
that of \emph{chain recurrent}\index{subject}{chain recurrent system}
dynamical systems. This direction had been explored previously by
\cite{Poon1989}.\index{names}{Poon, Yiu Tung}

A word is in order here on the notation $K^0(X,T)$. The letter $K$
is used in reference to $K$-theory. In this theory,
an abelian group is associated to a ring, based on the structure of idempotents
in the algebra of matrices over this ring. These algebras,
actually $C^*$-algebras,
will be introduced in Chapter~\ref{chapterBratteli}. In the
present case, the algebra is the $K_0$ group of the cross
product $C^*$-algebra arising from a dynamical system
(see~\cite{BoyleHandelman1996}).

The proof of Lemma~\ref{lemmaTorus} 
follows the lines of the proof given in \cite{Giordano&Handelman&Hosseini:2017}.
\index{names}{Giordano, Thierry}\index{names}{Handelman, David}\index{names}{Hosseini, Maryam}%

\subsection{Ordered group of a recurrent shift space}
Example \ref{exampleKimRoushWilliams}
showing that the ordered cohomology group of a shift
of finite type may fail to be a dimension group
is from~\cite[Example 3.3]{KimRoushWilliams2001}.
\index{names}{Kim, Ki Hang}%
\index{names}{Roush, Fred W.}%
\index{names}{Williams, Susan G.}%

\subsection{Invariant measures and states}
For  an introduction to probability measures 
on topological dynamical systems, 
see~\cite{KatokHasselblatt1995,BertheRigo2010}).

For an introduction to the notion of integral of a measurable
function (Section~\ref{sectionInvariant}),
see~\cite{Halmos1974}.
The original reference for the Krylov-Bogolyubov Theoreom
(Theorem~\ref{theoremKrylovBogolioubov}) is \citep{KrylovBogolioubov1937}.
The notion of measure-theoretic dynamical systems is
the central object of ergodic theory. See \cite{KatokHasselblatt1995}
for a systematic exposition. The original reference
for Poincar\'e Recurrence Theorem is
\cite{Poincare1890}.
\index{names}{Poincar\'e, Henry}%

The original reference for Oxtoby's Theorem (Theorem~\ref{theoremOxtoby})
is~\citep{Oxtoby1952}.
\index{names}{Oxtoby, John C.}%
The example of a minimal non uniquely ergodic
shift (Exercise~\ref{exerciseNonUniquelyErgodic}) is from \cite{Oxtoby1952}.
It is the first ever constructed sequence with this property.

The ergodic theorem, due to Birkhoff,
\index{names}{Birkhoff, George D.}%
 is also called the \emph{pointwise ergodic theorem}
or \emph{strong ergodic theorem},
in contrast with the \emph{mean ergodic theorem},
due to von Neumann
\index{names}{Neu@von Neuman, John}%
 which states
the weaker convergence in mean in an $L^2$-space.
The pointwise ergodic theorem is sometimes also called the Birkhoff-Khinchin
Theorem.
\index{names}{Khinchin, Aleksandr Y.}%

The fact that the shift associated with a primitive substitution
is uniquely ergodic (Theorem~\ref{theoremMichel}) is due to \cite{Michel1974}.
\index{names}{Michel, Pierre}
We follow the proof of \cite{Queffelec2010}. A one-sided
version of Lemma \ref{lemmaQueffelec}
appears in \cite{DumontThomas1989}.
\index{names}{Dumont, Jean-Marie}\index{names}{Thomas, Alain}%
The computation of the invariant measure on the shift associated
with a primitive substitution is developped in \cite{Queffelec2010}.
\index{names}{Queff\'elec, Martine}
Formula~\eqref{eqInvariantMeasure} is from~\cite[Corollary 5.14]{Queffelec2010},

Theorem~\ref{propositionKerov} is~Theorem 5.5 in \citep{HermanPutnamSkau1992},
\index{names}{Herman, Richard H.}\index{names}{Putnam, Ian F.}%
\index{names}{Skau, Christian F.}
where it is credited to~\cite{Kerov}.
\index{names}{Kerov, Serguei}%
% For a proof of
%the Caratheodory Extension Theorem, 
%see~\cite{Halmos1974}\index{names}{Halmos, Paul}
% or~\cite{Billingsley1995}.
%\index{names}{Billingsley, Patrick P.}%
The computation of the cohomology group of Sturmian
shifts is classical. It appears in particular in~\cite{Dartnell&Durand&Maass:2000}.
\subsection{Exercises}
The equivalent definition of $H(X,T,\Z)$ in terms of \v{C}ech
cohomology
\index{names}{Cech@\v{C}ech, Eduard}%
 (Exercise~\ref{exerciseCechCohomology}) is taken from~\cite{ParryTuncel1982}.
\index{names}{Parry, William}\index{names}{Tuncel, Selim}%
 The definition
of the \v{C}ech cohomology is not the classical one
but is equivalent. According to \cite{ParryTuncel1982}, it
defines the \emph{Brushlinski group}.
\index{subject}{Bruschlinski group}\index{names}{Bruschlinski, N.}%

The variant of Gottschalk and Hedlund
Theorem for $L^2$ presented in Exercise \ref{exerciseGH3} is from~\cite[Proposition II.2.11]{ParryTuncel1982}. 
The 
Schauder-Tychonoff fixed-point Theorem can be found
as \cite[Theorem V.10.5]{DunfordSchwartz1988}.
\index{names}{Dunford, Nelson}\index{names}{Schwartz, Jacob T.}%
while
the Markov-Kakutani fixed point theorem appears as Theorem V.10.6 in
\citep{DunfordSchwartz1988}.
The unique ergodicity of irrational rotations (Exercise~\ref{exerciseErgodicityRotations}) is known as the Kronecker-Weyl Theorem.
\index{names}{Kronecker, Leopold}\index{names}{Weyl, Hermann}%

The invariant Borel probability measure $\mu$ defined by Formula \eqref{eqContextualProba}
is called the \emph{contextual probability} measure on $Y$
in~\cite{HanselPerrin1983} (see also \cite{BlanchardPerrin1980}).
\index{names}{Hansel, Georges}%
The following result is proved in \cite{HanselPerrin1983}.
Assume that $X=A^\Z$, $Y=B^\Z$ and that $\pi$ is a Bernoulli measure
on $B^\Z$ such that $\nu([w])=\pi([\varphi(w)])$ for every $w\in A^*$.
  Then the contextual probability measure $\mu$
  is equal to $\pi$.
  \index{names}{Blanchard, Fran\c{c}ois}%

 %%%%%%%%%%%%%%%%%%
%  chapter Dimension groups and partitions in towers
%%%%%%%%%%%%%%%%%
\chapter{Partitions in towers}
\label{chapterDimensionGroupsPartitions}

In this chapter, we first present the important notion
of partition in towers. It is the basis of many of the constructions
that will follow.
We prove the Theorem of Herman, Putnam and Skau asserting
that every minimal Cantor system can be represented 
by a sequence of partitions in towers (Theorem~\ref{theoremKRPartitions}).

We next show how partitions in towers
allow to compute  the ordered cohomology group of minimal Cantor systems. 
We first define in Section~\ref{sectionOrderedGroupPartition}
the ordered group associated to a partition. The definition
uses the notions on induction introduced in Section~\ref{sectionInduced}.
In Section~\ref{sectionOrderedGroupSequences}, we prove
that the ordered cohomology group $K^0(X,T)$ is the direct
limit of the sequence of ordered groups associated with a sequence of partitions
in towers. This allows us to prove, as a main result of this
chapter, the theorem
of Herman, Putnam and Skau asserting
 that the ordered cohomology group of a minimal invertible
Cantor system is a simple dimension group (Theorem~\ref{theoremDimensionGroup}).

We will  use  these results to determine the dimension
groups of some minimal shift spaces and relate them
to several notions such as return words or Rauzy graphs. We will in particular
consider episturmian shifts 
(Proposition~\ref{propositionDimensionEpisturmian}), which will
illustrate the use of return words.  Next, we will use Rauzy
graphs to give a description of dimension groups of substitution
shifts (Proposition~\ref{propositionDimensionGroupSubstitutions}).
These methods will appear again in a new light
in the next chapter when we consider
Bratteli diagrams.
%%%%%%%%%%%%%%%%%%%%%%%
%\section{Some short historical overview}\marginpar{FD: a modifier}
%%%%%%%%%%%%%%%%%%%%%%%

\medskip

In this chapter all the dynamical systems $(X,T)$ we consider are 
invertible, that is,
such that $T$ is a homeomorphism. 
%We  thus  work with   the  two-sided {\em orbit} of $x\in X$, that  is,
%$\{ T^n x \mid n\in \Z \}$.

%%%%%%%%%%%%%%%%%%
\section{Partitions in towers}\label{sectionPartitionTowers}
Let $(X,T)$ be an invertible  Cantor system.
Let $B_1 , \dots , B_m$ be a family of disjoint nonempty open sets and $h_1, \dots , h_m$ be positive integers.
Assume that the family
\begin{displaymath}
\mathfrak{P}=\{T^jB_i\mid 1\le i\le m,\ 0\le j<h_i\}.
\end{displaymath}
\index{symbols}{P@$\mathfrak{P}$}
is a partition of $X$.
It implies that each element of $\mathfrak{P}$ is a clopen set. 
We say that it is the \emph{clopen partition in towers}
\index{subject}{partition!in towers}%
(or \emph{Kakutani-Rokhlin partition}\index{subject}{Kakutani-Rokhlin partition}
 or \emph{KR-partition}\index{subject}{KR-partition}) of $(X,T)$ built on $B_1 , \dots , B_m $ with \emph{heights} $h_1, \dots , h_m$.
\index{names}{Kakutani, Shizuo}%
\index{names}{Rokhlin, Vladimir A.}%

The number of towers is $m$, $\{T^jB_i\mid 0\le j<h_i\}$ is the $i$-th 
\emph{tower}\index{subject}{tower!of partition},
$h_i$ is its {\em height} and $B_i$ its \emph{basis}\index{subject}{basis!of tower}.
The union
$B(\mathfrak{P})=\cup_iB_i$ is the \emph{basis}
\index{subject}{basis!of partition} of the partition $\mathfrak{P}$.

Since $T$ is bijective, the shift sends
the elements at the top of the towers back to the bottom,
that is $T^{h_i}B_i\subset B(\Pg)$ for $1\le i\le m$.
As a consequence, when such a partition exists, the integers $h_i$ are unique. 

\begin{figure}[hbt]
\centering
\tikzset{node/.style={draw,minimum size=1cm,inner sep=0pt}}
\tikzset{title/.style={minimum size=0.2cm,inner sep=0pt}}
\begin{tikzpicture}
\node[node](B1)at(0,0){$B_1$};
\node[node](B2)at(2,0){$B_2$};
\node[title](etc)at(4,0){$\ldots$};
\node[node](Bm)at(6,0){$B_m$};
\node[node](TB1)at(0,1){};
\node[node](TB2)at(2,1){};
\node[node](TBm)at(6,1){};
\node[node](T2B1)at(0,2){};

\node[title](B1m)at(0,0.2){};\node[title](TB1m)at(0,1){};
\draw[->](B1m)edge node{}(TB1m);
\node[title](B2m)at(2,0.2){};\node[title](TB2m)at(2,1){};
\draw[->](B2m) edge node{}(TB2m);
\node[title](Bmm)at(6,0.2){};\node[title](TBmm)at(6,1){};
\draw[->](Bmm)edge node{}(TBmm);
\node[title](T2B1m)at(0,2){};
\draw[->](TB1m)edge node{}(T2B1m);
\end{tikzpicture}
\caption{A partition in towers.}\label{figurePartition}
\end{figure}
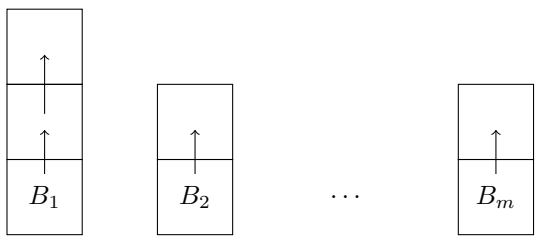

Thus, informally speaking, a partition in towers gives an approximate description of the action of $T$ on $X$. 
Each tower can be seen as a stack of clopen sets.
The transformation consists in climbing one step up the stack except at the top level where it goes back to the basis in some tower (see Figure~\ref{figurePartition}).

A fundamental observation is the following link between
partitions in towers and induction. Let $\Pg$ be a partition
in towers of $(X,T)$. Then the return time to the basis $B=B(\Pg)$
is bounded and thus the  system $(B,T_B)$ induced by $(X,T)$
on the clopen set $B$ is well defined. It is easy
to verify that $X$ is conjugate to the system
obtained from $(B,T_B)$ by the tower construction
using the function $f(x)=h_i$ if $x\in B_i$.

\subsection{Partitions and return words}
We give two examples of partitions in towers of a shift space.
The first example is related to return words. 
The partition relies on the location of two (possibly overlapping)
occurrences of a word $w$, with the second one occurring at a strictly
positive index.

\begin{proposition}\label{propositionPartitionReturn}
Let $(X,S)$ be a minimal subshift and let $w\in \mathcal{L}(X)$. For every $w\in \mathcal{L}(X)$,
the family
\begin{displaymath}
\Pg=\{S^j[vw]\mid v\in\RR'_X(w),0\le j<|v|\}
\end{displaymath}
is a partition in towers with basis the set of cylinders $[vw]$
for $v\in\RR'_X(w)$.
\end{proposition}

\begin{proof}
Let $x\in X$.
Since $(X,S)$ is minimal, there is a smallest
integer $i>0$ such that $x_{[i,i+|w|-1]}=w$. By definition of
a left return word, there is a unique integer $j\ge 0$ such that
$x_{[-j,i-1]}$ belongs to  $\RR'_X(w)$ (see Figure~\ref{figureReturn}).
\begin{figure}[hbt]
\centering
\tikzset{node/.style={circle,draw,minimum size=0.1cm,inner sep=0pt}}
\tikzset{title/.style={minimum size=0.2cm,inner sep=0pt}}
\begin{tikzpicture}
\node[title](x)at(-2,0){$x$};
\node[title](0)at(-0.5,0){$\ldots$};
\node[node](1)at(0,0){};
\node[title](1t)at(0,0.3){$-j$};
\node[node](2)at(2,0){};
\node[node](3)at(3,0){};
\node[title](3t)at(3,0.3){$0$};
\node[node](4)at(5,0){};
\node[title](4t)at(5,0.3){$i$};
\node[node](5)at(7,0){};
\node[title](6)at(7.5,0){$\ldots$};

\draw[above](1)edge node{$w$}(2);
\draw(2) edge node{}(3);
\draw(3)edge node{}(4);
\draw[above](4)edge node{$w$}(5);

\end{tikzpicture}
\caption{Two occurrences of $w$ in $x$}\label{figureReturn}
\end{figure}
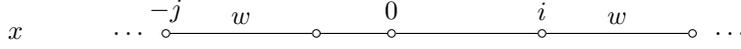

We set $v=x_{[-j,i-1]}\in\RR'_X(w)$ and observe that we have $0\le j<|v|$. Thus $x$ belongs to  $S^j[vw]$.
Since $v$ and $j$ are unique, this shows that $\Pg$ is a partition.
\end{proof}

\begin{example}\label{examplePartition}
Let $(X,S)$ be the two-sided Fibonacci shift on the alphabet $A=\{a,b\}$.
\index{subject}{Fibonacci!shift}\index{subject}{shift space!Fibonacci}%
We have $ab\in \mathcal{L}(X)$ and $\RR'_X(ab)=\{ab,aba\}$. The
corresponding partition in towers is represented in Figure~\ref{figureTower}.
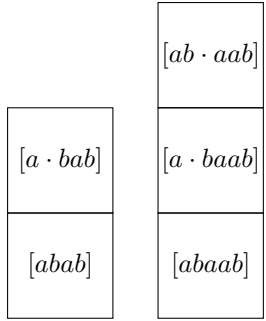
\begin{figure}[hbt]
\centering
\tikzset{node/.style={draw,minimum size=1.4cm,inner sep=0pt}}
\begin{tikzpicture}
\node[node](11)at(0,0){$[abab]$};
\node[node](12)at(0,1.4){$[a\cdot bab]$};
\node[node](21)at(2,0){$[abaab]$};
\node[node](22)at(2,1.4){$ [a\cdot baab]$};
\node[node](23)at(2,2.8){$ [ab\cdot aab]$};
\end{tikzpicture}
\caption{A partition in towers of the Fibonacci shift}\label{figureTower}
\end{figure}
\end{example}

\subsection{Partitions of substitution shifts}
The second example is a partition of a substitution shift.

\begin{proposition}\label{propositionPartitionSubstitution}
Let $\varphi:A^*\rightarrow A^*$ be a primitive substitution and let $X$
be the associated shift space. If $X$ is infinite, the family
\begin{displaymath}
\Pg(n)=\{S^j\varphi^n([a])\mid a\in A, 0\le j<|\varphi^n(a)|\}
\end{displaymath}
is, for every $n\ge 1$, a partition in towers with basis $\varphi^n(X)$.
\end{proposition}
\begin{proof}
Set $\psi=\varphi^n$. Since $\varphi$ is primitive and $X$ is infinite,
$\varphi$ is recognizable on $X$ and thus $\psi$ is recognizable on $X$.
This implies clearly that $\Pg(n)$ is a partition.
\end{proof}
The partition $\Pg(n)$ is called the \emph{partition associated to the substitution}\index{subject}{partition!in towers!associated to a substitution}\index{subject}{substitution!partition associated} $\varphi^n$.
We illustrate Proposition~\ref{propositionPartitionSubstitution} with the
following example.
\begin{example}\label{examplePartitionFibonacci}
Let $\varphi$ be the Fibonacci morphism and let $(X,S)$
be the Fibonacci shift.
\index{subject}{Fibonacci!shift}\index{subject}{shift space!Fibonacci}%
 The  partition $\Pg(n)$ corresponding
to Proposition~\ref{propositionPartitionSubstitution}
has two towers with basis
$\varphi^n([a])$ and $\varphi^n([b])$ respectively.
One has $\varphi^n([a])=[u_n]$ with $u_0=a$ and
where $u_n$ is, for $n\ge 1$, the
shortest right-special factor which has $\varphi(u_{n-1})$
as a prefix. A similar rule holds for $\varphi^n([b])$.
The partition in
towers for $n=2$ 
is identical with the partition
 represented in Figure~\ref{figureTower}.
\end{example}

\subsection{Sequences of partitions}

A partition $\alpha$ is a \emph{refinement} \index{subject}{partition!finer}\index{suject}{refinement!of partition}\index{subject}{partition!refinement} of
a partition $\rho$ if every element of $\alpha$ is a subset
of an element of $\rho$ or, equivalently,
if every element of $\rho$ is a union of elements of $\alpha$.
We also say that $\rho$ is \emph{coarser}\index{subject}{partition!coarser}
than $\alpha$. We denote $\alpha\ge \rho$.

Let $\mathcal{A}$ and $\mathcal{B}$ be two families of subsets of some set $X$.
We set
$$
\mathcal{A} \vee \mathcal{B} = \{A\cap B\mid A\in\mathcal{A} ,\ B\in\mathcal{B}\}.
$$
\index{symbols}{A@$\A\vee\B$}
For finitely many families $\mathcal{A}_1, \mathcal{A}_2, \dots , \mathcal{A}_n $ of subsets of $X$ we set 
$$
\bigvee_{i=1}^n \mathcal{A}_i = \mathcal{A}_1 \vee \mathcal{A}_2 \vee \dots  \vee \mathcal{A}_n 
$$
while $\A\cup\B$ is the union of  the families of $A\in \A$
and $B\in\B$.

We also say that a partition in towers $\Pg$ with basis $B$
is \emph{nested}
\index{subject}{partition!in towers!nested}%
\index{subject}{nested!partition}%
 in a partition in towers $\mathfrak P'$ with basis $B'$
if $B\subset B'$ and $\mathfrak P$ is a refinement, as a partition,
of $\mathfrak P'$. 

A sequence $(\Pg(n))$ of partitions in towers
is \emph{nested}
\index{subject}{nested!sequence of partitions}%
if $\Pg(n+1)$ is nested in $\Pg(n)$ for all $n\ge 1$.

The following statement shows that the sequence of  partitions
of a primitive substitution shift defined in Proposition \ref{propositionPartitionSubstitution} is nested.
\begin{proposition}\label{propositionPartitionSubstitutionNested}
  Let $\varphi:A^*\to A^*$ be a primitive substitution and let
  $X$ be the associated shift space. If $X$ is infinite, the
  sequence of partitions
  \begin{displaymath}
    \Pg(n)=\{S^j\varphi^n([a])\mid a\in A,0\le j<|\varphi^n(a)|\}
  \end{displaymath}
  is nested.
\end{proposition}
\begin{proof}
  If $\varphi(a)$ begins with $b$, we have $\varphi^{n+1}([a])\subset\varphi^n([b])$. Thus, the base of $\Pg(n+1)$ is contained in the base of $\Pg(n)$.
  Next, if $0\le j<|\varphi^{(n+1)}(a)|$,
  there is a factorization $\varphi(a)=xcy$ with $c\in A$ such that
  $|\varphi^n(x)|\le j<|\varphi^n(xc)|$. Set $k=j-|\varphi^n(x)|$.
  Then
  \begin{displaymath}
    S^j\varphi^{n+1}([a])\subset S^k\varphi^n([c])
  \end{displaymath}
  with $0\le k<|\varphi^n(c)|$. This shows that $\Pg(n+1)$ refines
  $\Pg(n)$.
  \end{proof}

We say that a sequence $(\Pg(n))_{n\ge 0}$ with
\begin{displaymath}
  \Pg(n)=\{T^jB_i(n)\mid 0\le j<h_i(n), 1\le i\le t(n)\}
\end{displaymath}
of KR-partitions of $X$
with bases $B(n)$ is a \emph{refining sequence}
\index{subject}{refining sequence}
if it satisfies the three following conditions.
\begin{itemize}
\item[(KR1)] $\cap_nB(n)=\{x\}$ for some $x\in X$, that is, the intersection of the bases consists in one point
  $x\in X$,
\item[(KR2)] the sequence $(\Pg(n))$ is nested,
\item[(KR3)] $\cup_{n\ge 1}\Pg(n)$ generates the topology of $X$
  (that is every open set is a union of elements of the partitions $\Pg(n)$).
\end{itemize}

Condition (KR3) can be expressed equivalently by the condition that
the sequence of partitions $(\Pg(n))$ tends to the point partition
(that is, for every $\varepsilon>0$, there is an $n$ such
that all elements of $\Pg(n)$ are contained in a ball of radius $\varepsilon$).

None of the conditions (KR1) or (KR3) implies the other one (see Exercises
\ref{exerciseKR1pasKR3}, \ref{exerciseKR3pasKR1}). It may be puzzling
that (KR1) does not imply (KR3). After all, the element $T^jB_i(n)$ of $\Pg(n)$
at height $j$ is the image by $T^j$ of the element $B_i(n)$ of $B(n)$ which tends to a point.
The reason is that $T$, although continuous, is not an isometry and that the ratio
of the diameters of $T^jB_i(n)$ and  $B_i(n)$ can be unbounded. This is
the phenomenon taking place in Exercise~\ref{exerciseKR1pasKR3}.

We will prove the following
important statement.

\begin{theorem}\label{theoremKRPartitions}
Let $(X,T)$ be a minimal Cantor system. There exists
a refining sequence of KR-partitions of $X$.
\end{theorem}
The proof uses the following technical result.

\begin{proposition}\label{propositionVersik}
Let $(X,T)$ be a minimal Cantor system.
Let ${\mathfrak Q}$ be a clopen partition of $X$ and $B$
be a clopen set. Then there exists a clopen partition
$B_1,\ldots,B_t$ of $B$ and integers $(h_i)_{1\le i\le t}$ such that
\begin{displaymath}
\Pg=\{T^jB_i\mid 0\le j<h_i,\ 1\le i\le t\}
\end{displaymath}
is finer than $\mathfrak{Q}$.
\end{proposition}
\begin{proof}
Let $r_B:X\rightarrow \Z$ be the first return map to $B$ defined by
$r_B(x)=\inf\{n>0\mid T^nx\in B\}$. Since $(X,T)$ is minimal,
$r_B$ is well-defined and continuous. Let $r_1,r_2,\ldots,r_{t'}$
be the set of values taken by $r_B$. For every
$i$ with $1\le i\le t'$, we define $B'_i=\{x\in B\mid r_B(x)=r_i\}$. Then
\begin{displaymath}
\Pg'=\{T^jB'_i\mid 0\le j<r_i,\ 1\le i\le t'\}
\end{displaymath}
is a clopen partition of $X$. Indeed, since $(B'_i)_{1\le i\le t'}$ is a partition of $B$,
the family $\Pg'$ is
formed of disjoint sets. Their union is a nonempty closed invariant subset
of $X$ and since $(X,T)$ is minimal, it is equal to $X$.
It is however not necessarily finer than $\mathfrak{Q}$.
Let $\mathfrak{Q}'=\{P'\cap Q\mid P'\in \Pg', Q\in\mathfrak{Q}\}$.
It suffices to find $\Pg$ finer than $\Qg'$ and which is a partition in towers.
Let $Q'$ be an atom of $\Qg'$. There exists a unique pair $(i_0,j_0)$
with $1\le i_0\le t'$
and $0\le j_0<r_{i_0}$ such that $Q'\subset T^{j_0}B'_{i_0}$.
Set $Q''=T^{j_0}B'_{i_0}\setminus Q'$.
We divide the tower $i_0$ into two new towers and obtain
a new KR-partition $\Pg''$ with $t'+1$ towers
\begin{eqnarray*}
\Pg''&=&\{T^jB'_i\mid 0\le j<r_i,1\le i\le t',i\ne i_0\}\\
     &&\cup \{T^jQ'\mid -j_0\le j<r_{i_0}-j_0\}\\
     && \cup \{T^j Q''\mid -j_0\le j<r_{i_0}-j_0\}.
\end{eqnarray*}
with a split of the $i_0$-tower propagating up and down the
split of $T^{j_0}B'_{i_0}$ in two parts, namely $Q'$
and $ Q''$.
We repeat this procedure for every atom of $\Qg'$. The result
is the desired KR-partition.
\end{proof}
\begin{proofof}{of Theorem~\ref{theoremKRPartitions}}
  Let $x\in X$.
We start choosing a decreasing sequence of clopen sets
$(B_n)_{n\ge 1}$ whose intersection is $\{x\}$ and an increasing
sequence of partitions $(\Pg'(n))_n$ generating the topology.
We apply Proposition~\ref{propositionVersik} to $\Qg=\Pg'(1)$
and $B=B_1$ to obtain $\Pg(1)$.

Applying Proposition~\ref{propositionVersik}
iteratively for $n\ge 2$ to $B=B_n$ and by setting now
\begin{displaymath}
\Qg=\Pg'(n)\vee\Pg(n-1),
\end{displaymath}
we obtain a partition $\Pg(n)$ with basis $B_n$ which is finer than $\Pg'(n)$ and $\Pg(n-1)$.

Condition (KR1) holds because for each $n$, the basis of $\Pg(n)$
is $B_n$ and $\cap_n B_n=\{x\}$ by hypothesis.

Condition (KR2) holds because, by construction, $\Pg(n)$ is nested in
 $\Qg$, which is nested in $\Pg(n-1)$.

Finally, condition (KR3) holds because  $\cup_{n\ge 1}\Pg(n)$ is finer
than $\cup_{n\ge 1}\Pg'(n)$, which generates the topology.
\end{proofof}
Note that, by definition,  
in a nested sequence $(\Pg(n))$ of partitions in towers, the sequence $B(\Pg(n))$ is decreasing.

We give two simple examples illustrating Theorem~\ref{theoremKRPartitions}.
The first one is the ring of $p$-adic integers (see
Section~\ref{sectionCantorSpaces}).
\begin{example}\label{examplepAdic}
We show that the odometer on the ring of $p$-adic integers
can be represented by a sequence of partitions in towers
with one tower.
For $n\ge 1$, let $B(n)=p^n\Z_p$, that is the ball of $\mathbb{Z}_p$ centered in $0$ of radius $p^n$. 
Then the family $\Pg(n)=\{T^jB(n)\mid 0\le j<p^n\}$ is, for
each $n$, a partition formed of one tower. It is easy to
verify that the sequence $(\Pg(n))$ satisfies the conditions
of Theorem~\ref{theoremKRPartitions}.
\end{example}
In the second example, we show how a nested sequence of
partitions can modified to become a refining sequence.
\begin{example}\label{exampleFibonacciPartitions}
Let $X$ be the two-sided Fibonacci shift.
\index{subject}{Fibonacci!shift}\index{subject}{shift space!Fibonacci}%
 For $n\ge 1$, let $\Pg(n)$
be the %\marginpar{FD: ce n est pas bien defini}
  partition $\Pg(n)=\{S^j\varphi^n([a]\mid a\in A, 0\le j<|\varphi^n(a)|\}$
(see Proposition~\ref{propositionPartitionSubstitution}).
Properties (KR1) and (KR3)
do not hold for this sequence of partitions. Indeed,
$\varphi^2:a\mapsto aba,b\mapsto ab$ has 
one right infinite fixed point $x$ (the \emph{Fibonacci word}\index{subject}{Fibonacci!word})
and two left infinite fixed points
$y,z$. The two points (actually fixed points) $y\cdot x$ and 
$z\cdot x$ belong to all $\varphi^{n}[a]$.

Using instead the substitution $\psi:a\mapsto baa,b\mapsto ba$
%\marginpar{FD: conjugate n est pas defini} 
(related to $\varphi^2$ by $a\psi(u)=\varphi^2(u)a$
for every word $u$), the shift defined remains the
same. Indeed, set $u_n=\varphi^{2n}(a)$ and $v_n=\psi^n(b)$. Then
$u_{n+1}a=a\psi(u_n)$ and $v_{n+1}=bu_nw_n$
with $w_n=a\psi(w_{n-1})$ (because $v_{n+1}=\psi(v_n)=\psi(bu_{n-1}w_{n-1})
=ba\psi(u_{n-1})\psi(w_{n-1})=bu_{n}a\psi(w_{n-1})=bu_nw_n$).
Thus $\psi^\omega(b)=b\varphi^\omega(a)$.
The  sequence of partitions $(\Pg(n))$ associated to the substitutions $\psi^n$ satisfies
the three conditions because there is a unique fixed point
(we shall see in Chapter~\ref{ch5:sec:examples} that this example describes a general situation).

\end{example}
%%%%%%%%%%%%%%%%%%%%%%%%
\section{Ordered group associated with a partition}
\label{sectionOrderedGroupPartition}
Let $(X,T)$ be a minimal Cantor system.
Let 
\begin{displaymath}
\mathfrak{P}=\{T^jB_i\mid 1\le i\le m,0\le j<h_i\}
\end{displaymath}
 be a KR-partition
 of $(X,T)$ built on  $B_1, \dots , B_m$ and with heights $h_1 , \dots , h_m$.
Denote $C(\Pg)$ \index{symbols}{C@$C(\Pg)$} the subgroup of $C(X,\Z)$ formed of the functions which are constant on every element of the partition $\Pg$ and 
$C^+(\Pg)=C(\Pg)\cap C(X,\Z_+)$. The triple $(C(\Pg),C^+(\Pg),\chi_X)$
is a unital ordered group.

Next, denote $G(\Pg)$ \index{symbols}{G@$G(\Pg)$} the subgroup of $C(B(\Pg),\Z)$ formed of the
functions constant on the basis $B_i$ of each tower,
denote $G^+(\Pg)=G(\Pg)\cap C(B(\Pg),\Z_+)$ and $\mathbf{1}_\Pg$ the function with value
$h_i$ on $B_i$.
We define the \emph{unital ordered group  associated to}
\index{subject}{unital!ordered group!associated to a partition} $\mathfrak{P}$
as the triple $(G(\Pg),G^+(\Pg),\mathbf{1}_\Pg)$.

Let $I(\Pg):C(\Pg)\rightarrow G(\Pg)$ \index{symbols}{I@$I(\Pg)$}
be the group morphism defined by
\begin{displaymath}
(I(\Pg)f)(x)=f^{(h_i)}(x)
\end{displaymath} 
for every $f\in C(\Pg)$ and $x\in B_i$. It is the restriction
to $C(\Pg)$ of the morphism $I_{B(\Pg)}$ introduced in Section
\ref{sectionInduced}. Thus, by Proposition~\ref{propositionIuRu}, we have
\begin{equation} \label{eqI(P)R_B(P)}
  R_{B(\Pg)}\circ I(\Pg)=\id_{C(\Pg)},\quad I(\Pg)\circ R_{B(\Pg)}=\id_{G(\Pg)}
     \end{equation}
where $R_{B(\Pg)}f$ is the map equal to $f$ on $B(\Pg)$ and equal to $0$
elsewhere.

By Proposition~\ref{propositionIuRu}, the kernel of $I(\Pg)$  consists in coboundaries
and there exists a morphism $\pi(\Pg):G(\Pg)\rightarrow H(X,T,\Z)$
\index{symbols}{pi@$\pi(\Pg)$}
which makes the diagram of Figure~\ref{figureDiagrampi(P)} commutative
(we denote by $\pi$ the canonical morphism from $C(X,\Z)$
onto $H(X,T,\Z)=C(X,\Z)/\partial_TC(X,\Z)$).
\begin{figure}[hbt]
\centering
\tikzset{node/.style={minimum size=0.2cm,inner sep=0.2pt}}
\begin{tikzpicture}(30,20)
\node[node](C)at(1.5,2){$C(\Pg)$};
\node[node](G)at(0,0){$G(\Pg)$};
\node[node](H)at(3,0){$H(X,T,\Z)$};

\draw[left,->](C)edge node{$I(\Pg$)}(G);
\draw[right,->](C)edge node{$\pi$}(H);
\draw[above,->](G)edge node{$\pi(\Pg$)}(H);
\end{tikzpicture}
\caption{The morphism $\pi(\Pg)$.}\label{figureDiagrampi(P)}
\end{figure}

\begin{proposition}\label{pi(P)}
The morphism $\pi(\Pg)$ defines a morphism of unital ordered groups
from $(G(\Pg),G^+(\Pg),\mathbf{1}_\Pg)$ to 
$K^0(X,T)$.
\end{proposition}
\begin{proof} It is clear that $\pi(\Pg)(G^+(\Pg))$ is included in $H^+(X,T,\Z)$.
Next, for any $x\in B_i$
\begin{displaymath}
(I(\Pg) \charac_X)(x)=h_i.
\end{displaymath}
Thus $\pi(\Pg)(\mathbf{1}_\Pg)=\pi(\chi_X)=\mathbf{1}_X$.
\end{proof}
Assume now that
\begin{displaymath}
\mathfrak{P}'=\{T^\ell B'_k\mid 1\le k\le m',\ 0\le \ell<h'_k\}
\end{displaymath}
is an other KR-partition in tower (with basis $B'_1, \dots , B'_{m'}$ and heights $h'_1 , \dots , h'_{m'}$)
which is nested in $\Pg$.
The morphism $R_{B(\Pg)}:C(B(\Pg),\Z)\rightarrow C(X,\Z)$ defined in Section~\ref{sectionInduced} maps $G(\Pg)$ to $C(\Pg)$, and thus in $C(\Pg')$.
The morphism 
\begin{displaymath}
I(\Pg',\Pg)=I(\Pg')\circ R_{B(\Pg)}
\end{displaymath}
\index{symbols}{I@$I(\Pg',\Pg)$}%
maps $G(\Pg)$ to $G(\Pg')$ (see the commutative diagram in Figure~\ref{figureI(P',P)}). It is clearly a morphism of ordered groups. We have
\begin{equation}
  I(\Pg',\Pg)\circ I(\Pg)=I(\Pg')\label{equationI(P',P)}
\end{equation}
where the right-hand side is actually
the restriction of $I(\Pg')$ to $C(\Pg)$.
Indeed,
$I(\Pg',\Pg)\circ I(\Pg)=I(\Pg')\circ R_{B(\Pg)}\circ I(\Pg)=I(\Pg')\circ\id_{C(\Pg)}$
by Equation~\eqref{eqI(P)R_B(P)}
Moreover,
\begin{equation}
\pi(\Pg')\circ I(\Pg',\Pg)=\pi(\Pg).
\end{equation}

\begin{figure}[hbt]
\centering
\tikzset{node/.style={circle,draw,minimum size=0.1cm,inner sep=0pt}}
\tikzset{title/.style={minimum size=0.2cm,inner sep=0pt}}
\begin{tikzpicture}(
\node[title](CP)at(0,3){$C(\Pg)$};\node[title](CP')at(3,3){$C(\Pg')$};
\node[title](GP)at(0,1.5){$G(\Pg)$};\node[title](GP')at(3,1.5){$G(\Pg')$};
\node[title](H)at(1.5,0){$H(X,T,\Z)$};

\draw[above](CP)edge node{$\subset$}(CP');
\draw[bend left, ->,right](CP)edge node{$I(\Pg)$}(GP);
\draw[bend left,->,left](GP)edge node{$R_{B(\Pg)}$}(CP);
\draw[->,right](CP')edge node{$I(\Pg')$}(GP');
\draw[->,above](GP)edge node{$I(\Pg',\Pg)$}(GP');
\draw[->,left](GP)edge node{$\pi(\Pg)$}(H);
\draw[->,right](GP')edge node{$\pi(\Pg')$}(H);
\end{tikzpicture}
\caption{The morphism $I(\Pg',\Pg)$}\label{figureI(P',P)}
\end{figure}

It can be useful to write the morphism $I(\Pg',\Pg)$ in matrix form.
We can make the following identification
\begin{displaymath}
G(\Pg)=\Z^m,\ G^+(\Pg)=\Z_+^m,\quad G(\Pg')=\Z^{m'}, \ G^+(\Pg')=\Z_+^{m'}.
\end{displaymath}
The units $\mathbf{1}_{\Pg}$ and $\mathbf{1}_{\Pg'}$ are identified
with the vectors of heights.
\begin{proposition}\label{propositionMatrixI(P,P')}
  The $m'\times m$-matrix $M(\Pg',\Pg)$
\index{symbols}{M@$M(\Pg',\Pg)$}%
of the morphism $I(\Pg',\Pg)$ is given by
\begin{equation}
M(\Pg',\Pg)_{k,i}=\Card\{\ell\mid 0\le\ell< h'_k, T^\ell B'_k\subset B_i\}
\label{equationMatrixPartitions}
\end{equation}
for $1\le k\le m'$ and $1\le i\le m$.
\end{proposition}
\begin{proof}
  Let
$f\in G(\Pg)$ be the function with value $\alpha_i$ on $B_i$
for $1\le i\le m$. Then,
for $x\in B'_k$, we have by \eqref{equationI(P',P)}
\begin{displaymath}
  I(\Pg',\Pg)f(x)=I(\Pg')\circ R_{B(\Pg)}f(x).
\end{displaymath}
Set $g=R_{B(\Pg)}f$. By definition of $R_{B(\Pg)}$, we have
\begin{displaymath}
  g(x)=\begin{cases}
  \alpha_i& \mbox{ if $x\in B_i$}\\0&\mbox{if $x\notin B(\Pg)$}\end{cases}
\end{displaymath}
Thus
\begin{displaymath}
  I(\Pg',\Pg)f(x)=I(\Pg')g(x)=g^{(h'_k)}(x)=\sum_{i=1}^m M(\Pg',\Pg)_{k,i}\alpha_i.
\end{displaymath}
\end{proof}
\begin{example}
Let $(X,S)$ be the two-sided Fibonacci shift.
\index{subject}{Fibonacci!shift}\index{subject}{shift space!Fibonacci}%
 Let $\Pg$ be 
the partition corresponding to the
return words on $a$. We have $\RR'_X(a)=\{a,ab\}$, and the
partition $\Pg$ is represented in Figure~\ref{figureTower2} on the left.
Let $\Pg'$ be
the partition corresponding to the return words
on $ab$ (see Example~\ref{examplePartition}) represented
in Figure~\ref{figureTower2} on the right.
\begin{figure}[hbt]
\centering
\tikzset{node/.style={draw,minimum size=1.2cm,inner sep=0pt}}
\begin{tikzpicture}

\node[node](11)at(0,0){$\cdot aa$};
\node[node][node](21)at(1.5,0){$\cdot aba$};
\node[node](22)at(1.5,1.2){$ a\cdot ba$};

\node[node](11)at(5,0){$\cdot abab$};\node[node](12)at(5,1.2){$a\cdot bab$};
\node[node](21)at(7,0){$\cdot abaab$};\node[node](22)at(7,1.2){$ a\cdot baab$};
\node[node](23)at(7,2.4){$ ab\cdot aab$};
\end{tikzpicture}

\caption{The partitions $\Pg$ and $\Pg'$.}\label{figureTower2}
\end{figure}

Since $a$ is a prefix of $ab$, the partition $\Pg'$ is nested
in $\Pg$.

The matrix $M(\Pg',\Pg)$ is
\begin{displaymath}
M(\Pg',\Pg)=\begin{bmatrix}0&1\\1&1\end{bmatrix}
\end{displaymath}
For example, the second row is $[1\ 1]$ because,
in the second tower of $\Pg'$, the lower
element  is contained in the 
basis of the second tower of $\Pg$ and the upper element in the basis
of the first one.
\end{example}
%%%%%%%%%%%%%%%%%%%%
\section{Ordered groups of sequences of partitions}
\label{sectionOrderedGroupSequences}
We now give a description of the ordered group of the minimal Cantor
system $(X,T)$ in terms of a refining sequence  $(\Pg(n))$ with
\begin{displaymath}
  \Pg(n)=\{T^j B_i(n)\mid 0\le j< h_i(n), 1\le i\le t(n)\}
\end{displaymath}
of $KR$-partitions. 
Denote for simplicity $B(n)$ the base of the partition $\Pg(n)$
and
\begin{eqnarray*}
    (G(n),G^+(n),\mathbf{1}_{n})&=&(G(\Pg(n)),\G^+(\Pg(n),\mathbf{1}_{\Pg(n)}),\\
    I(n+1,n)&=&I(\Pg(n+1),\Pg(n)).
\end{eqnarray*}
\begin{proposition}\label{propositionGroupSequencePartitions}
\label{sectionOrderdGroupSequencePartitions}
Let $(X,T)$ be an invertible minimal Cantor system 
%\marginpar{FD: homeo ?}. 
Let
$(\Pg(n))$ be a refining sequence of $KR$-partitions.
The  ordered group  $K^0(X,T)$
is the inductive limit 
of the unital ordered groups  $(G(n),G^+(n),\mathbf{1}_n)$
with the connecting morphisms $I(n+1,n)$.
\end{proposition}
The proof is given below. 
Before we need some lemmas. 

Let $(G,G^+,\mathbf{1})$ be the inductive limit of the sequence
of ordered groups 
$(G(n),G^+(n)),\mathbf{1}_{n})$ with the connecting
morphisms $I(n+1,n)=I(\Pg(n+1))$.
Let $i(n):G(n)\rightarrow G$ be the natural morphism.
Note that for $m>n$, we have $i(n)=i(m)\circ I(m,n)$.

\begin{figure}[hbt]
\centering

\begin{tikzpicture}(30,30)
\node(CP)at(1.5,3){$C(n)$};
\node(GP)at(0,1.5){$G(n)$};\node(H)at(3,1.5){$H(X,T,\Z)$};
\node(GS)at(1.5,0){$G$};

\draw[left,->](CP)edge node{$I(n)$}(GP);
\draw[right,->](CP)edge node{$\pi$}(H);
\draw[above,->](GP)edge node{$\pi(n)$}(H);
\draw[left,->](GP)edge node{$i(n)$}(GS);
\draw[right,->](GS)edge node{$\sigma$}(H);
\end{tikzpicture}
\caption{The morphism $\sigma$}\label{figureG(S)}
\end{figure}
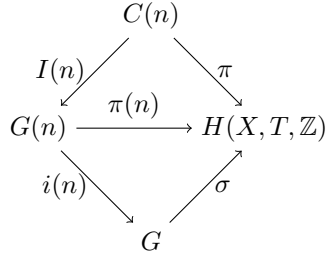

%The proof uses two lemmas.
\begin{lemma}\label{lemma1}
Let $(X,T)$ be an invertible minimal Cantor system.
Let $(\Pg(n))$ be a nested sequence of $KR$-partitions.
There is a unique morphism 
$\sigma:(G,G^+,\mathbf{1})\rightarrow K^0(X,T)$
such that $\sigma\circ i(n)=\pi(n)$ for every $n\ge 0$.
\end{lemma}
\begin{proof}
By Proposition~\ref{pi(P)}, the morphism $\pi(n):G(n)\rightarrow
H(X,T,\Z)$ is for every
$n\ge 0$ a morphism of unital ordered groups. 
By Proposition~\ref{propositionMorphismDirectLimit} there is
a unique morphism $\sigma:(G,G^+,\mathbf{1})\rightarrow K^0(X,T)$ such
that $\sigma\circ I(n)=\pi(n)$.
\end{proof}
We illustrate the proof in Figure~\ref{figureG(S)},
where the upper part reproduces Figure~\ref{figureDiagrampi(P)}
for $\Pg=\Pg(n)$ with $C(n)=C(\Pg(n))$, $I(n)=I(\Pg(n))$
 and $\pi(n)=\pi(\Pg(n))$).
\index{symbols}{C@$C(n)$}\index{symbols}{I@$I(n)$}\index{symbols}{pi@$\pi(n)$}%

The proof uses the following second lemma
(in which we do  not yet assume that $Y$ is reduced to one point,
as in Proposition~\ref{propositionGroupSequencePartitions}).
A sequence 
\begin{displaymath}
G_0\edge{\varphi_1}G_1\edge{\varphi_2}G_2\edge{}\cdots\edge{\varphi_n}G_n
\end{displaymath}
of group morphisms is \emph{exact}\index{subject}{exact sequence} if the image of each morphism
is the kernel of the next one.
\begin{lemma}\label{lemma2Partition}
Let $(X,T)$ be a minimal Cantor system with $T$ a homeomorphism. Let
$(\Pg(n))$ be a nested sequence of partitions
such that  $\Pg(n)$ converges to the point partition. 
Let $Y=\cap_{n\ge 0}B(n)$. 
There
is  a group morphism $r:C(Y,\Z)\rightarrow \ker(\sigma)$ such
that the sequence of group morphisms
\begin{equation}
0\rightarrow\Z\rightarrow C(Y,\Z)\stackrel{r}{\rightarrow}G\stackrel{\sigma}{\rightarrow}
H(X,T,\Z)\rightarrow 0 \label{eqSuiteExacte}
\end{equation}
is exact, with $\Z$ identified with the group of constant functions on $Y$
where $G$ is the inductive limit of the groups $G(n)$
and $\sigma$ is defined in Lemma~\ref{lemma1}.
%\marginpar{FD:X=Y}
\end{lemma}
\begin{proof}
  Recall that $C(n)$ is the set of functions
in $C(X,\Z)$ which are constant on every element of $\Pg(n)$.
  Since the sequence $(\Pg(n))$ tends to a partition in points,
  and since a continuous function is locally constant, we have 
$C(X,\Z)=\cup_nC(n)$.

This implies, as a first step, that $\sigma$ is surjective 
by Proposition~\ref{propositionMorphismDirectLimit}.

 As a second step, let us define the morphism $r$.
For this, let $u\in C(Y,\Z)$, let $h\in C(X,\Z)$ having $u$ as restriction to $Y$
and let $g=\partial_Th$.
For $n$ large enough, since $C(X,\Z)=\cup_n C(n)$, we have 
$g\in C(n)$. Let $f=I(n)( g)$ and let $\alpha=i(n)(f)$.

Since $\alpha$
is the image in $G$ of a coboundary, it belongs to the
kernel of $\sigma$. Indeed, we have (see Figure~\ref{figureG(S)})
\begin{eqnarray*}
\sigma (\alpha)&=&\sigma(i(n)f)=\pi(n)( f)\\
&=&\pi(n)\circ I(n)(g)=\pi(g)=0.
\end{eqnarray*}
Thus the map $r:u\mapsto \alpha$ is
a group morphism from $C(Y,\Z)$ onto $\ker(\sigma)$.

Let us finally prove that $\alpha=0$ if and only if $u$ is constant.

If $u$ is constant, $h$ is constant on $Y$ and thus on some neighborhood of $Y$.
For $m\ge n$ large enough, since the
sequence $B(n)$ is decreasing, $h$ is constant on $B(m)$.
We may  assume that $m$ is large enough
to have also $g=\partial_Th\in C(m)$. Then, by
Equation~\eqref{eq1}, for any $x\in B_i(m)$, since
$T^{h_i(m)}x$ belongs to $B(m)$, we have
\begin{equation}
I(m)g(x)=g^{(h_i(m))}(x)=h(T^{h_i(m)}x)-h(x)=0.\label{eqI(m)g(x)}
\end{equation}
Since $I(m)g=I(m,n)\circ I(n)g=I(m,n)f$,
we have also $I(m,n)f=0$ and thus $\alpha=i(n)f=i(m)\circ I(m,n)f=0$.

Conversely, if $\alpha=0$, there exists $m\ge n$ such that
$0=I(m,n)f=I(m)g$. Thus, by Equation~\ref{eqI(m)g(x)}, $h$ is constant
on a set which is dense in $B(m)$, which implies that it is constant
on $B(m)$ and therefore that $u$ is constant.

It follows that the kernel of $r$ is the group of constant
functions on $Y$.

This completes the proof that the sequence of Equation~\eqref{eqSuiteExacte}
is an exact sequence.
\end{proof}
\begin{proofof}{of Proposition~\ref{propositionGroupSequencePartitions}}
Since the intersection of the bases $B(n)$ consists in one point,
we have $\ker(r)=C(Y,\Z)$ and ${\rm Im}(r)=0$. Thus $\sigma$ is an isomorphism and
the exact sequence displayed in~\eqref{eqSuiteExacte} reduces to the isomorphism
of $G$ with $H(X,T,\Z)$.
\end{proofof}

It is often convenient to write the direct limit 
\begin{equation}
  G(0)\ledge{I(1,0)}G(1)\ledge{I(2,1)}G(2)\ledge{I(3,2)}\cdots\label{eqDirectLimitAbstract}
\end{equation}
defining $K^0(X,T)$ in terms of matrices. Each $G(n)$ is isomorphic to
$\Z^{t(n)}$
with the natural order and the morphism $I(n+1,n)$ is represented by the
matrix $M(n+1,n)=M(\Pg(n+1),\Pg(n))$
  defined by \eqref{equationMatrixPartitions}. Thus the direct
  limit \eqref{eqDirectLimitAbstract} can be written
  \begin{equation}
    Z^{t(0)}\ledge{M(1,0)}\Z^{t(1)}\ledge{M(2,1)}\Z^{t(2)}\ledge{M(3,2)}\cdots\label{eqDirectLimitConcrete}
  \end{equation}
  The matrices $M(n+1,n)$, called the \emph{connecting matrices},
  \index{subject}{connecting!matrix}\index{subject}{matrix!connecting}%
  will play an important role later
  when we introduce Bratteli diagrams.

We deduce from Proposition~\ref{propositionGroupSequencePartitions}
the following important result.
\begin{theorem}[Herman, Putnam, Skau]\label{theoremDimensionGroup}
For any minimal Cantor system $(X,T)$, with $T$ a homeomorphism, the ordered group $K^0(X,T)$
%\marginpar{FD: j ai ajoute T homeo}
is a simple dimension group.
\end{theorem}
\begin{proof}
By Theorem \ref{theoremKRPartitions}, there exists  a sequence
$(\Pg(n))$ of partitions in towers satisfying the hypotheses
of Proposition \ref{propositionGroupSequencePartitions}.
Thus $K^0(X,T)$ is the direct limit $G$ of the ordered groups
$(G(n))$. Since each $G(n)$ is isomorphic to some $\Z^m$
with the natural order, it follows that $K^0(X,T)$
is a dimension group.

For  $n\ge 1$, let $[i,n]$ be the class of the characteristic function $\charac_{B_i(n)}$
of  an element $B_i(n)$ of the basis 
%\marginpar{FD: le "n" est il au bon endroit ?} 
of $\Pg(n)$, let $I(i,n)$ be the ideal of $G$ formed
of the classes of the functions $f\in G$ such that
$-k[i,n]\le f\le k[i,n]$ for some positive integer $k$. 
It is
easy to see that $I(i,n)$ is an order ideal. Next, suppose that
$I$ is any nonzero order ideal in $G$. Let
$x$ be a nonzero positive element of $I$. It must
be represented by some strictly positive element in some
$G(n)$ and, if  $B_i(n)$ is an element of
$B(n)$ on which $x$ is positive, we have $0\le [i,n]\le x$
and hence $[i,n]$ is also in $I$. It follows that $I(i,n)\subset I$.

We have shown that if $G$ contains a nonzero
order ideal, then it contains one
of the form $I(i,n)$. On the other hand,
since $(X,T)$ is minimal, one has $I(i,n)=G$
for every element $B_i(n)$ of $B(n)$. Indeed,
consider the connecting matrix $M(m,n)=M(\Pg(m),\Pg(n))$.
Since $(X,T)$ is minimal, there
is an $m>n$ such that the matrix $M(m,n)$ has
all its elements positive. Thus $k[i,n]$ can exceed
any positive element $x$ of $G(m)$, which implies
$x\in I(i,n)$. Thus we conclude that $G$
is simple.
\end{proof}
In agreement with Theorem~\ref{theoremDimensionGroup}, the group
$K^0(X,T)$ for a minimal Cantor system $(X,T)$, with $T$ a homeomorphism, is called
the \emph{dimension group}\index{subject}{dimension group!of minimal dynamical system}
 of $(X,T)$.

We give two examples illustrating Theorem~\ref{theoremDimensionGroup}.
Many more examples will appear later.

\begin{example}\label{example2adic}
The dimension group of the odometer on the ring of $2$-adic integers
is the group of dyadic rationals. Indeed, we have seen 
in Example~\ref{examplepAdic} that it
can be represented by a sequence $(\Pg(n))$ of partitions with one
tower and the map $i_{n+1,n}$ is easily seen to be the multiplication
by $2$. Thus the dimension group $K^0(X,T)$ is the group
$\Z[1/2]$ (see Example~\ref{exampleDyadic}).
\end{example}

\begin{example}
We have seen before (Proposition~\ref{theoremDGSturm}) that the dimension
group of a Sturmian shift of slope $\alpha$ is $\Z+\Z\alpha$,
which is a simple dimension group.
\index{subject}{Sturmian!shift space}\index{subject}{shift space!Sturmian}%
\end{example}

%%%%%%%%%%%%%%%%%%%%%%%%%%%%%%
% section Dimension groups and return words
%%%%%%%%%%%%%%%%%%%%%%%%%%%%%%%%
\section{Dimension groups and return words}
\label{chapterReturnWords}

In this section, we show how to use return words
to compute the dimension group of a minimal shift space.
This method has the advantage of using in general abelian groups
of smaller dimension than with the cylinder functions
(as seen in Section~\ref{sectionGroupCylinder}).
 For example, in a Sturmian shift space,
the number of return words is constant while the word
complexity is linear.

In the first part,
we show that the group $K^0(X,T)$ is, for every
minimal shift space, the direct limit of a sequence
of groups associated with return words. In the
second part
we illustrate the use return words to compute the dimension groups 
of episturmian shifts.
%%%%%%%%%%%%%%%%%%%%
\subsection{Sequences of return words}\label{sectionReturnWords}
Let $X$ be a minimal shift space.
We have already introduced in Section~\ref{sectionReturnWords},
the sets $\RR_X(w)$ and $\RR'_X(w)$ of
right and left return words to $w\in \cL(X)$. 

We fix a point $x\in X$ and we denote $W_n(x)=\RR'_X(x_{[0,n-1]})$.
\index{symbols}{W@$W_n(X)$}
Let $G_n(x)$\index{symbols}{G@$G_n(x)$} be the group of maps from $W_n(x)$ to $\Z$, $G_n^+(x)$
the subset of nonnegative ones and $\mathbf{1}_n(x)$ the map
which associates to $v\in W_n(x)$ its length.

Since $x_{[0,n-1]}$ is a prefix of $x_{[0,n]}$, the set $W_{n+1}(x)$
is contained in the submonoid generated by the set $W_n(x)$
%\marginpar{FD:definir l etoile ?}
This means that for every $v\in W_{n+1}(x)$ there is a  decomposition
$v=w_1(v)w_2(v)\cdots w_{k(v)}(v)$ as a concatenation of
elements of $W_n(x)$. The decomposition is moreover easily seen to be unique.
For every
$\phi\in G_n(x)$, let $i_{n+1,n}\phi\in G_{n+1}(x)$ be defined by
\begin{displaymath}
(i_{n+1,n}\phi)(v)=\sum_{i=1}^{k(v)}\phi(w_i(v)).
\end{displaymath}
\index{symbols}{i@$i_{n+1,n}(x)$}%
It is clear that $i_{n+1,n}:(G_n(x),G_n^+(x),\mathbf{1}_n(x))
\rightarrow(G_{n+1}(x),G_{n+1}^+(x),\mathbf{1}_{n+1}(x))$ is a morphism
of unital ordered groups. Indeed, in particular, for $v\in W_{n+1}(x)$,
\begin{displaymath}
i_{n+1,n}(\mathbf{1}_n(x))(v)=\sum_{i=1}^{k(v)}\mathbf{1}_n(x))(w_i(v))=
\sum_{i=1}^{k(v)}|w_i(v)|=|v|=\mathbf{1}_{n+1}(x)(v).
\end{displaymath}
We will prove the following result using the partition in towers
associated with return words (see Proposition~\ref{propositionPartitionReturn}).
However, we will not be able to use Proposition~\ref{propositionGroupSequencePartitions} because the sequence of partitions in towers associated with
the sets $W_n(x)$ do not satisfy the hypotheses of
Proposition~\ref{propositionGroupSequencePartitions} (the intersection
of the bases does not consist in one point).

\begin{proposition}\label{propositionRecurrentSubshift}
Let $X$ be a minimal shift space. The group $K^0(X,S)$
is the inductive limit of the family $(G_n(x),G_n^+(x),\mathbf{1}_n(x))$
with the morphisms $i_{n+1,n}$.
\end{proposition}
\begin{proof}
By Proposition~\ref{propositionPartitionReturn}, for every $n> 0$,
the family $\Pg_n=\{S^j[vx_{[0,n-1]}]\mid v\in W_n(x), 0\le j<|v|\}$
is a partition in towers with basis $[x_{[0,n-1]}]$. We can identify
$G(n)=G(\Pg_n)$ with $G_n(x)$ so that
\begin{displaymath}
(G_n(x),G_n^+(x),\mathbf{1}_n(x))=(G(n),G^+(n),\mathbf{1}_{n})
\end{displaymath}
The morphisms $i_{n+1,n}$ are then identified to the morphisms
$I(n+1,n)$. Let 
$(G,G^+,\mathbf{1})$ be the direct limit
of the sequence of ordered groups $(G(n),G^+(n),\mathbf{1}_{n})$
with the morphisms $I(n+1,n)$. By Lemma~\ref{lemma1},
there is a unique morphism $\sigma:G\rightarrow H(X,T,\Z)$ such that
$\sigma\circ I(n)=\pi(n)$ for every $n\ge 1$.
We show that $\sigma$ is an isomorphism.

Let $y,z\in S^j[vx_{[0,n-1]}]$ for some $n\ge 1$, some $v\in W_n(x)$ and some
$j$ with $0\le j<|v|$. Then $y,z$ have the same prefix of length
$|v|-j+n$ and thus the same prefix of length $n$. Thus
$y_{[0,n-1]}=z_{[0,n-1]}$. It follows that the partition $\Pg_n$ is finer than
the partition in $n$-cylinders, and that $C(\Pg_n)$ contains
every cylinder function corresponding to some $\phi\in Z_n(X)$.
Consequently $\cup_{n\ge 1} C(\Pg_n)$ contains every cylinder function.
Since every $f\in C(X,\Z)$ is cohomologous to some cylinder function,
the morphism $\sigma$ is onto and maps $G^+$ onto
$H^+(X,S,\Z)$.

Assume now that $\alpha\in G$ is in the kernel of $\sigma$.
For some $n$, $\alpha$ is the image in $G$ of some $f\in C(\Pg_n)$
and $f$ is a coboundary. The function $g=R_{B(\Pg_n)}\circ I(n)f$
is, by Proposition~\ref{propositionIuRu},
 cohomologous to $f$ and thus is also a coboundary.
It is a cylinder function because it is constant on $[vx_{[0,n-1]}]$
for each $v\in W_n(x)$ and null outside $[x_{[0,n-1]}]$.
Thus, by Lemma~\ref{lemma4ii}, it is the coboundary of some
cylinder function $h$. For $m$ large enough, $h$
is constant on the basis $x_{[0,m-1]}$ of $\Pg_m$.
This implies that $I(m)f=I(m)g=0$ and the image $\alpha$
of $f$ in $G$ is $0$.
\end{proof}
\begin{example}\label{exampleSturmian2}
Let $X$ be the two-sided Fibonacci shift.
\index{subject}{Fibonacci!shift}\index{subject}{shift space!Fibonacci}%
 Let $x\in X$
be such that $x_0x_1\cdots$ is the Fibonacci word. We have
\begin{eqnarray*}
W_1(x)&=&\{a,ab\}\\
W_2(x)&=&\{ab,aba\}\\
W_3(x)&=&\{ab,aba\}\\
W_4(x)&=&\{aba,abaab\}
\end{eqnarray*}
In general, one has $W_{n+1}(x)=W_n(x)$ if $x_{[0,n-1]}$ is not right special
and otherwise $W_{n+1}(x)=\{v,vu\}$ if $W_n(x)=\{u,v\}$ with
$|u|<|v|$. Indeed, the right-special prefixes of $x$ are
its palindrome prefixes $u_n$ and $\RR'_X(u_n)=\{\varphi^n(a),\varphi^n(b)\}$
(this results from Equation~\ref{equationReturnPal}).
The set $W_n(x)$ has two elements for every $n\ge 1$
and thus $G_n(x)=\Z^2$ for every $n\ge 1$. 

The fact
that $W_n(x)$ has always two elements 
is true for every Sturmian shift, and thus
also that $G_n(x)=\Z^2$ for all $n\ge 1$, as we have already
seen in Example~\ref{exampleSturmian1}.
\end{example}
%%%%%%%%%%%%%%%%%%%
\subsection{Dimension groups of Arnoux-Rauzy shifts}\label{sectionDimensionGroupsSturmian}
We now use return words to describe the dimension group
of an arbitrary Arnoux-Rauzy shift. 
\index{subject}{Arnoux-Rauzy!shift}%
\index{subject}{shift space!Arnoux-Rauzy}%
\index{subject}{dimension group!of Arnoux-Rauzy shift}%
Recall from Section~\ref{sectionSturmianShifts}
that if $s$ is a standard episturmian sequence,
\index{subject}{episturmian!shift} there is
a word $x=a_0a_1\cdots$ called its directive word such that
$s=\Pal(x)$. Moreover, the words $u_n=\Pal(a_0\cdots a_{n-1})$
are the palindrome prefixes of $s$ and the set
of left return words to $u_n$ is, by Equation~\eqref{equationReturnPal},
\begin{displaymath}
\RR'_X(u_n)=\{L_{a_0\cdots a_{n-1}}(a)\mid a\in A\}
\end{displaymath}
Denote by $M_a$ the incidence matrix of the morphism $L_a$. Thus,
if $A=\{a,b\}$, we have
\begin{displaymath}
M_a=\begin{bmatrix}1&0\\1&1\end{bmatrix},\quad 
M_b=\begin{bmatrix}1&1\\0&1\end{bmatrix}.
\end{displaymath}
The following result allows us to compute the dimension group
of a strict episturmian shift, also called Arnoux-Rauzy shift.
\begin{proposition}\label{propositionDimensionEpisturmian}
Let $s$ be a strict and standard episturmian sequence on the alphabet $A$, let $x=a_0a_1\cdots$ be its directive
sequence  and let $X$
be the shift generated by $s$. The dimension group of $X$ is the direct
limit of the sequence
\begin{equation}
\Z^A\ledge{M_{a_0}}\Z^A\ledge{M_{a_1}}\Z^A\ledge{M_{a_2}}\ldots\label{eqDirectLimitAR}
\end{equation}
\end{proposition}
\begin{proof}
Set $u_n=\Pal(a_0\cdots a_{n-1})$ and $\alpha_n=|u_n|$. We set $W_n(s)=\RR'_X(s_{[0,n-1]})$
and we use the notation of the previous section with $s$ in place of $x$.
We identify $G_n(s)$ with $\Z^A$ via the bijection
 $a\mapsto L_{a_0\cdots a_{n-1}}(a)$ from $A$ onto $W_n(s)$.
It is enough to prove that the matrix of the map $i_{\alpha_{n+1},\alpha_n}$
is the matrix $M_{a_n}$. Note that for every $a\in A$,
since $L_{a_0\cdots a_n}(a)=L_{a_0\cdots a_{n-1}}(L_{a_n}(a))$, we have,
denoting $L_{[0,n]}=L_{a_0\cdots a_n}$,
\begin{equation}
L_{[0,n]}(a)=\begin{cases}
L_{[0,n-1]}(a_n)L_{[0,n-1]}(a)&\mbox{if $a\ne a_n$}\\
L_{[0,n-1]}(a)&\mbox{otherwise}.
\end{cases}\label{eqDecompw_i}
\end{equation}
This gives the decomposition of an element of $W_{\alpha_{n+1}}(s)$ 
as a product of elements of $W_{\alpha_n}(s)$.
Consider now $\phi\in G_{\alpha_n}(s)$. We consider $\phi$ as a column
vector with components $\phi_a$ for $a\in A$ (via the identification above).
Then, by \eqref{eqDecompw_i}, we have
\begin{displaymath}
(i_{\alpha_{n+1},\alpha_n}\phi)(b)=\begin{cases}
\phi_{a_n}+\phi_b&\mbox{if $b\ne a_n$}\\
\phi_b&\mbox{otherwise.}
\end{cases}
\end{displaymath}
This shows that $i_{\alpha_{n+1},\alpha_n}\phi=M_{a_n}\phi$ and completes the proof.
\end{proof}
The description of the direct limit $G$ of a sequence given in~\eqref{eqDirectLimitAR} with invertible matrices $M_{a_i}$ is similar to the case where all matrices $M_{a_i}$ are equal
(see Proposition~\ref{propositionDeltaM}). Indeed, set $M_n=M_{a_n}$. Let
\begin{equation}
  G=\{x\in \R^A\mid M_n\cdots M_0 x\in \Z^A\mbox{ for $n$ large enough}\}
  \label{eqDimensionGroupNonhomogeneous}
\end{equation}
and let
\begin{equation}
  G^+=\{x\in G\mid M_n\cdots M_0 x\in \Z_+^A\mbox{ for $n$ large enough }\}.
  \label{eqDimensionGroupNonhomogeneous+}
    \end{equation}
Then the direct limit of the sequence \eqref{eqDirectLimitAR} is
$(G,G^+,\1)$ where $\1$ is the vector with all components equal to $1$.
Indeed, let
\begin{displaymath}
  \Delta=\{(x_n)_{n\ge 0}\mid x_n\in\Z^A\mbox{ for all $n$}, x_{n+1}=M_{n}x_n
  \mbox{ for $n$ large enough}\},
\end{displaymath}
\begin{displaymath}
  \Delta^0=\{(x_n)_{n\ge 0}\in\Delta\mid x_n=0
  \mbox{ for $n$ large enough}\},
  \end{displaymath}
  and
   \begin{displaymath}
  \Delta^+=\{(x_n)_{n\ge 0}\in\Delta\mid  x_n\in\Z_+^A\mbox{ for $n$ large enough}\}.
  \end{displaymath}

Let $(x_n)_{n\ge 0}\in\Delta$ be such
that $x_{n+1}=M_{n}x_n$ for all $n\ge m$. Since $M_{m-1}\cdots M_0$
is invertible, there is a unique $x\in \R^A$ such that
$M_{m-1}\cdots M_0 x=x_m$. This defines a morphism from
$\Delta$ onto $G$ with kernel $\Delta^0$. Since this morphism
sends $\Delta^+$ onto $G^+$ and the vector $\1$ to itself,
this proves the claim.

We give two examples with the Fibonacci shift (which is Sturmian
and thus we already know its dimension group by Theorem~\ref{theoremDGSturm})
and the Tribonacci shift.
\begin{example}\label{exampleFibonacci4}
Let $s$ be the Fibonacci word, which generates the Fibonacci shift $X$
(see Example~\ref{exampleFibonacci0}). 
\index{subject}{Fibonacci!shift}\index{subject}{shift space!Fibonacci}%
The directive word of $s$
is $x=(ab)^\omega$. Indeed, one has by Justin's Formula
$x=L_{ab}(x)$ whence the result since
$L_{ab}=\varphi^2$ where $\varphi$ is the Fibonacci morphism.
It follows from Proposition~\ref{propositionDimensionEpisturmian}
that the dimension group of $X$ is the ordered group $\Delta_M$
of the matrix $M=M_bM_a=M(\varphi)^2$ where $\varphi$ is the Fibonacci
morphism. Thus we prove again that the dimension group of $X$ is the group
of algebraic integers
$\Z+\frac{1+\sqrt{5}}{2}\Z$, in agreement with Theorem~\ref{theoremDGSturm}.
\end{example}
\begin{example}\label{exampleTribonacci3}
Let now $s$ be the Tribonacci word which is the fixed point of
the morphism $\varphi:a\mapsto ab,b\mapsto ac,c\mapsto a$
(see Example~\ref{exampleTribonacci2}). Its directive word is,
as we have seen, $x=(abc)^\omega$. Thus the dimension group
of the Tribonacci shift $X$ generated by $s$ is the group $\Delta_M$ of the
incidence matrix $M$ of the morphism $\varphi$.
\index{subject}{Tribonacci!shift}\index{subject}{shift space!Tribonacci}%
We have
\begin{displaymath}
M=\begin{bmatrix}1&1&0\\1&0&1\\1&0&0\end{bmatrix}
\end{displaymath}
The dominant eigenvalue is the positive real number $\lambda$
such that $\lambda^3=\lambda^2+\lambda+1$. A corresponding row 
eigenvector is $[\lambda^2,\lambda,1]$. Thus 
$K^0(X,S)$ is isomorphic to $\Z[\lambda]$.
\end{example}
\subsection{Unique ergodicity of Arnoux-Rauzy shifts}
We have already seen that Sturmian shifts are uniquely ergodic
(Exercise~\ref{exerciseErgodicityRotations}). In fact,
the following holds more generally.\index{subject}{unique ergodicity! of Arnoux-Rauzy shifts}
\begin{theorem}\label{ARUniquelyErgodic}
Every Arnoux-Rauzy shift is uniquely ergodic.
\end{theorem}
%The proof uses the \emph{projective distance}\index{subject}{projective distance}
%defined on
%the set $S_n$
%of positive vectors of dimension $n$ with sum $1$ by
%\begin{displaymath}
%  d(x,y)=\log \frac{Q(x,y)}{q(x,y)}
%\end{displaymath}
%where $Q(x,y)=\max_{i=1}^n(x_i/y_i)$ and $q(x,y)=\min_{i=1}^n(x_i/y_i)$.
%It is easy to verify that $d$ is a distance on $S_n$
%(Exercise~\ref{exerciseProjectiveDistance}). This distance
%measures the 'angle' made by $x,y$ in the space $\R^n$.

%For a
%positive $n\times n$-matrix $M$, one has (see Exercise~\ref{exerciseContraction})
%for every $x,y\in S_n$
%\begin{equation}
%  d(Mx,My)\le d(x,y)
%  \end{equation}
%  expressing the fact that the action of $M$ 'contracts' the vectors in $S_n$.

 Let $s$ be a strict and standard  episturmian word and let $x=a_0a_1\cdots$
 be its directive word. Set $M_n=M_{a_n}^t$
 (thus $M_n$ is the composition matrix of $L_{a_n}$)
 and $M_{[0,n)}=M_{0}\cdots M_{n-1}$.
\begin{lemma}\label{lemmaAvilaDelecroix}
  There exists a sequence of matrices $\tilde{M}_n$ such that
  $\|\tilde{M}_{[0,n)}\|_\infty\le 1$ and
  \begin{equation}
    \tilde{M}_{[0,n)}x =M_{[0,n)}x\label{eqDHS}
    \end{equation}
        for every $x\in (M_{[0,n)})^{-1}\1^\perp$ where $\1$
          is the vector with all coefficients equal to $1$.
\end{lemma}
\begin{proof}
  
  Set
  \begin{displaymath}
v^{(n)}=\frac{^t\1 M_{[0,n)}}{\|^t\1M_{[0,n)}\|_1}.
  \end{displaymath}
 
  Define $\tilde{M}_n=M_n-V_n$ where $V_n$ is the matrix with all rows
  equal to $0$ except the row of index $a_n$, which is equal to
  $v^{(n+1)}/\|v^{(n+1)}\|_\infty$.

  Let us first prove that $\|\tilde{M}_n\|_\infty\le 1$
  (since $\tilde{M}_n$ has rows equal to a unit vector,
  we will actually have $\|\tilde{M}_{n}\|_\infty= 1$
  but this will imply $\|\tilde{M}_{[0,n)}\|_\infty\le 1$). For this, it is
  enough to prove that for every $n\ge 0$,
  \begin{equation}
    \|v^{(n)}\|_\infty\le \frac{1}{d-1}\label{eqnormevn}
  \end{equation}
  with $d=\Card(A)$. Indeed, the  the row $r$ of index $a_n$ of
  $\tilde{M}_n$ will be such that
  \begin{displaymath}
    \|r\|_1=d-v^{(n+1)}/\|v^{(n+1)}\|_\infty\le d-(d-1)=1.
    \end{displaymath}
  We prove \eqref{eqnormevn} by induction on $n$. It is true for $n=0$
  since $v^{(0)}=\begin{bmatrix}1/d&\ldots&1/d\end{bmatrix}$. Next, assume that it holds for $n$
  and thus that $v^{(n)}$ belongs to the simplex
    \begin{equation}
      S=\{x\in\R_+^A\mid \|x\|_1=1,\|x\|_\infty\le \frac{1}{d-1}\}.
      \label{eqSimplexAR}
    \end{equation}
    The extreme points of this simplex are the vectors
    with all coordinates equal to $1/(d-1)$ except one which is $0$.
    We may suppose that $a_n$ is the first index. The image by $M_n$ of
    the extreme points are the vectors
    \begin{displaymath}
      \begin{bmatrix}0&\frac{1}{d-1}&\ldots&\frac{1}{d-1}\end{bmatrix},
   \begin{bmatrix}\frac{1}{d-1}&\frac{1}{d-1}&\frac{2}{d-1}&\ldots\end{bmatrix},
     \begin{bmatrix}\frac{1}{d-1}&\frac{2}{d-1}&\frac{1}{d-1}&\ldots\end{bmatrix},\ldots
       \end{displaymath}
    All but the first have sum $2$. After normalizing to sum $1$ we obtain the vectors
    \begin{displaymath}
      \begin{bmatrix}0&\frac{1}{d-1}&\ldots&\frac{1}{d-1}\end{bmatrix},
   \begin{bmatrix}\frac{1}{2d-2}&\frac{1}{2d-2}&\frac{1}{d-1}&\ldots\end{bmatrix},
     \begin{bmatrix}\frac{1}{2d-2}&\frac{1}{d-1}&\frac{1}{2d-2}&\ldots\end{bmatrix},\ldots
    \end{displaymath}
    which all belong to $S$. This proves \eqref{eqnormevn}.

    Next, we also prove \eqref{eqDHS} by induction on $n$. It holds
    trivially for $n=0$.
  For $x\in (M_{[0,n+1)})^{(-1)}\1^\perp$, we have $\langle v^{(n+1)},x\rangle=0$
    and thus $\tilde{M}_nx=M_nx-V_nx=M_nx$. We have then
    \begin{displaymath}
      \tilde{M}_{[0,n+1)}x=\tilde{M}_{[0,n)}\tilde{M}_nx=\tilde{M}_{[0,n]}M_nx.
    \end{displaymath}
    Since $M_nx\in (M_{[0,n)})^{-1}\1^\perp$, we obtain, using the induction
    hypothesis,
    \begin{eqnarray*}
      \tilde{M}_{[0,n+1)}x&=&\tilde{M}_{[0,n]}M_nx\\
        &=&M_{[0,n)}M_nx=M_{[0,n+1)}x.
            \end{eqnarray*}
\end{proof}
\begin{example}
  Assume that $A=\{a,b,c\}$ and that $a_0=a$. Then 
  \begin{displaymath}
    M_0=\begin{bmatrix}1&1&1\\0&1&0\\0&0&1\end{bmatrix}\mbox{ and }
    v^{(1)}=\begin{bmatrix}\frac{1}{5}&\frac{2}{5}&\frac{2}{5}\end{bmatrix}
  \end{displaymath}
  Thus, we find
  \begin{displaymath}
    \tilde{M}_0=\begin{bmatrix}\frac{1}{2}&0&0\\0&1&0\\0&0&1\end{bmatrix}
    \end{displaymath}
\end{example}
Note that the simplex $S$ defined by \eqref{eqSimplexAR}
plays an important role for Arnoux-Rauzy shifts (see Exercise~\ref{exerciseSimplexAR}).

\begin{proofof}{of Theorem~\ref{ARUniquelyErgodic}}

 %Since every letter occurs infinitely often
%  in $x$, we can write $x=w_1w_2\cdots$ where every letter occurs
%  in each word $w_i$. Let $M_{w_i}$ be the matrix of the morphism $L_{w_i}$.
%  Since every letter appears in $w_i$, the matrix $M_{w_i}$ is positive.
%  By Theorem~\ref{theoremBirkhoffContraction}, the matrix
%  \begin{displaymath}
%M_n=M_{w_1}M_{w_2}\cdots M_{w_n}
%  \end{displaymath}
  %  tends to a matrix of rank  $1$.
  Let $v\in\cap_{n\ge 0}M_{[0,n)}\R^A_+$ with $\|v\|_1=1$.
    We will prove that $v$ is unique, which will imply
    by Equation~\ref{eqDimensionGroupNonhomogeneous+} that
  the dimension group of $x$ has a unique state, whence the
  conclusion that $X$ is uniquely ergodic by Theorem~\ref{propositionKerov}.

  Let $\pi_n$ be the projection on $ (M_{[0,n)})^{-1}\1^\perp$
    along $(M_{[0,n)})^{-1}v$. Note that $\pi_0M_{[0,n)}=M_{[0,n)}\pi_n$.
          Let $i_n\in\Z^A$ be the characteristic
  vector of $a_n$.
  By Lemma~\ref{lemmaAvilaDelecroix}, we have 
  \begin{eqnarray*}
    \|\pi_0 M_{[0,n)} i_n\|_\infty&=&\|M_{[0,n)} \pi_n(i_n)\|_\infty\\
    &=&\|\tilde{M}_{[0,n)} \pi_n(i_n)\|_\infty\le\|\pi_n(i_n)\|_\infty.
  \end{eqnarray*}
  Set $L_{[0,m)}=L_{a_0\cdots a_{m-1}}$. We can write for every $k\ge 1$
  \begin{displaymath}
    s_{[0,k)}=L_{[0,m-1)}(p_{m-1})L_{[0,m-2)}(p_{m-2})\cdots L_{[0,1)}(p_1)p_0
  \end{displaymath}
  where each $p_n$ is either empty or equal to $a_n$
  and where $m=m(k)$ is the least index such that $|L_{[0,m)}|>k$
    (this is the one-sided version of Lemma~\ref{lemmaQueffelec}).
  For a word $w$, denote by $\ell(w)\in\Z^A$ the vector 
  $(|w|_a)_{a\in A}$. Then
  \begin{displaymath}
    \ell(s_{[0,k)})=\sum_{n=0}^{m-1}\ell(L_{[0,n)}(p_n))=\sum_{n=0}^{m-1}M_{[0,n)}\ell(p_n)
  \end{displaymath}
  and therefore
  \begin{displaymath}
    \|\pi_0\ell(s_{[0,k]})\|_\infty=\|\sum_{n=0}^{m-1}\pi_0M_{[0,n)}\ell(p_n)\|_\infty\le m.
  \end{displaymath}
  Since $m(k)/\|\ell(s_{[0,k]})\|_\infty\to 0$ when $k\to\infty$,
  the direction of $\ell(s_{[0,k)})$ converges to that of $v$ and
    thus $v$ is determined by $s$.
\end{proofof}

\begin{example}
  Consider again the Tribonacci shift (see Example~\ref{exampleTribonacci3}).
  The unique vector $v\in \cap_{n\ge 0}M_{[0,n)}\R_+^A$
    is the vector $v=1/(1+\lambda^2)\begin{bmatrix}\lambda&1+\lambda&1\end{bmatrix}$ which is an eigenvector corresponding to the maximal
      eigenvalue $\lambda$ of the matrix $M$.
  \end{example}
%%%%%%%%%%%%%%%%%%%%%%
\section{Dimension groups and Rauzy graphs}
\label{chapterDimensionGroupsRauzyGraphs}
\label{sectionGroupsRauzyGraphs}
We now show how to use Rauzy graphs to compute
the dimension group of a minimal shift space.
We begin with considerations valid for all graphs.
We use the fundamental group
$G(\Gamma)$ of a connected graph $\Gamma$
to define the ordered cohomology group of
a connected graph. We then define the notion
of Rauzy graphs associated to a recurrent shift space $X$.
We show that the direct limit sequence of ordered cohomology groups
of the Rauzy graphs is the group $K^0(X,S)$ (Proposition~\ref{propositionDirectLimitGroupsRauzy}).
%%%%%%%%%%%%%%%%%%%%
\subsection{Group of a graph}\label{sectiongroupGraph}
Let $\Gamma=(V,E)$ be a finite oriented graph with $V$ as set of vertices and $E$
as set of edges. We have already met  basic notions concerning graphs such as
paths or cycles, but we recall them now for convenience
(see also Appendix~\ref{appendixGroups} where more details are given).
For an edge $e\in E$, we denote by $s(e)$ its \emph{source}\index{subject}{source!of edge}
\index{subject}{edge!source of}%
(also called its \emph{origin})\index{subject}{origin!of edge}\index{subject}{edge!origin of}
and by $r(e)$ its \emph{range}\index{subject}{range!of edge}\index{subject}{edge!range of}
\index{symbols}{alpha@$\alpha(e)$}\index{symbols}{omega@$\omega(e)$}%
(also called its \emph{end}).\index{subject}{end!of edge}\index{subject}{edge!end of}
There may be several edges
with the same source and range (thus, $G$ is actually a \emph{multigraph}\index{subject}{multigraph}).

Two edges $e,f$ are \emph{consecutive}\index{subject}{consecutive!edges}
\index{subject}{edge!consecutive} if the range of $e$ is the source of $f$.
A \emph{path}\index{subject}{path!in graph} in $\Gamma$ is a sequence
of consecutive edges. A \emph{cycle}\index{subject}{cycle!in graph}
is a path $(e_1,\ldots,e_n)$ such that the source of $e_1$ is the range of $e_n$.
To every path
$p=(e_1,\ldots,e_n)$ in $\Gamma$, we associate its
\emph{composition}\index{subject}{composition!of path}
 $\kappa(p)=e_1+\ldots+e_n$, which is an element
of the free abelian group $\Z(E)$ on the set $E$. 
The \emph{group of cycles}\index{subject}{group!of cycles}
of $\Gamma$, denoted $\Sigma(\Gamma)$\index{symbols}{Sigma@$\Sigma(\Gamma)$}, is the subgroup of $\Z(E)$ 
spanned by the compositions of the cycles
of $\Gamma$.

The elements of $\Z(E)$ can be represented by row vectors indexed by $E$.
Thus, the elements of $\Sigma(\Gamma)$ are also represented by
row vectors indexed by $E$.

We consider, for $v\in V$, the fundamental group $G(\Gamma,v)$
\index{symbols}{G@$G(\Gamma,v)$}%
of $\Gamma$ (see Appendix~\ref{appendixGroups}).
When $\Gamma$ is strongly connected, the group $\Sigma(\Gamma)$
is the abelianization of any of the groups $G(\Gamma,v)$. Indeed,
any cycle $p=p_{s(e)}ep_{r(e)}^{-1}$ can be written
$(p_{s(e)}eq)(p_{r(e)}q)^{-1}$ where $q$ is a path from $r(e)$ to $v$
and then $\kappa(p)=\kappa(p_{s(e)}eq)-\kappa(p_{r(e)}q)$. Moreover,
if $p$ is a cycle around $v'$, let $q$ be a path from 
$v$ to $v'$. Then $r=qpq^{-1}$ is a cycle around $v$ and
$\kappa(r)=\kappa(p)$. 

Let $C(\Gamma)=\Hom(\Sigma(\Gamma),\Z)$\index{symbols}{C@$C(\Gamma)$} and let
$C_+(\Gamma)$ be the submonoid of $C(\Gamma)$ formed
by the morphisms giving a nonnegative value to every
cycle of $\Gamma$.

It will be convenient, given a basis
$B$ of $\Sigma(\Gamma)$, to
represent an element   of $C(\Gamma)$ as a column
vector index by $B$. 

 The coboundary homomorphism $\partial:\Z^V\rightarrow \Z^E$
is defined by
\begin{displaymath}
(\partial\phi)(e)=\phi(r(e))-\phi(s(e))
\end{displaymath}
for every $\phi\in \Z^V$ and $e\in E$. We denote $H(\Gamma)=\Z^E/\partial\Z^V$
and $H_+(\Gamma)=\Z_+^E/\partial\Z^V$.

We identify $\Z^E$ with $\Hom(\Z(E),\Z)$ by duality. Thus, it will
be convenient to consider the elements of $\Z^E$ as column vectors
indexed by $E$.

\begin{example}\label{exampleGraph2}
Consider  the graph $\Gamma$
of Figure~\ref{figureGraph}.
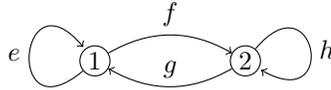
\begin{figure}[hbt]
\tikzset{node/.style={circle,draw,minimum size=0.4cm,inner sep=0pt}}
\tikzstyle{every loop}=[->,shorten >=1pt,looseness=12]
\tikzstyle{loop left}=[in=130,out=220,loop]
\tikzstyle{loop right}=[in=330,out=50,loop]
\centering
\begin{tikzpicture}(20,12)(0,-2)
\node[node](1) at (0,0){$1$};\node[node](2) at(2,0){$2$};

\draw[left](1)edge [loop left]node{$e$}(1);
\draw[above, bend left, ->](1) edge node{$f$}(2);
\draw[above, bend left, ->](2)edge node{$g$}(1);
\draw[right](2) edge [loop right]node{$h$}(2);
\end{tikzpicture}
\caption{A connected graph.}\label{figureGraph}
\end{figure}

 We use the basis $\{e,fg,fhf^{-1}\}$ of $G(\Gamma,1)$
and the corresponding basis $B=\{e,f+g,h\}$ of $\Sigma(\Gamma)$.
Thus $\Sigma(\Gamma)$ is formed of integer row vectors of size $3$
and is isomorphic to $\Z^3$, while $C(\Gamma)$ is formed of integer
column vectors
of the same size.
The matrix of the coboundary map is the  \emph{incidence
  matrix}\index{subject}{incidence!matrix!of graph}\index{subject}{graph!incidence matrix of}
of the graph $\Gamma$, which is the $V\times E$ matrix defined by
\begin{displaymath}
  D_{v,e}=\begin{cases}1&\mbox{ if $v=r(e)$ and $v\ne s(e)$}\\-1&\mbox{ if $v=s(e)$ and $v\ne r(e)$}\\
  0&\mbox{ otherwise.}\end{cases}
\end{displaymath}
In our example, we find
\begin{displaymath}
D=\kbordermatrix{&e&f&g&h\cr 1&0&-1&1&0\cr 2&0&1&-1&0}
\end{displaymath}
The group $\partial\Z^V$ is the group  generated by the rows of the matrix
$D$ and thus it is, in this example, isomorphic to $\Z$. Accordingly, the group $H(\Gamma)$
is isomorphic to $\Z^3$.
\end{example}
%The following is Lemma 5 in ~\cite{Host1995} where the proof uses
%a direct argument instead of the more standard use of a basis
%of the fundamental group of the graph as we do below.

Let $\Gamma=(V,E)$ be a strongly connected graph and let
$\rho:\Z^E\rightarrow C(\Gamma)$ be the morphism
assigning to an element of $\Hom(\Z(E),\Z)$ its restriction to $\Sigma(\Gamma)$. 
Then $\rho$ is a positive morphism such that
$\rho(\Z_+^E)=C_+(\Gamma)$. Given a basis $B$ of $\Sigma(\Gamma)$, the
matrix of the morphism $\rho$ is the matrix having as rows the
elements of $B$ (considered as row vectors indexed by $E$).

The following statement is just the dual of the classical statement
that the sequence
\begin{equation}
  0\to \Z(\Gamma)\edge{\kappa} \Z(E)\edge{\beta} \Z(V)\edge{\gamma}\Z
  \to 0,\label{eqSuiteExacte}
  \end{equation}
where $\kappa$ is the composition map, $\beta(e)=r(e)-s(e)$
and $\gamma(v)=1$ identically, is exact (Exercise~\ref{exerciseExactSequenceGraph}). We give, however, a direct proof for the sake of clarity.

\begin{proposition}\label{propositionGraph}
 For every strongly connected graph $\Gamma$, the sequence
\begin{displaymath}
0\rightarrow\Z\edge{\gamma} \Z^V\edge{\partial} \Z^E\edge{\rho} C(\Gamma)\rightarrow 0
\end{displaymath}
where $\gamma(i)$ is the constant map equal to $i$,
is exact and there is an isomorphism from $C(\Gamma)$ onto
$H(\Gamma)$
sending $C_+(\Gamma)$ onto $H_+(\Gamma)$.
\end{proposition}
\begin{proof}
  The equality $\im(\gamma)=\ker(\partial)$ results from the hypothesis that
  $\Gamma$ is strongly connected.
  
Let us now  show that $\im(\partial)=\ker(\rho)$.
If $\theta$ belongs to  $\Z^V$ and $\psi=\partial\theta$, then
$\psi(\kappa(p))=\theta(v')-\theta(v)$ for every path $p$ from $v$ to $v'$
and thus $\rho(\psi)=0$.
Conversely, suppose that $\rho(\psi)=0$. Let $T$ be a spanning
tree of $\Gamma$ rooted at $v$ and let $p_w$ be the path in $T$ from $v$
to $w$, for every $w\in V$. Let $\theta\in \Z^V$ be defined by
$\theta(w)=\psi(\kappa(p_w))$. 
Then it is easy to verify
that $\psi=\partial\theta$ and thus that $\psi$ belongs to $\im(\partial)$.

We now verify that $\rho$ is surjective. Let $T$
be a spanning tree of $\Gamma$ rooted at $v\in V$ and let
$B'$ be the corresponding basis of $G(\Gamma,v)$
given by \eqref{eqFondamental}.
Then the set $B=\{ \kappa(p) |p\in B'\}$ is a basis
of $C(\Gamma)$.
It is enough to show
that for each $\pi=\kappa(p)\in B$ there is some $\psi\in \Z^E$ such that
$\rho(\psi)$ is the unit vector $u_\pi$ (recall that the elements
of $C(\Gamma)$ are column vectors indexed by $B$).
Let $e\in E\setminus T$ be the
unique edge
such that $p=p_{s(e)}ep_{r(e)}^{-1}$.
Let $\psi\in \Z^E$ be the characteristic function of $e$.
Then clearly $\rho(\psi)=u_\pi$. 
Since $\psi\in \Z_+^E$,
this shows also that $\rho(\Z_+^E)=C_+(\Gamma)$.
\end{proof}
Let $\mathbf{1}_\Gamma$ be the function which associates to each
cycle its length. The triple $(C(\Gamma),C_+(\Gamma),\mathbf{1}_\Gamma)$ is an
 ordered group called the \emph{ordered cohomology
group}\index{subject}{ordered!cohomology group!of a graph}
\index{subject}{graph!ordered cohomology group of} associated to the graph $\Gamma=(V,E)$.
\begin{example}\label{exampleGraph3}
Consider again the graph $\Gamma$
of Figure~\ref{figureGraph}. As in Example~\ref{exampleGraph2}, we use the basis $\{e,fg,fhf^{-1}\}$ of $G(\Gamma,1)$
and the corresponding basis $B=\{e,f+g,h\}$ of $\Sigma(\Gamma)$.
Thus $\Sigma(\Gamma)$ is isomorphic to $\Z^3$.
The corresponding matrix of the morphism $\rho$ is
\begin{displaymath}
P=\begin{bmatrix}1&0&0&0\\0&1&1&0\\0&0&0&1\end{bmatrix}
\end{displaymath}
The rows of $P$ are the coefficients of the elements of $B$ in $\Z(E)$.
The  group $H(\Gamma)$ 
  is the group of linear maps from $\Sigma(\Gamma)$
to $\Z$ and thus it is isomorphic to $\Z^3$.
The submonoid $H_+(\Gamma)$ is formed by the non-negative
linear maps on the space generated by $B$. Thus
it is isomorphic to $\N^3$ and the ordered cohomology group
of $\Gamma$ is $\Z^3$
with the natural ordering. The order unit is the vector
$\begin{bmatrix}1&2&1\end{bmatrix}^t$.
\end{example}
We give below an example of a graph with an ordered cohomology
group which is not isomorphic to $\Z^n$ with the natural ordering.
\begin{example}\label{exampleNotNatural}
Let $\Gamma$ be the graph represented in Figure~\ref{figureGraph3}.
\begin{figure}[hbt]
\tikzset{node/.style={circle,draw,minimum size=0.4cm,inner sep=0pt}}

\centering
\begin{tikzpicture}(20,12)(0,-2)
\node[node](1) at (0,0){$1$};\node[node](2) at(2,0){$2$};

\draw[above, bend left=80, ->](1)edge node{$e$}(2);
\draw[above, bend left=30, ->](1) edge node{$f$}(2);
\draw[above, bend left=30, ->](2)edge node{$g$}(1);
\draw[above, bend left=80, ->](2) edge node{$h$}(1);
\end{tikzpicture}
\caption{The graph $\Gamma$.}\label{figureGraph3}
\end{figure}
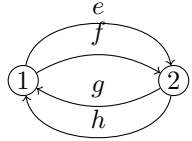

We take this time $\{e+g,f+g,f+h\}$ as basis of $\Sigma(\Gamma)$.
Thus the matrix $P$ is
\begin{displaymath}
P=\begin{bmatrix}1&0&1&0\\0&1&1&0\\0&1&0&1\end{bmatrix}
\end{displaymath}
The ordered cohomology group is again $\Z^3$ but this
time 
$H_+(\Gamma)=\{(\alpha,\beta,\gamma)\in\Z_+^3\mid \alpha+\gamma\ge \beta\}$.
Indeed, let $\theta=\begin{bmatrix}\alpha\\\beta\\\gamma\end{bmatrix}$
be an element of $H(\Gamma)$. Then $\theta\in H_+(\Gamma)$ if and only
if its value on any cycle of $\Gamma$ is nonnegative, that is,
if and only if $\alpha,\beta,\gamma\ge 0$ and
\begin{displaymath}
\theta(e+h)=\theta(e+g)-\theta(f+g)+\theta(f+h)=\alpha-\beta+\gamma
\end{displaymath}
is non-negative, that is, if $\alpha+\gamma\ge\beta$.
The order unit is $\begin{bmatrix}2&2&2\end{bmatrix}^t$. The ordered
group is not a Riesz group because, in $H(\Gamma)$, the sum of the
first two colums of $P$ is equal to the sum of the two last ones,
although none of these vectors can be written as a sum of
positive vectors.
\end{example}
%%%%%%%%%%%%%%%%%%%%%
\subsection{Rauzy graphs}\label{sectionRauzyGraphs}
Let $X$ be a shift space on the alphabet $A$.
Recall from Section~\ref{sectionChapter2RauzyGraphs} that 
 the \emph{Rauzy graph}\index{subject}{Rauzy!graph}\index{subject}{graph!Rauzy}
 of $X$ of order $n$ is the graph $\Gamma_n(X)$
with $\cL_{n-1}(X)$ as set of vertices 
 and $\cL_{n} (X)$ as set of edges. The edge $w$ goes from $u$ to $v$
if $w=ua=bv$ with $a,b\in A$.

\begin{example}
Let $X$ be the Fibonacci shift.
\index{subject}{Fibonacci!shift}\index{subject}{shift space!Fibonacci}%
 The Rauzy graphs of order
$n=1,2,3$ are represented in Figure~\ref{figureRauzy}
(with the edge from $w=ua$ labeled $a$).
\begin{figure}[hbt]
\centering
\tikzset{node/.style={circle,draw,minimum size=0.4cm,inner sep=0.2pt}}
\tikzstyle{every loop}=[->,shorten >=1pt,looseness=12]
\tikzstyle{loop left}=[in=130,out=220,loop]
\tikzstyle{loop right}=[in=330,out=50,loop]
\begin{tikzpicture}

\node[node](11) at (0,0){$\varepsilon$};

\draw[left](11) edge[loop left]node {$a$}(11);
\draw[right](11) edge[loop right]node {$b$}(11);

\node[node](21)at(3,0){$a$};\node[node](22)at(5,0){$b$};

\draw[left](21)edge[loop left]node{$a$}(21);
\draw[above, bend left, ->](21) edge node{$b$}(22);
\draw[below, bend left, ->](22)edge node{$a$}(21);

\node[node](31)at (6,0){$aa$};
\node[node](32)at (8,1){$ab$};
\node[node](33)at(8,-1){$ba$};

\draw[above, bend left, ->](31)edge node{$b$}(32);
\draw[right,bend left, ->](32)edge node{$a$}(33);
\draw[below, bend left, ->](33)edge node{$a$}(31);
\draw[left, bend left, ->](33)edge node{$b$}(32);

\end{tikzpicture}
\caption{The Rauzy graphs of order $n=1,2,3$ of the Fibonacci shift.}
\label{figureRauzy}
\end{figure}

\end{example}
There is a positive morphism from the group $C(\Gamma_n(X))$
to the group $C(\Gamma_{n+1}(X))$. Indeed, the prefix map
defines a projection from the graph $\Gamma_{n+1}(X)$ onto
the graph $\Gamma_n(X)$. The set of edges of the first one is mapped
onto the set of edges of the second one and the set of cycles
onto the set of cycles. It follows that this projection
defines a positive morphism from $\Sigma(\Gamma_{n+1}(X))$
to $\Sigma(\Gamma_n(X))$ and by duality from
$C(\Gamma_n(X))$ to $C(\Gamma_{n+1}(X))$.
\begin{proposition}\label{propositionDirectLimitGroupsRauzy}
For any recurrent shift space $X$,
the unital ordered group  $K^0(X,S)$ is the direct limit 
of the ordered groups associated with its Rauzy graphs.
\end{proposition}
\begin{proof}
Since the set of vertices of $\Gamma_n(X)$ is $\cL_{n-1}(X)$ and its set of edges
is $\cL_n(X)$, we can identify $H(\Gamma_n(X))$ to $G_n(X)$ (defined by \eqref{def:gnX}) and
$H_+(\Gamma_n(X))$ to $G_n^+(X)$.

The unital ordered group $(C(\Gamma_n(X)),C_+(\Gamma_n(X)),\1_{\Gamma_n(X)})$
 can be identified, by Proposition~\ref{propositionGraph},
with $(G_n(X),G_n^+(X),\1_n(X))$. Since
the morphism from $C(\Gamma_n(X))$ to $C(\Gamma_{n+1}(X))$
induced by taking the prefixes is the same as
the morphism $i_{n+1,n}$ from $G_n(X)$ to $G_{n+1}(X)$, the result follows from
Proposition~\ref{propositionRecurrentSubshift}.
\end{proof}
%\marginpar{DP: developper calcul pour Fibonacci}
\begin{example}
  We consider again the Fibonacci shift $X$. We already know
  that its dimension group is $\Z[\alpha]$ with $\alpha=(1+\sqrt{5})/2$
  (see Examples~\ref{exampleSturmian2}) and \ref{exampleFibonacci4}).
  We will not prove it again but
    our point  is to put in evidence the isomorphisms used
    in the proof of Proposition~\ref{propositionDirectLimitGroupsRauzy}.
    
  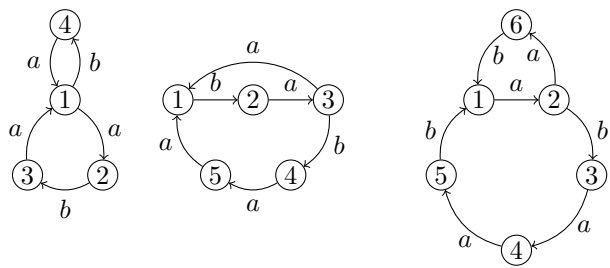
\begin{figure}[hbt]
    \centering
\tikzset{node/.style={circle,draw,minimum size=0.4cm,inner sep=0.2pt}}
\tikzstyle{every loop}=[->,shorten >=1pt,looseness=12]
\tikzstyle{loop left}=[in=130,out=220,loop]
\tikzstyle{loop right}=[in=330,out=50,loop]
\begin{tikzpicture}

%Gamma_4

  \node[node](1)at(.5,1){$1$};
  \node[node](2)at(1,0){$2$};
  \node[node](3)at(0,0){$3$};
\node[node](4)at(.5,2){$4$};

\draw[bend left,right,->](1)edge node{$a$}(2);
\draw[bend left,below,->](2)edge node{$b$}(3);
\draw[bend left,left,->](3)edge node{$a$}(1);
\draw[bend right,right,->](1)edge node{$b$}(4);
\draw[bend right,left,->](4)edge node{$a$}(1);
  %Gamma_5
\node[node](1)at(2,1){$1$};\node[node](2)at(3,1){$2$};
\node[node](3)at(4,1){$3$};
\node[node](4)at(3.5,0){$4$};
\node[node](5)at(2.5,0){$5$};

\draw[above,->](1)edge node{$b$}(2);
\draw[above,->](2)edge node{$a$}(3);
\draw[bend right=45,right,above,->](3)edge node{$a$}(1);
\draw[bend left,right,->](3)edge node{$b$}(4);
\draw[bend left,below,->](4)edge node{$a$}(5);
\draw[bend left,left,->](5)edge node{$a$}(1);

\node[node](1)at(6,1){$1$};
\node[node](2)at(7,1){$2$};
\node[node](3)at(7.5,0){$3$};
\node[node](4)at(6.5,-1){$4$};
\node[node](5)at(5.5,0){$5$};
\node[node](6)at(6.5,2){$6$};

\draw[->,above](1)edge node{$a$}(2);
\draw[bend left,right,->](2)edge node{$b$}(3);
\draw[bend left,right,->](3)edge node{$a$}(4);
\draw[bend left,below,->](4)edge node{$a$}(5);
\draw[bend left,left,->](5)edge node{$b$}(1);
\draw[bend right,left,->](2)edge node{$a$}(6);
\draw[bend right,right,->](6)edge node{$b$}(1);

\end{tikzpicture}
\caption{The Rauzy graphs of order $4,5,6$ of the Fibonacci shift.}
\label{figureRauzy2}
  \end{figure}
  The prefix $p_n$ of length $n$ of the Fibonacci word $x$ is the vertex
  labeled $1$ in $\Gamma_{n+1}(X)$ in Figures \ref{figureRauzy}
  and \ref{figureRauzy2}. Thus, the elementary cycles around $1$
  in $\Gamma_{n+1}(X)$ are labeled by the right return words to $p_n$
  (conjugate to the left return words forming the set denoted
  $W_n(x)$ in Example~\ref{exampleSturmian2}). For example,
  the vertex $1$ of $\Gamma_5(X)$ is the word $p_4=abaa$. The
  elementary cycles around $1$ are labeled by $baa$ and $babaa$,
  which are conjugate to the elements of $W_4(x)=\{aba,abaab\}$
  by the conjugacy $u\to (abaa)u(abaa)^{-1}$.

  The prefix
  map from $\Gamma_4(X)$ to $\Gamma_3(X)$ sends the vertices
  $2$ and $4$ to the vertex $ba$. Since $W_2(x)=W_3(x)$,
  the morphism from $\Sigma(\Gamma_3(X))$ to $\Sigma(\Gamma_4(X))$
  induced by the prefix map is the identity. In contrast,
  the morphism from $\Sigma(\Gamma_4(X))$ to $\sigma(\Gamma_5(X))$ is
  given by the matrix
  \begin{displaymath}
\kbordermatrix{ &ab&aba\\aba&0&1\\abaab&1&1}
  \end{displaymath}
  where the rows are indexed by the words of $W_4(X)$,
  in bijection with the (composition of the) cycles
  around $1$ in $\Gamma_5(X)$. The columns are indexed
  by the words of $W_3(X)$. Finally, the morphism from
  $C(\Gamma_4(X))$ to $C(\Gamma_5(X))$ is obtained by transposition
  of this matrix.
\end{example}
We will need several times (in the next section, and later in Chapter~\ref{chapterDendricShifts}) the following result relating return words
and Rauzy graphs.
\begin{proposition}\label{propositionJulien}
Let $X$ be a minimal shift space and let $u\in\cL(X)$.
There exists an $n\ge 1$ with the following property.
Let $x\in\cL_n(X)$ be a word ending with $u$ and let $S$
be the set of labels of paths from $x$ to itself in $\Gamma_{n+1}(X)$. Then
$S$ is contained in $\RR_X(u)^*$.
\end{proposition}
\begin{proof}
Let $n$ be the maximal length of the words in $u\RR_X(u)$.
Consider  $y\in S$. Since $y$ is the label of a path from
$x$ to $x$ in $\Gamma_{n+1}(X)$, the word $xy$ ends with $x$.
Thus there is a unique factorization $y=y_1y_2\cdots y_k$
in nonempty words $y_i$
where for each $i$ with $1\le i\le k$, the word $uy_i$ ends
with $u$ and has no other occurrence of $u$ except
as a prefix or as a suffix. But, by the choice of $n$,
the prefix of length $n$
of $uy_i$ has a factor $u$ other than as a prefix
and thus $|uy_i|\le n$. Now since
$uy_i$ is the label of a path of length at most $n$ in
$\Gamma_n(X)$, it is in $\cL(X)$. This implies that $y_i$
is in $\RR_X(u)$ and proves the claim.
\end{proof}

%%%%%%%%%%%%%%%%%%%%%%%%%%
\section{Dimension groups of substitution shifts}\label{chapterFixedPointMorphism}\label{sectionDGSubstitutionShifts}
We show in this section how to compute the dimension
group of a substitutive shift.

We begin by
establishing a connection between the two block presentation
of the shift space and its Rauzy graph of order two.

\subsection{Two-block presentation and Rauzy graph}
\label{sctionTwoBlockExtension}
Let $\varphi:A\rightarrow A^*$ be a substitution and let $X$ be the
shift space associated to $\varphi$.

We have seen that there is a matrix
$M$ associated to $\varphi$, called its composition matrix. Unfortunately,
 the composition (or incidence) matrix of the substitution
does not determine the ordered group of the substitution shift
associated to $\varphi$, as shown by the following example.
We show below that it is however determined by the incidence
matrix $M_2$ of the $2$-block presentation $\varphi_2$ of $\varphi$
(see Section~\ref{sectionSubstitutionSystems}
for the definition of the matrix $M_2$).
%A different and more complete approach has been developped by Forrest, see \cite{Forrest:97}.
%It applies to a wider family of shifts. 

\begin{example}\label{exampleMorse4}
Let $X$ be the Thue-Morse shift. The matrix $M$ associated to the
Thue-Morse morphism is
\begin{displaymath}
M=\begin{bmatrix}1&1\\1&1\end{bmatrix}
\end{displaymath}
and it is the same as for the substitution $a\mapsto ab,b\mapsto ab$.
The shift space corresponding to the second substitution
has two elements and its dimension group is $\Z$.
The dimension group of the Thue-Morse shift is not isomorphic to
$\Z$. Indeed (see Example~\ref{exampleMeasureMorse}), 
one has $\mu([aa])\ne\mu([ab])$ for the unique
invariant measure on $X$ and thus, 
by Proposition~\ref{propositionCoboundaryInt},
the difference $\charac_{[aa]}-\charac_{[ab]}$ of the characteristic
functions of the cylinders $[ab]$ and $[ba]$ is not a coboundary
(we will see shortly that actually the group $H(X,S,\Z)$ is isomorphic to $\Z [1/2] \times \Z$).
Thus the dimension groups are not the same for the two shift spaces.
\end{example}
We  prove a statement connecting the endomorphism of $R_2(X)$
%\marginpar{FD: pas defini} 
defined by the matrix $M_2$ with the fundamental group $\Sigma(\Gamma_2(X))$
of the Rauzy graph $\Gamma_2(X)$. Recall that $R_n(X)$ denotes the group of maps from $\cL_n(X)$
to $\R$) and that $M_2(X)$ operates on the left on the elements
of $R_2(X)$ considered as column vectors.

Recall  from Section~\ref{chapterOrderedGroupRecurrent} 
that the map $\partial_1:Z_1(X)\rightarrow Z_2(X)$ is the morphism
defined by $\partial_1\phi(ab)=\phi(b)-\phi(a)$. Since $\cL_1(X)$ (resp.
$\cL_2(X)$) is the set of vertices (resp. of edges) of the graph $\Gamma_2(X)$,
it coincides with the map $\partial$ of Section~\ref{sectiongroupGraph}
and we shall denote it simply $\partial$.

Let $\tau:A\rightarrow A$ be the map which sends each letter $a\in A$ to the first
letter of $\varphi(a)$. We define an endomorphism $I$ of $R_1(X)$ by
\begin{displaymath}
(I\phi)(a)=\phi(\tau(a))
\end{displaymath}
for every $\phi\in R_1(X)$ and $a\in A$.
%\marginpar{FD : $R_1$, $R_2$ ?}

Let $P$ be a matrix whose rows form a basis of the group
$\Sigma(\Gamma_2(X))$. The matrix $P$ is the matrix of the morphism $\rho$ in Proposition~\ref{propositionGraph}.

We may choose $P$ as a nonnegative matrix by choosing a basis
formed  of compositions of cycles of the graph.

Moreover these cycles may be chosen to correspond to return words
to a vertex $a\in A$, where $a$ is such that $\varphi(a)$
begins with $a$ (we may always find such $a$, up to replacing
$\varphi$ by some power). Indeed, by Proposition~\ref{propositionJulien}
applied with $u=a$,
there is an integer $n$ such that the set labels of paths
from $x$ to itself in $\Gamma_{n+1}(X)$, where
$x\in\cL_n(X)$ ends with $a$ is generated by $\RR_X(a)$.
This implies, by projection from $\Gamma_{n+1}(X)$
to $\Gamma_2(X)$, that  every cycle from $a$ to $a$ in $\Gamma_2(X)$,
is a product of return words to $a$ or their inverses.
Thus the compositions of return words to $a$
generate the
group $\Sigma(\Gamma_2(X))$.

We illustrate this on an example.

\begin{example}\label{exampleJulien2blocks}
Let $A=\{a,b,c,d\}$ and let $\varphi:a\to ab,b\to cda,c\to cd,d\to abc$
(we shall consider again this
substitution below and later in Example~\ref{exampleJulien}). We have
$\cL_2(X)=\{ab,ac,bc,ca,cd,da\}$ and the graph $\Gamma_2(X)$
is represented in Figure~\ref{figureGrapheGamma_2Julien}.
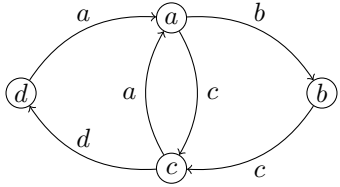
\begin{figure}[hbt]
\centering
\tikzset{node/.style={circle,draw,minimum size=0.4cm,inner sep=0.2pt}}
\begin{tikzpicture}
\node[node](a)at(2,2){$a$};
\node[node](b)at(4,1){$b$};
\node[node](c)at(2,0){$c$};
\node[node](d)at(0,1){$d$};

\draw[above,->,bend left](a)edge node{$b$}(b);
\draw[right,->,bend left](a)edge node{$c$}(c);
\draw[below,->,bend left](b)edge node{$c$}(c);
\draw[left,->,bend left](c)edge node{$a$}(a);
\draw[above,->,bend left](c)edge node{$d$}(d);
\draw[above,->,bend left](d)edge node{$a$}(a);
\end{tikzpicture}
\caption{The graph $\Gamma_2(X)$.}\label{figureGrapheGamma_2Julien}
\end{figure}

We have $\RR_X(a)=\{bca,bcda,cda\}$.
The choice for the matrix $P$ is then
\begin{displaymath}
  P=\kbordermatrix{
      &ab&ac&bc&ca&cd&da\\
   bca&1 &0 &1 &1 &0 &0\\
   bcda&1 &0 &1 &0 &1 &1\\
   cda &0 &1 &0 &0 &1 &1
  }
\end{displaymath}
For example, the first row corresponds to the return word $bca$
since it gives the sequence $(ab)(bc)(ca)$ of $2$-blocks.
\end{example}

\begin{proposition}\label{lemma0Substitutions}
We have $M_2\circ\partial=\partial\circ I$, so that $\partial(R_1(X))$ 
and $\partial(Z_1(X))$ are invariant by $M_2$. Moreover,
assuming that the rows of $P$ correspond to return words to a letter
$a\in A$ such that $\varphi(a)$ begins with $a$,
there is a unique nonnegative matrix $N_2$ such that
\begin{displaymath}
PM_2=N_2P.
\end{displaymath}
\end{proposition}
\begin{proof}
We first prove that $M_2\circ\partial=\partial\circ I$.
Let $\phi\in R_1(X)$. Then
\begin{eqnarray*}
(\partial\circ I\circ\phi)(ab)&=&(\partial\circ\phi\circ\tau)(ab)\\
&=&\phi\circ\tau(b)-\phi\circ\tau(a)=\phi(\tau(b))-\phi(\tau(a)).
\end{eqnarray*}
To evaluate $(M_2\circ\partial\phi)(ab)$, set $\varphi(a)=a_1a_2\cdots a_k$
and $\varphi(b)=b_1b_2\cdots b_\ell$. Then, identifying $A_2$ and $\cL_2(X)$,
we have
\begin{equation}
  \varphi_2(ab)=(a_1a_2)(a_2a_3)\cdots(a_kb_1)\label{eqvarphi_2}
\end{equation}
where $\varphi_2$ is the $2$-block presentation of $\varphi$.
Next, considering accordingly $M_2$ as an endomorphism of $R_2(X)$,  we have
for any $\psi\in R_2(X)$, using Equation~\eqref{eqvarphi_2},
\begin{displaymath}
(M_2\psi)(ab)=\psi(a_1a_2)+\psi(a_2a_3)+\ldots+\psi(a_kb_1).
\end{displaymath}
 We obtain
\begin{eqnarray*}
(M_2\partial\phi)(ab)&=&\phi(a_2)-\phi(a_1)+\phi(a_3)-\phi(a_2)+\ldots+\phi(b_1)-\phi(a_k)\\
&=&\phi(b_1)-\phi(a_1)=\phi(\tau(b))-\phi(\tau(a))
\end{eqnarray*}
and thus the conclusion.

Since $M_2\circ \partial=\partial\circ I$, the subgroups $\partial(R_1(X)$
and $\partial(Z_1(X))$ are invariant by $M_2$. 

Let $p=(a_1a_2)(a_2a_3)\cdots (a_{n-1}a_n)(a_na_1)$ be a cycle in $\Gamma_2(X)$
which is a row of $P$ and thus with $a_1=a$.
We can consider $p$ as a word on the alphabet $\cL_2(X)$ and compute its image
by the morphism $\varphi_2$. Since the rows
of $P$ correspond to return words to $a$,
the first letter of $\varphi(a_1)$ is $a$.
Then the first element of $\varphi_2(a_1a_2)$ is equal to $ab$ for some
$b\in A$, and the last element of $\varphi_2(a_na_1)$ is equal to
$ca$ for some $c\in A$ (as we have seen above computing $\varphi_2(ab)$).
It follows that $\varphi_2(p)$ is a cycle around $a$ in $\Gamma_2(X)$
and thus a composition of return words to $a$.
By the choice of $P$, it is a nonnegative combination of rows
of $P$.

This implies that $PM_2=N_2P$ where $N_2$ is the matrix of the
map sending the composition of a cycle $p$ on the composition
of $\varphi_2(p)$.
\end{proof}
Note that the two assertions of Proposition~\ref{lemma0Substitutions}
are related. Indeed, the matrix of the map $\partial$
is the transpose of the incidence matrix $D$ of the graph $\Gamma_2(X)$.
The rows of $D$ generate the orthogonal of the space generated by the
rows of $P$ (because \eqref{eqSuiteExacte} is an exact sequence). It is thus equivalent to say the
space generated by the rows of $P$ is invariant by $M_2$
and that the space $\partial(R_1(X))$ is invariant by $M_2$.
This duality is described in more detail in the following example.
\begin{example}
  Consider again the substitution shift of Example~\ref{exampleJulien2blocks}.
  Using the bijection $\{ab,ac,bc,ca,cd,da\}\to\{x,y,z,t,u,v\}$,
  we find $\varphi_2:x\to xz,y\to xz,z\to uvy, t\to uv, u\to uvx, v\to xzt$.
The matrices $M,M_2$ are
\begin{displaymath}
M=\begin{bmatrix}1&1&0&0\\1&0&1&1\\0&0&1&1\\1&1&1&0\end{bmatrix},\quad
M_2=\begin{bmatrix}1&0&1&0&0&0\\1&0&1&0&0&0\\0&1&0&0&1&1\\0&0&0&0&1&1\\
0&0&0&0&1&1\\1&0&1&1&0&0\end{bmatrix}
\end{displaymath}
The matrices $P$ and $N_2$ are then
\begin{displaymath}
P=\begin{bmatrix}1&0&1&1&0&0\\1&0&1&0&1&1\\0&1&0&0&1&1\end{bmatrix},\quad
N_2=\begin{bmatrix}0&1&1\\1&1&1\\1&1&0\end{bmatrix}
\end{displaymath}
so that $N_2$ has dimension less than $M$. The matrix $D$ is
\begin{displaymath}
  D=\kbordermatrix{ &ab&ac&bc&ca&cd&da\\a&-1&-1&0&1&0&1\\b&1&0&-1&0&0&0\\
    c&0&1&1&-1&-1&0\\d&0&0&0&0&1&-1}.
\end{displaymath}
One may easily verify that the rows of $D$
generate a vector space of dimension $3$ orthogonal to the
space generated by the rows
of $P$.
\end{example}
\begin{example}\label{exampleFibonacci5}

Let $\varphi:a\mapsto ab,b\mapsto a$ be the Fibonacci morphism and let $X$
be the Fibonacci shift.
\index{subject}{Fibonacci!shift}\index{subject}{shift space!Fibonacci}%
 Let $f:x\mapsto aa,y\mapsto ab,z\mapsto ba$
be a bijection from $A_2=\{x,y,z\}$
onto $L_2(X)=\{aa,ab,ba\}$. The morphism $\varphi_2$
is $x\mapsto yz,y\mapsto yz,z\mapsto x$.
The associated matrices $M$ and $M_2$ and the vector
$v$ spanning the group $\partial_1(R_1(X))$
are (see Example~\ref{exampleFibonacci4})
\begin{displaymath}
M=\begin{bmatrix}1&1\\1&0\end{bmatrix},\qquad
M_2=\begin{bmatrix}0&1&1\\0&1&1\\1&0&0\end{bmatrix},\qquad
v=\begin{bmatrix}0\\1\\-1\end{bmatrix}.
\end{displaymath}
The Rauzy graph $\Gamma_2(X)$ is represented in Figure~\ref{figureGamma_2Fibo}
with the edges labeled by the corresponding element of $A_2$.
\begin{figure}[hbt]
\centering
\tikzset{node/.style={circle,draw,minimum size=0.4cm,inner sep=0.2pt}}
\tikzstyle{every loop}=[->,shorten >=1pt,looseness=12]
\tikzstyle{loop left}=[in=130,out=220,loop]
\tikzstyle{loop right}=[in=330,out=50,loop]
\begin{tikzpicture}
\node[node](21)at(3,0){$a$};\node[node](22)at(5,0){$b$};

\draw[left](21)edge[loop left]node{$x$}(21);
\draw[above, bend left, ->](21) edge node{$y$}(22);
\draw[below, bend left, ->](22)edge node{$z$}(21);

\end{tikzpicture}

\caption{The Rauzy graph $\Gamma_2(X)$.}\label{figureGamma_2Fibo}
\end{figure}
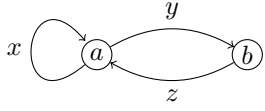

The matrices $P$ corresponding to the basis $\{x,yz\}$ and the
matrix $N_2$ associated with the action of $\varphi_2$ on the rows of $P$ are
\begin{displaymath}
P=\begin{bmatrix}0&1&1\\1&0&0\end{bmatrix},\qquad
N_2=\begin{bmatrix}1&1\\1&0\end{bmatrix}
\end{displaymath}
Thus we conclude that $N_2=M$ in this case.
\end{example}

%%%%%%%%%%%%%%%%%%%%%%%%%%%%%%
\subsection{Dimension group of a substitution shift}
\label{sectionDimensionGroupSubstitutive}
%We define the following unital ordered group associated 
%with  the matrix $M_2$. 
Recall  first from Section~\ref{sectionDirectLimits} that,
for a $d\times d$ matrix $M$, $\RR_{M}$ denotes the eventual range of $M$
and $\KK_{M}$ its eventual kernel. 
%Since $M_2$ is an $L_2(X)\times L_2(X)$-matrix, the elements 
%of $\RR_{M_2}$ and those of $\KK_{M_2}$
%can be considered as functions from $L_2(X)$ to $\R$, that is as elements
%of $R_2(X)$.

For $M$ nonnegative, we have defined a unital ordered group
$(\Delta_{M},\Delta^+_{M},\1_{M})$. By Equation~\eqref{eqDeltaM},
the group $\Delta_M$ is given by
\begin{displaymath}
\Delta_{M}=\{v\in\RR_{M}\mid \mbox{ for some } k\ge 1, M^kv\in \Z^d\},
\end{displaymath} 
%since it is the group formed of all $v\in\RR_{M_2}$ such that
%for some $k\ge 1$, the vector $M_2^kv$ has components in $\Z$,
%that is, belongs to $Z_2(X)$. 
with positive cone 
\begin{displaymath}
\Delta^+_{M}=\{v\in\RR_{M}\mid \mbox{ for some } k\ge 1, M^kv\in \Z_+^d\}
\end{displaymath}
and  order unit $\1_{M}$ equal to the projection on $\RR_M$ 
along $\KK_M$ of the vector
with all its components equal to $1$.

The following statement gives a surprisingly simple way to compute
the dimension group of a substitution shift.
Indeed, it shows that the dimension
group of a primitive substitution shift is
defined by the matrix $N_2$ of Proposition~\ref{lemma0Substitutions} and  the nonnegative matrix $P$ such that $PM_2=N_2P$. 
\begin{proposition}\label{propositionDimensionGroupSubstitutions}
Assume that $\varphi$ is a primitive substitution and that
the shift $X$ associated to $\varphi$ is infinite.
The dimension group $K^0(X,S)$ is isomorphic to
$(\Delta_{N_2}, \Delta_{N_2}^+,P 1_{M_2})$.
\end{proposition}
Note that the order unit is not the vector $\1_{N_2}$ but the projection
on $G_2(X)$ of the vector $1_{M_2}$ (see Example~\ref{exampleMorse5}).
The proof relies on several lemmas.

We first introduce a sequence of partitions in towers associated with
 the primitive morphism $\varphi$.
\begin{lemma}\label{lemma1Substitutions}
The sequence $(\Pg(n))_{n\ge 0}$ with
\begin{equation}
\Pg(n)=\{S^j\varphi^n([ab])\mid ab\in \cL_2(X), 0\le j<|\varphi^n(a)|\}.\label{sequencePartitionsSubst}
\end{equation}
is a nested sequence of partition in towers of $X$.
\end{lemma}
\begin{proof}
  Let $f:\cL_2(X)\to A_2$ be a bijection and let $X^{(2)}$
  be the second higher block presentation  of $X$.
We consider the sequence $(\Qg(n))$ of partitions in towers of the
second higher block presentation $X^{(2)}$ of $X$ associated with the
$2$-block presentation $\varphi_2^n$ of $\varphi^n$ as in 
Proposition~\ref{propositionPartitionSubstitution}. Thus
\begin{displaymath}
\Qg(n)=\{S^j\varphi_2^n([u])\mid u\in A_2, 0\le j<|\varphi_2^n(u)|\}.
\end{displaymath}
Let
$\pi:A_2\rightarrow A$ be the morphism assigning to $u\in A_2$
the first letter of $f(u)$. The extension of $\pi$ to $A_2^\Z$
 defines an isomorphism from $X^{(2)}$ onto $X$. Let $(\Pg(n))$ be the image
of the sequence $(\Qg(n))$ by the isomorphism $\pi$. 
Since $\pi([u])=[ab]$ when $f(u)=ab$, the partition
$\Pg(n)$ is given by Equation~\eqref{sequencePartitionsSubst}.
By Proposition~\ref{propositionPartitionSubstitutionNested},
the sequence $\Pg(n)$ is nested.
\end{proof}
Denote $G(n)=G(\Pg(n))$, $G^+(n)=G^+(\Pg(n))$ and
$1_n=1_{\Pg(n)}$. Let $(G,G^+,1)$ be the inductive limit of the sequence
$(G(n),G^+(n),1_n))$ with the morphisms $I(n+1,n)=I(\Pg(n+1),\Pg(n))$.
\begin{lemma}\label{lemma2Substitutions}
The map from $Z_2(X)$ to $C(X,\Z)$ sending $\phi\in Z_2(X)$ to
the map equal to $\phi(ab)$ on the cylinder $[ab]$ defines
an isomorphism of unital ordered groups from $(\Delta_{M_2},\Delta^+_{M_2},1_{M_2})$
onto $(G,G^+,1)$.
\end{lemma}
\begin{proof}
By associating to each $\phi\in Z_2(X)$ the function equal to $\phi(ab)$ on 
$\varphi^n([ab])$, we can identify $G(n)$ to $Z_2(X)$,
$G^+(n)$ to $Z_2^+(X)$ and $1_n$ to  the map
$ab\mapsto |\varphi^n(a)|$. Given $k$ with $0\le k<|\varphi^n(a)|$,
we have (see Figure~\ref{figurephin+1})
\begin{displaymath}
S^k(\varphi^{n+1}([ab])\subset\varphi^n([cd])\Leftrightarrow\left\{
\begin{array}{l}
\mbox{ there exists $\ell$ with $0\le \ell<|\varphi(a)|$ such that}\\
\mbox{$S^\ell([ab])\subset[cd]$ and $k=|\varphi^{n}(\varphi([ab])_{[0,\ell-1]})|$}.
\end{array}
\right.
\end{displaymath}

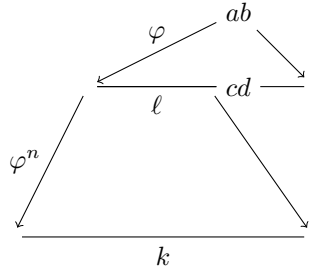
\begin{figure}[hbt]
\centering
\begin{tikzpicture}
\node(1)at (3,3){$ab$};
\node(2)at(1,2){};
\node(3bis)at(2.6,2){};\node(3)at(3,2){$cd$};\node(3ter)at(4,2){};
\node(4)at(0,0){};\node(5)at(4,0){};

\draw[below](2) edge node{$\ell$}(3);
\draw(3) edge node{}(2);
\draw[left, above,->](1)edge node{$\varphi$} (2);
\draw[->](1) edge node{}(3ter);
\draw(3) edge node{}(3ter);
\draw[left,->](2) edge node{$\varphi^n$}(4);
\draw[->](3bis) edge node{}(5);
\draw[below](4) edge node{$k$}(5);
\end{tikzpicture}
\caption{Representing $S^k(\varphi^{n+1}([ab])\subset\varphi^n([cd])$}
\label{figurephin+1}
\end{figure}

Therefore, for every $u,v\in A_2$ with $f(u)=ab$ and $f(v)=cd$, the number
$$
\Card\{k\mid0\le k<|\varphi^n(a)|,S^k\varphi^{n+1}([ab])\subset\varphi^n([cd])\}
$$
is the number of occurrences of $v$ in $\varphi_2(u)=(M_2)_{uv}$.
Consequently, we may identify the morphism $I(n+1,n))$ to the
morphism $M_2:Z_2(X)\rightarrow Z_2(X)$.

Thus, the inductive limit $(G,G^+,1)$ associated to the sequence
$(\Pg(n))$ of partitions in towers can be identified to 
$(\Delta_{M_2},\Delta^+_{M_2},1_{M_2})$. 
\end{proof}
Recall from Proposition~\ref{pi(P)} that
there is a morphism $\pi(n):G(n)\rightarrow H(X,S,\Z)$
making the diagram of Figure \ref{figureDiagrampi(Pn)} commutative
and from Lemma~\ref{lemma1} that there is a unique morphism $\sigma:(G,G^+,1)
\rightarrow K^0(X,S)$ such that $\sigma\circ I(n)=\pi(n)$ for every $n\ge 0$.

Recall also from Section~\ref{sectionGroupCylinder} that $\partial_1:R_1(X)\rightarrow R_2(X)$
is the morphism defined by $(\partial_1\phi)(ab)=\phi(b)-\phi(a)$.
\begin{lemma}\label{lemma3Substitutions}
The  morphism $\sigma:G\rightarrow H(X,S,\Z)$ is onto and its kernel
is $\Delta_{M_2}\cap \partial_1(R_1(X))$.
\end{lemma}
\begin{proof}
\noindent We first show that the morphism $\sigma$ is surjective. Since every function
in $C(X,\Z)$ is cohomologous to a cylinder function by Proposition~\ref{lemma4},
it is enough to consider  a cylinder function 
$\u(\phi)$ associated to $\phi\in Z_k(X)$. Choose $n$ so large that
$|\varphi^n(a)|>k$ for every letter $a\in A$. Since all elements of
the atoms of the partition $\Pg_n$ have the same prefix of length $k$,
the cylinder function $\u(\phi)$ is constant on every element
of the partition $\Pg_n$ and thus belongs to $C(n)$. Its
image by $I(n)$ is sent by $\pi(n)$ to the class of $\u(\phi)$
modulo the coboundaries (see Figure~\ref{figureDiagrampi(Pn)}).
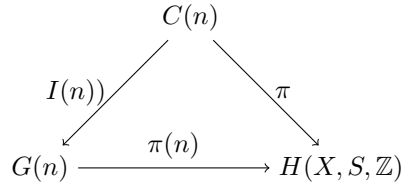
\begin{figure}[hbt]
\centering
\begin{tikzpicture}
\node(C)at(2,2){$C(n)$};
\node(G)at(0,0){$G(n)$};\node(H)at(4,0){$H(X,S,\Z)$};

\draw[left,->](C) edge node{$I(n)$)}(G);
\draw[right,->](C) edge node{$\pi$}(H);
\draw[above,->](G)edge node{$\pi(n)$}(H);
\end{tikzpicture}
\caption{The morphism $\pi(n)$.}\label{figureDiagrampi(Pn)}
\end{figure}

\noindent We now prove that the kernel of $\sigma$ is $\Delta_{M_2}\cap\partial_1(Z_1(X))$.
Let $\phi$ be in $G$, or equivalently in $\Delta_{M_2}$. In particular,
$\phi $ belongs to $R_2(X)$ and there exists $k\ge 0$ such that $M_2^k\phi\in Z_2(X)$.
Assume that $\sigma(\phi)=0$.
By definition this means that $\pi(\tilde{\phi})=0$ where
$\pi:C(n)\rightarrow H(X,S,\Z)$ is the natural projection
and $\tilde{\phi}$ is defined by
\begin{displaymath}
\tilde{\phi}(x)=\begin{cases}
(M_2^k\phi)(ab)&\mbox{if $x\in\varphi^k([ab])$ for some $ab\in \cL_2(X)$}\\
0&\mbox{otherwise}.
\end{cases}
\end{displaymath}
Note that, as in the proof of Proposition~\ref{lemma0Substitutions},
we identify $A_2$ with $\cL_2(X)$ and consequently consider
$M_2$ as an endomorphism of $R_2(X)$.

Let $g\in C(X,\Z)$ be such that $\tilde{\phi}=g\circ S-g$.
We claim that there exists $n\ge 0$ such that $g$ is constant
on the set $\varphi^n([ab])$ for every $ab\in \cL_2(X)$.
 Since $g$ is continuous, it
is locally constant and there
is an $m\ge 1$ be such that $g(x)$ depends only on $x_{[-m,m]}$. 
Choose $n>k$ so large
that $|\varphi^n(a)|>m$ for every letter $a$. 
Let $ab\in \cL_2(X)$ and $y,z\in \varphi^n([ab])$. 
Since $g$ depends only on $x_{[-m,m]}$, $g(S^mx)$
depends only on $x_{[0,2m]}$. Since $y,z\in\varphi^n([ab])$ and
$|\varphi^n(a)|,|\varphi^n(b)|>m$, $y$ and $z$ share the same
$2m$ first coordinates and therefore $g(S^mx)=g(S^my)$. On the other hand,
for all $0\le j<|\varphi^n(a)|$,
$S^jy$ and $S^jz$ are in the same atom of the partition $\Pg(n)$.
Since $m<|\varphi^n(a)|$ and since $\tilde{\phi}$ is constant on
the atoms of $\Pg(n)$, we obtain $\tilde{\phi}^{(m)}(y)=\tilde{\phi}^{(m)}(z)$.
Since finally 
$g=g\circ S^m-\tilde{\phi}^{(m)}$ by Equation~\eqref{eq1}, we conclude
that $g(y)=g(z)$.

 Let $\psi\in Z_2(X)$ be such that $\psi(ab)=g(x)$ for $x\in\varphi^n([ab])$. 
Then if
$x\in \varphi^n([ab])$ and $S^{|\varphi^n(a)|}x\in\varphi^n([bc])$, we have
with $\ell=|\varphi^n(a)|$,
\begin{eqnarray*}
\psi(bc)-\psi(ab)&=&g(S^{\ell}x)-g(x)\\
&=&\tilde{\phi}^{(\ell)}(x).
\end{eqnarray*}
Recall that
\begin{displaymath}
\tilde{\phi}^{(\ell)}(x)=\tilde{\phi}(x)+\tilde{\phi}\circ S(x)+\ldots
\tilde{\phi}(S^{\ell-1}(x))
\end{displaymath}
and note that the term $\tilde{\phi}\circ S^j(x)$ of
 this sum is equal to $(M^k\phi)(cd)$
if there is $cd\in \cL_2(X)$ such that $S^j\varphi^n([ab])\subset \varphi^k(|cd])$
and equal to $0$ otherwise. Thus
\begin{eqnarray*}
\tilde{\phi}^{(\ell)}(x)&=&\sum_{cd\in \cL_2(X)}\Card\{0\le j<\ell
\mid S^j\varphi^n([ab])\subset\varphi^k([cd])\}(M_2^k\phi)(cd)\\
&=&\sum_{cd\in \cL_2(X)}M_2^{n-k}(ab,cd)(M_2^k\phi)(cd)=(M_2^n\phi)(ab).
\end{eqnarray*}

Consequently, we have for every $abc\in \cL_3(X)$,
\begin{displaymath}
\psi(bc)-\psi(ab)=(M_2^n\phi)(ab).
\end{displaymath}
Choose $m$ so large that $|\varphi^m(a)|\ge 2$ for every $a\in A$ and define $\theta\in Z_1(X)$ by
\begin{displaymath}
\theta(a)=\psi(a_1a_2)
\end{displaymath}
for $a\in A$ if $\varphi^m(a)=a_1\cdots a_\ell$.

If $ab$ belongs to $ \cL_2(X)$ with $\varphi^m(a)=a_1\cdots a_r$ and $\varphi^m(b)=b_1\cdots b_s$,
we obtain
\begin{eqnarray*}
(M_2^{n+m}\phi)(ab)&=&(M_2^n\phi)(a_1a_2)+\ldots+(M_2^n\phi)(a_r b_1)\\
&=&\psi(b_1b_2)-\psi(a_1a_2)=\theta(b)-\theta(a)\\
&=&(\partial_1\theta)(ab).
\end{eqnarray*}
If follows that $M_2^{n+m}\phi$ belongs to $\partial_1(Z_1(X))$. 
 Choosing $m$ large enough, we may assume
that $M_2^{n+m}\phi$ is in $\RR_{M_2}$. Since $M_2$ defines an automorphism
of $\RR_{M_2}$, 
and since, by Proposition~\ref{lemma0Substitutions}, the subspace $\partial_1(R_1(X))$
is invariant by $M_2$, we conclude that $\phi$ is an element of $\Delta_{M_2}\cap\partial_1(R_1(X))$.
Thus the kernel of $\sigma$ is included
in $\Delta_{M_2}\cap\partial_1(R_1(X))$. Since the converse inclusion is obvious, the conclusion follows.
\end{proof}

\begin{proofof}{of Proposition~\ref{propositionDimensionGroupSubstitutions}}
By Lemma~\ref{lemma2Substitutions}, the ordered
groups $(\Delta_{M_2},\Delta_{M_2}^+,1_{M_2})$ and $(G,G^+,1)$
 can be identified.
By Lemma~\ref{lemma3Substitutions} the morphism $\sigma$ defines an isomorphism
 from $\Delta_{M_2}/(\Delta_{M_2}\cap\partial_1(R_1(X))$ onto $H(X,S,\Z)$.
But since $PM_2=N_2P$, we have also $PM_2^k=N_2^kP$ for every $k\ge 1$.
Thus the projection $v\mapsto Pv$ maps $\Delta_{M_2}$ onto
$\Delta_{N_2}$ and
we obtain 
\begin{displaymath}
H(X,S,\Z)\simeq \frac{\Delta_{M_2}}{\Delta_{M_2}\cap\partial_1(R_1(X))}\simeq\Delta_{N_2}.
\end{displaymath}
Similarly, we have $H^+(X,S,\Z)\simeq \Delta^+_{N_2}$. Finally, the map
$\sigma$ sends $1$ to $1_X$ and we conclude that
$K^0(X,S)$ is isomorphic to $(\Delta_{N_2},\Delta^+_{N_2},P 1_{M_2})$.
\end{proofof}
We give two examples of computation of the dimension group
of a substitution shift. Other examples are treated in the exercises.
\begin{example}\label{exampleFibonacci6}
Let $\varphi:a\mapsto ab,b\mapsto a$ be the Fibonacci morphism and let $X$
be the Fibonacci shift.
\index{subject}{Fibonacci!shift}\index{subject}{shift space!Fibonacci}%
 As seen in Example \ref{exampleFibonacci5},
we have $N_2=M$ and 
\begin{displaymath}
M=\begin{bmatrix}1&1\\1&0\end{bmatrix},\quad 
P=\begin{bmatrix}0&1&1\\1&0&0\end{bmatrix}.
\end{displaymath}
The maximal eigenvalue of $M$ is  $\lambda=(1+\sqrt{5})/2$
and
the vector $[\lambda\ 1]$  is a left
eigenvector of $M$.
Thus $H(X,T,\Z)=\Z^2$ and $H^+(X,T,\Z)=\{(\alpha,\beta)\mid \alpha\lambda+\beta\ge 0\}$. The order unit is 
\begin{displaymath}
P\begin{bmatrix}1\\1\\1\end{bmatrix}=\begin{bmatrix}2\\1\end{bmatrix}.
\end{displaymath}
Thus, using the map $(\alpha,\beta)\mapsto \alpha\lambda+\beta$,
we see that the dimension group of $X$ is isomorphic
the group of algebraic integers $\Z+\Z\lambda$
with the order induced by the reals and the order unit 
$2\lambda+1$.
(see Example~\ref{exampleGolden1} and Example~\ref{exampleSturmian2}).

To obtain a normalized subgroup with $1$ as unit, we
consider the automorphism of $\Z^2$ such that $(\alpha,\beta)\mapsto (\alpha-\beta,2\beta-\alpha)$. This is actually  the automorphism defined by the
matrix $M^{-2}$ since
\begin{displaymath}
M^{-2}\begin{bmatrix}\alpha\\\beta\end{bmatrix}=
\begin{bmatrix}1&-1\\-1&2\end{bmatrix}\begin{bmatrix}\alpha\\\beta\end{bmatrix}=
\begin{bmatrix}\alpha-\beta\\2\beta-\alpha\end{bmatrix}.
\end{displaymath}

Then $(2,1)$ maps to $(1,0)$. Since 
\begin{displaymath}
\alpha\lambda+\beta=\lambda^2(\lambda(\alpha-\beta)+
(2\beta-\alpha)),
\end{displaymath}
 we do not change the order. Thus we conclude that
the dimension group of $X$ is isomorphic with $\Z+\lambda\Z$
with the order induced by the reals and $1$ as order unit.
\end{example}
\begin{example}\label{exampleMorse5}
Let $\varphi:a\mapsto ab,b\mapsto ba$ be the Thue-Morse morphism
on the alphabet $A=\{a,b\}$
and let $X$ be the Thue-Morse shift. We have $\cL_2(X)=A^2$.
Set $A_2=\{x,y,z,t\}$ and let $f:A_2\rightarrow A^2$
be the bijection $x\mapsto aa,y\mapsto ab,z\mapsto ba,t\mapsto bb$.
Then $\varphi_2$ is the substitution
$x\mapsto yz,y\mapsto yt,z\mapsto zx,t\mapsto zy$.
The Rauzy graph $\Gamma_2(X)$ is represented in Figure~\ref{figureGamma_2Morse}.
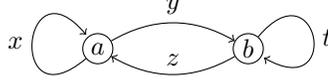
\begin{figure}[hbt]
\centering
\tikzset{node/.style={circle,draw,minimum size=0.4cm,inner sep=0pt}}
\tikzstyle{every loop}=[->,shorten >=1pt,looseness=12]
\tikzstyle{loop left}=[in=130,out=220,loop]
\tikzstyle{loop right}=[in=330,out=50,loop]
\begin{tikzpicture}
\node[node](1) at(0,0){$a$};
\node[node](2) at (2,0){$b$};

\draw[left](1) edge [loop left]node {$x$}(1);
%\draw[below,bend left, ->](1) edge node{$c$}(2);
\draw[above,bend left, ->](2) edge node{$z$}(1);
\draw[above,bend left, ->](1) edge node{$y$}(2);
%\draw[below,bend left, ->](2) edge node{$e$}(1);
\draw[right](2)edge[loop right]node{$t$}(2);
\end{tikzpicture}

\caption{The Rauzy graph $\Gamma_2(X)$.}\label{figureGamma_2Morse}
\end{figure}

The matrices $M$, $M_2$, $P$ and $N_2$
are 
\begin{displaymath}
M=\begin{bmatrix}1&1\\1&1\end{bmatrix},\quad
M_2=\begin{bmatrix}0&1&1&0\\0&1&0&1\\1&0&1&0\\0&1&1&0\end{bmatrix},\quad
P=\begin{bmatrix}1&0&0&0\\0&1&1&0\\0&0&0&1\end{bmatrix},\quad
N_2=\begin{bmatrix}0&1&0\\1&1&1\\0&1&0\end{bmatrix}
\end{displaymath}
Note that we already saw $M_2$ in Example \ref{exampleChristian}
where we have computed $\Delta_{M_2}$.
The eventual range $\RR_{N_2}$ of $N_2$ is generated by the vectors
\begin{displaymath}
v=\begin{bmatrix}1\\2\\1\end{bmatrix},\qquad 
w=\begin{bmatrix}1\\-1\\1\end{bmatrix},
\end{displaymath}
the first one being an eigenvector for the maximal eigenvalue $2$
and the second one for the eigenvalue $-1$.
Now, for $\alpha v+\beta w\in\RR_{N_2}$, since
\begin{displaymath}
  N_2^k(\alpha v+\beta w)=2^k\alpha v+(-1)^k\beta w
  =\begin{bmatrix}2^k\alpha+(-1)^k\beta\\2^{k+1}\alpha+(-1)^{k+1}\beta\\
  2^k\alpha+(-1)^k\beta\end{bmatrix},
\end{displaymath}
we have $N_2^k(\alpha v+\beta w)\in\Z^3$ if and only if
$\alpha=\frac{m}{3\cdot 2^k}$ and $\beta=\frac{n}{3}$
with $m,n\in\Z$ and $m+n\equiv 0\bmod 3$.
The order unit is 
\begin{displaymath}
P\begin{bmatrix}1\\1\\1\\1\end{bmatrix}=\begin{bmatrix}1\\2\\1\end{bmatrix}=v
\end{displaymath}

Thus, 
\begin{displaymath} 
\Delta_{N_2}\simeq\{(\alpha,\beta)\mid 3\alpha\in\Z[1/2],3\beta\in\Z, 3\alpha+3\beta\equiv 0\bmod 3\},
\end{displaymath}
with
 \begin{displaymath} 
\Delta_{N_2}^+\simeq\{(\alpha,\beta)\in\Delta_{N_2}\mid \alpha> 0\}\cup\{(0,0)\}.
\end{displaymath}
with order unit $(1,0)$.

%Thus,
%we obtain that the dimension group is isomorphic to the unital ordered group
%\begin{displaymath}
%(\Z[1/2]\times \Z,[(\Z_+[1/2]\setminus\{0\})\times\Z]\cup\{0\},(1,0))
%\end{displaymath}

The group $\Delta_{N_2}$ is actually
isomorphic to $\Z[1/2]\times\Z$, using the map $(\alpha,\beta)\to (\alpha+\beta,3\beta)$.
%the positive cone 
%$\Delta_{N_2}^+=\{(\gamma,\delta)\mid 3\gamma+\delta> 0\}\cup\{(0,0)\}$
%and the unit $(1,0)$.

To close the loop, let us express the  unique trace
on $K^0(X,S)$ which, by Proposition~\ref{propositionKerov} has the form
$\alpha_\mu$ where $\mu$ is the unique invariant probability
measure on $X$ (see Example~\ref{exampleMeasureMorse}).
We have
\begin{displaymath}
\alpha_\mu(\alpha v+\beta w)=\alpha
\end{displaymath}
since this map is a positive unital morphism from $K^0(X,S)$
to $(\R,\R_+,1)$. We find for example (in agreement with the
value given in Example~\ref{exampleMeasureMorse})
$\mu([aa])=1/6$ since the characteristic function of the cylinder $[aa]$
can be identified with the vector 
\begin{displaymath}
P\begin{bmatrix}1\\0\\0\\0\end{bmatrix}=\begin{bmatrix}1\\0\\0\end{bmatrix}
=\begin{bmatrix}1/2\\0\\1/2\end{bmatrix}+\begin{bmatrix}1/2\\0\\-1/2\end{bmatrix}
=1/6(v+2w)+\begin{bmatrix}1/2\\0\\-1/2\end{bmatrix}
\end{displaymath}
were the last expression is the decomposition in $\RR_{N_2}\oplus\KK_{N_2}$
\end{example}

%%%%%%%%%%%%%%%%%%%%%%%%%%%%
\section{Exercises}

\exosection{Section~\ref{sectionPartitionTowers}}
\begin{exercise}\label{exerciseRefinement}
Let $\Pg=\{T^jB_i\mid 1\le i\le m, 0\le j< h_i\}$ be
a partition in towers nested in a partition $\Pg'=\{T^hB'_k\mid 1\le k\le m',
0\le h\le h'_k\}$. Show that if
\begin{displaymath}
T^j B_i\subset T^h B'_k
\end{displaymath}
then
\begin{displaymath}
0\le j-h\le h_i-h'_k.
\end{displaymath}
\end{exercise}

\begin{exercise}\label{exerciseKR1pasKR3}
Let $\sigma$ be the substitution $\sigma:a\mapsto abb,b\mapsto aab$.
Let $\Pg(n)$ be the partition with basis
$\sigma^n(X)$ associated with $\sigma^n$. Let $\Pg'(n)$ be the
partition obtained by merging the two towers of $\Pg(n)$, that
is
\begin{displaymath}
\Pg'(n)=\{T^j\sigma^n(X)\mid 0\le j< 3^n\}
\end{displaymath}
Show that the sequence $(\Pg'(n))$ satisfies (KR1), (KR2) but not (KR3).
\end{exercise}

\begin{exercise}\label{exerciseKR3pasKR1}
Let $\sigma$ be the substitution $a\to acb,b\to bcb,c\to abb$.
Let $\Pg(n)$ be the partition associated with $\sigma^n$.
Show that the sequence $(\Pg(n))$ satisfies (KR2) and (KR3) but not (KR1).
\end{exercise}

\exosection{Section~\ref{chapterReturnWords}}
\begin{exercise}\label{exerciseSimplexAR}
  Let $S=\{x\in\R_+^A\mid \|x\|_1=1,\|x\|_\infty\le\frac{1}{d-1}$
  (see Figure~\ref{figureSimplexAR} for $d=3$).
  \begin{figure}[hbt]
    \centering
    \tikzset{node/.style={circle,draw,minimum size=0cm,inner sep=0pt}}
\tikzset{title/.style={minimum size=0cm,inner sep=0pt}}
    \begin{tikzpicture}
      \node[node](101)at(0,3){};\node[title]at(-.8,3){$[\frac 12,0,\frac 12]$};
      \node[node](112)at(2,3){};
      \node[node](011)at(4,3){};\node[title]at(4.6,3){$[0,\frac 12,\frac 12]$};
\node[node](211)at(1,1.5){};\node[node](121)at(3,1.5){};
      \node[node](110)at(2,0){};\node[title]at(2,-.4){$[\frac 12,\frac 12,0]$};
      \draw(101)--(011);\draw(011)--(110);\draw(110)--(101);
      \draw(211)--(112){};\draw(211)--(121){};\draw(112)--(121){};
      \fill[green](0,3)--(2,3)--(1,1.5);
      \fill[red](2,3)--(4,3)--(3,1.5);
      \fill[blue](1,1.5)--(3,1.5)--(2,0);
      \node[title]at(2,3.4){$[\frac 14,\frac 14,\frac 12]$};
      \node[title]at(.2,1.5){$[\frac 12,\frac 14,\frac 14]$};
      \node[title]at(3.7,1.5){$[\frac 14,\frac 12,\frac 14]$};
    \end{tikzpicture}
    \caption{The simplex $S$.}\label{figureSimplexAR}
  \end{figure}
  Show that for every $a\in A$ and $x\in S$, the vector $xM_a$
  is colinear to an element of $S$.
  \end{exercise}

\begin{exercise}\label{exerciseProjectiveDistance}
  For a positive vectors $x\in \R^n$, denote $\|x\|=\frac{\max x_i}{\min x_i}$.
  For positive vectors $x,y$ denote $xy=(x_iy_i)_{1\le i\le n}$
  and $x^{-1}=(x_i^{-1})_{1\le i\le n}$. Set
  $d(x,y)=\|xy^{-1}\|$ and
  \begin{equation}
    \delta(x,y)=\log d(x,y).
  \end{equation}
  Show that $\delta$ is a distance on the set $S_n$ of positive vectors $x$
  such that $\sum_{i=1}^nx_i=1$, called the \emph{projective distance}.
  \index{subject}{projective distance}%
\end{exercise}

\begin{exercise}\label{exerciseContraction}
  Show that if $M$ is a positive $m\times n$-matrix, then
  \begin{displaymath}
\delta(Mx,My)\le \delta(x,y)
  \end{displaymath}
  for every positive vectors $x,y\in\R^n$, where $\delta$ is the projective metric.
\end{exercise}
\begin{exercise}\label{exerciseBirkhoffContractionCoefficient1}
%The \emph{contraction coefficient}\index{subject}{contraction coefficient}
%of a positive matrix $M$ is
%\begin{equation}
%\tau(M)=\sup_{x,y\in S_n,x\ne y}\frac{\delta(Mx,My)}{\delta(x,y)}.
%\end{equation}
%This definition makes sense since $S_n$ is compact.
%By definition, we have
%\begin{equation}
% \tau(MN)\le\tau(M)\tau(N).
%  \end{equation}
%For a positive
%matrix $M$, set
%\begin{displaymath}
%  \phi(M)=\min_{i,j,r,s}\frac{M_{ij}M_{rs}}{M_{js}M_{ir}}
%\end{displaymath}

%  Show that for 
%  \begin{displaymath}
 %   M=\begin{bmatrix}1&2\\3&4\end{bmatrix},
%  \end{displaymath}
  % one has $\phi(M)=2/3$.
For an $m\times n$-matrix $M$, set $d(M)=\max_{i,j} d(m_i,m_j)$
where $m_i$ is the row of index $i$ of $M$. Set also
$\tau(r)=(\sqrt{r}-1)/(\sqrt{r}+1)$ for $r\ge 1$.

  A real number $c>0$ is called a \emph{contraction coefficient}
  \index{subject}{contraction coefficient}%
  \index{subject}{coefficient!contraction}%
  for a positive $m\times n$-matrix $M$ if
\begin{equation}
d(Mx,My)<d(x,y)^{c}\label{equationBirkhoff}
  \end{equation}
  whenever $d(M)>1$ and $d(x,y)>1$.

Prove, using the steps indicated below,
  that for every positive $m\times n$-matrix $M$, the real
  number $\tau(d(M))$, called \emph{Birkhoff contraction coefficient}
  \index{subject}{Birkhoff!contraction coefficient}%
  \index{names}{Birkhoff, Garrett}%
  is the smallest
  contraction coefficient of $M$.
%  \begin{equation}
%   \tau(M)=\frac{1-\sqrt{\phi(M)}}{1+\sqrt{\phi(M)}}.\label{equationBirkhoff}
%  \end{equation}
  \begin{enumerate}
  \item Show that it is enough to prove \eqref{equationBirkhoff} for $m=2$.
  \item Set
    \begin{displaymath}
      F(r,s,y)=\frac{\langle rs,y\rangle\langle \1,y\rangle}
      {\langle r,y\rangle\langle s,y\rangle}
    \end{displaymath}
    where $\langle x,y\rangle=\sum x_iy_i$ denotes the usual scalar product of $x,y$.
    Show that \eqref{equationBirkhoff} holds if and
    only if
    \begin{equation}
      F(r,s,y)<\|s\|^{\tau(\|r\|)}\label{equationF(r,s,y)}
    \end{equation}
    for every positive $r,s,y\in\R^n$ with $r,s$ nonconstant
    and where $\1$ is the vector with all components equal to $1$.
  \item Show that it is enough to prove \eqref{equationBirkhoff}
    for $n=m=2$. Hint: use the fundamental Theorem of Linear Programming
    \index{subject}{Fundamental Theorem!of Linear Programming}
    (see Appendix~\ref{appendixLinearAlgebra})
    which implies that among the set of nonnegative vectors $y$
    maximizing $\langle rs,y\rangle$
    under  the linear
    constraints $\langle r,y\rangle=1$ and $\langle s-\1,y\rangle=0$,
    there is one with at most two nonzero coordinates.
    \item Prove \eqref{equationF(r,s,y)} for $n=2$.
    \end{enumerate}
\end{exercise}

\exosection{Section~\ref{sectionGroupsRauzyGraphs}}
\begin{exercise}\label{exerciseExactSequenceGraph}
  Let $G$ be a strongly connected graph.
  Show that the sequence
  \begin{displaymath}
    0\to \Z(\Gamma)\edge{\kappa} \Z(E)\edge{\beta} \Z(V)\edge{\gamma}\Z
    \to 0,
    \end{displaymath}
where $\kappa$ is the composition map, $\beta(e)=r(e)-s(e)$
and $\gamma(v)=1$ identically, is exact.
  \end{exercise}
\exosection{Section~\protect{\ref{sectionDGSubstitutionShifts}}}

\begin{exercise}\label{exerciseChacon}
Consider the \emph{Chacon ternary substitution}
\index{subject}{Chacon!ternary!substitution}%
$\tau:a\mapsto aabc,b\mapsto bc,c\mapsto abc$. Show that
the dimension group of the associated shift space $X$ is
isomorphic to
$\Z[1/3]\times\Z$ with positive cone $\Z_+[1/3]\times\Z$
and unit $(3,-1)$.
\end{exercise}
\begin{exercise}\label{exerciseTribonacci}
Let $\varphi:a\mapsto ab,b\mapsto ac,c\mapsto a$
be the Tribonacci morphism and let $X$ be the corresponding
substitution shift.
\index{subject}{Tribonacci!shift}\index{subject}{shift space!Tribonacci}%
Show that the dimension group of $X,S)$ is the group
$\Z[\lambda]$ 
where $\lambda$ is the positive real solution of 
$\lambda^3=\lambda^2+\lambda+1$
(see also Example~\ref{exampleTribonacci3} where
we found the same result using return words).
\end{exercise}

%%%%%%%%%%%%%%%%%%%%%%%%%%
\section{Solutions}
\exosection{Section~\ref{sectionPartitionTowers}}
\begin{solution}{\ref{exerciseRefinement}}
Let $B,B'$ be the bases of $\Pg,\Pg'$. Since $\Pg$ is nested in
$\Pg'$, we have $B\subset B'$.
If $h>j$, then $T^{h-j}B'_k$ contains an element of $B\subset B'$,
a contradiction. Thus $0\le j-h$. Next, set $\ell=h_i-j$. Then 
 $T^{h+\ell}B'_k\subset B'$ implies $h+\ell\ge h'_k$. Thus
$j-h\le j-(h'_k-\ell)=h_i-h'_k$.
\end{solution}
\begin{solution}{\ref{exerciseKR1pasKR3}}
Since $\sigma$ is primitive and proper and $X(\sigma)$
is not periodic, the sequence $(\Pg(n))$ is a refining sequence of partitions
in towers by Lemma~\ref{lemmaPartitionProperSubstitution}
(to be proved in Chapter~\ref{ch5:sec:examples}). The sequence
$\Pg'(n)$ satisfies (KR1) (since $B'(n)=B(n)$) and (KR2) but it cannot
satisfy (KR3). Indeed, otherwise, $X(\sigma)$ would have a BV-representation
with one vertex at each level and would be an odometer 
(by Theorem~\ref{theoremBVrepresenationOdometers}).
\end{solution}

\begin{solution}{\ref{exerciseKR3pasKR1}}
Each $\Pg(n)$ is a partition in towers since $\sigma$ is primitive and $X(\sigma)$ is infinite.
The sequence satisfies (KR2) as any sequence of partitions built in this way.
It satisfies (KR3) because $\sigma^n[a]$ 
and $\sigma^n[c]$ tend both to $\sigma^\omega(b\cdot a)$ while
$\sigma^n[b]$ tends to $\sigma^\omega(b\cdot b)$. But condition (KR1) is not satisfied since there are two admissible fixed points.
\end{solution}
\exosection{Section~\ref{chapterReturnWords}}
\begin{solution}{\ref{exerciseSimplexAR}}
  Assume that $a$ is the first index and that $d=3$ for simplicity.
  The extreme points of $S$ are the vectors with one coefficient zero
  and all others equal to $1/(d-1)$.
   The image of the extreme
  points of $S$ are the vectors
  \begin{displaymath}
    \begin{bmatrix}0&\frac{1}{2}&\frac{1}{2}\end{bmatrix},
    \begin{bmatrix}\frac{1}{2}&\frac{1}{2}&1\end{bmatrix},
      \begin{bmatrix}\frac 12&1&\frac 12\end{bmatrix}
  \end{displaymath}
  After normalization to sum $1$, we obtain
  \begin{displaymath}
    \begin{bmatrix}0&\frac{1}{2}&\frac{1}{2}\end{bmatrix},
    \begin{bmatrix}\frac{1}{4}&\frac{1}{4}&\frac{1}{2}\end{bmatrix},
\begin{bmatrix}\frac{1}{4}&\frac{1}{2}&\frac{1}{4}\end{bmatrix},
  \end{displaymath}
  which all belong to $S$ (and actually to the red triangle of Figure~\ref{figureSimplexAR}).
  \end{solution}
\begin{solution}{\ref{exerciseProjectiveDistance}}
  It is clear that $\delta(x,y)\ge 0$ and, since the vectors have sum $1$,
  with equality if and only if $x=y$.
  We have to check the triangular inequality. Set $Q(x,y)=\max{x_i/y_i}$
  and $q(x,y)=\min x_i/y_i$.
  For every triple
  of positive vectors of sum $1$, we have
  \begin{displaymath}
x\le Q(x,y))y\le Q(x,y)Q(y,z)z
  \end{displaymath}
  and thus $x_i/z_i\le Q(x,y)Q(y,z)$ which implies $Q(x,z)\le Q(x,y)Q(y,z)$.
  Similarly $q(x,z)\ge q(x,y)q(y,z)$. Thus $\delta(x,z)\le \delta(x,y)+\delta(y,z)$.
 
\end{solution}
\begin{solution}{\ref{exerciseContraction}}
  Set $\hat{x}=Mx$ and $\hat{y}=My$. Since $M$ is positive,
  $\hat{x}$ and $\hat{y}$ are positive. We have
  \begin{displaymath}
\frac{\hat{x}_i}{\hat{y}_i}=\frac{a_{i1}x_1+\ldots+a_{in}x_n}{a_{i1}y_1+\ldots+a_{in}y_n}.
  \end{displaymath}
  Using the inequality
  \begin{displaymath}
    \min\frac{r_i}{s_i}\le\frac{r_1+\ldots+r_n}{s_1+\ldots+s_n}\le\max\frac{r_i}{s_i},
  \end{displaymath}
  we obtain $q(x,y)\le \hat{x}_i/\hat{y}_i\le Q(x,y)$ and thus
  the inequalities $Q(\hat{x},\hat{y})\le Q(x,y)$ and
  $q(\hat{x},\hat{y})\ge q(x,y)$ whence the inequality
  $d(Mx,My)\le d(x,y)$.
\end{solution}
\begin{solution}{\ref{exerciseBirkhoffContractionCoefficient1}}
  
  1. Assume the statement true for $m=2$.
    We have $d(Mx,My)=\frac{\langle m_i,x\rangle\langle m_i,y\rangle^{-1}}
    {\langle m_j,x\rangle\langle m_j,y\rangle^{-1}}$ for some indices $i,j$.
    Let $M'$ be the $2\times n$-matrix with rows $m_i,m_j$. Then
    $d(M'x,M'y)=d(Mx,My)>1$ and therefore $1<d(M')$. Then, by the hypothesis,
    \begin{eqnarray*}
      d(Mx,My)&=&d(M'x,M'y)<d(x,y)^{\tau(d(M'))}=d(x,y)^{\tau(d(m_i,m_j))}\\
      &<&d(x,y)^{\tau(M)}
    \end{eqnarray*}
    where the last inequality holds because $\tau$ is increasing.

    2. Consider a $2\times n$-matrix $M$ with rows $a,b$. We must
    prove that $\frac{\langle a,x\rangle\langle a,y\rangle^{(-1)}}
    {\langle b,x\rangle\langle b,y\rangle^{(-1)}}<d(x,y)^{\tau(d(a,b))}$
    or equivalently
    \begin{displaymath}
\frac{\langle a,b^{-1}x\rangle\langle a,b^{-1}y\rangle^{(-1)}}
    {\langle b,b^{-1}x\rangle\langle b,b^{-1}y\rangle^{(-1)}}<d(b^{-1}x,b^{-1}y)^{\tau(d(a,b)}.
    \end{displaymath}
    Set $r=ab^{-1}$ and $s=xy^{-1}$. Then we may rewrite the left-hand side
    above as
    \begin{eqnarray*}
     \frac{\langle ab^{-1},x\rangle\langle ab^{-1},y\rangle^{(-1)}}
          {\langle bb^{-1},x\rangle\langle bb^{-1},y\rangle^{(-1)}}
          &=&\frac{\langle ab^{-1},x\rangle\langle \1,y\rangle}
            {\langle ab^{-1},y\rangle\langle \1,x\rangle}\\
            &=&\frac{\langle r,sy\rangle\langle \1,y\rangle}
              {\langle r,y\rangle\langle \1,sy\rangle}=
              \frac{\langle rs,y\rangle\langle \1,y\rangle}
              {\langle r,y\rangle\langle s,y\rangle}=F(r,s,y),
    \end{eqnarray*}
    whence the desired formula.

    3. By scaling $r$ and $s$ suitably, we may assume that $\langle r,y\rangle=1$ and $\langle s-\1,y\rangle=0$. Then $F(r,s,y)=\langle rs,y\rangle$
    is linear in $y$. By the Fundamental Theorem of Linear Programming,
    there is a vector $y$ maximizing $\langle rs,y\rangle$ under
    the constraints $\langle r,y\rangle=1$ and $\langle s-\1,y\rangle=0$
    with at most two nonzero coordinates.
    Therefore, it is enough to prove \eqref{equationF(r,s,y)} for
    $r,s,y\in\R^2$.

    4. We may assume that $r,s,y$ have the form $(1,r),(1,s),(1,y)$
    for positive real numbers $r,s,y$
    with $r>1$ and $s\ne 1$. Then an easy computation gives
    \begin{displaymath}
      F(r,s,y)=\frac{1+rsy^2+(1+rs)y}{1+rsy^2+(r+s)y}
    \end{displaymath}
    If $s<1$, then $0>(r-1)(s-1)=(1+rs)-(r+s)$ so that
    $F(r,s,y)<1$. We may thus assume $s>1$.
    Since $1+rsy^2\ge 2y\sqrt{rs}$,
    we have $F(r,s,y)\le(\frac{1+\sqrt{rs}}{\sqrt{r}+\sqrt{s}})^2$,
    and it is sufficient to prove that for $r,s>1$,
    \begin{displaymath}
      \frac{1+\sqrt{rs}}{\sqrt{r}+\sqrt{s}}<(\sqrt{s})^{\tau(r)}
    \end{displaymath}
    Replacing $\sqrt{r},\sqrt{s}$ by $r,s$, we are left with proving
    that
    \begin{displaymath}
      \frac{1+rs}{r+s}<s^{\frac{r-1}{r+1}}
    \end{displaymath}
    Since both sides take the same value for $s=1$, it is enough to
    verify the corresponding inequality for their partial derivatives
    with respect to $s$, that is
    \begin{displaymath}
      \frac{r^2-1}{(r+s)^2}<\frac{r-1}{r+1}s^{\frac{-2}{r+1}}
    \end{displaymath}
    or equivalently
    \begin{displaymath}
      (\frac{r+1}{r+s})^2<s^{\frac{-2}{r+1}}
    \end{displaymath}
    The last inequality is equivalent to
    \begin{displaymath}
      (\frac{r+s}{r+1})>s^{\frac{1}{r+1}}.
    \end{displaymath}
    Again both sides take the same value for $s=1$. Taking again the
    partial devatives of both sides, we are left with the inequality
    \begin{displaymath}
      \frac{1}{r+1}>\frac{1}{r+1}s^{\frac{-r}{r+1}}.
    \end{displaymath}
    which is clearly true for $r>0$ and $s>1$.
\end{solution}

\exosection{Section~\ref{sectionGroupsRauzyGraphs}}
\begin{solution}{\ref{exerciseExactSequenceGraph}}
  Set $\Card(V)=n$ and $\Card(E)=m$. Clearly, $\kappa$ is injective
  and $\Z(\Gamma)=\im(\kappa)\subset\ker(\beta)$.
  Next $\im(\beta)=\ker(\gamma)$ and $\dim(\ker(\gamma)=n-1$.
  Thus
  \begin{displaymath}
m=\dim(\Z(E)=\dim(\ker(\beta))+\dim(\im(\beta))=\dim(\ker(\beta)+n-1.
    \end{displaymath}
  This implies that the dimension of $\ker(\beta)$ is $m-n+1$.
  On the other hand, since $G$ is strongly connected,
  any covering tree $T$ of $G$ has $n-1$ elements.
  Since $\Gamma$ has a basis of $\Card(E)-\Card(T)=m-n+1$
  elements (see Appendix~\ref{appendixGroups}),
  this implies that the dimension of $\Z(\Gamma)$ is $m-n+1$.
  We conclude that $\Z(\Gamma)=\ker(\beta)$.
  \end{solution}
\exosection{Section~\protect{\ref{sectionDGSubstitutionShifts}}}

\begin{solution}{\protect{\ref{exerciseChacon}}}
We have $\cL_2(X)=\{aa,ab,bc,ca,cb\}$. Set $A_2=\{x,y,z,t,u\}$
with the bijection taken in alphabetic order. The morphism
$\tau_2$ is $x\mapsto xyzt,y\mapsto xyzu,z\mapsto zt,t\mapsto yzt,u\mapsto yzy$.
The Rauzy graph $\Gamma_2(X)$ is represented in Figure~\ref{figureRauzyChacon}.
\begin{figure}[hbt]
\centering
\tikzset{node/.style={circle,draw,minimum size=0.4cm,inner sep=0.2pt}}

\begin{tikzpicture}
\node[node](t)at(0,0){$t$};\node[node](x)at(0,2){$x$};
\node[node](z)at(2,0){$z$};\node[node](y)at(2,2){$y$};
\node[node](u)at(4,0){$u$};

\draw[left,->](t)edge node{$a$}(x);\draw[left,above,->](t)edge node{$b$}(y);
\draw[above,->](x)edge node{$b$}(y);\draw[above,->](z)edge node{$a$}(t);
\draw[right,->](y)edge node{$c$}(z);
\draw[bend left,above,->](z)edge node{$b$}(u);
\draw[bend left,above,->](u)edge node{$c$}(z);
\end{tikzpicture}
\caption{The Rauzy graph $\Gamma_2(X)$.}\label{figureRauzyChacon}
\end{figure}
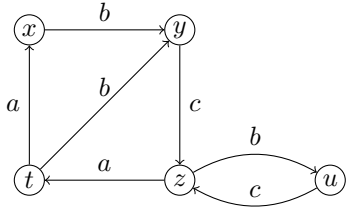

The matrices $M,M_2$, $P$
and $N_2$ are
\begin{displaymath}
M=\begin{bmatrix} 2&1&1\\0&1&1\\1&1&1\end{bmatrix},\quad
M_2=\begin{bmatrix}1&1&1&1&0\\0&0&1&1&0\\0&0&1&1&0\\0&1&1&1&0\\0&1&1&0&1
\end{bmatrix},\
\end{displaymath}
\begin{displaymath}
P=\begin{bmatrix}1&1&1&1&0\\0&1&1&1&0\\0&0&1&0&1\end{bmatrix},\quad
N_2=\begin{bmatrix}2&1&1\\1&1&1\\0&1&1\end{bmatrix}
\end{displaymath}
The eigenvalues of $N_2$ are $3,1,0$ and eigenvectors for $3,1$ are
\begin{displaymath}
v=\begin{bmatrix}3\\2\\1\end{bmatrix},\quad
w=\begin{bmatrix}1\\0\\-1\end{bmatrix}
\end{displaymath}
One has $N_2^n(\alpha v+\beta w)\in\Z^3$ if and only if $23^n\alpha\in\Z$
and $2\beta\in \Z$. Thus 
\begin{displaymath}
\Delta_{N_2}=\{(\alpha,\beta)\mid 2\alpha\in\Z[1/3], 2\beta\in\Z\},
\Delta^+_{N_2}=\{(\alpha,\beta)\in\Delta_{N_2}\mid \alpha>0\}\cup\{(0,0)\}.
\end{displaymath}
 The unit is
$P[1\ 1\ 1\ 1\ 1]^t=[4\ 3\ 1]^t=3/2v-1/2w$.
Using the isomorphism $(\alpha,\beta)\mapsto (2\alpha,2\beta)$,
we obtain the desired result.

\end{solution}
\begin{solution}{\protect{\ref{exerciseTribonacci}}}
Let $A_2=\{x,y,z,t,u\}$ be an alphabet in order preserving
bijection with
$\cL_2(X)=\{aa,ab,ac,ba,ca\}$. The morphism $\varphi_2:A_2^*\to A_2^*$
is $x\mapsto yt,y\mapsto yt,z\mapsto yt,t\to zu,u\mapsto x$.
The matrices $M$ and $M_2$ are
\begin{displaymath}
M=\begin{bmatrix}1&1&0\\1&0&1\\1&0&0\end{bmatrix},\quad
M_2=\begin{bmatrix}0&1&0&1&0\\0&1&0&1&0\\0&1&0&1&0\\0&0&1&0&1\\1&0&0&0&0
\end{bmatrix}.
\end{displaymath}
The Rauzy graph $\Gamma_2(X)$ is represented in Figure~\ref{figureRauzyTribo}.
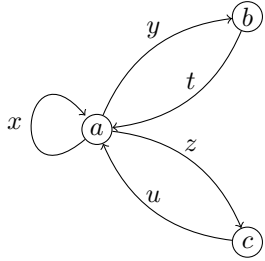
\begin{figure}[hbt]
\centering
\tikzset{node/.style={circle,draw,minimum size=0.4cm,inner sep=0.2pt}}
\tikzstyle{every loop}=[->,shorten >=1pt,looseness=12]
\tikzstyle{loop left}=[in=130,out=220,loop]
\begin{tikzpicture}
\node[node](a)at(0,1){$a$};
\node[node](b)at(2,2.5){$b$};\node[node](c)at(2,-0.5){$c$};

\draw[bend left,above,->](a)edge node{$y$}(b);
\draw[bend left,above,->](b)edge node{$t$}(a);
\draw[bend left,above,->](a)edge node{$z$}(c);
\draw[bend left,above,->](c)edge node{$u$}(a);
\draw[left,->](a)edge[loop left]node {$x$}(a);
\end{tikzpicture}
\caption{The Rauzy graph $\Gamma_2(X)$.}\label{figureRauzyTribo}
\end{figure}

Thus the matrix $P$ is
\begin{displaymath}
P=\begin{bmatrix}0&1&0&1&0\\0&0&1&0&1\\1&0&0&0&0\end{bmatrix},
\end{displaymath}
and $N_2=M$. The matrix
$M$ is invertible and its dominant eigenvalue
is the positive real number $\lambda$ such that $\lambda^3=\lambda^2+\lambda+1$.
A corresponding row eigenvector is $[\lambda^2,\lambda,1]$. Thus the dimension group is $\Z[\lambda]$.
\end{solution}

%\begin{solution}{\ref{exerciseHoltonZamboni}}
%\end{solution}
%%%%%%%%%%%%%%%%%%%%%%%ù
\section{Notes}

Kakutani-Rokhlin partitions owe their name to a result in ergodic
theory called \emph{Rokhlin's Lemma}\index{subject}{Rokhlin Lemma}
\index{subject}{Lemma!Rokhlin}\index{names}{Rokhlin, Vladimir A.}%
or \emph{Kakutani-Rokhlin Lemma}
\index{names}{Kakutani, Shizuo} which states that every aperiodic
measure-theoretic dynamical system can be represented by
an arbitrary high tower of measurable sets~\citep{Rokhlin1948,Kakutani1943}.

\subsection{Partitions in towers}
Proposition~\ref{propositionVersik}
 is~\cite[Lemma 3.1]{Putnam1989},
\index{names}{Putnam, Ian F.}%
where the credit of the construction is given to Vershik.
\index{names}{Vershik, Anatol M.}%
We follow the presentation of~\cite[Proposition 6.4.2]{Durand2010}.
\index{names}{Durand, Fabien}%
Theorem \ref{theoremKRPartitions} is
due to \cite{HermanPutnamSkau1992}.
\index{names}{Herman, Richard H.}%
\index{names}{Putnam, Ian F.}%
\index{names}{Skau, Christian F.}%
Theorem~\ref{theoremDimensionGroup} is also from \cite{HermanPutnamSkau1992}.
The proof of the fact that the dimension group is simple
is from \cite[Theorem 2.14]{Putnam2010}.

\subsection{Dimension groups and return words}

The results of Section~\ref{sectionReturnWords},
in particular Proposition \ref{propositionRecurrentSubshift}
      are from \cite{Host1995}
(see also~\cite{Host2000}). The unique ergodicity of Arnoux-Rauzy
shifts is a result of~\cite{DelecroixHejdaSteiner2013}.
Lemma~\ref{lemmaAvilaDelecroix} is from \cite{AvilaDelecroix2013}.
\index{names}{Avila, Artur}\index{names}{Delecroix, Vincent}%

The ergodic properties of products of sequences of matrices have
been investigated, especially in the context of nonhomogeneous Markov
chains defined by sequences of stochastic matrices
(instead of powers of a single stochastic matrix as in an ordinary
finite Markov chain).
Classical references are \cite{Hartfiel2002}\index{names}{Hartfiel, Darald J.}
and \cite{Seneta2006}.\index{names}{Seneta, Eugene}
The projective distance (Exercise~\ref{exerciseProjectiveDistance}), also called 
called \emph{Hilbert projective distance}\index{subject}{Hilbert!projective distance}\index{names}{Hilbert, David} can be found in \cite{Hartfiel2002}
(including the link with Hilbert original definition) and \cite{Seneta2006}.
Birkhoff contraction coefficient (Exercise~\ref{exerciseBirkhoffContractionCoefficient1}) was introduced in~\cite{Birkhoff1957} (see also~\cite[pp. 383-386]{Birkhoff1967}). The original proof by Birkhoff uses projective geometry.
The elementary proof presented here is from \cite{Carroll2004}.
\index{names}{Carrol, Joseph E.} 
\subsection{Dimension groups and Rauzy graphs}

The notion of fundamental group of a graph used
in Section~\ref{chapterDimensionGroupsRauzyGraphs}
 is classical in algebraic topology
(see~\cite{LyndonSchupp2001} for a
more detailed introduction to its direct definition
on a graph and its connection with spanning trees of the graph).
The ordered group $C(\Gamma)$ associated to a graph $\Gamma$
is called a \emph{graph group}\index{subject}{graph!group of}
in \cite{BoyleHandelman1996} while a direct limit of graph
groups as in Proposition~\ref{propositionDirectLimitGroupsRauzy}
is called a \emph{graphical group}.

\subsection{Dimension group of a substitution shift}
The results of this section, in particular Proposition~\ref{propositionDimensionGroupSubstitutions}, are from \cite{Host1995} (see also~\cite{Host2000}).

%%%%%%%%%%%%%%%%%%%%%%%%
%  Bratteli Diagrams %
%%%%%%%%%%%%%%%%%%%%%%%

\chapter{Bratteli diagrams}\label{chapterBratteliDiagrams}
We now introduce Bratteli diagrams. We will see that,
adding an order on the diagram and provided this order
is what we will call proper,
an ordered
Bratteli diagram defines in a natural way a topological dynamical system.
We prove the Bratteli-Vershik representation theorem:
any minimal topological dynamical system defined on a Cantor set can be obtained
in this way (Theorem~\ref{ch5:theo:BVmodel}). Thus every
minimal Cantor system $(X,T)$ can be represented by an 
ordered  Bratteli diagram,
called a BV-representation of $(X,T)$.

We will next introduce an equivalence on dynamical systems
called Kakutani equivalence and prove that it can be characterized
by a transformation on BV-representations.

We then prove one of the major results presented in this book,
namely the Strong Orbit Equivalence Theorem (Theorem~\ref{theoremStrongOrbitEquivalence}). This result shows that the dimension group is a complete
invariant for the so-called strong orbit equivalence.
As a complement, the Orbit Equivalence Theorem (Theorem~\ref{theoremOrbitEquivalence}) shows that the quotient of the dimension group by the
infinitesimal subgroup is a complete invariant for orbit equivalence.

In Section~\ref{sectionEtaleEquiv}, we develop a systematic study of equivalences
on Cantor spaces. We introduce the notion of \'etale equivalence relation
and prove that both the relations of orbit equivalence and of cofinality
in Bratteli diagrams are \'etale equivalences.

In the last section (Section~\ref{ch5:sec:entropy}), we discuss the
link with the notion of entropy (which had not been considered before
in this book).

%%%%%%%%%%%%%%%%%%%%%%%%%%%%%%
\section{Bratteli diagrams}\label{ch5:sec:2}
%%%%%%%%%%%%%%%%%%%%%%%%%%%%%%
\medskip
A {\em Bratteli diagram}\index{subject}{Bratteli diagram}\index{names}{Bratteli, Ola}
is an infinite directed graph $(V, E)$\index{symbols}{V@$(V,E)$}
 where the {\em vertex}
\index{subject}{vertex!of Bratteli diagram} set $V$ and the {\em edge}\index{subject}{edge}
set $E$ can be partitioned into nonempty finite sets
$$
V = V(0) \cup V(1) \cup V(2) \cup \cdots \ \ {\rm and} \ \ E = E(1) \cup E(2)
\cup \cdots
$$
with the following properties:
\begin{enumerate}
\item
$V(0) = \{ v(0) \}$ is a one-point set,
\item
$r(E(n)) \subseteq V(n), \; s(E(n)) \subseteq V (n-1), n = 1,2,
\ldots $,
\end{enumerate}
where $r:E\to V$\index{symbols}{r(e)@$r(e)$} is called the 
{\em range map}\index{subject}{range!map}\index{subject}{Bratteli diagram!range map}\index{subject}{map!range}
and $s:E\to V$\index{symbols}{s(e)@$s(e)$} the 
{\em source map}\index{subject}{source!map}\index{subject}{Bratteli diagram!source map}\index{subject}{map!source}. 
They satisfy $s^{-1} (v) \neq \emptyset$  for all $v \in V$ and
$r^{-1} (v) \neq \emptyset$  for all $v \in V \setminus V (0)$.

We use the terminology of graphs to handle Bratteli diagrams
(see also Appendix~\ref{appendixGroups}).
Note that Bratteli diagrams are actually  multigraphs
(there can be several edges between two vertices).
Since edges can only go from $V(n)$ to $V(n+1)$,
they are \emph{acyclic}, that is, there is no nontrivial cycle.
Some particular terms are of common use for acyclic graphs.
\index{subject}{acyclic graph}\index{subject}{graph!acyclic}
In particular, the vertex $v(0)$ is called the \emph{root}.
\index{subject}{root of Bratteli diagram}%
\index{subject}{Bratteli diagram!root of}%
A \emph{successor}\index{subject}{successor!of vertex}
 of a vertex $v\in V$ is a vertex
$w$ such that $v=s(e)$ and $w=r(e)$ for some $e\in E$.
Two edges $e,f\in E$ are \emph{consecutive}\index{subject}{consecutive!edges}
if $r(e)=s(f)$.
A \emph{path}\index{subject}{path!in Bratteli diagram} is
a sequence $(e_1,e_2,\ldots,e_n)$
of consecutive edges. The \emph{source}\index{subject}{source!of path}
of the path is $s(e_1)$ and its \emph{range}\index{subject}{range!of path}
is $r(e_n)$. A \emph{descendant}\index{subject}{descendant
of vertex} of a vertex $v$ is a vertex $w$ such that there is a path
from $v$ to $w$, that is a path with source $v$ and range $w$.

\begin{figure}[hbt]
\tikzset{node/.style={circle,draw,minimum size=0.1cm,inner sep=0pt}}
	\tikzset{title/.style={minimum size=0.5cm,inner sep=0pt}}
\begin{tikzpicture}
%\node[title]at (-4,5){level};
\node[title]at(-2,5){adjacency matrices};
%\node[title]at(-4,4){$n-1$};
\node[title]at(-2,3){$M(n)=\begin{bmatrix}2&0\\0&1\\1&1\\1&0\\0&1\end{bmatrix}$};
\node[node](v11)at(3,4){};
\node[node](v12)at(5,4){};
\node[title]at (7,4){$V(n-1)$};
\node[title]at (1,3){$E(n)$};
\node[node](v21)at(2,2){};\node[node](v22)at(3,2){};\node[node](v23)at(4,2){};
\node[node](v24)at(5,2){};\node[node](v25)at(6,2){};
\node[title]at(7,2){$V(n)$};
\node[title]at (-2,1){$M(n+1)=\begin{bmatrix}1&1&0&0&0\\0&0&1&1&1\end{bmatrix}$};
\node[title]at(1.5,1){$E(n+1)$};
\node[node](v31)at(3,0){};\node[node](v32)at(5,0){};\node[title]at(7,0){$V(n+1)$};

\draw[bend right,left](v11)edge node{}(v21);
\draw[right,below](v11)edge node{}(v21);
\draw[left,below](v11)edge node{}(v23);
\draw[left,below](v11)edge node{}(v24);
\draw[right](v12)edge node{}(v22);
\draw[right,below](v12)edge node{}(v23);
\draw[right](v12)edge node{}(v25);

\draw[left,below](v21)edge node{}(v31);
\draw[right,below,right](v22)edge node{}(v31);
\draw[left,below,left](v23)edge node{}(v32);
\draw[left,below,left](v24)edge node{}(v32);
\draw[right,below,right](v25)edge node{}(v32);
\end{tikzpicture}
\caption{Representation of a diagram between the levels $n-1$ and $n+1$.}\label{ch5:fig:ch5-diag1}
\end{figure}
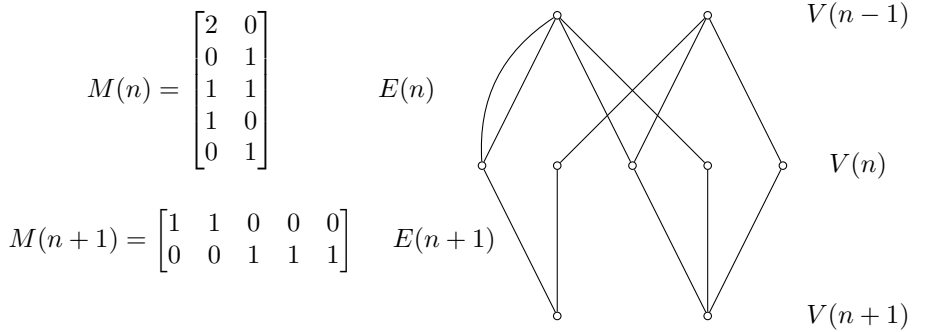

It is convenient to represent  the Bratteli
diagram by a picture with $V(n)$\index{symbols}{V@$V(n)$} the vertices at (horizontal) level $n$, and $E(n)$\index{symbols}{E@$E(n)$}
the edges (downward directed) connecting the vertices at level $n-1$
with those at level $n$.
Also, if $\Card(V(n-1))= t (n-1)$ and
$\Card(V(n))=t (n)$, then  $E(n)$ determines a $t (n) \times t (n-1)$ 
{\em adjacency matrix}\index{subject}{adjacency matrix!of  Bratteli diagram}\index{subject}{Bratteli diagram!adjacency matrix} $M(n)$\index{symbols}{M@$M(n)$}
defined by
\begin{equation}
M(n)_{r,s}=\Card\{e\in E(n)\mid r(e)=r,s(e)=s\}\label{eqIncidenceMatrixBD}
\end{equation}
(see Figure~\ref{ch5:fig:ch5-diag1}). Note that $M(n)$ is not
exactly the adjacency matrix of the graph $(V,E)$ because
 the edges are taken from range to source
(or, equivalently, the matrix is transposed) and moreover
because the adjacency matrix of $(V,E)$ is a
$V\times V$-matrix instead of a $V(n)\times V(n-1)$-matrix.

We say that  two Bratteli diagrams $(V, E)$
and $(V' , E' )$ are 
{\em isomorphic}\index{subject}{Bratteli diagram!isomorphic} 
whenever there exists a pair of bijections $f:V\to V'$, preserving the degrees, and $g:E \to E'$,
intertwining the respective source and range maps: $$s'\circ g = f\circ s \mbox{ and } r' \circ g = f\circ r \ . $$

Let $k , l \in \N$ with $1\le k < l$
and let $E_{k,l}$ denote the  set of paths
from $V(k)$ to $V(l)$.
Specifically,
$$
E_{k,l}=\{ ( e_{k} , \ldots , e_l ) \mid 
e_i \in E(i) ,  k  \leq i \leq  l ,
r(e_i ) = s(e_{i+1} ), k  \leq i \leq  l - 1 \}.
$$

Remark that  the adjacency matrix of $E_{k,l}$ is $M(l) \cdots M(k)$.
We define $r(e_{k} , \ldots , e_l ) = r(e_l)$ and $s(e_{k} ,
\ldots , e_l )  = s(e_{k})  .$

\subsection{Telescoping and simple diagrams}\label{sectionTelescopingSimple}
\begin{figure}
\centering
\tikzset{node/.style={circle,draw,minimum size=0.1cm,inner sep=0pt}}
	\tikzset{title/.style={minimum size=0.5cm,inner sep=0pt}}
\begin{tikzpicture}
\draw (-3, -1.5) node {$
\begin{bmatrix}
2&1\\
2&2
\end{bmatrix}
  $}  ;

\node[node](11)at(2,0){} ;
\node[node](12)at(5,0){} ;
\node[node](21)at(2,-3){} ;
\node[node](22)at(5,-3){} ;
\node[title]at(6,0){$V (n-1)$} ;
\node[title]at(6,-3){$V (n+1)$} ;
\draw (-.5, -1.5) node {$E_{n,n+1} $} ;
\draw (11) -- (21) ;
\draw (12) -- (22) ;
\draw (11) -- (22) ;
\draw (12) -- (21) ;
\draw (11) to[bend right] (21) ;
\draw (11) to[bend right] (22) ;
\draw (12) to[bend left] (22) ;
\end{tikzpicture}
\caption{Telescoping between the levels $n-1$ and $n+1$ in the diagram of Figure~\ref{ch5:fig:ch5-diag1}}\label{ch5:fig:ch5-diag2}
\end{figure}
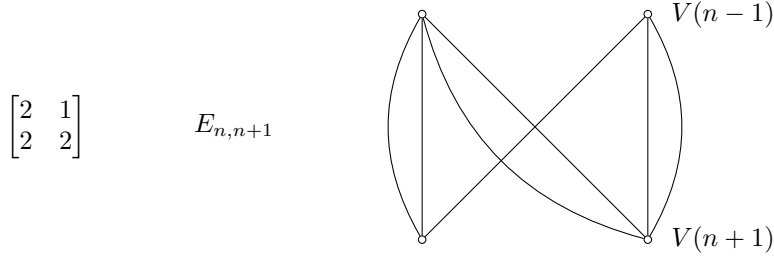

Given a Bratteli diagram $(V,E)$ and a sequence
$$
m_0 = 0 < m_1 < m_2 < \ldots
$$
in $\N$, we define the 
{\em telescoping}\index{subject}{telescoping!of Bratteli diagram}\index{subject}{Bratteli diagram!telescoping} 
of $(V,E)$ with respect to $\{ m_n \mid  n\in \N \}$
as the new Bratteli diagram $(V' , E' )$, where
$V' (n) = V (m_{n})$ and $E '(n) = E_{m_{n-1}+1, m_{n}}$
and the range and source maps are as above (see Figure~\ref{ch5:fig:ch5-diag2}).

We say that $(V,E)$ is a 
{\em simple Bratteli diagram}\index{subject}{Bratteli diagram!simple}\index{subject}{simple!Bratteli diagram} 
if there exists a
telescoping $(V',E')$ of $(V,E)$ such  that the adjacency matrices of
$(V',E')$ have only non-zero entries at each level.

We will use the following characterization of simple diagrams.
Let $(V,E)$ be a Bratteli diagram. A set $W\subset V$ is \emph{directed}
\index{subject}{directed!set} if  every edge having its source in $W$
has also its range in $W$. In symbols, for every $e\in E$
\begin{displaymath}
s(e)\in W\Rightarrow r(e)\in W.
\end{displaymath}
It is \emph{hereditary}\index{subject}{hereditary set} if
it satisfies for every $v\in V$ the following condition.
If every edge with source $v$
has its range in $W$, then $v$ itself is in $W$. In symbols, for every
$v\in V$ 
\begin{displaymath}
r(e)\in W\mbox{ for every edge $e$ such that $v=s(e)$}\Rightarrow v\in W.
\end{displaymath}

\begin{proposition}\label{propositionDirectedHereditary}
A Bratteli diagram is simple if and only if there
is no nonempty set both directed and hereditary 
other than $V$.
\end{proposition}
\begin{proof}
Assume first that $(V,E)$ is simple. Let $W\subset V$ be a
nonempty directed
and hereditary set. Since $W$ is nonempty, there is at least one $w$ in $W$.
Let $n$ be such that $w\in V(n)$.
Since $(V,E)$ is simple, there is an $m>n$ such that there is a path
from $w$ to every vertex in $V(m)$. Since $W$ is directed, this
implies $V(m)\subset W$. Since $W$ is hereditary, this implies
that all vertices of $V(n)$ for $n\le m$ are in $W$. Thus
$v(0)\in W$, which implies $V=W$.

Conversely, assume that $(V,E)$ is not simple. Let $v\in V(n)$ be such
that for every $m>n$ there is some $w\in V(m)$ which cannot be reached
from $v$. Consider the set $W$ of vertices $w\in V(m)$
for some $m\ge n$
for which there is an integer $p=p(w)>n$ such that all descendants of $w$
in $V(p)$ are descendants of $v$. It is a directed
set by definition. Suppose that some vertex $w\in V$ is such that
all its successors belong to $W$. Let $p$ be the supremum  of
the integers $p(u)$ for $u$ successor of $w$. Then all
descendants of $w$ in $V(p+1)$ are descendants of $v$ and thus
$w$ is in $W$. This shows that $W$ is hereditary.
Finally there is at least one vertex  in $V(n)$ 
which is not in $W$ since otherwise taking the supremum
of the integers $p(u)$ for $u\in V(n)$, we find that
all vertices in $V(p)$ are descendants of $v$.
Thus $W$ is a nonempty directed and hereditary set strictly
contained in $V$.
\end{proof}

\begin{example}\label{exampleDirectedHereditary}
Consider the Bratteli diagram represented in Figure~\ref{figureDirectedHereditary}.
\begin{figure}[hbt]
\centering
\tikzset{node/.style={circle,draw,minimum size=0.2cm,inner sep=0pt}}
	\tikzset{title/.style={minimum size=0.5cm,inner sep=0pt}}
\begin{tikzpicture}
\node[node](0)at(0,1){};
\node[node,fill=red](11)at(1,0){};\node[node](12)at(1,2){};
\node[node,fill=red](21)at(3,0){};\node[node](22)at(3,2){};
\node[node,fill=red](31)at(5,0){};\node[node](32)at(5,2){};
\node[node,fill=red](41)at(7,0){};\node[node](42)at(7,2){};
\node[title](51)at(8,0){$\cdots$};\node[title](52)at(8,2){$\cdots$};

\draw(0)edge node{}(11);\draw(0)edge node{}(12);
\draw(11)edge node{}(21);\draw(12)edge node{}(21);\draw(12)edge node{}(22);
\draw(21)edge node{}(31);\draw(22)edge node{}(31);\draw(22)edge node{}(32);
\draw(31)edge node{}(41);\draw(32)edge node{}(41);\draw(32)edge node{}(42);
\end{tikzpicture}
\caption{A non simple Bratteli diagram}\label{figureDirectedHereditary}
\end{figure}
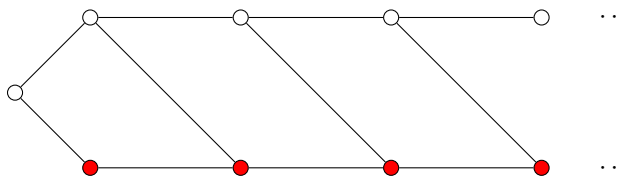

This diagram is not simple because the vertices of the lower level
can never reach the top level. Accordingly, the vertices at lower
level (excluding the root) form a directed and hereditary set.
\end{example}

We  denote by  $\sim$ the {\em telescoping equivalence}
\index{subject}{telescoping!equivalence}
  on Bratteli diagrams\index{subject}{Bratteli diagram!telescoping!equivalence}
 as the equivalence  relation
generated by isomorphism and telescoping.  It is not hard to show that
$(V^1,E^1) \sim (V^2,E^2)$ if and only if there exists a Bratteli diagram
$(V,E)$ such  that telescoping $(V,E)$ to odd levels $0 < 1 < 3 <
\ldots$ yields a telescoping of either $(V^1,E^1)$ or $(V^2,E^2)$, and
telescoping
$(V,E)$ to even levels $0 < 2 < 4 < \ldots$ yields a telescoping of the
other (Exercise~\ref{exerciseEquivalenceTelescoping}).
\begin{example}
Consider the Bratteli diagram of Figure~\ref{figureTelecopingEquivalence}
in the middle. Telescoping at even levels gives the diagram
on the left and telescoping at odd levels gives the diagram 
on the right. Thus the left and right diagrams are equivalent.
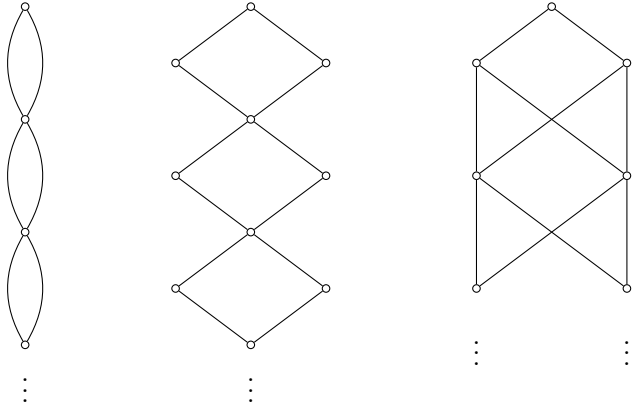
\begin{figure}[hbt]
\centering
\tikzset{node/.style={circle,draw,minimum size=0.1cm,inner sep=0pt}}
	\tikzset{title/.style={minimum size=0.5cm,inner sep=0pt}}
\begin{tikzpicture}
%odo
\node[node](0)at(1,4.5){};
\node[node](1)at(1,3){};
\node[node](2)at(1,1.5){};
\node[node](3)at(1,0){};
\node[title](4)at(1,-.5){$\vdots$};

\draw[bend right,left](0)edge node{}(1);\draw[bend left,right](0)edge node{}(1);
\draw[bend right,left](1)edge node{}(2);\draw[bend left,right](1)edge node{}(2);
\draw[bend right,left](2)edge node{}(3);\draw[bend left,right](2)edge node{}(3);
%odo2
\node[node](0)at(4,4.5){};
\node[node](11)at(3,3.75){};\node[node](12)at(5,3.75){};
\node[node](1)at(4,3){};
\node[node](21)at(3,2.25){};\node[node](22)at(5,2.25){};
\node[node](2)at(4,1.5){};
\node[node](31)at(3,0.75){};\node[node](32)at(5,0.75){};
\node[node](3)at(4,0){};
\node[title](41)at(4,-0.5){$\vdots$};

\draw(0)edge node{}(11);\draw(0)edge node{}(12);
\draw[left,near end](11)edge node{}(1);\draw[left,near end](12)edge node{}(1);
\draw[right,near end](1)edge node{}(21);\draw[right,near end](1)edge node{}(22);
\draw[left, near end](21)edge node{}(2);\draw[left, near end](22)edge node{}(2);
\draw[right,near end](2)edge node{}(31);\draw[right,near end](2)edge node{}(32);
\draw(31)edge node{}(3);\draw(32)edge node{}(3);

%Morse
\node[node](0)at(8,4.5){};
\node[node](11)at(7,3.75){};\node[node](12)at(9,3.75){};

\node[node](21)at(7,2.25){};\node[node](22)at(9,2.25){};
\node[node](31)at(7,0.75){};\node[node](32)at(9,0.75){};
\node[title](41)at(7,0){$\vdots$};\node[title](42)at(9,0){$\vdots$};

\draw(0)edge node{}(11);\draw(0)edge node{}(12);
\draw[left,near end](11)edge node{}(21);\draw[right,near end](11)edge node{}(22);
\draw[left,near end](12)edge node{}(21);\draw[right,near end](12)edge node{}(22);
\draw[left, near end](21)edge node{}(31);\draw[right,near end](21)edge node{}(32);
\draw[left,near end](22)edge node{}(31);\draw[right,near end](22)edge node{}(32);

\end{tikzpicture}
\caption{Three  equivalent Bratteli diagrams.}
\label{figureTelecopingEquivalence}
\end{figure}
\end{example}

\subsection{Dimension group of a Bratteli diagram}

Let $(V,E)$ be a Bratteli diagram. Let $V(n)=\{v_1,\ldots,v_{t(n)}\}$
and $G(n)=\Z^{t(n)}$. Let also $n_j$ be the number of paths
from $v(0)$ to $v_j\in V(n)$. We consider $G(n)$ as a unital ordered
group with the usual order and the unit 
$u(n)=\begin{bmatrix}n_1&\ldots&n_{t_n}\end{bmatrix}^t$. 
The \emph{dimension group}
\index{subject}{dimension group!of Bratteli diagram}%
\index{subject}{Bratteli diagram!dimension group}%
of $(V,E)$, denoted $\DG(V,E)$\index{symbols}{DG@$\DG(V,E)$}
 is the direct limit of the sequence
\begin{displaymath}
G(0)\edge{M(1)}G(1)\edge{M(2)}G(2)\ldots
\end{displaymath}
defined by the adjacency matrices $M(n)$.
%Thus by definition $\DG(V,E)=(G,G^+,1_G$) where $G=\Delta/\Delta^0$ with
%\begin{equation}
%\Delta=\{(x_n)\in\prod\Z^{t_n}\mid x_{n+1}=M(n)x_n\mbox{ for all large enough $n$}\}
%\label{equationDelta}
%\end{equation}
%and $\Delta^0=\{x\in\Delta\mid x_n=0\mbox{ for all $n$ large enough}\}$.

\begin{example}\label{exampleDBratteli}
The dimension group of the diagram represented in Figure~\ref{figureTelecopingEquivalence} on the left is the group $\Z[1/2]$ (see Example~\ref{exampleDyadic}).
\end{example}
The following result shows that the group $\DG(V,E)$ is a complete
invariant for the telescoping equivalence.
\begin{theorem}\label{theoremDGBratteli}
Two Bratteli diagrams $(V,E)$ and $(V',E')$ are telescoping
equivalent  if and only if the unital ordered groups $\DG(V,E)$
and $\DG(V',E')$ are isomorphic.
\end{theorem}
\begin{proof}
Taking a subsequence (starting at $0$) does not change the 
direct limit and thus telescoping does not change the dimension group.

Conversely, Let $(V, E)$ 
and $(V',E')$ be two Bratteli diagrams. 
Set $V(n)=\{v_1,\ldots,v_{t(n)}\}$
and $V'(n)=\{v'_1,\ldots,v'_{t'(n)}\}$. Set $G=\DG(V,E)$
and $G'=\DG(V',E')$.
We shall construct a Bratteli diagram $(W,F)$ that contracts to a contraction of $(V,E)$ on odd levels and to a contraction of $(V',E')$ on even levels.
It suffices to give the sets of vertices $W(n)$ and the incidence matrices $N(n)$ between consecutive levels.

We set $W(1) =  V (1) $ and $N (1) = M (1)$.
Looking at the canonical generators of $\Z^{t(1)}$ as elements of $G'$, we can consider  that 
they are elements of some $\Z^{t'(n_2)}$.
We set $W(2)  = V'(n_2)$,  and 
denote $N(2)$ the matrix of the
map it defines from $\Z^{t(1)}$ to $\Z^{t'(n_2)}$.
Again, the elements of $\Z^{t'_{n_2}}$ can be considered as elements of $G$,  
and thus 
belong to some $\Z^{t(n_3)}$.
We set $W(3)  = V(n_3)$ and 
we call $N(3)$ the map that  it defines from $\Z^{t'(n_2)}$ to $\Z^{t(n_3)}$.
Proceeding like this, we obtain the sequence 
$$
\Z 
\stackrel{N(1)}{\longrightarrow} 
\Z^{t(1)}
\stackrel{N(2)}{\longrightarrow}
\Z^{t'(n_2)}
\stackrel{N(3)}{\longrightarrow}
\Z^{t(n_3)}
\cdots
$$
that is sufficient to define the Bratteli diagram  we are looking for. 
\end{proof}
We illustrate this result with the following example.
\begin{example}
Consider the two diagrams of Figure~\ref{figureTelecopingEquivalence}
on the left and on the right. We have already seen 
in Example~\ref{exampleDBratteli}
that the dimension group of the first one is $\Z[1/2]$ obtained
as the direct limit of the sequence $\Z\edge{2}\Z\edge{2}\cdots$.
The dimension group of the second one is the direct limit
of the sequence $\Z^2\edge{M}\Z^2\edge{M}\cdots$
where $M$ is the matrix
\begin{displaymath}
M=\begin{bmatrix}1&1\\1&1\end{bmatrix}
\end{displaymath}
Since $\RR_M=\{\begin{bmatrix}x&x\end{bmatrix}^t\mid x\in\R\}$,
the isomorphism of the dimension groups is consequence of the
commutative diagram below.
\begin{displaymath}
\begin{CD}
x@> >>\begin{bmatrix}x&x\end{bmatrix}^t\\
@VV 2 V @VV M V\\
2x@> >>\begin{bmatrix}2x&2x\end{bmatrix}^t
\end{CD}
\end{displaymath}
\end{example}
Theorem~\ref{theoremDGBratteli} means that the properties
of a Bratteli diagram, or at least of its equivalence
class for telescoping should be read on its dimension
group. A first step in the direction is the following statement.
\begin{proposition}\label{propositionEquiSimpleDGDiagram}
A Bratteli diagram is simple if and only if its dimension group is simple.
\end{proposition}
\begin{proof}
Let $(V,E)$ be a Bratteli diagram and $G=\DG(V,E)$. Let
us first suppose that $(V,E)$ is simple. We have to show
that every nonzero element $g\in G^+$ is an order unit.
Let indeed $g\in G^+$ be nonzero and let $h\in G^+$.
We can choose $n\ge 1$ such that $g=i_n(x)$ 
and $h=i_n(y)$ with $x,y\in G(n)^+$.
Let $v\in V(n)$ be such that $x_v>0$. Since $(V,E)$ is simple,
we have $M(m)\cdots M(n+1)x>0$ for all $m$ large enough.
Then $M(m)\cdots M(n+1)y\le N M(m)\cdots M(n+1)x$ for $N$ large
enough.
This implies that $h<Ng$. Thus $G$ is simple.

Conversely, assume that $(V,E)$ is not simple. Let $W$ be a nonempty directed 
and hereditary set strictly contained in $V$. Let
$H$ be the set of $h\in D(V,E)$ which correspond to an 
$x=(x_n)$ with $x_n\in \Z^{V(n)}$ having the property that
for some $n\ge 1$
we have $x_{n,v}=0$ for every $v\notin W$
and $x_{m+1}=M(m+1)x(m)$ for all $m\ge n$.. Since $W$ is directed, 
every matrix $M(m)$ has the form
\begin{displaymath}
M(m)=\kbordermatrix{&W & \cr W&&\vrule & \cr \cline{2-4} &0&\vrule&}
\end{displaymath}
Thus, if $x$ satisfies
this property for $n$, it
holds for every $m\ge n$. Thus $H$ is a subgroup of $G$
which is clearly an ideal. Since $W\ne\emptyset$, we can
choose $w\in W\cap V(n)$ and $x$ such that $x_{n,w}>0$. 
Then $x_{m,v}>0$ for every descendant of $w$, which implies
that the class of $x$ is not $0$. Therefore, we have $H\ne\{0\}$.
Since $W$ is strictly contained in $V$ and since it is hereditary,
we have $H\ne G$. Thus $G$ is not simple.
\end{proof}

\begin{example}
Consider again the non simple Bratteli diagram of Example~\ref{exampleDirectedHereditary}. All matrices $M(n)$ for $n\ge 1$ are equal to
\begin{displaymath}
M=\begin{bmatrix}1&1\\0&1\end{bmatrix}
\end{displaymath}
The dimension group $G$ is $\Z^2$
with the lexicographic order, that is
$(\Z^2,\Z^+\times \Z\cup \{0\}\times \Z_+,0\}$ 
(see Example~\ref{exampleStates}). The set $\{0\}\times\Z$
is an order ideal.
\end{example}
\subsection{Ordered Bratteli diagrams}
\begin{figure}[hbt]
\tikzset{node/.style={circle,draw,minimum size=0.1cm,inner sep=0pt}}
	\tikzset{title/.style={minimum size=0.5cm,inner sep=0pt}}
\begin{tikzpicture}
\node[node](v11)at(3,4){};\node[title]at(3,4.4){$A$};
\node[node](v12)at(5,4){};\node[title]at(5,4.4){$B$};
\node[title]at (8,4){$V(n-1)$};
\node[title]at (0,3){$E(n)$};
\node[node](v21)at(2,2){};\node[title]at(1.6,2){$0$};
\node[node](v22)at(3,2){};\node[title]at(2.6,2){$1$};
\node[node](v23)at(4,2){};\node[title]at(3.6,2){$2$};
\node[node](v24)at(5,2){};\node[title]at(4.6,2){$3$};
\node[node](v25)at(6,2){};\node[title]at(5.6,2){$4$};
\node[title]at(8,2){$V(n)$};
\node[title]at(0,1){$E(n+1)$};
\node[node](v31)at(3,0){};\node[title]at(2.6,0){$a$};
\node[node](v32)at(5,0){};\node[title]at(4.6,0){$b$};
\node[title]at(8,0){$V(n+1)$};

\draw[bend right,left](v11)edge node{$1$}(v21);
\draw[right,below](v11)edge node{$2$}(v21);
\draw[left,below](v11)edge node{$2$}(v23);
\draw[left,below](v11)edge node{}(v24);
\draw[right](v12)edge node{}(v22);
\draw[right,below](v12)edge node{$1$}(v23);
\draw[right](v12)edge node{}(v25);

\draw[left,below](v21)edge node{$1$}(v31);
\draw[right,below,right](v22)edge node{$2$}(v31);
\draw[left,below,left](v23)edge node{$3$}(v32);
\draw[left,below,left](v24)edge node{$1$}(v32);
\draw[right,below,right](v25)edge node{$2$}(v32);
\end{tikzpicture}
\caption{Order on the diagram of Figure~\ref{ch5:fig:ch5-diag1}}\label{ch5:fig:ch5-diag1ordre}
\end{figure}
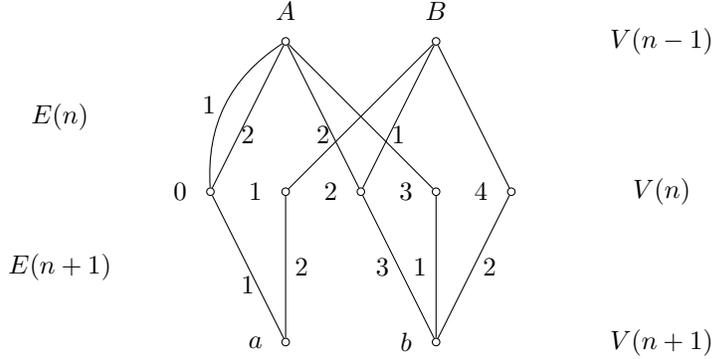

An 
{\em ordered Bratteli diagram}\index{subject}{Bratteli diagram!ordered}\index{subject}{ordered}\index{subject}{ordered!Bratteli diagram}
 $(V,E, \le )$\index{symbols}{V@$(V,E,\le)$}  is a Bratteli
diagram $(V,E)$ together with a partial order $\le$ on $E$ such that edges $e$,
$e'$
in $E$ are {\em  comparable}\index{subject}{comparable edges} if, and only if, $r(e) = r(e')$, in other words, we have
a linear order on each set $r^{-1} ( \{ v \} )$, where $v$ belongs to $V \setminus V(0)$ (see Figure~\ref{ch5:fig:ch5-diag1ordre}).

\medskip

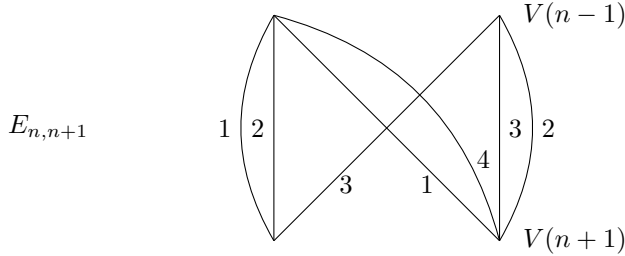
\begin{figure}[hbt]
\centering
\tikzset{node/.style={circle,draw,minimum size=0.1cm,inner sep=0pt}}
	\tikzset{title/.style={minimum size=0.5cm,inner sep=0pt}}
\begin{tikzpicture}
\draw (-1, -1.5) node {$E_{n,n+1} $} ;
\draw[left] (2,0) edge node{$2$}(2,-3) ;
\draw [bend right,left](2,0) edge node {$1$}(2,-3) ;
\draw [left,near end](2,0)edge node {$1$}(5,-3) ;
\draw [right,bend left,near end](2,0) edge node{$4$} (5,-3) ;
\draw [right,near end](5,0) edge node{$3$} (2,-3) ;

\draw [right](5,0)edge node{$3$} (5,-3) ;
\draw [bend left,right](5,0) edge node{$2$} (5,-3) ;

\draw (6,0) node {$V (n-1)$} ;
\draw (6,-3) node {$V (n+1)$} ;
\end{tikzpicture}
\caption{Telescoping of the diagram of Figure~\ref{ch5:fig:ch5-diag1ordre}.}
        \label{ch5:fig:ch5-diag2ordre}
\end{figure}

Note that if $(V, E, \le )$ is an ordered Bratteli diagram and $k < l$ in
${\Z}^+$, then the set $E_{k+1,l}$ of
paths from $V(k)$ to $V (l)$ may be given an induced ({\em lexicographic}) 
order\index{subject}{lexicographic!order}\index{subject}{order!lexicographic} as follows:
$$
(e_{k+1} , e_{k+2} , \ldots , e_l ) > (f_{k+1} , f_{k+2} , \ldots , f_l )
$$
if, and only if, for some $i$ with $k + 1 \leq i \leq l, \ e_j = f_j$
for $i < j \leq l$ and $e_i > f_i$.  
It is a simple observation that if
$(V, E, \le )$ is an ordered Bratteli diagram and $( V' , E' )$ is a
telescoping of $(V, E)$ as defined above, then with the induced order
$\le' , \; (V' , E' , \le' )$ is again an ordered Bratteli diagram.  
We say that
$(V' , E' , \le' )$ is a 
{\em telescoping}\index{subject}{telescoping!of ordered Bratteli diagram}\index{subject}{Bratteli diagram!telescoping}
 of $(V, E, \le )$ (see Figure~\ref{ch5:fig:ch5-diag2ordre}).

Again there is an obvious notion of isomorphism between ordered Bratteli
dia\-grams.
Let $\approx$ denote the equivalence relation on ordered Bratteli diagrams
generated by isomorphism and by telescoping.  
One can show that $G^1
\approx G^2$, where $G^1 = (V^1, E^1, \le^1), \; G^2 = (V^2, E^2,
\le^2)$, if, and only if, there exists an ordered Bratteli diagram $G
= (V,E, \le)$ such that telescoping $G$ to odd levels $0 < 1 < 3 <
\ldots$ yields a telescoping of either $G^1$ or $G^2$, and telescoping $G$
to even levels $0 < 2 < 4 < \ldots$ yields a telescoping of the
other (Exercise~\ref{exerciseEquivalenceTelescoping}).  This is analogous to the situation for the equivalence
relation $\sim$ on Bratteli diagrams as we discussed above.
%We write $G^1 \sim G^2$ to say $(V^1, E^1) \sim (V^2, E^2)$.
\begin{example}
The two first ordered diagrams of Figure~\ref{figureProperlyOrdered}
are equivalent. The third one is not
(although all three are equivalent as unordered Bratteli diagrams).
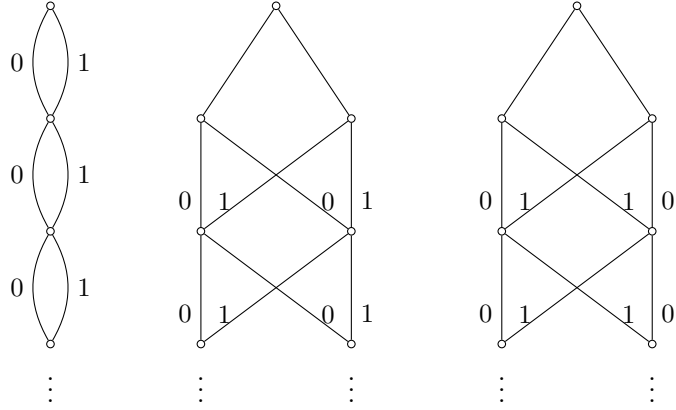
\begin{figure}[hbt]
\centering
\tikzset{node/.style={circle,draw,minimum size=0.1cm,inner sep=0pt}}
	\tikzset{title/.style={minimum size=0.5cm,inner sep=0pt}}
\begin{tikzpicture}
%odo
\node[node](0)at(1,4.5){};
\node[node](1)at(1,3){};
\node[node](2)at(1,1.5){};
\node[node](3)at(1,0){};
\node[title](4)at(1,-.5){$\vdots$};

\draw[bend right,left](0)edge node{$0$}(1);\draw[bend left,right](0)edge node{$1$}(1);
\draw[bend right,left](1)edge node{$0$}(2);\draw[bend left,right](1)edge node{$1$}(2);
\draw[bend right,left](2)edge node{$0$}(3);\draw[bend left,right](2)edge node{$1$}(3);
%odo2
\node[node](0)at(4,4.5){};
\node[node](11)at(3,3){};\node[node](12)at(5,3){};
\node[node](21)at(3,1.5){};\node[node](22)at(5,1.5){};
\node[node](31)at(3,0){};\node[node](32)at(5,0){};
\node[title](41)at(3,-.5){$\vdots$};\node[title](42)at(5,-.5){$\vdots$};

\draw(0)edge node{}(11);\draw(0)edge node{}(12);
\draw[left,near end](11)edge node{$0$}(21);\draw[right,near end](11)edge node{$0$}(22);
\draw[left,near end](12)edge node{$1$}(21);\draw[right,near end](12)edge node{$1$}(22);
\draw[left, near end](21)edge node{$0$}(31);\draw[right,near end](21)edge node{$0$}(32);
\draw[left,near end](22)edge node{$1$}(31);\draw[right,near end](22)edge node{$1$}(32);

%Morse
\node[node](0)at(8,4.5){};
\node[node](11)at(7,3){};\node[node](12)at(9,3){};
\node[node](21)at(7,1.5){};\node[node](22)at(9,1.5){};
\node[node](31)at(7,0){};\node[node](32)at(9,0){};
\node[title](41)at(7,-.5){$\vdots$};\node[title](42)at(9,-.5){$\vdots$};

\draw(0)edge node{}(11);\draw(0)edge node{}(12);
\draw[left,near end](11)edge node{$0$}(21);\draw[right,near end](11)edge node{$1$}(22);
\draw[left,near end](12)edge node{$1$}(21);\draw[right,near end](12)edge node{$0$}(22);
\draw[left, near end](21)edge node{$0$}(31);\draw[right,near end](21)edge node{$1$}(32);
\draw[left,near end](22)edge node{$1$}(31);\draw[right,near end](22)edge node{$0$}(32);

\end{tikzpicture}
\caption{Two properly ordered Bratteli diagrams and a non properly ordered one.}
\label{figureProperlyOrdered}
\end{figure}
\end{example}

The following notion will be important when we will deal with Bratteli diagrams and subshifts.
Fix $n\ge 1$ and let us consider $V(n-1)$ and $V (n)$ as alphabets.
For every letter $a\in V(n)$, consider the
ordered list $(e_1,\ldots ,e_k)$
of edges of $E(n)$ which range at $a$, 
and let $(a_1,\ldots ,a_k)$ be the ordered list of the labels of the sources
of these edges. 
This defines a morphism $\tau(n):a\mapsto a_1\cdots a_k$ from $V(n)^*$ to $V(n-1)^*$ we call {\em the morphism read on}\index{subject}{morphism!read on a Bratteli diagram}\index{subject}{Bratteli diagram!morphism read on}%
\index{subject}{substitution!read on Bratteli diagram} $E(n)$.
For example, in Figure~\ref{ch5:fig:ch5-diag1ordre}, the morphism we read on:
\begin{itemize}
\item
$E(n)$ is    $ \tau(n)    : \ 0\mapsto AA , \ 1 \mapsto B , \ 2 \mapsto BA , \ 3 \mapsto A , \ 4 \mapsto B,$
\item
$E(n+1)$ is  $ \tau(n+1) : \ a \mapsto 01 , \ b \mapsto 342 $,
\end{itemize}
and on Figure~\ref{ch5:fig:ch5-diag2ordre} the morphism we read on $E_{n,n+1}$ is $\sigma : a\mapsto AAB$, $b\mapsto ABBA$.  
We can check  of course  that we have $\sigma = \tau_{n} \circ \tau_{n+1}$.
Note that the matrix $M(n)$ is the incidence matrix of the morphism $\tau(n)$.
In our example, we find:
\begin{displaymath}
  M(n)=\begin{bmatrix}2&0\\0&1\\1&1\\0&1\end{bmatrix},\quad
  M(n+1)=\begin{bmatrix}1&1&0&0&0\\0&0&1&1&1\end{bmatrix}.
  \end{displaymath}
Note that the elements of the matrices $M(n)$ are integers (
since $M(n)_{r,s}$ is the number of edges from $s$ to $r$) 
whereas the label of an edge (which is also an integer)  does not represent
a multiplicity of the edge but its order with respect the
other edges with the same range.
%%%%%%%%%%%%%%%%%%%%%%%%%%%%%%%%%%%%%%%%%%%%%%%%%%%
\section{Dynamics for ordered Bratteli diagrams}\label{sectionDynamicsBratteli}
%%%%%%%%%%%%%%%%%%%%%%%%%%%%%%%%%%%%%%%%%%%%%%%%%%%
We shall see now how one can define a dynamics on the set of paths
in a Bratteli diagram.
\subsection{The Bratteli compactum}
Let $(V,E, \le )$ be an ordered Bratteli diagram.
Let $X_E$ denote the associated 
{\em infinite path space},\index{subject}{Bratteli diagram!infinite path}\index{subject}{path!in Bratteli diagram}\index{subject}{path!space}, that is,
$$
X_E = \{ (e_1 , e_2 , \ldots ) \mid  e_i \in E (i) , r (e_i ) =
s (e_{i+1} )  ,  i = 1, 2, \ldots \} .
$$
\index{symbols}{X@$X_E$}%
We exclude trivial cases and assume henceforth that $X_E$ is an infinite
set.  

Two paths in $X_E$ are said to be 
{\em cofinal}\index{subject}{cofinal paths}\index{subject}{path!cofinal}
if they have
the same tails, {\em i.e.}, the edges agree from a certain level on. 
We denote by $R_E$\index{symbols}{R@$R_E$} this equivalence, called the equivalence
of \emph{cofinality}
\index{subject}{cofinality equivalence}\index{subject}{equivalence!cofinality}%
on the set $X_E$. We will see that cofinality is an important notion
for Bratteli diagrams when we will define the dynamics on ordered
Bratteli diagrams.

The set $X_E$ is a closed subset of $\prod_{i\ge 1} E(i)$.
Since every $E(i)$ is finite, the product is compact and
thus $X_E$ is compact.
A basis for the topology is the family of {\em cylinder sets}\index{subject}{cylinder}
$$
[e_1 , e_2 , \ldots , e_k ]_E = \{ ( f_1 , f_2 , \ldots ) \in
X_E \mid  f_i = e_i , 1 \leq i \leq k \} \ .
$$

Each $[e_1 , \ldots , e_k ]_E$ is also closed, as is  easily seen.
When it will be clear from the context we will write $[e_1 , \ldots , e_k ]$ instead of $[e_1 , \ldots , e_k ]_E$.
Endowed with this topology, we call $X_E$ the
{\em Bratteli compactum}
\index{subject}{Bratteli compactum}\index{subject}{compactum}%
associated with $ (V , E , \le )$.
Let $d_E$ be the distance on $X_E$ defined by $d_E((e_n)_n , (f_n)_n) = \frac{1}{2^k}$ where $k=\inf \{ i  \mid  e_i \not = f_i \}$.
It clearly  defines  the topology of the cylinder sets.

If $(V,E)$ is a simple Bratteli diagram, then $X_E$ has no isolated
points, and so is a Cantor space
\index{subject}{Cantor!space}%
(recall that we assume $X_E$ to be infinite,
see Exercise~\ref{exerciseBratteliCantor}).
Moreover, each class of the equivalence $R_E$ (corresponding to cofinality)
is dense in $X_E$ (Exercise~\ref{exerciseCofinalDense}).

\medskip

Let $x = (e_1, e_2, \ldots)$ be an element of $X_E$.  We
will call $e_n$ the $n$th label of $x$ and denote it by $x(n)$.  We
let $X_E^{\max}$\index{symbols}{X@$X_E^{\max}$} denote those elements $x$ of $X_E$
\index{subject}{path!maximal}%
such that $x(n)$ is a maximal edge\index{subject}{edge!maximal} for all $n$ 
and $X_E^{\min}$\index{symbols}{X@$X_E^{\min}$} the analogous set for the minimal edges\index{subject}{edge!minimal}. 
\index{subject}{path!minimal}%

It is not difficult to show that $X_E^{\max}$ and
$X_E^{\min}$ are non-empty (see Exercise~\ref{exerciseMaxMin}).
Moreover, for every $v\in V$, the set of minimal edges
forms a \emph{spanning tree}\index{subject}{spanning tree}
of the graph $(V,E)$. This means that for every $v\in V$,
there is a unique path formed of minimal edges from $v$ to $v(0)$ (Exercise~\ref{exerciseMaxMin}). The same holds for maximal edges.

%\begin{definition}
The ordered Bratteli diagram $(V,E , \le)$ is {\em properly ordered}
\index{subject}{Bratteli diagram!properly ordered}%
\index{subject}{properly!ordered}\index{subject}{ordered!properly}%
if it is simple
    and if $X_E^{\max}$ and $X_E^{\min}$ both are a one point set:
    $X_E^{\max} = \{ x_{\max} \}$\index{symbols}{x@$x_{\max}$}
 and $X_E^{\min} = \{ x_{\min} \}$.\index{symbols}{x@$x_{\min}$}
 %\end{definition}

\begin{example}\label{exampleProperlyOrdered}
The Bratteli diagrams of Figure~\ref{figureProperlyOrdered} on the left 
and center are properly ordered
while the diagram on the right is not. Indeed, there are two
paths labeled with $0,0,0,\ldots$ and two paths labeled
$1,1,1,\ldots$.
\end{example}
Note that every simple Bratteli diagram can be properly
ordered (Exercise~\ref{exerciseProperlyOrdered}).

\subsection{The Vershik map}
We can now define, for a properly ordered Bratteli diagram $(V,E,\le)$,
 a map $T_E:X_E \rightarrow X_E$,
\index{symbols}{T@$T_E$}%
called the {\em Vershik map}\index{subject}{map!Vershik}\index{subject}{Vershik map}\index{names}{Vershik, Anatol M.} (or the {\em lexicographic map}\index{subject}{lexicographic!map}\index{subject}{map!lexicographic}),
associated with $(V,E, \le)$.

We set $T_E (x_{\max}) = x_{\min}$.  If $ x = (e_1, e_2, \ldots) \neq
x_{\max}$, let $k$ be the least integer such that $e_k$ is not a
maximal edge.  Let $f_k$ be the successor
\index{subject}{successor!of edge} of $e_k$ (relative to the
order $\le$ so that $r(e_k) =
r(f_k))$.  Define
\begin{displaymath}
  T_E (x) =  (f_1, \ldots , f_{k-1}, f_k, e_{k+1},
  e_{k+2} , \ldots),
\end{displaymath}
where $(f_1, \ldots, f_{k-1})$ is the
minimal edge in  $E_{1,k-1}$ with
range equal to $s(f_k)$.

Thus, the image by $T_E$ of a point $x\ne x_{\max}$ is its successor
in the lexicographic order. Moreover, $x$ and $T_E(x)$ are clearly cofinal.

The map $T_E$ is clearly continuous. It is moreover one-to-one (Exercise~\ref{exerciseBVhomeo}).
We call the resulting pair $(X_E, T_E)$\index{symbols}{V@$(V_E,T_E)$} a
{\em Bratteli-Vershik dynamical system}\index{subject}{Bratteli-Vershik!dynamical system}\index{subject}{dynamical system!Bratteli-Vershik}.
Since $X_E$ is a Cantor space, it
is a Cantor dynamical system\index{subject}{dynamical system!Cantor}\index{subject}{Cantor!dynamical system}.
\begin{proposition}\label{propositionBratteliCompactum}
  Let $(V,E,\le)$ be a properly ordered Bratteli diagram.
  The system $(X_E,T_E)$ is a minimal
  Cantor dynamical system.
\end{proposition}
The proof is left as an exercise (Exercise~\ref{exerciseBVMinimal}).
We will see below with Theorem \ref{ch5:theo:BVmodel} that the converse is also true.

In the sequel BV\index{subject}{Bratteli-Vershik!BV} will refer to {\em Bratteli-Vershik}.

\begin{example}\label{exampleBV1}
Consider the Bratteli diagram $(V,E)$ of Figure~\ref{figureProperlyOrdered}
on the left. The system $(X_E,T_E)$ is  isomorphic to the
odometer $(\Z_2,T)$
where $T(x)=x+1$. Indeed, as well-known, the addition of $1$
in base $2$ consists, on the representation
in base $2$ of numbers with a fixed number of digits, in taking the next sequence in the lexicographic order. Since the representation in base
$2$ of $2$-adic numbers is written with the least
significant digit on the left, addition of $1$
corresponds to the next element in the reverse
lexicographic order for right infinite sequences (see Figure \ref{figureAddition} below).
\begin{figure}[hbt]
\begin{displaymath}
\begin{array}{cccccc}
0&0&0&0&0&\cdots\\
\textcolor{red}{1}&0&0&0&0&\cdots\\
0&\textcolor{red}{1}&0&0&0&\cdots\\
\textcolor{red}{1}&1&0&0&0&\cdots\\
0&0&\textcolor{red}{1}&0&0&\cdots\\
\textcolor{red}{1}&0&1&0&0&\cdots\\
0&\textcolor{red}{1}&1&0&0&\cdots\\
\textcolor{red}{1}&1&1&0&0&\cdots
\end{array}
\end{displaymath}
\caption{The addition of $1$ in base $2$.}\label{figureAddition}
\end{figure}
\end{example}

%%%%%%%%%%%%%%%%%%%%%%%%%%%%%%%%%%%%%%%%%%%%%%%%%%%%%%%%%%%%%%%%%%%%%%%
\section{The Bratteli-Vershik model theorem}\label{sectionBVmodelTheorem}
%%%%%%%%%%%%%%%%%%%%%%%%%%%%%%%%%%%%%%%%%%%%%%%%%%%%%%%%%%%%%%%%%%%%%%%

Let $(X,T)$  be an invertible minimal Cantor dyna\-mical sys\-tem.
The properly ordered Bratteli
diagram $(V,E, \le)$ is a {\em BV-representation}\index{subject}{representation!BV} of $(X,T)$ if $(X_E,
T_E)$ is isomorphic to $(X,T)$.
We will show, as a main result of this chapter, that every  invertible
minimal Cantor system
has a BV-representation.
\subsection{From partitions in towers to Bratteli diagrams}
We will first show how to associate to any nested sequence of
partitions of an invertible system $(X,T)$ a
sequence of ordered Bratteli diagrams.

Let
\begin{displaymath}
\Pg(n)=\{T^j B_i(n) \mid  0\leq j< h_i(n), 1\leq i\leq t(n)\}
\end{displaymath}
 be a nested sequence of KR-partitions of $(X,T)$.
We may suppose that $\Pg (0) = \{ X \}$. Hence $t(0)=1$, $h_1(0)=1$ and $B_1 (0) = X$.

Let $V(n) = \{ (n,1) , \ldots , (n,t(n)) \}$, for $n\ge 0$. Thus the set
of vertices is, at each level $n$, the set  of towers of the partition $\Pg(n)$.

The set of edges records the inclusions of the
elements of $\Pg(n)$ in the elements of $\Pg(n-1)$. Specifically,
let $E(n)$ be the set  
of quadruples $(n,t',t,j)$ satisfying  
\begin{align}
\label{ch5:inclusion-def-DB}
T^{j} B_{t} (n) \subseteq B_{t'}(n-1) 
\end{align}
for $1\leq t' \leq  t(n-1)$, $1\leq t \leq t(n)$, $0\leq j\leq h_{t}(n)-1$ and $n\geq 1$. Note that, in particular, 
\begin{enumerate}
\item the index $j$ is such that $T^jB_t(n)$ is contained in the basis $B(n-1)$ of $\Pg(n-1)$,
\item the index $j$ is the return time of every element of $B_t(n)$
to $B(n-1)$.
\item  not all indices $j$ with $0\le j\le h_t(n)-1$ appear in these quadruples.
\end{enumerate}
The range and source maps are given by
\begin{align}
\label{ch5:def:sourcerange}
s (n,t',t,j) = (n-1,t') \hbox{ and } r(n,t',t,j)= (n,t) \ .
\end{align}

Two edges $e_1=(n_1,t'_1,t_1,j_1)$ and $e_2=(n_2,t'_2,t_2,j_2)$ are comparable whenever 
$n_1=n_2$ and $t_1=t_2$. 
In this case we define $e_1\geq e_2$ if $j_1\geq j_2$.
It is straightforward to verify  that $(V,E, \le )$ is an ordered Bratteli diagram.

It is useful to remark, from \eqref{ch5:def:sourcerange}, that $((n , t'_n ,t_n ,j_n))_n$ is an infinite path of 
$(V,E, \le )$ if, and only if, $t_{n-1}=t'_{n}$ for all $n\geq 1$.
Hence the paths of the Bratteli diagram have the form $((n , t_{n-1} ,t_{n} ,j_n))_n$ with $1\le t_{n-1}\le t(n-1)$, $1\le t_n\le t(n)$, $0\le j_n\le h_{t_n}(n)$
and
\begin{equation}
T^{j_n}B_{t_n}(n)\subset B_{t_{n-1}}(n-1)\label{eqDefinitionDB}
\end{equation}
Note that \eqref{eqDefinitionDB} implies (by Exercise~\ref{exerciseRefinement}) that
\begin{displaymath}
0\le j_n\le h_{t_n}(n) -h_{t_{n-1}}(n-1).
\end{displaymath}
Note also that $(n , t_{n-1} ,t_{n} ,j_n)$ is a minimal edge if and only if $j_n = 0$ and is maximal if and only if 
\begin{equation}
j_n=h_{t_n}(n) -h_{t_{n-1}}(n-1).\label{equationMaximalEdge}
\end{equation}
Additionally, if $(n,t_{n-1},t_n,j_n)$ is not a maximal edge, its
successor is an edge $(n,t'_{n-1},t_n,j'_n)$ with
\begin{equation}
j'_n=j_n+h_{t_{n-1}}(n-1)\label{eqSuccesseur}
\end{equation} since $j'_n$ is the least
integer such that $T^{j'_n-j_n}B_{t_n}(n)\subset B(n-1)$.

\begin{figure}[ht]
  \centering
   
\tikzset{node/.style={draw,minimum width=.6cm,minimum height=.3cm}}
\tikzset{title/.style={minimum size=0cm,inner sep=0pt}}
\begin{tikzpicture}
%premiere partition
\node[node](11)at(0,0){};\node[node](12)at(0,0.3){};
\node[node](21)at(1.5,0){};
\node[node](22)at(1.5,0.3){};\node[node](23)at(1.5,.6){};
\node[title]at(0,-.5){$A$};\node[title]at(1.5,-.5){$B$};

\node[title](11m)at(0,0){};\node[title](12m)at(0,.3){};
\draw[->](11m)edge node{}(12m);
\node[title](21m)at(1.5,0){};\node[title](22m)at(1.5,.3){};
\draw[->](21m)edge node{}(22m);
%deuxieme partition premiere tour
\node[node](11)at(4,0){};\node[node](12)at(4,.3){};
\draw[line width=1.5pt](3.7,.45)edge node{}(4.3,.45);
%\node[minimum width=.6cm,minimum height=.56cm,draw]at(5,.15){};
\node[title]at(3.3,.15){$A$};
\node[node](13)at(4,.6){};\node[node](14)at(4,.9){};
%\node[minimum width=.6cm,minimum height=.58cm,draw]at(5,.75){};
\draw[line width=1.5pt](3.7,1.05)edge node{}(4.3,1.05);
\node[title]at(3.3,.75){$A$};
\node[node](15)at(4,1.2){};\node[node](16)at(4,1.5){};
\node[node](17)at(4,1.8){};
\draw[line width=1.5pt](3.7,1.95)edge node{}(4.3,1.95);
\node[title]at(3.3,1.5){$B$};
\node[node](18)at(4,2.1){};\node[node](19)at(4,2.4){};
\node[title]at(3.3,2.25){$A$};

\node[title](211m)at(4,0){};\node[title](212m)at(4,.3){};
\draw[->](211m)edge node{}(212m);

%deuxieme partition deuxieme tour
\node[node](21)at(6,0){};\node[node](22)at(6,.3){};
\draw[line width=1.5pt](5.7,.45)edge node{}(6.3,.45);
\node[title]at(5.3,.15){$A$};
\node[node](23)at(6,.6){}; \node[node](24)at(6,.9){}; 
\node[title]at(5.3,.75){$A$};
\node[title](221m)at(6,0){};\node[title](212m)at(6,.3){};
\draw[->](221m)edge node{}(212m);

%trosieme tour
\node[node](31)at(8,0){};\node[node](32)at(8,.3){};\node[node](33)at(8,.6){};
\node[title]at(7.3,.3){$B$};
\node[title](31m)at(8,0){};\node[title](32m)at(8,.3){};
\draw[->](31m)edge node{}(32m);
\draw[line width=1.5pt](7.7,.75)edge node{}(8.3,.75);
\node[node](21)at(8,.9){};\node[node](22)at(8,1.2){};
%\draw[line width=1.5pt](7.7,1.05)edge node{}(8.3,1.05);
\node[title]at(7.3,1.05){$A$};
\node[node](21)at(8,1.5){};\node[node](22)at(8,1.8){};
\draw[line width=1.5pt](7.7,1.35)edge node{}(8.3,1.35);
\node[title]at(7.3,1.65){$A$};
\end{tikzpicture}          
  \caption{The partition ${\Pg} (1)$ consists of two towers called $A$ and $B$. 
        The dynamics $T$ acts vertically except for the last levels where it goes back to the base.
        The elements of the partition ${\Pg} (2)$ can be seen as the piling up of vertical pieces of the towers $A$ and $B$.} 
        \label{ch5:fig:ch5-towerexemple}
\end{figure}
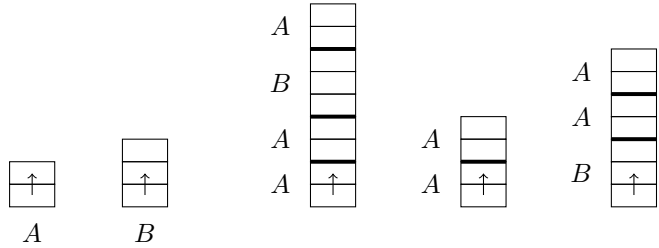

For example suppose  that $\Pg (n)$ is a refining sequence of KR-partitions  such that (see Figure~\ref{ch5:fig:ch5-towerexemple})
\begin{enumerate}
\item
${\Pg} (1) = \{ A(1) , T A(1) , B(1) , T B(1)  , T^2 B(1)  \}$ and
\item
${\Pg} (2) = \left\{ T^j B_i(2) \mid  0\leq j< h_i(2), 1\leq i\leq t(2) \right\}$ with 
\begin{enumerate}
\item $t(2) = 3$, $h_1 (2) = 9$, $h_2 (2) = 4$, $h_3 (2) = 7$,
\item $B_{1} (2) \subseteq A(1)$, $T^{2} B_{1} (2) \subseteq A(1)$, $T^{4} B_{1} (2) \subseteq B(1)$, $T^{7} B_{1} (2) \subseteq A(1)$,
\item $B_{2} (2) \subseteq A(1)$, $T^{2} B_{2} (2) \subseteq A(1)$,
\item $B_{3} (2) \subseteq B(1)$, $T^{3} B_{3} (2) \subseteq A(1)$, $T^{5} B_{3} (2) \subseteq A(1)$.
\end{enumerate}
\end{enumerate}
Note that every element of $\Pg(1)$ is a union of elements of $\Pg(2)$. In
particular, we have
\begin{displaymath}
  B(1)=T^4B_1(2)\cup B_3(2).
  \end{displaymath}

The corresponding Bratteli diagram is represented in Figure~\ref{ch5:fig:ch5-diagexemple}
with the value of $j$ indicated on the edge $(n,t',t,j)$.
Note that only the order of these indices matter for the
Bratteli diagram,
and not their particular value.

\begin{figure}[hbt]
\centering
\tikzset{node/.style={circle,draw,minimum size=0.1cm,inner sep=0pt}}
	\tikzset{title/.style={minimum size=0.5cm,inner sep=0pt}}
\begin{tikzpicture}
\node[node](v0)at(4,4){};\node[title]at(7,4){$V(0)$};
\node[title]at(0,3){$E(1)$};
\node[node](v11)at(2,2){};\node[node](v12)at(5,2){};\node[title]at(7,2){$V(1)$};
\node[title]at(0,1){$E(2)$};
\node[node](v21)at(1,0){};\node[node](v22)at(3,0){};\node[node](v23)at(6,0){};
\node[title]at(7,0){$V(2)$};

\draw[bend right,left](v0)edge node{$0$}(v11);\draw[right](v0)edge node{$1$}(v11);
\draw[right](v0)edge node{$1$}(v12);\draw[bend left,right](v0)edge node{$2$}(v12);
\draw[bend right,left](v0)edge node{$0$}(v12);
\draw[bend right=45,left,near end](v11)edge node{$0$}(v21);
\draw[left,near end](v11)edge node{$2$}(v21);
\draw[bend left,right](v11)edge node{$7$}(v21);
\draw[left,near end](v11)edge node{$0$}(v22);
\draw[bend left,right,near end](v11)edge node{$2$}(v22);
\draw[left,near end](v11)edge node{$3$}(v23);
\draw[bend left,right,near end](v11)edge node{$5$}(v23);
\draw[right,near end](v12)edge node{$4$}(v21);
\draw[right,near end](v12)edge node{$0$}(v23);
\end{tikzpicture}
\caption{Diagrammatic representation of ${\Pg} (0)$, ${\Pg} (1)$ and ${\Pg} (2)$.}
        \label{ch5:fig:ch5-diagexemple}
\end{figure}
The next two propositions give properties of the ordered Bratteli diagram
$(V,E,\le)$
associated as above with the nested sequence of partitions $(\Pg(n))$.
The first one relates the paths in the diagram and the nested
sequence of partitions.

\begin{proposition}
\label{ch5:lem:equivpathpart}
For every infinite path $(e_1,e_2,\ldots)$ with $e_n=(n , t_{n-1} ,t_{n} ,j_n)$ in $(V,E, \le )$,
the following assertions hold.
\begin{enumerate}
\item[\rm(i)] We have $0\leq \sum_{i=1}^n j_i \leq h_{t_n} (n) -1$ for every $n\ge 1$.
\item[\rm(ii)] For $1\le m\le n$, the waiting time for an element
of $B_{t_n}(n)$ to access $B_{t_m}(m)$ is $\sum_{i=m}^nj_i$.
\item[\rm(iii)]
The sequence 
\begin{equation}
C_n=T^{\sum_{i=1}^nj_i}B_{t_n}(n)\label{eqCn}
\end{equation}
is decreasing, that is
\begin{displaymath}
C_1\supset C_2\supset\ldots\supset C_n\supset\ldots
\end{displaymath}
and the map $\phi:(e_n)\to \cap_{n\ge 1}C_n$ sends distinct paths
to disjoint sets.
\item[\rm(iv)]
If the sequence $\Pg(n)$ generates the topology of $X$, then
$\Card(\cap_{n\geq 1} C_n) = 1$.
\end{enumerate}
\end{proposition}
\begin{proof}
(i) Let us prove the first assertion by induction on $n$. It is true for $n=1$
since $j_1\le h_{t_1}(1)-1$. Next, for $n\ge 2$, since 
$T^{j_n}B_{t_n}(n)\subset B_{t_{n-1}}(n-1)$, we have
$j_n\le h_{t_n}(n)-h_{t_{n-1}}(n-1)$ (see Exercise~\ref{exerciseRefinement}).
Since, by induction hypothesis, we have $j_{n-1}+\ldots+j_1\le h_{t_{n-1}}(n-1)$
we conclude that $j_n+\cdots+j_1\le h_{t_n}(n)$.

(ii) The statement is true for $n=m$. Next, we argue by induction on $n-m$.
Since $j_{n+1}$ is the return time from $B_{t_{n+1}}(n+1)$ to $B_{t_n}(n)$
and since, by induction hypothesis $\sum_{i=m}^n j_i$ is the return time
from $B_{t_n}(n)$ to $B_{t_m}(m)$, the conclusion follows.

(iii) By~\eqref{eqDefinitionDB}, we have for every $n\ge 1$
\begin{eqnarray*}
C_{n+1}&=&T^{\sum_{i=1}^{n+1} j_i} B_{t_{n+1}}(n+1)=T^{\sum_{i=1}^{n}j_i}T^{j_{n+1}}B_{t_{n+1}}(n+1)\\
&\subset&T^{\sum_{i=1}^{n}j_i}B_{t_n}(n)=C_n.
\end{eqnarray*}
Since $(C_n)_{n\ge 1}$ is a
decreasing sequence of closed sets and 
since $X$ is compact, its intersection is nonempty.
Since $0\leq \sum_{i=1}^n j_i \leq h_{t_n} (n) -1$ for every $n\ge 1$,
 $C_n$ is an element of the partition $\Pg(n)$. Let
 $(e_n)$ and $(e'_n)$ be distinct paths and let $C_n,C'_n$
be the associated sequences. We have $e_k\ne e'_k$ 
for some $k\ge 1$. Since $C_m$ and $C'_m$ are  distinct elements
of the partition $\Pg(m)$, they are disjoint. Thus $\phi(e_n)\cap\phi(e'_n)=\emptyset$.

(iv) If the sequence of partitions generates the topology,
the intersection of all $C_n$ is reduced to one point.
\end{proof}
In the second statement, we assume an hypothesis on $(\Pg(n))_{n\ge 1}$ mildly weaker than  that
of being a refining sequence (recall that there exist nested sequences of partitions
satisfying (KR1) but not (KR3), see Exercise~\ref{exerciseKR1pasKR3}).
\begin{proposition}\label{lemmaRefiningProperlyOrdered}
If  $(\Pg(n))_{n\ge 1}$ is a nested sequence satisfying (KR1), the Bratteli diagram
$(V,E,\le)$ is properly ordered.
\end{proposition}
\begin{proof}
First, the diagram $(V,E)$ is simple. Indeed, it is enough to prove
that there is an $n\ge 1$ such that there is a path
from every vertex at level $n$ to every vertex at level $1$.
For this, it enough to take $n$ such that all the $h_t(n)$
are larger than the maximum of the waiting times to access
an element of the partition $B(1)$.

Let us  show now that $X_E^{\rm min} $ consists of a single path. 
Let $(e_n)$ with $e_n=(n , t_{n-1} ,t_{n} ,j_n)$ be an infinite path of $X_E^{\rm min}$.
The edges comparable to $e_n$ are the edges of the form $(n,t,t_{n},j)$ for some $t$
and exactly one of them is of the form $(n,t,t_{n},0)$.
It is clearly a minimal edge. 
Hence $j_n = 0$ for all $n$.
But, by Lemma~\ref{ch5:lem:equivpathpart}
$$
\phi(e_n)=\cap_n T^0 B_{t_n} (n) \subseteq \cap_n B(n) 
$$
which consists of a single point by hypothesis.
Since $\phi$ sends distinct paths to disjoint sets, 
the path $(e_n)$ is the unique path of $X_E^{\rm min}$.

Similarly, if $(e_n)\in X_E^{\max}$, then $j_n=h_{t_n}(n)-h_{t_{n-1}(n-1)}$
by \eqref{equationMaximalEdge}. Then 
\begin{displaymath}
\phi(e_n)=\cap_n T^{h_{t_n}(n)-1}B_{t_n}(n)
\end{displaymath}
which is reduced to one point because $T(\phi(_n))$ belongs to the
intersection of the bases, whence the conclusion again.
\end{proof}
\subsection{The BV-representation theorem}
We can now state and prove the BV-representation theorem.
\begin{theorem}[Herman, Putnam, Skau]\label{ch5:theo:BVmodel}
\index{names}{Herman, Richard H.}\index{names}{Putnam, Ian F.}\index{names}{Skau, Christian F.}
For every invertible minimal Cantor system $(X,T)$,
there exists a properly ordered  Bratteli
diagram $(V,E, \le)$ such that $(X,T)$ is isomorphic to $(X_E,
T_E)$.

More precisely, for every refining sequence $(\Pg(n))$
of KR-partitions of $(X,T)$,
the Bratteli diagram $(V,E,\le)$ associated with $\Pg(n)$ is such that
$(X,T)$ is isomorphic to $(X_E,T_E)$.
\end{theorem}

\begin{proof}
By Theorem~\ref{theoremKRPartitions}, there exists a  be refining sequence of partitions $(\Pg (n))$
of $(X,T)$. By Proposition~\ref{lemmaRefiningProperlyOrdered},
the ordered Bratteli diagram $(V,E,\le )$  associated to $(\Pg(n))$
is properly ordered. This allows to consider
the  Cantor system $(X_E,T_E)$. It is minimal
by Proposition \ref{propositionBratteliCompactum}.

Consider the map $\phi : X_E \to X$ defined by 
\begin{equation}
\phi ((n , t_{n-1} ,t_{n} ,j_n)_n) =x \hbox{ where } \{x \}=\cap_{n\geq 1}  C_n
\label{eqDefinitionPhi}
\end{equation}
with $C_n=T^{\sum_{i=1}^n j_i} B_{t_{n}} (n)$.
It is well defined (Proposition~\ref{ch5:lem:equivpathpart}) and is a homeomorphism (see Exercise~\ref{ch5:ex:BVmodel}).
Note that $(\Pg(n))_n$ being a decreasing sequence of partitions, we also have $\{x \}=\cap_{n\geq N} C_n$ for all $N$.

There remains to show  that it commutes with the  dynamics.
Let $e=(e_n)_n$ be an infinite path of $X_E$ with $e_n = (n , t_{n-1} ,t_{n} ,j_n)$.
Suppose first that  $e$ is not the maximal path.
Then there exists $n_0$ such that $T_E (e) = e'_1 \cdots e'_{n_0-1}e'_{n_0}  e_{n_0+1} e_{n_0+2}\cdots $ where 
$e'_{n}=(n , t'_{n-1} ,t'_{n} ,j'_n)$, with $j'_n = 0$, $1\leq n \leq n_0-1$ and $e'_{n_0} = (n_0 , t'_{n_0-1} ,t_{n_0} ,j'_{n_0})$ is the successor of $e_{n_0}$.
Note that, for $1\leq n \leq n_0-1$, the edges $e_n$ being maximal we have $j_n =h_{t_n}(n) -h_{t_{n-1}}(n-1)$.

Since $(e'_{n_0})$ is the successor of $e_{n_0}$,
we have by~\eqref{eqSuccesseur}, $j'_{n_0} = j_{n_0} +h_{t_{n_0-1}} (n_0-1)$.
Hence 
$$
\sum_{1\leq n\leq n_0} j_{n} =j_{n_0} +  \sum_{1\leq n\leq n_0} h_{t_n}(n) -h_{t_{n-1}}(n-1) = j_{n_0} +h_{t_{n_0-1}} (n_0-1)-1 \ .
$$
while
\begin{eqnarray*}
\sum_{1\leq n\leq n_0} j'_{n} &=& j_{n_0} +h_{t_{n_0-1}}(n_0-1)\\
&=&\sum_{1\leq n\leq n_0} j_{n}+1.
\end{eqnarray*}

Set $C'_n=\cap_{n}T^{\sum_{i=1}^nj'_i}B_{t'_n}(n)$. Then
$\phi\circ T_E(e)=\phi(e')=\cap_nC'_n$. Since $\sum_{i=1}^{n_0}j'_i=
\sum_{i=1}^{n_0}j_i+1$
and $j'_n=j_n$ for $n>n_0$, we have $C'_n=TC_n$ for $n\ge n_0$.
This shows that $\phi\circ T_E(e)=T\phi(e)$ and thus the conclusion.

Suppose now  that $e$ is the maximal path.
Let $x_{\min }$ be the minimal path of $(V,E,\le)$. 
Then we have to prove that $\phi (x_{\min } ) = T (\phi (e))$.
But since ${\Pg} (0) = \{ X \}$, we have $h_{t_0} (0) = 1$ and consequently
\begin{align*}
T (\phi (e)) = & T \left( \bigcap_{n\geq 1} T^{\sum_{i=1}^{n } h_{t_i}(i) -h_{t_{i-1}}(i-1)} B_{t_{n}} (n)\right) \\
= & \bigcap_{n\geq 1} T^{h_{t_n}(n)} B_{t_{n}} (n) \subseteq \bigcap_{n\geq 1} \bigcup_{1\leq i\leq t(n)}B_{i} (n) = \{\phi (x_{\min }) \} \ .
\end{align*}
%Moreover, any contraction of $G$ yields a BV-representation of $(X,T)$ and some of  them have incidence matrices with positive entries. This
%shows that $G$ is simple.
\end{proof}

%As an illustration of Theorem~\ref{ch5:theo:BVmodel}, consider
%a strict episturmian word $s$ 
%\index{subject}{episturmian word} on the alphabet $A$. Let
%$x=a_1a_2\cdots$ be its directive sequence and let $X$ be the
%shift generated by $s$. Recall that for $a\in A$,
%we denote by $L_a$ the morphism which places $a$ before each
%letter $b\ne a$. Let $W_n=\RR_X(s_{[1,n]})$. 
%Recall from Section~\ref{sectionSturmianShifts}
%that  the words $u_n=\Pal(a_1\cdots a_{n})$
%are the palindrome prefixes of $s$ and the set
%of left return words to $u_n$ is
%\begin{displaymath}
%\RR'_X(u_n)=\{L_{a_1\cdots a_{n}}(a)\mid a\in A\}
%\end{displaymath}
%\begin{proposition}\label{propositionBVrepresentationEpisturmian}
%The ordered Bratteli diagram $(V,E)$ such that $L_{a_n}$
%is the sequence of morphisms read on $(V,E)$ is a BV-representation
%of the minimal shift $X$.
%\end{proposition}
%\begin{proof}
%Set $\alpha_n=|s_n|$ for $n\ge 1$ and $r_n(a)=L_{a_1\ldots a_n}(a)$
%for $a\in A$. We have
%\begin{equation}
%r_n(a)=\begin{cases}r_{n-1}(a_n)r_{n-1}(a)&\mbox{ if $a\ne a_n$}\\
%r_{n-1}(a_n)&\mbox{ otherwise}\end{cases}\label{equationr_n(a)}
%\end{equation}
%Let $\Pg(n)$ be the partition in towers ass
%\end{proof}

%We illustrate Theorem~\ref{ch5:theo:BVmodel} with the following
%example of a minimal shift space.
%\begin{example}

%\end{example}

\subsection{From Bratteli diagrams to partitions in towers}
\label{subsectionFromBDtoPT}
We now show that the construction used in Theorem~\ref{ch5:theo:BVmodel}
can done backwards, using the fact that the system $(X_E,T_E)$
has a natural sequence of partitions.
\begin{proposition}
  Let $(V,E,\le)$ be a properly ordered Bratteli diagram.
  For every $n\ge 1$, the family
  \begin{equation}
    \Pg(n)=\{[e_1,\ldots,e_n]\mid (e_1,\ldots,e_n)\in E_{1,n}\}\label{eqNaturalSequence}
  \end{equation}
  is a partition in towers of $(X_E,T_E)$ with basis the cylinders
  $[e_1,\ldots,e_n]$ formed of minimal edges and the sequence
  $(\Pg(n))$ is a refining sequence.
  \end{proposition}
%The system
%$(X_E,T_E)$  has a natural sequence of KR-partitions $(\Pg(n))$,
%\index{subject}{sequence!of partitions associated to a Bratteli diagram}%
%\index{subject}{Bratteli diagram!sequence of partitions associated to}%
%namely the partitions $[e_1,\ldots,e_n]$ for $(e_1,\ldots,e_n)\in E_{1,n}$.

\begin{proof}
The basis of the partition $\Pg(n)$ is the union of
the cylinders corresponding to minimal paths and can be identified
with $V(n)$. 
This sequence of partitions is a refining sequence.
Indeed, condition (KR1) is satisfied because $(V,E,\le)$
is properly ordered, and thus that there is a unique minimal
element in $X_E$. It is clear that the sequence $(\Pg(n))$
is nested and thus (KR2)  is satisfied. Finally, (KR3)
is also satisfied, by definition of the topology on $X_E$.
\end{proof}
The sequence of partitions $(\Pg(n))$ defined by \eqref{eqNaturalSequence}
is called the \emph{natural sequence of partitions}
associated to $(V,E,\le)$. 
\begin{proposition}\label{lemmaFormulaHeights}
  Each $\Pg(n)$ of the natural sequence of partitions is a partition in
$t(n)=\Card(V(n))$ towers.
  For each $v\in V(n)$, the height of tower $B_v(n)$ is
  \begin{equation}
h_v(n)=(M(n)\cdots M(1))_{1,v}.\label{eqh_t(n)}
\end{equation}
\end{proposition}
\begin{proof}
Since the number of paths from
$V(0)=\{1\}$ to $v\in V(n)$ is equal to $(M(n)\cdots M(1))_{1,v}$,
the height of the tower $B_v(n)$ with $v\in V(n)$ is given by \eqref{eqh_t(n)}.
\end{proof}
Note that $M(n)\cdots M(1)$  is actually a column vector so
that the right-hand side could also be denoted $(M(n)\cdots M(1))_v$.

For every minimal invertible Cantor system $(X,T)$,
the Bratteli diagram $G=(V,E)$ built from a refining sequence
$(\Pg(n)$ of KR-partitions
 is such that
$(X_E,T_E)$ is, by Theorem~\ref{ch5:theo:BVmodel}, conjugate to $(X,T)$. 
It is easy to verify that the natural sequence of partitions
and the original partition
coincide up to the conjugacy (Exercise~\ref{exerciseAllerRetour}).

The above correspondance between refining sequences of partitions
and properly ordered Bratteli diagrams invites to a translation
of the properties of each other. Thus,  telescoping a diagram
corresponds to taking a subsequence of the refining sequence.
We will see next another translation with the adjacency matrices
of a diagram.

\subsection{Dimension groups and BV-representation}
We have showed in Proposition~\ref{propositionGroupSequencePartitions}
 that for any minimal invertible Cantor
system $(X,T)$,  the dimension group $K^0(X,T)$ is a direct
limit of groups $G(n)$ associated with a sequence $(\Pg(n))$
of $KR$-partitions. We will now see how this group is
defined directly in terms of a BV representation of $(X,T)$.
\index{subject}{dimension group!of BV-dynamical system}%
\index{subject}{Bratteli-Vershik!dynamical system!dimension group}%

\begin{theorem}\label{theoremDGBratteliDiagram}
Let $(V,E,\le)$ be a properly ordered Bratteli diagram.
The group $K^0(X_E,T_E)$ is the dimension group $\DG(V,E)$
of the diagram $(V,E)$.
\end{theorem}
This relies on the important fact that,
as we shall see now, the adjacency matrices
$M(n)$ are precisely the connecting matrices $M(n+1,n)$
\index{subject}{connecting!matrix}%
introduced in Section~\ref{sectionOrderedGroupSequences}.
\begin{lemma}\label{lemmaConnectingMatrices}
  Let $(V,E,\le)$ be a properly ordered Bratteli diagram
  with adjacency matrices $M(n)$. Let $\Pg(n)$
  be the natural sequence of partitions associated to $(V,E,\le)$.
  Let $G(n)=G(\Pg(n))$ be the ordered group associated to the partition $\Pg(n)$
  and let $M(n+1,n)$ be the matrix of the connecting morphism $I(n+1,n)$.
  Then
  \begin{equation}
    M(n,n-1)=M(n)\label{eqConnectingMatrices}
    \end{equation}
\end{lemma}
\begin{proof}
  One has by Equation~\eqref{equationMatrixPartitions}
\begin{displaymath}
M(n,n-1)_{k,i}=
\Card\{j\mid  0\le j\le h_{k}(n), T_E^jB_{k}(n)\subset B_{i}(n-1)\}.
\end{displaymath}
for $1\le i\le t(n)$ and $0\le j <h_i(n)$.
Thus, by \eqref{eqDefinitionDB}, we have also
\begin{eqnarray*}
M(n,n-1)_{k,i}&=&\Card\{e\in E(n)\mid r(e)=k, s(e)=i\}\\
&=&M(n)_{k,i}.
\end{eqnarray*}
  \end{proof}

\begin{proofof}{of Theorem~\ref{theoremDGBratteliDiagram}}
Consider the natural sequence $(\Pg(n)$ of partitions associated to $(V,E,\le)$.
The basis of the partition $\Pg(n)$ is, as we have seen, the union of
the cylinders corresponding to minimal paths and can be identified
with $V(n)$. Set $t_n=\Card(V(n))$.
By Proposition~\ref{propositionGroupSequencePartitions}, the group
$K^0(X_E,T_E)$ is the direct limit of the groups $(G(n),G^+(n),\mathbf{1}_n)$
where $G(n)$ can be identified with $\Z^{t(n)}$ and with the
morphisms $I(n,n-1)$ defined (by Equation \eqref{equationMatrixPartitions})
by the connecting matrices $M(n-1,n)$.
Since $M(n-1,n)=M(n)$ by Lemma~\ref{lemmaConnectingMatrices},
this shows that $K^0(X,T)$ is the direct limit
$\Z^{t(0)}\edge{M(1)}\Z^{t(1)}\edge{M(2)}\cdots$, which is precisely $\DG(V,E)$.
\end{proofof}

\begin{example}
Let $(V,E)$ be the properly ordered Bratteli diagram represented in Figure~\ref{figureBratteliProperly}. 
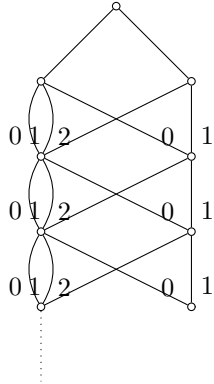
\begin{figure}[ht]
        \centering
\tikzset{node/.style={circle,draw,minimum size=0.1cm,inner sep=0cm}}
\tikzset{title/.style={minimum size=0cm,inner sep=0pt}}
 \begin{tikzpicture}
\node[node](0)at(1,6){};
\node[node](11)at(0,5){};\node[node](12)at(2,5){};
\node[node](21)at(0,4){};\node[node](22)at(2,4){};
\node[node](31)at(0,3){};\node[node](32)at(2,3){};
\node[node](41)at(0,2){};\node[node](42)at(2,2){};
\node[title](51)at(0,1){};\node[title](52)at(2,2){};
%E(1)
\draw(0)edge node{}(11);\draw(0)edge node{}(12);
%E(2)
\draw[left,near end,bend right](11)edge node{$0$}(21);
\draw[left,near end,bend left](11)edge node{$1$}(21);
\draw[right,near end](11)edge node{$0$}(22);
\draw[left,near end](12)edge node{$2$}(21);
\draw[right,near end](12)edge node{$1$}(22);
%E(3)
\draw[left,near end,bend right](21)edge node{$0$}(31);
\draw[left,near end,bend left](21)edge node{$1$}(31);
\draw[right,near end](21)edge node{$0$}(32);
\draw[left,near end](22)edge node{$2$}(31);
\draw[right,near end](22)edge node{$1$}(32);
%E(4)
\draw[left,near end,bend right](31)edge node{$0$}(41);
\draw[left,near end,bend left](31)edge node{$1$}(41);
\draw[right,near end](31)edge node{$0$}(42);
\draw[left,near end](32)edge node{$2$}(41);
\draw[right,near end](32)edge node{$1$}(42);
%etc
\draw[dotted](41)edge node{}(51);\draw[dotted](42)edge node{}(52);
\end{tikzpicture}
                
                \caption{A properly ordered Bratteli diagram.}
        \label{figureBratteliProperly}
\end{figure}
The adjacency matrix at each level is
\begin{displaymath}
M=\begin{bmatrix}2&1\\1&1\end{bmatrix}=\begin{bmatrix}1&1\\1&0\end{bmatrix}^2.
\end{displaymath}
Thus the dimension group $K^0(X_E,T_E)$ is $\Z+\lambda\Z$ with
$\lambda=(1+\sqrt{5})/2$ (see Example~\ref{exampleGolden1}).
\end{example}

%%%%%%%%%%%%%%%%%%%%%%%%%%%%%%%%%
\section{Kakutani equivalence}\label{sectionKakutani}
%%%%%%%%%%%%%%%%%%%%%%%%%%%%%%%%%

The minimal Cantor dynamical systems $(X,T)$ and $(Y,S)$ are {\em
Kakutani equivalent}\index{subject}{Kakutani equivalence} if they have (up to isomorphism) a common induced system, that is, 
there exist nonempty clopen sets $U\subseteq X$ and $V\subseteq Y$ such
that the induced systems $(X_U , T_U)$ and $(Y_V , S_V)$ are isomorphic
(see Exercise~\ref{exerciseKakutaniEquivalence} for a proof that
it is really an equivalence relation).

Two systems which are conjugate are Kakutani equivalent but the converse
is false. For example, The system with two points is Kakutani
equivalent with the system with one point but they are not conjugate.

Let us relate Kakutani equivalence to Bratteli diagrams.
If $(V,E,\le)$ is a properly ordered Bratteli diagram we may
change it into a new properly ordered Bratteli diagram $(V',E',\le')$
by making a finite change, that is, by  adding and/or removing any finite number of
edges (vertices), and then making arbitrary choices of linear orderings of
the edges meeting at the same vertex (for a finite number of vertices).  So
$(V,E,\le)$ and $(V',E',\le')$ are cofinally identical, that is, they only differ on finite
initial portions.  (Observe that this defines an equivalence relation on the
family of properly ordered Bratteli diagrams.)
We have the following nice characterization of the Kakutani equivalence.

\begin{theorem}[Giordano, Putnam, Skau]\index{names}{Giordano, Thierry}\index{names}{Putnam, Ian F.}\index{names}{Skau, Christian F.}\label{theoremKakutaniEquiv}
Let $(X_E,T_E)$ be the dynamical
system associated with the properly ordered Bratteli diagram $(V,E,\le)$.
Then the minimal Cantor dynamical system $(X,T )$ is Kakutani equivalent to $(X_E,T_E)$ if
and only if $( X,T)$ is isomorphic to $(X_{E'},T_{E'})$, where $(V',E',
\le')$ is obtained from $(V,E,\le)$ by a finite change as described above.
\end{theorem}

An interesting particular case of this result is the following. 
Let $U$ be a clopen set of $(X_E , T_E)$. 
It is a finite union of cylinder sets. 
We can suppose that  they all have the same length, that is, for some $n$, $U=\cup_{p\in P} [p]$ where  $P$ is a set of paths from level $n$ to level $0$. 
To obtain a BV-representation of the induced system on $U$ it suffices to take the properly ordered Bratteli diagram $(V',E',\le')$ which consists  of  all the paths starting with an element of $P$ endowed with the induced ordering. 
It is not too much work to prove that the induced system on $U$ is isomorphic to $(X_{E'} , T_{E'})$. 

Conversely, we will have several occasions to use the following
statement which describes the construction of a BV-representation
for a system from one of an induced system (see for example Proposition~\ref{theoremDGSubstitutionShifts}).

\begin{proposition}\label{propositionBVPrimitive}
Let $(X,T)$ be a minimal Cantor system and let $(U,T_U)$ be the induced
system on a clopen set $U\subset X$. Let $(X_E,T_E)$ be a BV-representation
of $(U,T_U)$ corresponding to a Bratteli diagram $(V,E)$
such that for every vertex $v\in V(1)$ the return time to $U$ is
constant and equal to $f(v)$ for every element of $[v]$.
Then $X$ has a BV-representation $(X_{E'},T_{E'})$ obtained
from $(X_E,T_E)$ by replacing each edge from $0$ to $v\in V(1)$
by $f(v)$ edges.
\end{proposition}
\begin{proof}
The system $(X_{E'},T_{E'})$ is isomorphic to the primitive
of $(X_E,T_E)$ relative to the function $x\mapsto f(v)$
when $x\in[v]$. Thus $(X,T)$ and $(X_{E'},T_{E'})$
are isomorphic.
\end{proof}
We illustrate Proposition~\ref{propositionBVPrimitive}
with the following simple example.
\begin{example}
Let $X$ be the set of integers of the form 
$x+3y$ where $x=0,1,2$ and $y\in \Z_2$
and $U$ be the set of those for which $x=0$. The systems
$(X,T)$ with $T$ being the addition of $1$ and the system induced on $U$
have the BV-representations shown in Figure~\ref{figureBVInduced}.
The diagram on the right is (up to the first level) the 
usual BV-representation of $\Z_2$. The diagram on the left
is the same except for the first level made of $3$ edges
in agreement with Proposition~\ref{propositionBVPrimitive}.
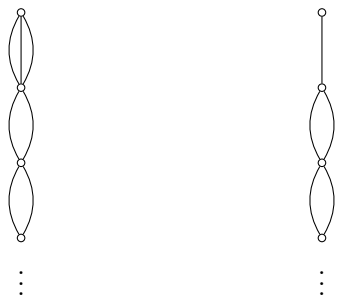
\begin{figure}[hbt]
\centering
\tikzset{node/.style={circle,draw,minimum size=0.1cm,inner sep=0cm}}
\tikzset{title/.style={minimum size=0cm,inner sep=0pt}}
\begin{tikzpicture}
\node[node](0)at(0,4){};
\node[node](1)at(0,3){};
\node[node](2)at(0,2){};
\node[node](3)at(0,1){};
\node[title]at(0,0.5){$\vdots$};

\draw[bend left](0)edge node{}(1);\draw[bend right](0)edge node{}(1);
\draw(0)edge node{}(1);
\draw[bend left](1)edge node{}(2);\draw[bend right](1)edge node{}(2);
\draw[bend left](2)edge node{}(3);\draw[bend right](2)edge node{}(3);

\node[node](0)at(4,4){};
\node[node](1)at(4,3){};
\node[node](2)at(4,2){};
\node[node](3)at(4,1){};
\node[title]at(4,0.5){$\vdots$};

\draw(0)edge node{}(1);
\draw[bend left](1)edge node{}(2);\draw[bend right](1)edge node{}(2);
\draw[bend left](2)edge node{}(3);\draw[bend right](2)edge node{}(3);
\end{tikzpicture}
\caption{The BV-representation of a system and of its derivative.}
\label{figureBVInduced}
\end{figure}
\end{example}
%%%%%%%%%%%%%%%%%%%%%%%%%%%%%%%%%%%%%%%%%%%%%%%
\section{The Strong Orbit Equivalence Theorem}\label{sectionSOE}
%%%%%%%%%%%%%%%%%%%%%%%%%%%%%%%%%%%%%%%%%%%%%%%

We say that  two dynamical systems $(X,T)$ and $(Y,S)$ are {\em orbit equivalent}
\index{subject}{orbit!equivalent}\index{subject}{equivalent!orbit}%
whenever there exists a homeomorphism $\phi : X \to Y$ sending orbits to orbits
$$
\phi \left(\{ T^n x \mid  n\in \Z \} \right) = \{ S^n \phi (x) \mid  n\in \Z \} \ ,
$$
for all $x\in X$.
This induces the existence of maps $\alpha : X\to \Z$ and $\beta : X \to \Z$ satisfying
for all $x\in X$
$$
\phi \circ T (x) = S^{\alpha (x)}\circ \phi (x) \hbox{ and } \phi \circ T^{\beta (x)} (x) = S\circ \phi (x) .
$$
These maps are called the \emph{orbit cocycles}
\index{subject}{orbit!cocycles}
associated to $\phi$.

When $\alpha$ and $\beta$ have at most one point of discontinuity, we say  that $(X,T)$ and $(Y,S)$ are {\em strongly orbit equivalent} (SOE).
\index{subject}{strong!orbit equivalence}%
\index{subject}{SOE|see{strongly orbit equivalent}}%
\index{subject}{orbit!equivalent!strongly} %
It is natural to consider such a definition because
one can show  that if $\alpha$ is continuous then $(X,T)$ is conjugate to $(Y,S)$ or to $(Y,S^{-1})$ (Exercise~\ref{exerciseBoyleOrbitEquiv}). 

The following result characterizes
strong orbit equivalence by means of Bratteli diagrams and dimension groups.

Recall from Section~\ref{sectionTelescopingSimple} that
two  Bratteli diagrams $(V,E)$ and 
$(V',E')$ are telescoping equivalent 
if they have a common
\emph{intertwining}\index{subject}{intertwining of Bratteli diagrams}
 \index{subject}{Bratteli diagram!intertwining}%
that is, a Bratteli diagram $(W,F)$ such that telescoping
to odd levels gives a telescoping of
$(V,E)$ and telescoping to even levels gives a telescoping of
$(V',E')$.

\begin{theorem}[Giordano,Putnam,Skau]\label{ch5:th:GPS}
\label{theoremStrongOrbitEquivalence}
\index{names}{Giordano, Thierry}%
\index{names}{Putnam, Ian F.}
\index{names}{Skau, Christian F.}
Let $(X,T)$ and $(X',T')$ be two minimal Cantor dynamical systems.
The following are equivalent:
\begin{enumerate}
\item
\label{ch5:th:equiv1}
There exist two BV-representations, $(V, E, \le)$ of $(X,T)$ and $(V',E', \le')$ of $(X',T')$, which have a common intertwining.

\item
\label{ch5:th:equiv2}
There exist two BV-representations, $(V, E, \le)$ of $(X,T)$ and $(V',E', \le')$ of $(X',T')$, 
and a homeomorphism $\psi:X_E\to X_{E'}$ such that $\psi(x)(n)$ depends only on $x(1) \dots x(n)$ and $\psi\left(x_{\rm u}\right) = x'_{\rm u}$, $u\in \{{\rm min}, {\max} \}$, and having the property that if $x$ and $y$ are cofinal from level $n$ on,  then $\psi(x)$ and $\psi(y)$ are cofinal from level $n+1$.
\item
\label{ch5:th:equiv3}
$(X,T)$ and $(X',T')$ are strong orbit equivalent.
\item
\label{ch5:th:equiv4}
The dimension groups $K^0(X,T)$ and $K^0(X',T')$ are isomorphic as unital ordered groups.
\end{enumerate}
\end{theorem}

\begin{proof}
Let us show that \ref{ch5:th:equiv1} implies \ref{ch5:th:equiv2}.
Let $(W,F)$ be a common intertwining of $(V,E)$ and $(V',E')$.
Since  a BV-representation is assumed to be simple, we can also suppose that all adjacency matrices have entries at least equal to two.
This means that every pair of vertices in consecutive levels has at least two connecting edges.

Let $x_{\rm min}$ and $x_{\rm max}$ be the minimal and maximal paths of $(V,E,\le)$, and  let $x'_{\rm min}$ and $x'_{\rm max}$ be those for $(V',E',\le')$.
There are unique paths $\tilde{x}_{\rm min}$ and $\tilde{x}'_{\rm min}$ in $(W,F)$ that contract respectively to $x_{\rm min}$ and $x'_{\rm min}$.
Choose a path $z_{\rm min}$ in $(W,F)$ passing through the same vertices as $\tilde{x}_{\rm min}$ does at odd levels and through the same vertices as $\tilde{x}'_{\rm min}$ at even levels.
We similarly construct a path $z_{\rm max}$ by taking care that  it does not share any common edge with $z_{\rm min}$. This is possible because the incidence matrices have entries larger than two.

Let us define two homeomorphisms $\phi : X_{F} \to X_{E}$ and $\phi' : X_{F} \to X_{E'}$.
In constructing $z_{\rm min}$, for each even $n$, we matched a pair of edges in $F(n)\circ F(n+1)$ with an edge in $E(n/2)$,  namely
$$
(z_{\rm min}(n) , z_{\rm min} (n+1) ) \to x_{\rm min} (n/2) \ ,
$$
and we matched a pair in $F(n+1)\circ F(n+2)$ with an edge in $E'((n+2)/2)$,
namely
$$
(z_{\rm min}(n+1) , z_{\rm min} (n+2) ) \to x'_{\rm min} ((n+2)/2).
$$

In the same way $(z_{\rm max}(n) , z_{\rm max} (n+1) )$ is matched with $x_{\rm max} (n/2)$ and 
$(z_{\rm max}(n+1) , z_{\rm max} (n+2) )$ with $x'_{\rm max} ((n+2)/2)$.
Now, for all even $n$, we  extend these matchings in an arbitrary way to bijections respecting the range and source maps from 
$F(n)\circ F(n+1)$ to $E(n/2)$ and from $F(n+1)\circ F(n+2)$ to $E'((n+2)/2)$.
This defines two homeomorphisms
$$
\phi : X_F \to X_{E} \hbox{ and } \phi' : X_F \to X_{E'} .
$$  
The homeomorphism $\psi=\phi'\circ \phi^{-1}$ has  the desired properties.

\medskip

Let us show that \ref{ch5:th:equiv2} implies \ref{ch5:th:equiv3}.
We will show that $(X_{E}, T_{E})$ and $(X_{E'}, T_{E'})$ are SOE.
In a minimal BV-representation, two points belong to the same orbit if and only if they are cofinal, except when it is the orbit of the minimal path.
This implies that $\psi$ maps orbits to orbits with the possible exception of the orbit of the minimal paths. 
But since $\psi\left(x_{\rm u}\right) = x'_{\rm u}$, $u\in \{{\rm min}, {\max} \}$, this is also true for the orbit of the minimal paths.
Consequently, there are maps $\alpha : X_{E}\to \Z$ and $\beta : X_{E} \to \Z$ uniquely defined by the relations
$$
\psi \circ T_{E} (x) = T_{E'}^{\alpha (x)}\circ \psi (x) \hbox{ and }
 \psi \circ T_{E}^{\beta(x)} (x) = T_{E'}\circ \psi (x) 
$$
for all $x\in X_{E}$.
It remains to prove that $\alpha$ and $\beta$ are continuous with the possible exception of $x_{\max }$ and $x_{\max }'$.
We do it for $\alpha$. 
It is similar for $\beta$.

Let $x=(x_n)_n \in X_{E}\setminus \{ x_{\max } \}$ and $k=\alpha (x)$.
Let $n_0$ be such that $(x_1, \ldots , x_{n_0})$ has a  non-maximal edge and the minimum number of 
paths from any vertex in $V_{n_0-1}$ to $V_0$ is greater than $k$.

Let $y$ belong to the cylinder $[x_1,\ldots , x_{n_0+1} ]$. 
It suffices to show that $\alpha (y) = k$.
The paths $T_{E} (x)$ and 
$T_{E} (y)$ start with the same $n_0+1$ first edges.
Thus, from the property of $\psi$, $\psi\circ T_{E} (x)$ and $\psi\circ T_{E} (y)$ start with the same $n_0+1$ first edges $f_1$, $f_2$, $\ldots$, $f_{n_0+1}$:
\begin{align*}
\psi\circ T_{E} (x) & = (f_1, f_2, \ldots, f_{n_0+1} , x'_{n_0+2} , \ldots ) \hbox{ and } \\
\psi\circ T_{E} (y) & = (f_1, f_2, \ldots, f_{n_0+1} , y'_{n_0+2} , \ldots ) .
\end{align*}

For the same reason, and because $x$ and $T_{E} (x)$, and, $y$ and $T_{E} (y)$ are cofinal from  $n_0+1$, $\psi(x)$ and $\psi(y)$ start with the same edges $g_1$, $g_2$, $\ldots$, $g_{n_0+1}$  and
\begin{align*}
\psi (x) & = (g_1, g_2, \ldots, g_{n_0+1} , x'_{n_0+2} , \ldots ) \hbox{ and } \\
\psi (y) & = (g_1, g_2, \ldots, g_{n_0+1} , y'_{n_0+2} , \ldots ) .
\end{align*}

But since there are at least $k$ paths from any vertex in $V_{n_0-1}$ to $V_0$, we deduce that 
$T^k_{E} ([g_1, g_2, \ldots, g_{n_0+1}]) = [f_1, f_2, \ldots, f_{n_0+1}]$ because
 $\psi \circ T_{E} (x) = T_{E'}^{k}\circ \psi (x)$.
Therefore $\psi \circ T_{E} (y) = T_{E'}^{k}\circ \psi (y)$ and $\alpha (y) = k$.

\medskip

Let us show that \ref{ch5:th:equiv3} implies \ref{ch5:th:equiv4}.
Let $\psi : (X,T) \to (X',T')$ be a SOE map.
Remark that $(X',T')$ is isomorphic to $(X,\psi^{-1}\circ T' \circ \psi)$. 
Hence we can suppose $X'=X$ and set $S=\psi\circ T'\circ\psi$
Then we have
$$
T(x) = S^{\alpha (x)} (x) \hbox{ and } S(x) = T^{\beta (x)} (x) 
$$
where $\alpha $ and $\beta $ are continuous everywhere, except
possibly on some $y$.

Let $A$ be a clopen set not containing $y$.
Since $\alpha$ is continuous on $A$, the set $\alpha (A)$ is compact and consequently finite:
there exist $n_1$, $\ldots$, $n_k$ such that $A= \cup_{1\leq i\leq k} A\cap \alpha^{-1} (\{ n_i \})$. 
Recall that   the 
 characteristic function
 of the set $A$ is denoted by $\charac_A$.
Hence $$TA = \bigcup_{1\leq i\leq k} S^{n_i}(A\cap \alpha^{-1} (\{ n_i \}))$$
 and thus, taking the characteristic function of each side,
$$\charac_{A}\circ T^{-1} = \sum_{1\leq i\leq k} \charac_{A\cap \alpha^{-1} (\{ n_i \})}\circ S^{-n_i} \ . $$
But since $f - f\circ S^{-n} = (\sum_{1\leq i\leq n} f\circ S^{-i})\circ S - (\sum_{1\leq i\leq n} f\circ S^{-i})$, 
we deduce that $\charac_A\circ T^{-1} - \charac_A$ belongs to $\partial_S (C(X,\Z))$.

Now suppose  that $A$ contains $y$. 
Remark that $\charac_A \circ T^{-1} -\charac_A =  \charac_{X\setminus A} - \charac_{X\setminus A} \circ T^{-1}$.
But since $y$ is not contained in $X\setminus A$ we deduce from the previous case that $\charac_A \circ T^{-1} -\charac_A$ belongs
to $\partial_S (C(X,\Z))$.
Since $T$ is invertible we proved that for every clopen set $U$, $\charac_U \circ T -\charac_U$ belongs to $\partial_{S} (C(X,\Z))$ and consequently that 
$\partial_T (C(X,\Z)) \subseteq \partial_{S} (C(X,\Z))$.
Proceeding similarly with the equality $S(x) = T^{\beta (x)} (x)$ we obtain $\partial_T (C(X,\Z)) = \partial_{S} (C(X,\Z))$.

This shows that $H(X,T,\Z)=H(X,S,\Z)$ and thus $K^0(X,T)=K^0(X,S)$.

\medskip

Finally the implication \ref{ch5:th:equiv4} implies \ref{ch5:th:equiv1}
results of Theorem~\ref{theoremDGBratteli}.

\end{proof}

We give below an example of strongly orbit equivalent minimal
Cantor systems which are not conjugate.
\begin{example}
Consider the primitive substitution $\sigma:a\mapsto ab,b\mapsto a^2b^2$.
As we shall see in Chapter~\ref{chapterSubstitutionShifts},
the shift $X(\sigma)$ is conjugate to $(X_E,T_E)$ where $(V,E)$
is the Bratteli diagram represented in Figure~\ref{figureSOE}
on the right (the morphism $\sigma$ is  read on $(V,E)$).
\begin{figure}[ht]
        \centering
\tikzset{node/.style={circle,draw,minimum size=0.1cm,inner sep=0cm}}
\tikzset{title/.style={minimum size=0cm,inner sep=0pt}}
\begin{tikzpicture}
%odometre
\node[node](0)at(0,5){};
\node[node](1)at(0,3){};
\node[node](2)at(0,1){};
\node[title](31)at(-0.3,0){};\node[title](32)at(0,0){};\node[title](33)at(.3,0){};

\draw[bend left](0)edge node{}(1);
%\draw(0)edge node{}(1);
\draw[bend right](0)edge node{}(1);
 \draw[bend left](1)edge node{}(2);
 \draw(1)edge node{}(2);
\draw[bend right](1)edge node{}(2); 
  \draw[in=90,out=-120](2)edge node{}(31);
 \draw(2)edge node{}(32);
\draw[in=90,out=-60](2)edge node{}(33);  
%intertwinning
\node[node](0)at(5,5){};
\node[node](11)at(4,4){};
\node[node](12)at(6,4){};
\node[node](2)at(5,3){};
\node[node](31)at(4,2){};
\node[node](32)at(6,2){};
\node[node](4)at(5,1){};
\node[node](51)at(4,0){};
\node[node](52)at(6,0){};

\draw(0)edge node{}(11);
%\draw[bend left](0)edge node{}(12);
\draw(0)edge node{}(12);
\draw(11)edge node{}(2);
\draw(12)edge node{}(2);
\draw(2)edge node{}(31);
\draw[bend left](2)edge node{}(32);
\draw[bend right](2)edge node{}(32);
\draw(31)edge node{}(4);
\draw(32)edge node{}(4);
\draw(4)edge node{}(51);
\draw[bend left](4)edge node{}(52);
\draw[bend right](4)edge node{}(52);
%substitution
\node[node](0)at(10,5){};
\node[node](11)at(9,4){};\node[node](12)at(11,4){};
\node[node](21)at(9,2){};\node[node](22)at(11,2){};
\node[node](31)at(9,0){};\node[node](32)at(11,0){};

\draw(0)edge node{}(11);\draw(0)edge node{}(12);
%\draw[bend right](0)edge node{}(12);
\draw[left,near end](11)edge node{$1$}(21);
\draw[bend right=15,left, near end](11)edge node{$1$}(22);
\draw[bend left=15,left,near end](11)edge node{$2$}(22);
\draw[left,near end](12)edge node{$2$}(21);
\draw[bend right,right,near end](12)edge node{$3$}(22);
\draw[bend left,right,near end](12)edge node{$4$}(22);
\draw(21)edge node{}(31){};
\draw[bend left=15](21)edge node{}(32){};\draw[bend right=15](21)edge node{}(32){};
\draw(22)edge node{}(31){};
\draw[bend left](22)edge node{}(32){};\draw[bend right](22)edge node{}(32){};
\end{tikzpicture}        
\caption{Three Bratteli diagrams.}\label{figureSOE}
\end{figure}
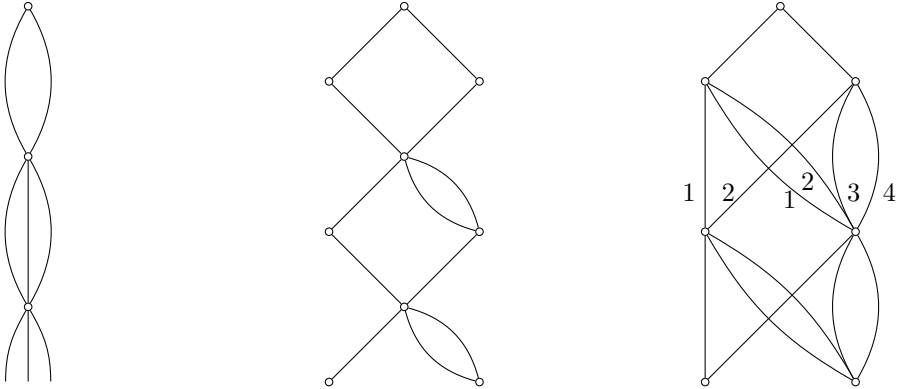
The incidence matrix of the diagram is 
\begin{displaymath}
M=\begin{bmatrix}1&1\\2&2\end{bmatrix}
\end{displaymath}
which has eigenvalues $0,3$. We claim that the dimension group is isomorphic to
$\Z[1/3]$ with the usual ordering but with order unit $2$,
that is, the group $\frac{1}{2}\Z[1/3]$. Indeed, in the
basis 
\begin{displaymath}
u=\begin{bmatrix}1\\2\end{bmatrix},v=\begin{bmatrix}1\\-1\end{bmatrix}
\end{displaymath}
we have 
\begin{displaymath}
\begin{bmatrix}1\\1\end{bmatrix}=\frac{2}{3}u+\frac{1}{3}v
\end{displaymath}
Since multiplication by $3$ is an automorphism, this proves the claim.
Thus $X(\sigma)$ has the same dimension group as the odometer
in base $p_n=2.3^{n-1}$. Indeed, the dimension group
of this odometer is, by Proposition~\ref{propositionDimensionGroupOdometer}
the subgroup of $\Q$ formed of the $p/q$
with $q$ dividing some $2.3^n$, which is
$\frac{1}{2}\Z[1/3]$.

 The shift $X(\sigma)$ and the odometer
in base $2.3^{n-1}$ are thus strong orbit equivalent by 
Theorem~\ref{ch5:th:GPS}. The BV-representation of this odometer
is represented in Figure~\ref{figureSOE} on the left.
An intertwining of the Bratteli diagrams is represented in the middle.
\end{example}

Concerning orbit equivalence, we have the following additional result
that we quote without proof.
Recall from Chapter~\ref{chapterDirectLimitsOrderedGroups} that we denote by $\Inf(G)$ the infinitesimal subgroup
\index{subject}{infinitesimal!subgroup} of a unital ordered group $G$.
\begin{theorem}\label{theoremOrbitEquivalence}
Let $(X,T)$ and $(X',T')$ be two Cantor minimal systems
and let $G=K^0(X,T)$, $G'=K^0(X',T')$ be their
dimension groups. The following conditions
are equivalent.
\begin{enumerate}
\item[\rm (i)] The systems $(X,T)$ and $(X',T')$ are orbit equivalent.
\item[\rm(ii)] The groups $G/\Inf(G)$ and $G'/\Inf(G')$ are
isomorphic.
\end{enumerate}
\end{theorem}
We give below an example of orbit equivalent shifts which are 
not strongly orbit equivalent.

\begin{example}
Let $\sigma$ be the primitive and proper substitution $a\to aab,b\to abb$.
The corresponding substitution shift is a Toeplitz shift
(in the same way as in Example~\ref{exampleToeplitz}].
A BV-representation of $X(\sigma)$ is given in Figure~\ref{figureOE}
on the right.
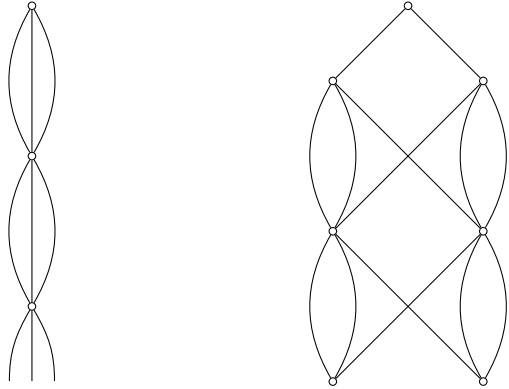
\begin{figure}[ht]
        \centering
\tikzset{node/.style={circle,draw,minimum size=0.1cm,inner sep=0cm}}
\tikzset{title/.style={minimum size=0cm,inner sep=0pt}}
\begin{tikzpicture}
%odometre
\node[node](0)at(0,5){};
\node[node](1)at(0,3){};
\node[node](2)at(0,1){};
\node[title](31)at(-0.3,0){};\node[title](32)at(0,0){};\node[title](33)at(.3,0){};

\draw[bend left](0)edge node{}(1);
\draw(0)edge node{}(1);
\draw[bend right](0)edge node{}(1);
 \draw[bend left](1)edge node{}(2);
 \draw(1)edge node{}(2);
\draw[bend right](1)edge node{}(2); 
  \draw[in=90,out=-120](2)edge node{}(31);
 \draw(2)edge node{}(32);
\draw[in=90,out=-60](2)edge node{}(33);  

%substitution
\node[node](0)at(5,5){};
\node[node](11)at(4,4){};\node[node](12)at(6,4){};
\node[node](21)at(4,2){};\node[node](22)at(6,2){};
\node[node](31)at(4,0){};\node[node](32)at(6,0){};

\draw(0)edge node{}(11);\draw(0)edge node{}(12);
%\draw[bend right](0)edge node{}(12);
\draw[bend left](11)edge node{}(21){};\draw[bend right](11)edge node{}(21){};
\draw(11)edge node{}(22){};
\draw(12)edge node{}(21){};
\draw[bend left](12)edge node{}(22){};\draw[bend right](12)edge node{}(22){};
\draw[bend left](21)edge node{}(31){};\draw[bend right](21)edge node{}(31){};
\draw(21)edge node{}(32){};
\draw(22)edge node{}(31){};
\draw[bend left](22)edge node{}(32){};\draw[bend right](22)edge node{}(32){};
\end{tikzpicture}        
\caption{ Bratteli diagrams of orbit equivalent systems.}\label{figureOE}
\end{figure}
The incidence matrix of this Bratteli diagram is
\begin{displaymath}
M=\begin{bmatrix}2&1\\1&2\end{bmatrix}
\end{displaymath}
It has eigenvalues $1,3$. The dimension group is
$G=\Z[1/3]\times\Z$ with $G^+=(\Z_+[1/3]\setminus\{0\})\times\Z\cup \{(0,0)\}$ and $1_G=(1,1)$.
Thus (see Example~\ref{exampleInfinitesimals}), the
quotient $G/\Inf(G)$ is isomorphic to $\Z[1/3]$. By Theorem~\ref{theoremOrbitEquivalence}, the shift $X(\sigma)$ is orbit equivalent to the odometer 
$Z_3$ whose BV-representation is shown on the left.
\end{example}

%%%%%%%%%%%%%%%%%%%%%%%%%%%%%%%%%%%ù
\section{Equivalences on Cantor spaces}\label{sectionEtaleEquiv}
%%%%%%%%%%%%%%%%%%%%%%%%%%%%%%%%%%
We systematically consider in this section pairs $(X,R)$
of a topological space $X$ and an equivalence
relation $R$ on $X$. After all, the
orbit equivalence suggests the idea of studying
this notion for its own sake. We will show that
this can be successfully realized.

We have
in mind two basic examples.

The first one is is the
cofinality equivalence\index{subject}{cofinality equivalence}
 $R_E$ on the space $X_E$ of paths in Bratteli
diagram $(V,E)$. Recall that 
\begin{displaymath}
R_E=\{(x,y)\in X_E\times X_E\mid \mbox{ for some $N\ge 1$, } x_n=y_n\mbox{ for all $n\ge N$}\}.
\end{displaymath}
Note that $X_E$ is actually a one-sided shift space and that $R_E$ is 
contained in the orbit equivalence on this shift.

The second one is the orbit equivalence\index{subject}{orbit!equivalence}
\begin{displaymath}
R_T=\{(x,y)\in X\times X\mid T^nx=y\mbox{ for some $n\in\Z$}\}
\end{displaymath}
when $(X,T)$ is a topological dynamical system.

 Let us give a first example of how
 one may study pairs $(X,R)$ of an equivalence
$R$ on a space $X$ on their own.

A relation $R$ on a space $X$ is \emph{minimal}\index{subject}{minimal!relation}
\index{subject}{relation!minimal} if for every $x\in X$, the equivalence
class of $x$ is dense in $X$. 

A set $Y\subset X$ is \emph{invariant}, or $R$-invariant \index{subject}{invariant!by relation} if it is saturated by $R$, that is, $Y$ is a union of classes
of $R$.

We have then the following statement.
\begin{proposition}\label{propositionMinimalEquivalence}
If $R$ is minimal, the only closed $R$-invariant sets in $X$ are $X$ and $\emptyset$.
\end{proposition}
\begin{proof}
Assume that $Y$ is a closed nonempty $R$-invariant subset of $X$. Then
$Y$ is dense in $X$ and thus $Y=X$.
\end{proof}
We will see below that the converse is also true under an additional condition.

The orbit equivalence $R_T$ of a  dynamical system $(X,T)$
is obviously minimal if and only if the system is minimal.

More, interestingly,
the cofinality equivalence $R_E$ on a  Bratteli diagram is  minimal
if and only if the diagram is simple.

Indeed, let $(V,E)$ be a simple Bratteli diagram
and let $x,y\in X_E$. For every $n\ge 1$ there is an $m\ge n$ such
that there is a path $(z_{n+1},\ldots,z_{m-1})$ from $r(y_n)$ to $s(x_m)$. Let
$z^{(n)}=(y_0,\ldots,y_n,z_{n+1}\ldots,z_{m-1},x_m,x_{m+1},\ldots)$. Then $x$
and $z^{(n)}$ are cofinal and $\lim z^{(n)}=y$. This shows that
the class of $x$ is dense in $X_E$.

Conversely, assume that $(V,E)$ is not simple. Then,
by Proposition~\ref{propositionDirectedHereditary}, we can
find a proper subset $W$ of $V$ which is
directed and hereditary. Let $Y$
be the set of infinite paths $y$ in $X_E$ which pass
by a vertex $w$ in $W$. By definition, the set $Y$ contains all
paths which coincide with $y$ until $w$. Thus $Y$ is open.
Its complement is a closed set $Z$. The set $Z$ is also invariant
because if $z\in Z$ is cofinal to $t$, then $t$ cannot
pass by a vertex in $W$ and thus $t$ is in $Z$.
Finally, the definition of $W$ implies that the set $Z$ is nonempty
and strictly included in $X_E$.
This shows that $R_E$ is not minimal.

\subsection{Local actions and \'etale equivalences}
Let $X,Y$  be topological spaces. If $U\subset X$
and $V\subset Y$ 
are clopen sets, a homeomorphism $\gamma:U\to V$ is called a \emph{partial
homeomorphism}\index{subject}{partial!homeomorphism}.

We denote by $s(\gamma)=U$ the \emph{source}\index{subject}{source}
\index{subject}{partial!homeomorphism!source}%
of $\gamma$ and by $r(\gamma)=V$ the
\emph{range}\index{subject}{range}\index{subject}{partial!homeomorphism!range}
of $\gamma$.

If $\gamma_1:U_1\to V_1$ and $\gamma_2:U_2\to V_2$ are such partial
homeomorphisms, we denote $\gamma_1\cap \gamma_2$ the map
equal to $\gamma_1$ on the set where the two functions agree.
Thus for the intersection to be a partial homeomorphism,
we require this set to be open.

A function $f:X\to Y$ is a \emph{local homeomorphism}
\index{subject}{local!homeomorphism} if for every $x\in X$,
 there is a clopen set $U\subset X$
containing $x$ such that $f(U)\subset Y$ is open and that $f|_U$ is a partial
homeomorphism.

For a relation $\rho$ on $X$, we may consider its \emph{inverse}
\index{subject}{inverse!of relation}\index{subject}{relation!inverse}%
which is the set
of pairs $(y,x)$ for $(x,y)\in\rho$. The \emph{composition}
\index{subject}{composition!of relations}\index{subject}{relation!composition}%
$\sigma\circ\rho$ of two relations $\sigma,\rho$
on $X$ is defined as usual by 
\begin{displaymath}
\sigma\circ\rho=\{(x,y)\in X\times X\mid (x,z)\in\sigma
\mbox{ and } (z,y)\in \rho \mbox{ for some $z\in X$}\}.
\end{displaymath}
Finally the \emph{intersection}\index{subject}{intersection of relations}
\index{subject}{relation!intersection}%
of two relations is well defined since a relation on $X$ is just a subset of $X\times X$.

We will find convenient to consider a map $\gamma:X\to X$
as a relation on $X$, via the identification of
$\gamma$ and its graph $\{(x,\gamma(x))\mid x\in X\}$. 
When $\gamma:U\to V$
is a partial homeomorphism, its inverse (considered as a function
or as a relation) $\gamma^{-1}:V\to U$
is again a partial homeomorphism. When $\gamma':U'\to V'$
is another partial homeomorphism, the composition
$\gamma\circ\gamma'$ always exists, even without the
requirement that $U'=V$ as usual for maps
(note that the notation of composition for relations reverses
the order of the factors). If $U'\cap V$ is empty,
then $\gamma\circ\gamma'$ is empty. 

We denote by $s,r$ the two canonical projections from $X\times X$
onto $X$ defined by $s(x,y)=x$ and $r(x,y)=y$. This is
consistent with the notation $s(\gamma)=U$ and $r(\gamma)=V$
for a partial homeomorphism $\gamma:U\to V$. For an open
set $U$, we denote $\id_U$ the identity map on $U$,
that is, $\id_U=\{(x,x)\mid x\in U\}$.

A collection $\Gamma$ of partial homeomorphisms of $X$ is a \emph{local action}
\index{subject}{local!action}%
if 
\begin{enumerate}
\item[(i)] The collection of sets $U\subset X$ such that $\id_U\in\Gamma$
forms a base of the topology.
\item[(ii)] The family $\Gamma$ is closed under taking inverses, composition
and intersection.
\end{enumerate}
Note that for $\gamma_1,\gamma_2\in\Gamma$ the intersection
$\gamma_1\cap\gamma_2$ is the map which is equal to $\gamma_1$
(and $\gamma_2$) on the set of points $x$ where $\gamma_1(x)=\gamma_2(x)$.
Thus the condition $\gamma_1\cap\gamma_2\in \Gamma$ implies
that if $\gamma_1$ and $\gamma_2$ agree on some point $x$, they
agree on a neighborhood of $x$.

Regarding a local action $\Gamma$ as a set of binary relations on $X$,
we can consider the union $\cup\Gamma$ of all its elements. Thus
$(x,y)\in\cup\Gamma$ if and only if $y=\gamma(x)$ for
some $\gamma\in\Gamma$.
\begin{proposition}\label{propositionLocalAction}
Let $\Gamma$ be a local action and let $R=\cup\Gamma$. Then
\begin{enumerate}
\item $R$ is
an equivalence relation.
\item $\Gamma$ is a basis for a topology on $R$.
\item With this topology, the source and range maps $s,r:\cup\Gamma\to X$
are local homeomorphisms.
\end{enumerate}
\end{proposition}
\begin{proof}
1. For every $x\in X$, there is by condition (i) a set $U\subset X$
containing $x$ such that $\id_U\in \Gamma$.
Thus  $R$
is reflexive. Since $\Gamma$ is closed by inverse, $R$ is symmetric.
Finally, since $\Gamma$ is closed by composition, the relation
$R$ is transitive (note that we did not use the closure
by intersection).

2. This results from the fact that the elements of $\Gamma$ cover $R$
and that $\Gamma$ is closed under intersection.

3. For $\gamma\in\Gamma$, denote $s_\gamma=s|_\gamma$. 
Then $s_\gamma:(x,\gamma(x))\mapsto x$ is a bijection
from $\gamma$ to $s(\gamma)$. We claim that $s_\gamma$
is a homeomorphism from $\gamma$ (as a subset of $R$
with the topology from $\Gamma$)
to $s(\gamma)$ (as a subset of $X$ with usual
topology). 

To show that $s_\gamma$ is continuous,
consider a clopen set $U\subset X$ such that
$\id|_U\in\Gamma$.  The set
\begin{displaymath}
s^{-1}_\gamma(U\cap s(\gamma))=\id_U\circ\gamma
\end{displaymath} is in $\Gamma$ since $\Gamma$ is closed by composition.
Thus it is an open set for the topology of $R$. Since
such sets $U$ generate the topology of $X$ by definition of a local action, 
we conclude that $s_\gamma$ is continuous. 

To show that $s^{-1}_\gamma$ is continuous, consider $\gamma'\in\Gamma$.
We have $s_\gamma(\gamma\cap\gamma')=s(\gamma\cap\gamma')$.
Since $\gamma\cap\gamma'\in\Gamma$, the set $s(\gamma\cap\gamma')$
is clopen. Thus, $s_\gamma^{-1}$ is also continuous.
This proves the claim. 

The proof concerning $r$ is symmetric.
\end{proof}

Let $X$ be a topological space.
An equivalence relation $R$ on  $X$ with its topology is
\emph{\'etale}\index{subject}{etale@\'etale!equivalence relation}
\index{subject}{equivalence!relation!\'etale} if its
topology arises from a local action $\Gamma$ on $X$ as in 
Proposition~\ref{propositionLocalAction}.

Note that the topology on $R$ defined by a local action
will in general not be the topology induced by $X$
on the product $X\times X$. 

The following result shows that the new notion can
be used to obtain a result  which is more precise than Proposition~\ref{propositionMinimalEquivalence}.
\begin{proposition}
An \'etale equivalence relation $R$ on a Cantor set $X$
is minimal if and only if the only closed $R$-invariant
subsets of $X$ are $X$ and $\emptyset$.
\end{proposition}
\begin{proof}
We have already seen that the condition is necessary.
To prove the converse, consider $x\in X$. To show
that its class $[x]_R$ is dense, consider its closure
$\overline{[x]_R}$. We will show that it is $R$-invariant.
For this, let $(y,z)\in R$ with $y\in\overline{[x]_R}$.
Choose a sequence $(y_n)$ in $[x]_R$ with limit $y$.
Since $R$ is \'etale, there is a partial homeomorphism $\gamma\in R$
such that $(y,z)\in\gamma$. By definition of the topology on $R$,
all $(y_n)$ are in $s(\gamma)$ for $n$ large enough. For such $n$,
we set $z_n=\gamma(y_n)$. Since $\gamma$ is continuous,
the sequence $(z_n)$ converges to $z$. But each $y_n$
is in $[x]_R$ and $(y_n,z_n)$ is also in $R$,
so that $z_n$ is in $[x]_R$. This shows that
$z$ is in $\overline{[x]_R}$ and thus that $\overline{[x]_R}$
is $R$-invariant. By hypothesis, this forces $\overline{[x]_R}=X$
and thus the class of $x$ is dense, whence the conclusion
that $R$ is minimal.
\end{proof}
We are now going to show that each of the two examples of relations
on a Cantor space given above are \'etale relations.
\subsection{The \'etale relation $R_E$}
Consider the cofinality relation $R_E$ on a Bratteli
diagram $(V,E)$. We first prove the following statement
which associates to it a local action. Let $(V,E)$ be a Bratteli diagram. For $n\ge 1$ and
$p,q\in E_{0,n}$ with $r(p)=r(q)$,
let 
\begin{displaymath}
\gamma(p,q)=\{(x,y)\in X_E\times X_E\mid x\in[p],y\in[q], x_k=y_k \ (k>n)\}.
\end{displaymath}
Then each $\gamma(p,q)$ is a partial homeomorphism from $[p]$ to $[q]$.
\begin{proposition}
The set $\Gamma_E$  of all partial homeomorphisms $\gamma(p,q)$
 is a local action.
\end{proposition}
\begin{proof}
Since $\gamma(p,p)=\id_{(p]}$, the family $\Gamma_E$ contains all
$\id_{[p]}$ and thus the family
of $U\subset X_E$ such that $\id_U$ is in $\Gamma_E$
is a basis of the topology of $X$. We have
$\gamma(p,q)^{-1}=\gamma(q,p)$ and thus $\Gamma_E$ is closed under
taking inverses. Next $\gamma(p,q)\circ\gamma(q,r)=\gamma(p,r)$
and thus $\Gamma_E$ is closed under composition. Finally,
$\gamma(p,q)\cap\gamma(p',q')$ is empty or equal
to $\gamma(p,q)$ if $p=p'$ and $q=q'$. Thus $\Gamma_E$ is closed under intersection.
\end{proof}

\begin{theorem}\label{theoremREetale}
Let $(V,E)$ be a Bratteli diagram. The
cofinality equivalence $R_E$ is an \'etale equivalence relation
with respect to the topology defined by $\Gamma_E$.
\end{theorem}
\begin{proof}
Let us show that $R_E=\cup\Gamma_E$. If $x,y$ are cofinal, there is an $N\ge 1$
such that $x_n=y_n$ for all $n\ge N$. Set $p=x_0x_1\cdots x_{N-1}$
and $q=y_0y_1\cdots y_{N-1}$. Then $(x,y)\in\gamma(p,q)$ and thus
$(x,y)\in \Gamma_E$. The converse is obvious.
\end{proof}

\subsection{The \'etale relation $R_T$}
We will prove that the orbit equivalence $R_T$ on a minimal
Cantor system is \'etale. We first prove that
one may associate to $R_T$ a local action.

\begin{proposition}
Let $(X,T)$ be a minimal Cantor system. Let $\Gamma_T$
be the set of all partial homeomorphisms of the form
$T^n|_U$ for $n\in \Z$ and $U\subset X$ clopen.
Then $\Gamma_T$ is a local action.
\end{proposition}
\begin{proof}
Since $X$ is a Cantor space, the clopen sets form a basis of the topology.
And for every clopen set $U$, we have $\id_U=T^0|_U$ and thus $\id_U$
is in $\Gamma_T$. Next, if $\gamma=T^n|_U$, then $\gamma^{-1}$
is the restriction of $T^{-n}$ to $T^nU$. Thus $\Gamma_T$ is closed by taking inverses. It is also closed under composition since for $\gamma=T^n|_U$
and $\gamma'=T^{n'}|_{U'}$ with $U'=T^nU$, we have $\gamma'\circ\gamma=T^{n+n'}|_U$. Finally, since $(X,T)$ is minimal $T^nx=T^{n'}x$ implies $n=n'$.
Thus $\gamma\cap\gamma'$ is either empty or equal to $\gamma$.
\end{proof}

\begin{theorem}\label{theoremRTEtale}
Let $(X,T)$ be a minimal Cantor dynamical system.
The equivalence $R_T$ on  $(X,T)$
defined by 
\begin{displaymath}
R_T=\{(x,T^nx)\mid x\in X,n\in\Z\}
\end{displaymath}
is \'etale with respect to the topology defined by the local action $\Gamma_T$.
\end{theorem}
\begin{proof}
This follows from the fact that $R_T=\cup\Gamma_T$.
\end{proof}

%%%%%%%%%%%%%%%%%%%%%%%%%%%%%%%%%%%%
\section{Entropy and Bratteli diagrams}
%%%%%%%%%%%%%%%%%%%%%%%%%%%%%%%%%%%%
\label{ch5:sec:entropy}

%It is known that every minimal Cantor dynamical system is SOE
%to a minimal Cantor dynamical system of entropy zero.
%In this section we give a  proof of this result.
%It suffices to prove that any Bratteli diagram can be telescoped to a diagram admitting an ordering such that 
%the associated Bratteli-Vershik dynamical system has entropy zero.

%Let us recall how it can be defined for subshifts $(X,S)$ and Bratteli-Vershik dynamical systems $(X_B ,V_B )$.

The \emph{topological entropy} of a dynamical system $(X,T)$ is defined
as
\begin{equation}
  h(T)=\sup_{\alpha}h(T,\alpha)
\end{equation}
where $\alpha$ runs over the open covers of $X$ (see Appendix~\ref{appendixTopologicalEntropy}).
The topological entropy
\index{subject}{topological!entropy}%
 of a shift space $(X,S)$, denoted by $h(X,S)$, or simply
 $h(S)$, is actually (see Theorem~\ref{corollaryEntropyShiftSpace})
 the growth rate of the number $p_n(X)$ of finite words of length $n$ occurring in 
elements of $X$, that is
\begin{equation}
h(S) = \lim_{n\to\infty} \frac{1}{n}\log p_n(X).
\end{equation}
\index{symbols}{h@$h(X,S)$}%
Thus $0\le h(S)\le\log\Card(A)$ for a shift on the alphabet $A$.
One can show easily that the  limit exists
and is equal to $\inf{\frac{1}{n}\log p_n(X)}$,
using the fact that $p_{n+m}(X)\le p_n(X)p_m(X)$
(see Exercise~\ref{exerciseEntropyLimit}).

As seen in Appendix~\ref{appendixTopologicalEntropy},
entropy is invariant under conjugacy
(see Exercise~\ref{exerciseEntropyInvariant} for the case of shift spaces).
It is a measure of the `size' of a shift space. For example,
an edge shift on a graph with adjacency matrix $M$ has
entropy $\log \lambda$ where $\lambda$ is the dominant eigenvalue of $M$
(Exercise~\ref{exerciseEntropyEdgeShift}).
One could expect that minimal shifts have entropy zero.
However, one can show that there are minimal shifts of arbitrary entropy.
Indeed, one has the following result.
\begin{theorem}[Grillenberger]\label{theoremGrillenberger}
  For every $k\ge 2$ and $0<h<\log k$, there exists a minimal shift space
  $(X,S)$ on $k$ symbols such that $h(X,S)=h$.
\end{theorem}
See Exercise~\ref{exerciseGrillenberger2} for the
existence of minimal shifts with  arbitrary large entropy.

We will show that
entropy is far from being invariant by strong orbit equivalence
and that actually
every minimal Cantor dynamical system is SOE
to a minimal Cantor dynamical system of entropy zero.

Let $(V,E,\le)$ be a properly ordered Bratteli diagram.
For $n\ge 1$, let $A_n$ be the set of paths $(e_1,\ldots,e_n)$
of length $n$
from $V(0)$ to $V(n)$. For $x=(e_1,e_2,\ldots) \in X_E$, denote
$\pi_n(x)=(e_1,\ldots,e_n)$.
Let $(X_n,S_n)$ be the shift space on the alphabet $A_n$
defined by
\begin{equation}
  X_n=\{(\pi_n(T_E^kx))_{k\in\Z}\mid x\in X\}\label{eqDefX_n}
\end{equation}
 the topological
entropy of $T_E$ is given by
\begin{equation}
  h(T_E)=\lim_{n\to\infty}h(S_n).
  \end{equation}
Indeed, let $\alpha_n$ be the partition of $X_E$
defined by the map $\pi_n$. Then $h(S_n)=h(T_E,\alpha_n)$
and, by Theorem~\ref{theoremWalters}, we have $h(T_E)=\lim_{n\to\infty}h(T_E,\alpha_n)$.

Hence, we need first to compute $h(S_n)$.
To this end we will need the following shift spaces.
When $W$ is a set of finite words, we denote by $W^\infty$ the set of all bi-infinite sequences formed by concatenation of  words belonging to $W$.
Let $S_W$ denote the shift map on $W^\infty$.
It is clear $(W^\infty , S_W )$ is a shift space.

\begin{lemma}
\label{ch5:lemma:majentropy}
Let $W$ be a set of $m$ finite words of length at least $l$. 
Then 
$$
h(S_W) \leq \frac{\log m}{l} .
$$
\end{lemma}

\begin{proof}
Let $k$ be the maximal length of the  words in $W$.
Let $w$ be a  word of length $n$ occurring in some sequence of $W^\infty$.
Then there exist $r$  words $u_1$, $\ldots$, $u_r$ of $W$, a prefix $s$ and suffix $p$ of some  words in $W$ such that
$w = su_1 \cdots u_r p$.
Since  $r \leq \frac{n}{l}$,  we deduce that  there are at most $k^2 m^{2+\frac{n}{l}}$  words 
in $\cL_n (W^\infty)$,  which ends the proof. 
\end{proof}

\begin{figure}[h]
        \centering
\tikzset{node/.style={circle,draw,minimum size=0.1cm,inner sep=0cm}}
\tikzset{title/.style={minimum size=0cm,inner sep=0pt}}
 \begin{tikzpicture}
\node[node](00)at(0,2){};\node[node](01)at(2,2){};\node[node](02)at(4,2){};
\node[node](1)at(2,0){};
\node[title]at(2,-.5){$v$};

\draw[bend right,near end,left](00)edge node{$0$}(1);
\draw[bend left,near end,left](00)edge node{$1$}(1);
\draw[left](01)edge node{$4$}(1);
\draw[bend right,near end,right](02)edge node{$2$}(1);
\draw[bend left,near end,right](02)edge node{$3$}(1);
\end{tikzpicture}
        \caption{An example of a consecutive order viewed from a vertex $v\in V(n+1)$ to $V(n)$.}
        \label{ch5:fig:consecutive}
\end{figure}
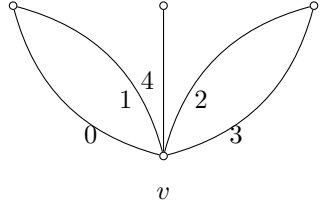

An order on a Bratteli diagram is a {\em consecutive order}
\index{subject}{consecutive!order}\index{subject}{order!consecutive} if 
whenever edges $e$, $f$ and $g$ have the same range, $e$ and $g$ have the same source and $e\leq f\leq g$, then $e$ and $f$ have the 
same source (see Figure~\ref{ch5:fig:consecutive}).

\begin{proposition}
\label{ch5:prop:entropyzero}
Let $(V,E,\le )$ be a properly ordered Bratteli diagram where $\le $ is a consecutive order.
Suppose that 
$$
\lim_{n\to +\infty} \frac{\log (\eta(n+1) \Card( V(n)) )}{\eta (n+1)} =0 ,
$$
where $\eta (n)$, $n\geq 1$, is the minimum number of edges from a vertex at level $n-1$ to a vertex at level $n$.
Then, $h(T_E) = 0$.
\end{proposition}

\begin{proof}
  Let $A_n$ be the set of paths
  of length $n$ in $(V,E)$
  starting at $V(0)$ and $X_n$ be  the shift space on $A_n$ defined
  by Equation~\eqref{eqDefX_n}.
Let $u$ be a vertex at level $n$, and let $p_1, \ldots , p_s$ be the paths from level 0 to $u$, listed in  increasing order.
We set $W(u)=p_1 \cdots  p_s$. 
We can consider  that it belongs to $A_n^s$.
Now, assume that $v$ is a vertex at level $n+1$, that $a_1,\ldots , a_t$ are the edges from $u$ to $v$,
listed in increasing order, 
and that $y$ is an infinite path in the Bratteli diagram such that $y_1 \cdots y_{n+1} = p_1 a_1$.
Then, 
$$
\left(
\pi_n
\left(
T_E^k (y)
\right)
\right)_{0\leq k \leq st-1}
=
W(u)^t .
$$

Note that $t$ is larger than $\eta (n+1)$.
Therefore, $X_n$ is included in $\mathcal{W}^\infty$ where
$$
\mathcal{W} = 
\left\{
W(u)^t \mid u\in V(n),  \ \eta (n+1) \leq t < 2\eta (n+1) 
\right\}  , 
$$
and consequently,  $h (S_n) \leq h (S_\mathcal{W} )$.
Since $\mathcal{W}$ consists of at most $\Card(V(n)) \eta (n+1)$ finite words of length at least $\eta (n+1 )$,
from Lemma~\ref{ch5:lemma:majentropy}
we obtain
$$
h (S_\mathcal{W} ) 
\leq 
\frac{\log \left(  \eta (n+1) \Card( V(n)) \right)}{\eta (n+1)} .
$$
This completes the proof.
\end{proof}

\begin{theorem}\label{theoremBoyleHandelmanEntropyZero}
Any minimal Cantor dynamical system is strongly orbit equivalent to a minimal Cantor dynamical system of entropy zero.
\end{theorem}

\begin{proof}
From Theorem~\ref{ch5:theo:BVmodel},  it suffices to consider a minimal Bratteli-Vershik dynamical system $(X_E, T_E)$.
Let $\eta (n)$, $n\geq 1$, be the minimum number of edges from a vertex at level $n-1$ to a vertex at level $n$. 

From Theorem~\ref{ch5:theo:BVmodel},  we know,  by contracting if needed,  that we can assume  the incidence matrices 
of $B=(V,E,\le )$  to have strictly positive entries.
Hence, contracting again if needed, we can suppose that 
$$
\lim_{n\to +\infty} \frac{\log (\eta(n+1) \Card( V(n)) )}{\eta (n+1)} =0 \ .
$$

Consider $(V',E',\le' )$ where $V'=V$,
$E'=E$ and $\le'$ is a consecutive order.
Then, from Proposition~\ref{ch5:prop:entropyzero}, $(X_{E'} , T_{E'} )$ has  zero entropy and, from Theorem~\ref{ch5:th:GPS}, is 
strongly orbit equivalent to $(X_E , T_E )$.
\end{proof}

Note this proof actually implies
 that all minimal BV-dynamical systems with a consecutive order have entropy zero.

We finally quote without proof the following more general statement.

\begin{theorem}[Sugizaki]\label{theoremSugisaki}
Let $\alpha \in [1, +\infty [$ and $(X,T)$ a minimal Cantor dynamical system.
There exists a minimal shift space of entropy $\log \alpha $ which is strongly orbit equivalent to $(X,T)$.
\end{theorem}

%%%%%%%%%%%%%%%%%%%%%%%%%%%%
\section{Exercises}

\exosection{Section~\ref{ch5:sec:2}}

\begin{exercise}\label{exerciseSimpleDiagram}
Show that the following conditions are equivalent for a Bratteli diagram
$(V,E)$.
\begin{enumerate}
\item[(i)] For each $m\ge 0$ and every vertex $v$ in $V(m)$, there exists
$n>m$ such that, for every $w\in V(n)$, there is a path from $v$ to $w$.
\item[(ii)] For each $m\ge 0$, there exists
$n>m$ such that, and every vertex $v$ in $V(m)$
and for every $w\in V(n)$, there is a path from $v$ to $w$.
\item[(iii)] The Bratteli diagram $(V,E)$ is simple.
\end{enumerate}
\end{exercise}

\begin{exercise}\label{exerciseEquivalenceTelescoping}
Show that the equivalence on Bratteli diagrams (ordered or not)
generated by telescoping
is given by $(V,E)\equiv(V',E')$ if there exists
a Bratteli diagram $(W,F)$ such that telescoping to odd levels
gives a telescoping of $(V,E)$ and telescoping to even levels
gives a telescoping of $(V',E')$.
\end{exercise}

\begin{exercise}\label{exerciseC*Equivalence}
Two nonnegative integer square matrices $M,N$ are said to be
\emph{$C^*$-equivalent}
\index{subject}{C*-equivalence@$C^*$-equivalence of matrices}%
if there are sequences $R_n,S_n$ of nonnegative matrices and $k_n,\ell_n$
of integers such that
\begin{equation}
M^{k_n}=R_nS_n,\quad N^{\ell_n}=S_nR_{n+1}.\label{eqC*Equiv}
\end{equation}
for all $n\ge 1$. Show that $C^*$-equivalence is an equivalence relation
containing shift equivalence. Hint: show that $M,N$ are $C^*$-equivalent 
if and only if the stationary Bratteli diagrams with matrices
$M$ and $N$ are equivalent modulo intertwining.

Let $M,N$ be the matrices
\begin{equation}
M=\begin{bmatrix}
1&1&0&0&0\\0&1&1&0&0\\0&0&1&1&0\\0&0&0&1&1\\1&0&0&0&1
\end{bmatrix}=I+P,
\quad
N=\begin{bmatrix}
0&1&1&0&0\\0&0&1&1&0\\1&0&0&0&1\\1&1&0&0&0\\0&0&0&1&1
\end{bmatrix}=QM
\label{eqMNC*EquivNotShiftEquiv}
\end{equation}
where $P$ and $Q$ denote the matrices of the permutations $(12345)$
and $(123541)$ respectively. Show that $M,N$ are $C^*$-equivalent
(it can be shown that they are 
 not shift equivalent, see the notes for a reference).
\end{exercise}

\exosection{Section~\ref{sectionDynamicsBratteli}}
\begin{exercise}\label{exerciseBratteliCantor}
Show that the Bratteli compactum $X_E$ of a simple Bratteli diagram $(V,E,\le)$
is a Cantor space.
\end{exercise}

\begin{exercise}\label{exerciseCofinalDense}
Show that a Bratteli diagram $(V,E)$ is simple if and only
if each class of the equivalence $R_E$ of cofinality is dense
in $R_E$.
\end{exercise}

\begin{exercise}\label{exerciseMaxMin}
Show that $X_E^{\rm max}$ and $X_E^{\rm min}$ are non-empty sets for
every ordered Bratteli diagram $(V,E,\le)$ and more precisely that
the set of maximal (resp. minimal) edges forms a spanning tree 
of $(V,E)$.
\end{exercise}

\begin{exercise}\label{exerciseProperlyOrdered}
Show that every simple Bratteli diagram can be properly ordered.
\end{exercise}

\begin{exercise}\label{exerciseBVhomeo}
Let $(V,E,\le)$ be properly ordered  Bratteli diagram. Show that
the Vershik map $T_E:X_E\to X_E$ is one-to-one.
\end{exercise}

\begin{exercise}\label{exerciseBVMinimal}
Show that for every properly ordered  Bratteli diagram $(V,E,\le)$, the system 
$(X_E,T_E)$ is minimal.
\end{exercise}

\begin{exercise}\label{exerciseTescopingConjugacy}
Show that a telescoping of Bratteli diagrams from $(V,E,\le)$ to $(V',E',\le')$
induces a conjugacy from $(X_E,T_E)$ to $(X_{E'},T_{E'})$.
\end{exercise}

\exosection{Section~\ref{sectionBVmodelTheorem}}

\begin{exercise}\label{ch5:ex:BVmodel}
Prove that the map $\phi$ in the proof of Theorem~\ref{ch5:theo:BVmodel} is a
homeomorphism.
\end{exercise}
\begin{exercise}\label{exerciseAllerRetour}
Let $(X,T)$ be a minimal invertible Cantor system.
Let $\Pg(n)$ be a refining sequence of partitions in towers
with bases $B(n)=\cup_{1\le t\le t(n)}B_t(n)$
 and let $(V,E)$
be the corresponding Bratteli diagram. Let $\phi:X_E\to X$
be the conjugacy from $(X_E,T_E)$ onto $(X,T)$ defined
by \eqref{eqDefinitionPhi}. Show that 
\begin{equation}
B_{t_n}(n)=\phi([e_1,\ldots,e_n])\label{equationAllerRetour}
\end{equation}
where $(e_1,\ldots,e_n)$ is the minimal path from $V(0)$
to $(n,t_n)\in V(v)$. Conclude that $\phi$ sends the partition
$[e_1,\ldots,e_n]$ for $(e_1,\ldots,e_n)\in E_{1,n}$ to
the partition $\Pg(n)$.
\end{exercise}
\exosection{Section~\ref{sectionKakutani}}
\begin{exercise}\label{exerciseDerivativePrimitive}
Let $(X_1,T_1)$ and $(X_2,T_2)$ be minimal topological
dynamical systems.
Say that $(X_1,T_1)$ is a \emph{derivative}
\index{subject}{derivative!of dynamical system}
of $(X_2,T_2)$ if $(X_1,T_1)$ is isomorphic to an induced transformation of $(X_2,T_2)$.
We also say in this case that $(X_2,T_2)$ is a \emph{primitive}
\index{subject}{primitive}%
of $(X_1,T_1)$. Show that two minimal transformations $S,T$ have a common
derivative if and only if they have a common primitive.
\end{exercise}
\begin{exercise}\label{exerciseKakutaniEquivalence}
Show that Kakutani equivalence is an equivalence relation.
\end{exercise}
\begin{exercise}\label{exerciseKakutani}
A Bratteli diagram is \emph{stationary}
\index{stationary!Bratteli diagram}%
 if all its
adjacency matrices are equal.
Prove that the family of systems isomorphic to a
stationary BV-dynamical system is stable under
Kakutani equivalence.
\end{exercise}

%\begin{exercise}\label{exerciseHoltonZamboni}
%Show that Theorem~\ref{ch5:th:holtonzamboni} does not hold for clopen
%sets instead of cylinders.\\
%\textbf{Hint.} Eigenvalues could help, see
%Section~\ref{ch5:sec:eigenvalues}.
%\end{exercise}
\exosection{Section~\ref{sectionSOE}}
\begin{exercise}\label{exerciseBoyleOrbitEquiv}
Let $(X,T)$ and $(Y,S)$ be two Cantor minimal systems
such that if there is homeomorphism $\phi:X\to Y$
which sends orbits to to orbits, that is such that there are two
maps $\alpha,\beta$
such that  $\phi\circ T(x)=S^{\alpha(x)}\circ\phi(x)$
and $\phi\circ T^{\beta(x)}(x)=S\circ\phi(x)$.
Our aim is to show that if $\alpha$ is continuous, then
$(X,T)$ is conjugate to $(Y,S)$ or to $(Y,S^{-1})$.
\begin{enumerate}
\item Show that, replacing $S$ by $\phi^{-1}\circ S\circ\phi$, we may 
assume that $X=Y$. Set $T^k(x)=S^{f_k(x)}(x)$. Show that $k\mapsto f_f(x)$
is a bijection from $\Z$ to $\Z$ which satisfies the
cohomological equation
\begin{equation}
f_k(T^jx))=f_{k+j}(x)-f_j(x).\label{eqCohomologicalExerciseSOE}
\end{equation}
\item Show that for every integer $M>0$ there is an integer  $\bar{M}$
such that $[-M,M]\subset\{f_k(x)\mid \bar{M}\le k\le\bar{M}\}$ for all $x\in X$.
\item For $m>0$, set
\begin{eqnarray*}
A_m&=&\{x\mid\forall n\ge m, f_n(x)>0\mbox{ and } f_{-n}(x)<0\},\\
B_m&=&\{x\mid\forall n\ge m, f_n(x)<0\mbox{ and } f_{-n}(x)>0\}.
\end{eqnarray*}
Show that there is $K>0$ such that $X=A_K$ or $X=B_K$. Assume that $X=A_K$
(the other case is symmetric).
\item Set 
\begin{eqnarray*}
P_m(x)&=&\Card\{f_i(x)\mid f_i(x)>0,|i|\le m\},\\
N_m(x)&=&\Card\{f_i(x)\mid f_i(x)<0,|i|\le m\},
\end{eqnarray*}
and $a(x)=(N_M(x)-P_M(x))/2$ (which is independent of $M$ for $M\ge K$.
Show that $a(Sx)=a(x)-j+1$ where $j$ is such that $f_j(x)=1$.
\item Show that $g(x)=T^{a(x)}x$ is a conjugacy from $(X,T)$ to $(X,S)$.
\end{enumerate}
Conclude that  $(X,T)$ is conjugate to $(Y,S)$ or to $(Y,S^{-1})$.
\end{exercise}

\exosection{Section~\ref{ch5:sec:entropy}}
\begin{exercise}\label{exerciseEntropyLimit}
Let $(u_n)$ be a sequence of real numbers such that $u_{n+m}\le u_n+u_m$
for all $n,m\ge 0$
(such a sequence is called \emph{subadditive}).
\index{subject}{subadditive!sequence}%
Show that the limit $\lim u_n/n$ exists 
and is equal to $\inf u_n/n$ (\emph{Fekete's Lemma}).
\index{names}{Fekete, Michael}%
\index{subject}{Fekete's Lemma}\index{subject}{Lemma!Fekete's}%
Conclude that the $\limsup$ in the definition of topological entropy
can be replaced by a limit.
\end{exercise}
\begin{exercise}\label{exerciseEntropyInvariant}
Show that entropy is invariant under conjugacy.
\end{exercise}
\begin{exercise}\label{exerciseEntropyEdgeShift}
Let $(X,S)$ be the edge shift on a graph $G$. Let $M$ be the
 adjacency matrix of $G$ and let $\lambda$ be its dominant eigenvalue.
Show that $h(X,S)=\log\lambda$.
\end{exercise}

\begin{exercise}\label{exerciseGrillenberger1}
  Let $m_j(k)$ be defined for $j,k\ge 1$ by $m_1(k)=1$
  and $m_j(k)=\prod_{i=0}^{j-2}k!^{(i)}$ for $k\ge 2$,
  where $k!^{(0)}=k$ and $k!^{(j+1)}=(k!^{(j)})!$.

  Let $\lambda_j(k)$  be defined for $j\ge 1$ by $m_{j+1}=m_je^{m_j\lambda_j}$
   or, equivalently, by $e^{m_{j+1}\lambda_{j+1}}=(e^{m_j\lambda_j})!$.
  Show that $\lambda_j(k)$ is nonincreasing
  with $j$ and that $\lambda(k)=\lim_{j\to\infty}\lambda_j(k)$ satisfies $\lim_{k\to\infty}(\log k-\lambda(k))=1$. Hint: use Stirling Formula
  \begin{equation}
1\le n!(2\pi n)^{-1/2}n^{-n}e^n\le e^\frac{1}{12n}.\label{eqStirling}
  \end{equation}
  \index{subject}{Stirling Formula}\index{subject}{Formula!Stirling}%
  \index{names}{Stirling, James}%
  \end{exercise}
\begin{exercise}\label{exerciseGrillenberger2}
 Let $U_1=\{0,1,\ldots,k-1\}$ and for $j\ge 2$,
\begin{displaymath}
U_{j+1}=\{ u_{\sigma(1)} \cdots u_{\sigma(\Card(U_j)}\mid u_i\in U_j, \sigma\in \Sg_{\Card(U_j)}\}
\end{displaymath}
where $\Sg_n$ is the \emph{symmetric group} on $\{1,2,\ldots,n\}$.
\index{subject}{symmetric!group}\index{symbols}{sigma@$\Sg_n$}%
Show that the words of $U_j$ have length $m_j$ and that $\Card(U_j)=e^{m_j\lambda_j}$ where $m_j,\lambda_j$ are as in Exercise~\ref{exerciseGrillenberger1}.
For each $j\ge 1$, let $\ell_j,r_j\in U_j$ be such that $\ell_j$ is a suffix of $\ell_{j+1}$ and $r_j$ is a prefix of $r_{j+1}$. Then define $x$ by $x_{[-m_j,m_j)}=\ell_jr_j$.
  
  Show that  $x$ is uniformly recurrent and
  such that $h(x)=\lambda(k)$. Show that for $k=2$, with $\ell_1=r_1=0$,
  the sequence $x$ is the
  Thue-Morse sequence.
\end{exercise}

%%%%%%%%%%%%%%%%%%%%%%%%%%
\section{Solutions}

\exosection{Section~\ref{ch5:sec:2}}

\begin{solution}{\ref{exerciseSimpleDiagram}}
(i)$\Rightarrow$(ii) is clear since every $V(m)$ is finite.
(ii)$\Rightarrow$(iii) Denote $f:\N\to\N$ the function defined
by $fm)=n$ if $n>m$ is the least integer such that (ii) holds.
Then the telescoping of $(V,E)$ with respect to this sequence
has the desired property.
(iii)$\Rightarrow$(i) is clear.
\end{solution}
\begin{solution}{\ref{exerciseEquivalenceTelescoping}}
We have to show that the relation is transitive. For this, consider
Bratteli diagrams $B_1,B_2,B_3$ such that $T_1$ is
an intertwining of $B_1,B_2$ and $T_2$ is an intertwining
of $B_2,B_3$. We will build an intertwining $T_3$ of $B_1,B_3$.
This will prove that the relation is transitive.

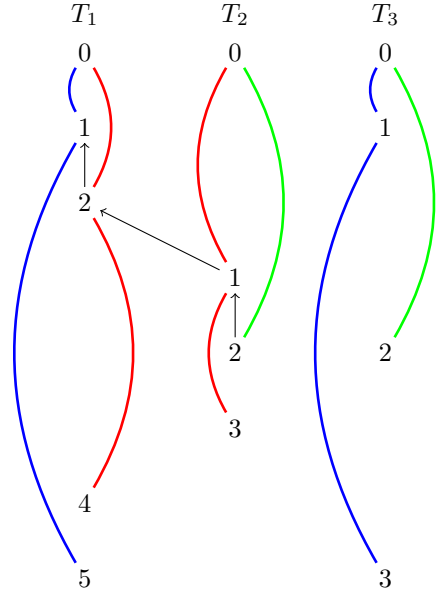
\begin{figure}[hbt]
\centering
\tikzset{node/.style={circle,draw,minimum size=0.4cm,inner sep=0.1cm}}
\tikzset{title/.style={minimum size=0.4cm,inner sep=0.4pt}}
\begin{tikzpicture}
\node[title]at(0,4.5){$T_1$};\node[title]at(2,4.5){$T_2$};
\node[title]at(4,4.5){$T_3$};
\node[title](0T1)at(0,4){$0$};\node[title](0T2)at(2,4){$0$};
\node[title](0T3)at(4,4){$0$};
\node[title](n1T1)at(0,3){$1$};\node[title](n1T3)at(4,3){$1$};
\node[title](n2)at(0,2){$2$};
\node[title](n3)at(2,1){$1$};
\node[title](n4T2)at(2,0){$2$};\node[title](n4T3)at(4,0){$2$};
\node[title](n5)at(2,-1){$3$};
\node[title](n6)at(0,-2){$4$};
\node[title](n7T1)at(0,-3){$5$};\node[title](n7T3)at(4,-3){$3$};

\draw[bend right,color=blue,line width=1pt](0T1)edge node{}(n1T1);
\draw[bend right,color=blue,line width=1pt](0T3)edge node{}(n1T3);
\draw[bend left,color=red,line width=1pt](0T1)edge node{}(n2);
\draw[bend right,color=red,line width=1pt](0T2)edge node{}(n3);
\draw[bend left,color=green,line width=1pt](0T2)edge node{}(n4T2);
\draw[bend left,color=green,line width=1pt](0T3)edge node{}(n4T3);
\draw[bend right,color=red,line width=1pt](n3)edge node{}(n5);
\draw[bend left,color=red,line width=1pt](n2)edge node{}(n6);
\draw[bend right,color=blue,line width=1pt](n1T1)edge node{}(n7T1);
\draw[bend right,color=blue,line width=1pt](n1T3)edge node{}(n7T3);

\draw[->](n4T2)edge node{}(n3);\draw[->](n3)edge node{}(n2);
\draw[->](n2)edge node{}(n1T1);
\end{tikzpicture}
\caption{The construction of the intertwining $T_3$.}\label{figureIntertwinning}
\end{figure}
Denote by $B_1(n,m)$ the matrix of $B_1$ between levels $n$ and $m$
and similarly for $B_2,B_3,T_1,T_2$. For every $n\ge 1$,
 there is an integer $\ell_1(n)$ such that 
\begin{displaymath}
T_1(n,0)=\begin{cases}B_1(\ell_1(n),0)&\mbox{ if $n$ is odd}\\
B_2(\ell_1(n),0)&\mbox{ if $n$ is even}.
\end{cases}
\end{displaymath}
and there is a similar integer $\ell_2(n)$ for $T_2$
related to $B_2$ and $B_3$. We start with $T_3(1,0)=T_1(1,0)$.
Next, by telescoping $T_2$, we can obtain $\ell_2(1)\ge\ell_1(2)$.
We define
\begin{displaymath}
T_3(2,1)=T_2(2,1)B_2(\ell_2(1),\ell_1(2))T_1(2,1).
\end{displaymath}
This corresponds to the path from level 2 in $T_2$
to level 1 in $T_1$ indicated in Figure~\ref{figureIntertwinning}.
We verify that, with this choice, $T_3(2,0)=B_3(\ell_2(2),0)$.
This follows by inspection of Figure~\ref{figureIntertwinning}
or by a patient verification as below.
\begin{eqnarray*}
T_3(2,0)&=&T_3(2,1)T_3(1,0)\\
&=&(T_2(2,1)B_2(\ell_2(1),\ell_1(2))T_1(2,1))T_3(1,0)\\
&=&(T_2(2,1)B_2(\ell_2(1),\ell_1(2))T_1(2,1))T_1(1,0)\\
&=&T_2(2,1)B_2(\ell_2(1),\ell_1(2))T_1(2,0)\\
&=&T_2(2,1)B_2(\ell_2(1),0)\\
&=&T_2(2,0)=B_3(\ell_2(2),0).
\end{eqnarray*}
We perform one more step to convince everybody that the
construction can continue in the same way forever.
We may assume, again by telescoping $T_2$ if necessary, that
$\ell_1(4)\ge \ell_2(3)$. Then we define
\begin{displaymath}
T_3(3,2)=T_1(5,4)B_2(\ell_1(4),\ell_2(3))T_2(3,2)
\end{displaymath}
and the reader may verify as above that $T_3(3,1)=B_1(\ell_1(5),\ell_1(1))$.
\end{solution}

\begin{solution}{\ref{exerciseC*Equivalence}}
The fact that $C^*$-equivalence is the same as
the existence of an intertwining of the stationary
Bratteli diagrams with matrices $M,N$ appears clearly
on the diagram of Figure~\ref{figureInterleavingStationary}.
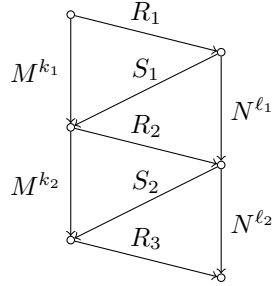
\begin{figure}[hbt]
\centering
\tikzset{node/.style={circle,draw,minimum size=0.1cm,inner sep=0cm}}
\tikzset{title/.style={minimum size=0.4cm,inner sep=0.4pt}}
\begin{tikzpicture}
\node[node](00)at(0,4.5){};\node[node](01)at(2,4){};
\node[node](10)at(0,3){};\node[node](11)at(2,2.5){};
\node[node](20)at(0,1.5){};\node[node](21)at(2,1){};

\draw[->,left](00)edge node{$M^{k_1}$}(10);
\draw[->,above](00)edge node{$R_1$}(01);
\draw[->,above](01)edge node{$S_1$}(10);
\draw[->,right](01)edge node{$N^{\ell_1}$}(11);
\draw[->,above](10)edge node{$R_2$}(11);
\draw[->,left](10)edge node{$M^{k_2}$}(20);
\draw[->,above](11)edge node{$S_2$}(20);
\draw[->,right](11)edge node{$N^{\ell_2}$}(21);
\draw[->,above](20)edge node{$R_3$}(21);
\end{tikzpicture}
\caption{Intertwining of stationary Bratteli digrams}
\label{figureInterleavingStationary}
\end{figure}
As a consequence, the $C^*$-equivalence is an equivalence relation
since the existence of an intertwining is an equivalence 
(Exercise~\ref{exerciseEquivalenceTelescoping}).

The matrices $M,N$ given by Equation~\eqref{eqMNC*EquivNotShiftEquiv}
have the same maximal eigenvalue $2$
and the same left and right eigenvectors 
$\begin{bmatrix}1&1&1&1&1\end{bmatrix}$. We may conjugate to
write 
\begin{displaymath}
M=\begin{bmatrix}2&0\\0&M_1\end{bmatrix},\quad 
N=\begin{bmatrix}2&0\\0&N_1\end{bmatrix}
\end{displaymath}
Let $\lambda$ be the maximal modulus of the eigenvalues of $M_1,M_2$
and let $\mu$ be the maximal modulus of eigenvalues of their inverses.
For every $k\ge 1$, we have
\begin{displaymath}
N^{-k}M^{ck}=\begin{bmatrix}2^{(c-1)k}&0\\0&0\end{bmatrix}
+\begin{bmatrix}0&0\\0&N_1^{-k}M_1^{ck}\end{bmatrix}
\end{displaymath}
After we conjugate back, the entries of the first term will
be at least $c_12^{(c-1)k}$ and those of the second term will
be at most $c_2\lambda^{ck}\mu^k$. Thus, if we choose
$c$ such that $2^{(c-1)k}>\lambda^c\mu$, all entries
of $N^{-k}M^{ck}$ will be positive. This shows
that for all large enough $k$, the matrix $N^{-k}M^n$ is positive
for all $n$ large enough. It is also integral because
$M,N$ have both determinant equal to $2$.

Consider now the construction of  sequences $(R_n,S_n)$ 
and $(k_n,\ell_n)$ such that
Equation~\eqref{eqC*Equiv} holds for all $n\ge 1$. Start with $S_1=I$.
Equations \eqref{eqC*Equiv} are equivalent to the equations
\begin{eqnarray*}
R_1&=&M^{k_1}\\
R_2&=&M^{-k_1}N^{\ell_1}\\
S_2&=&N^{-\ell_1}M^{k_1+k_2}\\
R_3&=&M^{-k_1-k_2}N^{\ell_1+\ell_2}
\end{eqnarray*}
and so on. By the preceding remark, we can successively choose
$k_1,\ell_1,k_2,\ell_2,\ldots$ in such a way that all $R_n,S_n$
are nonnegative and integral. Thus $M,N$ are $C^*$-equivalent.

%To see that they are not shift equivalent, we compute their
%characteristic plynomials. We find
%\begin{eqnarray*}
%\det(tI-M)&=&(t-2)(t^4-3t^3+4t^2-2t+1),\\
%\det(tI-N)&=&(t-2)(t^4+t^3+1).
%\end{eqnarray*}
%\marginpar{A TERMINER}
\end{solution}

\exosection{Section~\ref{sectionDynamicsBratteli}}

\begin{solution}{\ref{exerciseBratteliCantor}}
The  $X_E$ is compact.
If $(V',E')$ is obtained by telescoping $(V,E)$, the induced
map from $X_E$ to $X_{E'}$ is a homeomorphism. Thus, to show
that, if $(V,E)$ is simple then $X_E$ has no isolated points, we
can assume that there is an edge between any
vertices in $V(n-1)$ and $V(n)$. Assume that
$x=(e_1,e_2,\ldots)$ is an isolated point. Then
$[e_1,\ldots,e_n]=\{x\}$ for some $n\ge 1$. But then
$E(m)$ has to be reduced to $\{e_m\}$ for all $m>n$,
a contradiction with our hypothesis that $X_E$ is infinite.
\end{solution}

\begin{solution}{\ref{exerciseCofinalDense}}
Assume that $(V,E)$ is simple.
Let $e=(e_n)_{n\ge 1}$ and $f=(f_n)_{n\ge 1}$ be elements of $X_E$.
We show that $f$ belongs to the closure of the class of $e$.
Since $(V,E)$ is simple, for every $n\ge 1$, there is by
Exercise~\ref{exerciseSimpleDiagram}, an integer $m> n+1$ such
that there is a path $(g_{n+1},\ldots,g_{m-1})$
from $r(f_n)$ to $s(e_m)$. Then $h^{(n)}=(f_1,\ldots,f_n,g_{n+1},\ldots,g_{m-1},e_m,e_{m+1},\ldots)$ is a path in $X_E$ which is in the cofinality class of $e$.
Since the sequence $h^{(n)}$ tends to $f$ when $n\to\infty$, the
claim is proved.
\end{solution}
\begin{solution}{\ref{exerciseMaxMin}}
For each $n\ge 1$, set 
\begin{displaymath}
F_n=\{[e_1,\ldots,e_n]\in E_{1,n}
\mid \mbox{every $e_i$ is maximal}\}
\end{displaymath}
 An easy induction on $n$
shows that for every vertex $v\in V(n+1)$ there is an element
$[e_1,\ldots,e_n]$ in $F_n$ such that $r(e_n)=v$.
This proves that the set of maximal edges forms a spanning tree.
Next $F_n$ is closed.
 Thus, the set $X_E^{\rm max}$ is the intersection of the nonempty closed 
sets $F_n$. Since $X_E$ is compact, it is nonempty. The
argument for $X_E^{\rm min}$ is similar.
\end{solution}

\begin{solution}{\ref{exerciseProperlyOrdered}}
Let $(V,E)$ be a simple Bratteli diagram. Fix an
 order on $V(n)$ for each $n\ge 1$. We order the edges
$e\in E(n)$ in such a way that, for two edges $e,f$ with
the same range, one has $s(e)<s(f)\Rightarrow e<f$.
Then any long enough path made of minimal edges
leads to the minimal vertex.
\end{solution}

\begin{solution}{\ref{exerciseBVhomeo}}
The inverse of $T_E$ can be described as follows. 
First $T_E^{(-1)}(x_{\min})=x_{\max}$.
Next, for $e=(e_n)_{n\ge 1}$
distinct of $x_min$, let $k\ge 1$ be the minimal index
such that $e_k$ is not a  minimal edge. Then
$T_E^{-1}(e)=(f_1,\ldots,f_{k-1},f_k,e_{k+1},e_{k+2},\ldots)$ where
$f_k$ is the antecedent of $e_k$ and $(f_1,\ldots,f_{k-1})$ is
the maximal path from $v(0)$ to $s(f_k)$.
\end{solution}
\begin{solution}{\ref{exerciseBVMinimal}}
We claim that if
$e_m=f_m$ for $m> n$, and if $e_n< f_n$, then there is a $k\ge 0$
such that $T_E^k(e)=f$. This proves that the classes of cofinality 
are contained in the orbits of $T_E$ and thus implies
the statement since the classes of $R_E$ are dense by Exercise~\ref{exerciseCofinalDense}.

One proves the claim by induction on $n\ge 1$. It is clearly true
for $n=1$. Next, let $(g_1,\ldots,g_{n-1})$ be a path of maximal edges
from $v(0)$ to $s(e_n)$. Then there is an integer $k\ge 1$
and a path $(h_1,\ldots,h_{n-1})$ of minimal edges
such that $T_E^k(g_1,\ldots,g_{n-1},e_n,e_{n+1},\ldots)=(h_1,\ldots,h_{n-1},f_n,f_{n+1},\ldots)$. By induction hypothesis, there is an $\ell\ge 0$
such that $T_E^\ell(h_1,\ldots,h_{n-1},f_n,f_{n+1},\ldots)=f$. This proves the claim.
\end{solution}
\begin{solution}{\ref{exerciseTescopingConjugacy}}
Let $m_0=0<m_1<m_2<\ldots$ be the sequence defining the telescoping
from $(V,E,\le)$ to $(V',E',\le')$. Let $\varphi:X_E\to X_{E'}$
be the map defined by $y=\varphi(x)$ if $x=(e_1,e_2,\ldots)$ and
$y=(f_1,f_2,\ldots)$ with $f_n=(e_{m_n+1},\ldots,e_{m_{n+1}})$.
The map $\varphi$ is clearly a homeomorphism from $X_E$ onto
$X_{E'}$. To show that it is a conjugacy, we have to show that
$\varphi(T_E x)=T_{E'} \varphi(x)$. Consider first
the case of $x=x_{\max}$. Then $\varphi(T_E x)=\varphi(x_{\min})=y_{\min}$
while $T_{E'}\varphi(x)=T_{E'}y_{\max}=y_{\min}$. Next, let 
$x=(e_1,e_2,\ldots)\ne x_{\max}$
and $y=(e'_1,e'_2,\ldots)=\varphi(x)$. Let
$k$ be the least index $n$ such that $e_n$ is not maximal
and let $f_k$ be the successor of $e_k$. Let $n$ be such that
$m_n<k\le m_{n+1}$. Then 
\begin{displaymath}\varphi(T_E x)=\varphi(e_1,\ldots,f_k,\ldots)=
(e'_1,\ldots,f'_k,\ldots)
\end{displaymath}
 with $f'_k=(e_{m_n+1},\ldots,f_k,\ldots,e_{m_{n+1}})$
while 
\begin{displaymath}
T_{E'}\varphi(x)=T_{E'}(e'_1,e'_2,\ldots)=(e'_1,\ldots,f'_k,\ldots)
\end{displaymath}
since $f'_k$ is the successor of $e'_k$ in $E'_k$. Thus $\varphi$
is a conjugacy.
\end{solution}
\exosection{Section~\ref{sectionBVmodelTheorem}}

\begin{solution}{\ref{ch5:ex:BVmodel}}
Let $\psi:X\to X_E$ be defined as follows. For every $x\in X$ and $n\ge 1$,
since $\Pg(n)$ is a partition of $X$,
there is a unique pair $(k_n,t_n)$ with $1\le k_n\le h_{t_n}(n)-1$
and $1\le t_n\le t(n)$ such that $x\in T^{k_n}B_{t_n}(n)$.
We set $k_0=0$ and $t_0=1$.
Then $k_n\ge k_{n-1}$ for all $n\ge 1$.
The map $\psi(x)=(n,t_{n-1},t_n,k_n-k_{n-1})_{n\ge 1}$
is well defined and continuous. Since $\psi$
is obviously the inverse of $\varphi$, the result follows.
\end{solution}
\begin{solution}{\ref{exerciseAllerRetour}}
By definition of $(V,E)$, we have $e_i=(i,t_{i-1},t_i,j_i)$ for $1\le i\le n$.
Since each $e_i$ is minimal, we have $j_i=0$. Thus, by definition
of $\phi$,
\begin{eqnarray*}
\phi([e_1,\ldots,e_n])&=&T^{\sum_{i=1}^nj_i}B_{t_n}(n)\\
&=&B_{t_n}(n).
\end{eqnarray*}
\end{solution}
\exosection{Section~\ref{sectionKakutani}}
\begin{solution}{\ref{exerciseDerivativePrimitive}}
Let first $(Z,U)$ be a common primitive of $(X,S)$ and $(Y,T)$. Thus there are
nonempty clopen sets $A,B\subset Z$ such that $(X,S)\simeq(A,U_A)$ and $(Y,T)\simeq(B,U_B)$. For every $n\ge 0$,
$(A,U_A)$ is isomorphic to $(U^nA,U_{U^nA})$. Since $U$ is minimal, we may assume
that $C=U^nA\cap B\ne\emptyset$. Thus $(C,U_C)$ is a common derivative of $(X,S)$ and $(Y,T)$.

Conversely, let $(Z,U)$ be a common derivative of $(X,S)$ and $(Y,T)$.
Then, using the tower construction of Exercise~\ref{exerciseTowerConstruction},
we have $S\simeq(Z^g,U^g)$ and $T\simeq(Z^h,U^h)$. Let $(W,R)$ be the
result of the tower construction corresponding to
$g+h$. Then $(W,R)$ is isomorphic to $(X,S)^k$ with $k$ defined by
\begin{displaymath}
k(z,i)=\begin{cases}1&\mbox{ if $i<g(z)$}\\h(z)+1&\mbox{ if $i=g(z)$}
\end{cases}
\end{displaymath}
through the map 
\begin{displaymath}
\varphi(z,\ell)=\begin{cases}(z,\ell,1)&\mbox{ if $\ell<g(z)$}\\
(z,g(z),\ell+1)&\mbox{ otherwise}\end{cases}.
\end{displaymath}
Thus $(W,R)$ is a primitive of $(X,S)$.
In the same way, $(W,R)$ is a primitive of $(Y,T)$.
\end{solution}
\begin{solution}{\ref{exerciseKakutaniEquivalence}}
We have to show the transitivity of the relation $(X,S)\sim(Y,T)$
if $(X,S)$ and $(Y,T)$ are Kakutani equivalent. Suppose that $(X,S)\sim(Y,T)$
and $(Y,T)\sim(Z,U)$. Let $(X_1,S_1)$ be a common derivative of $(X,S)$
and $(Y,T)$ and let $(X_2,S_2)$ be a common derivative of $(Y,T)$
and $(Z,U)$ (see Figure~\ref{figureKakutaniEquiv}).
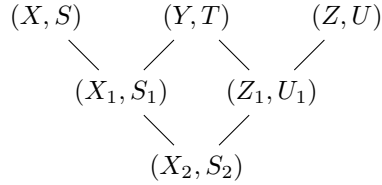
\begin{figure}[hbt]
\centering
\begin{tikzpicture}
\node(X)at(0,3){$(X,S)$};\node(Y)at(2,3){$(Y,T)$};\node(Z)at(4,3){$(Z,U)$};
\node(X1)at(1,2){$(X_1,S_1)$};\node(Z1)at(3,2){$(Z_1,U_1)$};
\node(X2)at(2,1){$(X_2,S_2)$};
\draw(X)edge node{}(X1);\draw(Y)edge node{}(X1);\draw(Y)edge node{}(Z1);
\draw(Z)edge node{}(Z1);
\draw(X1)edge node{}(X2);\draw(Z1)edge node{}(X2);
\end{tikzpicture}
\caption{Transitivity of Kakutani equivalence}\label{figureKakutaniEquiv}
\end{figure}
Since $(X_1,S_1)$ and $(Z_1,U_1)$ have a common primitive, namely $(Y,T)$
they have by Exercise~\ref{exerciseDerivativePrimitive}
a common derivative $(X_2,S_2)$. Thus $(X,S)\sim(Z,U)$.
\end{solution}
\begin{solution}{\ref{exerciseKakutani}}
Let $(V,E)$ be a stationary Bratteli diagram.
By Theorem~\ref{theoremKakutaniEquiv}, a system $(X,T)$ Kakutani equivalent
to $(X_E,T_E)$ is isomorphic to  a BV-system $(X_{E'},V_{E'})$ 
where $(V',E')$ is obtained from $(V,E)$ by a finite number of changes.
Then a telescoping of $(V',E')$ gives again a stationary Bratteli diagram.
Thus $(X,T)$ is conjugate to a stationary BV-system.
\end{solution}
%\begin{solution}{\ref{exerciseHoltonZamboni}}
%\end{solution}
\exosection{Section~\ref{sectionSOE}}
\begin{solution}{\ref{exerciseBoyleOrbitEquiv}}
\begin{enumerate}
\item Set $\tilde{S}=\phi^{-1}\circ S\circ\phi$. Then $T(x)=\tilde{S}^{\alpha(x)}(x)$.
Thus, replacing $S$ by $\tilde{S}$,  we may assume $X=Y$. Since $(X,T)$ 
is a Cantor system, there are no periodic points. Thus $k\mapsto f_k(x)$ 
is a bijection. It satisfies~\eqref{eqCohomologicalExerciseSOE} since
\begin{displaymath}
T^{k+j}x=S^{f_k(T^jx)}(T^jx)=S^{f_k(T^jx)+f_j(x)}x.
\end{displaymath}
\item Follows from compactness since the $f_k$ are continuous.
\item Set $n_0=\sup_{x\in X}|\alpha(x)|$ and choose $K$ such that
$[-n_0,n_0]\subset\{f_i(x)\mid -K\le i\le K\}$ for all $x\in X$.
Then $X=A_K\cup B_K$. Since
 $A_K,B_K$ are closed and invariant, the conclusion follows.
\item This follows by a counting argument since $f_{j+k}(x)=f_k(Sx)+1$
for every $k$
by~\eqref{eqCohomologicalExerciseSOE}.
\item We have 
\begin{eqnarray*}
g\circ S(x)&=&T^{a(Sx)}(Sx)=T^{a(x)-j+1}(T^j(x))\\
&=&T^{a(x)+1}(x)=T(T^{a(x)}(x))=T\circ g(x).
\end{eqnarray*}
\end{enumerate}
\end{solution}

\exosection{Section~\ref{ch5:sec:entropy}}
\begin{solution}{\ref{exerciseEntropyLimit}}
Set $\alpha=\inf u_n/n$.
We show that for every $\varepsilon>0$, 
we have $(u_n/n)-\alpha\le\varepsilon$ for all large enough $n$. 
Let $k$ be such that $u_k/k<\alpha+\varepsilon/2$. Then, for $0\le j<k$
and $m\ge 1$, we have
\begin{eqnarray*}
\frac{u_{mk+j}}{mk+j}&\le&\frac{u_{mk}}{mk+j}+\frac{u_j}{mk+j}\le\frac{u_{mk}}{mk}+\frac{u_j}{mk}\\
&\le&\frac{mu_k}{mk}+\frac{ju_1}{mk}\le\frac{u_k}{k}+\frac{u_1}{m}\le\alpha+\frac{\varepsilon}{2}+\frac{u_1}{m}.
\end{eqnarray*}
Hence, if $nmk+j$ is large enough so that $u_1/m<\varepsilon/2$, then
$u_n/n<\alpha+\varepsilon$.

 Set $v_n=p_n(X)$. Then $v_{n+m}\le v_nv_m$ since
a word of length $n+m$ is determined by its prefix of length $n$
and its suffix of length $m$. Thus the sequence $u_n=\log v_n$
is subadditive.
\end{solution}
\begin{solution}{\ref{exerciseEntropyInvariant}}
Let us show that if $\varphi:X\to Y$ is a surjective morphism 
from $(X,T)$ onto $(Y,S)$, then $h(Y)\le h(X)$. By Theorem~\ref{theoremCurtisHedlundLyndon}, $\varphi$ is the sliding block code associated to
some block map $f:\cL_{n+m+1}\to B$, extended to a map
from $\cL_{n+m+k}$ to $\cL_k(Y)$ for every $k\ge 1$.  Since $\varphi$
is surjective, $f$ is surjective. Thus $p_k(Y)\le p_{n+m+k}(X)$ for every
$k\ge 1$, which implies 
\begin{displaymath}
\frac{p_k(Y)}{k}\le\frac{p_k(X)}{k}+\frac{p_n(X)+p_m(X)}{k}.
\end{displaymath}
The second term of the right hand side tends to $0$ 
when $k\to\infty$ and thus $h(Y)\le h(X)$.
\end{solution}
\begin{solution}{\ref{exerciseEntropyEdgeShift}}
It is easy to show that
\begin{displaymath}
c\lambda^n\le\sum_{i,j}M^n_{i,j}\le d\lambda^n.
\end{displaymath}
for some constants $c,d$. Since $p_n(X)=\sum_{i,j}M^n_{i,j}$, the result follows.
\end{solution}

\begin{solution}{\ref{exerciseGrillenberger1}}
  Since $m_1(k)=1$ and $m_2(k)=k$, we have $\lambda_j(1)=0$. Next,
  $m_j(2)=2^{j-1}$ implies $\lambda_j(2)=2^{-j+1}\log 2$. Thus $\lambda(2)=0$.

  Next, for $k\ge 3$, let us first show by induction on $j$ that
  \begin{equation}
    e^{m_j\lambda_j}>m_j+1.\label{eqem_jlambda_j}
  \end{equation}
  Indeed, it is true for $j=1$ since $e^{m_1\lambda_1}=k>2$. Next, we have
  \begin{displaymath}
e^{m_{j+1}\lambda_{j+1}}=(e^{m_j\lambda_j})!>(e^{m_j\lambda_j}-2)e^{m_j\lambda_j}+1
  \end{displaymath}
  since $n!>(n-2)+1$ for $n\ge 3$. Thus, using the induction hypothesis,
  \begin{displaymath}
e^{m_{j+1}\lambda_{j+1}}>m_je^{m_j\lambda_j}+1=m_{j+1}+1.
  \end{displaymath}
  which proves \eqref{eqem_jlambda_j}.

  Since $n!<n^n$, we have
  \begin{displaymath}
    e^{m_{j+1}\lambda_{j+1}}<e^{m_j\lambda_je^{m_j\lambda_j}}=e^{m_{j+1}\lambda_{j+1}}
  \end{displaymath}
  showing that the sequence $(\lambda_j(k))_j$ is decreasing.

  Using now Stirling Formula
  with $n=e^{m_j\lambda_j}$ and taking the logarithm, we estimate
  \begin{displaymath}
    0\le m_{j+1}\lambda_{j+1}-\frac{1}{2}(\log2\pi+m_j\lambda_j)-m_{j+1}\lambda_j+ e^{m_j\lambda_j}\le
    \frac{1}{12}e^{-m_j\lambda_j}\le\frac{1}{36}
  \end{displaymath}
  where the last inequality follows from $e^{m_j\lambda_j}=k!^{(j-1)}\ge 3$.
  Dividing the leftmost and rightmost sides of this inequality by
  $m_{j+1}$, we obtain
  \begin{displaymath}
0\le\lambda_{j+1}-\lambda_j+\frac{1}{m_j}\le \frac{1}{2m_{j+1}}(\frac{1}{18}+\log2\pi+m_j\lambda_j)
  \end{displaymath}
  The right hand side can be bounded by $(2/m_{j+1})m_j\lambda_j<2\log k e^{-{m_j\lambda_j}}$. Summing up on $j$, we obtain
\[
 0\le\lambda(k)-\log k+\sum_{j\ge 1}\frac{1}{m_j}\le 2\log k \sum_{j\ge 1}e^{-m_j\lambda_j}.
 \]
 Now $1<\sum_{j\ge 1}\frac{1}{m_j}<\sum_{j\ge 0}1/k^j=k/(k-1)$ and thus
 $\lim_{k\to\infty}\sum_{j\ge 1}\frac{1}{m_j}=1$. Similarly,
   $\sum_{j\ge 1}e^{-m_j\lambda_j}<\sum_{j\ge 0}k^{-j-1}=1/(k-1)$
   and thus
   \begin{displaymath}
     0<\lim_{k\to\infty}(\lambda(k)-\log k -1)<\lim_{k\to\infty}2\frac{\log k}{k-1}=0.
     \end{displaymath}
\end{solution}
\begin{solution}{\ref{exerciseGrillenberger2}}
We have $\Card(U_j)=e^{m_j\lambda_j}$ and all words in $U_j$ have length $m_j$.

To show that $x$ is uniformly recurrent, consider $u\in\cL(x)$. Then
$u$ appears in some $\ell_jr_j$ with $\ell_j\ne r_j$. Then $\ell_j r_j$
appears in some $v\in U_{j+1}$ which appears itself in every
$w\in U_{j+2}$. Thus $u$ occurs in every $x_{[tm_j,(t+1)m_j)}$, which
    shows that $x$ is uniformly recurrent.

    %To show that $x$ is uniquely ergodic, we have to show that
    %$x$ has uniform frequencies (Proposition~\ref{propositionUniformFrequencies}. For this it is enough to show that every factor $u\in U_j$ has uniform
    %frequency. Set $c_j=\Card(U_j)$, $A_j=\{0,1,\ldots, c_j\}$ and
    %$U_j=\{u_1,\ldots,u_{c_j}\}$.
    %Let $\varphi:A_j^*\to U_1^*$ be the morphism defined by
    %$\varphi(i)=u_i$. Then $x=\varphi(y)$ where $y$ is a sequence
    %where every factor $y_{[tc_j,(t+1)c_j)}$ is a permutation of $A_j$.
      %Thus every letter of $y$ has uniform frequency $1/c_j$,
      %which implies that $u$ has frequency $1/m_jc_j$.

      There remains to show that $h(x)=\lambda(k)$. Set
       $c_j=\Card(U_j)=e^{m_j\lambda_j}$.
      Since $p_{m_j}(x)\ge\Card(U_j)$, we have $h(x)\ge \lambda(k)$. Next, consider $u\in\cL_{m_{j+1}}(x)$. There
      exist $u_1,\ldots,u_{c_j+1}\in U_j$ such that $u$ is a factor
      of $u_1\cdots u_{c_j+1}$. Thus
      \begin{displaymath}
p_{m_{j+1}}(x)\le m_jc_j^{c_j+1}=m_{j+1}c_j^{c_j}.
      \end{displaymath}
      This implies
      \begin{eqnarray*}
        h(x)&\le&\frac{1}{m_{j+1}}\log p_{m_{j+1}}(x)\le\frac{\log m_{j+1}}{m_{j+1}}+
        \frac{c_j}{m_{j+1}}\log c_j\\
        &\le&\frac{\log m_{j+1}}{m_{j+1}}+\lambda_j
      \end{eqnarray*}
      whence the conclusion since the right hand side tends to $\lambda(k)$
      when $j\to\infty$.

      For $k=2$, we have $U_1=\{0,1\}$ and $U_{j+1}=\{\tau^j(0),\tau^j(1)\}$
      where $\tau:0\mapsto 01,1\mapsto 10$ is the Thue-Morse morphism.
      Thus, with $\ell_1=r_1=0$, we obtain $x=\tau^\omega(0\cdot 0)$.
\end{solution}

%%%%%%%%%%%%%%%%%%%%%%%ù
\section{Notes}
Let us begin by a historical overview of the contents of this chapter.

In 1972 Ola Bratteli \citep{Bratteli1972}\index{names}{Bratteli, Ola} introduced special infinite graphs subsequently
called {\em Bratteli diagrams} which conveniently encoded the successive
embeddings of an ascending sequence $(\Ga_n)_{n\geq 0}$ of finite-dimensional semi-simple algebras over $\C$ (``multi-matrix algebras''). 
The sequence $(\Ga_n)_{n\geq 0}$ determines a so-called approximately finite-dimensional (AF) $C^*$-algebra (see Chapter~\ref{chapterBratteli}). 
Bratteli proved that the equivalence relation on Bratteli diagrams generated by the operation of telescoping is a complete isomorphism
invariant for AF-algebras. 

From a different direction came the extremely fruitful idea of A. M. Vershik \citep{Vershik:1985}\index{names}{Vershik, Anatol M.} to associate dynamics (called {\em adic transformations})\index{subject}{adic!transformation}\index{subject}{transformation!adic} with Bratteli diagrams
({\em Markov compacta})\index{subject}{Markov compactum} by introducing a lexicographic order on the infinite
paths of the diagram. 
By a careful refining of Vershik's construction,
R. H. Herman, I. F. Putnam and C. F. Skau \cite{HermanPutnamSkau1992}\index{names}{Herman, Richard H.}\index{names}{Putnam, Ian F.}\index{names}{Skau, Christian F.} succeeded in showing that every minimal Cantor dynamical
system is isomorphic to a Bratteli-Vershik dynamical system.

\subsection{Bratteli diagrams}
The BV representation theorem (Theorem~\ref{ch5:theo:BVmodel})
saying that $(X,T)$ can be topologically realized as a BV-dynamical system
is  the main result of \cite{HermanPutnamSkau1992}\index{names}{Herman, Richard H.}\index{names}{Putnam, Ian F.}\index{names}{Skau, Christian F.}. 
We recall  that A. M. Vershik obtained in \cite{Vershik:1985}\index{names}{Vershik, Anatol M.} such a result in a measure-theoretic context.

Theorem~\ref{theoremDGBratteli} is due to Elliott \citep{Elliott1976}
(see also \cite{Krieger1980b}).
\subsection{Kakutani equivalence}
Kakutani equivalence was introduced
in the context of measure-theoretic systems~\citep{Kakutani1943}.
\index{names}{Kakutani, Shizuo}
Theorem~\ref{theoremKakutaniEquiv} is  from~\cite{GiordanoPutnamSkau1995}.
Exercises~\ref{exerciseDerivativePrimitive} and~\ref{exerciseKakutani}
are from~\cite{OrnsteinRudolphWeiss1982}
\index{names}{Ornstein, Donald S.}\index{names}{Rudolph, Daniel J.}%
\index{names}{Weiss, Benjamin}%
 (the context is for measure-theoretic systems but the arguments transpose easily).
\subsection{Strong orbit equivalence}

The   Strong Orbit Equivalence theorem (Theorem~\ref{ch5:th:GPS}) is
due to \cite{GiordanoPutnamSkau1995}\index{names}{Giordano, Thierry}\index{names}{Putnam, Ian F.}\index{names}{Skau, Christian F.}.
We follow the proof proposed in \cite{GlasnerWeiss1995},
\index{names}{Glasner, Eli}\index{names}{Weiss, Benjamin}%
giving more details.
 
Exercise \ref{exerciseBoyleOrbitEquiv}
is due to   \cite{Boyle1983} (see~\cite[Theorem 2.4]{GiordanoPutnamSkau1995}).
\index{names}{Boyle, Michael}%
See also~\cite{BoyleTomiyama1998}.
\index{names}{Tomiyama, Jun}%
Two systems $(X,T)$ and $(Y,S)$ such that $(X,T)$ is conjugate to
$(Y,S)$ or to $(Y,S^{-1})$ are called \emph{flip equivalent}
\index{subject}{flip equivalence}.
\subsection{Equivalences on Cantor spaces}
The consideration of etale equivalence relations to
formulate results on orbit equivalence in topological
dynamical systems
is due to \cite{Putnam2010}. Etale
equivalence relations are a  particular case of the notion
of \emph{\'etale groupoid}\index{subject}{etale@\'etale!groupoid}
 introduced by \cite{Renault2014}.
\index{names}{Renault, Jean}%
We follow here the beautiful book of
\cite{Putnam2018}.
\subsection{Entropy}
We suggest \cite{Walters1982}\index{names}{Walters, Peter} or
\cite{LindMarcus1995} for an introduction to topological entropy.

The fact that there exist minimal systems of arbitrary entropy
(Theorem~\ref{theoremGrillenberger}) is due
to~\cite{Grillenberger1972}
\index{names}{Grillenberger, Christian}%
  (we give in Exercises \ref{exerciseGrillenberger1}
and \ref{exerciseGrillenberger2} the proof that the entropy
of a minimal system 
  can be arbitrary large).
The fact that every minimal Cantor dynamical system is SOE
to a minimal Cantor dynamical system of entropy zero
(Theorem~\ref{theoremBoyleHandelmanEntropyZero})
 is from \cite{Boyle&Handelman:1994}.
\index{names}{Boyle, Michael}\index{names}{Handelman, David}%

In \cite{Boyle&Handelman:1994}\index{names}{Boyle, Michael}\index{names}{Handelman, David} the authors use a different lemma
instead of Lemma~\ref{ch5:lemma:majentropy}.
They show that 
 $h(S_W) = \frac{\log m}{l}$, whenever $W$ is a set of $m$ distinct finite words of length $l$.

The original reference to Fekete's Lemma (Exercise~\ref{exerciseEntropyLimit})
is \citep{Fekete1923}.

Theorem \ref{theoremBoyleHandelmanEntropyZero}
 shows that there can be dynamical systems with different entropies   in a strong orbit equivalence class.
Hence it is natural to ask whether all entropies can be realised inside a given class.
M. Boyle\index{names}{Boyle, Michael} and D. Handelman\index{names}{Handelman, David} showed in \citep{Boyle&Handelman:1994}\index{names}{Boyle, Michael}\index{names}{Handelman, David}
that  it is true in the class of the odometer $(\Z_2 , x\mapsto x+1 )$.
Later, F. Sugisaki\index{names}{Sugisaki, Fumiani} proved in \citep{Sugisaki2003}\index{names}{Sugisaki, Fumiani}  that it is true in any strong orbit equivalence class.
Theorem \ref{theoremSugisaki}, showing that the realizations can be chosen
to be subshifts, is from  \cite{Sugisaki2007}.
\subsection{Exercises}
The solution of Exercise \ref{exerciseEquivalenceTelescoping} has
been kindly provided to us by Ian Putnam.
Exercise \ref{exerciseC*Equivalence} 
is from~\cite{BratteliJorgensenKimRoush2000}
where a proof that the matrices of Equation~\eqref{eqMNC*EquivNotShiftEquiv}
are not shift equivalent is given (we have not included it here
as it uses technical tools from algebraic number theory
which are out of our scope).
\index{names}{Bratteli, Ola}\index{names}{J{\o}rgensen, Palle E.}%
\index{names}{Kim, Ki Hang}\index{names}{Roush, Fred W.}%

The decidability of $C^*$-equivalence of matrices was proved
in~\cite{BratteliJorgensenKimRoush2001}. \index{subject}{decidability!of $C^*$-equivalence of matrices}
It proves
the decidability of the telescoping equivalence of stationary
diagrams. \index{subject}{decidability!of telescoping equivalence!of stationary diagrams}
Actually, the
equivalence of Bratteli diagrams is undecidable in general,
\index{subject}{undecidability!of telescoping equivalence} as shown
by \cite{MundiciPanti1993} (see also \cite{Mundici2003}).
\index{names}{Mundici, Daniele}\index{names}{Panti, Giovanni}%
The undecidability is proved by a reduction from the problem
of isomorphism of lattice ordered groups, proved undecidable
by \cite{GlassMadden1984}.
\index{names}{Glass, Andrew M.}\index{names}{Madden, James J.}%
The latter being itself obtained
by a reduction from the problem of 
piecewise linear
homeomorphism  for compact polyhedra
proved undecidable by \cite{Markov1958}.
\index{names}{Markov, Alexander A.}%

%%%%%%%%%%%%%%%%%%%%%%%%%%%%%%%
\chapter{Substitution shifts and  generalizations}
%%%%%%%%%%%%%%%%%%%%%%%%%%%%%%%
\label{ch5:sec:examples}\label{chapterSubstitutionShifts}

In this chapter, we describe the  BV-representations for substitution
shifts and their generalizations, namely linearly recurrent
shifts and $\Sa$-adic shifts. We  treat before
the case of odometers which will be needed after and show that
they can be represented by Bratteli diagrams
with one vertex at each level. 

We describe next the construction of
BV-representations for substitution shifts. The main
result is Theorem~\ref{ch5:subsec:Bratteli-substitution}
which states that
the family  of Bratteli-Vershik systems associated with stationary,
properly ordered Bratteli diagrams is (up to isomorphism) the disjoint union
of the family of infinite substitution minimal systems and the family  of stationary
odometer systems.

In the next sections, we treat the cases of linearly recurrent shifts (Section~\ref{ch5:subsec:LR})
and of $\Sa$-adic shifts (Section~\ref{sectionSadicShifts}).
In Section \ref{sectionDGSadic}, we give a description of the
dimension group of unimodular $\Sa$-adic shifts (Theorem~\ref{theo:dg}).
This will be applied, in Chapter~\ref{chapterDendricShifts},
to the case of dendric shifts. Finally, in Section~\ref{sectionDerivatives},
we prove a result characterizing substitutive sequences 
by a finiteness property of the set of their derivatives (Theorem~\ref{theoremCharacterisationSubstitutive}).

%%%%%%%%%%%%%%%%%%%%%%
\section{Odometers}
%%%%%%%%%%%%%%%%%%%%%%
\label{ch5:subsec:rep-odo}

Let $(p_n)_{n\ge 1}$ be a strictly increasing sequence 
of natural integers such that $p_n$ divides $p_{n+1}$ for all $n\ge 1$.
We endow the set $X=\prod_{n\ge 1} \Z / p_n \Z$ with the product topology of the discrete topologies. 
The set
$$
\Z_{(p_n)} = \{ (x_n)_{n\ge 1} \in X \mid x_{n} \equiv x_{n+1} \mod p_n \} 
$$
\index{symbols}{Z@$\Z_{(p_n)}$}%
is   a group for componentwise addition,
called the group of $(p_n)$-\emph{adic integers}.
\index{subject}{p-adic@$(p_n)$-adic integer}%
\index{subject}{integer!p-adic@$p_n$-adic}%
It is a compact topological group (see Exercise~\ref{ch5:ex:odom}).
A basis of the topology is given by the sets 
$$
[a_1 , a_2 , \ldots , a_m] := \{ (x_n) \in \Z_{(p_n)} \mid x_i = a_i , 1\leq i\leq m  \} \ .
$$
The neutral element is $0=(0,0,\ldots)$. 

As a variant of the definition, one starts with an arbitrary sequence 
$(q_n)_{n\ge 1}$
of natural integers $q_n\ge 2$ and considers the group
\begin{displaymath}
Y=\{(y_n)_{n\ge 0}\mid 0\le y_n\le q_{n+1}-1\}
\end{displaymath}
 of $(q_n)$-\emph{adic expansions}
\index{subject}{q-adic@$(q_n)$-adic expansions}%
 using for addition the sum with carry from
left to right as in an expansion of a real number with respect to an integer
basis. This is equivalent to the previous definition using $q_1=p_1$ and
$q_{n+1}=p_{n+1}/p_{n}$ for $n\ge 1$ (Exercise~\ref{exerciseExpansionBasispn}).

When $p_n = p^n$ for all $n\ge 1$, this defines the classical group
 of $p$-\emph{adic integers} $\Z_p$.
\index{subject}{p-adic@$p$-adic!integer}%
\index{subject}{integer!p-adic@$p$-adic}
The corresponding odometer has already been met several times
(see Example~\ref{examplepAdic}).

When $p_n=n!$ for all $n\ge 1$,
the group $\Z_{(n!)}$ is called the group of \emph{profinite integers}.
\index{subject}{profinite!integer}%
\index{subject}{integer!profinite}%
The corresponding expansion of integers denoted
$x=(\cdots c_3c_2c_1)_!$
if $x=c_1+c_22!+c_33!+\ldots$ with $0\le c_i<i$ is called
the \emph{factorial number system}\index{subject}{factorial!number system}.

Let $T : \Z_{(p_n)} \to \Z_{(p_n)}$ be the map $x\mapsto x+1$
where $1=(1,1,\ldots)$.
The pair  $(\Z_{(p_n)} , T)$ is called {\em odometer in base} $(p_n)$\index{subject}{odometer!base of}\index{subject}{base!odometer}\index{subject}{odometer}. It is a minimal dynamical system (indeed, the orbit of $0$ is dense
and $0$ is in the closure of the orbit of every point).

The odometer in base $n!$ is called the \emph{universal odometer}.
\index{subject}{universal!odometer}\index{subject}{odometer!universal}%

Let us compute the BV-representation of an odometer.
For all $n$, we set $B(n) = \{[0^{n-1}]\}$, $t(n)=1$, $h (n)=p_n$ and 
\begin{equation}
\Pg (n)=\{ T^j B (n) \mid 0\leq j \leq h(n) -1 \}.\label{equationPatitionsOdometer}
\end{equation}

Then  $(\Pg (n) )_n$ is a refining sequence of  KR-partitions  (see Exercise~\ref{ch5:ex:KR}).
Remark that $T^j B (n) = [j_0 j_1 \cdots j_{n-1}]$ where $j_i = j \mod p_i$.
The edges of the BV-representation of $(\Z_{(p_n)} , T)$ given in Section~\ref{sectionBVmodelTheorem} are of the form $(n,1,1,l)$, with $0\leq l \leq q_n -1= \frac{p_n}{p_{n-1}} -1$.

Thus an odometer has a BV-representation with one vertex at each level.
The converse is also clearly true. We can thus state the following
simple nice result.
\begin{theorem}\label{theoremBVrepresenationOdometers}
A Cantor dynamical system is an odometer if and only if it has
a BV-representation with one vertex at each level.
\end{theorem}

\begin{figure}[ht]
        \centering
\tikzset{node/.style={circle,draw,minimum size=0.1cm,inner sep=0cm}}
\tikzset{title/.style={minimum size=0cm,inner sep=0pt}}
\begin{tikzpicture}
\node[title]at(0,3){$V(0)$};\node[title]at(2,3){$V(1)$};
\node[title]at(4,3){$V(2)$};\node[title]at(6,3){$V(3)$};
\node[node](0)at((0,2){};\node[node](1)at((2,2){};
\node[node](2)at((4,2){};\node[node](3)at((6,2){};

\draw[bend right](0)edge node{}(1);\draw[bend left](0)edge node{}(1);
\draw[bend right=45](1)edge node{}(2);\draw[bend right=25](1)edge node{}(2);
\draw(1)edge node{}(2);\draw[bend left=25](1)edge node{}(2);
\draw[bend left=45](1)edge node{}(2);
\draw[bend left](2)edge node{}(3);\draw(2)edge node{}(3);
\draw[bend right](2)edge node{}(3);

\node[title]at(1,1){$E(1)$};\node[title]at(3,1){$E(2)$};
\node[title]at(5,1){$E(3)$};
\end{tikzpicture}
                
        \caption{The three first levels  of the BV representation of $(\Z_{(p_n)} , R)$.}
        \label{ch5:fig:ch5-diagodo}
\end{figure}
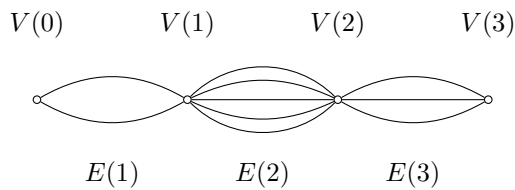

For example if $p_1 = 2$, $p_2 = 10$ and $p_3=30$, the first  three levels  are given in Figure~\ref{ch5:fig:ch5-diagodo}.

\begin{proposition}\label{propositionDimensionGroupOdometer}
The dimension group of the odometer in base $(p_n)$ is the subgroup
of $\Q$ formed of the $p/q$ with $q$ dividing some $p_n$.
\end{proposition}
\begin{proof}
The dimension group of $(\Z_{(p_n)},T)$ is the direct limit
of the sequence
\begin{displaymath}
\Z\edge{p_1}\Z\edge{p_2}\Z\cdots
\end{displaymath}
whence the statement. Note that the order unit is $1$.
\end{proof}
For example, the dimension group of the odometer in base $(2^n)$
is the group of dyadic rationals.

We will use in the next section
the following statement, which identifies 
 systems equivalent to an odometer in a non obvious way.

\begin{proposition}\label{propositionOdometerRankOne}
Let $(V,E,\le)$ be an ordered 
 Bratteli diagram such that
\begin{enumerate}
\item[\rm(i)] The vertex set for $n\ge 1$ is $V(n)=\{1,2,\ldots,t\}$.
\item[\rm(ii)] The incidence matrices $M(n)$
are all equal for $n\ge 2$ to a square matrix of rank $1$.
\item[\rm(iii)] The order is such that $e\le e'$ if $r(e)=r(e')$
and $s(e)\le s(e')$.
\end{enumerate}
 Then $(X_E,T_E)$ is topologically
conjugate to an odometer in base $(qp^n)_n$.
\end{proposition}
\begin{proof}
Set $M=M(n)$ for $n\ge 2$. Since $M$ has rank $1$, we
have  $M=xy$ with $x$ a nonnegative integer column vector and $y$ a 
nonnegative integer row vector. Consider the ordered Bratteli diagram
$(V',E',\le)$ with incidence matrices $M'(1)=M(1)$ and
$M'(2n)=x$, $M'(2n+1)=y$ for $n\ge 1$. By condition (iii), we can choose
the order on $(V',E',\le')$
 such that the telescoping  with respect to $0,1,3,5,\ldots$
gives $(V,E,\le)$. Now the telescoping of $G'$ with respect to $0,2,4,\ldots$
is the odometer in basis $(pq^n)_{n\ge 0}$ with $p=yM(1)$ and $q=yx$.
\end{proof}

\begin{figure}[hbt]
\centering
\tikzset{node/.style={circle,draw,minimum size=0.1cm,inner sep=0cm}}
\tikzset{title/.style={minimum size=0cm,inner sep=0pt}}
 \begin{tikzpicture}
\node[node](0)at(1,6){};
\node[node](11)at(0,5){};\node[node](12)at(2,5){};
\node[node](21)at(0,4){};\node[node](22)at(2,4){};
\node[node](31)at(0,3){};\node[node](32)at(2,3){};
\node[node](41)at(0,2){};\node[node](42)at(2,2){};
\node[title](51)at(0,1){};\node[title](52)at(2,2){};

\draw(0)edge node{}(11);\draw(0)edge node{}(12);
\draw[left,near end](11)edge node{$0$}(21);\draw[right,near end](11)edge node{$0$}(22);
\draw[left,near end](12)edge node{$1$}(21);\draw[right,near end](12)edge node{$1$}(22);
\draw[left,near end](22)edge node{$1$}(31);\draw[right,near end](22)edge node{$1$}(32);
\draw[left,near end](21)edge node{$0$}(31);\draw[right,near end](21)edge node{$0$}(32);
\draw[left,near end](32)edge node{$1$}(41);\draw[right,near end](32)edge node{$1$}(42);
\draw[left,near end](31)edge node{$0$}(41);\draw[right,near end](31)edge node{$0$}(42);
\draw[dotted](41)edge node{}(51);\draw[dotted](42)edge node{}(52);

\node[node](0')at(5,6){};
\node[node](1)at(5,4.5){};
\node[node](2)at(5,3.5){};
\node[node](3)at(5,2.5){};

\draw[bend left](0')edge node{}(1);\draw[bend right](0')edge node{}(1);
\draw[bend left](1)edge node{}(2);\draw[bend right](1)edge node{}(2);
\draw[bend left](2)edge node{}(3);\draw[bend right](2)edge node{}(3);
\end{tikzpicture}

                \caption{The Bratteli diagram with matrix $M$.}
        \label{figureOdometerRankOne}
\end{figure}
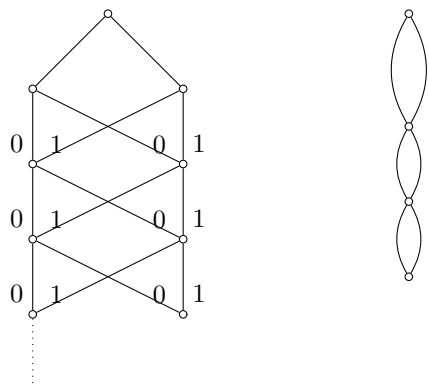

For example,  the odometer represented in Figure~\ref{figureOdometerRankOne}
has incidence matrices
\begin{displaymath}
M(1)=\begin{bmatrix}1\\1\end{bmatrix},\quad M(n)=\begin{bmatrix}1&1\\1&1\end{bmatrix}=\begin{bmatrix}1\\1\end{bmatrix}\begin{bmatrix}1&1\end{bmatrix}
\end{displaymath}
for $n\ge 2$.

Note finally the following simple but important observation.

\begin{proposition}\label{propositionOdometersSubstitutions}
The family of odometers is disjoint from the family of minimal shift spaces.
\end{proposition}
\begin{proof}
Any odometer $\Z_{(p_n)}$ is infinite and has a partition in towers with a basis
formed of one element, namely $\Z_{(p_n)}=\cup_{j=0}^{h-1} T^jB$
where $h=p_1$ and $B=[0]$. This is not possible for a minimal shift space $(X,S)$.
 Indeed, assume that $X=\cup_{j=0}^{h-1}T^jB$
where $B$ is a clopen set. Then there is an integer $n$ and a partition of $\cL_n(X)$
in $h$ sets $W_0,\ldots,W_{h-1}$ such that $T^jB=[W_j]$ for $0\le j<h$. Let $w\in W_0$.
Then the periodic point $\cdots ww\cdot ww\cdots$ belongs to $X$. Since $(X,S)$
is minimal, it must be periodic and thus finite, a contradiction.

\end{proof}
A different proof uses the notion of expansive system.
A dynamical system $(X,T)$, endowed with the distance $d$, is {\em expansive}
\index{subject}{expansive dynamical system}%
\index{subject}{dynamical system!expansive}%
 if there exists $\epsilon $ such that for all pairs of points $(x, y)$, $x\not = y$, there exists $n$ with $d(T^n x , T^n y)\geq \epsilon $.
We say that  $\epsilon $ is a {\em constant of expansivity}\index{subject}{constant!of expansivity} of $(X,T)$. Expansivity is a property
invariant by conjugacy (but the constant may not be the same).

A shift space is expansive while an odometer is not. Actually,
a topological dynamical system is a shift space if and only if it is expansive
(Exercise~\ref{exerciseExpansive2}).

At the opposite of expansive systems,
a topological dynamical system $(X,T)$, endowed with the distance $d$,
 is  said to be {\it equicontinuous}
 \index{subject}{equicontinuous!system}%
 \index{subject}{dynamical system!topological!equicontinuous}%
 whenever for every $\epsilon>0$, there exists $\delta>0$ such that
 $$ d(x,y)< \delta\Rightarrow  \sup_{n\in \Z }d(T^n x , T^n y)< \epsilon . $$

 For example, an odometer is equicontinuous (with $\delta=\epsilon$).

 A Bratteli diagram has the {\em equal path number property}
 \index{subject}{equal path number property} if for all $n\geq 1$ and $u,v\in V(n)$ we have $\Card(r^{-1} (u)) = \Card(r^{-1} (v))$.
Note that this implies that the number of paths from $u,v$
to $v(0)$ are the same.

Note that the equal path number property is preserved by telescoping.
 
Note also that the representation of odometers given in Section~\ref{ch5:subsec:rep-odo} shares  this property. We prove the following result
which connects all these notions together.

\begin{theorem}[Gjerde Johansen]\label{theoremEqualPathNumber}
A minimal shift is Toeplitz if and only if it has a
 BV-representation $(X_E , T_E)$ where  
$(V,E,\le )$ has the equal path number property.
\end{theorem}
\begin{proof}
Let $(X_E,T_E)$
be a BV-representation of the minimal shift $(X,S)$, associated to the sequence of refining KR-partition $(\Pg (n))$, which has
the equal path number property. 
Let $\phi:X_E\to X$ be a conjugacy.
Taking $(\Pg (n))_{n\geq n_0}$ instead of $(\Pg (n))$ preserves the equal path number property. 
Consequently, we may assume that the partitions $\Pg(n)$ are such
that all points $x,y$ in $\phi (B(n))$ satisfy $x_{[-n,n]} = y_{[-n,n]}$.
Then, the sequence $x=\phi(x_{min})$ is a Toeplitz sequence. Indeed,
one has $x_{n+kp}=x_n$ for all $k\in\Z$ with $k=\Card(r^{-1} (\{ u\}) )$ for all
$u\in V(n)$.

Conversely, let $(X,S)$ be a Toeplitz shift and let $x\in X$ be a Toeplitz sequence. 
For a clopen set $U\subset X$ and $y\in X$ we set 
\[
{\rm Per}_p(y,U)= \{ n \in \mathbb{Z} | S^{n+kp} (y) \in U \text{ for all } k\in \mathbb{Z} \} .
\] 
We first claim that for all $n\in \N$ there exist $p>0$ and a clopen partition $\{C, S^{-1}C, \dots , S^{-p+1} C\}$ such that $C$ is included in the cylinder set $U=[x_{[-n,-1]}.x_{[0,n]}]$.

The sequence $x$ being Toeplitz there exists $p>0$ such that ${\rm Per}_p(x,U)$ is nonempty.
We can suppose that $p$ is minimal,
that is, if $q$ is such that $0\leq q <p$ then ${\rm Per}_q(x,U)$ is empty.

The set $\{y\in X | {\rm Per}_p(y,U) \not = \emptyset \}$ being nonempty, closed and $S$-invariant, it is equal to $X$ by minimality.
Thus, the set ${\rm Per}_p (y , U)$ is nonempty for all $y\in X$.

Let $C$ be the closed subset  $\{y\in X | {\rm Per}_p (y , U) = {\rm Per}_p (x , U)\}$.
Observe that $\{C, S^{-1}C, \dots , S^{-p+1} C\}$ is a partition 
as $p$ is minimal. Since $S^{-p}C=C$, the nonempty
closed set $\cup_{i=0}^{p-1} S^{-i} C$ is
$S$-invariant set. Thus, it is equal to $X$.
The set $C$ being closed,  it is a clopen partition of $X$.
Moreover, because $0$ belongs to ${\rm Per}_p(x,U)$, it also belongs to ${\rm Per}_p(y,U)$ for all $y\in C$. 
Hence $y$ belongs to $U$ and $C$ is included in $U$. 
This proves the claim.

Let us now proceed as in Theorem \ref{theoremKRPartitions} with the exception that  the sequence of clopen sets $C_n$ will not be chosen at the initial
step but defined 
as adequate cylinder sets $[u.v]$ including $x$. 

We start choosing an increasing sequence of partitions $(\Pg'(n))_n$ generating the topology.
Let $U_1 = [x_{[-n_1 , n_1]}]$, $n_1 = 1$ and let $C_1$ be as in the claim (for the period $p_1=p$).

We apply Proposition~\ref{propositionVersik} to $\Qg=\Pg'(1)$
and $C=C_1$ to obtain the KR-partition $\Pg(1)$.
Observe that the height of each tower is $p_1$. 
The base is $C_1=\cup B_i \subset U_1$ where the $B_i$'s are atoms of  $\Pg(1)$.
One can suppose that $x$ belongs to $B_1$ and thus there exists some $n_2 > n_1$ such that the cylinder set $U_2=[x_{[-n_2 , -1]} . x_{[0,n_2]}]$ is included in $B_1$. 

The sequence $x$ being Toeplitz there exists $p_2>0$, taken minimal, such that $0$ belongs to ${\rm Per}_{p_2} (x, U_2 )$.
We apply the claim to obtain the clopen set $C_2$ and the period $p_2$.

Applying Proposition~\ref{propositionVersik}
iteratively for $n\ge 2$ to $C=C_n$ and by setting now
\begin{displaymath}
\Qg=\Pg'(n)\vee\Pg(n-1),
\end{displaymath}
we obtain KR-partition $\Pg(n)$ with basis $C_n$ which is finer than $\Pg'(n)$ and $\Pg(n-1)$ and whose heights are all equal to $p_n$ for some $p_n$.
Moreover, the sequence $(\Pg(n))_n$ of KR-partition is nested, generates the topology and the intersection of the bases is included in $\cap_n U_n = \{ x\}$. 

Thus the BV-representation corresponding to this partition has the equal path number property. 
\end{proof}

The class of BV-systems having the equal path number property is
larger than that of Toeplitz shifts.
Indeed, one can show that there exist BV-systems having the equal path number property that are neither expansive nor equicontinuous (see the Notes for a reference).

\begin{example}
Let $\sigma:0\to 01,1\to 00$ be the substitution
generating the period-doubling sequence (see Example~\ref{exampleToeplitz})
which is a Toeplitz shift.
A BV-representation of the corresponding shift is shown in Figure~\ref{figurePeriodDoubling}. It can be verified that it has the equal path number
property.
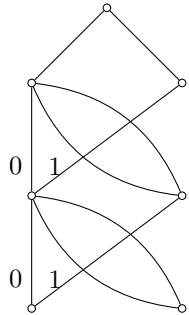
\begin{figure}[hbt]
\centering
\tikzset{node/.style={circle,draw,minimum size=0.1cm,inner sep=0cm}}
\tikzset{title/.style={minimum size=0cm,inner sep=0pt}}
\begin{tikzpicture}
\node[node](0)at(2,5){};
\node[node](11)at(1,4){};\node[node](12)at(3,4){};
\node[node](21)at(1,2.5){};\node[node](22)at(3,2.5){};
\node[node](31)at(1,1){};\node[node](32)at(3,1){};

\draw(0)edge node{}(11);\draw(0)edge node{}(12);
\draw[left,near end](11)edge node{$0$}(21);\draw[left,near end](12)edge node{$1$}(21);
\draw[bend right,left,near end](11)edge node{}(22);
\draw[bend left,right,near end](11)edge node{}(22);
\draw[left,near end](21)edge node{$0$}(31);\draw[left,near end](22)edge node{$1$}(31);
\draw[bend right,left,near end](21)edge node{}(32);
\draw[bend left,right,near end](21)edge node{}(32);
\end{tikzpicture}
\caption{The BV-representation of the period doubling shift.}
\label{figurePeriodDoubling}
\end{figure}
The orbit of $x_{\min}$ is shown in Figure~\ref{figurePeriodDoublingOrbit}.
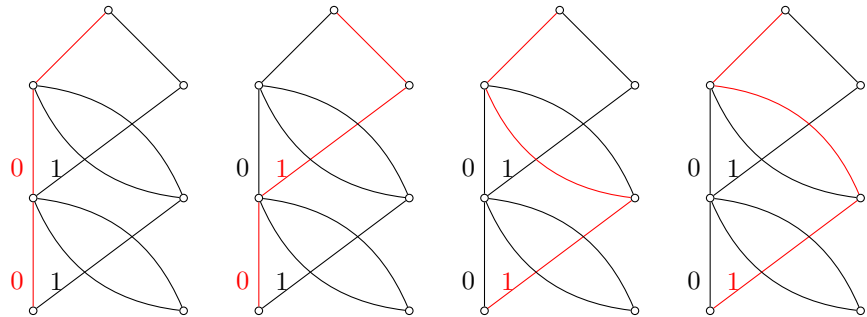
\begin{figure}[hbt]
\centering
\tikzset{node/.style={circle,draw,minimum size=0.1cm,inner sep=0cm}}
\tikzset{title/.style={minimum size=0cm,inner sep=0pt}}
\begin{tikzpicture}
%1
\node[node](0)at(2,5){};
\node[node](11)at(1,4){};\node[node](12)at(3,4){};
\node[node](21)at(1,2.5){};\node[node](22)at(3,2.5){};
\node[node](31)at(1,1){};\node[node](32)at(3,1){};

\draw[color=red](0)edge node{}(11);\draw(0)edge node{}(12);
\draw[left,near end,color=red](11)edge node{$0$}(21);\draw[left,near end](12)edge node{$1$}(21);
\draw[bend right,left,near end](11)edge node{}(22);
\draw[bend left,right,near end](11)edge node{}(22);
\draw[left,near end,color=red](21)edge node{$0$}(31);\draw[left,near end](22)edge node{$1$}(31);
\draw[bend right,left,near end](21)edge node{}(32);
\draw[bend left,right,near end](21)edge node{}(32);
%2
\node[node](0)at(5,5){};
\node[node](11)at(4,4){};\node[node](12)at(6,4){};
\node[node](21)at(4,2.5){};\node[node](22)at(6,2.5){};
\node[node](31)at(4,1){};\node[node](32)at(6,1){};

\draw(0)edge node{}(11);\draw[color=red](0)edge node{}(12);
\draw[left,near end](11)edge node{$0$}(21);\draw[left,near end,color=red](12)edge node{$1$}(21);
\draw[bend right,left,near end](11)edge node{}(22);
\draw[bend left,right,near end](11)edge node{}(22);
\draw[left,near end,color=red](21)edge node{$0$}(31);\draw[left,near end](22)edge node{$1$}(31);
\draw[bend right,left,near end](21)edge node{}(32);
\draw[bend left,right,near end](21)edge node{}(32);
%3
\node[node](0)at(8,5){};
\node[node](11)at(7,4){};\node[node](12)at(9,4){};
\node[node](21)at(7,2.5){};\node[node](22)at(9,2.5){};
\node[node](31)at(7,1){};\node[node](32)at(9,1){};

\draw[color=red](0)edge node{}(11);\draw(0)edge node{}(12);
\draw[left,near end](11)edge node{$0$}(21);\draw[left,near end,](12)edge node{$1$}(21);
\draw[bend right,left,near end,color=red](11)edge node{}(22);
\draw[bend left,right,near end](11)edge node{}(22);
\draw[left,near end](21)edge node{$0$}(31);\draw[left,near end,color=red](22)edge node{$1$}(31);
\draw[bend right,left,near end](21)edge node{}(32);
\draw[bend left,right,near end](21)edge node{}(32);
%4
\node[node](0)at(11,5){};
\node[node](11)at(10,4){};\node[node](12)at(12,4){};
\node[node](21)at(10,2.5){};\node[node](22)at(12,2.5){};
\node[node](31)at(10,1){};\node[node](32)at(12,1){};

\draw[color=red](0)edge node{}(11);\draw(0)edge node{}(12);
\draw[left,near end](11)edge node{$0$}(21);\draw[left,near end,](12)edge node{$1$}(21);
\draw[bend right,left,near end](11)edge node{}(22);
\draw[bend left,right,near end,color=red](11)edge node{}(22);
\draw[left,near end](21)edge node{$0$}(31);\draw[left,near end,color=red](22)edge node{$1$}(31);
\draw[bend right,left,near end](21)edge node{}(32);
\draw[bend left,right,near end](21)edge node{}(32);
\end{tikzpicture}
\caption{The orbit of $x_{\min}$.}
\label{figurePeriodDoublingOrbit}
\end{figure}
\end{example}
%%%%%%%%%%%%%%%%%%%%%%%%%%
\section{Substitutions}\label{sectionSubstitutionsBV}
%%%%%%%%%%%%%%%%%%%%%%%%%%

We will now consider the BV-representation of substitution shifts.
Let us note that we work  in all this section with bi-infinite
 words and two-sided shifts.
We  first need a new definition.

A Bratteli diagram $(V,E)$ is {\em stationary}
\index{subject}{Bratteli diagram!stationary}%
\index{subject}{stationary!Bratteli diagram}%
 if
there exists $k$ such that 
$k= \Card(V (n))$  for all $n$,  and if (by an appropriate labeling of the vertices)
the incidence matrices between level $n$ and $n+1$ are the same $k \times k$
matrix $M$ for all $n=1,2,\ldots $.  
In other words, beyond level $1$ the
diagram repeats itself. Clearly we may label the vertices in $V (n)$ as $v(n,a_1),
\cdots , v(n, a_k)$, where $A=\{ a_1, \ldots, a_k\}$ is a set of $k$
distinct symbols. The matrix $M$ is called the \emph{matrix of the stationary
diagram}.
\index{subject}{matrix!of a stationary diagram}%

The ordered Bratteli diagram $(V,E, \le)$ is  {\em stationary}\index{subject}{stationary!ordered Bratteli diagram} if $(V,E)$
is stationary, and the ordering on the edges with range $v(n,a_i)$ is the
same as the ordering on the edges with range $v(m,a_i)$ for $m,n=2,3,\ldots$
and $i=1,\ldots, k$.  In other words, beyond level $1$ the diagram with the
ordering repeats itself.

An odometer is \emph{stationary}\index{subject}{stationary!odometer} if it has
a stationary BV-representation. 
It is easy to see that the odometer in base $(p_n)$
is stationary if and only if $p_n=pq^n$ for some $p,q\ge 2$.

Let $(V,E,\le)$ be a stationary properly ordered Bratteli diagram.
The morphism read on $E(n)$ is constant from $n\ge 2$. 
We call it the {morphism read on}\index{subject}{morphism!read on a stationary Bratteli diagram}\index{subject}{Bratteli diagram!morphism read on} $(V,E,\le)$.

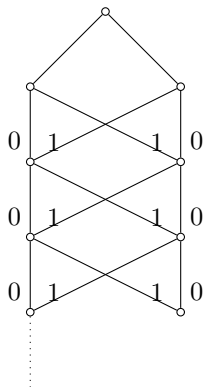
\begin{figure}[ht]
        \centering
\tikzset{node/.style={circle,draw,minimum size=0.1cm,inner sep=0cm}}
\tikzset{title/.style={minimum size=0cm,inner sep=0pt}}
 \begin{tikzpicture}
\node[node](0)at(1,6){};
\node[node](11)at(0,5){};\node[node](12)at(2,5){};
\node[node](21)at(0,4){};\node[node](22)at(2,4){};
\node[node](31)at(0,3){};\node[node](32)at(2,3){};
\node[node](41)at(0,2){};\node[node](42)at(2,2){};
\node[title](51)at(0,1){};\node[title](52)at(2,2){};

\draw(0)edge node{}(11);\draw(0)edge node{}(12);
\draw[left,near end](11)edge node{$0$}(21);\draw[right,near end](11)edge node{$1$}(22);
\draw[left,near end](12)edge node{$1$}(21);\draw[right,near end](12)edge node{$0$}(22);
\draw[left,near end](22)edge node{$1$}(31);\draw[right,near end](22)edge node{$0$}(32);
\draw[left,near end](21)edge node{$0$}(31);\draw[right,near end](21)edge node{$1$}(32);
\draw[left,near end](32)edge node{$1$}(41);\draw[right,near end](32)edge node{$0$}(42);
\draw[left,near end](31)edge node{$0$}(41);\draw[right,near end](31)edge node{$1$}(42);
\draw[dotted](41)edge node{}(51);\draw[dotted](42)edge node{}(52);
\end{tikzpicture}
                
                \caption{The Thue--Morse substitution read on
a Bratteli diagram.}
        \label{ch5:fig:ch5-diagmorse}
\end{figure}

Given a substitution $\sigma$, it is not possible in general to use the stationary
Bratteli diagram $(V,E,\le)$, with $\sigma$ read on $(V,E,\le)$, to represent
the shift space corresponding to $\sigma$.
For example, in the case of the Thue-Morse morphism
%\index{subject}{Thue-Morse!substitution}%
%\index{subject}{substitution!Thue-Morse}%
\index{names}{Morse, Marston} $a\mapsto ab$, $b\mapsto ba$ the Bratteli diagram  is given in Figure~\ref{ch5:fig:ch5-diagmorse}.
\index{subject}{Thue-Morse!Bratteli diagram}%

It is clear that it has two maximal and two minimal paths.
Hence this representation does not give a properly ordered Bratteli diagram.

\subsection{Main result}

We will show that one has however the following result. It implies in particular that every
infinite minimal substitution shift has a BV-representation.

\begin{theorem} 
\label{ch5:subsec:Bratteli-substitution}
The family $\mathcal{B}$ of Bratteli-Vershik systems associated with stationary,
properly ordered Bratteli diagrams is (up to isomorphism) the disjoint union
of the family of infinite substitution minimal systems and the family  of stationary
odometer systems. 
\end{theorem}

Furthermore, we will see that the correspondence in question is given by
an explicit and algorithmically effective construction.

A morphism $\sigma$ on the alphabet $A$ is
{\em proper}\index{subject}{substitution!proper}\index{subject}{proper!substitution}\index{subject}{morphism!proper} if there are two letters $r,l\in A$ such
that, for every $a\in A$, $r$ is the last letter of $\sigma(a)$ and $l$ is the first letter of $\sigma(a)$. It is called \emph{eventually proper}
\index{subject}{eventually!proper substitution}%
\index{subject}{substitution!eventually proper}%
if there is an integer $p\ge 1$ such that $\sigma^p$ is proper.
An eventually proper morphism $\sigma$ has exactly one fixed point which is
$\sigma^\omega(\ell\cdot r)$.

\begin{proposition}\label{propositionPrimitiveProper}
The morphism read on a stationary, properly ordered, Bratteli diagram
is primitive and eventually proper.
\end{proposition}
\begin{proof}
Let $\sigma$ be the  morphism read on $(V,E,\le)$.
Since $(V,E,\le)$ is properly ordered, it is simple. For every
$a,b\in A$, since $(V,E,\le)$ is simple, there is a path from $(1,b)$
to some $(n,a)$. Then $b$ occurs in $\sigma^n(a)$ showing that $\sigma$  is primitive.

Let $i(a)$ be the first letter of $\sigma(a)$.
For every $a\in A$ and $n\ge 1$, the source of 
the minimum edge with range $(n,a)$
is $(n-1,i(a))$. Thus, if the minimal rank of the maps $i^n:A\to A$
were larger than $1$, there would exist more than one minimal
infinite path, a contradiction with the hypothesis that $(V,E,\le)$
is properly ordered. This shows that there exists $n$ such that
$i^n(a)$ is the same for all $a\in A$. A symmetric argument holds
for the last letter. Thus $\sigma$ is eventually proper.
\end{proof}
We recall that a shift space $X$ is said to be {\em periodic}\index{subject}{periodic shift}
\index{subject}{shift space!periodic}%
if there exist $x \in X$ and  an integer $k$ such that $ X= \{ x , Sx , \ldots ,S^{k-1} x\}$. Thus a shift space $X$ is periodic if
and only if the dynamical system $(X,S)$ is periodic.
Likewise, it is said to be {\em aperiodic}\index{subject}{aperiodic!shift space}
\index{shift space!aperiodic}%
if the system is aperiodic, that is, does not contain any periodic point.
\index{subject}{subshift!periodic}\index{subject}{subshift!aperiodic}%
Thus a periodic shift is the same as a minimal finite shift space. Observe
that the property of being periodic is decidable (see Exercise~\ref{exerciseFixedPointPeriodic2}).
\subsection{Diagrams with simple hat}

We say  that a Bratteli diagram $(V,E)$ has a {\em simple hat}\index{subject}{simple!hat} whenever it has only simple edges between the top vertex and any vertex of
the first level. Note that the Bratteli diagram associated
with a nested sequence of partitions as in Section~\ref{sectionBVmodelTheorem}
has a simple hat.

The following result gives a proof of one direction of Theorem~\ref{ch5:subsec:Bratteli-substitution} in the particular case of diagrams with a simple hat.
\begin{proposition}\label{ch5:proposition:substitutionread}
Let $(V,E,\le)$ be a stationary, properly ordered Bratteli diagram with a simple hat, 
  let $\sigma: A^*\to A^*$  be the morphism read on $(V,E,\le)$,
and let $(X,S)$ be the substitution shift associated to $\sigma$.
\begin{enumerate}
\item
If $(X, S)$ is not periodic, then it is isomorphic to $(X_E,T_E)$.
\item
If $(X, S)$ is periodic,
then $(X_E,T_E)$ is isomorphic to an odometer in base $(qp^n)_n$, for some $p,q\ge 2$.
\end{enumerate}
\end{proposition}
We will use the following lemma.
\begin{lemma}\label{lemmaPartitionProperSubstitution}
Let $\sigma:A^*\to A^*$ be a primitive and eventually proper
morphism. If $X(\sigma)$ is infinite, the family
\begin{displaymath}
\Pg(n)=\{T^j\sigma^n([a])\mid a\in A,\ 0\le j<|\sigma^n(a)|\}
\end{displaymath}
is a refining sequence of partitions in towers.
\end{lemma}
\begin{proof}
By Proposition~\ref{propositionPartitionSubstitution} the 
partition
$\Pg(n)$ is for every $n\ge 1$ a KR-partition of $X(\sigma)$ with basis $\sigma^n(X)$.
The sequence $\Pg(n)$ is clearly nested. Since
$\sigma$ is eventually proper,
the intersection of the bases
is reduced to one point, namely the unique fixed point of $\sigma$. 
Finally, $\Pg(n)$ tends to the partition in points. Indeed, 
let $p\ge 1$ be such that all $\sigma^p(a)$ begin with $\ell$
and end with $r$. Then (see Figure~\ref{figurePartitionGenerates})
\begin{displaymath}
\sigma^n([a])=[\sigma^{n-p}(r)\cdot\sigma^n(a)\sigma^{n-p}(\ell)].
\end{displaymath}
Thus all words in $T^j\sigma^n([a])$ coincide on $[-m,m]$
for $m=\min_{b\in B}|\sigma^{n-p}(b)|$.
This shows that $\Pg(n)$ is a refining sequence of KR-partitions.

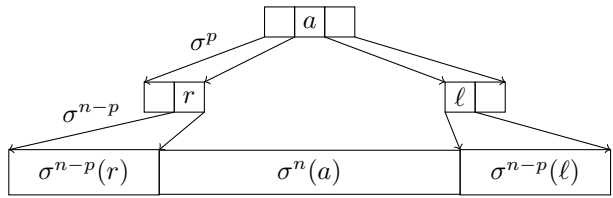
\begin{figure}[hbt]
\centering
\tikzset{node/.style={draw,minimum size=0.4cm,inner sep=0pt}}
	\tikzset{title/.style={minimum size=0.5cm,inner sep=0pt}}
\begin{tikzpicture}
\node[node]at(3.6,2){};\node[node]at(4,2){$a$};\node[node]at(4.4,2){};
\node[node]at(2,1){};\node[node]at(2.4,1){$r$};
\node[node]at(6,1){$\ell$};\node[node]at(6.4,1){};
\node[draw, minimum width=2cm]at(1,0){$\sigma^{n-p}(r)$};\node[draw,minimum width=4cm]at(4,0){$\sigma^{n}(a)$};
\node[draw, minimum width=2cm]at(7,0){$\sigma^{n-p}(\ell)$};

\draw[left,above,->](3.4,1.8)-- node{$\sigma^{p}$}(1.8,1.2);
\draw[->](3.8,1.8)-- node{}(2.6,1.2);
\draw[->](4.2,1.8)-- node{}(5.8,1.2);
\draw[->](4.6,1.8)-- node{}(6.6,1.2);
\draw[left,above,->](2.2,.8)-- node{$\sigma^{n-p}$}(0,.3);
\draw[->](2.6,.8)-- node{}(2,.3);
\draw[->](5.8,.8)-- node{}(6,.3);
\draw[->](6.2,.8)-- node{}(8,.3);
\end{tikzpicture}
\caption{The sequence $\Pg(n)$ generates the topology.}
\label{figurePartitionGenerates}
\end{figure}
\end{proof}
\begin{proofof}{of Proposition~\ref{ch5:proposition:substitutionread}}
Assume first that $(X,S)$ is aperiodic. By Proposition~\ref{propositionPrimitiveProper},
the morphism $\sigma$ is primitive and eventually proper.
By Lemma~\ref{lemmaPartitionProperSubstitution}, the family $\Pg(n)$
is a refining sequence of partitions in towers.

Since $(V,E,\le)$ has simple hat, the Bratteli diagram associated to the sequence of partitions $\Pg(n)$
is clearly equal to $(V,E,\le)$. Thus, by Theorem~\ref{ch5:theo:BVmodel},
$(X,S)$ is isomorphic with $(X_E,T_E)$.

Assume now that $(X,S)$ is periodic. Replacing $\sigma$
by some power does not modify $(X,S)$ (and replaces $(V,E,\le)$ by a periodic
telescoping) so that we may assume
that all words $\sigma(a)$ for $a\in A$ begin with
$\ell$
and end with $r$. The unique fixed point 
$x=\sigma^\omega(r\cdot\ell)$ of $\sigma$
is then periodic. Set $x=\cdots ww\cdot www\cdots$ with $w$ 
as short as possible. Then $w$ is a \emph{primitive} word,
\index{subject}{primitive!word}\index{subject}{word!primitive}%
that is, it is not a power of a shorter word. Let $n$ be large enough so that
$|\sigma^n(a)|\ge|w|$ for every $a\in A$. Then each
$\sigma^{2n}(a)$ is a word of period $|w|$ which
\begin{itemize}
\item begins with $w$ since it begins with $\sigma^n(\ell)$,
\item ends with $w$ since it ends with $\sigma^n(r)$.
\end{itemize}
Since $w$ is primitive, this forces each $\sigma^{2n}(a)$ to be a power
of $w$. Thus, replacing again $\sigma$ be one of its powers,
we may assume that every $\sigma(a)=w^{k_a}$ is a power of $w$.
In this case, the incidence matrix of the diagram $(V,E,\le)$
is such that for every $a,b\in A$
\begin{displaymath}
M_{a,b}=k_a|w|_b
\end{displaymath}
Thus the hypotheses of Proposition
\ref{propositionOdometerRankOne} are satisfied and $(X_E,T_E)$
is an odometer in basis $qp^n$ with $q=|w|$ and $p=\sum_{a,b\in A}k_a|w|_b$.
\end{proofof}
We illustrate the periodic case in the following example.
\begin{example}
Let $(V,E,\le)$ be the stationary Bratteli diagram represented in Figure~\ref{figureBratteliPeriodic}
on the left.
\begin{figure}[hbt]
\centering
\begin{tikzpicture}
\tikzset{node/.style={circle,draw,minimum size=0.4cm,inner sep=0pt}}
\tikzset{title/.style={minimum size=0.4cm,inner sep=0pt}}
\node[node](0)at(1,3){};
\node[node](11)at(0,2){$a$};\node[node](12)at(1,2){$b$};\node[node](13)at(2,2){$c$};
\node[node](21)at(0,1){$a$};\node[node](22)at(1,1){$b$};\node[node](23)at(2,1){$c$};
\node[node](31)at(0,0){$a$};\node[node](32)at(1,0){$b$};\node[node](33)at(2,0){$c$};
\node[title]at(1,-.5){$\vdots$};

\draw[left](0)edge node{}(11);\draw[left](0)edge node{}(12);\draw[left](0)edge node{}(13);
\draw[left](11)edge node{$1$}(21);\draw[left](11)edge node{$1$}(23);
\draw[left](12)edge node{$2$}(21);\draw[left](12)edge node{$2$}(23);
\draw[right](13)edge node{}(22);\draw[right](13)edge node{$3$}(23);
\draw[left](21)edge node{$1$}(31);\draw[left](21)edge node{$1$}(33);
\draw[left](22)edge node{$2$}(31);\draw[left](22)edge node{$2$}(33);
\draw[right](23)edge node{}(32);\draw[right](23)edge node{$3$}(33);
%%%
\tikzset{node/.style={circle,draw,minimum size=0.1cm,inner sep=0pt}}
\node[node](0)at(6,3){};\node[node](1)at(6,1.5){};\node[node](2)at(6,.5){};\node[node](3)at(6,-.5){};
\node[title]at(6,-1){$\vdots$};

\draw[bend left](0) edge node{}(1);\draw(0) edge node{}(1);\draw[bend right](0) edge node{}(1);
\draw[bend left=40](1) edge node{}(2);\draw[bend left=15](1) edge node{}(2);
\draw[bend right=15](1) edge node{}(2);\draw[bend right=40](1) edge node{}(2);
\draw[bend left=40](2) edge node{}(3);\draw[bend left=15](2) edge node{}(3);
\draw[bend right=15](2) edge node{}(3);\draw[bend right=40](2) edge node{}(3);
\end{tikzpicture}
\caption{The periodic case.}\label{figureBratteliPeriodic}
\end{figure}
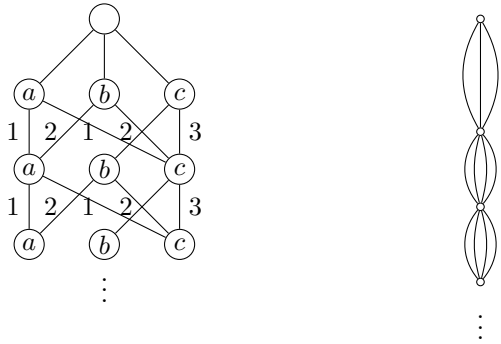
The morphism read on $(V,E,\le)$ is $\sigma:a\to ab,b\to c,c\to abc$. The substitution
shift defined by $\sigma$ is periodic since $\sigma(abc)=(abc)^2$. The corresponding
odometer is represented in Figure~\ref{figureBratteliPeriodic} on the right.
\end{example}

 A \emph{partition matrix} is a $t\times\ell$-matrix $P$
with coefficients $0,1$
\index{subject}{partition!matrix}\index{subject}{matrix!partition}%
 such that every column of $P$ has exactly one coefficient equal to $1$.
Such a matrix defines a partition $\theta$ of the set $\{1,2,\ldots,\ell\}$
of indices of the columns putting together $x,y$ if $P_{i,x}=P_{i,y}=1$.
It also defines a map $\pi$ from $\{1,2,\ldots,\ell\}$ onto $\{1,2,\ldots,k\}$
by $\pi(x)=i$ if $i$ is the index such that $P_{i,x}=1$.
\begin{example}
The matrix
\begin{displaymath}
P=\begin{bmatrix}1&0&0\\0&1&1\end{bmatrix},
\end{displaymath}
is a partition matrix. The corresponding partition $\theta$
is $\{1\},\{2,3\}$ and the map $\pi$ is $1\to 1,2\to 2,3\to 2$.
\end{example}

The following result will allow us to reduce to the case of a Bratteli diagram
with a simple hat.
\begin{proposition}\label{propositionLemma15Forrest}
For every stationary and properly ordered Bratteli diagram $(V,E,\le)$,
there is a stationary properly ordered Bratteli diagram $(V',E',\le')$
with a simple hat such that $(X_E,T_E)$ and $(X_{E'},T_{E'})$
are isomorphic.

More precisely, for every ordered
Bratteli diagram $B$, there is an ordered
Bratteli diagram $B'$ such that the following holds.
There exist two nonegative matrices $P,Q$,
with $P$ a partition matrix,
such that the matrices $M,M'$ of $B$ and $B'$ satisfy
\begin{displaymath}
M=PQ,\quad M'=QP.
\end{displaymath}
and  the vector $v=M(1)$ is such that $v=Pw$
where $w$ has all its components equal to $1$.
Moreover, if $B$ is properly ordered, then $B'$ is properly ordered.
\end{proposition}
\begin{proof}
Let $M$ be the $t\times t$-matrix equal to the incidence matrices $M(n)$
of the diagram $B$,
for all $n\ge 2$ and let $v=M(1)$. Set $\ell=\sum_{i=1}^tv_i$. 

Let
$P$ be the $t\times\ell$ partition matrix
 defined by 
\begin{displaymath}
P_{i,j}=\begin{cases}1&\mbox{ if $\sum_{k<i}v_k<j\le\sum_{k\le i}v_k$}\\
0&\mbox{ otherwise}.
\end{cases}
\end{displaymath}
Thus the rows of $P$ are the characteristic vectors of elements of the
partition of $\{1,2,\ldots,\ell\}$ into the $t$ sets 
$V_1=\{1,2,\ldots,v_1\}$, ...,$V_t=\{\ell-v_t+1,\ldots,\ell\}$.
We may assume, replacing if necessary $M$ by some power,
that, for every $i$, the set of edges with range $i$ has
at least $v_i$ elements.
We choose an $\ell\times t$ matrix $Q$ such that $M=PQ$.
This is equivalent to splitting the set of edges entering
 the vertex $i$ in $v_i$ nonempty subsets (the sum of the rows of $Q$
with index in $V_i$ is then the row of index $i$ of $M$).
The edges in each subset keep the order induced by the order on $B$.
Let $B'$ be the Bratteli diagram with incidence matrices 
$M'(1)=w=\begin{bmatrix}1\ 1\ldots 1\end{bmatrix}^t$ and
$M'(n)=QP$ for $n\ge 2$. We order the diagram $B'$
by the order induced by that of $B$. Since $v=Pw$,
$B$ and $B'$ can both be obtained by telescoping from the Bratteli 
diagram $C$ with incidence matrices $(w,P,Q,P,Q,\ldots)$.
If $B$ is properly ordered, $C$ is properly ordered and
consequently $B'$ also.
\end{proof}

We illustrate the construction in the following example.

\begin{example}
Let $(V,E,\le)$ be the Bratteli diagram represented in Figure~\ref{figureSplit}
on the left.
\begin{figure}[hbt]
\centering
\begin{tikzpicture}
\tikzset{node/.style={circle,draw,minimum size=0.1cm,inner sep=0pt}}
\tikzset{title/.style={minimum size=0.4cm,inner sep=0pt}}
\node[node](0)at(1,3){};
\node[node](11)at(0,2){};\node[node](12)at(2,2){};
\node[node](21)at(0,1){};\node[node](22)at(2,1){};
\node[node](31)at(0,0){};\node[node](32)at(2,0){};
\node[title]at(1,0){$\vdots$};

\draw(0)edge node{}(11);
\draw[bend left](0)edge node{}(12);\draw[bend right](0)edge node{}(12);
\draw(11)edge node{}(21);\draw(11)edge node{}(22);
\draw(12)edge node{}(21);\draw(12)edge node{}(22);
\draw(21)edge node{}(31);\draw(21)edge node{}(32);
\draw(22)edge node{}(31);\draw(22)edge node{}(32);
%%%%%%%%%%%%%%%%%%%%
\node[node](0)at(6,3){};
\node[node](11)at(5,2.5){};\node[node](12)at(6,2.5){};\node[node](13)at(7,2.5){};
\node[node](21)at(5.5,2){};\node[node](22)at(6.5,2){};
\node[node](31)at(5,1.5){};\node[node](32)at(6,1.5){};\node[node](33)at(7,1.5){};
\node[node](41)at(5.5,1){};\node[node](42)at(6.5,1){};
\node[node](51)at(5,.5){};\node[node](52)at(6,.5){};\node[node](53)at(7,.5){};
\node[node](61)at(5.5,0){};\node[node](62)at(6.5,0){};

\draw(0)edge node{}(11);\draw(0)edge node{}(12);\draw(0)edge node{}(13);
\draw(11)edge node{}(21);\draw(12)edge node{}(22);\draw(13)edge node{}(22);
\draw(21)edge node{}(31);\draw(21)edge node{}(32);
\draw(22)edge node{}(31);\draw(22)edge node{}(33);
\draw(31)edge node{}(41);\draw(32)edge node{}(42);\draw(33)edge node{}(42);
\draw(41)edge node{}(51);\draw(41)edge node{}(52);
\draw(42)edge node{}(51);\draw(42)edge node{}(53);
\draw(51)edge node{}(61);\draw(52)edge node{}(62);\draw(53)edge node{}(62);
%%%%%%%%%%
\node[node](0)at(10,3){};
\node[node](11)at(9,2.5){};\node[node](12)at(10,2.5){};\node[node](13)at(11,2.5){};
\node[node](21)at(9,1.5){};\node[node](22)at(10,1.5){};\node[node](23)at(11,1.5){};
\node[node](31)at(9,0.5){};\node[node](32)at(10,0.5){};\node[node](33)at(11,0.5){};

\draw(0)edge node{}(11);\draw(0)edge node{}(12);\draw(0)edge node{}(13);
\draw(11)edge node{}(21);\draw(11)edge node{}(22);
\draw(12)edge node{}(21);\draw(12)edge node{}(23);
\draw(13)edge node{}(21);\draw(13)edge node{}(23);
\draw(21)edge node{}(31);\draw(21)edge node{}(32);
\draw(22)edge node{}(31);\draw(22)edge node{}(33);
\draw(23)edge node{}(31);\draw(23)edge node{}(33);
\end{tikzpicture}
\caption{The transformation of $(V,E,\le)$ by splitting vertices.}
\label{figureSplit}
\end{figure}
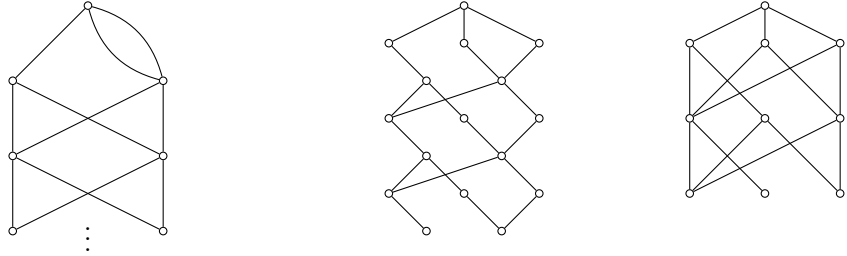
The matrices $M,P,Q,M'$ are
\begin{displaymath}
M=\begin{bmatrix}1&1\\1&1\end{bmatrix},\quad
P=\begin{bmatrix}1&0&0\\0&1&1\end{bmatrix},\quad
Q=\begin{bmatrix}1&1\\1&0\\0&1\end{bmatrix},\quad
M'=\begin{bmatrix}1&1&1\\1&0&0\\0&1&1\end{bmatrix}.
\end{displaymath}
The diagram with matrices $(w,P,Q,P,\ldots)$ is represented in the middle
of Figure~\ref{figureSplit} and the diagram with simple hat $(V',E',\le')$ 
with matrices $w,M',M',\ldots)$ on the right.
\end{example}

\subsection{A useful result}

The following result will be the key to build a BV-representation of
an aperiodic minimal substitution shift. It allows one
to replace 
a pair $(\tau,\phi)$ of morphisms with $\tau$ primitive 
 by a pair $(\zeta,\theta)$ with $\zeta$
primitive and $\theta$ letter-to-letter (see Figure~\ref{figureRauzyProper}).
\begin{figure}[hbt]
\centering
\tikzset{node/.style={circle,draw,minimum size=0.1cm,inner sep=0pt}}
\tikzstyle{every loop}=[->,shorten >=1pt]
\tikzstyle{loop above}=[in=40,out=130,loop]
\begin{tikzpicture}
\node(B)at(0,2){$B^\Z$};
\node(C)at(3,2){$C^\Z$};
\node(A)at(1.5,0){$A^\Z$};

\draw[above](B)edge[loop above]node{$\tau$}(B);
\draw[above](C)edge[loop above]node{$\zeta$}(C);
\draw[left,->](B)edge node{$\phi$}(A);
\draw[right,->](C)edge node{$\theta$}(A);
\draw[above,->](B)edge node{$\gamma$}(C);
\end{tikzpicture}
\caption{The eventually proper substitution $\zeta$}\label{figureRauzyProper}
\end{figure}
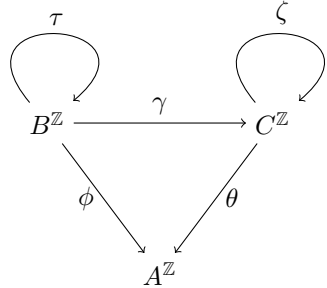
\begin{proposition}
\label{ch5:prop:rauzy}
Let $y\in B^\Z$ be an admissible fixed
point of a primitive substitution $\tau$ on the alphabet $B$ and let
 $\phi: B^*\to A^*$ be a  nonerasing morphism. Set $x=\phi(y)$ and let $X$ be the
subshift spanned by $x$. 

There exist a primitive substitution $\zeta$ on an alphabet $C$, an
admissible fixed point $z$ of $\zeta$, a morphism
$\gamma:B^*\to C^*$ and a map $\theta: C\to
A$ such that:
\begin{enumerate}
\item
\label{ch5:prop:rauzy:item1}
$\theta(z)=x$ and $\phi=\theta\circ\gamma$.
\item
\label{ch5:prop:rauzy:item2}
If  $\phi$ is injective and $\phi(B)$ is a circular code,
\index{subject}{circular code}%
 then $\theta$ is a conjugacy from $(X(\zeta),S)$ onto $X$,
\item
\label{ch5:prop:rauzy:item3}
If $\tau$ is eventually proper, then $\zeta$ is eventually proper.
\end{enumerate}
\end{proposition}

\begin{proof} 
The proof below is very simple, but the notation is, in an
unavoidable way, a bit heavy.  
By substituting  $\tau $   by  a power of  itself if needed, we can assume
that $|\tau(b)|\geq |\phi(b)|$ for all $b\in B$.  We define:
\begin{itemize}
\item
an alphabet $C$ by  $C =\left\{  b_p\mid b\in B, \, 
1\leq p\leq |\phi(b)| \right\} $,
\item
a map $\theta: C\to A$  by  $\theta(b_p)= \left(\phi(b)\right)_p $ ,
\item
a map $\gamma : B\to C^+ $  by  $\gamma (b)=b_1b_2\cdots b_{|\phi(b)|}$.
\end{itemize}
Clearly $\theta \circ \gamma = \phi$.
We define a substitution $\zeta$ on $C$  as follows.
For $b$ in $B$ and $1\leq p\leq |\phi(b)|$,  we set
$$
\zeta(b_p) =\begin{cases} 
\gamma \left( \left( \tau(b) \right)_p \right) & \hbox{ if }  1\leq p<|\phi(b)| \\
\gamma \left( \left( \tau(b) \right)_{[\phi(b),\tau(b)]} \right) & \hbox{ if }  p=|\phi(b)|  \ . 
\end{cases}
$$

Hence, for every $b\in B$, $\zeta \left( \gamma(b)\right) =
\zeta(b_1)\cdots\zeta(b_{|\phi(b)|})=\gamma \left( \tau(b)\right)$, {\em i.e.}, 
\begin{align}
\label{ch5:prop:rauzy:eq1}
\zeta\circ\gamma=\gamma \circ \tau
\end{align}
and it follows that  
$$
\zeta^n \circ \gamma = \gamma \circ \tau^n \hbox{ for all }n\geq 0 \ .
$$ 
We claim that $\zeta$ is primitive. 
Let $n$ be an
integer such that $b$ occurs in $\tau^n(a)$ for all $a,b\in B$. Let
$b_p$ and $c_q$ belong to $C$. 
By construction, $\zeta(b_p)$
contains $\gamma\left( \tau(b)_p\right)$ as a factor, thus
$\zeta^{n+1}(b_p)$ contains  $\zeta^n\left( \gamma\left( \tau(b)_p\right)\right)=
\gamma\left(\tau^n\left(\tau(b)_p\right)\right)$ as a factor. 
By the
choice of $n$, $c$ occurs in $\tau^n\left(\tau(b)_p\right)$, thus
$\gamma(c)$ is a factor of
$\gamma\left(\tau^n\left(\tau(b)_p\right)\right)$, and also of 
$\zeta^{n+1}(b_p)$. 
Since $(c_q)$ is a letter of $\gamma(c)$, $c_q$
occurs in  $\zeta^{n+1}(b_p)$ and our claim is proved.

\emph{Proof of 1.}
Let $z=\gamma(y)$. 
By (\ref{ch5:prop:rauzy:eq1})  we get $\zeta(z)=\gamma(\tau(y))=\gamma(y)=z$,
and $z$ is a fixed point of $\zeta$. By construction, $z$ is
uniformly recurrent, thus it is an admissible fixed point of
$\zeta$. Moreover, $\theta(z)=\theta(\gamma(y))=\phi(y)=x$, and
\ref{ch5:prop:rauzy:item1} is proved.

{\em Proof of \ref{ch5:prop:rauzy:item2}.}
Since $\theta$ commutes with the shift and maps $z$ to $x$,  and by minimality of the subshifts,  it maps $X(\zeta)$ onto $X$. 
There remains to prove that $\theta : X(\zeta) \to X$ is one-to-one. 
Let $\alpha\in X$. 
By definition of $X$, there exist $t\in X_\tau$ and an integer $p$,
with $0\leq p<|\phi(t_0)|$, such that $\alpha=S^p\phi(t)$. Let
$\beta$ be an element of $X(\zeta)$ with
$\theta(\beta)=\alpha$. 
By definition of $\gamma$, there exist some
$\delta\in X(\tau)$ and some integer $q$, with $0\leq q<|\gamma(\delta_0)|$,
such that $\beta=S^q\gamma(\delta)$.
It follows that
$S^q\phi(\delta)=\theta(\beta)=\alpha=S^p\phi(t)$.
Since $0\leq q<|\gamma (\delta_0)|=|\phi(\delta_0)|$ by construction of $\gamma$, 
since $\phi(B)$ is a circular code and since $\phi$ is injective, it follows that $\delta=t$ and $q=p$,
thus $\beta=S^p\gamma(t)$: $\beta$ is uniquely determined by
$\alpha$, and $\theta$ is one-to-one.

{\em Proof of \ref{ch5:prop:rauzy:item3}. }
Let $l\in B$ be  the first letter of $\tau(b)$ for every $b\in
R$. 
Let $b_p\in C$, and $c=\tau(b)_p$. 
By definition of $\zeta$,
the first letter of $\zeta(b_p)$ is $c_1$, and the first letter of
$\zeta^2(b_p)$ is the first letter of $\zeta(c_1)$, that is, $l_1$. 
By the same method, if $r$ is the last letter of $\tau(b)$ for every
$b\in B$, then the last letter of $\zeta^2(b_p)$ is $r_{|\phi(r)|}$ for
every $b_p\in C$.
\end{proof}

Note that the proof of Proposition~\ref{ch5:prop:rauzy} is very close to that of
Proposition \ref{propositionLemma15Forrest}
(see Exercise~\ref{exerciseForrestRauzy}).

%Given a primitive substitution $\sigma$ on the alphabet $A$ such that $(X_\sigma ,S)$ is aperiodic, let us describe how to construct a primitive proper substitution $\zeta$ such that
%$(X_\sigma , S)$ is isomorphic to $(X_\zeta , S)$.
%With the techniques used below, we do not need to suppose primitivity. 
%We only need $\sigma$ to generate  a minimal subshift.
%For example, this includes the Chacon substitution\index{subject}{Chacon!substitution} $0\mapsto 0010$ and $1\mapsto 1$.
%An illustration of the construction   given below is  given  in  Section~\ref{ch5:subsubsec:Chacon}.

%%%%%%%%%%%%%%%%%%%%%%%%%%%%
%\subsection{Return words}
%%%%%%%%%%%%%%%%%%%%%%%%%%%%

%In order to find $\zeta$ we need to introduce the notion of return words. For more details,
% the reader is referred {\em e.g.}, \cite{Durand:1998b}, \cite{Durand:2000}, \cite{Durand:2003}.
%We define an {\em occurrence of $u.v$ in  the  bi-infinite word $x$}\index{subject}{occurrence} to be an integer $n$ such that
%$x_{[n-|u|,n+|v|-1]}=uv$.
%A finite word $w$ on $A$ is a {\em return word to $u.v$}\index{subject}{return word}
%in $x$ if there exist two consecutive occurrences $j,k$ of $u.v$ in $x$
%such that $w=x_{[j,k-1]}$.\marginpar{eint}

%We need one more lemma concerning return words.
%\begin{lemma}\label{lemmaReturn}
%Let $\sigma:A\to A^*$ be a substitution genrerating an aperiodic minimal
%shift and let $x$ be an admissible fixed point of $\sigma$.
%\end{lemma}
Consider  a  morphism $\sigma:A\to A^*$ generating an aperiodic
minimal shift $X$. Replacing $\sigma$ by one of its powers, we can
choose an admissible fixed point $x$ of $\sigma$. Set $r=x_{-1}$ and $l=x_0$.
Note that since $x$ is a fixed point of $\sigma$, $\sigma(r)$ ends with
$r$ and $\sigma(l)$ begins with $l$.

Let $\RR_X(rl)$ be the set of right return words
\index{subject}{return!word}%
to $rl$. Every word in $\RR_X(rl)$ ends with $l$. Set 
\begin{displaymath}
\RR_X(r\cdot l)=l\RR_X(rl)l^{-1}
\end{displaymath}
\index{symbols}{R@$\RR_X(r\cdot l)$}%
which is the set of words of the form $lu$ for $ul\in\RR_X(rl)$.
Thus $w$ is in $\RR_X(r\cdot l)$ if and only if $rwl$ is in $\cL(X)$ and contains
exactly two occurrences of $rl$, one as a prefix and one as a suffix.
In particular all words in $\RR_X(r\cdot l)$ begin with $l$ and end with $r$.
Moreover, the set $\RR_X(r\cdot l)$ satisfies the
following properties.
\begin{enumerate}
\item[(i)] any word in $\cL(X)$ which begins with $l$ and ends with $r$
is a concatenation of words of $\RR_X(r\cdot l)$.
\item[(ii)] The set $\RR_X(r\cdot l)$ is a circular code and, more precisely, no word in $r\RR_X(r\cdot l) l$
overlaps non trivially a product of words of $\RR_X(r\cdot l)$.
\end{enumerate}
Since $X$ is minimal, $\RR_X(r\cdot l)$ is finite. Let $\phi:B\to \RR_X(r\cdot l)$ be
a coding morphism for $\RR_X(r\cdot l)$.

Since $x$ is uniformly recurrent, it has an infinite number
of occurrences of $rl$ at positive and at negative indices.
By the above properties of $\RR_X(r\cdot l)$, there is
a unique element $y$ of $B^\Z$ such that $\phi (y)=x$.

We now define as follows a morphism $\tau$ on the alphabet $B$.
For every $b\in B$, $\sigma\circ \phi(b)$ begins with $l$ and ends with $r$.
This implies that $\sigma\circ \phi(b)\in\RR_X(r\cdot l)^*$ and thus that
$\sigma\circ \phi(b)=\phi(w)$ for some unique $w\in B^*$.
We set $\tau(b)=w$. 

This defines a morphism $\tau$ on the alphabet $B$, characterised by
\begin{equation}
\phi\circ\tau=\sigma\circ\phi \ . \label{equationPhiTauSigma}
\end{equation}
It follows that $\phi\circ\tau^n=\sigma^n\circ\phi$ for each $n\geq 0$.

\begin{example}\label{exampleTauFibonacci}
Consider the Fibonacci morphism $\varphi:a\to ab,b\to a$
generating the Fibonacci shift $X(\varphi)$.
The sequence $x=\varphi^{2\omega}(a\cdot a)$ is an admissible fixed point
of $\varphi^2:a\to aba,b\to ab$. We have
\begin{displaymath}
\RR_X(a\cdot a)=\{aba,ababa\}.
\end{displaymath}
Set $B=\{a,b\}$ with $\phi(a)=aba$ and $\phi(b)=ababa$. Then
\begin{displaymath}
\phi\circ\tau(a)=\varphi^2\circ\phi(a)=\varphi^2(aba)=abaababa=\phi(ab)
\end{displaymath}
so that $\tau(a)=ab$. Similarly, we find $\tau(b)=abb$.
\end{example}
\begin{proposition}\label{propositiontauProperAperiodic}
The morphism $\tau$ 
defined by \eqref{equationPhiTauSigma}
is primitive, eventually proper and aperiodic.
\end{proposition}
 \begin{proof}
Let us show that $\tau$ is eventually proper. For this, let $n$ be such that $|\sigma^n(l)|>|\phi(y_0)|$.
Then for every $b\in B$, the first letter of  $\tau^n(b)$ , is $y_0$. Indeed,
$\phi\circ\tau^n(b)=\sigma^n\circ\phi(b)$ begins with $\sigma^n(l)$ 
(see Figure~\ref{figureTauProper}).
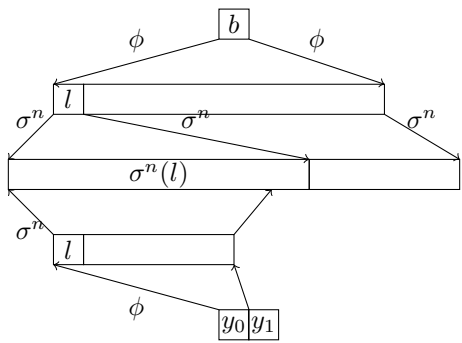
\begin{figure}[hbt]
\centering
\tikzset{node/.style={draw,minimum size=0.4cm,inner sep=0pt}}
	\tikzset{title/.style={minimum size=0.5cm,inner sep=0pt}}
\begin{tikzpicture}
\node[node](b)at (4,4){$b$};
\node[node,minimum width=.4cm](l)at(1.8,3){$l$};
\node[node,minimum width=4cm]at(4,3){};
\node[node,minimum width=4cm]at(3,2){$\sigma^n(l)$};
\node[node,minimum width=2cm]at(6,2){};
\node[node]at(1.8,1){$l$};
\node[node,minimum width=2cm]at(3,1){};
\node[node]at(4,0){$y_0$};\node[node,left]at(4.6,0){$y_1$};

\draw[->,above](3.8,3.8)-- node{$\phi$}(1.6,3.2);
\draw[->,above](4.2,3.8)-- node{$\phi$}(6,3.2);
\draw[->,above](1.6,2.8)-- node{$\sigma^n$}(1,2.2);
\draw[->,above](2,2.8)-- node{$\sigma^n$}(5,2.2);
\draw[->,above](6,2.8)-- node{$\sigma^n$}(7,2.2);
\draw[->,below](1.6,1.2)-- node{$\sigma^n$}(1,1.8);
\draw[->,below](4,1.2)-- node{}(4.5,1.8);
\draw[->,left,below](3.8,.2)-- node{$\phi$}(1.6,.8);
\draw[->,left,below](4.2,.2)-- node{}(4,.8);
\end{tikzpicture}
\caption{Proof that $\tau$ is eventually proper.}\label{figureTauProper}
\end{figure}
But  $\sigma^n(l)$ is a prefix of $x=\sigma^n(x)=\sigma^n\circ\phi(y)$. Thus, by the choice of $n$,
$\sigma^n\circ\phi(y_0)l$ is a prefix of $\sigma^n(l)$. This implies
that $\tau^n(b)$ begins with $y_0$. A similar argument shows that if $n$ is such
that $|\sigma^n(r)|>|\phi(y_{-1})|$, then every $\tau^n(b)$ end with $y_{-1}$.
Thus $\tau$ is eventually proper.

Let $k>0$ be an occurrence of $r.l$  large enough so that   every return word
$w\in\mathcal{R}$   appears  in the decomposition of $x_{[0,k)}$, that is,   every $b\in
B$   occurs in the finite word $u\in B^+$ defined by $\phi(u)=x_{[0,k)}$.
Let $n$ be so large that $|\sigma^n(l)|>k$. Let $b,c\in B$. As above,
$x_{[0,k)}l$ is a prefix of $\sigma^n(l)$, which is a prefix of
$\sigma^n(\phi(b))=\phi(\tau^n(b))$. Thus $u$ is a prefix of $\tau^n(b)$, and
$c$ occurs in $\tau^n(b)$, hence $\tau$ is primitive.

Moreover,
$$
\phi(\tau(y))=\sigma(\phi(y))=\sigma(x)=x=\phi(y) \ , 
$$
thus $\tau(y)=y$ since $\RR_X(r\cdot l)$ is circular, and $y$
is the unique fixed point of $\tau$. 
Since $\phi(y)=x$ is not periodic,
$y$ is not periodic. Thus $\tau$ is aperiodic.
\end{proof}
We are now ready to prove Theorem~\ref{ch5:subsec:Bratteli-substitution}.\\

\begin{proofof}{of Theorem~\ref{ch5:subsec:Bratteli-substitution}}
Let first $(V,E,\le)$ be a stationary properly ordered Bratteli diagram. By Proposition
\ref{propositionLemma15Forrest}, we may assume that $(V,E,\le)$ has a simple hat. Then,
using Proposition~\ref{ch5:proposition:substitutionread}, we conclude 
that $(X_E,T_E)$ is either an aperiodic minimal substitution shift or an odometer.
This proves the theorem in one direction.

Let us now prove the converse implication. If $(X,T)$ is an odometer, we have seen (Section~\ref{ch5:subsec:rep-odo})
that $(X,T)$ has a BV-representation with  a stationary, properly ordered, Bratteli diagram.
Moreover, $(X,T)$ cannot be at the same time a minimal substitution shift
by Proposition~\ref{propositionOdometersSubstitutions}.

Let finally $\sigma:A^*\to A^*$ be a substitution generating
a minimal shift space $X$. Let $\tau:B^*\to B^*$ be the
primitive, eventually proper and aperiodic
substitution defined by Proposition~\ref{propositiontauProperAperiodic}.

Let $\zeta:C\to C^*$ be the substitution given by Proposition~\ref{ch5:prop:rauzy}
with an admissible fixed point $z\in C^\Z$ and a map $\theta:C\to A$. 
Since $\phi$
is injective and $\RR_X(r\cdot l)=\phi(B)$ is a circular code, by assertion 2,
 $\theta$
is an isomorphism from $(X_\zeta,S)$ onto $X$. 
Moreover, by assertion 3, since
$\tau$ is eventually proper, $\zeta$ is eventually proper. 
Let $(V,E,\le)$ be the stationary properly ordered Bratteli diagram
such that $\zeta$ is the morphism read on $(V,E,\le)$. By Proposition~\ref{ch5:proposition:substitutionread},
the systems $(X_\zeta,S)$ and $(X_E,T_E)$ are isomorphic. This concludes the proof.

\end{proofof}
%%%%%%%%%%%%%%%%%%%%%%%%%%%%%%%%%%%%%%%%%%%%%%%%%%%%%
%\subsection{An example: the Chacon substitution}\label{ch5:subsubsec:Chacon}
%%%%%%%%%%%%%%%%%%%%%%%%%%%%%%%%%%%%%%%%%%%%%%%%%%%%%

Note that  the preceding results imply the
following property of substitution shifts.
\begin{corollary}\label{corollaryMinimalPrimitive}
 Every infinite minimal substitution shift 
is isomorphic to a primitive eventually proper substitution shift.
\end{corollary}
\begin{proof}
 A minimal substitution shift is, by Theorem~\ref{ch5:subsec:Bratteli-substitution},
isomorphic to $(X_E,T_E)$ with $(V,E,\le)$ a 
stationary properly ordered Bratteli diagram. 
By Proposition~\ref{ch5:proposition:substitutionread}, the system $(X_E,T_E)$ is isomorphic
to the substitution shift $(X_\sigma,S)$ where $\sigma$ is the morphism read on $(V,E,\le)$.
But by Proposition~\ref{propositionPrimitiveProper}, the morphism $\sigma$ is primitive
and eventually proper,
which proves the statement.
\end{proof}
Let us illustrate this result on the  case of the binary
Chacon substitution $\sigma:0\to 0010,1\to 1$.
\index{subject}{Chacon!binary!substitution}%
\index{subject}{substitution!binary Chacon}%
\index{names}{Chacon, Rafael V.}%
The substitution is not primitive but the corresponding shift space
is minimal (see Exercise~\ref{exerciseChaconMinimal}).

Let $x=\sigma^\omega (0.0)$. 
This is the {\em Chacon binary sequence}.
\index{subject}{Chacon!binary!sequence}%
\index{subject}{sequence!Chacon binary}%
\index{names}{Chacon, Rafael V.}%
Using the return words to $0.0$ we see that  $x=\phi (y)$ where $\phi : \{ a,b,c \}^* \rightarrow \{ 0,1\}^*$ is defined by $\phi (a) = 0$, $\phi (b) = 010$ and $\phi (c) = 01010$, and $y=\tau^\omega (b.a)$ where $\tau$ is defined by $\tau (a) = ab$, $\tau (b) = acb$ and $\tau (c) = accb$. 
According to the  proof of Proposition~\ref{ch5:prop:rauzy} we need to take $\tau^2$ instead of $\tau$ and we take
\begin{enumerate}
\item
$C = \{ a_1, b_1,b_2,b_3, c_1, c_2,c_3,c_4,c_5 \}$,
\item
$\theta : C \to \{ 0,1 \}$ given by the following table
$$
\hskip -1cm
\begin{array}{|l|l|l|l|l|l|l|l|l|l|}
\hline
\alpha         & a_1 & b_1& b_2 & b_3 & c_1& c_2 & c_3 & c_4 & c_5\\
\hline
\theta (\alpha ) & 0     & 0    & 1     & 0    & 0    & 1     & 0     & 1     & 0  \\
\hline
\end{array}
$$
\item
$\gamma : \{ a,b,c \} \to C^+$ defined by 
$$
\gamma (a) = a_1 , \ \gamma (b) = b_1 b_2 b_3 , \ \gamma (c) = c_1 \cdots  c_5  \  , 
$$
\item
the substitution $\zeta : C^* \to C^*$  defined by the following tables
\begin{align*}
& \begin{array}{|l|l|l|l|l|l|l|l|l|l|}
\hline
\alpha          & a_1             & b_1      & b_2      & b_3            \\
\hline
\zeta (\alpha ) & \gamma (abacb) & \gamma (a) & \gamma (b) & \gamma (accbacb) \\ 
\hline
\end{array} \\
& \begin{array}{|l|l|l|l|l|l|l|l|l|l|}
\hline
\alpha          & c_1      & c_2      & c_3      & c_4      & c_5\\
\hline
\zeta (\alpha ) & \gamma (a) & \gamma (b) & \gamma (a) & \gamma (c) & \gamma (cbaccbacb)   \\
\hline
\end{array}
\end{align*}
\end{enumerate}

\begin{figure}[hbt]
\centering
\tikzset{node/.style={circle, draw,minimum size=0.1cm,inner sep=0pt}}
	\tikzset{title/.style={minimum size=0.5cm,inner sep=0pt}}
\begin{tikzpicture}
\node[node](0)at(5,6){};

\node[node](11)at(0,5){};\node[node](12)at(1.5,5){};\node[node](13)at(3,5){};\node[node](14)at(4.5,5){};
\node[node](15)at(6,5){};\node[node](16)at(7.5,5){};\node[node](17)at(9,5){};
\node[node](18)at(10.5,5){};\node[node](19)at(12,5){};

\node[node](21)at(3,4){};\node[node](22)at(5,4){};\node[node](23)at(7,4){};

\node[node](31)at(0,2){};\node[node](32)at(1.5,2){};\node[node](33)at(3,2){};\node[node](34)at(4.5,2){};
\node[node](35)at(6,2){};\node[node](36)at(7.5,2){};\node[node](37)at(9,2){};
\node[node](38)at(10.5,2){};\node[node](39)at(12,2){};

\draw(0)edge node{}(11);\draw(0)edge node{}(12);
\draw(0)edge node{}(13);\draw(0)edge node{}(14);
\draw(0)edge node{}(15);\draw(0)edge node{}(16);
\draw(0)edge node{}(17);\draw(0)edge node{}(18);
\draw(0)edge node{}(19);

\draw(11)edge node{}(21);
\draw(12)edge node{}(22);\draw(13)edge node{}(22);\draw(14)edge node{}(22);
\draw(15)edge node{}(23);\draw(16)edge node{}(23);\draw(17)edge node{}(23);
\draw(18)edge node{}(23);\draw(19)edge node{}(23);
%Q1
\draw[bend left=5](21)edge node{}(31);\draw[bend right=5](21)edge node{}(31);
\draw[bend left=5](22)edge node{}(31);\draw[bend right=5](22)edge node{}(31);
\draw(23)edge node{}(31);
%Q2
\draw(21)edge node{}(32);
%Q3
\draw(22)edge node{}(33);
%Q4
\draw[bend left=5](21)edge node{}(34);\draw[bend right=5](21)edge node{}(34);
\draw[bend left=5](22)edge node{}(34);\draw[bend right=5](22)edge node{}(34);
\draw(23)edge node{}(34);
%Q5678
\draw(21)edge node{}(35);\draw(22)edge node{}(36);
\draw(21)edge node{}(37);\draw(23)edge node{}(38);
%Q9
\draw[bend left=5](21)edge node{}(39);\draw[bend right=5](21)edge node{}(39);
\draw[bend left=5](21)edge node{}(39);
\draw[bend right=5](22)edge node{}(39);\draw(22)edge node{}(39);
\draw[bend left=5](22)edge node{}(39);
\draw[bend left=10](23)edge node{}(39);\draw[bend left=5](23)edge node{}(39);
\draw[bend right=10](23)edge node{}(39);\draw[bend right=5](23)edge node{}(39);
\end{tikzpicture}
\caption{The BV-representation of the Chacon binary shift}
\label{figureBVrepresentationChacon}
\end{figure}
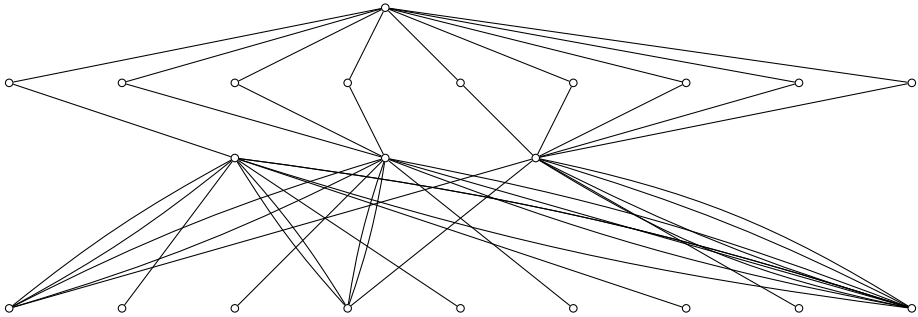

A BV-representation of the Chacon shift (that is, the subshift generated by $x$) is isomorphic to $(X_E , T_E)$ where $(V,E,\le)$ a stationary properly ordered Bratteli diagram with a simple hat such that $\zeta$ is the morphism  read on it. This diagram is obtained by telescoping at odd levels the
diagram of Figure~\ref{figureBVrepresentationChacon}.

As already mentioned (and developed in Exercise~\ref{exerciseForrestRauzy}),
an alternative route to obtain a Bratteli diagram with a simple hat
is to use Proposition~\ref{propositionLemma15Forrest}. In the
present example, we would obtain
\begin{displaymath}
M=\begin{bmatrix}2&2&1\\3&3&3\\4&4&5\end{bmatrix},\quad
P=\begin{bmatrix}1&0&0&0&0&0&0&0&0\\0&1&1&1&0&0&0&0&0\\0&0&0&0&1&1&1&1&1\end{bmatrix},\quad
Q=\begin{bmatrix}2&2&1\\1&0&0\\0&1&0\\2&2&3\\1&0&0\\0&1&0\\1&0&0\\0&0&1\\
2&3&4\end{bmatrix},
\end{displaymath}and thus
\begin{displaymath}
M'=\begin{bmatrix}
2&2&2&2&1&1&1&1&1\\
1&0&0&0&0&0&0&0&0\\
0&1&1&1&0&0&0&0&0\\
2&2&2&2&3&3&3&3&3\\
1&0&0&0&0&0&0&0&0\\
0&1&1&1&0&0&0&0&0\\
1&0&0&0&0&0&0&0&0\\
0&0&0&0&1&1&1&1&1\\
2&3&3&3&4&4&4&4&4
\end{bmatrix}
\end{displaymath}

%%%%%%%%%%%%%%%%%%%%%%%%%%%%%%%%%%%%%%%%%%%%%%%%
%\subsubsection{Generality in the Chacon example}
%%%%%%%%%%%%%%%%%%%%%%%%%%%%%%%%%%%%%%%%%%%%%%%%

%Let $\sigma : A \to A^*$ be a substitution and consider the sets 
%\begin{align*}
%L (\sigma ) & = \{ w\in A^* \mid w \hbox{ is a factor of some } \sigma^n (a) , a\in A , n\in \N \} \hbox { and }\\
%\Omega (\sigma ) & = \{ x\in A^\Z \mid x_i x_{i+1} \cdots x_j \in L (\sigma ) ,i,j \in \Z \} \ .
%\end{align*}
%For $\Omega (\sigma )$ to be non-empty, it is necessary and sufficient that there exists $e\in A$ with 
%\begin{align}
%\label{ch5:condition1}
%\lim_{n\to +\infty} \vert \sigma^n (e)\vert = +\infty \ .
%\end{align}
%Without lost of generality we suppose that all letters of $A$ have an occurrence in some element of $\Omega (\sigma )$ and that 
%\begin{align}
%\label{ch5:condition2}
%\forall a\in A, \  \exists n \in \N, \  \vert \sigma^n (a) \vert_a \geq 1 \ .
%\end{align}T

\subsection{Dimension groups and BV representation of substitution shifts}
We now derive from the previous results a description
of the dimension group of substitution shifts.

Let $\sigma:A^*\to A^*$ be a substitution generating
an aperiodic minimal  shift $X$. As in the proof of Theorem
\ref{ch5:subsec:Bratteli-substitution}, changing if necessary $\sigma$ for some power of $\sigma$, let $x\in A^\Z$ be an admissible
fixed point of  $\sigma$ and set $r=x_{-1}$ and $\ell=x_0$.
Let $\RR_X(r\cdot \ell)=\ell\RR_X(r\ell)\ell^{-1}$ and let
$\phi:B\to\RR_X(r\cdot\ell)$ be a coding morphism for
$\RR_X(r\cdot\ell)$. Let $\tau:B^*\to B^*$ be the morphism
such that $\phi\circ\tau=\sigma\circ\phi$.
\begin{theorem}\label{theoremDGSubstitutionShifts}
Let $\sigma$ be a substitution generating
an aperiodic minimal  shift $X$ and let $\tau:B^*\to B^*$ be as
above. Let $M$ be the composition matrix
of $\tau$. Then $K^0(X,S)=(\Delta_M,\Delta_M^+,v)$
where $v$ is the image in $\Delta_M$ of the
vector with components $|\phi(b)|$ for $b\in B$.
\end{theorem}
\begin{proof}
By Proposition~\ref{propositiontauProperAperiodic},
the substitution $\tau$ is primitive, eventually proper and aperiodic.
Thus, by Proposition~\ref{ch5:proposition:substitutionread},
the system $(X_\tau,S)$ is isomorphic to $(X_E,T_E)$
where $B=(V,E,\le)$ is the Bratteli diagram with $\tau$
read on $B$. By Theorem \ref{theoremDGBratteliDiagram},
the dimension group of $(X_\tau,S)$
is isomorphic with $(\Delta_M,\Delta_M^+,1_M)$. Since $(X_\tau,S)$
is the system induced by $X$ on the clopen set $[r\cdot\ell]$,
the result follows from Proposition~\ref{propositionIuRu}.
\end{proof}
The result can of course also be deduced from 
Proposition~\ref{propositionBVPrimitive}.

We will use the following simple result (note that
it is a particular case of Theorem~\ref{theoremKakutaniEquiv} on Kakutani equivalence).

Let $\tau:B^*\to B^*$ be a primitive eventually proper 
nonperiodic substitution and
let $X(\tau)$ be the associated substitution shift.
Let $(V,E,\le)$ be a properly ordered  simple Bratteli diagram
with simple hat
such that $\tau$ is the morphism read on $(V,E,\le)$.
By Proposition~\ref{ch5:proposition:substitutionread},
the shift $X(\tau)$ is conjugate to $(X_E,T_E)$.
We identify the set $V\setminus\{0\}$ to $B\times\N$.
Let $\phi:B^*\to A^*$ be  a  morphism recognizable on $X(\tau)$
\index{subject}{recognizable!morphism}
and let $Y=X(\tau,\phi)$ be the corresponding substitutive shift.

\begin{proposition}\label{propositionBratteliPrimitive}
The shift $Y$ is conjugate to $(X_{E'},T_{E'})$ where $(V',E')$
is the Bratteli diagram obtained from $(X,E)$ by replacing
each edge from $0$ to $(b,1)$ by $|\phi(b)|$ edges
$(0,b,i)$ with $0\le i < |\phi(b)|$.
\end{proposition}
\begin{proof}
Let $U$ be the clopen subset of $X_{E'}$ formed of the paths
with a first edge of the form $(0,b,0)$ for some $b\in B$.
Ths system induced by $(X_{E'},T_{E'})$ on $U$ is clearly
$X_E$. Thus $X_{E'}$ is the primitive of $X_E$
relative to the function $f(x)=|\phi(x_0)|$.
On the other hand, since $\phi$ is recognizable
on $X$, by Proposition\ref{propositionRecognizableInduced},
$Y$ is the primitive
of $X$ relative to the function $f(x)=|\phi(x)|$.
Thus $X_{E'}$ and $Y$ are conjugate.
\end{proof}
We will see several examples of application of this result below.

\subsection{Some examples}

We will illustrate the preceding results on some classical
examples. 
\begin{example}\label{exampleBVrepresentationFibonacci}
Let $X$ be the Fibonacci shift
generated by $\sigma:a\to ab,b\to a$.
\index{subject}{Fibonacci!shift space!BV-representation}%
\index{subject}{BV-representation!of Fibonacci shift}%
 Consider the fixed point $\sigma^{2\omega}(a\cdot a)$. Let $\phi:\{u,v\}\to\{a,b\}^*$ 
be a coding morphism
for $\RR_X(a\cdot a)=\{aba,ababa\}$. The morphism $\tau:B^*\to B^*$
such that $\phi\circ \tau=\sigma^2\circ\phi$ is
$\tau:u\to uv,v\to uvv$ (see Example~\ref{exampleTauFibonacci}). The Bratteli diagram 
with simple hat corresponding to $\tau$ is shown in Figure~\ref{figureBVrepresentationFibonacci} on the left. The Bratteli
diagram of $X$ (according to Proposition~\ref{propositionBratteliPrimitive}) is 
shown on the right.
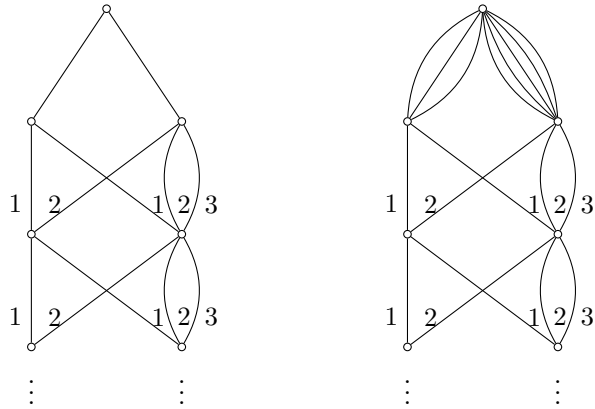
\begin{figure}[hbt]
\centering
\tikzset{node/.style={circle,draw,minimum size=0.1cm,inner sep=0pt}}
\tikzset{title/.style={minimum size=0.4cm,inner sep=0pt}}
\begin{tikzpicture}
%X(tau)
\node[node](0)at(1,4.5){};
\node[node](11)at(0,3){};\node[node](12)at(2,3){};
\node[node](21)at(0,1.5){};\node[node](22)at(2,1.5){};
\node[node](31)at(0,0){};\node[node](32)at(2,0){};
\node[title](41)at(0,-.5){$\vdots$};\node[title](42)at(2,-.5){$\vdots$};

\draw(0)edge node{}(11);\draw(0)edge node{}(12);
\draw[left,near end](11)edge node{$1$}(21);\draw[right,near end](11)edge node{$1$}(22);
\draw[left,near end](12)edge node{$2$}(21);
\draw[bend left,right,near end](12)edge node{$3$}(22);
\draw[bend right,right,near end](12)edge node{$2$}(22);
\draw[left,near end](21)edge node{$1$}(31);\draw[right,near end](21)edge node{$1$}(32);
\draw[left,near end](22)edge node{$2$}(31);
\draw[bend left,right,near end](22)edge node{$3$}(32);
\draw[bend right,right,near end](22)edge node{$2$}(32);

%X(sigma)
\node[node](0)at(6,4.5){};
\node[node](11)at(5,3){};\node[node](12)at(7,3){};
\node[node](21)at(5,1.5){};\node[node](22)at(7,1.5){};
\node[node](31)at(5,0){};\node[node](32)at(7,0){};
\node[title](41)at(5,-.5){$\vdots$};\node[title](42)at(7,-.5){$\vdots$};

\draw[bend left](0)edge node{}(11);\draw(0)edge node{}(11);
\draw[bend right](0)edge node{}(11);
\draw[bend left=30](0)edge node{}(12);\draw[bend left=15](0)edge node{}(12);
\draw(0)edge node{}(12);
\draw[bend right=15](0)edge node{}(12);\draw[bend right=30](0)edge node{}(12);
\draw[left,near end](11)edge node{$1$}(21);\draw[right,near end](11)edge node{$1$}(22);
\draw[left,near end](12)edge node{$2$}(21);
\draw[bend left,right,near end](12)edge node{$3$}(22);
\draw[bend right,right,near end](12)edge node{$2$}(22);
\draw[left,near end](21)edge node{$1$}(31);\draw[right,near end](21)edge node{$1$}(32);
\draw[left,near end](22)edge node{$2$}(31);
\draw[bend left,right,near end](22)edge node{$3$}(32);
\draw[bend right,right,near end](22)edge node{$2$}(32);
\end{tikzpicture}
\caption{The BV-representation of the Fibonacci shift}\label{figureBVrepresentationFibonacci}
\end{figure}
According to Theorem~\ref{theoremDGBratteliDiagram}, the
dimension group is the group $(\Delta_M,\Delta_M^+,v)$
with 
\begin{displaymath}
M=\begin{bmatrix}2&1\\1&1\end{bmatrix},\quad v=\begin{bmatrix}5\\3\end{bmatrix},
\end{displaymath}
(exchanging the letters $u,v$ to write $M$ and $v$). This is consistent
with the fact that $K^0(X,S)=\Z[\lambda]$ with $\lambda=(1+\sqrt{5})/2$
(see Examples \ref{exampleK0Fibonacci} and \ref{exampleFibonacci6}) because
\begin{displaymath}
M=\begin{bmatrix}1&1\\1&0\end{bmatrix}^2 \mbox{ and } v=M^2\begin{bmatrix}1\\0\end{bmatrix}.
\end{displaymath}
\end{example}
\begin{example}\label{exampleBVrepresentationThueMorse}
Let now $X$ be the Thue-Morse shift
generated by$\sigma:0\to 01,1\to 10$.
\index{subject}{Thue-Morse!shift space!BV-representation}%
\index{subject}{BV-representation!of Thue-Morse shift}%
Let $x=\sigma^{2\omega}(0\cdot 0)$.
Set $B=\{1,2,3,4\}=$ and let $\phi:B^*\to \{0,1\}^*$ be the
coding morphism for $\RR_X(0\cdot 0)=\{01101,0110,01011010,010110\}$.
Let $\tau:B^*\to B^*$ be the substitution such that
$\phi\circ\tau=\sigma^2\circ\phi$. We have
$\tau:1\to 1234,2\to 124,3\to 13234,4\to 1324$.
Using again Proposition~\ref{propositionBratteliPrimitive}, we obtain
the BV-representation of Figure~\ref{figureBVrepresentationMorse}.

\begin{figure}[hbt]
\centering
\tikzset{node/.style={circle,draw,minimum size=0.1cm,inner sep=0pt}}
\tikzset{title/.style={minimum size=0.4cm,inner sep=0pt}}
\begin{tikzpicture}
\node[node](0)at(4.5,6){};
\node[node](11)at(0,4){};\node[node](12)at(3,4){};
\node[node](13)at(6,4){};\node[node](14)at(9,4){};
\node[node](21)at(0,2){};\node[node](22)at(3,2){};
\node[node](23)at(6,2){};\node[node](24)at(9,2){};
\node[node](31)at(0,0){};\node[node](32)at(3,0){};
\node[node](33)at(6,0){};\node[node](34)at(9,0){};
\node[title](41)at(0,-.5){$\vdots$};\node[title](42)at(3,-.5){$\vdots$};
\node[title](43)at(6,-.5){$\vdots$};\node[title](44)at(9,-.5){$\vdots$};

\draw[bend left=15](0)edge node{}(11);\draw[bend left=10](0)edge node{}(11);
\draw[bend left=5](0)edge node{}(11);\draw[bend right=5](0)edge node{}(11);
\draw[bend right=10](0)edge node{}(11);\draw[bend right=15](0)edge node{}(11);

\draw[bend left=10](0)edge node{}(12);\draw[bend left=5](0)edge node{}(12);
\draw[bend right=5](0)edge node{}(12);\draw[bend right=10](0)edge node{}(12);

\draw[bend left=20](0)edge node{}(13);\draw[bend left=15](0)edge node{}(13);
\draw[bend left=10](0)edge node{}(13);\draw[bend left=5](0)edge node{}(13);
\draw[bend right=5](0)edge node{}(13);\draw[bend right=10](0)edge node{}(13);
\draw[bend right=15](0)edge node{}(13);\draw[bend right=20](0)edge node{}(13);

\draw[bend left=15](0)edge node{}(14);\draw[bend left=10](0)edge node{}(14);
\draw[bend left=5](0)edge node{}(14);\draw[bend right=5](0)edge node{}(14);
\draw[bend right=10](0)edge node{}(14);\draw[bend right=15](0)edge node{}(14);

\draw[left,very near end](11)edge node{$1$}(21);
\draw[left,very near end](12)edge node{$2$}(21);
\draw[left,very near end](13)edge node{$3$}(21);
\draw[left,very near end](14)edge node{$4$}(21);

\draw[left,very near end](11)edge node{$1$}(22);
\draw[left,very near end](12)edge node{$2$}(22);
\draw[left,very near end](14)edge node{$3$}(22);

\draw[left,very near end](11)edge node{$1$}(23);
\draw[left,very near end](12)edge node{$3$}(23);
\draw[left,very near end,bend right](13)edge node{$2$}(23);
\draw[right,very near end,bend left](13)edge node{$4$}(23);
\draw[right,very near end](14)edge node{$5$}(23);

\draw[right,very near end](11)edge node{$1$}(24);
\draw[right,very near end](12)edge node{$3$}(24);
\draw[right,very near end](13)edge node{$2$}(24);
\draw[right,very near end](14)edge node{$4$}(24);

\draw[left,very near end](21)edge node{$1$}(31);
\draw[left,very near end](22)edge node{$2$}(31);
\draw[left,very near end](23)edge node{$3$}(31);
\draw[left,very near end](24)edge node{$4$}(31);

\draw[left,very near end](21)edge node{$1$}(32);
\draw[left,very near end](22)edge node{$2$}(32);
\draw[left,very near end](24)edge node{$3$}(32);

\draw[left,very near end](21)edge node{$1$}(33);
\draw[left,very near end](22)edge node{$3$}(33);
\draw[left,very near end,bend right](23)edge node{$2$}(33);
\draw[right,very near end,bend left](23)edge node{$4$}(33);
\draw[right,very near end](24)edge node{$5$}(33);

\draw[right,very near end](21)edge node{$1$}(34);
\draw[right,very near end](22)edge node{$3$}(34);
\draw[right,very near end](23)edge node{$2$}(34);
\draw[right,very near end](24)edge node{$4$}(34);
\end{tikzpicture}
\caption{The BV-representation of the Thue-Morse shift}\label{figureBVrepresentationMorse}
\end{figure}
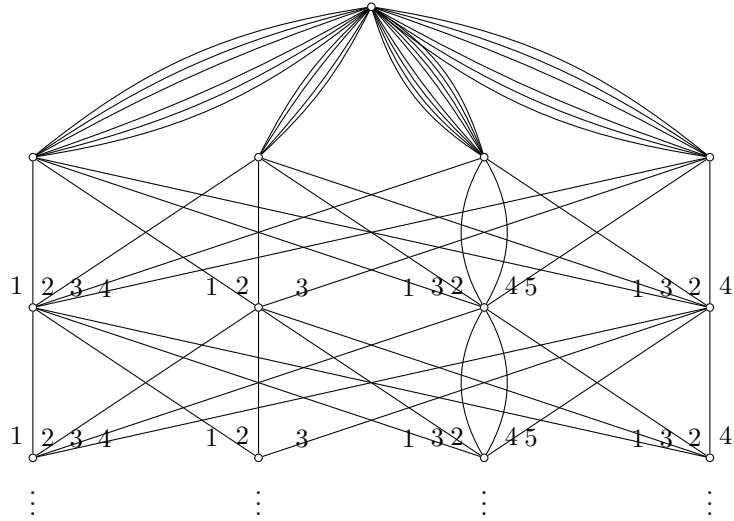
The dimension group $K^0(X,S)$ is thus $(\Delta_M,\Delta_M^+,v)$
with
\begin{displaymath}
M=\begin{bmatrix}
1&1&1&1\\1&1&0&1\\1&1&2&1\\1&1&1&1\end{bmatrix},\mbox{ and }
v=\begin{bmatrix}4\\3\\5\\4\end{bmatrix}
\end{displaymath}
 This is consistent with the result found in Example~\ref{exampleMorse5}
where we identified $K^0(X,S)$ to $\Z[1/2]\times\Z$ with the natural order
of $\R^2$ and unit $(1,0)$.
Indeed, the eventual range of $M$ is generated by the vectors
\begin{displaymath}
x=\begin{bmatrix}1\\1\\1\\1\end{bmatrix},\quad
y=\begin{bmatrix}1\\0\\2\\1\end{bmatrix}.
\end{displaymath}
In this basis, the matrix $M$ takes the form
\begin{displaymath}
N=\begin{bmatrix}3&1\\2&2\end{bmatrix}
\end{displaymath}
with eigenvectors $z=x+y$ for the eigenvalue $4$ and $t=x-2y$ for the eigenvalue $1$. This allows us to identify $K^0(X,S)$ via the map
$\alpha x+\beta y\mapsto(\alpha,\beta)$ with the pairs
$(\alpha,\beta)$ such that $2^n\alpha,\beta\in \Z$ for some $n\ge 1$. 
The order unit is
$v=3x+y\mapsto(3,1)$. Since $v=Nx$, we can normalize the unit to $(1,0)$.

In Example~\ref{exampleMorse5}
we obtained the dimension group as the group of a matrix of dimension $3$
instead of $4$ as here. Actually the method used in Example~\ref{exampleMorse5}
can be used also to build a BV-representation of the Thue-Morse 
shift with $3$ vertices at each level $n\ge 1$ 
(Exercise \ref{exerciseThueMorseDim3}).
\end{example}
%%%%%%%%%%%%%%%%%%%%%%%%%%%%%%%%%%%%%%%%%
\section{Linearly recurrent shifts}
%%%%%%%%%%%%%%%%%%%%%%%%%%%%%%%%%%%%%%%%%
\label{ch5:subsec:LR}

In this section, we study the case of linearly recurrent shifts,
already introduced in Chapter~\ref{chapterTopologicalDynamicalSystems}.
This family contains primitive substitution shifts 
(Proposition~\ref{propositionPrimitiveSubstitutiveIsLR}). The main result
is a characterization of the BV-representation of aperiodic linearly
recurrent shifts (Theorem~\ref{ch5:theorem:LRBVrepresentation}).

\subsection{Properties of linearly recurrent shifts}
%\begin{definition}
\label{ch5:def:LR}

In what follows we show that a subshift is linearly recurrent if and only if
 it has a BV-representation where the incidence matrices have positive
 entries and belong to a finite set.

The class of linearly recurrent shift spaces is clearly closed under conjugacy
(Exercise~\ref{exerciseLRClosed}).

A simple example of a minimal shift which is not LR is given below.
\begin{example}\label{exampleNotLR}
Let $u_n$ be the sequence of words defined by $u_{-1}=b$, $u_0=a$
and $u_{n+1}=u_n^nu_{n-1}$ for $n\ge 1$. The subshift generated
by the infinite word $u$ having all $u_n$ as prefixes is not LR.
This  results from the property of LR shift spaces with constant $K$
to be $(K+1)$-th power free (Proposition~\ref{ch5:proposition:firstpropLR} below).
\end{example}

We will use Theorem~\ref{ch5:subsec:Bratteli-substitution} to prove the following result, which
gives an easily verifiable condition for a substitution shift to be LR
(the condition is, for example, visibly satisfied by
the Chacon binary substitution).
%(the definition of linearly recurrent subshift is given in Section~\ref{ch5:subsec:LR}).
It proves that the property of being LR is decidable
\index{subject}{decidability!of LR for substitution shift} for a substitution shift.
\begin{theorem}[Damanik, Lenz]\label{ch5:th:DamanikLenz2006}
\index{names}{Damanik, David}\index{names}{Lenz, Daniel H.}%
Let $\sigma:A\to A^*$ be a substitution
and let $(X , S)$ be the corresponding substitution subshift. 
Let $e\in A$ be a letter such that $\lim_{n\to\infty}|\sigma^n(e)|=\infty$.
Then the following are equivalent.
\begin{enumerate}
\item[\rm(i)]
The letter $e$
 occurs with bounded gaps in $\cL(X)$
and for every $a\in A$, one has $|\sigma^n(e)|_a\ge 1$ for some $n\ge 1$.
\item[\rm(ii)]
\label{ch5:condmin}
$(X , S)$ is minimal.
\item[\rm(iii)]
\label{ch5:condLR}
$(X  , S)$ is linearly recurrent.
\end{enumerate}
\end{theorem}
\begin{proof}
(i)$\Rightarrow$ (ii) We have to show that for every $n\ge 1$,
the word $\sigma^n(e)$ occurs in every long enough word of $\cL(X)$.
Since $e$ occurs with bounded gaps,
there is $k$  such that every word of $\cL_k(X)$ contains $e$.
Choose $w$ of length $|w|\ge (3+k)\max_{a\in A}|\sigma^n(a)|$.
Then $w$ contains a factor $\sigma^n(u)$ for some word $u\in\cL_k(X)$ 
and thus contains $\sigma^n(e)$.

(ii)$\Rightarrow$ (iii) By Corollary~\ref{corollaryMinimalPrimitive}
$X$ is isomorphic to a primitive substitutive shift. Thus the
assertion results from Proposition~\ref{propositionPrimitiveSubstitutiveIsLR}.

(iii)$\Rightarrow$ (i) is obvious.
\end{proof}

A fourth equivalent statement could be added,
namely that $(X  , S)$ is uniquely ergodic (see the notes for a reference).

From Theorem~\ref{ch5:th:DamanikLenz2006}, we know that all minimal substitution shifts are linearly recurrent.
We will prove now some important properties of LR shifts.

\begin{proposition}
\label{ch5:proposition:firstpropLR}
Let $X$ be a non-periodic shift space and suppose that  it is linearly recurrent with constant $K$. 
Then:  
\begin{enumerate}
\item Every word of $\cL_n(X)$ appears in every word of $\cL_{(K+1)n-1}(X)$.
\item
The number of distinct factors of length $n$ of
$X$ is at most $Kn$. 
\item
$\cL(X)$ is $(K+1)$-power free, that is, for every nonempty word $u$,
 $u^n\in\cL(X)$ implies $n\le K$.
\item 
For all $ u\in\cL(X)$ and for all $w \in {\mathcal R}_X(u)$ we have $\frac{1}{K}|u| < |w| \le K|u|$.
\item 
If $u\in\cL(X)$, then  $\Card({\mathcal R}_X(u) ) \leq K(K+1)^2$.
\end{enumerate}
\end{proposition}
\begin{proof}
The first three assertions are already in Proposition~\ref{propositionLRhasLinearComplexity}.

4. Suppose that for $u\in\cL(X)$ and $w\in\RR_X(u)$, we have $|u|\ge K|w|$.
Since $u$ is a suffix of $uw$, the length of $w$ is a period of $u$
and thus $w^{K+1}$ is a suffix of $uw$, a contradiction with assertion 3.

5. Let $u,v\in\cL(X)$ with $|v|=(K+1)^2|u|$. By assertion 1, every element 
of $\cL(X)$ of length $(K+1)|u|$ appears in $v$ and thus all words of $u\RR_X(u)$
are factors of $v$. By assertion 4, we have
\begin{equation}
\Card(\RR_X(u))|u|/K\le \sum_{w\in \RR_X(u)}|w|.\label{eqChristophe1}
\end{equation}
Now,  every word in $u\RR_X(u)$
occurs in  $v$.  Considering all
successive occurrences of $u$ in $v$,
we have   a factor  of $v$
of the form $ux$ with $x=x_1x_2\cdots x_k$ where all $x_i$ are in $\RR_X(u)$
and every occurrence of $u$ in $v$ is either as a prefix of $ux$
or as a suffix of some $ux_i$.
Since the $x_i$ are all the factors of $v$ which belong to $\RR_X(u)$,
 every word in $\RR_X(u)$ appears as some $x_i$. This implies that
\begin{equation}
\sum_{w\in \RR_X(u)}|w|\le |v|=(K+1)^2|u|,\label{eqChristophe2}
\end{equation}
%\marginpar{Il y a un probleme la}
whence the conclusion $\Card(\RR_X(u))\le K(K+1)^2$.
\end{proof}

Note the following corollary of Proposition~\ref{ch5:proposition:firstpropLR}.

\begin{corollary}\label{corollaryComplexityPrimitive}
The factor complexity of a minimal substitution shift is at most linear.
\end{corollary}
Indeed, by Theorem~\ref{ch5:th:DamanikLenz2006},
a minimal substitution shift is linearly recurrent.
This statement is not true for a substitution shift which is not primitive
(Exercise~\ref{exerciseExampleComplexityn^2}).

We now come back to return words.
\index{subject}{return!word}%
Let $X$ be a minimal shift space. For $u,v$ such that $uv\in\cL(X)$,
we consider as in the proof of Theorem \ref{ch5:subsec:Bratteli-substitution},
the set $\RR_X(uv)$ of right return words to $uv$ and the set
\begin{displaymath}
\RR_X(u\cdot v)=v\RR_X(uv)v^{-1},
\end{displaymath}
\index{symbols}{R@$\RR_X(u\cdot v)$}%
called the set of \emph{return words} to $u\cdot v$.
As we will manipulate return words,  it is important to observe that a finite word $w\in A^+$ is a return word to $u.v$ in $X$ if, and
only if, 
\begin{enumerate}
\item
$uwv$ is in $\cL(X)$, and
\item
$v$ is a prefix of $wv$ and $u$ is a suffix of $uw$, and
\item
the finite word $uwv$ contains exactly two occurrences of the finite word $uv$.
\end{enumerate}

If $v$ is the empty word, we have $\RR_X(u\cdot v)=\RR_X(u)$ and if $u$
is the empty word, then $\RR(u\cdot v)=\RR'_X(v)$, the set of left return words to $v$.
 When it will be clear from the context   we will use to ${\mathcal R}(u\cdot v)$ in place of ${\mathcal R}_X(u\cdot v)$.

%We enumerate the elements $w$ of ${\mathcal R}_{x,u.v}$ in the order of the
%first appearance of $uwv$ in $x_{[-|u|,+\infty )}$. 
Let $B_{u\cdot v}$ be a finite alphabet in bijection with $\RR(u\cdot v)$ by $
\theta_{u\cdot v}:B_{u\cdot v}\to\RR(u\cdot v)$.
The map $\theta_{u\cdot v}$ extends to
a morphism from $B_{u\cdot v}^*$ to $A^*$ and the set $\theta_{u\cdot v}(B_{u\cdot v} ^*)$ 
consists of all concatenations of return words to $u\cdot v$.

 The set $\RR(u\cdot v)$ is a circular code and, more precisely, no word
of $\RR(u\cdot v)$ overlaps no trivially a product of words of $\RR(u\cdot v)$.
In particular,
the map $\theta_{u\cdot v}: B_{u\cdot v}^*\to A^*$ is one-to-one.

\subsection{BV-representation of linearly recurrent shifts}
The following result characterizes linearly recurrent shift spaces in terms of BV-representations. 
Recall that a dynamical system $(X,T)$, endowed with the distance $d$, is {\em expansive}
\index{subject}{expansive dynamical system}%
\index{subject}{dynamical system!expansive}%
 if there exists $\epsilon $ such that for all pairs of points $(x, y)$, $x\not = y$, there exists $n$ with $d(T^n x , T^n y)\geq \epsilon $.
We say that  $\epsilon $ is a {\em constant of expansivity}\index{subject}{constant!of expansivity} of $(X,T)$.
It is easy to see that the shift spaces are expansive but that odometers are not.

\begin{theorem}
\label{ch5:theorem:LRBVrepresentation}
An aperiodic shift space is linearly recurrent if, and only if, it has a
BV-repre\-sen\-tation satisfying:
\begin{enumerate}
\item
\label{ch5:proposition:cond1LR}
its incidence matrices have positive entries and belong to a finite set of matrices,
\item
\label{ch5:proposition:cond2LR}
for all $n\geq 1$ the morphism read on $E(n)$ is eventually proper.
\end{enumerate}
\end{theorem}
\index{subject}{BV-representation!of linearly recurrent shifts}
\begin{proof}
Let $X$ be an aperiodic LR shift space. 
It suffices to construct a sequence of KR-partitions having the desired properties.

From Proposition~\ref{ch5:proposition:firstpropLR} there exists an integer $K$ such that 
for all $u$ occurring in some $x\in X$ and all $w\in {\mathcal R}(u)$,  we have
$$
\frac{|u|}{K} 
\leq 
|w| 
\leq 
K |u|.
$$
Moreover, every word of $\cL_n(X)$ occurs in every word of
$\cL_{(K+1)n-1}(X)$.

We set $\alpha = (K+1)^2$. 
Let $x= (x_n )_n$ be an element of $X$. 
For each non-negative
integer $n$, we set $u_n = x_{-\alpha^n} \cdots x_{-2} x_{-1}$,
$v_n = x_0x_1\cdots x_{\alpha^n-1}$
and ${\mathcal R}_n = {\mathcal R}(u_n.v_n)$.
%$B_n =B_{u_n\cdot v_n}$ and $\theta_n = \theta_{u_n\cdot v_n}$. 

By the choice of $K$, every word of $\RR_X(u_n\cdot v_n)$ appears
in $u_{n+1}v_{n+1}$. Indeed, $\alpha= (K+1)^2$
implies $2\alpha^{n+1}= (K+1)^22\alpha^n$ and thus
$|u_{n+1}v_{n+1}|= (K+1)^2|u_nv_n|$.

Now define for all $n\ge 0$
$$
{\mathcal P} (n) = 
\left\{
S^j [u_n .w v_n ] \mid w\in {\mathcal R}_n ,  \ 0\leq j< |w| 
\right\}.
$$

The verification that $({\mathcal P} (n))_n $ is a sequence of KR-partitions having the desired properties is left to the reader.

Let now $(V,E,\le)$ be a BV-representation of an aperiodic
shift space satisfying \ref{ch5:proposition:cond1LR} and \ref{ch5:proposition:cond2LR}. Since $(X_E,T_E)$ is conjugate to a shift space, it is expansive.
Let $\epsilon $ be a constant of expansivity of $(X_E,T_E)$.
Let $n_0$ be a level of $(V,E,\le)$ such that all cylinders $[e_1 , \dots , e_{n_0}]$ are included in a ball of radius $\epsilon/2$.
For any vertex $v\in V(n_0)$ let $h_v$ denote the number of paths from $V(0)$
to $v$.
Now consider the alphabet
$$
A = \{ (v,j) \mid v\in V(n_0), \ 0\leq j< h_v \},
$$
the map $C :  X_E  \to  A$ defined by
\begin{displaymath}
    C(e(n))_n = (r(e_{n_0}), j)
\end{displaymath}
if $(e(n))_{1\leq n\leq n_0}$ is the $j$-th finite path in $(V,E,\le)$ from $V(0)$
to $r(e_{n_0})$  with respect to the order on $(V,E,\le)$,
and finally define $\varphi :  X_E  \to  A^{\Z} $ by
$$
 \varphi(x)=\left(C\circ T_E^n (x)\right)_{n\in \Z }.
$$

 By the choice of $\epsilon$, the system
 $(X_E , T_E)$ is isomorphic to $(\Omega , S)$ where $\Omega = \varphi (X_E)$ and $S$ is the shift on $A$.
There remains to show that  $(\Omega , S)$ is LR.

Let $K = \sup_{n\geq n_0}\max_{v\in V(n)} \sum_{v'\in V(n-1)} { M}(n)_{v,v'}$. 
Condition~\ref{ch5:proposition:cond1LR} implies  that $K$ is finite.
Let $L = \max_{v,v'\in V(n_0-1)} (h_v/h_{v'})$.

Let $v\in V(n_0-1)$ and let $w$ be a return word to $u=(v,0) (v, 1) \cdots , (v,h_v-1)$.
Due to Condition~\ref{ch5:proposition:cond1LR} and Condition~\ref{ch5:proposition:cond2LR} we have
\begin{equation}
\label{ch5:inequality:DBLR}
|w| \leq 2K \max_{v'\in V(n_0-1)} h_{v'} \leq 2K L |u|.
\end{equation}

Now for all $n\geq n_0$, let $\tau_n : V(n) \to V(n-1)^*$ be the morphism read on $E(n)$.
We set $W= \{ (v,0) (v, 1) \cdots , (v,h_v-1) \mid v\in V(n_0-1) \}$ and we define the morphism $\sigma : V(n_0-1) \to A^*$ by
$\sigma (v) = (v,0) (v, 1) \cdots , (v,h_v-1)$.

It is clear that all the elements of $\Omega$ are concatenations of finite words belonging to $W$.
They are also concatenations of finite words belonging to $\sigma \circ  \tau_{n_0} (V (n_0 ))$, and more
generally, of finite words belonging to $\sigma \circ \tau_{n_0} \circ  \cdots \circ \tau_{n} (V (n) )$ for all $n\geq n_0$.
As for \eqref{ch5:inequality:DBLR}, we can prove that all return words to some elements of $\sigma \circ \tau_{n_0} \circ \cdots \circ \tau_{n} (V(n))$ satisfy the same inequality. 

Now let $u$ be any non-empty finite word appearing in some word of $\Omega$ and $w$ be a return word to $u$.
There exists $n$ such that 
$$
\max_{v\in V(n)} |\sigma \circ \tau_{n_0} \circ  \cdots \circ  \tau_{n} (v)|\leq |u| < \max_{v\in V(n+1)} |\sigma \circ \tau_{n_0} \cdots \circ \tau_{n+1} (v)|.
$$

Then $u$ is a factor of some $\sigma \circ \tau_{n_0} \circ \cdots \circ \tau_{n+1} (vv')$, $v$ and $v'$ belonging to $V(n+1)$.
From Condition~\ref{ch5:proposition:cond1LR} and Condition~\ref{ch5:proposition:cond2LR}, 
we deduce that $vv'$ is a factor of some $\tau_{n+2} \circ \tau_{n+3} (v'')$, $v'' \in V(n+3)$. 
Then, $u$ is a factor of $\sigma \circ \tau_{n_0}\circ  \cdots \circ \tau_{n+3} (v'')$.
Consequently
\begin{align*}
|w| \leq 2K L |\sigma \circ \tau_{n_0} \circ  \cdots \circ \tau_{n+3} (v'')| \leq 2K L^4 |u|
\end{align*}
and $(\Omega , S)$ is LR.
\end{proof}

Let us call {\em linearly recurrent} any Cantor dynamical systems
having a BV-representation \index{subject}{linearly!recurrent!Cantor dynamical system}\index{subject}{Cantor!dynamical system!linearly recurrent} satisfying
\ref{ch5:proposition:cond1LR} and \ref{ch5:proposition:cond2LR} in
Theorem \ref{ch5:theorem:LRBVrepresentation}. 
By Proposition \ref{propositionBratteliCompactum} linearly recurrent systems are minimal. 
Actually, more can be said about these dynamical systems but we need the following theorem 
which can be seen as an extension of Proposition~\ref{ch5:proposition:substitutionread}.

We say  that a minimal Cantor dynamical system $(X,T)$ has {\em topological rank}
\index{subject}{rank!topological}\index{subject}{topological!rank} $k$ if $k$ is the smallest integer such that 
$(X,T)$ has a BV-representation $(X_E , T_E)$ where
 the sequence of number of vertices $(\Card(V (n)))_n$ bounded by $k$. 
When such a $k$ does not exist,  we say that  it has {\em  infinite topological rank}.
\index{subject}{rank!topological!infinite}%
Of course, linearly recurrent BV-dynamical systems have finite topological rank. 
Using the notion of equicontinuous system introduced in Section~\ref{ch5:subsec:rep-odo}, we have the following result.

\begin{theorem}[Downarowicz, Maass]\label{theoremDownarowiczMaass}
Let $(X,T)$ be a minimal Cantor dynamical system with topological rank $k\in \N$.
Then, $(X , T)$ is expansive if and only if $k\geq 2$.
Otherwise it is equicontinuous.
\end{theorem}
\index{names}{Downarovicz, Tomasz}\index{names}{Maass, Alejandro}%
Let us take the notation of the proof of Theorem~\ref{ch5:theorem:LRBVrepresentation}. We call $\sigma_{n_0}$ the morphism $\sigma$ and $A_{n_0}$ the alphabet $A$.
Let $X_{n_0}$ be the subset of $A_{n_0}^\Z$ consisting of all the sequences
$x$ such that for all $i,j$,  $x_i x_{i+1} \cdots x_j$ is a factor of $\sigma_{n_0}\circ \tau_{n_0} \cdots \circ \tau_{n} (v)$ for some $n\in \N$ and $v\in V (n)$.
It can be checked that $(X_{n_0} , S)$ is a minimal subshift.

\begin{corollary}
Let $(X_E , T_E)$ be a BV-dynamical system with finite topological rank.
Then, $(X_E , T_E)$ is expansive if and only if there exists $n_0$ such that $(X_{n_0} , S)$ is not periodic.

Moreover, if the cylinders $[e_1 , \ldots , e_{n_0 -1}] $ of $X_E$ are all included in balls of radius $\frac{\epsilon}{2}$, $\epsilon$ being a constant of  expansivity, then $(X_E , T_E)$ is isomorphic to $(X_{n_0} , S)$.
\end{corollary}

Note that once some $(X_{n_0} , S)$ is not periodic, then $(X_{n} , S)$ is aperiodic for all  $n\geq n_0$.

\subsection{Unique ergodicity of linearly recurrent systems}

We will prove the following property of linearly recurrent Cantor systems.
It generalizes the property of unique ergodicity for
primitive substitution shifts (Theorem~\ref{theoremMichel}).
\begin{theorem}\label{theoremLRisUE}
Every linearly recurrent Cantor system is uniquely ergodic.
\end{theorem}\index{subject}{unique ergodicity!of linearly recurrent systems}

Let $(X,T)$ be a linearly recurrent system.
By definition, there
is a BV-representation $(V,E,\le)$  satisfying
the conditions of Theorem~\ref{ch5:theorem:LRBVrepresentation},
and thus such that
the incidence matrices have positive entries and belong to a finite set
of matrices.

Let  $(\Pg(n)=\{T^{-j}B_{k}(n); 1\leq k\leq t(n), 0\le j<h_{k}(n) \};n\ge 0 )$ be the refining sequence of partitions associated to
$(X_E,V_E)$ as in Section~\ref{subsectionFromBDtoPT}.
Let $( M(n) = (m_{l,k}(n); 1\leq l \leq t(n) , 1\leq k \leq t(n-1)) ; n\geq 1 )$ be the associated sequence of matrices.

The following lemma expresses a property of linearly recurrent systems,
which could actually be taken for definition (Exercise~\ref{exerciseAltDefLR}).
\begin{lemma}\label{lemmaCondh_k(n)}
  There is an integer $L\ge 1$ such
  that
  \begin{equation}
h_i(n)\le Lh_j(n-1)\label{equationLR}
  \end{equation}
  for every $n\ge 1$, $1\le i\le t(n)$ and $1\le j\le t(n-1)$.
\end{lemma}
\begin{proof}
  Let $h(n)$ be the vector with components $h_i(n)$ for $1\le i\le t(n)$.
  By~\eqref{eqh_t(n)}, we have
  \begin{displaymath}
    h(n)=M(n)h(n-1).
  \end{displaymath}
  Since there is a finite number of distinct matrices $M(n)$, we can define
  $K=\max\|M(n)\|_\infty$. Let $i_0$ be such that $h_{i_0}(n-1)=\max h_i(n-1)$.
  Then
  \begin{displaymath}
    h_{i_0}(n-1)\le h_j(n)\le K h_{i_0}(n-1)
\end{displaymath}
  where the first inequality holds because $M(n)$ is a positive integer matrix.
  This implies that for every $j,k$ and $n\ge 1$, we have
  \begin{displaymath}
    \frac{1}{K}\le\frac{h_j(n)}{h_k(n)}\le K
  \end{displaymath}
  and thus
  \begin{displaymath}
h_j(n)\le K h_{i_0}(n-1)\le K^2h_i(n-1).
    \end{displaymath}
  Consequently \eqref{equationLR}
  will hold with $L=K^2$.
  \end{proof}

We notice that for each $T$--invariant probability measure $\mu$ and for every $n \ge 1$ and $1\leq k\leq t(n-1)$ we have
\begin{equation}
\label{B1}
\mu(B_k(n-1))=\sum_{l=1}^{t(n)}m_{l,k}(n)\,\mu(B_l(n)) \end{equation}
and
\begin{equation}
\label{B2}
\sum_{k=1}^{t(n)}h_k(n)\mu(B_k(n))=1.
\end{equation}

Let us observe that the real numbers $\mu(B_k(n-1))$ should be positive. 
Indeed, from Equation \eqref{B2} there is at least one index $k$ for which $\mu(B_k(n))$ is not zero. 
Hence, the coefficients $m_{l,k}(n)$ being positive, it is also the case for $\mu(B_k(n-1))$.

To prove that linearly recurrent systems are uniquely ergodic we need the following lemma in which $L$ is the integer given by Lemma~\ref{lemmaCondh_k(n)}.

\begin{lemma}
\label{binf}
Let $\mu$ be an invariant measure of $(X,T)$. Then, for all $n\ge 0$ and $1\leq k\leq t(n)$ we have $$
\frac{1}{L}  \leq h_{k}(n)\mu(B_{k}(n))\leq 1 .
$$
\end{lemma}
\begin{proof}
The inequality on the right is obvious from Equation \eqref{B2}.

Fix $k$ with $1\leq k\leq t(n)$. By Equation~\eqref{B1}, since all the entries of $M(n+1)$ are positive, we get $$
\mu(B_k(n))\geq\sum_{\ell=1}^{t(n+1)}\mu(B_l(n+1))\ .
$$
By Lemma~\ref{lemmaCondh_k(n)}, for every $l$ we have 
\begin{equation}\label{B3}
h_k(n)\geq h_l(n+1)/L
\end{equation}
thus
$$
h_k(n)\mu (B_k(n))\geq \sum _{\ell=1}^{t(n+1)}
\frac{h_l(n+1)}{L}\mu(B_l(n+1))=\frac 1L.
$$
\end{proof}

\begin{proofof}{of Theorem~\ref{theoremLRisUE}}
Let $(X,T)$ be a linearly recurrent system.
Given a $T$--invariant probability measure $\mu$, we define the numbers
$$
\mu_{n,k}=\mu(B_{k}(n)),\ n\geq 0,\; 1\leq k\leq t(n)\ .
$$
These nonnegative numbers satisfy the relations
\begin{equation}
\label{eq:rel-mu}
\mu_{0,1}=1\text{ and }\mu_{n-1,k}=\sum_{l=1}^{t(n)}\mu_{n,l}m_{l,k}(n)
\text{ for }
n\ge 1,\ 1\leq k\leq t (n-1).
\end{equation}
In a matrix form: with $V(n)=(\mu_{n,1},\dots,\mu_{n,t(n)})$, we have
$V(n-1)=V(n)M(n)$.
Conversely, let the nonnegative numbers $(\mu_{n,k} ; n\geq 0,\ 1\leq k\leq t(n)$) satisfy these conditions. Since the partitions $\Pg (n)$ are clopen and span the topology of $X$, it is immediate to check that there exists a unique invariant probability measure $\mu$ on $X$ with
$\mu_{n,k}=\mu(B_{k}(n))$ for every $n\ge 0$ and $k\in \{ 1, \dots , t(n)\}$.

By Lemma~\ref{binf} and Equation \eqref{B3},
for the constant  constant $\delta = 1/L^2$ one has
$$
\mu_{n,i}
\geq
\frac{1}{Lh_i (n)} 
\geq
\frac{1}{L^2h_k (n-1)} 
\geq
\frac{1}{L^2} \mu_{n-1,k}
=
\delta\mu_{n-1,k}
$$
for every $n\geq 1$ and $(i,k)\in \{ 1, \dots , t(n)\}\times \{ 1 , \dots , t(n-1) \}$,
and every invariant measure $\mu$. Without loss of generality we can assume $\delta<1/2$.
Let $\mu,\mu'$ be two invariant measures, and $\mu_{n,k},\mu'_{n,k}$ be defined as above. We define
$$
S_n=\max_k\frac{\mu'_{n,k}}{\mu_{n,k}} =\frac{\mu'_{n,i}}{\mu_{n,i}} \ ,\ s_n=\min _k\frac{\mu'_{n,k}}{\mu_{n,k}}=\frac{\mu'_{n,j}}{\mu_{n,j}} \  \hbox{ and } r_n=\frac{S_n}{s_n}
$$
for some $i,j$.
For every $k\in \{ 1, \dots , t(n-1) \}$ we have
\begin{align*}
\mu'_{n-1,k}
&=\sum_{l\neq j}\mu'_{n,l}m_{l,k}(n)+\mu'_{n,j}m_{j,k}(n)\\
&\leq S_n\sum_{l\neq j}\mu_{n,l}m_{l,k}(n)+s_n\mu_{n,j}m_{j,k}(n)\\
&=S_n\mu_{n-1,k}-(S_n-s_n)\mu_{n,j}m_{j,k}(n)
\leq S_n\mu_{n-1,k}-(S_n-s_n)\mu_{n,j}\\
&\leq \mu_{n-1,k}s_n\bigl(r_n(1-\delta)+\delta\bigr)\ .
\end{align*}

And in similar way, for every $k\in \{ 1, \dots , t(n-1) \}$ we have

$$
\mu'_{n-1,k} \geq \mu_{n-1,k}s_n\bigl(\delta r_n+(1-\delta)\bigr).
$$
We deduce that
$$
r_{n-1}\leq \phi(r_{n})\text{ where }
\phi(x)=\frac {(1-\delta)x+\delta}{\delta x+(1-\delta)}.
$$

The function $\phi$ is increasing on $[0,+\infty)$ and tends to $(1-\delta)/\delta$ at $+\infty$. Writing $\phi^m=\phi\circ\dots\circ\phi$ ($m$ times), for every $n,m \ge 0$ we have
$1\leq r_n\leq\phi^m(r_{n+m})\leq \phi^{m-1}((1-\delta)/\delta)$. Taking the limit with $m\to+\infty$, we get $r_n=1$ and thus $\mu=\mu'$.
\end{proofof}

Another proof of Theorem~\ref{theoremLRisUE}
is proposed in Exercise~\ref{exerciseDelecroixBoshernitzan}.

%%%%%%%%%%%%%%%%%%%%%%%%%%%%%%%
\section{$\mathcal{S}$-adic representations}\label{sectionSadicShifts}
%%%%%%%%%%%%%%%%%%%%%%%%%%%%%%%
We introduce now the notion of an $\Sa$-adic representation
of a shift, which generalizes the representation of shifts as
substitution shifts. We will replace the iteration of a morphism by 
the application of an arbitrary sequence
of morphisms. The concept is of course close to that of a BV-representation
but more flexible. For example, every stationary Bratteli diagram
defines  a substitution shift (with the substitution read
on the diagram), but not conversely since, as we have seen,
one has to use a conjugacy to obtain a stationary Bratteli
diagram for a substitution shift (when the substitution is not proper).
The same situation occurs for $\Sa$-adic representations.

In this section, all morphisms are assumed to be nonerasing
and alphabets are always assumed to have cardinality at least 2.
%Using concatenation, we extend $\sigma$ to~$\mathcal{A}^\mathbb{N}$ and~$\mathcal{A}^\mathbb{Z}$. 
%With a morphism $\tau :\A^* \to \B^*$, where $\A$ and $\B$ are finite alphabets, we classically associate an {\em incidence matrix} $M_\tau$ indexed by $\B \times \A$ such that  for every $(b,a)\in \B \times \A$, its entry at position $(b,a)$ is the number of occurrences of $b$ in $\tau(a)$.

\subsection{Directive sequence of morphisms}
For a morphism $\sigma:A^*\to B^*$, we denote as usual 
$|\sigma|=\max_{a\in A}|\sigma(a)|$ and $\langle\sigma\rangle=\min_{a\in A}|\sigma(a)|$.

Let $\Sa$ be a family of morphisms.
Let $(A_n)_{n\ge 0}$ be a sequence of finite alphabets and
let $\tau = (\tau_n)_{n \geq 0}$  be a sequence of morphisms with
$\tau_n :A_{n+1}^*\to A_n^*$ and $\tau_n\in\Sa$.
We thus have an infinite sequence of morphisms
\begin{displaymath}
  \cdots \edge{\tau_2}A_2^*\edge{\tau_1}A_1\edge{\tau_0}A_0^*
  \end{displaymath}
% such that  $|\tau_1 \circ \cdots \circ \tau_{n-1}|$
%goes to infinity when $n$ increases. 
For $0\leq n< N$, we define $\tau_{[n,N)} = \tau_n \circ \tau_{n+1} \circ \dots \circ \tau_{N-1}$
\index{symbols}{tau@$\tau_{[n,N)}$}%
%\begin{displaymath}
%\tau_{[n,N)} = \begin{cases}\tau_n \circ \tau_{n+1} \circ \dots \circ \tau_{N-1}
%&\mbox{ if $n<N$}\\{\rm Id}& \mbox{ otherwise}\end{cases}
%\end{displaymath}
 and $\tau_{[n,N]} = \tau_n \circ \tau_{n+1} \circ \dots \circ \tau_{N}$. 
For $n\geq 0$, the \emph{language $\mathcal{L}^{(n)}({\tau})$ of level $n$}
 associated with $\tau$ is defined~by 
\[
\mathcal{L}^{(n)}(\tau) = \big\{w \in A_n^* \mid \mbox{$w$ occurs in $\tau_{[n,N)}(a)$ for some $a \in A_N$ and $N>n$}\big\}.
\]
\index{symbols}{L@$\cL^{(n)}(\tau)$}%
%Since $|\tau_{[1,n)}|$ goes to infinity when $n$ increases,
Thus  $\cL^{(n)}(\tau)$ is the set of factors of $\tau_n(\cL^{(n+1)}(\tau))$ for $n\ge 0$ and 
\begin{displaymath}
  \cdots \edge{\tau_1}\cL^{(1)}(\tau)\edge{\tau_0}\cL^{(0)}(\tau).
\end{displaymath}
More generally, we have $\cL^{(n)}(\tau)$ is the set of factors of $\tau_{[n,n+m)}(\cL^{(n+m)}(\tau))$.

The language $\mathcal{L}^{(n)}(\tau)$ defines a shift space
$X^{(n)}(\tau)$  called the {\em shift generated by $\mathcal{L}^{(n)}(\tau)$}.
More precisely, $X^{(n)}(\tau)$ is the set of points $x \in A_n^\mathbb{Z}$ such that $\mathcal{L} (x) \subseteq \mathcal{L}^{(n)}(\tau)$.  
\index{symbols}{X@$X^{(n)}(\tau)$}%
We have then
\begin{displaymath}
  \cdots \edge{\tau_1}X^{(1)}(\tau)\edge{\tau_0}X^{(0)}(\tau).
    \end{displaymath}

Note that it may happen that $\mathcal{L}(X^{(n)}(\tau))$ is strictly contained in $\mathcal{L}^{(n)}(\tau)$ (if, for example, all $\tau_n$
are equal to a morphism which is not a substitution) or even empty
(if $|\tau_{[0,n)}|$ is bounded).

  We say
that $\tau$ is a \emph{directive sequence}
if for every $n$,
\begin{displaymath}
  \cL(X^{(n)}(\tau))=\cL^{(n)}(\tau).
  \end{displaymath}
 This
implies that $|\tau_{[0,n)}|$ tends to infinity with $n$
and that the language  $\cL^{(n)}(\tau)$ is extendable for every $n\ge 0$.
%We also say that $(\tau_n)$ is a $\Sa$-adic \emph{system}
%\index{subject}{S-adic@$\Sa$-adic!system}
Thus,  a constant sequence $\tau=(\sigma,\sigma,\ldots)$
is a directive sequence if and only if $\sigma$ is a substitution.

When $\tau$ is directive sequence,
we set $X(\tau) = X^{(0)}(\tau)$ and $\cL(\tau)=\cL^{(0)}(\tau)$.
We call $(X(\tau),S)$ the \emph{$\Sa$-adic shift} 
\index{subject}{S-adic@$\Sa$-adic!shift}
with \emph{directive sequence}~$\tau$.
\index{subject}{directive!sequence of morphisms}\index{subject}{S-adic@$\Sa$-adic!shift!directive sequence}%
We also say that $\tau$ is an $\Sa$-adic \emph{representation} of
$X=X(\tau)$.\index{subject}{S-adic@$\Sa$-adic!representation}

As a first example, let us see that $\Sa$-adic representations capture
the notion of substitutive shift.
\begin{example}
  Every substitutive shift $X(\sigma,\phi)$, where
  $\sigma:B^*\to B^*$ is a substitution
  and $\phi:B^*\to A^*$ is nonerasing, has an $\Sa$-adic representation
with $\Sa=\{\sigma,\phi\}$. Indeed, $X(\sigma,\phi)$
has the $\Sa$-adic representation $(\phi,\sigma,\sigma,\ldots)$.
\end{example}
As a second example, we recover the Arnoux-Rauzy shifts.
\begin{example}
  Let $s$ be a strict standard episturmian sequence with directive sequence $x=a_0a_1\cdots$ (see Section~\ref{sectionSturmianShifts}).
  The shift generated by $s$ has the $\Sa$-adic representation
  $\tau=(\tau_n)_{n\ge 0}$ where 
  $\tau_n=L_{a_n}$ and where each $L_a:A^*\to A^*$ is the elementary automorphism
  defined for $b\in A$ by
  $L_a(b)=ab$ for $b\ne a$ and $L_a(a)=a$. Indeed, every word
  in $\cL(s)$ is a factor of some $L_{a_0\cdots a_{n-1}}(a_n)$
  which is in $\cL^0(\tau)$. Conversely, since $s$ is strict,
  every letter appears infinitely often in $x$. This
  implies that all letters appear in $L_{a_0\cdots a_{n-1}}^{-1}(s)$
  and thus that $L_{a_0\cdots a_{n-1}}(a)$ appears in $s$ for every $a\in A$.
  This shows that $\cL^0(\tau)=\cL(s)$.
  \end{example}

As for BV-representations, we have the notion of \emph{telescoping}
\index{subject}{telescoping!of $\Sa$-adic sequence} of a
directive sequence. Given a  sequence of morphisms $\tau=(\tau_n)$
and a sequence $(n_m)$ of integers with $n_0=0<n_1<n_2<\ldots$,
a telescoping of $\tau$ with respect to $(n_m)$ is the
 sequence of morphisms $\tau'=(\tau_{[n_m,n_{m+1})})_{m\ge 0}$. If $\tau$ is a directive
sequence, then $\tau'$ is a directive sequence. Moreover
$\tau$ and $\tau'$ define the same shift $X(\tau)=X(\tau')$.

As for substitution shifts, in which  fixed points play
an important role, we have for directive sequences, the notion
of \emph{limit point}.
\index{subject}{limit!point of directive sequence}%
\index{subject}{point!limit}%
Let
$\tau=(\tau_n)_{n\ge 0}$ be an $\Sa$-adic system with $\tau_n:A_{n+1}^*\to A_n^*$.
A sequence $x\in A_0^\N$ is called a limit point of $\tau$
if there is a sequence $(w^{(n)})$ of sequences $w^{(n)}\in A_n^\N$
such that $w^{(n)}=\tau_n(w^{(n+1)})$ with $x=w^{(0)}$.

\begin{example}
Let $\sigma$ be a substitution with fixed point $x=\sigma^\omega(a)$.
Then $x$ is a limit point of the directive sequence
$(\sigma,\sigma,\ldots)$. Indeed, the constant sequence $w^{(n)}=\sigma^\omega(a)$
is such that $w^{(n+1)}=\sigma(w^{(n)})$ and $x=w^{(0)}$.
\end{example}
Observe that, contrary to the case of substitutive sequences,
every sequence has an $\Sa$-adic representation (Exercise~\ref{ExerciseSadicCassaigne}).

\subsection{Primitive directive sequences}

We say that the sequence $\tau$ is {\em primitive}
\index{subject}{primitive!sequence of morphisms}%
 if, for any $n\geq 1$, there exists $N>n$ such that for all $a \in A_N$, $\tau_{[n,N)}(a)$ contains occurrences of all letters of $A_{n}$.

Set $M_{[n,N)}=M_{\tau_n} M_{ \tau_{n+1}}\cdots  M_{\tau_{N-1}}$
\index{symbols}{M@$M_{[n,N)}$} where
$M_{\tau_n}$ is the composition matrix of $\tau_n$. Then
$\tau$ is primitive if for any $n\ge 1$, there exists $N>n$ such that
$M_{[n,N)} >0$, 

Observe that if $\tau$ is primitive, then 
$\langle\tau_{[1,n)}\rangle$ goes to infinity when $n$ increases.

When $\tau$ is primitive,
we can use  alternative definitions
of a limit point $x$.
\begin{enumerate}
\item[(i)] There is  a sequence $a_n\in A_n$ of letters such that $x$
is the limit of $\tau_{[0,n)}(a_n^\omega)$
\item[(ii)] There is  a sequence $a_n\in A_n$ of letters such that 
$\{x\}=\cap_{n\ge 0}[\tau_{[0,n)}(a_n)]$ where
$[w]$ is the cylinder defined by $w$.
\end{enumerate}
(See Exercise~\ref{exerciseLimitPoint}).

The following result generalizes the fact that every growing
morphism has a power with a fixed point (Proposition \ref{propositionGrowing}).
\begin{proposition}\label{propositionLimitPoint}
Every primitive $\Sa$-adic sequence has a limit point.
\end{proposition}
\begin{proof}
Let $\tau=(\tau_n)$ be a primitive directive sequence. For $a\in A_{n+1}$,
denote by $f_n(a)$ the first letter of $\tau_n(a)$. Since $A_0$
is finite, there is a letter $a$ in the
decreasing sequence of sets $\cap_nf_{[0,n)}(A_n)$. By construction,
there is a sequence $(a_n)$ such that $a_0=a$ and $a_n=f_n(a_{n+1})$.
Set $u_n=\tau_{[0,n)}(a_n)$.
By construction, every $u_n$ is a prefix of $u_{n+1}$
and since $\tau$ is primitive, their lengths
tend to infinity. Let $u^{(0)}$
be the unique one-sided sequence having all $u_n$ as prefixes.
By shifting the sequences $(\tau_n)$ and $(a_n)$
(at the first step, replace $(\tau_n)_{n\ge 0}$
by $(\tau_n)_{n\ge 1}$ and so on), we define
a sequence $u^{(k)}$ of sequences such that $u^{(k)}=\tau_k(u^{(k+1)})$.
Thus $u^{(0)}$ is a limit point of $\tau$.
\end{proof}

The following proposition generalizes the fact that a primitive
substitution shift is minimal (Proposition~\ref{propositionPrimitiveMinimal}).
\begin{proposition}\label{propositionSadicPrimitive}
If $\tau$ is primitive,  
it is a directive sequence  and $X(\tau)$ is  minimal. 
\end{proposition}
\begin{proof}
Since $\langle\tau_{[n,N)}\rangle\to\infty$ with $N$, the language $\cL^{(n)}(\tau)$
is extendable and thus $\tau$ is a directive sequence.

Let us show that $X(\tau)$ is uniformly recurrent. 
For this, let $u\in\cL(\tau)$. By definition of $\cL(\tau)$, there is an $n>1$
and a letter $a\in A_n$ such that $u$ appears in $\tau_{[1,n)}(a)$.
Since $\tau$ is primitive, there is an $N>n$ such that for
all $b\in A_N$, the letter $a$ appears in $\tau_{[n,N)}(b)$. Let $b,c\in A_N$
be such that $bc\in\cL^{(N)}(\tau)$. Then $u$ appears in 
$\tau_{[1,N)}(b)$ and in $\tau_{[1,N)}(c)$ at a distance at 
most equal to $2|\tau_{[1,N)}|$ (see Figure~\ref{figureSadicPrimitive}). This shows that $X(\tau)$ 
is uniformly recurrent.
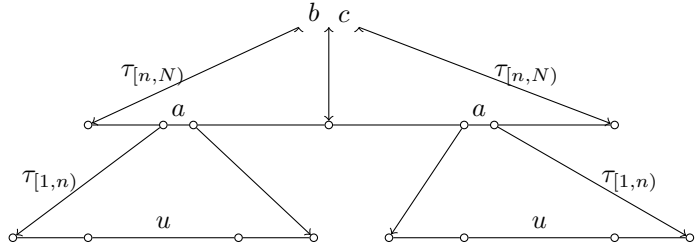
\begin{figure}[hbt]
\centering
\tikzset{node/.style={circle,draw,minimum size=0.1cm,inner sep=0pt}}
\tikzset{title/.style={circle,minimum size=0.1cm,inner sep=0pt}}
\begin{tikzpicture}
\node[title](b)at(4,3.5){$b$};\node[title](c)at(4.4,3.45){$c$};
\node[node](21)at(1,2){};\node[node](22)at(2,2){};\node[node](23)at(2.4,2){};
\node[node](24)at(4.2,2){};\node[node](25)at(6,2){};\node[node](26)at(6.4,2){};
\node[node](27)at(8,2){};
\node[node](31)at(0,0.5){};\node[node](32)at(1,0.5){};\node[node](33)at(3,0.5){};
\node[node](34)at(4,0.5){};\node[node](35)at(5,0.5){};\node[node](36)at(6,0.5){};
\node[node](37)at(8,0.5){};\node[node](38)at(9,0.5){};

\draw[left,->](3.8,3.3)edge node{$\tau_{[n,N)}$}(21);
\draw[left,->](4.2,3.3)edge node{}(24);
\draw[right,->](4.6,3.3)edge node{$\tau_{[n,N)}$}(27);
\draw[above](21)edge node{}(22);\draw[above](22)edge node{$a$}(23);
\draw[above](23)edge node{}(24);\draw[above](24)edge node{}(25);
\draw[above](25)edge node{$a$}(26);
\draw[above](26)edge node{}(27);
\draw[left,->](22)edge node{$\tau_{[1,n)}$}(31);
\draw[left,->](23)edge node{}(34);
\draw[left,->](25)edge node{}(35);
\draw[right,->](26)edge node{$\tau_{[1,n)}$}(38);
\draw[above](31)edge node{}(32);\draw[above](32)edge node{$u$}(33);
\draw[above](33)edge node{}(34);\draw[above](35)edge node{}(36);
\draw[above](36)edge node{$u$}(37);
\draw[above](37)edge node{}(38);
\end{tikzpicture}
\caption{The shift $X(\tau)$ is minimal.}\label{figureSadicPrimitive}
\end{figure}
\end{proof}
A morphism $\sigma:A^*\to B^*$ is said to be {\em left proper}
\index{subject}{left!proper!morphism}\index{subject}{morphism!left proper}%
 (resp. {\em right proper})
\index{subject}{right!proper!morphism}\index{subject}{morphism!right proper}%
 when there exist a letter $b \in B$ such that 
for all $a \in A$, $\sigma(a)$ starts with $b$ (resp., ends with $b$).
Thus it is proper if it is both left and right proper.

We say that  sequence $\tau=(\tau_n)$
of morphisms is {\em left proper} (resp.
\emph{ right proper}, resp. {\em proper})
\index{subject}{left!proper!sequence of morphisms}%
\index{subject}{right!proper!sequence of morphisms}%
\index{subject}{proper!sequence of morphisms}%
 whenever each morphism $\tau_n$ is 
left proper (resp.  right proper, resp. proper).
We also say that a shift is a left proper (resp.  right proper, resp.
primitive) $\Sa$-adic shift if there exists a left proper (resp.
  right proper, resp. primitive) sequence of morphisms
 $\tau$ such that $X = X(\tau)$.

Let us give another way to define $X(\tau)$ when $\tau$
 is primitive and proper. 
For a non erasing morphism $\sigma:A^*\to B^*$, let $\Omega (\sigma )$
\index{symbols}{Omega@$\Omega(\sigma)$} be the
closure of $\cup_{k\in \Z} S^k \sigma ( A^\Z )$.

\begin{lemma}
\label{lemma:omega}
Let $\tau = (\tau_n : A_{n+1}^*\to A_n^*)_{n \geq 0}$ be a primitive and proper sequence of morphisms.
Then, 
$$
X(\tau) = \cap_{n\ge 0} \Omega (\tau_{[0,n]} ) .
$$
\end{lemma}
\begin{proof}
It is equivalent to prove that
\begin{displaymath}
\cL(\tau)=\cap_{n\ge 1}\Fac(\tau_{[0,n)}(A_n^*))
\end{displaymath}
where $\Fac(L)$ denotes the set of factors of the words in $L$.
\index{symbols}{fac@$\Fac(L)$}%
Consider first $w\in \cL(\tau)$. By definition of $\cL(\tau)$, there
is an $n\ge 2$ and a letter $a\in A_n$
such that $w$ appears in $\tau_{[0,n)}(a)$. Since $\tau$ is primitive,
 for all $N>n$ large enough and every a letter $b\in A_N$, the letter
$a$ is a factor of $\tau_{[n,N)}(b)$. This proves that
$w\in\Fac(\tau_{[0,N)}(A_N)$ for arbitrary large $N$ and thus
that $w\in\cap_{n\ge 1}\Fac(\tau_{[0,n)}(A_n^*))$.

Conversely, let $\ell_n\in A_n$ (resp. $r_n\in A_n$) 
be the first letter (resp. last letter)
of all $\tau_n(a)$ for $a\in A_{n+1}$. Let $w\in\cap_{n\ge 1}\Fac(\tau_{[0,n)}(A_n^*))$. Since $\tau$ is primitive we may assume, telescoping if necessary
the sequence $\tau_n$, that $\langle \tau_n\rangle\ge 2$. For
each $n\ge 1$, there is some $w_n\in A_n^*$ such that $w$ is a factor
of $\tau_{[0,n)}(w_n)$ and we can choose $w_n$ of minimal length.
Then $|w|\ge 2^n(|w_n|-2)$. Thus there is an $n\ge 1$ such that
$w_n$ has length at most $2$. If $|w_n|=1$ we obtain the conclusion $w\in \cL(\tau)$. Assume that $w_n$ has length $2$. If $w_n\ne r_n\ell_n$, then
$w_n$ is a factor of some $\tau_{n+1}(a)$ for $a\in A_{n+1}$
and we are done. Otherwise, consider a letter $c\in A_{n+2}$ and 
 two consecutive letters $d,e$ of $\tau_{n+1}(c)$. Then $\tau_n\tau_{n+1}(c)$
has a factor $r_n\ell_n$ and thus, we can choose $w_{n+2}=c$.

\end{proof}

Observe that  the hypotheses that $\mathbf{\tau}$ is both primitive 
 and proper cannot be dropped in the previous statement.
Without these hypotheses, the inclusion $X(\tau) \subset \cap_{n\in \N} \Omega (\tau_{[0,n]} )$  still  holds (under the mild assumption that
$A_n\subset \Fac(\tau_n(A_{n+1}))$, but not  the  reverse inclusion, as shown by the
following examples.
\begin{example}
Take  for  the  directive sequence  $\tau$   the constant  sequence  equal to   $ \tau$, defined by $\tau (0) = 0010$ and $\tau (1)=1$
(this is the Chacon binary substitution). 
The directive sequence $\tau$ is not primitive, and $\cdots 111\cdots $ belongs to $\cap_{n\in \N} \Omega (\tau_{[0,n]} ) $ but not to $ X(\tau)$.
\end{example}
\begin{example}
In the case of  the non-proper constant  sequence  given by  $ \tau$ with $\tau (0) = 0100$ and $\tau (1) = 101$,  the fixed point $\tau^\omega(1\cdot 1)$ belongs to $\cap_{n\in \N} \Omega (\tau_{[0,n]} )$,  but not to $X(\tau)$, since $11$   appears in no element       of $X(\tau)$.
\end{example}

With a left proper morphism $\sigma:A^* \to B^*$ such that $b \in B$ is the first letter of all images $\sigma(a)$, $a \in A$, we associate the right proper morphism $\overline{\sigma}:A^* \to B^*$ defined by 
$b\overline{\sigma}(a) = \sigma(a)b$ for all $a \in A$. 

\begin{lemma}
\label{lemma:proper}
Let $X$ be an $\Sa$-adic shift  generated by the
primitive and left proper directive sequence $\tau = (\tau_n:A_{n+1}^* \to A_n^*)_{n \geq 1}$.
Then $X$  is also generated by 
the primitive and proper directive sequence $\tilde{\tau} = (\tilde{\tau}_n)_{n \geq 1}$, where for all $n$, $\tilde{\tau}_n = \tau_{2n-1} \overline{\tau}_{2n}$.
In particular,  if $\tau$ is unimodular, then so is $\tilde{\tau}$.
\end{lemma}
\begin{proof}
Each morphism $\tilde{\tau}_n$ is trivially proper.
It is also clear that the unimodularity of $\boldsymbol{\tau}$ is preserved in this process.
Now, let $\sigma:A^* \to B^*$ be a proper morphism such that $b \in B$ is the first letter of all images $\sigma(a)$, $a \in A$;
 for all $x\in A^\Z$, one has 
$
\overline{\sigma} (x) = S{\sigma}(x) .
$
Together with Lemma \ref{lemma:omega} this   ends the proof.
\end{proof}

\subsection{Unimodular $\Sa$-adic shifts.}
We  also say that  sequence $\tau=(\tau_n)$ of morphisms
is {\em unimodular}\index{subject}{unimodular!sequence of morphisms}
whenever, for all $n \geq 1$, $A_{n} = A$ and the matrix $M(\tau_n)$
is unimodular, that is, has determinant of absolute value 1.
We say that a shift is a \emph{unimodular} (resp. \emph{proper},
resp. \emph{primitive}) \index{subject}{proper!S-adic@$\Sa$-adic shift}
\index{subject}{S-adic@$\Sa$-adic!shift!proper}%
\index{subject}{primitive!S-adic@$\Sa$-adic shift}%
\index{subject}{S-adic@$\Sa$-adic!shift!primitive}%
 $\Sa$-adic shift \index{subject}{unimodular!$\Sa$-adic shift}
\index{subject}{S-adic@$\Sa$-adic!shift!unimodular}if there exists a unimodular (resp. proper, resp. primitive) sequence of morphisms
 $\tau$ such that $X = X(\tau)$.

\begin{lemma}
\label{lemma:aperiodicity}
All   primitive unimodular proper $\Sa$-adic shifts are aperiodic. 
\end{lemma}
\begin{proof}
Let $\tau$ be a primitive unimodular proper directive sequence on the alphabet $A$ of cardinality $d \geq 2$. Suppose that the shift $X(\tau)$
is periodic, that is $X(\tau)=\{x,Tx,\ldots,T^{p-1}x\}$.

Since $\tau$ is primitive,
there is some $n\ge 1$ such that $|\tau_{[1,n]}(a)|\ge p$ for all $a\in A$.

Let $b\in A$ and set $\tau_{n+1} (b) = b_0 b_1 \cdots b_k$. 
Since the directive sequence $\tau$ is proper,  $b_0 b_1 \cdots b_k b_0$ is also a word in $\mathcal{L}^{(n+1)}  (\tau)$. Then 
$w=\tau_{[1,n]}(b_0b_1\cdots b_kb_0)$ is a word which has both
period $p$ (as a factor of $x$) and period $|\tau_{[1,n+1]}(b)|$.
By Fine-Wilf Theorem (Exercise~\ref{exerciseFineWilf}), 
we have $|\tau_{[1,n+1]}(b)|\equiv0\bmod p$. But then
every row of the matrix $M(\tau_{[1,n+1]})$ has a sum
divisible by $p$ and thus its
 determinant   is a multiple of $p$,
which contradicts the  unimodularity of  $\tau$.
\end{proof}

\subsection{Recognizability and unimodular $\Sa$-adic shifts.}

We have seen in Section~\ref{sectionSubstitutionSystems}
the definition of recognizability of morphisms.
Let us recall this definition here.
 Let $\varphi:A^*\to B^*$ is a nonerasing morphism. Let $X$
with $X\subset A^\Z$
be a shift space and let $Y$ be the closure of $\varphi(X)$ under the
shift. Every $y\in Y$ has a representation as $y=S^k\varphi(x)$
with $x\in X$ and $0\le k<|\varphi(x_0)|$.
We say that 
$\varphi$ is \emph{recognizable} on the shift $X$ for the point $y$
if $y$ has only one such representation.

We say that $\varphi$ is \emph{recognizable}\index{subject}{recognizable!morphism} on $X$ if it is
recognizable on $X$ for every point $y\in Y$. We also
say that $\varphi$ is \emph{recognizable} on $X$ for aperiodic points
if $\varphi$ is recognizable on $X$ for every aperiodic point $y\in Y$.
Finally, we say that $\varphi$ is \emph{recognizable for aperiodic points}
if it is recognizable on the full shift for aperiodic points. 

Recall from Section~\ref{sectionRecognizable} that recognizability
can expressed using the system
$X^\varphi=\{(x,i)\mid x\in X, 0\le i<|\varphi(x_0)|\}$
and the morphism $\widehat{\varphi}:X^\varphi\to Y$ defined by
$\widehat{\varphi}(x,i)=S^i\varphi(x)$. Then $\varphi$
is recognizable on $X$ for $y$ if
$\widehat{\varphi}^{-1}(y)$ has only one element.

This formulation  will be helpful in the proof of the following
statement.

\begin{proposition}\label{propositionCompositionRecognizable}
Let $\sigma:A^*\to B^*$ and $\tau:B^*\to C^*$ be morphisms.
Let $X\subset A^\Z$ be a shift space and let $Y$ be the subshift
generated by $\sigma(X)$. Then 
\begin{enumerate}
\item $\tau\circ\sigma$ is
recognizable on $X$ if and only if $\sigma$ is recognizable
on  $X$ and $\tau$ is recognizable on $Y$.
\item If $\sigma$ is recognizable on $X$ for aperiodic points
and $\tau$ is recognizable on $Y$ for aperiodic points,
then $\tau\circ\sigma$ is recognizable on $X$ for aperiodic points.
\item If $\tau\circ\sigma$ is recognizable on $X$ for aperiodic points,
then  $\tau$ is recognizable on $Y$ for aperiodic points.
\end{enumerate}
\end{proposition}
\begin{proof}
Let $Z$ be the subshift of $C^\Z$ generated by $\tau(Y)$.
Set $\rho=\tau\circ\sigma$. We have $\hat{\rho}=\hat{\tau}\circ\alpha$
(see Figure~\ref{figureAlpha})
where $\alpha: X^\rho\to Y^\tau$ is the following map.
\begin{figure}[hbt]
\centering
\begin{tikzpicture}
\node(X)at(0,2){$X^\rho$};\node(Z)at(2,2){$Z$};
\node(Y)at(1,0){$Y^\tau$};
\draw[above,->](X)edge node{$\hat{\rho}$}(Z);
\draw[left,->](X)edge node{$\alpha$}(Y);
\draw[right,->](Y)edge node{$\hat{\tau}$}(Z);
\end{tikzpicture}
\caption{The map $\alpha$.}\label{figureAlpha}
\end{figure}
For each $(x,k)\in X^\rho$, there is a unique pair $(i,j)$ such that
$k=|\tau(b_0\cdots b_{i-1})|+j$ with $j<|\tau(b_i)|$
with $\sigma(x_0)=b_0\cdots b_{|\sigma(x_0)|}$. Then
\begin{displaymath}
\alpha(x,k)=(\hat{\sigma}(x,i),j)
\end{displaymath}

1. It is clear that
 $\alpha$ is injective if and only if $\hat{\sigma}$ is injective.
 Since $\alpha$ is surjective, it follows that $\hat{\rho}$
is injective if and only $\hat{\sigma}$ and $\hat{\tau}$
are injective.

2 and 3 follow easily since $\hat{\tau}$ sends periodic points
to periodic points.
\end{proof}
Recall (see Exercise~\ref{exerciseElementarySubstitution}) that
a morphism $\varphi:A^*\to B^*$ is called \emph{elementary}
\index{subject}{elementary!morphism}
if for every decomposition $\varphi=\beta\circ\alpha$
with $\alpha:A^*\to C^*$ and $\beta:C^*\to B^*$, one has
$\Card(C)\ge\Card(A)$. In particular, an elementary
morphism is nonerasing.
\begin{proposition}\label{propositionElementaryMorphism}
If $\varphi:A^*\to B^*$ be  is elementary,
 the extension of $\varphi$ to $A^\N$ is  injective.
\end{proposition}
The proof is Exercise~\ref{exerciseElementarySubstitution}.
%\begin{proof}
%We argue by induction on $\|\varphi\|=\sum_{a\in A}|\varphi(a)|$.
%The statement is true if $\|\varphi\|=\Card(A)$ since in this case
%an elementary morphism is a bijection from $A$ to $B$.
%Assume that the extension of $\varphi$ to $A^\N$ is not injective.
%Let $x,x'\in A^\N$ be such that $\varphi(x)=\varphi(x')$.
%We may assume that $x$ begins with $a\in A$ and $x'$ with $a'\ne a$.
%If $\varphi(a)=\varphi(a')$, then $\varphi$ is clearly not elementary.
%Otherwise, we may assume that $|\varphi(a')|>|\varphi(a)|$. Set
%$\varphi(a')=\varphi(a)u$. Let $C=(A\cup a'')\setminus a'$
%where $a''$ is a new letter. Let $\alpha:A^*\to C^*$
% be defined by $\alpha(a')=aa''$ and $\alpha(b)=b$ for $b\ne a'$.
%Then $\varphi=\beta\circ\alpha$ with 
%\begin{displaymath}
%\beta(c)=\begin{cases}u&\mbox{ if $c=a''$}\\\varphi(c)&\mbox{ otherwise}
%\end{cases}
%\end{displaymath}
%Then $\|\beta\|<\|\varphi\|$ and thus $\beta$ is not elementary by induction
%hypothesis. Therefore $\varphi$ is not elementary.
%\end{proof}

We will now prove the following result. It holds in
particular when $M(\sigma)$ is unimodular.
\begin{theorem}\label{theo:BSTY}
Let $\sigma:A^*\to B^*$ be an elementary morphism 
%such that $M(\sigma)$ has rank $\Card(A)$. 
Then $\sigma$ is recognizable for aperiodic points.
\end{theorem}
\begin{proof}
We proceed by induction on $|\sigma|=\max_{a\in A}|\sigma(a)|$. 

If $|\sigma(a)|=1$
for all $a\in A$, then  $\sigma$ is a permutation of $A$, whence recognizable
in the full shift.

Assume now that $\sigma$ is not recognizable on aperiodic points.
Since $M(\sigma)$ has rank $\Card(A)$, the morphism $\sigma$
is elementary. Thus, by Proposition~\ref{propositionElementaryMorphism},
$\sigma$ is injective on $A^\N$.
Thus there exist $x,x'\in A^\N$ and $w$ with $0<|w|<|\sigma(x_0)|$
such that $\sigma(x)=w\sigma(x')$ for some proper suffix $w$ of $\sigma(a')$.
Set $\sigma(a')=vw$. We can then write $\sigma=\sigma_1\circ\tau_1$
whith $\tau_1:A^*\to A_1^*$ and $\sigma_1:A_1^*\to B^*$
and $A_1=A\cup\{a''\}$ where $a''$ is a new letter. We have
$\tau_1(a')=a'a''$ and $\tau_1(a)=a$ otherwise. Next
$\sigma_1(a')=v$, $\sigma_1(a'')=w$ and $\sigma_1(a)=\sigma(a)$ otherwise.
Note that $|\sigma_1|=|\sigma|$.
Since $\tau_1$ is injective on $A^\N$,
$\sigma_1$ is not injective on $A_1^\N$. By Proposition~\ref{propositionElementaryMorphism} again, $\sigma_1$ is not elementary. Thus $\sigma_1=\sigma_2\circ\tau_2$
with $\tau_2:A_1^*\to A_2^*$, $\sigma_2:A_2^*\to B^*$ and
$\Card(A_2)<\Card(A_1)=\Card(A)+1$. On the other hand, since $\sigma=(\sigma_2\circ\tau_2)\circ \tau_1$
is elementary, we have $\Card(A)\le\Card(A_2)$. Thus
\begin{displaymath}
\Card(A)\le\Card(A_2)\le\Card(A)
\end{displaymath}
and consequently $\Card(A)=\Card(A_2)$.
Since $|\sigma_2|<|\sigma_1|=|\sigma|$, 
we may apply the induction hypothesis to obtain that
$\sigma_2$ is recognizable for aperiodic points. Since $\tau_1$
and $\tau_2$ are recognizable on the full shift, we obtain
the conclusion by Proposition \ref{propositionCompositionRecognizable}.
\end{proof}
We will use the following corollary.
\begin{corollary}\label{corollaryBSTY}
Let $\sigma:A^*\to B^*$ be a morphism 
such that $M(\sigma)$ has rank $\Card(A)$. 
Then $\sigma$ is recognizable for aperiodic points.
\end{corollary}
This follows directly from Theorem~\ref{theo:BSTY}.
Indeed, a morphism $\sigma:A^*\to B^*$ such that $M(\sigma)$ has rank $\Card(A)$
is elementary.
\begin{example}
The matrix of the morphism $\sigma:a\to ab,b\to aa$ is
\begin{displaymath}
M(\sigma)=\begin{bmatrix}1&2\\1&0\end{bmatrix}
\end{displaymath}
which has rank $2$. The morphism is not recognizable for
$a^\infty$ but for every other point $y$, the occurrence of a $b$
is enough to determine the representation of
$y$ as $y=S^k\sigma(x)$. For example, $y_0=b$ forces
$k=1$ and $x_0=a$. All the other values of $x_n$ are then determined.
\end{example}

We use Corollary~\ref{corollaryBSTY} to prove the following.
\begin{proposition}
\label{prop:minimalSSIprimitif}
Let $X$ be an $\Sa$-adic shift with unimodular  and proper
directive sequence $\tau$.
Then $X$ is aperiodic and minimal if and only if $\tau$ is primitive.
\end{proposition}

\begin{proof}
Recall that any $\Sa$-adic subshift with a primitive directive sequence is minimal (Proposition~\ref{propositionSadicPrimitive}) and that aperiodicity is proved in Lemma \ref{lemma:aperiodicity}.

We only have to show that the condition is necessary.
We assume that $X$ is aperiodic and minimal.
For all $n \geq 1$, $X^{(n)}(\tau)$ is trivially aperiodic. 
Let us show that it is minimal.

Assume by contradiction that for some $n \geq 1$, $X^{(n)}(\tau)$ is minimal, but not $X^{(n+1)}(\tau)$.
There exist $u \in \mathcal{L}(X^{(n+1)}(\tau))$ and $x \in X^{(n+1)}(\tau)$ such that $u$ does not occur in $x$. Since $\tau$ is unimodular,
by Corollary~\ref{corollaryBSTY}, $\tau_n$
is recognizable for aperiodic points
and thus $\{\tau_n([v]) \mid v \in \La(X^{(n+1)}(\tau)) \cap A^{|u|}\}$ is a finite clopen  partition of 
$\tau_n(X^{(n+1)}(\tau))$.
Thus, considering $y = \tau_n(x)$, by minimality of $X^{(n)}(\tau)$, there exists $k \geq 0$ such that $S^k y$ is in $\tau_n([u])$.
Take $z \in [u]$ such that $S^ky = \tau_n(z)$.
Since $y$ is aperiodic and since we also have $S^k y = S^{k'} \tau_n(S^\ell x)$ for some $\ell \in \N$ and $0 \leq k' < |\tau_n(x_\ell)|$,
we  obtain  that $\tau_n(z) = S^{k'} \tau_n(S^\ell x)$ with $z \in [u]$, $S^\ell x \notin [u]$ and $0 \leq k' < |\tau_n(x_\ell)|$.
This contradicts the fact that $\tau_n$ is recognizable for aperiodic points.

We now show that $\lim_{n \to +\infty} \langle\tau_{[0,n)}\rangle = +\infty$.
We again proceed by contradiction, assuming that $(\langle\tau_{[0,n)}\rangle)_{n \geq 1}$ is bounded.
Then there exists $N > 0$ and a sequence $(a_n)_{n \geq N}$ of letters in $A$ such that for all $n \geq N$, $\tau_n(a_{n+1}) = a_n$.
We claim that there are arbitrary long words of the form $a_N^k$ in $\La(X^{(N)}(\tau))$ which contradicts the fact that $X^{(N)}(\tau)$ is minimal and aperiodic.
Since $\tau$ is proper, for all $n \geq N$ and all $b \in A$, $\tau_n(b)$ starts and ends with $a_n$.
Since $|\tau_{[0,n)}|$ goes to infinity, there exists a sequence $(b_n)_{n \geq N}$ of letters in $A$ such that $|\tau_{[N,n)}(b_n)|$ goes to infinity and for all $n \geq N$, $b_n$ occurs in $\tau_n(b_{n+1})$.
This implies that there exists $M \geq N$ such that for all $n \geq M$, $b_n \neq a_n$ and, consequently, that $\tau_n(b_{n+1}) = a_n u_n$ for some word $u_n$ containing $b_n$.  
It is then easily seen that, for all $k \geq 1$, $a_M^k$ is a prefix of $\tau_{[M,M+k)}(b_{M+k})$, which proves the claim.

We finally show that $\tau$ is primitive.
If not, there exist $N \geq 1$ and a sequence $(a_n)_{n \geq N}$ of letters in $A$ such that for all $n > N$, $a_N$ does not occur in $\tau_{[N,n)}(a_n)$.
Since $(|\tau_{[N,n)}(a_n)|)_n$ goes to infinity, this shows that there are arbitrarily long words in $\La(X^{(N)}(\tau))$ in which $a_N$ does not occur.
%Since $\tau$ is unimodular, there is also a sequence $(a'_n)_{n \geq N}$ of letters in $A$ such that $a_N = a_N'$ and for all $n \geq N$, $a_n'$ occurs in $\tau_n(a_{n+1}')$.
%Again using the fact that $|\tau_{[N,n)}(a'_n)|$ goes to infinity, this shows that $a_N$ belongs to $\La(X^{(N)}(\tau))$.
We conclude that $X^{(N)}(\tau)$ is not minimal, a contradiction.
\end{proof}

%%%%%%%%%%%%%%%%%%%%%
\section{Dimension groups of unimodular $\Sa$-adic shifts.}
\label{sectionDGSadic}

In this section we first  prove a key result of this chapter, namely Theorem \ref{theo:cohoword},  which   states that,
for primitive unimodular $\Sa$-adic shifts,
the group $K^0(X,T)$ is generated, as an additive group,
  by the classes  of the   characteristic functions of letter cylinders. 
We then  deduce a simple expression for  the dimension group of primitive unimodular proper $\Sa$-adic subshifts.

\subsection{From letters  to  factors}
We recall  that $\chi_U$ stands for the characteristic function of the set $U$.

\begin{theorem}
\label{theo:cohoword}
Let $X$ be a primitive unimodular proper $\Sa$-adic shift.
Any function  $f\in C (X, \mathbb{Z})$ is cohomologous to some integer linear combination of  the form  $\sum_{a\in A} \alpha_a \chi_{[a]}\in C (X,\mathbb{Z})$. Moreover, the classes  $[\chi_{[a]}]$, $a\in A$, are $\mathbb{Q}$-independent.
\end{theorem}

\begin{proof}
Let $\boldsymbol{\tau} = (\tau_n:A^* \to A^*)_{n\geq 0}$ be a primitive unimodular proper  directive sequence of $X$. Hence $X = X(\boldsymbol{\tau})$.
Using Proposition~\ref{prop:minimalSSIprimitif}, all shift spaces $X_{\boldsymbol{\tau}}^{(n)}$ are minimal and aperiodic and $\min_{a \in A} |\tau_{[0,n)}(a)|$ goes to infinity when $n$ increases.

Let us show that the group $H(X,S )=C(X, {\mathbb Z})/ \partial C(X,{\mathbb Z}) $ is spanned by the set   of  classes  of characteristic functions of  letter
cylinders
$\{ [\chi_{[a]}] \mid a \in A \}$.
By Corollary~\ref{corollaryBSTY} and using the fact that $X$ is minimal and aperiodic, one has,   for all  positive  integer   $n$, that  
$$
\mathcal{P}_n = \{ S^k \tau_{[1,n]} ([a]) \mid 0\leq k < |\tau_{[0,n]} (a)| , a \in A \}
$$
is a finite  partition of $X$ into clopen sets.  This provides a   family of nested Kakutani-Rokhlin  partitions in towers.

We first claim that  $H(X,S ) $ is spanned by the set of classes     $\cup_n \Omega_n $, where
$$\Omega_n =  \{ [\chi_{\tau_{[0,n]} ([a])}] \mid a \in A \} \quad n \geq 1.$$
In other words,  $H(X,S)$ is spanned  by the set of  classes  of  characteristic functions of  bases   of the  sequence of   partitions  $({\mathcal P})_n$.
It suffices to check that, for all $u^- .u^+\in \La (X)$, the class $[\chi_{[u^-.u^+]}]$ is a linear integer combination of elements belonging
 to some $\Omega_n$.

Let  us check this assertion. Let $u^-.u^+$ be in $\La (X)$. 
Since $\min_{a \in \A} |\tau_{[0,n)}(a)|$ goes to infinity, there exists $n$ such that  $|u^-|,|u^+| < \min_{a\in \A} |\tau_{[0,n)} (a) |$.
The directive sequence $ \boldsymbol{\tau}$ being proper, there exist words $w,w'$ with  respective  lengths $|w|=|u^-|$ and $|w'|=|u^+|$
 such that all images $\tau_{[0,n]} (a)$ start with $w$ and end with $w'$.

Let $x \in [u^-. u^+]$. 
Let $a \in A$ and $k \in \mathbb{N}$, $0\leq k < |\tau_{[0,n]} (a)|$, such that $x$ belongs to the atom $S^k \tau_{[0,n]} ([a])$.
Observing that $\tau_{[0,n]} ([a])$ is included in $[w'.\tau_{[0,n]} (a)w]$, this implies that the full  atom $S^k \tau_{[0,n]} ([a])$ is included in $[u^-. u^+]$.
Consequently $[u^-. u^+]$ is a finite union of atoms in $\mathcal{P}_n$.
But each characteristic function  of   an atom  of the form   $ S^k \tau_{[0,n]} ([a]) $   is   cohomologous to    $\chi_{\tau_{[0,n]} ([a])}$. 
  This thus    proves the claim. 

Now we claim that each element of $\Omega_n $ is a linear integer combination of elements in $\{ [\chi_{[a]}] \mid a \in A \}$.
More precisely,  let us  show that $\chi_{\tau_{[0,n]}([b])}$ is cohomologous to 
$$
\sum_{a\in A}  (M_{[0,n]}^{-1})_{b,a}\chi_{[a]} .
$$

Let $a\in A$ and $n\geq 0$.
One has $[a] = \cup_{B \in \mathcal{P}_n } (B \cap [a]) $ and thus $\chi_{[a]} $ is cohomologous to the map 
$$
\sum_{b\in A}  (M_{[0,n]})_{a,b}\chi_{\tau_{[0,n]}([b])} ,
$$
by using  the fact  that the  maps $  \chi_{S^k \tau_{[0,n]} ([a]) }$   are   cohomologous to    $\chi_{\tau_{[0,n]} ([a])}$.
This means that for $U = ([\chi_{[a]} ] ) _{a\in A} \in H (X,S)^A$ and $V = ([\chi_{\tau_{[0,n]}([a])}])_{a\in A} \in H (X,S)^A$, one has 
$$
U = M_{[0,n]}V 
$$
and as a consequence $V = M_{[0,n]}^{-1}U$.
This proves the claim and the first part of the theorem.
 
To show the independence, suppose that  there exists some row vector $\alpha = (\alpha_a )_{a\in A} \in \mathbb{Z}^A$ such that $\sum_a \alpha_a   [\chi_{[a]}] = 0$.
Hence there is some $f\in C (X, \mathbb{Z} )$ such that $\sum_a \alpha_a   \chi_{[a]} = f\circ S - f$.
The morphisms of the directive sequence $\boldsymbol{\tau}$ being proper, for all $n$, there are letters $a_n, b_n$ such that all images $\tau_n (c)$, $c \in A$, start with $a_n$ and end with $b_n$.
From this, it is classical to check that $(\mathcal{P}_n)_n$ generates the topology of $X$ (the proof is the same as in the proof
of Lemma \ref{lemmaPartitionProperSubstitution} which  is concerned with  the particular case $\tau_{n+1}=\tau_n$ for all $n$).
 
We now  fix some $n$ for which  $f$ is constant on each atom of $\mathcal{P}_n$.
Observe that for all $x\in X$ and all $k \in \N$, one has $f (S^k x ) - f(x) = \sum_{j=0}^{k-1} \alpha_{x_j}$. 
Let $c\in A$ and $x\in \tau_{[0,n+1]} ([c])$.
Then, $x$ and $S^{ |\tau_{[0,n+1]} (c)|} (x)$ belong to  $\tau_{[0,n]} ([a_{n+1}])$.
Hence, $f (S^{ |\tau_{[0,n+1]} (c)|} x ) - f(x) = 0 $, and thus
$$
(\alpha M_{[0,n]})_c  = \sum_{j=0}^{|\tau_{[0,n+1]} (c)|-1} \alpha_{x_j} = 0 .
$$ 
This holds for all $c$, hence $\alpha M_{[0,n]} = 0$, which yields $\alpha = 0$, by invertibility of the matrix $ M_{[0,n]}$.
\end{proof}

Observe that in the previous result, we can relax the assumption of minimality. Indeed, 
one checks that the same proof works if we assume that $X$ is  aperiodic  (recognizability then holds by Corollary~\ref{corollaryBSTY})
 and that  $\min_{a \in A} |\tau_{[1,n)}(a)|$ goes to infinity.

We now derive  two corollaries  from Theorem~\ref{theo:cohoword}   dealing respectively with   invariant measures  and with the image  subgroup.

\begin{corollary}\label{coro:measures}
Let $X$ be a primitive unimodular proper $\Sa$-adic shift over the alphabet $A$ and let $\mu, \mu' \in \mathcal M X$.
If $\mu$ and $\mu'$ coincide on  the  letters, then they are equal, 
that is, if $\mu([a]) = \mu'([a])$ for all a in $A$, then $\mu(U) = \mu'(U)$, for any clopen subset $U$ of $X$.
\end{corollary}

\begin{corollary}\label{coro:freq}
Let $X$ be a primitive  unimodular proper $\Sa$-adic shift over  the alphabet $ A$.
The image subgroup of $X$  satisfies
$$I(X,S) = \bigcap_{\mu\in\mathcal{M}(X,S)}\left\lbrace\sum_{a\in{A}}\mathbb{Z}\mu([a])\right\rbrace.$$
\end{corollary}

\subsection{An explicit description of the dimension group}
Theorem~\ref{theo:cohoword} allows us
to give a precise description of the dimension group of primitive unimodular proper $\Sa$-adic shifts.

\begin{theorem}\label{theo:dg}
Let $X$ be a primitive unimodular and proper $\Sa$-adic shift. 
The linear map $\Phi : H (X,S) \to \mathbb{Z}^A$ defined by $\Phi (  [\chi_{[a]}] ) = e_a$, where $\{ e_a \mid a \in A \}$ is the canonical base of $ \mathbb{Z}^A$,
defines an isomorphism of dimension groups from $K^{0}(X,S)$ onto
\begin{align}
\label{align:ZdGD}
\left( {\mathbb Z} ^A, \,  \{ {\bf x} \in {\mathbb Z}^A \mid \langle {\bf x},  \boldsymbol{\mu} \rangle > 0 \mbox{ for all }  
 \mu \in  {\mathcal M} (X,S) \}\cup \{{\mathbf 0}\},\,  \bf{1}\right ),
 \end{align}
where the entries of $\bf{1}$ are equal to $1$
and $\boldsymbol{\mu}$ is the vector $(\mu([a]))_{a\in A}$.
\end{theorem}  
\index{subject}{dimension group!of unimodular shift}
\begin{proof}
From Theorem \ref{theo:cohoword}, $\Phi$ is well defined and is a group isomorphism from $H (X,S)$ onto $\Z^d$.
We obviously have $\Phi (  [1] ) = \Phi(\sum_{a \in A} [\chi_{[a]}]) = \bf{1}$ and it remains to show that 
\[
\Phi(H^+(X,S)) 
= 
\{ {\bf x} \in {\mathbb Z}^A \mid \langle {\bf x},  \boldsymbol{ \mu} \rangle > 0 \mbox{ for all }  
 \mu \in  {\mathcal M} (X,S) \}
 \cup 
 \{{\mathbf 0}\}.
\]
Any element of $H^+(X,S)$ is of the form $[f]$ for some $f \in C(X,\N)$.
From Theorem \ref{theo:cohoword}, there exists a unique vector ${\bf x} = (x_a)_{a \in A}$ such that $[f] = \sum_{a \in A} x_a [\chi_{[a]}]$.
Since $f$ is non-negative, we have, for any $\mu \in  {\mathcal M} (X,S)$, 
\[
	\langle	\Phi([f]),\boldsymbol{\mu} \rangle
	=
	\sum_{a \in A} x_a \mu([a])
	=
	\int f d\mu
	\geq 0,
\]
with equality if and only if $f=0$ (in which case ${\bf x} = {\bf 0}$).

For the other inclusion, assume that ${\bf x} = (x_a)_{a \in A} \in {\mathbb Z}^A$ satisfies $\langle {\bf x},  \boldsymbol{ \mu} \rangle > 0 \mbox{ for all }  
 \mu \in  {\mathcal M} (X,S)$ (the case ${\bf x} = {\bf 0}$ is trivial).
We consider the function $f = \sum_{a \in A} x_a \chi_{[a]}$. According to 
  Proposition~\ref{lemma2},  the existence of  $f' \in [f]$ such that $f'$ is non-negative is  equivalent to  the   existence of   a   lower bound for  ergodic sums $f^{(n)}$.
  Assume by contradiction that there exists a  point $x \in X$ such that the sequence $\left( \sum_{k=0}^n f \circ S^k(x)\right)_{n \geq 0}$ is not bounded from below.

There is a an increasing sequence of positive integers $(n_i)_{i \geq 0}$ such that 
\[
	\lim_{i \to +\infty} \sum_{k=0}^{n_i-1} f \circ S^k(x) = -\infty.
        \]
         Consider, as in the proof of Theorem~\ref{theoremKrylovBogolioubov},
  the sequence
  \begin{displaymath}
    \mu_n=\frac{1}{n}\sum_{i=0}^{n-1}\delta_{S^ix}
  \end{displaymath}
  where $\delta_x$ is the Dirac measure at $x$. 
  Extracting a subsequence $(m_i)_{i \geq 0}$ of $(n_i)_{i \geq 0}$,
  we find a subsequence of $(\mu_n)$ converging
  to an invariant measure
  $\mu$. Thus there exists $\mu \in \mathcal{M}(X,S)$ satisfying
\[
	\langle {\bf x}, \boldsymbol{\mu} \rangle = \int f d\mu = \lim_{i \to +\infty} \frac{1}{m_i} \sum_{k=0}^{m_i-1} f \circ S^k(x) \leq 0,
\]
which contradicts our hypothesis.
The sequence $\left( \sum_{k=0}^n f \circ S^k(x)\right)_{n \geq 0}$ is thus bounded from below and we conclude  by using Lemma~\ref{lemma2}.
\end{proof}

\begin{remark}
We cannot remove the hypothesis of being left or right proper in Theorem \ref{theo:dg}.
Consider indeed  the subshift $X$ defined by the  primitive unimodular non-proper substitution $\tau$  defined over $\{a,b\}^*$ as  $\tau \colon   a   \mapsto aab,  b \mapsto ba$.
  The dimension group  of $X$  is  isomorphic to 
\begin{displaymath}
\left( {\mathbb Z}^3 , \left\{ {\mathbf x} \in {\mathbb Z}^3  :    \langle {\mathbf x},  {\mathbf v}  \rangle  > 0 \right\} , (1,2,1) \right)
\end{displaymath}
where ${\mathbf v}=  [ 2\lambda,\lambda+1,1]$ with
$\lambda=(1+\sqrt{5})/2$ (Exercise~\ref{exercise(aab,ba)}).
Thus, although the matrix of $\tau$ is unimodular, the dimension
group is $\Z^3$ and not $\Z^2$.

\end{remark}

%%%%%%%%%%%%%%%%%%%%%%%%%%%%%%%%%%%%%
% section derivatives
%%%%%%%%%%%%%%%%%%%%%%%%%%%%%
\section{Derivatives of substitutive sequences}\label{sectionDerivatives}
In this section, we will prove a finiteness result characterizing
substitutive shifts. It can be viewed as a counterpart
of Boshernitzan-Carroll Theorem concerning interval exchange transformations
(Theorem~\ref{theoremBoschernitzanCarroll}). \marginpar{A expliquer mieux}

Let $X$ be a minimal shift space and let $u$ be a word of $\cL(X)$. Let 
$\varphi_u:A_u^*\to A^*$ 
be a coding morphism for the set $\RR'_X(u)$ of left return words
to $u$. We normalize the alphabet $A_u$\index{symbols}{B@$B_u$} by setting
$A_u=\{0,1,\ldots,n-1\}$ with $n=\Card(\RR'_X(u))$.
Recall from Section~\ref{sectionInduced}
that the \emph{derivative shift}\index{subject}{derivative!shift}
 of $X$ with respect to the cylinder $[u]$ is
the shift  $Y=\{y\in A_u^\Z\mid \varphi_u(y)\in X\}$.
%The map $\varphi_u$ is injective since $\RR'_X(u)$ is a circular code
%by Proposition~\ref{propositionReturnCircular}.
For $x\in X$ and $u\in\cL(x)$, the \emph{derivative sequence}\index{subject}{derivative!sequence}
 of $x$ with
respect to $u$ or $u$-derivative of $x$, denoted $\D_u(x)$
\index{symbols}{D@$\D_u(x)$} is 
defined as follows. Let $i\ge 0$ be the least integer
such that $(T^ix)^+$ begins with $u$, that is,
such that $T^ix$ is in $[u]$. Then $\D_u(x)$
is the unique $y\in Y$
such that $T^ix=\varphi_u(y)$. 
In particular, when $u$ is a prefix of $x^+$, we have $\D_u(x)=y$
where $y$ is the unique word of $Y$ such that $x=\varphi_u(y)$.
Thus, by definition, we have in this case
\begin{equation}
x=\varphi_u(\D_u(x)).\label{eqD_u(x)}
\end{equation}
%In view of comparing the different $\D_u(x)$,
%we assume that the alphabet $B_u$ is normalized as 
%$B_u=\{0,1,\ldots,\Card(\RR_X(u))-1\}$.
\subsection{A finiteness condition}
We will prove the following characterization of primitive substitutive sequences
by a finiteness condition.
\begin{theorem}
\label{theoremCharacterisationSubstitutive}
Let $x$ be a uniformly recurrent two-sided sequence. The following conditions
are equivalent
\begin{enumerate}
\item[\rm(i)]
The sequence $x$ is  primitive substitutive.
\item[\rm(ii)]
 The set of  $u$-derivative sequences of $x$ is finite, $u$ being a word of $\cL(x)$.
\item[\rm(iii)]
The set of its $u$-derivative sequences is finite, $u$ being a prefix of $x^+$.
\end{enumerate}
\end{theorem}

Note that, as a corollary of 
Theorem~\ref{theoremCharacterisationSubstitutive}, we obtain that
every shift of a substitutive sequence is substitutive. Indeed,
if $y=T^nx$ is a shift of a substitutive sequence $x$,
we have $\cL(x)=\cL(y)$ and
the set of $u$-derivatives of $x$ and $y$ are the same.
Thus $y$ is substitutive.
A direct proof is given in Exercise~\ref{exerciseShiftDerivative}.

Note also that a uniformly recurrent substitutive two-sided sequence 
is actually primitive substitutive (see Corollary~\ref{corollaryMinimalPrimitive}). Thus we can replace condition (i) by
\begin{enumerate}
\item[\rm(i')] \emph{The sequence $x$ is substitutive.}
\end{enumerate}

Observe finally that in the definition of a substitutive sequence $x$,
one may relax all conditions on the pair $(\sigma,\phi)$
of morphisms defining $x$ provided $x$ is a two-sided
infinite sequence (see Exercise~\ref{exerciseCobham}).

We begin by considering the easy case where $x$ is periodic.
\begin{proposition}\label{propositionPeriodicCase}
Every periodic two-sided sequence is substitutive and has a finite
number of derivatives.
\end{proposition}
\begin{proof}
Let $x=v^\infty$ (recall that $u^\infty=\cdots vv\cdot vv\cdots$)
with $v=a_0a_1\cdots a_{n-1}$. Let $B=\{b_0,\ldots,b_{n-1}\}$
and let $\sigma:B^*\to B^*$ be defined by $\sigma(b_i)=b_jb_{j+1}$
with $j=2i\bmod n$. Set $w=b_0b_1\cdots b_{n-1}$.
Then $\sigma(w)=w^2$
and thus $y=w^\infty$ is an admissible fixed point of $\sigma$.
We have $x=\phi(y)$ where $\phi:B^*\to A^*$ is the morphism
defined by $\phi(b_i)=a_i$ for $0\le i<n$. We conclude that
$x$ is substitutive.

Since $X$ has only finitely many different cylinders $[u]$, the number
of its derivatives is also finite.
\end{proof}

\subsection{Sufficiency of the condition}
The implication (iii)$\Rightarrow$ (i) will result of 
Proposition~\ref{proposition(iii)implies(i)}. 
\begin{proposition}\label{proposition(iii)implies(i)}
Let $x$ be an aperiodic uniformly recurrent 
two-sided sequence such that the set
of derivatives $\D_u(x)$, for $u$ prefix of $x^+$, is finite.
Then $x$ is primitive substitutive.
\end{proposition}
\begin{proof}
 There exists a sequence of prefixes $(u_n)_{n\ge 1}$
 of $x^+$ such that $|u_n|<|u_{n+1}|$ and $\D_{u_n} (x) = \D_{u_{n+1}} (x)$,
 for all $n\ge 1$. Clearly this implies, for all $n\ge 0$, 
 that $u_n$ is a prefix of $u_{n+1}$. Take $u = u_1$.

We denote by $\RR'(u)$
the set $\RR'_X(u)$ where $X$ is the shift
generated by $x$.
The sequence $x$ being uniformly recurrent we can choose $N$ so large that
 every factor of length $N$ of $x$ has an occurrence of each 
$ru$ for $r\in\RR'(u)$. Since $x$ is not periodic, the minimal length
of a word in $\RR'(u_n)$ cannot be bounded independently
of $n$. Otherwise
there exists a word $v$ such that $vu_n$ begins with $u_n$
for all $n$ large enough. As a consequence,
 $|v|$ is a period of all $u_n$ and thus of $x^+$.
Since $x$ is uniformly recurrent, this  forces
the period of all words in $\cL(x)$ to be bounded and thus $x$ is
 periodic, a contradiction.
Thus there exists  some $w=u_l$ such that
 $|rw|>N$ for all $r\in \RR'(w)$. %We have $R(u) = R(v)$, we set $R=R(u)$.
 %The set $\RR(v)$ is included in $\RR(u)^+$.

Let $\varphi_u:A_u^*\to A^*$ be a coding morphism for $\RR'(u)$.
Since $\RR'(u)$ is a suffix code,
the map $\varphi_u$ is one to one. Since $\D_u(x)=\D_w(x)$
we have $A_u=A_w$. Set $B=A_u=A_w$.
 Since $\RR'(w)\subset\RR'(u)^*$,
 we can define a morphism $\tau : B^* \rightarrow B^*$
 such that $\varphi_u \circ\tau = \varphi_w$. By the choice of $u$ and $w$,
 for every $i,j$ in $B$ the word $\varphi_u (j)w $ appears
 in $\varphi_w (i)$. By definition of return words,
this implies that $j$ appears in $\tau(i)$.
This means that $\tau$ is a primitive substitution.
 We have
$$
\varphi_u \circ\tau (\D_u (x)) = \varphi_w (\D_w (x)) = x=\varphi_u(\D_u(x)) .
$$
Since $\varphi_u$ is injective
by Proposition~\ref{propositionCircularCodesInjective},
we obtain  $\tau (\D_u (x)) = \D_u (x)$. Hence $\D_u (x)$ is a fixed point of $\tau$. Since $x = \varphi_u (\D_u (x))$,
 Proposition \ref{ch5:prop:rauzy} implies that there is a primitive
substitution $\zeta$ on an alphabet $C$,
an admissible fixed point $z$ of $\zeta$ and a map $\theta:C\to A$
such that $\theta(z)=x$. Thus $x$ is primitive substitutive.
\end{proof}
\bigskip

%Let $\sigma : A^*\rightarrow A^*$ be a primitive substitution with fixed point $y$ and dominant eigenvalue $\alpha$, and $\phi$ be a letter to letter morphism from $A$ to $B$ such that $x = \phi (y)$. We say in this case that $x$ 
%and $y$ are $\alpha$-substitutive.

\subsection{Necessity of the condition}
We now consider a primitive substitutive sequence $x$. By definition,
there is a primitive substitution $\sigma:B^*\to B^*$, a fixed point
$y$ of $\sigma$ and a letter-to-letter morphism $\phi:B^*\to A^*$
such that $x=\phi(y)$.
We first prove the following preliminary result.
\begin{proposition}
\label{derptf}
Let $y$ be a two-sided fixed point of a primitive
substitution $\sigma:B^*\to B^*$.
Let $u$ be a non-empty prefix of $y^+$. 
The derived sequence $\D_u (y)$ is the fixed point of a primitive substitution $\sigma_u : B_u^* \rightarrow B_u^*$ which satisfies 
\begin{equation}
\varphi_u \circ \sigma_u = \sigma \circ \varphi_u .\label{equationsigma_u}
\end{equation}
\end{proposition}

\begin{proof} For every $i\in B_u$,
the word $u$ is a prefix of $\sigma (u)$ and $\varphi_u (i)u$.
Next
\begin{displaymath}
\varphi_u(i)u\in\cL(X)\Rightarrow \sigma(\varphi_u(i))u\in\cL(X).
\end{displaymath}
hence $\sigma ( \varphi_u (i))$  belongs to $\RR'_X(u)^+$.
 We can therefore define a morphism $\sigma_u : B_u \rightarrow B_u^+$ by
$$
\varphi_u \circ \sigma_u = \sigma \circ \varphi_u .
$$ 
For all $n\ge 1$ we have 
$\varphi_u \circ \sigma_u^n = \sigma^n \circ \varphi_u $.
Let $i,j\in B_u$.
 For $n$ large enough, since $y$ is uniformly recurrent,
the word  $\varphi_u(j)u$ appears in $\sigma^n(\varphi_u(i))$.
By definition of return words, this implies that $j$ is a factor
of $\sigma_u^n(i)$. Thus $\sigma_u$ is primitive. Moreover 
\begin{equation}
\varphi_u \circ \sigma_u (\D_u (y)) =  \sigma \circ \varphi_u (\D_u (y)) = \sigma (y) = y=\varphi_u(\D_u(y)).\label{eqphi_usigma_u}
\end{equation}
Since $\varphi_u$ is injective by Proposition~\ref{propositionCircularCodesInjective}, it follows that 
$\sigma_u (\D_u (y)) = \D_u (y)$. Thus $\D_u(y)$
is a fixed point of $\sigma_u$.
%Let $\beta$ be the domininant eigenvalue of $M(\sigma_u)$
%and let $v$ be a corresponding nonnegative eigenvector.
%Set $w=M(\varphi_u)v$. Then, using \eqref{equationsigma_u},
%\begin{displaymath}
%M(\sigma)w=M(\varphi_u)M(\sigma_u)v=\beta M(\varphi_u)v=\beta w
%\end{displaymath}
%which implies $\alpha=\beta$. We conclude that $\D_u (y)$ is $\alpha$-substitu\-tive.
\end{proof}
\medskip

The substitution $\sigma_u : B_u^* \rightarrow B_u^{*}$
\index{symbols}{sigma@$\sigma_u$} will be called 
a {\it return  substitution}\index{subject}{return!substitution}.

\medskip

\begin{example}
  If $\sigma : B^*\rightarrow B^{*}$ is the substitution defined by
  $\sigma (a) = aba$ and $\sigma (b) = aa$. Then $\RR'(a)=\{a,ab\}$.
Set $B_a=\{1,2\}$ with $\varphi_a(1)=a$ and $\varphi_a(2)=ab$.
 Then the return substitution $\sigma_a$ is given
 by $\sigma_a (1) = 21$ and $\sigma_a ( 2 ) = 2111$.
\end{example}
\begin{proposition}
\label{suitedersub}
Let $\phi:B^*\to A^*$ be a letter-to-letter morphism.
Let $y\in B^\Z$ be uniformly recurrent, let $x=\phi(y)$,
and let $u$ be a prefix of $x^+$. Then,  there exists a prefix $v$ of $y^+$
and a morphism $\lambda_u:B_v^*\to A_u^*$   such that  
$\varphi_u \circ \lambda_u = \phi\circ\varphi_v$ 
and $\lambda_u (\D_v (y)) = \D_u (x)$.
\end{proposition}
\begin{proof} Let $v$ be the unique prefix of $y^+$ such that $\phi (v) = u$. 
If $w$ is a return word to $v$ then $\phi (w)$ is a concatenation of 
return words to $u$. The morphism $\varphi_u$ being one to one,
 we can define a morphism $\lambda_u : B_v^* \rightarrow A_u^*$ by
 $\varphi_u \circ\lambda_u = \phi\circ\varphi_v$. We have,
as in Equation~\eqref{eqphi_usigma_u},
\begin{displaymath}
\varphi_u\circ\lambda_u(\D_v(y))=\phi\circ\varphi_v(\D_v(y))
=\phi(y)=x=\varphi_u(\D_u(x)).
\end{displaymath}
This implies that $\lambda_u (\D_v (y) )= \D_u (x)$. 
\end{proof}

Consequently to prove that (i) implies (iii) it suffices to prove that the sets $\{
 \sigma_v \mid v \hbox{ prefix of } y^+ \}$ and $\{ \lambda_u \mid u \hbox{ prefix of } x^+ \}$ are finite. 

\begin{proposition}
\label{morphfini}
Let $y\in B^\Z$ be a fixed point of a primitive
substitution $\sigma$, let $\phi:B^*\to A^*$ be a letter-to-letter
substitution and let $x=\phi(y)$.
The sets $\{  \sigma_v \mid v \hbox{ prefix of } y^+ \}$ and $\{ \lambda_u \mid u \hbox{ prefix of } x^+ \}$ are finite.
\end{proposition}

\begin{proof}
The periodic case is easy to check hence we suppose that $y$ is non-periodic.
By Proposition~\ref{propositionPrimitiveSubstitutiveIsLR},
the sequence $y$ is linearly recurrent, say with constant $K$.
We start by proving that the set
\begin{displaymath}
\{  \sigma_v : B_v^*  \rightarrow B_v^{*}  \mid v \hbox{ prefix of } y^+ \}
\end{displaymath}
 is finite. For this, it suffices to prove that $\Card(B_v)$ and $|\sigma_v (i)|$ are bounded independently of $v$ and $i\in B_v$. 
Set as usual $|\sigma| = \sup \{|\sigma (a)| \mid a\in A  \}$
and $\langle\sigma\rangle = \inf \{|\sigma (a)| \mid a\in A  \}$.

Let $v$ be a non-empty prefix of $y^+$, let $i$ be an element of $B_{v}$ and 
let $w =\varphi_{v}(i)$ be a left return word to $v$.
 Since $y$
is not periodic, we have $| w| \leq K| v|$ and thus
$| \sigma (w) | \leq |\sigma| K | v|$. 
The length of each element of $\RR'(v)$ is larger than $| v|/K$ 
by assertion 4 of Proposition \ref{ch5:proposition:firstpropLR}.
Now, since $wv\in \cL(y)$ and since $wv$ begins with $v$,
we have that the word $\sigma(w)v$ is in $\cL(y)$ and begins
with $v$. Thus we can write $\sigma(w)=x_1x_2\cdots x_k$
with $k\ge 1$ and $x_i\in\RR'(v)$ for $1\le i\le k$. Since
$|x_i|\ge [v|/K$, we have
\begin{displaymath}
[v||\sigma|K\ge|\sigma(w)|\ge k|v|/K
\end{displaymath}
and we conclude that $k\le|\sigma|K^2$.
Therefore $\sigma_v(i)=i_1i_2\cdots i_k$ with $x_j=\varphi_v(i_j)$.
Thus we obtain finally 
that $| \sigma_{v}(i) | \leq |\sigma|K^2$. 
Moreover we know from assertion 5 of Proposition \ref{ch5:proposition:firstpropLR} that $\Card(B_v) \leq K (K+1)^2$. This ends the first part of the proof.

Let now $u$ be a factor of $x$ and $v$ be a factor of $y$ such that $\phi (v) = u$. The length of a return word to $u$ in $x$ is bounded
by the length of a return word to $v$ in $y$ and thus $x$ is linearly
recurrent with constant $K$. Consequently, since $x$ is non periodic,
$\langle\varphi_u\rangle\ge|u|/K$
by Proposition~\ref{ch5:proposition:firstpropLR} again.
We have then for every prefix $u$ of $x^+$,
\begin{eqnarray*}
|\lambda_u (i) | (1/K) |u| 
&\leq&  |\lambda_u (i) |\langle\varphi_u\rangle 
\leq  | \varphi_{u}(\lambda_u (i))|\\
&=&  |\phi(\varphi_v(i))|=| \varphi_{v}( i)| 
\leq |\varphi_v| 
\leq K|v| = K|u|.
\end{eqnarray*}
Hence $|\lambda_u (i) |\leq K^2$. This completes the proof.
\end{proof}
\subsection{End of the proof}
We now conclude the proof of the main result.

\begin{proofof}{of Theorem\ref{theoremCharacterisationSubstitutive}}
We have already proved that (iii) implies (i) (Proposition~\ref{proposition(iii)implies(i)}). We have also proved that (i) implies (iii) (Proposition~\ref{morphfini}).

Since (ii) clearly implies (iii),
 there only remains to prove that (i) implies (ii).
We start with some notation. Let $t$ be a word with prefix  $s$. By $s^{-1}t$ we mean the word $r$ such that $t = sr$. In this way we have $s s^{-1} t = t$.

Since the sequence $x$
is primitive substitutive, there is a primitive
substitution $\sigma$, an admissible fixed point $y$ of $\sigma$
and a letter-to-letter morphism $\phi$ such that $x=\phi(y)$.

Since $y$ is a fixed point of a primitive substitution,
by Proposition \ref{propositionPrimitiveSubstitutiveIsLR}
it is linearly recurrent, say with constant $K$. It
follows that $x$ is linearly recurrent with the same constant $K$.

Let $u$ be a word of $\cL(x)$ and $v$ be
such that $vu$ is a prefix of $x^+$ and $u$ has exactly one occurrence in $vu$. 
Since $v$ is a suffix of a word in $\RR'_X(u)$,
we have $|v|\le K|u|$ and thus $|vu| \leq (K+1) |u|$.

If $w$ is a left return word to $vu$ then $u$ is a prefix of $v^{-1} w vu$ and 
$v^{-1} w vu=(v^{-1}wv)u$ is a word of $\cL(x)$. Hence $v^{-1} w v$  is a concatenation of return words to $u$. Thus, we can define $\phi_{v,u} : B_{vu}^* \rightarrow B_u^*$ by
\begin{equation}
\varphi_u \circ\phi_{v,u} (i) = v^{-1} \varphi_{vu} (i) v,\label{eqphi_v,u}
\end{equation}
for all $i\in B_{vu}$. We have 
\begin{equation}
\phi_{v,u} (\D_{vu} (x)) = \D_u (x).\label{eqphi_v,uD_vu}
\end{equation} 
Indeed, by definition of $v$, the integer $i=|v|$
is the least integer $i\ge 0$ such that $T^ix$ begins with $u$.
Thus, by definition of $\D_u(x)$, we have 
\begin{equation}
T^{-|v|}x=\varphi_u(\D_u(x)).\label{eqv-1x}
\end{equation}
Now, by definition of the morphism $\phi_{u,v}$, if $\varphi_{vu}(y)=x$, we have
using iteratively \eqref{eqphi_v,u} on the prefixes of $y$,
\begin{equation}
\varphi_u\circ\phi_{v,u}(y)=T^{-|v|}\varphi_{vu}(y)=T^{-|v|}x.\label{eqvarphi_uphi_vu}
\end{equation}
We conclude using Equations \eqref{eqv-1x} and \eqref{eqvarphi_uphi_vu}
and the equality $y=\D_{vu}(x)$
that
\begin{displaymath}
\varphi_u(\phi_{v,u}(y))=\varphi_u(\D_u(x))
\end{displaymath}
whence \eqref{eqphi_v,uD_vu} since $\varphi_u$ is injective.

The set $\{ \D_u (x) ; u=x_{[0,n]}, n\geq 0 \}$ being finite, it suffices to prove that the set
$$
H = \{ \phi_{v,u} : B_{vu} \rightarrow B_u ; vu=x_{[0,n]}, |vu|\leq (K+1)|u| , n\geq 0 \} 
$$
 is finite to conclude. By Proposition
\ref{ch5:proposition:firstpropLR} the length of every word in $\RR_X(u)$
is at least $|u|/K$. Thus for every $i\in B_{vu}$, we have
$|\varphi_u(\phi_{v,u}(i))|\ge |\phi_{v,u}(i)||u|/K$.
Therefore, or all $i\in B_{vu}$ we have, using \eqref{eqphi_v,u},
$$
| \phi_{v,u} (i)| \leq \frac{|\varphi_u(\phi_{v,u} (i))|}{|u|/K}
=\frac{|\varphi_u(i)|}{|u|/K} \leq K^2.
$$
Moreover $\Card(B_s)\leq K(K+1)^2$ for all words $s\in \cL(x)$ hence $H$ is finite.\hfill 
\end{proofof}

\begin{example}\label{exampleFibonacciSubstitutive}
Let $\varphi:a\to ab,b\to a$
be the Fibonacci morphism and let $x=\varphi^\omega(a)$
be the Fibonacci word. Let $u_1,u_2,\ldots$ be the palindrome prefixes of $x$.
Since the directive word of $x$
is $(ab)^\omega$, we have by Equation \eqref{equationReturnPal} 
$\varphi_{u_{2n}}=(L_{ab})^n$ and $\varphi_{u_{2n+1}}=(L_{ab})^nL_a$.
Let $\overline{x}$ be the result of exchanging $a,b$ in $x$.
We have $x=L_{a}(\overline{x})$  because $x=\Pal((ab)^\omega)=L_a(\Pal(ba)^\omega))=
L_a(\overline{x})$. Next $x=L_{ab }(x)$ by a similar computation
(or also because $L_{ab}=\varphi^2$). Thus $\D_{u_{2n}}(x)=x$
and $\D_{u_{2n+1}}=\overline{x}$.

If $u$ is any factor of $x$, we have
$\RR'_x(u)=v^{-1}\RR'_X(u_n)v$ where $v$ is the shortest word
such that $vu$ is a prefix of some palindrome prefix $u_n$.
Consequently, $\D_u(x)=\D_{u_n}(x)$.
Thus the derivatives of $x$
with respect to a factor of $x$ are $x$ and  $\bar{x}$.
\end{example}
\subsection{A variant of the main result}
We have considered  in the first part of this section
substitutive sequences. We now consider substitutive shifts.

First of all, we have the following result which shows
that the class of minimal primitive substitutive shifts is closed under
conjugacy.

\begin{theorem}
  A minimal shift is primitive substitutive if and only if it is conjugate
  to a primitive substitution shift.
\end{theorem}
\begin{proof}
  Let first $\gamma:X\to Y$
  be a conjugacy from a primitive substitution shift $X=X(\sigma)$
  where $\sigma$ is primitive to a shift space $Y$. Every conjugacy is a composition
  of a $k$-block code $\gamma_k:X\to X^{(k)}$ and a letter-to-letter morphism
  $\phi:X^{(k)}\to Y$. We have $X^{(k)}=X(\sigma_k)$ where $\sigma_k$ is the $k$-block
  presentation of $\sigma$, which is primitive by Proposition~\ref{propositionphi_kPrimitive}.
  Thus $X^{(k)}$ is a primitive substitution shift and consequently, $Y$ is primitive
  substitutive.

  Conversely, if $Y$ is primitive substitutive, by Proposition \ref{ch5:prop:rauzy},
  it is conjugate to a primitive substitution shift.
\end{proof}
%\marginpar{vrai pour minimal substitutive?oui, par Holton/Zamboni. Autre preuve?}

The following variant of Theorem~\ref{theoremCharacterisationSubstitutive}
characterizes substitutive shifts. In the following statement,
we consider two shifts as different if they cannot be identified
by  renaming  the alphabet.

\begin{theorem}\label{theoremHoltonZamboni}
A minimal shift space $X$ is primitive substitutive if and only if there
is a finite number of different derived systems
on clopen sets $[u]$ for $u\in\cL(X)$.
\end{theorem}
\begin{proof}
Consider $X=X(\sigma,\phi)$ with
$\sigma:B^*\to B^*$ primitive and $\phi:B^*\to A^*$ letter-to-letter.
Let $y$ be an admissible fixed point of $\sigma$ and let $x=\phi(y)$.
Let $u\cL(X)$, let $U=[u]$ and let $(X_U,S_U)$ be the shift  
induced by $X$ on $U$. Since $X$ is minimal, we have $u\in\cL(x)$ and $\D_u(x)\in X_U$.
Since $X_U$ is also minimal, it is generated by $\D_u(x)$. By Theorem
\ref{theoremCharacterisationSubstitutive}, since $x$ is
uniformly recurrent and primitive substitutive, there is a finite
number of derived sequences $\D_u(x)$. Thus
there is only a finite number of derived shifts $(X_U,T_U)$.

Conversely, suppose that there are a finite number of 
systems induced by $X$ on clopen sets $U=[u]$ for $u\in\cL(X)$.
If $X$ is finite, then it is substitutive by Proposition~\ref{propositionPeriodicCase}.
Otherwise, fix some $x\in X$. Since $X$ is minimal
and infinite, $x$ is aperiodic and uniformly recurrent.
We argue as in the proof
of Proposition~\ref{proposition(iii)implies(i)}.
There exists an infinite sequence of prefixes $u_n$ of $x^+$ with $[u_n|<|u_{n+1}|$
such that $X_{[u_n]}=X_{[u_{n+1}]}$. Set $u=u_1$. Since $X$ is minimal,
there is an $N$ such that every word of $\RR'_X(u)u$ appears
in every word of $\cL_N(X)$. Since $x$ is aperiodic, the words
of $\RR'_X(u_n)$ cannot be of bounded length (the argument is the same
as in the proof of Proposition~\ref{proposition(iii)implies(i)}.
Thus there is some $w=u_l$ such that $|rw|>N$ for all $r\in\RR'_X(w)$.
By hypothesis, we can identify $B_u$ and $B_w$ to an alphabet $B$
and set $Y=X_{[u]}=X_{[w]}$. Since $\RR'_X(w)\subset\RR'_X(u)$, we can define
a morphism $\tau:B^*\to B^*$ such that $\varphi_u\circ\tau=\varphi_w$.
By the choice of $u$ and $w$, the morphism $\tau$ is primitive.
Then $\tau(Y)=Y$. Let $y\in B^*$ be an admissible fixed point of $\tau$
and let $z=\varphi_u(y)$.
By Proposition~\ref{ch5:prop:rauzy}, there is a primitive
substitution $\zeta:C^*\to C^*$, an admissible fixed point $t$
of $\tau$ and a letter-to-letter morphism $phi:C^*\to A^*$
such that $z=\phi(t)$.
Thus $X$ is primitive substitutive.
\end{proof}

\begin{example}
Consider again the Fibonacci shift $X=X(\varphi)$ as in Example~\ref{exampleFibonacciSubstitutive}. There are only two different shifts induced by $X$ on clopen sets
$[u]$ for $u\in\cL(X)$, namely $X$ and its image $\overline{X}$ by exchange
of $a,b$.
\end{example}

%%%%%%%%%%%%%%%%%%%%
\section{Exercises}
\exosection{Section~\ref{ch5:subsec:rep-odo}}
\begin{exercise}\label{ch5:ex:odom}
Prove that the set $\Z_{(p_n)}$ of $(p_n)$-adic integers is a compact topological group and actually a topological ring.
\end{exercise}
\begin{exercise}\label{exerciseZdense}
Show that the set of $x\in\Z_{(p_n)}$ such that
$x_n$ is eventually constant is a dense subgroup isomorphic to $\Z$.
\end{exercise}
\begin{exercise}\label{exerciseTopologicalConjugacyOdometers}
Show that two odometers $(X,T)$
and $(X',T')$ are topologically conjugate if and only
if the groups $X$ and $X'$ are isomorphic.
\end{exercise}
\begin{exercise}\label{exerciseExpansionBasispn}
Let $(q_n)_{n\ge 1}$ be an increasing sequence of natural integers
and let $Y=\{(y_n)_{n\ge 0}\mid 0\le y_n<q_{n+1}\}$ be
the group of $(q_n)$-adic expansions.
Set $p_1=q_1$ and $p_{n+1}=p_nq_{n+1}$ for $n\ge 1$.
Show that the map $\varphi:x\mapsto y$ defined by $y_0=x_1$ and
$y_n=(x_{n+1}-x_n)/p_n$ for $n\ge 1$ is homeomorphism from $\Z_{(p_n)}$
onto $Y$. 
\end{exercise}
\begin{exercise}\label{exerciseFactorial}
Show that every integer $n$ has a unique factorial expansion
\begin{displaymath}
x=c_1+c_22!+\ldots
\end{displaymath}
with $0\le c_i\le i$.
\end{exercise}
\begin{exercise}\label{exerciseFactorial2}
Show that $-1=(\cdots321)_!$.
\end{exercise}
\begin{exercise}\label{exerciseInverseLimit}
Let $(G_n)_{n\ge 1}$ be a sequence of groups and let $\varphi_n:G_{n+1}\to G_n$
be a sequence of morphisms. Set $X=\prod_{n\ge 1}G_n$.
The \emph{inverse limit}
\index{subject}{inverse!limit}%
of the sequence $(G_n,\varphi_n)$ is the set 
\begin{displaymath}
X=\{(x_n)_{n\ge 1}\in X\mid \varphi_n(x_{n+1})=x_n\}
\end{displaymath}
Show that the group $\Z_{(p_n)}$ is the inverse limit of the
groups $\Z/p_n\Z$ with $\varphi_n(x_{n+1})=x_{n+1}\bmod p_n$.
\end{exercise}
\begin{exercise}\label{ch5:ex:KR}
Show that the sequence of partitions $\Pg(n)$ defined by Equation~\eqref{equationPatitionsOdometer} is a refining sequence.
\end{exercise}
\begin{exercise}\label{exerciseSuperNatural}
A \emph{supernatural number}
\index{subject}{supernatural number}\index{number!supernatural}%
is a formal product $\prod_{p}p^{n_p}$ where for each prime number
$p$ we have  $n_p\in\N\cup{\infty}$.
If $(p_n)_{n\ge 1}$ is a sequence of integers with $p_n|p_{n+1}$
for all $n\ge 1$ and $(X,T)$ is the associated odometer,
we associate to $X$ the supernatural integer $\sigma(X)=\prod_p p^{n_p}$
where  $p^{n_p}$ is the maximal power of $p$ which divides 
some $p_n$ (and thus all $p_m$ for $m\ge n$).
Show that  two odometers
$(X,T)$ and $(X',T')$ are topologically conjugate 
if and only if $\sigma(X)=\sigma(X')$.
\end{exercise}
\begin{exercise}\label{exerciseProfiniteGroups}
A \emph{profinite group}
\index{subject}{profinite!group}%
\index{subject}{group!profinite}%
is an inverse limit of finite groups. Show that
for every odometer $(X,T)$, the group $X$ is profinite.
\end{exercise}

\begin{exercise}\label{exerciseExpansive}
Let $(X,T)$ be a invertible topological dynamical system which is expansive.
Let $\varepsilon$ be the expansivity constant. Let $\gamma$
be a finite cover of $X$ by sets $C_1,C_2,\ldots,C_r$
of diameter at most $\varepsilon$. 

Show that
 the diameter of the
elements of the cover
$\vee_{j=-n}^nT^{-j}\gamma$ converges to $0$ when $j\to\infty$
(if $\alpha,\beta$ are covers of $X$, then $\alpha\vee\beta$ is formed
of the $A\cap B$ for $A\in\alpha$ and $B\in\beta$ and $T^{-j}\alpha$
by the $T^{-j}A$ for $A\in\alpha$).
\end{exercise}

\begin{exercise}\label{exerciseExpansive2}
Let $(X,T)$ be an invertible topological dynamical system which is expansive
with constant $\varepsilon$. We aim at proving that $(X,T)$
is conjugate to a shift space.

For this, let $\{B_0,B_1,\ldots,B_{k-1}\}$ be a cover by open balls with radius 
$\varepsilon$. Let $C_0=\overline{B_0}$ and for $n\ge 1$, let
$C_n=\overline{B_n}\setminus(B_0\cup\cdots\cup B_{n-1})$. 

1. Show that  $\gamma=\{C_0,\ldots,C_{k-1}\}$ is a closed cover
of $X$ with $\diam(C_i)<\varepsilon$ for each $i$, 
$C_i\cap C_j=\partial C_i\cap \partial C_j$ if $i\ne j$ and
$\partial C_i$ having no interior (we denote
 $\partial C=C\setminus \interior(C)$ where $\interior(C)$ is the interior of $C$).

2. Set $D=\cup_{i=0}^{k-1}\partial C_i$
and $D_\infty=\cup_{n\in\Z}T^n D$. For $x\in X\setminus D_\infty$,
let $y=\psi(x)$ be the sequence in $A^\Z$ with $A=\{0,,\ldots,k-1\}$
defined by $y_n=i$ if $T^n(x)\in C_i$. Show that $\psi$
extends to a conjugacy from $X$ to a subshift of $A^\Z$.
\end{exercise}

\exosection{Section~\ref{sectionSubstitutionsBV}}
\begin{exercise}\label{exerciseSubstitutiveShift}
Suppose  that $X$ is a minimal substitution shift on the alphabet
$A$.
Let $x\in X$, $\varphi : A^* \to B^*$ be a non erasing morphism and $y=\varphi (x)$.
Consider $(Y, S)$ the subshift generated by $y$.
Prove that  it is isomorphic to a primitive substitution shift.
\end{exercise}

\begin{exercise}\label{exerciseChaconTernaryBinary}
The \emph{Chacon ternary substitution}
\index{subject}{Chacon!ternary!substitution}%
\index{subject}{substitution!ternary Chacon}%
is the primitive substitution $\tau:0\to 0012,1\to 12,2\to 012$
(see also Exercise~\ref{exerciseChacon}). Show that 
$w_n=\tau^n(0)$ satisfies the recurrence relation $w_{n+1}=w_nw_n1w'_n$
where $w'_n$ is obtained from $w_n$ by changing the initial letter $0$ into a $2$.
Deduce from this that the $1$-block map  $\theta:0,2\to 0,1\to 1$ defines a conjugacy
from the substitution shift $X(\tau)$ defined by $\tau$ to the substitution shift
$X(\sigma)$ defined by the Chacon binary substitution $\sigma:0\to 0010,1\to 1$.
\end{exercise}
\begin{exercise}\label{exerciseComplexityChacon}
Show that the factor complexity of the Chacon ternary shift is $p_n(X)=2n+1$
(hint: show that the bispecial words are $0$ and the words $\alpha^n(012),\alpha^n(120)$
for $n\ge 0$ where $\alpha(w)=012\tau(w)$ and where $\tau$ is the ternary Chacon substitution).
\end{exercise}
\begin{exercise}\label{exerciseForrestRauzy}
Use Proposition~\ref{propositionLemma15Forrest}
to prove Proposition~\ref{ch5:prop:rauzy}. Develop the construction
on the example of $\tau:0\to 01,1\to 10$ being the Thue-Morse morphism
and $\phi:0\to ab,1\to a$.
\end{exercise}
\begin{exercise}\label{exerciseThueMorseDim3}
Let $X$ be the Thue-Morse shift generated by $\varphi:a\to ab,b\to ba$.
\index{subject}{Thue-Morse!shift space!BV-representation}%
\index{subject}{BV-representation!of Thue-Morse shift}%
Let $f:\{aa,ab,ba,bb\}\to A_2=\{x,y,z,t\}$
 and let $\varphi_2:x\to yz,y\to yt,z\to zx,t\to zy$
be the $2$-block presentation of $\varphi$.
Let $B=\{u,v,w\}$
and let $\phi:B^*\to A_2^*$ be a coding morphism 
for $f(a\RR_X(a))=f(\{aa,aba,abba\})=\{x,yz,ytz\}$.
Show that there is a morphism $\tau$ such that 
$\varphi_2\circ\phi=\phi\circ\tau$ and derive a BV-representation
of $X$ with $3$ vertices at each level $n\ge 1$.
\end{exercise}
\exosection{Section~\ref{ch5:subsec:LR}}
\begin{exercise}\label{exerciseLRClosed}
Show that if $X$ is LR with constant $K$, then its $k$-th block presentation
$(X^{(k)},S)$ is LR with $K\le K(k-1)+1$. Conclude that the class
of LR shifts is closed under conjugacy.
\end{exercise}

\begin{exercise}\label{exerciseExampleComplexityn^2}
Consider the morphism $\sigma:a\to abd,b\to bb,c\to c, d\to dc$.
Show that the factor complexity $p_n(X)$ of the shift $X=X(\sigma)$ is not linear.
\end{exercise}

\begin{exercise}\label{exerciseAltDefLR}
  Show that a Cantor system $(X,T)$ is linearly recurrent if and only if there
  is a refining sequence $(\Pg(n))$ of partitions in towers
  with positive matrices $M(n)$ and 
  with heights $(h_i(n))_{1\le i\le t(n)}$
  satisfying \eqref{equationLR}
  for some $L\ge 1$.
  \end{exercise}
\begin{exercise}\label{exerciseDelecroixBoshernitzan}
  For a probability measure $\mu$ on a shift space $X$,
  denote
  \begin{displaymath}
    \varepsilon_n(\mu)=\min\{\mu(u)\mid u\in\cL_n(X)\}.
    \end{displaymath}
  Show that a shift space $X$ is linearly recurrent if and only if
  there is an invariant probability measure $\mu$ on $X$ such
  that $\inf n\varepsilon_n(\mu)>0$. Derive from this an alternative
  proof that linearly recurrent shifts are uniquely ergodic.
  \end{exercise}
\exosection{Section~\ref{sectionSadicShifts}}
\begin{exercise}\label{exerciseLimitPoint}
 Show that
the following conditions are equivalent for a one-sided
sequence $x$.
\begin{itemize}
\item[(i)] $x$ is a limit point of a primitive $\Sa$-adic system
$\tau$.
\item[(ii)] There is a primitive $\Sa$-adic
system $\tau$ and a sequence $(a_n)$ of letters $a_n\in A_n$
such that $x=\lim\tau_{[0,n)}(a_n^\omega)$.
\item[(iii)] There is a primitive $\Sa$-adic system
$\tau$ and a sequence $(a_n)$ of letters $a_n\in A_n$
such that $\{x\}=\cap_n[\tau_{[0,n)}(a_n)]$ where $[w]$ denotes the
cylinder $\{y\in A_0^\N\mid y_{[0,|w|)}=w\}$.
\end{itemize}
\end{exercise}
\begin{exercise}\label{ExerciseSadicCassaigne}
Show that for every sequence $x\in A^\N$, there are 
morphisms $(\sigma_a)_{a\in A}:(A\cup\#)^*\to (A\cup\#)^*$ 
and $\phi:(A\cup \#)^*\to A^*$ where $\#$ is a letter not in $A$
such that $x=\lim \phi\circ\sigma_0\circ\sigma_1\cdots\sigma_n(\#)$
and thus that $(\varphi,\sigma_0,\sigma_1,\ldots)$
is an $\Sa$-adic representation of $x$.
\end{exercise}

\begin{exercise}\label{exercise(aab,ba)}
Let $X$ be the shift generated by the substitution $\tau:a\to aab,b\to ba$.
Show that the dimension group of $X$ is
isomorphic to 
\begin{displaymath}
\left( {\mathbb Z}^3 , \left\{ {\mathbf x} \in {\mathbb Z}^3  :    \langle {\mathbf x},  {\mathbf v}  \rangle  > 0 \right\}\cup\{0\} , (1,2,1) \right)
\end{displaymath}
where ${\mathbf v}=  [ 2\lambda,\lambda+1,1]$ with
$\lambda=(1+\sqrt{5})/2$.
\end{exercise}
\begin{exercise}\label{exerciseM(u)}
Let $\varphi:A^*\to B^*$ be a nonerasing morphism. Set $U=\varphi(A)$.
For $u\in B^*$, set
\begin{displaymath}
M_\varphi(u)=\{v\in U^*\mid uv\in U^*u\}
\end{displaymath}
The set $M_\varphi(u)$ is for every $u\in B^*$ a submonoid of $U^*$.

Prove that if
 a morphism $\varphi:A^*\to B^*$ is recognizable on a subshift $X$ of $A^\Z$
for aperiodic points, then for every $u\in \cL(X)$
which is not in $U^*$,
the period of words in $M_\varphi(u)\cap\cL(X)$ 
is bounded. 

Prove that the converse is true if $\varphi(A)$ is a prefix code
and $\varphi$ is injective on $A$.
\end{exercise}
\begin{exercise}\label{exerciseLeftPermutative}
A morphism $\varphi:A^*\to B^*$ is \emph{left permutative}\index{subject}{left!permutative morphism}\index{subject}{morphism!left permutative} if
every word $\varphi(a)$ for $a\in A$ begins with a distinct letter.
In particular, $\varphi$ is injective on $A$ and $\varphi(A)$
is a prefix code.

Prove that
if $\varphi:A^*\to B^*$ is left permutative, then it is recognizable at
aperiodic points. Hint: use Exercise~\ref{exerciseM(u)}
\end{exercise}
\exosection{Section~\ref{sectionDerivatives}}
\begin{exercise}\label{exerciseErasing}
Show that the sequence $001^\omega$ is substitutive
but not purely substitutive.
\end{exercise}
\begin{exercise}\label{exerciseShiftDerivative}
Show that every shift of a substitutive sequence is substitutive.
Hint: consider a $k$-th higher block presentation of the
shift generated by the sequence.
\end{exercise}
\begin{exercise}\label{exerciseCobham}
Let $\tau:B^*\to B^*$ be a morphism prolongable on $a\in B$
and let $y=\tau^\omega(a)$.
Let $\phi:B^*\to A^*$ be a morphism such that $x=\phi(y)$ is an infinite
word. 
Show that
the sequence $x$ is substitutive. Hint: adapt the proof
of Proposition~\ref{ch5:prop:rauzy} to the case where $\tau$
is not primitive and $\phi$ is possibly erasing.
\end{exercise}
\begin{exercise}\label{exerciseNonerasing}
Let $\sigma:\{0,1,2\}^*\to\{0,1,2\}^*$ be the substitution

\begin{displaymath}
  0\to 01222,1\to 10222,2\to\varepsilon
\end{displaymath}
and let
$x=\sigma^\omega(0)$. Show that $x$ is not the fixed point
of a nonerasing substitution. Hint: show that erasing the
letter $2$ in $x$ gives the Thue-Morse sequence.
\end{exercise}
%%%%%%%%%%%%%%%%%%%%%
\section{Solutions}
\exosection{Section~\ref{ch5:subsec:rep-odo}}
\begin{solution}{\protect{\ref{ch5:ex:odom}}}
Each $\Z/p_n\Z$ is a compact topological group
for the discrete topology.
The direct product of compact topological groups is
a compact topological group for the product topology by Tychonov theorem.
\index{subject}{Tychonov Theorem}\index{subject}{Theorem!Tychonov}%
\index{names}{Tychonov, Andrei N.}%
Since $\Z_{(p_n)}$ is a closed subgroup of the direct product
of the $\Z/p_n\Z$, the result follows.
Each $\Z/p_n\Z$ is also a finite ring and $\Z_{(p_n)}$ is a subring
of their direct product which is also a topological ring.
\end{solution}
\begin{solution}{\ref{exerciseZdense}}
Let $G$ be the set of $x=(x_n)\in\Z_{(p_n)}$ such that $x_n$
is eventually constant. The map $\varphi:G\to \Z$ such that $\varphi(x)$
is the value of all  $x_n$ for $n$ large enough is an injective morphism from
$G$ into $\Z$. Since
the sequence $(p_n)$ is strictly increasing, it is onto. It is
clear that $G$ is dense in $\Z_{(p_n)}$.
\end{solution}
\begin{solution}{\ref{exerciseTopologicalConjugacyOdometers}}
Assume first that $\varphi:X\to X'$ is a topological conjugacy from
$(X,T)$ onto $(X',T')$. Set $\alpha(x)=\varphi(x)-\varphi(0)$. Then
$\alpha$ is another conjugacy from $(X,T)$ onto $(X',T')$. It satisfies
$\alpha(0)=0$ and $\alpha(1)=\alpha(T0)=T'\alpha(0)=T'0=1$. 
 Since $\Z$ is dense in $X$ and $X'$,
 it is an isomorphism. Conversely, assume that $\varphi:X\to X'$
is a group isomorphism. Since $\varphi(1)$ generates $\Z$,
we have  $\varphi(1)=1$ or $-1$. In the first case,
$\varphi$ is a conjugacy. In the second case, the map $\psi(x)=-\varphi(x)$
is a conjugacy.
\end{solution}
\begin{solution}{\protect{\ref{exerciseExpansionBasispn}}}
The inverse of $\varphi$ is the map from $X$ to $\Z_{(p_n)}$ defined by
\begin{displaymath}
x_n=y_0+y_1p_1+\ldots+y_{n-1}p_{n-1}.
\end{displaymath}

\end{solution}
\begin{solution}{\protect{\ref{exerciseFactorial}}}
Set $x_n=x\bmod(n+1)!$ with $x_0=0$. Then
\begin{displaymath}
c_n=(x_{n}-x_{n-1})/n!
\end{displaymath}
with $n\ge 1$ gives the unique factorial expansion of $x$.
\end{solution}
\begin{solution}{\ref{exerciseFactorial2}}
This holds because $1+2.2!+\ldots+n.n!=(n+1)!-1$, as one may verify by induction on $n$.
\end{solution}
\begin{solution}{\protect{\ref{exerciseInverseLimit}}}
The verification is easy.
\end{solution}
\begin{solution}{\protect{\ref{ch5:ex:KR}}}
Condition (KR1) is satisfied since $\cap_{n\ge 1}B(n)=\{0\}$. Condition (KR2)
is also satisfied because $B(n+1)\subset B(n)$ and
for every $j\le p_{n+1}$,
\begin{displaymath}
B^j(n+1)=\cup_{k\equiv j\bmod p_n}B^kB(n).
\end{displaymath}
Finally, since
\begin{displaymath}
\cap_{n\ge 1}T^{x_n}B(n)=\{x\}
\end{displaymath}
condition (KR3) is also satisfied.
\end{solution}
\begin{solution}{\ref{exerciseSuperNatural}}
If $(X,T)$ and $(X',T')$ are conjugate, then
$\sigma(X)=\sigma(X')$. Indeed, for every prime $p$, $p^n$ divides some $p_n$
if and only if the group $X$ has elements of order $p^n$.

Conversely, let $\sigma=\prod_p p^{n_p}$ be a supernatural number.
Consider the odometer $(Y,S)$
corresponding to the sequence $(q_1,q_2,\ldots)$ defined by
\begin{displaymath}
q_n=\prod_{n_p=\infty,p\le n}p^{n}\times\prod_{n_p<\infty,p\le n}p^{n_p}.
\end{displaymath}
By construction, $\sigma(Y)=\sigma$.

Let now $(X,T)$ be an arbitrary odometer corresponding to $(p_1,p_2,\ldots)$
such that $\sigma(X)=\sigma$. Note that for each $x\in X$
and each integer  $q\ge 1$ the sequence $(x_n\bmod q)_{n\ge 1}$ is nondecreasing
and eventually equal to an integer denoted $\max (x\bmod q)$.

Consider the continuous map $y\mapsto x$ from $Y$ to $X$ defined by 
\begin{displaymath}
x_n=\max(y\bmod p_n).
\end{displaymath}
 Its inverse is the map $x\mapsto y$ from $X$ to $Y$ where $y_n$
is the unique integer $<q_n$ such that
\begin{displaymath}
y_n\equiv \max x\bmod p^m
\end{displaymath}
for all $p<n$ with $m=n$ if $n_p=\infty$ and $m=n_p$ otherwise.
 The existence and uniqueness of $y_n$ are guaranteed
by the Chinese Remainder Theorem. This shows that 
 $(X,T)$ and $(Y,S)$ are topologically conjugate. Thus
all  odometers such that $\sigma(X)=\sigma(Y)$ are conjugate.
\end{solution}
\begin{solution}{\ref{exerciseProfiniteGroups}}
We have seen in Exercise~\ref{exerciseInverseLimit} that
the odometer $(X,T)$ associated with the sequence $(p_1,p_2,\ldots)$
the inverse limit of the finite groups $\Z/p_n\Z$.
\end{solution}

\begin{solution}{\ref{exerciseExpansive}}
Arguing by contradiction, assume that there is an
$\varepsilon_0>0$, a stricly increasing
 sequence $n_i$  and a sequence $x_i,y_i$ of points
such that $d(x_i,y_i)\ge\varepsilon_0$ and
$x_i,y_i\in \vee_{j=-n_i}^{j=n}T^{-j}C_{i,j}$ where $C_{i,j}\in\gamma$.
Choose a subsequence $i_k$ such that $x_{i_k}\to x$ and
$y_{i_k}\to y$. Then $d(x,y)\ge\varepsilon_0$. For every
$j\in\Z$, there is an $\ell_j$ such that $x_{i_k},y_{i_k}\in T^{-j}C_{\ell_j}$
for an infinity of $k$. Thus $x,y\in T^{-j}\overline{C_{\ell_j}}$
and therefore $d(x,y)\le\varepsilon$. This implies
$x=y$, a contradiction.
\end{solution}
\begin{solution}{\ref{exerciseExpansive2}}
1. It is clear that $\gamma$ is a closed cover of $X$
by sets with radius at most $\varepsilon$. Next, for $i<j$,
we have $C_i\cap C_j=C_i\cap \partial C_j$ because
 $\interior(C_j)=B_j\setminus\overline{(B_0\cup\ldots\cup B_{j-1})}$
and thus $C_i\cap C_j=C_i\cap \partial C_j=\partial C_i\cap\partial C_j$ 
because $\partial C_i\cap\interior(C_j)\subset B_i\setminus(B_0\cup\ldots\cup B_{j-1})=\emptyset$.

By Baire Category Theorem, the set $X\setminus D_\infty$ is dense.
The map $\psi:X\setminus D_\infty$ is injective. Indeed,
suppose that $\psi(x)= \psi(y)$. 
For every $\alpha>0$, by Exercise~\ref{exerciseExpansive},
we can choose $N$ such that $\diam(\vee_{|n|\le N}T^n\gamma)<\alpha$.
Then $\psi(x)_n=\psi(y)_n$ implies that $x,y$ are in the same 
element of $\vee_{|n|\le N}T^n\gamma$ and thus $d(x,y)<\alpha$.
This shows that $x=y$. The extension of $\psi$
to $X$ is then a conjugacy from $X$ to the closure
of $\psi(X)$ which is a subshift of $A^\Z$.
\end{solution}

\exosection{Section~\ref{sectionSubstitutionsBV}}
\begin{solution}{\ref{exerciseSubstitutiveShift}}
Let $\sigma$ be a substitution such that $X=X(\sigma)$.
Replacing if necessary $\sigma$ by one of its powers,
let $x$ be an admissible fixed point of $\sigma$. Set $r=x_{-1}$
and $\ell=x_0$. Let $\psi:B\to \RR_X(r\cdot\ell)$ be a coding
morphism and let $\tau:B^*\to B^*$ be the substitution
such that $\psi\circ\tau=\sigma\circ\psi$. By
Proposition~\ref{propositiontauProperAperiodic}, the morphism
$\tau$ is primitive. Thus by Proposition~\ref{ch5:prop:rauzy},
applied to the substitution $\tau$ and the morphism $\phi=\varphi\circ\psi$,
the shift space $Y$ is conjugate to a primitive substitution shift
$X(\zeta)$
(see Figure~\ref{figureConjugacySubstitutive}).
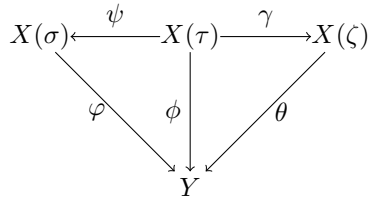
\begin{figure}[hbt]
\centering
\tikzset{node/.style={circle,draw,minimum size=0.1cm,inner sep=0pt}}
\tikzset{title/.style={minimum size=0.4cm,inner sep=0pt}}
\begin{tikzpicture}
\node[title](X)at(0,2){$X(\sigma)$};
\node[title](Z)at(2,2){$X(\tau)$};
\node[title](Xzeta)at(4,2){$X(\zeta)$};
\node[title](Y)at(2,0){$Y$};
\draw[->,above](Z)edge node{$\psi$}(X);
\draw[->,above](Z)edge node{$\gamma$}(Xzeta);
\draw[->,left](X)edge node{$\varphi$}(Y);
\draw[->,left](Z)edge node{$\phi$}(Y);
\draw[->,right](Xzeta)edge node{$\theta$}(Y);
\end{tikzpicture}
\caption{The conjugacy from $X(\zeta)$ onto $Y$.}\label{figureConjugacySubstitutive}
\end{figure}
\end{solution}

\begin{solution}{\ref{exerciseChaconTernaryBinary}}
Set $w_n=0t_n$ and thus $w'_n=2t_n$ for $n\ge 0$. We have then
\begin{displaymath}
w_{n+1}=\tau(w_n)=0012\tau(t_n)=0\tau(2t_n)=0\tau(w'_n)
\end{displaymath}
showing that $t_{n+1}=\tau(w'_n)$ for $n\ge 0$. Thus
\begin{eqnarray*}
w_{n+1}&=&\tau(w_{n-1}w_{n-1}1w'_{n-1})=w_{n}w_{n}12\tau(w'_{n-1})\\
&=&w_{n}w_{n}12t_{n}=w_{n}w_{n}1w'_{n}.
\end{eqnarray*}
The map $\theta$ sends the infinite word $\tau^\omega(0)$ to $\sigma^\omega(0)$
and thus maps $(X(\tau),S)$ to $(X(\sigma),S)$. Its inverse is the map which replaces
$0$ by $2$ when there is a $1$ just before. Thus $\theta$ is a conjugacy.
\end{solution}
\begin{solution}{\ref{exerciseComplexityChacon}}
The bispecial words of length at most $3$
are $\varepsilon$, $0$, $012$ and $120$.
Their extension graphs are shown in Figure~\ref{figureChacon}.

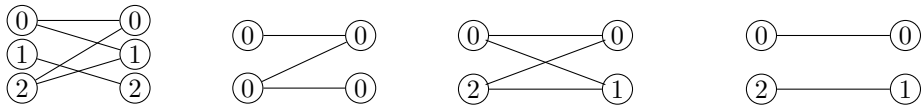
\begin{figure}[hbt]
\centering
%%% IMAGE EN TIKZ
	\tikzset{node/.style={circle,draw,minimum size=.4cm,inner sep=0.1pt}}
	\begin{tikzpicture}
% epsilon
          \node[node](0l)at(-2,2.2){$0$};\node[node](0r)at(-.5,2.2){$0$};
          \node[node](1l)at(-2,1.75){$1$};\node[node](1r)at(-.5,1.75){$1$};
          \node[node](2l)at(-2,1.3){$2$};\node[node](2r)at(-.5,1.3){$2$};
          \draw(0l)edge node{}(0r);\draw(0l)edge node{}(1r);
          \draw(1l)edge node{}(2r);
          \draw(2l)edge node{}(0r);\draw(2l)edge node{}(1r);
% 0    
      \node[node](0l)at (1,2){$0$};\node[node](1r)at (2.5,2){$0$};
          \node[node](2l)at (1,1.3){$0$};\node[node](0r)at (2.5,1.3){$0$};
          \draw(0l)edge node{}(1r);
          \draw(2l)edge node{}(1r);
          \draw(2l)edge node{}(0r);
		\node[node](abcal) at (4,2){$0$};
		\node[node](abccl) [below= 0.3cm of abcal] {$2$};
		\node[node](abcar) [right= 1.5cm of abcal] {$0$};
		\node[node](abcbr) [right= 1.5cm of abccl] {$1$};
		\path[draw]
			(abcal) edge node {} (abcar)
			(abccl) edge node {} (abcbr);
		\path[draw, shorten <=0 -1pt, shorten >=-1pt]
			(abcal) edge node {} (abcbr)
			(abccl) edge node {} (abcar);
		\node[node](bcaal) [right= 1.5 of abcar] {$0$};
		\node[node](bcacl) [below= 0.3cm of bcaal] {$2$};
		\node[node](bcaar) [right= 1.5cm of bcaal] {$0$};
		\node[node](bcabr) [right= 1.5cm of bcacl] {$1$};
		\path[draw]
			(bcaal) edge node {} (bcaar)
			(bcacl) edge node {} (bcabr);
	\end{tikzpicture}

 \caption{The extension graphs of $\varepsilon$, $0$, $012$ and $120$.}
 \label{figureChacon}
\end{figure}

Let now $\alpha$ be the map  defined by $\alpha(x) = 012 \tau(x)$.
A bispecial word $y$ of length at least $4$ begins and ends with $012$ (see the trees
of left and right special words in Figure~\ref{figureChacon2}).

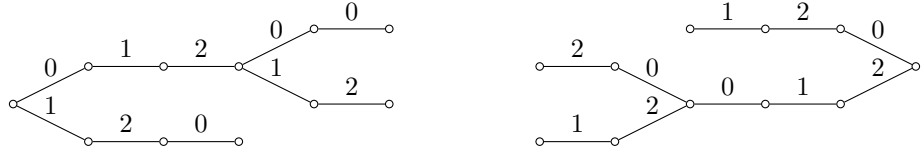
\begin{figure}[hbt]
\centering
\tikzset{node/.style={circle,draw,minimum size=0.1cm,inner sep=0pt}}
\begin{tikzpicture}
%left special
\node[node](e)at(0,0.5){};
\node[node](0)at(1,1){};\node[node](1)at(1,0){};
\node[node](01)at(2,1){};\node[node](12)at(2,0){};
\node[node](012)at(3,1){};\node[node](120)at(3,0){};
\node[node](0120)at(4,1.5){};\node[node](0121)at(4,0.5){};
\node[node](01200)at(5,1.5){};\node[node](01212)at(5,0.5){};

\draw[above](e)edge node{$0$}(0);\draw[above](e)edge node{$1$}(1);
\draw[above](0)edge node{$1$}(01);\draw[above](1)edge node{$2$}(12);
\draw[above](01)edge node{$2$}(012);\draw[above](12)edge node{$0$}(120);
\draw[above](012)edge node{$0$}(0120);\draw[above](012)edge node{$1$}(0121);
\draw[above](0120)edge node{$0$}(01200);\draw[above](0121)edge node{$2$}(01212);
%right special
\node[node](e)at(12,1){};
\node[node](0)at(11,1.5){};\node[node](2)at(11,0.5){};
\node[node](20)at(10,1.5){};\node[node](12)at(10,0.5){};
\node[node](120)at(9,1.5){};\node[node](012)at(9,0.5){};
\node[node](0012)at(8,1){};\node[node](2012)at(8,0){};
\node[node](20012)at(7,1){};\node[node](12012)at(7,0){};

\draw[above](e)edge node{$0$}(0);\draw[above](e)edge node{$2$}(2);
\draw[above](0)edge node{$2$}(20);\draw[above](2)edge node{$1$}(12);
\draw[above](20)edge node{$1$}(120);\draw[above](12)edge node{$0$}(012);
\draw[above](012)edge node{$0$}(0012);\draw[above](012)edge node{$2$}(2012);
\draw[above](0012)edge node{$2$}(20012);\draw[above](2012)edge node{$1$}(12012);

\end{tikzpicture}
\caption{The trees of left and right special words.}\label{figureChacon2}
\end{figure}
Thus $y=012\tau(x)=\alpha(x)$ where $x$ is a bispecial word.

Let us verify that if the extension graph of $x$ is the graph $\E(012)$ (see Figure~\ref{figureChacon}),
 the same holds for the extension graph of $y = \alpha(x)$.
Indeed, since $0x0 \in \cL(X)$, the word $\tau(0x0) = 0012 \tau(x) 0012 = 0y0012$ 
is also in $\cL(X)$ and thus $(0,0) \in \E(y)$.
Since $2x0\in \cL(X)$ and since a letter $2$ is always preceded by a letter $1$, we have $bcxa \in \cL(X)$.
Thus $\tau(12x0) = 12y0012 \in \cL(X)$ and thus $(2,0) \in \E(y)$.
The proof of the other cases is similar.
The same property holds for a word $x$ with the extension graph on the right of Figure~\ref{figureChacon}.

We conclude that $b_n(X)=0$ for every $n\ge 0$. Indeed, this is true for $n=0,1$ since
$m(\varepsilon)=m(0)=0$ and for $n=3$ since $m(012)=1$, $m(120)=-1$. Let $n\ge 4$
be such that there are bispecial words of length $n$. These
bispecial words are of the form $\alpha^k(012),\alpha^k(120)$ 
(note that these two words have the same length) and thus 
\begin{displaymath}
b_n=m(\alpha^k(012))+m(\alpha^k(120))=m(012)+m(120)=1-1=0.
\end{displaymath}
\end{solution}
\begin{solution}{\ref{exerciseForrestRauzy}}
Let $M$ be the composition matrix of $\tau$ and let $v$
be the vector $(|\phi(b)|)_{b\in B}$. Let $B=(V,E,\le)$
be the stationary diagram with matrices $(v,M,M,\ldots)$.
Let $P,Q$ be the nonnegative matrices, with $P$ a partition matrix,
such that $M=PQ$ defined by Proposition~\ref{propositionLemma15Forrest}.
Set $M'=QP$ and $w=\begin{bmatrix}1&1&\ldots&1\end{bmatrix}^t$.
Let $B'=(V',E',\le')$ be the ordered Bratteli diagram with
matrices $(w,M',M',\ldots)$ with the order induced from $(V,E,\le)$.

Let $C$ be the index of the columns of $P$ identified
with $\{b_p\mid 1\le p\le|\phi(b)|\}$. Let $\zeta$
be the morphism read on $B'$. Let $\gamma:B\to C^*$ be the morphism
defined by $\gamma(b)=b_1\cdots b_{|\phi(b)|}$. Let
finally $\theta:C\to A$ be defined by $\theta(b_p)=\phi(b)_p$.

It is clear that $\zeta\circ\gamma=\gamma\circ\tau$ since
both are equal to the morphism read on $PQP$ (with the order
induced by $B$). The equality $\phi=\theta\circ\gamma$
holds by definition of $\gamma$ and $\theta$.
This implies assertion 1.

Assertion 2 is clear since $\phi$ is injective from $X(\tau)$ to
$X$.

Assume finally that $\tau$ is eventually proper. We can also assume
that $X(\tau)$ is aperiodic since otherwise the result is trivial.
By Proposition~\ref{propositionBratteliPrimitive} the
system $(X_E,T_E)$ is isomorphic to $(X(\tau),S)$. Since $B$
is properly ordered, $B'$ is also properly ordered and thus
$\zeta$ is eventually proper by Proposition~\ref{propositionPrimitiveProper}.

Consider $\tau:0\to 01,1\to 10$ and $\phi:0\to 01,1\to 0$.
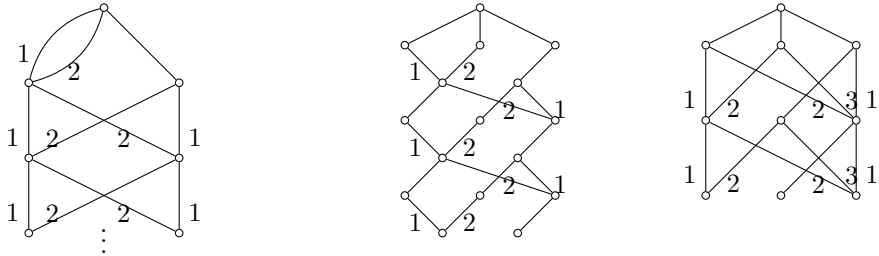
\begin{figure}[hbt]
\centering
\begin{tikzpicture}
\tikzset{node/.style={circle,draw,minimum size=0.1cm,inner sep=0pt}}
\tikzset{title/.style={minimum size=0.4cm,inner sep=0pt}}
%B
\node[node](0)at(1,3){};
\node[node](11)at(0,2){};\node[node](12)at(2,2){};
\node[node](21)at(0,1){};\node[node](22)at(2,1){};
\node[node](31)at(0,0){};\node[node](32)at(2,0){};
\node[title]at(1,0){$\vdots$};

\draw[bend left,right,near end](0)edge node{$2$}(11);
\draw[bend right,left,near end](0)edge node{$1$}(11);
\draw(0)edge node{}(12);
\draw[left,near end](11)edge node{$1$}(21);
\draw[left,near end](11)edge node{$2$}(22);
\draw[left,near end](12)edge node{$2$}(21);
\draw[right,near end](12)edge node{$1$}(22);
\draw[left,near end](21)edge node{$1$}(31);
\draw[left,near end](21)edge node{$2$}(32);
\draw[left,near end](22)edge node{$2$}(31);
\draw[right,near end](22)edge node{$1$}(32);
%%%%%%%%%%%%%%%%%%%%
\node[node](0)at(6,3){};
\node[node](11)at(5,2.5){};\node[node](12)at(6,2.5){};\node[node](13)at(7,2.5){};
\node[node](21)at(5.5,2){};\node[node](22)at(6.5,2){};
\node[node](31)at(5,1.5){};\node[node](32)at(6,1.5){};\node[node](33)at(7,1.5){};
\node[node](41)at(5.5,1){};\node[node](42)at(6.5,1){};
\node[node](51)at(5,.5){};\node[node](52)at(6,.5){};\node[node](53)at(7,.5){};
\node[node](61)at(5.5,0){};\node[node](62)at(6.5,0){};

\draw(0)edge node{}(11);\draw(0)edge node{}(12);\draw(0)edge node{}(13);
\draw[left,near end](11)edge node{$1$}(21);
\draw[right,near end](12)edge node{$2$}(21);\draw(13)edge node{}(22);
\draw(21)edge node{}(31);
\draw[left,near end](21)edge node{$2$}(33);
\draw(22)edge node{}(32);
\draw[right,near end](22)edge node{$1$}(33);
\draw[left,near end](31)edge node{$1$}(41);
\draw[right,near end](32)edge node{$2$}(41);
\draw(33)edge node{}(42);
\draw(41)edge node{}(51);
\draw[left,near end](41)edge node{$2$}(53);
\draw(42)edge node{}(52);
\draw[right,near end](42)edge node{$1$}(53);
\draw[left,near end](51)edge node{$1$}(61);
\draw[right,near end](52)edge node{$2$}(61);
\draw(53)edge node{}(62);
%%%%%%%%%%
\node[node](0)at(10,3){};
\node[node](11)at(9,2.5){};\node[node](12)at(10,2.5){};\node[node](13)at(11,2.5){};
\node[node](21)at(9,1.5){};\node[node](22)at(10,1.5){};\node[node](23)at(11,1.5){};
\node[node](31)at(9,0.5){};\node[node](32)at(10,0.5){};\node[node](33)at(11,0.5){};

\draw(0)edge node{}(11);\draw(0)edge node{}(12);\draw(0)edge node{}(13);
\draw[left,near end](11)edge node{$1$}(21);
\draw[left,very near end](11)edge node{$2$}(23);
\draw[right,very near end](12)edge node{$2$}(21);
\draw[right,near end](12)edge node{$3$}(23);
\draw(13)edge node{}(22);
\draw[right,near end](13)edge node{$1$}(23);
\draw[left,near end](21)edge node{$1$}(31);
\draw[left,very near end](21)edge node{$2$}(33);
\draw[right,very near end](22)edge node{$2$}(31);
\draw[right,near end](22)edge node{$3$}(33);
\draw(23)edge node{}(32);
\draw[right,near end](23)edge node{$1$}(33);
\end{tikzpicture}
\caption{The diagrams $B$ and $B'$.}
\label{figureSplit2}
\end{figure}
The matrices $M,P,Q,M'$ are
\begin{displaymath}
M=\begin{bmatrix}1&1\\1&1\end{bmatrix},\quad
P=\begin{bmatrix}1&1&0\\0&0&1\end{bmatrix},\quad
Q=\begin{bmatrix}1&0\\0&1\\1&1\end{bmatrix},\quad
M'=\begin{bmatrix}1&1&0\\0&0&1\\1&1&1\end{bmatrix}.
\end{displaymath}
The morphisms $\gamma,\zeta$ and $\theta$ are
\begin{eqnarray*}
\gamma&:&0\to ab,1\to c\\
\zeta&:&a\to ab,b\to c,c\to cab\\
\theta&:&a\to a,b\to b,c\to a
\end{eqnarray*}
\end{solution}
\begin{solution}{\ref{exerciseThueMorseDim3}}
%Let $x=a_1a_2\cdots a_k$ be in $a\RR_X(a)$. Thus $a_1=a_k=a$.
%Then, with a slight
%abuse of notation,
%$f(x)=(a_1a_2)(a_2a_3)\cdots(a_{k-1}a_k)$. Since $\varphi(a)$ begins with $a$,
%the first element of $\varphi_2(a_1,a_2)$ is $(ab)$ and
%the last element of $\varphi_2(a_{k-1}a_k)$ is of the form $(ca)$
%with $c=a$ or $c=b$. Thus
%\begin{eqnarray*}
%\varphi_2\circ f(x)&=&\varphi_2(a_1a_2)\cdots \varphi_2(a_{k-1}a_k)\\
%&=&(ab)\cdots (ca)
%\end{eqnarray*}

 %This implies
%\begin{displaymath}
%\varphi_2\circ \phi(B)= \varphi_2\circ f(a\RR_X(a))
%\subset f(a\RR_X(a)^*)\subset\phi(B^*)
%\end{displaymath}
%and proves the existence of $\tau$.
Set $\tau:u\to v,v\to wu,w\to wvu$. Then it is easy
to verify that $\varphi_2\circ \phi=\phi\circ\tau$
(note that the existence of such $\tau$ is guaranteed for every
primitive morphism $\varphi$,
 see the proof of Proposition~\ref{lemma0Substitutions}).
The substitution $\tau$ is eventually proper and aperiodic
and $\phi(B)$ is circular. 
Thus, by Proposition~\ref{propositionBratteliPrimitive},
a BV-representation of $X$ is shown in Figure \ref{figureThueMorseDim3}.
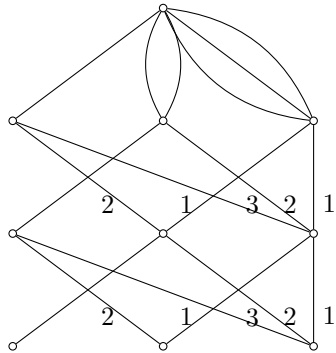
\begin{figure}[hbt]
\centering
\tikzset{node/.style={circle,draw,minimum size=0.1cm,inner sep=0pt}}
\tikzset{title/.style={minimum size=0.4cm,inner sep=0pt}}
\begin{tikzpicture}
\node[node](0)at(2,6){};
\node[node](11)at(0,4.5){};\node[node](12)at(2,4.5){};\node[node](13)at(4,4.5){};
\node[node](21)at(0,3){};\node[node](22)at(2,3){};\node[node](23)at(4,3){};
\node[node](31)at(0,1.5){};\node[node](32)at(2,1.5){};\node[node](33)at(4,1.5){};

\draw(0)edge node{}(11);
\draw[bend left](0)edge node{}(12);\draw[bend right](0)edge node{}(12);
\draw[bend left](0)edge node{}(13);\draw(0)edge node{}(13);
\draw[bend right](0)edge node{}(13);
\draw[left,near end](11)edge node{$2$}(22);
\draw[right,near end](11)edge node{$3$}(23);
\draw[left,near end](12)edge node{}(21);
\draw[right,near end](12)edge node{$2$}(23);
\draw[left,near end](13)edge node{$1$}(22);
\draw[right,near end](13)edge node{$1$}(23);
\draw[left,near end](21)edge node{$2$}(32);
\draw[right,near end](21)edge node{$3$}(33);
\draw[left,near end](22)edge node{}(31);
\draw[right,near end](22)edge node{$2$}(33);
\draw[left,near end](23)edge node{$1$}(32);
\draw[right,near end](23)edge node{$1$}(33);
\end{tikzpicture}
\caption{A BV-representation of the Thue-Morse shift.}
\label{figureThueMorseDim3}
\end{figure}
\end{solution}
\exosection{Section~\ref{ch5:subsec:LR}}
\begin{solution}{\ref{exerciseLRClosed}}
Let $f:\cL_k(X)\to A_k$ be a bijection extended as usual
to a map $f:\cL_{n+k-1}(X)\to \cL_n(X^{(k)})$. For every $w\in\cL_n(X^{(k)})$, we have
\begin{displaymath}
\RR_{X^{(k)}}(w)=f(\RR_X(u))
\end{displaymath}
where $w=f(u)$. Thus we have for every $x\in\RR_X(u)$
\begin{eqnarray*}
|f(x)|&=&|x|-k+1\le K|u|-k+1\le K(n+k-1)-k+1\\
&\le& K(k-1)+1.
\end{eqnarray*}
\end{solution}
\begin{solution}{\ref{exerciseExampleComplexityn^2}}
Let $x=\sigma^\omega(a)$.
We have $x=a\prod_{i\ge 0}b^{2^i}dc^i$ as one may verify
by computing $\sigma(x)$. For $0\le i\le j\le (n-2)/2$, the word
$b^idc^jbw$ in in $\cL(x)$ for some $w\in\{b,c,d\}^{n-i-j+2}$.
This shows that $p_n(x)$ grows like $n^2$.
\end{solution}
\begin{solution}{\ref{exerciseAltDefLR}}
The necessary condition is proven with Lemma \ref{lemmaCondh_k(n)}.

Let us show it is sufficient.
One has $h_j (n-1) \leq h_i (n)\leq L h_j (n-1)$, for all $i,j,n$.
Let $M(n+1,n)$ be the connecting matrices associated to the sequence of partitions $(\Pg (n))$.

and $M(n) (h_i()$ ... en attente de definition de matrices de partitions emboitees

\end{solution}
\begin{solution}{\ref{exerciseDelecroixBoshernitzan}}
  Assume first that $X$ is linearly recurrent with constant $K$.
  Let $\mu$ be an invariant probability measure on $X$.
  For $u\in \cL_n(X)$, the family $\Pg=\{S^j[vu]\mid v\in\RR'_X(u),0\le j<|v|\}$
  is a partition of $X$ (see Proposition~\ref{propositionPartitionReturn}).
  Thus, we have
  \begin{displaymath}
    \sum_{v\in\RR'_X(u)}|v|\mu([vu])=1
  \end{displaymath}
  Since $\mu$ is invariant, we have $\sum_{v\in\RR'_X(u)}\mu([vu])=\mu([u])$.
  Thus
  we obtain
  \begin{displaymath}
    1\le Kn\sum_{v\in\RR'_X(u)}\mu([u])= Kn\mu([u])
  \end{displaymath}
  and finally $n\varepsilon_n(\mu)\ge 1/K$.

  Conversely, assume $\inf n\varepsilon_n( \mu)>\varepsilon>0$
  where $\mu$ is some $S$-invariant measure on $X$. Let $u\in \cL_n (X)$ and
  $w=w_1w_2\cdots w_{N}\in\RR'_X(u)$ be of length $N$. We need to bound $N/n$.
  %We may assume $k\ge n$. 
  Consider, for $n\le k\le N$, the set $W_k $ of words in $\cL_N(X)$
  which end with $w_1\cdots w_k$. 
 We set $U_k = \cup_{u\in W_k} [u]$. 
  From the assumption we have
  $\mu(U_k)=\mu([w_1\cdots w_k])\ge \varepsilon/k$.
  Since the sets $U_k$ are disjoint
  we have
  \begin{displaymath}
    1\ge \sum_{k=n}^N\frac{\varepsilon}{k}\ge\varepsilon\int_{n}^{N+1}\frac{dx}{x}
    \ge\varepsilon\log((N+1)/n)
  \end{displaymath}
  and therefore we obtain the bound $N/n\le\exp(1/\varepsilon)$.
  This shows the linear recurrence for the constant $K = \exp(1/\varepsilon)$. 

Let $u\in \cL_n (X)$. From Proposition \ref{ch5:proposition:firstpropLR} for all $v\in \RR'_X (u)$ one has $(1/K)|u|\leq |v| \leq K |u|$.
Consequently,
$$
|u|\mu ([u])=|u|\sum_{v\in \RR_X'(u)}  \mu ([vu])\leq K\sum_{v\in \RR_X'(u)} |v| \mu ([vu])=K.
$$
Let $\mu,\nu$ be two distinct ergodic measures on $X$.
Since $\_mu,\_nu$ are mutually singular (Theorem \ref{propositionMutuallySingular})
the ratio  $\nu([u])/\mu([u])$ should be unbounded when $|u|$ goes to infinity
(indeed, otherwise, the ratio $\mu(U)/\nu(U)$ for $U,V$ clopen would also be bounded,
a contradiction with the fact that there exist clopen sets $U_0,V_0$ such that
$\mu(U\Delta U_0),\nu(U\Delta U_0)$ are arbitrarily small). 
But one has
$$
\frac{\nu([u])}{\mu([u])} = \frac{|u|\nu([u])}{|u|\mu([u])}\leq \frac{K}{\varepsilon} .
$$ 
Thus $X$ is uniquely ergodic.
  \end{solution}
\exosection{Section~\ref{sectionSadicShifts}}
\begin{solution}{\ref{exerciseLimitPoint}}
For a primitive $\Sa$-adic system, we
have $\langle\tau{[0,n)}\rangle\to\infty$. Thus (ii) and (iii) are equivalent
with the same system $\tau$.

(i) $\Rightarrow$ (iii) Let $a_n$ be the first letter of $w^{(n)}$.
Since $\tau_{[0,n)}(a_n)$ is a prefix of $x$ with length
which tends to $\infty$, this proves (iii).

(iii) $\Rightarrow$ (i) Set $w^{(0)}=x$.
Let $y_n$ be in $\tau_{[1,n)}(a_n)$ for each $n\ge 1$.
Let $y_{n_i}$ be a  subsequence of the $y_i$ converging to $y\in A_1^\N$. We replace $\tau$
by its telescoping with respect to the subsequence $n_i$
and the sequence $(a_n)$ by the sequence $(a_{n_i})$. Set $w^{(1)}=y$.
Then $x=\tau_0(w^{(1)})$
and $\{w^{(1)}\}=\cap_{n\ge 1}\tau[1,n)(a_n)$. Continuing in this
way, we build a sequence $w^{(n)}$ with the required properties.

\end{solution}
\begin{solution}{\ref{ExerciseSadicCassaigne}}
Define $\sigma_a$ by 
\begin{displaymath}
\sigma_a(b)=\begin{cases}a\#&\mbox{ if $b=a$}\\b&\mbox{ otherwise}
\end{cases}
\end{displaymath}
\end{solution}

\begin{solution}{\ref{exercise(aab,ba)}}
The set $\cL_2(X)$ is $\{aa,ab,ba,bb\}$ put in bijection
with $\{x,y,z,t\}$. The morphism $\tau_2$ is $x\to xyz$,
$y\to xyt$, $z\to zx$, $t\to zy$.
The  composition matrices  $M$ and $M_2$ are
\begin{displaymath}
M=\begin{bmatrix}2&1\\1&1\end{bmatrix},
\quad
M_2=\begin{bmatrix}1&1&1&0\\1&1&0&1\\1&0&1&0\\0&1&1&0\end{bmatrix}
\end{displaymath}
The graph $\Gamma_2(X)$ is the complete graph on two vertices and thus
a choice for the matrix $P$ with rows  a basis of its cycles is
\begin{displaymath}
P=\begin{bmatrix}1&0&0&0\\0&1&1&0\\0&0&0&1
\end{bmatrix},\quad N=\begin{bmatrix}1&1&0\\2&1&1\\0&1&0\end{bmatrix}
\end{displaymath}
The corresponding matrix $N$ such that $PM_2=NP$ is shown on the right.
A left eigenvector for the maximal eigenvalue $\lambda^2=(3+\sqrt{5})/2$
is $\begin{bmatrix}2\lambda&\lambda+1&1\end{bmatrix}$. We conclude
by Proposition~\ref{propositionDimensionGroupSubstitutions} that the dimension group has the form
indicated (note that $N$ is the matrix of Example~\ref{examplePrimitiveUnimodular}).
\end{solution}
\begin{solution}{\ref{exerciseM(u)}}
Assume  that $M_\varphi(u) \cap\cL(X)$
contains words with arbitrary large period. Then there is an
aperiodic point in $X$ of the form
$\cdots v_{-1}v_0v_1\cdots= \cdots w_{-1}w_0w_1\cdots$
where $uv_i=w_iu$, $v_i\in M_\varphi(u)$, $w_i\in U^*$
 and $u\in\cL(X)\setminus U^*$.
Thus $\varphi$ is not recognizable on $X$
for aperiodic points.

Assume now that $\varphi(B)$ is a prefix code and that
$\varphi$ is injective on $B$. This implies that
if $\varphi(x)=\varphi(x')$,
then $x^+=x'^+$. 

Suppose that $\varphi$ is not recognizable on $X$
at an aperiodic point $y=\varphi(x)$ and
let $(x',k)\ne(x,0)$ be such that 
$y=S^k\varphi(x')$ with $0\le k<|\varphi(x_0)|$.
We cannot have $x=x'$ since $y$ is aperiodic. Thus
$k\ne 0$.
For an infinity
of  $n<0$ there is an $m<0$ and $k_n$ with $0\le k_n<|\varphi(x_n)|$
such that $\varphi(S^nx)=S^{k_n}\varphi(S^mx')$. By the preceding remark,
we have $k_n\ne 0$. Let $u_n$ be the word of length $k_n$
such that $\varphi(S^nx)=u_n\varphi(S^mx')$. By shifting $y$
is necessary, we can assume that
all $u_n$ are equal to $u$ (see Figure~\ref{figureM(u)}).

 Then all words $v_m=\varphi(x'_{m}\cdots x'_{-1})$
are in $M_\varphi(u)$, which thus contains  words 
of arbitrary large period (see Figure~\ref{figureM(u)}).
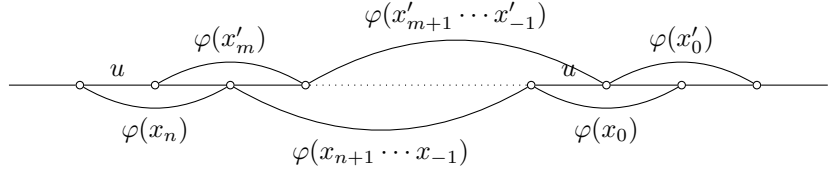
\begin{figure}[hbt]
\centering
\tikzset{node/.style={circle,draw,minimum size=0.1cm,inner sep=0pt}}
\tikzset{title/.style={circle,minimum size=0.1cm,inner sep=0pt}}
\begin{tikzpicture}
\node[title](-infty)at(-1,2){};
\node[node](n)at(0,2){};\node[node](m)at(1,2){};
\node[node](n+1)at(2,2){};\node[node](m+1)at(3,2){};

\node[node](0)at(6,2){};\node[node](1)at(7,2){};
\node[node](2)at(8,2){};\node[node](3)at(9,2){};
\node[title](infty)at(10,2){};

\draw[above](-infty)edge node{}(n);\draw[above](n)edge node{$u$}(m);
\draw[above](m)edge node{}(n+1);\draw[above](n+1)edge node{}(m+1);
\draw[below,bend right](n)edge node{$\varphi(x_{n})$}(n+1);
\draw[above,bend left](m)edge node{$\varphi(x'_{m})$}(m+1);
\draw[dotted](m+1)edge node{}(0);
\draw[above](0)edge node{$u$}(1);\draw[above](1)edge node{}(2);
\draw[below,bend right](0)edge node{$\varphi(x_{0})$}(2);
\draw[above,bend left](1)edge node{$\varphi(x'_{0})$}(3);
\draw(2)edge node{}(3);
\draw(3)edge node{}(infty);
\draw[above,bend left](m+1)edge node{$\varphi(x'_{m+1}\cdots x'_{-1})$}(1);
\draw[below,bend right](n+1)edge node{$\varphi(x_{n+1}\cdots x_{-1})$}(0);
\end{tikzpicture}
\caption{The factors of $y$ in $M_\varphi(u)$.}\label{figureM(u)}
\end{figure}
\end{solution}
\begin{solution}{\ref{exerciseLeftPermutative}}
Let $P$ be the set of prefixes of the words of $U=\varphi(A)$.
Since $\varphi$ is left permutative, there is for every $p\in P$
exactly one $a\in A$ such that $pa\in P\cup U$.

Assume that $M_\varphi(u)$ contains words arbitrary large period.
Let $v,w\in M_\varphi(u)$ be two words with distinct periods. Replacing if
necessary $v,w$ by some power, we may assume that none is a prefix
of the other. Let $r$ be the longest common prefix of $v,w$
and let $p\in P$ be such that $r=xp$ with $x\in U^*$.
Then $p$ has two right extensions by distinct letters, a contradiction.
Thus $\varphi$ is recognizable at aperiodic points
by Exercise~\ref{exerciseM(u)}.
\end{solution}
\exosection{Section~\ref{sectionDerivatives}}
\begin{solution}{\ref{exerciseErasing}}
Set $x=001^\omega$.
Let $\varphi:0\to 01,1\to 1$ and $\phi:0\to 00,1\to 1$.
Then $\varphi^\omega(0)=01^\omega$ and $\phi\circ\varphi^\omega(0)=x$.
Thus $001^\omega$ is substitutive. It cannot be purely
substitutive because if $\psi(x)=x$, we have $\psi(0)=001^n$
with $n\ge 0$ and thus $\psi(x)$ begins with $001^n001^n$.
\end{solution}
\begin{solution}{\ref{exerciseShiftDerivative}}
Let $\varphi:A^*\to A^*$ be a substitution and let $y$
be a fixed point of $\varphi$. Let $\phi:A^*\to B^*$
be a letter-to-letter morphism and let $x=\phi(y)$. 
Let $X=X(\varphi)$ be the shift generated by $y$.
  For $k\ge 1$, let
$f:\cL_{k}(X)\to A_k$ be a bijection from $\cL_{k}(X)$
onto an alphabet $A_k$ and let $\gamma_{k}:X\to A_k^\Z$
be the corresponding higher block code. The shift
$X^{(k)}$ is the substitution shift generated by
the $k$-th block presentation $\varphi_k$ of $\varphi$
and $z=\gamma_k(y)$ is a fixed point of $\varphi_k$. 
Define $\theta:A_k\to B$ by 
\begin{displaymath}
\theta\circ f(y_n\cdots y_{n+k-1})=\phi(y_{n+k-1}).
\end{displaymath}
Then $\theta(z)=\phi(T^{k-1}x)$ which shows that the sequence
$T^{k-1}x$ is substitutive.
\end{solution}
\begin{solution}{\ref{exerciseCobham}}
We first prove that we can modify the pair $(\tau,\phi)$
in such a way that 
\begin{equation}
|\phi\circ\tau(b)|\ge\phi(b)\label{eqphiCircvarphi}
\end{equation}
 for every $b\in B$ and with strict inequality when $b=a$..

Since $\lim|\tau^n(a)|=\infty$, there are $1\le j<k$ such that
\begin{displaymath}
|\phi\circ\tau^j(b)|\le|\phi\circ\tau^k(b)|
\end{displaymath}
for every $b\in B$ with strict inequality if $b=a$.
Set  $\tau'=\tau^{k-j}$, $\phi'=\phi\circ\tau^j$.
Then 
\begin{eqnarray*}
|\phi'\circ\tau'(b)|&=&|\phi\circ\tau^j\circ\tau{k-j}(b)|=|\phi\circ\tau^k(b)|\\
&\ge&|\phi\circ\tau^j(b)|=|\phi'(b)|
\end{eqnarray*}
for every $b\in B$ with strict inequality when $b=a$.

We now assume that $(\tau,\phi)$ satisfies Equation \eqref{eqphiCircvarphi}.
Proceeding as in the proof of Proposition~\ref{ch5:prop:rauzy},
we define an alphabet $C=\{b_b\mid b\in B, 1\le p\le\phi(b)\}$,
a map $\theta:C\to A$ by $\theta(b_p)=(\phi(b))_p$ and a map
$\gamma:B\to C^+$ by $\gamma(b)=b_1b_2\cdots b_{|\phi(b)|}$.
In this way, we have $\theta\circ\gamma=\phi$.

Finally, we define the substitution $\zeta$
essentially as in the proof of Proposition~\ref{ch5:prop:rauzy}
(with the difference that this time the inequality $|\tau(b)|\ge|\phi(b)|$
is replaced by the weaker inequality \eqref{eqphiCircvarphi}).
For every $b\in B$, we have 
\begin{displaymath}
|\gamma\circ\tau(b)|=|\phi\circ\tau(b)|\ge|\phi(b)|
\end{displaymath}
by \eqref{eqphiCircvarphi}. Thus, we can define words $w_1,w_2,\ldots,w_{|\phi(b)|}\in C^*$ such that
\begin{displaymath}
\gamma\circ\tau(b)=w_1w_2\ldots w_{|\phi(b)|}
\end{displaymath}
with $|w_1|>1$ when $b=a$. Then we define the morphism $\zeta:C^*\to C^*$
by 
\begin{displaymath}
\zeta(b_p)=w_p
\end{displaymath}
We then have by construction $\zeta\circ\gamma=\gamma\circ\tau$
and thus $x=\theta\circ\zeta^\omega(a_1)$.

\end{solution}
\begin{solution}{\ref{exerciseNonerasing}}
Let $\chi:0\to 0,1\to 1,2\to\varepsilon$ be the morphism erasing $0$.
Let $\mu:0\to 01,1\to 10$ be the Thue-Morse morphism
and let $t=\mu^\omega(0)$. Since $\mu\circ\chi=\chi\circ\sigma$,
we have $t=\chi(x)$. 

Let $\tau$ be a nonerasing substitution
such that $x=\tau^\omega(0)$. Then $\chi(\tau(2)^3)=\chi(\tau(2))^3$
is a factor of $t$ which is a cube and thus $\chi(\tau(2))=\varepsilon$.
Since the factors of $x$ in $2^*$ are $\varepsilon,2,22,222$, this
forces $\tau(2)=2$. Now $\tau(0)$ is a prefix of $x$ and thus
$\tau(0)=01u$ for some $u\in\{0,1,2\}^*$. Next, since
$\tau(1222)=\tau(1)222$ cannot end with $2^4$, we have
$\tau(1)=ya$ with $y\in \{0,1,2\}^*$ and $a\in\{0,1\}$.
Finally, $\tau(10)=ya01u$ has a factor of length $3$
in $\{0,1\}^*$ while the factors of $x$ in $\{0,1\}^*$
are of length at most $2$, a contradiction.
\end{solution}
%%%%%%%%%%%%%%%%%%%%%%%%%
\section{Notes}
\section{Odometers}
Odometers, also called \emph{solenoids}~\citep{KatokHasselblatt1995}
\index{subject}{solenoid}%
or \emph{adding machines}
\index{adding machine}%
\citep{Brown1976},
are a classical object in dynamical systems theory.

The ring $\Z_p$ of $p$-adic integers is contained in a field, called
the field of $p$-\emph{adic numbers}.
\index{subject}{p-adic@$p$-adic!number}\index{subject}{numbers!p-adic@$p$-adic}%
A classical reference to $p$-adic numbers and $p$-adic analysis
is~\cite{Koblitz1984}.
\index{names}{Koblitz, Neal}%
The factorial representation of integers (Exercise~\ref{exerciseFactorial})
is described in~\cite{Knuth1998}. Supernatural numbers 
(Exercise~\ref{exerciseSuperNatural}),
also called 
\emph{generalized natural numbers}
\index{subject}{generalized!natural numbers}%
 or \emph{Steinitz numbers}
\index{subject}{Steinitz number}\index{names}{Steinitz, Ernst}%
are used to define orders of profinite groups.
The fact that they give a complete invariant for odometers
(Exercise~\ref{exerciseSuperNatural}) is a result
proved by \cite{Glimm1960}
\index{names}{Glimm, James G.}%
 in the context of
\emph{uniformly hyperfinite algebras}
\index{subject}{uniformly!hyperfinite $C^*$-algebra}%
\index{subject}{algebra!uniformly hyperfinite}%
of \emph{UHF-algebras}
\index{subject}{UHF-algebra}\index{subject}{algebra!UHF}%
(see Chapter~\ref{chapterBratteli}).

The notion of expansive system is classical in topological dynamics (see~\cite{KatokHasselblatt1995}
for example). Exercises \ref{exerciseExpansive}
and  ~\ref{exerciseExpansive2}
are from \cite[Theorem 5.24]{Walters1982}
\index{names}{Walters, Peter}%
 (see also \cite{Kurka2003}\index{names}{K{\r{u}}rka, Petr}).

The notion of equicontinuity is classical in analysis
for a family of functions.\index{subject}{equicontinuous!family of functions}
The condition
defining an equicontinuous
dynamical system is equivalent to the equicontinuity of the family $T^n$
of maps from $X$ to itself.

The  equal path number property 
and Theorem \ref{theoremEqualPathNumber} are from  \cite{Gjerde&Johansen:2000}\index{names}{Gjerde, Richard}\index{names}{Johansen, Orjan}.

%substitutions
\subsection{Substitution shifts}
In \cite{VershikLivshits1992}\index{names}{Vershik, Anatol M.}\index{names}{Livshits, Anatol N.} the authors showed that when $\sigma$ is a primitive substitution then the subshift it generates can be represented (in a measure-theoretic sense) by an ordered Bratteli diagram $B$ where $\sigma$ is the substitution we read on $B$. 

Theorem \ref{ch5:subsec:Bratteli-substitution} 
  was first proven in \cite{Forrest:97}\index{names}{Forrest, Alan H.}.
The proofs given in that paper are mostly of existential nature and do
not state a method to compute effectively the BV-representation associated with substitution systems.
Another proof was given in \cite{Durand&Host&Skau:1999}\index{names}{Durand, Fabien}\index{names}{Host, Bernard}\index{names}{Skau, Christian F.} that provides such an algorithm.
Proposition \ref{ch5:proposition:substitutionread}
is from \cite{Durand&Host&Skau:1999}.

Proposition~\ref{propositionLemma15Forrest} is  \cite[lemma 15]{Forrest:97}. As mentioned
in this paper, the proof uses an important technique called \emph{state-splitting}
or \emph{symbol splitting}.
\index{subject}{state!splitting}\index{subject}{symbol splitting}%
We actually use in the proof of Proposition~\ref{propositionLemma15Forrest}
an \emph{output split}.
\index{subject}{output!split}%
 See \citep{LindMarcus1995} for a systematic presentation of state splitting.

Proposition \ref{ch5:prop:rauzy} is also
 from \cite{Durand&Host&Skau:1999}. It  is a modification
 of an unpublished result of Rauzy.\index{names}{Rauzy, G\'erard}.

 Corollary~\ref{corollaryMinimalPrimitive} appears (with
 another proof) in \citep{MaloneyRust2018}.
 \index{names}{Maloney, Gregory R.}\index{names}{Rust, Dan}%
 A similar result concerning the more general class
 of minimal substitutive shifts appears in \cite{Durand2013b},
 namely that every uniformly recurrent substitutive sequence (resp. shift)
 is primitive substitutive.

For more  details about the Chacon substitution, see \citep{Ferenczi1995}
or~\citep{PytheasFogg2002}.
\index{names}{Ferenczi, S\'ebastien}

%Theorem~\ref{ch5:th:DamanikLenz2006} is from
%\cite{Damanik&Lenz:2006}.
%\index{names}{Damanik, David}\index{names}{Lenz, Daniel H.}%

%LR shifts

%Proposition~\ref{ch5:proposition:onetoone}
%is  proved in \cite{Durand&Host&Skau:1999}\index{names}{Durand, Fabien}\index{names}{Host, Bernard}\index{names}{Skau, Christian F.}.

\subsection{Linearly recurrent shifts}
Linearly recurrent shifts, also called 
 \emph{linearly repetitive}
\index{subject}{linearly!repetitive}%
 shifts were introduced in \cite{Durand&Host&Skau:1999}.

Theorem \ref{ch5:th:DamanikLenz2006}
\index{names}{Damanik, David}\index{names}{Lenz, Daniel H.}%
is from \cite{Damanik&Lenz:2006}
(see also \cite{Shimomura2019}).\index{names}{Shimomura, Takashi}
The fact that the equivalent conditions of Theorem~\ref{ch5:th:DamanikLenz2006}
are also equivalent to  unique ergodicity is from \cite{Durand2000}.

Proposition~\ref{ch5:proposition:firstpropLR}
is proved in \cite{Durand&Host&Skau:1999}\index{names}{Durand, Fabien}\index{names}{Host, Bernard}\index{names}{Skau, Christian F.}.
The corollary asserting that the factor
complexity of a primitive substitution shift is at most
linear can be found in \cite{Michel1976} (see also \cite{Pansiot1984}).
Exercise~\ref{exerciseExampleComplexityn^2} is
from \cite[Example 10.4.1]{AlloucheShallit2003}.
The factor complexity of substitutive shifts has been
extensively studied. By a result
of \cite{EhrenfeuchtLeeRozenberg1975},
\index{names}{Ehrenfeucht, Andrew}\index{names}{Lee, Kwok Pun}%
\index{Rozenberg, Gregorz} one has always $p_n(X)= O(n^2)$.
See \cite{AlloucheShallit2003}
\index{names}{Allouche, Jean-Paul}\index{names}{Shallit, Jeffrey O.} for a survey of this question.

Theorem
\ref{ch5:theorem:LRBVrepresentation} is from 
 \cite{Durand2003}.
\index{names}{Durand, Fabien}%
Theorem \ref{theoremDownarowiczMaass} is from 
\cite{DownarowiczMaass2008}
\index{names}{Maass, Alejandro}\index{names}{Downarowicz, Tomasz}.

The unique ergodicity of linearly recurrent Cantor systems
(Theorem~\ref{theoremLRisUE})
is from \cite{Cortez&Durand&Host&Maass:2003}.
\index{names}{Cortez, Mar{\'i}a Isabel}\index{names}{Host, Bernard}
This result has been generalized by
\cite[Corollary4.14]{BezuglyiKwiatkowskiMedynetsSolomiak2013}
\index{names}{Bezuglyi, Sergey}\index{names}{Kwiatkowski, Jan}\index{names}{Medynets, Konstantin}\index{names}{Solomiak, Boris}%
to the case
of an infinite set of positive matrices $M(n)$ such that
$\|M(n)\|_1\le Cn$ for some constant $C$. The proof
uses Birkhoff contraction coefficient (Exercise~\ref{exerciseBirkhoffContractionCoefficient1}).
\index{subject}{Birkhoff!contraction coefficient}

%S-adic
\subsection{$\Sa$-adic shifts}
The notion of $\Sa$-adic shift
was introduced in~\cite{Ferenczi1996},
\index{names}{Ferenczi, S\'ebastien}%
 using a terminology initiated by Vershik and coined out by Bernard Host.
\index{names}{Host, Bernard}%
 For more information, see \cite{PytheasFogg2002} or
  \cite{BertheSteinerThuswaldnerYassawi2019} 
\index{names}{Berth\'e, Val\'erie}\index{names}{Steiner, Wolfgang}%
\index{names}{Thuswaldner, J\"org M.}\index{names}{Yassawi, Reem}%
or \cite{BertheDelecroix2014}.
\index{names}{Delecroix, Vincent}%
See also \cite{Thuswaldner2020} for a recent survey on $\Sa$-adic systems.

Proposition~\ref{propositionSadicPrimitive}
is from~\cite[Lemma 7]{Durand2000}). 

Lemma~\ref{lemma:proper} is a weaker version of~\cite[Corollary 2.3]{DurandLeroy2012}.
\index{names}{Durand, Fabien}\index{names}{Leroy, Julien}%
Theorem \ref{theo:BSTY} is 
\cite[Theorem 3.1]{BertheSteinerThuswaldnerYassawi2019}.
\index{names}{Steiner, Wolfgang}\index{names}{Thuswaldner, J\"org M.}
Proposition~\ref{propositionLimitPoint} is from
\cite{ArnouxMizutaniSellami2014}.
\index{names}{Arnoux, Pierre}\index{names}{Mizutani, Masahiro}%
\index{names}{Sellami, Tarek}%

The original reference for the Theorem of Fine-Wilf 
(Exercise~\ref{exerciseFineWilf}) is \citep{FineWilf1965}.

Proposition~\ref{propositionElementaryMorphism}
is~\cite[Exercise 5.1.5]{BerstelPerrinReutenauer2009}
where it is proved as a variant of a result called the
\emph{Defect Theorem}\index{subject}{Defect Theorem} (see~\cite{Lothaire1983}).

The notion of recognizability for aperiodic points
has been introduced by \cite{BezuglyiKwiatkowskiMedynets2009}
who have used it to establish a generalization
of Mosse's Theorem to non primitive morphisms.
Theorem~\ref{theo:BSTY} is from \cite[Theorem 5.1]{BertheSteinerThuswaldnerYassawi2019}. Exercise~\ref{exerciseLeftPermutative}
is from~\cite[Lemma 3.3]{BertheSteinerThuswaldnerYassawi2019}.
It improves some results of \cite{BezuglyiKwiatkowskiMedynets2009}.
The argument on which the solution of Exercise~\ref{exerciseLeftPermutative}
is based, is used in a different context in~\cite{PerrinRindone2003}.
\index{names}{Perrin, Dominique}\index{names}{Rindone, Giuseppina}

Theorem \ref{theo:cohoword} is from 
\cite{BertheCecchiDurandLeroyPerrinPetite2020}.
\index{names}{Leroy, Julien}\index{names}{Cecchi Bernales, Paulina}%
\index{names}{Petite, Samuel}\index{names}{Perrin, Dominique}%

Corollary \ref{coro:measures}   extends a statement initially proved for interval exchanges by \cite{FerencziZamboni2008}.
\index{names}{Ferenczi, S\'ebastien}\index{names}{Zamboni, Luca Q.}%
In Corollaries \ref{coro:measures} and \ref{coro:freq},
   the assumption of being proper can be dropped. 
The   proof  then  uses the measure-theoretical 
Bratteli-Vershik representation of 
primitive unimodular $\Sa$-adic subshift given in \cite[Theorem 6.5]{BertheSteinerThuswaldnerYassawi2019}. 

Note that  one  recovers, with the description
of dimension groups of primitive unimodular proper $\Sa$-adic shifts
of Theorem~\ref{theo:dg},    the results   obtained in the case of interval exchanges,in   \cite{Putnam1989}. See also \cite{Putnam1992,Gjerde&Johansen:2002}.

\subsection{Derivatives of substitutive sequences}
Substitutive sequences have been considered early by 
\cite{Cobham1968}\index{names}{Cobham, Alan} who has proved the statement that
every infinite sequence $x=\phi(\tau^\omega(a))$
is substitutive, whatever be the morphisms $\phi:B^*\to A^*$
and $\tau:B^*\to B^*$ (Exercise~\ref{exerciseCobham}). This
result was proved independently by \cite{Pansiot1983}
\index{names}{Pansiot, Jean-Jacques}%
(see the presentation in \cite{AlloucheShallit2003}).
\index{names}{Allouche, Jean-Paul}\index{names}{Shallit, Jeffrey O.}%
We follow the proof given in \cite{CassaigneNicolas2003}
\index{names}{Cassaigne, Julien}\index{names}{Nicolas, Fran\c{c}ois}%
The effective computability of the representation
as a substitutive sequence was proved by \cite{Honkala2009}
\index{names}{Honkala, Juha}%
and \cite{Durand2013}.

Theorem~\ref{theoremCharacterisationSubstitutive}
is from \cite{Durand1998}.
\index{names}{Durand, Fabien}%
 It is closely related
with the results of \cite{HoltonZamboni1999}
\index{names}{Holton, Charles}\index{names}{Zamboni, Luca Q.}%
 who proved independently
that conditions (i) and (ii) are equivalent.
Theorem~\ref{theoremHoltonZamboni} is from \cite[Theorem 1.3]{HoltonZamboni1999}

\section{Exercises}

The notion of profinite group (Exercise~\ref{exerciseProfiniteGroups})
generalizes the profinite integers of Section~\ref{ch5:subsec:rep-odo}.
See \cite{AlmeidaCostaKyriakoglouPerrin2020} for a description of the link between profinite groups or semigroups
and symbolic dynamics.
\index{names}{Almeida, Jorge}\index{names}{Costa, Alfredo}%
\index{names}{Kyriakoglou, Revekka}\index{names}{Perrin, Dominique}%

The Chacon ternary substitution
\index{subject}{Chacon!ternary!substitution}%
(Exercises~\ref{exerciseChaconTernaryBinary}
and \ref{exerciseComplexityChacon})
is from \cite{PytheasFogg2002}.

Exercise~\ref{exerciseDelecroixBoshernitzan} is
from \cite{Delecroix2015}
\index{names}{Delecroix, Vincent}%
(see also~\cite{BesbesBoshernitzanLenz2013})
\index{names}{Besbes, Adnene}\index{names}{Boshernitzan, Michael}%
\index{names}{Lenz, Daniel H.} where the result is credited to
Boshernitzan (private communication). The condition
$\inf n\varepsilon_n(X)>0$ is a reinforcement of the
condition $\limsup n\varepsilon_n(X)>0$ which is
proved to be equivalent to unique ergodicity in
\cite{Boshernitzan1992}.

%%%%%%%%%%%%%%%%%%%%%%%%%%%%
\chapter{Dendric shifts}  %
%%%%%%%%%%%%%%%%%%%%%%%%%%%%
\label{chapterDendricShifts}
In this chapter, we define the important class of dendric shifts,
which are defined by a condition on the possible extensions of a word
in their language. We prove a striking property of the sets
of return words, namely that for every minimal dendric 
shift, the set of return words forms a basis of the free group
(Theorem~\ref{theoremReturn}).
We show that they have a finite $\Sa$-adic representation (Theorem~\ref{theoremSadicDendric}).
We  illustrate these results on the
class of Sturmian shifts (Section~\ref{sectionChapter6Sturmian})
which are a particular case of interval
exchange shifts considered in the next chapter.
We next present the class of specular shifts (Section~\ref{sectionSpecular})
  which is
build to generalize the class of linear involutions,
itself a natural generalization of interval exchanges, and also presented in the
next chapter.
%%%%%%%%%%%%%%%%%%%%%%%%%%%%%%%
\section{Dendric shifts}\label{sectionDendric}
%%%%%%%%%%%%%%%%%%%%%%%%%%%%%%%
Let $X$ be a shift space on the alphabet $A$.
We will assume in this chapter that $A\subset \cL(X)$.
For $w \in \cL(X)$ and $n \ge 1$, we denote
\begin{eqnarray*}
L_X(w) & = & \{ a \in A \mid aw \in \cL(X) \} \\
R_X(w) & = & \{ b\in A \mid wb \in \cL(X) \} \\
E_X(w) & = & \{ (a,b) \in L_X(w) \times R_X(w) \mid awb \in \cL(X) \}
\end{eqnarray*}
\index{symbols}{L@$L_X(w)$}\index{symbols}{R@$R_X(w)$}%
\index{symbols}{E@$E_X(w)$}%
The \emph{extension graph}
\index{subject}{extension!graph}\index{subject}{graph!extension}%
  of $w$, denoted $\E_X(w)$,\index{symbols}{E@$\E_X(w)$}
 is the undirected bipartite graph whose set of vertices is the disjoint
 union of $L_X(w)$ and $R_X(w)$ and whose 
 edges are the elements of $E_X(w)$.

When the context is clear, we denote $L(w), R(w), E(w)$ and $\E(w)$ instead of $L_X(w), R_X(w), E_X(w)$ and $\E_X(w)$.

When in need to distinguish the
disjoint copies of $L(w)$ and $R(w)$
forming the vertices of the extension graph $\E(w)$,
we denote them  by $1\otimes L(w)$
\index{symbols}{one@$1\otimes L(w)$} and $R(w)\otimes 1$.
\index{symbols}{R@$R(w)\otimes 1$}%
A path in an undirected graph is \emph{reduced}
\index{subject}{reduced path}\index{path!reduced}%
if it does not contain successive equal edges.
For any $w \in \cL(X)$, since any vertex of $L_X(w)$ is connected
by an edge to at least one vertex of $R_X(w)$, the bipartite graph $\E_X(w)$ is a tree if and only if there is a unique reduced path between every pair of vertices of $L_X(w)$ (resp. $R_X(w)$).

The shift $X$ is said to be \emph{eventually dendric}
\index{subject}{eventually!dendric shift}\index{subject}{shift space!eventually dendric}
 with \emph{threshold}
\index{subject}{threshold of eventually dendric shift}%
\index{subject}{eventually!dendric shift!threshold}%
 $m \ge 0$ if $\E_X(w)$ is a tree for every word $w \in \cL_{\ge m}(X)$.
It is said to be \emph{dendric}
\index{subject}{dendric!shift space}\index{subject}{shift space!dendric}%
 if we can choose $m = 0$.
Thus, a shift $X$ is dendric if and only if $\E_X(w)$ is a tree for every word $w \in \cL(X)$.

When $X$ is a dendric shift (resp. eventually dendric shift), we also say that 
$\cL(X)$ is a \emph{dendric set} (resp. eventually dendric set).
\index{subject}{dendric!set}

An important observation is that, in any shift space,
 for a word $w\in\cL(X)$ which is not bispecial,
the graph $\E(w)$ is always a tree. Indeed, if $w$ is not left-special,
all vertices of $R(w)$ are connected to the unique vertex of $L(w)$
and thus $\E(w)$ is a tree. Also, if $w$ is bispecial and such that
$\E(w)$ is a tree, then $w$ is neutral. Indeed, in a tree
$m(w)=e(w)-\ell(w)-r(w)+1=0$ by a well known property
of trees.
We begin with an example of a non minimal dendric shift.
\begin{example}
Let $X$ be the shift space such that $\cL(X)=a^*ba^*$ (we denote $a^*=\{a^n\mid n\ge 0\}$).
The bispecial words are the words in $a^*$.
Their extension graph is the tree represented in Figure~\ref{figurea*ba*}.
Thus $X$ is a dendric shift.
\begin{figure}[hbt]
	\centering
	%%% IMAGE EN TIKZ
	\tikzset{node/.style={rectangle,draw,rounded corners=1.2ex}}
	\begin{tikzpicture}
	\node[node](a1al) {$a$};
	\node[node](a1bl) [below= 0.3cm of a1al] {$b$};
	\node[node](a1ar) [right= 1.5cm of a1al] {$a$};
	\node[node](a1br) [below= 0.3cm of a1ar] {$b$};
	\path[draw,thick, shorten <=0 -1pt, shorten >=-1pt]
	(a1al) edge node {} (a1br)
	(a1bl) edge node {} (a1ar);
	\path[draw,thick, shorten <=0 pt, shorten >=-0pt]
	(a1al) edge node {} (a1ar);
	
	\end{tikzpicture}

 \caption{The graph $\E(a)$.}
 \label{figurea*ba*}
\end{figure}
\end{example}
An important example of minimal dendric shifts is formed by 
\emph{strict episturmian shifts}
\index{subject}{strict episturmian!shift}%
 (also called \emph{Arnoux-Rauzy shifts}).
\index{subject}{Arnoux-Rauzy!shift}%
\index{subject}{shift space!Arnoux Rauzy}%
\begin{proposition}\label{propositionARisDendric}
Every Arnoux-Rauzy shift is dendric. In particular, every Sturmian shift
is dendric.
\end{proposition}
\begin{proof}
Let $X$ be an Arnoux-Rauzy shift. For every bispecial
word $w\in \cL(X)$, there is exactly one letter $\ell$ such that $\ell w$ is
right-special and one letter $r$ such that $wr$ is left-special.
Thus the extension
graph $\E(w)$ has exactly two vertices $\ell,r$ which have degree more than
one, with $\ell\in L(w)$ and  $r\in R(w)$. Any vertex 
distinct from $\ell,r$ is connected
to either $\ell$ or $r$ by an edge but not to both and $\ell,r$
are connected by an edge. Thus $\E(w)$ is a tree.
\end{proof}

\begin{example}
\label{exampleFibo}
Let $X$ be the \emph{Fibonacci shift},
\index{subject}{Fibonacci!shift}\index{subject}{shift space!Fibonacci}%
 which is generated by the morphism $a \mapsto ab, b \mapsto a$.
It is a Sturmian shift (Example \ref{exampleFibonacci0bis}).
The graph $\E(a)$ is shown in Figure~\ref{figureFibo}.

	\begin{figure}[hbt]
	\centering
	%%% IMAGE EN TIKZ
	\tikzset{node/.style={rectangle,draw,rounded corners=1.2ex}}
	\begin{tikzpicture}
	\node[node](a1al) {$a$};
	\node[node](a1bl) [below= 0.3cm of a1al] {$b$};
	\node[node](a1ar) [right= 1.5cm of a1al] {$a$};
	\node[node](a1br) [below= 0.3cm of a1ar] {$b$};
	\path[draw,thick, shorten <=0 -1pt, shorten >=-1pt]
	(a1al) edge node {} (a1br)
	(a1bl) edge node {} (a1ar);
	\path[draw,thick, shorten <=0 pt, shorten >=-0pt]
	(a1bl) edge node {} (a1br);
	
	\end{tikzpicture}

 \caption{The graph $\E(a)$.}
 \label{figureFibo}
\end{figure}
\end{example}

A shift space $X$ is said to be a \emph{eventually
dendric of characteristic}
\index{subject}{eventually!dendric!of characteristic $c$}%
\index{subject}{characteristic!of eventually dendric shift}%
 $c$ if
\begin{enumerate}
\item for
 any $w \in \cL_{\ge 1}(X)$, the extension graph $\E(w)$
 is a tree and 
\item the graph $\E(\varepsilon)$ is a disjoint union of $c$ trees.
\end{enumerate}
Eventually dendric shifts of characteristic $c \ge 1$ are eventually dendric.
 Indeed, since the extension graphs of all nonempty words are trees, the shift space is eventually dendric with threshold $1$.
It is a dendric shift if $c=1$.
\begin{example}
\label{exampleJulien}
Let $X$ be the shift generated by the morphism $a \mapsto ab, b \mapsto cda, c \mapsto cd, d\mapsto abc$.
It is dendric of characteristic $2$ (Exercise~\ref{exerciseExampleJulienCassaigne}).
The extension graph $\E(\varepsilon)$ is shown in Figure~\ref{figureSpecular}.

\begin{figure}[hbt]
	\centering
	\tikzset{node/.style={rectangle,draw,rounded corners=1.4ex}}
	\begin{tikzpicture}
	\node[node](eal) at(0,1){$a$};
	\node[node](ebl) at(0,0){$b$};
	\node[node](ebr) at(2,1){$b$};
	\node[node](ecr) at(2,0) {$c$};
	\path[draw,thick]
	(eal) edge node {} (ebr)
	(ebl) edge node {} (ecr);
	\path[draw,thick, shorten <=0 -1pt, shorten >=-1pt]
	(eal) edge node {} (ecr);
	\node[node](ecl) at(4,1) {$c$};
	\node[node](edl) at(4,0) {$d$};
	\node[node](edr) at(6,1) {$d$};
	\node[node](ear) at(6,0) {$a$};
	\path[draw,thick]
	(ecl) edge node {} (edr)
	(edl) edge node {} (ear);
	\path[draw,thick, shorten <=0 -1pt, shorten >=-1pt]
	(ecl) edge node {} (ear);
	\end{tikzpicture}
	\caption{The extension graph $\E(\varepsilon)$.}
	\label{figureSpecular}
\end{figure}
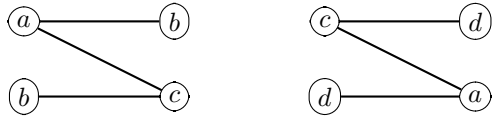

\end{example}

\begin{example}\label{ex:Tribo}
Let $X$ be the \emph{Tribonacci shift},
\index{subject}{Tribonacci!shift}\index{subject}{shift space!Tribonacci}%
 which is the substitution shift generated by the \emph{Tribonacci substitution}
\index{Tribonacci!substitution}%
  $\sigma: a \mapsto ab, b \mapsto ac, c \mapsto a$.
		It is an Arnoux-Rauzy shift
(see Example~\ref{exampleTribonacci}) and thus a dendric shift.
\end{example}

The following statement shows that eventually dendric shifts have at most
linear complexity. Recall from Chapter~\ref{chapterTopologicalDynamicalSystems}
that we denote $p_n(X)=\Card(\cL_n(X))$ and that $s_n(X)=p_{n+1}(X)-p_n(X)$.
\begin{proposition}
\label{propositionComplexity}
Let $X$ be an eventually dendric shift  on the alphabet $A$.
Then $X$ has at most linear complexity, that is,
 $p_n(X)\le Kn$ for some constant $K$. If $X$ is dendric, then
\begin{equation}
p_n(X)=(\Card(A)-1)n+1\label{eqComplexityDendric}
\end{equation}
\end{proposition}
\begin{proof}
Let $b_n(X)=s_{n+1}(X)-s_n(X)$.
Since $X$ is eventually dendric, there is
 $n \ge 1$  such that the extension graph of every word in $\cL_{\ge n}(X)$ is a tree. 
Then $b_p(X) = 0$ for every $p \ge n$. 
Indeed, by Proposition~\ref{propositionCassaigne}, we have $b_p(X)=\sum_{w\in \cL_p(X)}m(w)$. Since all words of length $p$ in $\cL(X)$ are neutral,
the conclusion follows.
Thus $s_p(X) = s_{p+1}(X)$ for every $p \ge n$, whence our conclusion.

If $X$ is dendric, then $b_n(X)=0$ for all $n\ge 1$ and thus
$s_n(X)$ is constant. Since $s_1(X)=\Card(A)-1$
by the assumption $A\subset\cL(X)$, this implies
$p_n(X)=(\Card(A)-1)n+1$.
\end{proof}
A more precise formula can be given for the complexity of
eventually dendric shifts (Exercise~\ref{exerciseComplexityEentuallyDendric}).

\begin{corollary}
The minimal dendric shifts on two letters are the Sturmian shifts.
\end{corollary}
\begin{proof}
We have already seen that a Sturmian shift is dendric 
(Proposition~\ref{propositionARisDendric}). Conversely, if $X$ is a dendric
shift on two letters, its factor complexity is $p_n(X)=n+1$
by Proposition~\ref{propositionComplexity} and thus $X$ is Sturmian.
\end{proof}
On more than two letters, the class of dendric shifts is larger
than the class of Arnoux-Rauzy shifts since it contains, as we shall
see in Chapter~\ref{chapterIET}, the class of interval exchange shifts.

The converse of Proposition~\ref{propositionComplexity} is not true, as shown by the following example.

\begin{example}
The \emph{Chacon ternary shift}
 \index{subject}{Chacon!ternary!shift}%
 is the substitution shift $X$ on the alphabet
$A=\{0,1,2\}$ generated 
by the Chacon ternary substitution \index{subject}{Chacon!ternary!substitution}
 $\tau:0\to 0012,1\to 12,2\to 012$.
Its complexity is $p_n(X)=2n+1$ (see Exercise~\ref{exerciseComplexityChacon}).
It is not eventually dendric (Exercise~\ref{exerciseChaconNotEventuallyDendric}).
\end{example}
\subsection{Generalized extension graphs}
We will need to consider extension graphs which correspond to extensions
by words instead of letters.
Let $X$ be a shift space. For $w\in \cL(X)$, and $U,V\subset \cL(X)$,
let 
$L_U(w)=\{\ell\in U\mid \ell w\in \cL(X)\}$,
 let $R_V(w)=\{r\in V\mid wr\in \cL(X)\}$ and 
let $E_{U,V}(w)=\{(\ell,r)\in U\times V\mid \ell w r\in \cL(X)\}$.
\index{symbols}{L@$L_U(w)$}\index{symbols}{R@$R_V(w)$}%
\index{symbols}{E@$E_{U,V}(w)$}%
The \emph{generalized extension graph}
\index{subject}{generalized!extension graph}\index{subject}{extension!graph!generalized}%
 of $w$ relative to
$U,V$ is the following undirected graph $\E_{U,V}(w)$. The set of vertices is
made of two disjoint copies of $L_U(w)$ and $R_V(w)$.
 The edges are the elements of $E_{U,V}(w)$. The extension graph $\E(w)$ defined previously
corresponds
to the case where $U,V=A$.

\begin{example}
Let $X$ be the Fibonacci shift.
\index{subject}{Fibonacci!shift}\index{subject}{shift space!Fibonacci}%
 Let $w=a$, $U=\{aa,ba,b\}$ and let
$V=\{aa,ab,b\}$. The graph $\E_{U,V}(w)$ is represented in
Figure~\ref{figureStrongTree}.
\begin{figure}[hbt]
\centering
\tikzset{node/.style={circle,draw,minimum size=0.4cm,inner sep=0.2pt}}
\begin{tikzpicture}
\node[node](b)at(0,0){$b$};\node[node](ba)at(0,.7){$ba$};
\node[node](ab)at(2,0){$ab$};\node[node](b')at(2,.7){$b$};

\draw(ba)edge node{}(b');\draw(b)edge node{}(b');
\draw(b)edge node{}(ab);
\end{tikzpicture}
\caption{The graph $\E_{U,V}(w)$.}\label{figureStrongTree}
\end{figure}
\end{example}

The following property shows that in a dendric shift, not only
the extension graphs but, under appropriate hypotheses,
 all generalized extension graphs
are acyclic.
\begin{proposition}\label{PropStrongTreeCondition}
Let $X$ be a shift space and $n\ge 1$ be such
that for every $w\in\cL_{\ge n}(X)$, the graph $\E(w)$ is acyclic.
  Then, for any $w\in \cL_{\ge n}(X)$, any
 finite 
suffix code  $U$ and any  finite  prefix code $V$,
 the generalized extension graph $\E_{U,V}(w)$ is  acyclic.
\end{proposition}
The proof uses the following lemma.
\begin{lemma}\label{lemmaTree}
Let $X$ be a  shift space. Let $w\in \cL(X)$ and let $U,V,T\subset \cL(X)$.
Let
 $\ell\in \cL(X)\setminus U$  be such that $\ell w\in \cL(X)$. Set
 $U'=(U\setminus T\ell
)\cup\ell$. If the graphs $\E_{U',V}(w)$
and $\E_{T,V}(\ell w)$
are acyclic   then $\E_{U,V}(w)$ is  acyclic.
\end{lemma}
\begin{proof}
 Assume  that
$\E_{U,V}(w)$ contains a cycle $C$. If the cycle does not use any
of the vertices
in $U'$, it defines a cycle in the graph $\E_{T,V}(\ell w)$
obtained by replacing each vertex $t\ell$ for $t\in T$ by a vertex $t$.
Since $\E_{T,V}(\ell w)$ is acyclic, this is impossible.
If it uses a vertex of $U'$ it defines a cycle of
the graph $\E_{U',V}(w)$ obtained by replacing each possible vertex $t\ell$
by $\ell$ (and suppressing the possible identical
successive edges created by the identification). This is impossible since $\E_{U',V}(w)$ is acyclic.
Thus $\E_{U,V}(w)$ is  acyclic.
\end{proof}

\begin{proofof}{of Proposition~\ref{PropStrongTreeCondition}}
We show by induction on the sum of the lengths of the words in $U,V$
that for any $w\in \cL_{\ge n}(X)$, the graph $\E_{U,V}(w)$ is  acyclic.

Let $w\in \cL_{\ge n}(X)$. We may assume that $U=L_U(w)$ and $V=R_V(w)$
and also that $U,V\ne\emptyset$.
If $U,V\subset A$, the property is true. 

Otherwise, assume for example that $U$ contains words of length at
least $2$. 
Let $u\in U$ be of maximal length. Set $u=a\ell$ with $a\in A$.
Let $T=\{b\in A\mid b\ell\in U\}$.
 Then $U'=(U\setminus T\ell
)\cup\ell$ is a suffix code and $\ell w\in \cL(X)$ since $U=U(w)$. 

By induction hypothesis, the graphs $\E_{U',V}(w)$
and $\E_{T,V}(\ell w)$
are acyclic. 
By lemma~\ref{lemmaTree}, the graph  $\E_{U,V}(w)$ is acyclic.
\end{proofof}
We prove now a similar statement concerning connectedness. 
For $w\in\cL(X)$, we say that a suffix code $U\subset \cL(X)$ is
\emph{$(X,w)$-maximal} if every word $v$ such that $vw\in\cL(X)$
is comparable with a word in $U$ for the suffix order.
For $w=\varepsilon$, we say $X$-maximal instead of $(X,\varepsilon)$-maximal.
Thus a suffix code is $X$-maximal if it is not strictly contained
in any suffix code $U'\subset\cL(X)$ 
\index{subject}{prefix!code!X-maximal@$X$-maximal}%
\index{subject}{X-maximal@$X$-maximal!prefix code}%
\index{subject}{suffix!code!X-maximal@$X$-maximal}%
\index{subject}{X-maximal@$X$-maximal!suffix code}%
 The same definitions hold symmetrically for prefix codes.
Thus a prefix code $V\subset \cL(X)$ is $(X,w)$-maximal if
every word $v$ such that $wv\in\cL(X)$ is comparable
for the prefix order with a word of $V$.

For example, when $X$ is recurrent, the set $\RR_X(w)$
is an $(X,w)$-maximal prefix code and $\RR'_X(w)$
is an $(X,w)$-maximal suffix code.
\begin{proposition}\label{propStrongTreeConditionBis}
Let $X$ be an eventually dendric shift with threshold $n$.
  For any $w\in \cL_{\ge n}(X)$, any
 finite 
$(X,w)$-maximal suffix code  $U\subset \cL(X)$ 
and any  finite  $(X,w)$-maximal prefix code $V\subset \cL(X)$,
 the generalized extension graph $\E_{U,V}(w)$ is  a tree.
\end{proposition}

For a shift space $X$, two finite sets $U,V\subset\cL(X)$  and $w\in\cL(X)$, denote
$\ell_U(w)=\Card(L_U(w))$, $r_V(w)=\Card(R_V(w))$ and
$e_{U,V}(w)=\Card(\E_{U,V}(w))$. Next, we define
\begin{displaymath}
m_{U,V}(w)=e_{U,V}(w)-\ell_{U,V}(w)-r_{U,V}(w)+1.
\end{displaymath}
Thus, for $U=V=A$, the integer $m_{U,V}(w)$ is the multiplicity $m(w)$ of $w$.
\begin{lemma}\label{lemmaExtendedMultiplicity}
Let $X$ be a shift space and $n\ge 0$ be such that $m(w)=0$
for every $w\in\cL_{\ge n+1}(X)$. Then for every $w\in\cL_{\ge n}(X)$, every finite
$(X,w)$-maximal suffix code $U$ and every finite 
$(X,w)$-maximal prefix code $V$,
we have $m_{U,V}(w)=m(w)$.
\end{lemma}
\begin{proof}
We use an induction on the sum of the lengths of the words in $U$ and in $V$.
We may assume that $Uw,wV\subset\cL(X)$.

If $U,V$ contain only words of length $1$,
 since $U$ (resp. $V$) is an $(X,w)$-maximal suffix (resp. prefix) code,
 we have 
$U = L(w)$ and $V = R(w)$
and there is nothing to prove.
Assume next that one of them, say $V$, contains words of length at least $2$.
 Let $p$ be a nonempty proper prefix of $V$. Set $V'=(V\setminus pA)\cup\{p\}$.
If $wp\notin \cL(X)$, then $m_{U,V}( w ) = m_{U,V'}( w )$
 and the conclusion follows by induction
hypothesis. Thus we may assume that $wp \in\cL(X)$. Then
\begin{displaymath}
m_{U,V'} ( w ) - m_{U,V} ( w ) =  e_{U , A} ( wp ) - \ell_U ( wp ) - r_A ( wp ) + 1 = m_{U,A} ( wp ).
\end{displaymath}
By induction hypothesis, we have $m_{U,V} ( w ) = m ( w )$. But $m_{U,A} ( wp ) = 0$
since $|wp|\ge n+1$, whence the conclusion.
\end{proof}

\begin{proofof}{of Proposition~\ref{propStrongTreeConditionBis}}
Let $w\in\cL_{\ge n}(X)$. By Proposition~\ref{PropStrongTreeCondition},
the graph $\E_{U,V}(w)$ is acyclic. Since, by Lemma~\ref{lemmaExtendedMultiplicity}, we have $m_{U,V}(w)=0$, it follows that $\E_{U,V}(w)$ is a tree.
\end{proofof}

\subsection{Return Theorem}

\index{subject}{Return Theorem}\index{subject}{Theorem!Return}%

We will now prove the following result (called
the \emph{Return Theorem}). Recall that we assume in this chapter that
$A\subset\cL(X)$
for a shift space $X$ on the alphabet $A$. 
\begin{theorem}\label{theoremReturn}
Let $X$ be a minimal dendric shift on the alphabet $A$.
For every $u\in\cL(X)$, the set $\RR_X(u)$
is a basis of the free group on $A$.
\end{theorem}

Note that that, in the particular case of an episturmian shift $X$, the property
results directly from Equation 
\eqref{equationReturnPal}. Indeed, the set of return words $\RR'_X(u)$
is for every $u\in\cL(X)$ conjugate to a set of the form $\alpha(A)$
where $\alpha$ is an automorphism of the free group on $A$.

A shift space $X$ is \emph{neutral}\index{subject}{neutral!shift space}
\index{subject}{shift space!neutral}%
if every word $u$ in $\cL(X)$ is neutral. A dendric shift is of course
neutral but the converse is false (see Exercise~\ref{exerciseNeutral}).
The first step of the proof of Theorem~\ref{theoremReturn} is the following statement. It shows
in particular that, under the hypotheses below, the
cardinality of sets of return words in constant. We had
already met this property in the case of strict episturmian shifts
(see Equation~\ref{equationReturnPal})
\begin{theorem}\label{theoremCardinality}
If $X$ is a recurrent neutral shift such that $A\subset\cL(X)$, then 
for every $u\in\cL(X)$, one has $\Card(\RR_X(u))=\Card(A)$.
\end{theorem}

Note the surprising consequence that every recurrent neutral shift
is minimal. Indeed, if $X$ is recurrent, all sets of return words are finite
by Theorem~\ref{theoremCardinality}. Thus
the shift is minimal. Thus, we could weaken the hypothesis in Theorem~\ref{theoremReturn}
to require $X$ to be only recurrent.

The proof of Theorem~\ref{theoremCardinality} uses the following lemma.
\begin{lemma}\label{lemmaLeftproba}
Let $X$ be a neutral shift. For every $v\in\cL(X)$, set $\rho(v)=r_X(v)-1$.
\index{symbols}{rho@$\rho_X(v)$}%
Then one has
\begin{equation}
\sum_{a\in L(v)}\rho(av)=\rho(v).\label{eqLeftProba}
\end{equation}
\end{lemma}
\begin{proof}
Since $v$ neutral, we have $e_X(v)-\ell_X(v)-r_X(v)+1=0$. Thus
\begin{eqnarray*}
\sum_{a\in L(v)}\rho(av)&=& \sum_{a\in L(v)}(r_X(av)-1)=e_X(v)-\ell_X(v)\\
&=&r_X(v)-1=\rho(v)
\end{eqnarray*}.
\end{proof}

\begin{proofof}{of Theorem~\ref{theoremCardinality}}
Let $U$ be the set of proper prefixes of $u\RR_X(u)$
which are not proper prefixes of $u$. We claim that $U$ is
an $X$-maximal suffix code. Indeed, assume first that $v,v'\in U$
with $v$ is a proper suffix of $v'$. Then $u$ is a proper prefix of $v$
and  thus $u$ appears as a factor of $v'$ otherwise than a suffix (see
Figure~\ref{figureSuffixCode}), a contradiction.

\begin{figure}[hbt]
\centering
\tikzset{node/.style={circle,draw,minimum size=0.1cm, inner sep=0cm}}
\begin{tikzpicture}

\node[node](v')at (0,2){};\node[node](v'r)at(4,2){};
\node[node](v)at (1,1.5){};\node[node](vr)at (4,1.5){};
\node[node](w)at(1,1){};\node[node](wr)at (3,1){};

\draw[above](v')edge node{$v'$}(v'r);
\draw[above](v)edge node{$v$}(vr);
\draw[above](w)edge node{$u$}(wr);
\draw[dotted](1,2)edge node{}(1,1);
\end{tikzpicture}
\caption{The set $U$ is a suffix code.}\label{figureSuffixCode}
\end{figure}
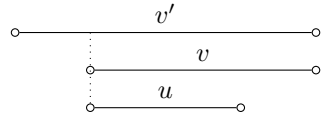

Now, since $X$ is recurrent,
 every long enough word in $\cL(X)$ has a factor equal to $u$.
Thus it has a suffix which begins with $u$ and has no other
factor equal to $u$. This suffix is in $U$. This proves the claim
concerning $U$.

Set, as in Lemma~\ref{lemmaLeftproba},
$\rho(v)=r_X(v)-1$ for every $v\in \cL(X)$.

Consider first the tree formed by the set $P$ of  prefixes of $u\RR_X(u)$. The children
of $p\in P$ are the $pa\in P$ for $a\in A$. Since $u\RR_X(u)$
is a prefix code, the leaves are the elements of $u\RR_X(u)$.
The internal nodes are the elements of $Q=P\setminus u\RR_X(u)$
 As in any finite
tree, the number of leaves minus 1 is equal to the sum over the internal nodes $v$
of the integers $d(v)-1$, where $d(v)$ is the number of children of $v$.
For $v\in Q$, since $\RR_X(u)$ is an $(X,u)$-maximal prefix code, we have 
\begin{displaymath}d(v)=\begin{cases}r_X(v)&\mbox{ if $v\in U$}\\1&\mbox{otherwise.}
\end{cases}
\end{displaymath} 
Thus
\begin{equation}
  \Card(\RR_X(u))-1=
\Card(u\RR_X(u))-1=\sum_{v\in Q}(d(v)-1)=\sum_{v\in U}\rho(v).\label{eqSommeDegres}
\end{equation}

Consider now the  tree  formed by the  set $S$ of suffixes of $U$. The root is $\varepsilon$
and the children
of a word $v\in S$ are the words $av\in S$ with $a\in A$. Since $U$ is an $X$-maximal suffix code,
the leaves of the tree $S$
are the elements of $U$ and the children of  $v\in S\setminus U$ are all the $av$ for $a\in L_X(v)$. 
Since  for every internal node of $S$, the sum
of the $\rho(av)$ taken over the children of $v$ is equal to $\rho(v)$
(by Lemma~\ref{lemmaLeftproba}), we have
\begin{equation}
\sum_{v\in U}\rho(v)=\rho(\varepsilon)=\Card(A)-1\label{eqProbaSuf}
\end{equation}
with the last equality resulting from the hypothesis $A\subset \cL(X)$.
Comparing \eqref{eqSommeDegres} and \eqref{eqProbaSuf}, we obtain the
desired equality.
\end{proofof}

We illustrate the proof with the following example.
\begin{example}
Let $X$ be the Fibonacci shift and let $u=aa$.
\index{subject}{Fibonacci!shift}\index{subject}{shift space!Fibonacci}%
 We have
$\RR_X(u)=\{baa,babaa\}$ and thus $U=\{aa,aab,aaba,aabab,aababa\}$.
The tree formed by the set $\RR_X(u)$ is represented in Figure~\ref{figureTreeT}
on the left. The tree $S$ is represented on the right
 with the value of $\rho$ indicated on the leaves.
\begin{figure}[hbt]
\centering
\tikzset{node/.style={circle,draw,minimum size=0.1cm,inner sep=0pt}}
\tikzset{leaf/.style={draw,minimum size=0.1cm}}
\begin{tikzpicture}
\node[node](1)at(-6.5,1){};\node[node](b)at(-5.5,1){};
\node[node](ba)at(-4.5,1){};\node[leaf](baa)at(-3.5,1.5){};\node[node](bab)at(-3.5,.5){};
\node[node](baba)at(-2.5,.5){};\node[leaf](babaa)at(-1.5,.5){};

\draw[above](1)edge node{$b$}(b);\draw[above](b)edge node{$a$}(ba);
\draw[above](ba)edge node{$a$}(baa);\draw[above](ba)edge node{$b$}(bab);
\draw[above](bab)edge node{$a$}(baba);\draw[above](baba)edge node{$a$}(babaa);

\node[leaf](aababa)at(-1,1){$0$};
\node[node](ababa)at(0,1){};\node[leaf](aabab)at(0,0){$0$};
\node[leaf](aaba)at(1,2){$1$};\node[node](baba)at(1,1){};\node[node](abab)at(1,0){};
\node[node](aba)at(2,1.5){};\node[leaf](aab)at(2,1){$0$};\node[node](bab)at(2,0){};
\node[leaf](aa)at(3,2.5){$0$};\node[node](ba)at(3,1.5){};\node[node](ab)at(3,1){};
\node[node](a)at(4,2){};\node[node](b)at(4,1){};
\node[node](1)at(5,1.5){};

\draw[above](aababa)edge node{$a$}(ababa);
\draw[above](ababa)edge node{$a$}(baba);\draw[above](aabab)edge node{$a$}(abab);
\draw[above](aaba)edge node{$a$}(aba);\draw[above](baba)edge node{$b$}(aba);
\draw[above](abab)edge node{$a$}(bab);
\draw[above](aba)edge node{$a$}(ba);\draw[above](aab)edge node{$a$}(ab);
\draw[above](bab)edge node{$b$}(ab);
\draw[above](aa)edge node{$a$}(a);
\draw[above](ba)edge node{$b$}(a);
\draw[above](ab)edge node{$a$}(b);
\draw[above](a)edge node{$a$}(1);
\draw[above](b)edge node{$b$}(1);

\end{tikzpicture}
\caption{The trees $\RR_X(u)$ and  $S$.}\label{figureTreeT}
\end{figure}
The unique leaf  of $S$ with a nonzero value of $\rho$ is the
unique right-special word which belongs to $U$, namely $aaba$.
\end{example}
Observe that the $X$-maximal suffix code
$U$ used in the proof is closely related with the
partition in towers built from the set $\RR_X(u)$
as in  Proposition~\ref{propositionPartitionReturn}
(see Exercise~\ref{exercisePartitionsCodes}).

We now come to the second part of the proof of Theorem~\ref{theoremReturn}.

In a graph $G=(V,E)$ labeled by an alphabet $A$,
 we consider for every edge $e$ from $v$
to $w$ with label $a$, an \emph{inverse edge}
\index{subject}{inverse!edge} $e^{-1}$ which goes from $w$ to $v$
and is labeled $a^{-1}$.
A \emph{generalized path}\index{subject}{generalized!path}
 in  $G$ is a sequence formed
of consecutive  edges or their inverses. 
The label of a generalized path is the reduced word which is the reduction of the label of the
path. Thus it is an element of the free group on $A$.

\begin{lemma}\label{lemmaRauzyGroup}
Let $X$ be a dendric shift such that $A\subset \cL(X)$. For every $n\ge 1$, the group defined
by the Rauzy graph $\Gamma_n(X)$ with respect to one of its vertices
is the free group on $A$.
\end{lemma}
\begin{proof}
We will show that a sequence of Stallings foldings reduces any
Rauzy graph $\Gamma_{n+1}(X)$ to $\Gamma_n(X)$. 

Consider two
vertices $ax,bx$ of $\Gamma_{n+1}(X)$ differing only by the first letter.
Since the extension graph of $x$ is a tree, there is a path
$a_0,b_1,\ldots,a_{k-1},b_k,a_k$ in $\E(x)$ such that $a=a_0$
and $b=a_k$. The successive Stallings foldings at $xb_1,\ldots,xb_k$
identify the vertices $a_0x,\ldots,a_kx$. In this way,
$\Gamma_{n+1}(X)$ is mapped onto $\Gamma_n(X)$. 

Thus the groups defined by the Rauzy graphs $\Gamma_n(X),\Gamma_{n-1}(X),\ldots,
\Gamma_1(X)$ are all identical. Since $A\subset \cL(X)$,
the graph $\Gamma_1(X)$ defines the free group on $A$,
and thus the same is true for $\Gamma_n(X)$.
\end{proof}
We illustrate the proof of Lemma~\ref{lemmaRauzyGroup} with the following example.
\begin{example}
Consider the Fibonacci shift
\index{subject}{Fibonacci!shift}\index{subject}{shift space!Fibonacci}%
 and the Rauzy graphs $\Gamma_n(X)$
for $n=1,2,3$ represented in Figure~\ref{figureLemmaRauzy}.
\begin{figure}[hbt]
\centering
\tikzset{node/.style={circle,draw,minimum size=0.4cm,inner sep=0.2pt}}
\tikzstyle{every loop}=[->,shorten >=1pt,looseness=12]
\tikzstyle{loop left}=[in=130,out=220,loop]
\tikzstyle{loop right}=[in=330,out=50,loop]
\begin{tikzpicture}

\node[node](11) at (0,0){$\varepsilon$};

\draw[left](11) edge[loop left]node {$a$}(11);
\draw[right](11) edge[loop right]node {$b$}(11);

\node[node](21)at(3,0){$a$};\node[node](22)at(5,0){$b$};

\draw[left](21)edge[loop left]node{$a$}(21);
\draw[above, bend left, ->](21) edge node{$b$}(22);
\draw[below, bend left, ->](22)edge node{$a$}(21);

\node[node](31)at (6,0){$aa$};
\node[node](32)at (8,.7){$ab$};
\node[node](33)at(8,-.7){$ba$};

\draw[above, bend left, ->](31)edge node{$b$}(32);
\draw[right,bend left, ->](32)edge node{$a$}(33);
\draw[below, bend left, ->](33)edge node{$a$}(31);
\draw[left, bend left, ->](33)edge node{$b$}(32);

\end{tikzpicture}
\caption{The Rauzy graphs of order $n=1,2,3$ of the Fibonacci shift.}
\label{figureLemmaRauzy}
\end{figure}
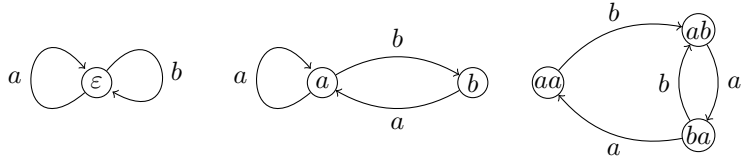
Since there are edges labeled $b$ in $\Gamma_3(X)$ from $aa$ and $ba$
  to $ab$, a Stallings folding merges $aa$ and $ba$.
The result is $\Gamma_2(X)$. Similarly, since there are edges labeled $a$
from the vertices $a$ and $b$ to vertex $a$ in $\Gamma_2(X)$,
 we merge the vertices $a$ and $b$. The result
is $\Gamma_1(X)$.
\end{example}

\begin{figure}[hbt]
\centering
\tikzset{node/.style={circle,draw,minimum size=0.1cm,inner sep=0.1pt}}
\begin{tikzpicture}
\node[node](x)at(0,1){$x$};\node[node](p)at(2,2){$p$};
\node[node](q)at(2,0){$q$};

\draw[above,bend left,->](x)edge node{$u$}(p);
\draw[above,bend left=20,->](p)edge node{$t$}(x);
\draw[right,bend right,->](q)edge node{$a$}(p);
\draw[above,bend left=20,->](x)edge node{$w$}(q);
\draw[above,bend left,->](q)edge node{$v$}(x);
\end{tikzpicture}
\caption{The decomposition of the path $\pi$.}\label{figureGroupDefined}
\end{figure}
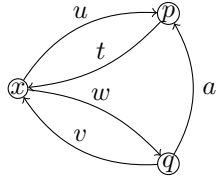

\begin{lemma}\label{lemmaGroupDefined}
Let $G$ be a  strongly connected labeled graph and $x$ be a vertex. The group defined
by $G$ with respect to $x$ is generated by the set $S$ of labels
of paths from $x$ to $x$ in $G$.
\end{lemma}
\begin{proof}
Consider a generalized path $\pi$ from $x$ to $x$ labeled $y$.
We have to prove that $y$ belongs to the subgroup $\langle S\rangle$
generated by $S$.
We use an induction on the number $r$ of inverse edges used in 
the path $\pi$. If $r=0$, then $y$ is in $S$. Otherwise, we can
write $y=ua^{-1}v$ where $x\edge{u}p\edge{a^{-1}}q\edge{v}x$
is a factorization of the path $\pi$. Since $G$ is strongly
connected, there are (ordinary) paths $p\edge{t}x$ and $x\edge{w}q$.
Then 
\begin{displaymath}
y=utt^{-1}a^{-1}w^{-1}wv=(ut)(wat)^{-1}wv.
\end{displaymath}
By definition we have $wat\in S$ (see Figure~\ref{figureGroupDefined})
and
by induction hypothesis, we have $ut,wv\in \langle S\rangle$. This shows
that $y\in \langle S\rangle$ and concludes the proof.
\end{proof}

We are now ready for the proof of Theorem~\ref{theoremReturn}.\\

\begin{proofof}{of Theorem~\ref{theoremReturn}}
 Let $n\ge 1$ be such that the property of Proposition \ref{propositionJulien}
holds for $x\in\cL_n(X)$.

 The inclusion $S\subset \RR_X(u)^*$ implies the inclusion
$\langle S\rangle\subset \langle \RR_X(u)\rangle$. But, by Lemma~\ref{lemmaGroupDefined}, $S$ generates the group
defined by $\Gamma_{n+1}(X)$. By Lemma~\ref{lemmaRauzyGroup}, this group is the whole free group on $A$. Thus $\RR_X(u)$ generates the free group on $A$. Since
any generating set of $F(A)$ having $\Card(A)$
elements is a basis, and since $\RR_X(u)$ has $\Card(A)$ elements
by Theorem~\ref{theoremCardinality},
this implies our conclusion.

\end{proofof}
\begin{example}\label{exampleReturnTheorem}
Let $A=\{u,v,w\}$, and let $X=X(\sigma)$ be the shift generated
by the substitution $\sigma:u\to vuwwv,v\to vuww,w\to vuwv$. The shift
$X$ is minimal since $\sigma$ is primitive. It is also dendric.
Consider indeed the morphism $\phi:u\to aa,v\to ab,w\to ba$. Then
we have $\phi\circ \sigma=\varphi^3\circ\phi$ where $\varphi$
is the Fibonacci morphism. Indeed, we have 
\begin{displaymath}
\phi\circ\sigma(u)=\phi(vuwwv)=abaababaab=\varphi^3(aa)
\end{displaymath}
and similarly for $v,w$. This shows that $X$ is obtained by
reading the Fibonacci shift with non-overlapping blocks of length $2$
and thus that $X$ is dendric (this is actually a particular case of
Theorem~\ref{theoremBifixDecoding}). We have
\begin{displaymath}
\RR_X(u)=\{wwvu,wwvvu,wvvu\}
\end{displaymath}
which is a basis of the free group on $\{u,v,w\}$. Note that
$\phi(\RR_X(u))$ is a basis of a subgroup of index $2$
of the free group on $\{a,b\}$.
\end{example}
\subsection{Derivatives of minimal dendric shifts}

Let $X$ be a minimal shift space and let $u\in\cL(X)$. Let $\varphi$
be a bijection from an alphabet $B$ onto the set $\RR_X(u)$
extended as usual to a morphism from $B^*$ into $A^*$.
The shift space $Y=\varphi^{-1}(X)$ is called the \emph{derivative}
\index{subject}{derivative!of dendric shift}%
of $X$ with respect to $u$. Actually, $Y$ is the
derivative system of $X$ on the clopen set $[u]$
(see Section~\ref{sectionInduced}).
We will use (in the proof of Theorem~\ref{theoremSadicDendric}) the following closure property
of the family of minimal dendric shifts.

\begin{theorem}\label{propositionReturns}
Any derivative  of a minimal dendric shift
 is a minimal dendric shift on the same number of letters.
\end{theorem}
\begin{proof}
Let $X$ be a minimal dendric shift on the alphabet $A$
containing 
$A$, let $u\in \cL(X)$
and let $\varphi$ be a bijection from an alphabet $B$ onto
 $U=\RR_X(u)$.
By Theorem~\ref{theoremCardinality}, the set $\RR_X(u)$ has $\Card(A)$ elements.
Thus we may choose $B=A$.

Set $Y=\varphi^{-1}(X)$. Since $Y$ is an induced system,
 it is
minimal.

Consider $y\in \cL(Y)$ and set $x=\varphi(y)$. Let
$\varphi'$
be the bijection from $A$ onto $U'=\RR'_X(u)$ such that $u\varphi(b)=\varphi'(b)u$ for every $b\in B$.
For $a,b\in B$, we have
\begin{displaymath}
(a,b)\in \E(y)\Leftrightarrow (\varphi'(a),\varphi(b))\in \E_{U',U}(ux),
\end{displaymath}
where $\E_{U',U}(ux)$ denotes the generalized extension graph of $ux$
relative to $U',U$.
Indeed, 
\begin{displaymath}
ayb\in \cL(Y)
\Leftrightarrow u\varphi(a)x\varphi(b)\in \cL(X)
\Leftrightarrow
\varphi'(a)ux\varphi(b)\in \cL(X).
\end{displaymath}
The set $U'$ is a $(X,u)$-maximal suffix code and the set $U$ is a
$(X,u)$-maximal prefix code. By
Proposition~\ref{propStrongTreeConditionBis}
the generalized extension graph $\E_{U',U}(ux)$ is a tree. Thus the
graph $\E(y)$
is a tree. This shows that $Y$ is a dendric shift.
\end{proof}

\begin{example}\label{exampleDerivativeTribo}
Let $X$ be the Tribonacci shift (see Example~\ref{ex:Tribo}).
\index{subject}{Tribonacci!shift}\index{subject}{shift space!Tribonacci}%
It is the shift generated by the substitution $\sigma$ defined by $\sigma(a)=ab$,
$\sigma(b)=ac$, $\sigma(c)=a$.
We have $\RR_X(a)=\{a,ba,ca\}$.  Let  $\varphi:A\to \RR_X(a)$ be the morphism
 defined by  $\varphi(a)=a$, $\varphi(b)=ba$,
$\varphi(c)=ca$ and let $\varphi':A\to\RR'_X(a)$ be such that
$a\varphi(x)=\varphi'(x)a$ for all $x\in A$. We have $\sigma=\varphi'\circ \pi$ where
$\pi$ is the circular permutation $\pi=(abc)$.
Let $x=\varphi^\omega(a)$.
Set  $z=\varphi'^{-1}(x)$. Since $\varphi'\pi(x)=x$, we have $z=\pi(x)$. 
Thus the derivative  of $X$ with respect
to $a$ is the shift $\pi(X)$.
\end{example}
\subsection{Bifix codes in dendric shifts}
A \emph{bifix code}
\index{subject}{bifix code}\index{subject}{code!bifix}%
on the alphabet $A$
is a set $U$ of words on $A$ which is both a prefix code and a suffix code.
For example, for every $n\ge 1$, a set of words of length $n$ is a bifix code.

%A bifix code $U\subset \cL(X)$ is $X$-\emph{maximal}

Let $X$ be a  shift space and let $U\subset \cL(X)$ be a finite
bifix code which is an $X$-maximal prefix and
suffix code. Let $f$ be a coding morphism for $U$. 
Then $f^{-1}(\cL(X))$ is factorial and extendable.
The shift space
$Y$ such that  $\cL(Y)=f^{-1}(\cL(X))$ 
is called the \emph{decoding}\index{subject}{decoding of a shift}
of $X$ by $U$. We denote $Y=f^{-1}(X)$.

As a particular case, for $n\ge 1$, the set $U=\cL_n(X)$ is a bifix code
which is both an $X$-maximal prefix code and an $X$-maximal
suffix code. The corresponding decoding $Y=f^{-1}(X)$ is called
the coding of $X$ by \emph{non overlapping blocks} of length $n$.
\index{subject}{coding!by non overlapping blocks}%
The map $f$ is then a conjugacy from $Y$ onto the system
$(X,S^n)$ (see Example~\ref{exampleFiboDecoding}).

The following result expresses a closure property of the family of dendric shifts.
\begin{theorem}\label{theoremBifixDecoding}
The decoding of a dendric shift by a finite bifix code which
is an $X$-maximal prefix and suffix code is a dendric shift.
\end{theorem}
\begin{proof}
Let $U\subset\cL(X)$ be a finite bifix code which is an $X$-maximal
prefix and suffix code. Let $f:B^*\to A^*$ be a coding morphism for $U$
and let $Y=f^{-1}(X)$. For $w\in\cL(Y)$ and $a,b\in B$, we have
\begin{displaymath}
(a,b)\in\E_Y(w)\Leftrightarrow (f(a),f(b))\in \E_{U,U}(f(w))
\end{displaymath}
and thus $E_Y(w)$ is a tree by Proposition~\ref{propStrongTreeConditionBis}.
This shows that $Y$ is dendric.
\end{proof}

\begin{example}\label{exampleFiboDecoding}
Let $X$ be the Fibonacci shift on $\{a,b\}$ and consider the bifix code
$U=\{aa,ab,ba\}$. The corresponding decoding of $X$ is the shift
of Example~\ref{exampleReturnTheorem}.
\end{example}

We will prove  the following  result.

\begin{theorem}\label{theoremBifixBasisFreeGroup}
Let $X$ be a dendric shift on the alphabet $A$. A finite
bifix code $U\subset\cL(X)$ which
is a basis of the free group on $A$ is equal to $A$.
\end{theorem}
%Let $X$ be a shift space. A bifix code $U\subset \cL(X)$ is $X$-\emph{maximal}
%\index{subject}{maximal bifix code}\index{subject}{bifix code!maximal}%
%if it is not properly contained in another bifix code $V\subset \cL(X)$.
%\begin{theorem}\label{theoremFiniteIndexBasis}
%Let $X$ be a dendric shift. A finite bifix code is $X$-maximal 
%if and only if it is a basis of a subgroup of finite index
%of the free group on $A$.
%\end{theorem}
The following example shows that the hypothesis that $X$ is dendric
is necessary in Theorem~\ref{theoremBifixBasisFreeGroup}.

\begin{example}\label{exampleWen1}
Let $A=\{a,b,c\}$ and $U=\{ab,acb,acc\}$. The set $U$ is a bifix code. It
is also a basis of the free group on $A$. Indeed, we have
$accb=(acb)(ab)^{-1}(acb)$ and $b=(acc)^{-1}(accb)$. 
Next, $a=(ab)b^{-1}$ and $c=a^{-1}(acb)b^{-1}$. Thus $a,b,c$ belong
to the group generated by $U$. Observe that we can verify directly that no dendric shift $X$
can be such that $U\subset \cL(X)$. Indeed, this would force
$ab,cb,cc,ac\in\cL(X)$ and thus
the extension graph of $\varepsilon$ to contain a cycle (see Figure~\ref{figureWen}).
\end{example}
%Note that Theorem~\ref{theoremFiniteIndexBasis} gives a combinatorial
%interpretation of the fact that the complexity of a dendric shift $X$
%is $p_n(X)=(\Card(A)-1)n+1$ (Equation~\eqref{eqComplexityDendric}).
 %Indeed, let $U=\cL_n(X)$. Then
%$U$ is clearly an $X$-maximal bifix code. By 
%Theorem~\ref{theoremFiniteIndexBasis}, it is a basis of a subgroup of
%finite index $k$. By \emph{Schreier's Formula}
%\index{subject}{Schreier's Formula}\index{subject}{Formula!Schreier}%
%\index{names}{Schreier, Oscar}%
%this implies that $\Card(U)=(\Card(A)-1)k+1$. By \eqref{eqComplexityDendric}
%we have $k=n$. Thus $U$ is a basis of the kernel of the morphism
%$a\to 0$ from the free group onto $\Z/n\Z$.

To prove Theorem \ref{theoremBifixBasisFreeGroup},
%and \ref{theoremFiniteIndexBasis}, 
we introduce the following notion. Let $U$ be a bifix code and let $P$ (resp. $S$)
be the set of  proper prefixes (resp. suffixes)
of the words of $U$. The \emph{incidence graph}\index{subject}{incidence!graph
  of bifix code}\index{subject}{graph!incidence of bifix code}%
\index{subject}{bifix code!incidence graph of}%
  of $U$ is the
following undirected graph. Its set of vertices is the disjoint union
of $P$ and $S$. The edges are $(\varepsilon,\varepsilon)$
and the pairs $(p,s)\in P\times S$ such that
$ps\in U$.

\begin{example}
Let $X$ be the Fibonacci shift
\index{subject}{Fibonacci!shift}\index{subject}{shift space!Fibonacci}%
 and let $U=\cL_3(X)$. The incidence
graph of $U$ is represented in Figure~\ref{figureIncidenceGraph}
(in each of the three parts, the vertices on the left
are in $P$ and those on the right in $S$).
\begin{figure}[hbt]
\centering
\tikzset{node/.style={circle,draw,minimum size=0.4cm,inner sep=0.3pt}}
\begin{tikzpicture}
\node[node](1l)at(-3,1){$\varepsilon$};\node[node](1r)at(-2,1){$\varepsilon$};
\draw(1l)edge node{}(1r);

\node[node](al)at(0,1.4){$a$};\node[node](bl)at(0,.8){$b$};
\node[node](bar)at(1,1.7){$ba$};\node[node](abr)at(1,1.1){$ab$};
\node[node](aar)at(1,.5){$aa$};

\draw(al)edge node{}(bar);\draw(al)edge node{}(abr);
\draw(bl)edge node{}(abr);\draw(bl)edge node{}(aar);

\node[node](aal)at(3,1.7){$aa$};\node[node](bal)at(3,1.1){$ba$};
\node[node](abl)at(3,.5){$ab$};
\node[node](br)at(4,1.4){$b$};\node[node](ar)at(4,.8){$a$};

\draw(aal)edge node{}(br);\draw(bal)edge node{}(br);
\draw(bal)edge node{}(ar);\draw(abl)edge node{}(ar);
\end{tikzpicture}
\caption{The incidence graph of $U=\cL_3(X)$.}\label{figureIncidenceGraph}
\end{figure}
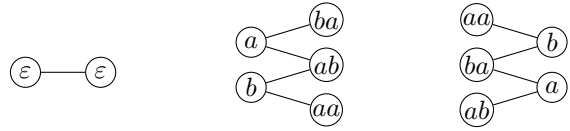
\end{example}

\begin{proposition}\label{newLemma633}
Let $X$ be an dendric
shift and let $U\subset \cL(X)$ be a 
bifix code. Let $P$ (resp. $S$) be the set of proper prefixes
(resp. proper suffixes) of $U$  and let $G$ be the incidence graph of $U$. Then the
following
assertions hold.
\begin{enumerate}
\item[\rm(i)]
The graph $G$ is
acyclic.
\item[\rm(ii)]
The intersection of $P'=P\setminus\{\varepsilon\}$ (resp. $S'=S\setminus\{\varepsilon\}$) with each connected
component
of $G$ is a suffix (resp. prefix) code.
\item[\rm(iii)]
 For every reduced path $(v_1,u_1,\ldots,u_n,v_{n+1})$ in $G$ with
$u_1,\ldots,u_n\in P'$ and $v_1,\ldots,v_{n+1}$ in $S'$,
the longest common prefix of $v_1,v_{n+1}$ is a proper prefix of all
$v_1,\ldots,v_n,v_{n+1}$. 
\item[\rm(iv)]
Symmetrically, for every reduced path $(u_1,v_1,\ldots,v_n,u_{n+1})$ in $G$
with $u_1,\ldots,$ $u_{n+1}\in P'$ and $v_1,\ldots,v_n\in S'$,
the longest common suffix of $u_1,u_{n+1}$ is a proper
suffix of $u_1,u_2,\ldots,u_{n+1}$.
\end{enumerate}
\end{proposition}

\begin{proof}
Assertions (iii) and (iv) imply Assertions (i) and (ii). Indeed, assume that
(iii) holds. Consider
a  reduced path $(v_1,u_1,\ldots,u_n,v_{n+1})$ in $G$ with
$u_1,\ldots,u_n\in P'$ and $v_1,\ldots,v_{n+1}$ in $S'$. If
$v_1=v_{n+1}$,
then $v_1$ is a prefix of all $v_i$ and in particular of
$v_2$, a contradiction since $U$
is a bifix code. Thus $G$ is acyclic and (i) holds. Next, if
$v_1$, $v_{n+1}$ are comparable for the prefix order, their longest
common prefix is one of them, a contradiction with (iii) again.
The assertion on $P'$ is proved in an analogous way using  assertion (iv).

We prove  simultaneously (iii) and (iv)
by induction on $n\ge 1$.

The assertions holds for $n=1$. Indeed, if $u_1v_1,u_1v_2\in U$
and if $v_1\in \cL(X) $ is a prefix of
$v_2\in S'$,
then $u_1v_1$ is a prefix of $u_1v_2$, a contradiction with the
hypothesis
that $U$ is a prefix code. The same
holds symmetrically for $u_1v_1,u_2v_1\in U$ since $U$ is a suffix code.

Let $n\ge 2$ and assume that the assertions hold for any path
of length at most $2n-2$. We treat the case of a path
$(v_1,u_1,\ldots,u_n,v_{n+1})$ in $G$ with
$u_1,\ldots,u_n\in P'$ and $v_1,\ldots,v_{n+1}$ in $S'$. The other
case is symmetric.

Let $p$ be the longest common prefix of $v_1$ and $v_{n+1}$.
We may assume that $p$ is nonempty since otherwise the
statement is obviously true.
  Any two elements of the
set $L=\{u_1,\ldots,u_n\}$ are connected by a path of length
at most $2n-2$ (using elements of $\{v_2,\ldots v_n\}$).
Thus, by induction hypothesis, 
$L$ is a suffix code. Similarly, any two elements of the set 
$R=\{v_1,\ldots,v_n\}$ are connected by a path of length
at most $2n-2$ (using elements of $\{u_1,\ldots u_{n-1}\}$).
Thus $R$ is a prefix code.  We cannot have $v_1=p$ since
otherwise, using the fact that $u_np$ is a prefix of $u_nv_{n+1}$
and thus in $\cL(X)$, the generalized extension
graph $\E_{L,R}(\varepsilon)$ would have
the cycle $(p,u_1,v_2,\ldots,u_n,p)$, a contradiction since $\E_{L,R}(\varepsilon)$
is acyclic by Proposition~\ref{PropStrongTreeCondition}. Similarly, we
cannot
have $v_{n+1}=p$.

Set $W = p^{-1}R$ and $R' = (R \setminus pW) \cup p$. Since
$R$ is a prefix code and since
$p$ is a proper prefix of $R$, the set $R'$ is
a prefix code.
Suppose that $p$ is not a proper prefix of all $v_2,\dots,v_n$. 
Then there exist $i,j$ with $1\le i<j\le n+1$ such that
$p$ is a proper prefix of $v_i,v_j$ but not of any $v_{i+1},\ldots,v_{j-1}$.
Then $v_{i+1},\dots,v_{j-1} \in R'$ and there is the cycle 
$(p,u_i,v_{i+1},u_{i+1},\dots,v_{j-1},u_{j-1},p)$ in the graph
$\E_{U,V'}(\varepsilon)$. 
This is in contradiction with Proposition~\ref{PropStrongTreeCondition} because, 
$V'$ being a prefix code, $\E_{L,R'}(\varepsilon)$ is acyclic. 
Thus $p$ is a proper prefix of all $v_2,\dots,v_n$.
\end{proof}
Let $X$ be a dendric shift and let $U\subset\cL(X)$ be a bifix code.
Let $P$ be the set of proper  prefixes of the words of $U$.
Let $\theta_U$\index{symbols}{theta@$\theta_U$} be the equivalence on $P$ 
defined by $p\equiv q\bmod\theta_U$ if $p,q$ are in the
same connected component of the incidence graph of $U$.
Note that, since $U$ is bifix, the class of $\varepsilon$ is reduced
to $\varepsilon$.
The \emph{coset graph}\index{subject}{coset!graph}\index{subject}{graph!coset}
of $U$ is
the following labeled graph. The set vertices is the set $R$
 of  classes of $\theta_U$.
 There is an edge labeled $a$ from the class of $p$ to the class of $q$ 
in each of the following cases
\begin{enumerate}
\item[(i)]   $q=pa$,
\item[(ii)] $q=\varepsilon$ and  $pa\in U$.
\end{enumerate}

\begin{example}
Let $X$ be the Fibonacci shift. The coset
graph of $\{a,bab\}$ is shown in Figure~\ref{figureCosetAutomaton}
on the left and the coset graph of $\{a,bab,baab\}$ on the right.
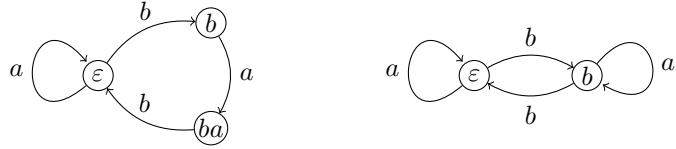
\begin{figure}[hbt]
\centering
\tikzset{node/.style={circle,draw,minimum size=0.4cm,inner sep=0.3pt}}
\tikzstyle{every loop}=[->,shorten >=1pt,looseness=12]
\tikzstyle{loop left}=[in=130,out=220,loop]
\tikzstyle{loop right}=[in=330,out=50,loop]
\begin{tikzpicture}
\node[node](1)at(0,0){$\varepsilon$};\node[node](b)at(1.5,.7){$b$};
\node[node](ba)at(1.5,-.7){$ba$};

\draw[left,->](1)edge[loop left] node{$a$}(1);
\draw[above,->,bend left](1)edge node{$b$}(b);
\draw[right,->,bend left](b)edge node{$a$}(ba);
\draw[above,->,bend left](ba)edge node{$b$}(1);

\node[node](1)at(5,0){$\varepsilon$};\node[node](b)at(6.5,0){$b$};

\draw[left,->](1)edge[loop left] node{$a$}(1);
\draw[above,->,bend left](1)edge node{$b$}(b);
\draw[below,->,bend left](b)edge node{$b$}(1);
\draw[right,->](b)edge[loop right] node{$a$}(b);
\end{tikzpicture}
\caption{The coset graphs of $\{a,bab\}$ and of $\{a,bab,baab\}$.}\label{figureCosetAutomaton}
\end{figure}
\end{example}
A \emph{simple path}\index{subject}{simple!path}\index{subject}{path!simple}
from a vertex $v$ to itself in a graph is a path which is
not a concatenation of two nonempty paths from $v$ to itself.
%A \emph{right coset} \index{subject}{right!coset} of a subgroup $H$
%in a group $G$ is a set of the form $Hg=\{hg\mid h\in H\}$ for
%some $g\in G$. Two cosets are equal or disjoint.
\begin{proposition}\label{propositionCosetGraph}
Let $X$ be a dendric shift and let $U\subset\cL(X)$
 be a finite bifix code. Let $P$ be the set of
proper prefixes of $U$ and let $H=\langle U\rangle$ \index{symbols}{U@$\langle U\rangle$}
be the subgroup generated by $U$. Let also $C$
be the coset graph of $U$.
\begin{enumerate}
\item  For every $p,q\in P$, $p\equiv q\bmod \theta_U$ implies $Hp=Hq$.
\item If $p\equiv p'\bmod\theta_U$ and if $p\edge{a}q$, $p'\edge{a}q'$
are edges in $C$, then $q\equiv q'\bmod\theta_U$.
\item Every $u\in U$ is the label of a simple path from $\varepsilon$ to
itself in $C$.
\item The graph $C$ is the
Stallings graph of the subgroup $\langle U\rangle$ generated by $U$.
\end{enumerate}
\end{proposition}
\begin{proof}
Let $C=(R,E)$ be the coset graph of $U$ and let $G$ be its incidence graph.

1. The first assertion is clear since $\theta_U$ is the equivalence
on $P$ generated by the pairs $p,q$ such that there is an $s$ with
$ps,qs\in U$ and thus $p,q\in Hs^{-1}$.

2. We assume that $q=pa$ and $q'=p'a$. The other cases are similar.
Let $s,s'$ be such that $qs,q's'\in U$. 
Let $p=u_0,v_1,u_1,\ldots,v_n,u_n=p'$ be a path from $p$
to $p'$ in the incidence graph $G$. Set $v_0=as$ and $v_{n+1}=as'$.
Then $(v_0,u_0,\ldots,u_n,v_{n+1})$ is a path in $G$. But since the letter
$a$ is a common prefix of $v_0$ and $v_{n+1}$, by Proposition~\ref{newLemma633}
it is also a common prefix of all $v_i$. Set $v_i=av'_{i}$ for 
$0\le i\le n+1$. Then $(u_0a,v'_1,u_1a,\ldots,v'_n,u_na)$
is a path from $q$ to $q'$ in the coset graph $C$ and thus $q\equiv q'\bmod\theta_U$. 

3. This follows from the fact that $C$ can be obtained by Stallings
foldings from the graph on $P$ with edges $p\edge{a}q$
if either $pa=q$ or $q=\varepsilon$ and $pa\in U$.

4. Let  $K$ be the group
defined by the coset graph $C$. Let us show that $K$ is equal to $H$. By construction,
we have $U\subset K$ and thus $H\subset K$. The converse
follows easily from Assertion 1.

Let us finally show that $C$ is Stallings reduced.
Assume that $p,q\in P$ are such that there are edges with the same label $a$
from the class $\bar{p},\bar{q}$ of $p,q$
to the same vertex $\bar{r}$. 
Let $v$ be the label of a path from $\bar{r}$ to $\varepsilon$ which does not
pass by $\varepsilon$ before. 
Then $pav,qav\in U$ and thus $p\equiv q\bmod\theta_U$,
which implies that $\bar{p}=\bar{q}$. This shows that $C$ is Stallings reduced.
\end{proof}
We are now ready to prove Theorem~\ref{theoremBifixBasisFreeGroup}.
%and \ref{theoremFiniteIndexBasis}.

\begin{proofof}{of Theorem~\ref{theoremBifixBasisFreeGroup}}
Let $X$ be a dendric shift on the alphabet $A$.
Let $U\subset \cL(X)$ be a bifix code which is a basis of the
free group on $A$. By Proposition~\ref{propositionCosetGraph},
the coset graph $C$ of $U$
has only one vertex $\varepsilon$ and loops $\varepsilon\edge{a}\varepsilon$
for every $a\in A$. Since, by Proposition~\ref{propositionCosetGraph}
again, every word of $U$ is the
label of a simple path from $\varepsilon$
to itself in the coset graph of $U$, we have $U\subset A$ and thus $U=A$.
\end{proofof}

%\begin{proofof}{of Theorem~\ref{theoremFiniteIndexBasis}}
%Let $C$ be the coset graph of $U$. Let $V$ be the set
%of labels of simple paths from $\varepsilon$ to itself in $C$.
%It is a bifix code which, by Proposition~\ref{propositionCosetGraph}
%contains $U$. This implies that $U$ is an $X$-maximal prefix code.
%Indeed, if $w\in\cL(X)$ is such that $U\cup\{w\}$ is a prefix code,
%then $w$ 
%\end{proofof}

\subsection{Tame automorphisms}
An automorphism $\alpha$ of the free group on $A$ is \emph{positive}
\index{subject}{automorphism!positive}\index{subject}{positive!automorphism}%
if $\alpha(a)\in A^+$ for every $a\in A$.
We say that a positive automorphism of the free group on $A$ is \emph{tame}
\index{subject}{tame!automorphism}\index{subject}{automorphism!tame}%
if it belongs to the submonoid generated by the permutations of $A$ and
the automorphisms $\alpha_{a,b}$, $\tilde{\alpha}_{a,b}$ defined 
for $a,b\in A$ with $a\ne b$ by
\begin{displaymath}
\alpha_{a,b}(c)=\begin{cases}ab&\text{ if $c=a$,}\\c&\text{
    otherwise} \end{cases}
\quad
\text{ and }
\quad
\tilde{\alpha}_{a,b}(c)=\begin{cases}ba&\text{ if $c=a$,}\\c&\text{
    otherwise.} \end{cases}
\end{displaymath}
\index{symbols}{alpha@$\alpha_{a,b}$}%
Thus $\alpha_{a,b}$ places a letter $b$ after each $a$ and
 $\tilde{\alpha}_{a,b}$ places a letter $b$ before each $a$.
The above automorphisms and the permutations of $A$
are called the \emph{elementary}
\index{subject}{elementary!automorphism}\index{subject}{automorphism!elementary}%
positive automorphisms on $A$.
The monoid of positive automorphisms is not finitely generated
as soon as the alphabet has at least three generators
(see the Notes Section).

A basis $U$ of the free group is \emph{positive}\index{subject}{positive!basis}
 if $U\subset A^+$.
A positive basis $U$ of the free group is \emph{tame}\index{subject}{tame!basis} if there exists
a tame automorphism $\alpha$ such that $U=\alpha(A)$.

\begin{example}
The set $U=\{ba,cba,cca\}$ is a tame basis of the free group on $\{a,b,c\}$.
Indeed, one has the following sequence of elementary automorphisms.
\begin{displaymath}
(b,c,a)\edge{\alpha_{c,b}}(b,cb,a)\edge{\tilde{\alpha}_{a,c}^2}(b,cb,cca)
\edge{\alpha_{b,a}}(ba,cba,cca).
\end{displaymath}
The fact that $U$ is a basis can be checked directly since
$(cba)(ba)^{-1}=c$, $c^{-2}(cca)=a$ and finally $(ba)a^{-1}=b$.
\end{example}
The following result will play a key role in the proof of the
main result of this section (Theorem~\ref{theoremTame}).
\begin{proposition}\label{propAuxiliary}
A set $U\subset A^+$ is a tame basis of the free group on $A$
if and only if $U=A$ or there is a tame basis $V$ of the free
group on $A$ and $u,v\in V$ such that
$U=(V\setminus v)\cup uv$ or $U=(V\setminus u)\cup uv$.
\end{proposition}
\begin{proof}
Assume first that $U$ is a tame basis of the free group on $A$.
Then $U=\alpha(A)$ where $\alpha$ is a tame automorphism
of $\langle A\rangle$. Then $\alpha=\alpha_1\alpha_2\cdots\alpha_n$
where the $\alpha_i$ are elementary positive automorphisms.
We use an induction on $n$. If $n=0$, then $U=A$. 
If $\alpha_n$ is a permutation of $A$, then
$U=\alpha_1\alpha_2\cdots\alpha_{n-1}(A)$ and the result holds by
induction
hypothesis. Otherwise, set $\beta=\alpha_1\cdots\alpha_{n-1}$ and
 $V=\beta(A)$. By induction hypothesis, $V$ is tame.
If $\alpha_n=\alpha_{a,b}$, set $u=\beta(a)$ and $v=\beta(b)=\alpha(b)$.
Then 
\begin{eqnarray*}
U&=&\alpha(A\setminus a)\cup\alpha(a)=\beta(A\setminus a)\cup\beta(ab)\\
&=&(V\setminus u)\cup uv
\end{eqnarray*} 
and thus the condition is satisfied.
The case were $\alpha_n=\tilde{\alpha}_{a,b}$ is symmetrical.

Conversely, assume that $V$ is a tame basis and that $u,v\in V$
are such that $U=(V\setminus u)\cup uv$. Then, there is a tame automorphism
$\beta$ of $F(A)$ such that $V=\beta(A)$. Set $a=\beta^{-1}(u)$
and $b=\beta^{-1}(v)$. Then $U=\beta\circ\alpha_{a,b}(A)$ and thus $U$ is a tame
basis.
\end{proof}
We note the following corollary. 
\begin{corollary}\label{corollaryTame}
A tame basis of the free group which is a bifix code is the alphabet.
\end{corollary}
\begin{proof}
Assume that $U$ is a tame basis which is not the alphabet.
By Proposition~\ref{propAuxiliary} there is a tame basis $V$ and $u,v\in V$
such that $U=(V\setminus v)\cup uv$ or  $U=(V\setminus u)\cup uv$.
In the first case, $U$ is not prefix. In the second one, it is not suffix.
\end{proof}

\begin{example}\label{exampleWen}
The set $U=\{ab,acb,acc\}$ is a basis of the free group
on $\{a,b,c\}$ (see Example~\ref{exampleWen1}).
 The set $U$ is bifix and thus it is not a tame basis by
Corollary~\ref{corollaryTame}.
\end{example}
The following result is a remarkable consequence of Theorem~\ref{theoremBifixBasisFreeGroup}.
\begin{theorem}\label{theoremTame}
Any basis of the free group included in the language 
of a minimal dendric shift is tame.
\end{theorem}
\begin{proof}
Let $X$ be a minimal dendric shift. Let $U\subset \cL(X)$
be a basis of the free group on $A$.  We use an induction
on the sum $\lambda(U)$ of the lengths of the words of $U$. If $U$
is bifix, by Theorem~\ref{theoremBifixBasisFreeGroup}, we have $U=A$.
Next assume, for example that $U$ is not prefix. Then there
are nonempty words $u,v$ such that $u,uv\in U$. Let
$V=(U\setminus uv)\cup v$, so that $V\subset\cL(X)$. Then $V$ is a basis of the free
group and $\lambda(V)<\lambda(U)$. By induction
hypothesis, $V$ is tame. Since $U=(V\setminus v)\cup uv$,
$U$ is tame
by Proposition~\ref{propAuxiliary}.
\end{proof}
\begin{example}
The set $U=\{ab,acb,acc\}$ is a basis of the free group which is not
tame (see Example~\ref{exampleWen}). Accordingly, the extension
graph $\E(\varepsilon)$ relative to the set of factors
of $U$ is not a tree (see Figure~\ref{figureWen}).
\begin{figure}[hbt]
\centering
\tikzset{node/.style={circle,draw,minimum size=0.4cm,inner sep=0.3pt}}
\begin{tikzpicture}(20,10)
\node[node](al)at(0,1){$a$};\node[node](cl)at(0,0){$c$};
\node[node](br)at(2,1){$b$};\node[node](cr)at(2,0){$c$};
\draw(al)edge node{}(br);\draw(al)edge node{}(cr);
\draw(cl)edge node{}(br);\draw(cl)edge node{}(cr);
\end{tikzpicture}
\caption{The graph $\E(\varepsilon)$.}\label{figureWen}
\end{figure}
\end{example}
%%%%%%%%%%%%%%%%%%
%\subsection{$S$-adic representations}\label{sectionSadic}
\subsection{$\Sa$-adic representation of dendric shifts}
We now study the $\Sa$-adic representations of dendric shifts.
We first recall a general construction allowing to build $S$-adic representations of any minimal aperiodic shift (Proposition~\ref{prop: S-adic UR set}) which is based on return words.
Using Theorem~\ref{theoremTame}, we show that this construction actually 
provides $\mathcal{\Sa}_e$-representations of minimal dendric
shifts (Theorem~\ref{theoremSadicDendric}), where $\mathcal{\Sa}_e$ is the set of elementary positive automorphisms
 of the free group on $A$.

%\subsection{$S$-adic representations}
Let $\Sa$ be a set of morphisms and $\tau = (\tau_n)_{n \ge 1}$ 
be a directive sequence of morphisms in $\Sa$ with $\tau_n:A_{n+1}^* \to A_n^*$ and $A_1=A$. When $\tau$ is primitive and proper, 
we have by Lemma~\ref{lemma:omega}
 \begin{equation}
\cL(\tau)=\bigcap_{n \ge 1} \Fac(\tau_{[0,n]}(A_{n+1}^*))
\label{eqSadic}
\end{equation}
where $\Fac(L)$ denotes the set of factors of the words in $L$.

%A minimal shift space $X$ is said to be {\em branching}
%\index{subject}{branching!shift space}\index{subject}{shift space!branching}%
 %if $\cL(X)$ contains at least one right-special factor of 
%each length. Equivalently, $X$ is branching if its factor complexity
%$p_n(X)$ is such that $p_{n+1}(X)-p_n(X)\ge 1$ for all $n\ge 1$.
%A branching shift is aperiodic.

The next  proposition provides a general construction to get a primitive proper $\Sa$-adic representation 
of any aperiodic minimal shift space $X$.

\begin{proposition}
\label{prop: S-adic UR set}
An aperiodic shift $X$ is minimal if and only if it has a primitive 
proper $\Sa$-adic 
representation for some (possibly infinite) set $\Sa$ of 
morphisms.
\end{proposition}
\begin{proof}
The direct implication follows from Proposition~\ref{prop:minimalSSIprimitif}.

Let us prove the converse.
Since $X$ is aperiodic, we have $p_{n+1}>p_n(X)$ for all $n\ge 1$
by Theorem~\ref{theoremCovenHedlund}.
Thus there is for each $n\ge 1$,  a 
 right-special word $u_n$ of length $n$ in $\cL(X)$
 such that $u_n$ is a suffix of $u_{n+1}$.
By assumption, $X$ is minimal so that
 $\mathcal{R}_X(u_{n+1})$ is finite for all $n$. 
Since $u_{n+1}$ is right-special, $\mathcal{R}_X(u_{n+1})$ has cardinality at least 2 for all $n$.
For all $n$, let $A_n = \{0,\dots,\Card(\mathcal{R}_X(u_n))-1\}$ and let $\alpha_n: A_n^* \to A^*$ be a 
coding morphism for $\mathcal{R}_X(u_n)$.
The word $u_n$ being suffix of $u_{n+1}$, we have $\alpha_{n+1}(A_{n+1}) \subset \alpha_n(A_n^+)$.
Since $\alpha_n(A_n) = \mathcal{R}_X(u_n)$ is a prefix code, there is a unique 
morphism $\tau_n: A_{n+1}^* \to A_n^*$ such that 
$\alpha_n \circ\tau_n = \alpha_{n+1}$.
For all $n$ we get 
$\mathcal{R}_X(u_n) = \alpha_0 \circ\tau_0 \circ\tau_1 \circ\cdots\circ \tau_{n-1}(A_n)$ 
and $\cL(X) = \bigcap_{n \in \N} \Fac(\alpha_0 \circ\tau_0 \circ\cdots\circ \tau_n(A_{n+1}^*))$.
Without loss of generality, we can suppose that $u_0 = \varepsilon$ and $A_0 = A$.
In this case we get $\alpha_0 = \id$ and the shift space $X$ thus has an $\Sa$-adic 
representation with $\Sa = \{\tau_n \mid n \in \N\}$.

By construction, each morphism $\sigma_n$ is right proper. By Lemma
\ref{lemma:proper}, we can modify the morphisms $\sigma_n$
to make them proper.

Finally, by Proposition~\ref{prop:minimalSSIprimitif} again, the $\Sa$-adic representation
is primitive.
\end{proof}

%\subsection{$S$-adic representation of tree sets}
Even for minimal shifts with linear factor complexity, the set of morphisms $\Sa=\{\tau_n \mid n \ge 0\}$ 
considered in Proposition~\ref{prop: S-adic UR set} is usually  infinite as well as the sequence of alphabets 
$(A_n)_{n \ge 0}$ is usually  unbounded.
% (see~\cite{Durand-Leroy-Richomme}).
For dendric shifts, the next theorem significantly improves the only if part of 
Proposition~\ref{prop: S-adic UR set}. Indeed, for such sets, the set $\Sa$ can be replaced by the set $\mathcal{S}_e$ 
of elementary positive automorphisms. In particular, $A_n$ is equal to $A$ for all $n$.
\index{symbols}{Se@$\Sa_e$}%

\begin{theorem}\label{theoremSadicDendric}
Every minimal dendric shift  has 
\begin{enumerate}
\item a primitive proper unimodular $\Sa$-adic representation and also
\item a primitive $\mathcal{S}_e$-adic representation.
\end{enumerate}
\end{theorem}

\begin{proof}
For any non-ultimately periodic sequence
 $(u_n)_{n\ge 0}$ of words of $\cL(X)$
such that $u_0 = \varepsilon$ and $u_n$ is suffix of $u_{n+1}$,
 the sequence of morphisms $(\tau_n)_{n \ge 0}$ built in the proof 
of Proposition~\ref{prop: S-adic UR set} is a primitive proper
 $\Sa$-adic representation of $X$ with $\Sa = \{\tau_n \mid n \ge 0\}$.
By  Theorem~\ref{theoremReturn}, the set $\mathcal{R}_X(u_1)$ is a basis of the free group on $A$. This implies that the matrix $M(\tau_n)$ is unimodular.
This proves the first assertion.

To prove assertion 2, all we need to do is to consider such a sequence $(u_n)_{n \ge 0}$ such that $\tau_n$ is tame for all $n$.

Let $u_1 = a^{(0)}$ be a letter in $A$. Since $X$ is dendric, the set
$\RR_X(u_1)$ has $\Card(A)$ elements by Theorem~\ref{theoremCardinality}.
Let  $\tau_0: A \to\mathcal{R}_X(u_1) $
be a bijection.
By Theorem~\ref{theoremReturn} again, since
the set $\mathcal{R}_X(u_1)$ is a basis of the free group on $A$,
by Theorem~\ref{theoremTame},
it is a tame basis. Thus the morphism $\tau_0: A^* \to A^*$ is 
a tame automorphism.
Let $a^{(1)} \in A$ be a letter and set $u_2 = \tau_0(a^{(1)})$. 
Thus $u_2 \in \mathcal{R}_X(u_1)$ and $u_1$ is a suffix of $u_2$.
By Theorem~\ref{propositionReturns}, the derived shift
 $X^{(1)} = \tau_0^{-1}(X)$ is a minimal dendric shift on the alphabet $A$. 
We thus reiterate the process with $a^{(1)}$ and we conclude by induction with $u_n = \tau_0 \cdots \tau_{n-2}(a^{(n-1)})$ for all $n \geq 2$.
\end{proof}
%%%%%%%%%%%%%%%%%%
We illustrate Theorem~\ref{theoremSadicDendric} by the following example.
\begin{example}\label{exampleSadicTree}
Let $A=\{a,b,c\}$, let $\sigma$ be the substitution
defined by $\sigma(a)=ac$, $\sigma(b)=bac$, $\sigma(c)=cbac$
and let $X$ be the substitution shift generated by $\sigma$.
It can be shown that $X$ is dendric (Exercise~\ref{exerciseExampleJulien}).
We have $\sigma=\alpha_{a,c}\alpha_{b,a}\alpha_{c,b}$. Thus 
$S$ has the $\mathcal{S}_e$-adic representation $(\sigma_n)_{n\ge 0}$ given
by the periodic sequence $\sigma_{3n}=\alpha_{a,c}$, $\sigma_{3n+1}=\alpha_{b,a}$, $\sigma_{3n+2}=\alpha_{c,b}$.
\end{example}
The converse of Theorem~\ref{theoremSadicDendric} is not true, as shown by Example~\ref{exampleNotTree} below (see also Exercise~\ref{exerciseBrun}). 
\begin{example}\label{exampleNotTree}
Let $A=\{a,b,c\}$ and let $\sigma:a\mapsto ac, b\mapsto bac, c\mapsto cb$.
The substitution
shift  generated by $\sigma$ (it
is generated by the fixed point $\sigma^\omega(a)$) is not dendric
since $bb,bc,cb,cc\in \cL(X)$ and thus $\E(\varepsilon)$ has a cycle, although 
$\sigma$ is a tame automorphism since $\sigma=\alpha_{a,c}\alpha_{c,b}\alpha_{b,a}$.
\end{example}

\subsection{Dimension groups of dendric shifts.}
Recall that $\M(X,S)$ denotes the set of invariant measures
on a shift space $X$.
\begin{theorem}\label{theoremDGDendric}
The dimension group of a minimal dendric shift $X$ on 
the alphabet $A$
is $(G,G^+,1_G)$ with $G=\Z^A$,
$G^+=\{x\in\Z^A\mid \langle x,\boldsymbol{\mu}\rangle>0,\mu\in\M(X,S)\}\cup \0$ and $1_G=\1$
where $\1$ is the vector with all components equal to $1$
and $\boldsymbol{\mu}$ is the vector $(\mu([a])_{a\in A}$.
\end{theorem}
\begin{proof}
By Theorem \ref{theoremSadicDendric}, $X$ has a primitive proper
and unimodular $\Sa$-adic representation. Thus
the form of the dimension group is given by Theorem~\ref{theo:dg}.
\end{proof}

\begin{example}
Consider again the dendric shift $X$ generated by the unimodular substitution
$\sigma:a\to ac,b\to bac,c\to cbac$ of Example~\ref{exampleSadicTree}.
Since every word in the image of $\sigma$
ends with $ac$, the morphism 
$\tau:x\mapsto c\sigma(x)c^{-1}$ (where $c^{-1}$ is the inverse
of $c$ in the free group) is  proper and
 $X=X(\tau)$. Since $\tau$ is proper, the shift $X$
has a stationary BV-representation with matrix $M(\sigma)=M(\tau)$.
This implies by Theorem~\ref{theoremDGBratteliDiagram} that the dimension group of $X$
is the group of the matrix $M(\sigma)$.
The matrix $M(\sigma)$ is
\begin{displaymath}
M(\sigma)=\begin{bmatrix}1&0&1\\1&1&1\\1&1&2\end{bmatrix}.
\end{displaymath}
The dominant eigenvalue is the largest root $\lambda$ of
$\lambda^3-4\lambda^2+5\lambda-1=0$. The vector
$w=\begin{bmatrix}\lambda&\lambda-1&(\lambda-1)^2\end{bmatrix}$
is a corresponding eigenvector. The map $x\mapsto\langle x,w\rangle$
sends $K^0(X,S)$ onto $\Z[\lambda]$. The image of the unit
vector $\1_M$ is $\lambda^2$ which a unit of $\Z[\lambda]$. Thus the dimension
group of $X$ is isomorphic to $\Z[\lambda]$.
\end{example}
%%%%%%%%%%%%%%%%%%%%%%%%%%%%%%%
\section{Sturmian shifts}\label{sectionChapter6Sturmian}
%%%%%%%%%%%%%%%%%%%%%%%%%%%%%%%
We illustrate the preceding results on the family of Sturmian shifts.
We have seen that Sturmian shifts are dendric (Proposition~\ref{propositionARisDendric}).
In the particular case of dendric shifts, all the general results
concerning dendric shifts can be formulated more precisely.
We have already seen that for the Return Theorem.
We give below the complete descripition of the $\Sa$-adic representations of Sturmian
shifts.
\subsection{BV-representation of Sturmian shifts}

We define the  morphisms  \( \rho _{n} \)  and \( \gamma _{n} \), $n\ge 1$  
from \( \left\{ 0,1\right\}  \)  to  \( \left\{ 0,1\right\} ^{*} \)    by 

\[
\begin{array}{l}
\rho _{n}\left( 0\right) =01^{n}\\
\rho _{n}\left( 1\right) =01^{n+1}
\end{array}\textrm{ and }\begin{array}{l}
\gamma _{n}\left( 0\right) =10^{n+1}\\
\gamma _{n}\left( 1\right) =10^{n} \   .
\end{array}
\]
  
The morphisms $\rho_n,\gamma_n$ are related as follows to the
elementary automorphisms $L_0,L_1$ introduced in Section~\ref{sectionSturmianShifts}. We have for every $n\ge 1$ and every $u\in\{0,1\}^*$,
\begin{displaymath}
1^n\rho_n(u)=L_{1^n0}(u)1^n,\quad 0^n\gamma_n(u)=L_{0^n1}(u)0^n,
\end{displaymath}
as one verifies easily for $u=0,1$, which implies the identities for all $u$.

The following result will be used to
give a KR-representation of Sturmian shifts.

\begin{proposition}
\label{ch5:theorem:morsehedlund}  
Let  $X$  be a Sturmian shift of slope $\alpha=[0,1+d_1,d_2,\ldots]$.  
There is a Sturmian shift space $Y$ and an $n\ge 1$
such that either $X=\rho_n(Y)$ or $X=\gamma_n(Y)$.
More precisely, if $d_1>0$, then $n=d_1$, $X=\gamma_n(Y)$
and $Y$ is the Sturmian shift of slope $[0,d_2,d_3,\ldots]$. If $d_1=0$,
then $n=d_2$, $X=\rho_n(Y)$ and $Y$ is the Sturmian shift of slope $[0,d_3,d_4,\ldots]$. 
\end{proposition}
\begin{proof}
Let $X^+$ be the one-sided shift space associated to $X$
and let $x$ be the standard Stumian word which belongs to $X^+$.
By Theorem~\ref{standardEpisturmianTheorem},
we have $x=\Pal(\Delta)$ where
 $\Delta=0^{d_1}1^{d_2}\cdots$ is the directive word of $x$. 
Assume that $d_1>0$ and set $\Delta=0^n1\Delta'$ with $n=d_1$.
Then, by Justin Formula \eqref{eqJustin}, we have
$x=L_{0^n1}(\Pal(\Delta'))$. We have seen that $L_{0^n1}(u)0^n=0^n\gamma_n(u)$
for every $u\in\{0,1\}^*$. Thus $x=0^nx'$, where
$x'$ is a concatenation of words
in $\{10^{n+1},10^n\}$. The word $y=\Pal(\Delta')$ is Sturmian
and $x'=\gamma_n(y)$. This shows that $X^+=\gamma_n(Y^+)$
where $Y^+$ is the one-sided shift generated by $y$.
Finally, we obtain $X=\gamma_n(Y)$ where $Y$
is the two-sided shift associated to $Y^+$.
The case where $\Delta$ begins with $1$
is analogous.
\end{proof}

\begin{example}\label{exampleBVSturm}
Let $X=X(\sigma)$ with $\sigma:0\to 01,1\to 010$. 
Actually, $X$ is the Sturmian shift of slope $\alpha=\sqrt{2}-1$.
Indeed, we have $\alpha=[0,2,2,\ldots]$ and thus the directive word of $c_\alpha$ is $011001100\cdots=(0110)^\omega$. The standard word
with slope $\alpha$ is thus the fixed point of the morphism
$L_{0110}:0\to 01010,1\to 0101001$ which is $\sigma^2$.
We have
$\sigma(0)0=0\gamma_1(1)$ and $\sigma(1)0=0\gamma_1(0)$. Thus
$X=\gamma_1(Y)$ where $Y$ is obtained from $X$ by exchanging
$0$ and $1$.
\end{example}

Let  
\( \left( X, S \right)  \)  
be  a  Sturmian  shift. For
\( a\in \left\{ 0,1\right\}  \),  
recall that 
$[ a] =\{( x_{i}) _{i\in \Z }\in X\mid x_{0}=a\}  $.   

\begin{theorem}
\label{ch5:corollary:sturmKR}  
Let $X$ be a Sturmian shift on $\{0,1 \}$.
\index{subject}{Sturmian!shift space}
There  exists  a  sequence  \( \left( \zeta _{n}\right) _{n\in \N } \)  taking  values  in  
\( \left\{ \rho _{1},\gamma _{1},\rho _{2},\gamma _{2},\ldots \right\}  \)  
such  that $({\mathcal P} (n))_n$ is a refining
sequence of KR-partitions with ${\mathcal P} (1) = \{[0] , [1] \}$ and, for $n\geq 2$,
$$  
{\mathcal P} (n)
=
\left\{ 
S^{k}\zeta _{1}\cdots \zeta _{n}\left( \left[ a\right] \right) \mid 0\leq k<\left| \zeta _{1}\cdots \zeta _{n-1}\left( a\right) \right| , \  a\in \left\{ 0,1\right\} \right\}  \ .
$$
\end{theorem}
 \begin{proof}
We apply iteratively Proposition~\ref{ch5:theorem:morsehedlund}
to obtain the sequence $(\zeta_n)$.
To see that condition (KR1) (the intersection of the bases is
reduced to a point) is satisfied,  note  that, if for example $\zeta_n = \rho_i$, then 
$$
\zeta _{1}\cdots \zeta _{n}\left( \left[ a\right] \right) 
\subseteq
\left[
\zeta _{1}\cdots \zeta _{n-1} (1) . \zeta _{1}\cdots \zeta _{n}\left( a \right) \zeta _{1}\cdots \zeta _{n-1} (0) 
\right] \ .
$$
\end{proof}
Let  $ \left( X,S \right)$ be  a Sturmian  shift  
and $({\mathcal P} (n))_{n}$ be the  
sequence  of partitions given  by  Corollary~\ref{ch5:corollary:sturmKR}.  
With  such  a  sequence  is  associated  an  ordered  Bratteli-Vershik  diagram  
$B=\left( V,E,\le \right)$ which  can  be  described  as  follows.  
For  all  $n\geq 1$,  $V_{n}$ consists  of  two  vertices, the substitution read on $E_{n+1} $  is $\zeta _{n}$, with 
$E(1)$ consisting  of a simple hat. We have thus proved the following statement.
\begin{corollary}
A shift space is Sturmian if and only if it has a BV-representation with
simple hat, with
two vertices at every level and such that the substitution read on $E_{n+1}$
is some $\rho_i$ or $\gamma_i$.
\end{corollary}
\index{subject}{BV-representation!of Sturmian shifts}

\begin{example}
Let $X=X(\sigma)$ with $\sigma$ as in Example~\ref{exampleBVSturm}.
The sequence $(\zeta_n)$ can be chosen to be $\zeta_n=\tau$
where $\tau:0\to 10,1\to 100$ (the composition of $\gamma_1$
with the exchange of $0,1$). The corresponding Bratteli diagram
is shown in Figure~\ref{figureBVSturm}.
\begin{figure}[hbt]
\centering
\tikzset{node/.style={circle,draw,minimum size=0.1cm,inner sep=0pt}}

\begin{tikzpicture}
\node[node](0)at(0,1){};
\node[node](11)at(1,0.3){};\node[node](12)at(1,1.7){};
\node[node](21)at(3,0.3){};\node[node](22)at(3,1.7){};
\node[node](31)at(5,0.3){};\node[node](32)at(5,1.7){};

\draw(0)edge node{}(11);\draw(0)edge node{}(12);
\draw[below,near end](11)edge node{$2$}(21);
\draw[bend right=20,above,near end](11)edge node{$3$}(22);\draw[bend left=20,above,near end](11)edge node{$2$}(22);
\draw[below,near end](12)edge node{$1$}(21);\draw[above,near end](12)edge node{$1$}(22);
\draw[below,near end](21)edge node{$2$}(31);
\draw[bend right=20,above,near end](21)edge node{$3$}(32);\draw[bend left=20,above,near end](21)edge node{$2$}(32);
\draw[below,near end](22)edge node{$1$}(31);\draw[above,near end](22)edge node{$1$}(32);
\node at(5.4,0){$\cdots$};\node at(5.4,2){$\cdots$};
\end{tikzpicture}
\caption{The BV-representation of the Sturmian shift $X$.}\label{figureBVSturm}
\end{figure}
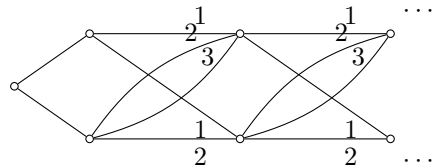
\end{example}

\subsection{Linearly recurrent Sturmian shifts}

We obtain as a corollary of Theorem~\ref{ch5:corollary:sturmKR}, the following result.
\begin{corollary}\label{corollarySturmianLR}
A Sturmian shift $X$ of slope $\alpha=[a_0,a_1,\ldots]$
 is linearly recurrent if and only if
the coefficients $a_i$ are bounded.
\end{corollary}
\begin{proof}
Assume first that the coefficients are not bounded. By Exercise
\ref{exerciseSturmLR}, the shift $X$ is not linearly recurrent.

Conversely, if the coefficients $a_i$ are bounded, the sequence
$(\zeta_n)$ has a finite number of terms and $X$ is linearly
recurrent by Theorem~\ref{ch5:corollary:sturmKR}
\end{proof} 

We can add one more equivalent condition in Corollary~\ref{corollarySturmianLR},
namely: $\cL(X)$ is $K$-power free for some $K\ge 2$
(Exercise~\ref{exerciseSturmPowerFree}).

We note the following additional result.

\begin{theorem}\label{theoremDartnellDurandMaas}
Let $X$ be a Sturmian shift with slope $\alpha$.
Then $X$ is a substitutive shift if and only if $\alpha$ is quadratic.
\end{theorem}
\begin{proof}
  Assume first that $X$ is Sturmian with a slope $\alpha$ which is quadratic.
  By Lagrange Theorem (see Appendix~\ref{appendixAlgebraicNumberTheory}),
  the continued fraction expansion of $\alpha=[0,1+d_1,d_2,\ldots]$ is eventually periodic.
  The directive
  word $x=0^{d_1}1^{d-2}\cdots$ of the standard word $s=c_\alpha$ is eventually periodic.
  Set $x=uy$ with $y=v^\omega$. Then $s=L_u(y)$ and $y$ is a fixed point of $L_v$.
  This shows that $s$ is substitutive and thus that $X$ is substitutive.

  Conversely, \marginpar{A completer}
\end{proof}

\subsection{Derivatives of episturmian shifts}

We have seen that a derivative of a minimal dendric shift
is a minimal dendric shift on the same alphabet
(Theorem~\ref{propositionReturns}).
This implies, since the Sturmian shifts
are the minimal dendric shifts on two letters,
that a derivative of a Sturmian shift is Sturmian.
We prove the following more general statement.

\begin{theorem}\label{theoremDerivativeEpisturmian}
  Any derivative of an episturmian (resp. Arnoux-Rauzy) shift
  is episturmian (resp. Arnoux-Rauzy) on the same alphabet.
\end{theorem}
\begin{proof}
  Let $s=\Pal(x)$ be the standard episturmian sequence with directive
  sequence $x=a_0a_1\cdots$. Every derivative of the shift space $X$
  generated by $s$ is conjugate to the derivative of $X$
  with respect to some palindrome prefix $u_n=\Pal(a_0a_1\cdots a_{n-1})$
  of $s$. By \eqref{equationReturnPal}, the set of left return
  words to $u_n$ is $\RR'_X(u_n)=\{L_{a_0\cdots a_{n-1}}(a)\mid a\in A\}$.
  This shows that the derivative $X'$ of $X$ with respect to $u_n$
  is generated by the standard episturmian word $s'=\Pal(x')$
  with $x'=a_na_{n+1}\cdots$.
  Thus $X'$ is episturmian. If $X$ is Arnoux-Rauzy, then
  every letter appears infinitely often in $x$
  and thus also in $x'$, which implies that $X'$
  is also an Arnoux-Rauzy shift.
  \end{proof}

\begin{example}
  Let $X$ be the Tribonacci shift (see Example~\ref{exampleDerivativeTribo}).
  It is the shift generated by the standard episturmian sequence
  $s=\Pal(abc)^\omega$. The derivative of $X$ with respect to $a$
  is generated by $s'=\Pal(bca)^\omega$.
  \end{example}

%%%%%%%%%%%%%%%%%%%%%%%%%%
\section{Specular shifts}\label{sectionSpecular}
%%%%%%%%%%%%%%%%%%%%%%%%%
We end this chapter with the description of a family of shifts which
generalizes dendric shifts and is build as an abstract
model for the transformations
called linear involutions described in the next chapter.
\subsection{Specular groups}
We begin with the definition of a class of groups which generalizes free groups.
We consider an alphabet $A$ with an involution $\theta:A\rightarrow A$,
possibly with some fixed points. We also consider the
group $G_\theta$\index{symbols}{G@$G_{\theta}$} generated by $A$ with the relations $a\theta(a)=1$
for every $a\in A$. Thus $\theta(a)=a^{-1}$ for $a\in A$.
The set $A$ is called a \emph{natural} set
of generators of $G_\theta$. 
\index{subject}{natural!set of generators}

When $\theta$ has no fixed point, we can set $A=B\cup B^{-1}$
by choosing a set of representatives of the orbits of $\theta$
for the set $B$. The group $G_\theta$ is then the free
group on $B$, denoted $F_B$. 

In general, the group $G_\theta$ is a free product of
a free group and a finite number of copies of $\Z/2\Z$, that
is $G_\theta=\Z^{*i}*(\Z/2\Z)^{*j}$ where $i$  is the number
of orbits of $\theta$ with two elements  and $j$ the number of its fixed points. 
Such a group will be called a \emph{specular
group} \index{subject}{specular!group}\index{subject}{group!specular} of type $(i,j)$. These groups are very close to free groups, as we will see.
The integer $\Card(A)=2i+j$ is called the \emph{symmetric rank}
\index{subject}{symmetric!rank} of the specular group
$\Z^{*i}*(\Z/2\Z)^{*j}$. Two specular groups are isomorphic if and only if they have the same type. Indeed, the commutative image of a group of type
$(i,j)$ is $\Z^i\times (\Z/2\Z)^j$ and the uniqueness of $i,j$ follows
from the fundamental theorem of finitely generated abelian groups
(see Appendix~\ref{appendixGroups}).
\begin{example}\label{exampleSpecularGroup}
Let $A=\{a,b,c,d\}$ and let $\theta$ be the involution which exchanges $b,d$
and fixes $a,c$. Then $G_\theta=\Z*(\Z/2\Z)^2$ is a specular group of symmetric
rank
$4$.
\end{example}
The Cayley graph of a specular group $G_\theta$ with respect to the set
of natural generators $A$ is a regular tree where each vertex has degree
$\Card(A)$.
The specular groups are actually
characterized by this property.

By the Kurosh Subgroup Theorem,
\index{subject}{Kurosh Subgroup Theorem}\index{subject}{Theorem!Kurosh}%
\index{names}{Kurosh, Aleksandr G.}%
 any subgroup of a free product $G_1*G_2*\cdots *G_n$
is itself a free product of a free group and of groups conjugate
to subgroups of the $G_i$ (see Appendix~\ref{appendixGroups}). Thus, we have,
replacing  the Nielsen-Schreier Theorem of free groups, the following
result.
\begin{theorem}\label{theoremKurosh}
Any subgroup of a specular group is specular.
\end{theorem}
It also follows from the Kurosh Theorem that the
elements of order $2$ in a specular group $G_\theta$ are the conjugates of
the $j$ fixed points of $\theta$ and  this number is thus the number
of conjugacy classes of elements of order $2$. Indeed, an element
of order $2$ generates a subgroup conjugate to one the subgroups
generated by the letters.

A word on the alphabet $A$ is $\theta$-\emph{reduced} (or simply  reduced)
if it has no factor of the
form $a\theta(a)$ for $a\in A$. It is clear that any element of a specular
group is represented by a unique reduced word.

A subset of a group $G$ is called \emph{symmetric}
\index{subject}{symmetric!set}
if it is closed under taking inverses. A set $X$ in a specular group $G$
is called a \emph{monoidal basis}\index{subject}{monoidal basis} of $G$ if it is symmetric,  if the monoid that it generates is
$G$  and if any product
$x_1x_2\cdots x_m$ of elements of $X$
such that $x_kx_{k+1}\ne 1$ for $1\le k\le m-1$
is distinct of $1$. The alphabet $A$ is a monoidal basis of $G_\theta$
and the symmetric rank of a specular group is the cardinality
of any monoidal basis (two monoidal bases have the same cardinality
since the type is invariant by isomorphism).

Let $H$ be a subgroup of a specular group $G$.
 Let $Q$ be a set of reduced words on $A$
which is a prefix-closed set of representatives of the right cosets $Hg$
of $H$.
Such a set is traditionally called a \emph{Schreier transversal}
\index{subject}{Schreier!transversal} for $H$
(the proof of its existence is classical in the free group and it
is the same in any specular group).

Let
\begin{equation}
U=\{paq^{-1}\mid a\in A, p,q\in Q, pa\not\in Q, pa\in Hq\}.\label{setSchreier}
\end{equation}
Each word $x$ of $U$ has a unique factorization $paq^{-1}$ with $p,q\in Q$
and $a\in A$. The letter $a$ is called the \emph{central part}
\index{subject}{central part of word} of $x$.
The set $U$
is a monoidal basis of $H$, called the \emph{Schreier basis}
\index{subject}{Schreier!basis}%
relative to $Q$ (the proof is the same as in the free group, see Appendix
\ref{appendixGroups}).

One can deduce directly Theorem~\ref{theoremKurosh} from these
properties of $U$. Indeed, let $\varphi:B\rightarrow U$
be a bijection from a set $B$ onto $U$ which extends to a morphism
from $B^*$ onto $H$. Let $\sigma:B\rightarrow B$
be the involution sending each $b$ to $c$ where $\varphi(c)=\varphi(b)^{-1}$.
Since the central parts never cancel, if a nonempty word
$w\in B^*$  is $\sigma$-reduced
then $\varphi(w)\ne 1$. This shows
that  $H$ is isomorphic to the group $G_\sigma$. Thus $H$ is specular.

If $H$ is a subgroup of index $n$ of a specular group $G$ of symmetric rank $r$,
the symmetric rank $s$ of $H$ is
\begin{equation}
s=n(r-2)+2. \label{SchreierSpecular}
\end{equation}
This formula replaces Schreier's Formula (which corresponds to the case $j=0$). It can be proved as follows.
Let $Q$ be a Schreier transversal for $H$ and let $U$ be the
corresponding Schreier basis.
 The number
of elements of $U$ is $nr-2(n-1)$. Indeed, this
is the number of pairs $(p,a)\in Q\times A$
 minus the $2(n-1)$ pairs $(p,a)$ such that $pa\in Q$ with $pa$
reduced or $pa\in Q$ with $pa$ not reduced.  This gives Formula
\eqref{SchreierSpecular}.
\begin{example}\label{exampleEvenLength}
Let $G$ be the specular group of Example~\ref{exampleSpecularGroup}.
Let $H$ be the subgroup formed by the elements 
 represented by a reduced word of even length. The
set $Q=\{1,a\}$ is a prefix-closed
set of representatives of the two cosets of $H$. The representation of
$G$ by permutations on the cosets of $H$ is represented in Figure~\ref{figureSchreier}.
\begin{figure}[hbt]
\centering
\tikzset{node/.style={circle,draw,minimum size=0.4cm,inner sep=0.4pt}}
\tikzset{title/.style={circle,minimum size=0.1cm,inner sep=0pt}}
\tikzstyle{loop left}=[in=130,out=220,loop]
\tikzstyle{loop right}=[in=330,out=50,loop]
\begin{tikzpicture}
\node[node](1)at(0,0){$1$};\node[node](a)at(2,0){$a$};

\draw[bend left,->,above](1)edge node{$a,b,c,d$}(a);
\draw[bend left,->,below](a)edge node{$a,b,c,d$}(1);
\end{tikzpicture}
\caption{The representation of $G$ by permutations on the cosets of $H$.}
\label{figureSchreier}
\end{figure}
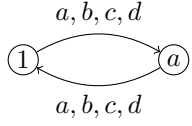
The monoidal basis corresponding to Formula \eqref{setSchreier} is 
$U=\{ab,ac,ad,ba,ca,da\}$.
The symmetric rank of $H$ is $6$, in agreement with Formula~\eqref{SchreierSpecular} and $H$ is a free group of rank $3$.
\end{example}

\begin{example}
Let again $G$ be the specular group of Example~\ref{exampleSpecularGroup}.
Consider now the subgroup $K$ stabilizing $1$ in the representation of
$G$ by permutations on the set $\{1,2\}$ of Figure~\ref{figureSchreier2}.
\begin{figure}[hbt]
\centering
\tikzset{node/.style={circle,draw,minimum size=0.4cm,inner sep=0.4pt}}
\tikzset{title/.style={circle,minimum size=0.1cm,inner sep=0pt}}
\tikzstyle{loop left}=[in=130,out=220,loop]
\tikzstyle{loop right}=[in=330,out=50,loop]
\begin{tikzpicture}
\node[node](1)at(0,0){$1$};\node[node](2)at(2,0){$2$};

\draw[left,->](1)edge[loop left]node{$a,c$}(1);
\draw[bend left,above](1)edge node{$b,d$}(2);
\draw[above,->,bend left](2)edge node{$b,d$}(1);
\draw[right,->](2)edge[loop right]node{$a,c$}(2);
\end{tikzpicture}
\caption{The representation of $G$ by permutations on the cosets of $K$.}
\label{figureSchreier2}
\end{figure}
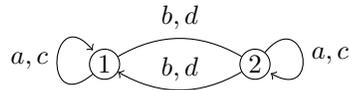
We choose $Q=\{1,b\}$. The set $U$ corresponding to Formula \eqref{setSchreier} is  $U=\{a,bad,bb,bcd,c,dd\}$. The group $K$ is isomorphic to $\Z*(\Z/2\Z)^{*4}$.
\end{example}

Any specular group $G=G_\theta$ has a free subgroup of index $2$. Indeed, let
$H$ be the subgroup formed of the reduced words of even length. It has
clearly index $2$. It is free because it does not contain any element
of order $2$ (such an element is conjugate of a fixed point of $\theta$
and thus is of odd length).

We will need two more properties of specular groups. Both are well-known
to hold for free groups (see Appendix~\ref{appendixGroups}).

A group $G$ is called \emph{residually finite}\index{subject}{residually finite group}\index{subject}{group!residually finite}  if for every element
$g\ne 1$ of $G$, there is a morphism $\varphi$ from $G$
onto a finite group such that $\varphi(g)\ne 1$.

\begin{proposition}\label{propResiduallyFinite}
Any specular group is residually finite.
\end{proposition}
\begin{proof}
Let $K$ be a free subgroup of index $2$ in the specular group $G$.
Let $g\ne 1$ be in $G$. If $g\notin K$, then the image of $g$ in $G/K$
is nontrivial. Assume $g\in K$. Since $K$ is free, it is residually finite.
Let $N$ be a normal subgroup of finite index
of $K$ such that $g\notin N$. Consider the representation of $G$
on the right cosets of $N$. Since $g\notin N$, the image of $g$
in this finite group is nontrivial.
\end{proof}

A group $G$ is said to be \emph{Hopfian}\index{subject}{Hopfian group}
\index{subject}{group!Hopfian} if every surjective morphism from
$G$ onto $G$ is also injective.
By a result of Malcev, any finitely generated residually finite group
is Hopfian.
We thus deduce from Proposition~\ref{propResiduallyFinite} that
any specular group is Hopfian.
As a consequence, we have the following result, which will be used later.

\begin{proposition}\label{propSymBasis}
Let $G$ be a specular group of type $(i,j)$ and let $U\subset G$ be a symmetric set with $2i+j$ elements. If $U$
generates $G$, it is a monoidal basis of $G$.
\end{proposition}
\begin{proof}
Let $A$ be a set of natural generators of $G$.
Considering the commutative image of $G$, we obtain that $U$
contains $j$ elements of order $2$. Thus there is a bijection
$\varphi$ from $A$ onto $U$ such that $\varphi(a^{-1})=\varphi(a)^{-1}$
for every $a\in A$. The map $\varphi$ extends to a morphism
from $G$ to $G$ which is surjective since $U$ generates $G$. Then $\varphi$ being surjective, it
is also injective
since $G$ is Hopfian, and thus $U$ is a monoidal basis of $G$.
\end{proof}
\subsection{Specular shifts}

We assume given an involution $\theta$ on the alphabet $A$ generating
the specular group $G_\theta$.

A symmetric, factorial and extendable 
set $S$ of reduced words on the alphabet $A$
is called
a \emph{laminary set}\index{subject}{laminary set} on $A$ relative to $\theta$. Thus the elements of a laminary set $S$
are elements of the specular group $G_\theta$ and the set $S$ is contained
in $G_\theta$.

A \emph{specular shift}\index{subject}{specular!shift} is a shift space $X$
such that $\cL(X)$ is a laminary set on $A$ which is also dendric of characteristic $2$. Thus, in a specular shift, the extension graph of every nonempty word
is a tree and the extension graph of the empty word is a union of two
disjoint trees. If $X$ is a specular shift, we also say that
$\cL(X)$ is a \emph{specular set}.\index{subject}{specular!set}

The following is a very simple example of a specular shift.
\begin{example}\label{exampleabab}
Let $A=\{a,b\}$ and let $\theta$ be the identity on $A$. Then the periodic
shift formed of the infinite repetitions of $ab$  is a specular shift.
\end{example}
As a second example, we find every dendric shift giving rise to
a  specular shift (that we may consider as degenerate).
\begin{example}
Let $A=B\cup B^{-1}$ be a symmetric alphabet and let $\theta:b\to b^{-1}$.
For every dendric shift $Y$ on the alphabet $B$, the set $L=\cL(Y)\cup\cL(Y)^{(-1)}$
is a laminary set which is dendric of characteristic $2$. Thus
the shift $X$ such that $\cL(X)=L$ is specular.
\end{example}
The next example is more interesting.
\begin{example}\label{exampleJulienSpecular}
Let $A=\{a,b,c,d\}$ and let $X=X(\sigma)$ 
where $\sigma:A^*\to A^*$ is defined by
\begin{displaymath}
\sigma(a)=ab,\quad \sigma(b)=cda,\quad \sigma(c)=cd,\quad \sigma(d)=abc.
\end{displaymath}
We have already seen that  $X$ is eventually dendric of characteristic $2$
(Example~\ref{exampleJulien}). We will see
later (Example~\ref{exampleJulien2}) that $X$ is specular  relative to the involution $\theta=(bd)$.
\end{example}
The following result shows in particular
 that in a specular shift the two trees forming $\E(\varepsilon)$
are isomorphic since they are exchanged by the bijection $(a,b)\rightarrow (b^{-1},a^{-1})$. To distinguish the
disjoint copies of $L(w)$ and $R(w)$
forming the vertices of the extension graph $\E(w)$,
we denote them  by $1\otimes L(w)$ and $R(w)\otimes 1$.
\begin{proposition}\label{propositionClelia}
Let $X$ be a specular shift.
Let $\cT_0,\cT_1$ be the two trees such that $\E_X(\varepsilon)=\cT_0\cup \cT_1$.
For any $a,b\in A$ and $i=0,1$, one has $(1\otimes a,b\otimes 1)\in \cT_i$
if and only if  $(1\otimes b^{-1},a^{-1}\otimes 1)\in \cT_{1-i}$
\end{proposition}
\begin{proof}
Assume that $(1\otimes a,b\otimes 1)$
and  $(1\otimes b^{-1},a^{-1}\otimes 1)$ are both in $\cT_0$. Since $\cT_0$
is a tree, there is a path from $1\otimes a$ to $a^{-1}\otimes 1$. 
We may assume that this path is reduced, that is, does not use
consecutively twice the same edge. Since this path is of odd
length, it has the form $(u_0,v_1,u_1,
\ldots,u_p,v_p)$ with $u_0=1\otimes a$ and $v_p=a^{-1}\otimes 1$.
Since $\cL(X)$ is symmetric, we also have a reduced path
$(v_p^{-1},u_p^{-1},\cdots,u_1^{-1},u_0^{-1})$ which is in $\E(\varepsilon)$
 (for $u_i=1\otimes a_i$, we denote $u_i^{-1}=a_i^{-1}\otimes 1$
 and similarly for $v_i^{-1}$) and thus in $\cT_0$ since $\cT_0,\cT_1$
are disjoint.
Since $v_p^{-1}=u_0$, these two paths have the same origin and end.
But if a path of odd length is its own inverse, its central
edge has the form $(x,y)$ with $x=y^{-1}$, as one verifies easily
by induction on the length of the path. This is a contradiction with
the fact that the words of $\cL(X)$ are reduced.
Thus the two paths  are distinct.
This implies that $\E(\varepsilon)$ has
a cycle, a contradiction.
\end{proof}

We say
that a laminary set $S$ is \emph{orientable}\index{subject}{orientable!laminary set}\index{subject}{laminary set!orientable} if there exist two factorial
sets $S_+,S_-$ such that $S=S_+\cup S_-$ with $S_+\cap S_-=\{\varepsilon\}$
and for any $x\in S$, one has $x\in S_-$ if and only if $x^{-1}\in S_+$
(where $x^{-1}$ is the inverse of $x$ in $G_\theta$).

The following result shows in particular that for 
any dendric shift $X$ on the alphabet $B$, the set $\cL(X)\cup \cL(X)^{-1}$
is a specular set on the alphabet $A=B\cup B^{-1}$.
\begin{theorem}
Let $X$ be a specular shift on the alphabet $A$.
Then, $\cL(X)$ is orientable if and only if there is a partition
$A = A_{+} \cup A_{-}$ of the alphabet $A$ and a dendric shift $Y$ on the
alphabet $B = A_{+}$ such that $\cL(X) = \cL(Y) \cup \cL(Y)^{-1}$.
\end{theorem}

\begin{proof}
Let $X$ be a specular shift on the alphabet $A$ which is orientable.
Let $(S_+,S_-)$ be the corresponding pair of subsets of $S=\cL(X)$.
 The sets $S_+,S_-$ are biextendable, since $S$ is.
Set $A_+=A\cap S_+$ and $A_-=A\cap S_-$. Then $A=A_+\cup A_-$ is a partition
of $A$ and, since $S_-,S_+$ are factorial, we have $S_+\subset A_+^*$
and $S_-\subset A_-^*$.
Let $\cT_0,\cT_1$ be the two trees such that $\E(\varepsilon)=\cT_0\cup \cT_1$.
Assume that a vertex of $\cT_0$ is in $A_+$. Then all vertices
of $\cT_0$ are in $A_+$ and all vertices of $\cT_1$ are in $A_-$.
Moreover, $\E_{S_+}(\varepsilon)=\cT_0$ and $\E_{S_-}(\varepsilon)=\cT_1$.
Thus $S_+=\cL(Y)$ with $Y$ a dendric shift and $S_-=\cL(Y)^{-1}$.
\end{proof} 

The following result follows easily from Proposition~\ref{propositionCassaigne}.

\begin{proposition}\label{propComplexity}
The factor complexity  of a specular shift  containing the alphabet
  $A$ is
$p_n=n(k-2)+2$ for $n\ge 1$ with $k=\Card(A)$. 
\end{proposition}

\subsection{Doubling maps}
We now introduce a construction which allows one to build specular shifts.

A \emph{transducer}\index{subject}{transducer} is a labeled graph with vertices
in a set $Q$ and edges labeled in $\Sigma\times A$. The set $Q$ is called the set of
\emph{states}\index{subject}{state!of transducer}, the set $\Sigma$ is called the \emph{input alphabet}
\index{subject}{input!alphabet}\index{subject}{alphabet!input}%
\index{subject}{transducer!input alphabet}%
and $A$ is called the \emph{output alphabet}.
\index{subject}{output!alphabet}\index{subject}{alphabet!output}%
\index{subject}{transducer!output alphabet}%
The graph obtained by erasing the output
letters is called the \emph{input automaton}. \index{subject}{input!automaton}
\index{subject}{automaton!input} Similarly, the \emph{output automaton}
\index{subject}{output!automaton}\index{subject}{automaton!output}%
 is obtained by erasing the input letters.

Let $\A$ be a transducer with set of states $Q=\{0,1\}$ on the input alphabet
$\Sigma$ and the output alphabet $A$. We assume that
\begin{enumerate}
\item Every
letter of $\Sigma$ acts on $Q$ as a permutation.
\item the output labels of the edges are all distinct.
\end{enumerate}
  We define two maps $\delta_0,\delta_1:\Sigma^*\rightarrow A^*$ corresponding to 
the choice of
$0$ and $1$ respectively as initial vertices.\index{symbols}{delta@$(\delta_0,\delta_1)$}
 Thus $\delta_0(u)=v$ (resp. $\delta_1(u)=v$)
if the path starting at state $0$ (resp. $1$) with input label $u$
has output $v$.
The pair $\delta=(\delta_0,\delta_1)$ is called a \emph{doubling map}
\index{subject}{doubling!map}%
and the transducer $\A$ a \emph{doubling transducer}.
\index{subject}{doubling!transducer}
The \emph{image} of a set $T$ on the alphabet $\Sigma$
by the doubling map $\delta$ is the set $S=\delta_0(T)\cup \delta_1(T)$.

If $\A$ is a doubling transducer, we define an involution $\theta_\A$ on $A$
\index{symbols}{theta@$\theta_\A$}%
as follows. For any $a\in A$, let $(i,\alpha,a,j)$ be the edge with input
label $\alpha$ and output label $a$. We define $\theta_\A(a)$ as the output
label of the edge starting at $1-j$ with input label $\alpha$. Thus, $\theta_\A(a)=\delta_i(\alpha)=a$
if $i+j=1$ and $\theta_\A(a)=\delta_{1-i}(\alpha)\ne a$ if $i=j$.

Recall that the \emph{reversal}\index{subject}{reversal!of word}
 of a word $w=a_1a_2\cdots a_n$ is the word
$\tilde{w}=a_n\cdots a_2a_1$.\index{symbols}{w@$\tilde{w}$} A set $S$ of words is closed under
reversal if $w\in S$ implies $\tilde{w}\in S$ for every $w\in S$.
%As is well known, any Sturmian set is closed under reversal
%(see~\cite{BerstelDeFelicePerrinReutenauerRindone2012}).

\begin{theorem}\label{theoremDoubling}
For any dendric shift $X$ 
on the alphabet $\Sigma$, such that $\cL(X)$ is closed under reversal
and any doubling map $\delta$, the image of $\cL(X)$ by $\delta$
  is a specular set relative to the
involution $\theta_\A$.
\end{theorem}
\begin{proof}
Set $T=\cL(X)$ and  $S=\delta_0(T)\cup \delta_1(T)$. The set $S$ is clearly
biextendable.
Assume that  $x=\delta_i(y)$ for $i\in\{0,1\}$ and $y\in T$.
Let $j$ be the end of the path starting at $i$ and with input
label $y$. Since each letter acts on the two elements of $Q$ as the identity
or as a transposition,
there is a path labeled $\tilde{y}$ from $j$ to
$i$ and  also a path labeled $\tilde{y}$ from $1-j$ to $1-i$.
 Thus $x^{-1}=\delta_{1-j}(\tilde{y})$. Since $T$  is closed under reversal,
$x^{-1}\in\delta_{1-j}(T)$.
This shows that $S$ is symmetric and that it is laminary.

Next, for any nonempty word $x=\delta_i(y)$, the graph $\E_S(x)$ is isomorphic to the graph $\E_T(y)$. Indeed, let $j$ be the end of the path with origin $i$
and input label $y$. For $a,b\in A$, one has $axb\in S$ if and only if
$cyd\in T$ where $c$ (resp. $d$) is the input label of the edge with output label
$a$ (resp. $b$) ending in $i$ (resp. with origin $j$).

Finally, the graph $\E_S(\varepsilon)$ is the union
of two trees isomorphic to $\E_T(\varepsilon)$. Indeed, consider the 
map $\pi$ from $S\cap A^2$ onto $\{0,1\}$ which assigns to $ab\in S\cap A^2$
the state $i$ which is the end of the edge of $\A$ with output 
label $a$ (and the
origin of the edge with output label $b$). Set $S_i=\pi^{-1}(i)$. We have a partition 
$S\cap A^2=S_0\cup S_1$ such that each $S_i$ is isomorphic to $\E_T(\varepsilon)$
and moreover $ab\in S_i$ if and only if $(ab)^{-1}\in S_{1-i}$.
Thus $S$ is specular.
\end{proof}

We now give several examples of specular shifts obtained by
a doubling map. The first one is obtained by doubling the Fibonacci
shift.
\begin{example}\label{exampleFiboDouble}
Let $\Sigma=\{\alpha,\beta\}$ and let $X$ be the Fibonacci shift.
\index{subject}{Fibonacci!shift}\index{subject}{shift space!Fibonacci}%
 Let $\delta$ be
the doubling map given by the transducer of Figure~\ref{figureFiboDouble}
on the left.
\begin{figure}[hbt]
\centering
\tikzset{node/.style={circle,draw,minimum size=0.4cm,inner sep=0.4pt}}
\tikzset{title/.style={circle,minimum size=0.1cm,inner sep=0pt}}
\tikzstyle{loop left}=[in=130,out=220,loop]
\tikzstyle{loop right}=[in=330,out=50,loop]
\begin{tikzpicture}

\node[node](0)at(0,.5){$0$};\node[node](1)at(2,.5){$1$};

\draw[left,->](0)edge[loop left] node{$\beta\mid d$}(0);
\draw[above,bend left,->](0)edge node{$\alpha\mid a$}(1);
\draw[below,bend left,->](1)edge node{$\alpha\mid c$}(0);
\draw[right, ->](1)edge[loop right] node{$\beta\mid b$}(1);

\node[node](bl)at(4,0){$b$};\node[node](al)at(4,1){$a$};
\node[node](cr)at(6,0){$c$};\node[node](br)at(6,1){$b$};

\draw(bl)edge node{}(cr);\draw(al)edge node{}(cr);
\draw(al)edge node{}(br);

\node[node](dl)at(8,0){$d$};\node[node](cl)at(8,1){$c$};
\node[node](ar)at(10,0){$a$};\node[node](dr)at(10,1){$d$};

\draw(dl)edge node{}(ar);\draw(cl)edge node{}(ar);
\draw(cl)edge node{}(dr);

\end{tikzpicture}
\caption{A doubling transducer and the extension graph $\E_S(\varepsilon)$.}\label{figureFiboDouble}
\end{figure}
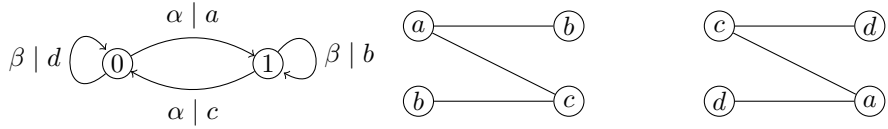

 Then $\theta_\A$ is the involution $\theta$ of Example~\ref{exampleSpecularGroup} and the image of $\cL(X)$ by $\delta$ is a specular set $S$
on the alphabet $A=\{a,b,c,d\}$.
The graph $\E_S(\varepsilon)$ is represented
in Figure~\ref{figureFiboDouble} on the right.

%The letters $a,c$ are odd and $b,d$ are even.

Note that $S$ is the set of factors of the fixed point $g^\omega(a)$
of the morphism
\begin{displaymath}
g: a\mapsto abcab,\quad b\mapsto cda,\quad c\mapsto cdacd,\quad d\mapsto abc.
\end{displaymath}
The morphism $g$ is obtained by applying the doubling map
to the cube $f^3$ of the Fibonacci morphism $f$ in such a way that
$g^\omega(a)=\delta_0(f^\omega(\alpha))$.

\end{example}
In the next example (due to Julien Cassaigne), the specular set
is obtained using a morphism of smaller size.
\begin{example}\label{exampleJulien2}
Let $A=\{a,b,c,d\}$.
Let $T$ be the set of factors of the fixed point $x=f^\omega(\alpha)$ of
the morphism $f:\alpha\mapsto \alpha\beta,\beta\mapsto \alpha\beta\alpha$. It is a Sturmian set.
Indeed, $x$ is the characteristic sequence of slope $-1+\sqrt{2}$
(see Example~\ref{exampleBVSturm}).
%(see~\cite{Lothaire2002}).
The sequence $s_n=f^n(\alpha)$ satisfies $s_n=s_{n-1}^2s_{n-2}$ for $n\ge 2$.
The image $S$ of $T$ by
the doubling  automaton of Figure~\ref{figureFiboDouble} is
the set of factors of the fixed point $\sigma^\omega(a)$
of the morphism $\sigma$ from $A^*$ into
itself defined by 
\begin{displaymath}
\sigma(a)=ab,\quad \sigma(b)=cda,\quad \sigma(c)=cd,\quad \sigma(d)=abc.
\end{displaymath}
Thus the set $S$ is the same as that of Example~\ref{exampleJulien}.
\end{example}

\subsection{Odd and even words}\label{sectionOddEven}
%\subsection{Odd and even letters.}
We introduce a notion which plays, as we shall see,
 an important role in the study
of specular shifts.
Let $X$ be a specular shift on the alphabet $A$ with $A\subset \cL(X)$.
 Any letter $a\in A$ occurs exactly twice as a vertex of $\E(\varepsilon)$,
one as an element of $L(\varepsilon)$ and one as an element of
$R(\varepsilon)$.
 A letter $a\in A$ is said to be \emph{even}
\index{subject}{even!letter}\index{subject}{letter!even}%
if its two occurrences appear in the same tree. Otherwise, it is
said to be \emph{odd}. 
\index{subject}{odd!letter}\index{subject}{letter!odd}%
Observe that if a specular shift $X$ is recurrent, there is at
least one odd letter. 

\begin{example}
Let $X$ be the shift of period $ab$ as in Example~\ref{exampleabab}.
Then $a$ and $b$ are odd.
\end{example}

A word $w\in S$ is said to be \emph{even}\index{subject}{even!word}
\index{subject}{word!even} if it has an even number
of odd letters. Otherwise it is said to be \emph{odd}.\index{subject}{odd!word}
\index{subject}{word!odd}%
The set of even words has the form $U^*\cap S$ where $U\subset S$
 is a bifix code, called the \emph{even code}\index{subject}{even!code}. The set $U$
is the set of even words without a nonempty even
prefix (or suffix). 
\begin{proposition}\label{propositionEvenCode}
 Let $X$ be a minimal specular shift.
The even code is a finite bifix code
which is an $X$-maximal prefix and suffix code.
\end{proposition}
\begin{proof} 
The even code $U$ is bifix by definition.
To prove that it is an $X$-maximal prefix code,
let us verify that
any $w\in S$ is comparable for the prefix order with an element of
the even code
$X$. If $w$ is even, it is in $U^*$. Otherwise, since $S$ is recurrent,
there is a word $u$ such that $wuw\in S$. If $u$ is even, then
$wuw$ is even and thus $wuw\in U^*$. Otherwise $wu$ is even and
thus $wu\in U^*$. This shows that $U$ is $X$-maximal. 
A word of the form $awb$ with $a,b$ odd and $w$ even
cannot be an internal factor of $U$. Indeed, if
$pawbq$ is even, either $p,q$ are even and then $p,awb,q$ are
in $U^*$ or $p,q$ are odd and $pa,w,bq$ are in $U^*$.
Since $X$ is minimal, this implies that $U$ is finite.
\end{proof}
\begin{example}\label{exampleEvenCode}
Let $X$ be the specular shift of Example~\ref{exampleJulien}. The letters
$b,d$ are even and the letters $a,c$ are odd. The even code is
\begin{displaymath}
U=\{abc,ac,b,ca,cda,d\}.
\end{displaymath}
\end{example}

Denote by $\cT_0,\cT_1$\index{symbols}{T@$\cT_0,\cT_1$} the two trees such that $\E(\varepsilon)=\cT_0\cup \cT_1$.
We consider the directed graph $\G$ with vertices $0,1$ and edges all the triples
$(i,a,j)$ for $0\le i,j\le 1$ and $a\in A$ such that $(1\otimes b,a\otimes 1)\in \cT_i$ and
$(1\otimes a,c\otimes 1)\in \cT_j$  for some $b,c\in A$. 
The graph $\G$ is called the
\emph{parity graph}\index{subject}{parity graph}\index{subject}{graph!parity}
 of $S$. Observe that for every letter $a\in A$ there
is exactly one edge labeled $a$ because $a$ appears exactly once as a left
(resp. right) vertex in $\E(\varepsilon)$.

Note that, when $X$ is a specular shift obtained by a doubling map
using a transducer $\A$,
the parity graph of $X$ is the output automaton of $\A$.

\begin{example}
Let $X$ be the specular shift of Example~\ref{exampleFiboDouble}.
The parity graph of $X$ is represented in Figure~\ref{figureParityGraph}.
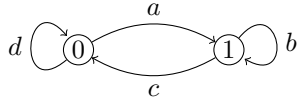
\begin{figure}[hbt]
\centering
\tikzset{node/.style={circle,draw,minimum size=0.4cm,inner sep=0.4pt}}
\tikzset{title/.style={circle,minimum size=0.1cm,inner sep=0pt}}
\tikzstyle{loop left}=[in=130,out=220,loop]
\tikzstyle{loop right}=[in=330,out=50,loop]
\begin{tikzpicture}(20,10)
\node[node](0)at(0,5){$0$};\node[node](1)at(2,5){$1$};

\draw[left,->](0)edge[loop left]node{$d$}(0);
\draw[above,bend left,->](0)edge node{$a$}(1);
\draw[below,bend left,->](1)edge node{$c$}(0);
\draw[right,->](1)edge[loop right]node{$b$}(1);
\end{tikzpicture}
\caption{The parity graph.}\label{figureParityGraph}
\end{figure}
It is the output automaton of the doubling transducer of Figure~\ref{figureFiboDouble}.
\end{example}

\begin{proposition}\label{propositionPartition}
Let $X$ be a  specular shift and let $\G$ be its parity graph.
Let $S_{i,j}$\index{symbols}{S@$S_{i,j}$} be the set of words in $S=\cL(X)$ which are the label of a path
from $i$ to $j$ in the graph $\G$. 
\begin{enumerate}
\item[(1)] The family $(S_{i,j}\setminus\{\varepsilon\})_{0\le i,j\le 1}$ is a partition of $S\setminus \{\varepsilon\}$.
\item[(2)] For $u\in S_{i,j}\setminus\{\varepsilon\}$ and $v\in S_{k,\ell}\setminus\{\varepsilon\}$, if $uv\in S$,
then $j=k$. 
\item[(3)] $S_{0,0}\cup S_{1,1}$ is the set of even words.
\item[(4)] $S_{i,j}^{-1}=S_{1-j,1-i}$.
\end{enumerate}
\end{proposition}
\begin{proof}
We first note that for $a,b\in A$ such that $ab\in S$, there is a path
in $\G$ labeled $ab$. Since $(a,b)\in E(\varepsilon)$,
there is a $k$  such that $(1\otimes a,b\otimes 1)\in\T_k$. Then we have
$a\in S_{i,k}$ and $b\in S_{k,j}$ for some $i,j\in\{0,1\}$. This shows
that $ab$ is the label of a path from $i$ to $j$ in $\G$.

Let us prove by induction on the length of a nonempty word $w\in S$ 
that there exists a unique pair $i,j$ such that $w\in S_{i,j}$.
The property is true for a letter, by definition of the extension
graph $\E(\varepsilon)$ and for words of length $2$
by the above argument. Let next $w=ax$ be in $S$ with $a\in A$
and $x$ nonempty. By induction hypothesis, there is a unique pair
$(k,j)$ such that $x\in S_{k,j}$. Let $b$ be the first letter
of $x$. Then the edge of $\G$ with label $b$ starts in $k$. Since
$ab$ is the label of a path, we have $a\in S_{i,k}$ for some
$i$ and thus $ax\in S_{i,j}$.
The other assertions follow easily
(Assertion (4), follows from Proposition~\ref{propositionClelia}).
\end{proof}
Note that Assertion (4) implies that no nonempty even word is its own
inverse. Indeed, $S_{0,0}^{-1}=S_{1,1}$ and $S_{1,1}^{-1}=S_{0,0}$.

\begin{proposition}\label{propositionOddReturns}
Let $X$ be a specular shift and let $S=\cL(X)$. 
If $x,y\in S$ are nonempty words such that $xyx^{-1}\in S$, then $y$ is odd.
\end{proposition}
\begin{proof}
Let $i,j$ be such that $x\in S_{i,j}$. Then $x^{-1}\in S_{1-j,1-i}$ by 
Assertion (4) of
Proposition~\ref{propositionPartition} and thus $y\in S_{j,1-j}$ by Assertion
 (2).
 Thus $y$ is odd by Assertion (3).
\end{proof}

Recall that for a shift space  $X$, a finite
bifix code $U\subset \cL(X)$ which is an $X$-maximal prefix and
suffix code and  a coding morphism  $f$ for $U$,
the shift space
$Y$ such that  $\cL(Y)=f^{-1}(\cL(X))$ 
is called a \emph{decoding}\index{subject}{decoding of a shift}
 of $X$ by $U$. We denote $Y=f^{-1}(X)$.

\begin{theorem}\label{theoremDecodingEven}
The decoding of a minimal
specular shift by the even code is  a union of two minimal
dendric shifts. More precisely, let $X$ be
a minimal
specular shift and let $f$ be a coding morphism for the even code. The
 shifts $Y_0,Y_1$ 
such that
$\cL(Y_0)=f^{-1}(S_{0,0})$ and $\cL(Y_1)=f^{-1}(S_{1,1})$
are isomorphic minimal dendric shifts.
\end{theorem}
\begin{proof}
We show that the shift $Y_0$ such that
$\cL(Y_0)=f^{-1}(S_{0,0})$ is a minimal  dendric shift. The
proof for $f^{-1}(S_{1,1})$ is the same. 
Set $T_0=\cL(Y_0)$.

Set $S=\cL(X)$.
Since
$X$ is minimal, for every $u\in S$, there exists $n\ge 1$
such that $u$ is a factor of any word $w$ in $S$ of length $n$.
But if $u,w\in S_{0,0}$ are such that $w=\ell ur$, then $\ell,r\in S_{0,0}$.
Thus $Y_0$ is minimal.

We now show that $Y_0$ is dendric.
Let $U$ be the even code.
 Set $U_0=U\cap S_{0,0}$, $U_1=U\cap S_{1,1}$
 and $\E_0(w)=\E_{U_0,U_0}(w)$.  It is enough
to show that the graph $\E_0(w)$ is a tree for any $w\in S_{0,0}$.

Assume first that $w$ is nonempty.
Note first that $\E_0(w)=\E_{U,U}(w)$. Indeed, 
if $x,y\in U$ are such that $xwy\in S$, 
one has $x,y\in U_0$ and $xwy\in S_{0,0}$.
 But  the graph $\E_{U,U}(w)$ is a tree by Proposition
\ref{propStrongTreeConditionBis}.

Suppose now that $w=\varepsilon$. First,
since $\E(\varepsilon)$ is a union of two trees, it  is acyclic, and thus
the graph $\E_0(\varepsilon)$ is acyclic by Proposition~\ref{PropStrongTreeCondition}. Next, since every nonempty
word in $S$ is neutral, by Lemma~\ref{lemmaExtendedMultiplicity},
we have $m_{U,U}(\varepsilon)=m(\varepsilon)=-1$. This implies that
$\E_{U,U}(\varepsilon)$ is a union of two trees. Since $\E_{U,U}(\varepsilon)$
is the disjoint union of $\E_0(\varepsilon)$ and $\E_{U_1,U_1}(\varepsilon)$, 
this
implies that each one is a tree.

Clearly,
$Y_0$ and $Y_1$ are isomorphic. Indeed, let $f:B^*\to A^*$
be a coding morphism for $U$. Set $B_0=f^{-1}(U_0)$ and $B_1=f^{-1}(U_1)$.
Then $\alpha:B_0\to B_1$ defined by $\alpha(b)=f^{-1}(b^{-1})$
defines an isomorphism from $Y_0$ onto $Y_1$.
\end{proof}
Note that the decoding of a dendric shift $X$ by an $X$-maximal prefix
and suffix code is again dendric (Exercise~\ref{exerciseMaximalBifixDecoding}).
\begin{example}
Let $X$ be the shift space of Example~\ref{exampleJulien}. 
We have seen that it generated by the morphism
\begin{displaymath}
\sigma:a\mapsto ab,\quad b\mapsto cda,\quad c\mapsto cd,\quad d\mapsto abc.
\end{displaymath}
The even code $U$
is given in Example~\ref{exampleEvenCode}. Let $\Sigma=\{a,b,c,d,e,f\}$
and let $g$ be the
coding morphism for $X$ given by 
\begin{displaymath}
a\mapsto abc,\quad b\mapsto ac,\quad c\mapsto b,\quad d\mapsto ca,\quad e\mapsto cda,\quad f\mapsto d.
\end{displaymath}
The decoding of $X$
by $U$ is a union of two dendric shifts which generated by the two morphisms
\begin{displaymath}
a\mapsto afbf,\ b\mapsto af,\ f\mapsto a
\end{displaymath}
and 
\begin{displaymath}
c\mapsto e,\ d\mapsto ec,\ e\mapsto ecdc
\end{displaymath}
These two morphisms are actually the restrictions to $\{a,b,f\}$ and $\{c,d,e\}$
of
the morphism $g^{-1}\sigma g$.
\end{example}

\subsection{Complete return words}
Let $X$ be a shift space and let  $U\subset \cL(X)$
be a bifix code. An \emph{internal factor}\index{subject}{internal factor}
of a word $u$ is a word $v$ such that $u=pvs\in U$ for nonempty words $p,s$.
A  \emph{complete return word}
\index{subject}{complete!return word}%
to $U$ is a word of $\cL(X)$ with
a proper prefix in $U$, a proper suffix in $U$
but no internal factor in $U$. We denote by $\CR_X(U)$ the set
of complete return words to $U$.\index{symbols}{CR@$\CR_X(U)$}

  The set $\CR_X(U)$ is a bifix code.
If $X$ is minimal, $\CR_X(U)$ is finite for any finite set
 $U$. For $x\in \cL(X)$, we denote $\CR_X(x)$\index{symbols}{CR@$\CR_X(w)$}
 instead of $\CR_X(\{x\})$.

\begin{example}\label{exampleCompleteFiboDouble}
Let $X$ be the specular shift of Example~\ref{exampleFiboDouble}. One has
\begin{eqnarray*}
\CR_X(a)=\{abca,abcda,acda\},
\CR_X(b)=\{bcab,bcdacdab,bcdacdacdab\}\\
\CR_X(c)=\{cabc,cdabc,cdac\},
\CR_X(d)=\{dabcabcabcd,dabcabcd,dacd\}.
\end{eqnarray*}
\end{example}
The \emph{kernel}\index{subject}{kernel of bifix code}
of a bifix code $U$ is the set of words of $U$ which are factor
of another word in $U$.
The following result is a generalization of Theorem~\ref{theoremCardinality}.
The proof is very similar (Exercise~\ref{exerciseCardCompleteReturn}).

\begin{theorem}\label{propositionNew}
Let $X$ be a minimal specular shift  containing the alphabet $A$.
 For any  finite nonempty bifix code $U\subset S$ with empty kernel, we have
\begin{equation}
 \Card(\CR_X(U))=\Card(U)+\Card(A)-2.\label{formulaComplete}
\end{equation}
\end{theorem}

%As a consequence of Theorem~\ref{propositionNew}, one has the following
%statement. 
%\begin{corollary}\label{corollaryCardCompleteReturns}
%Let $X$ be a minimal specular shift on the
%alphabet $A$ with $A\subset \cL(X)$. For any finite nonempty bifix code $U\subset \cL(X)$ with empty kernel, one has
%\begin{displaymath}
%\Card(\CR_X(U))=\Card(U)+\Card(A)-2.
%\end{displaymath}
%\end{corollary}
The following example illustrates Theorem~\ref{propositionNew}.
\begin{example}\label{exampleJulien3}
Let $X$ be the specular shift on the alphabet $A=\{a,b,c,d\}$
of Example~\ref{exampleJulien}. 
We have
\begin{displaymath}
\CR_X(\{a,b\})=\{ab,acda,bca,bcda\}.
\end{displaymath}
It has four elements in agreement with Theorem~\ref{propositionNew}.
\end{example}

\subsection{Right  return words}

We now come to right return words in specular shifts.
Note that when $S$ is a laminary set $\RR_S(x)^{-1}=\RR'_S(x^{-1})$.

\begin{proposition}\label{propositionReturnsEven}
Let $X$ be a specular shift and let $u\in \cL(X)$ be a nonempty word.
All the words of $\RR_X(x)$ are even.
\end{proposition}
\begin{proof}
If $w\in\RR_X(u)$, we have $uw=vu$ for some $v\in \cL(X)$. If $u$
is odd, assume that $u\in S_{0,1}$. Then $w\in S_{1,1}$.
Thus $w$ is even. If $u$ is even, assume that $u\in S_{0,0}$. Then
$w\in S_{0,0}$ and $w$ is even again.
\end{proof}
We now establish the following result, which replaces,
for specular shifts, Theorem~\ref{theoremCardinality} for dendric shifts.
\begin{theorem}\label{theoremCardRightReturns}
Let $X$ be a minimal specular shift. For any $u\in \cL(X)$,
the set $\RR_X(u)$ has $\Card(A)-1$ elements.
\end{theorem}
\begin{proof}
This follows directly from Theorem~\ref{propositionNew}
with $U=\{u\}$
since $\Card(\RR_X(u))=\Card(\CR_X(u))$.
\end{proof}

\begin{example}\label{exampleFiboDoubleFirst}
Let $X$ be the specular shift of Example~\ref{exampleFiboDouble}. We have
\begin{eqnarray*}
\RR_X(a)&=&\{bca,bcda,cda\},\\
\RR_X(b)&=&\{cab,cdacdab,cdacdacdab\},\\
\RR_X(c)&=&\{abc,dabc,dac\},\\
\RR_X(d)&=&\{abcabcd,abcabcabcd,acd\}.
\end{eqnarray*}
\end{example}

\subsection{Mixed  return words}
Let $S$ be a laminary set.
For $w\in S$ such that $w\ne w^{-1}$, we consider  complete return words to the
set $X=\{w,w^{-1}\}$. 
\begin{example}\label{exampleInvolution3bis}
Let $X$ be the substitution shift
generated by the substitution
$f:a\to cb^{-1},b\to c,c\to ab^{-1}$.
We shall verify later that $X$ is specular (Example~\ref{exampleInvolution3}).
We have 
\begin{eqnarray*}
\CR_S(\{a,a^{-1}\})&=&\{ab^{-1}cba^{-1},ab^{-1}cbc^{-1}a,a^{-1}cb^{-1}c^{-1}a,\\
&&\qquad  ab^{-1}c^{-1}ba^{-1},a^{-1}cbc^{-1}a,a^{-1}cb^{-1}c^{-1}ba^{-1}\}\\
\CR_S(\{b,b^{-1}\})&=&\{ba^{-1}cb,ba^{-1}cb^{-1},bc^{-1}ab^{-1},b^{-1}cb,
b^{-1}c^{-1}ab^{-1},b^{-1}c^{-1}b\},\\
\CR_S(\{c,c^{-1}\})&=&\{cba^{-1}c,cbc^{-1},cb^{-1}c^{-1},c^{-1}ab^{-1}c,
c^{-1}ab^{-1}c^{-1},c^{-1}ba^{-1}c\}.
\end{eqnarray*}

\end{example}
The following result shows that, at the cost of taking return
words to a set of two words, we recover a situation similar
to that of dendric shifts. 
\begin{theorem}\label{theoremCardReturns}
Let $X$ be a minimal specular shift
on the alphabet $A$ such that $A\subset \cL(X)$.
 For any $w\in \cL(X)$ such that $w\ne w^{-1}$, the set of  complete 
return  words to $\{w,w^{-1}\}$ has $\Card(A)$ elements.
\end{theorem}
\begin{proof}
The statement results
directly of Theorem~\ref{propositionNew}. 
\end{proof}
\begin{example}\label{exampleCompleteFiboDoubleMixed}
Let $X$ be the specular shift of Example~\ref{exampleFiboDouble}. In view of
the values of $\CR_S(b)$ and $\CR_S(d)$
given in Example~\ref{exampleCompleteFiboDouble},
we have
\begin{displaymath}
\CR_S(\{b,d\})=\{bcab,bcd,dab,dacd\}.
\end{displaymath}
\end{example}

Two words $u,v$
are said to \emph{overlap}\index{subject}{overlap of words}
 if a nonempty suffix of one of them
is a prefix of the other. In particular a nonempty word overlaps with
itself.

We now consider  the return words to $\{w,w^{-1}\}$ with $w$
such that $w$ and $w^{-1}$ do not overlap. 
This is true for every $w$ in
a laminary set $S$ where the involution $\theta$
has no fixed point (in particular when $X$ is the natural coding of a linear
involution, as we shall see). In this case, the group $G_\theta$ is free and for any $w\in S$, the words $w$ and $w^{-1}$ do not overlap. 

With a complete return word $u$ to $\{w,w^{-1}\}$, we associate a word $N(u)$ obtained as follows. If $u$ has  $w$ as prefix, we erase it
and if $u$ has a suffix $w^{-1}$, we also erase it. Note that these two operations can be made in any order since
$w$ and $w^{-1}$ cannot overlap.

The \emph{mixed  return words to $w$}\index{subject}{mixed return word}
\index{subject}{return!word!mixed} are the words $N(u)$ associated with
 complete return words $u$ to $\{w,w^{-1}\}$. We denote by $\MR_X(w)$ the set
of  mixed return words to $w$ in $X$. \index{symbols}{MR@$\MR_X(w)$}

Note that $\MR_X(w)$ is symmetric and that $w\MR_X(w)w^{-1}=\MR_S(w^{-1})$. Note
also that if $S$ is orientable, then
\begin{displaymath}
\MR_X(w)=\RR_X(w)\cup \RR_X(w)^{-1}=\RR_X(w)\cup \RR'_X(w^{-1}).
\end{displaymath}

\begin{example}\label{exampleMixedReturn}
Let $X$ be the substitution shift generated by the morphism
$f:a\to cb^{-1},b\to c,c\to ab^{-1}$ extended to an automorphism
of the free group on $\{a,b,c\}$. We shall see
later that $X$ is actually specular (Example~\ref{exampleInvolution3}).
We have 
\begin{eqnarray*}
\MR_S(a)&=&\{b^{-1}cb,b^{-1}cbc^{-1}a,a^{-1}cb^{-1}c^{-1}a
,b^{-1}c^{-1}b,a^{-1}cbc^{-1}a,a^{-1}cb^{-1}c^{-1}b\}\\
\MR_S(b)&=&\{a^{-1}cb,a^{-1}c,c^{-1}a,b^{-1}cb,
b^{-1}c^{-1}a,b^{-1}c^{-1}b\},\\
\MR_S(c)&=&\{ba^{-1}c,b,b^{-1},c^{-1}ab^{-1}c,
c^{-1}ab^{-1},c^{-1}ba^{-1}c\}.
\end{eqnarray*}
\end{example}

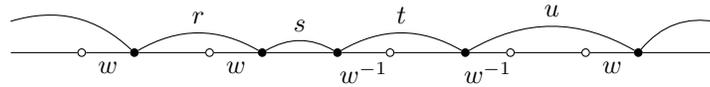
\begin{figure}[hbt]
\centering
\tikzset{node/.style={circle,draw,minimum size=0.1cm,inner sep=0.4pt}}
\tikzset{title/.style={circle,minimum size=0.1cm,inner sep=0pt}}
\begin{tikzpicture}

\node[title](0)at(.5,0){};\node[title](1b)at(.5,.4){};
\node[node](1)at(1.5,0){};\node[node,fill=black](2)at(2.2,0){};
\node[node](3)at(3.2,0){};\node[node,fill=black](4)at(3.9,0){};
\node[node,fill=black](5)at(4.9,0){};\node[node](6)at(5.6,0){};\node[node,fill=black](7)at(6.6,0){};
\node[node](8)at(7.2,0){};\node[node](9)at(8.2,0){};
\node[node,fill=black](10)at(8.9,0){};\node[title](11)at(10.0,0){};
\node[title](11b)at(10.0,.4){};

\draw(0)edge node{}(1);\draw[below](1)edge node[below]{$w$}(2);
\draw[bend left](1b)edge node{}(2);
\draw(2)edge node{}(3);
\draw[below](3)edge node{$w$}(4);\draw[bend left,above](2)edge node{$r$}(4);
\draw(4)edge node{}(5);\draw[below](5)edge node{$w^{-1}$}(6);
\draw[bend left,above](5)edge node{$t$}(7);
\draw[bend left,above](4)edge node{$s$}(5);
\draw(6)edge node{}(7);
\draw[below](7)edge node{$w^{-1}$}(8);\draw[above](8)edge node{}(9);
\draw[below](9)edge node{$w$}(10);
\draw[bend left,above](7)edge node{$u$}(10);
\draw[bend left,above](10)edge node{}(11b);
\draw(10)edge node{}(11);
\end{tikzpicture}
\caption{A uniformly recurrent infinite word factorized as an infinite product
$\cdots rstu\cdots$ of mixed   return words to $w$.}\label{figureFact}
\end{figure}

Observe  that any uniformly recurrent 
bi-infinite word $x$ such that $F(x)=S$ can be uniquely written 
as a concatenation of  mixed return words (see Figure~\ref{figureFact}).
Note that successive occurrences of $w$ may overlap but that
successive occurrences of $w$ and $w^{-1}$ cannot.

We have the following cardinality result. 
\begin{theorem}\label{theoremCardFirstMixed}
Let $X$ be a minimal specular shift on the
alphabet $A$ such that $A\subset\cL(X)$. For any $w\in S$ such that
$w,w^{-1}$ do not overlap, the set $\MR_X(w)$ has $\Card(A)$ elements.
\end{theorem}
\begin{proof}
This is a direct consequence of Theorem~\ref{theoremCardReturns} since
$\Card(\MR_X(w))=\Card(\CR_X(\{w,w^{-1}\})$ when $w$ and $w^{-1}$ do not overlap.
\end{proof}
Note that the bijection between $\CR_X(\{w,w^{-1}\})$ and $\MR_X(w)$
is illustrated in Figure~\ref{figureFact}.
\begin{example}\label{exampleCardFirstMixed}
Let $X$ be the specular shift of Example~\ref{exampleFiboDouble}. The
value of $\CR_X(\{b,d\})$ is given in Example~\ref{exampleCompleteFiboDoubleMixed}.
Since $b,d$ do not overlap, the set
\begin{displaymath}
\MR_X(b)=\{cab,c,dac,dab\}
\end{displaymath}
has four elements in agreement with Theorem~\ref{theoremCardFirstMixed}.
\end{example}
%As a corollary, we obtain the following result.
%\begin{corollary}
%Let $X$ be the natural coding of a linear involution without
%connection on the alphabet $A$. For any $w\in S$, the set
%$\MR_S(w)$ has $\Card(A)$ elements.
%\end{corollary}

\subsection{The Return Theorem for specular shifts}

By Theorem~\ref{theoremReturn},
the set of right   return words to a given word
in a minimal dendric shift on the alphabet $A$
such that $A\subset\cL(X)$ is a basis of
the free group on $A$. We will see a counterpart of this result 
for  specular shifts.

Let $S$ be a specular set.
The \emph{even subgroup}\index{subject}{even!subgroup}
 is the group formed by the even words.
It is a subgroup of index $2$ of $G_\theta$ with
symmetric rank $2(\Card(A)-1)$ by \eqref{SchreierSpecular}
generated by the even code.
 Since no even word is its own
inverse (by Proposition~\ref{propositionPartition}), it is a free group. Thus its rank is $\Card(A)-1$.

The following result replaces, for specular shifts, the Return Theorem
of dendric sifts (Theorem~\ref{theoremReturn}).

\begin{theorem}\label{theoremReturns}
Let $X$ be a minimal specular shift. 
For any $w\in \cL(X)$, the set of  right  return words to $w$ is a  basis
of the even subgroup.
\end{theorem}
\begin{proof}
Set $S=\cL(X)$.
We first consider the case where $w$ is even. Let $f:B^*\rightarrow A^*$
be a coding morphism for the even code $U\subset S$. 
Consider the partition $(S_{i,j})$, as in Proposition~\ref{propositionPartition},
and set $U_0=U\cap S_{0,0}$, $U_1=U\cap S_{1,1}$.
By Theorem~\ref{theoremDecodingEven}, the shift $f^{-1}(X)$ is the union of
 two minimal
dendric shifts, $Y_0$ and $Y_1$
on the alphabets $B_0=f^{-1}(U_0)$ and $B_1=f^{-1}(U_1)$ respectively. We may assume that
$w\in S_{0,0}$. Then $\RR_X(w)$ is the image by $f$ of the set $R=\RR_{Y_0}(f^{-1}(w))$. By Theorem~\ref{theoremReturn}, the set $R$ is
a basis of the free group on $B_0$. Thus $\RR_S(w)$ is a basis of the
image of $F_{B_0}$ by $f$, which is the even subgroup.

Suppose now that $w$ is odd. Since the even code
is an $X$-maximal bifix code, there exists  an odd word $u$ such that $uw\in S$.
Then $\RR_X(uw)\subset \RR_X(w)^*$. By what precedes, the set $\RR_X(uw)$
generates the even subgroup and thus the group generated by $\RR_X(w)$
contains the even subgroup.  Since all words in $\RR_X(w)$ are even, the
group generated by $\RR_X(w)$ is contained in the even subgroup, whence the
equality. We conclude by Theorem~\ref{theoremCardRightReturns}.
\end{proof}
\begin{example}
Let $X$ be the specular shift of Example~\ref{exampleFiboDouble}. The sets
of  right return words to $a,b,c,d$ are given in Example~\ref{exampleFiboDoubleFirst}. Each one is a basis of the even subgroup.
\end{example}
Concerning mixed  return words, we have the following statement.
\begin{theorem}\label{theoremFirstMixed}
Let $X$ be a minimal specular shift.
For any $w\in \cL(X)$ such that $w,w^{-1}$ do not overlap, the set $\MR_X(w)$
is a monoidal basis of the group $G_\theta$.
\end{theorem}
\begin{proof}
Since $w$ and $w^{-1}$ do not overlap, we have $\RR_X(w)\subset \MR_X(w)^*$.
Thus, by Theorem~\ref{theoremReturns}, the group $\langle\MR_X(w)\rangle$
contains the even subgroup. But $\MR_X(w)$ always contains odd
words. Indeed, assume that $w\in S_{i,j}$. Then $w^{-1}\in S_{1-j,1-i}$
and thus any $u\in\MR_X(w)$ such that $wuw^{-1}\in S$ is odd.
Since the even group is a maximal subgroup of $G_\theta$, 
this implies that $\MR_X(w)$ generates
the group $G_\theta$. Finally since $\MR_S(w)$ has $\Card(A)$ elements
by Theorem~\ref{theoremCardFirstMixed}, we obtain the conclusion
by Proposition~\ref{propSymBasis}.
\end{proof}
\begin{example}
Let $X$ be the specular shift of Example~\ref{exampleFiboDouble}. We have
seen in Example~\ref{exampleCardFirstMixed} that
$\MR_X(b)=\{c,cab,dab,dac\}$.
This set is a monoidal basis of $G_\theta$ in agreement with Theorem~\ref{theoremFirstMixed}.
\end{example}
%Since, in the free group, a reduced word $w$ and its inverse do not overlap, we have the following corollary of Theorem~\ref{theoremFirstMixed}
%in the case where the involution $\theta$ has no fixed points.

%\begin{corollary}
%Let $X$ be the natural coding of a linear involution without
%connection  on the alphabet $A=B\cup B^{-1}$.
%For any $w\in S$, the set
%$\MR_S(w)$
%is a monoidal basis of $F_B$.
%\end{corollary}
%\begin{example}
%Let $T$ be the linear involution of Example~\ref{exampleInvolution3}.
%We have seen in Example~\ref{exampleMixedReturn} that
%$\MR_S(b)=\{a^{-1}cb,a^{-1}c,c^{-1}a,b^{-1}cb,b^{-1}c^{-1}a,b^{-1}c^{-1}b\}$.
%It is a monoidal basis of the free group on $\{a,b,c\}$.
%\end{example}

\subsection{Dimension groups of specular shifts.}
We have the following description of dimension groups
of minimal specular shifts. It shows that they
are dimension groups of dendric shifts, except possibly
for the order unit.
\begin{theorem}\label{theoremDGSpecular}
Let $X$ be a minimal specular shift on a $k$ letter alphabet $A$.
The dimension group of $X$ is, as ordered group, isomorphic to 
the dimension group of a minimal dendric shift on $k-1$
letters.
\end{theorem}
\begin{proof}
Let $w\in\cL(X)$. By Theorem~\ref{theoremReturns}, the
set $\RR_X(w)$ is a basis of the even group. Let $Y$
be the shift space induced by $X$ on $[w]$.
Assume that $w$ is even and, for instance that $w\in S_{0,0}$.
Let $f$ is a coding morphism for the even code $U$
with $U_0=U\cap S_{0,0}$. Let 
$Y_0$ be the shift such that $\cL(Y_0)=f^{-1}(S_{0,0})$
Then $\RR_X(w)\subset U_0^*$ and thus $Y\subset Y_0$.
 Since $Y_0$ is a  minimal dendric shift
by Proposition~\ref{theoremDecodingEven}, we have $Y=Y_0$
and we obtain the conclusion that $Y$ is a minimal
dendric shift. Since $\Card(\RR_X(w))=k-1$, this completes the proof
by Proposition~\ref{propositionBVPrimitive}.
\end{proof}

\begin{example}
Let $X$ be the specular shift  generated by the morphism
$\varphi:a\to ab,b\to cda,c\to cd, d\to abc$ (see Example
\ref{exampleJulienSpecular}). The set of return words to $a$
is $\RR_X(a)=\{bca,bcda,cda\}$. It is a basis of the even
group, itself generated by the even code $U=\{abc,ac,b,ca,cda,d\}$.
Let $f:\{ab,ac,bc,ca,cd,da\}\to A_2=\{u,v,w,x,y,z\}$ 
and let $\varphi_2:u\to uw,v\to uw,w\to yzv,x\to yz,y\to yz,z\to uwx$
be the $2$-block presentation of $\varphi$. 
Let $B=\{r,s,t\}$ and let $\phi:B^*\to A_2^*$
be a coding morphism for $f(a\RR_X(a))=\{uwx,uwyz,vyz\}$.
The morphism $\tau:r\to st,s\to str,t\to sr$ is such that
$\varphi_2\circ\phi=\phi\circ\tau$. The matrix $M(\tau)$
is
\begin{displaymath}
M(\tau)=\begin{bmatrix}0&1&1\\1&1&1\\1&1&0\end{bmatrix}
\end{displaymath}

The matrix $M(\tau)$ has eigenvalues $-1,1-\sqrt{2}$
and $\lambda=1+\sqrt{2}$ which is its dominant eigenvalue. A row eigenvector
corresponding to $\lambda$ is $w=\begin{bmatrix}1&\sqrt{2}&1\end{bmatrix}$.
The dimension group is thus $G=\Z^3$ with 
$G^+=\{(a_r,a_s,a_t)\in\Z^3\mid 
a_r+a_s\sqrt{2}+a_t>0\}\cup\{\0\}$ and unit $u=\begin{bmatrix}3&4&3\end{bmatrix}^t$ (the unit is given by the lengths of the words of $\phi(B)$).
The infinitesimal group is generated by the eigenvector 
$\begin{bmatrix}1&0&-1\end{bmatrix}^t$ corresponding
to the eigenvalue $-1$. The quotient is the image
of $G$ by the map $v\to w\cdot v$. It is isomorphic to
$\frac{1}{2}\Z[\sqrt{2}]$ (the unit is sent by  this map to $6+4\sqrt{2}=2(1+\sqrt{2})^2$).

It is interesting to make the following observation. We
have seen before (Example~\ref{exampleJulien2}) that the shift $X$ is 
obtained by a doubling map from the Sturmian shift $Y$
generated by the morphism $\sigma:a\to ab,b\to aba$.
Since the map sending $a,c$ to $a$ and $b,d$ to $b$
is a morphism from $X$ onto $Y$, we know from
Proposition~\ref{proposition4.7.1} that
there is a natural embedding of $K^0(Y,S)$ in
$K^0(X,T)$. Let us look in more detail
how this is related to the doubling map.

The morphism $\sigma$ is eventually proper and thus the
dimension group $K^0(Y,S)$ is the group of the matrix
\begin{displaymath}
M(\sigma)=\begin{bmatrix}1&1\\2&1\end{bmatrix}.
\end{displaymath}
The group is found to be $\Z[\sqrt{2}]$. Thus, up to
the unit, $K^0(Y,S)$ is the same as the quotient of $K^0(X,T)$
by the infinitesimal group. This can be verified directly as follows.
We  have in $X$, first writing  down all possible extensions of $c$
and next using  the form of $\RR_X(a)$
\begin{eqnarray*}
\charac_{[c]}&=&\charac_{[ab\cdot ca]}+\charac_{[a\cdot cda]}+\charac_{ab\cdot cda}\\
&\sim& \charac_{[abca]}+\charac_{[acda]}+\charac_{[abcda]}\\
&=&\charac_{[a]}
\end{eqnarray*}
where $\sim$ denotes the cohomology equivalence. Since 
\begin{displaymath}
M(\varphi)=\begin{bmatrix}
1&1&0&0\\1&0&1&1\\0&0&1&1\\1&1&1&0
\end{bmatrix}
\end{displaymath}
the values of the invariant probability measure $\mu$ of $X$ 
on $a,b,c,d$ are proportional to the left eigenvector of $M(\varphi)$
for $\lambda=1+\sqrt{2}$, which is
$[\sqrt{2}, 1, \sqrt{2}, 1]$.
Thus, using as basis the characteristic
functions of $[a],[b],[d]$,
 we find the group $K^0(X,S)$ as $G=\{(\alpha,\beta,\delta)\in\Z^3
\mid \alpha\sqrt{2}+\beta+\delta>0\}\cup\{0\}$. 
The natural embedding of  $H(Y,S,\Z)$ in $H(X,S,\Z)$
is induced by the map $(\alpha,\beta)\mapsto (2\alpha,\beta,\beta)$.
The image is as expected embedded in $\frac{1}{2}\Z[\sqrt{2}]$.
\end{example}

%%%%%%%%%%%%%%%%%%%%
\section{Exercises}

\exosection{Section~\ref{sectionDendric}}
\begin{exercise}\label{exerciseExampleJulienCassaigne}
Let $X$ be the substitution shift generated by the substitution
$a\to ab,b\to cda,c\to cd,d\to abc$.
Show that the graph of every nonempty word in $\cL(X)$
is a tree.
\end{exercise}
\begin{exercise}\label{exerciseComplexityEentuallyDendric}
  Show that if $X$ is eventually dendric with threshold $m$, its
  complexity is,  for $n\ge m$,
  given by $p_n(X)=Kn+L$ with $K=s_m(X)$ and $L=p_m(X)-ms_m(X)$.
  \end{exercise}
\begin{exercise}\label{exerciseChaconNotEventuallyDendric}
Show that the Chacon ternary shift is not eventually dendric.
\end{exercise}
\begin{exercise}\label{exerciseNeutral}
Let $B=\{1,2,3\}$ and $A=\{a,b,c,d\}$. Let $\tau:A^*\to B^*$
 be the morphism $a\to 12,b\to 2,c\to 3$ and $d\to 13$. Let 
$X$ be the shift on $A$ generated by the morphism $a\to ab,b\to cda,c\to cd,d\to abc$
of Example~\ref{exampleJulien}.
Show that the shift $Y=\tau(X)$ is neutral.
\end{exercise}

\begin{exercise}\label{exerciseMaximalBifix}
Let $X$ be a shift space. A bifix code $U\subset \cL(X)$ is $X$-\emph{maximal}
\index{subject}{bifix code!X-max@$X$-maximal}%
\index{subject}{X-maximal@$X$-maximal!bifix code}%
if it is not strictly contained in another bifix code $V\subset\cL(X)$.
Show that, if $X$ recurrent, every finite $X$-maximal bifix code is
also $X$-maximal as a prefix code (and, symmetrically as a suffix code).
Hint: show that if $U$ is neither $X$-maximal as a prefix code
and as a suffix code, it is not $X$-maximal as a bifix code.
\end{exercise}

\begin{exercise}\label{exercisePartitionsCodes}
  Let $X$ be a minimal shift space. For $u\in\cL(X)$,
  let $\Pg=\{S^i[vu]\mid v\in\RR'_X(u), 0\le i<|v|\}$
  be the partition in towers of Proposition \ref{propositionPartitionReturn}.
  Show that the set
  \begin{displaymath}
    T=\{t\in\cL(X)\mid  [s\cdot t]\in\Pg \mbox{ for some $s$}\}
  \end{displaymath}
  is an $X$-maximal prefix code.
  \end{exercise}

\begin{exercise}\label{exerciseSaturation}
Let $X$ be a minimal dendric shift on the alphabet $A$.
Denote by $F(A)$ the free group on $A$. Let $U\subset \cL(X)$
be a finite $X$-maximal bifix code.
Let $H$ be the subgroup of $F(A)$ generated by $U$.
Show that $H\cap \cL(X)=U^*\cap \cL(X)$.
Hint: consider the coset graph of $U$ and the set $V$ of labels
of simple paths from $\varepsilon$ to $\varepsilon$.
\end{exercise}
\begin{exercise}\label{exerciseDegree1}
%Let $X$ be a shift space.
%A bifix code
%$U\subset \cL(X)$ is called $X$-\emph{complete}
%\index{subject}{X@$X$-complete bifix code}%
%\index{subject}{bifix code!X-complete@$X$-complete}%
%if $U$ is both an $X$-maximal prefix code and an $X$-maximal
%suffix code (in particular, $U$ is an $X$-maximal bifix code).
%\index{subject}{bifix code!X@$X$-maximal}%
Let$X$ be a recurrent shift space and let $U\subset\cL(X)$ be 
a finite $X$-maximal bifix code.
A \emph{parse}
\index{subject}{parse of a word}%
 of a word $w$  is a triple $(s,x,p)$ such that
$w=sxp$ with $s$ a proper suffix of a word in $U$, $x\in U^*$
and $p$ a proper prefix of a word of $U$. 
Let $d_U(w)$\index{symbols}{d@$d_U(w)$} be the number of
parses of $w$. Show that for every $u\in\cL(X)$ and $a\in R_X(u)$ 
\begin{equation}
d_U(ua)=\begin{cases}d_U(u)&\mbox{ if $ua$ has a suffix in $U$}\\
d_U(u)+1&\mbox{otherwise}\end{cases}\label{eqDegree1}
\end{equation}
Show that 
one has for $u,v,w\in\cL(X)$, the inequality
\begin{equation}
d_U(v)\le d_U(uvw)\label{eqDegree}
\end{equation}
with equality if $v$ is not a factor of a word in $U$.
\end{exercise}
\begin{exercise}\label{exerciseDegree2}
Let $X$ be a recurrent shift space and 
let $U$ be a finite $X$-maximal  bifix code.
 The $X$-\emph{degree}
\index{subject}{X-degree@$X$-degree of a bifix code}%
\index{subject}{bifix code!X-degree@$X$-degree of}%
of $U$, denoted $d_U(X)$\index{symbols}{d@$d_U(X)$} is the maximum of the numbers $d_U(w)$
for  $w\in\cL(X)$.
Show that every word in $\cL(X)$ which is not a factor of a word
in $U$ has $d_U(X)$ parses.
\end{exercise}
\begin{exercise}\label{exerciseDegree3}
Show that the $X$-degree of $\cL_n(X)$ is equal to $n$.
\end{exercise}
\begin{exercise}\label{exerciseSuffixesMaximalBifixCode}
Let $X$ be a recurrent shift space. Let $U\subset\cL(X)$ be a finite
$X$-maximal bifix code. Set $d=d_U(X)$. Show that the set $S$ of nonempty proper
suffixes of $U$ is a disjoint union of $d-1$ $X$-maximal prefix codes.
Hint: consider for $2\le i\le d$ the set $S_i$ of proper suffixes $s$
of $U$ such that $d_U(s)=i$.
\end{exercise}
\begin{exercise}\label{exerciseCardinalityTheorem}
Let $X$ be a minimal dendric shift on the alphabet $A$
such that $A\subset\cL(X)$. Show that for every finite
 $X$-maximal bifix code $U$, one has 
\begin{equation}
\Card(U)=(\Card(A)-1)d_U(X)+1.\label{eqCardinalityBifix}
\end{equation}
Hint: use Exercise~\ref{exerciseSuffixesMaximalBifixCode}.
\end{exercise}
\begin{exercise}\label{exerciseDegree4}
Let
$X$ be a minimal dendric shift on the alphabet $A$
such that $A\subset \cL(X)$. Show that 
a finite bifix code $U\subset\cL(X)$ is   $X$-maximal with $X$-degree $d$ if and only if
it is a basis of a subgroup of index $d$ of the free group $F(A)$.
Hint: consider a word $w\in\cL(X)$ which is not a factor of $U$.
Let $Q$ be the set of suffixes of $w$ which are proper prefixes
of $U$. Show that $K=\{v\in F(A)\mid Qv\subset HQ\}$ is 
equal to the free group $F(A)$.

Hint: Consider the set $Q$ of suffixes of a word $w$ having $d_U(X)$ parses.
Show the the set $V$ of words $v$ such that $Qv\subset \langle U\rangle Q$
is a subgroup of the free group containing $\RR_X(w)$.
This implies that $U$ generates a subgroup of index $d=d_U(X)$.
Conclude using Exercise~\ref{exerciseCardinalityTheorem}
and \emph{Schreier's Formula}
\index{subject}{Schreier!Formula}\index{subject}{Formula!Schreier}%
\index{names}{Schreier, Otto}%
asserting that a basis $U$ of a subgroup of index $d$ of the free
group on $A$ has $\Card(U)=d(\Card(U)-1)+1$ elements.

\end{exercise}
\begin{exercise}\label{exerciseExampleJulien}
Show that the substitution
shift generated by $\sigma:a\to ac, b\to bac, c\to cbac$
of Example~\ref{exampleSadicTree}
is dendric.
\end{exercise}

\begin{exercise}\label{exerciseRBC}
Let $X$ be a shift space. A bispecial word $w\in\cL(X)$ is called
\emph{regular}\index{subject}{regular!bispecial word}
\index{subject}{bispecial!word!regular} if there are
unique $a\in L(w)$ and $b\in R(w)$ such that
$aw$ is right special and $wb$ is left special. In other words,
the graph $\E(w)$ is a tree with paths of length at most $3$.
The shift $X$ is said to satisfy the \emph{regular bispecial condition}
\index{subject}{regular!bispecial condition}\index{subject}{bispecial!regular condition}%
if there is an $n\ge 1$ such that every bispecial word in $\cL_n(X)$
is regular.
Show that $X$ is eventually dendric if and only if
satisfies the regular bispecial condition.
\end{exercise}

\begin{exercise}\label{exerciseBrun}
A primitive $\Sa$-adic shift  $X$ is called a \emph{Brun shift}
\index{subject}{Brun shift}\index{subject}{shift space!Brun}%
if it is generated by a directive sequence $\tau=(\tau_n)$ such that
every $\tau_n$ is an elementary automorphism 
$\beta_{ab}=\tilde{\alpha}_{ba}$
(which places $a$ before $b$)
and if for all $n$, we have that $\tau_n\circ\tau_{n+1}$ is
either equal to $\beta_{ab}^2$ or to $\beta_{ab}\alpha_{ca}$
for some $a,b,c\in A$ with $a\ne b$ and $a\ne c$.

1. Prove that the morphism $\sigma:a\to cbccba, b\to cbccb, c\to cbccbacbc$
generates a Brun shift which is not dendric.

2. Show that a Brun shift is a proper unimodular $\Sa$-adic shift.
\end{exercise}

\exosection{Section~\ref{sectionChapter6Sturmian}}

\begin{exercise}\label{exerciseSturmPowerFree}
Let $X$ be a Sturmian shift of slope $\alpha=[a_0,a_1,\ldots]$. Show
that the following conditions are equivalent.
\begin{itemize}
\item[(i)] $X$ is linearly recurrent.
\item[(ii)] The coefficients $a_i$ are bounded.
\item[(iii)] $\cL(X)$ is $K$-power free for some $K\ge 1$.
\end{itemize}
\end{exercise}

\exosection{Section~\ref{sectionSpecular}}
\begin{exercise}\label{exerciseEvenCode}
Let $X$ be a recurrent specular shift. Show that the even
code is an $X$-maximal bifix code of $X$-degree $2$.
\end{exercise}
\begin{exercise}\label{exerciseMaximalBifixDecoding}
Prove that the decoding of dendric shift $X$ by an $X$-maximal
prefix and suffix code is dendric.
\end{exercise}
\begin{exercise}\label{exerciseCardCompleteReturn}
Prove Theorem~\ref{propositionNew}. Hint: Proceed as for the proof
of Theorem~\ref{theoremCardinality}.
\end{exercise}
%%%%%%%%%%%%%%%%%%%%%
\section{Solutions}

\exosection{Section~\ref{sectionDendric}}
\begin{solution}{\ref{exerciseExampleJulienCassaigne}}
Let $x$ be the fixed point $\sigma^\omega(a)$.
 Let $\pi$ be the morphism from $A^*$ onto $\{a,b\}^*$
defined by $\pi(a)=\pi(c)=a$ and $\pi(b)=\pi(d)=b$. The image of $x$
by $\pi$ is the Sturmian word $y$ which is the fixed point of the morphism
$\tau:a\mapsto ab,\ b\mapsto aba$. The word $x$ can be obtained back from
$y$ by changing one every other letter $a$  into a $c$ and any letter $b$
after a $c$ into a $d$ (see Figure~\ref{figureDouble}). Thus every word in $\cL(y)$  gives rise to $2$ words in $\cL(x)$.
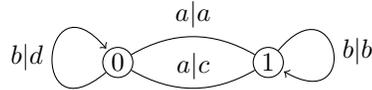
\begin{figure}[hbt]
\centering
\tikzset{node/.style={circle,draw,minimum size=0.4cm,inner sep=0.1pt}}
\tikzstyle{every loop}=[->,shorten >=1pt,looseness=12]
\tikzstyle{loop left}=[in=130,out=220,loop]
\tikzstyle{loop right}=[in=330,out=50,loop]
\begin{tikzpicture}
\node[node](0)at(0,0){$0$};
\node[node](1)at(2,0){$1$};

\draw[left](0)edge[loop left]node{$b|d$}(0);
\draw[above,bend left](0)edge node{$a|a$}(1);\draw[above,bend left](1)edge node{$a|c$}(0);
\draw[right](1)edge[loop right]node{$b|b$}(1);
\end{tikzpicture}
\caption{The inverse of the map $\pi$.}\label{figureDouble}
\end{figure}
In this way every bispecial word $w$ of $\cL(y)$ gives two bispecial words
$w',w''$ of $\cL(x)$ and their extension graphs  are isomorphic to $\E(w)$.
This proves the claim.
\end{solution}
\begin{solution}{\ref{exerciseComplexityEentuallyDendric}}
  We have for $n\ge m$, $p_n(X)-p_m(X)=s_{n-1}+\cdots+s_m=(n-m)s_m$.
  Thus $p_n(X)=ns_m+(p_m-ms_m)$.
  \end{solution}
\begin{solution}{\ref{exerciseChaconNotEventuallyDendric}}
This follows from the fact that, for every $n$, the word $\alpha^n(012)$
is bispecial not neutral by Exercise~\ref{exerciseComplexityChacon}.
\end{solution}
\begin{solution}{\ref{exerciseNeutral}}
Let $g:\{a,c\}A^*\cap A^*\{a,c\}\rightarrow B^*$ be the map defined by
\begin{displaymath}
g(w)=\begin{cases}3\tau(w)&\text{if $w$ begins and ends with $a$}\\
3\tau(w)1&\text{if $w$ begins with $a$ and ends with $c$}\\
2\tau(w)&\text{if $w$ begins with $c$ and ends with $a$}\\
2\tau(w)1&\text{if $w$ begins with $c$ and ends with $c$}\\
\end{cases}
\end{displaymath}
It can be verified that the set of bispecial words of
$\cL(Y)$ is  the union of $\{\varepsilon, 2, 31\}$ and of the
 images by $g$  of nonempty bispecial words of $\cL(X)$ 
(described in the solution of Exercise~\ref{exerciseExampleJulienCassaigne}).
One may verify that these words are neutral.
 Since the  words $\varepsilon$, $2$, $31$ are also neutral,
the shift space $X$ is neutral. 
\end{solution}

\begin{solution}{\ref{exerciseMaximalBifix}}
Consider a word $u\in\cL(X)$ which is not a factor of a word in $U$. Since $X$
is recurrent, there is a word $v$ such that $uvu\in \cL(X)$. Define
a relation $\rho$ on the set $P$ of proper prefixes of $u$ by $(p,q)\in\rho$
if one of the following conditions is satisfied (see Figure~\ref{figureRelationRho}).
\begin{enumerate}
\item[(i)]$q\in pU$.
\item[(ii)]  $u=ps=qt$ and $svq=xyz$ with $x,z\in U$, $y\in U^*$,
$s$ a proper prefix of $x$, $q$ a proper suffix of $z$.
\end{enumerate}
\begin{figure}[hbt]
\centering
\tikzset{node/.style={circle,draw,minimum size=0.1cm,inner sep=0pt}}
\tikzset{title/.style={minimum size=0cm,inner sep=0pt}}
\begin{tikzpicture}
\node[node](p)at(0,0){};\node[node](s)at(1,0){};\node[node](v)at(2,0){};
\node[title](xr)at(2.5,0){};\node[title](zl)at(3.5,0){};
\node[node](q)at(4,0){};\node[node](t)at(5,0){};\node[node](r)at(6,0){};

\draw[below](p)edge node{$p$}(s);\draw[below](s)edge node{$s$}(v);
\draw[below](v)edge node{$v$}(q);\draw[below](q)edge node{$q$}(t);
\draw[below](t)edge node{$t$}(r);
\draw[bend left,above](s)edge node{$x$}(xr);\draw[bend left,above](xr)edge node{$y$}(zl);
\draw[bend left,above](zl)edge node{$z$}(t);
\end{tikzpicture}
\caption{The relation $\rho$.}\label{figureRelationRho}
\end{figure}
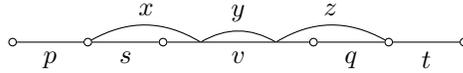
Since $U$ is bifix, the relation $\rho$ is a partial bijection from $P$
to itself.
Assume first that $U$ is $X$-maximal as a suffix code. Then the partial map $\rho$
is onto. This implies that it is a bijection and thus $u$ has a prefix in $U$.
Since this is true for every long enough $u\in\cL(X)$, it implies that $U$ is $X$-maximal
as a prefix code.

Assume now that $U$ is neither $X$-maximal as  a suffix code
nor as a prefix code. Let $y,z\in\cL(X)$ be such that $U\cup y$ is a prefix
code and $U\cup z$ is a suffix code. Since $X$ is recurrent, there is
a $v$ such that $yvz\in\cL(X)$. Then $U\cup yvz$ is a bifix code, a
contradiction.
\end{solution}

\begin{solution}{\ref{exercisePartitionsCodes}}
  This follows easily from the fact that $T$ is also
  the set of suffixes of $\RR'_X(u)u$ which are not suffix of $u$.
  \end{solution}
\begin{solution}{\ref{exerciseSaturation}}
One has clearly $U^*\cap\cL(X)\subset \langle U\rangle\cap\cL(X)$.
Conversely,
consider the coset graph $C$ of $U$ and let $V$ be the set of labels
of simple paths from $\varepsilon$ to itself in $C$. By
Proposition~\ref{propositionCosetGraph}, we have $U\subset V$. 
Since $C$ is Stallings reduced, $V$ is a bifix code and since $U$ is $X$-maximal,
this implies $U=V\cap \cL(X)$. Thus
\begin{displaymath}
\langle U\rangle\cap\cL(X)\subset \langle V\rangle\cap \cL(X)
=V^*\cap\cL(X)=U^*\cap\cL(X)
\end{displaymath}
whence the conclusion.
\end{solution}
\begin{solution}{\ref{exerciseDegree1}}
If $ua$ has a suffix in $X$, the number of parses of $ua$ and $u$
are the same. Otherwise, since $U$ is a maximal suffix code,
$ua$ is a suffix of a word in $U$ and thus
$ua$ has one more parse than $u$, namely $(ua,\varepsilon,\varepsilon)$.
This proves \eqref{eqDegree1}.

Next $d_U(v)\le d_U(uvw)$ since $U$ is $X$-maximal as a prefix code
and as a suffix code.
Indeed, every parse of $v$ extends
to a parse of $uvw$ . Next, if $v$ is not a factor of a word in $U$,
 let $(s,x,p)$ be a parse of $uvw$. 
Since $v$ is not a factor
of a word of $U$, it cannot be a factor of any of $s,x$ or $p$. Thus
there is a parse $(q,y,r)$ of $v$ and a factorization
$x=zyt$ with $z,y,t\in U^*$ such that $sz=uq$ and $rw=zp$. This
shows that every parse of $uvw$ is an extension of a parse of $v$
and thus that $d_U(uvw) =d_U(v)$. This proves ~\eqref{eqDegree}.
\end{solution}
\begin{solution}{\ref{exerciseDegree2}}
Let $w\in \cL(X)$ be such that $d_U(w)=d_U(X)$. Let $u\in \cL(X)$
 not a factor of a word in $U$. Since $X$ is recurrent, there
is a $v\in\cL(X)$ such that $uvw\in\cL(X)$. By Equation~\eqref{eqDegree},
we have $d_U(u)=d_U(uvw)\ge d_U(w)$. Thus $d_U(u)=d_U(X)$.
\end{solution}
\begin{solution}{\ref{exerciseDegree3}}
Every word of length at least $n$ has clearly $n$ parses.
\end{solution}
\begin{solution}{\ref{exerciseSuffixesMaximalBifixCode}}
For each $i$ with $2\le i\le d_U(X)$, 
let $S_i$ be the set of proper suffixes $s$ of $U$ such that $d_U(s)=i$.
Then $S_i$ is a prefix code. Indeed, if $s,t\in S_i$ and if $s$ is a proper
prefix of $t$, then $d_U(s)=d_U(t)$ implies that $t$ has a suffix
in $U$ by Equation \eqref{eqDegree1} of Exercise~\ref{exerciseDegree1}, a contradiction. Next, $S_i$ is
an $X$-maximal prefix code. Indeed, let $w\in\cL(X)$ be long enough
so that $d_U(w)=d$. Then $w$ has nonempty prefixes $s_2,\ldots,s_d$
such that $s_i\in S_i$ for $2\le i\le d$. We conclude that
the set $S$ of nonempty proper prefixes of $U$ is a disjoint union
of $d-1$ $X$-maximal prefix codes.
\end{solution}
\begin{solution}{\ref{exerciseCardinalityTheorem}}
Set $d=d_U(X)$.
Let $P$ be the set of proper prefixes of $U$. By the well-known
formula relating the number of leaves of a tree to the number
of children of its interior nodes, we have
\begin{equation}
\Card(U)-1=\sum_{p\in P}(r_X(p)-1).\label{eqCardUSumRho}
\end{equation}
By (the dual of) Exercise~\ref{exerciseSuffixesMaximalBifixCode}, the set $P\setminus\{\varepsilon\}$
is a disjoint union of $d-1$ $X$-maximal suffix codes $V_1,\ldots,V_{d-1}$. 
Set $\rho(u)=r_X(u)-1$ and, for $V\subset\cL(X)$,
denote $\rho(V)=\sum_{v\in V}\rho(v)$.
By Lemma~\ref{lemmaLeftproba}, we have $\rho(u)=\sum_{a\in L(u)}\rho(au)$. This implies
that for any $X$-maximal suffix code $V$, one has 
\begin{displaymath}
\rho(V)=\rho(\varepsilon)=\Card(A)-1
\end{displaymath}
where the last equality results of the hypothesis $A\subset\cL(X)$.
Thus, we have $\Card(U)-1=\rho(P)=\rho(\varepsilon)+\sum_{i=1}^{d-1}\rho(V_i)
=(\Card(A)-1)d$.
\end{solution}

\begin{solution}{\ref{exerciseDegree4}}
Suppose first that $U$ is a finite $X$-maximal bifix code of $X$-degree
 $d=d_U(X)$ and let $H=\langle U\rangle$ be the subgroup generated by $U$.
Let $w\in\cL(X)$ be a word which not a factor of a word in $U$.
Let $Q$ be the set of suffixes of $w$ which are proper prefixes of $U$.
Then, by Exercise~\ref{exerciseDegree2}, $w$ has $d$ parses
and thus $\Card(Q)=d$. 

Moreover, we claim that it follows from Exercise~\ref{exerciseSaturation} that the cosets $Hq$ for $q\in Q$ are distinct. Indeed,
let $p,q\in Q$ be such that $Hp=Hq$. Since $p,q$ are suffixes
of $w$, one is a suffix of the other. Assume that $q=tp$.
Then $Hp=Htp$ implies $Ht=H$ and thus $t\in H$. By Exercise~\ref{exerciseSaturation}, this implies $t\in U^*$. Since $t$
is a proper prefix of $U$, we conclude that $t=\varepsilon$
and thus $p=q$, which establishes the claim.

 Consider the set
 $K=\{v\in F(A)\mid Qv\subset HQ\}$.
It is  a subgroup of $F(A)$. Indeed, by what precedes, the map $p\mapsto q$
if $pv\in Hq$ is a permutation of $Q$ for every $v\in K$.

Next, we have $\RR_X(w)\subset K$. In fact, consider $v\in \RR_X(w)$.
For every $p\in Q$, since $U$ is $X$-complete, there is
some $x\in U^*$ and some proper prefix $q$ of $U$ such that
$pv=xq$. But since $v$ is in $\RR_X(w)$, $pv$ ends with $w$
and thus $q\in Q$.

Now, by Theorem \ref{theoremReturn}, $\RR_X(w)$ generates $F(A)$
and thus $K=F(A)$. We conclude that $F(A)\subset HQ$ and thus that
$Q$ is a set of representatives of the cosets of $H$. Thus
$H$ has index $d$. Since $U$ generates a subgroup of index $d$
and since $\Card(U)-1=d(\Card(A)-1)$, we conclude by Schreier's Formula
that $U$ is a basis of $H$.

Conversely, if the bifix code $U\subset \cL(X)$ is a basis of a subgroup $H$ of index $d$, let $C$ be the coset graph of $U$. 
By Proposition~\ref{propositionCosetGraph}, $C$ is the Stallings
graph of a subgroup of index $d$. Moreover $U$ is contained
in the set $V$ of labels of simple paths from $\varepsilon$ to $\varepsilon$.
Then $W=V\cap \cL(X)$ is an $X$-maximal bifix code of $X$-degree at most $d$.
By Exercise~\ref{exerciseCardinalityTheorem}, we have
$\Card(W)\le d(\Card(A)-1)$
and thus, by Schreier's Formula $\Card(W)\le \Card(U)$. 
This forces $U=W$ and concludes the proof.
\end{solution}
\begin{solution}{\ref{exerciseExampleJulien}}
The right-special words
are the suffixes of the words $\sigma^n(c)$ for $n\ge 1$ and the left-special words
are the prefixes of the words $\sigma^n(a)$ or $\sigma^n(c)$ for $n\ge 1$, 
as one may verify. 
   Let us show by induction on the length of $w$
that for any bispecial word $w\in \cL(X)$, the graph $\E(w)$ is a tree.
It is true for $w=c$ and $w=ac$. Assume that $|w|\ge 2$. Either
$w$ begins with $a$ or with $c$. Assume the first case. Then $w$ begins
and ends with $ac$. We must have $w=ac\sigma(u)$ where $u$ is a bispecial word
beginning and ending with $c$.  In the second case, $w$ begins with
$cbac$ and ends with $ac$. We must have $w=cbac\sigma(u)$ where
$u$ is a bispecial word beginning with $a$. In both cases,
by induction hypothesis, $\E(u)$ is a tree
and thus $\E(w)$ is a tree.
\end{solution}
\begin{solution}{\ref{exerciseRBC}}
Denote by $LS(X)$ (resp. $LS_{\ge n}(X)$) the set of left-special words in $\cL(X)$
(resp. $\cL_{\ge n}(X)$).

Assume first that $X$ is eventually dendric with threshold $m$.
Then any  word $w$ in $LS_{\ge m}(X)$ has at least one 
 right extension in $LS(X)$.
Indeed, since $R_1(w)$ has at least two elements and since the graph $\E_1(w)$ is connected, there is at least one element $r$ of $R_1(w)$ which is connected by an edge to more than one element of $L_1(w)$ and thus that $rw\in LS(X)$.

Next, the symmetric of Equation~\eqref{eqLeftProba} shows 
that for any  $w \in LS_{\ge m}(X)$ which has more than one 
 right extension in $LS(X)$, one has $\ell(wb) < \ell(w)$ for each such extension. 
Thus the number of  words in $LS_{\ge m}(X)$ which are
prefix of one another, and which have more than one right extension,
 is bounded by $\Card(A)$.
This proves that there exists an $n \ge m$ such that for any $w \in LS_{\ge n}(X)$ there is exactly one $b \in A$ for which $wb \in LS(X)$ is left-special. Moreover, one has then $\ell(wb) = \ell(w)$ by the symmetric
of Equation~\eqref{eqLeftProba}. This proves the uniqueness of $r$.
The proof for left extensions of right-special words is symmetric.

Conversely, assume that the regular bispecial 
condition is satisfied for some integer $n$.
For any word $w$ in $\cL_{\ge n}(X)$, the graph $\E_1(w)$ is acyclic since all vertices in $R_1(w)$ except at most one have degree $1$.
Let $w\in LS_{\ge n}(X)$. If $w$ is not bispecial,
 there is exactly one $b$ such that $w\in\cL(X)$
and then $wb\in LS(X)$. If it is bispecial,
it is regular and there is exactly one $b$ such that
$wb\in LS(X)$. Thus, in both cases, there is exactly 
one $b$ such that $wb\in LS(X)$. One has then $L(wb)\subset L(w)$.
Thus there exists an $N\ge n$ such that for every $w\in LS_{\ge N}(X)$,
one has $wb\in LS(X)$ and moreover $L(wb)=L(w)$. But for such a
$w$, the extension graph $\E(w)$ is connected and thus it is a tree.
This shows that $X$ is eventually dendric.
\end{solution}
\begin{solution}{\ref{exerciseBrun}}
1. We have $\sigma=\beta_{cb}\circ\beta_{bc}\circ\beta_{cb}\circ\beta_{ba}\circ\beta_{ac}$ 
since for example
\begin{displaymath}
a\edge{\beta_{ac}}a\edge{\beta_{ba}}ba\edge{\beta_{cb}}cba\edge{\beta_{bc}}bcba\edge{\beta_{cb}}cbccba.
\end{displaymath}
 The extension graph of the word $w=cbccb$ is shown in Figure~\ref{figureE(cbccb)}.
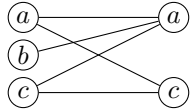
\begin{figure}[hbt]
\centering
\tikzset{node/.style={circle,draw,minimum size=0.4cm,inner sep=0pt}}
\tikzset{title/.style={minimum size=0cm,inner sep=0pt}}
\begin{tikzpicture}
\node[node](al)at(0,3){$a$};\node[node](ar)at(2,3){$a$};
\node[node](bl)at(0,2.5){$b$};
\node[node](cl)at(0,2){$c$};\node[node](cr)at(2,2){$c$};

\draw(al)edge node{}(ar);\draw(al)edge node{}(cr);
\draw(bl)edge node{}(ar);
\draw(cl)edge node{}(ar);\draw(cl)edge node{}(cr);
\end{tikzpicture}
\caption{The extension graph of $w$.}\label{figureE(cbccb)}
\end{figure}
Since this graph has a cycle, the shift $X(\sigma)$ is not dendric.

2. Let $\alpha=(\beta_{a_nb_n})$ be the directive sequence of morphisms
defining the Brun shift $X$. Since $\alpha$ is primitive, there
exists an increasing sequence $(n_k)$ of integers such
that the set $\{a_{n_i}\mid n_k\le n_i<n_{k+1}\}$ is equal to $A$.
Then all morphisms $\beta_{n_k}\circ\cdots\circ\beta_{n_{n+1}-1}$ are left
proper. Indeed, let $\sigma=(\beta_{a_1a_2})^{i_1}\circ(\beta_{a_2a_3})^{i_2}
\circ\cdots\circ(\beta_{a_na_{n+1}})^{i_n}$ with $i_j=1$ or $2$.
Then all words $\sigma(a_1),\ldots,\sigma(a_{n+1})$ begin with $a_1$.
Using finally Lemma~\ref{lemma:proper}, we obtain the conclusion.
\end{solution}

\exosection{Section~\ref{sectionChapter6Sturmian}}

\begin{solution}{\ref{exerciseSturmPowerFree}}
The equivalence of (i) and (ii) is Corollary~\ref{corollarySturmianLR}.
Next, the implication (i) $\Rightarrow$ (iii)
results from Proposition~\ref{propositionLRhasLinearComplexity}. 
Finally we have seen in Exercise~\ref{exercisecfrac2}
that (iii) $\Rightarrow$ (ii).

\end{solution}

\exosection{Section~\ref{sectionSpecular}}
\begin{solution}{\ref{exerciseEvenCode}}
This follows from the fact that any 
 word of the even code  $U$ of length at least $2$ is not an internal factor 
of $U$ and has two parses.
\end{solution}
\begin{solution}{\ref{exerciseMaximalBifixDecoding}}
Let $X$ be a dendric shift on the alphabet $A$. Let
$U\subset \cL(X)$ be a bifix code which is both
an $X$-maximal prefix and suffix code.
Let $f;B^*\to A^*$ be a coding morphism for $U$.
Let $Y$ be the shift space such that $\cL(Y)=f^{-1}(\cL(X))$.
For every $w\in\cL(Y)$, the generalized extension graph $\E_{U,U}(f(w))$
is a tree by Proposition~\ref{propStrongTreeConditionBis}.
Since this graph is isomorphic to $\E_Y(w)$, this proves that
$Y$ is dendric.
\end{solution}
\begin{solution}{\ref{exerciseCardCompleteReturn}}
Set $S=\cL(X)$.
Let $P$ be the set 
of proper prefixes of $\CR_X( U )$ . For $q \in P$ , we define 
$\alpha ( q ) = \Card \{ a \in A \mid qa \in P \cup \CR_X( U )\} - 1$. For
$P'\subset P$, we set $\alpha ( P ') = \sum_{p \in P'} \alpha ( p )$ .
Since $\CR_X( U )$ is a finite nonempty prefix code, we have, by a well-known property of trees, $\Card ( CR_X( U )) = 1 + \alpha ( P )$.
Let $P'$ be the set of words in $P$ which are proper prefixes of $U$
 and let $Y = P \setminus P'$ . Since $P'$ is the set of proper prefixes
of $U$ , we have $\alpha ( P') = \Card ( U ) - 1$.

For $u\in\cL(X)$, set 
\begin{displaymath}
\rho(u)=\begin{cases}r_X(u)-1&\mbox{if $u\ne\varepsilon$}\\
\Card(A)-2&\mbox{otherwise}.\end{cases}
\end{displaymath}
In this way, we have $\rho(u)=\sum_{a\in L(u)}(au)$ for every $u\in\cL(X)$.
 Thus, if $Y$ is an $X$-maximal suffix code, we have
$\sum_{y\in Y}\rho(y)=\rho(\varepsilon)=\Card(A)-2$ .

Since $P \cup \CR_X( U ) \subset S$, one has $\alpha ( q ) \le \rho( q )$
 for any nonempty $q \in P$ . Moreover, since $X$ is recurrent, and since $U$ has empty kernel,
any word of $S$ with a prefix in $U$ is comparable for the prefix order with a word of $\CR_X( U )$ . This implies that for any $q \in Y$
and any $b \in \RR( q )$ , one has $qb \in P \cup \CR( U )$. Consequently, we have $\alpha ( q ) = \rho( q )$ for any $q \in Y$ . Thus we have shown that
$\Card ( \CR_X( U )) = 1 + \alpha ( P' ) + \rho ( Y ) = \Card ( U) + \rho ( Y )$. 

Let us show that $Y$ is an $X$-maximal suffix code. This will imply
our conclusion. Suppose that $q , uq \in Y$ with $u$ nonempty. Since $q$ is in $Y$ , it has a proper prefix in $U$ . But
this implies that $uq$ has an internal factor in $U$, a contradiction. Thus $Y$ is a suffix code. 
Consider $w \in S$. Then, for any $x \in U$ , there is some $u \in S$
 such that $xu w \in S$. Let $y$ be the shortest suffix of $xuw$ which has
a proper prefix in $U$ . Then $y \in Y$ . This shows that $Y$ is an $X$-maximal suffix code.
\end{solution}

%%%%%%%%%%%%%%%%%%%%%%%%%
\section{Notes}

%dendric
\subsection{Dendric shifts}
The languages of dendric shifts were
 introduced in~\cite{BertheDeFeliceDolceLeroyPerrinReutenauerRindone2015}
\index{names}{Berth\'e, Val\'erie}\index{names}{De Felice, Clelia}%
\index{names}{Dolce, Francesco}\index{names}{Leroy, Julien}%
\index{names}{Perrin, Dominique}\index{names}{Reutenauer, Christophe}%
\index{names}{Rindone, Giuseppina}% 
under the name of \emph{tree sets}.
\index{subject}{tree set}%
The language of the shift of Example~\ref{exampleJulien}
is a tree set of characteristic $2$ (\cite[Example 4.2]{BertheDeFeliceDelecroixDolceLeroyPerrinReutenauerRindone2017}) 
and it is actually a specular set.
Tree sets of characteristic $c \ge 1$ 
were introduced in~\cite{BertheDeFeliceDelecroixDolceLeroyPerrinReutenauerRindone2017}
(see also \cite{DolcePerrin2017}).
The shift spaces of
Example \ref{exampleJulien} and of Exercise \ref{exerciseNeutral} (a neutral shift which is not dendric)
are due to Julien Cassaigne~\citep{Cassaigne2015}.

It can be proved that the class of eventually dendric shifts is closed
under conjugacy (see~\cite{DolcePerrin2019}).

Theorem~\ref{theoremReturn} is from~\cite{BertheDeFeliceDolceLeroyPerrinReutenauerRindone2015}. Theorem~\ref{theoremCardinality} was proved earlier
in \cite{BalkovaPelantovaSteiner2008}
\index{names}{Balkov\'{a}, L'ubom\'{\i}ra}%
\index{names}{Pelantov\'{a}, Edita}%
\index{names}{Steiner, Wolfgang}%
 and
was proved even earlier for episturmian shifts in
\cite{JustinVuillon2000}.
\index{names}{Justin, Jacques}\index{names}{Vuillon, Laurent}%

Theorem~\ref{propositionReturns} is from~\cite[Theorem 6.5.19]{BertheDeFeliceDolceLeroyPerrinReutenauerRindone2015b}.
It generalizes the fact that the derived word of a Sturmian
word
is Sturmian (see~\cite{JustinVuillon2000}).

Theorem~\ref{theoremBifixBasisFreeGroup} is from~\cite{BertheDeFeliceDolceLeroyPerrinReutenauerRindone2015b}. It is obtained as a corollary
of a result (called the Finite Index Basis Theorem) proved
in \cite{BertheDeFeliceDolceLeroyPerrinReutenauerRindone2015c}
(see Exercise~\ref{exerciseDegree4}).

The fact that the group of positive automorphisms of a free
group on three letters is not finitely generated is from
\cite{TanWenZhang2004}.
\index{names}{Tan, Bo}\index{names}{Wen, Zhi Xiong}\index{names}{Zhang, Yi Ping}%
The word `tame' used for tame
automorphisms (as opposed to wild) is used here on analogy
with its use in ring theory (see~\cite{Cohn1985}).
\index{names}{Cohn, Paul M.}%
 The tame
automorphisms as introduced here should, strictly speaking, be called
positive tame automophisms since the group of all automorphisms,
positive or not, is tame in the sense that it is generated by the
elementary automorphisms.

Theorem~\ref{theoremSadicDendric} is from \cite{BertheDeFeliceDolceLeroyPerrinReutenauerRindone2015b}.

In the case of a ternary alphabet, a characterization
of tree sets by their $S$-adic representation can be proved,
showing that there exists a  \emph{B\"uchi automaton}\index{subject}{B\"uchi automaton}
\index{subject}{automaton!B\"uchi} \index{names}{Bu@B\"uchi, Richard J.} on the alphabet $\mathcal{S}_e$
recognizing
the set of  $\Sa_e$-adic representations
of uniformly recurrent tree sets~\citep{Leroy2014,Leroy2014bis}.

The Brun shifts \index{subject}{Brun shift} of Exercise~\ref{exerciseBrun}
arise in the generalization of continued fractions expansion to
triples of integers instead of pairs \citep{Brun1958}.
\index{names}{Brun, Viggo}%
The study of these algorithms as dynamical system was initiated
in \cite{BertheSteinerThuswaldner2019}.
\index{names}{Berth\'e, Val\'erie}\index{names}{Steiner, Wolfgang}%
\index{names}{Thuswaldner, J\"org M.} The relation
of these shifts with dendric shifts was studied in
\cite{LabbeLeroy2016}.\index{names}{Labb\'e S\'ebastien}%
\index{names}{Leroy, Julien}%
 It is shown in \citep{AvilaDelecroix2013},
based on an analysis similar to that
done for Arnoux-Rauzy shifts, that Brun shifts (on a ternary alphabet)
are uniquely ergodic.

The regular bispecial condition 
(Exercise \ref{exerciseRBC}) was introduced by~\cite{DamronFickenscher2019}.
\index{names}{Damron, Michael}\index{names}{Fickenscher, Jon}%
They proved that for an irreducible shift satisfying this condition,
the number of ergodic measures is at most $(K+1)/2$
where $K$ is the limiting value of $p_{n+1}(X)-p_n(X)$. This
generalizes a result proved independently by  \cite{Katok1973}
\index{names}{Katok, Anatole B.}%
and \cite{Veech1978}
\index{names}{Veech, William A.}%
 concerning interval exchange transformations.

\subsection{Sturmian shifts}
Theorem
\ref{ch5:theorem:morsehedlund}   is already in
\index{names}{Morse, Marston}\index{names}{Hedlund, Gustav A.} \cite{MorseHedlund1940}. 
Assertion~1  follows  from  Theorem~7.1  
in  \cite{MorseHedlund1940} and  Assertion~2  is  Theorem~8.1.

Corollary \ref{corollarySturmianLR} is from \cite{Durand&Host&Skau:1999}.
It is of course related to the fact that a Sturmian word of slope $\alpha$ is
$k$-power free if and only if the coefficients of the expansion of
$\alpha$ as a continued fraction are bounded (Exercise~\ref{exercisecfrac2}).

Theorem \ref{theoremDartnellDurandMaas} is from
 \cite{Dartnell&Durand&Maass:2000}.\index{names}{Dartnell, Pablo R.}\index{names}{Durand, Fabien}\index{names}{Maass, Alejandro}
 It is proved as a consequence
of the fact that a Sturmian shift
 is a substitution shift if and only if the sequence
 $(\zeta_n )_n$ is ultimately periodic
(see also \cite{Kurka2003} and \cite{AraujoBruyere2005}).
\index{names}{K{\r{u}}rka, Petr}\index{names}{Araujo, Isabel}%
\index{names}{Bruy\`ere, V\'eronique}%
For similar theorems   characterizing   purely substitutive  Sturmian words,  see the references in \cite{Lothaire2002}\index{names}{Lothaire, M.}.
In particular, it is proved in \cite{Allauzen1998} that,
\index{names}{Allauzen, Cyril}
for $0<\alpha<1$
the caracteristic sequence $c_\alpha$ 
of slope $\alpha$ is purely substitutive if
and only if $\alpha$ is quadratic and its conjugate is $>1$.

It is not surprising that the condition on $\alpha$
for $c_\alpha$ to be purely substitutive is more
restrictive than for the shift of slope $\alpha$ to 
be substitutive.
Indeed, if $\alpha$ is quadratic, the Sturmian shift of slope
$\alpha$ contains a purely substitutive sequence, but it need
not be the characterstic sequence $c_\alpha$. 
\marginpar{A revoir. Le shift de pente alpha est-il purement
substitutif si et seulement si $c_\alpha$ est purement substitutif?}

%Specular
\subsection{Specular shifts}
The notion of specular shift was introduced in \cite{BertheDeFeliceDelecroixDolceLeroyPerrinReutenauerRindone2017}. The idea of considering
laminary sets is from~\cite{HilionCoulboisLustig2008}
(see also~\cite{LopezNarbel2013}).
\index{names}{Lopez, Michel}\index{names}{Narbel, Pierre}%
\index{names}{Hilion, Arnaud}\index{names}{Coulbois, Thierry}%
\index{names}{Lustig, Martin}%

For a proof of Kurosh subgroup Theorem, see \cite{MagnusKarrassSolitar2004}.
\index{names}{Magnus, Wilhelm}\index{names}{Karrass, Abraham}%
\index{names}{Solitar, Donald}%
Specular groups are 
characterized by the property
of their  Cayley graphs to be regular (see~\cite{Harpe2000}).
\index{names}{Harpe@de la Harpe, Pierre}%
See also \cite{Harpe2000} concerning the notion of virtually free group.
Actually, specular groups can be studied as groups acting on trees as developed in the Bass-Serre theory \citep{Serre2003}.
\index{names}{Serre, Jean-Pierre}\index{names}{Bass, Hyman}

A group having a free subgroup of finite index is called \emph{virtually
free}.\index{subject}{virtually free group}\index{subject}{group!virtually free}
 On the other hand, a finitely generated group is
said to be \emph{context-free}\index{subject}{context-free!group}
 if, for some presentation, the set
 of words equivalent to $1$ is a context-free language.
 \index{subject}{context-free!language}%
By Muller and Schupp's theorem \citep{MullerSchupp1983},
a finitely generated group is virtually free if and only if
it is context-free. 
Thus a specular group is context-free. One may
verify this directly as follows. A
context-free grammar\index{subject}{context-free!grammar}
generating the words equivalent to $1$ for the
natural presentation of a specular group $G=G_\theta$ is the
grammar with one nonterminal symbol $\sigma$ and the rules
\begin{equation}
\sigma\rightarrow a\sigma a^{-1}\sigma\quad (a\in A),\quad \sigma\rightarrow 1.
\label{eqGrammarContextFree}
\end{equation}

The proof that the grammar given by Equation~\eqref{eqGrammarContextFree}
generates the set of words
equivalent to $1$ is similar to that used in~\cite{Berstel1979}
\index{names}{Berstel, Jean}%
for the so-called Dyck-like languages.

Theorem \ref{theoremDecodingEven} is the counterpart for minimal specular shifts
of the main result of~\cite[Theorem 6.1]{BertheDeFeliceDolceLeroyPerrinReutenauerRindone2015b} asserting that the family of uniformly recurrent tree sets 
of characteristic $1$ is closed
under maximal bifix decoding.

Theorem~\ref{propositionNew} is proved in~\cite[Theorem 3]{DolcePerrin2017}.

The definition of mixed return words
comes from the fact that, when $S$ is the natural coding
of a linear involution, we are interested in the transformation induced
on $I_w\cup\sigma_2(I_w)$ (see~\cite{BertheDelecroixDolcePerrinReutenauerRindone2017b}).
The natural coding of a point in $I_w$ begins
with $w$ while the natural coding of a point $z$ in $\sigma_2(I_w)$ `ends'
with $w^{-1}$ in the sense that the natural coding of  $T^{-|w|}(z)$ 
begins with $w^{-1}$.

A geometric proof and interpretation  of Theorem~\ref{theoremFirstMixed} is
 given in~\cite{BertheDelecroixDolcePerrinReutenauerRindone2017b}.
It is shown that the set of mixed return words 
are a symmetric basis of a fundamental group corresponding
to a surface built above the linear involution.

%%%%%%%%%%%%%%%%%%%%%%%%%%%%%%%%%%%%%%%ù
\chapter{Interval exchange transformations}
%%%%%%%%%%%%%%%%%%%%%%%%%%%%%%%%%
\label{chapterIET}
In the chapter, we study a class of dynamical systems
obtained by iterating simple geometric transformations
on an interval. These transformations, called interval
exchange form a classical family of dynamical systems.
They can be seen as a generalization of the rotations
of the circle, already met several times.
They define by a natural coding shift spaces spaces
which turn out to be dendric shifts. Thus
we obtain, after Arnoux-Rauzy shifts a large
class of dendric shifts on which we will describe
more precisely return words and $\Sa$-adic representations.

In the first part, we introduce interval exchange transformations.
We define the natural representation of an interval exchange transformation,
which is shown to be dendric shifts (Proposition~\ref{propositionPlanarDendric}). We introduce the notion of regular interval exchange transformation
and prove that it implies the minimality of the transformation 
(Theorem~\ref{theoremIDOC}).
We develop the notion of Rauzy induction and characterize the
subintervals reached by iterating the transformation (Theorem~\ref{theo:birauzy1}).
We generalize Rauzy induction to a two-sided version and
characterize the intervals reached by this more general transformation
(Theorem~\ref{theo:birauzy2}). We link these transformations
with automorphisms of the free group (Theorem~\ref{theo:inductionbi}).
We also relate these results with the theorem of Boshernizan and Carroll
giving a finiteness condition on the systems induced by an interval
exchange when the lengths of the intervals belong to a quadratic field
(Theorem \ref{theo:quadratic}).

In a second part, we present linear involutions, which form a larger
family by adding the possibility of symmetries in addition
to translations. We introduce the notion of connection
and relate it to the minimality of the
transformation (Proposition~\ref{propBL}). We show that the natural coding of a linear
involution without connection is specular (Theorem~\ref{theoremInvolutionSpecular}).

%%%%%%%%%%%%%%%%%%%%%%%%%%%%%%%%%%%%%%%%%%%%%%
\section{Interval exchange transformations}
%%%%%%%%%%%%%%%%%%%%%%%%%%%%%%%%%%%%%%%%%%%%%%
\label{ch5:subsection:iet}

A \emph{semi-interval}\index{subject}{semi-interval} is a nonempty
subset of the real line of the form $[\ell,r)=\{z\in\R\mid\ell\le z<r\}$.
Thus it a left-closed and right-open interval. For two semi-intervals
$\Delta,\Gamma$, we denote $\Delta<\Gamma$ if $x<y$ for
every $x\in\Delta$ and $y\in\Gamma$. A partition
$\{I_1,\ldots,I_k\}$ of the semi-interval $I$ indexed
by $1,\ldots,k$ is \emph{ordered}
\index{subject}{ordered!partition} if $\Delta_1<\ldots<\Delta_k$.

Let $\{\Delta_1, \ldots , \Delta_k \}$
\index{symbols}{Delta@$(\Delta_1, \ldots , \Delta_k )$} be an
ordered partition of the semi-interval $I=[\ell,r)$ into $k\geq 2$ disjoint semi-intervals.
A {\em $k$-interval exchange transformation}\index{subject}{interval exchange!transformation}\index{subject}{transformation!interval exchange} on $I$ is an onto map $T : [\ell,r) \to [\ell,r)$ where $T:\Delta_i \to [\ell,r)$ is a translation.

For every $i$ with $1\le i\le k$, there is a real number $\alpha_i$
such that
\begin{displaymath}
Tx=x+\alpha_i
\end{displaymath}
for every $x\in\Delta_i$. The numbers $\alpha_i$ are called the
\emph{translation values}\index{subject}{translation values} of $T$.

\begin{figure}[hbt]
\centering
\tikzset{node/.style={circle,draw,minimum size=0.2cm,inner sep=0pt}}
\tikzset{title/.style={circle,minimum size=0.1cm,inner sep=0pt}}
\begin{tikzpicture}

\node[node](0H)at(0,1){};\node[title]at(0,1.4){$d_1=\ell$};
\node[node](1H)at(2,1){};\node[title]at(2,1.4){$d_2$};
\node[node](2H)at(3,1){};\node[title]at(3,1.4){};
\node[node](21H)at(4,1){};\node[title]at(4,1.4){$d_{\pi(1)}$};
\node[node](22H)at(6.5,1){};\node[title]at(6.5,1.4){};
\node[node](3H)at(7.5,1){};\node[title]at(7.5,1.4){$d_k$};
\node[node](4H)at(10,1){};\node[title]at(10,1.4){$r$};

\node[node](0B)at(0,0){};
\node[node](1B)at(2.5,0){};
\node[node](2B)at(4.5,0){};
\node[node](3B)at(8,0){};
\node[node](4B)at(10,0){};

\node[title]at(0,.5){$\vdots$};\node[title]at(10,.5){$\vdots$};

\draw[above,line width=1pt](0H)edge node{$\Delta_1$}(1H);
\draw[above,line width=1pt](1H)edge node{$\Delta_2$}(2H);
\draw[dotted](2H)edge node{}(21H);
\draw[above,line width=1pt,color=red](21H)edge node{$\Delta_{\pi(1)}$}(22H);
\draw[dotted](22H)edge node{}(3H);
\draw[above,line width=1pt](3H)edge node{$\Delta_k$}(4H);
\draw[below,line width=1pt,color=red](0B)edge node{$T\Delta_{\pi(1)}$}(1B);
\draw[below,line width=1pt](1B)edge node{$T\Delta_{\pi(2)}$}(2B);
\draw[dotted](2B)edge node{}(3B);
\draw[below,line width=1pt](3B)edge node{$T\Delta_{\pi(k)}$}(4B);
\draw[->](5,1)edge node{}(1.5,0){};%\draw[->](4,1)edge node{}(1,0);
%\draw[->](5.75,1)edge node{}(2.75,0);
\end{tikzpicture}
\caption{A $k$-interval exchange transformation.}\label{figurekinterval}
\end{figure}
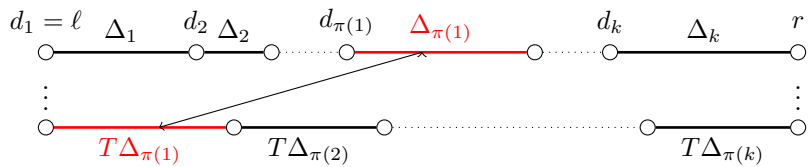

We denote by $\pi$ the permutation of $\{ 1, \ldots , k \}$ such that 
  $T\Delta_{\pi(1)},\ldots,T\Delta_{\pi(k)}$
is an ordered partition (see Figure~\ref{figurekinterval}). 

If $\lambda_i$ denotes the length of the interval $\Delta_i$ and $\lambda=(\lambda_1,\ldots,\lambda_k)$, the transformation is defined by $\lambda$ and the permutation 
$\pi$. We
denote $T=T_{\lambda,\pi}$.

One may imagine an interval exchange transformation as the following
`physical' operation. Break the semi-interval $I$ into $k$ pieces
and rearrange them in a different order.

Such a transformation is not continuous and thus does not fit well
in the framework of topological dynamical systems. It is actually
a measure preserving transformation and thus an interval exchange
transformation is a measure theoretic dynamical system.

We shall see below how
the interval can be modified to obtain a topological dynamical system.

For $k=2$, an interval exchange transformation 
on $I=[0,1)$ is a rotation.
We will often consider the semi-interval interval $[0,1)$ to simplify
the notation. Indeed,
for $0<\alpha<1$,
set $\Delta_1=[0,1-\alpha)$ and $\Delta_2=[1-\alpha,1)$. The
corresponding interval exchange transformation $T$ is
$T(x)=x+\alpha\bmod 1$. In this particular case, the
transformation is continuous provided one identifies the
two endpoints of the interval (see Example~\ref{exampleRotations}).

An interval exchange transformation is invertible
since a translation is one-to-one. Its inverse
is again an interval exchange transformation,
on the intervals $\Delta'_i=T\Delta_i$.

A power $T^n$ of an interval exchange transformation is again an interval
exchange transformation. Since we have seen that $T^{-1}$ is also 
an interval exchange transformation, we may assume $n>0$.
Consider  the nonempty sets of the form
\begin{displaymath}
\Delta_{i_0,\ldots,i_{n-1}}^{(n)}=\Delta_{i_0}\cap T^{-1}\Delta_{i_1}\cap\ldots\cap T^{-n+1}\Delta_{i_{n-1}}.
\end{displaymath}
These sets are  semi-intervals because for every $i,j$, the set
\begin{displaymath}
\Delta_i\cap T^{-1}\Delta_j=T^{-1}(\Delta_j\cap T\Delta_i)=T^{-1}(\Delta_j\cap \Delta'_i)
\end{displaymath}
is a semi-interval or empty and next
\begin{eqnarray*}
\Delta_{i_0}\cap T^{-1}\Delta_{i_1}\cap\ldots\cap T^{-n+1}\Delta_{i_n}&=&\\
\Delta_{i_0}\cap T^{-1}(\Delta_{i_1}\cap T^{-1}(\Delta_{i_2}\cap&\cdots&\cap
T^{-1}(\Delta_{i_{n-2}}\cap T^{-1}\Delta_{i_{n-1}})\cdots).
\end{eqnarray*}
Thus $T^n$ is an exchange of the nonempty intervals $\Delta_{i_0,\ldots,i_{n-1}}^{(n)}$.

\begin{example}
A $3$-interval exchange transformation is represented in
Figure~\ref{figure3interval}. 
The associated permutation is the cycle $\pi=(123)$.
%Let $\alpha,\beta,\gamma$ be positive real numbers such that
%$\alpha+\beta+\gamma=1$.
\begin{figure}[hbt]
\centering
\centering
\tikzset{node/.style={circle,draw,minimum size=0.2cm,inner sep=0pt}}
\begin{tikzpicture}

\node[node,fill=red](0H)at(0,1){};
\node[node,fill=blue](1H)at(3,1){};
\node[fill=green,node](2H)at(5,1){};
\node[node](3H)at(6.5,1){};

\node[node,fill=blue](0B)at(0,0){};
\node[node,fill=green](1B)at(2,0){};
\node[node,fill=red](2B)at(3.5,0){};
\node[node](3B)at(6.5,0){};

\draw[red,line width=2pt](0H)edge node{}(1H);
\draw[blue,line width=2pt](1H)edge node{}(2H);
\draw[green,line width=2pt](2H)edge node{}(3H);
\draw[blue,line width=2pt](0B)edge node{}(1B);
\draw[green,line width=2pt](1B)edge node{}(2B);
\draw[red,line width=2pt](2B)edge node{}(3B);
\draw[->](1.5,1)edge node{}(5,0){};\draw[->](4,1)edge node{}(1,0);
\draw[->](5.75,1)edge node{}(2.75,0);
\end{tikzpicture}
\caption{A $3$-interval exchange transformation.}\label{figure3interval}
\end{figure}
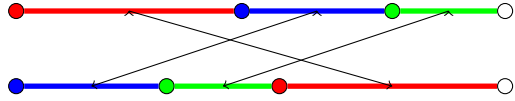
\end{example}
\subsection{Minimal interval exchange transformations}
The \emph{separation points}\index{subject}{separation points}
 of an interval exchange are the
left end points $d_1=0,d_2,\ldots,d_k$
of the semi-intervals $\Delta_i$ for $1\le i\le k$.
We denote by $\Sep(T)$ the set of separation points of $T$.
\index{symbols}{Sep@$\Sep(T)$}

An interval exchange transformation $T$ on $I=[0,1)$ is called
\emph{regular}
\index{subject}{regular!interval exchange}%
\index{subject}{interval exchange!regular}%
 if the orbits of the nonzero separation
points  are infinite and disjoint.
Note that the orbit of $0$ cannot be disjoint of the others,
since one has $T(d_i)=0$ for some $i$ with $1\le i\le k$,
but that the orbit of $0$ is infinite when $T$ is regular.

Equivalently, an interval exchange transformation is regular
if there is no triple $(i,j,n)$ for $2\le i,j\le k$, 
and $n> 0$ such that $T^nd_i=d_j$. Such a triple
is called a \emph{connection}.
\index{subject}{connection!of interval exchange}%
Indeed, a connection $(i,j,n)$ with $i=j$ 
corresponds to a finite orbit and with $i\ne j$
to intersecting orbits.

As an example, a rotation of angle $\alpha$ is regular
if and only if $\alpha$ is irrational. Indeed, $\alpha$ is
rational if and only if the orbit of $\alpha$ is finite
(in which case all orbits are finite).

We say that an interval exchange is \emph{minimal}
\index{subject}{minimal!interval exchange}%
\index{subject}{interval exchange!minimal} if the orbit of every point
is dense.
\begin{theorem}[Keane]\label{theoremIDOC}
A regular  interval exchange transformation is minimal. 
\end{theorem}
The converse is not true. Indeed, consider the the rotation of 
irrational angle
$\alpha$ on $I=[0,1)$ as a $3$-interval exchange with the partition $[0,1-2\alpha),
[1-2\alpha,1-\alpha),[1-\alpha,0)$ (see Figure~\ref{figure3interval2}). The transformation is minimal,
as any rotation of irrational angle but $T(1-2\alpha)=1-\alpha$
and thus there is a connection.

We first prove the following result. Let $T$
be an interval exchange on $I=[0,1)$. We shall consider the transformation $T_1$
induced by  $T$ on a 
semi-interval $I_1=[a,b)\subset[0,1)$. It is defined by $T_1(z)=T^{n(x)}(z)$
where $n(z)$ is the least integer $n\ge 1$ such that $T^n(z)\in I_1$.
The return time $n(z)$ exists, even if $T$ is not minimal,
 by the Poincar\'e Recurrence Theorem,
because $T$ preserves the Lebesgue measure.
\index{subject}{Lebesgue!measure}%
\index{subject}{Poincar\'e Recurrence Theorem}%
\index{subject}{Theorem!Poincar\'e Recurrence}%
\index{names}{Poincar\'e, Henry}%
Such an induction on  a semi-interval is called a \emph{Rauzy induction}
and we will have more to say on this later.
\index{subject}{Rauzy!induction}\index{names}{Rauzy, G\'erard}%
\begin{theorem}[Rauzy]\label{lemmaRauzyInduction}
Let $T$ be a $k$-interval exchange transformation 
on $I=[0,1)$ and let $X_1=[a,b)\subset[0,1)$
be a semi-interval.  The transformation $T_1$ induced by $T$
on $X_1$ is a $k_1$-interval exchange transformation with $k_1\le k+2$.
\end{theorem}
\begin{proof}
Consider the set $Y=\{a,b,d_2,\ldots,d_k\}$ where $d_2,\ldots,d_k$ are the
nonzero separation points. For $y\in Y$, let $s(y)$ be the least integer
$s\ge 0$ if it exists such that $T^{-s(y)}(y)\in(a,b)$. The points $T^{-s(y)}$
divide the semi-interval $[a,b)$ in $k_1\le k+2$ semi-intervals
$\Delta'_1,\ldots,\Delta'_{k_1}$. For each semi-interval $\Delta'_i$
consider the number $n_i$ which is the least $n\ge 1$ such 
that $T^n\Delta'_i\cap [a,b)\ne\emptyset$. 

Then, for $1\le p\le n_i$ the transformations $T^p$ are continuous
on $\Delta'_i$ and we have $T^{n_i}\Delta'_i\subset [a,b)$.
Indeed, otherwise, for some $p$ with $1\le p\le n_i-1$, the semi-interval
$T^p\Delta'_i$ contains one of the points $y\in Y$ and then $s(y)=p$
so that the point $T^{-p}y$ would lie within $\Delta'_i$, a contradiction
with the definition of the semi-intervals $\Delta'_i$. This
shows that $n_i$ is the return time of all $x\in\Delta'_i$
and thus that $T_1$ is the interval exchange of $[a,b)$
corresponding to the intervals $\Delta'_i$.
\end{proof}
\begin{proofof}{of Theorem~\ref{theoremIDOC}}
Let us first show that $T$ has no periodic point. Assume that $T^n$ has
a fixed point, that is a point $x\in I$ such that 
$T^nx=x$. Since $T^n$ is an interval exchange transformation
on the intervals $\Delta_{i_0,\ldots,i_{n-1}}^{(n)}$, one of these intervals is formed
of fixed points. But the left end of this interval is in the
orbit of some of the separation points $d_i$. This contradicts the hypothesis
that the orbit of all nonzero separation points is infinite unless $i=1$
and $n=1$.
But in this case, the interval $\Delta_1$ is fixed by $T$, which implies
that $Td_j=d_2$ for some $j\ge 2$, a contradiction again with the
hypothesis.

Suppose now that the orbit $O(z)$ of some $z\in[0,1)$ is not dense. We can find
a semi-interval $[a,b)$ disjoint  from $O(z)$. Let $T_1$
be the transformation induced by $T$ on $[a,b)$. By Theorem~\ref{lemmaRauzyInduction}, $T_1$ is an interval exchange on semi-intervals $\Delta'_j$. Let
$n_j$ denote the return time to $\Delta'_j$. Set
\begin{displaymath}
F=\bigcup_j\bigcup_{n=0}^{n_j-1}T^n\Delta'_j.
\end{displaymath}
The set $F$ can be written as a union of a finite number of non intersecting semi-intervals.
Hence its connected components are also semi-intervals, say $F_s$
and their number is finite. Let $G$ be the set of the left end
points of all the semi-intervals $F_s$. By definition, the set
$F$ is invariant by $T$ and therefore, for every $g\in G$, we have either
$Tg\in G$ or the point $g$ is a point of discontinuity of $T$
or $g=0$, that is $x=d_i$ for some $i$ with $1\le i\le k$.
Since $G$ is finite and since $T$ has no periodic points, there
is for every $g\in G$ some $n\ge 0$ such that $T^nx=d_i$ with $1\le i\le k$.
Similarly, there is for every $x\in G$ an integer $m> 0$ 
such that $T^{-m}x=d_j$ with $2\le i\le k$. Then, we have 
\begin{displaymath}
T^nx=T^{n+m}d_j=d_i.
\end{displaymath}
Since the orbits of $d_2,\ldots,d_k$ are disjoint, the only possibility
is $i=1$. Then $d_j=T^{-1}d_1$ and thus $m=1$, $n=0$ and $x=0$. Thus
$G=\{0\}$, which implies $F=[0,1)$, a contradiction with
the fact that, by construction $F\cap O(x)=\emptyset$.
\end{proofof}
%The converse of Theorem~\ref{theoremIDOC} is false. For example,
% the rotation of irrational angle $\alpha$ on $[0,1)$
% is a $3$-interval exchange on the intervals 
%$[0,1-2\alpha,1-\alpha,1)$. It is minimal since $\alpha$ is irrational
%but the orbits of the separation points are not disjoint since
%$T(1-2\alpha)=1-\alpha$.

The following necessary condition for minimality of an interval exchange transformation is useful.
A permutation $\pi$ of an ordered set $A$ is called \emph{decomposable}
\index{subject}{decomposable permutation}%
\index{subject}{permutation!decomposable} if there exists an element $b \in A$ such that the set $B$ of elements strictly less than $b$ is nonempty and such that $\pi(B) = B$.
Otherwise it is called \emph{indecomposable}.
\index{subject}{indecomposable permutation}%
\index{subject}{permutation!indecomposable}%
If an interval exchange transformation $T=T_{\lambda, \pi}$ is minimal, the permutation $\pi$ is indecomposable.
Indeed, if $B$ is a set as above, the set of orbits of the points in the set $S = \cup_{a \in B}I_a$ is closed and strictly included in $[\ell,r[$.
The following example shows that the indecomposablity of $\pi$ is not sufficient for $T$ to be minimal.

\begin{example}
Let $A = \{a, b, c\}$ and $\lambda$ be such that $\lambda_a = \lambda_c$.
Let $\pi$ be the transposition $(ac)$.
Then $\pi$ is indecomposable but $T_{\lambda, \pi}$ is not minimal since it is the identity on $I_b$.
\end{example}

The iteration of a $k$-interval exchange transformation is, in general, an interval exchange transformation operating on a larger number of semi-intervals.

\begin{proposition}
\label{pro:tn}
Let $T$ be a regular $k$-interval exchange transformation.
Then, for any $n \geq 1$, $T^n$ is a regular $n(k-1)+1$-interval exchange transformation.
\end{proposition}
\begin{proof}
Since $T$ is regular, the set $\cup_{i=0}^{n-1} T^{-i}(d)$ where $d$ runs over the set of $s-1$ nonzero separation points of $T$ has $n(k-1)$ elements.
These points partition the interval $[\ell,r[$ in $n(k-1)+1$ semi-intervals on which $T$ is a translation.
\end{proof}

We close this subsection with a lemma which will be useful in Section~\ref{subsec:morphic}.

\begin{lemma}
\label{lem:distance}
Let $(T,I)$ be a minimal interval exchange transformation. For every $N > 0$ there exists an $\varepsilon > 0$ such that for every $z \in I$ and for every $n>0$, one has
$$ \left| T^n(z) - z \right| < \varepsilon \quad \Longrightarrow \quad n \geq N.$$
\end{lemma}
\begin{proof}
Let $\alpha_1, \alpha_2, \ldots, \alpha_s$ be the translation values of $T$.
For every $N>0$ it is sufficient to choose
$$\varepsilon = \min \left\{ \left| \textstyle{\sum_{i_j = 1}^M \alpha_{i_j}} \right|  \; \mid \;\;\; 1 \leq i_j \leq s \; \mbox{ and } \; M \leq N \right\}.$$
\end{proof}

\subsection{Natural coding}\label{subsec:nc}
To every interval exchange transformation $T$ on $I=[\ell,r)$, we may associate a shift
space called its \emph{natural representation}
\index{subject}{natural!representation!of interval exchange}%
\index{subject}{interval exchange!natural representation}%
 defined as follows.

 Set $A=\{1,2,\ldots,k\}$
and let $\gamma:I\to A^\Z$ be the map defined by $x=\gamma(z)$ 
with $x=(x_n)$ and
\begin{displaymath}
x_n=a\quad\mbox{ if\ $T^n(z)\in \Delta_a$.}
\end{displaymath}
The two-sided sequence $y=\gamma(z)$ is called the \emph{natural coding}
\index{subject}{natural!coding!of interval exchange} of $z$.
The natural representation of $T$, denoted $X(T)$,
 is the closure in $A^\Z$ of $\gamma([\ell,r))$. 

An \emph{interval exchange shift}\index{subject}{interval exchange!shift}
\index{subject}{shift space!interval exchange}
is the natural representation $X(T)$ of an interval exchange $T$.
It is said to be regular if the interval exchange is regular.
We denote  $\cL(T)=\cL(X(T))$.

The map $\gamma$ satisfies $\gamma\circ T=S\circ\gamma$ where $S$ denotes
as usual the shift. Thus the natural representation of a minimal interval
exchange transformation is a minimal shift space and the
converse is also true.
Note that if $T$ is minimal, the map $\gamma$ is injective. Indeed,
if $z\ne z'$, there is for every $a$
an $n$ such that $T^nz\in\Delta_a$ but $T^n z'\notin\Delta_a$.

\begin{example}
\label{ex:t}\label{ex:returns}
\label{exampleRotation2alpha}\label{ex:2alpha}
Consider the interval exchange transformation $T$ represented in Figure
\ref{figure3interval2} with $\alpha=(3-\sqrt{5})/2$.
\begin{figure}[hbt]
\centering
\tikzset{node/.style={circle,draw,minimum size=0.1cm,inner sep=0pt}}
\tikzset{title/.style={circle,minimum size=0.1cm,inner sep=0pt}}

\begin{tikzpicture}

\node[node,color=red](0h)at(0,1){};\node[title]at(0,1.4){$0$};
\node[node,color=blue](1-2alpha)at(2.36,1){};\node[title]at(2.36,1.4){$1-2\alpha$};
\node[node,color=green](1-alpha)at(6.18,1){};\node[title]at(6.18,1.4){$1-\alpha$};
\node[node](1h)at(10,1){};\node[title]at(10,1.4){$1$};
\draw[node,color=red,line width=1,above](0h)edge node{$a$}(1-2alpha);
\draw[node,color=blue,line width=1,above](1-2alpha)edge node{$b$}(1-alpha);
\draw[node,color=green,line width=1,above](1-alpha)edge node{$c$}(1h);

\node[node,color=blue](0b)at(0,0){};\node[title]at(0,0.4){$0$};
\node[node,color=green](alpha)at(3.82,0){};\node[title]at(3.82,0.4){$\alpha$};
\node[node,color=red](2alpha)at(7.64,0){};\node[title]at(7.64,0.4){$2\alpha$};
\node[node](1b)at(10,0){};\node[title]at(10,0.4){$1$};
\draw[color=blue,line width=1,below](0b)edge node{$b$}(alpha);
\draw[color=green,line width=1,below](alpha)edge node{$c$}(2alpha);
\draw[color=red,line width=1,below](2alpha)edge node{$a$}(1b);
\end{tikzpicture}
\caption{The rotation of angle $2\alpha$.}\label{figure3interval2}
\end{figure}
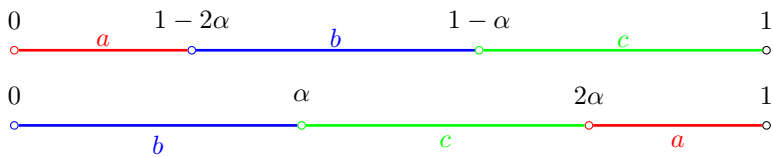
It is the rotation $x\mapsto x+2\alpha$ represented as a regular interval
exchange on 
three semi-intervals $\Delta_a=[0,1-2\alpha)$,
$\Delta_b=[1-2\alpha,1-\alpha)$, and $\Delta_c=[1-\alpha,1)$.

Since $T$ is minimal, the set $\cL(T)$ is uniformly recurrent.
%In Subsection~\ref{subsec:morphic} we will show that the set $F = F(T)$ is the set of factors of the fixpoint of a primitive morphism. 
The words of length at most $6$ of the set $\cL(T)$ are represented in Figure~\ref{fig:setf}.

\begin{figure}[hbt] 
\centering
\tikzset{node/.style={circle,draw,minimum size=0.2cm,inner sep=0pt}}
\tikzset{title/.style={circle,minimum size=0.1cm,inner sep=0pt}}
\begin{tikzpicture}
\node[node](1)at(0,3){};
\node[node](a)at(1,4){};
\node[node](b)at(1,3){};
\node[node](c)at(1,2){};

\node[node](ac)at(2,4.5){};
\node[node](ba)at(2,3.5){};
\node[node](bb)at(2,2.5){};
\node[node](cb)at(2,2){};
\node[node](cc)at(2,1){};

\node[node](acb)at(3,5){};
\node[node](acc)at(3,4.5){};
\node[node](bac)at(3,3.5){};
\node[node](bba)at(3,2.5){};
\node[node](cba)at(3,2){};
\node[node](cbb)at(3,1.5){};
\node[node](ccb)at(3,1){};

\node[node](acbb)at(4,5.5){};
\node[node](accb)at(4,4.5){};
\node[node](bacb)at(4,4){};
\node[node](bacc)at(4,3.5){};
\node[node](bbac)at(4,2.5){};
\node[node](cbac)at(4,2){};
\node[node](cbba)at(4,1.5){};
\node[node](ccba)at(4,1){};
\node[node](ccbb)at(4,.5){};

\node[node](acbba)at(5,5.5){};
\node[node](accba)at(5,5.1){};
\node[node](accbb)at(5,4.5){};
\node[node](bacbb)at(5,4){};
\node[node](baccb)at(5,3.5){};
\node[node](bbacb)at(5,3){};
\node[node](bbacc)at(5,2.5){};
\node[node](cbacc)at(5,2){};
\node[node](cbbac)at(5,1.5){};
\node[node](ccbac)at(5,1){};
\node[node](ccbba)at(5,.5){};

\node[node](acbbac)at(6,5.5){};
\node[node](accbac)at(6,5.1){};
\node[node](accbba)at(6,4.5){};
\node[node](bacbba)at(6,4){};
\node[node](baccba)at(6,3.7){};
\node[node](baccbb)at(6,3.2){};
\node[node](bbacbb)at(6,3){};
\node[node](bbaccb)at(6,2.5){};
\node[node](cbaccb)at(6,2){};
\node[node](cbbacb)at(6,1.7){};
\node[node](cbbacc)at(6,1.3){};
\node[node](ccbacc)at(6,1){};
\node[node](ccbbac)at(6,.5){};

\draw[above](1)edge node{$a$}(a);
\draw[above](1)edge node{$b$}(b);
\draw[above](1)edge node{$c$}(c);
\draw[above](a)edge node{$c$}(ac);
\draw[above](b)edge node{$a$}(ba);
\draw[above](b)edge node{$b$}(bb);
\draw[above](c)edge node{$b$}(cb);
\draw[above](c)edge node{$c$}(cc);
\draw[above](ac)edge node{$b$}(acb);
\draw[above](ac)edge node{$c$}(acc);
\draw[above](ba)edge node{$c$}(bac);
\draw[above](bb)edge node{$a$}(bba);
\draw[above](cb)edge node{$a$}(cba);
\draw[above,near end](cb)edge node{$b$}(cbb);
\draw[above](cc)edge node{$b$}(ccb);
\draw[above](acb)edge node{$b$}(acbb);
\draw[above](acc)edge node{$b$}(accb);
\draw[above](bac)edge node{$b$}(bacb);
\draw[above,near end](bac)edge node{$c$}(bacc);
\draw[above](bba)edge node{$c$}(bbac);
\draw[above](cba)edge node{$c$}(cbac);
\draw[above](cbb)edge node{$a$}(cbba);
\draw[above](cbba)edge node{$c$}(cbbac);
\draw[above](ccb)edge node{$a$}(ccba);
\draw[above,near end](ccb)edge node{$b$}(ccbb);
\draw[above](acbb)edge node{$a$}(acbba);
\draw[above](accb)edge node{$a$}(accba);
\draw[above,near end](accb)edge node{$b$}(accbb);
\draw[above](bacb)edge node{$b$}(bacbb);
\draw[above](bacc)edge node{$b$}(baccb);
\draw[above](bbac)edge node{$b$}(bbacb);
\draw[above,near end](bbac)edge node{$c$}(bbacc);
\draw[above](cbac)edge node{$c$}(cbacc);
\draw[above](cbba)edge node{$c$}(cbbac);
\draw[above](ccba)edge node{$c$}(ccbac);
\draw[above](ccbb)edge node{$c$}(ccbba);
\draw[above](acbba)edge node{$c$}(acbbac);
\draw[above](accba)edge node{$c$}(accbac);
\draw[above](accbb)edge node{$a$}(accbba);
\draw[above](bacbb)edge node{$a$}(bacbba);
\draw[above,near end](baccb)edge node{$a$}(baccba);
\draw[above,near end](baccb)edge node{$b$}(baccbb);
\draw[above](bbacb)edge node{$b$}(bbacbb);
\draw[above](bbacc)edge node{$b$}(bbaccb);
\draw[above](cbacc)edge node{$b$}(cbaccb);
\draw[above](cbbac)edge node{$b$}(cbbacb);
\draw[above,near end](cbbac)edge node{$c$}(cbbacc);
\draw[above](ccbac)edge node{$c$}(ccbacc);
\draw[above](ccbba)edge node{$c$}(ccbbac);
\end{tikzpicture}
\caption{The words of length $ \leq 6$ of the set $\cL(T)$.}
\label{fig:setf}
\end{figure}
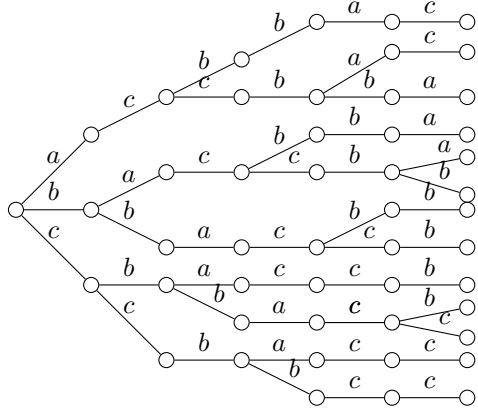

The shift $X=X(T)$ is actually a primitive substitution shift.
Indeed, since $T$ is the rotation of angle $2\alpha$,
the shift $X(T)$ is the coding by non overlapping blocks
of length $2$ of the Fibonacci shift, which is the natural
coding of the rotation of angle $\alpha$. But
the Fibonacci morphism $\varphi:a\to ab, b\to a$ is such that
$\varphi^3:a\to abaab,b\to aba$ sends words of even length
to words of even length. Using the coding $a\to aa,b\to ab,c\to ba$,
the action of $\varphi^3$ on $\{a,b,c\}$ is
$\psi:a\to baccb,b\to bacc,c\to bacb$. Thus $X=X(\psi)$. 
\end{example}

\subsection{Cantor version of interval exchange}

The natural coding can be made  continuous at the cost
of modifying the space $I=[0,1]$ as described below.

Suppose that $T$ is minimal, that is, all its orbits are dense in $[0,1)$.
Let $\Sep(T) = \{ d_1,\ldots , d_k \}$ be the set of separation points
 and set $\mathcal{O} (T) = \{ T^j d \mid j\in \Z , \  d\in \Sep (T) \}$.
We define
$$
X = \left( [0,1) \setminus \mathcal{O} (T) \right) \bigcup \left\{ z^{-}, z^{+} \mid  z\in \mathcal{O} (T) \right\}  
$$
where $0^- = 1$. 
Defining $x<z^- < z^+<y$ for all $z\in \mathcal{O} (T)$ 
and $x,y\in [0,1)\setminus\mathcal{O}(T)$ such that $x<z<y$
(with the exception of $0^-\geq x$ for all $x\in X$), this 
extends the natural order on $[0,1)$ to $X$.
Endowed with the topology of intervals, $X$ is a Cantor space because $\mathcal{O} (T)$ is dense in $[0,1)$.

Let $F : X\to X$ defined by $F (x) = T (x)$ if $x\in [0,1) \setminus \mathcal{O} (T)$ and $F (z^\epsilon) = T (z)^\epsilon$ if $z\in \mathcal{O} (T)$ where $\epsilon \in \{ +,-\}$.
The pair $(X,F)$ is a minimal Cantor dynamical system, we will refer to as the {\em Cantor version of the interval exchange} $T$.\index{subject}{interval exchange!Cantor version}\index{subject}{Cantor!version of interval exchange}

Let $\phi : X \to [0,1)$  be defined by $\phi (z^+)=\phi (z^-)=z$ and $\phi(z)=z$ when $z\not \in \mathcal{O} (T)$. 
This is an onto continuous map.
It is one-to-one everywhere except on a countable set of points.
Moreover,   $\phi \circ F = T\circ \phi$.

Let $T$ be a $k$-interval exchange transformation corresponding to
semi-intervals $\Delta_i$. Let $A=\{1,2,\ldots,k\}$ and let
$X$ be the natural representation of $T$. For $w\in A^*$,  denote by $I(w)$
\index{symbols}{I@$I(w)$}%
the set defined by $I(\varepsilon)=I$  and by
\begin{equation}
I(au)=\Delta_a\cap T^{-1}(I(u)).\label{eqI(au)}
\end{equation}
It can be verified that every nonempty $I(w)$ is a semi-interval
(Exercise~\ref{exerciseIw}) and that $I(a)=\Delta_a$.
Note that $I(ua)=I(u)\cap T^{-|u|}\Delta_a$. Note also that
\begin{equation}
 b\in R(w) \mbox{ if and only if }I(w)\cap T^{-|w|}\Delta_b\ne\emptyset
\label{eqConditionR(w)}
\end{equation}
and 
\begin{equation}
a\in L(w)\mbox{ if and only if }T\Delta_a\cap
I(w)\ne\emptyset
\label{eqConditionL(w)}
\end{equation}
where $R(w)=\{a\in A\mid wa\in\cL(T)\}$ and $L(w)=\{a\in A\mid aw\in\cL(T)\}$.
Condition \eqref{eqConditionR(w)} holds because $I(wb)=I(w)\cap T^{-|w|}\Delta_b$ and condition
\eqref{eqConditionL(w)} because $I(aw)=\Delta_a\cap T^{-1}(I(w))$, which implies
$T(I(aw))=T\Delta_a\cap I(w)$.
%In particular, (i) implies that $(I(wb))_{b\in R(w)}$ is an ordered
%partition of $I(w)$ with respect to $<_1$.
\begin{proposition}
One has $w\in\cL(X)$ if and only if $I(w)\ne\emptyset$.
\end{proposition}
\begin{proof}
The statement follows from the fact that $w\in\cL(X)$ if and only
if $[w]\ne\emptyset$ and that an easy induction shows that
$[w]\ne\emptyset$ if and only if $I(w)\ne\emptyset$.
\end{proof}

Besides the semi-intervals $I(w)$, we will also need
symmetrically the sets $J(w)$ defined, for $w\in A^*$, by $J(\varepsilon)=I$
 and by
\begin{equation}
J(ua)=TJ(u)\cap T\Delta_a\label{eqJua}
\end{equation}
for $a\in A$ and $u\in A^*$. As for the $I(w)$, the nonempty sets
$J(w)$ are semi-intervals (Exercise~\ref{exerciseJw}).

\subsection{Planar dendric shifts}

A shift space is called \emph{planar dendric} if it is dendric and if there are
two orders $\le_1$ and $\le_2$ on the alphabet $A$ such
that for every $w\in\cL(X)$, the tree $\E(w)$
is compatible with these orders, that is,  for every $(a,b),(c,d)\in \E(w)$,
one has $a\le_1 c$ if and only if $b\le_2 d$.

Thus, placing the
vertices of $L(w)$ on a vertical line in the order given by  $<_1$ and those of $R(w)$
on a parallel line in the order given by $<_2$,
the tree $\E(w)$ becomes planar. Note that the orders
$<_1$ and $<_2$ do not depend on $w$.
\begin{proposition}\label{propositionPlanarDendric}
A regular interval exchange transformation shift is
a minimal and planar dendric shift.
\end{proposition}
\begin{proof}
Assume that $T$ is a regular interval exchange transformation relative
to $(\Delta_a)_{a\in A}$. Let $X$ be the natural representation of $T$.
 We consider, on the alphabet $A$, the two orders
defined by 
\begin{enumerate}
\item $a<_{\rm top}b$ if $\Delta_a$ is to the left of $\Delta_b$,
\item $a<_{\rm bottom}b$ if $T\Delta_a$ is to the left of $T\Delta_b$.
\end{enumerate}
Note that, if $wb,wc$ are in $\cL(X)$ then $b<_{\rm top}c$ if and only if
$I(wb)$ is to the left of $I(wc)$. Indeed, $I(wb),I(wc)\subset I(w)$
and $T^{|w|}$ is a translation on $I(w)$.

For $a,a'\in L(w)$ with $a<_{\rm bottom} a'$, there is a unique reduced
path in $\E(w)$ from $a$ to
$a'$ which is the sequence
 $a_1,b_1,\ldots a_n$ with $a_1=a$ and $a_n=a'$ with
$a_1<_{\rm bottom}a_2<_{\rm bottom}\cdots<_{\rm bottom}a_n$, $b_1<_{\rm top}b_2<_{\rm top}\cdots <_{\rm top} b_{n-1}$ and
$T\Delta_{a_{i}}\cap I_{wb_i}\ne\emptyset$, $T\Delta_{a_{i+1}}\cap
I_{wb_{i}}\ne\emptyset$
for $1\le i\le n-1$ (see Figure~\ref{figG(w)}). Note that
the hypothesis that $T$ is regular is needed here since otherwise
the right boundary of $T\Delta_{a_i}$ could be the left boundary of
$I(wb_i)$. Thus $\E(w)$ is a tree.
It is compatible with the orders $<_{\rm bottom},<_{\rm top}$ since the above shows that $a<_{\rm bottom} a'$ 
implies that the letters $b_1,b_{n-1}$ such that
$(a,b_1),(a',b_{n-1})\in E(w)$
satisfy $b_1\le_{\rm top} b_{n-1}$.
\end{proof}
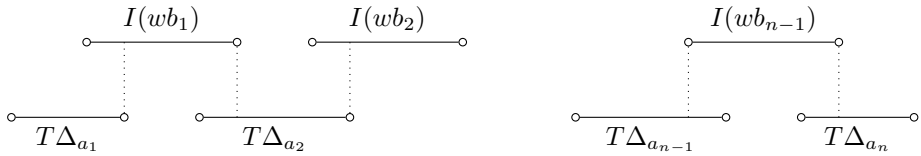
\begin{figure}[hbt]
\centering
\tikzset{node/.style={circle,draw,minimum size=0.1cm,inner sep=0pt}}
\begin{tikzpicture}
\node[node](H1)at(0,1){};
\node[node](H2)at(2,1){};\node[node](H3)at(3,1){};\node[node](H4)at(5,1){};
\node[node](H5)at(8,1){};\node[node](H6)at(10,1){};

\node[node](B1)at(-1,0){};\node[node](B2)at(.5,0){};\node[node](B3)at(1.5,0){};
\node[node](B4)at(3.5,0){};
\node[node](B5)at(6.5,0){};\node[node](B6)at(8.5,0){};
\node[node](B7)at(9.5,0){};
\node[node](B8)at(11,0){};

\draw[above](H1)edge node{$I(wb_1)$}(H2);
\draw[above](H3)edge node{$I(wb_2)$}(H4);
\draw[above](H5)edge node{$I(wb_{n-1})$}(H6);

\draw[below](B1)edge node{$T\Delta_{a_1}$}(B2);
\draw[below](B3)edge node{$T\Delta_{a_2}$}(B4);
\draw[below](B5)edge node{$T\Delta_{a_{n-1}}$}(B6);
\draw[below](B7)edge node{$T\Delta_{a_n}$}(B8);

\draw[dotted](B2)edge node{}(.5,1){};\draw[dotted](2,0)edge node{}(H2){};
\draw[dotted]{}(B4)edge node{}(3.5,1);
\draw[dotted](8,0)edge node{}(H5);\draw[dotted](10,0)edge node{}(H6);
\end{tikzpicture}
\caption{A path from $a_1$ to $a_n$ in $\E(w)$.}\label{figG(w)}
\end{figure}
The converse of Proposition~\ref{propositionPlanarDendric} is also true
(see the Notes Section).
\begin{example}
  Let $X$ be the interval exchange of Example~\ref{exampleRotation2alpha}.
The extension graph of $\varepsilon$ is represented in Figure~\ref{figureExtensionGraphRotation2alpha}.
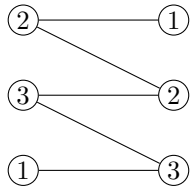
\begin{figure}[hbt]
\centering
\tikzset{node/.style={circle,draw,minimum size=0.4cm,inner sep=0.2pt}}
\begin{tikzpicture}
\node[node](2l)at(0,3){$2$};
\node[node](3l)at(0,2){$3$};
\node[node](1l)at(0,1){$1$};
\node[node](1r)at(2,3){$1$};
\node[node](2r)at(2,2){$2$};
\node[node](3r)at(2,1){$3$};

\draw(2l)edge node{}(1r);\draw(2l)edge node{}(2r);
\draw(3l)edge node{}(2r);\draw(3l)edge node{}(3r);
\draw(1l)edge node{}(3r);
\end{tikzpicture}
\caption{The extension graph $\E(\varepsilon)$.}\label{figureExtensionGraphRotation2alpha}
\end{figure}
The orders $<_{top}$ and $<_{bottom}$ are
\begin{displaymath}
1<_{top}2<_{top}3\quad\mbox{and}\quad 2<_{bottom}<3<_{bottom}1
\end{displaymath}
\end{example}

The following example shows that the Tribonacci shift is not a planar dendric shift
and thus, that the Tribonacci shift is not a regular interval exchange shift
(this can also be proved directly, see Exercise~\ref{exerciseTribonacciNotIET}).

\begin{example}\label{exampleTribonacciPasPlanaire}
Let $X$ be the Tribonacci shift (see Example~\ref{exampleTribonacci}).
\index{subject}{Tribonacci!shift}\index{subject}{shift space!Tribonacci}%
The words $a,aba$ and $abacaba$ are bispecial. Thus
 the words $ba,caba$ are right-special and the
words $ab, abac$ are left-special.
The graphs $\E(\varepsilon),\E(a)$ and $\E(aba)$ are shown in
Figure~\ref{figureTribonacci}.
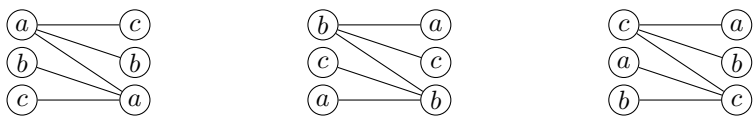
\begin{figure}[hbt]
\centering
\tikzset{node/.style={circle,draw,minimum size=0.4cm,inner sep=0.2pt}}
\begin{tikzpicture}

\node[node](a)at(0,1){$a$};\node[node](b)at(0,.5){$b$};\node[node](c)at(0,0){$c$};
\node[node](c')at(1.5,1){$c$};\node[node](b')at(1.5,.5){$b$};\node[node](a')at(1.5,0){$a$};

\draw(a)edge node{}(a');\draw(a)edge node{}(b');\draw(a)edge node{}(c');
\draw(b)edge node{}(a');\draw(c)edge node{}(a');

\node[node](b)at(4,1){$b$};\node[node](c)at(4,.5){$c$};\node[node](a)at(4,0){$a$};
\node[node](a')at(5.5,1){$a$};\node[node](c')at(5.5,.5){$c$};\node[node](b')at(5.5,0){$b$};

\draw(b)edge node{}(a');\draw(b)edge node{}(b');\draw(b)edge node{}(c');
\draw(c)edge node{}(b');\draw(a)edge node{}(b');

\node[node](c)at(8,1){$c$};\node[node](a)at(8,.5){$a$};\node[node](b)at(8,0){$b$};
\node[node](a')at(9.5,1){$a$};\node[node](b')at(9.5,.5){$b$};\node[node](c')at(9.5,0){$c$};

\draw(c)edge node{}(a');\draw(c)edge node{}(b');\draw(c)edge node{}(c');
\draw(a)edge node{}(c');\draw(b)edge node{}(c');

\end{tikzpicture}
\caption{The graphs $\E(\varepsilon),\E(a)$ and $\E(aba)$ in the Tribonacci
  set.}
\label{figureTribonacci}
\end{figure}
One sees easily that it not possible to find two orders on $A$ making
the representation of the three graphs simultaneously planar. 
\end{example}

\subsection{Rauzy induction}

Since the natural representation of a minimal interval exchange $T$
is dendric, it has by Theorem~\ref{theoremSadicDendric}
a primitive $\Sa_e$-adic representation with a directive
sequence $\tau=(\tau_n)$ formed of elementary morphisms $\tau_n\in\Sa_e$.
We will now give a method, called \emph{Rauzy induction}
\index{subject}{Rauzy!induction} which allows to
build directly this representation from the interval exchange $T$.

\begin{proposition}\label{propositionRauzyInduction}
Let $T$ be the $k$-interval exchange transformation on $[0,1)$
 on the intervals $\Delta_i$ of lengths $\lambda_i$
with the permutation $\pi$. Let $r=\min\{\lambda_k,\lambda_{\pi(k)}\}$.
If $T$ is minimal, the transformation induced by $T$ on the interval $[0,r)$
is again a  minimal $k$-interval exchange transformation $T'$.
\end{proposition}
\begin{proof} 
Set $\ell=\pi(k)$.
Assume first that $\lambda_{\ell}<\lambda_k$ (see Figure~\ref{figureRauzyInduction}).
\begin{figure}[hbt]
\centering
\tikzset{node/.style={circle,draw,minimum size=0.1cm,inner sep=0pt}}
\tikzset{title/.style={circle,minimum size=0.1cm,inner sep=0pt}}
\begin{tikzpicture}
%bas
\node[node](0H)at(0,0){};\node[title]at(0,.4){$0$};
\node[node](1H)at(1.5,0){};
\node[node,fill=black](2H)at(2.5,0){};\node[title]at(2.5,.4){$T\Delta_{k}$};
\node[node](3H)at(3.5,0){};
\node[node,](4H)at(5,0){};\node[title]at(5,-.4){$r$};
\node[node](5H)at(6,0){};\node[title]at(6,.4){$1$};
%haut
\node[node](0B)at(0,1){};\node[title]at(0,1.4){$0$};
\node[node](01B)at(2,1){};\node[node](02B)at(3,1){};
\node[node](1B)at(4,1){};
\node[node,fill=black](2B)at(5,1){};\node[title]at(5,1.4){$\Delta_{k}$};
\node[node](3B)at(6,1){};\node[title]at(6,1.4){$1$};

\draw[dotted](0H)edge node{}(1H);
\draw[below,line width=1pt](2H)edge node{$T\Delta''_{k}$}(3H);
\draw[below,line width=1pt,color=red](1H)edge node{$T\Delta'_{k}$}(2H);
\draw[dotted](3H)edge node{}(4H);
\draw[above,line width=1pt](4H)edge node{$T\Delta_\ell$}(5H);
\draw[dotted](0B)edge node{}(01B);
\draw[above,line width=1pt](01B)edge node{$\Delta_\ell$}(02B);
\draw[dotted](02B)edge node{}(1B);
\draw[below,line width=1pt](2B)edge node{$\Delta''_{k}$}(3B);
\draw[below,color=red,line width=1pt](1B)edge node{$\Delta'_{k}$}(2B);
\draw[dotted](2B)edge node{}(4H);
%\draw[->](5,1)edge node{}(1.5,0){};%\draw[->](4,1)edge node{}(1,0);
%\draw[->](5.75,1)edge node{}(2.75,0);
\end{tikzpicture}
\caption{The  first case of Rauzy induction.}\label{figureRauzyInduction}
\end{figure}
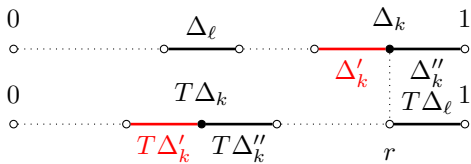
 Let $\Delta''_{k}=\Delta_{k}\cap T\Delta_\ell$ and
let
$\Delta_k=\Delta'_k\cup \Delta''_{k}$ be a partition of $\Delta_k$.
 The transformation $T'$ induced
by $T$ on $[0,r)$ is such that for every $x\in[0,r)$,
\begin{displaymath}
T'x=\begin{cases}Tx&\mbox{ if $Tx\notin \Delta_{\ell}$}\\T^2x&\mbox{ otherwise.}
\end{cases}
\end{displaymath}
Thus $T'$ is an interval exchange on the $k$ intervals 
\begin{displaymath}
\Delta_1,\ldots,\Delta_{k-1},\Delta'_k.
\end{displaymath}
In the case $\lambda_k<\lambda_{\pi(k)}$, we replace $T$ by $T^{-1}$
and we find the first case again. The case $\lambda_{\pi(k)}=\lambda_k$
cannot occur because $T$ is minimal.
\end{proof}

The natural representation of $T'$ is the shift space
$X'=\sigma^{-1}(X)$ where $X$ is the natural coding of $X$
and $\sigma$ is the elementary morphism
\begin{displaymath}
\alpha_{\ell,k}(i)=\begin{cases}\ell k&\mbox{ if $i=\ell$}\\i&\mbox{ otherwise}\end{cases}
\end{displaymath}
which places a $k$ after each $\ell$.

In the case $\lambda_k<\lambda_{\pi(k)}$, we obtain the morphism
$\tilde{\alpha}_{k,\ell}$ which places an $\ell$
before each $k$.
\begin{example}\label{exampleRotation2alphaRauzy1}
Consider the $3$-interval exchange of
Example~\ref{exampleRotation2alpha} which is a rotation of angle $2\alpha$.
\begin{figure}[hbt]
\centering
\tikzset{node/.style={circle,draw,minimum size=0.1cm,inner sep=0pt}}
\tikzset{title/.style={circle,minimum size=0.1cm,inner sep=0pt}}

\begin{tikzpicture}
%haut
\node[node,color=red](0h)at(0,1){};\node[title]at(0,1.4){$0$};
\node[node,color=blue](1-2alpha)at(2.36,1){};\node[title]at(2.36,1.4){$1-2\alpha$};
\node[node,color=green](1-alpha)at(6.18,1){};\node[title]at(6.18,1.4){$1-\alpha$};
\node[node](1h)at(7.64,1){};\node[title]at(7.64,1.4){$2\alpha$};
\draw[color=red,line width=1.5,above](0h)edge node{$a$}(1-2alpha);
\draw[color=blue,line width=1.5,above](1-2alpha)edge node{$b$}(1-alpha);
\draw[color=green,line width=1.5,above](1-alpha)edge node{$c$}(1h);
%bas
\node[node,color=blue](0b)at(0,0){};\node[title]at(0,0.4){$0$};
\node[node,color=green](alpha)at(3.82,0){};\node[title]at(3.82,0.4){$\alpha$};
\node[node](4alpha-1)at(5.24,0){};\node[title]at(5.24,0.4){$4\alpha-1$};
\node[node,color=red](2alpha)at(7.64,0){};\node[title]at(7.64,0.4){$2\alpha$};
\draw[color=blue,line width=1.5,below](0b)edge node{$b$}(alpha);
\draw[color=green,line width=1.5,below](alpha)edge node{$c$}(4alpha-1);
\draw[color=red,line width=1.5,below](4alpha-1)edge node{$a$}(2alpha);

\end{tikzpicture}
\caption{The effect of Rauzy induction on the rotation of angle $2\alpha$.}\label{figure3intervalRauzy}
\end{figure}
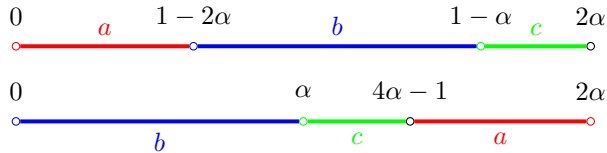
The effect of Rauzy induction is represented in 
Figure~\ref{figure3intervalRauzy}. The corresponding morphism
is $\alpha_{a,c}:a\to ac,b\to b,c\to c$.
\end{example}

Iterating Rauzy induction, one may obtain a transformation
defined on $\Delta_1\cup\ldots\cup\Delta_{k-1}$. This allows
to obtain a BV-representation described as follows.
\begin{theorem}\label{theoremBVrepresentationIET}
Let $(X,S)$ be the natural representation of a minimal $k$-interval exchange.
Then, $(X,S )$ has a BV-representation $(X_E , \varphi_E)$ where $(V,E,\le )$ is such that
\begin{enumerate}
\item
$\Card(V (1)) = k$ and $\Card(V(i)) - \Card(V(i+1)) \in \{ 0,1 \}$ for all $i\geq 1 $,
\item
for all $i\geq 1$ when $V(i-1)=V(i)$ the incidence matrix ${M}(i)$ of $E (i)$ has the following form
$$
{M}(i) =
\begin{bmatrix}
1       & 0      & \dots & 0      & 0      & \dots & 0 & s_1\\
0       & 1      & \dots & 0      & 0      & \dots & 0 & s_2\\
\vdots & \vdots &       & \vdots & \vdots &        & \vdots & \vdots \\
0       & 0      & \dots & 1      & 0      & \dots & 0 & s_l\\
0       & 0      & \dots & 1      & 0      & \dots & 0 & s_{l+1}\\
0       & 0      & \dots & 0      & 1      & \dots & 0 & s_{l+2}\\
\vdots & \vdots &       & \vdots & \vdots &        & \vdots & \vdots \\
0       & 0      & \dots & 0      & 0      & \dots & 1 & s_{k}
\end{bmatrix} \ ,
$$
where $s_i \in \{ 0,m,m+1 \}$, $s_l = m$ and $s_{l+1} = m+1$ for some $m\geq 0$.
When $\Card(V(i)) - \Card(V(i+1)) = 1$, the line $l+1$ does not exist. 
All the entries of ${M}(1)$ are equal to 1.
\end{enumerate}
\end{theorem}

\begin{example}
Let us consider again the interval exchange which is the
rotation of angle $2\alpha$ of Example~\ref{exampleRotation2alpha}.
A first step of Rauzy induction gives the interval exchange of
Example~\ref{exampleRotation2alphaRauzy1}. A second step
gives the interval exchange represented in Figure~\ref{figure3intervalRauzy2}.
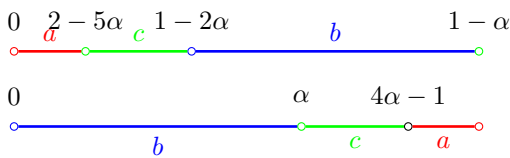
\begin{figure}[hbt]
\centering
\tikzset{node/.style={circle,draw,minimum size=0.1cm,inner sep=0pt}}
\tikzset{title/.style={circle,minimum size=0.1cm,inner sep=0pt}}

\begin{tikzpicture}
%haut
\node[node,color=red](0h)at(0,1){};\node[title]at(0,1.4){$0$};
\node[node,color=green](2-5alpha)at(.95,1){};\node[title]at(.95,1.4){$2-5\alpha$};
\node[node,color=blue](1-2alpha)at(2.36,1){};\node[title]at(2.36,1.4){$1-2\alpha$};
\node[node,color=green](1-alpha)at(6.18,1){};\node[title]at(6.18,1.4){$1-\alpha$};

\draw[color=red,line width=1,above](0h)edge node{$a$}(2-5alpha);
\draw[color=green,line width=1,above](2-5alpha)edge node{$c$}(1-2alpha);
\draw[color=blue,line width=1,above](1-2alpha)edge node{$b$}(1-alpha);

%bas
\node[node,color=blue](0b)at(0,0){};\node[title]at(0,0.4){$0$};
\node[node,color=green](alpha)at(3.82,0){};\node[title]at(3.82,0.4){$\alpha$};
\node[node](4alpha-1)at(5.24,0){};\node[title]at(5.24,0.4){$4\alpha-1$};
\node[node,color=red](1-alphab)at(6.18,0){};\node[title]at(6.18,0.4){};
\draw[color=blue,line width=1,below](0b)edge node{$b$}(alpha);
\draw[color=green,line width=1,below](alpha)edge node{$c$}(4alpha-1);
\draw[color=red,line width=1,below](4alpha-1)edge node{$a$}(1-alphab);

\end{tikzpicture}
\caption{A second step of Rauzy induction on the rotation of angle $2\alpha$.}\label{figure3intervalRauzy2}
\end{figure}
The combination of the two steps corresponds to the morphism
$\varphi:a\to ac,b\to b, c\to cac$. The composition matrix of this 
morphism is of the form given in Theorem~\ref{theoremBVrepresentationIET}
with the order $b<a<c$. Indeed, with this order
\begin{displaymath}
M(\varphi)=\begin{bmatrix}1&0&0\\0&1&1\\0&1&2\end{bmatrix}
\end{displaymath}
which is of form given in case 2 with $s_1=0$, $s_2=1$, $s_3=2$
and $\ell=1$, $m=1$.
\end{example}

\subsection{Induced transformations and admissible semi-intervals}
Let $T$ be a minimal interval exchange transformation on $I=[\ell,r)$.
Let $J\subset [\ell, r)$ be a  semi-interval.
Since $T$ is minimal, for each $z \in [\ell, r)$ there is an integer $n>0$ such that $T^n(z) \in J$.

As seen before, the transformation induced
 by $T$ on $J$ is the transformation $S : J\rightarrow J$ defined for $z \in J$ by $S(z) = T^n(z)$ with $n = \min\{ i>0 \mid T^i(z) \in J\}$. We also say that $S$ is the \emph{first return map}
 (of $T$) on $J$.
The semi-interval $J$ is called the \emph{domain} of $S$, denoted $D(S)$.

\begin{example}
Let $T$ be the transformation of Example~\ref{exampleRotation2alpha}.
Let $J = [0, 2\alpha)$.
The transformation induced by $T$ on $J$ is
$$
S(z) =
\begin{cases}
T^2(z)	&	\text{ if $0 \leq z< 1-2\alpha$} \\
T(z)	&	\text{ otherwise}.
\end{cases}
$$
\end{example}

Let $T$ be an interval exchange transformation relative to  $(\Delta_a)_{a \in A}$. Denote by $d_a$ for $a\in A$ the separation points.
For  $\ell < t < r$, the semi-interval $[\ell, t)$  is \emph{right admissible} for $T$ if there is a $n \in \Z$ such that $t = T^n(d_a)$ for some $a\in A$ and
\begin{enumerate}
\item[(i)] if $n > 0$, then $t < T^h(d_a)$ for all $h$ such that $0 < h < n$,
\item[(ii)] if $n \leq 0$, then $t < T^h(d_a)$ for all $h$ such that $n < h \leq 0$.
\end{enumerate}
We also say that $t$ itself is right admissible.
Note that all semi-intervals $[\ell, d_a)$ with $\ell < d_a$ are right admissible.
Similarly, all semi-intervals $[\ell, Td_a)$ with $\ell < Td_a$ are right admissible.

\begin{example}
\label{ex:division}
Let $T$ be the interval exchange transformation of Example~\ref{exampleRotation2alpha}.
The semi-interval $[0, t)$ for $t = 1-2\alpha$ or $t = 1-\alpha$ is right admissible since $1-2\alpha = d_b$ and $1-\alpha = d_c$.
On the contrary, for $t = 2-3\alpha$, it is not right admissible because $t = T^{-1}(d_c)$ but $d_c < t$ contradicting (ii).
\end{example}

The following result is the basis for the definition of Rauzy induction.

\begin{theorem}[Rauzy]
\label{theo:rauzy1}
Let $T$ be a regular $k$-interval exchange transformation and let $J$ be a right admissible interval for $T$.
The transformation induced by $T$ on $J$ is a regular $k$-interval exchange
transformation.
\end{theorem}

Note that  a $k$-interval exchange transformation $T$ on $[\ell, r)$,
the transformation induced by $T$ 
on any semi-interval included in $[\ell, r)$ is always an interval exchange transformation on at most $k+2$ intervals by Theorem~\ref{lemmaRauzyInduction}.

\begin{example}
\label{ex:induced1}
Consider again the transformation of Example~\ref{exampleRotation2alpha}.
The transformation induced by $T$ on the semi-interval $I = [0,2\alpha)$ 
is the $3$-interval exchange transformation represented in Figure~\ref{fig:3ietind}.
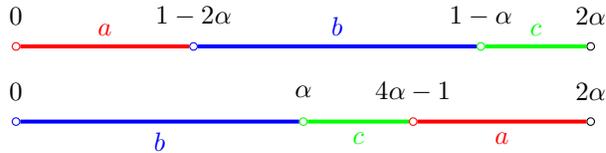
\begin{figure}[hbt]
\centering
\tikzset{node/.style={circle,draw,minimum size=0.1cm,inner sep=0.2pt}}
\tikzset{title/.style={circle,minimum size=0.1cm,inner sep=0pt}}
\begin{tikzpicture}
\node[node,color=red](0h)at(0,1){};\node[title]at(0,1.4){$0$};
\node[node,color=blue](1-2alpha)at(2.36,1){};\node[title]at(2.36,1.4){$1-2\alpha$};
\node[node,color=green](1-alpha)at(6.18,1){};\node[title]at(6.18,1.4){$1-\alpha$};
\node[node](2alphah)at(7.64,1){};\node[title]at(7.64,1.4){$2\alpha$};
\draw[color=red,line width=1.5,above](0h)edge node{$a$}(1-2alpha);
\draw[color=blue,line width=1.5,above](1-2alpha)edge node{$b$}(1-alpha);
\draw[color=green,line width=1.5,above](1-alpha)edge node{$c$}(2alphah);

\node[node,color=blue](0b)at(0,0){};\node[title]at(0,.4){$0$};
\node[node,color=green](alpha)at(3.82,0){};\node[title]at(3.82,.4){$\alpha$};
\node[node,color=red](4alpha-1)at(5.28,0){};\node[title]at(5.28,.4){$4\alpha-1$};
\node[node](2alpha)at(7.64,0){};\node[title]at(7.64,.4){$2\alpha$};
\draw[color=blue,line width=1.5,below](0b)edge node{$b$}(alpha);
\draw[color=green,line width=1.5,below](alpha)edge node{$c$}(4alpha-1);
\draw[color=red,line width=1.5,below](4alpha-1)edge node{$a$}(2alpha);
\end{tikzpicture}
\caption{The transformation induced on $I$.}
\label{fig:3ietind}
\end{figure}
\end{example}

The notion of left admissible interval is symmetrical to that of right admissible.
For $\ell < t < r$, the semi-interval $[t, r)$ is \emph{left admissible} for $T$ if there is a $k \in \Z$ such that $t = T^k(d_a)$ for some $a \in A$ and
\begin{enumerate}
\item[(i)] if $k > 0$, then $T^h(d_a) < t$ for all $h$ such that $0 < h < k$,
\item[(ii)] if $k \leq 0$, then $T^h(d_a) < t$ for all $h$ such that $k <h \leq 0$.
\end{enumerate}
We also say that $t$ itself is left admissible.
Note that, as for right induction, the semi-intervals $[d_a, r)$ and $[Td_a, r)$ are left admissible.
The symmetrical statements of Theorem~\ref{theo:rauzy1} also hold for left admissible intervals.

Let us now generalize the notion of admissibility to a two-sided version.
For a semi-interval $J = [u,v) \ \subset [\ell, r)$, we define the following functions on $[\ell, r)$:
$$
\rho^+_{J,T}(z) = \min\{ n > 0 \mid T^n(z) \in\ (u,v) \}, \quad \rho^-_{I,T}(z) = \min\{ n \geq 0 \mid T^{-n}(z) \in\ (u,v)\}.
$$

We then define for every $z\in[\ell,r)$ three sets.
First, let $E_{J,T}(z)$ be the following set of indices.
$$
E_{J,T}(z) = \{k \mid -\rho_{J,T}^-(z) \leq k < \rho_{J,T}^+(z)\}.
$$
Next, the set of \emph{neighbors} of $z$ with respect to $J$ and $T$ is
$$
N_{J,T}(z) = \{T^k(z) \mid k \in E_{J,T}(z)\}.
$$
Finally, the set of \emph{division points} of $I$ with respect to $T$ is the finite set
$$
\Div(J,T) = \bigcup_{i=1}^s N_{J,T}(d_i).
$$

We now formulate the following definition.
For $\ell \leq u < v \leq r$, we say that the semi-interval $J = [u,v)$ is \emph{admissible} for $T$ if $u,v \in \Div(J,T) \cup \{r\}$.

Note that a semi-interval $[\ell, v)$ is right admissible if and only if  it is admissible and that a semi-interval $[u, r)$ is left admissible if and only if it is
admissible.
Note also that $[\ell, r)$ is admissible.

Note also that for a regular interval exchange transformation relative to a partition $(\Delta_a)_{a \in A}$, each of the semi-intervals $\Delta_a$ (or $T\Delta_a$) is admissible although only the first one is right admissible (and the last one is left admissible).
Actually, we will prove that for every word $w$, the semi-intervals $I(w)$ and $J(w)$ are admissible.
In order to do that, we need the following lemma.

\begin{lemma}
\label{lem:boundary}
Let $T$ be a $k$-interval exchange transformation on the semi-interval $[\ell,r)$.
For any $n \geq 1$, the set  $P_{n} = \{T^h(d_i) \mid 1 \leq i \leq k,\ 1 \leq h \leq n\}$ is the set of $(k-1)n+1$ left boundaries of
 the semi-intervals $J(y)$ for
all words $y$ with $|y| = n$.
\end{lemma}
\begin{proof}
Let $Q_n$ be the set of left boundaries of the intervals $J(y)$ for $|y| = n$.
Since $\Card(\cL(T) \cap A^n) = (k-1)n+1$ by Proposition~\ref{pro:tn}, we have $\Card(Q_n) = (k-1)n+1$.
Since $T$ is regular the set $R_n = \{T^h(d_i) \mid 2 \leq i \leq k,\ 1 \leq h \leq n\}$ is made of $(k-1)n$ distinct points.
Moreover, since 
$$
d_1 = T(d_{\pi(1)}),\ T(d_1) = T^2(d_{\pi(1)}), \ldots, T^{n-1}(d_1) = T^n(d_{\pi(1)}),
$$
we have $P_n = R_n\cup \{T^n(d_1)\}$.
This implies $\Card(P_n) \leq (k-1)n+1$.
On the other hand, if $y = b_0 \cdots b_{n-1}$, then $J(y) = \cap_{i=0}^{n-1} T^{n-i}I(b_{i})$.
Thus the left boundary of each $J(y)$ is the left boundary of some $T^hI(a)$ for some $h$ with $1 \leq h \leq n$ and some $a \in A$.
Consequently $Q_n \subset P_n$.
This proves that $\Card(P_n) = (k-1)n+1$ and that consequently $P_n = Q_n$.
\end{proof}

A dual statement holds for the semi-intervals $I(y)$.

\begin{proposition}
\label{pro:jadm}
Let $T$ be a $k$-interval exchange transformation on the semi-interval $[\ell,r)$.
For any $w\in F(T)$, the semi-interval $J(w)$ is admissible.
\end{proposition}
\begin{proof}
Set $|w| = n$ and $J(w) = [u,v)$.
By Lemma~\ref{lem:boundary}, we have $u = T^g(d_i)$ for $1 \leq i \leq k$ and $1 \leq g \leq n$.
Similarly, we have $v = r$ or $v = T^d(d_j)$ for $1 \leq j \leq k$ and $1 \leq d \leq n$.

For  $1 < h < g$, the point $T^h(d_i)$ is the left boundary of some semi-interval $J_y$ with $|y| = n$ and thus $T^h(d_i) \notin J(w)$.
This shows that $g \in E_{J(w),T}(d_i)$ and thus that $u \in \Div(J(w),T)$.

If $v = r$, then $v \in \Div(J(w),T)$.
Otherwise, one shows in the same way as above that $v \in \Div(J(w),T)$.
Thus $J(w)$ is admissible.
\end{proof}
Note that the same statement holds for the semi-intervals $I(w)$ instead of the semi-intervals $J(w)$ (using the dual statement of Lemma~\ref{lem:boundary}).

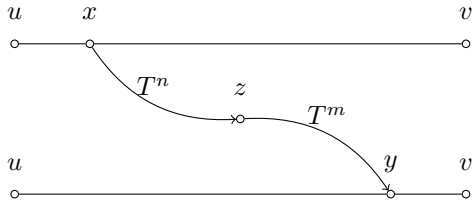
\begin{figure}[hbt]
\centering
\tikzset{node/.style={circle,draw,minimum size=0.1cm,inner sep=0pt}}
\tikzset{title/.style={circle,minimum size=0.1cm,inner sep=0pt}}
\begin{tikzpicture}
\node[node](u1)at(0,2){};\node[title]at(0,2.4){$u$};
\node[node](x)at(1,2){};\node[title]at(1,2.4){$x$};
\node[node](v1)at(6,2){};\node[title]at(6,2.4){$v$};
\node[node](z)at(3,1){};\node[title]at(3,1.4){$z$};
\node[node](u3)at(0,0){};\node[title]at(0,0.4){$u$};
\node[node](y)at(5,0){};\node[title]at(5,.4){$y$};
\node[node](v3)at(6,0){};\node[title]at(6,.4){$v$};
\node[title](d1)at(8,2){};
\node[title](d2)at(8,1){};
\node[title](d3)at(8,0){};

\draw(u1)edge node{}(x);
\draw(x)edge node{}(v1);
\draw(u3)edge node{}(y);
\draw(y)edge node{}(v3);
\draw[bend right,above,->](x)edge node{$T^n$}(z);
\draw[bend left,above,->](z)edge node{$T^m$}(y);
\end{tikzpicture}
\caption{The neighbors of $z$ with respect to $I = [u,v[$.}
\label{fig:adm}
\end{figure}

It can be useful to reformulate the definition of a division point and of an admissible pair using the terminology of graphs.
Consider the graph with vertex set $[\ell, r)$ and edges the pairs $(z, T(z))$ for $z \in [\ell, r[$.
Then, if $T$ is minimal and $I$ is a semi-interval, for any $z \in [\ell, r)$, there is a path $P_{I,T}(z)$ such that its origin $x$ and its end $y$ are in $I$, $z$ is on the path, $z \neq y$ and no vertex of the path except $x, y$ are in $I$ (actually $x = T^{-n}(z)$ with $n = \rho^-_{I,T}(z)$ and $y = T^m(z)$ with $m = \rho^+_{I,T}(z)$).
Then the division points of $I$ are the vertices which are on a path $P_{I,T}(d_i)$ but not at its end (see Figure~\ref{fig:adm}).

The following is a generalization of Theorem~\ref{theo:rauzy1}.
Recall that $\Sep(T)$ denotes the set of separation points of $T$, i.e. the points $d_1 = 0, d_2, \ldots, d_k$ (which are the left boundaries of the semi-intervals $\Delta_1, \ldots, \Delta_k$).

\begin{theorem}
\label{theo:birauzy1}
Let $T$ be a regular $k$-interval exchange transformation on $[\ell, r)$.
For any admissible semi-interval $I = [u,v)$, the transformation $S$ induced by $T$ on $I$ is a regular $k$-interval exchange transformation with separation points $\Sep(S) = \Div(I,T)\cap I$.
\end{theorem}
\begin{proof}
Since $T$ is regular, it is minimal.
Thus for each $i \in \{2, \ldots, k\}$ there are points $x_i, y_i \in (u,v)$ such that there is a path from $x_i$ to $y_i$ passing by $d_i$ but not containing any point of $I$ except at its origin and its end.
Since $T$ is regular, the $x_i$ are all distinct and the $y_i$ are all distinct.

Since $I$ is admissible, there exist $\lambda,\rho \in \{1, \ldots, k\}$ such that $u \in N_{I,T}(d_\lambda)$ and $v \in N_{I,T}(d_\rho)$.
Moreover, since $u$ is a neighbor of $d_\lambda$ with respect to $I$, $u$ is on the path from $x_\lambda$ to $y_\lambda$ (it can be either before or after $d_\lambda$).
Similarly, $v$ is on the path from $x_\rho$ to $y_\rho$ (see Figure~\ref{fig:biinduced} where $u$ is before $d_\lambda$ and $v$ is after $d_\rho$).

\begin{figure}[hbt]
\centering
\tikzset{node/.style={circle,draw,minimum size=0.1cm,inner sep=0pt}}
\tikzset{title/.style={circle,minimum size=0.1cm,inner sep=0pt}}
\begin{tikzpicture}

\node[node,color=blue](uh)at(0,2){};\node[title]at(0,2.4){$u$};
\node[node](x_{i_2})at(.7,2){};
\node[node,color=green](x_g)at(2.5,2){};\node[title]at(2.5,2.4){$x_\lambda$};
\node[node](x_{i_{k+1}})at(3.2,2){};
\node[node,color=red](x_{i_j})at(4,2){};\node[title]at(4,2.4){$x_j$};
\node[node](x_{i_{j+1}})at(5.5,2){};
\node[node,color=yellow](x_d)at(7,2){};\node[title]at(7,2.4){$x_\rho$};
\node[node](x_e)at(7.8,2){};
\node[node](vh)at(10,2){};\node[title]at(10,2.4){$v$};

\node[node](gamma_g)at(3,1){};\node[title]at(3,1.4){$d_\lambda$};
\node[node](gamma_{i_j})at(5,1){};\node[title]at(5,1.4){$d_j$};
\node[node](gamma_d)at(8.2,1){};\node[title]at(8,1.4){$d_\rho$};

\node[node,color=green](ub)at(0,0){};\node[title]at(0,.4){$u$};
\node[node](y_{i_2})at(.7,0){};
\node[node,color=blue](y_g)at(4,0){};\node[title]at(4,.4){$y_\lambda$};
\node[node](y_{i_{k+1}})at(4.7,0){};
\node[node,color=red](y_{i_j})at(6,0){};\node[title]at(6,.4){$y_j$};
\node[node](y_{i_{j+1}})at(7.5,0){};
%\node[node](x_d)at(7,2){};\node[title]at(7,2.4){$x_d$};
\node[node,color=yellow](y_d)at(8.6,0){};\node[title]at(8.6,.4){$y_\rho$};
\node[node](y_e)at(9.4,0){};
\node[node](vb)at(10,0){};\node[title]at(10,.4){$v$};

\draw[color=blue,line width=1](uh)edge node{}(x_{i_2}){};
\draw(x_{i_2})edge node{}(x_g);
\draw[color=green,line width=1](x_g)edge node{}(x_{i_{k+1}});
\draw(x_{i_{k+1}})edge node{}(x_{i_j});
\draw[color=red,line width=1](x_{i_j})edge node{}(x_{i_{j+1}});
\draw[](x_{i_{j+1}})edge node{}(x_d);
\draw[color=yellow,line width=1](x_d)edge node{}(x_e);
\draw(x_e)edge node{}(vh);
\draw[bend right,->](x_g)edge node{}(uh);
\draw[bend right=15,->](uh)edge node{}(gamma_g);
\draw[bend left,->](gamma_g)edge node{}(y_g);
\draw[bend right,->](x_{i_j})edge node{}(gamma_{i_j});
\draw[bend left,->](gamma_{i_j})edge node{}(y_{i_j});
\draw[bend right,->](x_d)edge node{}(gamma_d);
\draw[bend right,->](gamma_d)edge node{}(vh);
\draw[bend left,->](vh)edge node{}(y_d);
\draw[color=green,line width=1](ub)edge node{}(y_{i_2});
\draw(y_{i_2})edge node{}(y_g);
\draw[color=blue,line width=1](y_g)edge node{}(y_{i_{k+1}});
\draw(y_{i_{k+1}})edge node{}(y_{i_j});
\draw[color=red,line width=1](y_{i_j})edge node{}(y_{i_{j+1}});
\draw[](y_{i_{j+1}})edge node{}(y_d);
\draw[color=yellow,line width=1](y_d)edge node{}(y_e);
\draw(y_e)edge node{}(vb);
\end{tikzpicture}
\caption{The transformation induced on $[u, v)$.}
\label{fig:biinduced}
\end{figure}
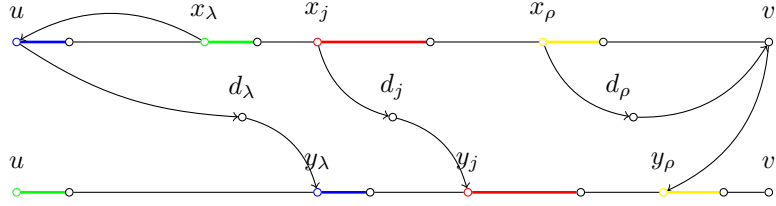

Set $x_1 = y_1 = u$.
Let $(I_j)_{1 \leq j \leq k}$ be the partition of $I$ in semi-intervals such that $x_j$ is the left boundary of $I_j$ for $1 \leq j \leq k$.
Let $J_j$ be the partition of $I$ such that $y_j$ is the left boundary of $J_j$ for $1 \leq j \leq k$.
We will prove that
$$
S(I_j) =
\begin{cases}
J_j	&	\text{if $j \neq 1,\lambda$} \\
J_1	&	\text{if $j = \lambda$} \\
J_\lambda	&	\text{if $j = 1$}
\end{cases}
$$
and that the restriction of $S$ to $I_j$ is a translation.

Assume first that $j \neq 1,\lambda$.
Then $S(x_j) = y_j$.
Let $n$ be such that $y_j = T^n(x_j)$ and denote $I'_j = I_j \setminus x_j$.
We will prove by induction on $h$ that for $0 \leq h \leq k-1$, the set $T^h(I'_j)$ does not contain $u, v$ or any $x_i$.
It is true for $h = 0$.
Assume that it holds up to $h < n-1$.

For any $h'$ with $0 \leq h' \leq h$, the set $T^{h'}(I_j')$ does not contain any $\gamma_i$.
Indeed, otherwise there would exist $h''$ with $0 \leq h'' \leq h'$ such that $x_i \in T^{h''}(I'_j)$, a contradiction.
Thus $T$ is a translation on $T^{h'}(I_j)$.
This implies that $T^h$ is a translation on $I_j$.
Note also that $T^h(I'_j) \cap I = \emptyset$.
Assume the contrary.
We first observe that we cannot have $T^h(x_j) \in I$.
Indeed, $h < n$ implies that $T^h(x_j) \notin (u, v)$.
And we cannot have $T^h(x_j) = u$ since $j \neq \lambda$.
Thus $T^h(I'_j) \cap I \neq \emptyset$ implies that $u \in T^h(I'_j)$, a contradiction.

Suppose that $u = T^{h+1}(z)$ for some $z \in I'_j$.
Since $u$ is on the path from $x_\lambda$ to $y_\lambda$, it implies that for some $h'$ with $0 \leq h' \leq h$ we have $x_\lambda = T^{h'}(z)$, a contradiction with the induction hypothesis.
A similar proof (using the fact that $v$ is on the path from $x_\rho$ to $y_\rho$) shows that $T^{h+1}(I'_j)$ does not contain $v$.
Finally suppose that some $x_i$ is in $T^{h+1}(I'_j)$.
Since the restriction of $T^h$ to $I_j$ is a translation, $T^h(I_j)$ is a semi-interval.
Since $T^{h+1}(x_j)$ is not in $I$ the fact that $T^{h+1}(I_j) \cap I$ is not empty implies that $u \in T^h(I_j)$, a contradiction.

This shows that $T^n$ is continuous at each point of $I'_j$ and that $S = T^n(x)$ for all $x \in I_j$.
This implies that the restriction of $S$ to $I_j$ is a translation into $J_j$.

If $j=1$, then $S(x_1) = S(u) = y_\lambda$.
The same argument as above proves that the restriction of $S$ to $I_1$ is a translation form $I_1$ into $J_\lambda$.
Finally if $j = \lambda$, then $S(x_\lambda) = x_1 = u$ and, similarly, we obtain that the restriction of $S$ to $I_\lambda$ is a translation into $I_1$.

Since $S$ is the transformation induced by the transformation $T$ which is one to one, it is also one to one.
This implies that the restriction of $S$ to each of the semi-intervals $I_j$ is a bijection onto the corresponding interval $J_j,J_1$ or $J_\lambda$ according to the value of $j$.

This shows that $S$ is an $s$-interval exchange transformation.
Since the orbits of the points $x_2, \cdots, x_k$ relative to $S$ are included in the orbits of $d_2, \ldots, d_k$, they are infinite and disjoint.
Thus $S$ is regular.

Let us  finally show that $\Sep(S) = \Div(I,T) \cap I$.
We have $\Sep(S) = \{x_1, x_2, \ldots, x_k\}$ and $x_i \in N_{I,T}(d_i)$.
Thus $\Sep(S) \subset \Div(I,T) \cap I$.
Conversely, let $x \in \Div(I,T) \cap I$.
Then $x \in N_{I,T}(d_i) \cap I$ for some $1 \leq i \leq k$.
If $i\ne 1,\lambda$, then $x = x_i$.
If $i = 1$, then either $x = u$ (if $u = \ell$) or $x = x_{\pi(1)}$ since $d_1 = T(d_{\pi(1)})$.
Finally, if $i = \lambda$ then $x = u$ or $x = x_\lambda$.
Thus $x \in \Sep(S)$ in all cases.
\end{proof}

\subsection{A closure property}
We will prove that the family of regular 
interval exchange shifts is closed under derivation.
The same property holds for minimal dendric shifts
by Theorem~\ref{propositionReturns}. Thus regular interval
exchange shifts form a subfamily closed under derivation
of minimal dendric shifts.
\begin{theorem}
\label{theo:returns}
Any derived  shift of a regular $k$-interval exchange shift $X$
with respect to some $w\in \cL(X)$ is a regular $k$-interval exchange shift.
\end{theorem}
We first prove the following lemma.
\begin{lemma}
\label{lem:returnsind}
Let $T$ be a regular interval exchange transformation and let $(X,S)$
be its natural representation. 
For $w \in \cL(X)$, let $T'$ be the transformation induced by $T$ on $J_w$.
One has $x\in \RR_X(w)$ if and only if
$$
\Sigma_T(z) = x\Sigma_T(T'(z))
$$
for some $z \in J_w$.
\end{lemma}
\begin{proof}
Assume first that $x\in \RR_X(w)$.
Then for any $z \in J_{w} \cap I_x$, we have $T'(z) = T^{|x|}(z)$ and
$$
\Sigma_T(z) = x\Sigma_T(T^{|x|}(z)) = x\Sigma_T(T'(z)).
$$
Conversely, assume that $\Sigma_T(z) = x\Sigma_T(T'(z))$ for some $z \in J_w$.
Then $T^{|x|}(z) \in J_w$ and thus $wx \in A^*w$ which implies that $x \in \Gamma_X(w)$.
Moreover $x$ does not have a proper prefix in $\Gamma_X(w)$ and thus $x\in \RR_X(w)$.
\end{proof}

Since a regular interval exchange shift is recurrent, the previous lemma says that the natural coding of a point in $J_w$ is a concatenation of first return words to $w$.
Moreover, note also that $T^n(z) \in J_w$ if and only if the prefix of length $n$ of $\Sigma_T(z)$ is a return word to $w$.

\begin{proofof}{of Theorem~\ref{theo:returns}}
Let $T$ be a regular $s$-interval exchange transformation and let $X$
be its natural representation.

Let $w \in \cL(X)$.
Since the semi-interval $J_w$ is admissible according to Proposition~\ref{pro:jadm}, the transformation $T'$ induced by $T$ on $J_w$ is, by Theorem~\ref{theo:birauzy1}, a $k$-interval exchange transformation.
The corresponding partition of $J_w$ is the family $\left( J_{wx} \right)_{x \in \RR_X(w)}$.

Using Lemma~\ref{lem:returnsind} and the observation following it, it is clear that $\Sigma_T(z) = \varphi(\Sigma_S(z))$, where $z$ is a point 
of $J_w$ and $\varphi : A^* \to \RR_X(w)^*$ is a coding morphism for $\RR_X(w)$.

Set $x = \Sigma_T(T^{-|w|}(z))$ and $y = \Sigma_T(z)$.
Then $x = wy$ and thus $\Sigma_S(z) = \D_w(x)$.
This shows that the derived shift of $X$ with respect to $w$ is the natural representation of $T'$.
\end{proofof}
\begin{example}
  Consider  the interval exchange of Example~\ref{exampleRotation2alpha}. The
  set of return words to $c$ is $\RR_X(c)=\{bac,bbac,c\}$. The
  natural coding of
  transformation induced on the interval $\Delta_c$ is the interval exchange
  represented in Figure~\ref{figureExampleInduced}.
  The intervals of the upper level are labeled by $\RR'_X(c)$
  and those of the lower level by $\RR_X(c)$
  (we shall describe this transformation in more detail in Examples
  \ref{ex:simphi6T} and~\ref{exampleAutomorphisms}).
  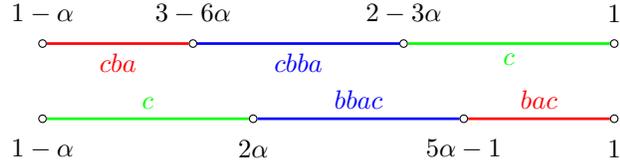
\begin{figure}[hbt]
\centering
\tikzset{node/.style={circle,draw,minimum size=0.1cm,inner sep=0pt}}
\tikzset{title/.style={circle,minimum size=0.1cm,inner sep=0pt}}
\begin{tikzpicture}
  \node[node](alpha)at(0,1){};\node[title]at(0,1.4){$1-\alpha$};
  \node[node](2-4alpha)at(2,1){};\node[title]at(2,1.4){$3-6\alpha$};
  \node[node](1-alpha)at(4.8,1){};\node[title]at(4.8,1.4){$2-3\alpha$};
  \node[node](2alpha)at(7.6,1){};\node[title]at(7.6,1.4){$1$};
  \node[node](alphab)at(0,0){};\node[title]at(0,-.4){$1-\alpha$};
  \node[node](4alpha-1)at(2.8,0){};\node[title]at(2.8,-.4){$2\alpha$};
  \node[node](7alpha-2)at(5.6,0){};\node[title]at(5.6,-.4){$5\alpha-1$};
  \node[node](2alphab)at(7.6,0){};\node[title]at(7.6,-.4){$1$};

  \draw[color=red,line width=1,below](alpha)edge node{$cba$}(2-4alpha);
  \draw[color=blue,line width=1,below](2-4alpha)edge node{$cbba$}(1-alpha);
  \draw[color=green,line width=1,below](1-alpha)edge node{$c$}(2alpha);
  \draw[color=green,line width=1,above](alphab)edge node{$c$}(4alpha-1);
  \draw[color=blue,line width=1,above](4alpha-1)edge node{$bbac$}(7alpha-2);
  \draw[color=red,line width=1,above](7alpha-2)edge node{$bac$}(2alphab);
\end{tikzpicture}
\caption{The transformation induced on $\Delta_c$}\label{figureExampleInduced}
    \end{figure}
  \end{example}

\section{Branching Rauzy induction}
\label{sec:rauzy}

In this section we introduce a branching version of Rauzy
induction and generalize Rauzy's results to the two-sided case (Theorems~\ref{theo:birauzy1} and~\ref{theo:birauzy2}).
In particular we characterize in Theorem~\ref{theo:birauzy2} the admissible semi-intervals for an interval exchange transformation.

Let $T = T_{\lambda,\pi}$ be a regular $k$-interval exchange transformation on $[\ell, r)$.
Set $Z(T) = [\ell, \max\{d_{k}, Td_{\pi(k)}\})$.

Note that $Z(T)$ is the largest semi-interval which is right-admissible for $T$.
We denote by $\psi(T)$ the transformation induced by $T$ on $Z(T)$.

%The following result is Theorem 23 in~\cite{Rauzy1979}.

\begin{theorem}[Rauzy]
\label{theo:rauzy2}
Let $T$ be a regular interval exchange transformation.
A semi-interval $I$ is right admissible for $T$ if and only if there is an integer $n \geq 0$ such that $I = Z(\psi^n(T))$.
In this case, the transformation induced by $T$ on $I$ is $\psi^{n+1}(T)$.
\end{theorem}

The map $T \mapsto \psi(T)$ is called the \emph{right Rauzy induction}.
There are actually, as we have seen,
 two cases according to $d_{k} < Td_{\pi(k)}$ (Case 0) or $d_{k} > Td_{\pi(k)}$ (Case 1).
We cannot have $d_{k} = Td_{\pi(k)}$ since $T$ is regular.

In Case 0, we have $Z(T) = [\ell, Td_{\pi(k)})$ and for any $z\in Z(T)$,
$$
S(z) =
\begin{cases}
T^2(z)	&	\text{if $z \in I(a_{\pi(k)})$} \\ 
T(z)	&	\text{otherwise}.
\end{cases}
$$
The transformation $S$ is the interval exchange transformation relative to $(K_a)_{a \in A}$ with $K(a) = I(a)\cap Z(T)$ for all $a \in A$.
Note that $K(a) = I(a)$ for $a \neq a_k$.
The translation values $\beta_a$ are defined as follows, denoting $\alpha_i, \beta_i$ instead of $\alpha_{a_i}, \beta_{a_i}$,
$$
\beta_i =
\begin{cases}
\alpha_{\pi(k)} + \alpha_k	&	\text{if $i = \pi(k)$} \\
\alpha_i			&	\text{otherwise.}
\end{cases}
$$
In summary, in Case 0, the semi-interval $J(a_\pi(k))$ is suppressed, the semi-interval $J(a_k)$ is split into $SK(a_k)$ and $SK(a_{\pi(k)})$.
The left boundaries of the semi-intervals $K(a)$ are the left boundaries of the semi-intervals $I(a)$.
The transformation is represented  in Figure~\ref{fig:rauzyind0}, in which the left boundary of the semi-interval $S(K(a_{\pi(k)})$ is  $Sd_{\pi(k)}$.

\begin{figure}[hbt]
\centering
\tikzset{node/.style={circle,draw,minimum size=0.1cm,inner sep=0pt}}
\tikzset{title/.style={circle,minimum size=0.1cm,inner sep=0pt}}
\begin{tikzpicture}

\node[node](0h)at(0,4){};\node[title]at(0,4.4){$\ell$};
\node[node,fill=blue](gamma_{pi(k)})at(2,4){};\node[title]at(2,4.4){$d_{\pi(k)}$};
\node[node](mu_{pi(k)})at(3,4){};
\node[node,fill=red](gamma_k)at(8,4){};\node[title]at(8,4.4){$d_k$};
\node[node](1h)at(10,4){};\node[title]at(10,4.4){$r$};
\draw(0h)edge node {}(gamma_{pi(k)});
\draw[color=blue,line width=1.5,below](gamma_{pi(k)})edge node{$a_{\pi(k)}$}(mu_{pi(k)});
\draw(mu_{pi(k)})edge node{}(gamma_k);
\draw[color=red,line width=1.5,below](gamma_k)edge node{$a_s$}(1h);

\node[node](0b)at(0,3){};
\node[node,fill=red](delta_k)at(4,3){};\node[title]at(4,3.4){$Td_k$};
\node[node](nu_k)at(6,3){};
\node[node,fill=blue](delta_{pi(k)})at(9,3){};\node[title]at(9,3.4){$Td_{\pi(k)}$};
\node[node](1b)at(10,3){};
\draw(0b)edge node{}(delta_k);
\draw[color=red,line width=1.5,below](delta_k)edge node{$a_k$}(nu_k);
\draw(nu_k)edge node{}(delta_{pi(k)});
\draw[color=blue,line width=1.5,below](delta_{pi(k)})edge node{$a_{\pi(k)}$}(1b);

\node[title]at(5,2){\huge $\downarrow$};\node[title]at(9,2){$\vdots$};

\node[node](0h)at(0,1){};\node[title]at(0,1.4){$\ell$};
\node[node,fill=blue](gamma_{pi(k)})at(2,1){};\node[title]at(2,1.4){$d_{\pi(k)}$};
\node[node](mu_{pi(k)})at(3,1){};
\node[node,fill=red](gamma_k)at(8,1){};\node[title]at(8,1.4){$d_k$};
\node[node](1h)at(9,1){};
\draw(0h)edge node{}(gamma_{pi(k)});
\draw[color=blue,line width=1.5,below](gamma_{pi(k)})edge node{$a_{\pi(k)}$}(mu_{pi(k)});
\draw(mu_{pi(k)})edge node{}(gamma_k);
\draw[color=red,line width=1,below](gamma_k)edge node{$a_k$}(1h);

\node[node](0b)at(0,0){};
\node[node,fill=red](delta_{k,S})at(4,0){};\node[title]at(4,.4){$Sd_{k}$};
\node[node,fill=blue](delta_{pi(k),S})at(5,0){};\node[title]at(5,.4){$Sd_{\pi(k)}$};
\node[node](nu_{pi(k),S})at(6,0){};
\node[node](last)at(9,0){};
\draw(0b)edge node{}(delta_{k,S});
\draw[color=red,line width=1.5,below](delta_{k,S})edge node{$a_k$}(delta_{pi(k),S});
\draw[color=blue,line width=1.5,below](delta_{pi(k),S})edge node{$a_{\pi(k)}$}(nu_{pi(k),S});
\draw(nu_{pi(k),S})edge node{}(last);

\end{tikzpicture}
\caption{Case 0 in Rauzy induction.}
\label{fig:rauzyind0}
\end{figure}

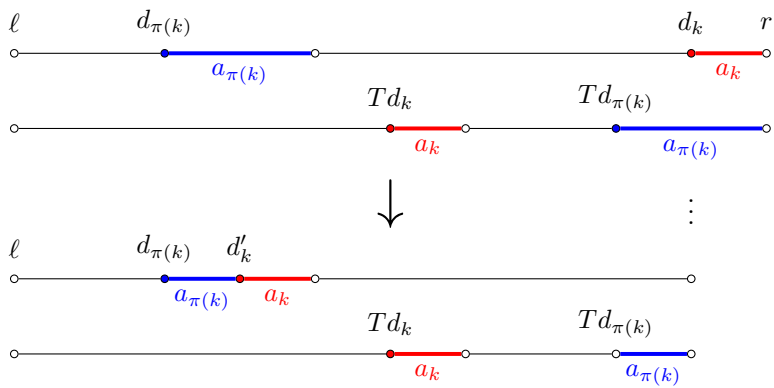
\begin{figure}[hbt]
\centering
\tikzset{node/.style={circle,draw,minimum size=0.1cm,inner sep=0pt}}
\tikzset{title/.style={circle,minimum size=0.1cm,inner sep=0pt}}
\begin{tikzpicture}

\node[node](0h)at(0,4){};\node[title]at(0,4.4){$\ell$};
\node[node,fill=blue](gamma_{pi(s)})at(2,4){};\node[title]at(2,4.4){$d_{\pi(k)}$};
\node[node](mu_{pi(s)})at(4,4){};
\node[node,fill=red](gamma_s)at(9,4){};\node[title]at(9,4.4){$d_k$};
\node[node](1h)at(10,4){};\node[title]at(10,4.4){$r$};
\draw(0h)edge node{}(gamma_{pi(s)});
\draw[color=blue,line width=1.5,below](gamma_{pi(s)})edge node{$a_{\pi(k)}$}(mu_{pi(s)});
\draw(mu_{pi(s)})edge node{}(gamma_s);
\draw[color=red,line width=1.5,below](gamma_s)edge node{$a_k$}(1h);

\node[node](0b)at(0,3){};
\node[node,fill=red](delta_s)at(5,3){};\node[title]at(5,3.4){$Td_k$};
\node[node](nu_s)at(6,3){};
\node[node,fill=blue](delta_{pi(s)})at(8,3){};\node[title]at(8,3.4){$Td_{\pi(k)}$};
\node[node](1b)at(10,3){};
\draw(0b)edge node{}(delta_s);
\draw[color=red,line width=1.5,below](delta_s)edge node{$a_k$}(nu_s);
\draw(nu_s)edge node{}(delta_{pi(s)});
\draw[color=blue,line width=1.5,below](delta_{pi(s)})edge node{$a_{\pi(k)}$}(1b);

\node[title]at(5,2){\huge $\downarrow$};\node[title]at(9,2){$\vdots$};

\node[node](0h)at(0,1){};\node[title]at(0,1.4){$\ell$};
\node[node,fill=blue](gamma_{pi(s)})at(2,1){};\node[title]at(2,1.4){$d_{\pi(k)}$};
\node[node,fill=red](gamma_{s,S})at(3,1){};\node[title]at(3,1.4){$d'_{k}$};
\node[node](mu_{pi(s)})at(4,1){};
\node[node](1hb)at(9,1){};

\draw(0h)edge node{}(gamma_{pi(s)});
\draw[color=blue,line width=1.5,below](gamma_{pi(s)})edge node{$a_{\pi(k)}$}(gamma_{s,S});
\draw[color=red,line width=1.5,below](gamma_{s,S})edge node{$a_k$}(mu_{pi(s)});
\draw(mu_{pi(s)})edge node{}(1hb);

\node[node](0b)at(0,0){};
\node[node,fill=red](delta_sb)at(5,0){};\node[title]at(5,.4){$Td_{k}$};
\node[node](nu_{s})at(6,0){};
\node[node](delta_{pi(s)})at(8,0){};\node[title]at(8,.4){$Td_{\pi(k)}$};
\node[node](1b)at(9,0){};

\draw(0b)edge node{}(delta_sb);
\draw[color=red,line width=1.5,below](delta_sb)edge node{$a_k$}(nu_{s});
\draw(nu_{s})edge node{}(delta_{pi(s)});
\draw[color=blue,line width=1.5,below](delta_{pi(s)})edge node{$a_{\pi(k)}$}(1b);

\end{tikzpicture}
\caption{Case 1 in Rauzy induction.}
\label{fig:rauzyind1}
\end{figure}

In Case 1, we have $Z(T) = [\ell, d_{k})$ and for any $z \in Z(T)$,
$$
S(z) =
\begin{cases}
T^2(z)	&	\text{if $z \in T^{-1}(I_{a_k})$} \\
T(z)	&	\text{otherwise}.
\end{cases}
$$
The transformation $S$ is the interval exchange transformation relative to  $(K(a))_{a\in A}$ with
$$
K(a) =
\begin{cases}
T^{-1}I(a)		&	\text{if $a = a_k$} \\
T^{-1}(TI(a)\cap Z(T))	&	\text{otherwise.}
\end{cases}
$$
Note that $K(a) = I(a)$ for $a \neq a_k$ and $a \neq a_{\pi(k)}$.
Moreover $K(a) = S^{-1}(TI(a) \cap Z(T))$ in all cases.
The translation values $\beta_i$ are defined by
$$
\beta_i =
\begin{cases}
\alpha_{\pi(k)} + \alpha_k	&	\text{if $i = k$} \\
\alpha_i			&	\text{otherwise.}
\end{cases}
$$
In summary, in Case 1, the semi-interval $I(a_s)$ is suppressed, the semi-interval $I(a_{\pi(k)})$ is split into $K(a_{\pi(k)})$ and $K(a_k)$.
The left boundaries of the semi-intervals $SK(a)$ are the left boundaries of the semi-intervals $J(a)$.
The transformation is represented in Figure~\ref{fig:rauzyind1}, where the left boundary of the semi-interval $K(a_k)$ is  $d'_{k}$.

\begin{example}
\label{ex:induced2}
Consider again the transformation $T$ of Example~\ref{exampleRotation2alpha}.
Since $Z(T) = [0, 2\alpha)$, the transformation $\psi(T)$ is represented in Figure~\ref{fig:3ietind}.
The transformation $\psi^2(T)$ is represented in Figure~\ref{fig:3ietind2}.

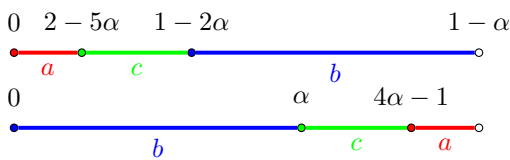
\begin{figure}[hbt]
\centering
\tikzset{node/.style={circle,draw,minimum size=0.1cm,inner sep=0pt}}
\tikzset{title/.style={circle,minimum size=0.1cm,inner sep=0pt}}
\begin{tikzpicture}
\node[node,fill=red](0h)at(0,1){};\node[title]at(0,1.4){$0$};
\node[node,fill=green](2-5alpha)at(.9,1){};\node[title]at(.9,1.4){$2-5\alpha$};
\node[node,fill=blue](1-2alpha)at(2.36,1){};\node[title]at(2.36,1.4){$1-2\alpha$};
\node[node](1-alphah)at(6.18,1){};\node[title]at(6.18,1.4){$1-\alpha$};
\draw[color=red,line width=1.5,below](0h)edge node{$a$}(2-5alpha);
\draw[color=green,line width=1.5,below](2-5alpha)edge node{$c$}(1-2alpha);
\draw[color=blue,line width=1.5,below](1-2alpha)edge node{$b$}(1-alphah);

\node[node,fill=blue](0b)(0,0){};\node[title]at(0,.4){$0$};
\node[node,fill=green](alpha)at(3.82,0){};\node[title]at(3.82,.4){$\alpha$};
\node[node,fill=red](4alpha-1)at(5.28,0){};\node[title]at(5.28,.4){$4\alpha-1$};\node[node](1-alpha)at(6.18,0){};
\draw[color=blue,line width=1.5,below](0b)edge node{$b$}(alpha);
\draw[color=green,line width=1.5,below](alpha)edge node{$c$}(4alpha-1);
\draw[color=red,line width=1.5,below](4alpha-1)edge node{$a$}(1-alpha);
\end{tikzpicture}
\caption{The transformation  $\psi^2(T)$.}
\label{fig:3ietind2}
\end{figure}
\end{example}

The symmetrical notion of \emph{left Rauzy induction} is defined similarly
as follows.

Let $T = T_{\lambda, \pi}$ be a regular $s$-interval exchange transformation on $[\ell, r)$. %Let $f_i$ be the right boundary of the semi-interval $\Delta_i$.
Set $Y(T) = [\min\{d_2, Td_{\pi(2)}\}, r)$.
We denote by $\varphi(T)$ the transformation induced by $T$ on $Y(T)$.
The map $T \mapsto \varphi(T)$ is called the \emph{left Rauzy induction}. 

%Note that one has also $Y(T) = [\min\{d_2, Td_{\pi(2)}\}, r)$.

The symmetrical assertions of Theorem~\ref{theo:rauzy2} also hold for left admissible intervals.

\subsection{Branching induction}
The following is an addition  to Theorem~\ref{theo:rauzy1}.

\begin{theorem}
\label{theo:birauzy2}
Let $T$ be a regular $k$-interval exchange transformation on $[\ell,r)$.
A semi-interval $I$ is admissible for $T$ if and only if there is a sequence $\chi \in \{\varphi, \psi\}^*$ such that $I$ is the domain of $\chi(T)$.
In this case, the transformation induced by $T$ on $I$ is $\chi(T)$.
\end{theorem}

We first prove the following lemmas, in which we assume that $T$ is a regular $s$-interval exchange transformation on $[\ell, r)$.
Recall that $Y(T), Z(T)$ are the domains of $\varphi(T), \psi(T)$ respectively.

\begin{lemma}
\label{lem:init}
If a semi-interval  $I$ strictly included in $[\ell, r)$ is admissible for $T$, then either $I \subset Y(T)$ or $I\subset Z(T)$.
\end{lemma}
\begin{proof}
Set $I = [u,v)$.
Since $I$ is strictly included in $[\ell, r)$, we have either $\ell < u$ or $v < r$.
Set $Y(T) = [y, r)$ and $Z(T) = [\ell, z)$.

Assume that $v < r$.
If $y \leq u$, then $I \subset Y(T)$.
Otherwise, let us show that $v \leq z$.
Assume the contrary.
Since $I$ is admissible, we have $v = T^j(d_i)$ with $j \in E_{I,T}(d_i)$ for some $i$ with $1 \leq i \leq k$.
But $j > 0$ is impossible since $u < T(d_i) < v$ implies $T(d_i) \in\ (u,v)$, in contradiction with the fact that $j < \rho_I^+(d_i)$.
Similarly, $j \leq 0$ is impossible since $u < d_i < v$ implies $d_i \in\ ]u,v)$.
Thus $I\subset Z(T)$.

The proof in the case $\ell < u$ is symmetric.
\end{proof}

%The next lemma is  the two-sided version of Lemma 22 in~\cite{Rauzy1979}.

\begin{lemma}
\label{lem:bi22}
Let $T$ be a regular $k$-interval exchange transformation on $[\ell, r)$.
Let $J$ be an admissible semi-interval for $T$ and let $S$ be the transformation induced by $T$ on $J$.
A semi-interval $I \subset J$ is admissible for $T$ if and only if it is admissible for $S$.
Moreover $\Div(J,T) \subset \Div(I,T)$.
\end{lemma}
\begin{proof}
Set $J = [t,w)$ and $I = [u,v)$.
Since $J$ is admissible for $T$, the transformation $S$ is a regular $k$-interval exchange transformation by Theorem~\ref{theo:birauzy1}.

Suppose first that $I$ is admissible for $T$.
Then $u = T^g(d_i)$ with $g \in E_{I,T}(d_{i})$ for some $1 \leq i \leq k$, and $v = T^{d}(d_{j})$ with $d \in E_{I,T}(d_{j})$ for some $1 \leq j \leq k$ or $v = r$.

Since $S$ is the transformation induced by $T$ on $J$, there is a separation point $x$ of $S$ of the form $x = T^m(d_{i})$ with $m = -\rho^-_{J,T}(d_i)$ 
and thus $m \in E_{J,T}(d_i)$.
Thus $u = T^{g-m}(x)$.

Assume first that $g-m > 0$.
Since $u, x \in J$, there is an integer $n$ with $0 < n \leq g-m$ such that $u = S^n(x)$.

Let us show that $n \in E_{I,S}(x)$.
Assume by contradiction that $\rho_ {I,S}^+(x) \leq n$.
Then there is some $j$ with $0 < j \leq n$ such that $S^j(x) \in (u, v)$.
But we cannot have $j = n$ since $u \notin\ (u,v)$.
Thus $j < n$.

Next, there is $h$ with $0 < h< g-m$ such that $T^h(x) = S^j(x)$.
Indeed, setting $y = S^j(x)$, we have $u = T^{g-m-h}(y) = S^{n-j}(y)$ and thus $h < g-m$.
If $0 < h \leq -m$, then $T^h(x) = T^{m+h}(d_i) \in I\subset J$ contradicting the hypothesis that $m \in E_{J,T}(d_i)$.
If $-m < h < g-m$, then $T^h(x) = T^{m+h}(d_i) \in I$, contradicting the fact that $g \in E_{I,T}(d_i)$.
This shows that $n \in E_{I,S}(x)$ and thus that $u \in \Div(I,S)$.

Assume next that $g-m \leq 0$.
There is an integer $n$ with $g-m \leq n \leq 0$ such that $u = S^n(x)$.
Let us show that $n \in E_{I,S}(x)$.
Assume by contradiction that $n < -\rho^-_{I,S}(x)$.
Then there is some $j$ with $n < j < 0$ such that $S^j(x) = T^h(x)$.
Then $T^h(x) = T^{h+m}(d_i) \in I$ with $g < h+m < m$, in contradiction with the hypothesis that $m \in E_{I,T}(d_i)$.

We have proved that $u \in \Div(I,S)$.
If $v = r$, the proof that $I$ is admissible for $S$ is complete.
Otherwise, the proof that $v \in \Div(I,S)$ is similar to the proof for $u$.

Conversely, if  $I$ is admissible for $S$, there is some $x \in \Sep(S)$ and $g \in E_{I,S}(x)$ such that $u = S^g(x)$.
But $x = T^m(d_i)$ and since $u, x\in J$ there is some $n$ such that $u = T^n(d_i)$.

Assume for instance that $n > 0$ and suppose that there exists $k$ with $0 < k < n$ such that $T^k(d_i) \in (u, v)$.
Then, since $I \subset J$, $T^k(d_i)$ is of the form $S^h(x)$ with $0 < h < g$ which contradicts the fact that $g \in E_{I,S}(x)$.
Thus $n \in E_{I,T}(d_i)$ and $u \in \Div(I,T)$.

The proof is similar in the case $n \leq 0$.

If $v = r$, we have proved that $I$ is admissible for $T$.
Otherwise, the proof that $v \in \Div(I,T)$ is similar.

Finally, assume that $I$ is admissible for $T$ (and thus for $S$).
For any $d_i \in \Sep(T)$, one has
$$
\rho^-_{I,T}(d_i) \geq \rho^-_{J,T}(d_i)
\quad \text{ and } \quad
\rho^+_{I,T}(d_i) \geq \rho^+_{J,T}(d_i)
$$
showing that $\Div(J,T) \subset \Div(I,T)$.
\end{proof}

The last lemma is the key argument to prove Theorem~\ref{theo:birauzy2}.
%It is a tree version of the argument used by Rauzy in~\cite{Rauzy1979}.

\begin{lemma}
\label{lem:finiteness}
For any admissible interval $I \subset[\ell, r)$, the set $\F$ of sequences $\chi \in \{\varphi, \psi\}^*$ such that $I \subset D(\chi(T))$ is finite.
\end{lemma}
\begin{proof}
The set $\F$ is suffix-closed.
Indeed it contains the empty word because $[\ell, r)$ is admissible.
Moreover, for any $\xi, \chi \in \{\varphi,\psi\}^*$, one has $D(\xi \chi(T)) \subset D(\chi(T))$ and thus $\xi \chi \in \F$ implies $\chi \in \F$. 

The set $\F$ is finite.
Indeed, by Lemma~\ref{lem:bi22}, applied to $J = D(\chi(T))$, for any $\chi \in \F$, one has $\Div(D(\chi(T)),T) \subset \Div(I,T)$.
In particular, the boundaries of $D(\chi(T))$ belong to $\Div(I,T)$.
Since $\Div(I,T)$ is a finite set, this implies that there is a finite number of possible semi-intervals $D(\chi(T))$.
Thus there is  no infinite word with all its suffixes in $\F$.
Since the sequences $\chi$ are binary, this implies that $\F$ is finite.
\end{proof}

\begin{proofof}{ of Theorem~\ref{theo:birauzy2}}
  We consider $\chi$ as a word on the alphabet $\varphi,\psi$.
We first prove by induction on the length of $\chi$ that the domain $I$ of $\chi(T)$ is admissible and that the transformation induced by $T$ on $I$ is $\chi(T)$.
It is true for $|\chi| = 0$ since $[\ell, r)$ is admissible and $\chi(T) = T$.
Next, assume that $J = D(\chi(T))$ is admissible and that the transformation induced by $T$ on $J$ is $\chi(T)$.
Then $D(\varphi \chi(T))$ is admissible for $\chi(T)$ since $D(\varphi \chi(T)) = Y(\chi(T))$.
Thus $I = D(\varphi\chi(T))$ is admissible for $T$ by Lemma~\ref{lem:bi22} and the transformation induced by $T$ on $I$ is $\varphi \chi(T)$.
The same proof holds for $\psi \chi$.

Conversely, assume that $I$ is admissible.
By Lemma~\ref{lem:finiteness}, the set $\F$ of sequences $\chi \in \{\varphi, \psi\}^*$ such that $I\subset D(\chi(T))$ is finite. 

Thus there is some $\chi \in \F$ such that $\varphi \chi, \psi \chi \notin \F$.
If $I$ is strictly included in $D(\chi(T))$, then by Lemma~\ref{lem:init} applied to $\chi(T)$, we have $I\subset Y(\chi(T)) = D(\varphi\chi(T))$ or $I\subset Z(\chi(T))=D(\psi\chi(T))$, a contradiction.
Thus $I = D(\chi(T))$.
\end{proofof}

We close this subsection with a result concerning the dynamics of the branching induction.

\begin{theorem}
\label{theo:length}
For any sequence $(T_n)_{n \geq 0}$ of regular interval exchange transformations such that $T_{n+1} = \varphi(T_n)$ or $T_{n+1} = \psi(T_n)$ for all $n \geq 0$, the length of the domain of $T_n$ tends to $0$ when $n \rightarrow \infty$.
\end{theorem}
\begin{proof}
Assume the contrary and let $I$ be an open interval included in the domain of $T_n$ for all $n \geq 0$.
The set $\Div(I,T) \cap I$ is formed of $s$ points.
For any pair $u, v$ of consecutive elements of this set, the semi-interval $[u, v)$ is admissible.
By Lemma~\ref{lem:finiteness}, there is an integer $n$ such that the domain of $T_n$ does not contain $[u, v)$, a contradiction.
\end{proof}

\subsection{Equivalent transformations}
Let $[\ell_1, r_1)$, $[\ell_2, r_2)$ be two semi-intervals of the real line.
Let $T_1 = T_{\lambda_1, \pi_1}$ be an $s$-interval exchange transformation relative to a partition of $[\ell_1, r_1)$ and $T_2 = T_{\lambda_2, \pi_2}$ another $s$-interval exchange transformations relative to $[\ell_2, r_2)$.
We say that $T_1$ and $T_2$ are \emph{equivalent} if $\pi_1 = \pi_2$ and $\lambda_1 = c \lambda_2$ for some $c > 0$.
Thus, two interval exchange transformations are equivalent if we can obtain the second from the first by a rescaling following by a translation.
We denote by $\left)T_{\lambda, \pi}\right]$ the equivalence class of $T_{\lambda, \pi}$.

\begin{example}
\label{ex:simphi6T}
Let $S = T_{\mu,\pi}$ be the $3$-interval exchange transformation on a partition of the semi-interval $)2\alpha, 1)$, with $\alpha= (3-\sqrt{5})/2$, represented in Figure~\ref{fig:simphi6T}.
$S$ is equivalent to the transformation $T = T_{\lambda,\pi}$ of Example~\ref{ex:2alpha}, with length vector $\lambda = \left( 1-2\alpha, \alpha, \alpha \right)$ and permutation the cycle $\pi = (132)$.
Indeed the length vector $\mu = \left( 8\alpha-3, 2-5\alpha, 2-5\alpha \right)$ satisfies $\mu = \frac{2-5\alpha}{\alpha} \lambda$.

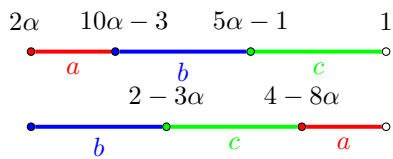
\begin{figure}[hbt]
\centering
\tikzset{node/.style={circle,draw,minimum size=0.1cm,inner sep=0pt}}
\tikzset{title/.style={circle,minimum size=0.1cm,inner sep=0pt}}
\begin{tikzpicture}
\node[node,fill=red](2alphah)at(0,1){};\node[title]at(0,1.4){$2\alpha \, \, \,$};
\node[node,fill=blue](10alpha-3h)at(1.12,1){};\node[title]at(1.12,1.4){$\, \, \, \, 10\alpha-3$};
\node[node,fill=green](5alpha-1h)at(2.92,1){};\node[title]at(2.92,1.4){$5\alpha-1$};
\node[node](1h)at(4.72,1){};\node[title]at(4.72,1.4){$1$};
\draw[color=red,line width=1.5,below](2alphah)edge node{$a$}(10alpha-3h);
\draw[color=blue,line width=1.5,below](10alpha-3h)edge node{$b$}(5alpha-1h);
\draw[color=green,line width=1.5,below](5alpha-1h)edge node{$c$}(1h);

\node[node,fill=blue](2alphab)at(0,0){};
\node[node,fill=green](2-3alphab)at(1.8,0){};\node[title]at(1.8,.4){$2-3\alpha$};
\node[node,fill=red](4-8alphab)at(3.6,0){};\node[title]at(3.6,.4){$4-8\alpha$};
\node[node](1b)at(4.72,0){};
\draw[color=blue,line width=1.5,below](2alphab)edge node{$b$}(2-3alphab);
\draw[color=green,line width=1.5,below](2-3alphab)edge node{$c$}(4-8alphab);
\draw[color=red,line width=1.5,below](4-8alphab)edge node{$a$}(1b);
\end{tikzpicture}
\caption{The transformation $S$.}
\label{fig:simphi6T}
\end{figure}
\end{example}

Note that if $T$ is a minimal (resp. regular) interval exchange transformation and $[S] = [T]$, then $S$ is also minimal (resp. regular).

For an interval exchange transformation $T$ we consider the directed labeled graph $\mathcal{G}(T)$, called the \emph{induction graph} of $T$, defined as follows.
The vertices are the equivalence classes of transformations obtained starting from $T$ and applying all possible $\chi \in \left\{ \psi, \varphi \right\}^*$.
There is an edge labeled $\psi$ (resp. $\varphi$) from a vertex $[S]$ to a vertex $[U]$ if and only if $U = \psi(S)$ (resp $\varphi(S)$) for two transformations $S \in [S]$ and $U \in [U]$.

\begin{example}
\label{ex:equivalencegraph}
Let $\alpha = \frac{3-\sqrt{5}}{2}$ and $R$ be a rotation of angle $\alpha$. 
It is a $2$-interval exchange transformation on $[0,1)$ relative to the partition $[0, 1-\alpha)$, $[1-\alpha, 1)$.
The induction graph $\mathcal{G}(R)$ of the transformation is represented in the left of Figure~\ref{fig:graphrotation}.
\end{example}

Note that for a $2$-interval exchange transformation $T$, one has $[\psi(T)] = [\varphi(T)]$, whereas in general the two transformations are not equivalent.

The induction graph of an interval exchange transformation can be infinite.
A sufficient condition for the induction graph to be finite is given in Section~\ref{sec:quadratic}.
$$ $$
Let us now introduce a variant of this equivalence relation (and of the related graph).
%We will now consider  two transformation as equivalent up to reflection (and up to the separation points).
%In Section~\ref{subsec:nc} we will justify this choise in terms of natural coding of two specular points.

For a $k$-interval exchange transformation $T = T_{\lambda, \pi}$, with length vector $\lambda = \left( \lambda_1, \lambda_2, \ldots, \lambda_k \right)$, we define the \emph{mirror transformation} $\widetilde{T} = T_{\widetilde{\lambda}, \tau \circ \pi}$ of $T$, where $\widetilde{\lambda} = \left( \lambda_k, \lambda_{k-1}, \ldots, \lambda_1 \right)$ and $\tau : i \mapsto (k-i+1)$ is the permutation that reverses the names of the semi-intervals.

Given two interval exchange transformations $T_1$ and $T_2$ on the same alphabet relative to two partitions of two semi-intervals $[\ell_1, r_1)$ and $[\ell_2, r_2)$ respectively, we say that $T_1$ and $T_2$ are \emph{similar} either if $[T_1] = [T_2]$ or $[T_1] = [\widetilde{T_2}]$.
Clearly, similarity is also an equivalent relation.
We denote by $\langle T \rangle$ the class of transformations similar to $T$.

\begin{example}
\label{ex:phi6T}
Let $T$ be the interval exchange transformation of Example~\ref{ex:2alpha}.
The transformation $U = \varphi^6(T)$ is represented in Figure~\ref{fig:phi6T} (see also Example~\ref{ex:u}).
It is easy to verify that $U$ is similar to the transformation $S$ of Example~\ref{ex:simphi6T}.
Indeed, we can obtain the second transformation (up to the separation points and the end points) by taking the mirror image of the domain.

Note that the order of the labels, i.e. the order of the letters of the alphabet, may be different from the order of the original transformation.

\begin{figure}[hbt]
\centering
\tikzset{node/.style={circle,draw,minimum size=0.1cm,inner sep=0pt}}
\tikzset{title/.style={circle,minimum size=0.1cm,inner sep=0pt}}
\begin{tikzpicture}
\node[node,fill=blue](2alphah)at(0,1){};\node[title]at(0,1.4){$2\alpha$};
\node[node,fill=red](2-3alphah)at(1.8,1.0){};\node[title]at(1.8,1.4){$2-3\alpha$};
\node[node,fill=green](4-8alphah)at(3.6,1.0){};\node[title]at(3.6,1.4){$4-8\alpha$};
\node[node](1h)at(4.72,1){};\node[title]at(4.72,1.4){$1$};
\draw[color=blue,line width=1.5,below](2alphah)edge node{$b$}(2-3alphah);
\draw[color=red,line width=1.5,below](2-3alphah)edge node{$a$}(4-8alphah);
\draw[color=green,line width=1.5,below](4-8alphah)edge node{$c$}(1h);

\node[node,fill=green](2alphab)at(0,0){};
\node[node,fill=blue](10alpha-3b)at(1.12,0){};\node[title]at(1.12,.4){$10\alpha-3$};
\node[node,fill=red](5alpha-1b)at(2.92,0){};\node[title]at(2.92,.4){$5\alpha-1$};
\node[node](1b)(4.72,0){};
\draw[color=green,line width=1.5,below](2alphab)edge node{$c$}(10alpha-3b);
\draw[color=blue,line width=1.5,below](10alpha-3b)edge node{$b$}(5alpha-1b);
\draw[color=red,line width=1.5,below](5alpha-1b)edge node{$a$}(1b);
\end{tikzpicture}
\caption{The transformation $U$.}
\label{fig:phi6T}
\end{figure}
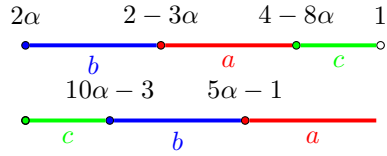
\end{example}

As  the equivalence relation, similarity  also preserves minimality and regularity.

Let $T$ be an interval exchange transformation.
We denote by
$$\mathcal{S}(T) = \bigcup_{n \in \Z} T^n \big(\Sep(T)\big)$$
the union of the orbits of the separation points.
Let $S$ be an interval exchange transformation similar to $T$.
Thus, there exists a bijection $f : D(T) \setminus \mathcal{S}(T) \to D(S) \setminus \mathcal{S}(S)$. This bijection is given by an affine transformation, namely a rescaling following by a translation if $T$ and $S$ are equivalent and a rescaling following by a translation and a reflection otherwise.
By the previous remark, if $T$ is a minimal exchange interval transformation and $S$ is similar to $T$, then the two interval exchange sets $\cL(T)$ and $\cL(S)$ are equal up to permutation, that is there exists a permutation $\pi$ such that one for every $w=a_0 a_1 \cdots a_{n-1} \in \cL(T)$ there exists a unique word $v = b_0 b_1 \cdots b_{n-1} \in \cL(S)$ such that $b_i = \pi(a_i)$ for all $i = 1, 2, \ldots n-1$.

In a similar way as before, we can use the similarity in order to construct a graph.
For an interval exchange transformation $T$ we define $\widetilde{\mathcal{G}}(T)$ the \emph{modified induction graph} of $T$ as the directed (unlabeled) graph with vertices the similar classes of transformations obtained starting from $T$ and applying all possible $\chi \in \left\{ \psi, \varphi \right\}^*$ and an edge from $\langle S \rangle$ to $\langle U \rangle$ if $U=\psi(S)$ or $U=\varphi(S)$ for two transformations $S \in \langle S \rangle$ and $U \in \langle U \rangle$.

Note that this variant appears naturally when considering the Rauzy induction of a $2$-interval exchange transformation as a continued fraction expansion.
There exists a natural bijection between the closed interval $[0,1]$ of the real line and the set of $2$-interval exchange transformation given by the map $x \mapsto T_{\lambda, \pi}$ where $\pi = (12)$ and $\lambda = \left( \lambda_1, \lambda_2 \right)$ is the length vector such that $x = \frac{\lambda_1}{\lambda_2}$.

In this view, the Rauzy induction corresponds to the Euclidean algorithm, i.e. the map $\mathcal{E} : \R_+^2 \to \R_+^2$ given by
\begin{displaymath}
\mathcal{E}(\lambda_1, \lambda_2) =
\begin{cases}
\left( \lambda_1 - \lambda_2, \lambda_2 \right)
&\text{ if $\lambda_1 \geq \lambda_2$}\\
\left( \lambda_1, \lambda_2 - \lambda_2 \right)
&\text{otherwise}.
\end{cases}
\end{displaymath}

Applying iteratively the Rauzy induction starting from $T$ corresponds then to the continued fraction expansion of $x$.

\begin{example}
Let $\alpha$ and $R$ be as in Example~\ref{ex:equivalencegraph}.
The modified induction graph $\widetilde{\mathcal{G}}(R)$ of the transformation is represented on the right of Figure~\ref{fig:graphrotation}.
Note that the ratio of the two lengths of the semi-intervals exchanged by $T$ is
\begin{displaymath}
  \frac{1-\alpha}{\alpha} = \frac{1+\sqrt{5}}{2} = \phi = 1 + \frac{1}{1+\frac{1}{1 + \cdots}}.
  \end{displaymath}

\begin{figure}[hbt]
\centering
\tikzset{node/.style={circle,draw,minimum size=0.5cm,inner sep=0.1pt}}
\tikzset{title/.style={circle,minimum size=0.1cm,inner sep=0pt}}
\tikzstyle{loop right}=[in=-30,out=30,loop]
\begin{tikzpicture}
\node[node](R1)at(0,0){$[R]$};
\node[node](R2)at(2,0){};
\draw[bend left,->,above](R1)edge node{$\psi, \varphi$}(R2);
\draw[bend left,->,below](R2)edge node{$\psi, \varphi$}(R1);

\node[node](R)at(5.5,0){$\langle R \rangle$};
\draw(R)edge[loop right,->](R){};
\end{tikzpicture}
\caption{Induction graph and modified induction graph of the rotation $R$ of angle $\alpha$= $(3-\sqrt{5})/2$.}
\label{fig:graphrotation}
\end{figure}
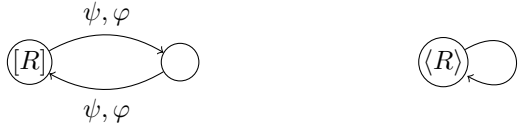
\end{example}

\subsection{Induction and automorphisms}
Let $T = T_{\lambda,\pi}$ be a regular interval exchange on $[\ell,r)$ relative to $(I(a))_{a \in A}$.
Set $A=\{a_1,\ldots,a_k\}$.
Recall now from Subsection~\ref{subsec:nc} that for any $z \in [\ell,r)$, the natural coding of $T$ relative to $z$ is the infinite word $\Sigma_T(z) = b_0 b_1 \cdots$ on the alphabet $A$ with $b_n \in A$ defined for $n \geq 0$ by $b_n = a$ if $T^n(z) \in I_a$.

Denote by $\theta_1$, $\theta_2$ the morphisms from $A^*$ into itself defined by
$$
\theta_1(a) =
\begin{cases}
a_{\pi(k)}a_k	&	\text{if $a = a_{\pi(k)}$} \\
a		&	\text{otherwise}
\end{cases},
\qquad
\theta_2(a) =
\begin{cases}
a_{\pi(k)}a_k	&	\text{if $a=a_k$} \\
a		&	\text{otherwise}
\end{cases}.
$$
The morphisms $\theta_1, \theta_2$ extend to automorphisms of the free group on $A$.

\begin{proposition}
\label{pro:autoelem}
Let $T$ be a regular interval exchange transformation on the alphabet $A$ and let $S = \psi(T)$, $I = Z(T)$.
There exists an automorphism $\theta$ of the free group on $A$ such that $\Sigma_T(z) = \theta(\Sigma_S(z))$ for any $z \in I$.
\end{proposition}
\begin{proof}
Assume first that $\gamma_{k} < \delta_{\pi(k)}$ (Case 0).
We have $Z(T) = [\ell, \delta_{\pi(k)})$ and for any $x \in Z(T)$,
$$
S(z) =
\begin{cases}
T^2(z)	&	\text{if $z \in K(a_{\pi(k)}) = I(a_{\pi(k)})$} \\ 
T(z)	&	\text{otherwise}.
\end{cases}
$$
We will prove by induction on the length of $w$ that for any $z \in I$, $\Sigma_S(z) \in wA^*$ if and only if $\Sigma_T(z) \in \theta_1(w)A^*$.
The property is true if $w$ is the empty word.
Assume next that $w = av$ with $a \in A$ and thus that $z \in I(a)$.
If $a \neq a_{\pi(k)}$, then $\theta_1(a) = a$, $S(z) = T(z)$ and
$$
\Sigma_S(z) \in avA^* \Leftrightarrow \Sigma_S(S(z)) \in vA^* \Leftrightarrow \Sigma_T(T(z) )\in \theta_1(v)A^* \Leftrightarrow \Sigma_T(z) \in \theta_1(w)A^*.
$$
Otherwise, $\theta_1(a) = a_{\pi(k)}a_k$, $S(z) = T^2(z)$.
Moreover, $\Sigma_T(z) = a_{\pi(k)}a_k\Sigma_T(T^2(z))$ and thus
$$
\Sigma_S(z) \in avA^* \Leftrightarrow \Sigma_S(S(z)) \in vA^* \Leftrightarrow \Sigma_T(T^2(z)) \in \theta_1(v)A^* \Leftrightarrow \Sigma_T(z) \in \theta_1(w)A^*.
$$
If $Td_{\pi(k)} < d_{k}$ (Case 1), we have $Z(T) = [\ell,d_{k})$ and for any $z \in Z(T)$,
$$
S(z) =
\begin{cases}
T^2(z)	&	\text{if $z \in K(a_k) = T^{-1}I(a_k)$} \\
T(z)	&	\text{otherwise}.
\end{cases}
$$
As in Case 0, we will prove by induction on the length of $w$ that for any $z \in I$, $\Sigma_S(z) \in wA^*$ if and only if $\Sigma_T(z) \in \theta_2(w)A^*$.

The property is true if $w$ is empty.
Assume next that $w = av$ with $a \in A$.
If $a \neq a_k$, then $\theta_2(a) = a$, $S(z) = T(z)$ and $z \in K(a) \subset I(a)$.
Thus
$$
\Sigma_S(z) \in avA^* \Leftrightarrow \Sigma_S(S(z)) \in vA^* \Leftrightarrow \Sigma_T(T(z)) \in \theta_2(v)A^* \Leftrightarrow \Sigma_T(z) \in \theta_2(w)A^*.
$$
Next, if $a = a_k$, then $\theta_2(a) = a_{\pi(k)}a_k$, $S(z) = T^2(z)$ and $z \in K_{a_k} = T^{-1}(I_{a_k}) \subset I(a_{\pi(k)})$.
Thus
$$
\Sigma_S(z) \in avA^*\Leftrightarrow \Sigma_S(S(z)) \in vA^* \Leftrightarrow \Sigma_T(T^2(z)) \in \theta_2(v)A^* \Leftrightarrow \Sigma_T(z) \in \theta_2(w)A^*.
$$
where the last equivalence results from the fact that $\Sigma_T(z) \in a_{\pi(k)}a_kA^*$.
This proves that $\Sigma_T(z) = \theta_2(\Sigma_S(z))$.
\end{proof}

\begin{example}
Let $T$ be the transformation of Example~\ref{ex:2alpha}.
The automorphism $\theta_1$ is defined by
$$
\theta_1(a) = ac, \quad \theta_1(b) = b, \quad \theta_1(c) = c.
$$
The right Rauzy induction gives the transformation $S = \psi(T)$ computed in Example~\ref{ex:induced1}.
One has $\Sigma_S(\alpha) = bacba \cdots$ and $\Sigma_T(\alpha) = baccbac \cdots = \theta_1(\Sigma_S(\alpha))$.
\end{example}

We state the symmetrical version of Proposition~\ref{pro:autoelem} for left Rauzy induction.
The proof is analogous.

\begin{proposition}
\label{pro:autoelemsym}
Let $T$ be a regular interval exchange transformation on the alphabet $A$ and let $S = \varphi(T)$, $I = Y(T)$.
There exists an automorphism $\theta$ of the free group on $A$ such that $\Sigma_T(z) = \theta(\Sigma_S(z))$ for any $z \in I$.
\end{proposition}

Combining Propositions~\ref{pro:autoelem} and~\ref{pro:autoelemsym}, we obtain the following statement.

\begin{theorem}
\label{theo:inductionbi}
Let $T$ be a regular interval exchange transformation.
For $\chi \in \{\varphi,\psi\}^*$, let $S = \chi(T)$ and let $I$ be the domain of $S$.
There exists an automorphism $\theta$ of the free group on $A$ such that $\Sigma_T(z) = \theta(\Sigma_S(z))$ for all $z \in I$.
\end{theorem}
\begin{proof}
The proof follows easily by induction on the length of $\chi$ using Propositions~\ref{pro:autoelem} and~\ref{pro:autoelemsym}.
\end{proof}

Note that if the transformations $T$ and $S = \chi(T)$, with $\chi \in \left\{\psi, \varphi\right\}^*$ , are equivalent, then there exists a point $z_0 \in D(S) \subseteq D(T)$ such that $z_0$ is a fixed point of the isometry that transforms $D(S)$ into $D(T)$ (if $\chi$ is different from the identity map, this point is unique).
In that case one has $\Sigma_S (z_0) = \Sigma_T (z_0) = \theta\left(\Sigma_S (z_0)\right)$ for an appropriate automorphism $\theta$, i.e. $\Sigma_T (z_0)$ is a fixed point
of an appropriate automorphism.

We now prove the following statement, which gives for regular
interval exchanges, a direct proof of the Return Theorem
(Theorem~\ref{theoremReturn}).

\begin{corollary}
Let $T$ be a regular interval exchange transformation and let $X=X(T)$.
For $w \in \cL(T)$, the set $\RR_X(w)$ is a basis of the free group on $A$.
\end{corollary}
\begin{proof}
By Proposition~\ref{pro:jadm}, the semi-interval $J(w)$ is admissible.
By Theorem~\ref{theo:birauzy2} there is a sequence $\chi \in \{\varphi,\psi\}^*$ such that $D(\chi(T)) = J(w)$.
Moreover, the transformation $S = \chi(T)$ is the transformation induced by $T$ on $J(w)$.
By Theorem~\ref{theo:inductionbi} there is an automorphism $\theta$ of the free group on $A$ such that $\Sigma_T(z) = \theta(\Sigma_S(z))$ for any $z \in J(w)$.

By Lemma~\ref{lem:returnsind}, we have $x \in \RR_X(w)$ if and only if $\Sigma_T(z) = x\Sigma_T(S(z)))$ for some $z \in J(w)$.
This implies that $\RR_X(w) = \theta(A)$.
Indeed, for any $z \in J(w)$, let $a$ is the first letter of $\Sigma_S(z)$.
Then
$$
\Sigma_T(z) = \theta(\Sigma_S(z)) = \theta(a\Sigma_S(S(z))) = \theta(a)\theta(\Sigma_S(Sz)) = \theta(a)\Sigma_T(S(z)).
$$
Thus $x \in \RR_X(w)$ if and only if there is $a\in A$ such that $x = \theta(a)$.
This proves that the set $\RR_X(w)$ is a basis of the free group on $A$.
\end{proof}

We illustrate  this result with the following examples.

\begin{example}\label{exampleAutomorphisms}
We consider again the transformation $T$ of Example~\ref{ex:2alpha} and $X =X(T)$.
We have $\RR_X(c) = \{bac,bbac,c\}$ (see Example~\ref{ex:returns}).
We represent in Figure~\ref{fig:illustrate} the sequence $\chi$ of Rauzy inductions such that $J(c)$ is the domain of $\chi(T)$.

\begin{figure}[hbt]
\centering
\tikzset{node/.style={circle,draw,minimum size=0.1cm,inner sep=0pt}}
\tikzset{title/.style={circle,minimum size=0.1cm,inner sep=0pt}}
\begin{tikzpicture}

\node[node,fill=red](0h)at(0,0){};
\node[node,fill=blue](1-2alpha)at(0,.78){};
\node[node,fill=green](1-alpha)at(0,2.06){};
\node[node](1h)at(0,3.3){};
\draw[color=red,line width=1.5,left](0h)edge node{$a$}(1-2alpha);
\draw[color=blue,line width=1.5,left](1-2alpha)edge node{$b$}(1-alpha);
\draw[color=green,line width=1.5,left](1-alpha)edge node{$c$}(1h);

\node[node,fill=blue](0b)at(1,0){};
\node[node,fill=green](alpha)at(1,1.2){};
\node[node,fill=red](2alpha)at(1,2.54){};\node[node](1b)at(1,3.3){};
\draw[color=blue,line width=1.5,left](0b)edge node{$b$}(alpha);
\draw[color=green,line width=1.5,left](alpha)edge node{$c$}(2alpha);
\draw[color=red,line width=1.5,left](2alpha)edge node{$a$}(1b);

\node[title]at(2,2){\LARGE$\edge{a\mapsto ac}$};%\put(18,0){$\ldots$}

\node[node,fill=red](0h)at(3,0){};
\node[node,fill=blue](1-2alpha)at(3,.78){};
\node[node,fill=green](1-alpha)at(3,2.06){};
\node[node](2alphah)at(3,2.54){};
\draw[color=red,line width=1.5,left](0h)edge node{$a$}(1-2alpha);
\draw[color=blue,line width=1.5,left](1-2alpha)edge node{$b$}(1-alpha);
\draw[color=green,line width=1.5,left](1-alpha)edge node{$c$}(2alphah);

\node[node,fill=blue](0b)at(4,0){};
\node[node,fill=green](alpha)at(4,1.2){};
\node[node,fill=red](4alpha-1)at(4,1.7){};
\node[node](2alpha)at(4,2.54){};
\draw[color=blue,line width=1.5,left](0b)edge node{$b$}(alpha);
\draw[color=green,line width=1.5,left](alpha)edge node{$c$}(4alpha-1);
\draw[color=red,line width=1.5,left](4alpha-1)edge node{$a$}(2alpha);

\node[title]at(5,2){\LARGE$\edge{a\mapsto ba}$};%\put(48,25.4){$\ldots$}

\node[node,fill=red](0h)at(6,.78){};
\node[node,fill=blue](1-2alpha)at(6,1.57){};
\node[node,fill=green](1-alpha)at(6,2.06){};
\node[node](2alpha)at(6,2.54){};
\draw[color=red,line width=1.5,left](0h)edge node{$a$}(1-2alpha);
\draw[color=blue,line width=1.5,left](1-2alpha)edge node{$b$}(1-alpha);
\draw[color=green,line width=1.5,left](1-alpha)edge node{$c$}(2alpha);

\node[node,fill=blue](0b)at(7,.78){};
\node[node,fill=green](alpha)at(7,1.2){};
\node[node,fill=red](4alpha-1)at(7,1.7){};
\node[node](2alpha)at(7,2.54){};
\draw[color=blue,line width=1.5,left](0b)edge node{$b$}(alpha);
\draw[color=green,line width=1.5,left](alpha)edge node{$c$}(4alpha-1);
\draw[color=red,line width=1.5,left](4alpha-1)edge node{$a$}(2alpha);

\node[title]at(8,2){\LARGE$\edge{b\mapsto ba}$};%\put(78,25.4){$\ldots$}

\node[node,fill=red](0h)at(9,1.2){};
\node[node,fill=blue](1-2alpha)at(9,1.57){};
\node[node,fill=green](1-alpha)at(9,2.06){};
\node[node](2alpha)at(9,2.54){};
\draw[color=red,line width=1.5,left](0h)edge node{$a$}(1-2alpha);
\draw[color=blue,line width=1.5,left](1-2alpha)edge node{$b$}(1-alpha);
\draw[color=green,line width=1.5,left](1-alpha)edge node{$w$}(2alpha);

\node[node,fill=green](alpha)at(10,1.2){};
\node[node,fill=blue](4alpha-1)at(10,1.7){};
\node[node,fill=red](7alpha-2)at(10,2.2){};
\node[node](2alpha)at(10,2.54){};
\draw[color=green,line width=1.5,left](alpha)edge node{$c$}(4alpha-1);
\draw[color=blue,line width=1.5,left](4alpha-1)edge node{$b$}(7alpha-2);
\draw[color=red,line width=1.5,left](7alpha-2)edge node{$a$}(2alpha);

\end{tikzpicture}
\caption{The sequence $\chi \in \{\varphi,\psi\}^*$}
\label{fig:illustrate}
\end{figure}
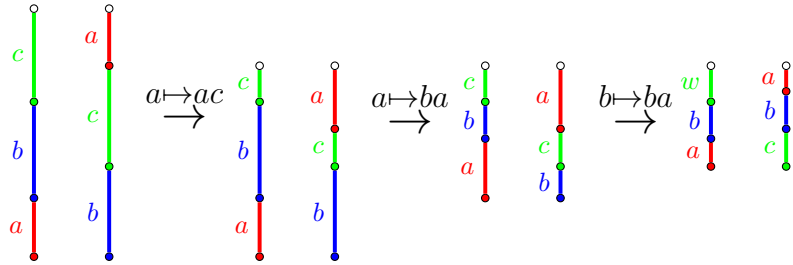

The sequence is composed of a right induction followed by two left inductions.
We have indicated on each edge the associated automorphism (indicating only the image of the letter which is modified).
We have $\chi = \varphi^2\psi$ and the resulting composition $\theta$ of automorphisms gives
$$
\theta(a) = bac, \quad \theta(b) = bbac, \quad \theta(c)=c.
$$
Thus $\RR_X(c) = \theta(A)$.
\end{example}

\begin{example}
\label{ex:u}
Let $T$ and $X$ be as in the preceding example.
Let $U$ be the transformation induced by $T$ on $J_a$.
We have $U = \varphi^6(T)$ and a computation shows that for any $z \in J_a$, $\Sigma_T(z) = \theta(\Sigma_U(z))$ where $\theta$ is the automorphism of the free group on $A = \{a,b,c\}$ which is the coding morphism for $\RR_X(a)$ defined by:
$$
\theta(a) = ccba, \quad \theta(b) = cbba, \quad \theta(c)=ccbba.
$$
One can verify that $\cL(U) = \cL(S)$, where $S$ is the transformation obtain from $T$ by permuting the labels of the intervals according to the permutation $\pi = (acb)$.

Note that $\cL(U) = \cL(S)$ although $S$ and $U$ are not identical, even up to rescaling the intervals.
Actually, the rescaling of $U$ to a transformation on $[0,1)$ corresponds to the mirror image of $S$, obtained by taking the image of the intervals by a symmetry centered at $1/2$.
\end{example}

Note that in the above examples, all lengths of the intervals belong to the quadratic number field $\Q[\sqrt{5}]$.

In the next Section we will prove that if a regular interval exchange transformation $T$ is defined over a quadratic field, then the family of transformations obtained from $T$ by the Rauzy inductions contains finitely many distinct transformations up to rescaling.

\section{Interval exchange over a quadratic field}
\label{sec:quadratic}

An interval exchange transformation is said to be defined 
over a number field $K \subset \R$ if the lengths of all exchanged
semi-intervals belong to $K$.
 Let $T$ be a minimal interval exchange transformation on semi-intervals defined
over a quadratic number field. Let $(T_n )_{n\ge 0}$ be a sequence of interval exchange transformation such that $T_0 = T$
and $T_n+1$ is the transformation induced by $T_n$ on one of its exchanged semi-intervals $I_n$ . 
\begin{theorem}[Boshernitzan, Carroll]\label{theoremBoshernitzanCarroll}
If $T$ is defined over a quadratic number field, up to rescaling
all semi-intervals $I_n$ to the same length, the sequence $(T_n )$ contains finitely many distinct transformations.
\end{theorem}

It is possible to generalize this result and prove that, under the above hypothesis on the lengths of the semi-intervals and up to rescaling and translation, there are finitely many transformations obtained by the branching Rauzy induction defined in Section~\ref{sec:rauzy}.

\begin{theorem}
\label{theo:quadratic}
Let $T$ be a regular interval exchange transformation defined over a quadratic field.
The family of all induced transformation of $T$ over an admissible semi-interval contains finitely many distinct transformations up to equivalence.
\end{theorem}

An immediate corollary of Theorem~\ref{theo:quadratic} is the following.

\begin{corollary}
\label{cor:graph}
Let $T$ be a regular interval exchange transformation defined over a quadratic field. Then the induction graph $\mathcal{G}(T)$ and the modified induction graph $\widetilde{\mathcal{G}}(T)$ are finite.
\end{corollary}

\begin{example}
Let $T$ be the regular interval exchange transformation of Example~\ref{ex:2alpha}.
The modified induction graph $\widetilde{\mathcal{G}}(T)$ is represented in Figure~\ref{fig:graph}.
The transformation $T$ belongs to the similarity class $\langle T_1 \rangle$ as well as transformations $S$ of Example~\ref{ex:simphi6T} and $U$ of Example~\ref{ex:phi6T}.
The transformations $\psi(T)$ and $\psi^2(T)$ of Example~\ref{ex:induced2} belongs respectively to classes $\langle T_2 \rangle$ and $\langle T_4 \rangle$, while the two last transformations of Figure~\ref{fig:illustrate}, namely $\varphi \psi(T)$ and $\varphi^2 \psi(T)$, belongs respectively to $\langle T_5 \rangle$ and $\langle T_7 \rangle$.
Finally, the left Rauzy induction sequence from $T$ to $U = \varphi^6 (T)$ corresponds to the loop $\langle T_1 \rangle \to \langle T_3 \rangle \to \langle T_4 \rangle \to \langle T_6 \rangle \to \langle T_7 \rangle \to \langle T_8 \rangle \to \langle T_1 \rangle$ in $\widetilde{\mathcal{G}}(T)$ (indicated by thick edges in Figure~\ref{fig:graph}).

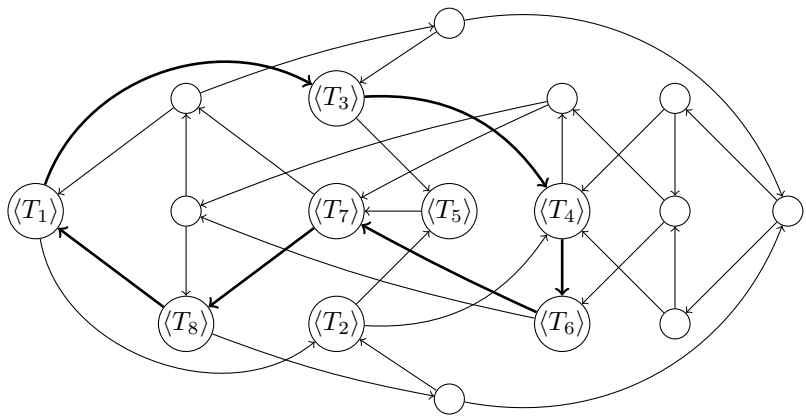
\begin{figure}[hbt]
\centering
\tikzset{node/.style={circle,draw,minimum size=0.4cm,inner sep=0.1pt}}
\tikzset{title/.style={circle,minimum size=0.1cm,inner sep=0pt}}
\begin{tikzpicture}

\node[node](1)at(-5,0){$\langle T_1 \rangle$};
\node[node](2)at(-1,1.5){$\langle T_3 \rangle$};
\node[node](3)at(-1,-1.5){$\langle T_2 \rangle$};
\node[node](4)at(2,0){$\langle T_4 \rangle$};
\node[node](5)at(.5,0){$\langle T_5 \rangle$};
\node[node](6)at(2,1.5){};
\node[node](7)at(2,-1.5){$\langle T_6 \rangle$};
\node[node](8)at(-1,0){$\langle T_7 \rangle$};
\node[node](9)at(-3,0){};
\node[node](10)at(-3,1.5){};
\node[node](11)at(-3,-1.5){$\langle T_8 \rangle$};
\node[node](13)at(.5,2.5){};
\node[node](14)at(.5,-2.5){};
\node[node](15)at(5,0){};
\node[node](16)at(3.5,1.5){};
\node[node](17)at(3.5,-1.5){};
\node[node](18)at(3.5,0){};

\draw[bend left=50,line width=1,->](1)edge node{}(2);
\draw[bend right=60,->](1)edge node{}(3);
\draw[bend left,line width=1,->](2)edge node{}(4);
\draw[->](2)edge node{}(5);
\draw[bend right,->](3)edge node{}(4);
\draw[->](3)edge node{}(5);
\draw[->](4)edge node{}(6);
\draw[line width=1,->](4)edge node{}(7);
\draw[->](5)edge node{}(8);
\draw[bend right=2,->](6)edge node{}(8);
\draw[bend right=5,->](6)edge node{}(9);
\draw[bend left=2,line width=1,->](7)edge node{}(8);
\draw[bend left=5,->](7)edge node{}(9);
\draw[->](8)edge node{}(10);
\draw[line width=1,->](8)edge node{}(11);
\draw[->](9)edge node{}(10);
\draw[->](9)edge node{}(11);
\draw[->](10)edge node{}(1);
\draw[bend left=5,->](10)edge node{}(13);
\draw[->,line width=1](11)edge node{}(1);
\draw[bend right=5,->](11)edge node{}(14);
\draw[->](13)edge node{}(2);
\draw[bend left=40,->](13)edge node{}(15);
\draw[->](14)edge node{}(3);
\draw[bend right=40,->](14)edge node{}(15);
\draw[->](15)edge node{}(16);
\draw[->](15)edge node{}(17);
\draw[->](16)edge node{}(4);
\draw[->](16)edge node{}(18);
\draw[->](17)edge node{}(4);
\draw[->](17)edge node{}(18);
\draw[->](18)edge node{}(6);
\draw[->](18)edge node{}(7);

\end{tikzpicture}
\caption{Modified induction graph of the transformation $T$.}
\label{fig:graph}
\end{figure}
\end{example}

\subsection{Primitive substitution shifts}
\label{subsec:morphic}

In this section we prove
 an important property of interval exchange transformations
 defined over a quadratic field, namely that the related interval exchange 
shifts are primitive substitution shifts.

\begin{theorem}
\label{theo:morphic}
Let $T$ be a regular interval exchange transformation defined over a quadratic field.
The interval exchange shift $X(T)$ is a primitive substitution shift.
\end{theorem}

\begin{example}
Let $T = T_{\lambda, \pi}$ be the transformation of Example~\ref{ex:2alpha} (see also~\ref{ex:t}).
The shift $X(T)$ is a primitive substitution shift, as we have seen.
This can be obtained as a consequence of Theorem~\ref{theo:morphic}.
Indeed the transformation $T$ is regular and the length vector $\lambda = (1-2\alpha, \alpha, \alpha)$ belongs to $\Q\left[\sqrt{5}\right]^3$.
\end{example}

In order to prove Theorem~\ref{theo:morphic} we need some preliminary results.

\begin{proposition}
\label{pro:primitive}
Let $T, \chi(T)$ be two equivalent regular interval exchange transformations with $\chi \in \left\{ \varphi, \psi \right\}^*$.
There exists a primitive morphism $\theta$ and a point $z \in D(T)$ such that the natural coding of $T$ relative to $z$ is a fixed point of $\theta$.
\end{proposition}
\begin{proof}
Since $T$ is regular, it is minimal and thus the set $\cL(T)$ is uniformly recurrent.
Thus, there exists a positive integer $N$ such that every letter of the alphabet appears in every word of length $N$ of $\cL(T)$.
Moreover, by Theorem~\ref{theo:length}, applying iteratively the Rauzy induction, the length of the domains tends to zero.

Consider $T' = \chi^m (T)$, for a positive integer $m$, such that $D(T')~<~\varepsilon$, where $\varepsilon$ is the positive real number for which, by Lemma~\ref{lem:distance}, the first return map for every point of the domain is ``longer'' than $N$, i.e. $T'(z) = T^{n(z)}(z)$, with $n(z) \geq N$, for every $z \in D(T')$.

By Theorem~\ref{theo:inductionbi} and the remark following it, there exists an automorphism $\theta$ of the free group and a point $z \in D(T') \subseteq D(T)$ such that the natural coding of $T$ relative to $z$ is a fixed point of $\theta$, that is $\Sigma_T (z) = \theta\left( \Sigma_{T} (z) \right)$.

By the previous argument, the image of every letter by $\theta$ is longer than $N$, hence it contains every letter of the alphabet as a factor.
Therefore, $\theta$ is a primitive morphism.
\end{proof}

Using the previous results we can finally prove Theorem~\ref{theo:morphic}.

\begin{proofof}{ of Theorem~\ref{theo:morphic}}
By Theorem~\ref{theo:quadratic} there exists a regular interval transformation $S$ such that we can find in the induction graph $\mathcal{G}(T)$ a path from $[T]$ to $[S]$ followed by a cycle on $[S]$.
Thus, by Theorem~\ref{theo:inductionbi} there exist a point $z \in D(S)$ and two automorphisms $\theta, \eta$ of the free group such that $\Sigma_T(z) = \theta \left( \Sigma_S (z) \right)$, with $\Sigma_S (z)$ a fixed point of $\eta$.

By Proposition~\ref{pro:primitive} we can suppose, without loss of generality, that $\eta$ is primitive.
Therefore, $X(T)$ is a primitive substitution shift.
\end{proofof}

%%%%%%%%%%%%%%%%%%%%%%%%%%%%%%%%%%%%%%%%%%%%%%%%%%%%%%%%%%%%%%%%%%%
%section Linear Involutions
%%%%%%%%%%%%%%%%%%%%%%%%%%%%%%%%%%%%%%%%%%%%%%%%%%%%%%%%%%%%%%%%%%%
\section{Linear involutions}\label{sectionInvolutions}

Let $A$ be an alphabet of cardinality $k$ with an involution $\theta$
and the corresponding specular group $G_\theta$. Note that we allow
$\theta$ to have fixed points. Recall that, in the group $G_\theta$,
we have $\theta(a)=a^{-1}$. Thus, when $\theta$ has no fixed points,
the alphabet $A$ can be identified with $B\cup B^{-1}$
in such a way that $G\theta$ is the free group on $B$.

We consider two copies $I\times \{0\}$ and
$I\times \{1\}$ of an open interval $I$ of the real line and denote
 $\hat{I}=I\times \{0,1\}$.\index{symbols}{I@$\hat{I}$}
We call the sets $I\times \{0\}$ and $I\times \{1\}$ the two
\emph{components} of $\hat{I}$. We consider each component as an open interval.

A \emph{generalized permutation}\index{subject}{generalized!permutation}
 on $A$ of type $(\ell,m)$, with $\ell+m=k$,
  is a bijection $\pi:\{1,2,\ldots,k\}\rightarrow A$.
We represent it by a two line array
\begin{displaymath}
\pi=\begin{pmatrix} \pi(1)\ \pi(2)\ \ldots \pi(\ell)\\
\pi(\ell +1)\ \ldots \pi(\ell+m)
\end{pmatrix}
\end{displaymath}
A \emph{length data}\index{subject}{length!data}
 associated with $(\ell,m,\pi)$\index{symbols}{l@$(\ell,m,\pi)$} is a nonnegative
vector $\lambda\in \R_+^{A}=\R_+^{k}$ such that
\begin{displaymath}
\lambda_{\pi(1)}+\ldots+\lambda_{\pi(\ell)}=
\lambda_{\pi(\ell+1)}+\ldots+\lambda_{\pi(k)}
\text{ and }\lambda_a=\lambda_{a^{-1}}\text{ for all }a\in A.
\end{displaymath}

We consider a partition of $I\times \{0\}$ (minus $\ell-1$
points) in $\ell$ open intervals
$I_{\pi(1)},\ldots,I_{\pi(\ell)}$ of lengths $\lambda_{\pi(1)},\ldots,\lambda_{\pi(\ell)}$
and a partition of $I\times \{1\}$ (minus $m-1$ points) in $m$ open intervals
$I_{\pi(\ell+1)},\ldots,I_{\pi(\ell+m)}$ of lengths $\lambda_{\pi(\ell+1)},\ldots,\lambda_{\pi(\ell+m)}$. Let $\Sigma$ be the set of $k-2$ \emph{division points}
\index{subject}{division points} separating
the intervals $I_a$ for $a\in A$.

The \emph{linear involution}\index{subject}{linear involution} on $I$ relative to these data is the
map $T=\sigma_2\circ\sigma_1$ defined on the set
$\hat{I}\setminus\Sigma$, formed of $\hat{I}$ minus the $k-2$ division points,
and
which is
 the composition
of two involutions defined as follows. 
\begin{enumerate}
\item[(i)]The first involution $\sigma_1$\index{symbols}{sigma@$\sigma_1$}
 is defined on $\hat{I}\setminus\Sigma$.
It is such that for each $a\in A$, its restriction to $I_a$
is either a translation or a symmetry from $I_a$ onto $I_{a^{-1}}$.

\item[(ii)]The second involution $\sigma_2$\index{symbols}{sigma@$\sigma_2$}
exchanges the two components of
  $\hat{I}$.
It  is defined, for $(x,\delta)\in \hat{I}$,
by $\sigma_2(x,\delta)=(x,1-\delta)$. The image of $z$ by $\sigma_2$
is called the \emph{mirror image} of $z$.
\end{enumerate}
We also say that $T$ is a linear involution on $I$ and relative to
 the alphabet $A$
or that it is a $k$-linear involution to express the fact
that the alphabet $A$ has $k$ elements.

\begin{example}\label{exampleLinear}
Let $A=\{a,b,c,d,a^{-1},b^{-1},c^{-1},d^{-1}\}$ and
\begin{displaymath}
\pi=\begin{pmatrix}a&b&a^{-1}&c\\c^{-1}&d^{-1}&b^{-1}&d
\end{pmatrix}
\end{displaymath}
Let $T$ be the $8$-linear involution corresponding to the length data
represented in Figure~\ref{figureLinear} (we represent $I\times\{0\}$
above $I\times \{1\}$) with the assumption that the restriction
of $\sigma_1$ to $I_a$ and $I_d$ is a symmetry while its restriction
to $I_b,I_c$ is a translation.

\begin{figure}[hbt]
\centering
\tikzset{node/.style={circle,draw,minimum size=0.1cm,inner sep=0pt}}
\tikzset{title/.style={circle,minimum size=0.1cm,inner sep=0pt}}
\begin{tikzpicture}
%--interval nodes--
\node[node](h0)at(0,1){};\node[node](h1)at(2,1){};\node[node](h2)at(5,1){};
\node[node](h3)at(7,1){};\node[node](h4)at(11,1){};\node[title](0)at(12,1){$I\times\{0\}$};
\node[node](b0)at(0,0){};\node[node](b1)at(4,0){};
\node[node](b2)at(6,0){};\node[node](b3)at(9,0){};
\node[node](b4)at(11,0){};\node[title](1)at(12,0){$I\times\{1\}$};
%translation nodes
\node[title](z)at(.5,.8){$z$};\node[title](Tz)at(6.5,-.2){$Tz$};
\node[title](T2z)at(2.5,-.2){$T^2z$};

%interval lines
\draw[above,line width=2pt,color=red](h0)edge node{$a$}(h1);
\draw[above,line width=2pt,color=blue](h1)edge node{$b$}(h2);
\draw[above,line width=2pt,color=red](h2)edge node{$a^{-1}$}(h3);
\draw[above,line width=2pt,color=green](h3)edge node{$c$}(h4);
\draw[below,line width=2pt,color=green](b0)edge node{$c^{-1}$}(b1);
\draw[below,line width=2pt,color=violet](b1)edge node{$d^{-1}$}(b2);
\draw[below,line width=2pt,color=blue](b2)edge node{$b^{-1}$}(b3);
\draw[below,line width=2pt,color=violet](b3)edge node{$d$}(b4);

%translation lines
\draw[bend left,->](.5,1)edge node{}(6.5,1);\draw[->](6.5,1)edge node{}(6.5,0);
\draw[->](6.5,0)edge node{}(2.5,1);\draw[->](2.5,1)edge node{}(2.5,0);
\end{tikzpicture}
\caption{A linear involution.}\label{figureLinear}
\end{figure}
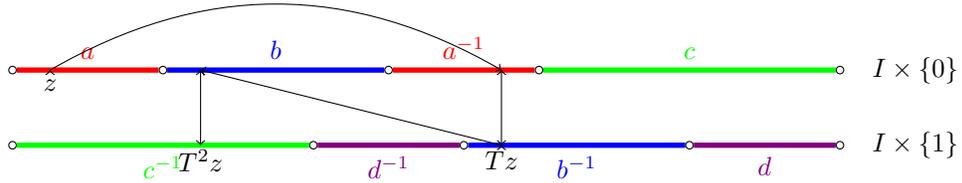
We indicate on the figure the effect of the transformation $T$ on a point
$z$ located in the left part of the interval $I_a$. The point
$\sigma_1(z)$ is located in the right part of $I_{a^{-1}}$ and the point
$T(z)=\sigma_2\sigma_1(z)$ is just below on the left of $I_{b^{-1}}$.
 Next, the point $\sigma_1T(z)$ is located on the left part of $I_b$
and the point $T^2(z)$ just below.
\end{example}
Thus the notion of linear involution is an extension of the notion 
of  interval exchange transformation in the following sense.
Assume that 
\begin{enumerate}
\item[(i)] $\ell=m$, 
\item[(ii)] for each letter $a\in A$, the interval
$I_a$ belongs to $I\times\{0\}$ if and only if $I_{a^{-1}}$ belongs
to $I\times\{1\}$,  
\item[(iii)] the restriction of $\sigma_1$
to each subinterval is a translation. 
\end{enumerate}
Then, the restriction of $T$
to $I\times \{0\}$ is an interval exchange (and so is its restriction to 
$I\times \{1\}$ which is the inverse of the first one). Thus,
in this case, $T$ is a pair of mutually inverse interval exchange transformations.

Note that we consider here interval exchange transformations
defined by a partition of an open interval  minus
$\ell-1$ points in $\ell$
 open intervals. The usual notion of interval exchange transformation
uses a partition of a semi-interval in a finite number of semi-intervals.
One recovers the usual notion of interval exchange transformation
 on a semi-interval
 by attaching to each open interval its left endpoint.

A linear involution $T$ is a bijection from $\hat{I}\setminus\Sigma$
onto $\hat{I}\setminus\sigma_2(\Sigma)$.
Since $\sigma_1,\sigma_2$ are involutions and $T=\sigma_2\circ\sigma_1$,
 the inverse of $T$
is $T^{-1}=\sigma_1\circ\sigma_2$.

The set $\Sigma$ of division points is also the set of singular points
of $T$ and their mirror images are the singular points of $T^{-1}$
(which are the points where $T$ (resp. $T^{-1}$) is not defined).
Note that these singular points $z$ may be `false' singularities, in the sense
that $T$ can have a continuous extension to an open neighborhood of $z$.

Two particular cases of linear involutions deserve attention.

A linear involution $T$ on the alphabet $A$ 
relative to a generalized permutation $\pi$ of type
$(\ell,m)$ 
is said to be \emph{nonorientable}\index{subject}{nonorientable linear involution}\index{subject}{linear involution!nonorientable} if there are indices $i,j\le \ell$ such that
$\pi(i)=\pi(j)^{-1}$  (and thus indices $i,j\ge \ell+1$
such that $\pi(i)=\pi(j)^{-1}$).  In other words,  there is  some $a\in A$ for which  $I_a$
and $I_{a^{-1}}$ belong to  the same component  of $\hat{I}$. Otherwise
$T$ is said to be \emph{orientable}.  \index{subject}{orientable!linear involution}\index{subject}{linear involution!orientable}

A linear involution $T=\sigma_2\circ \sigma_1$ on $I$
relative to the alphabet $A$
is said to be \emph{coherent}\index{subject}{coherent linear involution}
\index{subject}{linear involution!coherent} if, for each $a\in A$, the restriction
of $\sigma_1$ to $I_a$ is a translation if and only if $I_a$
and $I_{a^{-1}}$ belong to distinct components of $\hat{I}$. 

\begin{example}
The linear involution of Example~\ref{exampleLinear} is coherent.
\end{example}

Linear
involutions which are orientable and coherent
correspond to interval exchange transformations, 
whereas  orientable but non coherent  linear
involutions are called \emph{interval exchanges with flip}.
\index{subject}{interval exchange!with flip}

%\subsubsection{Orientation covering}

%Let $T$ be a linear involution on $I$ relative to $A$. Let
%$\varphi:\hat{I}\rightarrow \Z/2\Z$ be the map from $\hat{I}=I\times\{0,1\}$
%into the additive group $\Z/2\Z$ defined for $z\in \hat{I}$ by
%\begin{displaymath}
%\varphi(z)=\begin{cases}0&\text{if $z\in I_a$ with $a$ even}\\
%1&\text{otherwise}
%\end{cases}
%\end{displaymath}
%Let $S$ be the skew product of $T$ with $\Z/2\Z$ relative to $\varphi$.
%By definition $S$ is the transformation on $\hat{I}\times\{0,1\}$
%defined for $(z,i)\in\hat{I}\times\{0,1\}$ by
%\begin{displaymath}
%S(z,i)=(T(z),i+\varphi(z)).
%\end{displaymath}
%The transformation $S$ is called the \emph{orientation covering} of $T$
%(see\cite{BertheDelecroixDolcePerrinReutenauerRindone2014} for the justification of this term). The orientation covering is actually an interval exchange transformation with flip on each of the sets $I_0=I\times\{0,0\}\cup I\times\{1,1\}$
%and $I_1=I\times\{1,0\}\cup I\times\{0,1\}$. Indeed, 

%\subsubsection{Minimal involutions}

A \emph{connection}\index{subject}{connection!of linear involution}
\index{subject}{linear involution!connection} of a linear involution $T$ is a triple $(x,y,n)$ 
where $x$ is a singularity of $T^{-1}$, $y$ is a singularity of $T$,
$n\ge 0$ and $T^n x = y$. 

\begin{example}\label{exampleConnexion}
Let us consider the linear involution $T$ which is the same as in Example~\ref{exampleLinear} but such that the restriction of $\sigma_1$ to $I_c$ 
is a symmetry. Thus $T$ is not coherent. We assume that $I=]0,1[$,
that $\lambda_a=\lambda _d$. Let $x=(1-\lambda_d,0)$ and $y=(\lambda_a,0)$.

Then $x$ is a singularity of $T^{-1}$ ($\sigma_2(x)$ is the left endpoint
of $I_d$), $y$ is a singularity of $T$ (it is the right endpoint of $I_a$)
and $T(x)=y$. Thus $(x,1,y)$ is a connection.
\end{example}

\begin{example}\label{exampleInvolution3}
Let $T$ be the  linear involution on $I=]0,1[$
represented in Figure~\ref{figureLinear3}. We assume that the
restriction of $\sigma_1$ to $I_a$ is a translation
whereas the restriction to $I_b$ and $I_c$ is a symmetry.
We choose $(3-\sqrt{5})/2$ for the length of the interval $I_c$
(or $I_b$). With this choice, $T$ has no connections.
\begin{figure}[hbt]
\centering
\tikzset{node/.style={circle,draw,minimum size=0.1cm,inner sep=0pt}}
\tikzset{title/.style={circle,minimum size=0.1cm,inner sep=0pt}}
\begin{tikzpicture}
\node[node](h0)at(0,1){};\node[node](b)at(2,1){};\node[node](bbar)at(6,1){};\node[node](h1)at(10,1){};
\node[node](b0)at(0,0){};\node[node](cbar)at(4,0){};\node[node](abar)at(8,0){};\node[node](b1)at(10,0){};

\draw[above,line width=2pt,color=red](h0)edge node{$a$}(b);
\draw[above,line width=2pt,color=blue](b)edge node{$b$}(bbar);
\draw[above,line width=2pt,color=blue](bbar)edge node{$b^{-1}$}(h1);
\draw[below,line width=2pt,color=green](b0)edge node{$c$}(cbar);
\draw[below,line width=2pt,color=green](cbar)edge node{$c^{-1}$}(abar);
\draw[below,line width=2pt,color=red](abar)edge node{$a^{-1}$}(b1);
\end{tikzpicture}
\caption{A linear involution without connections.}\label{figureLinear3}
\end{figure}
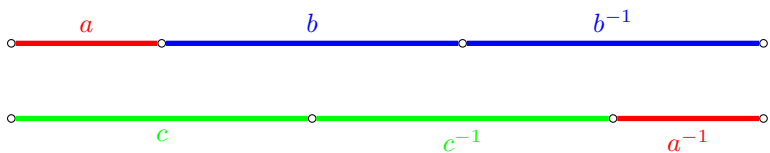
\end{example}

Let $T$ be a  linear involution without connections. Let 
\begin{equation}
O=\bigcup_{n\ge 0}T^{-n}(\Sigma) \quad \text{and }\quad \hat{O}=O\cup \sigma_2(O)\label{eqO}
\end{equation}
be respectively the negative orbit of the singular points and its
closure under mirror image.
Then $T$ is a bijection from $\hat{I}\setminus \hat{O}$ onto itself.
Indeed, assume that $T(z)\in \hat{O}$. If $T(z)\in O$ then
$z\in O$. Next if $T(z)\in \sigma_2(O)$, then
$T(z)\in \sigma_2(T^{-n}(\Sigma))=T^n(\sigma_2(\Sigma))$ for some $n\ge 0$. We cannot
have $n=0$ since $\sigma_2(\Sigma)$ is not in the image of $T$.
Thus $z\in T^{n-1}(\sigma_2(\Sigma))=\sigma_2(T^{-n+1}(\Sigma))\subset
\sigma_2(O)$. Therefore in both cases $z\in \hat{O}$. The converse
implication
is proved in the same way.
\subsection{Minimal linear involutions}
 A linear involution $T$ 
on $I$ without connections
is  \emph{minimal}\index{subject}{minimal!linear involution}
\index{subject}{linear involution!minimal} if for any point $z\in\hat{I}\setminus\hat{O}$ 
the nonnegative orbit of $z$  is dense in $\hat{I}$. 

Note that when a linear involution is  orientable, that is,  when it is a pair
of interval exchange transformations (with or without flips),   the interval exchange transformations
can be minimal although the linear involution is not since each component
of $\hat{I}$ is stable by the action of $T$.

\begin{example}\label{exampleNonCoherent}
Let us consider the non coherent linear involution $T$ which is the same as in Example~\ref{exampleLinear} but such that the restriction of $\sigma_1$ to $I_c$ 
is a symmetry, as in Example~\ref{exampleConnexion}. We assume that $I=]0,1[$,
that $\lambda_a=\lambda _d$, that $1/4<\lambda_c<1/2$
 and that $\lambda_a+\lambda_b<1/2$.
 Let $z=1/2+\lambda_c$ (see Figure~\ref{figureLinearNonCoherent}).
We have then $T^3(z)=z$, showing that $T$ is not minimal. Indeed,
since $z\in I_c$, we have
$T(z)=1-z=1/2-\lambda_c$. Since
$T(z)\in I_a$ we have $T^2(z)=(\lambda_a+\lambda_b)+(\lambda_a-1+z)=z-\lambda_c=1/2$.
Finally, since $T^2(z)\in I_{d^{-1}}$, we obtain $1-T^3(z)=T^2(z)-\lambda_c=1-z$
and thus $T^3(z)=z$.

\begin{figure}[hbt]
\centering
\tikzset{node/.style={circle,draw,minimum size=0.1cm,inner sep=0pt}}
\tikzset{title/.style={circle,minimum size=0.1cm,inner sep=0pt}}
\begin{tikzpicture}
%--interval nodes--
\node[node](h0)at(0,1){};\node[node](h1)at(2,1){};\node[node](h2)at(5,1){};
\node[node](h3)at(7,1){};\node[node](h4)at(11,1){};\node[title](0)at(12,1){$I\times\{0\}$};
\node[node](b0)at(0,0){};\node[node](b1)at(4,0){};
\node[node](b2)at(6,0){};\node[node](b3)at(9,0){};
\node[node](b4)at(11,0){};\node[title](1)at(12,0){$I\times\{1\}$};
%translation nodes
\node[title](z)at(9.5,.8){$z=T^3z$};
\node[title](Tz)at(1.5,.8){$Tz$};\node[title](T2z)at(5.5,-.2){$T^2z$};

%interval lines
\draw[above,line width=2pt,color=red](h0)edge node{$a$}(h1);
\draw[above,line width=2pt,color=blue](h1)edge node{$b$}(h2);
\draw[above,line width=2pt,color=red](h2)edge node{$a^{-1}$}(h3);
\draw[above,line width=2pt,color=green](h3)edge node{$c$}(h4);
\draw[below,line width=2pt,color=green](b0)edge node{$c^{-1}$}(b1);
\draw[below,line width=2pt,color=violet](b1)edge node{$d^{-1}$}(b2);
\draw[below,line width=2pt,color=blue](b2)edge node{$b^{-1}$}(b3);
\draw[below,line width=2pt,color=violet](b3)edge node{$d$}(b4);

%translation lines
\draw[bend left,->](1.5,1)edge node{}(5.5,1);\draw[->](5.5,1)edge node{}(5.5,0);
\draw[->](5.5,0)edge node{}(9.5,0);\draw[->](9.5,0)edge node{}(9.5,1);
\draw[->](9.5,1)edge node{}(1.5,0);\draw[bend right,->](5.5,0)edge node{}(9.5,0);
\draw[->](1.5,0)edge node{}(1.5,1);
\end{tikzpicture}
\caption{A non coherent linear involution.}\label{figureLinearNonCoherent}
\end{figure}
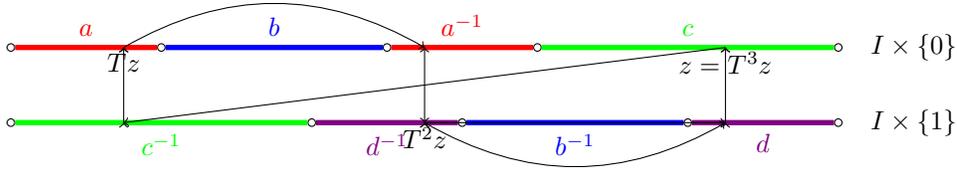
\end{example}

We quote without proof the following result, analogous
to Keane Theorem (Theorem~\ref{theoremIDOC}) for interval exchange transformations. 
\begin{proposition}\label{propBL}
Let $T$ be a linear involution without connections on $I$. If $T$
is nonorientable, it is minimal. Otherwise, its restriction to each component
of $\hat{I}$ is minimal. 
\end{proposition}

\subsection{Natural coding}
Let $T$ be a linear involution on $I$, let $\hat{I}=I\times\{0,1\}$ and let
$\hat{O}$ be the set defined by Equation~\eqref{eqO}.

Given $z\in \hat{I}\setminus \hat{O}$, the \emph{infinite natural
  coding}\index{subject}{infinite natural coding} of $T$ relative to $z$ is the infinite word
$\Sigma_T(z)=(a_n)_{n\in\Z}$
on the alphabet $A$ defined by
\begin{displaymath}
a_n=a\quad\text{ if }\quad T^n(z)\in I_a.
\end{displaymath}
We first observe that the factors of $\Sigma_T(z)$ are
reduced. Indeed, assume that $a_n=a$ and $a_{n+1}=a^{-1}$ with $a\in
A$. Set $x=T^n(z)$ and $y=T(x)=T^{n+1}(z)$. Then $x\in I_a$
and $y\in I_{a^{-1}}$. But $y=\sigma_2(u)$ with $u=\sigma_1(x)$. Since $x\in
I_a$, we have $u\in I_{a^{-1}}$. This implies that
$y=\sigma_2(u)$ and $u$ belong to the same component of $\hat{I}$, a contradiction.

We denote by $X(T)$ the set of  infinite natural codings
of $T$. We say that $X(T)$ is 
the \emph{natural coding} of $T$. We also denote $\cL(T)=\cL(X(T))$.
%It is also called the \emph{lamination set} of $T$ in \cite{HilionCoulboisLustig2008}.

\begin{example}
Let $T$ be the linear involution of Example~\ref{exampleInvolution3}.
The words of length at most $3$ of $S=\cL(T)$ are represented 
in Figure~\ref{figureSetS}.
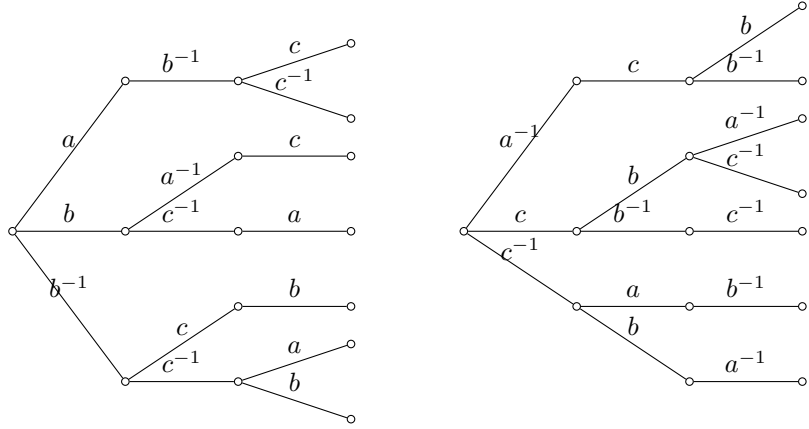
\begin{figure}[hbt]
\centering
\tikzset{node/.style={circle,draw,minimum size=0.1cm,inner sep=0pt}}
\tikzset{title/.style={circle,minimum size=0.1cm,inner sep=0pt}}
\begin{tikzpicture}

%gauche
\node[node](1)at(0,2){};
\node[node](a)at(1.5,4){};
\node[node](b)at(1.5,2){};\node[node](bbar)at(1.5,0){};
\node[node](abbar)at(3,4){};
\node[node](babar)at(3,3){};\node[node](bcbar)at(3,2){};
\node[node](bbarc)at(3,1){};\node[node](bbarcbar)at(3,0){};
\node[node](abbarc)at(4.5,4.5){};\node[node](abbarcbar)at(4.5,3.5){};
\node[node](babarc)at(4.5,3){};\node[node](bcbara)at(4.5,2){};
\node[node](bbarcb)at(4.5,1){};
\node[node](bbarcbara)at(4.5,.5){};\node[node](bbarcbarb)at(4.5,-.5){};

\draw[above](1)edge node{$a$}(a);\draw[above](1)edge node{$b$}(b);\draw[above](1)edge node{$b^{-1}$}(bbar);
\draw[above](a)edge node{$b^{-1}$}(abbar);
\draw[above](b)edge node{$c^{-1}$}(bcbar);\draw[above](b)edge node{$a^{-1}$}(babar);
\draw[above](bbar)edge node{$c$}(bbarc);\draw[above](bbar)edge node{$c^{-1}$}(bbarcbar);
\draw[above](abbar)edge node{$c$}(abbarc);\draw[above](abbar)edge node{$c^{-1}$}(abbarcbar);
\draw[above](babar)edge node{$c$}(babarc);\draw[above](bcbar)edge node{$a$}(bcbara);
\draw[above](bbarc)edge node{$b$}(bbarcb);
\draw[above](bbarcbar)edge node{$a$}(bbarcbara);\draw[above](bbarcbar)edge node{$b$}(bbarcbarb);

%droit
\node[node](1)at(6,2){};
\node[node](abar)at(7.5,4){};\node[node](c)at(7.5,2){};\node[node](cbar)at(7.5,1){};
\node[node](abarc)at(9,4){};
\node[node](cb)at(9,3){};\node[node](cbbar)at(9,2){};
\node[node](cbara)at(9,1){};\node[node](cbarb)at(9,0){};
\node[node](abarcb)at(10.5,5){};\node[node](abarcbbar)at(10.5,4){};
\node[node](cbabar)at(10.5,3.5){};\node[node](cbcbar)at(10.5,2.5){};
\node[node](cbbarcbar)at(10.5,2){};
\node[node](cbarabbar)at(10.5,1){};\node[node](cbarbabar)at(10.5,0){};

\draw[above](1)edge node{$a^{-1}$}(abar);\draw[above](1)edge node{$c$}(c);\draw[above](1)edge node{$c^{-1}$}(cbar);
\draw[above](abar)edge node{$c$}(abarc);
\draw[above](c)edge node{$b$}(cb);\draw[above](c)edge node{$b^{-1}$}(cbbar);
\draw[above](cbar)edge node{$a$}(cbara);\draw[above](cbar)edge node{$b$}(cbarb);
\draw[above](abarc)edge node{$b$}(abarcb);\draw[above](abarc)edge node{$b^{-1}$}(abarcbbar);
\draw[above](cb)edge node{$a^{-1}$}(cbabar);\draw[above](cb)edge node{$c^{-1}$}(cbcbar);
\draw[above](cbbar)edge node{$c^{-1}$}(cbbarcbar);
\draw[above](cbara)edge node{$b^{-1}$}(cbarabbar);\draw[above](cbarb)edge node{$a^{-1}$}(cbarbabar);

\end{tikzpicture}
\caption{The words of length at most $3$ of $\cL(X)$.}\label{figureSetS}
\end{figure}

The set $S$ can  actually be defined directly as the set of factors
of the substitution
\begin{displaymath}
f:a\mapsto cb^{-1},\quad b\mapsto c,\quad c\mapsto ab^{-1}.
\end{displaymath}
which extends to an automorphism of the free group on $\{a,b,c\}$.
Indeed, the application twice of Rauzy induction on $T$
gives a linear involution which is the same as $T$ (with the
two copies of $[0,1]$ interchanged). This gives the
explanation of
 why the substitution shift of Example~\ref{exampleInvolution3bis}
is specular.
%(see~\cite{BertheDelecroixDolcePerrinReutenauerRindone2014}).

\end{example}

%The following is Proposition 5.3 in~\cite{BertheDelecroixDolcePerrinReutenauerRindone2014}.
Define, as for an interval exchange transformation, for a word
$w=a_0a_1\cdots a_n$ with $a_i\in A\cup A^{-1}$,
\begin{displaymath}
I_w=I_{a_0}\cap T^{-1}(I_{a_1})\cap\ldots\cap T^{-n}(I_{a_n}).
\end{displaymath}
It is clear that
\begin{equation}
w\in \cL(T)\Leftrightarrow I_w\ne\emptyset.\label{equationI_w}
\end{equation}
\begin{proposition} \label{prop:inverse}
Let $T$ be a linear involution.
The set $\cL(T)$ is a laminary set.
\end{proposition}
\begin{proof}
Set $T=\sigma_2\circ\sigma_1$.
 We claim that for any nonempty word  $u\in {\mathcal L } (T)$, one has $I_{u^{-1}}=\sigma_1T^{|u|-1}(I_u)$.

To prove the  claim, we use an induction on the
length of $u$. The property holds for $|u|=1$ by definition of $\sigma_1$.
 Next, consider
 $u\in {\mathcal L } (T)$ and $a\in A\cup A^{-1}$ such that $ua\in {\mathcal L } (T)$.

Since $T^{-1}=\sigma_1\circ\sigma_2$, we have, using the induction hypothesis,
\begin{eqnarray*}
\sigma_1T^{|u|}(I_{ua})&=&\sigma_1T^{|u|}(I_u\cap T^{-|u|}(I_a))
=\sigma_1T^{|u|}(I_u)\cap \sigma_1(I_a)\\
&=&\sigma_1\sigma_2\sigma_1T^{|u|-1}(I_u)\cap \sigma_1(I_a)
=\sigma_1\sigma_2(I_{u^{-1}})\cap I_{a^{-1}}=I_{a^{-1}u^{-1}}
\end{eqnarray*}
where the last equality results from
$I_{a^{-1}u^{-1}}= T^{-1} I_{u^{-1}}  \cap  I_{a^{-1}}  $.

We  easily  deduce from the claim that   the set ${\mathcal L } (T)$ is closed under taking inverses. 
 Furthermore it  is a   factorial subset of the group $G_\theta$. It  is  thus a laminary set.
\end{proof}

We prove the following result.
\begin{theorem}\label{theoremInvolutionSpecular}
The natural coding of a linear involution without connections
 is a 
 specular
shift.
\end{theorem}

We first prove the following lemma, which replaces for linear
involutions Conditions \eqref{eqConditionR(w)}
and \eqref{eqConditionL(w)}.
\begin{lemma}\label{lemma9.6} Let $T$ be a linear involution. For every nonempty word $w$ and letter $a \in A$, one has
\begin{enumerate}
\item[\rm(i)] $a \in L ( w ) \Leftrightarrow \sigma_2 ( I_{a^{-1} }) \cap I_w\ne \emptyset$ ,
\item[\rm(ii)] $a \in R ( w )\Leftrightarrow  \sigma_2 ( I_a)\cap I_{w^{-1}}\ne \emptyset$.
\end{enumerate}
\end{lemma}
\begin{proof} By \eqref{equationI_w}, we have $a \in L ( w )$
 if and only if $I_{aw}\ne\emptyset$,
 which is also equivalent to $T ( I_{aw} )\ne\emptyset$ . 
By definition of $I_{aw}$, we have $T ( I_{aw} ) = T ( I_a ) \cap I_w$. 
Since $T = \sigma_2 \circ \sigma_1$ and since $\sigma_1 ( I_a ) = I_{a^{-1}}$ , 
we have $a\in L ( w )$ if and only if $\sigma_2 ( I_{a^{-1}} ) \cap I_w \ne\emptyset$ .
Next, since $\cL ( T )$ is closed under taking inverses by Proposition~\ref{prop:inverse}, we have $aw \in \cL(T)$ if and only if $w^{-1} a^{-1} \in \cL(T)$. 
Thus $a \in R ( w )$
if and only if $a^{-1} \in L ( w^{-1} )$ , whence the second equivalence
\end{proof}
Given a linear involution $T$ on $I$ , we introduce two orders on $\cL ( T )$
 as follows. For any $u , v \in \cL ( T )$ , one has
\begin{enumerate}
\item[(i)] $u <_R v$ if and only if $I_u < I_v$ ,
\item[(ii)] $u <_L v$ if and only if $I_{u^{-1}} < I_{v^{-1}}$ .
\end{enumerate}
\begin{lemma}\label{lemma9.7} Let $T$ be a linear involutions on $I$ without connexion. Let $w \in \cL ( T )$ and $a , a' \in L ( w )$
 (resp. $b, b' \in R ( w )$). Then $1\otimes a$,
$1\otimes a'$ (resp. $b\otimes 1$ , $b'\otimes 1$) are in the same connected component of $\E ( w )$ if and only if $I_{a^{-1}}$ , $I_{a'^{-1}}$ (resp. $I_b$,
$ I_{b'}$ ) are in the same
component of $\hat{I}$ .
\end{lemma}
\begin{proof}
 If $( 1\otimes a , b\otimes 1 ) \in \E ( w $) , then $\sigma_2 ( I_{a^{-1}} )\cap I wb \ne\emptyset$ . Thus $I_{a^{-1}}$ and $I_{wb}$
 belong to distinct components of $\hat{I}$ . Consequently,
if $ a , a' \in L ( w )$ (resp. $R ( w )$) belong to the same connected component of $\E ( w )$, then $I_{a^{-1}} , I_{a'^{-1}}$ (resp. $I_{wa}$ , $I_{wa'}$ ) belong to
the same component of $\hat{I}$ .
Conversely, let $a , a' \in L ( w )$ be such that $a , a'$ belong to the same component of $\hat{I}$ . We may assume that $a <_L a'$. There
is a reduced path (i.e., it does not use twice consecutively the same edge)
 in $\E ( w )$ from $a$ to $a'$ which is the sequence
$a_1 , b_1 , . . . , b_{n-1} , a_n$ with $a_1 = a$ and $a_n = a'$
 with $a_1 <_L a_2 <_L · · · <_L a_n$ , $wb_1 <_R wb_2 <_R · · · <_R wb_{n-1}$ and
$\sigma_2 ( I_{a_i^{-1}} )\cap I_{wb_i}\ne\emptyset$ ,
$\sigma_2 ( I_{a_{i+1}^{-1}} )\cap I_{wb_i} \ne\emptyset$ for $1 \le i \le n -1$
(see Figure~\ref{figurePathInvolution} for an illustration).
\begin{figure}[hbt]
\centering
\tikzset{node/.style={circle,draw,minimum size=0.1cm,inner sep=0pt}}
\begin{tikzpicture}
\node[node](H1)at(0,1){};
\node[node](H2)at(2,1){};\node[node](H3)at(3,1){};\node[node](H4)at(5,1){};
\node[node](H5)at(8,1){};\node[node](H6)at(10,1){};

\node[node](B1)at(-1,0){};\node[node](B2)at(.5,0){};\node[node](B3)at(1.5,0){};
\node[node](B4)at(3.5,0){};
\node[node](B5)at(6.5,0){};\node[node](B6)at(8.5,0){};
\node[node](B7)at(9.5,0){};
\node[node](B8)at(11,0){};

\draw[above](H1)edge node{$wb_1$}(H2);
\draw[above](H3)edge node{$wb_2$}(H4);
\draw[above](H5)edge node{$wb_{n-1}$}(H6);

\draw[below](B1)edge node{$a_1^{-1}$}(B2);
\draw[below](B3)edge node{$a_2^{-1}$}(B4);
\draw[below](B5)edge node{$a_{n-1}^{-1}$}(B6);
\draw[below](B7)edge node{$a_n^{-1}$}(B8);

\draw[dotted](B2)edge node{}(.5,1){};\draw[dotted](2,0)edge node{}(H2){};
\draw[dotted]{}(B4)edge node{}(3.5,1);
\draw[dotted](8,0)edge node{}(H5);\draw[dotted](10,0)edge node{}(H6);
\end{tikzpicture}
\caption{A path from $a_1$ to $a_n$ in $\E(w)$.}\label{figurePathInvolution}
\end{figure}
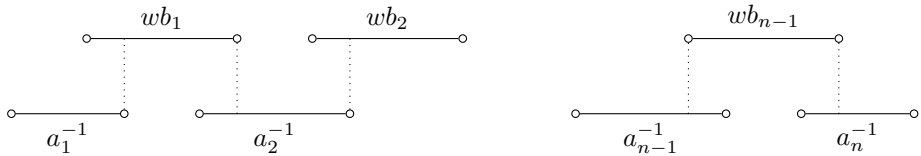

Note that the hypothesis that $T$ is without connection is needed since otherwise the right boundary of $\sigma_2 ( I_{a^{-1}} )$ could
be the left boundary of $I_{wb_i}$.
The assertion concerning $b , b' \in R ( w )$ is a consequence 
of the first one since $b , b' \in R ( w )$ if and only if 
$b^{-1} , b'^{-1} \in
L ( w^{-1} $). 
\end{proof}

\begin{proofof}{of Theorem~\ref{theoremInvolutionSpecular}}
Let $T$ be a linear involution without connections. By Proposition~\ref{prop:inverse}, the set $\cL(T)$ is symmetric. Since it is by definition extendable
and formed of reduced words,
it is a laminary set.  

There remains to show that $X(T)$ is dendric of characteristic $2$.
Let us first prove that
for any $w \in \cL ( T )$, the graph $\E ( w )$ is acyclic. Assume that 
$( 1\otimes a_1 , b_1\otimes 1 ,\ldots, 1\otimes a_n , b_n\otimes 1 )$
 is a path in $\E ( w )$ with
$(a_1 , . . . , a_n) \in L ( w )$
 and $b_1 , . . . , b_n \in R ( w )$. We may assume that the path is reduced, that $n \ge 2$ and also that $a_1 <_L a_2$. It
follows that $a_1 <_L \ldots <_L a_n$ and 
$wb_1 <_R \ldots <_R wb_n$ (see Figure~\ref{figurePathInvolution}). 
Thus it is not possible to have an edge $( a_1 , b_n )$, which
shows that $\E ( w )$ is acyclic.

Let $a , a' \in A$. If $I_{a^{-1}}$ and $I_{a'^{-1}}$ are in the same component of
 $\hat{I}$, then $1\otimes a$ and $1\otimes a'$ are in the same connected component
of $\E ( \varepsilon)$. Thus $\E ( \varepsilon )$ is a union of two trees.

Next, if $w \in \cL(T)$ is nonempty and $1\otimes a , 1\otimes a' \in L ( w )$,
 then $I_{a^{-1}}$ and $I_{a'^{-1}}$ are in the same component of $\hat{I}$ 
(by Lemma\ref{lemma9.6}), and
thus $1\otimes a , 1\otimes a'$ are in the same connected component 
of $\E ( w )$. Thus $\E ( w )$ is a tree.
\end{proofof}
We now present an example of a linear involution on an alphabet $A$
where the involution $\theta$ has fixed points.
\begin{example}\label{exampleFiboDoubleInvolution}
Let $A=\{a,b,c,d\}$ be as in Example~\ref{exampleSpecularGroup}
(in particular, $d=b^{-1}$, $a=a^{-1}$, $c=c^{-1}$).
\begin{figure}[hbt]
\centering
\tikzset{node/.style={circle,draw,minimum size=0.1cm,inner sep=0.4pt}}
\tikzset{title/.style={circle,minimum size=0.1cm,inner sep=0pt}}
\begin{tikzpicture}(100,15)
\node[node](h0)at(0,1){};\node[node](bbar)at(6.18,1){};
\node[node](h1)at(10,1){};
\node[node](b0)at(0,0){};\node[node](c)at(3.82,0){};
\node[node](b1)at(10,0){};

\draw[above,color=red,line width=2pt](h0)edge node{$a$}(bbar);
\draw[above,color=yellow,line width=2pt](bbar)edge node{$d$}(h1);
\draw[below,color=blue,line width=2pt](b0)edge node{$b$}(c);
\draw[below,color=green,line width=2pt](c)edge node{$c$}(b1);
\end{tikzpicture}
\caption{A linear involution on $A=\{a,b,c,d\}$.}\label{figureLinearFiboDouble}
\end{figure}
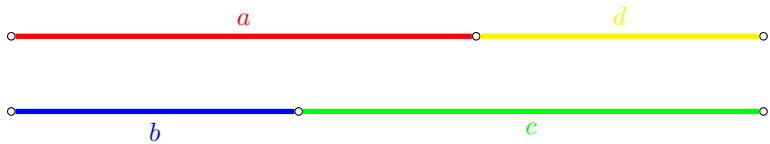
Let $T$ be the linear involution represented in Figure~\ref{figureLinearFiboDouble} with $\sigma_1$ being a translation on $I_b$ and a symmetry on $I_a,I_c$.
Choosing $(3-\sqrt{5})/2$ for the length of $I_b$, the involution is without connections. Thus $\cL(T)$ is a specular set.

Note that the natural coding of the linear
involution $T$ is equal to the
set  of factors of the shift of Example~\ref{exampleFiboDouble}. Indeed,
consider the interval exchange $V$ on the interval
$Y=]0,1[$ represented
in Figure~\ref{figureFibonacciDouble} on the right,
which is obtained by using two copies of the 
interval exchange $U$ defining the Fibonacci set (represented in Figure~\ref{figureFibonacciDouble}
on the left).
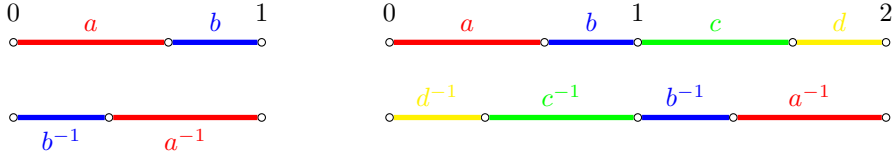
\begin{figure}[hbt]
\centering
\tikzset{node/.style={circle,draw,minimum size=0.1cm,inner sep=0.4pt}}
\tikzset{title/.style={circle,minimum size=0.1cm,inner sep=0pt}}
\begin{tikzpicture}
%gauche
\node[node](h0)at(0,1){};\node[title]at(0,1.4){$0$};
\node[node](b)at(2.06,1){};
\node[node](h1)at(3.3,1){};\node[title]at(3.3,1.4){$1$};
\node[node](b0)at(0,0){};
\node[node](a)at(1.27,0){};\node[node](b1)at(3.3,0){};

\draw[above,line width=2pt,color=red](h0)edge node{$a$}(b);
\draw[above,line width=2pt,color=blue](b)edge node{$b$}(h1);
\draw[below,line width=2pt,color=blue](b0)edge node{$b^{-1}$}(a);
\draw[below,line width=2pt,color=red](a)edge node{$a^{-1}$}(b1);

%droit
\node[node](h0)at(5,1){};\node[title]at(5,1.4){$0$};
\node[node](b)at(7.06,1){};\node[node](h1)at(8.3,1){};\node[title]at(8.3,1.4){$1$};
\node[node](d)at(10.36,1){};\node[node](h2)at(11.6,1){};
\node[title]at(11.6,1.4){$2$};

\node[node](b0)at(5,0){};\node[node](cbar)at(6.27,0){};\node[node](b1)at(8.3,0){};
\node[node](abar)at(9.57,0){};\node[node](b2)at(11.6,0){};

\draw[above,line width=2pt,color=red](h0)edge node{$a$}(b);
\draw[above,line width=2pt,color=blue](b)edge node{$b$}(h1);
\draw[above,line width=2pt,color=green](h1)edge node{$c$}(d);
\draw[above,line width=2pt,color=yellow](d)edge node{$d$}(h2);

\draw[above,line width=2pt,color=yellow](b0)edge node{$d^{-1}$}(cbar);
\draw[above,line width=2pt,color=green](cbar)edge node{$c^{-1}$}(b1);
\draw[above,line width=2pt,color=blue](b1)edge node{$b^{-1}$}(abar);
\draw[above,line width=2pt,color=red](abar)edge node{$a^{-1}$}(b2);

\end{tikzpicture}
\caption{Interval exchanges $U$ and $V$ for the Fibonacci set and its doubling.}
\label{figureFibonacciDouble}
\end{figure}
Let $X=]0,1[\times\{0,1\}$ and let $\alpha:Y\rightarrow X$ be the
map defined by
\begin{displaymath}
\alpha(z)=\begin{cases}(z,0)&\text{if $z\in]0,1[$}\\(2-z,1)&\text{otherwise.}
\end{cases}
\end{displaymath}
Then $\alpha\circ V=T\circ\alpha$ and thus $\cL(V)=\cL(T)$.
%The interval exchange $U$ is actually the orientation covering
%of the linear involution $T$.
%(see~\cite{BertheDelecroixDolcePerrinReutenauerRindone2014}).
\end{example}

%%%%%%%%%%%%%%%%%%%%%%%%%%%
\section{Exercises}
\exosection{Section~\ref{ch5:subsection:iet}}
\begin{exercise}\label{exerciseTribonacciNotIET}
Show that an Arnoux-Rauzy shift is not an interval exchange transformation.
  \end{exercise}
\begin{exercise}\label{exerciseInvariantMeasureIET}
  Let $T=T_{\lambda,\pi}$ be a $k$-interval exchange
  transformation on $I=[0,r)$ exchanging
    the intervals $\Delta_1,\ldots,\Delta_k$. Let $\mu$
    be an invariant  measure on $(I,T)$ and set $\mu_i=\mu(\Delta_i)$ for $1\le i\le k$.
    Let
    $\varphi:I\to I$ be the map defined by
    \begin{displaymath}
      \varphi(x)=\mu([0,x)).
    \end{displaymath}
    Show that if $T$ is minimal, then $\varphi$ is a continuous isomorphism
    of dynamical systems from $(I,T)$ onto
    $(I,T_{\mu,\pi})$ where $T_{\mu,\pi}$ is the
    interval exchange transformation   on $[0,\mu(I))$ defined
    by the vector $\mu=(\mu_1,\ldots,\mu_k)$ and the permutation $\pi$.
  \end{exercise}

\begin{exercise}\label{exerciseMatriceAntisym}
Let $T=T_{\lambda,\pi}$ be as above.
Set $Tx=x+\alpha_i$ for $x\in\Delta_i$. Show that
there is an antisymmetric  matrix $M_\pi$ such that the column vector
$\alpha$ with coordinates $\alpha_i$ is defined by
\begin{equation}
\alpha=M_\pi\lambda.\label{eqAlphaLambda}
\end{equation}
Compute the matrix $M_\pi$ and $M_{\pi'}$ for $\pi=(132)$ and $\pi'=(13)$.
\end{exercise}

\begin{exercise}\label{exerciseOrthogonal}
Let $M_\pi$ be the antisymmetric matrix of Exercise~\ref{exerciseMatriceAntisym}.
Show that if $\mu,\nu$ are invariant measures on $(I,T)$, then
\begin{displaymath}
  \mu^t M_\pi\nu=0
\end{displaymath}
where $\mu,\nu$ denote the vectors
$(\mu(\Delta_i)), (\nu(\Delta_i))$.
\end{exercise}
\begin{exercise}\label{exerciseLemma2.11}
  Let $T=T_{\lambda,\pi}$ be minimal and such that $\pi$ has the property
  that $\pi(j)\ne\pi(j+1)$ for all $j$ or, equivalently,
  that $T_{\lambda,\pi}$ has exactly $k-1$ points of discontinuity.
  Let $V_{\lambda,\pi}$ be the vector space generated by the cone $\I(I,T)$
  \index{symbols}{I(I,T)@$\I(I,T)$}%
  of invariant  measures on $(I,T)$.
  Show that $V_{\lambda,\pi}\cap\ker M_\pi=\{0\}$.
\end{exercise}
\begin{exercise}\label{exerciseTheorem2.12}
  Show that if $T=T_{\lambda,\mu}$ is minimal, then $V_{\lambda,\mu}$
  has dimension at most $1/2\rank(M_\pi)$
  and consequently that $(I,T)$ admits at most
  $1/2\rank(M_\pi)$ ergodic measures.\index{subject}{interval exchange!ergodic measures}
  \end{exercise}
\begin{exercise}\label{exerciseIw}
  Define for every word $w\in A^*$ a subset $I(w)$ of $I$ as follows.
Set $I(\varepsilon)=I$   and $I(au)=\Delta_a\cap T^{-1}(I(u))$
for $u\in A^*$ and $a\in A$.
Show that every nonempty $I(w)$ is a semi-interval.
\end{exercise}
\begin{exercise}\label{exerciseJw}
Let $J(w)$ be defined by $J(\varepsilon)=I$
 and by
\begin{equation}
J(ua)=TJ(u)\cap T\Delta_a\label{eqJua}
\end{equation}
for $a\in A$ and $u\in A^*$. Show that the nonempty sets $J(w)$
are semi-intervals. 
\end{exercise}
%%%%%%%%%%%%%%%%%%%%%%%
\section{Solutions}
\exosection{Section~\ref{ch5:subsection:iet}}
\begin{solution}{\ref{exerciseTribonacciNotIET}}
  Let $(T,I)$ be an interval exchange transformation. Suppose that
  its natural coding $X=X(T)$ is an Arnoux-Rauzy shift. Since
  $X(T)$ is minimal, $T$ is also minimal and thus the
  natural coding $\gamma:I\to X$ is injective. For $w\in\cL(X)$
  long enough, the length of the interval $I(w)$ is arbitrary small
  and thus $\Card(R(w))\le 2$. This is implies that $\Card(A)=2$and that $X$
  is a Sturmian shift.
  \end{solution}
\begin{solution}{\ref{exerciseInvariantMeasureIET}}
  If $\varphi$ is not continuous, there is some $x\in[0,1)$
    such that $\mu(\{x\})>0$. This is not possible since $T$
    is minimal. For the same reason, $\varphi$ is strictly increasing
    and thus bijective. Thus $\varphi$ is a homeomorphism.

    Let us show that $\varphi$ is a morphism of dynamical systems.
    For $1\le i\le k$, set
    \begin{displaymath}
      \Delta_i=[\beta_{i-1},\beta_i), \quad T\Delta_{i}=[\gamma_{\pi(i)-1},\gamma_{\pi(i)})
    \end{displaymath}
    and $\delta_i=\varphi(\beta_i),\varepsilon_j=\varphi(\gamma_j)$.
    \begin{figure}[hbt]
      \centering
      \tikzset{node/.style={circle,draw,minimum size=0.1cm,inner sep=0.4pt}}
\tikzset{title/.style={circle,minimum size=0.1cm,inner sep=0pt}}
\begin{tikzpicture}
  \node[node](betai-1)at(.5,1){};\node[title](betai-1n)at(.5,1.5){$\beta_{i-1}$};
  \node[node](betai)at(2.5,1){};\node[title](betain)at(2.5,1.5){$\beta_{i}$};
\node[node](gammapi-1)at(0,0){};\node[title](gammapi-1n)at(0,-0.5){$\gamma_{\pi(i)-1}$};
\node[node](gammapi)at(2,0){};\node[title](gammapin)at(2,-.5){$\gamma_{\pi(i)}$};

\draw[above](betai-1)-- node{$\Delta_i$}(betai);
\draw[left,->](betai-1)-- node{$T$}(gammapi-1);
\draw[left,->](betai)-- node{$T$}(gammapi);
\draw[above](gammapi-1)-- node{$T\Delta_i$}(gammapi);

% partie Tmu
  \node[node](deltai-1)at(4.5,1){};\node[title](deltai-1n)at(4.5,1.5){$\delta_{i-1}$};
  \node[node](deltai)at(7.5,1){};\node[title](deltain)at(7.5,1.5){$\delta_{i}$};
\node[node](epsilonpi-1)at(4,0){};\node[title](epsilonpi-1n)at(4,-0.5){$\varepsilon_{\pi(i)-1}$};
\node[node](epsilonpi)at(7,0){};\node[title](epsilonpin)at(7,-.5){$\varepsilon_{\pi(i)}$};

\draw[above](deltai-1)-- node{$\varphi(\Delta_i)$}(deltai);
\draw[left,->](deltai-1)-- node{$T'$}(epsilonpi-1);
\draw[left,->](deltai)-- node{$T'$}(epsilonpi);
\draw[above](epsilonpi-1)-- node{$\varphi(T\Delta_i)$}(epsilonpi);
      \end{tikzpicture}
      \caption{The action of $T$ and $T'$.}\label{figureTmu}
    \end{figure}
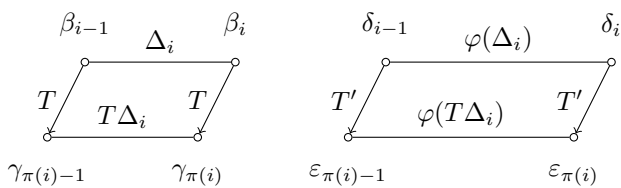
    
    Let $T'=\varphi\circ T\circ \varphi^{-1}$.
    Since $\varphi$ is increasing, we have $T'(\varphi(\Delta_i))\subset [\varepsilon_{\pi (i)-1},\varepsilon_{\pi(i)})$ (see Figure~\ref{figureTmu}). Since
      $\mu$ is invariant, $T'$ is a translation from $\varphi(\Delta_i)$
      onto $[\varepsilon_{\pi(i)-1},\varepsilon_{\pi(i)})$. Thus $T'=T_{\mu,\pi}$.

  \end{solution}

\begin{solution}{\ref{exerciseMatriceAntisym}}
We have 
\begin{equation}
\alpha_i=\sum_{\pi(j)<\pi(i)}\lambda_{j}-\sum_{j<i}\lambda_j.
\label{eqAlpha_i}
\end{equation}
Indeed, the first term of the right hand side is the
left boundary of $T\Delta_i$ and the second one is the left boundary
of $\Delta_i$. A pair $(i,j)$ of integers with $1\le i<j\le k$
is an \emph{inversion}\index{subject}{inversion of permutation}
\index{subject}{permutation!inversion of}%
of a permutation $\sigma\in\Sg_k$
if $\sigma(i)>\sigma(j)$. Setting
\begin{displaymath}
M_{ij}=\begin{cases}1&\mbox{ if $(i,j)$ is an inversion of $\pi$}\\
-1&\mbox{ if $(j,i)$ is an inversion of $\pi$}\\
0&\mbox{ otherwise}
\end{cases}
\end{displaymath}
we obtain the desired symmetric  matrix.

For $\pi=(123)$ and $\pi'=(13)$ the matrices are respectively 
\begin{displaymath}
  M_\pi=\begin{bmatrix}0&0&1\\0&0&1\\-1&-1&0\end{bmatrix},\quad
  M_{\pi'}=\begin{bmatrix}0&1&1\\-1&0&1\\-1&-1&0\end{bmatrix}.
  \end{displaymath}
\end{solution}
\begin{solution}{\ref{exerciseOrthogonal}}
  Let $\I(I,T)$ denote the family of invariant  measures
  on $(I,T)$.
  We have, by Exercise~\ref{exerciseInvariantMeasureIET}
  \begin{displaymath}
    \mu\in\I(I,T_{\nu,\pi})\Leftrightarrow \nu\in\I(I,T_{\mu,\pi}).
  \end{displaymath}
  Thus, since $M_\pi$ depends only on $\pi$,
  we only need to prove that, for $T=T_{\lambda,\pi}$
  and $\mu\in\I(I,T)$, one has $\mu^tM_\pi\lambda=0$.

  We have indeed, since $\mu$ is invariant,
  \begin{eqnarray*}
    0&=&\int(Tx-x)d\mu\\
    &=&\sum_{i=1}^k(M_\pi\lambda)_i\mu(\Delta_i)=\mu^t M_\pi\lambda.
    \end{eqnarray*}

\end{solution}

\begin{solution}{\ref{exerciseLemma2.11}}
  Let $\lambda',\lambda''$ be invariant measures on $(I,T)$
  and suppose that $\lambda'-\lambda''\in \ker M_\pi$. Since
  the interval exchange transformations $T_{\lambda',\pi}$ and
  $T_{\lambda'',\pi}$ are continuously isomorphic, we have for all $n\ge 1$
  and $1\le i\le k$
  \begin{equation}
\charac_{\Delta'_i}^{(n)}(0)=\charac_{\Delta''_i}^{(n)}(0)\label{eqM(n,pi,0)}
  \end{equation}
  where $\Delta'_i,\Delta''_i$ are the intervals exchanged
  by $T_{\lambda',\pi}$ and $T_{\lambda'',\pi}$ respectively.
  Indeed, set $\lambda_t=t\lambda'+(1-t)\lambda''$ for $0\le t\le 1$
  and let $\Delta_{t,i}=[\beta_{t,i-1},\beta_{t,i}]$ be the intervals
  exchanged by $T_{\lambda_t,\pi}$.
  For all $n,i$ either $T^n_{\lambda_t,\pi}0>\beta_{t,i-1}$ for all $t$,
  or $T^n_{\lambda_t,\pi}0=\beta_{t,i-1}$ for all $t$, or
  $T^n_{\lambda_t,\pi}0<\beta_{t,i-1}$ for all $t$.

  Now, using \eqref{eqM(n,pi,0)} and the fact that
  $M_\pi\lambda'=M_\pi\lambda''$, we obtain
  \begin{eqnarray*}
    T^n_{\lambda',\pi}0&=&\sum_{i=0}^k\charac_{\Delta'_i}^{(n)}(0)M_\pi\lambda'\\
    &=&\sum_{i=0}^k\charac_{\Delta''_i}^{(n)}(0)M_\pi\lambda''\\
    &=&T^n_{\lambda'',\pi}0
  \end{eqnarray*}
  Since $T_{\lambda',\pi}$ and $T_{\lambda'',\pi}$ are minimal, we conclude
  that $T_{\lambda',\pi}=T_{\lambda'',\pi}$ and thus that $\lambda'=\lambda''$.
\end{solution}
\begin{solution}{\ref{exerciseTheorem2.12}}
  We may assume that $T_{\lambda,\pi}$ has $k-1$ points of discontinuity
  (otherwise, we can reduce $k$ without changing
  the rank of $M_\pi$). Let $W=\R^k/(V_{\lambda,\pi}+\ker(M_\pi))$
  and let $\alpha:V_{\lambda,_pi}\to W^*$ be the linear map
  defined for $v\in V_{\lambda,\pi}$ and $u\in\R^k$ by
  \begin{displaymath}
    (\alpha u)(v)=u^tM_\pi v.
  \end{displaymath}
  Exercise \ref{exerciseLemma2.11} shows that this map is well defined
  and that it is injective. Since $V_{\lambda,\pi}\cap \ker (M_\pi)=\{0\}$,
  we have
  \begin{equation}
    \dim W=k-\dim V_{\lambda,\pi}-\dim\ker(M_\pi)=\rank M_\pi-\dim V_{\lambda,\pi}.\label{eqdimW}
  \end{equation}
  Since $\alpha$ is injective, we have also $\dim V_{\lambda,\pi}\le \dim W$
  and thus, using \eqref{eqdimW}, $\dim V_{\lambda,\pi}\le\rank M_\pi-\dim V_{\lambda,\pi}$ whence the conclusion $2\dim V_{\lambda,\pi}\le \rank M_\pi$.
  \end{solution}
\begin{solution}{\ref{exerciseIw}}
We use an induction on the length of $w$. The property is true
if $w$ is the empty word. Next, assume that $I(w)$ is a semi-interval
and let $a$ be a letter. Then $T(I(aw))=T(\Delta_a)\cap I(w)$ is a semi-interval
since $T(\Delta_a)$ is a semi-interval and also $I(w)$ by induction hypothesis.
Since $I(aw)\subset \Delta_a$, the set $T(I(aw))$ is a translation
of $I(aw)$, which is therefore also a semi-interval.
\end{solution}
\begin{solution}{\ref{exerciseJw}}
The proof is symmetrical to the proof for $I(w)$, using
this time the fact that $T^{-1}J(wa)=J(w)\cap \Delta_a$
is a semi-interval and thus that $J(wa)$ is a semi-interval
since $J(wa)\subset T \Delta(a)$.
\end{solution}
%%%%%%%%%%%%%%%%%%%%%%%%
\section{Notes}
% IET
\subsection{Interval exchange transformations}
Interval exchange transformations were introduced by \cite{Keane1975}
\index{names}{Keane, Michael}%
who proved Theorem~\ref{theoremIDOC}. The condition defining
regular interval exchange transformations is also called
the \emph{infinite disjoint orbit condition} or \emph{idoc}
\index{subect}{idoc}.

The idea of Rauzy
induction and Theorem~\ref{lemmaRauzyInduction} are
from~\cite{Rauzy1979}.
\index{names}{Rauzy, G\'erard}%
For more details on interval exchange transformations we refer to \cite{CornfeldFominSinai1982}, 
that we follow closely for the proof of Theorems
\ref{lemmaRauzyInduction} and ~\ref{theoremIDOC}.
\index{names}{Cornfeld, Isaac P.}\index{names}{Fomin, Sergei V.}\index{names}{Sinai, Yakov V.}%
The Cantor version of interval exchange
transformations was introduced by \cite{Keane1975}.

The notion of planar dendric shifts and Proposition~\ref{propositionPlanarDendric} are from \cite{BertheDeFeliceDolceLeroyPerrinReutenauerRindone2015d}.
The converse of Proposition \ref{propositionPlanarDendric},
characterizing the languages of regular interval exchange transformations,
is proved in \cite{FerencziZamboni2008}.

Theorem~\ref{theoremBVrepresentationIET}
is due to \cite{Gjerde&Johansen:2002}.
\index{names}{Gjerde, Richard}\index{names}{Johansen, Orjan}% 
They also showed that there are BV-dynamical systems satisfying the hypothesis of the theorem that are not isomorphic to a Cantor version of an interval exchange transformation.

The BV-representation of non-minimal Cantor systems can also be considered.
In \cite{Medynets2006}\index{names}{Medynets, Konstantin} the author shows that for Cantor dynamical systems without periodic points (but not  necessarily minimal)
a BV-representation can also be given.
It is applied in \cite{BezuglyiKwiatkowskiMedynets2009}
\index{names}{Bezuglyi, Sergey}\index{names}{Kwiatkowski, Jan}\index{names}{Medynets, Konstantin} to subshifts generated by non-primitive substitutions.
The authors show that  they have stationary BV-representations as in the minimal case.

Theorem \ref{theo:rauzy1} is Theorem 14 in~\cite{Rauzy1979}
while Theorem~\ref{theo:birauzy2} is Theorem 23.
Lemma \ref{lem:bi22} is  the two-sided version of Lemma 22 in~\cite{Rauzy1979}.

We have  noted that for any $s$-interval exchange transformation on $[\ell, r)$ and any semi-interval $I$ of $[\ell, r[$, the transformation $S$ induced by $T$ on $I$ is an interval exchange transformation on at most $s+2$-intervals
(Lemma~\ref{lemmaRauzyInduction}).
Actually, it follows from the proof of Lemma 2, page 128 in~\cite{CornfeldFominSinai1982} that, if $T$ is regular and $S$ is an $s$-interval exchange transformation with separation points $\Sep(S) = \Div(I,T)\cap I$, then $I$ is admissible.
Thus the converse of Theorem~\ref{theo:birauzy1} is also true.

Branching Rauzy induction is introduced in \cite{DolcePerrin2017r}
where Theorem~\ref{theo:birauzy2} appears. Actually, left
Rauzy induction is already considered in \cite{Veech1990}
and Theorem~\ref{theo:birauzy2} appears independently
in~\cite{Fickenscher2017}.

Proposition \ref{pro:autoelem} appears in~\cite{Jullian2013}.
\index{names}{Jullian, Yann}%

 On the relation between Rauzy induction
and continued fractions, see~\cite{MiernowskiNogueira2013} for more details.
\index{names}{Miernowski, Tomasz}\index{names}{Nogueira, Arnaldo}%

Theorem~\ref{theoremBoshernitzanCarroll} is from
\cite{BoshernitzanCarroll1997}.
\index{names}{Boshernitzan, Michael}\index{names}{Carroll, C.R.}%
In the same paper, an extension to right Rauzy induction is suggested (but not complety developed). Theorem~\ref{theo:quadratic} is
from \cite{DolcePerrin2017r}.

It was conjectured by Keane in \citep{Keane1975}
\index{names}{Keane, Michael} that a minimal
interval exchange transformation is uniquely ergodic. This
was disproved by \cite{Keynes1976}\index{names}{Keynes, Harvey B.}
\index{names}{Newton, Dan} with a $5$-interval exchange
transformation. Later Keane found examples with $k=4$
\cite{Keane1977} and conjectured that almost all minimal
interval exchange transformations are uniquely ergodic.
The conjecture was proved independently by \cite{Veech1982}
\index{names}{Veech, William A.}
and by \cite{Masur1982}.\index{names}{Masur, Howard}
 Exercises \ref{exerciseInvariantMeasureIET} to \ref{exerciseTheorem2.12}
are from \cite{Veech1978}.

\subsection{Linear involutions}
Linear involutions were introduced in \cite{DanthonyNogueira1990}.
\index{names}{Danthony, Claude}\index{names}{Nogueira, Arnaldo}%
It is also an extension of the notion of interval exchange with flip
\cite{Nogueira1989,NogueiraPiresTroubetzkoy2013}.
\index{names}{Pires, Benito}\index{names}{Troubetzkoy, Serge}%
Our definition is somewhat more general than
the one used in~\cite{DanthonyNogueira1990}
and also that of \cite{BertheDelecroixDolcePerrinReutenauerRindone2017b}.
\index{names}{Berth\'e, Val\'erie}\index{names}{Delecroix, Vincent}%
\index{names}{Dolce, Francesco}\index{names}{Perrin, Dominique}%
\index{names}{Reutenauer, Christophe}\index{names}{Rindone, Giuseppina}%
Orientable linear involutions correspond to orientable laminations, whereas 
 coherent
linear involutions correspond to orientable surfaces. Thus coherent
nonorientable involutions correspond to nonorientable laminations
on orientable surfaces.

Contrary to what happens with interval exchanges,
it is shown in~\cite{DanthonyNogueira1990} that non coherent linear 
involutions are almost surely not minimal.

Proposition \ref{propBL} (already  proved in~\cite[Proposition 4.2]{BoissyLanneau2009}
\index{names}{Boissy, Corentin}\index{names}{Lanneau, Erwan}%
 for the class of coherent involutions)
is~\cite[Proposition 3.7]{BertheDelecroixDolcePerrinReutenauerRindone2017b}.
The proof uses Keane's theorem proving that an
interval exchange transformation without connections is minimal~\citep{Keane1975}. 

The interval exchange $U$ built
in Example~\ref{exampleFiboDoubleInvolution} is actually the orientation covering
of the linear involution $T$
(see~\cite{BertheDelecroixDolcePerrinReutenauerRindone2017b}).

%%%%%%%%%%%%%%%%%%%%%%%%%%%%%
% chapter Bratteli Diagrams and C* algebras
%%%%%%%%%%%%%%%%%%%%%%%%%%%%%

\chapter{Bratteli diagrams and $C^*$-algebras}\label{chapterBratteli}

In this chapter, we give a short introduction to 
the connection between Bratteli diagrams and $C^*$-algebras. 
These notions are closely related to our subject. 
%We refer the reader to~\cite{Davidson1996} for a more detailed exposition.
Indeed, one may associate to every ordered Bratteli diagram both
a $C^*$-algebra and 
a Cantor  system. Minimal Cantor systems correspond to the so-called properly
ordered Bratteli diagrams and the dimension groups of the Cantor system
and of the algebra are the same 
%(see the beautiful book~\cite{Putnam2018}).

We define approximately finite dimensional algebras (AF algebras) as follows.
A $C^*$-algebra $\Ga$ is an AF algebra if it is the closure of an increasing sequence of finite dimensional subalgebras $(\Ga_k)_{k\ge 0}$.
Let $t_k$ be the dimension of $\Ga_k$.
To such an algebra $\Ga$ we associate a dimension group as a direct limit of abelian groups $\Z^{t_k}$. 

The main object of this chapter is to prove the theorem of Elliott (Theorem \ref{theoremElliott}).
It asserts that two unital approximately finite algebras are $*$-isomorphic if and only if their dimension groups are isomorphic. 
In this way the difficult problem of isomorphism of $C^*$-algebras is, in the case of AF algebras, reduced to the easier problem of isomorphism of dimension groups. 

In the first section, we show how Bratteli diagrams 
can describe  embeddings of sequences of finite dimensional
algebras. In the next section (Section~\ref{sectionC*algebras},
we give a brief introduction to $C^*$-algebras. 
In Section \ref{sectionDirectLimitsC*} we introduce enveloping
$C^*$-algebras and use them to define direct limits
of $C^*$-algebras. Elliott's Theorem (Theorem \ref{theoremElliott})
is presented in Section~\ref{sectionApproximatelyFinite}

%%%%%%%%%%%%%%%%%%%%%
\section{Bratteli diagrams}\label{sectionBratteli}
Let us first give a modified version of
Bratteli diagrams adapted to the purpose of this chapter
and which essentially differs by the addition of
integer labels to the vertices.
A Bratteli diagram, as defined in Chapter~\ref{chapterBratteliDiagrams},
 is an infinite labeled multigraph.
Its set of nodes of level $k$ consists of
pairs $(k,j)$ with $p\ge 0$ and $1\le j\le t_k$
with  $t_0=1$ (there is only one vertex at level $0$).
There are $a_{ij}^{(k)}\ge 0$ arrows from $(k,j)$ to $(k+1,i)$.
We denote by $M_k$ the $t_{k+1}\times t_{k}$-matrix with
coefficients $a_{ij}^{(k)}$.
We label the vertex $(k+1,i)$ of the diagram by the integer 
\begin{equation}
n(k+1,i)=\sum_{j=1}^{t_p}a_{ij}^{(k)}n(k,j)\label{equationA}
\end{equation}
with $n(0,1)=1$.
Conversely, any sequence $(M_k)_{k\geq 1}$ of nonnegative integer matrices defines a Bratteli diagram as above. By convention, we assume that 
$M_0=\begin{bmatrix}1&1&\ldots&1\end{bmatrix}^t$ and 
that $n(0,1)=1$. Thus $n(1,i)=1$ for $1\le i\le t_1$.
\begin{example}\label{exampleBratteli1}
The graph represented in Figure~\ref{figureBratteli1}
with all matrices $M_k$ for $k\ge 1$ equal to the matrix $\begin{bmatrix}1&1\\1&0\end{bmatrix}$ is a Bratteli diagram. The levels $k=0,1,2,\ldots$ are growing horizontally from left to right
  and on each level the vertices are numbered vertically 
  from bottom to top.
The integer labeling each node $(k,j)$ is $n(k,j)$.
\begin{figure}[hbt]
\centering

\begin{tikzpicture}(70,10)
\node(0)at(0,.5){$1$};
\node(10)at(1,1){$1$};\node(11)at(1,0){$1$};
\node(20)at(2,1){$1$};\node(21)at(2,0){$2$};
\node(30)at(3,1){$2$};\node(31)at(3,0){$3$};
\node(40)at(4,1){$3$};\node(41)at(4,0){$5$};
\node(50)at(5,1){$5$};\node(51)at(5,0){$8$};
\node(60)at(6,0){$\ldots$};\node(61)at(6,0){$\ldots$};

\draw[->](0)edge node{}(10);
\draw[->](0)edge node{}(11);
\draw[->](10)edge node{}(21);\draw[->](11)edge node{}(20);
\draw[->](11)edge node{}(21);
\draw[->](20)edge node{}(31);\draw[->](21)edge node{}(30);
\draw[->](21)edge node{}(31);
\draw[->](30)edge node{}(41);\draw[->](31)edge node{}(40);
\draw[->](31)edge node{}(41);
\draw[->](40)edge node{}(51);\draw[->](41)edge node{}(50);
\draw[->](41)edge node{}(51);
\draw[->](50)edge node{}(61);\draw[->](51)edge node{}(60);
\draw[->](51)edge node{}(61);
\end{tikzpicture}
\caption{A Bratteli diagram.}\label{figureBratteli1}
\end{figure}
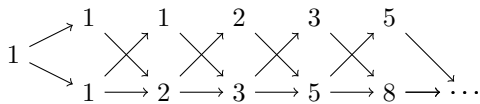

\end{example}
\begin{example}\label{exampleBratteli2}
The graph represented in Figure~\ref{figureBratteli2} is a Bratteli
diagram which corresponds to all matrices $M_n$
equal to $\begin{bmatrix}1&1\\0&1\end{bmatrix}$.
\begin{figure}[hbt]
\centering

\begin{tikzpicture}(70,10)
\node(0)at(0,.5){$1$};
\node(10)at(1,0){$1$};\node(11)at(1,1){$1$};
\node(20)at(2,0){$2$};\node(21)at(2,1){$1$};
\node(30)at(3,0){$3$};\node(31)at(3,1){$1$};
\node(40)at(4,0){$4$};\node(41)at(4,1){$1$};
\node(50)at(5,0){$5$};\node(51)at(5,1){$1$};
\node(60)at(6,0){$\ldots$};\node(61)at(6,1){$\ldots$};

\draw(0)edge node{}(10);\draw[->](0)edge node{}(11);
\draw[->](10)edge node{}(20);\draw[->](11)edge node{}(20);
\draw[->](11)edge node{}(21);
\draw[->](20)edge node{}(30);\draw[->](21)edge node{}(30);
\draw[->](21)edge node{}(31);
\draw[->](30)edge node{}(40);\draw[->](31)edge node{}(40);
\draw[->](31)edge node{}(41);
\draw[->](40)edge node{}(50);\draw[->](41)edge node{}(50);
\draw[->](41)edge node{}(51);
\draw[->](50)edge node{}(60);\draw[->](51)edge node{}(60);
\draw[->](51)edge node{}(61);
\end{tikzpicture}
\caption{A second Bratteli diagram.}\label{figureBratteli2}
\end{figure}
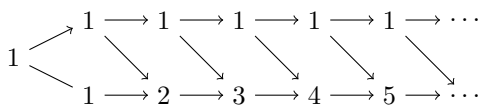
\end{example}
We will now explain how a Bratteli diagram describes
a sequence \begin{displaymath}
\Ga_1\stackrel{\pi_1}{\rightarrow}\Ga_2\stackrel{\pi_2}{\rightarrow}
\Ga_3\stackrel{\pi_3}{\rightarrow}\cdots
\end{displaymath}
of algebras $\Ga_k$ and morphisms $\pi_k:\Ga_k\rightarrow \Ga_{k+1}$.

Let $M$ be an $\ell\times k$-matrix with coefficients $a_{ij}\ge 0$,
let $(m_1,\ldots,m_t)$ be integers and let $(n_1,\ldots,n_s)$
be such that
\begin{displaymath}
n_i=\sum_{i=1}^ta_{ij}m_j\quad (1\le i\le \ell).
\end{displaymath}

Let $\Ga_1=\M_{m_1}\oplus\ldots\oplus \M_{m_t}$ and
$\Ga_2=\M_{n_1}\oplus\ldots\oplus \M_{n_s}$ where
$\M_n$\index{symbols}{M@$\M_n$} denotes the algebra of $n\times n$-matrices
with complex coefficients.
We associate to an $s\times t$-matrix $M$ with integer nonnegative coefficients $a_{ij}$ 
the morphism $\varphi:\Ga_1\rightarrow \Ga_2$ defined as
$\varphi=\varphi_1\oplus\ldots\oplus \varphi_t$ with
\begin{equation}
\varphi_i=\id_{m_1 }^{(a_{i1})}\oplus\ldots\oplus\id_{m_t}^{(a_{it})}\label{eqPartialMult}
\end{equation}
where $\id_n$\index{symbols}{id@$\id_n$} is the identity matrix of $\M_n$ and $\id_n^{(e)}:\M_n\rightarrow \M_{en}$
is the morphism 
$x\mapsto (x,\ldots,x)$ ($e$ times).\index{symbols}{id@$\id_n^{(e)}$}
The matrix $M$ is called the matrix of \emph{partial multiplicities}
\index{subject}{matrix!of partial multiplicities}\index{subject}{partial!multiplicities!matrix of}
associated with $\varphi$ and $\varphi$ is the morphism determined by $M$.

We associate to a Bratteli diagram the sequence
\begin{equation}
\Ga_1\stackrel{\pi_1}{\rightarrow}\Ga_2\stackrel{\pi_2}{\rightarrow}
\Ga_3\stackrel{\pi_3}{\rightarrow}\cdots \label{eqAF}
\end{equation}
of algebras $\Ga_k$ with algebra morphisms $\pi_k:\Ga_k\rightarrow \Ga_{k+1}$
where 
\begin{equation}
\Ga_k=\M_{n(k,1)}\oplus\ldots\oplus \M_{n(k,t_k)}\label{eqCstarFinDim}
\end{equation}
and where $\pi_k:\Ga_k\rightarrow \Ga_{k+1}$ is associated with the matrix
$M_k$ as above.
\begin{example}\label{exampleFibo1}
The sequence of algebras $\Ga_k$ associated to the Bratteli diagram of 
Example~\ref{exampleBratteli1} is $\Ga_k=\M_{n_k}\oplus\M_{n_k-1}$
where $(n_k)$ is the \emph{Fibonacci sequence}
\index{subject}{Fibonacci!sequence of numbers}
 defined by $n_1=n_2=1$ and $n_{k+1}=n_k+n_{k-1}$
for $k\ge 2$. The morphism $\pi_k$ is defined by
$\pi_k(x,y)=(x\oplus y,x)$.
\end{example}
\begin{example}\label{exampleK1}
The sequence of algebras $\Ga_k$ associated to the Bratteli diagram
of Figure~\ref{figureBratteli2} is $\Ga_k=\M_k\oplus\C$. The 
corresponding morphism $\pi_k:\Ga_k\rightarrow \Ga_{k+1}$
is defined by  $\pi_k(x,\lambda)=(x\oplus\lambda,\lambda)$.
\end{example}
As a more general notion of Bratteli diagram, one may consider \emph{partial
diagrams}\index{subject}{Bratteli diagram!partial}\index{subject}{partial!diagram}
 in which the labels $n(p,i)$ satisfy the inequalities
\begin{equation}
n(k+1,j)\ge\sum_{i=1}^{t_p}a_{ji}^{(k)}n(k,i).\label{equationAbis}
\end{equation}
instead of an equality.
\begin{example}
The diagram represented in Figure~\ref{figureBratteli2bis}
with all matrices $M_k$ equal to $[1]$
is a partial Bratteli diagram.
\begin{figure}[hbt]
\centering

\begin{tikzpicture}

\node(20)at(2,0){$1$};
\node(30)at(3,0){$2$};
\node(40)at(4,0){$3$};
\node(50)at(5,0){$4$};
\node(60)at(6,0){$\ldots$};

\draw[->](20)edge node{}(30);
\draw[->](30)edge node{}(40);
\draw[->](40)edge node{}(50);
\draw[->](50)edge node{}(60);
\end{tikzpicture}
\caption{A partial Bratteli diagram.}\label{figureBratteli2bis}
\end{figure}

\end{example}
A partial Bratteli diagram represents as above a sequence of $*$-morphisms
between finite dimensional algebras, but the morphisms
are this time not unital.
\begin{example}\label{exampleK1bis}
The Bratteli diagram of Figure~\ref{figureBratteli2bis}
corresponds to the sequence of algebras $(\M_k)_{k\ge 1}$
with the morphisms $\pi_k:\M_k\rightarrow \M_{k+1}$ defined by 
$\pi_k(x)=x\oplus 0$.
\end{example}
%%%%%%%%%%%%%%%
\section{$C^*$-algebras}\label{sectionC*algebras}
%Recall that a \emph{Banach algebra} is a normed algebra $\Ga$ which is complete
%and with the norm such that $\norm{AB}\le \norm{A}\norm{ B}$ for all $A,B\in \Ga$.

A \emph{$*$-algebra}\index{subject}{algebra@$*$-algebra} $\Ga$ is a complex  algebra with an idempotent map (called the \emph{adjoint}\index{subject}{adjoint!of element})
$A\mapsto A^*$ related to the algebra
structure by 
\begin{equation}
(A+B)^*=A^*+B^*,\quad (AB)^*=B^*A^*,\mbox{ and }
(\lambda A)^*=\bar{\lambda}A^*\label{eqStar1}
\end{equation}
for every $A,B\in \Ga$ and $\lambda\in\C$. It is usual to denote
with capitals the elements of a $C^*$-algebra. The notation $A$
should not lead the reader into a confusion with an alphabet,
nor the notation $A^*$ with the free monoid on $A$.

A \emph{$C^*$-algebra}\index{subject}{Cstar-algebra@$C^*$-algebra} $\Ga$ is a Banach
$*$-algebra  such that for all $A\in\Ga$, the identity
\begin{equation}
\norm{A^*A}=\norm{A}^2\label{eqStar2}
\end{equation}
\index{symbols}{A@$\norm{A}$}%
is satisfied. Equation~\ref{eqStar2} is called the $C^*$-\emph{identity}.
\index{subject}{Cstar-identity@$C^*$-identity} It implies
that $\norm{A}=\norm{A^*}$, that is, the adjoint
map $A\mapsto A^*$ is an isometry
(see Exercise~\ref{exerciseAltCstarIdentity}).
Thus \eqref{eqStar2} implies also
\begin{equation}
  \norm{A^*A}=\norm{A}\norm{A^*}\label{eqStar2bis}
\end{equation}

A \emph{morphism}\index{subject}{morphism!of $C^*$-algebras}
\index{subject}{Cstar-algebra@$C^*$-algebra!morphism of}
from a $C^*$-algebra $\Ga$ to a $C^*$-algebra $\Gb$
(also called a   $*$-\emph{morphism}\index{subject}{morphism@$*$-morphism})
is an algebra morphism $\pi$ such that $\pi(A^*)=\pi(A)^*$
for all $A\in \Ga$. It can be shown that this implies
$\norm{\pi(A)}\le\norm{A}$ for every $A\in\Ga$ 
(Exercise~\ref{exerciseNormPi(A)}).

A $C^*$-algebra is \emph{unital}\index{subject}{unital!C@$C^*$-algebra}
\index{subject}{Cstar-algebra@$C^*$-algebra!unital} if it has an identity element $I$.
It follows from \eqref{eqStar1} that $I^*=I$ and from \eqref{eqStar2}
that $\norm{I}=1$.

 Each algebra $\M_k$ is a $C^*$-algebra,
using the usual conjugate  transpose  $A^*$ of $A$
as adjoint of $A$ and using the matrix norm induced by the Hermitian norm.  
Indeed, Equation~\eqref{eqStar2} is a consequence of the Cauchy-Schwartz
inequality $|\langle x,y\rangle|\le\norm{x}\norm{y}$ which implies
\begin{displaymath}
\norm{A^*A}=\sup_{\norm{x}=\norm{y}=1}\langle A^*Ax,y\rangle
=\sup_{\norm{x}=\norm{y}=1}\langle Ax,Ay\rangle=
\norm{A}^2.
\end{displaymath}
More generally,  the algebra $\Gb(\HH)$\index{symbols}{BH@$\Gb(\HH)$} of bounded linear operators
on a Hilbert space $\HH$
 is a $C^*$-algebra. The adjoint of $A$ is
defined by $\langle A^*x,y\rangle=\langle x,Ay\rangle$ for all
$x,y\in \HH$.

As a second fundamental example, the space $C(X,\C)$ of continuous
complex valued functions on a compact space $X$ is a $C^*$-algebra
with complex conjugation as adjoint operation. We have indeed
\begin{displaymath}
\|\overline{f}f\|=\sup_{x\in X}|\overline{f(x)}f(x)|=\sup_{x\in X}|f(x)|^2=\|f\|^2
\end{displaymath}

An element $A$ of a $C^*$-algebra is \emph{self-adjoint}\index{subject}{element!self-adjoint}\index{subject}{self-adjoint!element}
if $A=A^*$. It is
\emph{normal}\index{subject}{normal!element of $C^*$-algebra} if $AA^*=A^*A$.
Recall that $\spr(A)$ denotes the spectral radius of $A$
(see Appendix~\ref{appendixTopo}). The following
property is clear when $\Ga=\M_k$ (as well as the necessity of the hypothesis
that $A$ is normal).
\begin{proposition}\label{propositionNormal}
In a $C^*$-algebra, one has $\|A\|=\spr(A)$ for every normal
element of $A$.\index{symbols}{spec@$\spr$}
\end{proposition}
\begin{proof}
Let first $A$ be self-adjoint. Then by repeated use of Equation
\eqref{eqStar2}, we have
\begin{equation}
\spr(A)=\lim\|A^{2^n}\|^{2^{-n}}=\|A\|
\end{equation}
When $A$ is normal, using the preceding case, we have
\begin{eqnarray*}
\spr(A)^2&\le&\|A\|^2=\|A^*A\|=\lim\|(A^*A)^n\|^{1/n}\\
&\le&\lim(\|(A^*)^n\|\|(A^n)\|)^{1/n}=\spr(A)^2.
\end{eqnarray*}
\end{proof}
One can deduce from this property that the norm is unique in a $C^*$-algebra
(Exercise \ref{exerciseNormUnique}).

A subalgebra $\Gb$ of a $C^*$-algebra $\Ga$ is \emph{self-adjoint}\index{subject}{self-adjoint!subalgebra}\index{subject}{subalgebra!self-adjoint} if $A^*$ is in $\Gb$ for every
$A\in\Gb$. It is easy to verify that a closed self-adjoint subalgebra of
a $C^*$-algebra is again a $C^*$-algebra.

The direct sum $\Ga\oplus\Gb$\index{symbols}{A@$\Ga\oplus\Gb$}
of two $C^*$-algebras $\Ga,\Gb$ is itself a $C^*$-algebra, using
on the set $\Ga\times\Gb$ the componentwise sum, product and $*$, and defining
$\norm{(A,B)}=\max(\norm{A},\norm{B})$.

\subsection{Ideals in $C^*$-algebras}
An \emph{ideal}\index{subject}{ideal!of $C^*$-algebra} in a $C^*$-algebra 
$\Ga$ is a norm-closed
two-sided ideal of the algebra $\Ga$. 

It can be shown that every ideal $\cal J$ is self-ajoint,
that is, such that $J^*$ is in $\cal J$ for every $J\in\cal J$.

The quotient $\Ga/\Id$ \index{symbols}{A@$\Ga/\Id$} of a $C^*$-algebra $\Ga$ by an ideal $\Id$
is the set of cosets \index{subject}{coset!of ideal}
\index{subject}{ideal!coset of} $A+\Id=\{A+J\mid J\in \Id\}$ with
the adjoint defined by $(A+\Id)^*=A^*+\Id$ and the norm
defined by $\|A+\Id\|=\inf_{J\in\Id}\|A+J\|$. 

Since $\Id$
is self-adjoint, we have $\|(A+\Id)^*\|=\|A+\Id\|$. Actually,
the $C^*$ norm condition is also satisfied and thus one has the
following statement.

\begin{theorem}\label{theoremQuotientAlgebra}
For every ideal $\cal J$ in a $C^*$-algebra
$\Ga$, the quotient algebra $\Ga/\cal J$ is a $C^*$-algebra.
\end{theorem}

The following now shows that the basic isomorphism theorem
for algebras still holds for $C^*$-algebras.
\begin{theorem}\label{theoremSeconIsomorphismThm}
If $\cal J$ is an ideal in $C^*$-algebra $\Ga$ and that
$\Gb$ is a $C^*$-subalgebra of $\Ga$. Then $\Gb+{\cal J}$
is a $C^*$-algebra and
\begin{displaymath}
\Gb/(\Gb+{\cal J})\simeq (\Gb+{\cal J})/{\cal J}.
\end{displaymath}
\end{theorem}
The following result is classical.
\begin{theorem}\label{theoremC*semisimple}
Every $C^*$-algebra is semisimple.
\end{theorem}
We give the proof for a finite dimensional $C^*$-algebra. Consider a $C^*$-algebra $\Ga$
of matrices contained in $\M_n$. Suppose that $I$ is a nilpotent ideal of $\Ga$. Up
to a change of basis, all matrices in $I$ are upper diagonal with zeroes on the diagonal.
But since $I$ is self adjoint, this forces $I=\{0\}$.
%%%%%%%%%%%%%%%%%%
\subsection{Finite dimensional $C^*$-algebras}

As a direct sum of $C^*$-algebras,
any algebra of the form \eqref{eqCstarFinDim} is itself a $C^*$-algebra.
Actually, any $C^*$-algebra which is finite dimensional (as a vector space) is of this form.
\begin{theorem}\label{theoremFiniteDimensionalC*}
Every finite dimensional $C^*$-algebra is $*$-isomorphic to a direct sum of full
matrix algebras 
\begin{displaymath}
\Ga=\M_{n_1}\oplus\M_{n_2}\oplus\ldots\oplus \M_{n_k}.
\end{displaymath}
\end{theorem}
\begin{proof}
  Since a $C^*$-algebra is semisimple by Theorem~\ref{theoremC*semisimple},
  the result  follows from Wedderburn Theorem
  (see Appendix~\ref{appendixGroups}).
  \end{proof}

%Each algebra $\M_n$ can be considered as a subalgebra of $\M_{n+1}$
%by sending $A$ to $A\oplus 0$. Moreover, all algebras $\M_n$
%can be considered as subalgebras of the $C^*$-algebra of norm
%continuous operators on the Hilbert space $\ell^2$ of
%sequences of complex numbers $x=(x_n)_{n\ge 0}$
%with the norm $\norm{x}^2=\sum_{n\ge 0}|x_n|^2$. One can
%for example consider $\M_n$ as the set of operators on the
%subspace of sequences $x$ such that $x_m=0$ for $m\ge n$.

 An element $U$ of
a $C^*$-algebra 
is \emph{unitary}\index{subject}{unitary!element}\index{subject}{element!unitary}
if $UU^*=U^*U=I$. Two elements $A,B$ are
\emph{unitarily equivalent}\index{subject}{unitary!equivalence} if 
$A=UBU^*$ where $U$ is a unitary element. Two morphisms $\varphi,\psi:\Ga\rightarrow \Gb$
are unitarily equivalent if 
there is a unitary element $U$ in $\Gb$
such that $\varphi=\Ad ( U) \circ \psi $ where
$\Ad ( U)(B)=UBU^*$ for every $B\in\Gb$.\index{symbols}{Ad@$\Ad(U)$}

The following result shows that the theory
of finite dimensional $C^*$-algebras is similar to that
of ordinary semisimple algebras.

\begin{proposition}\label{propSemiSimple}
Suppose that $\varphi$ is a unital $*$-morphism of  finite
dimensional $C^*$-algebras  from $\Ga=\M_{m_1}\oplus\ldots\oplus \M_{m_t}$
into $\Gb=\M_{n_1}\oplus\ldots\oplus\M_{n_s}$. Then $\varphi$ is determined
up to unitarily equivalence by its $s\times t$-matrix of
partial multiplicities.
\end{proposition}
\begin{proof}
  By the analysis made in Appendix \ref{appendixGroups}, there is an
  invertible matrix $U$ such that $\Ad (U)\circ \varphi=\varphi_1\oplus\ldots\oplus\varphi_t$
  is of the form $\varphi_i=\id_{m_1}^{(a_{i1})}\oplus\ldots\oplus\id^{(a_{it})}$ for $1\le i\le t$.

  There remains to show that if $\varphi:M\mapsto UMU^{-1}$ is a $C^*$-algebra morphism from $\M_n$
  onto itself, then $U$ is unitary. Consider the elementary matrices $E_{ij}\in\M_n$
  having all entries equal to $0$ except the entry $i,j$ which is $1$. We have
  \begin{equation}
    \varphi(E_{ij})=\ell_ir_j
  \end{equation}
  where $\ell_i$ is the column of index $i$ of $U$ and $r_j$ is the row of index $j$ of $U^{-1}$. Such a decomposition with
  vectors $r_i,\ell_i$ such that $r_i\ell_i=1$ is unique.
  Since $E_{ij}^*=E_{ji}$ and since $\varphi(E_{ij})^*=E_{ij}^*$, we have $r_j^*=\ell_j$
  and $\ell_i^*=r_i$. Thus $U^{-1}=U^*$ and $U$ is unitary.
  \end{proof}

Thus, any sequence
 $\Ga_1\stackrel{\varphi_1}{\rightarrow}\Ga_2\stackrel{\varphi_1}{\rightarrow}\Ga_3\ldots$  of $*$-morphisms between finite dimensional $C^*$-algebras
is determined by a Bratteli diagram. This fact plays a fundamental role in the sequel.
%%%%%%%%%%%%%%
\subsection{Direct limits of $C^*$-algebras}\label{sectionDirectLimitsC*}
Let $\Ga$ be a $*$-algebra. A $C^*$-\emph{seminorm}\index{subject}{seminorm}
on $\Ga$ is a seminorm $\rho$ on $\Ga$ such that for all $A,B\in\Ga$,
\begin{displaymath}
\rho(AB)\le\rho(A)\rho(B),\quad \rho(A^*)=\rho(A),\ \mbox{ and }
\rho(A^*A)=\rho(A)^2.
\end{displaymath}
If additionally, $\rho$ is in fact a norm, it is called a $C^*$-norm.

If $\varphi:\Ga\rightarrow \Gb$ is a $*$-morphism, the map
from $\Ga$ into $\R_+$ defined by
$\rho(A)=\norm{\varphi(A)}$ is a $*$-seminorm on $\Ga$.

If $\rho$ is a $C^*$-seminorm on a $*$-algebra $\Ga$, the set
${\cal N}=\rho^{-1}(0)$ is self adjoint ideal of $\Ga$ and we get a $C^*$-norm
on the quotient $*$-algebra $\Ga/{\cal N}$ by setting $\norm{A+{\cal N}}=\rho(A)$.
If $\Gb$ denotes the Banach space completion of $\Ga/\cal N$ with this
norm, it can be checked that the operations on $\Ga/\cal N$ extend
uniquely to $\Gb$ to make a $C^*$-algebra, called the
\emph{enveloping}\index{subject}{enveloping $C^*$-algebra} $C^*$-algebra of
$\Ga$ with respect to $\rho$.

Let now $(\Ga_k)$ be a sequence of $C^*$-algebras with morphisms
$\varphi_{k,k+1}:\Ga_k\rightarrow \Ga_{k+1}$. The set
\begin{displaymath}
\Gb=\{(A_k)_{k\ge 0}\mid A_k\in\Ga_k, A_{k+1}=\varphi_{k,k+1}(A_k) \text{ for every
$k$ large enough}\}
\end{displaymath}
is a $*$-subalgebra of the direct product $\prod_{k\ge 0}\Ga_k$. Since
the $\varphi_{k,k+1}$ are $C^*$-algebra morphisms, they are
norm decreasing and we can define a semi-norm on $\Gb$ by setting
$\rho(A)=\lim\norm{A_n}$ for $A=(A_n)_{n\ge 0}$. Note that
$\rho(A)=0$ for every $A=(A_n)_{n\ge 0}$ such that $A_n=0$ for $n$
large enough.

The enveloping $C^*$-algebra $\Ga$ of $\Gb$ with respect to $\rho$
is called the \emph{direct limit} \index{subject}{direct limit!of $C^*$-algebras}
\index{subject}{Cstar-algebra@$C^*$-algebra!direct limit of}%
of the sequence $(\Ga_n)$, denoted
$\displaystyle{\lim_\rightarrow}\Ga_k$.\index{symbols}{A@$\lim_\rightarrow\Ga_k$}

As in the case of a direct limit of groups, there is a natural
morphism $\varphi_k$ from each $\Ga_k$ into $\Ga$ which sends $A\in\Ga_k$
to the class of any sequence $(A_\ell)_{\ell\ge 0}$ such that $A_k=A$
and $A_{\ell+1}=\varphi_{\ell,\ell+1}(A_\ell)$ for all $\ell\ge k$.

We will see examples of direct limits of $C^*$-algebras in the next section.

\subsection{Positive elements}
An element $A$ of a $C^*$-algebra $\Ga$ is \emph{positive}\index{subject}{positive!element in $C^*$-algebra}\index{subject}{element!positive}
if $A=A^*$ (that is, $A$ is self-adjoint)
and its spectrum $\sigma(A)$ is contained in $\R_+$.

A useful property of positive elements is the existence,
for every positive element $A\in \Ga$ of a unique positive square root
denoted $A^{1/2}$.

It can be shown that if $A,B\in\Ga$ are positive, then $A+B$ is positive.
As a consequence, the positive elements determine an order on the set of self-adjoint
elements by $A\le B$ if $B-A$ is positive. Indeed, if $A\le B$
and $B\le C$, then $C-A=(C-B)+(B-A)$ which is positive and thus $A\le C$.

\begin{theorem}
For every $A\in \Ga$, the element $A^*A$ is positive.
\end{theorem}

%A \emph{multiplicative linear functional}\index{subject}{multiplicative
%linear functional} on a Banach algebra $\Ga$ is
%nonzero homomorphism from $\Ga$ into $\C$.
A \emph{positive linear functional}\index{subject}{positive!linear functional}
on $C^*$-algebra is a 
linear functional $f:\Ga\to \C$ such that
$f(A)\ge 0$ whenever $A\ge 0$. A \emph{state}\index{subject}{state!in $C^*$-algebra} is a positive linear functional of norm $1$.

For example, let $\HH$ be a Hilbert space
and let $\pi$ is a $*$-morphism from $\Ga$ into $\B(\HH)$. Then
\begin{displaymath}
f(A)=\langle\pi(A)x,x\rangle
\end{displaymath}
is a positive linear functional. Indeed, if $A\ge 0$, using
the square root $A^{1/2}$ of $A$, we have
\begin{displaymath}
f(A)=\langle\pi(A^{1/2})^2x,x\rangle=\|\pi(A^{1/2})x\|^2\ge 0.
\end{displaymath}
It is a state if $\Ga$ is unital and $\|x\|=1$.
%%%%%%%%%%%%%%%%%%
\section{Approximately finite algebras}\label{sectionApproximatelyFinite}
A $C^*$-algebra is \emph{approximately finite dimensional}\index{subject}{Cstar-algebra@$C^*$-algebra!approximately finite}\index{subject}{approximately finite algebra}\index{subject}{algebra!approximately finite}
 (or an \emph{AF} algebra)\index{subject}{algebra!AF} if it is 
the closure of an increasing union of finite dimensional subalgebras $(\Ga_k)_{k\geq 0}$.
When $\Ga$ is unital, we further stipulate that $\Ga_0$ consists
of scalar multiples of the identity element $1$.

As an equivalent definition, an \emph{AF} algebra is a direct limit
of finite dimensional $C^*$-algebras. Indeed, if $\varphi_{k,k+1}:\Ga_k\rightarrow \Ga_{k+1}$ are morphisms, then $\Ga=\displaystyle{\lim_\rightarrow} \Ga_k$
is an \emph{AF} algebra since it is the closure of the
images $\varphi_n(\Ga_k)$ of the finite dimensional $C^*$-algebras $\Ga_k$ by
the natural morphisms $\varphi_k$.

Let $\Ga=\overline{\cup_{k\ge 1}\Ga_k}$ be an AF algebra
with embeddings $\alpha_k$ from $\Ga_k$ into $\Ga_{k+1}$
We associate to the sequence $(\Ga_k)$ the Bratteli diagram
defined by the sequence of matrices of partial multiplicities
of the morphisms $\alpha_k$.

An AF algebra is separable \index{subject}{separable!space}.
Indeed, every finite dimensional $C^*$-algebra is separable.
The following statement, which we will admit, gives
an equivalent definition of AF algebras which does not depend
on a sequence of subalgebras.

\begin{theorem}\label{theoremCharacAF}
A $C^*$-algebra is AF if and only if it is separable and
for every $n\ge 1$, every $A_1,\ldots,A_n$ in $\Ga$ and 
every $\varepsilon>0$ there is a finite dimensional
$C^*$-subalgebra $\Gb$ of $\Ga$ and $B_1,\ldots,B_n$
such that $\|A_i,B_i\|\le\varepsilon$ for $1\le i\le n$.
\end{theorem}

As a simple example of a $C^*$-algebra
which is not an AF algebra,
the algebra  $C([0,1])$ is not an AF algebra (Exercise~\ref{exerciseNotAF}).

\begin{example}\label{exampleFibo2}
Consider again the sequence  of algebras $\Ga_k=\M_{n_k}\oplus\M_{n_{k-1}}$
with the morphisms $\pi_k:\Ga_k\rightarrow \Ga_{k+1}$
where $n_k$ is the Fibonacci sequence
(Example~\ref{exampleFibo1}). The direct limit of the sequence
of $C^*$-algebras $\Ga_k$
is the \emph{Fibonacci algebra}\index{subject}{Fibonacci!algebra}.
\end{example}
\begin{example}\label{exampleK2}
  Consider the sequence $\Ga_k=\M_k\oplus\C$ of algebras 
 with the
  embeddings $\pi_k(x,\lambda)=(x\oplus\lambda,\lambda)$, that is
   \begin{displaymath}
\pi_k(A,\lambda)=(\begin{bmatrix}A&0\\0&\lambda\end{bmatrix},\lambda)
  \end{displaymath}
  of Example~\ref{exampleK1}. The corresponding AF algebra is
 the $C^*$-algebra $\Gk+\C I$ where $\Gk$
is the $C^*$-algebra of  compact operators
\index{subject}{compact!operator} on the Hilbert space $\ell^2(\C)$
of sequences $x=(x_n)$ such that $\|x\|_2<\infty$. 
The algebra $\Gk\oplus\C I $ is the direct limit of the algebras $\overline{\cup_{k\ge 1}\Ga_k}$.
Thus $\Gk\oplus\C I$ is a unital AF algebra. The $C^*$-algebra
$\Gk$ itself (which is not unital) is obtained as the closure of the sequence
of subalgebras $\M_k$ as in Example~\ref{exampleK1bis}.
\end{example}

\begin{example}\label{exampleCAR}
Consider $\Ga_k=\M_{2^k}$ with the embedding of $\Ga_k$ into
$\Ga_{k+1}$ being
\begin{displaymath}
\varphi_k(A)=\begin{bmatrix}A&0\\0&A\end{bmatrix}
\end{displaymath}
The corresponding AF algebra is called the \emph{CAR algebra}.
\index{subject}{CAR algebra}\index{subject}{algebra!CAR}%
The corresponding Bratteli diagram is represented in Figure~\ref{figureCARAlgebra}.
\begin{figure}[hbt]
\centering
\tikzset{node/.style={circle,minimum size=0.4cm,inner sep=0pt}}
\tikzset{title/.style={minimum size=0cm,inner sep=0pt}}
\begin{tikzpicture}
\node[node](0)at(0,0){$1$};\node[node](1)at(2,0){$2$};
\node[node](2)at(4,0){$4$};
\node[node](3)at(6,0){$8$};
\node[title](4h)at(7,.3){};\node[title](4b)at(7,-.3){};
\node[title]at(7.5,0){$\cdots$};

\draw[bend left,->](0)edge node{}(1);\draw[bend right](0)edge node{}(1);
\draw[bend left,->](1)edge node{}(2);\draw[bend right](1)edge node{}(2);
\draw[bend left,->](2)edge node{}(3);\draw[bend right](2)edge node{}(3);
\draw[bend left=15](3)edge node{}(4h);\draw[bend right=15](3)edge node{}(4b);
\end{tikzpicture}
\caption{A Bratteli diagram for the CAR algebra.}\label{figureCARAlgebra}
\end{figure}
We recognize of course the BV-representation of the $(2^n)$ odometer.
Let $\Gb_k=\M_{2^k}\oplus \M_{2^k}$. Since
\begin{displaymath}
\varphi_{k+1}(\begin{bmatrix}B_1&0\\0&B_2\end{bmatrix})=
\left[
\begin{array}{cc|cc}
B_1&0&0&B_2\\
0&B_2&0&0\\ \hline
0&0&B_1&0\\
0&0&0&B_2
\end{array}
\right]
\end{displaymath}
the morphism $\varphi_{k+1}$ embeds $\Gb_k$ into $\Gb_{k+1}$.
This shows that $\Ga$ is also the algebra defined
by the diagram of Figure~\ref{figureCARAlgebra}.
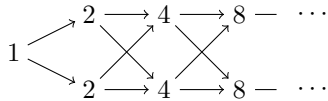
\begin{figure}[hbt]
\centering
\tikzset{node/.style={circle,minimum size=0.4cm,inner sep=0pt}}
\tikzset{title/.style={minimum size=0cm,inner sep=0pt}}
\begin{tikzpicture}
\node[node](0)at(0,.5){$1$};
\node[node](11)at(1,0){$2$};\node[node](12)at(1,1){$2$};
\node[node](21)at(2,0){$4$};\node[node](22)at(2,1){$4$};
\node[node](31)at(3,0){$8$};\node[node](32)at(3,1){$8$};
\node[title](4h)at(3.5,0){};\node[title](4b)at(3.5,1){};
\node[title]at(4,0){$\cdots$};\node[title]at(4,1){$\cdots$};

\draw[->](0)edge node{}(11);\draw[->](0)edge node{}(12);
\draw[->](11)edge node{}(21);\draw[->](11)edge node{}(22);
\draw[->](12)edge node{}(21);\draw[->](12)edge node{}(22);
\draw[->](21)edge node{}(31);\draw[->](21)edge node{}(32);
\draw[->](22)edge node{}(31);\draw[->](22)edge node{}(32);
\draw[](31)edge node{}(4h);\draw[](32)edge node{}(4b);
\end{tikzpicture}
\caption{Another Bratteli diagram for the CAR algebra.}\label{figureCARAlgebra2}
\end{figure}
\end{example}
We now give two examples of AF-algebras with non stationary diagram.
\begin{example}\label{exampleC(A^N)}
Let $X=\{0,1\}^\N$ be the one-sided full shift on $\{0,1\}$.
The space $\Ga=C(X)$ of continuous complex valued functions on $X$
is a commutative $C^*$-algebra. Let $\Ga_k$ be the subalgebra
of functions which are constant on each cylinder $[w]$ with
$w$ of length $k$. Then $\Ga=\overline{\cup\Ga_k}$
and thus $\Ga$ is an AF algebra. A Bratteli diagram for $\Ga$
is represented in Figure~\ref{figureBratteliC(X)}.
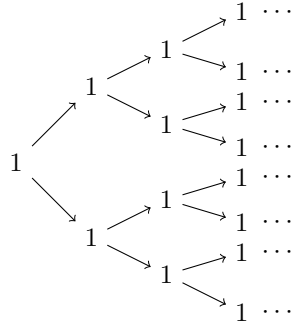
\begin{figure}[hbt]
\centering
\tikzset{node/.style={circle,minimum size=0.4cm,inner sep=0pt}}
\tikzset{title/.style={minimum size=0cm,inner sep=0pt}}
\begin{tikzpicture}
\node(0)at(0,0){$1$};
\node(11)at(1,-1){$1$};\node(12)at(1,1){$1$};
\node(21)at(2,-1.5){$1$};\node(22)at(2,-.5){$1$};\node(23)at(2,.5){$1$};\node(24)at(2,1.5){$1$};
\node(31)at(3,-2){$1$};\node(32)at(3,-1.20){$1$};\node(33)at(3,-.8){$1$};\node(34)at(3,-.2){$1$};
\node(35)at(3,.2){$1$};\node(36)at(3,.8){$1$};\node(37)at(3,1.2){$1$};\node(38)at(3,2){$1$};
\node(41)at(3.5,-2){$\cdots$};\node(42)at(3.5,-1.20){$\cdots$};\node(43)at(3.5,-.8){$\cdots$};\node(44)at(3.5,-.2){$\cdots$};
\node(45)at(3.5,.2){$\cdots$};\node(46)at(3.5,.8){$\cdots$};\node(47)at(3.5,1.2){$\cdots$};\node(48)at(3.5,2){$\cdots$};

\draw[->](0)edge node{}(11);\draw[->](0)edge node{}(12);
\draw[->](11)edge node{}(21);\draw[->](11)edge node{}(22);
\draw[->](12)edge node{}(23);\draw[->](12)edge node{}(24);
\draw[->](21)edge node{}(31);\draw[->](21)edge node{}(32);
\draw[->](22)edge node{}(33);\draw[->](22)edge node{}(34);
\draw[->](23)edge node{}(35);\draw[->](23)edge node{}(36);
\draw[->](24)edge node{}(37);\draw[->](24)edge node{}(38);
\end{tikzpicture}
\caption{A diagram for the algebra $C(\{0,1\}^\N)$.}\label{figureBratteliC(X)}
\end{figure}
\end{example}

\begin{example}\label{exampleGICAR}
Consider the Bratteli diagram represented in Figure~\ref{figurePascalTriangle}.
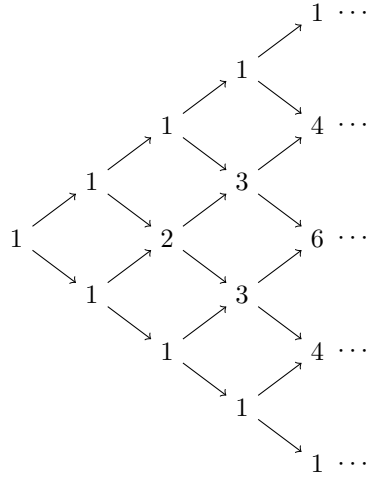
\begin{figure}[hbt]
\centering
\tikzset{node/.style={circle,minimum size=0.4cm,inner sep=0pt}}
\tikzset{title/.style={minimum size=0cm,inner sep=0pt}}
\begin{tikzpicture}
\node(0)at(0,0){$1$};
\node(11)at(1,-.75){$1$};\node(12)at(1,.75){$1$};
\node(21)at(2,-1.5){$1$};\node(22)at(2,0){$2$};\node(23)at(2,1.5){$1$};;
\node(31)at(3,-2.25){$1$};\node(32)at(3,-.75){$3$};\node(33)at(3,.75){$3$};\node(34)at(3,2.25){$1$};
\node(41)at(4,-3){$1$};\node(42)at(4,-1.5){$4$};\node(43)at(4,0){$6$};\node(44)at(4,1.5){$4$};
\node(45)at(4,3){$1$};
\node(51)at(4.5,-3){$\cdots$};\node(52)at(4.5,-1.5){$\cdots$};\node(53)at(4.5,0){$\cdots$};\node(54)at(4.5,1.5){$\cdots$};
\node(55)at(4.5,3){$\cdots$};

\draw[->](0)edge node{}(11);\draw[->](0)edge node{}(12);
\draw[->](11)edge node{}(21);\draw[->](11)edge node{}(22);
\draw[->](12)edge node{}(22);\draw[->](12)edge node{}(23);
\draw[->](21)edge node{}(31);\draw[->](21)edge node{}(32);
\draw[->](22)edge node{}(32);\draw[->](22)edge node{}(33);
\draw[->](23)edge node{}(33);\draw[->](23)edge node{}(34);
\draw[->](31)edge node{}(41);\draw[->](31)edge node{}(42);
\draw[->](32)edge node{}(42);\draw[->](32)edge node{}(43);
\draw[->](33)edge node{}(43);\draw[->](33)edge node{}(44);
\draw[->](34)edge node{}(44);\draw[->](34)edge node{}(45);
\end{tikzpicture}
\caption{The Pascal triangle.}\label{figurePascalTriangle}
\end{figure}
The corresponding AF algebra $\Ga$ is called the \emph{GICAR} algebra.
\index{subject}{GICAR algebra}\index{subject}{algebra!GICAR}%
We have $\Ga=\displaystyle{\lim_\rightarrow}\Ga_k$ where $\Ga_k$ is
\begin{displaymath}
\Ga_k=\M_{k\choose 0}\oplus\M_{k\choose 1}\oplus\ldots\oplus\M_{k\choose k}
\end{displaymath}
where 
\begin{displaymath}
{k\choose p}=\frac{k!}{(k-p)!p!}
\end{displaymath}
 is the binomial coefficient.
The embedding of $\Ga_k$ into $\Ga_{k+1}$ is
\begin{displaymath}
\varphi_k(A_0,A_1,\ldots,A_k)=(A_0,A_0+A_1,\ldots,A_{k-1}+A_k,A_k).
\end{displaymath}
\end{example}

The following result, which we will admit, is of fundamental importance.

\begin{theorem}[Bratteli]\label{theoremBratteli}
Let $\Ga=\overline{\cup\Ga_m}=\overline{\cup \Gb_n}$ be an AF
algebra obtained from two chains of finite dimensional $C^*$-algebras
$\Ga_m,\Gb_n$. There exists subsequences $m_k,n_k$ 
and an automorphism $\alpha$ of $\Ga$ such that
\begin{displaymath}
\Ga_{m_k}\subset\alpha(\Gb_{n_k})\subset\Ga_{m_{k+1}}
\end{displaymath}
for all $k\ge 1$.
\end{theorem}
We deduce the following characterisation of Bratteli diagrams
defining isomorphic AF algebras. It uses the notion
of intertwining of diagrams introduced in Chapter~\ref{chapterBratteliDiagrams}
and playing an essential role in the Strong Orbit Equivalence Theorem
(Theorem~\ref{theoremStrongOrbitEquivalence}).
\begin{corollary}\label{corollaryBratteli}
Two Bratteli diagrams define isomorphic AF algebras if and
only if they have a common intertwining.
\end{corollary}
\begin{proof}
This is a direct consequence of Theorem~\ref{theoremBratteli}.
Indeed, The Bratteli diagram corresponding to the sequence
\begin{displaymath}
\Ga_{m_1}\subset \alpha(\Gb_{n_1})\subset \Ga_{m_2}\subset\ldots
\end{displaymath}
is an intertwining of the diagrams corresponding to the sequences
$(\Ga_m)$ and $(\alpha(\Gb_n))$, the latter being the same
as the diagram corresponding to the sequence $(\Gb_n)$.
\end{proof}

\subsection{Simple AF algebras}
As for ordinary algebras, a $C^*$-algebra is \emph{simple}
\index{subject}{simple!C*-algebra@$C^*$-algebra}\index{subject}{C*-algebra@$C^*$-algebra!simple}%
if it has no nontrivial ideals. As well known every full matrix
algebra is simple.
The following result explains why the term of simple Bratteli 
diagram was chosen.
\begin{theorem}\label{theoremSimpleAF}
The AF algebra defined by a Bratteli diagram is simple
if and only if the diagram is simple.
\end{theorem}
We first prove the following lemma. The first statement holds
for general $C^*$-algebras because the hypothesis that
the algebras $\Ga_n$ are finite dimensional is not used in the proof.
\begin{lemma}\label{lemmaIdealAFalgebra}
If $\Id$ is an ideal of an AF algebra $\Ga=\overline{\cup_{k\ge 1}\Ga_k}$, then
\begin{displaymath}
\Id=\overline{\cup_{k\ge 1}(\Id\cap\Ga_k)}=\overline{\Id\cap\cup_{k\ge 1}\Ga_k}.
\end{displaymath}
In particular $\Id$ is an AF algebra.
\end{lemma}
\begin{proof}
Consider the commutative diagram
\begin{displaymath}
\begin{CD}
\Id\cap\Ga_n@>  >>\Ga_n@>\pi>>\Ga_n/(\Id\cap \Ga_n)\\
@VV        V @VV V @VV\alpha V\\
\Id@> >>\Ga_n+\Id@>\pi>>(\Ga_n+\Id)/\Id
\end{CD}
\end{displaymath}
where the unlabeled arrows are the canonical injections
and $\pi$ is the canonical quotient map from $\Ga$
to $\Ga/\Id$.
The map $\alpha$ is an isomorphism by Theorem~\ref{theoremSeconIsomorphismThm}.
If $J\in\Id$, since $\A=\overline{\cup_{k\ge 1}\Ga_k}$,
there is for every $\varepsilon>0$ an element $A$ of 
some $\Ga_k$ such that $\|J-A\|<\varepsilon$. Then, by definition of the
norm in the quotient,
$\|A+\Id\|<\varepsilon$. Since $\alpha$
is an isomorphism, we have also $\|A+(\Id\cap \Ga_k)\|<\varepsilon$
and thus there is $J'\in \Id\cap \Ga_k$ such that
$\|A-J'\|<\varepsilon$. Since $\|J-J'\|<2\varepsilon$,
we conclude that $\Id=\overline{\cap_{k\ge 1}(\Id\cap \Ga_k)}$.
\end{proof}
We can now give the proof of Theorem~\ref{theoremSimpleAF}.

\begin{proofof}{of Theorem~\ref{theoremSimpleAF}}
Let $(V,E)$ be a Bratteli diagram and let $\Ga$ be the
corresponding AF algebra.
We shall establish a one-to-one correspondence between
the ideals of $\Ga$ and the sets $W$ of vertices 
which are both directed \index{subject}{directed!set}  and
hereditary \index{subject}{hereditary set}

By Proposition~\ref{propositionDirectedHereditary},  $(V,E)$ is simple if and only if
there is no nontrivial directed and hereditary subset of $V$.
Thus the above correspondence will prove the theorem.

Let first $\Id$ be an ideal of $\Ga$. For each $p\ge 1$, the
set $\Id_p=\Id\cap \Ga_p$ is an ideal of $\Ga_p=\M_{(p,1)}\oplus\ldots\oplus \M_{(p,k)}$. Since every summand $\M_{(p,i)}$ is simple, the ideal
$\Id_p$ corresponds to a set $W_p$ of vertices $(p,i)$.
Let $W=\cup_{p\ge 1}W_p$.

Suppose that $v=(p,i)$ is in $W_p$ and that there is an edge
$(p,i)\to(p+1,j)$. Then, denoting $\alpha_p$ the injection of $\Ga_p$
into $\Ga_{p+1}$, we have
\begin{equation}
\Id\cap\M_{(p+1,i)}\supset\alpha_p(\M_{(p,i)})\cap\M_{(p+1,j)}\ne\emptyset.
\end{equation}
Hence $\Id$ contains $\M_{(p+1,j)}$ and thus $(p+1,j)$ is in $W_{p+1}$.
This shows that $W$ is directed. 

Next, let $(p,i)$ be a vertex and let $J=\{j\mid (p,i)\to(p+1,j)\}$.
Assume that $W$ contains all $(p+1,j)$ such that
$\in J$. Then
\begin{equation}
\alpha_p(\M_{(p,i)})\subset \sum_{j\in J}\M_{(p+1,j)}\subset \Id
\label{eqHereditary}
\end{equation}
which implies that $(p,i)\in W$. Thus $W$ is also hereditary.
This allows us to associate
to every ideal of $\Id$ a directed and hereditary set $W=\alpha(\Id)$ of vertices
of the diagram. By Lemma~\ref{lemmaIdealAFalgebra}, the ideal
$\Id$ can be recovered from $W$ and thus the map$\alpha$  is injective.

Conversely, let $W$ be a directed and hereditary subset of $W$.
Let $\Id_p$ be the ideal of $\Ga_p$ corresponding to the vertices
of $W$ at level $p$. Since $W$ is directed, the sequence $\Id_p$
is increasing. The closure of the union is an ideal $\Id$ of
$\Ga$. Since $W$ is hereditary, we have $\Id_p=\Ga_p\cap\Id$.
Indeed, if $\M_{(p,i)}$ is in $\Id$, there is a $q\ge p$
such that $\M_{(p,i)}\in \Ga_p$. All descendants of $(p,i)$
at level $q$ are then in $W_q$ and thus $(p,i)$ is in $W$
because $W$ is hereditary.
Thus $W=\alpha(\Id)$ and this completes the proof.
\end{proofof}

\begin{example}
The algebra $\C I+\Gk$ of Example~\ref{exampleK2} is not simple
since $\G_k$ is a proper ideal. Accordingly,
the Bratteli diagram of Figure~\ref{figureBratteli2}
is not simple. The vertices of the lower level form
a directed hereditary set of vertices. The restriction of the diagram
to this set of vertices is the partial Bratteli
diagram of Figure~\ref{figureBratteli2bis},
which represents the $C^*$-algebra $\Gk$.
\end{example}

\subsection{Dimension groups of AF algebras}
We   define the dimension group of an AF algebra as follows.
Let $\Ga=\overline{\cup_{n\ge 1} \Ga_n}$ be an AF algebra where
the sequence $(\Ga_n)_{n\geq 1}$ is defined by the sequence of 
$k_{n+1}\times k_n$-matrices $A_n$. The dimension group $K_0(\Ga)$
\index{symbols}{K@$K_0(\Ga)$} of $\Ga$
is the direct limit of the sequence
\begin{displaymath}
\Z^{k_1}\stackrel{A_1}{\rightarrow}\Z^{k_2}\stackrel{A_2}{\rightarrow}
\Z^{k_3}\stackrel{A_3}{\rightarrow}\cdots
\end{displaymath}
Thus, the dimension group of an AF algebra is the dimension
group of its Bratteli diagram (although the diagram
is not unique, the group is well defined, see below).
Note that if $\Ga=\M_{n_1}\oplus\ldots\oplus\M_{n_k}$ is a finite dimensional
AF algebra, one obtains $K_0(\Ga)=(\Z^k,\Z_+^k,(n_1,\ldots,n_k)^t)$.
Thus one can also write 
\begin{displaymath}
K_0(\Ga)=\lim_{\rightarrow} K_0(\Ga_n)
\end{displaymath}
with the connecting morphisms given by the matrices $A_n$.

It follows from Corollary~\ref{corollaryBratteli} that
the definition of the dimension group is independent
of the sequence $\Ga_n$ such that $\Ga=\overline{\cup\Ga_n}$.
Indeed, by Theorem \ref{theoremDGBratteli}, the dimension
groups of Bratteli diagrams euivalent by intertwining
are isomorphic.

\begin{example}\label{exampleFibo3}
Consider the Fibonacci algebra $\Ga$ of Example~\ref{exampleFibo2}.
The dimension group $K_0(\Ga)$ is the group $\Z[\alpha]$
where $\alpha=(1+\sqrt{5})/2$. Indeed, it is the direct limit
$\Z^2\stackrel{A}{\rightarrow}\Z^2
\stackrel{A}{\rightarrow}\Z^2\cdots$
with $A=\begin{bmatrix}1&1\\1&0\end{bmatrix}$
and thus the assertion results from Example~\ref{exampleGolden1}.
\end{example}
\begin{example}
The dimension group of the algebra $\C I\oplus\Gk$ (Example~\ref{exampleK2})
is the ordered group $\Z^2$ with positive cone
$\{(x,y)\mid y>0\}\cup\{(x,0)\mid x\ge 0\}$
(see Example~\ref{exampleTriangular}).
\end{example}
\begin{example}
The dimension group of the CAR algebra is $\Z[1/2]$ (see Example
\ref{exampleDyadic}).
\end{example}
The next example introduces  an ordered group not seen before.
\begin{example}
Consider the GICAR algebra (Example~\ref{exampleGICAR}).
We have $K_0(\Ga_n)=\Z^{n+1}$ with morphisms
\begin{displaymath}
\varphi_n(a_0,a_1,\ldots,a_n)=(a_0,a_0+a_1,\ldots,a_n).
\end{displaymath}
Let us represent the group $K_0(\Ga_n)$ as the set of polynomials
with integer coefficients
of degree at most $n$ though the map 
\begin{displaymath}
(a_0,a_1,\ldots,a_n)\mapsto a_0+a_1x+\ldots+a_nx^n.
\end{displaymath}
Since
\begin{displaymath}
(1+x)(a_0+a_1x+\ldots+a_nx^n)=a_0+(a_0+a_1)x+\ldots+(a_{n-1}+a_nx^n+a_nx^{n+1},
\end{displaymath}
the morphisms $\varphi_n$ are now replaced by the multiplication
by $1+x$. Thus, the dimension group of the GICAR algebra is
the group
\begin{displaymath}
G=\{(1+x)^{-n}p(x)\mid p\in\Z[x]\mbox{ of degree at most }n, n\ge 0\}.
\end{displaymath}
The positive cone $G_+$ corresponds to the polynomials $p$
such that $(1+x)^Np(x)>0$ for some $N\ge 1$.
One can show that a polynomial is of this form if and only
if $p(x)>0$ for $x\in ]0,\infty[$ 
(Exercise~\ref{exercisePositivePolynomials}).
\end{example}

A \emph{trace}\index{subject}{trace!in $C^*$-algebra} in a $C^*$-algebra
$\Ga$ is
a state $\tau$ such that $\tau(AB)=\tau(BA)$ for all $A,B\in\Ga$.
When $\Ga=\M_n$, there is a unique trace on $\Ga$ which is
\begin{displaymath}
\tau(A)=\frac{1}{n}\Tr(A)
\end{displaymath}
where $\Tr(A)$ is the usual trace of the matrix $A$.
More generally, when $\Ga=\M_{n_1}\oplus\ldots\oplus\M_{n_k}$
is a finite dimensional $C^*$-algebra, the traces  on $\Ga$
are the maps $\tau$ such that
\begin{equation}
\tau(A_1\oplus\ldots\oplus A_k)=\frac{t_1}{n_1}\Tr(A_1)+\ldots+\frac{t_k}{n_k}\Tr(A_k)\label{eqTraceFiniteDim}
\end{equation}
where $t_1,\ldots,t_k\ge 0$ are such that $\sum_{j=1}^kt_k=1$.

\begin{theorem}
The traces on an AF algebra $\Ga$ are in one-to-one correspondence
with the states on $K_0(\Ga)$.
\end{theorem}
\begin{proof}
Let first $\tau$ be a trace on the AF algebra $\Ga=\overline{\cup\Ga_n}$
where $\Ga_n=\M_{n_1}\oplus\ldots\oplus\M_{n_{k_n}}$. The
restriction of $\tau$ to $\Ga_n$ is a trace on $\Ga_n$ and thus it
is given by Equation~\eqref{eqTraceFiniteDim}. 
Set $t_j=\tau(I_j)$ where $I_j$ is the identity of the $\M_{n_j}$.
Let $\tau_n$
be the map on $K_0(\Ga_n)$  defined by
\begin{displaymath}
\tau_n(x_1,\ldots,x_{k_n})=t_1x_1+\ldots+t_{k_n}x_{k_n}
\end{displaymath}
Then $\tau_n$ is a state on $(\Z^{k_n},\Z_+^{k_n},\begin{bmatrix}n_1\\ \vdots
\\n_{k_n}\end{bmatrix})$.

Conversely, for every state $\sigma$ on $K_0(\Ga)$, we define
a state on $\Ga_n=\M_{n_1}\oplus\ldots\oplus\M_{n_k}$ by \eqref{eqTraceFiniteDim} with $k_i=\sigma(e_i)$ and $e_i$ the $i$-th basis vector of $\Z^{k_n}$. 
\end{proof}

\subsection{Elliott Theorem}
The following classification result is due to Elliott.

\begin{theorem}[Elliott]\label{theoremElliott}
Two unital AF algebras $\Ga$ and $\Gb$ are $*$-isomorphic if and only
if their dimension groups are isomorphic.
\end{theorem}
We first prove two lemmas which treat particular cases.

Let $\varphi$ be unital $*$-morphism from
$\Ga=\M_{n_1}\oplus\ldots\oplus\M_{n_k}$ into
$\Gb=\M_{m_1}\oplus\ldots\oplus\M_{m_\ell}$
as in Proposition~\ref{propSemiSimple}. We denote
by $\varphi_*$ the morphism from $K_0(\Ga)=(\Z^k,\Z_+^k,(n_1,\ldots,n_k)^t)$
 to $K_0(\Gb)=(\Z^\ell,\Z_+^\ell,(m_1,\ldots,m_\ell)^t)$
which corresponds to the multiplication by the matrix of partial
multiplicities $A$ of $\varphi$. Since $A$
has nonnegative coefficients the morphism $\varphi_*$
 is positive and since $A$
satisfies~\eqref{equationA}, it is unital. Thus
$\varphi_*$ is a unital morphism of ordered groups.

The next result shows how to recover morphisms
of finite dimensional AF algebras
from morphisms of their dimension groups.
\begin{lemma}\label{lemmaInduced}
Suppose that $\Ga,\Gb$ are finite dimensional $C*$-algebras
and that $\psi$ is a unital morphism of ordered groups from $K_0(\Ga)$ into $K_0(\Gb)$.
Then there is a unital $*$-morphism $\varphi:\Ga\rightarrow \Gb$ such that
$\varphi_*=\psi$ and $\varphi$ is unique up to unitary equivalence.
\end{lemma}
\begin{proof}
Set $\Ga=\M_{n_1}\oplus\ldots\oplus\M_{n_k}$
and $\Gb=\M_{m_1}\oplus\ldots\oplus\M_{m_\ell}$. Then 
$K_0(\Ga)=(\Z^k,\Z_+^k,(n_1,\ldots,n_k)^t)$
 and $K_0(\Gb)=(\Z^\ell,\Z_+^\ell,(m_1,\ldots,m_\ell)^t)$.
Let $A$ be the matrix of the morphism $\psi$, which is such that
$\psi(v)=Av$ for every $v\in\Z^k$. Since $\psi$ is positive, the coefficients
of $A$ are nonnegative. And since $\psi$ is unital,
 $\eqref{equationA}$ is satisfied. 
Thus the result follows from Proposition~\ref{propSemiSimple}.
\end{proof}
We now improve the last result by replacing $\Gb$ by an AF algebra.
For this, we introduce a notation. Let
$\Gb=\overline{\cup_{m\ge 1}\Gb_m}$ be an AF algebra 
with embeddings $\beta_m:\Gb_m\rightarrow \Gb$.
Then $\displaystyle{K_0(\Gb)=\lim_\rightarrow K_0(\Gb_m)}$ with connecting
morphisms $\beta_{m,n}$. We denote by $\beta_{m*}$
the natural morphism from $K_0(\Gb_m)$ into $K_0(\Gb)$
which corresponds to the connecting morphisms $\beta_{mn*}$.
Thus $\beta_{m*}(x)$ is, for $x\in \Gb_m$, the class of sequences
$(x_k)\in \prod_{k\ge 1}K_0(\Gb_k)$ such that $x_n=x$
and $x_{m+1}=\beta_{m+1,m*}(x_m)$ for all $m\ge n$.
\begin{lemma}\label{lemmaStar}
Let $\Ga$ be a finite dimensional $C^*$-algebra. Let
$\Gb=\overline{\cup_{m\ge 1}\Gb_m}$ be an AF algebra 
with the embeddings $\beta_m:\Gb_m\rightarrow \Gb$.
 Let $\psi:K_0(\Ga)\rightarrow K_0(\Gb)$
 be a unital morphism of ordered groups. Then there is
an integer $m$ and a $*$-morphism $\varphi$ from $\Ga$ into
$\Gb_m$ such that $\beta_{m*}\varphi_*=\psi$. Moreover, $\varphi$
is unique up to unitary equivalence.
\end{lemma}
\begin{proof}
Set $\Ga=\M_{n_1}\oplus\ldots\oplus\M_{n_k}$. Then, as we have seen,
$K_0(\Ga)$ is the unital ordered group $(\Z^k,\Z_+^k,(n_1,\ldots,n_k)^t)$.
Let $e_j\in \Z^k$ be the $j$-th basis vector of $\Z^k$. For each $j=1,\ldots,k$,
 since $K_0(\Gb)$ is the direct limit of the ordered groups $K_0(\Gb_m)$
with natural morphisms $\beta_{m*}$,
there is, by definition of the direct limit,
  an integer $m$ such that $\psi(e_j)$ is in $\beta_{m*}(K_0^+(\Gb_m))$. Taking
the maximum of these integers, we may assume that this holds
with the same $m$
for all $j=1,\ldots,k$. Set $\psi(e_j)=\beta_{m*}(v_j)$ for some $v_j\in K_0^+(\Gb_m)$.

We define a morphism $\rho:\Z^k\rightarrow K_0(\Gb_m)$ by $\rho(e_j)=v_j$
for $j=1,\ldots,k$. Then $\psi=\beta_{m*}\rho$
and we have the commutative diagram below
\begin{figure}[hbt]
\centering
\tikzset{edge/.style={->}}
\begin{tikzpicture}
\node(A)at(0,2){$K_0(\Ga)$};\node(Bm)at(2,2){$K_0(\Gb_m)$};
\node(B)at(2,0){$K_0(\Gb)$};

\draw[above](A)edge node{$\rho$}(Bm);
\draw[below](A)edge node{$\psi$}(B);
\draw[right](Bm)edge node{$\beta_{m*}$}(B);
\end{tikzpicture}
\end{figure}

The morphism $\rho$ is actually a morphism of unital ordered groups.
Indeed, $\rho$ is positive since $\rho(e_j)$ belongs to  $K_0^+(\Gb_m)$,
and $\rho$ is unital since 
\begin{eqnarray*}
\rho(n_1,\ldots,n_k)^t&=&\sum_{j=1}^kn_j\rho(e_j)=\sum_{j=1}^kn_j\psi(e_j)\\
&=&\psi(n_1,\ldots,n_k)^t=\1_{K_0(\Gb)}
\end{eqnarray*}
We now apply Lemma~\ref{lemmaInduced} to obtain a unital $*$-morphism
$\varphi:\Ga\rightarrow \Gb_m$ such that $\varphi_*=\rho$.

Assume that $\varphi':\Ga\rightarrow \Gb_{m'}$ is another map with the same properties.
Taking the maximum of $m,m'$,
we may assume that $m=m'$. Again by definition of the direct limit,
there is an integer $p\ge m$ such that $\beta_{p,m*}(\varphi_*)=\beta_{p,m*}(\varphi'_*)$.
We replace $\varphi,\varphi'$ by $\beta_{p,m}(\varphi )$ and
$\beta_{p,m}\varphi'$. Since $\beta_{p,m}(\varphi_*)=\beta_{p,m*}(\varphi_*)$,
applying Lemma~\ref{lemmaInduced}, we obtain that
$\varphi,\varphi'$ are unitarily equivalent.
\end{proof}

\begin{proofof}{of Theorem~\ref{theoremElliott}}
We will prove that given a unital isomorphism 
$\rho:K_0(\Ga)\rightarrow K_0(\Gb)$, there is a $*$-isomorphism
$\varphi:\Ga\rightarrow \Gb$ such that $\varphi_*=\rho$.

Let $\Ga=\lim\Ga_m$ and $\Gb=\lim\Gb_n$ 
with natural maps $\alpha_m$ and $\beta_n$ respectively. We will
build a commutative diagram with increasing sequences
$m_1<m_2<\ldots$ and $n_1<n_2<\ldots$ as represented below
so that $\varphi_*=\rho$ and $\psi_*=\rho^{-1}$.
 %\begin{equation}
%\begin{CD}
%\Ga_{m_1} @>\alpha_{m_1m_2}>> \Ga_{m_2}@>\alpha_{m_2m_3}>>\Ga_{m_3}@>\alpha_{m_%3m_4}
%>>\cdots@> >>\Ga\\
%@VV{\varphi_1}V               @VV{\varphi_2}V              @VV\varphi_3V
%@VV V       @VV{\varphi}V\\
%\Gb_{n_1} @>\beta_{n_1n_2}>> \Gb_{n_2}@>\beta_{n_2n_3}>>\Gb_{n_3}@>>\beta_{n_3n_4}
%>>\cdots@> >>\Gb
%\end{CD}\label{DiagramElliott}
%\end{equation}
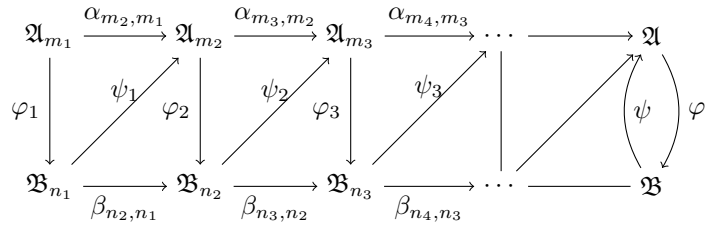
\begin{figure}[hbt]
\centering
\tikzset{edge/.style={->}}
\begin{tikzpicture}
\node(A1)at(0,2){$\Ga_{m_1}$};\node(A2)at(2,2){$\Ga_{m_2}$};
\node(A3)at(4,2){$\Ga_{m_3}$};\node(A4)at(6,2){$\cdots$};\node(A)at(8,2){$\Ga$};
\node(B1)at(0,0){$\Gb_{n_1}$};\node(B2)at(2,0){$\Gb_{n_2}$};
\node(B3)at(4,0){$\Gb_{n_3}$};\node(B4)at(6,0){$\cdots$};\node(B)at(8,0){$\Gb$};

\draw[above,->](A1)edge node{$\alpha_{m_2,m_1}$}(A2);
\draw[above,->](A2)edge node{$\alpha_{m_3,m_2}$}(A3);
\draw[above,->](A3)edge node{$\alpha_{m_4,m_3}$}(A4);\draw[->](A4)edge node{}(A);
\draw[left,->](A1)edge node{$\varphi_1$}(B1);
\draw[left,->](A2)edge node{$\varphi_2$}(B2);
\draw[left,->](A3)edge node{$\varphi_3$}(B3);\draw(A4)edge node{}(B4);
\draw[bend left,right,->](A)edge node{$\varphi$}(B);
\draw[below,->](B1)edge node{$\beta_{n_2,n_1}$}(B2);\draw[below,->](B2)edge node{$\beta_{n_3,n_2}$}(B3);
\draw[below,->](B3)edge node{$\beta_{n_4,n_3}$}(B4);\draw(B4)edge node{}(B);

\draw[above,->](B1)edge node{$\psi_1$}(A2);\draw[above,->](B2)edge node{$\psi_2$}(A3);
\draw[above,->](B3)edge node{$\psi_3$}(A4);
\draw[->](B4)edge node{}(A);
\draw[bend left,right,->](B)edge node{$\psi$}(A);
\end{tikzpicture}
\caption{The sequences $m_1,m_2\ldots$ and $n_1,n_2,\ldots$.}\label{diagramGene}
\end{figure}

We start with $m_1=1$. Apply Lemma~\ref{lemmaStar} to the
map $\rho\circ \alpha_{m_1*}:K_0(\Ga_{m_1})\rightarrow K_0(\Gb)$
to obtain an integer $n_1$ and a morphism 
$\varphi_1:\Ga_{m_1}\rightarrow \Gb_{n_1}$ such that $\beta_{n_1*}\varphi_{1*}=
\rho\alpha_{m_1*}$ (see Figure~\ref{figurephipsi} on the left).
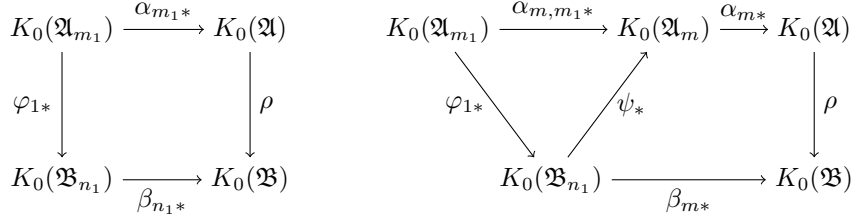
\begin{figure}[hbt]
\centering
\begin{tikzpicture}

\node(Am)at(0,2){$K_0(\Ga_{m_1})$};
\node(A)at(2.5,2){$K_0(\Ga)$};
\node(Bm)at(0,0){$K_0(\Gb_{n_1})$};\node(B)at(2.5,0){$K_0(\Gb)$};

\draw[above,->](Am)edge node{$\alpha_{m_1*}$}(A);
\draw[left,->](Am)edge node{$\varphi_{1*}$}(Bm);
\draw[right,->](A)edge node{$\rho$}(B);
\draw[below,->](Bm)edge node{$\beta_{n_1*}$}(B);

\node(Am1)at(5,2){$K_0(\Ga_{m_1})$};\node(Am2)at(8,2){$K_0(\Ga_m)$};
\node(A2)at(10,2){$K_0(\Ga)$};
\node(Bn)at(6.5,0){$K_0(\Gb_{n_1})$};\node(B2)at(10,0){$K_0(\Gb)$};

\draw[above,->](Am1)edge node{$\alpha_{m,m_1*}$}(Am2);
\draw[above,->](Am2)edge node{$\alpha_{m*}$}(A2);
\draw[left,->](Am1)edge node{$\varphi_{1*}$}(Bn);
\draw[right,->](Bn)edge node{$\psi_*$}(Am2);
\draw[below,->](Bn)edge node{$\beta_{m*}$}(B2);
\draw[right,->](A2)edge node{$\rho$}(B2);

\end{tikzpicture}
\caption{The construction of $\varphi_1$ and $\psi_1$}\label{figurephipsi}
\end{figure}

Now apply Lemma~\ref{lemmaStar} to the map
 $\rho^{-1}\circ \beta_{n_1*}:K_0(\Gb_{n_1})\rightarrow K_0(\Ga)$ to obtain a positive integer $m$ and a morphism $\psi$ of $\Gb_{n_1}$ into $\Ga_m$ such that
$\beta_{m*}\circ \psi_*=\rho^{-1} \circ \beta_{m_1*}$ (see Figure~\ref{figurephipsi}
on the right). Then, since 
\begin{displaymath}
\alpha_{m*}\circ \psi_* \circ \varphi_{1*}=\alpha_{m*}\circ \alpha_{m_1,m*}
\end{displaymath}
there is an integer $m_2\ge m$ such that
\begin{displaymath}
\alpha_{m_2,m*} \circ \psi_* \circ \varphi_{1*}=\alpha_{m_2 ,m*}\circ \alpha_{m_1m*}.
\end{displaymath}
By Lemma~\ref{lemmaInduced}, this implies that
$\alpha_{m_2m}\psi\varphi_1$ and $\alpha_{m_2m}$ are unitarily equivalent
and thus that there is a unitary element $U\in \Ga_{m_2}$ such that
$\alpha_{m_2m}=\Ad (U)\circ \alpha_{m_2m}\circ \psi \circ \varphi_1$. We set 
$\psi_1=\Ad (U)\circ \alpha_{m_2m}\circ \psi$. In this way, we have completed the
upper left triangle of Figure~\ref{diagramGene}.

We proceed in the same way to obtain the sequence $m_1<m_2<\ldots$,
the sequence $n_1<n_2<\ldots$ and the morphisms $\varphi_k,\psi_k$ as in the diagram
of Figure~\ref{diagramGene}.

Since $\Ga$ is the closure of $\cup_{m\ge 1}\Ga_m$
and similarly for $\Gb$, there is a unique 
continuous map $\varphi:\Ga\rightarrow\Gb$ such that its
restriction to $\Ga_k$ is equal to
$\varphi_k$. Moreover, $\varphi$ is a unital $*$-morphism.
Similarly, there is a unique continuous map $\psi:\Gb\rightarrow \Ga$
whose restriction to $\Gb_k$ is $\psi_k$. Since 
$\psi_k\circ \varphi_k$
is the embedding $\alpha_{m_{k+1} , m_k}$ of $\Ga_k$ into $\Ga_{k+1}$,
its restriction to $\Ga_k$ is the identity. This implies
that $\psi \circ \varphi=\id_\Ga$. Similarly $\varphi \circ \psi=\id_\Gb$.

Thus $\varphi,\psi$ are mutually inverse isomorphisms from $\Ga$
to $\Gb$. Finally, in view of the diagram of Figure~\ref{figurephipsi},
we have the commutative diagrams below and thus $\varphi_*=\rho$ as asserted.
\begin{displaymath}
\begin{CD}
K_0(\Ga_{m_k})@>{\alpha_{m_k*}}>>K_0(\Ga)\\
@VV\varphi_{k*}V           @VV\varphi_*V\\
K_0(\Gb_{n_k})@>{\beta_{n_k*}}>>K_0(\Gb)
\end{CD}
\quad
\mbox{ and }
\quad
\begin{CD}
K_0(\Ga_{m_k})@>{\alpha_{m_k*}}>>K_0(\Ga)\\
@VV\varphi_{k*}V           @VV\rho V\\
K_0(\Gb_{n_k})@>{\beta_{n_k*}}>>K_0(\Gb)
\end{CD}
\end{displaymath}
\end{proofof}
\section{Exercises}
\exosection{Section \ref{sectionC*algebras}}
\begin{exercise}\label{exerciseAltCstarIdentity}
  Show that the $C^*$-identity \eqref{eqStar2}
  implies \eqref{eqStar2bis} $\norm{A}=\norm{A^*}$.
  \end{exercise}
\begin{exercise}\label{exerciseNormUnique}
Show that the norm in a $C^*$-algebra is unique.
\end{exercise}
\begin{exercise}\label{exerciseNormPi(A)}
Show that if $\pi:\Ga\to \Gb$ is a $*$-morphism, then $\|\pi(A)\|\le\|A\|$
for all $A\in\Ga$.
\end{exercise}
\exosection{Section \ref{sectionApproximatelyFinite}}
\begin{exercise}\label{exerciseIdBratteli}
Show that if $\Ga=\overline{\cup_{n\ge 1}\Ga_n}$ and
$\Gb=\overline{\cup_{n\ge 1}\Gb_n}$ have the
same Bratteli diagram, they are isomorphic.
\end{exercise}
\begin{exercise}\label{exerciseUHFAlgebra}
A $C^*$-algebra is called \emph{uniformly hyperfinite} of UHF
\index{subject}{uniformly!hyperfinite $C^*$-algebra}%
\index{subject}{Cstar@$C^*$-algebra!UHF}%
if it is the increasing union of full matrix algebras $\M_{k_n}$
with $k_1|k_2|\ldots$. The
\emph{supernatural number}
\index{subject}{supernatural number}\index{subject}{number!supernatural}
associated to such $C^*$-algebra $\Ga$
is $\delta(\Ga)=\prod_{p}p^{n_p}$
where the product is over all prime numbers $p$ and $n_p\in\N\cup\infty$
is the supremum of the exponents of powers of
$p$ which divide $k_n$. Show that two UHF algebras are isomorphic
if and only if $\delta(\Ga)=\delta(\Gb)$.
\end{exercise}
\begin{exercise}\label{exerciseNotAF}
Show that $C([0,1])$ is not an AF algebra.
\end{exercise}
\begin{exercise}\label{exercisePositivePolynomials}
Show that for every $(a,b)\in\R^2$ with $b\ne 0$
there is an $N\ge 1$ such that $(1+x)^N(x^2-2ax+a^2+b^2)$
has positive coefficients. Conclude that
a polynomial $p$ is such that $(1+x)^Np(x)$ has positive coefficients
for some $N\ge 1$ if and only if $p(x)>0$ for every $x\in]0,\infty[$.
\end{exercise}
\section{Solutions}
\exosection{Section \ref{sectionC*algebras}}
\begin{solution}{\ref{exerciseAltCstarIdentity}}
  Since $\norm{A}^2=\norm{A^*A}\le\norm{A^*}\norm{A}$, we have
  $\norm{A}\le \norm{A^*}\le\norm{A^{**}}=\norm{A}$ and thus $\norm{A}=\norm{A^*}$. 

  Conversely, 
  \end{solution}
\begin{solution}{\ref{exerciseNormUnique}}
If $A$ is self-adjoint, then $\|A\|=\spr(A)$ by Proposition\ref{propositionNormal}. In the general case, we have
\begin{displaymath}
\|A\|^2=\|A^*A\|=\spr(A^*A)
\end{displaymath}
since $A^*A$ is self-adjoint.
\end{solution}
\begin{solution}{\ref{exerciseNormPi(A)}}
The  spectrum of $\pi(A)$ is contained in the spectrum of $A$. Thus,
if $A$ is self-adjoint, $\pi(A)$ is also
self-adjoint and we have by Proposition~\ref{propositionNormal}
\begin{displaymath}
\|\pi(A)\|=\spr(A)\le spr(A)=\|A\|.
\end{displaymath}
In the general case, 
\begin{displaymath}
\|\pi(A)\|^2=\|\pi(A)^*\pi(A)\|=\|\pi(A^*A)\|\le\|A^*A\|=\|A\|^2.
\end{displaymath}
\end{solution}
\exosection{Section \ref{sectionApproximatelyFinite}}
\begin{solution}{\ref{exerciseIdBratteli}}
Denote $\alpha_n$ the embedding of $\Ga_n$ into $\Ga_{n+1}$
and by $\beta_n$ the embedding of $\Gb_n$ into $\Gb_{n+1}$.
Each $\A_n$ is isomorphic to $\B_n$, as one can prove easily
by induction on $n$. Let $\varphi_n:\Ga_n\to \Gb_n$
be such isomorphism. Then $\varphi_{n+1}\circ\alpha_n$
and $\beta_n\circ\varphi_n$ are embeddings of $\Ga_n$
into $\Gb_{n+1}$ with the same partial multiplicities. 
By Proposition~\ref{propSemiSimple}, there is a unitary element $U_{n+1}$
of $\Gb_{n+1}$ such that 
\begin{displaymath}
\beta_n\circ\varphi_n=\Ad(U_{n+1})(\varphi_{n+1}\circ\alpha_n).
\end{displaymath}
Define recursively
$V_n\in \Gb_n$ and $\psi_n:\Ga_n\to\Gb_n$ by
 $\psi_1=\varphi_1$ and $V_1=I$ and 
\begin{displaymath}
V_{n+1}=\beta_n(V_n)U_{n+1}\mbox{ and } \psi_{n+1}=\Ad(V_{n+1})\varphi_{n+1}
\end{displaymath}
for $n\ge 1$. We then have
\begin{eqnarray*}
\beta_n\circ\psi_n&=&\beta_n\circ\Ad(V_n)\varphi_n=\Ad(\beta_n(V_n))\beta_n\circ\varphi_n\\
&=&\Ad(\beta_n(V_n))\Ad(U_{n+1})(\varphi_{n+1}\circ\alpha_n)\\
&=&\Ad(\beta_n(V_n)U_{n+1})\varphi_{n+1}\circ\alpha_n=\psi_{n+1}\circ\alpha_n.
\end{eqnarray*}
where in the first line, we have used the identity
$\beta\circ\Ad(V)\varphi=\Ad(\beta(V))\beta\circ\varphi$
for $*$-morphisms $\varphi:\Ga\to\Gb$ and $\varphi:\Gb\to\Gb'$.
We therefore have the commutative diagram
\begin{displaymath}
\begin{CD}
\Ga_1@>\alpha_1>>\Ga_2@>\alpha_2>>\Ga_3@>\alpha_3>>\cdots\\
@VV\psi_1    V @VV\psi_2V       @VV\psi_3V\\
\Gb_1@>\beta_1>>\Gb_2@>\beta_2>>\Gb_3@>\beta_3>>\cdots
\end{CD}
\end{displaymath}
showing that there is a map $\psi:\Ga\to \Gb$ which extends the maps $\psi_n$
and is a $*$-isomorphism from $\cup_{n\ge 1}\Ga_n$ to $\cup_{ge 1}\Gb_n$.
Since $\psi$ is an isometry, it extends to a $*$-isomorphism
from $\Ga$ onto $\Gb$.
\end{solution}
\begin{solution}{\ref{exerciseUHFAlgebra}}
This is a consequence of Elliott Theorem. In fact, two UHF algebras
$\Ga,\Gb$ are such that $\delta(\Ga)=\delta(\Gb)$ if and only if
the odometers associated to the Bratteli diagrams of $\Ga$ and $\Gb$
have the same supernatural number and thus have isomorphic
dimension groups (as seen in Exercise~\ref{exerciseSuperNatural}).
\end{solution}
\begin{solution}{\ref{exerciseNotAF}}
Assume that $C([0,1])=\overline{\cup\Ga_n}$ with
$\Ga_n$ finite dimensional. Let $\Gb$ be the subalgebra of
polynomials. Then $\B=\overline{\cup \Ga_n\cap\Gb)}$
with $\Ga_n\cap\Gb$ being a finite dimensional subalgebra
of $\Gb$. Bu the only finite dimensional subalgebra of $\Gb$
is formed by the constant functions, a contradiction.
\end{solution}
\begin{solution}{\ref{exercisePositivePolynomials}}
Set $\alpha=2a$ and $\beta=a^2+b^2$. We may assume that $\alpha>0$.
We have
\begin{eqnarray*}
(1+x)^N(x^2-\alpha x+\beta)&=&(x^2-\alpha x+\beta)\sum_{k=0}^N{N\choose k}x^k\\
&=&\sum_{k=0}^{N+2}(\beta{N\choose k}-\alpha{N\choose k-1}+{N\choose k-2})x^k\\
&=&\sum_{k=0}^{N+2}\frac{N!}{k!(N+2-k)!}a_{N,k}x^k
\end{eqnarray*}
where
\begin{eqnarray*}
a_{N,k}&=&\beta(N+2-k)(N+1-k)-\alpha k(N+2-k)+k(k-1)\\
&=&(1+\alpha+\beta)k^2-(2\beta+\alpha)Nk-(3\beta+2\alpha+1)k+\beta(N^2+3N+2)\\
&=&(1+\alpha+\beta)(k-\frac{\beta+\alpha/2}{1+\alpha+\beta}N)^2+
\beta(3N-3k+2)-(2\alpha+1)k\\
&&+(1+\alpha+\beta)^{-1}N^2((1+\alpha+\beta)\beta-(\beta+\alpha/2)^2).
\end{eqnarray*}
Since $1+\alpha+\beta>0$ and $k\le N$, we obtain
\begin{displaymath}
a_{N,k}\ge(1+\alpha+\beta)^{-1}N^2(\beta-\frac{\alpha^2}{4})+2\beta-(2\alpha+1)N
\end{displaymath}
which is positive for $N$ large enough, given that $4\beta>\alpha^2$.

The condition is clearly necessary. Conversely, a polynomial
satisfying the condition has positive leading coefficient
and no  positive real root. Thus
it factors as
\begin{displaymath}
p(x)=c\prod_i(x+\lambda_i)\prod_j(x^2-2a_jx+a_j^2+b_j^2)
\end{displaymath}
with $\lambda_i\ge 0$ and $b_j>0$. By what we have seen previously,
there are integers $N_j$ such that $(1+x)^{N_j}(x^2-2a_jx+a_j^2+b_j^2)$
has positive coefficients. Thus $N=\sum_jN_j$ is a solution.
\end{solution}
\section{Notes}
The reader is referred
for a more detailed presentation to the numerous monographs on the
subject, including \cite{Davidson1996}\index{names}{Davidson, Kenneth R.}
which we follow here
and also the edition of \cite{EilersOlesen2018}
\index{names}{Eilers, Soren}\index{names}{Olesen, Dorte}%
\index{names}{Pedersen, Gert K.} edited by Soren Eilers
and Dorte Oleson which we occasionally follow.

The proof that for every ideal $\cal J$ in a $C^*$-algebra
$\Ga$, the quotient algebra $\Ga/\cal J$ is a $C^*$-algebra
(Theorem~\ref{theoremQuotientAlgebra})
is in
\cite[Theorem I.5.4]{Davidson1996}.

For a proof of Proposition \ref{propSemiSimple},
see~\cite[Corollary III.2.1]{Davidson1996}.

Theorem~\ref{theoremCharacAF} appears in \cite[Theorem 2.2]{Bratteli1972}
and actually before in~\cite[Theorem 1.13]{Glimm1960}.
\index{names}{Glimm, James G.}%

The definition of the dimension group of an AF algebra given
here is not the usual one. The standard presentation
involves a development of the K-theory of AF algebras.
Indeed, $K_0$ is a functor which assigns an ordered abelian
group to each ring based on the structure of idempotents
in the matrix algebra over the ring. It occurs that
for AF algebras, this group coincides with the dimension
group of a Bratteli diagram defining the algebra,
a result due to ~\cite{Elliott1976}.
\index{names}{Elliott, George A.}%

The CAR algebra (Example~\ref{exampleCAR})
comes from quantum mechanics (see \cite{BratteliRobinson1987}).
\index{names}{Bratteli, Ola}\index{names}{Robinson, Derek W.} It is named for the 
\emph{Canonical Anticommutation
Relations algebra}. If $V$ is a vector space with a nonsingular
 symmetric bilinear form, the unital $*$-algebra generated
by the elements of $V$ subject to the relations 
\begin{eqnarray*}
fg+gf&=&\langle f,g\rangle\\
f^*=f
\end{eqnarray*}
for every $f,g\in V$ is the CAR algebra on $V$. It can be shown that
the CAR algebra as defined in Example~\ref{exampleCAR} is isomorphic the
CAR algebra on a separable Hilbert space (see \cite{Davidson1996}).

The name of the GICAR algebra 
(Example\ref{exampleGICAR}) stands for the \emph{Gauge Invariant}
CAR and is also called the \emph{current algebra} (see~\cite{Davidson1996}).
On positive polynomials (Exercise~\ref{exercisePositivePolynomials}), see 
\cite{Handelman1985}.

 Elliott's theorem (Theorem \ref{theoremElliott})
appeared in \cite{Elliott1976}.
The proof follows that of~\cite[Theorem IV.5.3]{Davidson1996}.

Let us finally mention perspectives on aspects which have
not been treated here.
We have been able to associate to every minimal Cantor system
a BV-representation (Theorem \ref{ch5:theo:BVmodel}) and thus
a $C^*$-algebra (by the results of this chapter)
in a way that conjugate systems define
isomorphic algebras. This has
been pursued using different constructions
for general shift spaces as an extension of the early
constructions of Krieger for shifts of finite type
(see \cite{CuntzKrieger1980}, \cite{Matsumoto1999}
and \cite{KriegerMatsumoto2002}). The link of this approach with substitutive
shifts was described in~\cite{CarlsenEilers2004}.
\index{names}{Krieger, Wolfgang}\index{names}{Cuntz, Joachim}%
\index{names}{Matsumoto, Kengo}\index{names}{Carlsen, Toke M.}%
\index{names}{Eilers, Soren}

The question of the computability of a Bratteli diagram for an AF-algebra
has been investigated by \cite{Mundici2020}. In particular,
the diagram is computable for AF-algebras $\Ga$ such the
dimension group $K^0(\Ga)$ is lattice ordered \citep{Mundici2020}.

\appendix
\chapter{Topological, metric and normed spaces}\label{appendixTopo}
We give in this appendix a list of the main notions and results
of general topology and elementary analysis used in the book.
This appendix (so as the following ones)
is intended to serve as a memento for notions
that some of the readers may have missed in part or forgotten long ago.
\section{Topological spaces}
A \emph{topological space}\index{subject}{topological!space}
\index{subject}{space!topological}%
 is a set $X$ with a family of
subsets, called \emph{open sets},\index{subject}{open!set}
containing $X$ and $\emptyset$ and closed by union and finite
intersection.
The complement of an open set is a \emph{closed set}\index{subject}{closed set}.
The sets that are both open and closed are called \emph{clopen sets}.
\index{subject}{clopen set}%

An element of a topological space is commonly called a \emph{point}
\index{subject}{point!in topological space}.

The topology for which all subsets are open is called the
\emph{discrete topology}.\index{subject}{discrete topology}
\index{subject}{topology!discrete}

The \emph{closure}\index{subject}{closure in topological space}
of a subset $S$ of a topological space $X$ is the
intersection of the closed sets containing $S$. 
It is also the smallest closed set containing $S$.

The \emph{interior}\index{subject}{interior of set}
 of a set $S$ is the union of all open sets
contained in $S$.

A \emph{neighborhood}\index{subject}{neighborhood} of a point
$x$ is a set $U$ containing an open set $V$ which contains $x$.
Thus $U$ is a neighborhood of $x$ if $x$ belongs to the
interior of $U$.

Any subset $S$ of a topological space inherits the topology of $X$
by considering as open sets the family of $S\cap U$ for $U\subset X$ open.
 This is called the \emph{induced topology}
\index{subject}{induced!topology}\index{subject}{topology!induced} on $S$ and $U$ is called
a \emph{subspace} of $S$.\index{subject}{subspace of topological space}

A topological space $X$ is a \emph{Hausdorff space}
\index{subject}{Hausdorff space}\index{subject}{space!Hausdorff}%
\index{names}{Hausdorff, Felix}%
 if for every distinct $x,y\in X$ there are disjoint open sets
$U,V$ such that $x\in U$ and $y\in V$.

An \emph{isolated point}\index{subject}{isolated point}
\index{subject}{point!isolated}%
in a topological space is a point $x$ such that $\{x\}$ is open.

A sequence $(x_n)$ of points in $X$ is \emph{convergent}\index{subject}{convergent sequence} to a limit $x\in X$ if  every open set $U$ containing $x$
contains all $x_n$ for $n$ large enough.
When $X$ is Hausdorff, the limit is unique. We denote $\lim x_n$
the limit.

A family $\F$ of subsets of a topological space $X$ is a \emph{basis}
\index{subject}{basis!of topology}\index{subject}{topology!basis of} of the topology if every open set
is a union of elements of $\F$. One also says that $\F$
\emph{generates}\index{subject}{generated topology}
\index{subject}{topology!generated}%
the topology of $X$. For  $\F$ to be a basis of some topology,
it is necessary and sufficient that it satisfies the two following conditions.
\begin{enumerate}
\item[(i)] Every point belongs to some $U\in\F$.
\item[(ii)] For every $U,V\in\F$, there is some $W\subset U\cap V$
such that $W\in\F$.
\end{enumerate}

Any direct product $\prod_{i\in I}X_i$
of topological spaces $X_i$ is a topological space
for the topology (called the \emph{product topology}\index{subject}{product topology}) with basis of open sets formed by the sets
$\prod_{i\in I}U_i$ with $U_i\subset X_i$ open sets such that
  $U_i=X_i$ for all but a finite number
of indices $i\in I$.

As a particular case, the set $Y^X$ of functions $f:X\to Y$
between topological spaces $X,Y$ is a subset
of the product $\prod_{x\in X}Y_x$. The topology induced by the product
 topology on $Y^X$ is the topology of \emph{pointwise convergence}.
\index{subject}{pointwise convergence} A sequence $(f_n)$
converges in this topology to a function $f$ 
if $\lim f_n(x)=f(x)$ for every $x\in X$.
 This notion is weaker than the notion of uniform convergence (see below).

A function $f:X\to Y$ between topological spaces $X,Y$ is \emph{continuous}
\index{subject}{continuous!function}\index{subject}{function!continuous}
 if $f^{-1}(U)$ is open
for every open set $U\subset Y$. A continuous invertible function
with a continuous inverse is called a \emph{homeomorphism}.
\index{subject}{homeomorphism}

A topological space is \emph{separable}\index{subject}{separable!space}
\index{subject}{space!separable}%
if it has a countable dense subset. For example, the real line
is separable with the usual topology but not separable with
the discrete topology. 

\section{Metric spaces}

A \emph{distance}\index{subject}{distance} on a set $X$
is a map $d:X\to \R_+$ such that 
\begin{enumerate}
\item[(i)] $d(x,y)=d(y,x)$, for every $x,y\in X$,
\item[(ii)] $d(x,y)=0$ if and only if $x=y$,
\item[(iii)] $d(x,z)\le d(x,y)+d(y,z)$ for every $x,y,z\in X$
(\emph{triangular inequality}).\index{subject}{triangular inequality}
\end{enumerate}
A set $X$ with a distance $d$ is called a \emph{metric space}.
\index{subject}{metric space}\index{subject}{space!metric}

For $x\in X$
and $\varepsilon>0$,
the \emph{open ball}\index{subject}{open!ball} centered at $x$
with radius $\varepsilon$ is
\begin{displaymath}
B(x,\varepsilon)=\{y\in X\mid d(x,y)<\varepsilon\}.
\end{displaymath}
\index{symbols}{B@$B(x,\varepsilon)$}
The topology with basis of open sets the open balls is
the \emph{topology defined}\index{subject}{topology!defined by a distance}
 by the distance $d$. In this topology, a sequence $(x_n)$
converges to a limit $x$ if for every $\varepsilon$ there is
an $N$ such that $d(x_n,x)\le\varepsilon$ for all $n\ge N$.

A subset $S$ of a metric space $X$ is \emph{bounded}
\index{subject}{bounded!set} if it is contained in some 
open ball or, equivalently if the set of distances $d(x,y)$
for $x,y\in S$ is bounded.

The \emph{diameter}
\index{subject}{diameter of set}
of a set $C$ in a metric space $X$, denoted $\diam(C)$,
\index{symbols}{diam@$\diam(C)$}%
is the maximal distance of two points in the set. 

A topological space is \emph{metrizable} \index{subject}{metrizable space}
\index{subject}{space!metrizable}%
if its topology can be defined by a distance. 
Not all topological spaces are metrizable,
even if all topological spaces  met in this book are metrizable spaces.
It is however useful to use the terminology of general
topological spaces, even when working with metric spaces.

A metric space is Hausdorff. Indeed, if $x\ne y$,
the balls $B(x,\varepsilon)$ and $B(y,\varepsilon)$ are
disjoint if $\varepsilon<d(x,y)/2$ by the triangular inequality.

In a metric space, many notions become easier to handle.
For example, a set $U$ in a metric space $X$ is closed
if and only if $U$ contains the limits of all converging sequences
of elements of $X$.

A sequence of functions $f_n:X\to Y$ between metric
spaces $X,Y$ is \emph{uniformly convergent}
\index{subject}{uniform!convergence} to $f:X\to Y$ if for every $\epsilon>0$ there
is an $N\ge 1$ such that for all $n\ge N$, one has
\begin{displaymath}
d(f_n(x),f(x))\le\varepsilon
\end{displaymath}
for every $x\in X$. Uniform convergence implies pointwise convergence
but the converse is false.

The limit of a uniformly convergent sequence of continuous function
is continuous.

A sequence $(x_n)$ in a metric space $(X,d)$ is a \emph{Cauchy sequence}
\index{subject}{Cauchy sequence}\index{names}{Cauchy, Augustin Louis}%
if for every $\varepsilon>0$ there is an integer $n\ge 1$
such that for all $i,j\ge n$ one has $d(x_i,x_j)\le\varepsilon$.
Every convergent sequence is a
Cauchy sequence but the converse is not true in general.

A metric space $(X,d)$ is \emph{complete}\index{subject}{complete!metric space}
\index{subject}{space!complete}%
if every Cauchy sequence converges. Every metric space $X$
can be embedded in a complete metric space $\hat{X}$,
called its \emph{completion}\index{subject}{completion!of metric space}
 such that
$X$ is dense in $\hat{X}$. The space $\hat{X}$ is build
as the set of classes of Cauchy sequences which are
equivalent, in the sense that $\lim d(x_n,y_n)=0$. The distance
on $\hat{X}$ is defined by continuity.
\section{Normed spaces}
When $X$ is a complex vector space, a 
\emph{norm}\index{subject}{norm!on vector space} on $X$ is a map $x\to ||x\|$
from $X$ to $\R_+$ such that
\begin{enumerate}
\item[(i)] $\|x\|=0$ if and only if $x=0$,
\item[(ii)] $\|\alpha x\|=|\alpha|\|x\|$ for $\alpha\in\C$ and $x\in X$,
and 
\item[(iii)] $\|x+y\|\le\|x\|+\|y\|$ for every $x,y\in X$.
\end{enumerate}
The space $X$ is called a \emph{normed space}.
\index{subject}{normed!space}\index{subject}{space!normed}%
If the separation property $\|x\|=0\Rightarrow x=0$ is missing,
one speaks of a \emph{seminorm}\index{subject}{seminorm}.
A norm defines a distance by $d(x,y)=\|x-y\|$. Indeed, using (iii),
we have
\begin{displaymath}
d(x,y)+d(y,z)=|x-y|+|y-z|\le|x-z|=d(x,z).
\end{displaymath}

The \emph{completion}\index{subject}{completion!of normed space}
 of a normed space $X$ is again a normed space.
The addition $x+y$ is defined as follows. If $(x_n)$ and $(y_n)$
are Cauchy sequences in $X$, then $(x_n+y_n)$ is a Cauchy
sequence and if $(x'_n)$ is equivalent to $(x_n)$, then
$(x'_n+y_n)$ is equivalent to $(x_n+y_n)$. Thus the addition
is well-defined on the classes.

A  complex algebra $\Ga$   is a \emph{normed algebra}
\index{subject}{normed!algebra} if
it has a norm $\|\ \|$ such that $\|AB\|\le\|A\|\|B\|$
for all $A,B\in\Ga$. The \emph{completion}\index{subject}{completion!of normed algebra} of a normed algebra $X$ is again a normed
algebra. One defines the product of two Cauchy sequences $(x_n)$
and $(y_n)$ as the sequence $(x_ny_n)$. The result is, as for the sum,
compatible with the equivalence.

If $X,Y$ are normed spaces, a linear map $f:X\to Y$ 
(also called a \emph{linear operator})\index{subject}{linear!operator}
\index{subject}{operator}%
 is continuous
if and only if it is \emph{bounded}\index{subject}{bounded!operator}
\index{subject}{operator!bounded}, that is, if
there is a constant $c\ge 0$ such that 
$\|f(x)\|/\|x\|\le c$ for all $x\ne 0$. The least such $c$ is called
the \emph{norm}\index{subject}{norm!of operator} of $f$
\emph{induced}\index{subject}{norm!induced} by the norms of $X,Y$,
denoted $\|f\|$.
This turns the space $L(X,Y)$ of bounded linear maps
from $X$ to $Y$ into a normed space and the algebra
$\Gb(X)$ of bounded linear maps from $X$ into $X$
into a normed algebra.

For every $p\ge 1$, the formula
\begin{displaymath}
\|x\|_p=(\sum_{i=1}^n |x_i|^p)^{1/p}
\end{displaymath}
\index{symbols}{x_p@$\lvert x\rvert_p$}%
defines a norm on $\R^n$, called the $L^p$-norm.
\index{subject}{L^p@$L^p$-norm}\index{subject}{norm!L@$L^p$}%
The $L^2$-norm is called the \emph{Euclidean norm}\index{subject}{Euclidean norm}
\index{subject}{norm!Euclidean} and the associated distance the
\emph{Euclidean distance}.

Additionally
\begin{displaymath}
\|x\|_\infty=\max\{|x_1|,|x_2|,\ldots,|x_n|\}
\end{displaymath}
\index{symbols}{x_infty@$\lvert x\rvert_\infty$}%
is called the $L^\infty$-norm.
\index{subject}{L^infty@$L^\infty$-norm}\index{subject}{norm!L@$L^\infty$}%
The matrix norm induced by the $L^\infty$-norm is
\begin{displaymath}
\|M\|_\infty=\max \|M_i\|_1
\end{displaymath}
where $M_i$ is the row of index $i$ of $M$.

An \emph{inner product}\index{subject}{inner product} on a 
vector space $X$ over $\C$ is 
 map $(x,y)\in X\times X\mapsto \langle x,y\rangle\in\C$ such that 
\begin{enumerate}
\item[(i)] $\langle x,y\rangle=\overline{\langle y,x\rangle}$ (conjugate symmetry)
\item[(ii)] $\langle x+y,z\rangle=\langle x,y\rangle+\langle y,z\rangle$
and  $\langle \alpha x,y\rangle=\alpha\langle x,y\rangle$
(sesquilinearity)
\item[(iii)] $\langle x,x\rangle\ge 0$
and ($\langle x,x\rangle=0$ only if $x=0$) (definite positivity)
\end{enumerate}
\index{symbols}{x@$\langle x,y\rangle$}%
for every $x,y,z\in X$ and $\alpha\in \C$. By (iii), we can
define a nonnegative real number $\|x\|$ by $\|x\|^2=\langle x,x\rangle$.
By the \emph{Cauchy-Schwarz inequality}\index{subject}{Cauchy-Schwarz inequality},
\index{names}{Schwarz, Hermann}\index{names}{Cauchy, Augustin Louis}%
we have
\begin{displaymath}
|\langle x,y\rangle|\le\|x\|\|y\|
\end{displaymath}
It follows from this inequality that $\|x\|$ satisfies
the triangular inequality and thus is a norm.

A \emph{Hilbert space}\index{subject}{Hilbert!space}
\index{names}{Hilbert, David} is 
a complex vector space with an inner product which 
is complete for the topology induced by the norm.

The space $\C^n$ is a Hilbert space for the inner product
\begin{displaymath}
\langle x,y\rangle=\sum_{i=1}^nx_i\overline{y_i}
\end{displaymath}
where $\bar{y}$ denotes the complex conjugate of $y$.
The corresponding norm $\|x\|_2$ is called the \emph{Hermitian norm}
\index{subject}{Hermitian norm}%
This extends to the sequences $x=(x_i)$ such that $\|x\|_2<\infty$
which form the Hilbert space $\ell_2(\C)$.

A \emph{Banach algebra}
\index{subject}{Banach algebra}\index{subject}{algebra!Banach}%
\index{names}{Banach, Stefan}%
is a complete normed algebra. The algebra $\M_n$ of $n\times n$-matrices
with complex elements is a Banach algebra (for the norm
induced by some norm on $\C^n$). 

If $\Ga$ is a Banach algebra,
the \emph{spectrum}\index{subject}{spectrum of element} of $A\in\Ga$
is
\begin{displaymath}
\sigma(A)=\{\lambda\in\C\mid \lambda I-A \mbox{ not invertible}\}
\end{displaymath}
\index{symbols}{sigma@$\sigma(A)$}%
It is a nonempty closed set. The \emph{spectral radius}\index{subject}{spectral radius} of $A\in\Ga$ is
\begin{displaymath}
\spr(A)=\sup_{\lambda\in \sigma(A)}|\lambda|
\end{displaymath}
The spectral radius satisfies $\spr(A)\le\|A\|$
and it is given by the formula
\begin{equation}
\spr(A)=\lim\|A^n\|^{1/n}.
\end{equation}

The set of bounded linear operators from a Hilbert space $H$ into itself,
denoted $\Gb(H)$, is a Banach algebra. The \emph{adjoint}\index{subject}{adjoint!of linear operator}\index{subject}{operator!adjoint}
of $T\in\Gb(H)$ is the unique operator $T^*$ such that
\begin{displaymath}
\langle Tx,y\rangle=\langle x,T^*y\rangle
\end{displaymath}
for every $x,y\in H$. For $H=\C^n$, the matrix $M$ of $T^*$
is the \emph{conjugate transpose}
\index{subject}{conjugate!transpose of matrix}\index{subject}{matrix!conjugate transpose}%
$M^*$ of $M$ defined by $M^*_{i,j}=\overline{M_{j,i}}$.

The operator $T\in\Gb(H)$ is said to be \emph{unitary}\index{subject}{unitary!operator}
if $TT^*=T^*T=I$, that is, if $T$ and $T^*$ are mutually inverse.

%%%%%%%%%%%%%%%%%%%%%%%%%%%%%%%%%
\section{Compact spaces}

A Hausdorff space $X$ is \emph{compact}\index{subject}{compact!space}
\index{subject}{space!compact} if for every family $(U_i)_{i\in I}$
of open sets with union $X$ there is a finite subfamily with
union $X$.

As a reformulation of the definition, say that a family
$(U_i)_{i\in I}$ of subsets of $X$ is an \emph{cover}
\index{subject}{cover of space}
of $X$ if  $X=\cup_{i\in I}U_i$ and an \emph{open cover}
\index{subject}{open!cover}\index{subject}{cover of space!open} if additionally
the $U_i$ are open. Then $X$ is compact if and only if
from every open cover one may extract a finite cover.

Here a few important properties of compact spaces.
\begin{enumerate}
\item A subset of a compact space is compact (for the induced topology)
if and only if it is closed.
\item A metric space is compact if and only if every sequence has
a converging subsequence.
\item The image of a compact space by a continuous map is compact.
\end{enumerate}
It follows easily from the second property that every compact metric space
is complete.

A subset of a topological space is \emph{relatively compact}
\index{subject}{relatively compact} if its closure is compact.
Thus any subset of a compact space is relatively compact.

\index{subject}{Tychonov Theorem}\index{subject}{Theorem!Tychonov}%
\index{names}{Tychonov, Andrei N.}%
\begin{theorem}[Tychonov]
Any product of compact spaces is compact.
\end{theorem}
This implies in particular that the set $A^\Z$
of sequences on a finite set $A$ is compact. This particular case
can of course be proved directly using \emph{K\"onig's Lemma}:
\index{subject}{K\"onig Lemma}\index{subject}{Lemma!K\"onig}%
\index{names}{K\"onig, D\'enes}%
every  sequence of infinite words on a finite alphabet has a converging
subsequence.

We state below (in one of its equivalent forms)
\emph{Baire Category Theorem}.\index{subject}{Baire Category Theorem}
\index{subject}{Theorem!Baire Category}%
\index{names}{Baire, Ren\'e-Louis}%
\begin{theorem}[Baire]\label{theoremBaireCategory}
In compact metric space, a countable intersection of dense
open sets is dense.
\end{theorem}
In particular, in a nonempty compact metric space, a countable
intersection of dense open sets is nonempty. Taking the complements,
the space is not  the countable union of closed sets
with empty interior.

A function $f:X\to Y$ between metric spaces $X,Y$ is said
to be \emph{uniformly continuous}
\index{subject}{uniformly!continuous}\index{subject}{uniform!continuity}
if for every $\varepsilon>0$ there exists $\delta>0$
such that for all $x,y\in X$, $d(x,y)<\varepsilon$ implies 
$d(f(x),f(y))<\delta$. Any continuous function on a compact
metric space is uniformly continuous (Heine-Cantor Theorem).
\index{subject}{Theorem!Heine-Cantor}\index{subject}{Heine-Cantor Theorem}
\index{names}{Heine, Eduard}\index{Cantor, Georg}%

The space $C(X,\R)$ of continuous real valued functions on a compact
metric space $X$ is a metric space for the norm
$\|f\|=\sup_{x\in X}|f(x)|$. The corresponding topology is the
topology of uniform convergence.

A linear map from a Banach space $X$ to a Banach space $Y$
is \emph{compact}\index{subject}{compact!operator}
\index{subject}{operator!compact} if the image of a bounded subset
is relatively compact. A compact operator is continuous.

\section{Connected spaces}
A topological space is \emph{disconnected}\index{subject}{disconnected!space}
\index{subject}{space!disconnected} if it is  a union
of disjoint nonempty open sets. Otherwise, the space is
\emph{connected}.\index{subject}{connected!space}
\index{subject}{space!connected}
The following are important properties of connected sets.
 1. A continuous image of a connected set is connected and
2. Any union of connected sets is connected.

Every topological space decomposes as a union of disjoint connected subsets,
called its \emph{connected components}.\index{subject}{connected!component}
 The connected component of $x$ is the union of
all connected sets containing $x$.

A topological space is \emph{totally disconnected}
\index{subject}{totally disconnected}
 if its connected components are
singletons. 
Any product of totally disconnected spaces is totally disconnected
 and every subspace of a totally disconnected space is totally
disconnected.

A topological space is \emph{zero-dimensional}
\index{subject}{zero-dimensional space}\index{subject}{space!zero-dimensional}%
 if it admits a basis consisting of
clopen sets.
A compact space is zero-dimensional if and only if it is totally disconnected

\section{Notes}

These notes follow mostly \cite{Willard2004} but a similar content
can be found in most textbooks on topology. Most authors 
(but all of them, see for example \cite{Rudin1987})
assume as we do a compact space to be Hausdorff.

\chapter{Measure and Integration}
\label{appendixMeasureIntegration}
We give below a short review of some definitions and concepts
concerning measure theory.
\section{Measures}
A family $\F$ of subsets of a set $X$
is a $\sigma$-\emph{algebra}\index{subject}{sigma-algebra@$\sigma$-algebra}
if it contains $X$ and is closed under complement and countable union.

Let $X$ be a topological space. The family of \emph{Borel subsets}
\index{subject}{Borel!subset}%
\index{names}{Borel, Emile}%
of $X$ is the closure of the family of open sets under
countable union  and complement.
It is thus the $\sigma$-algebra generated by the family of open sets.

A  \emph{measure} on a $\sigma$-algebra $\F$ is a map
$\mu:\F\to\R_+\cup\{\infty\}$ such that $\mu(\emptyset)=0$ and
which is countably additive, that is
such that $\mu(\cup_nX_n)=\sum_n\mu(X_n)$ for any
family $X_n$ of pairwise disjoint elements of $\F$
(note that we only consider measures with nonnegative values,
also called \emph{positive measures}). \index{subject}{positive!measure}\index{subject}{measure!positive}
The measure $\mu$ is \emph{finite}\index{subject}{finite!measure}%
\index{subject}{measure!finite} if $\mu(X)<\infty$.

The set $X$ is then called a 
\emph{measurable space}\index{subject}{measurable!space}
\index{subject}{space!measurable}%
and the elements of $\F$ are the \emph{measurable subsets}
\index{subject}{measurable!subset} of $X$.
We also say that $(X,\F,\mu)$ is a \emph{measure space}.\index{subject}{measure!space}\index{symbols}{XF@$(X,\F,\mu)$}

A property of points in a measurable space  is said to hold \emph{almost everywhere}
\index{subject}{almost everywhere (a.e.)} (often abbreviated a.e.)
if the set of points for which it is false has measure $0$.

A \emph{probability measure}
\index{subject}{probability measure}\index{subject}{measure!probability}%
%\index{subject}{Borel!probability measure}%
 on $X$ is a  measure $\mu$
such that $\mu(X)=1$.

Note the following important properties of measures.
\begin{enumerate}
\item A measure is \emph{monotone},\index{subject}{monotone measure}
  \index{subject}{measure!monotone}%
that is, if $U\subset V$, then $\mu(U)\le\mu(V)$.
\item A measure is \emph{subadditive}\index{subject}{subadditive!measure},
that is,
$\mu(\cup_n U_n)\le\sum_n\mu(U_n)$ for every family $X_n$ of measurable sets.
\end{enumerate}

%It is \emph{positive}
%\index{subject}{positive!measure}%
%\index{subject}{measure!positive}%
%if $\mu(U)\ge 0$ for every set $U\in\F$.

A  function $f:X\to Y$ between  measurable spaces $X,Y$
is \emph{measurable}\index{subject}{measurable!function}
\index{subject}{function!measurable}
if $f^{-1}(U)$ is a measurable set for every measurable subset $U$ in $Y$.

A Borel measure
\index{subject}{measure!Borel}%
\index{subject}{Borel!measure}%
on a topological space $X$ is a measure
on the family of Borel sets of $X$.

As an example, in a Hausdorff space, 
for every $x\in X$, the \emph{Dirac measure}
\index{subject}{Dirac measure}%
\index{subject}{measure!Dirac}%
\index{names}{Dirac, Paul}%
$\delta_x$, defined by $\delta_x(U)=1$ if $x\in U$
and $0$ otherwise, is a Borel probability measure on $X$.
\index{symbols}{delta@$\delta_x(U)$}%

As another example, the usual Borel measure on $\R^n$ is
called the \emph{Lebesgue measure}.\index{subject}{Lebesgue!measure}
\index{names}{Lebesgue, Henri}%

The Carath\'eodory Extension Theorem below
\index{subject}{Carath\'eodory Extension Theorem}%
\index{subject}{Theorem!Carath\'eodory Extension}%
\index{names}{Carath\'eodory, Constantin}%
shows that one may extend measures defined on a Boolean algebra
to one on a $\sigma$-algebra. Recall that a nonempty family $R$ of subsets
of a set $X$ is a \emph{Boolean algebra}\index{subject}{Boolean algebra}
\index{subject}{algebra!Boolean} if for every $U,V\in R$
one has $U\cup V,U\cap V,X\setminus U\in R$. A map
$\mu:R\to [0,1]$ is a probability measure on $R$
if whenever the union of disjoint sets $U_n\in R$
is in $R$, then $\mu(\cup_n U_n)=\sum_n\mu(U_n)$.
%The following result holds for any Boolean algebra
%of Borel sets instead of the algebra of clopen sets
%but we will not need this more general statement.
\begin{theorem}[Carath\'eodory]\label{theoremCaratheodoryExtension}
  Let $X$ be a topological space.
 Any probability measure $\mu$ on a Boolean algebra $R$ 
 has a unique extension $\mu^*$ to a  probability measure on the $\sigma$-algebra $\F$ generated
 by $R$.
 One has for every  set $U\in\F$,
\begin{equation}
  \mu^*(U)=\inf\sum_{n\ge 0}\mu(U_n)
\end{equation}
for $U_n\in R$ such that $U\subset \cup_{n\ge 0}U_n$.
\end{theorem}
This result holds in particular when $R$ is the Boolean algebra of clopen
spaces, $\F$ the $\sigma$-algebra of Borel sets and $X$ is a Cantor space.
It shows in particular that any Borel measure on a Cantor space can be approximated
by its value on clopen sets. Indeed, it implies that for every
Borel set $U$ and every $\varepsilon>0$, there is a clopen set $U_0$ such that
$\mu(U\Delta U_0)<\varepsilon$, where $U\Delta U_0=(U\cup U_0)\setminus (U\cap U_0)$
is the \emph{symmetric difference} \index{subject}{symmetric!difference}
of $U,U_0$.

\section{Integration}

Let $\mu$ be a positive measure on a measurable space $X$.
For $A\subset X$, we denote by $\charac_A$ the \emph{characteristic function}
\index{subject}{characteristic!function}%
\index{symbols}{chi@$\charac_A$}%
of $A$ defined by $\charac_A(x)=1$ if $x\in A$ and $0$ otherwise.
 A \emph{simple
function}\index{simple function}  is of the form
\begin{displaymath}
s=\sum_{i=1}^n\alpha_i\chi_{A_i}
\end{displaymath}
where the $A_i$ are measurable sets,
 and $\alpha_i\in\R$.
When $s$ is a simple function, we define
for a measurable set $U$,
\begin{displaymath}
\int_U s d\mu=\sum_{i=1}^n\alpha_i\mu(A_i\cap U).
\end{displaymath}
with the convention $0\cdot\infty=0$.

If $f:X\to [0,\infty]$ is measurable, the
\emph{integral} (or \emph{Lebesgue integral})
\index{subject}{integral}\index{subject}{Lebesgue!integral}%
 of $f$ over $U$ is defined as
\begin{displaymath}
\int_U f d\mu=\sup\int_U s d\mu
\end{displaymath}
the supremum being taken over the simple measurable functions
$s$ such that $0\le s\le f$. For $U=X$, we omit the subscript $U$.

For an arbitrary measurable function $f:X\to \R$, we define
\begin{displaymath}
\int f d\mu=\int f^+ d\mu -\int f^-d\mu
\end{displaymath}
\index{symbols}{int@$\int f d\mu$}%
where $f^+=\max\{f,0\}$ and $f^-=-\min\{f,0\}$.
A measurable function $f:X\to \R$  is \emph{integrable}
\index{subject}{integrable function} if
\begin{displaymath}
\int|f|d\mu<\infty.
\end{displaymath}
The set of integrable functions is denoted $L^1(X,\mu)$.

The following result is known as the \emph{Monotone Convergence Theorem}.
\begin{theorem}\label{theoremMonotoneConvergence}
Let $(f_n)$ be a sequence of measurable functions on $X$ such that
$0\le f_1(x)\le f_2(x)\le\cdots$  and
$f_n(x)\to f(x)$ for every $x\in X$. Then $f$ is measurable and
\begin{displaymath}
\int f_n d\mu\to\int f d\mu
\end{displaymath}
as $n\to\infty$.
\end{theorem}
It follows from Theorem~\ref{theoremMonotoneConvergence}
that if $f_n:X\to [0,\infty]$ is measurable for $n\ge 1$,
then
\begin{displaymath}
\int\sum_{n\ge 1}f_n d\mu=\sum_{n\ge 1}\int f_n d\mu
\end{displaymath}

The next result is known as the \emph{Dominated Convergence Theorem}.
\index{subject}{Dominated Convergence Theorem}
\begin{theorem}[Lebesgue]\label{theoremDominatedConvergence}
Suppose $(f_n)$ is a sequence of measurable functions from $X$
to $\R$ such that $f(x)=\lim_{n\to\infty}f_n(x)$ exists for every $x\in X$.
If there is a function $g\in L^1(X)$ such that $|f_n(x)|\le g(x)$
for all $x\in X$, then
\begin{displaymath}
\lim_{n\to\infty}\int f_n d\mu=\int f d\mu.
\end{displaymath}
\end{theorem}

Let $(X,\F,\mu)$ be a measure space and let $T:X\to X$ be
a measurable map. The map $\mu\circ T^{-1}$ is a measure
on $\F$. The following formula is called the \emph{change of variable formula}.
\index{subject}{change of variable formula}\index{subject}{Formula!change of variable}%
Its name relates to the usual change of variable
used to compute integrals of real functions of one variable.
We have, for every real-valued measurable function $f$ on $X$,
\begin{equation}
\int (f\circ T) d\mu=\int fd(\mu\circ T^{-1}).\label{eqChangeVariables}
\end{equation}
 Formula
\eqref{eqChangeVariables} is easy to verify
in the case of a simple function $s=\sum_{i=1}^n\alpha_i\chi_{A_i}$
Indeed,  we have
$\int (s\circ T) d\mu=\sum_{i=1}^n\alpha_i\mu(T^{-1}A_i)=\int s d(\mu\circ T^{-1})$.

Every Borel measure $\mu$ on a topological space $X$ defines
a bounded linear map from $C(X,\R)$ to $\R$ by $f\mapsto \int fd\mu$.
By the \emph{Riesz representation Theorem}
\index{subject}{Riesz!Representation Theorem}%
\index{subject}{Theorem!Riesz Representation}%
\index{names}{Riesz, Frigyes}%
the converse is true when $X$ is a compact metric space.
Thus the space of measures on a compact metric space $X$
can be identified with the (topological) dual of the space
$C(X,\R)$.

Let $\mu,\nu$ be two  Borel probability  measures on $X$.
Then $\nu$ is \emph{absolutely continuous}
\index{subject}{absolutely continuous measure}
\index{subject}{measure!absolutely continuous} with respect to $\mu$,
denoted $\nu\ll\mu$, if for every Borel set $U\subset X$
such that $\mu(U)=0$ one has $\nu(U)=0$.

\index{subject}{Radon-Nikodym Theorem}\index{subject}{Theorem!Radon-Nikodym}%
\index{names}{Radon, Johann}\index{names}{Nikodym, Otto}%
\begin{theorem}[Radon-Nikodym]
If $\nu\ll\mu$, there is a nonnegative $\mu$-integrable  function $f$ such that 
$\nu(U)=\int_U f d\mu$ for every measurable set $U\subset X$.
\end{theorem}

%As we have seen,
%every Borel probability measure defines a linear form on the space
%$C(X,\R)$  by $f\to \int f d\mu$. This allows to
%consider the set of Borel probability measure as 
%a subspace of the dual space of $C(X,\R)$.

Two measures $\mu,\nu$ on $X$ are \emph{mutually singular}
\index{subject}{mutually singular} if there is a measurable set $U$
such that $\mu(U)=\nu(X\setminus U)=0$.

A subset $S$ of a real vector space is \emph{convex}\index{subject}{convex!set}
if for every $x,y\in S$ and $t\in [0,1]$, one has
$tx+(1-t)y\in S$.  The \emph{convex hull} of $S$ is
the intersection of all convex sets containing $S$.\index{subject}{convex!hull}
The convex hull of a finite set of points in $\R^d$
is called a \emph{simplex}\index{subject}{simplex}.

The set ${\cal M}(X)$ of Borel  probability measures on 
 a compact metric
space $X$ is a convex subspace of the dual space of
$C(X,\R)$.  The \emph{weak-star topology}
\index{subject}{weak-star topology} on ${\cal M}(X)$
is  the topology for which
a sequence $(\mu_n)$ converges to $\mu\in{\cal M}(X)$ if
for all $f\in C(X,\R)$
\begin{displaymath}
\int f d\mu_n\to\int fd\mu.
\end{displaymath}
\begin{theorem}[Banach-Alaoglu]\label{theoremBanachAlaoglu}
\index{names}{Banach, Stefan}\index{names}{Alaoglu, Leonidas}%
For any compact metric space $X$, the space $\M(X)$
is metrizable and compact for the weak-star topology.
\end{theorem}
\index{subject}{Banach-Alaoglu Theorem}\index{subject}{Theorem!Banach-Alaoglu}%

An \emph{extreme point}\index{subject}{extreme point} of a convex set
$K$ is a point which does not belong to any open line segment in $K$.
The closed convex hull of a set $K$ is the closure
of the intersection of all closed convex subsets containing $K$.
\begin{theorem}[Krein-Milman]
Every compact convex set in a normed space is the closed
 convex hull of its extreme points.
\end{theorem}

\section{Notes}
We have mostly followed the classical \cite{Rudin1987} for this summary of notions
on measures and integration.
\index{names}{Rudin, Walter}%
We restricted the presentation to positive measures to simplify
the picture.
For the Carath\'eodory Extension Theorem
(also known as Kolmogorov extension Theorem), 
\index{names}{Kolmogorov, Andrei}
see \cite{Halmos1974}.
On the weak-star topology (or weak*-topology)
and the Banach-Alaoglu Theorem, see \cite{Rudin1991}.
The Krein-Milman Theorem is proved in~\cite[Theorem 3.23]{Rudin1991}.

\chapter{Algebraic Number Theory}\label{appendixAlgebraicNumberTheory}
In this appendix, we give an introduction to
the basic concepts and results of algebraic number theory.
We assume the reader to know the basic concepts of algebra,
as the notion of \emph{ring},\index{subject}{ring}
of \emph{ideal} \index{subject}{ideal!of ring}%
in a ring and of a \emph{principal ring}.\index{subject}{ring!principal}
\section{Algebraic numbers}
An \emph{algebraic number}
\index{subject}{algebraic!number}\index{subject}{number!algebraic}%
is a complex number $x$ solution of an equation
\begin{equation}
x^n+a_{n-1}x^n+\ldots+a_1x+a_0=0\label{eqAlgebraic}
\end{equation}
with coefficients in $\Q$. It is an \emph{algebraic integer}
\index{subject}{algebraic!integer}\index{subject}{integer!algebraic}%
if the coefficients are in $\Z$.
We denote by $\Q[x]$ (resp. $\Z[x]$) the subring generated by $\Q$
(resp. $\Z$) and $x$.
\index{symbols}{Q@$\Q[x]$}\index{symbols}{Z@$\Z[x]$}

\begin{theorem}\label{theoremAlgebraicModule}
For  $x\in\C$, the following conditions are equivalent.
\item[\rm(i)] $x$ is an algebraic number (resp. and algebraic integer).
\item[\rm(ii)] The $\Q$-vector space $\Q[x]$  (resp. the 
$\Z$-module $\Z[x]$) is finitely generated.
\end{theorem}

As a consequence, all elements of $\Q[x]$ (resp. $\Z[x)$)
are algebraic numbers (resp. algebraic integers). Moreover,
$\Q[x]$ is a field.

Let $K$ be a field containing $\Q$. It is called 
 an \emph{algebraic extension}
\index{subject}{algebraic!extension}\index{subject}{extension!algebraic}%
 of $\Q$, or also a \emph{number field},\index{subject}{number!field}
if all its elements are algebraic over $\Q$.

When $x$ is an algebraic number, the set of polynomials $p(X)$
such that $p(x)=0$ is a nonzero ideal of the ring $\Q[X]$.
Since $\Q[X]$ is a principal ring, this ideal is generated
by a unique polynomial of the form \eqref{eqAlgebraic}
(that is, with leading coefficient $1$). This polynomial
is called the \emph{minimal polynomial}\index{subject}{minimal!polynomial}
of $x$.

The  \emph{degree} \index{subject}{degree of extension}  over
$\Q$ of an extension $K$, denoted $[K:\Q]$
\index{symbols}{K@$[K:\Q]$}%
 is the dimension of the $\Q$-vector space $K$.
By Theorem \ref{theoremAlgebraicModule} every extension of finite
degree is algebraic (the converse is not true).

\section{Quadratic fields}
A \emph{quadratic field}\index{number field!quadratic}%
\index{subject}{quadratic!number field} is an extension of degree $2$ of $\Q$.
Every quadratic field is of the form $\Q[\sqrt{d}]$ where
$d$ is an integer without square factor.

\begin{theorem}
Let $K=\Q(\sqrt{d})$ be a quadratic field.
\begin{enumerate}
\item If $d\equiv 1,3 \bmod 4$ the ring of integers of $K$
is $\Z+\sqrt{d}\Z$.
\item If $d\equiv 1 \bmod 4$, the ring of integers of $K$ is
$\Z+\lambda \Z$ with $\lambda=(1+\sqrt{d})/2$.
\end{enumerate}
\end{theorem}

\paragraph{Norm and trace}
Let $K$ be a number field. For $x\in K$, 
the multiplication by $x$ in $K$ is a $\Q$-linear map $\rho(x)$.
The \emph{norm}\index{subject}{norm!in number field}
of $x$, denoted $N(x)$,\index{symbols}{N@$N(x)$} is the determinant of $\rho(x)$.
Its \emph{trace},\index{subject}{trace!in number field}
denoted $\Tr(x)$\index{symbols}{Tr@$\Tr(x)$} is the trace of $\rho(x)$.

If $x$ is an algebraic integer, then $N(x)$ and $\Tr(x)$
are integers.

For example, if $K=\Q(\sqrt{d}$) and $x=a+b\sqrt{d}$, then $N(x)=a^2-b^2d$
and $\Tr(x)=2a$.

\paragraph{Discriminant}
Let $K$ be a number field of degree $n$ over $\Q$. For 
$x_1,x_2,\ldots,x_n\in K$, the \emph{discriminant}
\index{subject}{discriminant!of basis}
of $(x_1,x_2,\ldots,x_n)$ is
\begin{displaymath}
D(x_1,x_2,\ldots,x_n)=\det(\Tr(x_ix_j)).
\end{displaymath}
\index{symbols}{D@$D(x_1,x_2,\ldots,x_n)$}%
One has $D(x_1,\ldots,x_n)=0$ if and only if $x_1,x_2,\ldots,x_n$ is a basis
of $K$.

If $y_1,y_2,\ldots,y_n$ is such that $y_i=\sum_{j=1}^na_{ij}x_j$, then
$D(y_1,y_2,\ldots,y_n)=\det(a_{ij})^2D(x_1,\ldots,x_n)$.
Thus the set of discriminants of the bases of the $\Z$-module
of algebraic integers of $K$
is an ideal of $\Z$. The generator of this ideal is called
the \emph{discriminant}\index{subject}{discriminant!of number field}
 of $K$.

For example, the discriminant of $\Q(\sqrt{d})$ is $4d$
if $d\equiv 2,3 \bmod 4$ and $d$ if $d\equiv 1 \bmod 4$.

\section{Classes of ideals}
Let $K$ be a number field and let $A$ be its ring of algebraic integers.
 Two ideals $I,J$ of  $A$ 
 are \emph{equivalent}\index{subject}{equivalent!ideals}
if there are nonzero $\alpha,\beta\in A$ such that $\alpha I=\beta J$. 

The equivalence classes of ideals form a group with respect
to the product, called
the \emph{class group}\index{subject}{group!of ideal classes}
of $F$. The neutral element is the class of the ideal $I=A$,
which can be shown to be formed of the principal ideals of $F$.
Thus, when $A$ is a principal ring, the class group has only one element.

The ring of integers of a number field may fail to be principal
and thus to have unique factorization. For example, in $\Q[\sqrt{-5}]$,
we have
\begin{displaymath}
(1+\sqrt{5})(1-\sqrt{5})=2\cdot 3
\end{displaymath}
although $1+\sqrt{5}$ has no nontrivial divisor. However,
one has the following result.
\begin{theorem}[Dirichlet]
For every number field, the class group is finite.
\end{theorem}
\index{subject}{Dirichlet!class group Theorem}
\index{names}{Dirichlet, Johann}
\paragraph{Units}

A \emph{unit}\index{subject}{unit!in number field} in a number
field $K$ is an invertible element of the ring $A$ of integers of $K$.
The set of units forms a multiplicative group, called
the \emph{group of units}\index{subject}{group!of units} of $A$.

An algebraic integer in a number field $K$ is a unit if and
only if $N(x)=\pm 1$.

For example, the group of units in a real quadratic field $\Q[\sqrt{d}]$
with $d\ge 2$, is formed of the $a+b\sqrt{b}$ with 
$a,b\in\Q$ solution of
\begin{equation}
a^2-db^2=\pm 1, \label{eqPell}
\end{equation}
which is known as \emph{Pell's equation}.\index{subject}{Pell's equation}
\index{subject}{equation!Pell}%

The following result is known as \emph{Dirichlet Unit Theorem}.
\index{subject}{Dirichlet!Unit Theorem}%
\index{subject}{Theorem!Dirichlet Unit}%
\begin{theorem}[Dirichlet]
The group of units of a number field $K$ is of the form
$\Z^r\times G$ where $r\ge 0$ and $G$ is a finite
cyclic group formed by the roots of unity contained in $K$.
\end{theorem}
The integer $r$ is $r=r_1+r_2-1$ where $r_1$ is the number 
of real embeddings of $K$ in $\C$ and $2r_2$ the number of complex embeddings.
One has $n=r_1+2r_2$ where $n$ is the degree of $K$ over $\Q$.

For example, in a real quadratic field $\Q[\sqrt{d}]$, the group of units is
$\Z\times\Z/2\Z$. Every unit is of the form $\pm u^n$
where $u$ is a \emph{fundamental unit}\index{subject}{fundamental!unit}.

For $d=2$, the positive fundamental unit is $1+\sqrt{2}$.
For $d=5$, it is $(1+\sqrt{5})/2$.

As another example, if $K=\Q[i]$, one has $r=0$. The group of roots of unity is of
order $4$ and is generated by $-i$. The corresponding ring
of integers is called the ring of \emph{Gaussian integers}.
\index{subject}{Gaussian integers}\index{integer!Gaussian}%

\section{Continued fractions}

Every irrational real number $\alpha>0$ has a
 unique \emph{continued fraction}
\index{subject}{continued fraction} expansion of $\alpha$
\begin{displaymath}
\alpha=a_0+\cfrac{1}{a_1+\cfrac{1}{a_2+\cfrac{1}{\ddots}}}
\end{displaymath}
where $a_0,a_1,\ldots$ are integers with $a_0\ge 0$ and $a_n>0$ for $n\ge 1$. We denote $\alpha=[a_0;a_1,\ldots]$.
\index{symbols}{a@$[a_0;a_1,\ldots]$}%
The integers $a_i$ are the \emph{coefficients}\index{subject}{coefficient!of continued fraction}\index{subject}{continued fraction!coefficient}
of the continued fraction.

\begin{theorem}[Lagrange]
A real number $\alpha>0$ is quadratic
if and only if its continued fraction expansion is eventually periodic .
\end{theorem}
\index{subject}{Lagrange Theorem}\index{subject}{Theorem!Lagrange}%
\index{names}{Lagrange, Joseph-Louis}%
For example, we have $[0;1,1,1,\ldots]=\frac{\sqrt{5}-1}{2}$
and $[0;2,2,\ldots]=\sqrt{2}-1$.
\paragraph{Notes}
For this brief introduction to algebraic number theory, we have
followed \cite{HardyWright2008}
\index{names}{Hardy, Godfrey}\index{names}{Wright, Edward M.}%
 and \cite{Samuel1970}
\index{names}{Samuel, Pierre}%
to which the reader is referred for a complete presentation.

\chapter{Groups, graphs and algebras}
\label{appendixGroups}
We recall in this appendix some of the basic notions of algebra
concerning free groups, fundamental graphs and simple algebras. We begin with some general notions on groups.
\section{Groups}
A \emph{group}\index{subject}{group}
is a set $G$ with an associative operation, a neutral element $1_G$ and an inverse $g^{-1}$
for every element $g\in G$. A \emph{subgroup}\index{subject}{subgroup} of a group $G$
is a subset $H$ closed for the operation and containing the inverses of its elements.
A subgroup $H$ of a group $G$ is \emph{normal}\index{subject}{normal!subgroup}
\index{subject}{subgroup!normal} if $gHg^{-1}=H$ for every $g\in G$.

A \emph{morphism}\index{subject}{morphism!of groups}\index{subject}{group!morphism}
from a group $G$ to a group $H$ is a map $\varphi:G\to H$ such that $\varphi(g)\varphi(h)=\varphi(gh)$
for every $g,h\in G$. The \emph{kernel}\index{subject}{kernel!of morphism}
\index{subject}{morphism!kernel} of a morphism $\varphi:G\to H$ is the set $\varphi^{-1}(1_H)$.
It is a normal subgroup of $G$.

Let $H$ be a subgroup of a group $G$. A \emph{right coset}
\index{subject}{right!coset}\index{subject}{coset!right}\index{subject}{subgroup!coset of}%
of $H$ in $G$ is a set
of the form $Hg=\{hg\mid h\in H\}$ for some $g\in G$. Two
cosets are equal or disjoint. 
The \emph{index}\index{subject}{index of subgroup} of $H$
in $G$, denoted $[G:H]$,\index{symbols}{G@$[G:H]$}
is the number of distinct cosets of $H$ in $G$.

The \emph{free abelian group}\index{subject}{free!abelian group} on a set $A$, denoted $\Z(A)$, is the additive
group formed of linear combinations $\sum_{a\in A}n_aa$ with $n_a\in\Z$.
When $A$ is finite with $n$ elements, it is isomorphic with $\Z^n$.

The following statement
is known as the \emph{Fundamental Theorem of abelian groups}.
\index{subject}{Fundamental Theorem!of abelian groups}%
A group is \emph{cyclic}\index{subject}{cyclic group}\index{subject}{group!cyclic}
A finite cyclic group is \emph{primary}\index{subject}{primary cyclic group}
\index{subject}{cyclic group!primary}
if its order is a power of a prime.
\begin{theorem}
Every finitely generated abelian group, is in a unique way
a direct product of primary cyclic groups and infinite
cyclic groups.
\end{theorem}
Thus every finitely generated abelian group $G$  can be written uniquely
as 
\begin{displaymath}
G=\Z^n\times\Z/q_1\Z\times\cdots\times\Z/q_m\Z
\end{displaymath}
with $n,m\ge 0$ and all $q_i$ powers of primes.

\section{Free groups}
 Recall that the \emph{free group}
\index{subject}{free!group}\index{subject}{group!free}%
on the alphabet $A$ is formed of the words on $A\cup A^{-1}$ which are
reduced (that is, contain no $aa^{-1}$ or $a^{-1}a$ for $a\in A$). The
product of two reduced words $u,v$ is the unique reduced word
obtained by reduction of $uv$. 

We denote by $F(A)$\index{symbols}{F@$F(A)$} the free group on $A$.
For $U\subset F(A)$, we denote by $\langle U\rangle$ the subgroup generated by $U$.

The following result is known as the Nielsen-Schreier Theorem,
\index{subject}{Nielsen-Schreier Theorem}\index{subject}{Theorem!Nielsen-Schreier}%
\index{names}{Nielsen, Jacob}\index{names}{Schreier, Otto}%
\begin{theorem}[Nielsen-Schreier]
 Every subgroup of the free group is free. 
\end{theorem}
Thus, every subgroup $H$ of the free
group has a generating set $U$, called a \emph{basis} of $H$,
\index{subject}{basis!of subgroup} such that $H$ is isomorphic to the free group on $U$. Two basis of a subgroup have the same number of elements,
called the \emph{rank}\index{subject}{rank!of subgroup}
of the subgroup.

Let $H$ be a subgroup of the free group $F(A)$.
 Let $Q$ be a set of reduced words on $A$
which is a prefix-closed set of representatives of the right cosets $Hg$
of $H$.
Such a set is traditionally called a \emph{Schreier transversal}
\index{subject}{Schreier!transversal} for $H$. Let
\begin{equation}
U=\{paq^{-1}\mid a\in A, p,q\in Q, pa\not\in Q, pa\in Hq\}.\label{SchreierBasis}
\end{equation}
Each word $x$ of $U$ has a unique factorization $paq^{-1}$ with $p,q\in Q$
and $a\in A$. The letter $a$ is called the \emph{central part}
\index{subject}{central part of word} of $x$.
The set $U$
is a  basis of $H$, called the \emph{Schreier basis}
\index{subject}{Schreier!basis}%
relative to $Q$. It is clear that $U$ is closed by taking
inverses. Indeed, if $x=paq^{-1}\in U$, then  $x^{-1}=qa^{-1}p^{-1}$.
We cannot have $qa^{-1}\in Q$ since otherwise $p\in Hqa^{-1}$ implies
$p=qa^{-1}$ by uniqueness of the coset representative and finally $pa\in Q$.
It generates $H$ as a monoid because if $x=a_1a_2\cdots a_m\in H$
with $a_i\in A$,
then $x=(a_1p_1^{-1})(p_1a_2p_2^{-1})\cdots(p_{m-1}a_m)$ with $a_1\cdots a_k\in Hp_k$
for $1\le k\le m-1$ is a factorization of $x$ in elements of $X\cup \{1\}$.
Finally, if a product $x_1x_2\cdots x_m$ of elements of $U$ 
is equal to $1$, then $x_kx_{k+1}=1$ for some index $k$
since the central part $a$ never
cancels in a product of two elements of $U$.

If $A$ has $k$ elements and $H$ is a subgroup of finite index $n$
in the free group on $A$, the rank $r$ of $H$ satisfies the equality
\begin{displaymath}
r-1=n(k-1)
\end{displaymath}
called \emph{Schreier's Formula}.\index{subject}{Schreier!Formula}
\index{subject}{Formula!Schreier}

A group $G$ is called \emph{residually finite}\index{subject}{residually finite group}\index{subject}{group!residually finite}  if for every element
$g\ne 1$ of $G$, there is a morphism $\varphi$ from $G$
onto a finite group such that $\varphi(g)\ne 1$.

The following can be proved directly without much difficulty
by associating to $g\ne 1$ a map $\varphi:F(A)$ into
the symmetric group on $\{1,2,\ldots,n\}$ such that $\varphi(g)\ne 1$. 
\begin{theorem}\label{theoremFGResiduallyFinite}
A finitely generated free group is residually finite.
\end{theorem}

A group $G$ is said to be \emph{Hopfian}\index{subject}{Hopfian group}
\index{subject}{group!Hopfian} if any surjective morphism from
$G$ onto $G$ is also injective.
\index{names}{Mal'tsev, Anatolii I.}%
\begin{theorem}[Malcev]\label{theoremMalcev}
 Any finitely generated residually finite group
is Hopfian.
\end{theorem}
In particular, by Theorem~\ref{theoremFGResiduallyFinite},
a finitely generated free group is Hopfian. This
can be proved directly as follows.
Let $\alpha:F(A)\to F(A)$ be a surjective morphism. Since
$\alpha(A)$ generates $F(A)$, it cannot have less than $\Card(A)$ elements.
Thus $\Card(\alpha(A))=\Card(A)$ and consequently $\alpha(A)$ is a basis
of $F(A)$. This implies that $\alpha$ is an automorphism.
\section{Free products}

Given two groups $G$ and $H$, the \emph{free product}\index{subject}{free!product of groups} of $G$
and $H$, denoted $G*H$\index{symbols}{G@$G*H$} is the set of all $g_1h_1\cdots g_nh_n$
with $n\ge 1$, $g_i\in G$ and $h_i\in H$ for $1\le i\le n$.
Thus the free group on $A$ is the free product of the infinite
cyclic groups generated by the $a\in A$.

The following result is known as the Kurosh Subgroup Theorem,
\index{subject}{Kurosh Subgroup Theorem}\index{subject}{Theorem!Kurosh}%
\index{names}{Kurosh, Aleksandr G.}%
\begin{theorem}[Kurosh]
 Any subgroup of a free product $G_1*G_2*\cdots *G_n$
is itself a free product of a free group and of groups conjugate
to subgroups of the $G_i$.
\end{theorem}

\section{Graphs}
Let $\Gamma=(V,E)$ be a graph
with $V$ as set of vertices and $E$ as set of edges. Each edge
$e\in E$ has an \emph{origin} $\alpha(e)$
\index{subject}{origin!of edge}\index{subject}{edge!origin of}%
\index{symbols}{alpha(e)@$\alpha(e)$}%
and an \emph{end} $\omega(e)$ \index{subject}{end!of edge}\index{subject}{edge!end of}
\index{symbols}{omega@$\omega(e)$}%
(also called its source and range). 
We allow multiple edges and thus we are considering \emph{directed multigraphs}.
\index{subject}{directed!graph}\index{subject}{graph!directed}
One also considers \emph{undirected}\index{subject}{undirected graph}
\index{subject}{graph!undirected} graphs in which an edge
is an unordered pair of vertices (or, equivalently, there is
an edge from $v$ to $w$ if and only if there is one from $w$ to $v$).
It is also possible to use undirected multigraphs, associating
to an edge a set of two vertices called its ends.

A \emph{labeled graph}\index{subject}{labeled graph}
\index{subject}{graph!labeled} is a  graph $G$ with a label on every edge,
which a letter of some finite alphabet $A$. We usually
denote $p\edge{a}q$ an edge from $p$ to $q$ labeled $a$.
%The label of a path $p_1\edge{a_1}p_2\cdots\edge{a_n}p_{n+1}$
%is the word $a_1a_2\cdots a_n$.

As a transition with what precedes, to every group $G$ with a given set $S$ of generators
is associated a labeled graph called its \emph{Cayley graph}.
\index{subject}{Cayley graph}\index{subject}{graph!Cayley}%
\index{names}{Cayley, Arthur}%
Its set of vertices is $G$ and there is an edge from $g$ to $h$
labeled $s$ whenever $gs=h$.

Two edges $e$ and $f$ are {\em consecutive}\index{subject}{edge!consecutive} whenever $\omega (e) = \alpha (f)$.
A {\em path}\index{subject}{path} in the graph is, as usual, a sequence $(e_1 , \ldots  , e_n)$
of consecutive edges. 
Its {\em origin}\index{subject}{path!origin}\index{subject}{origin!of path}
is $\alpha (e_1)$, its {\em end}\index{subject}{path!end}
\index{subject}{end!of path}$\omega (e_n)$ and the integer $n$ is its {\em length}.\index{subject}{length!of path}
By convention, there is, for each vertex $v$ a path of length $0$ of origin $v$ and end $v$. 
The path is a \emph{cycle}\index{subject}{cycle!in graph}
if its origin and end are equal. 
If $(e_1 , \dots , e_n )$ and $(f_1 , \dots , f_m)$ are two paths such that $\omega (e_n) = \alpha (f_1)$, then the concatenation of these paths is the path $(e_1 , \dots , e_n , f_1 , \dots , f_m)$.

For each
edge $e\in E$, we consider an inverse edge $e^{-1}$ from
$\omega(e)$ to $\alpha(e)$. 
By convention $(e^{-1})^{-1}$ is $e$.
A \emph{generalized
path} in $G$ is a sequence $(e_1,\ldots,e_n)$ of edges or their inverses such that, for all $i\in [1, n-1]$, the edges $e_i$ and $e_{i+1}$  are consecutive.
The generalized path is {\em reduced} whenever $e_i$ is different from $e_{i+1}^{-1}$ for all $i$.
Two paths are  {\em equivalent} whenever they can be obtain one from an other by a sequence of insertions or deletions of a sequence $(e , e^{-1})$. 
Every generalized path is equivalent to a unique reduced generalized path.
It is a \emph{generalized cycle}\index{subject}{cycle!generalized} if its origin and end are equal.

Let $v\in V$ be a vertex of the graph $\Gamma=(V,E)$. 
The \emph{fundamental group}\index{subject}{fundamental!group}
 $G(\Gamma,v)$\index{symbols}{G@$G(\Gamma,v)$} is the group formed by the generalized
cycles around $v$. When $G$ is connected, its isomorphism class does not depend
on $v$ and we denote $G(\Gamma)$ the fundamental group of $\Gamma$.
%\marginpar{FD: $\Sigma$?}

As well known, the fundamental group of a connected graph $G$
is a free group and a basis can be obtained as follows.
 Let $T$ be a \emph{spanning tree}\index{subject}{spanning tree}
of $\Gamma$ rooted at $v$, that is, a set of edges such that
for every $w\in V$ there is a unique path $p_w$ from $v$ to $w$
using the edges in $T$. Then, the set
\begin{equation}
\{p_{\alpha(e)}e p_{\omega(e)}^{-1}\mid e\in E\setminus T\}\label{eqFondamental}
\end{equation}
is a basis of $G(\Gamma,v)$.
\begin{example}\label{exampleGraph}
Let $\Gamma$ be the graph represented in Figure~\ref{figureGraphAppendix}
with $V=\{1,2\}$ and $E=\{e,f,g,h\}$.
\begin{figure}[hbt]
\tikzset{node/.style={circle,draw,minimum size=0.4cm,inner sep=0pt}}
\tikzstyle{every loop}=[->,shorten >=1pt,looseness=12]
\tikzstyle{loop left}=[in=130,out=220,loop]
\tikzstyle{loop right}=[in=330,out=50,loop]
\centering
\begin{tikzpicture}(20,12)(0,-2)
\node[node](1) at (0,0){$1$};\node[node](2) at(2,0){$2$};

\draw[left](1)edge [loop left]node{$e$}(1);
\draw[above, bend left, ->](1) edge node{$f$}(2);
\draw[above, bend left, ->](2)edge node{$g$}(1);
\draw[right](2) edge [loop right]node{$h$}(2);
\end{tikzpicture}
\caption{A connected graph.}\label{figureGraphAppendix}
\end{figure}

The set $T=\{f\}$ is a spanning tree  rooted at $1$. The corresponding
basis of $G(\Gamma,1)$ is $\{e,fg,fhf^{-1}\}$.
\end{example}
\section{Stallings graph}

For a graph $\Gamma=(V,E)$ labeled by an alphabet $A$, the \emph{group defined}
\index{subject}{group!defined by a labeled graph}%
by $\Gamma$ with respect to a vertex $v$ is the subgroup of the free group on $A$
formed by the labels of generalized paths
from $v$ to itself. Thus, if all labels are distinct, the
graph defined by $G$ is the fundamental graph of $G$.

A \emph{Stallings folding}\index{subject}{Stallings!folding}
\index{names}{Stallings, John R.}%
  of a labeled graph is the following transformation. Suppose that
two vertices $p,p'$ of a graph $G$   have  edges 
(or inverse edges) going to $q$
 with the same label $a$. Then we change $G$ to $G'$
by merging $p$ and $p'$.
A Stallings folding does not change the subgroup defined
by the graph. Indeed,  any generalized path in $G'$ can be obtained from
a generalized path in $G$ by insertion of paths of length $2$ labeled $aa^{-1}$.
A graph on which no Stallings folding can be performed is
called \emph{Stallings reduced}.
\index{subject}{Stallings!reduced graph}%

Given a finitely generated subgroup $H$ of the free group, there
is a unique Stallings reduced
 graph which defines $H$. This graph is called the \emph{Stallings graph}
\index{subject}{Stallings!graph of a subgroup} of the subgroup $H$.
\begin{example}
Let $A=\{a,b\}$ and let $H$ be subgroup of $F(A)$
generated by $\{a,bab,bb\}$. 

\begin{figure}[hbt]
\tikzset{node/.style={circle,draw,minimum size=0.4cm,inner sep=0pt}}
\tikzstyle{every loop}=[->,shorten >=1pt,looseness=12]
\tikzstyle{loop left}=[in=130,out=220,loop]
\tikzstyle{loop right}=[in=330,out=50,loop]
\centering
\begin{tikzpicture}
%1
\node[node](1) at (0,0){$1$};\node[node](2) at(2,0){$2$};
\node[node](3)at(1,-1){$3$};

\draw[left,->](1)edge [loop left]node{$a$}(1);
\draw[above, bend left, ->](1) edge node{$b$}(2);
\draw[above, bend left, ->](2)edge node{$b$}(1);
\draw[right,bend left,->](2) edge node{$a$}(3);
\draw[right,bend left,->](3) edge node{$b$}(1);

%2
\node[node](1) at (4,0){$1$};\node[node](2) at(6,0){$2$};

\draw[left](1)edge [loop left]node{$a$}(1);
\draw[above, bend left, ->](1) edge node{$b$}(2);
\draw[above, bend left, ->](2)edge node{$b$}(1);
\draw[right](2) edge [loop right]node{$a$}(2);
\end{tikzpicture}
\caption{A Stallings folding.}\label{figureGraphAppendix2}
\end{figure}
A graph defining $H$ is represented in Figure~\ref{figureGraphAppendix2}
on the left. The Stallings folding merging $2$ and $3$ gives
the graph represented on the right, which is the Stallings graph
of $H$.
\end{example}

\section{Semisimple algebras}

Let  $\Ga$ be an algebra over the field $\C$ of complex numbers.
It is said to be \emph{simple}\index{subject}{simple!algebra}
\index{subject}{algebra!simple} if it has no proper   nonzero two-sided ideal.
As an equivalent definition, an algebra is semisimple if does
not contain nonzero nilpotent ideals.
The algebra $\M_n$ of $n\times n$-matrices with coefficients in $\C$ is simple.
Any automorphism of $\M_n$ is an inner automorphism,
that is, of the form $M\mapsto AMA^{-1}$
for some invertible matrix $A\in\M_n$.

A direct sum $\Ga_1\oplus\ldots\oplus\Ga_k$ of simple algebras $\Ga_i$ is called a \emph{semisimple} algebra.
\index{subject}{semisimple algebra}\index{subject}{algebra!semisimple}%
\begin{theorem}[Wedderburn]
  A finite dimensional algebra $\Ga$ over $\C$
  is semisimple if and only if $\Ga=\M_{n_1}\oplus \ldots \oplus\M_{n_t}$.
\end{theorem}
\index{subject}{Wedderburn Theorem}\index{subject}{Theorem!Wedderburn}%
\index{names}{Wedderburn, Joseph}%

Moreover, if $\Ga=\M_{n_1}^{(a_1)}\oplus\cdots\oplus\M_{n_k}^{(a_k)}$
where the $n_i$ are distinct and each $\M_{n_i}$ is repeated $a_i$ times,
the $n_i$ and the $a_i$ are determined uniquely up to a permutation.
Consequently any embedding of an algebra $\Ga_1=\M_{m_1}\oplus\ldots\oplus\M_{m_t}$
into an algebra $\Ga_1=\M_{n_1}\oplus\cdots\oplus\M_{n_s}$ is, up to conjugacy,
of the type $\varphi=\varphi_1\oplus\ldots\oplus\varphi_t$ with
\begin{displaymath}
  \varphi_i=\id_{n_1}^{(a_{i1})}\oplus\ldots\oplus\id^{(a_{it})}
  \end{displaymath}
where $\id_n$ is the identity of $\M_n$ and $\id_n^{(a)}:\M_n\to\M_{an}$ is the morphism
$x\mapsto (x,\ldots,x)$ ($a$ times). In other terms, each morphism $\varphi_i:\Ga_1\to \M_{n_i}$
has the form

  \begin{displaymath}
  \varphi_i(M_1,\ldots,M_t)=\kbordermatrix{
   &a_{i1}  &  &a_{it}&  \cr
   a_{i1} &  \begin{array}{ccc|}M_1& & \\& \ddots&\\ & &M_1\\ \hline\end{array}&\cr
      & & \ddots\cr 
   a_{it}  &         &    &  \begin{array}{|ccc}\hline M_t & & \\ &\ddots& \\ & &M_t\end{array}
  }
  \end{displaymath}  
\section{Notes}

We have only briefly recalled the basic definitions and properties of free groups.
For a systematic exposition (in particular of the Nielsen-Schreier Theorem), see
\citep{LyndonSchupp2001} or \citep{MagnusKarrassSolitar2004}.
\index{names}{Schupp, Paul}%

See  \cite{LyndonSchupp2001} (p.197)
for a proof of  Malcev Theorem~\ref{theoremMalcev}.
The Stallings graphs are from \cite{Stallings1983}
\index{names}{Stallings, John R.}%
(see also \citep{KapovichMyasnikov2002}).
\index{names}{Kapovich, Ilya}\index{names}{Myasnikov, Alexei}%

The elementary properties of semisimple algebras (serving as a preparation
of the properties of finite dimensional $C^*$-algebras in Chapter~\ref{chapterBratteli}),
in particular Wedderburn Theorem,
can be found in \cite{Lang2002}.\index{names}{Lang, Serge}

%\appendix
\chapter{Linear algebra}\label{appendixLinearAlgebra}

In this appendix, we present some results concerning matrices
playing an important role in the book.

\section{Perron-Frobenius Theorem}
A $P\times Q$-matrix with coefficients in a ring  $K$ is a map $M:P\times Q\to K$.
The image of $(p,q)$, denoted $M_{p,q}$, is called the coefficient of $M$ at row $p$
and column $q$. The sum of two $P\times Q$-matrices $M,N$ is the $P\times Q$-matrix $M+N$ defined
by $(M+N)_{p,q}=M_{p,q}+N_{p,q}$. The product of a $P\times Q$-matrix $M$ and
a $Q\times R$ matrix $N$ is the $P\times R$-matrix $MN$  defined by
\begin{displaymath}
  (MN)_{p,r}=\sum_{q\in Q}M_{p,q}N_{q,r}.
  \end{displaymath}

We denote by $M^t$ the \emph{transpose}\index{subject}{transpose of matrix}
\index{subject}{matrix!transpose} of $M$, which is defined by $M_{p,q}^t=M_{q,p}$.
The matrix $M$ is \emph{symmetric}\index{subject}{symmetric!matrix}
if $M^t=M$ and \emph{antisymmetric} \index{subject}{antisymmetric matrix}
\index{subject}{matrix!antisymmetric} if $M^t=-M$.

Coming to a more familiar notation,
an $m\times n$-matrix is a $P\times Q$-matrix with $P=\{1,2,\ldots,m\}$
and $Q=\{1,2,\ldots,n\}$.

A matrix with real coefficients
is \emph{nonnegative}\index{subject}{matrix!nonnegative}
\index{subject}{nonnegative!matrix} (resp. \emph{positive})
\index{subject}{positive!matrix}\index{subject}{matrix!positive}%
if all its elements are nonnegative (resp. positive).
The same terms are used for vectors.
\index{subject}{vector!nonnegative}\index{subject}{nonnegative!vector}\index{positive!vector}%
\index{subject}{vector!positive}%

A nonnegative square matrix $M$ is \emph{irreducible}
\index{subject}{irreducible!matrix}\index{subject}{matrix!irreducible}%
if for every pair $i,j$ of indices, there is an integer $n\ge 1$
such that $M^n_{i,j}>0$.

A nonnegative square matrix $M$ is \emph{primitive}\index{subject}{primitive! matrix}\index{subject}{matrix!primitive},
if there is some
$n\ge 1$ such that all entries of $M^n$ are positive. A primitive matrix is irreducible
but not conversely.

To every finite graph $G=(V,E)$ is associated a nonnegative matrix called its
\emph{adjacency matrix}.\index{subject}{adjacency matrix!of graph}
\index{subject}{graph!adjacency matrix}\index{subject}{matrix!adjacency}%
It is the $V\times V$-matrix with coefficients
\begin{displaymath}
  M_{v,w}=\Card\{e\in E\mid \alpha(e)=v,\omega(e)=w\}.
\end{displaymath}
The matrix is irreducible if and only if $G$ is strongly connected.
It is primitive if and only if
\begin{enumerate}
\item[(i)] $G$ is strongly connected.
\item[(ii)] For every $v\in V$ the $\gcd$ of the lengths
  of the cycles around $v$ is $1$.
  \end{enumerate}

The following result\index{subject}{Perron-Frobenius Theorem}\index{subject}{Theorem!Perron-Frobenius} 
is well-known.
\begin{theorem}[Perron, Frobenius]
\label{theorem:perron}\index{names}{Perron, Oskar}\index{names}{Frobenius, Georg F.}
Let $M$ be a nonnegative   real square matrix. 
Then
\begin{enumerate}
\item[\rm(i)] $M$ has a 
 positive  eigenvalue $\lambda_M$ such that $|\mu|\le\lambda_M$ for
every  eigenvalue $\mu$ of $M$.
\item[\rm(ii)] There corresponds to $\lambda$
  a nonnegative eigenvector $v$.
  \end{enumerate}
If $M$ is irreducible (resp. primitive), then $\lambda_M$ is simple, one has $|\mu|<\lambda_M$
for every other eigenvalue and there corresponds to $\lambda_M$
  a positive eigenvector $v$. Every nonnegative eigenvector is colinear to $v$.
  One has additionally
  \begin{enumerate}
\item[\rm(iii)] The sequence $(M^n/\lambda_M^n)$
converges in mean (resp. converges) at geometric rate 
 to the matrix $wv$ where $v,w$ are  positive
left and right eigenvectors such that $Mw=\lambda_Mw$, $vM=\lambda_Mv$
and $vw=1$.
\end{enumerate}
\end{theorem}

The theorem expresses in particular that if a matrix $M$ is primitive, its spectral radius $\rho(M) = \max\{|\lambda| \mid \lambda \in \mathrm{Spec}(M)\}$ is an eigenvalue of $M$ which is algebraically simple.
Furthermore, any eigenvalue of $M$ different from $\rho(M)$ has modulus less than $\rho(M)$.
We call $\rho(M)$ the {\em dominant eigenvalue of $M$}. 
\index{subject}{dominant eigenvalue}
By abuse of language, when $M$ is the composition matrix of a primitive endomorphism $\sigma$, we call $\rho(M)$ the {\em dominant eigenvalue of $\sigma$}.

The term \emph{geometric rate}\index{subject}{geometric rate of convergence}
 of convergence used in assertion (iii) means
that there is a constant $c>0$ and a real number $r<1$ such that
for all $n\ge 0$
\begin{equation}
\left|\left|\frac{1}{\lambda_M^n}M^n-wv\right| \right|\le c\ r^n.\label{eqGeometricRate}
\end{equation}
We can choose for $r$ the quotient $r=\mu/\rho(M)$ where
$\mu$ is the maximum of the $|\lambda|$ for $\lambda$
an eigenvalue of $M$ other than $\rho(M)$.

As an important particular case, a real nonnegative $A\times A$-matrix $M$ is \emph{stochastic}
\index{subject}{stochastic matrix}\index{subject}{matrix!stochastic}%
if its rows have sum $1$, that is,
\begin{displaymath}
  \sum_{b\in A}M_{a,b}=1
  \end{displaymath}
for every $a\in A$. It is easy to verify that the dominant eigenvalue
of a stochastic matrix is $1$. The column vector with all coefficients $1$
is a corresponding eigenvector. When $M$ is irreducible, by Theorem~\ref{theorem:perron},
the sequence $\frac{1}{n}\sum_{i=0}^{n-1}M^i$ converges to the matrix with all rows equal
to the row eigenvector $v$ of sum $1$ for the eigenvalue $1$, that is, such that $vM=v$.
\section{Linear programming}

A \emph{linear program}\index{subject}{linear!program} is
given by a triple $(A,b,f)$
where  $f:\R^n\to\R$ is a linear map,  $A$ is a real $m\times n$-matrix
and  $b\in\R^m$ is a vector. A \emph{solution}\index{subject}{solution of linear
  program} is a vector $x\in\R_+^n$ such that $Ax=b$ and $f(x)\ge f(y)$
for all $y\in\R_+^n$ such that $Ay=b$. Thus,
a solution maximizes $f$ under the \emph{constraint} $Ax=b$.

A solution $x=(x_1,\ldots,x_n)$ of the equation $Ax=b$ is a \emph{basic solution}
if $x_i=0$
for for a set of $n-m$ indices $i$.

The Fundamental Theorem of Linear Programming is the following result.
\begin{theorem}
  If the linear program given by $(A,b,f)$ with $A$
  of rank $m$ has  a solution,
  it has a basic solution.
\end{theorem}
As a very simple example, if $A=\begin{bmatrix}1&1&1\end{bmatrix}$,
$b=1$
and $f(x)=x_1+x_2+x_3$, then $\begin{bmatrix}1/2&1/2&0\end{bmatrix}$
  is a basic solution.
  
The proof is closely related to the \emph{simplex algorithm}
\index{subject}{simplex!algorithm}%
which consists in building a sequence of basic solutions
of the equation $Ax=b$ increasing successively the values
of the function $f$.

\section{Notes}

The Perron-Frobenius Theorem is a classical and very useful result.
A proof can be found in \cite{Queffelec2010} or in the
classical \cite{Gantmacher1959}.

See \cite{Chvatal1983}
\index{names}{Chvatal, Va\v{s}ek}%
or \cite{LuenbergerYe2016}
\index{names}{Luenberger, David G.}\index{names}{Ye, Yinyu}%
for an introduction to Linear Programming and a proof of the
Fundamental Theorem.

\chapter{Topological entropy}\label{appendixTopologicalEntropy}
We give a short introduction to the notion of topological entropy.
\section{Entropy and covers}
Let $(X,T)$ be a topological dynamical system.

If $\alpha,\beta$ are covers of $X$, then
\emph{join} of $\alpha,\beta$, denoted $\alpha\vee\beta$ is
the cover formed
of the $A\cap B$ for $A\in\alpha$ 
and $B\in\beta$. This extends to the join $\vee_{i=0}^{n-1}\alpha_i$
of $n$ covers
$\alpha_i$ for $0\le i\le n-1$. We also denote  $T^{-j}\alpha$
the cover formed 
by the $T^{-j}A$ for $A\in\alpha$.

If $\alpha$ is an open cover of $X$, define the \emph{entropy}
\index{subject}{entropy!of cover} of $\alpha$
by $H(\alpha)=\log N(\alpha)$ where $N(\alpha)$ is the number of elements
of a finite subcover of $\alpha$ with minimal cardinality.

The \emph{entropy of $T$ relative to}
\index{subject}{entropy!relative to cover} $\alpha$ is the possibly infinite
real number
\begin{displaymath}
  h(T,\alpha)=\lim_{n\to\infty}\frac{1}{n}H(\vee_{i=0}^{n-1}T^{-i}\alpha)
\end{displaymath}
(the limit can be shown to exist).

The \emph{topological entropy}\index{subject}{topological!entropy}
\index{subject}{entropy!topological} of $T$ is
\begin{displaymath}
  h(T)=\sup_\alpha h(T,\alpha)
  \end{displaymath}
where $\alpha$ ranges over all open covers of $X$. It can be shown that
entropy is an invariant under conjugacy.
\section{Computation of the entropy}
The following result is useful to compute the topological entropy.
We denote by $\diam(\alpha)=\sup_{A\in \alpha}\diam(A)$ the diameter of
an open cover $\alpha$ where $\diam(A)$ is the diameter of the set $A$.
\begin{theorem}\label{theoremWalters}
  If $(\alpha_n)_{n\ge 0}$ is a sequence of open covers such that
  $\diam(\alpha_n)\to 0$, then $h(T)=\lim_{n\to\infty}h(T,\alpha_n)$.
  \end{theorem}

 The topological entropy of a shift
space $X$ has a simple expression in terms of its factor complexity $p_n(X)$.

\begin{theorem}\label{corollaryEntropyShiftSpace}
  The topological entropy of a shift space $(X,S)$ is
  $h(S)=\lim_{n\to\infty}\frac{1}{n}\log p_n(X)$.
\end{theorem}
Actually, the topological entropy of  a shift space
on the alphabet $A$ is the
entropy of $S$ relative to the natural partition 
$\{[a]\mid a\in A\}$.

\section{Notes}
Topological entropy is a counterpart for topological
dynamical systems of the notion of entropy for measure-theoretic system.
Our presentation follows \cite{Walters1975}. Theorem~\ref{theoremWalters}
is Theorem 7.6 and Theorem~\ref{corollaryEntropyShiftSpace}
is Theorem 7.13.

%\include{chapter10}

%-------------------- les references
\renewcommand{\bibname}{References}
%\phantomsection  
\renewcommand{\bibsection}{%
  \chapter*{\bibname}\addcontentsline{toc}{chapter}{\bibname}} 
\setlength{\bibhang}{\parindent}
\markboth{\bibname}{\bibname}\markright{\bibname}
\bibliographystyle{plainnat}
\bibliography{dimensionGroups}

\printindex{names}{\bf Name index}
\printindex{subject}{\bf Subject index}
\printindex{symbols}{\bf Index of symbols}
\end{document}